\documentclass[11pt,reqno]{amsart}
\usepackage{amssymb,mathrsfs,color}
\allowdisplaybreaks
\usepackage{txfonts}

\usepackage{bbm}
\usepackage{enumerate}
\usepackage{graphicx}
\usepackage{xcolor} 
\usepackage{geometry}
\usepackage{amsfonts,amssymb,amsmath,mathtools}
\usepackage{slashed}
\usepackage{footnote}

\usepackage{ mathrsfs }
\usepackage[title]{appendix}
\usepackage{wrapfig}

\usepackage{cancel}

\usepackage{cite}

\definecolor{green}{rgb}{0,0.8,0} 

\renewcommand{\Re}{\mathrm{Re}}
\renewcommand{\Im}{\mathrm{Im}}

\newcommand{\0}{\emptyset}

\newcommand{\eps}{\epsilon}
\newcommand{\veps}{\varepsilon}

\newcommand{\bff}{{\bf f}}

\newcommand{\bfv}{{\bf v}}

\newcommand{\bbC}{\mathbb C}

\newcommand{\bbN}{\mathbb N}

\newcommand{\bbR}{\mathbb R}
\newcommand{\bbS}{\mathbb S}

\newcommand{\bbZ}{\mathbb Z}

\newcommand{\calD}{\mathcal D}

\newcommand{\calF}{\mathcal F}

\newcommand{\calH}{\mathcal H}

\newcommand{\calK}{\mathcal K}
\newcommand{\calL}{\mathcal L}

\newcommand{\calQ}{\mathcal Q}
\newcommand{\calR}{\mathcal R}
\newcommand{\calS}{\mathcal S}
\newcommand{\calT}{\mathcal T}

\newcommand{\calV}{\mathcal V}

\definecolor{deepgreen}{cmyk}{1,0,1,0.5}

\newcommand{\C}{\mathbb{C}}

\newcommand{\R}{\mathbb{R}}

\newcommand{\Z}{\mathbb{Z}}

\newcommand{\fy}{\varphi}

\newcommand{\x}{\xi}

\makeatletter

\newcommand{\Rmnum}[1]{\expandafter\@slowromancap\romannumeral #1@}
\makeatother

\renewcommand{\bar}{\underline}

\newcommand{\EQ}[1]{\begin{equation}\begin{split} #1 \end{split}\end{equation}}

\newcommand{\Del}[1]{}

\numberwithin{equation}{section}

\newtheorem{theorem}{Theorem}[section]
\newtheorem{corollary}[theorem]{Corollary}
\newtheorem{lemma}[theorem]{Lemma}
\newtheorem{proposition}[theorem]{Proposition}
\newtheorem{remark}[theorem]{Remark}
\newtheorem{definition}[theorem]{Definition}

\newcommand{\sign}{\operatorname{sign}}

\renewcommand\Re{\mathrm{Re}\,}
\renewcommand\Im{\mathrm{Im}\,}

\renewcommand{\div}{\mathrm{div}\,}

\newcommand{\zetab}{\overline{\zeta}}

\newcommand{\frakE}{\mathfrak{E}}
\newcommand{\frakD}{\mathfrak{D}}

\newcommand{\frakL}{\mathfrak{L}}
\newcommand{\frakN}{\mathfrak{N}}
\newcommand{\xb}{\overline{x}}
\newcommand{\yb}{\overline{y}}
\newcommand{\zb}{\overline{z}}
\newcommand{\trho}{\widetilde{\rho}}

\newcommand{\vphi}{\varphi}
\newcommand{\pphi}{\Pi_{\Phi^{\perp}}\varphi}
\newcommand{\hvphi}{\widehat{\varphi}}

\newcommand{\hfrakN}{\widehat{\mathfrak{N}}}

\newcommand{\frakA}{\mathfrak{A}}
\newcommand{\tV}{\widetilde{V}}
\newcommand{\Ai}{\textrm{Ai}}
\newcommand{\Bi}{\textrm{Bi}}

\newcommand{\tf}{\widetilde{f}}
\newcommand{\ttau}{\widetilde{\tau}}

\newcommand{\tbff}{\widetilde{\bff}}

\newcommand{\thbar}{\widetilde{\hbar}}
\newcommand{\tQ}{\widetilde{Q}}
\newcommand{\xbar}{\underline{x}}
\newcommand{\alphab}{\underline{\alpha}}
\newcommand{\taub}{\overline{\tau}}
\newcommand{\qb}{\overline{q}}
\newcommand{\wb}{\overline{w}}

\newcommand{\tVb}{\underline{\widetilde{V}}}

\newcommand{\tiltau}{\tilde{\tau}}

\newcommand{\xh}{\overline{x}^{\hbar}}
\newcommand{\xhb}{\overline{x}^{\widetilde{\hbar}}}
\newcommand{\txh}{\widetilde{x}^{\hbar}}

\newcommand{\Fh}{\mathcal{F}^{\hbar}}
\newcommand{\Fhb}{\mathcal{F}^{\widetilde{\hbar}}}
\newcommand{\dftau}{\frac{\lambda'(\tau)}{\lambda(\tau)}}

\newcommand{\Sh}{S^{\hbar}}
\newcommand{\ddtau}{-\left(\left(\partial_{\tau}+\frac{\lambda'}{\lambda}R\partial_{R}\right)^{2}+\frac{\lambda'}{\lambda}\left(\partial_{\tau}+\frac{\lambda'}{\lambda}R\partial_{R}\right)\right)}

\newcommand{\ddDtau}{-\left(\left(\partial_{\tau}+\frac{\lambda'}{\lambda}R\partial_{R}\right)^{2}+3\frac{\lambda'}{\lambda}\left(\partial_{\tau}+\frac{\lambda'}{\lambda}R\partial_{R}\right)\right)}
\newcommand{\tH}{\tilde{H}}
\newcommand{\Shj}{S^{\hbar_{j}}}
\newcommand{\Shk}{S^{\hbar_{k}}}
\newcommand{\Sho}{S^{\hbar_{1}}}
\newcommand{\Sht}{S^{\hbar_{2}}}
\newcommand{\Shth}{S^{\hbar_{3}}}

\def\hbzwd{\hbar^{\frac23}}

\newtheorem{defi}{Definition}

\def\what{\widehat}
\def\calR{\mathcal{R}}
\def\one{\mathbbm{1}}

\def\nn{\nonumber}

\begin{document}

\title{A stability theory beyond the co-rotational setting for critical Wave Maps blow up}
\author{Joachim Krieger
\and Shuang Miao 
\and Wilhelm Schlag}


\begin{abstract}
	We exhibit non-equivariant perturbations of the  blowup solutions constructed in \cite{KST} for energy critical wave maps into $\bbS^2$. Our admissible class of perturbations is an open set in some sufficiently smooth topology and vanishes near the light cone. We show that the blowup solutions from \cite{KST}  are rigid under such perturbations, including the space-time location of blowup. As blowup is approached, the dynamics agree with the classification obtained in \cite{DJKM}, and all six symmetry parameters converge to limiting values. Compared to the previous work \cite{KMiao} in which the rigidity of the blowup solutions from \cite{KST}  under equivariant perturbations was proved, the class of perturbations considered in the present work does not impose any symmetry restrictions. Separation of variables and decomposing into angular Fourier modes leads to an infinite system of coupled nonlinear equations, which we solve for small admissible data. The nonlinear analysis is based on the distorted Fourier transform, associated with an infinite family of Bessel type Schr\"odinger operators on the half-line indexed by the angular momentum~$n$. A semi-classical WKB-type spectral analysis relative to the parameter $\hbar=\frac{1}{n+1}$ for large $|n|$ allows us to effectively determine the distorted Fourier basis for the entire infinite family. Our linear analysis is based on the global Liouville-Green transform as in the earlier works \cite{CSST, CDST}.  
\end{abstract}

\maketitle


\section{Introduction}
This work is dedicated to the energy critical Wave Maps equation $u: \R^{2+1}\longrightarrow \bbS^2$, and more specifically, to developing a theory which allows to establish the stability (and in fact, rigidity) of a certain class of its finite time blow up solutions {\it{without any symmetry restrictions}}. Specifically, the final result established shows that the finite time blow up solutions constructed in \cite{KST, GaoK} {\it{in the co-rotational setting}}, and with blow up scaling $\lambda(t) = t^{-1-\nu}$ sufficiently close to the self similar rate, i.e.,  $0<\nu<\nu_*$ and $\nu_*>0$ sufficiently small, are stable under {\it{arbitrary}} sufficiently smooth and small perturbations, {\it{without any symmetry restrictions}}. The precise result, stated below, will show that the perturbations result in slightly modulated profiles in perfect agreement with the general classification result due to Duyckaerts-Jia-Kenig-Merle~\cite{DJKM}, and it appears to be the first result which shows the full stability of a non-scattering solution in the setting of Wave Maps (as usual,  the open set of perturbations is relative to a stronger topology than energy). 
A peculiarity of the solutions constructed in \cite{KST, GaoK} is the fact that even before blow up, they are only of finite regularity (in fact $H^{1+\nu-}$-regularity, with $\nu$ as before), on account of a shock they display on the light cone centered at the singularity. This shock appears to endow the solutions with a strong rigidity, provided $\nu$ is sufficiently small, which may be interpreted as the reason why the space-time location of the singularity for the perturbed solution is un-changed, a feature which would be clearly false for blow up solutions which are smooth before the blow up time. It is important to keep in mind that the small perturbations added to the original data are much smoother than the original data. It is by now known that many different types of hyperbolic equations admit blow up solutions of similar character, and it is hoped that the present work may provide techniques to analyze their stability properties as well. In fact, the asymptotic Fourier methods developed here may be of use for much more general problems, beyond the narrow context of highly specific blow up solutions considered in this work. 
Before stating the main theorem, we recall the basic setup:

\subsection{Blow up in the co-rotational setting} 
Recall that a Wave Map  $u: \R^{2+1}\longrightarrow \bbS^2\hookrightarrow \R^3$ is a vector valued solution to the system
\begin{equation}\label{eq:BasicWMS2target1}
	\Box u = \left(|u_t|^2 - |\nabla_x u|^2\right)u,\quad \Box = -\partial_t^2 + \Delta_x
\end{equation}
This model is known to admit a rich class of static solutions of finite energy, due to the large class of finite energy harmonic maps from $\R^2$ to $\bbS^2$. Amongst these, the {\it{stereographic projection map}} $\mathcal{Q}: \R^2\longrightarrow \bbS^2$ is the one of minimal energy amongst the non-trivial harmonic maps. It is an example of a so-called {\it{co-rotational harmonic map}}, in that it admits the representation 
\[
\mathcal{Q}(r,\theta) = \left(\begin{array}{c}\cos\theta \sin Q\\ \sin\theta \sin Q\\ \cos Q\end{array}\right), 
\]
where $Q$ is only dependent on $r$, and in fact for the stereographic projection we have $Q(r) = 2\arctan r$. Making the co-rotational ansatz for \eqref{eq:BasicWMS2target1} leads to a remarkable simplification of the equation, removing the derivatives. In fact, making the ansatz 
\[
u(t, r) =  \left(\begin{array}{c}\cos\theta \sin v\\ \sin\theta \sin v\\ \cos v\end{array}\right),
\]
where now $v = v(t, r)$ is a {\it{scalar unknown}}, one derives the equation
\begin{equation}\label{eq:WMcorotational}
	\Box v = \frac{\sin 2v}{2r^2},\quad \Box = -\partial_t^2 + \partial_r^2 + \frac{1}{r}\partial_r. 
\end{equation}
The fact that $v(t,r) = 2\arctan r$ is a static solution of this problem is then easily verified by direct calculation. The simplified model \eqref{eq:WMcorotational} has been the subject of intense numerical and theoretical investigations, see for example \cite{STZ1, STZ2, ShStruwe}. In particular, the crucial role that harmonic maps play in the formation of singularities for this model was first stressed in~\cite{Struwe}, and numerical investigations showing the bubbling-off of co-rotational harmonic maps were done in~\cite{BCT}. 
The first rigorous construction of blow-up solutions within the co-rotational setting is due to~\cite{KST}, later followed by a different blow up construction in~\cite{RaphRod}. While the former constructed a {\it{continuum of blow up solutions}} of the form 
\[
v(t, r) = Q(\lambda(t)r) +\epsilon(t, r)
\]
where $\lambda(t) = t^{-1-\nu}$, $\nu>\frac12$ (later refined to $\nu>0$ in~\cite{GaoK}, see also~\cite{KS_FullRange}), and moreover with $\epsilon(t, \cdot)\in H^{1+\nu-}$, the solutions in~\cite{RaphRod} are $C^\infty$ before formation of the singularity, and were shown to be {\it{stable within the co-rotational class}} (as usual,  the open set of perturbations is relative to a stronger topology than energy), a feature missing in~\cite{KST, GaoK}. Co-rotational stability of the solutions in \cite{KST, GaoK} for $\nu$ small enough was proved recently in~\cite{KMiao}. It is to be noted that the technique for showing stability in~\cite{KMiao} are completely different than the ones giving stability in~\cite{RaphRod}. In particular, the work \cite{KMiao} relies heavily on the use of the distorted Fourier transform. In the present work, Fourier techniques will also play an all-encompassing role.

\subsection{Beyond the co-rotational setting}

As is clear from comparison of \eqref{eq:WMcorotational}, \eqref{eq:BasicWMS2target1}, the stability problem outside of the co-rotational context is of a quite different character than the co-rotational one: the nonlinear source terms display derivatives, and in particular they lose the smoothing property which the purely semilinear source terms in the co-rotational setting possess. In fact, the classical local well-posedness theory as developed for example in \cite{KlMach, KlSel1, KlSel2}, reveals that in order to get strong local well-posedness all the way to the critical scaling level (which is essentially required for us since the solutions we perturb are only of regularity $H^{1+\nu-}$, with $\dot{H}^1$ the scaling critical space), a special algebraic cancellation feature in the nonlinearity is generally required (so-called null-structure); the Wave Maps problem has such a structure. Problems with this structure can be handled by means of the $H^{s,\,\delta}$ (also referred to as $X^{s,b}$-spaces) functional framework, which, however, appears intimately adapted to the free wave equation. Perturbing around the co-rotational blow ups leads to a very different kind of operator than the standard d'Alembertian. 
While it would be very interesting to develop abstract function spaces similar to the $H^{s,\delta}$-spaces to handle more general hyperbolic problems encompassing in particular our setup, our strategy in this work is to pursue a reductive ansatz, employing suitable coordinates and a separation of variables to arrive at a countable family of $1+1$-dimensional wave operators indexed by the angular frequency, and which we analyze asymptotically to reduce the analysis to a setup similar to the one in~\cite{KMiao}. 
The asymptotic analysis of the infinite family of wave operators uses techniques developed in~\cite{CSST, CDST}. This allows us to develop sufficiently precise parametrices for the wave operators in order to analyze the propagation of the shock on the light cone in the various nonlinear interactions and at various angular frequencies. An abstract $H^{s,\delta}$ framework is then substituted for in this work by a precise concept of  {\em admissible singular functions}, see Definition~\ref{defi:xsingulartermsngeq2adm}, which is both compatible with the family of wave propagators at the various angular frequencies, as well as the structure of the nonlinear source terms (since, being formulated on the distorted Fourier side, it allows a natural translation to the physical side). 
In addition to handling an infinite family of wave propagators at different angular frequencies, one also needs to deal with the inherent instabilities at three exceptional angular momenta, namely $n = 0,\pm1$. The Fourier representations at these angular momenta are exceptional, due to the presence of resonances/eigenvectors at frequency zero. These cause degeneracies in the corresponding spectral measures, and require a separate analysis of the {\em continuous spectral part} and the {\em unstable (discrete) spectral part}.  We then rely on the fact that all these wave operators admit a super-symmetric partner, owing to a remarkable product structure already taken advantage of in the co-rotational setting in \cite{RaphRod, KMiao}, which reduces an analysis of the wave equations at exceptional angular momenta to an analysis of the super-symmetric partner equation (which no longer exhibits a singular spectral measure), as well as the ODE governing the evolution of the instability. 
\\
It turns out that the conjunction of the shock on the light cone and the nonlinear source terms involving derivatives produce source terms for the evolution of the unstable modes which would lead to divergence. Handling this issue requires the use of some form of modulation theory, which, however, differs from the more standard one in that it is the {\it{action of internal symmetries on the singular part of the original co-rotational blow up}} which is being used to counteract other singular terms arising from the interaction of the perturbation and the original co-rotational singularity (i.e., a shock along the light cone). The only relevant symmetries here are scaling, three types of rotations on the target, as well as Lorentz transforms, which are in fact all symmetries leaving the light cone invariant.   

\subsection{Statement of the main theorem}\label{subsec:maintheorem}
The following is the result of the developments of this work, proved at the very end. Recall that the co-rotational finite time blow up solutions in \cite{KST, GaoK} are of the form (in polar coordinates)
\[
\Phi(t,\theta,r) = \left(\begin{array}{c}\cos\theta\sin U\\ \sin\theta\sin U\\ \cos U\end{array}\right)
\]
where $U = U(t,r) = Q(\lambda(t)r) + \epsilon(t,r)$ where $\lambda(t) = t^{-1-\nu}$, $Q(r) = 2\arctan r$, and the error term $ \epsilon(t,r)$ vanishes asymptotically in a suitable sense in the light cone centred at the singularity $(t,r) = (0,0)$ as $t$ approaches $0$. The following theorem will be rendered completely precise in the ensuing analysis, with a delicate description of the radiation term $\delta\Psi$ in the light cone. 
\begin{theorem}\label{thm:Main} Let $\Phi = \left(\begin{array}{c}\cos\theta\sin U\\ \sin\theta\sin U\\ \cos U\end{array}\right)$, $U = U(t,r)$ be one of the finite time co-rotational blow up solutions constructed in \cite{KST, GaoK}, with $\nu>0$ sufficiently small, and restricted to the space time slab $(0,t_0]\times \R^2$, where $t_0 = t_0(\nu)$ is sufficiently small. Then there are $\delta_* = \delta_*(\nu, t_0)>0$, $C_1 = C_1(\nu)\gg 1$ with the following property: let $(\delta\Phi_1, \delta\Phi_2): \R^2\longrightarrow \R^3$ with the property that 
	\[
	\left\|(\delta\Phi_1, \delta\Phi_2)\right\|_{H^{C_1}(\R^2)}<\delta_*,
	\]
	and such that 
	\[
	\Psi[t_0]: = \left(\Phi(t_0,\cdot)+\delta\Phi_1, \partial_t\Phi(t_0,\cdot) + \delta\Phi_2\right): \R^2\longrightarrow \bbS^2\times T\bbS^2
	\]
	constitutes a data set for \eqref{eq:BasicWMS2target1}, and furthermore $\delta\Phi_j$, $j = 1,2$, is supported inside the disk $r = |x|<\frac{t_0}{2}$. Then the solution for \eqref{eq:BasicWMS2target1}, with data $\Psi[t_0]$ at time $t = t_0$ admits a solution $\Psi$ on the space-time slab $(0,t_0]\times \R^2$ which can be written in the form 
	\[
	\Psi= \mathcal{L}_{v(t)}\mathcal{R}_{h(t)}^{\alpha(t),\beta(t)}\mathcal{S}_{c(t)}\left( \Phi + \delta\Psi\right),
	\]
	where $ \mathcal{L}_{v(t)}$ is a time dependent Lorentz transform, $\mathcal{R}_{h(t)}^{\alpha(t),\beta(t)}$ a general rotation on the target $\bbS^2$ parametrized as the product of three elementary rotations,  and $\mathcal{S}_{c(t)}$ a continuously varying change of scale. The radiation part $\delta\Psi$ vanishes asymptotically in the light cone near the singularity in the energy topology, and the modulation parameters $v(t), \alpha(t),\beta(t), h(t),c(t)$ approach finite limiting values as $t\rightarrow 0$.  
	
\end{theorem}

\section*{Acknowledgements} The authors would like to express their gratitude to Roland Donninger and Sergiu Klainerman for their interests in this work. Part of this work was done when the second author was conducting postdoctoral research at the \'Ecole Polytechnique F\'ed\'erale de Lausanne, and he thanks that institution for its hospitality and support.  J.~Krieger was supported by the Swiss National Fund via Consolidator Grant BSCGI0 157694, S.~Miao was supported by National Key R \& D Program of China 2021YFA1001700, and  NSFC grants 12071360, 12221001, \& 12326344. W.~Schlag was supported by the NSF grant DMS-1902691 and DMS-2054841.

\section{Summary of the sections of the paper}

Very roughly speaking, the paper is divided into an elaborate linear part, which provides the basis for the distorted Fourier analysis at arbitrary angular frequencies, as well as a nonlinear part, where a machinery is developed to handle the various nonlinear source terms, in particular the null-form terms, and with irregular ``inputs". The nonlinear part also encompasses the analysis of the unstable modes arising at angular frequencies $n = 0,\pm 1$, which forces us to employ the modulation techniques. Here we give a brief description of the linear as well as the nonlinear part.
To get things started, we of course need the right formulation of the perturbation problem, which is accomplished in section~\ref{sec:perturb formulation} following an ansatz by Davila, DelPino, Wei~\cite{DavilaDelPWei}. Specifically, letting as before $\Phi$ denote the unperturbed co-rotational blow up solution, and introducing the frame (which constitutes an orthonormal basis for $T_{\Phi}\bbS^2$)
\[
E_1 = \left(\begin{array}{c}\cos\theta \cos U\\ \sin\theta\cos U\\ -\sin U\end{array}\right),\quad E_2 = \left(\begin{array}{c}-\sin\theta\\ \cos\theta\\ 0\end{array}\right),
\]
we make the ansatz 
\begin{align*}
	\Psi = \mathcal{L}_{v(t)}\mathcal{R}_{h(t)}^{\alpha(t),\beta(t)}\mathcal{S}_{c(t)}\left(\Phi + \sum_{j=1,2}\varphi_j E_j + \frak{q}\right),
\end{align*}
where $\frak{q}$ is a multiple of $\Phi$, of quadratic size in $(\varphi_1,\varphi_2)$ (provided the latter are small). 
\\
In section~\ref{sec:perturb formulation}, in a first approximation, we ignore the modulations, and derive the leading order wave equations for the $\varphi_j$. For this, we decompose these functions into discrete Fourier series in relation to the angular parameter $\theta$, i.e., we write $\varphi_j = \sum_{n\in Z}\varphi_j(n)e^{i n\theta}$. In order to extract wave operators whose spatial part is time independent, we change variables to $\tau = \int_t^\infty \lambda(s)\,ds$, $R = \lambda(t)\cdot r$, $\lambda(t) = t^{-1-\nu}$. Furthermore, we introduce the variables 
\[
\varepsilon_{\pm}(n) = \varphi_1(n) \mp i\varphi_2(n),\quad n\in \bbZ, 
\]
which satisfy wave equations with elliptic part given by the family of operators $H_n^{\pm}$ in \eqref{def: Hn pm}. 
This is then the family of operators whose spectral theory is analyzed asymptotically for large values of $n$ (the small values can be studied by more standard methods) in the linear part of the paper.

	\subsection{The linear part} This encompasses sections \ref{subsec:linear start} to \ref{subsec:linear end}. 
	The operators governing the linearized flow, restricted to a fixed angular frequency $n$,  are given by 
	\begin{align}\label{def: Hn pm} 
		H^{+}_{n}:=\partial_{R}^{2}+\frac{1}{R}\partial_{R}-f_{n}(R)+g_{n}(R),\quad H^{-}_{n}:=\partial_{R}^{2}+\frac{1}{R}\partial_{R}-f_{n}(R)-g_{n}(R).
	\end{align}
	where
	\begin{align*}
		f_{n}(R) &=\frac{n^{2}+1}{R^{2}}-\frac{8}{(R^{2}+1)^{2}},\quad g_{n}(R)=\frac{2n}{R^{2}}-\frac{4n}{R^{2}+1},\\
		-f_{n}(R)+g_{n}(R) &=-\frac{(n-1)^{2}}{R^{2}}-\frac{4n}{R^{2}+1}+\frac{8}{(R^{2}+1)^{2}},\\
		-f_{n}(R)-g_{n}(R) &=-\frac{(n+1)^{2}}{R^{2}}+\frac{4n}{R^{2}+1}+\frac{8}{(R^{2}+1)^{2}}
	\end{align*}
	These operators are symmetric in $L^2(R\, dR)$ with domain consisting of $C^2((0,\infty))$ functions of compact support in $(0,\infty)$. 
	Conjugating $H^{\pm}_n$ by the weights~$R^{\frac12}$ yields the  Schr\"odinger operators
	\begin{equation}
		\nn 
		\begin{split}
			\calH^{+}_{n}=R^{\frac12}  H^{+}_{n} R^{-\frac12} = &\partial^{2}_{R}+\frac{1}{4R^{2}}-\frac{(n-1)^{2}}{R^{2}}-\frac{4n}{R^{2}+1}+\frac{8}{(R^{2}+1)^{2}},\\
			\calH^{-}_{n}=R^{\frac12}  H^{-}_{n} R^{-\frac12}  =&\partial^{2}_{R}+\frac{1}{4R^{2}}-\frac{(n+1)^{2}}{R^{2}}+\frac{4n}{R^{2}+1}+\frac{8}{(R^{2}+1)^{2}}.
		\end{split}
	\end{equation}
	which are symmetric in $L^2(dR)$ over the same domain. The potentials here are strongly singular, meaning that they are not integrable on $(0,1)$. The nonlinear analysis will be based on the (distorted) Fourier transform associated with these operators, see Gesztesy, Zinchenko~\cite{GZ} for the general theory of this Fourier transform,  as well as Section~\ref{sec:DFT} below.  All operators $-\calH^{\pm}_{n}$ are nonnegative, and for $|n|\ge2$ they do not possess  any $0$ energy modes. However, $\calH^{+}_{n}=\calH^{-}_{-n}$ for $n=0,\pm1$ each exhibit a $0$ energy state, which is either an eigenvalue or resonance. In Section~\ref{subsec:modulation} below we relate these $0$ modes to six modes of the matrix operator $\frakL$ of~\eqref{real linear operator} which are generated by all available symmetries: translations of $\R^2$, Lorentz transforms of Minkowski space $\R^{1+2}_{t,x}$, dilations, and rotations of the sphere~$\bbS^2$. This is an $8$-parameter family, but translations and Lorentz transforms lead to the same modes.
	This is an essential feature of  our main theorem since translations do not appear in the asymptotic  description of the blowup solutions. 
	Even though separating variables  might seem like a heavy-handed approach, it appears necessary to proceed in this fashion since 
	the physical symmetries do not affect any frequencies $|n|\ge2$ in the angle.  
	
	In order to effectively diagonalize each of the  linearized operators (for fixed angular frequency) by means of its associated distorted Fourier transform, we need to gain precise control on the solutions of the equations $-\calH^{\pm}_{n}\phi = E^2\phi$, in the entire range of all three variables $n,E,R$. We treat $|n|\le N$ separately and by different methods from the regime $|n|\ge N$, where $N$ is some large absolute constant.  Small $|n|$ are handled by the perturbative techniques of~\cite{KST} and \cite{SSS1, SSS2}. The large $|n|$ regime, which is intrinsically a singular perturbation problem, is reduced to a semi-classical spectral problem $\left(-\hbar^2\partial_R^2 + V(R;\hbar)\right)f=(\hbar E)^2 f$ with $\hbar = (n+1)^{-1}$ and $V(R;\hbar)=R^{-2}+\hbar W(R;\hbar)$. A change of variables $x=\hbar E R$ moves the turning point $V(R;\hbar)-(\hbar E)^2=0$ to the vicinity of $x=1$.  As in~\cite{CSST}, we need to carry out  the so-called Langer correction by adding $\hbar^2/(4x^2)$ to the potential.  Otherwise our WKB type approximations will diverge as $E\to0$ (see~\cite{CSST} for more on this).  The resulting problem is of the form, see \eqref{Schlag eq alt} and~\eqref{def Q0}, 
	\[
	-\hbar^{2}\tf''(x)+Q_0(x)\tf(x)=0,\quad Q_0(x):=x^{-2}-1+\text{corrections}.
	\]
	To this we apply the global Liouville-Green transform~\eqref{def tau}, viz.
	\[
	\tau(x,\alpha;\hbar):=\sign\left(x-x_{t}(\alpha;\hbar)\right)\left|\frac{3}{2}\int_{x_{t}(\alpha;,\hbar)}^{x}\sqrt{\left|Q_{0}(u,\alpha;\hbar)\right|}\, du\right|^{\frac{2}{3}}
	\]
	also known as Langer transform, which turns the previous equation into a perturbed Airy equation in the variable~$\tau\in(-\infty,\infty)$, see~\eqref{perturbed Airy}. 
	In Section~\ref{sec:fund sys} we solve this perturbed Airy equation so as to obtain a fundamental system in~$\tau$. Our strategy is inspired by the one of~\cite{CSST, CDST}. It is absolutely essential for us that our approximation is accurate up to a {\em multiplicative correction} of the form $1+\hbar a(\tau;E,\hbar)$ where $a$ is uniformly bounded in the entire three-parameter regime. In addition, we control the derivatives in~$E$ and $\tau$. In particular, the functions $a$ cannot contain any rapidly oscillating factors, and all oscillations have to be included in the main term. This appears to be in stark contrast to the vast body of the semi-classical literature in which the errors are additive, or only the limit $\hbar\to0$ is considered. Anything but the asymptotic descriptions of the fundamental system and the spectral measure as we derive them below uniformly in all three parameters, would be insufficient for our nonlinear analysis.  One surprising feature of our analysis, which does not appear in~\cite{CSST, CDST} is a parameter regime in which $a(\tau;E,\hbar)$ can form arbitrarily long plateaus without any decay, see Lemma~\ref{lem:a0fine2}.

\subsection{The nonlinear part} Roughly speaking, the nonlinear part is divided into three parts. In the first part, section~\ref{sec:bilinlargensmooth}, we develop the basic bilinear and multilinear estimates allowing us to control the source terms {\it{away from the light cone}}, i.e.,  in the region where all source terms have a regularity like that of the perturbation of the initial data. Here the main difficulty consists in proving estimates which experience no losses in the angular frequency, and with respect to a space which is consistent with the wave parametrix and the {\it{transference operator}} which comes up in the solution of the linear homogeneous problem associated to $H_n^{\pm}$. We shall refer this space as the {\it{smooth space}}, since elements in it have a high degree of differentiability, but on the flip side they have little structure beyond the smoothness. 
In terms of the distorted Fourier variables $\xb(\tau, \xi)$, these norms are of the form 
\[
\left\|\left(\hbar^2\xi\right)^{1-\frac{\delta}{2}}\left\langle \hbar^2\xi\right\rangle^{\delta + \frac32}\xb(\tau, \xi)\right\|_{L^2_{d\xi}},
\]
where $\hbar = \frac{1}{|n|+1}$ is in terms of the angular frequency $n$, and $\left\langle x\right\rangle = \sqrt{1+x^2}$. It is important that these norms include the weight $\left(\hbar^2\xi\right)^{1-\frac{\delta}{2}}$ for small frequencies, and the  different weight $\left(\hbar^2\xi\right)^{\frac52 + \frac{\delta}{2}}$ for very large frequencies. The fact that one uses the variable $\hbar^2\xi$ rather than $\xi$ comes from the fine structure of the eigenfunctions and the transference operator at angular frequency $n$, and more specifically the fact that the separation point between the non-oscillatory and oscillatory regimes (so-called turning point) is defined in terms of $R\xi^{\frac12}\hbar$. 
We observe that other than some technical but non-essential issues arising near the spatial origin $R = 0$, the small angular frequencies and in particular the exceptional ones can be handled very similarly to the large ones in this setting. \\
In the second part, covering section~\ref{sec:multilinestimatesnearcone}, section~\ref{sec:ngeq2close}, we develop a functional framework which encompasses the shock singularity along the light cone, and which allows us to handle the null-form nonlinear source terms, as well as the remaining terms, for angular frequencies $|n|\geq 2$. Due to the complexity of the function space, we do this in several stages, starting with a basic {\em prototype singular term}  (Definition~\ref{defi:xsingulartermsngeq2proto} ), which gets refined to the more general {\em admissible singular function} concept in Definition~\ref{defi:xsingulartermsngeq2adm}, and this gets combined with the {\em smooth space} concept from before in the final Definition~\ref{defi:goodfourierrepngeq2} for angular frequencies $|n|\geq 2$. This gets complemented by an analogous definition for the exceptional modes, which also includes the desired bounds for the unstable modes. 
Completing the required bounds for the angular frequencies $|n|\geq 2$ then amounts to proving estimates for all the source terms as in Prop.~\ref{prop:ngeq2finalsourcetermestimatesingoodspaces}, as well as solving the inhomogeneous wave equation in its Fourier realization. The latter is accomplished in section~\ref{sec:ngeq2close}, where the same method as in~\cite{KST2} is used to iterate away the linear non-local source terms due to the transference operator. The reiteration method used here relies crucially on the asymptotic bounds for the transference operator kernel established in section \ref{sec:transference}.  
\\
It then remains to handle the exceptional modes. Here we trade in the difficulty involved in handling large angular frequencies (these only occur in pairs in the corresponding source terms, meaning they become harmless due to the rapid decay in terms of the angular frequency) for the difficulty controlling the unstable part. Precisely, it is here we have to take advantage of the freedom to continuously modulate in the symmetries in order to neutralize certain top order singular terms, the details being carried out in section \ref{sec:exceptional modes}.
Finally, the iteration scheme leading to the desired solution is outlined in the section \ref{sec:final}. 

\section{Description of the key components of the strategy}\label{sec:strategydetails}

\subsection{The distinction between the angular frequencies $|n|\geq 2$ and the exceptional ones $n\in \{0, \pm 1\}$.}

In order to control the perturbation, we shall resort to a decomposition into angular Fourier modes, \eqref{Fourier expansion}, which upon passage to the variables $\veps_{\pm}: = \varphi_1\mp i\varphi$ and their angular Fourier modes leads to the infinite system of coupled equations \eqref{eq diag phys}. 
In order to control each such angular mode $\veps_{\pm}(n)$, which is a function of the rescaled variables $\tau, R$, see \eqref{tau R variable def}, we use distorted Fourier representations \eqref{eq:hatfn}. The corresponding wave operators at angular frequency $n$ and expressed on the distorted Fourier side are then given by \eqref{eq tau xi hbar D}, \eqref{def calD}; the important {\it{transference operator}} $\mathcal{K}_{\hbar}^{(0)}$ appearing there will be analyzed carefully in Prop.~\ref{prop: K operator}, Prop.~\ref{prop: K operator negative} for the ``generic'' angular momenta $|n|\geq 2$, and separately for the exceptional angular momenta and their corresponding supersymmetric elliptic parts. \\
The leading order wave propagators at angular momentum $n$ and expressed on the distorted Fourier side are then given by \eqref{parametrix hbar homo}, \eqref{parametrix hbar inhomo}. In order to control the nonlinear interactions, it will be crucial to deduce {\it{estimates with good enough temporal decay}} for these wave propagators. This decay will not be the result of dispersive type estimates which are common in problems involving scattering toward temporal infinity, but rather they result from {\it{suitable weights}} in the norms and the scaling factors in the wave propagator. These norms, which will be naturally expressed in terms of the distorted Fourier coefficients, need to also lead to sufficient control of the underlying physical functions. 
\\

Here an important distinction appears between the angular frequencies $|n|\geq 2$ on the one, and the exceptional ones $n\in \{0, \pm 1\}$ on the other hand. To wit, the products
\begin{equation}\label{eq:key product}
	\sqrt{\rho_n(\xi)}\cdot \phi_n(R,\xi), 
\end{equation}
which are canonical in the sense that they correspond to a Fourier basis with corresponding spectral measure $1$, have the property that they vanish at $\xi = 0$ precisely if $|n|\geq 2$. It is this feature which allows us to introduce the norms \eqref{def S0 hbar} when $|n|\geq 2$, which do enjoy temporal decay upon application of the principal wave propagator as witnessed by Prop.~\ref{prop: homo para decay hbar}, and at the same time lead to good weighted $L^\infty$-bounds for the ``physical realization'' of a function as in Lemma~\ref{lem:derLinfty}.  
\\

Concerning the exceptional modes $n\in \{0, \pm 1\}$, there is a trick to remedy the failure of the vanishing of \eqref{eq:key product} at $\xi = 0$, namely the passage to the supersymmetric partner operator. This procedure, which is explained in subsection~\ref{subsec:exceptionalnsupersymmetry}, eliminates the root mode/resonance at $\xi = 0$, and leads to a spectral measure $\rho(\xi)$ with an extra factor $\xi$. This latter observation, which appears first in \cite{BKT}, and was also used in \cite{KMiao}, allows one to introduce analogous norms for the continuous spectral part of the exceptional modes and such that one has good control over the non-root/resonant part in \eqref{eq:levelup-rep}, \eqref{eq:levelup1rep} and \eqref{eq:levelup0rep}. The key issue for controlling the exceptional angular momenta then resides in controlling the resonant/root parts $c_n\cdot \phi_n(R)$, $n = 0, \pm 1$, which is where the modulation theory will play a crucial role. We observe that the coefficients $c_n$ solve the ordinary differential equations \eqref{eq:c+evolutioneqn}, \eqref{eq:c0evolutioneqn} and \eqref{eq:c-1evolutioneqn}.

\subsection{Controlling the unstructured smooth part and the structured singular/shock part}

The norm \eqref{def S0 hbar} and the related one $\big\|\cdot\big\|_{S_1^{\hbar}}$ used to control the distorted Fourier transform of derivative terms, all for angular momenta $|n|\geq 2$, and their analogues for the exceptional angular momenta defined in subsection~\ref{subsec:nexcspaces}, are sufficient to control the perturbation {\it{away from the light cone}}, and more precisely in the region $r\leq \frac{|t|}{2}$. Near the light cone, we need to introduce a function space for the distorted Fourier transforms of the various angular momentum components $\veps_{\pm}(n)$ and which reflect the symbolic expansion of the shock of the unperturbed solution $\Phi$, and more specifically the function $U$ describing it. 
\\
Precisely, we shall want this function space (which needs to be specified for each angular frequency $n$!), characterizing the distorted Fourier coefficients, to have the following properties:
\begin{itemize}
	\item It contains the functions arising as distorted Fourier transforms of the shock localized around the light cone, at fixed time. More specifically, the shock is described by linear combinations of the functions $(t-r)^{\frac12+k\nu}\cdot \log^l(t-r)$ with suitable time dependent factors. 
	\item For any of the multilinear source terms, assuming the factors (at various angular frequencies $n_j$) have distorted Fourier transforms in this space, then the wave propagator applied to such a source term leads to a function whose Fourier transform is again in this space, up to better errors in $S_0^{\hbar}$. 
	\item The complicated nature of the wave propagator on the Fourier side, which in particular requires application of the transference operator infinitely many times (see the proof of Lemma~\ref{lem:offdiagonaltransferencegainngeq2}), necessitates preservation of the space under the transference operator (of the same angular frequency). 
	\item The principal singular part needs to transform in a simple manner when applying any of the six symmetries which preserve the light cone. For this last aspect, we in fact introduce a somewhat more special function space, which we call of `restricted type'. 
\end{itemize}

We note that while there is a uniform way to define the space for angular momenta $|n|\geq 2$, we need to modify the definition specifically for each exceptional mode $n\in \{0, \pm 1\}$. 
\\

The final concept of ``good functions'', which for angular momenta $|n|\geq 2$ is introduced in Definition~\ref{defi:xsingulartermsngeq2adm} and Definition~\ref{defi:goodfourierrepngeq2}, results from the need to address each of the above points. More precisely, we first introduce the concept of prototypical singular part in Definition~\ref{defi:xsingulartermsngeq2proto}, which naturally arises when expressing functions involving a shock which are a linear combination of the functions $(t-r)^{\frac12+k\nu}\cdot \log^l(t-r)$ (with suitable weights depending on $t$) on the distorted Fourier side, as detailed in Lemma~\ref{lem:singPhysicaltoFourierngeq2}. This class of functions is unfortunately not compatible with the action of the wave propagator \eqref{parametrix hbar homo}, \eqref{parametrix hbar inhomo}, which forces us to introduce a more sophisticated concept of `admissibly singular functions' to describe the distorted Fourier transform of the $\veps_{\pm}(n)$ in Definition~\ref{defi:xsingulartermsngeq2adm}. The source terms will be essentially describable in terms of prototypical expansions. An important point here is the fact that we are free to modify source terms outside the light cone, which facilitates the transition from functions in physical space to functions on the Fourier side. \\
We then need to show that the action of the transference operator is compatible with this concept of function, which is accomplished in Proposition~\ref{prop:transferenceonsingularngeq2} (in case of angular frequencies $|n|\geq 2$). Since the function space is characterised in part by sufficient smoothness of certain coefficient functions with symbol type bounds, we need high order `diagonal type' derivative bounds for the kernel of the transference operator to preserve this structure. This is reflected by the delicate estimate\footnote{A key point of this estimate is that it is uniform in $\hbar$.} \eqref{eq:trace der}, whose proof relies on the refined asymptotics for the Fourier basis $\phi_n(R,\xi)$ in Prop.~\ref{prop:DFT nlarge}.

Furthermore, we need to show that all the source terms can be placed into these function spaces provided all the factors are in such a space, which is eventually achieved in Prop.~\ref{prop:preliminarytrilinearnullformgeneralinputsngeq2} (again for angular frequencies $|n|\geq 2$, with a suitable analogue for the exceptional angular modes), up to re-iterating the equations arbitrarily many times. This does not cause difficulties since we shall construct the solutions at the end via an iterative process, and re-iterating the equation then means going to the level of an earlier iterate, whence in fact there are only finitely many re-iterations of the equation then. 

\subsection{Controlling the evolution of the coefficients $c_n$, $n = 0,\,\pm 1$. Modulation theory}

The chief difficulty for controlling the evolution of the root mode coefficient as $c_+$, whose evolution equation is given by \eqref{eq:c+evolutioneqn}, consists in bounding the contribution of the source term $\lim_{R\rightarrow 0}H_1^{+}\varepsilon_+^1$, and analogously for the other evolution equations \eqref{eq:c0evolutioneqn} and \eqref{eq:c-1evolutioneqn}. Taking advantage of the formula \eqref{eq:hardsourceforcplus} for this term, and substituting the incoming principal singular part (in the sense of Definition~\ref{defi:xsingulartermsnless2proto}), one arrives at the formula \eqref{eq:c_+sourcetermobstruction}. Here the sum on the right represents the obstruction to obtaining decaying bounds for $c_+$, and is precisely due to the contribution of the principal incoming part of the shock on the light cone, which affects points on the symmetry axis close to the singularity. 
\\
This obstruction disappears if one imposes a finite number of vanishing conditions \eqref{eq:n=1vanishingrelations}. We note that analogous considerations apply to the evolution equations for the other instabilities. 
\\
In section~\ref{subsec:modulation}, the effect of the modulations, and more specifically, the {\it{primary modulations}} on the singular part is studied, and it is shown that suitable choice of the parameters $\alpha(t), \beta(t), h(t), c(t),\nu_1(t), \nu_2(t)$ allows one to force the vanishing conditions from before. This is exemplified by Lemma~\ref{lem:alphamodulationn=1singulartermspropagated}, which explains how to choose the parameter $\alpha$, with similar lemmas applying to the other parameters. 
\\
Unfortunately, for technical reasons it appears that even after the primary modulations which partly `de-singularize' the effect of the shock, one still does not have quite enough decay for the parameter $c_0(\tau)$ controlling the $n=0$ resonance to close the estimates, and a further {\it{secondary modulation step}}, which however, does not affect the shock part anymore (and neither the terminal values for $\alpha, \beta$ etc), seems necessary. Specifically, we need to modify the bulk part $\Phi$ as in \eqref{eq:secondarymod}. This allows us to effectively completely eliminate $c_0(\tau)$. 

\section{Perturbing equivariant blowup solutions}\label{sec:perturb formulation}
We consider   critical wave maps from $\bbR^{1+2}_{t,x} \to \bbS^{2}$, not necessarily equivariant. We recall the extrinsic description of wave maps of this type.  A smooth map $u:=(u^{1},u^{2},u^{3}):\bbR^{1+2}\rightarrow \bbS^{2}\hookrightarrow\bbR^{3}$ is a wave map if at every point $z=(t,x)\in\bbR^{1+2}$  
\begin{align}\label{orthogonal condition}
 \Box u\perp T_{u(z)}\bbS^{2},\quad \Box u(z)=\left(\Box u^{1}(z),\Box u^{2}(z),\Box u^{3}(z)\right).
\end{align}
Here $\Box$ denotes the operator $\Box:=\partial^{\alpha}\partial_{\alpha}=-\partial_{t}^{2}+\partial_{1}^{2}+\partial_{2}^{2}$. 
Formally speaking, \eqref{orthogonal condition} is the Euler-Lagrange equation of critical points of the action
\[
\calL(\partial u) = \frac12 \int_{\R^{1+2}_{t,x}}  \left(-  |u_t|_g^2  + |u_{x_1}|_g^2 + |u_{x_2}|_g^2\right)\, dtdx
\]
and $|\cdot|_g$ is the metric on the target manifold. 
For $\bbS^{2}$, note that $u$ is a normal vector field, therefore there is a scalar function $\lambda: \bbR^{1+2}\rightarrow\bbR$ such that
\begin{align}\label{orthogonal S2}
 \Box u=\lambda u.
\end{align}
On the other hand, $u$ takes values on $\bbS^{2}$, which implies $\langle u,u\rangle=1$. Here $\langle\cdot\rangle$ denotes the Euclidean inner product in $\bbR^{3}$. Hence 
\begin{align*}
 \lambda=\langle\Box u,u\rangle=\partial^{\alpha}\langle\partial_{\alpha}u,u\rangle-\langle\partial_{\alpha}u,\partial^{\alpha}u\rangle=-\langle\partial_{\alpha}u,\partial^{\alpha}u\rangle=|\partial_{t}u|^{2}-|\nabla u|^{2}.
\end{align*}
Therefore the extrinsic description of a wave map from $\bbR^{1+2}$ to $\bbS^{2}$ is  
\begin{align}\label{extrinsic eq}
 \Box u=\left(|\partial_{t}u|^{2}-|\nabla u|^{2}\right)u,\quad \textrm{or}\quad S(u):=-\partial^{2}_{t}u+\Delta u+\left(|\nabla u|^{2}-|\partial_{t}u|^{2}\right)u=0.
\end{align}
For more general surfaces than $\bbS^2$ the nonlinearity is replaced by the second fundamental form contracted with the Minkowski tensor, i.e., $A(\partial_\alpha,\partial_\beta)\eta^{\alpha\beta}$. 
As a wave equation,  \eqref{extrinsic eq}   conserves   energy. 

\begin{proposition}\label{prop: L2 conservation}
 Let $u:\bbR^{2}\times[0,T) \rightarrow \bbS^{2}\hookrightarrow\bbR^{3}$ be a smooth solution to \eqref{extrinsic eq} such that $(u(0,\cdot), u_t(0,\cdot))\in \dot H^1(\R^2)\times L^2(\R^2)$. Then the energy
 \begin{align}\label{L2 energy}
  E(t):=\frac{1}{2}\int_{\bbR^{3}}\left(|\partial_{t}u|^{2}+|\nabla u|^{2}\right)\, dx
 \end{align}
is constant for all times $t\in[0,T)$.
\end{proposition}
\begin{proof}
The proof proceeds by inner product with $\partial_{t}u$ on both sides of \eqref{extrinsic eq}, and integrating by parts. The right-hand side contributes nothing since 
\begin{align*}
 \left(|\partial_{t}u|^{2}-|\nabla u|^{2}\right)u\cdot\partial_{t}u=\frac{1}{2} \left(|\partial_{t}u|^{2}-|\nabla u|^{2}\right)\partial_{t}|u|^{2}=0,
\end{align*}
due to $|u|\equiv1$.
\end{proof}

Let $\Phi(t,r,\theta)$ be a smooth $1$-equivariant  wave map written in spherical coordinates. 
\begin{align}\label{background wave map}
 \Phi(t,r,\theta):=(\cos\theta\sin U(t,r),\sin\theta\sin U(t,r), \cos U(t,r)).
\end{align}
We introduce the orthonormal frame
\EQ{ 
E_1 & =\partial_U = (\cos\theta\cos U,\sin\theta\cos U,-\sin U) \\
E_2 & = (\sin U)^{-1}  \partial_\theta = (-\sin\theta,\cos\theta,0) 
}
on the sphere (treating the azimuth angle $U$ as an independent variable). We write non-equivariant  perturbations of $\Phi$ in the form, cf.~\cite{DavilaDelPWei},  
\begin{align}\label{perturbed sol}
 \Psi:=\Phi+\pphi+\mathfrak{q}
\end{align}
We expand the tangential vector field  $\pphi$  in the form 
\begin{align}\label{pphi expand}
\begin{split}
 \pphi(t,r,\theta):=&\vphi_{1}(t,r,\theta)(\cos\theta\cos U,\sin\theta\cos U,-\sin U)+\vphi_{2}(t,r,\theta)(-\sin\theta,\cos\theta,0)\\
 :=&\vphi_{1}(t,r,\theta)E_{1}+\vphi_{2}(t,r,\theta)E_{2}.
 \end{split}
\end{align}
The term $\mathfrak{q}$ in \eqref{perturbed sol} is a correction which insures that $|\Psi|=1$. We will specify it later. It is quadratic in the size of $\pphi$, which we assume small. 
Both $\Phi$ and $\Psi$ satisfy~\eqref{extrinsic eq}. Subtracting them yields 
\EQ{\nn 
& \Box \Psi - \Box \Phi = \Box \pphi + \Box \mathfrak{q} = \left(|\partial_{t}\Psi|^{2}-|\nabla \Psi|^{2}\right)\Psi  - \left(|\partial_{t}\Phi|^{2}-|\nabla \Phi|^{2}\right)\Phi   \\
& = \left(|\partial_{t}\Phi+\partial_{t}\pphi+\partial_{t}\mathfrak{q} |^{2}-|\nabla \Phi+\nabla\pphi+\nabla\mathfrak{q}|^{2}\right)(\Phi+\pphi+\mathfrak{q})  - \left(|\partial_{t}\Phi|^{2}-|\nabla \Phi|^{2}\right)\Phi \\
& = 2 ( \partial_{t}\Phi\cdot \partial_{t}\pphi - \nabla \Phi \cdot \nabla\pphi) \Phi + 2 ( \partial_{t}\Phi\cdot \partial_{t}\mathfrak{q} - \nabla \Phi \cdot \nabla\mathfrak{q}) \Phi \\
& \quad + ( | \partial_{t}\pphi+\partial_{t}\mathfrak{q} |^{2}-| \nabla\pphi+\nabla\mathfrak{q}|^{2} )\Phi +  \left(|\partial_{t}\Phi+\partial_{t}\pphi+\partial_{t}\mathfrak{q} |^{2}-|\nabla \Phi+\nabla\pphi+\nabla\mathfrak{q}|^{2}\right)(\pphi+\mathfrak{q}) 
}
Moving all terms which depend linearly on $\pphi$ to the left yields
\EQ{\label{eq:lin pphi}
&  \Box \pphi+\left(|\nabla\Phi|^{2}-|\Phi_{t}|^{2}\right)\pphi+2\left(\nabla\Phi\cdot\nabla\pphi-\Phi_{t}\cdot(\pphi)_{t}\right)\Phi \\
 & = - \Box \mathfrak{q} + 2 ( \partial_{t}\Phi\cdot \partial_{t}\mathfrak{q} - \nabla \Phi \cdot \nabla\mathfrak{q}) \Phi  + ( | \partial_{t}\pphi+\partial_{t}\mathfrak{q} |^{2}-| \nabla\pphi+\nabla\mathfrak{q}|^{2} )\Phi \\
 & \quad+  \left(|\partial_{t}\Phi+\partial_{t}\pphi+\partial_{t}\mathfrak{q} |^{2} -|\Phi_{t}|^{2} + |\nabla\Phi|^{2} -|\nabla \Phi+\nabla\pphi+\nabla\mathfrak{q}|^{2}\right)\pphi \\
 & \quad + \left(|\partial_{t}\Phi+\partial_{t}\pphi+\partial_{t}\mathfrak{q} |^{2}-|\nabla \Phi+\nabla\pphi+\nabla\mathfrak{q}|^{2}\right)\mathfrak{q} =: \mathfrak{Q}
} 
Here $\mathfrak{Q}$ are terms of second and  higher orders in $\vphi_1,\vphi_2$. 
From the definitions we compute  the following first order derivatives: 
\EQ{\nn
\partial_r E_1 &= -U_r \Phi,\quad \partial_r \Phi = U_r E_1,\quad \partial_r E_2=0 \text{\ \ (same for $\partial_t$)}\\
\partial_\theta E_1 &= \cos U E_2,\quad \partial_\theta E_2 = -\cos U E_1 -\sin U \Phi, \quad \partial_\theta \Phi = \sin U E_2
}
Further, with  $\nabla= \what{e_r} \partial_r + r^{-1}\what{e_\theta}\partial_\theta$  the gradient in polar coordinates: 
\begin{align}
 \nabla\Phi=&  \what{e_r} U_{r}E_{1} + \what{e_\theta}  \frac{\sin U}{r}E_{2},\quad 
             \Phi_{t}= U_{t}E_{1}, \nn \\
             |\nabla\Phi|^{2}-|\Phi_{t}|^{2}=&U^{2}_{r}-U^{2}_{t}+\frac{1}{r^{2}}\sin^{2}U,\nn\\
             (\pphi)_{t} 
             =&\vphi_{1,t}E_{1}-\vphi_{1}U_{t}\Phi+\vphi_{2,t}E_{2},\quad
             \Phi_{t}\cdot(\pphi)_{t}=U_{t}\vphi_{1,t},  \label{pphi 1st deri} \\
             \nabla\pphi =& \big(\vphi_{1,r}E_{1}-\vphi_{1}U_{r}\Phi+\vphi_{2,r}E_{2}\big)\what{e_r} + 
            \frac{1}{r} \big ( {\vphi_{1,\theta}} E_{1}+ \vphi_{1}\cos UE_{2}+ {\vphi_{2,\theta}} E_{2}- {\vphi_{2}}  (\cos U E_{1}+\sin U\Phi)\big) \what{e_\theta},\nn\\
             \nabla\pphi\cdot\nabla\Phi=&\vphi_{1,r}U_{r}+ {\frac{1}{2r^{2}}\sin(2U)\vphi_{1}}+\frac{1}{r^{2}}\sin U\vphi_{2,\theta}.\nn
\end{align}
The d'Alembertian $\Box\pphi$ is obtained as follows (using $\Delta=\div \nabla$ in polar coordinates): 
\begin{align*}
 (\pphi)_{tt}
 &= (\vphi_{1,tt} -\vphi_{1}U_{t}^{2}) E_{1}+\vphi_{2,tt}E_{2}-(2\vphi_{1,t}U_{t}+\vphi_{1}U_{tt})\Phi,\\
 \Delta\pphi &=  \frac{1}{r} \big( r\big(\vphi_{1,r}E_{1}+\vphi_{2,r}E_{2}-\vphi_{1}U_{r}\Phi\big)\big)_r  + \frac{1}{r^2}  \big ( ({\vphi_{1,\theta}} - {\vphi_{2}}  \cos U) E_{1}+ ({\vphi_{2,\theta}}+\vphi_{1}\cos U )E_{2}- {\vphi_{2}}  \sin U\Phi \big)_\theta \\
 & =    \big(\vphi_{1,rr}+\frac{1}{r}\vphi_{1,r}+\frac{1}{r^{2}}\vphi_{1,\theta\theta} -\frac{2\vphi_{2,\theta}}{r^{2}}\cos U -\vphi_{1}U^{2}_{r}  -\frac{\vphi_{1}\cos^2 U}{r^{2}}\big) E_1 \\
 & \qquad +\big(\vphi_{2,rr}+\frac{1}{r}\vphi_{2,r}+\frac{1}{r^{2}}\vphi_{2,\theta\theta} +\frac{2\cos U}{r^{2}}\vphi_{1,\theta} -\frac{\vphi_{2}}{r^{2}}\big)E_2 \\ 
 & \qquad + \big( -2\vphi_{1,r}U_{r}-\vphi_{1}U_{rr}-\frac{1}{r}\vphi_{1}U_{r}-\frac{\sin 2U}{2r^{2}}\vphi_{1}-\frac{2\sin U}{r^{2}}\vphi_{2,\theta}\big)\Phi
 \end{align*}
Now we are ready to plug in all the above formulas into the left hand side of \eqref{eq:lin pphi}. We have:
\begin{align}\label{Box pphi}
\begin{split}
 \Box\pphi=&\left(-\vphi_{1,tt}+\vphi_{1,rr}+\frac{1}{r}\vphi_{1,r}+\frac{1}{r^{2}}\vphi_{1,\theta\theta}+\vphi_{1}(U^{2}_{t}-U^{2}_{r})-\frac{\cos^{2}U}{r^{2}}\vphi_{1}-\frac{2\cos U}{r^{2}}\vphi_{2,\theta}\right)E_{1}\\
 &+\left(-\vphi_{2,tt}+\vphi_{2,rr}+\frac{1}{r}\vphi_{2,r}+\frac{1}{r^{2}}\vphi_{2,\theta\theta}-\frac{1}{r^{2}}\vphi_{2}+\frac{2\cos U}{r^{2}}\vphi_{1,\theta}\right)E_{2}\\
  &+\left(2\vphi_{1,t}U_{t}+\vphi_{1}U_{tt}-2\vphi_{1,r}U_{r}-\vphi_{1}U_{rr}-\frac{1}{r}\vphi_{1}U_{r}-\frac{\sin 2U}{2r^{2}}\vphi_{1}-\frac{2\sin U}{r^{2}}\vphi_{2,\theta}\right)\Phi
 \end{split}
\end{align}
Now we substitute \eqref{pphi 1st deri} and \eqref{Box pphi} into the left-hand side of \eqref{eq:lin pphi}. The coefficient for $E_{1}$ gives
\begin{align}\label{coe eq rough 1}
 \Box\vphi_{1}-\frac{\cos 2U}{r^{2}}\vphi_{1}-\frac{2\cos U}{r^{2}}\vphi_{2,\theta}=\mathfrak{Q}\cdot E_1
 \end{align}
The coefficient for $E_{2}$ is
\begin{align}\label{coe eq rough 2}
 \Box \vphi_{2}-\frac{\cos^{2}U}{r^{2}}\vphi_{2}+(U^{2}_{r}-U^{2}_{t})\vphi_{2}+\frac{2\cos U}{r^{2}}\vphi_{1,\theta}= \mathfrak{Q}\cdot E_2
\end{align}
The coefficient for $\Phi$ on the left-hand side of \eqref{eq:lin pphi} is 
\[
\Big(U_{tt} - U_{rr} - \frac{1}{r} U_r + \frac{\sin(2U)}{2r^2}\Big)\vphi_1 =0
\]
due to the equivariant wave map equation satisfied by $U$.  Therefore, the right-hand side of \eqref{eq:lin pphi} also vanishes, $\mathfrak{Q}\cdot \Phi=0$. 

\subsection{The linear operators governed by the rescaled harmonic map, separation of variables}\label{subsec:linear start}
We let $U$ be the equivariant blowup solution from~\cite{KST,KS_FullRange} which takes the form
\begin{align}\label{equivariant sol}
 U(t,r)=Q(\lambda(t)r)+\eps(t,r)=Q(R)+\eps(t,r)=2\arctan R+\eps(t,r), \qquad   \lambda(t)= t^{-1-\nu},\nu>0
\end{align}
Here $Q$ is the unique $1$-equivariant harmonic map, and the energy of $\eps$ interior to a light cone $0\le r\le t$ vanishes as $t\to0+$: 
\[
\lim_{t\to0+}\int_{r\le t}\big( \eps_t^2 + \eps_r^2\big)(t,r)\, rdr =0
\]
Later we will need much finer properties of $\eps$, but for now we simplify the left-hand sides of~\eqref{coe eq rough 1} and~\eqref{coe eq rough 2} by replacing all occurrences of $U$ with~$Q$. 
The errors thus created will be treated perturbatively in the main nonlinear analysis. 
We introduce the new independent variables
\begin{align}\label{tau R variable def}
 R:=\lambda(t)r=t^{-1-\nu}r,\quad \tau:=\nu^{-1}t^{-\nu}=\int_t^\infty \lambda(s)\, ds.
\end{align}
The spatial Laplacian is given by
\begin{align}\label{spatial laplace}
 \partial_{r}^{2}+\frac{1}{r}\partial_{r}+\frac{1}{r^{2}}\partial^{2}_{\theta}=\lambda^{2}\left(\partial_{R}^{2}+\frac{1}{R}\partial_{R}+\frac{1}{R^{2}}\partial^{2}_{\theta}\right), 
\end{align}
and we have with $Q=Q(R)$
\begin{align}\label{linear coefficients 1}
 \begin{split}
  \cos Q = \frac{1-R^{2}}{1+R^{2}},\quad \cos^{2}Q = \frac{1-2R^{2}+R^{4}}{1+2R^{2}+R^{4}},\quad \cos (2Q) = \frac{R^{4}-6R^{2}+1}{R^{4}+2R^{2}+1}.
 \end{split}
\end{align}
On the other hand, we have (with $\dot\lambda(t)=\frac{d\lambda(t)}{dt}$) 
\begin{align}\label{linear coefficients 2}
\begin{split}
& \partial_{r}\, Q(R)= \frac{2\lambda(t)}{1+R^{2}},\quad (\partial_r \, Q(R))^{2}=\frac{4\lambda(t)^{2}}{R^{4}+2R^{2}+1},\\
&\partial_{t}\, Q(R)  = \frac{2R}{1+R^{2}}\frac{\dot\lambda(t)}{\lambda(t)} ,\quad (\partial_{t}\, Q(R))^{2} = \frac{4R^{2}}{R^{4}+2R^{2}+1}\frac{(1+\nu)^{2}}{t^{2}}=\lambda^{2}\frac{4R^{2}}{R^{4}+2R^{2}+1}\frac{(1+\nu)^{2}}{\nu^{2}\tau^{2}}
 \end{split}
\end{align}
In these new variables and after dividing by $\lambda^2$, we arrive at the new spatial linear operator acting on $\vphi=(\vphi_1,\vphi_2)$, 
\begin{align}\label{linear elliptic operator}
 \begin{split}
  \calL\,\vphi:=\left(\begin{array}{c}
                     \left(\partial_{R}^{2}+\frac{1}{R}\partial_{R}+\frac{1}{R^{2}}\partial_{\theta}^{2}-\frac{R^{4}-6R^{2}+1}{R^{2}(R^{4}+2R^{2}+1)}\right)\vphi_{1}-\frac{2-2R^{2}}{R^{2}(1+R^{2})}\partial_{\theta}\vphi_{2}\\
                     \left(\partial_{R}^{2}+\frac{1}{R}\partial_{R}+\frac{1}{R^{2}}\partial_{\theta}^{2}-\frac{R^{4}-6R^{2}+1}{R^{2}(R^{4}+2R^{2}+1)}-\frac{(1+\nu)^{2}}{\nu^{2}\tau^{2}}\frac{4R^{2}}{1+2R^{2}+R^{4}}\right)\vphi_{2}+\frac{2-2R^{2}}{R^{2}(1+R^{2})}\partial_{\theta}\vphi_{1}
                    \end{array}\right).
 \end{split}
\end{align}
The term $-\frac{(1+\nu)^{2}}{\nu^{2}\tau^{2}}\frac{4R^{2}}{1+2R^{2}+R^{4}}\vphi_{2}$ derives from the term $-U_{t}^{2}\vphi_{2}$ in~\eqref{coe eq rough 2}, and  would be absent if $\Phi$ were a harmonic map in the variable $R$. On the other hand, due to the decay factor $\tau^{-2}$ we will be able to place this term on the right-hand side of the equation and treat it perturbatively. Therefore,  the linear operator driving the entire 
analysis is
\begin{align}\label{real linear operator}
 \begin{split}
  \frakL\, \vphi:=\left(\begin{array}{c}
                     \left(\partial_{R}^{2}+\frac{1}{R}\partial_{R}+\frac{1}{R^{2}}\partial_{\theta}^{2}-\frac{R^{4}-6R^{2}+1}{R^{2}(R^{4}+2R^{2}+1)}\right)\vphi_{1}-\frac{2-2R^{2}}{R^{2}(1+R^{2})}\partial_{\theta}\vphi_{2}\\
                     \left(\partial_{R}^{2}+\frac{1}{R}\partial_{R}+\frac{1}{R^{2}}\partial_{\theta}^{2}-\frac{R^{4}-6R^{2}+1}{R^{2}(R^{4}+2R^{2}+1)}\right)\vphi_{2}+\frac{2-2R^{2}}{R^{2}(1+R^{2})}\partial_{\theta}\vphi_{1}
                    \end{array}\right).
 \end{split}
\end{align}
The change of variables in time is given by
\EQ{\nn 
\partial_t &= \lambda(\partial_\tau + \beta R\partial_R),\qquad \beta(\tau)=\frac{\lambda'(\tau)}{\lambda(\tau)}=(1+1/\nu)\tau^{-1},\quad \lambda'=\frac{d\lambda}{d\tau} \\
\partial_t^2 &= \lambda^2\big[ (\partial_\tau + \beta R\partial_R)^2 + \beta (\partial_\tau + \beta R\partial_R)\big] 
}
whence the full dynamical problem in $(\tau,R)$ is of the form 
\EQ{\label{eq:lin dyn}
- \Big[  (\partial_\tau + \beta R\partial_R)^2 + \beta (\partial_\tau + \beta R\partial_R) \Big] \vphi +  \frakL\,\vphi  = \lambda^{-2} \binom{ \mathfrak{Q}\cdot E_1}{\mathfrak{Q}\cdot E_2} +\frakE + \binom{0}{\frac{(1+\nu)^{2}}{\nu^{2}\tau^{2}}\frac{4R^{2}}{1+2R^{2}+R^{4}}\vphi_{2}}
}
where  $\frakE$ is the error term obtained by replacing $U$ with $Q$ in \eqref{coe eq rough 1} and~\eqref{coe eq rough 2}. 
The precise form of the right-hand side of~\eqref{eq:lin dyn} will be determined in Section~\ref{sec:nonlin}.  To analyze  the linear dynamics,  we (formally) expand $\vphi_{1}$ and $\vphi_{2}$  into Fourier series
\begin{align}\label{Fourier expansion}
\begin{split}
 \vphi_{1}(t,R,\theta)=&\sum_{n}\hvphi_{1}(n,t,R)e^{in\theta},\quad \hvphi_{1}(n,t,R):=\int_{0}^{2\pi}\vphi_{1}(t,R,\theta)e^{- in\theta}\, \frac{d\theta}{2\pi}, \quad n\in\bbZ,\\
 \vphi_{2}(t,R,\theta)=&\sum_{n}\hvphi_{2}(n,t,R)e^{ in\theta},\quad \hvphi_{2}(n,t,R):=\int_{0}^{2\pi}\vphi_{2}(t,R,\theta)e^{- in\theta}\, \frac{d\theta}{2\pi}, \quad n\in\bbZ.
 \end{split}
\end{align}
For fixed $n\in\bbZ$, the operator applied on the Fourier coefficients is given by
\begin{align}\label{operator on Fourier}
 \begin{split}
  \frakL_{n}\,\hvphi(n):=\left(\begin{array}{c}
                     \left(\partial_{R}^{2}+\frac{1}{R}\partial_{R}-\frac{n^{2}}{R^{2}}-\frac{R^{4}-6R^{2}+1}{R^{2}(R^{4}+2R^{2}+1)}\right)\hvphi_{1}(n)- in\frac{2-2R^{2}}{R^{2}(1+R^{2})}\hvphi_{2}(n)\\
                     \left(\partial_{R}^{2}+\frac{1}{R}\partial_{R}-\frac{n^{2}}{R^{2}}-\frac{R^{4}-6R^{2}+1}{R^{2}(R^{4}+2R^{2}+1)}\right)\hvphi_{2}(n)+ in\frac{2-2R^{2}}{R^{2}(1+R^{2})}\hvphi_{1}(n)
                    \end{array}\right).
 \end{split}
\end{align}
We can write this operator as a $2\times 2$ matrix
\begin{align}
\label{operator Fourier matrix}
\begin{split}
 \frakL_{n}\,\hvphi(n):=&\frakA_{n}\;\left(\begin{array}{c}
                                                \hvphi_{1}(n,t,R)\\\hvphi_{2}(n,t,R)
                                               \end{array}\right)\\
\frakA_{n}:=&                                               \left(\begin{array}{cc}
                        \partial_{R}^{2}+\frac{1}{R}\partial_{R}-\frac{n^{2}}{R^{2}}-\frac{R^{4}-6R^{2}+1}{R^{2}(R^{4}+2R^{2}+1)}&- in\frac{2-2R^{2}}{R^{2}(1+R^{2})}\\
                         in\frac{2-2R^{2}}{R^{2}(1+R^{2})}&\partial_{R}^{2}+\frac{1}{R}\partial_{R}-\frac{n^{2}}{R^{2}}-\frac{R^{4}-6R^{2}+1}{R^{2}(R^{4}+2R^{2}+1)}
                       \end{array}\right)
\end{split}
                       \end{align}
\subsection{Symmetries and zero modes}

We now exhibit nonzero solutions of $\frakL\,\vphi=0$ in relation to the symmetries of the nonlinear problem.  First, the dilation symmetry $s\mapsto Q(e^s R)$ yields 
$$\partial_{s=0}\, Q(e^s R) = \frac{R}{1+R^2}=:\psi(R)$$
which leads to two modes 
\EQ{\label{eq:dilmodes}
\frakL \binom{\psi}{0} = \binom{0}{0}, \qquad \frakL \binom{0}{\psi} = \binom{0}{0}
}
This is the resonance mode of the linearized equivariant model~\cite{KST,KMiao}, i.e., of $\frakA_{0}$. 

To determine the zero modes obtained from space translations,  we first write the ground state harmonic map $\Xi(R,\theta)$ in the form
\EQ{\label{Phi in R theta}
 \Xi(R,\theta) &= \big( \cos\theta \sin Q, \sin\theta \sin Q, \cos Q  \big) \\
 &= \big(\cos\theta\frac{2R}{1+R^{2}},\sin\theta\frac{2R}{1+R^{2}},\frac{1-R^{2}}{1+R^{2}}\big).
}
The two translations are given by
\begin{align}\label{space translations}
 \partial_{1}=\cos\theta\, \partial_{R}-\frac{\sin\theta}{R}\partial_{\theta},\quad \partial_{2}=\sin\theta\, \partial_{R}+\frac{\cos\theta}{R}\partial_{\theta}.
\end{align}
We have the following derivatives on $\Xi$ (setting $Q(R)=2\arctan R$)
\begin{align}\label{Phi deri}
 \begin{split}
  \partial_{R}\Xi=&\left(\cos\theta\frac{2-2R^{2}}{(1+R^{2})^{2}},\sin\theta\frac{2-2R^{2}}{(1+R^{2})^{2}},-\frac{4R}{(1+R^{2})^{2}}\right)\\
  =&\frac{2}{1+R^{2}}\left(\cos\theta\cos Q,\sin\theta\cos Q,-\sin Q\right)=\frac{2}{1+R^{2}}E_{1},\\
  \partial_{\theta}\Xi=&\left(-\sin\theta\frac{2R}{1+R^{2}},\cos\theta\frac{2R}{1+R^{2}},0\right)= \frac{2R}{1+R^{2}} E_2
 \end{split}
\end{align}
and thus
\begin{align}\label{translations on Phi}
 \begin{split}
  \partial_{1}\Xi=&\cos\theta\, \partial_{R}\Xi-\frac{\sin\theta}{R}\partial_{\theta}\Xi
  =\frac{2}{1+R^{2}}\cos\theta E_{1}-\sin\theta\frac{2}{1+R^{2}}E_{2},\\
  \partial_{2}\Xi=&\sin\theta\, \partial_{R}\Xi+\frac{\cos\theta}{R}\partial_{\theta}\Xi
  =\frac{2}{1+R^{2}}\sin\theta E_{1}+\frac{2}{1+R^{2}}\cos\theta E_{2}.
 \end{split}
\end{align}
Therefore, two linearly independent solutions of $\calL \vphi=0$ corresponding to all space translations are given by
\begin{align}\label{null: space trans}
\begin{split}
 \left(\vphi_{1}(R),\vphi_{2}(R)\right)=&\left(\cos\theta\frac{2}{1+R^{2}},-\sin\theta\frac{2}{1+R^{2}}\right),\\
 \left(\vphi_{1}(R),\vphi_{2}(R)\right)=&\left(\frac{2}{1+R^{2}}\sin\theta,\frac{2}{1+R^{2}}\cos\theta\right).
 \end{split}
\end{align}
On the level of fixed angular momenta these solutions generate the following zero modes
\EQ{\nn
\frakA_{1} \binom{(1+R^2)^{-1}}{i(1+R^2)^{-1}}=\binom{0}{0},\qquad \frakA_{-1} \binom{(1+R^2)^{-1}}{-i(1+R^2)^{-1}}=\binom{0}{0}
}

Next, we turn our attention to Lorentz transforms which  are of the form 
\[
L_\alpha:\quad(t,x_1,x_2)\mapsto (t\cosh \alpha + x_1\sinh\alpha, t\sinh\alpha + x_1\cosh\alpha,x_2),\quad \alpha\in\R
\]
whence $Q(R)$ transforms into $$Q_\alpha(t,x_1,x_2)=Q(R_\alpha(t,x_1,x_2)),\quad R_\alpha(t,x_1,x_2):=\sqrt{( t\sinh\alpha + x_1\cosh\alpha)^2+x_2^2}$$ 
with derivative 
\[
\partial_\alpha\Big|_{\alpha=0} Q_\alpha(t,x_1,x_2) = \frac{2x_1 t}{R(1+R^2)} = \frac{2t\cos\theta}{1+R^2} 
\]
The stationary equivariant wave map associated with azimuth angle $Q(R)$ is 
\EQ{\label{eq:uWM}
u(x_1,x_2) & = (\cos\theta\sin Q(R), \sin\theta\sin Q(R), \cos Q(R)) \\
&= \left( \frac{x_1}{R} \sin Q(R), \frac{x_2}{R} \sin Q(R), \cos Q(R)\right)
}
which transforms into 
\EQ{\nn 
u_\alpha (t, x_1,x_2) 
&= \left( \frac{t\sinh\alpha + x_1\cosh\alpha}{R_\alpha} \sin Q(R_\alpha ), \frac{x_2}{R_\alpha} \sin Q(R_\alpha), \cos Q(R_\alpha)\right)
}
Differentiating in $\alpha$ we obtain 
\EQ{\nn
\partial_\alpha\Big|_{\alpha=0} u_\alpha(t,x_1,x_2)  & = t\left( \frac{x_2^2}{R^3} \sin Q + \frac{2x_1^2}{R^2(1+R^2)} \cos Q, -\frac{x_1x_2}{R^3}\sin Q + \frac{2x_1x_2}{R^2(1+R^2)}\cos Q, -\frac{2x_1}{R(1+R^2)} \sin Q\right) \\
&= t\left( \frac{\sin^2\theta}{R}\sin Q + \frac{2\cos^2\theta}{1+R^2}\cos Q, -\frac{\sin\theta\cos\theta}{R}\sin Q + \frac{2\sin\theta\cos\theta}{1+R^2}\cos Q, -\frac{2\cos\theta}{1+R^2}\sin Q\right)\\
&= t \frac{2\cos\theta}{1+R^2}E_1  - t \frac{\sin\theta\sin Q}{R} E_2 =  t \frac{2\cos\theta}{1+R^2}E_1  - t \frac{2\sin\theta}{1+R^2} E_2
}
This is exactly $t\vphi$ where $\vphi$ agrees with the first line of~\eqref{null: space trans}. The second line of~\eqref{null: space trans} is obtained by the Lorentz transform in $(t,x_2)$.  
Thus, we do not obtain any new solutions to $\calL\vphi=0$ from Lorentz transforms. 

Next we consider the rotations acting on the target. The linearized kernels corresponding to the three rotations are given by 
\begin{align}\label{three rotations}
 \begin{split}
  \left(-\sin\theta\sin Q,\cos\theta\sin Q,0\right),\quad \left(0, -\cos Q,\sin\theta\sin Q\right),\quad
  \left(\cos Q, 0, -\cos\theta\sin Q\right).
 \end{split}
\end{align}
These are obtained by acting on $u$ as in  \eqref{eq:uWM} with the infinitesimal rotation matrices
\[
\left[\begin{matrix} 
0 & -1 & 0 \\  1 & 0 & 0 \\ 0 & 0 & 0 
\end{matrix} \right], \quad \left[\begin{matrix} 
0 &   0  & 0 \\  0 & 0 & -1 \\ 0 & 1 & 0 
\end{matrix} \right], \quad \left[\begin{matrix} 
0 & 0 & -1 \\  0 & 0 & 0 \\ 1 & 0 & 0 
\end{matrix} \right],
\]
respectively. 
We expand these tangent maps in the basis $(E_{1},E_{2},\Phi)$. We have
\begin{align}\label{rotations expansion}
 \begin{split}
  \left(-\sin\theta\sin Q,\cos\theta\sin Q,0\right)=&\sin Q E_{2} = \frac{2R}{1+R^2} E_2,\\
  \left(0, -\cos Q,\sin\theta\sin Q\right)=&-\sin\theta E_{1}-\cos\theta\cos Q E_{2} = -\sin\theta E_{1}-\cos\theta\frac{1-R^2}{1+R^2} E_{2},\\
  \left(-\cos Q, 0, \cos\theta\sin Q\right)=&\cos\theta E_{1}-\sin\theta\cos Q E_{2} = \cos\theta E_{1}-\sin\theta \frac{1-R^2}{1+R^2} E_{2}.
 \end{split}
\end{align}
The first line is one of the two modes from~\eqref{eq:dilmodes}. The others are new $0$ modes which we write as 
\EQ{\label{eq:res2}
\vphi = \left(\sin\theta, \cos\theta \frac{1-R^2}{1+R^2}\right),\quad \vphi=\left( \cos\theta, -\sin\theta \frac{1-R^2}{1+R^2}\right) 
}
Finally we consider the space rotation acting on the domain $\bbR^{2}_{x}\times \bbR^{1}_{t}$. Since we are only concerned with the space rotation, we consider the stationary wave map \eqref{eq:uWM}. If we represent the rotation $r_{\alpha}$ by the matrix
\begin{align*}
	r_{\alpha}=\left(\begin{array}{cc}
	\cos\alpha&-\sin\alpha\\\sin\alpha&\cos\alpha
	\end{array}\right),
\end{align*}
then 
\begin{align*}
	&\left(r_{\alpha}u\right)(x_{1},x_{2})=\left(\frac{\cos\alpha x_{1}-\sin\alpha x_{2}}{R}\sin Q(R),\frac{\sin\alpha x_{1}+\cos\alpha x_{2}}{R}\sin Q(R),\cos Q(R)\right)\\
	\Rightarrow\quad &\frac{d}{d\alpha}(r_{\alpha}u)|_{\alpha=0}=\left(-\frac{x_{2}}{R}\sin Q(R),\frac{x_{1}}{R}\sin Q(R), 0\right)=\sin Q E_{2},
\end{align*}
which was already considered in \eqref{rotations expansion}. Therefore we do not obtain any new solutions to $\calL\varphi=0$ from the space rotation in $\bbR^{2}_{x}\times \bbR^{1}_{t}$.

Next, we diagonalize  the operator from \eqref{operator Fourier matrix} and thus relate the $0$ modes from \eqref{eq:dilmodes}, \eqref{null: space trans}, and~\eqref{eq:res2} to eigenvalues or resonances of certain scalar self-adjoint Schr\"odinger operators in $L^2(R\,dR)$. In this way we shall moreover see that there are no other $0$ modes than those generated from symmetries.

\subsection{Diagonalizing the linear operators}
 Let
\begin{align}\label{fn gn expression}
f_{n}(R):=\frac{n^{2}}{R^{2}}+\frac{R^{4}-6R^{2}+1}{R^{2}(R^{4}+2R^{2}+1)},\quad g_{n}(R):=n\frac{2-2R^{2}}{R^{2}(R^{2}+1)}
\end{align} 
Writing the matrix \eqref{operator Fourier matrix} in the form
 \begin{align}\label{formal matrix}
 \frakA_{n}:=\left(\begin{array}{cc}
 a&-ib\\ib&a
 \end{array}\right),\quad a:=\partial_{R}^{2}+\frac{1}{R}\partial_{R}-f_{n}(R),\quad b:=g_{n}(R), 
 \end{align}
we change variables to reduce it to a diagonal matrix: 
 \begin{align*}
 \frakD_{n}:=\left(\begin{array}{cc}
 a+b&0\\0&a-b
 \end{array}\right)=P\frakA_{n} P^{-1},\quad P=\left(\begin{array}{cc}
 1&-i\\1&i
 \end{array}\right),\quad P^{-1}=\left(\begin{array}{cc}
 \frac{1}{2}&\frac{1}{2}\\\frac{i}{2}&-\frac{i}{2}
 \end{array}\right)
 \end{align*}
 Coming back to \eqref{operator Fourier matrix}, we transform  
 $(\hvphi_{1},\hvphi_{2})$ into 
 \begin{align}\label{diag unknown}
 \veps:=(\veps_{+},\veps_{-})^{T}:=(\hvphi_{1}-i\hvphi_{2},\hvphi_{1}+i\hvphi_{2})^{T},
 \end{align}
 so that  the operator acting on $\veps$ is given by
 \begin{align}\label{diag operator}
 \begin{split}
 P\frakA_{n}  \binom{\hvphi_{1}}{\hvphi_{2}} &=\left(\begin{array}{cc}
 \partial_{R}^{2}+\frac{1}{R}\partial_{R}-f_{n}(R)+g_{n}(R)&0\\
 0&\partial_{R}^{2}+\frac{1}{R}\partial_{R}-f_{n}(R)-g_{n}(R)
 \end{array}\right)\left(\begin{array}{c}
 \veps_{+}\\\veps_{-}
 \end{array}\right)
 \end{split}
 \end{align}
 We denote the operators on the diagonal 
  \begin{align}\label{def Hn}
  H^{+}_{n}:=\partial_{R}^{2}+\frac{1}{R}\partial_{R}-f_{n}(R)+g_{n}(R),\quad H^{-}_{n}:=\partial_{R}^{2}+\frac{1}{R}\partial_{R}-f_{n}(R)-g_{n}(R).
 \end{align}
Next, we transform the $0$ modes from \eqref{eq:dilmodes}, \eqref{null: space trans}, and~\eqref{eq:res2} which we found above  according to~\eqref{diag unknown}. This yields
\EQ{\label{eq:Hmodes}
H_0^{\pm} \frac{R}{1+R^2}=0, \quad H_1^{+} \frac{1}{1+R^2} = H_{-1}^{-} \frac{1}{1+R^2}=0,\quad H_1^{-} \frac{R^2}{1+R^2} = H_{-1}^{+} \frac{R^2}{1+R^2}=0
}
These relations will play an important role in our analysis, especially the middle one due to $ \frac{1}{1+R^2} \in L^2(R\, dR)$. In other words, $0$ is an eigenvalue of $H_1^+=H_{-1}^{-}$. In contrast, for $H_0^{\pm}$ and $H_{-1}^+ = H_1^{-}$ we only have a zero energy resonance since the associated functions are not in $L^2(R\,dR)$. While these resonances affect the spectral theory, cf.~Section~\ref{sec:LAN},  they do not appear explicitly in the dynamical analysis. 
 
 \subsection{Determination of perturbative corrections and the precise nonlinearity}\label{sec:nonlin}
 In this section we compute the precise nonlinearity for the equation of $(\vphi_{1},\vphi_{2})$.  Recall~\eqref{perturbed sol}, viz. 
 \begin{align}\nn
 \Psi:=\Phi+\pphi+\mathfrak{q}
\end{align}
The term $\mathfrak{q}$ is not uniquely determined from this ansatz. We now choose it to be parallel to~$\Phi$ itself. This can always be done provided $\vphi$ is small. Thus, we set
 \begin{align}\label{perturbed sol precise}
 \begin{split}
\Psi:=\Phi+\pphi+a\left(\pphi\right)\Phi  ,\quad a\left(\pphi\right)=\sqrt{1-|\pphi|^{2}}-1.
 \end{split}
 \end{align}
Equation \eqref{eq:lin pphi} then transforms into the following one 
\begin{align}\label{pphi eq precise 1}
 \begin{split}
  &\Box \pphi+\left(|\nabla\Phi|^{2}-|\Phi_{t}|^{2}\right)\pphi+2\left(\nabla\Phi\cdot\nabla\pphi-\Phi_{t}\cdot(\pphi)_{t}\right)\Phi:=N_{p},
 \end{split}
\end{align}
with right-hand side 
\begin{align}\label{pphi eq precise 2}
 \begin{split}
  N_{p}=&\left(|\partial_{t}\Psi|^{2}-|\nabla\Psi|^{2}\right)a\left(\pphi\right)\Phi+\left(|\partial_{t}\left(a(\pphi)\Phi\right)|^{2}-|\nabla\left(a(\pphi)\Phi\right)|^{2}\right)\Phi\\
  &+\left(\partial_{t}\left(\Phi+\pphi\right)\cdot\partial_{t}\left(a(\pphi)\Phi\right)-\nabla\left(\Phi+\pphi\right)\cdot\nabla\left(a(\pphi)\Phi\right)\right)\Phi\\
  &+\left(|\partial_{t}\Psi|^{2}-|\partial_{t}\Phi|^{2}-|\nabla\Psi|^{2}+|\nabla\Phi|^{2}\right)\pphi-\Box\left(a(\pphi)\Phi\right).
 \end{split}
\end{align}
Next,  we expand $N_{p}$ in the frame $\Phi, E_{1}, E_{2}$ and determine the coefficients for $E_{1}$ and $E_{2}$. For the latter it suffices to retain the terms in the last two lines in \eqref{pphi eq precise 2}, as the others are parallel to~$\Phi$. In fact, of these two final terms we can discard another term parallel to $\Phi$, and expand the following expression:
\begin{align}\label{non linear pre}
 \begin{split}
  \left(|\partial_{t}\Psi|^{2}-|\partial_{t}\Phi|^{2}-|\nabla\Psi|^{2}+|\nabla\Phi|^{2}\right)\pphi-2\nabla\left(a(\pphi)\right)\cdot\nabla\Phi+2\left(a(\pphi)\right)_{t}\Phi_{t}-a(\pphi)\Box\Phi.
 \end{split}
\end{align}
Now we compute each term in \eqref{non linear pre}:
\begin{align}\label{non linear 1a}
 \begin{split}
  \Psi_{t}=&U_{t}E_{1}+\vphi_{1,t}E_{1}-\vphi_{1}U_{t}\Phi+\vphi_{2,t}E_{2}+a(\pphi)U_{t}E_{1}+\left(a(\pphi)\right)_{t}\Phi,
 \end{split}
\end{align}
which yields 
\begin{align}\label{non linear 1b}
 \begin{split}
  |\partial_{t}\Psi|^{2}-|\partial_{t}\Phi|^{2}=&\left(\left(a(\pphi)\right)_{t}-U_{t}\vphi_{1}\right)^{2}+\vphi_{2,t}^{2}\\
  &+\left(\vphi_{1,t}+a(\pphi)U_{t}\right)^{2}+2U_{t}\cdot\left(\vphi_{1,t}+a(\pphi)U_{t}\right).
 \end{split}
\end{align}
For $\nabla\Psi$, we obtain 
\begin{align}\label{non linear 2a}
 \begin{split}
  \nabla\Psi=&\nabla\Phi+\nabla\left(\pphi\right)+\nabla\left(a(\pphi)\Phi\right)\\
  =&\left(\begin{array}{c}
           U_{r}E_{1}\\\frac{\sin U}{r}E_{2}
          \end{array}\right)+\left(\begin{array}{c}
          \partial_{r}\left(a(\pphi)\right)\Phi\\
          \frac{1}{r}\partial_{\theta}\left(a(\pphi)\right)\Phi
          \end{array}\right)+a(\pphi)\left(\begin{array}{c}
           U_{r}E_{1}\\\frac{\sin U}{r}E_{2}
          \end{array}\right)\\
          &+\left(\begin{array}{c}
                                    \vphi_{1,r}E_{1}-\vphi_{1}U_{r}\Phi+\vphi_{2,r}E_{2}\\
                                    \frac{1}{r}\vphi_{1,\theta}E_{1}+\frac{\cos U}{r}\vphi_{1}E_{2}+\frac{1}{r}\vphi_{2,\theta}E_{2}-\frac{1}{r}\vphi_{2}\left(\sin U\Phi+\cos U E_{1}\right)
                                   \end{array}\right),
                                   \end{split}
\end{align}
which gives
\begin{align}\label{non linear 2b}
 \begin{split}
  &-|\nabla\Psi|^{2}+|\nabla\Phi|^{2}\\
  =&-\left(2a(\pphi)+\left(a(\pphi)\right)^{2}\right)\left(U_{r}^{2}+\frac{\sin^{2}U}{r^{2}}\right)-\left(\left(a(\pphi)\right)^{2}_{r}+\frac{1}{r^{2}}\left(a(\pphi)\right)^{2}_{\theta}\right)\\
  &-\vphi_{1,r}^{2}-U_{r}^{2}\vphi_{1}^{2}-\vphi^{2}_{2,r}-\frac{1}{r^{2}}\left(\vphi_{1,\theta}-\cos U\vphi_{2}\right)^{2}-\frac{1}{r^{2}}\left(\cos U\vphi_{1}+\vphi_{2,\theta}\right)^{2}-\frac{\sin ^{2}U}{r^{2}}\vphi^{2}_{2}\\
  &-2U_{r}\vphi_{1,r}-\frac{2\sin U}{r^{2}}\left(\cos U\vphi_{1}+\vphi_{2,\theta}\right)+2\left(a(\pphi)\right)_{r}U_{r}\vphi_{1}+\frac{2}{r^{2}}\left(a(\pphi)\right)_{\theta}\sin U\vphi_{2}\\
  &-2a(\pphi)U_{r}\vphi_{1,r}-2\frac{a(\pphi)}{r^{2}}\sin U\left(\cos U\vphi_{1}+\vphi_{2,\theta}\right).
 \end{split}
\end{align}
Combining \eqref{non linear 1b} and \eqref{non linear 2b}, we obtain
\begin{align}\label{non linear 1 2}
 \begin{split}
  &|\partial_{t}\Psi|^{2}-|\partial_{t}\Phi|^{2}-|\nabla\Psi|^{2}+|\nabla\Phi|^{2}\\
  =&-\left(2a(\pphi)+\left(a(\pphi)\right)^{2}\right)\left(U_{r}^{2}-U^{2}_{t}+\frac{\sin^{2}U}{r^{2}}\right)+\left(a(\pphi)\right)_{t}^{2}-\left(a(\pphi)\right)^{2}_{r}-\frac{1}{r^{2}}\left(a(\pphi)\right)^{2}_{\theta}\\
  &+\vphi_{1,t}^{2}-\vphi_{1,r}^{2}-\frac{1}{r^{2}}\vphi_{1,\theta}^{2}+\vphi_{2,t}^{2}-\vphi_{2,r}^{2}-\frac{1}{r^{2}}\vphi_{2,\theta}^{2}-\frac{1}{r^{2}}\vphi_{2}^{2}  +2\vphi_{1}\left(\left(a(\pphi)\right)_{r}U_{r}-\left(a(\pphi)\right)_{t}U_{t}\right)\\
  &+\left(U^{2}_{t}-U^{2}_{r}-\frac{\cos^{2}U}{r^{2}}\right)\vphi_{1}^{2}
  +2\left(1+a(\pphi)\right)\left(U_{t}\vphi_{1,t}-U_{r}\vphi_{1,r}-\frac{\sin U}{r^{2}}\cos U\vphi_{1}\right)\\
  &+\frac{2\cos U}{r^{2}}\left(\vphi_{1,\theta}\vphi_{2}-\vphi_{1}\vphi_{2,\theta}\right)-\frac{2\sin U}{r^{2}}\vphi_{2,\theta}
  +\frac{2\sin U}{r^{2}}\left(\left(a(\pphi)\right)_{\theta}\vphi_{2}-a(\pphi)\vphi_{2,\theta}\right).
 \end{split}
\end{align}
We write $-2\nabla\left(a(\pphi)\right)\cdot\nabla\Phi$ in the form 
\begin{align}\label{non linear 3}
 \begin{split}
  &-2\nabla\left(a(\pphi)\right)\cdot\nabla\Phi\\
  =&\frac{2\left(\pphi\cdot\nabla\pphi\right)\cdot\nabla\Phi}{\sqrt{1-|\pphi|^{2}}}\\
  =&\frac{2}{\sqrt{1-|\pphi|^{2}}}\left(U_{r}\left(\vphi_{1}\vphi_{1,r}+\vphi_{2}\vphi_{2,r}\right)E_{1}+\frac{\sin U}{r^{2}}\left(\vphi_{1}\vphi_{1,\theta}+\vphi_{2}\vphi_{2,\theta}\right)E_{2}\right)
 \end{split}
\end{align}
and $2\left(a(\pphi)\right)_{t}\Phi_{t}$ turns into 
\begin{align}\label{non linear 4}
 \begin{split}
  2\left(a(\pphi)\right)_{t}\Phi_{t}=&-2\frac{\pphi\cdot\left(\pphi\right)_{t}\Phi_{t}}{\sqrt{1-|\pphi|^{2}}}
  =-\frac{2}{\sqrt{1-|\pphi|^{2}}}\left(\vphi_{1}\vphi_{1,t}+\vphi_{2}\vphi_{2,t}\right)U_{t}E_{1}.
 \end{split}
\end{align}
The final term in \eqref{non linear pre} can be written as
\begin{align}\label{non linear 5}
 \begin{split}
  a(\pphi)\left(|\nabla\Phi|^{2}-|\Phi_{t}|^{2}\right)\Phi=&a(\pphi)\left(U^{2}_{r}-U^{2}_{t}+\frac{\sin^{2} U}{r^{2}}\right)\Phi,
 \end{split}
\end{align}
since $\Phi$ is a wave map. Now we can write the equations \eqref{coe eq rough 1}--\eqref{coe eq rough 2} in the form
\begin{align}\label{coe eq precise 1}
\begin{split}
& \Box\vphi_{1}-\frac{\cos 2U}{r^{2}}\vphi_{1}-\frac{2\cos U}{r^{2}}\vphi_{2,\theta}\\=&\left(|\Psi_{t}|^{2}-|\Phi_{t}^{2}|-|\nabla\Psi|^{2}+|\nabla\Phi|^{2}\right)\vphi_{1}\\
&+\frac{2}{\sqrt{1-|\pphi|^{2}}}\left(U_{r}(\vphi_{1}\vphi_{1,r}+\vphi_{2}\vphi_{2,r})-U_{t}(\vphi_{1}\vphi_{1,t}+\vphi_{2}\vphi_{2,t})\right):=N(\vphi_{1}),
\end{split}
 \end{align}
and 
\begin{align}\label{coe eq precise 2}
 \begin{split}
  &\Box\vphi_{2}-\frac{\cos^{2}U}{r^{2}}\vphi_{2}+(U^{2}_{r}-U^{2}_{t})\vphi_{2}+\frac{2\cos U}{r^{2}}\vphi_{1,\theta}\\
  =&\left(|\Psi_{t}|^{2}-|\Phi_{t}^{2}|-|\nabla\Psi|^{2}+|\nabla\Phi|^{2}\right)\vphi_{2}\\
  &+\frac{2\sin U}{r^{2}\sqrt{1-|\pphi|^{2}}}\left(\vphi_{1}\vphi_{1,\theta}+\vphi_{2}\vphi_{2,\theta}\right):=N(\vphi_{2}).
 \end{split}
\end{align}
In \eqref{coe eq precise 1}-\eqref{coe eq precise 2}, $|\Psi_{t}|^{2}-|\Phi_{t}^{2}|-|\nabla\Psi|^{2}+|\nabla\Phi|^{2}$ is given by \eqref{non linear 1 2}. To write down the exact equation in $(\tau,R)$-variable, we need to put the error between $U$ and $2\arctan R$ as well as the term $-\frac{(1+\nu)^{2}}{\nu^{2}\tau^{2}}\frac{4R^{2}}{R^{4}+2R^{2}+1}$ on the right hand side. Let the background equivariant solution be given by $U=Q+\eps$. We start by computing 
\begin{align}\label{difference equivariant 1}
 \begin{split}
  \cos 2U-\cos 2Q=&-2\sin\left(U+Q\right)\sin\eps=-2\sin\left(2Q+\eps\right)\sin\eps\\
 \cos^{2}U-\cos^{2}Q=&\frac{1}{2}\left(\cos 2U-\cos 2Q\right)=-\sin\left(2Q+\eps\right)\sin\eps\\
 2\cos U-2\cos Q=&-4\sin\left(\frac{2Q+\eps}{2}\right)\sin\frac{\eps}{2}
 \end{split}
\end{align}
For the difference from the term $U^{2}_{r}-U^{2}_{t}$, we arrive at
\begin{align}\label{difference equivariant 2}
 \begin{split}
  &U^{2}_{r}-U^{2}_{t}-Q^{2}_{r}+Q^{2}_{t}\\
  =&(Q+\eps)^{2}_{r}-Q^{2}_{r}-(Q+\eps)^{2}_{t}+Q^{2}_{t}  =2Q_{r}\eps_{r}+\eps^{2}_{r}-2Q_{t}\eps_{t}-\eps^{2}_{t}\\
  =&\frac{2\lambda^{2}}{1+R^{2}}\partial_{R}\eps+\lambda^{2}\left(\partial_{R}\eps\right)^{2}\\
  &-2\lambda^{2}\left(\partial_{\tau}Q+\frac{\lambda'(\tau)}{\lambda}R\partial_{R}Q\right)\left(\partial_{\tau}\eps+\frac{\lambda'(\tau)}{\lambda}R\partial_{R}\eps\right)-\lambda^{2}\left(\partial_{\tau}\eps+\frac{\lambda'(\tau)}{\lambda}R\partial_{R}\eps\right)^{2}\\
  =&\lambda^{2}\left(\frac{2\partial_{R}\eps}{1+R^{2}}+\left(\partial_{R}\eps\right)^{2}-\frac{\lambda'(\tau)}{\lambda}\frac{8R}{1+R^{2}}\left(\partial_{\tau}\eps+\frac{\lambda'(\tau)}{\lambda}R\partial_{R}\eps\right)-\left(\partial_{\tau}\eps+\frac{\lambda'(\tau)}{\lambda}R\partial_{R}\eps\right)^{2}\right).
 \end{split}
\end{align}
Therefore, equations \eqref{coe eq precise 1}--\eqref{coe eq precise 2} take the form
\begin{align}\label{coe eq more precise 1}
 \begin{split}
  &\Box\vphi_{1}-\frac{\cos 2Q}{r^{2}}\vphi_{1}-\frac{2\cos Q}{r^{2}}\vphi_{2,\theta}\\
  =&N(\vphi_{1})-\frac{2\sin(2Q+\eps)\sin\eps}{r^{2}}\vphi_{1}-\frac{4\sin\left(\frac{2Q+\eps}{2}\right)\sin\frac{\eps}{2}}{r^{2}}\vphi_{2,\theta}
  =:\frakN(\vphi_{1})
 \end{split}
\end{align}
and 
\begin{align}
  &\Box\vphi_{2}-\frac{\cos^{2}Q}{r^{2}}\vphi_{2}+\left(Q^{2}_{r}-Q^{2}_{t}\right)\vphi_{2}+\frac{2\cos Q}{r^{2}}\vphi_{1,\theta} \nn \\
  =&N(\vphi_{2})-\frac{\sin(2Q+\eps)\sin\eps}{r^{2}}\vphi_{2}+\frac{4\sin\left(\frac{2Q+\eps}{2}\right)\sin\frac{\eps}{2}}{r^{2}}\vphi_{1,\theta} \label{coe eq more precise 2} \\
  &-\lambda^{2}\left(\frac{2\partial_{R}\eps}{1+R^{2}}+\left(\partial_{R}\eps\right)^{2}-\frac{\lambda'(\tau)}{\lambda}\frac{8R}{1+R^{2}}\left(\partial_{\tau}\eps+\frac{\lambda'(\tau)}{\lambda}R\partial_{R}\eps\right)-\left(\partial_{\tau}\eps+\frac{\lambda'(\tau)}{\lambda}R\partial_{R}\eps\right)^{2}\right)\vphi_{2}
  =:\frakN(\vphi_{2}).\nn 
\end{align}
Passing  to $(\tau,R)$-variables, we obtain for $\vphi=:(\vphi_{1},\vphi_{2})^{T}$, 
\begin{align}\label{coe eq tau R}
 \begin{split}
  &-\left(\left(\partial_{\tau}+\frac{\lambda'}{\lambda}R\partial_{R}\right)^{2}+\frac{\lambda'}{\lambda}\left(\partial_{\tau}+\frac{\lambda'}{\lambda}R\partial_{R}\right)\right)\vphi+\frakL\vphi=\frakN(\vphi),\\
  &\frakN(\vphi)=:\lambda^{-2}\left(\begin{array}{c}
                                                        \frakN(\vphi_{1})\\\frakN(\vphi_{2})+\frac{(1+\nu)^{2}\lambda^{2}}{\nu^{2}\tau^{2}}\frac{4R^{2}}{1+2R^{2}+R^{4}}\vphi_{2}
                                                       \end{array}\right).
 \end{split}
\end{align}
The nonlinearity takes the form 
\begin{align}
  \lambda^{-2}\frakN(\vphi_{1})=&\lambda^{-2}N(\vphi_{1})-\frac{2\sin(2Q+\eps)\sin\eps}{R^{2}}\vphi_{1}-\frac{4\sin\left(\frac{2Q+\eps}{2}\right)\sin\frac{\eps}{2}}{R^{2}}\vphi_{2,\theta}\nonumber \\
  \lambda^{-2}\frakN(\vphi_{2})=&\lambda^{-2}N(\vphi_{2})-\frac{\sin(2Q+\eps)\sin\eps}{R^{2}}\vphi_{2}+\frac{4\sin\left(\frac{2Q+\eps}{2}\right)\sin\frac{\eps}{2}}{R^{2}}\vphi_{1,\theta}\label{nonlinearity tau R}\\
  &-\left(\frac{2\partial_{R}\eps}{1+R^{2}}+\left(\partial_{R}\eps\right)^{2}-\frac{\lambda'(\tau)}{\lambda}\frac{8R}{1+R^{2}}\left(\partial_{\tau}\eps+\frac{\lambda'(\tau)}{\lambda}R\partial_{R}\eps\right)-\left(\partial_{\tau}\eps+\frac{\lambda'(\tau)}{\lambda}R\partial_{R}\eps\right)^{2}\right)\vphi_{2}.\nonumber
\end{align}
 We will effectively solve the equations \eqref{coe eq tau R} in the backward light cone centered at the singularity by means of an iterative scheme, solving linear inhomogeneous problems consecutively. Assuming the data to be $C^\infty$ in the light cone (which is not really necessary), the iterates $\Pi_{\Phi^{\perp}}\vphi^{(j)}$, say, will also be $C^\infty$ when interpreted as functions on $I\times \R^2$. Now write this as 
\begin{align*}
\Pi_{\Phi^{\perp}}\vphi = \vphi_1\left(\begin{array}{c}\cos\theta\cos U\\ \sin\theta\cos U\\ -\sin U\end{array}\right) + \vphi_2\left(\begin{array}{c}-\sin\theta\\ \cos\theta\\0\end{array}\right)
\end{align*}
Expanding the functions $\vphi_j = \sum_{n}\hat{\vphi}_j(n)e^{in\theta}$, $j = 1, 2$, we infer the conditions that the functions 
\[
\hat{\vphi}_1(n)\sin U\cdot e^{in\theta}\in C^\infty,\quad\left(\hat{\vphi}_1(n)\cos U + i\hat{\vphi}_2(n)\right)e^{i(n+1)\theta}\in C^\infty,
\]
where we have summarized the condition for the first two components by reverting to complex notation. These imply that $\hat{\vphi}_1(n)$ has to vanish at least of order $|n|-1$ at the origin and can be expanded in a power series in $R$ with all powers congruent $2$ modulo $|n|-1$, while $\hat{\vphi}_1(n) + i\hat{\vphi}_2(n) = \veps_{-}(n)$ has to vanish to order at least $|n+1|$, with power series around $R = 0$ only containing powers congruent $2$ modulo $|n+1|$. Using the fact that $\overline{\hat{\vphi}_{1}(n)-i \hat{\vphi}_{2}(n)}=\hat{\vphi}_{1}(-n)+i\hat{\vphi}_{2}(-n)$ we know that $\veps_{+}(n)$ may not vanish at $R=0$ only for $n=1$.
This of course implies that the source terms for the equations of $\veps_{+}(n), \veps_{-}(n)$ need to have the Taylor expansions around $R = 0$ with the same properties (since the operators $H_n^{\pm}$ kill the lowest order term $R^{n\mp 1}$). Thus recalling \eqref{diag unknown}, \eqref{diag operator} and writing the equations for $\veps_{+}(n), \veps_{-}(n)$ in the form 
\begin{align}\label{eq diag phys}
\left(-\left(\partial_{\tau} + \frac{\lambda_{\tau}}{\lambda}R\partial_R\right)^2 - \frac{\lambda_{\tau}}{\lambda}\left(\partial_{\tau} + \frac{\lambda_{\tau}}{\lambda}R\partial_R\right) + H_n^{\pm}\right)\veps_{\pm}(n) = F_{\pm},\,
\end{align}
the right-hand sides will have Taylor expansions around $R = 0$ with the same properties. Crucially there is no reason at all for the individual terms that we derive later on in the formula for $F_{\pm}$ to have these properties at $R = 0$, but we do know a priori that {\it{their sum}} will have it.  

\section{Spectral analysis of the linear operators}\label{sec:LAN}
 This section is devoted entirely to the spectral analysis of the linear operators on the diagonal on the left-hand side of~\eqref{diag operator}. This refers to the determination of the spectral measures and the generalized Fourier transform of each of these operators for all angular momenta~$n\in\Z$, the case $n=0$ already having been dealt with  in the equivariant case~\cite{KST,KMiao}. 
 The most demanding aspect will be analyzing  the asymptotic behavior as $n\to\pm\infty$. As in \cite{CSST, CDST} this will be accomplished by means of a Liouville-Green transform applied to an equivalent semiclassical problem. However, we will need to go considerably further than these aforementioned references. We begin with elementary properties of the shape of the spectrum and the behavior at zero energy, followed by the spectral analysis of $n=\pm1$, and finally we present the more delicate turning point analysis for $|n|\ge2$, viewed as a semi-classical problem.  
 \subsection{Selfadjointness and resolvent}
The operators arising in~\eqref{diag operator} are 
 \begin{align}\label{diag operator redefine}
  H^{+}_{n}:=\partial_{R}^{2}+\frac{1}{R}\partial_{R}-f_{n}(R)+g_{n}(R),\quad H^{-}_{n}:=\partial_{R}^{2}+\frac{1}{R}\partial_{R}-f_{n}(R)-g_{n}(R).
 \end{align}
 where
 \begin{align*}
 f_{n}(R) &=\frac{n^{2}+1}{R^{2}}-\frac{8}{(R^{2}+1)^{2}},\quad g_{n}(R)=\frac{2n}{R^{2}}-\frac{4n}{R^{2}+1},\\
-f_{n}(R)+g_{n}(R) &=-\frac{(n-1)^{2}}{R^{2}}-\frac{4n}{R^{2}+1}+\frac{8}{(R^{2}+1)^{2}},\\
 -f_{n}(R)-g_{n}(R) &=-\frac{(n+1)^{2}}{R^{2}}+\frac{4n}{R^{2}+1}+\frac{8}{(R^{2}+1)^{2}}
\end{align*}
If $n=0$ (linearized equivariant operator) we will simply write $H_0$ in place of $H_0^+=H_0^{-}$. 
 These operators are symmetric in $L^2(R\, dR)$ with domain given by $C^2((0,\infty))$ functions of compact support in $(0,\infty)$. 
 In the following we will use standard terminology and techniques from spectral theory, see \cite{ReedSimon2,Teschl,GZ,Fulton}.  For that purpose it is convenient to also consider the 
 conjugations of $H^{\pm}_n$ by the weights~$R^{\frac12}$, viz. 
 \begin{equation}
 \label{eq:calHn}
 \begin{split}
  \calH^{+}_{n}=R^{\frac12}  H^{+}_{n} R^{-\frac12} = &\partial^{2}_{R}+\frac{1}{4R^{2}}-\frac{(n-1)^{2}}{R^{2}}-\frac{4n}{R^{2}+1}+\frac{8}{(R^{2}+1)^{2}},\\
  \calH^{-}_{n}=R^{\frac12}  H^{-}_{n} R^{-\frac12}  =&\partial^{2}_{R}+\frac{1}{4R^{2}}-\frac{(n+1)^{2}}{R^{2}}+\frac{4n}{R^{2}+1}+\frac{8}{(R^{2}+1)^{2}}.
 \end{split}
 \end{equation}
These operators are clearly symmetric in $L^2(dR)$ over the same dense family as before.

\begin{proposition}
\label{prop: no negative spectrum}
 $\calH^{+}_{n}$ and $\calH^{-}_{n}=\calH^{+}_{-n}$ with $n\in \bbN$ are nonpositive symmetric operators in the Hilbert space $L^{2}((0,\infty))$, with domain  given by $C^2$ functions compactly supported in $(0,\infty)$.  Each $\calH^{+}_{n}$ is strongly singular at $R=0$ and in the limit point case at $R=\infty$. Moreover,  for all $n\ne1$ these operators are in the limit point case at  $R=0$, whereas $\calH^+_1$ is limit circle at the origin. Thus, $\calH_n^+$ is essentially self-adjoint for $n\ne1$, while $\calH^+_1$ is not. We denote by $\calH^{+}_{n}=\calH^{-}_{-n}$, the self-adjoint operator in $L^{2}(dR)$ given by the  Friedrichs extensions.

  The spectrum of each $\calH_n^+$ is $(-\infty,0]$, which is also the essential spectrum. At the threshold $0$ one has the following properties: $\calH_0^+$ has a resonance with state $R^{\frac32}/(1+R^2)$, $\calH_1^+=\calH_{-1}^{-}$ has an eigenvalue with eigenfunction $\sqrt{R}/(1+R^2)$, and $\calH_{-1}^+=\calH_{1}^-$  has a resonance with state ${R^{\frac52}}/(1+R^2)$.  These are precisely the modes originating from symmetries, see~\eqref{eq:Hmodes}. 
  
 Finally, the  resolvent kernel of $(-H_n^+ + 1)^{-1}$ satisfies,  for any $f\in C([0,\infty))$ with compact support, 
\begin{align}\label{resolvent kernel Hn plus}
((-H_n^+ +1)^{-1}f)(R) = -R^{-\frac12} \int_0^\infty G_n(R,R';-1) \sqrt{R'} f(R')\, dR'  = cR^{|n-1|} (1+O(R^2)) 
\end{align}
as $R\to 0$. Here $G_n$ is the resolvent kernel of $(\calH_n^+ - 1)^{-1}$. In particular,  unless $n=1$, the left-hand side vanishes at $R=0$. 
\end{proposition}
\begin{proof}
 Let $f\in C^{2}_{0}((0,\infty))$ with compact support and note that  
 \begin{align*}
  \langle H_{n}^{+}f,f\rangle=&\int_{0}^{\infty}\partial_{R}\left(R\partial_{R}f\right)f\, dR-\int_{0}^{\infty}\frac{(n-1)^{2}}{R}f^{2}\, dR+\int_{0}^{\infty}\frac{8R}{(R^{2}+1)^{2}}f^{2}\, dR-\int_{0}^{\infty}\frac{4nR}{R^{2}+1}f^{2}\, dR\\
  \leq &-\int_{0}^{\infty}(f'(R))^{2}R\, dR-\int_{0}^{\infty}\frac{(n-1)^{2}}{R}f^{2}\, dR<0
 \end{align*}
 for all $n\ge2$. 
Similarly, for $H^{-}_{n}$ with $n\ge2$ we compute that 
\[
 \langle H_{n}^{+}f,f\rangle= -\int_{0}^{\infty}(f'(R))^{2}R\, dR+\int_{0}^{\infty} U(R) f^2(R)R\, dR <0 
\]
since 
\begin{align*}
 U(R) = &-\frac{(n+1)^{2}}{R^{2}}+\frac{4n}{R^{2}+1}+\frac{8}{(R^{2}+1)^{2}}\\
 =&\frac{-(n+1)^{2}R^{4}-\left(2(n+1)^{2}-4n(R^{2}+1)-8\right)R^{2}-(n+1)^{2}}{R^{2}(R^{2}+1)^{2}}\\
 =&\frac{-(n-1)^{2}R^{4}+(-2n^{2}+6)R^{2}-(n+1)^{2}}{R^{2}(R^{2}+1)^{2}}<0, 
\end{align*}
if $n\ge2$. Hence, the quadratic forms of $H_n^{\pm}$ over $L^2(R\,dR)$ are nonpositive for all $|n|\ge2$ on a dense family.  So we can pass to the standard  Friedrichs extensions for these operators (and equivalently, for their conjugates $\calH_n^{\pm}$ over $L^2(dR)$).  The equivariant case $n=0$ was treated in~\cite{KST}, and one has a negative operator here as well with a $0$ energy resonance. 

One checks that $\calH_1^+ \big(\sqrt{R}/(1+R^2)\big)=0$. Since $\sqrt{R}/(1+R^2)\in L^2((0,\infty))$ is a nonnegative eigenfunction, it must be the ground state and $\calH_1^+$ has no discrete spectrum. Thus, $\calH_1^+\le 0$ as a symmetric operator. 

On the other hand $\calH_{-1}^+ \big({R^{\frac52}}/(1+R^2)\big)=0$. The function $R^{\frac52}/(1+R^2)$ fails to lie in $L^2(dR)$ and is referred to as a zero energy resonance. The meaning of this is that its asymptotic behavior is {\em subordinate} at both $R=0$ and $R=\infty$. This refers to the fact that near $R=\infty$, $\calH_{-1}^+$ has a fundamental system with asymptotic behavior $\sqrt{R}$ and $\sqrt{R}\log R$, respectively. On the other hand, near $R=0$ the fundamental basis behaves as $R^{\frac52}$ and $R^{-\frac32}$, resp.  Using this function we conclude that 
\[
\calH_{-1}^+ = - D^* D, \quad D=\frac{d}{dR}+U, \quad U(R):= -\frac{5}{2R} + \frac{2R}{1+R^2}
\]
from which it again follows that $\calH_{-1}^+\le0$. 

The limit point behavior at $R=\infty$ of all $\calH_n^+$ follows from Theorem~X.8 in~\cite{ReedSimon2}.  The equation $\calH_n^+ f = zf$, with $z\in\C$ has a fundamental system $f_1(R,z), f_2(R,z)$ with the asymptotic behaviors  $f_1(R,z)\sim R^{n-\frac12}$ and $f_2(R,z)\sim R^{-n+\frac32}$, respectively,  as $R\to0+$ provided $n\ne1$. Note that  one of these is not in~$L^2((0,1))$, which we discard. We are left with one which we call the regular solution. 
So if $n\ne1$, then $\calH_{n}^{+}$ is in the limit point case at both $R=0$ and $R=\infty$. Therefore by Theorem~X.7 in~\cite{ReedSimon2}, the closure of $\calH_{n}^{+}$ for all $n\ne1$ is self-adjoint. 
In particular, this unique self-adjoint extension must agree with the Friedrichs extension. 
If $n=1$, then these solutions satisfy $f_1(R,z)\sim R^{\frac12}$ and $f_2(R,z)\sim R^{\frac12}\log R$, respectively, as $R\to0+$ (for any $z\in\C$).  We now show that the Friedrichs extension of $\calH_{1}^{+}$ does not allow a logarithmic singularity as $R\rightarrow0_{+}$. From the fact $\calH_{1}^{+}\left(\frac{\sqrt{R}}{1+R^{2}}\right)=0$ we conclude that 
\begin{align*}
\calH_{1}^{+}=-\calD^{*}\calD,\quad \calD=\frac{d}{dR}+\calV(R),\quad \calV(R):=\frac{3R^{2}-1}{2R(1+R^{2})}.
\end{align*}
By its construction, the domain of the Friedrichs extension of $\calH_{1}^{+}$ is a subspace of the completion of $C_{0}^{2}\left((0,\infty)\right)$ under the norm
\begin{align}\label{norm Firedrichs Ex}
\begin{split}
	\left\langle-\calH_{1}^{+}f,f\right\rangle=\left\|\calD f\right\|_{L^{2}(0,\infty)}^{2}
	=\int_{0}^{\infty}\left|f^{\prime}(R)+\frac{3R^{2}-1}{2R(1+R^{2})}f(R)\right|^{2}\,dR
	\end{split}
\end{align}
By explicit calculation, one verifies that for any function with the property that $f_{2}(R)\sim \sqrt{R}\log R$ as $R\rightarrow0_{+}$, with the corresponding asymptotic behavior of the derivative, 
this integral over $0<R<1$ is infinite. This means that  only the solution $f_1(R,z)\sim R^{\frac12}$ as $R\to0+$ is regular for $\calH_{1}^{+}$. This uniquely determines the resolvent kernel. Indeed, 
denoting the regular solution by~$\phi_n(R,z)$, and the Weyl solution $\psi_n(R,z)\in L^2((1,\infty))$ for $z\in\C\setminus\R$, we have the resolvent kernel of $(\calH_n^+-z)^{-1}$ satisfying
 \EQ{\label{eq:Green}
 G_n(R,R';z)= \frac{\phi_n(R,z)\psi_n(R',z)}{W_n(z)},\qquad 0<R<R'
 }
with $W_n(z)$ the Wronskian of $\phi_n,\psi_n$.  The asymptotic behavior~\eqref{resolvent kernel Hn plus} follows from the preceding. 
\end{proof}

The importance of \eqref{resolvent kernel Hn plus} for $n=1$ lies with the fact that the resolvent kernel $\left(\calH_{1}^{+}-1\right)^{-1}(R,R')$ has no logarithmic singularity as $R\rightarrow0_+$.  
As the previous proof shows, this is due to the  Friedrichs extension of $\calH_{1}^{+}$. Note that one cannot naively impose boundary conditions at $R=0$ due to the strong singularity of the potential. The operator $H_1^+$ plays a special role in our analysis, and is the only one with a zero energy eigenfunction.

Henceforth, we can consider both $H_n^{\pm}$ and $\calH_n^{\pm}$ as selfadjoint operators on $L^2(R\,dR)$, respectively~$L^2(dR)$, without further mention.  
The relevance of this lies with the distorted Fourier transform which these operators therefore possess~\cite{GZ,Fulton} and which will play a decisive role in our entire nonlinear analysis. The main point of this section is to exhibit the Fourier bases and the spectral measures of each~$\calH^{\pm}_n$. As already noted, the limit $n\to\pm\infty$ here is particularly delicate and will be treated as a semiclassical problem. 
We will see that the spectrum of all $\calH_n^+$ is purely absolutely continuous if $n\ne 1$, and $\calH_1^+$ is  absolutely continuous with a simple eigenvalue at~$0$. Finally, one can easily prove by truncation and the quadratic form computations appearing in the proof of Proposition~\ref{prop: no negative spectrum} that no $\calH_n^+$ with $|n|\ge2$ exhibits a $0$ energy resonance. But we have no need for this result so we skip it. 
  
  \subsection{The distorted Fourier transform}\label{sec:DFT}
We close this introductory discussion by recalling the distorted Fourier transform associated with an operator 
\EQ{\label{eq:T op}
T=-\frac{d^2}{dR^2}+V(R) \text{\ \ on \ \ } L^2((0,\infty)), 
} 
see~\cite{GZ,Fulton,EvKalf}. Here $V$ is real-valued, continuous, decaying as $R\to\infty$, and strongly singular at $R=0$. 
In this paper we only encounter potentials $V$ which are of the form $V(R)=R^{-2}(a_0+a_1 R^2 + O(R^4))$ as $R\to0$, with an analytic function $O(R^4)$. And $V(R)=O(R^{-2})$ as $R\to\infty$. 
While we are not necessarily assuming that $T$ is in the limit point case at $R=0$,  we suppose  so first for simplicity. Then let $\phi(\cdot,z), \theta(\cdot,z)$ be a 
fundamental system of $Tf=zf$, analytic near the real line, and real-valued on it. Assume the Wronskian $W(\theta(\cdot,z),\phi(\cdot,z))=1$ and that $\phi(\cdot,z)\in L^2((0,1))$.
Note that $\theta$ is unique only up to addition of a multiple of $k(z)\phi(\cdot,z)$ with a constant $k(z)$ that is analytic and real on~$\R$. And $\phi(\cdot,z)$ is unique only up to multiplication by a nonzero {\em analytic function} near $\R$ which is real on~$\R$. 
 Furthermore, 
let $\psi(R,z)$ be a Weyl-Titchmarsh solution for $\Im z>0$, i.e., it lies in  $L^2((1,\infty))$ for $\Im z>0$. 
The normalization of $\psi(R,z)$ is not too relevant, but it can be convenient to assume that $\psi$ has 
 asymptotic behavior $z^{-\frac14} e^{iz^{\frac12} R}$ as $R\to\infty$. Because of the decay of $V$ this asymptotic behavior can be achieved. This normalization implies $W(\overline{\psi(\cdot,z)}, \psi(\cdot,z))=2i$ provided $\Im z=0$.   The (generalized) Weyl-Titchmarsh $m$ function is then defined as 
\EQ{\label{eq:Cpsi}
C\,\psi(\cdot,z) =  \theta(\cdot,z) + m(z) \phi(\cdot,z),\quad C\ne0
}
with some constant $C$. 
Therefore, 
\EQ{\label{eq:gen m}
m(z) =\frac{W(\theta(\cdot,z),\psi(\cdot,z))}{W(\psi(\cdot,z),\phi(\cdot,z))}
}
This does not depend on the normalization of the Weyl-Titchmarsh solution. 
In view of the degrees of freedom we have in defining $\phi, \theta$ the $m$-function is far from unique. 
A spectral measure of $T$ is obtained as the limit
\EQ{\label{eq:spec meas}
\rho((\lambda_1,\lambda_2]) = \frac{1}{\pi} \lim_{\delta\to0+} \lim_{\epsilon\to0+} \int_{\lambda_1+\delta}^{\lambda_2+\delta} \Im m(\lambda+i\epsilon)\,d\lambda
}
The distorted Fourier transform of $f\in C((0,\infty))$ with compact support, 
\EQ{\label{eq:dist FT}
\hat{f}(\xi) = \int_0^\infty \phi(R,\xi) f(R) \, dR
}
is a unitary transformation $L^2((0,\infty))\to L^2(\R,\rho)$, with inverse transform 
\[
f(R) = \int_\R \phi(R,\xi) \hat{f}(\xi) \,\rho(d\xi),
\]
see~\cite{GZ} for the details and the existence of these integrals. Note that no assumption is made on the measure. In fact, if $T$ has eigenvalues, then $\rho$ will have atoms those points and the unitarity will contain the projections onto the eigenfunctions. 
Once this representation is obtained, we can derive equivalent ones as follows: if $h\in C(\R)$ is positive, then replacing 
$\phi(R,\xi) $ with $h(\xi)\phi(R,\xi)$ and $\rho(d\xi)$ with $h(\xi)^{-2}\rho(d\xi)$ leads to another -- equivalent -- distorted Fourier transform. 
The spectral measures in our case will be purely absolutely continuous, with the possible exception of an atom at~$0$ (as for $\calH_1^+$). But there cannot be any other atoms,
 or a singular continuous part.

In practice we will compute the density of the spectral measure as follows. First, we connect the solutions $\phi$ with the Weyl-Titchmarsh solutions, i.e., 
\EQ{\label{eq:phiapsi}
\phi(\cdot,\xi) = a(\xi) \psi(\cdot,\xi) + \overline{a(\xi) \psi(\cdot,\xi)},\quad \xi>0
}
Recall that we normalize $W(\theta,\phi)=1$ but we do not assume a normalization of the Weyl-Titchmarsh solutions $\psi$. Then $C(\xi)\psi(\cdot,\xi)=\theta(\cdot,\xi)+m(\xi)\phi(\cdot,\xi)$, see~\eqref{eq:Cpsi}, and 
\EQ{\nn 
CW(\psi(\cdot,\xi), \phi(\cdot,\xi)) & = 1 \\
|C(\xi)|^2 W(\psi(\cdot,\xi),\overline{\psi(\cdot,\xi)}) &= W(\theta(\cdot,\xi)+m(\xi)\phi(\cdot,\xi), \theta(\cdot,\xi)+\overline{m}(\xi)\phi(\cdot,\xi))  \\
&= \overline{m}(\xi)-m(\xi)=-2i\Im m(\xi)
}
In view of \eqref{eq:spec meas}, the density of the spectral measure is 
\EQ{\label{eq:a rep rho}
\frac{d\rho}{d\xi}(\xi) &= \frac{1}{\pi} \Im m = \frac{ W(\psi(\cdot,\xi), \overline{\psi(\cdot,\xi)})}{-2i\pi |W(\psi(\cdot,\xi),\phi(\cdot,\xi))|^2} \\
&= \frac{1 }{2i\pi |a(\xi)|^2  W(\psi(\cdot,\xi), \overline{\psi(\cdot,\xi)})} 
}
This relation was essential in \cite{KST, KMiao} since we can easily find $W(\psi(\cdot,\xi), \overline{\psi(\cdot,\xi)})$, say by computing the Wronskian at $R=\infty$, and $a(\xi)$ is found by matching representations of $\phi$ (such as~\eqref{eq:seriesJ1}) with an expansion such 
as~\eqref{eq:Hosc}. More precisely, we obtain an upper bound on $|a(\xi)|$ in this way by matching at $\xi^{\frac12}R\simeq1$, and the lower bound follows directly from~\eqref{eq:phiapsi} and a lower bound on $\phi$ since $| \phi(\cdot,\xi)|\le 2|a(\xi)||\psi(\cdot,\xi) |$.

Finally, if $T$ is in the limit-circle case at $R=0$, then we need to select a selfadjoint extension from an infinite family of possibilities, for example by means of a boundary condition. In the strongly singular case this cannot be done naively. In Proposition \ref{prop: no negative spectrum} we found that the Friedrichs extension turned out to be the correct choice since it guarantees the regularity property \eqref{resolvent kernel Hn plus}. 
The remainder of this section is devoted to the determination of a distorted Fourier transform for each of the angular momenta $n\in\Z$. In this regard the low modes $n=0,\pm 1$ play a special role due to the appearance of a threshold resonance or eigenvalue. 

We will freely switch between the self-adjoint operators $H_n^+$ in $L^2(R\, dR)$ and their versions $\calH_n^+$ which are self-adjoint in~$L^2(dR)$. 
The former are more natural  as they arise in the linearized wave map. In the context of~\eqref{eq:T op}, with $(Mf)(R):= R^{-\frac12} f(R)$, 
we have $T = M^{-1}\circ S\circ M$ where 
\[
S= -\frac{d^2}{dR^2} - \frac{1}{R}\frac{d}{dR} + \frac{1}{4R^2}  +V(R)
\]
If $\phi(\cdot,\xi)$ is the Fourier basis from~\eqref{eq:dist FT} relative to $T$, then $\tilde \phi(R,\xi):=R^{-\frac12} \phi(R,\xi)$ is the Fourier basis relative to~$S$ with the same spectral measure, and 
$S(\tilde \phi(\cdot,\xi))=\xi\tilde \phi(\cdot,\xi)$.  In fact, one has 
\EQ{
\label{eq:DFT2}
\hat{f}(\xi) &= \int_0^\infty \tilde \phi(R,\xi) f(R)R\, dR =  \int_0^\infty   \phi(R,\xi) f(R)R^{\frac12} \, dR \\
R^{\frac12} f(R) & = \int_0^\infty \phi(R,\xi) \hat{f}(\xi)\, \rho(d\xi), \qquad f(R)  = \int_0^\infty \tilde\phi(R,\xi) \hat{f}(\xi)\, \rho(d\xi)
} 
The unitarity in this context means that
$
\| f\|_{L^2(R\,dR)} = \|\hat{f}\|_{L^2(d\rho)}. 
$
We will formulate conditions below under which the integrals in~\eqref{eq:DFT2} converge. 

\subsection{Linearized operators at low angular momenta $n=0,\pm1$}\label{subsec:exceptionalnsupersymmetry}
In this section, we determine the Fourier transform associated with the nonpositive operators  
\begin{align}\label{linear n1}
\begin{split}
H^{+}_{1}=&\partial_{R}^{2}+\frac{1}{R}\partial_{R}-\frac{4}{R^{2}+1}+\frac{8}{(R^{2}+1)^{2}}=H^{-}_{-1}\\
H^{-}_{1}=&\partial_{R}^{2}+\frac{1}{R}\partial_{R}-\frac{4}{R^{2}}+\frac{4}{R^{2}+1}+\frac{8}{(R^{2}+1)^{2}}=H^{+}_{-1}.
\end{split} 
\end{align}
As noted in Proposition~\ref{prop: no negative spectrum}, $0$ is an eigenvalue of $H^+_{1}=H^{-}_{-1}$ with eigenfunction $1/(1+R^2)$, and a resonance of $H^+_{-1}=H^{-}_{1}$ with resonance function $R^2/(1+R^2)$.  These explicit solutions permit us to factorize the operators as follows.  From 
\begin{align}\label{D plus}
 \calD_{+}\left(\frac{1}{1+R^{2}}\right):=\left(\partial_{R}+\frac{2R}{1+R^{2}}\right)\left(\frac{1}{1+R^{2}}\right)=0,
\end{align}
we infer that $H^{+}_{1}=-\calD^{*}_{+}\calD_{+}$ where
\begin{align}\label{D star plus}
 \calD^{*}_{+}:=-\partial_{R}+\frac{R^{2}-1}{R(1+R^{2})}.
\end{align}
Similarly for $H^{-}_{1}$, we have
\begin{align}\label{D minus and star}
 H^{-}_{1}=-\calD^{*}_{-}\calD_{-},\quad \calD_{-}=\partial_{R}-\frac{2}{R}+\frac{2R}{1+R^{2}},\quad \calD_{-}^{*}=-\partial_{R}-\frac{3}{R}+\frac{2R}{1+R^{2}}.
\end{align}
We now associate with these second order operators their ``super-symmetric" versions which do not exhibit zero energy modes, viz. 
\begin{align}\label{H1 star}
 \tilde H^{+}_{1}:=\calD_{+}\calD^{*}_{+}=-\partial_{R}^{2}-\frac{1}{R}\partial_{R}+\frac{1}{R^{2}},\quad \tilde H^{-}_{1}:=\calD_{-}\calD^{*}_{-}=-\partial_{R}^{2}-\frac{1}{R}\partial_{R}+\frac{9}{R^{2}}-\frac{8}{R^{2}+1}.
\end{align}
We will develop the Fourier transform for these operators, and not for the original ones. This appears to be essential in our technique since the spectral measures of the super-symmetric versions are much less singular at threshold energies. The super-symmetric cousin of the linearized equivariant operator $H_0$ already played a key role in~\cite{KMiao}, see in particular their Section~4. The dynamical implications of the switching will not be apparent until much later when we solve the dynamical problem. 
$\tilde H^{+}_{1}$ is a pure Bessel operator whose Fourier transform is explicitly known. We begin with that operator.  The analysis of~$n=0$ was already carried out in~\cite[Section 5]{KST}. See also~\cite[Section 4]{KMiao} for the supersymmetric treatment of~$\calH_0$.

\subsubsection{$n=1$ mode}
 We recall the following facts about the Bessel operator $$\calH=-\frac{d^2}{dR^2} + \frac{3}{4R^2}=R^{\frac12} \tilde H^{+}_{1}R^{-\frac12}$$ from \cite[Section 4]{GZ}. 
A fundamental system of solutions of
\begin{align}\label{eigen eq}
 \calH f =zf
\end{align}
is given by
\begin{align}\label{funda sys}
 R^{\frac12} J_1(z^{\frac{1}{2}}R),\quad  R^{\frac12} Y_1(z^{\frac{1}{2}}R),\quad z\in\bbC\setminus\{0\},\quad R\in(0,\infty),
\end{align}
with $J_1(\cdot)$ and $Y_1(\cdot)$ the usual Bessel functions of order $1$.  A fundamental system of the type described in Section~\ref{sec:DFT} is given by
\begin{align}\label{modified funda sys}
\begin{split}
 \Phi_{1}(R,z)=&\frac{\pi}{2}z^{-\frac{1}{2}} R^{\frac12} J_{1}(z^{\frac{1}{2}}R),\\
 \Theta_{1}(R,z)=&z^{\frac{1}{2}} R^{\frac12} (-Y_{1}(z^{\frac{1}{2}}R)+\pi^{-1}\log(z)J_{1}(z^{\frac{1}{2}}R)).
\end{split}
 \end{align}
These two functions in \eqref{modified funda sys} extend to entire functions with respect to $z\in\bbC$ for fixed $R\ne0$ and they are real for $z\in\bbR$. Moreover, the normalizations are such that 
\begin{align}\label{Wronskian}
 W\left(\Theta_{1}(\cdot,z),\Phi_{1}(\cdot,z)\right)=\frac{\pi}{2} W(R^{\frac12}J_1(z^{\frac12}R), R^{\frac12}Y_1(z^{\frac12}R))=1,\quad z\in\bbC.
\end{align}
A Weyl-Titchmarsh solution to \eqref{eigen eq}    is given by
\begin{align}\label{L2 funda sys}
 \begin{split}
 R^{\frac12} H^{(1)}_{1}(z^{\frac{1}{2}}R)= R^{\frac12}\big( J_{1}(z^{\frac{1}{2}}R)+iY_{1}(z^{\frac{1}{2}}R)\big),\quad z\in\bbC\setminus[0,\infty),\quad R\in(0,\infty),
 \end{split}
\end{align}
with $H^{(1)}_{1}$ the Hankel function of order $1$. To account for the normalizations in \eqref{funda sys}, we modify \eqref{L2 funda sys} to
\begin{align}\label{modified L2 funda sys}
 \begin{split}
 \psi(R,z)=&z^{\frac{1}{2}} i  R^{\frac12}H^{(1)}_{1}(z^{\frac{1}{2}}R)\\
 =&z^{\frac{1}{2}}  R^{\frac12}\big(-Y_{1}(z^{\frac{1}{2}}R)+iJ_{1}(z^{\frac{1}{2}}R)\big)\\
 =&\Theta_{1}(z,R)+m(z)\Phi_{1}(z,R),\quad z\in\bbC\setminus[0,\infty),\quad R\in(0,\infty).
 \end{split}
\end{align}
In particular, the generalized Weyl-Titchmarsh function is 
\begin{align}\label{m function}
 m(z)=\frac{2}{\pi}z(i-\frac{1}{\pi}\log z),\quad z\in\bbC\setminus[0,\infty).
\end{align}
The branch of $\log z$ is such that it is real for $z>0$. In particular, $m(\xi+i0)=0$ if $\xi<0$ and 
the spectral measure $\rho(\xi)$ is absolutely continuous and of the explicit form
\EQ{\label{eq:n1 spec meas}
\rho(d\xi)=\frac{1}{\pi}\Im m(\xi+i0)\one_{[\xi>0]}\, d\xi = \frac{2}{\pi^2} \, \xi\,\one_{[\xi>0]}\, d\xi 
}
We have therefore obtained the following representation of the distorted Fourier transform associated with the super-symmetric cousin  $\tilde H^{+}_{1}$ of~$H^{+}_{1}$. 
As noted in Section~\ref{sec:DFT} the Fourier representation is not unique, and can be normalized in infinitely many ways. 

The power series representation for $J_{1}(u)$ is, see~\cite{AS},  
\EQ{\label{eq:seriesJ1}
J_{1}(u)=\sum_{m=0}^{\infty}\frac{(-1)^{m}}{m!\Gamma(m+2)}\left(\frac{u}{2}\right)^{2m+1} = \frac{u}{2} - \frac{u^3}{16} + O(u^5) \qquad u\to0
}
This implies that the Fourier basis of Lemma~\ref{lem:H1p} admits the following convergent expansion 
\begin{align}\label{phi1 small R}
\frac{\pi}{2} \xi^{-\frac12} J_1(\xi^{\frac12} R) =  \frac{\pi}{4} R + R \sum_{\ell=1}^\infty b_\ell (\xi R^2)^\ell
\end{align}
for all arguments. However,  it is only useful for small argument $R\xi^{\frac{1}{2}}\ll1$. 
For large arguments the solutions exhibit  oscillatory behavior. In fact, the Weyl-Titchmarsh solution from~\eqref{modified L2 funda sys}
admit the asymptotic series expansion, see~\cite{AS}, 
\EQ{
\label{eq:Hosc}
 \psi(R,\xi)= \xi^{\frac{1}{2}} i  R^{\frac12}H^{(1)}_{1}(\xi^{\frac{1}{2}}R) \sim  e^{-\frac{\pi i}{4}}\, \xi^{\frac{1}{4}} e^{i\xi^{\frac{1}{2}}R}\Big( \sqrt{\frac{2}{\pi}}
 +\frac{a_1}{\xi^{\frac{1}{2}}R} +\frac{a_2}{(\xi^{\frac{1}{2}}R)^2} +\ldots\Big)
}
as $\xi R^2\to\infty$. This means that $\psi(R,\xi) =  \xi^{\frac{1}{4}} e^{i\xi^{\frac{1}{2}}R}\sigma(\xi^{\frac{1}{2}}R)$ where $\sigma(q)$ is a smooth function of $q\ge1$ so that
for all $\ell\ge0$, $m\ge1$, 
\begin{align}\label{eq: def sigma}
\Big| (q\partial_q)^\ell \Big( e^{\frac{\pi i}{4}}\,\sigma(q) - \sqrt{\frac{2}{\pi}}
 - \frac{a_1}{q} - \frac{a_2}{q^2}  - \ldots - \frac{a_m}{q^m}\Big) \Big|\le C_{m,\ell} \, q^{-m-\ell-1}
\end{align}
We now connect $\Phi_{1}(R,\xi)$ and $\psi(R,\xi)$ as in \eqref{eq:phiapsi}, i.e., 
\EQ{\label{eq:phian1}
\Phi_{1}(\cdot,\xi) = a(\xi) \psi(\cdot,\xi) + \overline{a(\xi) \psi(\cdot,\xi)},\quad \xi>0
}
Using \eqref{eq:a rep rho}   we find   that 
\EQ{\label{eq:an1}
W(\psi(\cdot,\xi), \overline{\psi(\cdot,\xi)}) = -\frac{4i\xi}{\pi}, \quad a(\xi)=-\frac{\pi}{4i\xi},\quad \frac{d\rho}{d\xi}(\xi)  = \frac{1}{8\xi|a(\xi)|^2} = \frac{2\xi}{\pi^2} \one_{[\xi>0]} 
}
which agrees with \eqref{eq:n1 spec meas}.

\begin{lemma}
\label{lem:H1p} 
The distorted Fourier transform associated with $\calH=-\frac{d^2}{dR^2} + \frac{3}{4R^2}$  has the following property:  for any   $f\in C^2((0,\infty))$ with
\[
\int_0^\infty \Big( R^{-1} |f(R)| + |f'(R)| + R|f''(R)|\Big)\, dR \le M<\infty
\]
the Fourier transform 
\EQ{\label{eq:hatf exists}
\hat{f}(\xi) = \lim_{L\to\infty} \int_0^L \Phi_{1}(R,\xi) f(R)\, dR
}
with $\phi(R,\xi)$ as in~\eqref{modified funda sys}, 
exists for all $\xi>0$ and 
\EQ{\label{eq:abs xi int}
\int_0^\infty |\hat{f}(\xi)| |\Phi_{1}(R,\xi)|\, \xi\, d\xi\lesssim M
}
Thus, 
\EQ{\label{eq:H1p}
f(R) &  =   \frac{2}{\pi^2} \int_0^\infty  \hat{f}(\xi)\Phi_{1}(R,\xi)  \xi\,   d\xi, \quad \forall\; R>0
}
with absolutely convergent integrals. 
\end{lemma}
\begin{proof}
In view of \eqref{modified funda sys}, \eqref{eq:seriesJ1}, \eqref{eq:Hosc}, \eqref{eq:phian1}, and~\eqref{eq:an1}, 
\begin{align}\label{phi 1 bound}
|\Phi_{1}(R,\xi)|\lesssim \min(R^{\frac32}, \xi^{-\frac34}) \qquad \forall\; \xi>0,\; R>0
\end{align}
Thus, the left-hand side of \eqref{eq:abs xi int} is bounded by $\int_0^\infty |\hat{f}(\xi)|\, \xi^{\frac14}\, d\xi$ which we now proceed to estimate. 
Consider a partition of unity
  \EQ{\label{eq:POU}
  1= \chi_0(u) + \sum_{j=0}^\infty \chi(2^{-j} u) \qquad\forall\; u\ge0
  }
 with smooth functions $\chi_0$ supported on $[0,1]$ and $\chi$ supported on $[\frac12,2]$, respectively.   
 For any $\xi>0$ we define
  \EQ{
  \label{eq:xi splitn1} 
  A(\xi) &= \int_0^\infty \chi_0(R^2\xi) \Phi_{1}(R,\xi)\,   f(R) \, dR  \\
  B_j(\xi) &= \int_0^\infty \chi(2^{-j} R^2\xi) \Phi_{1}(R,\xi)\,  f(R) \, dR 
  }
 so that at least formally 
 \EQ{\label{eq:hatf sum}
 \hat{f}(\xi) = A(\xi) +\sum_{j=0}^\infty B_j(\xi)
 }
 Then 
 \EQ{\label{eq:Abdn1}
 |A(\xi)| &\lesssim \int_0^{\xi^{-\frac12}} R^{\frac32} |f(R)| \, dR.
 }
and 
 \EQ{\label{eq:A xiint}
  \int_0^\infty  |A(\xi)| \xi^{\frac14}  d\xi & \lesssim \int_0^\infty \int_0^{R^{-2}}  \xi^{\frac14}  \, d\xi  \; |f(R)|R^{\frac32} \, dR 
\lesssim \int_0^\infty R^{-1}|f(R)|\, dR
 }
 To bound $B_j$ we integrate by parts twice. To do so, we write  
 \EQ{\nn 
 B_j(\xi) &= \int_0^\infty \chi(2^{-j} R^2\xi) \, \xi^{\frac14} e^{i\xi^{\frac12}R} \sigma(R\xi^{\frac12})a(\xi)     f(R) \, dR  \\
 &\quad + \int_0^\infty \chi(2^{-j} R^2\xi) \, \xi^{\frac14} e^{-i\xi^{\frac12}R} \; \overline{\sigma(R\xi^{\frac12})a(\xi) }    f(R) \, dR 
 }
 It suffices to deal with the first line. Thus, by \eqref{eq:an1}, 
 \EQ{\label{eq:Bj bdn1} 
 |B_j(\xi)| &\lesssim \xi^{-\frac74} \int_0^\infty |\partial_R^2 \big( \chi(2^{-j} R^2\xi)  \sigma(R\xi^{\frac12})   f(R)  \big)|\, dR\\
 &\lesssim  \xi^{-\frac74}  \int_0^\infty \one_{[R^2\xi\simeq 2^j]} \big(   |f''(R)| + R^{-1} |f'(R)|+R^{-2} |f(R)|\big)\, dR
 }
 Summing in $j\ge0$ and inserting this bound into $\int_0^\infty |\hat{f}(\xi)|\, \xi^{\frac14}\, d\xi$ yields
 \EQ{\nn
\int_0^\infty \sum_{j=0}^\infty |B_j(\xi)| \xi^{\frac14}  \, d\xi &\lesssim\int_0^\infty \sum_{j=0}^\infty  \xi^{-\frac32} \int_0^\infty \one_{[R^2\xi\simeq 2^j]} \big(   |f''(R)| + R^{-1} |f'(R)|+R^{-2} |f(R)|\big)\, dR\,d\xi \\
&\lesssim \int_0^\infty  \big( R |f''(R)| +   |f'(R)|+R^{-1} |f(R)|\big)\, dR 
 } 
 In combination with \eqref{eq:A xiint} this proves~\eqref{eq:abs xi int}. 
 
The existence of the limit \eqref{eq:hatf exists} pointwise in $\xi>0$ follows from the convergence of the series~\eqref{eq:hatf sum}  for fixed $\xi>0$. Indeed, summing in~\eqref{eq:Bj bdn1} one obtains 
 \[
 \sum_{j=0}^\infty  |B_j(\xi)| \lesssim \xi^{-\frac74} \int_0^\infty\one_{[R^2\xi\gtrsim 1]} \big( |f''(R)| + R^{-1} |f'(R)|+R^{-2} |f(R)| \big)\, dR<\infty
 \]
 The finiteness in last line follows from the fact that $R\gtrsim \xi^{-\frac12}>0$ is bounded away from $0$, and from our assumption on~$f$.  
\end{proof}

\subsubsection{$n=-1$ mode}
None of the linearized operators associated with angular momenta other than $n=1$ can be explicitly reduced to an exact Bessel equation. However, they are perturbations of  exact Bessel operators. We use the technique of~\cite{KST} to obtain the asymptotic form of the Fourier basis and the spectral measure. The goal will be to obtain expansions similar to~\eqref{phi1 small R}, \eqref{eq:Hosc}
for the operator 
\begin{align*}
-\calH^{+}_{-1}=-\partial_{R}^{2}+\frac{15}{4R^{2}}-\frac{4}{R^{2}+1}-\frac{8}{(R^{2}+1)^{2}}.
\end{align*}
We exhibit a fundamental system of $\calH^{+}_{-1} f=0$. From Proposition~\ref{prop: no negative spectrum} one solution is $\Phi^{-1}_{0}(R):=\frac{R^{\frac{5}{2}}}{1+R^{2}}$. We seek another one with $W(\Theta^{-1}_{0}, \Phi^{-1}_{0})=1$.  This leads to the ODE 
\begin{align*}
 -\partial_{R}\Theta^{-1}_{0}+\left(\frac{5}{2R}-\frac{2R}{1+R^{2}}\right)\Theta^{-1}_{0}=\frac{1+R^{2}}{R^{\frac{5}{2}}}
\end{align*}
with general solution  
\begin{align*}
\Theta^{-1}_{0}(R)&=\frac{1 + 4 R^2 + C R^4  - 4 R^4 \log(R) }{ 4 R^{\frac32} (1 + R^2) }
\end{align*}
We fix the constant $C=-5$ so that $\Theta^{-1}_{0}(1)=0$: 
\[
\Theta^{-1}_{0}(R)=\Phi^{-1}_{0}(R)\left(\frac{4R^{2}+1}{4R^{4}}-\log R-\frac{5}{4}\right). 
\]
Following \cite{KST}, we use this fundamental system at energy $\xi=0$ to construct the Fourier basis $\Phi^{-1}(R,\xi)$  of the operator $-\calH^{+}_{-1}$ perturbatively.

\begin{proposition}\label{prop: basis n-1 small}
For all $R\ge0, \xi\ge0$, we have the following expansion for $\Phi^{-1}(R,\xi)$
\begin{align*}
\Phi^{-1}(R,\xi)=\Phi^{-1}_{0}(R)+R^{\frac{1}{2}}\sum_{j=1}^{\infty}(-R^{2}\xi)^{j}\Phi_{j}(R^{2})
\end{align*}
which converges absolutely. It converges uniformly if $R\xi^{\frac{1}{2}}$ remains bounded. 
Here $\Phi_{j}(u)\ge0$ are smooth functions of $u\ge0$ satisfying 
\EQ{\label{eq:Phijupper}
\Phi_{j}(u)\leq \frac{1}{j!}\frac{u}{1+u},\quad \textrm{for all}\quad u\geq 0, \; j\ge1
}
and $\Phi_1(u)\ge c_1 \frac{u}{1+u}$ for all $u\ge0$ with some absolute constant $c_1>0$. 
\end{proposition}
\begin{proof}
We make the ansatz 
\begin{align*}
\Phi^{-1}(R,\xi)=R^{\frac{1}{2}}\sum_{j=0}^{\infty}(-\xi)^{j}\bff_{j}(R),\qquad \bff_{0}(R):=\frac{R^{2}}{1+R^{2}}
\end{align*}
with
\begin{align*}
\calH^{+}_{-1}\left(R^{\frac{1}{2}}\bff_{j}\right)=R^{\frac{1}{2}}\bff_{j-1},\quad \bff_{-1}:=0.
\end{align*}
This leads to the recursion
\begin{align*}
\bff_{j}(R) &= \int_{0}^{R}R^{-\frac{1}{2}}s^{\frac{1}{2}}\left(\Phi^{-1}_{0}(R)\Theta^{-1}_{0}(s)-\Theta^{-1}_{0}(R)\Phi^{-1}_{0}(s)\right)\bff_{j-1}(s)\, ds\\ 
&= \int_{0}^{R}\frac{-s^{4}R^{4}\log s^{4}+s^{4}R^{4}\log R^{4}-s^{4}(4R^{2}+1)+R^{4}(4s^{2}+1)}{4sR^{2}(1+R^{2})(1+s^{2})}\bff_{j-1}(s)\, ds.
\end{align*}
Introducing the new variables $v:=s^2, u:=R^{2}$, as well as the new functions $\tbff_{j}(s^{2}):=\bff_{j}(s)$, we obtain
\EQ{\label{eq:tilf rec}
\tbff_{j}(u)= \int_{0}^{u}\frac{-u^{2}v^{2}\log v^{2}+u^{2}v^{2}\log u^{2}-v^{2}(4u+1)+u^{2}(4v+1)}{8uv(1+u)(1+v)}\tbff_{j-1}(v)\, dv.
}
Substituting $v=tu$ with $0<t\le 1$ yields 
\[
\tbff_{j}(u)= u\int_{0}^{1}\frac{4tu+1 - 2u^{2}t^{2}\log t- t^{2}(4u+1)}{8t(1+u)(1+tu)}\tbff_{j-1}(tu)\, dt
\]
Setting $h_j(u):=\tbff_{j}(u)/u^{j+1}, u>0, j\ge0 $ one has for all $k\ge0$, 
\EQ{\label{eq:hk rec}
h_{k+1}(u)= \int_{0}^{1}\frac{1-t^2 + 4tu (1-t)  - 2u^{2}t^{2}\log t }{8(1+u)(1+tu)}t^{k} h_{k}(tu)\, dt,\quad h_0(v)=\frac{1}{1+v}
}
Inductively this shows that $h_j(u)\in C^\infty([0,\infty))$. Since the kernel is positive and $h_0\ge0$, we also have $h_j\ge0$. Returning to the original variables, we infer that
\[
\Phi_{j}(R^2) = R^{-2j} f_j(R) = R^2 h_j(R^2),\quad \Phi_j(u) = uh_j(u)\ge0 
\]
is smooth in $u\ge0$.  To bound $h_j$ from above we make the following claim
\EQ{\label{eq:hk upper}
0\le \int_0^1  \frac{1-t^2 + 4ut (1-t)  - 2u^{2}t^{2}\log t }{(1+tu)^2} t^k \, dt \le \frac{7}{k+1},\qquad \forall\; k\ge0
} 
for all $u\ge0$. Assuming the claim we obtain 
inductively from \eqref{eq:hk rec} that
\[
h_k(u)\le \frac{1}{k! (1+u)}, \quad k\ge 0, \; u>0
\]
which gives us \eqref{eq:Phijupper}.   To prove the claim we make the following estimates
\EQ{\nn 
 \int_0^1  \frac{1-t^2 + 4ut (1-t)  - 2u^{2}t^{2}\log t }{(1+tu)^2} t^k \, dt &\le  \int_0^1  5 t^k \, dt -2 \int_0^1 t^k \log t\, dt = \frac{5}{k+1} + \frac{2}{(k+1)^2} \le \frac{7}{k+1}. 
}
For the lower bound on $h_1(u)$ we compute 
\EQ{\nn 
h_{1}(u)&= \int_{0}^{1}\frac{1-t^2 + 4tu (1-t)  - 2u^{2}t^{2}\log t }{8(1+u)(1+tu)^2}\, dt
}
In particular,  $h_1(0)=\frac{1}{12}$ and for any $u\ge0$ and with some absolute positive constants $c_0$, $c_1$
\EQ{\nn 
h_{1}(u)&\ge  \int_{\frac14}^{\frac34}\frac{1-t^2 + 4tu (1-t)  - 2u^{2}t^{2}\log t }{8(1+u)(1+tu)^2}\, dt \ge c_0 \int_{\frac14}^{\frac34}\frac{1 + u +u^{2} }{(1+u)(1+u)^2}\, dt 
\ge \frac{c_1}{1+u}
}
so that $\Phi_1(u)\ge c_1\frac{u}{1+u}$ as claimed. 
\end{proof}

Inspection of the proof shows that we can obtain similar lower bounds on all $\Phi_j$, but we have no need for them. 
Now we turn to the case when $R\xi^{\frac{1}{2}}\gtrsim 1$. Our goal is to obtain an asymptotic expansion similar to the one for the Hankel function~\eqref{eq:Hosc}. 

\begin{proposition}\label{prop: basis n-1 large}
For $R^{2}\xi\gtrsim1$, $\xi>0$ a Weyl-Titchmarsh function of  $\calH^{+}_{-1}$ is given by 
\EQ{\label{eq:PsiPlus-1}
\Psi^{+}_{-1}(R,\xi):=\xi^{-\frac{1}{4}}e^{iR\xi^{\frac{1}{2}}}\sigma_{-1}(R\xi^{\frac{1}{2}},R)
}
where $\sigma_{-1}$ is smooth in $q\gtrsim 1, R>0$ and admits the asymptotic approximation
\EQ{\nn 
\sigma_{-1}(q,R)\sim\sum_{j=0}^{\infty}q^{-j}\Psi^{+}_{j}(R),
}
in the sense that for all large integers $j_{0}$, and all indices $\alpha,\beta$
\begin{align*}
\sup_{R>0}\Big|(R\partial_{R})^{\alpha}(q\partial_{q})^{\beta}\big[\sigma_{-1}(q,R)-\sum_{j=0}^{j_{0}}q^{-j}\Psi^{+}_{j}(R)\big]\Big|\leq C_{j_{0},\alpha,\beta}\, q^{-j_{0}-1}
\end{align*}
for all $q\gtrsim 1$.  One has 
\EQ{\label{eq:sigma-1 exp}
 \Psi^{+}_{0}(R)=1,\quad \Psi^{+}_{1}(R)=\frac{15i}{8}-R\int_{R}^{\infty}\frac{2i}{1+s^{2}}\, ds=-\frac{i}{8}+ {O\left(\frac{1}{1+R^{2}}\right)}
}
as $R\to\infty$. For all $j\ge0$ the coefficient functions  $\Psi^{+}_{j}(R)$ 
are zero order symbols, i.e., 
\begin{align*}
\sup_{R>0}\left|(R\partial_{R})^{k}\Psi^{+}_{j}(R)\right|<\infty
\end{align*}
and  they are analytic at infinity. 
\end{proposition}
\begin{proof}
The argument is basically a verbatim repetition of the proof of~\cite[Proposition 5.6]{KST}.  From \eqref{eq:PsiPlus-1}, 
\[
\sigma_{-1}(R\xi^{\frac{1}{2}},R)= \Psi^{+}_{-1}(R,\xi)\xi^{\frac{1}{4}}e^{-iR\xi^{\frac{1}{2}}}
\]
which satisfies the conjugated equation
\[
\left( -\partial_R^2 -2i\xi^{\frac12}\partial_R + \frac{15}{4R^{2}}-\frac{4}{R^{2}+1}-\frac{8}{(R^{2}+1)^{2}}\right) \sigma_{-1}(R\xi^{\frac{1}{2}},R)=0
\]
This replaces eq.~(5.8) in~\cite{KST}. The ansatz 
\[
 \sigma_{-1}(R\xi^{\frac{1}{2}},R) = \sum_{j=0}^\infty \xi^{-\frac{j}{2}} f_j(r)
\]
leads to the recursion 
\[
2i\partial_{R} f_j(R) = \left( -\partial_R^2  + \frac{15}{4R^{2}}-\frac{4}{R^{2}+1}-\frac{8}{(R^{2}+1)^{2}}\right) f_{j-1}
\]
for all $j\ge1$ with $f_0\equiv 1$. Solving from $R=\infty$ where $f_j$ vanishes for $j\ge1$ we conclude that 
\[
f_j(R) = \frac{i}{2} \partial_{R}f_{j-1}(R) + \frac{i}{2} \int_R^\infty  \Big(   \frac{15}{4s^{2}}-\frac{4}{s^{2}+1}-\frac{8}{(s^{2}+1)^{2}}\Big)\, f_{j-1}(s)\, ds
\]
For $j=1$ we obtain the formula for $\Psi_1^+(R):=Rf_1(R)$ stated in \eqref{eq:sigma-1 exp}. 
The remainder of the proof of~\cite[Proposition 5.6]{KST} does not depend on the specific form of the potential and carries over verbatim to the case at hand. 
\end{proof}

To find the spectral measure 
 associated with $\calH^{+}_{-1}=\calH^{-}_{1}$ we now follow~\cite[Proposition 5.7]{KST}. 

\begin{proposition}
\label{prop: spectral measure n-1}
We have 
\begin{align*}
\Phi^{-1}(R,\xi)=2\Re\left(a_{-1}(\xi)\Psi^{+}_{-1}(R,\xi)\right)
\end{align*}
where $a(\xi)$ is smooth, always nonzero, and is of  size
\begin{align*}
|a_{-1}(\xi)|\simeq &\; \langle \xi\rangle^{-1}
\end{align*}
Moreover, it satisfies the symbol-type bounds
\begin{align*}
\left|(\xi\partial_{\xi})^{k}a_{-1}(\xi)\right|\leq c_{k}\quad \textrm{uniformly in }\quad \xi>0.
\end{align*}
The spectral measure of $-\calH_{-1}^+$ is absolutely continuous on $\xi\ge0$ with density
\begin{align*}
\frac{d\rho_{-1}(\xi)}{d\xi}=\frac{1}{4\pi}|a_{-1}(\xi)|^{-2}\simeq \langle\xi\rangle^{2}
\end{align*}
with the same zero order symbol bounds. 
\end{proposition}
\begin{proof}
We use the formalism introduced in \eqref{eq:phiapsi}--\eqref{eq:a rep rho}. I.e., we set 
\begin{align}\label{Phi Psi relation}
\Phi^{-1}(R,\xi)=a_{-1}(\xi)\Psi^{+}_{-1}(R,\xi)+\overline{a_{-1}(\xi)\Psi^{+}_{-1}(R,\xi)}
\end{align}
whence
\EQ{\label{eq:a-1W}
a_{-1}(\xi)=\frac{i}{2}W(\Phi^{-1},\overline{\Psi^{+}_{-1}}).
}
We evaluate this Wronskian in the regime where both the $\Phi^{-1}(R,\xi)$ and $\Psi^{+}_{-1}(R,\xi)$ asymptotics are useful, namely, $R\xi^{\frac{1}{2}}\simeq 1$. By 
Proposition~\ref{prop: basis n-1 small},  $|\Phi^{-1}(R,\xi)|\lesssim \xi^{-\frac{1}{4}}$ if~$R\simeq \xi^{-\frac{1}{2}}\ge1$ and $|\Phi^{-1}(R,\xi)|\lesssim \xi^{-\frac{5}{4}}$ if~$R\simeq \xi^{-\frac{1}{2}}\le1$.  From the equation $\calH_{-1}^+ \Phi^{-1}(R,\xi) = -\xi \Phi^{-1}(R,\xi)$ whence 
\[
|\partial_R^2\Phi^{-1}(R,\xi)|\lesssim (\xi+R^{-2}) |\Phi^{-1}(R,\xi)|  
\]
and thus $|\partial_R^2\Phi^{-1}(R,\xi)|\lesssim \xi |\Phi^{-1}(R,\xi)|$ for all $R\simeq \xi^{-\frac12}$. By calculus, $|\partial_R\Phi^{-1}(R,\xi)|\lesssim \xi^{\frac{1}{4}}$ if~$R\simeq \xi^{-\frac{1}{2}}\ge1$ and $|\partial_R\Phi^{-1}(R,\xi)|\lesssim \xi^{-\frac{3}{4}}$ if~$R\simeq \xi^{-\frac{1}{2}}\le1$. 
On the other hand, by Proposition \ref{prop: basis n-1 large},
\[
 |\Psi^{+}_{-1}(\xi^{-\frac{1}{2}},\xi)|\lesssim \xi^{-\frac14},\qquad |\partial_R\Psi^{+}_{-1}(\xi^{-\frac{1}{2}},\xi)|\lesssim 
\xi^{\frac14}
\]
whence 
\begin{align*}
|a_{-1}(\xi)|=&\left|\Phi^{-1}(\xi^{-\frac{1}{2}},\xi)\partial_{R}\Psi^{+}_{-1}(\xi^{-\frac{1}{2}},\xi)-\partial_{R}\Phi^{-1}(\xi^{-\frac{1}{2}},\xi)\Psi^{+}_{-1}(\xi^{-\frac{1}{2}},\xi)\right|
\lesssim \langle \xi\rangle^{-1}
\end{align*}
which gives the upper bound on $|a_{-1}(\xi)|$. The derivative bounds are obtained by differentiating the Wronskian~\eqref{eq:a-1W} in~$\xi$. For $ \Psi^{+}_{-1}$ one uses 
the symbol bounds from Proposition~\ref{prop: basis n-1 large}, and for $\Phi^{-1}(R,\xi)$ we use Proposition~\ref{prop: basis n-1 small} as well as the equation $\calH_{-1}^+ \Phi^{-1}(R,\xi)=-\xi \Phi^{-1}(R,\xi)$.

For the lower bound on $|a_{-1}(\xi)|$, we use a simpler argument than the one in~\cite{KST}. Directly from~\eqref{Phi Psi relation}, 
\[
|\Phi^{-1}(R,\xi)|\le 2 |a_{-1}(\xi)| |\Psi^{+}_{-1}(R,\xi)| ,\qquad |a_{-1}(\xi)| \ge\frac{|\Phi^{-1}(R,\xi)|}{2|\Psi^{+}_{-1}(R,\xi)|}
\]
for all $R>0$. On the one hand, $|\Psi^{+}_{-1}(R,\xi)|\le C\xi^{-\frac14}$ for all $\xi>0$ and $R^2\xi\ge \frac{1}{2}$ (say). 
From Proposition~\ref{prop: basis n-1 small} we have the lower bound 
\[
|\Phi^{-1}(R,\xi)| \ge \Phi^{-1}_{0}(R)\left( 1- \sum_{j=1}^{\infty}\frac{(R^{2}\xi)^{j}}{j!}\right) \ge \Phi^{-1}_{0}(R)\Big( 2- e^{R^{2}\xi}\Big) 
\]
So if, $R^2\xi=\frac12$, then 
\[
|a_{-1}(\xi)| \gtrsim \xi^{\frac14} \Phi^{-1}_{0}(1/\sqrt{2\xi})\simeq \langle\xi\rangle^{-1}
\]
matching the upper bound. For the spectral measure, see~\eqref{eq:a rep rho}.
\end{proof}

Given $f\in L^2(R\, dR)$, let $\tf(R)=R^{\frac{1}{2}}f(R)\in L^2(dR)$. Computing the distorted Fourier transform of $\tf$ relative to $-\calH_{1}^{-}$ gives,   at least formally 
\EQ{\label{eq:tfxxi}
\tf(R)=&\int_{0}^{\infty}x_{-1}(\xi)\Phi^{-1}(R,\xi)\, \rho_{-1}(d\xi),\qquad 
x_{-1}(\xi)=\int_{0}^{\infty}\tf(R)\Phi^{-1}(R,\xi)\, dR
}
The convergence of the first integral holds if $f$ is smooth and compactly supported in $(0,\infty)$ since then $x_{-1}(\xi)$ is rapidly decaying in~$\xi$, and bounded for small $\xi>0$. 
But in fact much less is needed.

\begin{lemma}
\label{lem:H-1conv}
Let $f\in L^2(R\,dR)\cap C^2((0,\infty))$ so that
\EQ{\label{eq:f f'' ass}
 \int_0^\infty  \big( R^{\frac32} |f''(R)| + R^{\frac12} |f'(R)|+R^{-\frac12} |f(R)|\big)\, dR =M<\infty
}
Then for all $\xi>0$ the limit 
\EQ{\label{eq:x-1 exists}
x_{-1}(\xi)=\lim_{L\to\infty} \int_{0}^{L} f(R)\Phi^{-1}(R,\xi)\, R^{\frac12}\, dR
}
exists and satisfies
\EQ{\label{eq:in xi abs}
\int_0^\infty |x_{-1}(\xi)| |\Phi^{-1}(R,\xi)|\, \rho_{-1}(d\xi) \lesssim M
}
In particular, the first integral in \eqref{eq:tfxxi} converges absolutely to $R^{\frac12}\,f(R)$ pointwise in $R>0$. 
\end{lemma}
\begin{proof}
By Propositions~\ref{prop: basis n-1 small}, \ref{prop: basis n-1 large}, and~\ref{prop: spectral measure n-1} we have
\EQ{\label{eq:Phi-1up}
\sup_{R>0} |\Phi^{-1}(R,\xi)|\lesssim \xi^{-\frac14}\langle\xi\rangle^{-1}
}
for all $\xi>0$. Indeed, Propositions~\ref{prop: basis n-1 small} implies that
  $$|\Phi^{-1}(R,\xi)|\lesssim \frac{R^{\frac52}}{1+R^2}\lesssim \frac{\xi^{-\frac54}}{1+\xi} \simeq \xi^{-\frac14}\langle\xi\rangle^{-1} $$
  for $R^2\xi\lesssim 1$, whereas Propositions~\ref{prop: basis n-1 large}, and~\ref{prop: spectral measure n-1} immediately imply the bound~\eqref{eq:Phi-1up} 
  in the range $R^2\xi\gtrsim 1$. It follows that
  \EQ{\label{eq:xi int}
  \int_0^\infty |x_{-1}(\xi)| |\Phi^{-1}(R,\xi)|\, \rho_{-1}(d\xi) &  \lesssim  \int_0^\infty |x_{-1}(\xi)| \xi^{-\frac14} \langle\xi\rangle\, d\xi
  }
  For any $\xi>0$ we define
  \EQ{
  \label{eq:xi split} 
  A(\xi) &= \int_0^\infty \chi_0(R^2\xi) \Phi^{-1}(R,\xi)\, R^{\frac12} f(R) \, dR  \\
  B_j(\xi) &= \int_0^\infty \chi(2^{-j} R^2\xi) \Phi^{-1}(R,\xi)\, R^{\frac12} f(R) \, dR 
  }
  relative to a partition of unity as in~\eqref{eq:POU}. Thus, at least formally, 
  \[
   x_{-1}(\xi) = A(\xi) + \sum_{j=0}^\infty B_j(\xi) \\
  \]
  We will show that the series does converge absolutely pointwise in $\xi>0$. 
 Then 
 \EQ{\label{eq:Abd}
 |A(\xi)| &\lesssim \int_0^{\xi^{-\frac12}} R^{\frac52}\langle R\rangle^{-2} |f(R)|R^{\frac12}\, dR.
 }
 Inserting the right-hand side into the upper bound of~\eqref{eq:xi int}  yields 
 \EQ{\nn
  \int_0^\infty |A(\xi)| \xi^{-\frac14} \langle\xi\rangle\, d\xi & \lesssim \int_0^\infty \int_0^{R^{-2}}  \xi^{-\frac14} \langle\xi\rangle\, d\xi R^{3}\langle R\rangle^{-2} |f(R)|\, dR \\
&\lesssim \int_0^\infty R^{-\frac12}|f(R)|\, dR
 }
 In the oscillatory regime we use the previous two propositions and integrate by parts twice. First, 
 \EQ{\nn 
 B_j(\xi) &= \int_0^\infty \chi(2^{-j} R^2\xi) \, \xi^{-\frac14} e^{i\xi^{\frac12}R} \sigma_{-1}(R\xi^{\frac12},R)a_{-1}(\xi)    R^{\frac12} f(R) \, dR  \\
 &\quad + \int_0^\infty \chi(2^{-j} R^2\xi) \, \xi^{-\frac14} e^{-i\xi^{\frac12}R} \; \overline{\sigma_{-1}(R\xi^{\frac12},R)a_{-1}(\xi) } \;  R^{\frac12} f(R) \, dR 
 }
 It suffices to deal with the first line. Thus, 
 \EQ{\label{eq:Bj bd} 
 |B_j(\xi)| &\lesssim \xi^{-\frac54}\langle\xi\rangle^{-1}\int_0^\infty |\partial_R^2 \big( \chi(2^{-j} R^2\xi)  \sigma_{-1}(R\xi^{\frac12},R)   R^{\frac12} f(R)  \big)|\, dR\\
 &\lesssim  \xi^{-\frac54}\langle\xi\rangle^{-1} \int_0^\infty \one_{[R^2\xi\simeq 2^j]} \big( R^{\frac12} |f''(R)| + R^{-\frac12} |f'(R)|+R^{-\frac32} |f(R)|\big)\, dR
 }
 Summing in $j\ge0$ and inserting this bound into~\eqref{eq:xi int} we obtain 
 \EQ{\nn
\int_0^\infty \sum_{j=0}^\infty |B_j(\xi)| \xi^{-\frac14} \langle\xi\rangle\, d\xi &\lesssim\int_0^\infty \sum_{j=0}^\infty  \xi^{-\frac32} \int_0^\infty \one_{[R^2\xi\simeq 2^j]} \big( R^{\frac12} |f''(R)| + R^{-\frac12} |f'(R)|+R^{-\frac32} |f(R)|\big)\, dR\,d\xi \\
&\lesssim \int_0^\infty  \big( R^{\frac32} |f''(R)| + R^{\frac12} |f'(R)|+R^{-\frac12} |f(R)|\big)\, dR 
 } 
 In summary, \eqref{eq:in xi abs} holds. On the other hand, the existence of the limit \eqref{eq:x-1 exists} is implicit in the argument, and in fact follows from the convergence of the series~\eqref{eq:xi split}  for fixed $\xi>0$. Indeed, summing in~\eqref{eq:Bj bd} one obtains 
 \[
 \sum_{j=0}^\infty  |B_j(\xi)| \lesssim \xi^{-\frac54}\langle\xi\rangle^{-1} \int_0^\infty\one_{[R^2\xi\gtrsim 1]} \big( R^{\frac12} |f''(R)| + R^{-\frac12} |f'(R)|+R^{-\frac32} |f(R)|\big)\, dR<\infty
 \]
 The finiteness in last line follows from the fact that $R\gtrsim \xi^{-\frac12}>0$ is bounded away from $0$, and from~\eqref{eq:f f'' ass}. 
\end{proof} 

We now derive the spectral representation of $\tilde H^{-}_{1}=\calD_{-} \calD_{-}^*$, cf.~\eqref{H1 star}. 
The generalized eigenbasis of $H^{-}_{1}$ is $R^{-\frac12}\Phi^{-1}(R,\xi)$, viz. 
\[
H^{-}_{1}\left(R^{-\frac{1}{2}}\Phi^{-1}(R,\xi)\right)= \;\; - \xi R^{-\frac{1}{2}}\Phi^{-1}(R,\xi)
\]
Returning to $f\in L^2(R\,dR)$ as in Lemma~\ref{lem:H-1conv},  we have in the sense stipulated by the lemma, 
\begin{align*}
f(R)=\int_{0}^{\infty}x_{-1}(\xi)R^{-\frac{1}{2}}\Phi^{-1}(R,\xi)\, \rho_{-1}(d\xi),\quad x_{-1}(\xi)=\left\langle R^{-\frac{1}{2}}\Phi^{-1}(R,\xi),f(R)\right\rangle_{L^{2}_{RdR}}.
\end{align*}
Applying $\calD_{-}$ to  the first of the above relations, we obtain at least formally
\begin{align*}
 \calD_{-}f(R)=&\int_{0}^{\infty}x_{-1}(\xi)\calD_{-}\left(R^{-\frac{1}{2}}\Phi^{-1}(R,\xi)\right)\, \rho_{-1}(d\xi) \\
 =&\int_{0}^{\infty}x_{-1}(\xi)\xi^{-1}\calD_{-}\left(R^{-\frac{1}{2}}\Phi^{-1}(R,\xi)\right)\, \trho_{-1}(d\xi) 
\end{align*}
where $\trho_{-1}(d\xi):=\xi\rho_{-1}(d\xi)$. In view of Proposition~\ref{prop: basis n-1 small} and the fact $\calD_{-}(R^{-\frac{1}{2}}\Phi^{-1}_0(R))=0$
it is natural to introduce the new basis 
\begin{align}\label{basis n-1 diff}
 \phi_{-1}(R,\xi):=\xi^{-1}\calD_{-}\left(R^{-\frac{1}{2}}\Phi^{-1}(R,\xi)\right),
\end{align}
which satisfies $\lim_{\xi\to0+}\phi_{-1}(R,\xi)=-\calD_{-}(R^2\Phi_1(R^2))$ and 
\begin{align*}
 \tilde H^{+}_{-1}\phi_{-1}(R,\xi)=\xi\phi_{-1}(R,\xi).
\end{align*}
So we arrive at the  Fourier representation formula for $\calD_{-}f$ relative to the operator $\tilde H_1^{-}=\tilde H_{-1}^+$: 
\EQ{\label{eq:D_f}
 \calD_{-}f(R)&=\int_{0}^{\infty}x_{-1}(\xi)\phi_{-1}(R,\xi)\, \trho_{-1}(d\xi),\quad 
 x_{-1}(\xi)=\left\langle\calD_{-}f,\phi_{-1}(\cdot,\xi)\right\rangle_{L^{2}_{RdR}}.
}
In the following  proposition we will state conditions on $f$ under which this Fourier transform converges. 
To obtain the second identity in~\eqref{eq:D_f}, we compute 
\begin{align*}
 \left\langle\calD_{-}f,\phi_{-1}(\cdot,\xi)\right\rangle_{L^{2}_{RdR}}=&\xi^{-1}\left\langle\calD_{-}f,\calD_{-}\left(R^{-\frac{1}{2}}\Phi^{-1}(\cdot,\xi)\right)\right\rangle_{L^{2}_{RdR}}\\
 =&-\xi^{-1}\left\langle f,H^{-}_{1}\left(R^{-\frac{1}{2}}\Phi^{-1}(\cdot,\xi)\right)\right\rangle_{L^{2}_{RdR}}
 =\left\langle f,R^{-\frac{1}{2}}\Phi^{-1}(\cdot,\xi)\right\rangle_{L^{2}_{RdR}}=x_{-1}(\xi).
\end{align*}
which is justified if $f\in C^1((0,\infty))$ has compact support. 

\begin{proposition}
\label{prop: spectral measure -1 diff}
Let $f\in  C^3((0,\infty))$ so that
\EQ{\label{eq:f f'' ass3}
 \int_0^\infty  \big( R  |f'''(R)| +  |f''(R)|+R^{-1} |f'(R)| + R^{-2}|f(R)|\big)\, dR =M<\infty
}
Then the  Fourier coefficient 
\[
y(\xi):=\left\langle\calD_{-}f,\phi_{-1}(\cdot,\xi)\right\rangle_{L^{2}_{RdR}} = \lim_{L\to\infty} \int_0^L \calD_{-}f(R) \phi_{-1}(R,\xi) R\, dR
\]
exists pointwise in $\xi>0$, and 
\[
\calD_{-}f(R)=\int_{0}^{\infty}y(\xi)\phi_{-1}(R,\xi)\, \trho_{-1}(d\xi)
\]
converges absolutely,  pointwise in $R>0$. Here  $\phi_{-1}(R,\xi)$ as in \eqref{basis n-1 diff} satisfies $ \tilde H^{+}_{-1}\phi_{-1}(R,\xi)=\xi\phi_{-1}(R,\xi)$. 
 The spectral measure $\trho_{-1}(\xi)$ associated to $\phi_{-1}(R,\xi)$ satisfies 
 $ \frac{d\trho_{-1}(\xi)}{d\xi}\simeq \xi \langle\xi\rangle^{2}$ for all $\xi>0$.  Moreover, we have the unitarity 
 $$\|\calD_{-}f\|_{L^2(R\,dR)}=\big\| \left\langle\calD_{-}f,\phi_{-1}(\cdot,\xi)\right\rangle_{L^{2}_{RdR}} \big\|_{L^2(\tilde\rho_{-1})}$$
 provided the left-hand side is finite. 
\end{proposition}
\begin{proof}
The properties of the spectral measure follow directly from Proposition \ref{prop: spectral measure n-1}. Recall from \eqref{D minus and star} that
\[\calD_{-}f(R)=f'(R)-\frac{2}{R}f(R)+\frac{2R}{1+R^{2}}f(R)
\]
and $|\calD_{-}f(R)|\lesssim |f'(R)|+R^{-1}|f(R)|$. From Proposition~\ref{prop: basis n-1 small} and $\calD_{-}(R^{-\frac{1}{2}}\Phi^{-1}_0(R))=0$, we infer that 
$$|\phi_{-1}(R,\xi)|\lesssim R^3\langle R^2\rangle^{-1}\lesssim \xi^{-\frac12}\langle\xi\rangle^{-1}\text{\ \ if\ \ }R^2\xi\le1$$
Similarly, from  Propositions~\ref{prop: spectral measure n-1} and~\ref{prop: basis n-1 large}, we obtain the same bound in the regime $R^2\xi\ge1$. 
The remainder of the proof now is analogous to that of Lemma~\ref{lem:H-1conv} and we do not write it out. The unitarity is proved in~\cite{GZ}. 
\end{proof} 

\subsection{Turning point analysis and the spectral measure for $n\geq 2$}\label{sec: ngeq2 linear}
The goal of this subsection is to find the asymptotic behavior of the spectral measures for small and large energies. This will be done by means of a careful asymptotic analysis of a fundamental system for all energies. We begin with $\calH_{n}^{+}$, and will present the needed modifications for $\calH^{-}_{n}=\calH^+_{-n}$ later. We will use a perturbation argument developed in the papers \cite{CSST} and \cite{CDST}. As in those papers, we write the eigenvalue problem $-\calH^{+}_{n}f=E^2f$ as follows:
\begin{align}\label{Schlag eq}
\begin{split}
 -\frac{1}{(n+1)^{2}}\partial_{R}^{2}f+V(R)f=&\frac{E^{2}}{(n+1)^{2}} f,\\ V(R):=&\frac{1}{R^{2}}-\frac{1}{4(n+1)^{2}R^{2}}-\frac{4n}{(n+1)^{2}}\frac{1}{R^{2}(R^{2}+1)}-\frac{8}{(n+1)^{2}}\frac{1}{(R^{2}+1)^{2}}.
 \end{split}
\end{align}
Switching to semiclassical notation, we introduce  $\hbar:=\frac{1}{n+1}$ and write $V(R)=V_n(R)=V(R;\hbar)$ as
\begin{align}\label{Schlag eq potential alter}
 \begin{split}
  V(R)=V(R;\hbar)=&\frac{1}{R^{2}}\left(1-\frac{\hbar^{2}}{4}-\frac{4\hbar}{R^{2}+1}+\frac{4\hbar^{2}(1-R^{2})}{(R^{2}+1)^{2}}\right)\\
  :=&\frac{1}{R^{2}}\left(1+\frac{15\hbar^{2}}{4}-4\hbar\right)+\veps(R^{2};\hbar)
 \end{split}
\end{align}
Here

\begin{align*}
 \veps(R^{2};\hbar):=\frac{4\hbar}{R^{2}+1}-\frac{4\hbar^{2}(R^{2}+3)}{(R^{2}+1)^{2}}
\end{align*}
is a bounded smooth function on $[0,\infty)$. Henceforth, it will be understood that $\hbar\in(0, \frac{1}{3}]$. 

We will construct a fundamental system for \eqref{Schlag eq} on $R>0$ for all $E>0$. We first scale $E$ out by introducing $x:=\hbar E R$. If we define $\tf(x):=f(R)$, then \eqref{Schlag eq} becomes

\begin{align}\label{Schlag eq alt}
 \begin{split}
  -\hbar^{2}\tf''(x)+Q(x)\tf(x)=0,\quad Q(x):=\hbar^{-2}E^{-2}V\left(\frac{x}{\hbar E}\right)-1.
 \end{split}
\end{align}
More precisely, we have, with $\alpha:=\hbar E$,
\begin{align*}
 Q(x,\alpha;\hbar)=&x^{-2}\left(1+\frac{15\hbar^{2}}{4}-4\hbar\right)+\alpha^{-2}\veps\left(\frac{x^{2}}{\alpha^{2}};\hbar\right)-1.
\end{align*}
As usual, see for example \cite{CSST}, we need to modify the potential by adding the \emph{Langer correction}: 
\begin{align}\label{def Q0}
 Q_{0}(x;\alpha,\hbar):=Q(x;\alpha,\hbar)+\frac{\hbar^{2}}{4x^{2}}=x^{-2}\left(1-2\hbar\right)^2+\alpha^{-2}\veps\left(\frac{x^{2}}{\alpha^{2}};\hbar\right)-1.
\end{align}
This modification is crucial in order to obtain the correct asymptotics as $E\to0+$ in the ensuing WKB analysis.  In the following lemma we consider $Q_0$ as a function of $x^2$ rather than~$x$, and $\alpha^2$ instead of~$\alpha$. This lemma is sharp in the sense that it does not hold in the full range $\hbar\in(0,\frac12 ]$. 

\begin{lemma}
\label{lem:Q0diff}
The function $Q_1(x^2;\alpha^2,\hbar)=Q_0(x;\alpha,\hbar)$ satisfies
 \begin{align*}
 - (\partial_y Q_1)(y;a,\hbar)  &\simeq y^{-2}, \quad  (\partial_a Q_1)(y;a,\hbar) <0 \\
 |\partial_a Q_1(y;a,\hbar)| &\lesssim  \hbar (y+a)^{-2}
 \end{align*}
 uniformly in $y>0, a>0$ and $\hbar\in(0,\frac13]$. The implied constants are absolute.  The higher order derivatives satisfy, uniformly in  $\hbar\in(0,\frac13]$
  \begin{align*}
 |(\partial_y^k Q_1)(y;a,\hbar) | &\le C_k\,  y^{-k-1} \text{\ \ for all \ \ } y>0, a>0, k\ge 2 \end{align*}
 and 
  \begin{align*}
 |\partial_y^k \partial_a^\ell  Q_1(y;a,\hbar)| &\le C_{k,\ell}\,\hbar\,  (y+a)^{-k-\ell-1}
\end{align*}
 for all $y>0, a>0$, $k\ge 0, \ell\ge1$. 
\end{lemma}
\begin{proof}
With 
 \begin{align*}
 Q_{1}(y,a)=&(2\hbar-1)^{2}y^{-1}+\frac{4\hbar(1-3\hbar)}{y+a}+\frac{8\hbar^{2}y}{(y+a)^{2}}-1
 \end{align*}
 we have, with $\partial_{u}Q_{1}$ denoting the partial derivative of $Q_{1}$ with respect to its first variable,
 \begin{equation}
 \label{eq:pyQ}
 \begin{aligned}
 -a^{2} (\partial_u Q_1)(ay,a) &= (2\hbar-1)^{2}y^{-2}+\frac{4\hbar(1-3\hbar)}{(1+y)^2}+\frac{8\hbar^{2}(y-1)}{(1+y)^{3}} \\
&= y^{-2} \Big[(2\hbar-1)^{2} + 4\hbar y^2(1+y)^{-3}\big( (1-3\hbar)(1+y) +  2\hbar   (y-1)  \big) \Big] \\
&= y^{-2} \Big[(2\hbar-1)^{2} + 4\hbar y^2(1+y)^{-3}\big( 1-5\hbar +  (1-\hbar )y \big) \Big] 
\end{aligned}
\end{equation}
We claim that the term in brackets on the right-hand side is $\simeq 1$ uniformly in $y>0$ and $\hbar\in(0, \frac{1}{3}]$. 
 By inspection, this is the case if $y\ge1$ or $y>0$ and $\hbar\in (0,\frac15]$. It thus suffices to check that the polynomial 
 \[
 P_\hbar(y) = (1-2\hbar)^{2}(1+y)^{3} + 4\hbar y^2\big( 1-5\hbar +  (1-\hbar )y  \big) \simeq 1
 \]
 uniformly on the rectangle $(y,\hbar)\in [0,1]\times [\frac15,\frac13]$. We have
 \[
 P_\hbar'(y) = 3 (1 + y)^2 - 4 \hbar^2 (-3 + 4 y) - 4 \hbar (3 + 4 y)
 \]
 The discriminant of the quadratic polynomial $P_h'$ is 
 \[
 D=-3 \hbar - 5 \hbar^2 + 32 \hbar^3 + 16 \hbar^4
 \]
 On checks that $D<0$ for all $\hbar\in  [\frac15,\frac13]$. Thus $P_\hbar'(y)$ does not change sign in $y\in \R$ for all such~$\hbar$. 
 Clearly, $P_h'(y)>0$ for large $y$  whence $P_\hbar(y)$ is increasing for all $\hbar\in  [\frac15,\frac13]$. In summary, for all $y\in [0,1]$
 \[
 P_\hbar(0)=(1-2\hbar)^{2} \le  P_\hbar(y) \le  P_\hbar(1)=8(1-2\hbar)^{2}+8\hbar \big( 1-3\hbar)
 \]
 Thus the claim above holds and from \eqref{eq:pyQ} we have 
 \[
 - (\partial_y Q_1)(y,a)  \simeq y^{-2} 
 \]
 uniformly in $y>0, a>0$ and $\hbar\in(0,\frac13]$. 
 
For the $a$-derivative we compute
\[
0<- (\partial_a Q_1)(y,a)  =  \frac{4\hbar(1-3\hbar)}{(y+a)^2}+\frac{16\hbar^{2}y}{(y+a)^{3}} \le 6\hbar  (y+a)^{-2}
\]
uniformly in $y>0, a>0$ and $\hbar\in(0,\frac13]$.  The estimates for the higher order derivatives follow by using the fact that $Q_{1}(y;a,\hbar)$ is a rational function for both $y$ and $a$.
\end{proof}

In the original variables the previous lemma takes the following form. 

\begin{corollary}
\label{cor:Q0diff}
One has 
 \EQ{\label{eq:palQ0}
 - (\partial_x Q_0)(x; \alpha,\hbar)  &\simeq x^{-3}, \quad  (\partial_\alpha Q_0)(x;\alpha,\hbar) <0 \\
 |\partial_\alpha Q_0(x;\alpha,\hbar)| &\lesssim  \hbar\, \alpha (x^2+\alpha^2)^{-2}
}
 uniformly in $x>0,\alpha >0$ and $\hbar\in(0,\frac13]$. The implied constants are absolute.  The higher order derivatives satisfy, uniformly in  $\hbar\in(0,\frac13]$
  \begin{align*}
 |(\partial_x^k Q_0)(x;\alpha,\hbar) | &\le C_k\,  x^{-k-2} \text{\ \ for all \ \ } x>0, \alpha >0, k\ge 2 \end{align*}
 and 
  \EQ{\label{eq:Q0kell}
 |\partial_x^k \partial_\alpha^\ell  Q_0(x;\alpha,\hbar)| &\le C_{k,\ell}\,\hbar\,   (x+\alpha)^{-k-\ell-2}
}
 for all $x>0, \alpha>0$, $k\ge 0, \ell\ge1$. 
\end{corollary}
\begin{proof}
We set $a=\alpha^2$ and $y=x^2$. By the chain rule 
\[
(2\alpha)^{-1} \partial_\alpha = \partial_a,\qquad  (2x)^{-1} \partial_x =\partial_y
\]
and  the corollary follows from the previous lemma. For example, the second resp.\ third  derivatives are 
\[
 \partial_\alpha^2 = 2\partial_a + 4a \partial_a^2,\quad  \partial_\alpha^3 = \alpha(12\partial_a^2 + 8\alpha^2\partial_a^3)
\]
and similarly in $x$. 
\end{proof}

We remark that the estimate in \eqref{eq:Q0kell} is not optimal if $k$ or $\ell$ are odd. Indeed, in that case one has vanishing at $x=0$, respectively $\alpha=0$ as in~\eqref{eq:palQ0}.
The following lemma introduces the unique turning point $x_{t}(\alpha,\hbar)$, i.e., the root of $Q_{0}(x;\alpha,\hbar)=0$. The same non-optimality remark again applies to~\eqref{eq:xt ell} below for odd~$\ell$. 

\begin{lemma}\label{lem: monotonicity of root in alpha}
 $Q_{0}(x;\alpha,\hbar)=0$ has a unique root $x_{t}(\alpha;\hbar)\in(1-2\hbar,1)$. It is strictly monotone decreasing and smooth in $\alpha>0$. Moreover, 
 \begin{align}
-\partial_\alpha x_t(\alpha;\hbar) & = |\partial_\alpha x_t(\alpha;\hbar)| \lesssim  \hbar \, \alpha (1+\alpha)^{-4} \nonumber \\
|\partial_\alpha^\ell x_t(\alpha;\hbar)| & \lesssim  \hbar\,   (1+\alpha)^{-\ell-2},\qquad \ell\ge1 \label{eq:xt ell}
\end{align}
for all $\alpha>0$ with a uniform constant in $\hbar\in (0,\frac13]$. 
\end{lemma}
\begin{proof}
By the previous lemma 
$
Q_{0}(x; \alpha,\hbar)$ is strictly monotone  decreasing in both $x$ and $\alpha$.  Thus, $x_t$ exists uniquely and is strictly decreasing in $\alpha$. 
 Thus, to determine the range of $x_t$ as a function of $\alpha$, it suffices to consider the limits $\alpha\rightarrow0^{+}$ and $\alpha\rightarrow\infty$. 
 When $\alpha=0$, we have
 \begin{align*}
  Q_{0}(x;0,\hbar)=x^{-2}-1,\quad \quad x_{t}(0;\hbar)=1.
 \end{align*}
When $\alpha=\infty$, we have 
\begin{align*}
 Q_{0}(x;\infty,\hbar)=(2\hbar-1)^{2}x^{-2}-1,\quad \quad x_{t}(\infty;\hbar)=1-2\hbar.
\end{align*}
Taking an $\alpha$ derivative of $$Q_0(x_t(\alpha;\hbar);\alpha,\hbar)=Q_1(x_t(\alpha;\hbar)^2;\alpha^2,\hbar)=0$$
yields 
\begin{equation}\label{eq:dxt} 
\partial_\alpha x_t(\alpha;\hbar) (\partial_x Q_0)(x_t(\alpha;\hbar);\alpha,\hbar) +  \partial_\alpha Q_0(x_t(\alpha;\hbar);\alpha,\hbar)=0
\end{equation}
whence 
\[
-\partial_\alpha x_t(\alpha;\hbar) = |\partial_\alpha x_t(\alpha;\hbar)| \lesssim \hbar  \alpha (1+\alpha^2)^{-2}
\]
Taking another $\alpha$ derivative in \eqref{eq:dxt} yields 
\[
|\partial_\alpha^2 x_t(\alpha;\hbar)| \lesssim  \hbar (1+\alpha)^{-4}
\]
as claimed. By Leibnitz's rule, $\partial_\alpha^\ell x_t(\alpha;\hbar)$  is a linear combination of all terms of the form 
\[
\partial_\alpha^{\nu_1} x_t \cdot\ldots\cdot \partial_\alpha^{\nu_r} x_t \;  \partial_x^r \partial_\alpha^m Q_0
\]
evaluated at $x=x_t(\alpha;\hbar)$, where $\nu_1+\ldots+\nu_r+m=\ell$. By induction, if \eqref{eq:xt ell} holds, then if $\nu_i<\ell$ for all $i$, and if $r\ge1$, 
\begin{align*}
 |\partial_\alpha^{\nu_1} x_t \cdot\ldots\cdot \partial_\alpha^{\nu_r} x_t \;  \partial_x^r \partial_\alpha^m Q_0(x_t;\alpha,\hbar)| 
&\lesssim \hbar^{r} (1+\alpha)^{-(\nu_1+\ldots+\nu_r+2r)} \,  (1+\alpha)^{-r-m-2} \\
&\lesssim \hbar^{r}  (1+\alpha)^{-\ell -3r -2} 
\end{align*}
which is largest when $r=1$ (if $m>0$, then one has $\hbar^{r+1}$ and not $\hbar^{r}$). Thus, 
\begin{align*}
| \partial_\alpha^{\ell} x_t (\alpha,\hbar) \partial_x Q_0(x_t(\alpha;\hbar);\alpha,\hbar)| &\lesssim \hbar  (1+\alpha)^{-\ell -5 } + |\partial_\alpha^\ell Q_0(x_t(\alpha;\hbar);\alpha,\hbar)|  \\
&\lesssim \hbar  (1+\alpha)^{-\ell-2}.
\end{align*}
By induction, \eqref{eq:xt ell} holds for all $\ell$ as claimed. 
 \end{proof}

\subsection{Liouville-Green transform, reduction to Airy's equation}\label{sec:transfAiry}

We will now apply a global Liouville-Green transform \footnote{{For more background on this Langer transform and turning point theory in general, see \cite{MillerP} and \cite{Wasow}.}} to the equation~\eqref{Schlag eq alt}, cf.~\cite[Chapter 6]{Olver}. As usual, this refers to a change of both the independent and dependent variables. 
The new independent variable $\tau=\tau(x,\alpha;\hbar)$ is given by 
\begin{align}\label{def tau}
 \tau(x,\alpha;\hbar):=\sign\left(x-x_{t}(\alpha;\hbar)\right)\left|\frac{3}{2}\int_{x_{t}(\alpha;,\hbar)}^{x}\sqrt{\left|Q_{0}(u,\alpha;\hbar)\right|}\, du\right|^{\frac{2}{3}}.
\end{align}
Before we rigorously analyze $\tau$, we first take a look at how the equation \eqref{Schlag eq alt} transforms under  this change of variables. As we have seen, for $\hbar\le \frac{1}{3}$, $Q_{0}(x)$ is strictly monotone for $x>0$, so is $\tau$ by its definition. Therefore the map $\tau(\cdot,\alpha,\hbar):(0,\infty)\rightarrow\bbR$ is injective. 
When $x>x_{t}(\alpha,\hbar)$, we have
\begin{align*}
 \tau(x;\alpha,\hbar)=&\left(\frac{3}{2}\int_{x_{t}(\alpha,\hbar)}^{x}\sqrt{-Q_{0}(u;\alpha,\hbar)}\, du\right)^{\frac{2}{3}}\\
 \frac{d\tau}{dx}=&\left(\frac{3}{2}\int_{x_{t}(\alpha,\hbar)}^{x}\sqrt{-Q_{0}(u;\alpha,\hbar)}\, du\right)^{-\frac{1}{3}}
 \cdot\sqrt{-Q_{0}(x;\alpha,\hbar)}
 \end{align*}
 and 
 \begin{align*}
 \frac{d^{2}\tau}{dx^{2}}=&\frac{1}{2}\left(\frac{3}{2}\int_{x_{t}(\alpha,\hbar)}^{x}\sqrt{-Q_{0}(u;\alpha,\hbar)}\, du\right)^{-\frac{4}{3}}Q_{0}(x;\alpha,\hbar)\\
 &-\left(\frac{3}{2}\int_{x_{t}(\alpha,\hbar)}^{x}\sqrt{-Q_{0}(u;\alpha,\hbar)}du\right)^{-\frac{1}{3}}\frac{1}{2\sqrt{-Q_{0}(x;\alpha,\hbar)}}\frac{dQ_{0}(x;\alpha,\hbar)}{dx}
\end{align*}
For simplicity, we will use ``$\prime$'' to denote the derivative with respect to $x$ and ``$\cdot$'' to denote the derivative with respect to $\tau$. Also, we will sometimes suppress the dependence on $\alpha$ and $\hbar$ when there is no confusion. Following \cite{CSST}, we define
\begin{equation}\label{eq:qdef}
 q:=-\frac{Q_{0}}{\tau} \text{\ \ so that\ \ } \frac{d\tau}{dx}=\tau'=\sqrt{q},\quad   \frac{d}{d\tau}=q^{-\frac{1}{2}}\frac{d}{dx}.
\end{equation}
The new dependent variable is given by $w:=\sqrt{\tau'}\tf$. Then, with $\dot{q}^{\frac{1}{4}}=\partial_\tau({q}^{\frac{1}{4}})$, $\ddot{q}^{\frac{1}{4}}=\partial_\tau^2({q}^{\frac{1}{4}})$ we have 
\begin{align*}
 &w:=\sqrt{\tau'}\tf=q^{\frac{1}{4}}\tf,\quad  \dot{w}=\dot{q}^{\frac{1}{4}}\tf+q^{-\frac{1}{4}}\tf'\\
 &\ddot{w}=\ddot{q}^{\frac{1}{4}}\tf+\dot{q}^{\frac{1}{4}}q^{-\frac{1}{2}}\tf'-q^{-\frac{1}{2}}\dot{q}^{\frac{1}{4}}\tf'+q^{-\frac{3}{4}}\tf''=q^{-\frac{3}{4}}\tf''+\ddot{q}^{\frac{1}{4}}\tf.
\end{align*}
 Using the equation \eqref{Schlag eq alt}, we obtain
\begin{align}\label{perturbed Airy pre}
 \begin{split}
  &\ddot{w}=q^{-\frac{3}{4}}\left(\frac{Q_{0}}{\hbar^{2}}-\frac{1}{4x^{2}}\right)\tf+\ddot{q}^{\frac{1}{4}}\tf=-\tau\hbar^{-2}w+\ddot{q}^{\frac{1}{4}}q^{-\frac{1}{4}}w-\frac{1}{4x^{2}}q^{-1}w\\
   &-\hbar^{2}\ddot{w}(\tau)=\tau w(\tau)-\hbar^{2}\ddot{q}^{\frac{1}{4}}q^{-\frac{1}{4}}w(\tau)+\hbar^{2}\frac{1}{4x^{2}q}w\\
  &-\hbar^{2}\ddot{w}(\tau)=:\tau w(\tau)-\hbar^{2}\tV(\tau;\alpha,\hbar) w(\tau)
 \end{split}
\end{align}
where 
\begin{equation}\label{eq:master}
\tV(\tau;\alpha,\hbar) = \ddot{q}^{\frac{1}{4}}q^{-\frac{1}{4}} - \frac{1}{4x^{2}q}.
\end{equation}
In order to use the results in \cite{CDST}, we still need to modify the last equation in \eqref{perturbed Airy pre} such that it is consistent with the equation (D.9) in \cite{CDST}. To this end, we introduce the new variable $\zeta:=-\tau$ and the last equation in \eqref{perturbed Airy pre} becomes
\begin{align}\label{perturbed Airy}
 \hbar^{2}\ddot{w}(\zeta)=\zeta w(\zeta)+\hbar^{2}\tV(-\zeta,\alpha;\hbar)w(\zeta).
\end{align}
The case for $x<x_{t}(\alpha,\hbar)$ can be handled similarly and results in the same equation. Below  we will analyze the behavior of $\tau$ in terms of $x$ in different regimes. In particular, there are three regimes to consider: 

\begin{itemize}
 \item When $x$ is close to $x_{t}(\alpha,\hbar)$. More precisely, $|x-x_{t}(\alpha,\hbar)|\le \frac{1}{2}x_t(\alpha;\hbar)$.
 \item When $x\rightarrow 0^{+}$. More precisely, $0<x\leq \frac12 x_{t}(\alpha;\hbar)$.
 \item When $x\rightarrow\infty$. More precisely, $x\ge \frac32 x_t(\alpha;\hbar)$.
\end{itemize}
On the other hand, the size of $\alpha=\hbar E$ affects the position of $x_{t}(\alpha,\hbar)$. Sometimes we distinguish the two different cases: $\alpha\ll1$ and $\alpha\gtrsim 1$. This will play a crucial role when $x\in(0,\frac{1}{2}]$.

Now we can start our analysis on the behavior of $\tau$ and the potential $\tV$ in the perturbed Airy equation \eqref{perturbed Airy}.

\begin{lemma}
\label{lem: Lemma 3.2 CDST}
For $x\in\left[\frac{1}{2}x_{t}(\alpha;\hbar), \frac{3}{2}x_{t}(\alpha;\hbar) \right]=:J_1(\alpha;\hbar)$,  the function $\tau$ defined in \eqref{def tau} satisfies
\begin{align*}
 \tau(x;\alpha,\hbar)   = (x-x_t(\alpha;\hbar)) \Phi(x;\alpha,\hbar) 
 \end{align*}
 where $\Phi(x;\alpha,\hbar)\simeq 1$ uniformly in $x\in J_1$, $\alpha>0$ and $\hbar\in(0,\frac13]$.  Moreover, uniformly in that range of parameters, 
  \begin{equation}\label{eq:kell g0}
 |\partial_x^k  \Phi(x;\alpha,\hbar)| \lesssim 1,\quad |\partial_x^k  \partial_\alpha^\ell \Phi(x;\alpha,\hbar)| \lesssim  \hbar\,  (1+\alpha)^{-\ell-2}
 \end{equation}
 for all $k\ge0$, $\ell\ge1$. Finally, $q=q(x;\alpha,\hbar)$ as in~\eqref{eq:qdef} satisfies the exact same properties. 
\end{lemma}
\begin{proof}
 We have
 \EQ{\label{eq:Q0 zero}
- Q_0(x;\alpha,\hbar)   &=  Q_0(x_t(\alpha;\hbar);\alpha,\hbar) - Q_0(x;\alpha,\hbar)  \\
 &= (x-x_t(\alpha;\hbar))\int_0^1(-\partial_x Q_0)(x_t(\alpha;\hbar)+s(x-x_t(\alpha;\hbar)) ;\alpha,\hbar)  \, ds \\
 & = (x-x_t(\alpha;\hbar)) g(x;\alpha,\hbar)
}
 By Corollary~\ref{cor:Q0diff} and Lemma~\ref{lem: monotonicity of root in alpha} we have $g(x;\alpha,\hbar)\simeq 1$ uniformly in $x\in J_1$, $\alpha>0$ and $\hbar\in(0,\frac13]$.  Moreover, uniformly in that range of parameters, 
 \begin{equation}\label{eq:kell g}
 |\partial_x^k  g(x;\alpha,\hbar)| \lesssim 1,\quad |\partial_x^k  \partial_\alpha^\ell g(x;\alpha,\hbar)| \lesssim  \hbar\, (1+\alpha)^{-\ell-2}
 \end{equation}
 for all $k\ge0$, $\ell\ge1$.  We next show how the estimates for higher order derivatives of $g$ are derived, in the case of pure $x$ or $\alpha$ derivatives, the case of mixed derivatives following similarly: 
 \begin{align*}
&	\partial_{x}g(x;\alpha,\hbar)=\int_{0}^{1}(-s\partial_{x}^{2}Q_{0})(x_{t}(\alpha,\hbar)+s(x-x_{t}(\alpha,\hbar));\alpha,\hbar)ds\\
\Rightarrow\quad &\partial_{x}^{k}g(x;\alpha,\hbar)=\int_{0}^{1}(-s^{k}\partial_{x}^{k+1}Q_{0})(x_{t}(\alpha,\hbar)+s(x-x_{t}(\alpha,\hbar));\alpha,\hbar)ds.
\end{align*} 
The estimate on $|\partial_{x}^{k}g(x;\alpha,\hbar)|$ then follows from Corollary~\ref{cor:Q0diff}. The $\alpha$-derivatives are given as follows:
\begin{align*}
	\partial_{\alpha}g(x;\alpha,\hbar)=&\int_{0}^{1}(-(1-s)\partial_{x}^{2}Q_{0})(x_{t}(\alpha,\hbar)+s(x-x_{t}(\alpha,\hbar));\alpha,\hbar)(\partial_{\alpha}x_{t})(\alpha,\hbar)ds\\&-\int_{0}^{1}(\partial_{x}\partial_{\alpha}Q_{0})(x_{t}(\alpha,\hbar)+s(x-x_{t}(\alpha,\hbar));\alpha,\hbar)ds,
\end{align*}
and more generally, using Fa\`a di Bruno's formula as well as Leibniz's rule
\begin{align*}
	&\partial^{\ell}_{\alpha}g(x;\alpha,\hbar)\\
	=&-\sum_{0\leq p'\leq p\leq \ell,\atop \sum p_j=p
	}C_{\{p_j\},p,\ell}\int_{0}^{1}\left((1-s)^{p'}\partial_{x}^{p'+1}\partial_{\alpha}^{\ell-p}Q_{0}\right)(x_{t}(\alpha,\hbar)+s(x-x_{t}(\alpha,\hbar));\alpha,\hbar)\prod_{j=1}^{p'}(\partial_{\alpha}^{p_j}x_{t})(\alpha,\hbar)ds
\end{align*}
Therefore the desired estimates follow from Corollary~\ref{cor:Q0diff} as well as Lemma~\ref{lem: monotonicity of root in alpha}.
  Inserting \eqref{eq:Q0 zero} into \eqref{def tau} yields 
 \begin{align}\label{def Phi}
 \begin{split}
 \tau(x;\alpha,\hbar) &= (x-x_t(\alpha;\hbar)) \Big|\frac32 \int_0^1 \sqrt{sg(x_t (\alpha;\hbar)+ s(x-x_t (\alpha;\hbar));\alpha,\hbar)}\, ds\Big|^{\frac23}\\
 &=: (x-x_t(\alpha;\hbar)) \Phi(x;\alpha,\hbar)
 \end{split}
 \end{align}
 where $\Phi$ satisfies the same estimates as $g$, namely $\Phi\simeq 1$ uniformly in the parameter range above and the derivatives satisfy~\eqref{eq:kell g}. 
 As a result, 
 \[
 q(x;\alpha,\hbar) = g(x;\alpha,\hbar)/\Phi(x;\alpha,\hbar),\text{\ \ and\ \ } \tau'(x;\alpha,\hbar) = \sqrt{q}
 \]
 also satisfies these exact same properties. 
 \end{proof}
 
For the potential function
\begin{align}\label{tV def}
 \tV=-\frac{1}{4x^{2}q}+q^{-\frac{1}{4}}\ddot{q}^{\frac{1}{4}},
\end{align}
see \eqref{perturbed Airy} and \eqref{eq:master} we have the following immediate corollary.

\begin{proposition}\label{prop: tV close xt}
Uniformly in the parameter range of  Lemma~\ref{lem: Lemma 3.2 CDST} we have
 \begin{equation}
   |\partial_\tau^k   \tV(\tau;\alpha,\hbar)| \lesssim 1,\quad |\partial_\tau^k  \partial_\alpha^\ell  \tV(\tau;\alpha,\hbar)| \lesssim  \hbar\,  (1+\alpha)^{-\ell-2}
 \end{equation}
 for all $k\ge0$, $\ell\ge1$. 
\end{proposition}
\begin{proof}
This follows from the properties of $q$ stated in the previous result. 
\end{proof}

Next we turn to the case when $x\in J_2:=[\frac32 x_t(\alpha;\hbar),\infty)$. 

\begin{lemma}
\label{lem: Lemma 3.4 CDST}
For $x\in J_2$ we write  $\tau = \Big(\frac32\ttau\Big)^{\frac23}$. Then $\ttau>0$  satisfies
\EQ{\nn 
\ttau(x;\alpha,\hbar)  & = x-y(\alpha;\hbar) + \rho(x;\alpha,\hbar) \\
y(\alpha;\hbar) & \simeq1,\quad |\partial_\alpha^\ell y(\alpha;\hbar)|\lesssim \hbar \, (1+\alpha)^{-\ell-2},\quad\ell\ge1 \\
\rho(x;\alpha,\hbar)&\simeq x^{-1},\quad |\partial_x^k  \rho(x;\alpha,\hbar) |\lesssim x^{-k-1}, \quad k\ge 0 \\
|\partial_x^k \partial_\alpha^\ell \rho(x;\alpha,\hbar) | &\lesssim \hbar \, (x+\alpha)^{-\ell-k-1}, \quad k\ge 0, \ell\ge 1
}
uniformly in  $ \alpha>0,0<\hbar\leq\frac{1}{3}$.
\end{lemma}
\begin{proof}
We have 
\EQ{\label{eq:Q0darst}
-Q_0(x;\alpha,\hbar) = 1 - x^{-2}(1-2\hbar)^2 - \frac{4\hbar}{x^2+\alpha^2} + \frac{4\hbar^2(x^2+3\alpha^2)}{(x^2+\alpha^2)^2}
}
whence
\begin{align*}
\ttau(x;\alpha,\hbar) &:= \int_{x_t(\alpha;\hbar)}^x \sqrt{-Q_0(u;\alpha,\hbar)}\, du = x-x_t(\alpha;\hbar) + \int_{x_t(\alpha;\hbar)}^x \Big( \sqrt{-Q_0(u;\alpha,\hbar)}-1\Big)\, du  \\
& = x-x_t(\alpha;\hbar) - \kappa(\alpha;\hbar)+ \rho(x;\alpha,\hbar)
\end{align*}
where
\EQ{\label{eq:rho def}
\rho(x;\alpha,\hbar):=\int_{x}^\infty \frac{u^{-2}(1-2\hbar)^2 + \frac{4\hbar(1-3\hbar)}{u^2+\alpha^2} + \frac{8\hbar^2u^2}{(u^2+\alpha^2)^2}}{1+\sqrt{-Q_0(u;\alpha,\hbar)}}\, du 
}
and 
\EQ{\nn
\kappa(\alpha;\hbar)&:= \int_{x_t(\alpha;\hbar)}^\infty \frac{u^{-2}(1-2\hbar)^2 + \frac{4\hbar(1-3\hbar)}{u^2+\alpha^2} + \frac{8\hbar^2 u^2}{(u^2+\alpha^2)^2}}{1+\sqrt{-Q_0(u;\alpha,\hbar)}}\, du \\
&= \frac{1}{x_t(\alpha;\hbar)} \int_{1}^\infty \frac{ v^{-2}(1-2\hbar)^2 + \frac{4\hbar(1-3\hbar)}{ v^2+\beta^2} + \frac{8\hbar^2 v^2}{(v^2+\beta^2)^2}}{1+\sqrt{-Q_0(x_t(\alpha;\hbar)v ;\alpha,\hbar)}}\, dv
}
with $\beta(\alpha;\hbar)=\frac{\alpha}{x_t(\alpha;\hbar)}$.    By the results of Section~\ref{sec:transfAiry}, uniformly in $\alpha>0$ and $\hbar\in(0,\frac13]$,
\[
\kappa(\alpha;\hbar)\simeq1,\quad |\partial_\alpha^\ell \kappa(\alpha;\hbar)|\lesssim \hbar \, (1+\alpha)^{-\ell-2},\quad\ell\ge1
\]
and
\EQ{\nn
\rho(x;\alpha,\hbar)&\simeq x^{-1},\quad |\partial_x^k  \rho(x;\alpha,\hbar) |\lesssim x^{-k-1}, \\
|\partial_x^k \partial_\alpha^\ell \rho(x;\alpha,\hbar) | &\lesssim \hbar \, (x+\alpha)^{-\ell-k-1}
}
as claimed. 
\end{proof}

\begin{remark}\label{rem:largerJ}
The previous analysis covers a larger interval, not just $J_1$. For example, we can set $J_1:=[x_0,x_1]$ where $0\ll x_0\ll 1 \ll x_1$ are fixed. 
\end{remark}

Based on Lemma~\ref{lem: Lemma 3.4 CDST}, we can now describe the behavior of the potential $\tV$ in the same regime of parameters. 

\begin{proposition}\label{prop: tV x large}
For all $x\in J_2$, $\alpha>0$, and $\hbar\in(0,\frac{1}{3}]$ 
 \EQ{ \label{eq:tVJ2}
  |\partial_{\tau}^{k} \tV(\tau; \alpha,\hbar)| & \leq C_{k}\, \tau^{-2-k}\\ 
  |\partial_{\tau}^{k}\partial^{\ell}_{\alpha}\tV(\tau; \alpha,\hbar)| & \leq C_{k,\ell}\, \hbar\, \tau^{-2-k}(1+\alpha)^{-\ell-1}
}
 for all $k\ge0$, $\ell\ge1$. 
\end{proposition}
\begin{proof}
 By Lemma~\ref{lem: Lemma 3.4 CDST}, with some constant $c$, 
 \EQ{\nn
 \tau^3  &= c (x-y(\alpha;\hbar) + \rho(x;\alpha,\hbar))^2  = (\xi + \tilde{\rho}(\xi;\alpha,\hbar))^2
 } 
 where $\xi = c^{\frac12}(x-y(\alpha;\hbar))$. By the previous remark it suffices to consider the case $\xi\ge1$. Then $\tilde\rho$ satisfies the same estimates relative to~$\xi$ as $\rho$ does as a function of~$x$. 
 It is convenient to introduce the new variable 
 \[
 \eta=\xi+ \tilde\rho(\xi;\alpha,\hbar)
 \]
 Thus, $\tau^3=\eta^2$ and 
 \begin{align}\label{eq:qtau4}
 \begin{split}
3 \tau^2 \frac{\partial\tau}{\partial\eta} &=  2 \eta,\quad \Big(\frac{\partial\tau}{\partial\eta} \Big)^2 = \frac49 \tau^{-1},\quad \frac{\partial\eta}{\partial\tau}=\frac32 \tau^{\frac12}  \\
q(\tau;\alpha,\hbar) &= \Big(\frac{\partial\tau}{\partial x} \Big)^2 =   \tau^{-1}(1+\rho'(x;\alpha,\hbar))^2  \\
x^2 q(\tau;\alpha,\hbar) &=  (c^{-\frac12}\xi+y(\alpha;\hbar))^2 (1+\rho'(x;\alpha,\hbar))^2 \tau^{-1}\\
& = \tau^2 \, \big(c^{-\frac12} (1+\xi^{-1}\tilde\rho(\xi;\alpha,\hbar))^{-1}+y(\alpha;\hbar)/\eta\big)^2 (1+\rho'(x;\alpha,\hbar))^2
\end{split}
 \end{align}
 In view of \eqref{tV def}, and with $\sigma:=q^{-1}\dot q$, 
 \EQ{\label{eq:tV tau}
 \tV(\tau;\alpha,\hbar) &= -\frac{1}{4} \tau^{-2} \big(c^{-\frac12}(1+\xi^{-1}\tilde\rho(\xi;\alpha,\hbar))^{-1}+y(\alpha;\hbar)/\eta\big)^{-2} (1+\rho'(x;\alpha,\hbar))^{-2} +\frac14 \dot\sigma  + \frac{1}{16} \sigma^2 
 }
On the one hand, 
\EQ{\nn
\xi^{-1}\tilde\rho(\xi;\alpha,\hbar) &= O(\xi^{-2})=O(\tau^{-3}) \\
y(\alpha;\hbar)/\eta &= O(\tau^{-\frac32}) \\
\rho'(x;\alpha,\hbar) &= O(x^{-2}) = O(\tau^{-3})
}
On the other hand, by \eqref{eq:qdef},  $\tau'(x)=\sqrt{q}=O(\tau^{-\frac12})$ (recall ${\ }'$ refers to $x$ derivatives)
\EQ{\nn
\sigma &= -\tau^{-1}  +2(1+\rho'(x;\alpha,\hbar))^{-1} q^{-\frac12} \rho'' ,\qquad q^{-\frac12} \rho''= O(\tau^{-4})\\
\dot\sigma &= \tau^{-2}  +2q^{-\frac12}\big((1+\rho'(x;\alpha,\hbar))^{-1} q^{-\frac12} \rho'' )' = \tau^{-2} + O(\tau^{-5})
}
as $\tau\to\infty$.  Inserting these estimates into~\eqref{eq:tV tau} implies the first line of~\eqref{eq:tVJ2}, at least for $k=0$. As for the derivatives in $\tau$, we have
\[
\frac{\partial}{\partial\tau}=q^{-\frac12}\frac{\partial}{\partial x}= (1+\rho'(x;\alpha,\hbar))^{-1} \tau^{\frac12}  \frac{\partial}{\partial x}
\]
Each derivative in $x$ gains a power of $x\simeq \tau^{\frac32}$, but we then lose a $\tau^{\frac12}$ factor by the previous line, resulting in a total gain of a~$\tau^{-1}$ factor for each application of $\frac{\partial}{\partial\tau}$ to~$\tV$. 

For the $\alpha$ derivatives we compute 
\EQ{\nn 
|\partial_\alpha q(\tau;\alpha,\hbar)| &= 2  \, \tau^{-1}|(1+\rho'(x;\alpha,\hbar))  \partial_x\partial_\alpha\rho(x;\alpha,\hbar) | \lesssim \hbar \tau^{-1} (1+\alpha)^{-3} \\
q(\tau;\alpha,\hbar)^{-1}|\partial_\alpha q(\tau;\alpha,\hbar)| &\lesssim \hbar  (1+\alpha)^{-3}
}
whence
\EQ{\nn
|\partial_\alpha \sigma| &\lesssim |\partial_\alpha \rho'(x;\alpha,\hbar)| q^{-\frac12} |\rho''| + q^{-\frac12} |\rho''| q^{-1} |\partial_\alpha q| + q^{-\frac12} |\partial_x^2\partial_\alpha\rho(x;\alpha,\hbar) |  \\
&\lesssim \hbar (1+\alpha)^{-3} \tau^{-4} + \hbar \tau^{\frac12} (x+\alpha)^{-4} \lesssim \hbar \tau^{-1} (1+\alpha)^{-3} \\
|\partial_\alpha \, \sigma^2| & \lesssim \hbar \tau^{-2}(1+\alpha)^{-3}
}
Similarly, one checks that 
$
|\partial_\alpha \, \dot\sigma|  \lesssim \hbar \tau^{-2} (1+\alpha)^{-3}. 
$
However, if we take $\partial_\alpha$ of the first term in~\eqref{eq:tV tau}, then we obtain $\tau^{-2}$ multiplied by each of these three terms 
\[
 \xi^{-1}\partial_\alpha \tilde\rho(\xi;\alpha,\hbar),\quad \partial_\alpha (y(\alpha;\hbar)/\eta),\quad \partial_\alpha\rho'(x;\alpha,\hbar)
\]
which are on the order of, respectively, 
\[
\hbar \tau^{-\frac32} (1+\alpha)^{-2},\quad \hbar \tau^{-\frac32} (1+\alpha)^{-3},\quad \hbar (x+\alpha)^{-3}=O(\hbar \tau^{-\frac32} (1+\alpha)^{-2})
\]
uniformly in the range of parameters under consideration. In view of all these contributions, 
\[
|\partial_\alpha \tV(\tau;\alpha,\hbar)|\lesssim \hbar \tau^{-2}(1+\alpha)^{-2},
\]
as claimed. The higher derivatives are controlled similarly. 
 \end{proof}
 
 It remains to analyze the interval $J_0:=(0,\frac12x_t(\alpha;\hbar)]$.

\begin{lemma}
\label{lem:tilVxtau}
For $\tau<0$ one has the representation 
\EQ{
\label{eq:tilVxtau}
\tV(\tau;\alpha,\hbar) &= \frac{5}{16\tau^2} - \tau \varphi(x;\alpha,\hbar) \\
\varphi(x;\alpha,\hbar) &= \frac{1}{4Q_2(x;\alpha,\hbar)} \Big( x\mu'(x;\alpha,\hbar)-  \frac14\mu^2(x;\alpha,\hbar)\Big)\\
&= \frac{xQ'_{2}(x)}{4Q_{2}(x)^{2}}+\frac{x^{2}Q''_{2}(x)}{4Q_{2}(x)^{2}}-\frac{5x^{2}Q'_{2}(x)^{2}}{16Q_{2}(x)^{3}}.
}
where $Q_2(x;\alpha,\hbar):=x^2 Q_0(x;\alpha,\hbar)$ and $\mu(x;\alpha,\hbar)= xQ_2'(x;\alpha,\hbar)/Q_2(x;\alpha,\hbar)$. 
\end{lemma}
\begin{proof} 
By \eqref{eq:tV tau},  and suppressing the $\alpha,\hbar$ dependence from the notation, 
\[
\tV(x) = -\frac{1}{4x^2 q(x)} + \frac14\dot\sigma(x) + \frac{1}{16}\sigma^2(x),\qquad \sigma=\dot q/q
\]
Now, 
\EQ{\nn 
\sigma &=  \frac{-\tau}{Q_0(x)} \left( \frac{Q_0(x)}{\tau^2} - \frac{Q_0'(x)}{\tau}q^{-\frac12}\right) = -\frac{1}{\tau} + (-\tau)^{\frac12} \frac{Q_0'(x)}{Q_0^{\frac32}(x)} \\
\sigma^2 &=  \frac{1}{\tau^2} +2 (-\tau)^{-\frac12} \frac{Q_0'(x)}{Q_0^{\frac32}(x)} -\tau \frac{Q_0'(x)^2}{Q_0(x)^{3}} \\
\dot\sigma &=\frac{1}{\tau^2} -\frac12(-\tau)^{-\frac12}  \frac{Q_0'(x)}{Q_0^{\frac32}(x)} -\tau Q_0(x)^{-\frac12} \left( \frac{Q_0'(x)}{Q_0^{\frac32}(x)}\right)' \\
& = \frac{1}{\tau^2} -\frac12(-\tau)^{-\frac12}  \frac{Q_0'(x)}{Q_0^{\frac32}(x)}  - \tau\left(\frac{Q_0''(x)}{Q_0(x)^2} - \frac32 \frac{Q_0'(x)^2}{Q_0(x)^3}\right) 
}
whence
\[
\frac14\dot\sigma(x) + \frac{1}{16}\sigma^2(x) = \frac{5}{16 \tau^2} + \tau\left(-\frac14\frac{Q_0''(x)}{Q_0(x)^2} + \frac{5}{16}\frac{Q_0'(x)^2}{Q_0(x)^3}\right)
\]
and
\EQ{\label{eq:tVQ1Q0}
\tV(x) = \frac{5}{16 \tau^2} + \frac{ \tau}{ 4Q_2(x)}\left(1 -\frac{x^4Q_0''(x)}{Q_2(x)} + \frac{5}{4}\frac{(x^3Q_0'(x))^2}{Q_2(x)^2}\right) 
}
Inserting
\EQ{\nn
Q_0(x) & = x^{-2}Q_2(x),\quad x^3 Q_0'(x) = xQ_2'(x)-2Q_2(x) \\
x^4 Q_0''(x) &= 6Q_2(x) -4xQ_2'(x) + x^2 Q_2''(x)
}
into \eqref{eq:tVQ1Q0} yields 
\[
\tV(x) = \frac{5}{16 \tau^2} - \frac{ \tau}{ 4Q_2(x)} \frac{4xQ_2(x)Q_2'(x)-5x^2Q_2'(x)^2+ 4x^2Q_2(x) Q_2''(x)}{4Q_2(x)^2}
\]
Setting $\mu(x)= xQ_2'(x)/Q_2(x)$ we have 
\EQ{\nn 
x\mu' (x) =  x^2Q_2''(x)/Q_2(x) +\mu(x) -\mu(x)^2
}
and thus
\[
\tV(x) = \frac{5}{16 \tau^2} - \frac{ \tau}{ 4Q_2(x)}(  x\mu'(x)- \mu(x)^2/4 )
\]
as claimed. 
\end{proof}
 
 From \eqref{eq:Q0darst}, 
 \EQ{\label{eq:Q2 form}
Q_2(x;\alpha,\hbar) &=  (1-2\hbar)^2 + \frac{4\hbar x^2}{x^2+\alpha^2} - \frac{4\hbar^2 x^2(x^2+3\alpha^2)}{(x^2+\alpha^2)^2} -x^2  \\
xQ_2'(x;\alpha,\hbar) &= 8\hbar \alpha^2  \frac{x^2}{(x^2+\alpha^2)^2} + 8\hbar^2 \alpha^2 \frac{ x^2 (x^2-3\alpha^2)}{(x^2+\alpha^2)^3}-2x^2
}
 
We begin the analysis of $\tau$ and $\fy(x;\alpha,\hbar)$ in case  $\alpha\gtrsim1$. 

\begin{lemma}\label{lem: Lemma 3.3 CDST}
If $0<x\leq \frac{1}{2}x_t(\alpha,\hbar)$, then the function $\tau$ defined by \eqref{def tau} has the form
 \begin{align*}
  \frac{2}{3}\left(-\tau(x;\alpha, \hbar)\right)^{\frac{3}{2}}
  =&-(1-2\hbar)\log x+\veps_{1}(x; \alpha, \hbar)
 \end{align*}
with $\veps_{1}$ satisfying for all $k\ge0$, $\ell\ge1$,  uniformly in $x\in J_0(\alpha,\hbar)$
\begin{align*}
 |\partial^{k}_{x}\veps_{1}(x;\alpha,\hbar)|\leq C_{k,\ell},\qquad 
 |\partial^{\ell}_{\alpha}\partial^{k}_{x}\veps_{1}(x;\alpha,\hbar)|\leq C_{k,\ell}\, \hbar\, \langle\alpha\rangle^{-2-\ell}
\end{align*}
for all $\alpha\gtrsim1$ and $0<\hbar\leq\frac{1}{3}$. Moreover, $\veps_1$ is analytic as a function of $x$ in a neighborhood of~$0$ and $\partial_x \veps_1(0;\alpha,\hbar)=0$. 
\end{lemma}
\begin{proof}
 In this case we have
 \EQ{\nn 
\sqrt{Q_0(u;\alpha,\hbar)} &= u^{-1} \sqrt{Q_2(u;\alpha,\hbar)} \\
&= u^{-1}\sqrt{ (1-2\hbar)^2 + \frac{4\hbar u^2}{u^2+\alpha^2} - \frac{4\hbar^2 u^2(u^2+3\alpha^2)}{(u^2+\alpha^2)^2} -u^2 }\\
&   =\frac{1-2\hbar}{u}+f(u;\alpha,\hbar),
}
where $f$ is analytic in small $u$, and bounded uniformly (with all derivatives in~$u$) in $\alpha\gtrsim 1$. Moreover, $f(0;\alpha,\hbar)=0$ and  
\begin{align*}
 |\partial^{\ell}_{\alpha}\partial^{k}_{u}f(u,\alpha;\hbar)|\leq C_{k,\ell}\,\hbar \langle\alpha\rangle^{-2-\ell},\quad \ell\ge1,\; k\ge0.
\end{align*}
Note that $f$ does not decay in $\alpha$, since 
\[
\lim_{\alpha\to\infty} \sqrt{Q_0(u;\alpha,\hbar)} = u^{-1}\sqrt{ (1-2\hbar)^2 - u^2}
\]
but  derivatives of $f$ relative to $\alpha$ do decay in $\alpha$. 
We split the integral in the definition of $\tau$ in the form, with $0<x_0\ll 1$ fixed independently of $\alpha,\hbar$, 
\begin{align*}
 \int_{x}^{x_{t}(\alpha, \hbar)}=\int_{x}^{x_0}+\int_{x_0}^{x_{t}(\alpha, \hbar)}
\end{align*}
and obtain
\EQ{\label{eq:Q0 sum}
 \int_{x}^{x_t}\sqrt{Q_0(u;\alpha,\hbar)}\, du & = -(1-2\hbar)\log x+ (1-2\hbar)\log x_0  + \int_0^{x_0} f(u;\alpha,\hbar)  \, du  - \int_0^{x} f(u;\alpha,\hbar)  \, du \\
 &\quad\qquad\qquad + \int_{x_0}^{x_{t}(\alpha, \hbar)} \sqrt{Q_0(u;\alpha,\hbar)}\, du   \\
 & = -(1-2\hbar)\log x+ \veps_{1}(x; \alpha, \hbar)
}
The final integral here is governed by \eqref{eq:Q0 zero} and \eqref{eq:kell g}. In fact, from \eqref{eq:Q0 zero}
\EQ{\nn 
\int_{x_0}^{x_{t}(\alpha, \hbar)} \sqrt{Q_0(u;\alpha,\hbar)}\, du = \int_0^{x_t(\alpha;\hbar)-x_0} \sqrt{g(x_t(\alpha;\hbar)-v;\alpha,\hbar)}\, \sqrt{v}\, dv
}
which obeys the desired bounds as a function of $\alpha$ by Lemma~\ref{lem: monotonicity of root in alpha} and  \eqref{eq:kell g}. The term
\[
\int_0^{x_0} f(u;\alpha,\hbar)  \, du
\]
satisfies the same properties, as does 
\[
\int_0^{x} f(u;\alpha,\hbar)  \, du = O(x^2) \text{\ \ as\ \ } x\to0
\]
Adding up these different contributions in \eqref{eq:Q0 sum} concludes the proof. 
\end{proof}

To express $x$ as a function of $\tau$ we invert the function in the previous proposition. 

\begin{corollary}
\label{cor:x tau} 
There exists a constant $x_0\in J_0(\alpha;\hbar)$ so that 
for all $x\in (0,x_0]$, with the parameters as in the previous proposition,  there is a representation
\[
x=X\left(\exp\left( - \frac{2}{3(1-2\hbar)}\left(-\tau(x;\alpha, \hbar)\right)^{\frac{3}{2}}\right);\alpha, \hbar\right)
\]
where $\tau\le \tau_0(\alpha,\hbar)$, and $\tau_{0}(\alpha,\hbar)$ is the value of the map $\tau=\tau(x;\alpha,\hbar)$ at $x=x_{0}$. Here $X$ is a diffeomorphism  $[0,y_0(\alpha;\hbar)]\to [0,x_{0}]$, where $y_{0}(\alpha,\hbar)$ is defined such that $X(y_{0}(\alpha,\hbar);\alpha,\hbar)=x_{0}$,  with $$X'(0;\alpha, \hbar)=\exp\left(  \frac{1}{1-2\hbar} \veps_1(0;\alpha,\hbar) \right).$$
We have $X'(y;\alpha, \hbar)\gtrsim 1$ for all $y\in[0,y_0]$ and $\alpha\gtrsim 1$ and $\hbar\in(0,\frac13]$, as well as 
\EQ{\label{eq:diff bds}
|\partial_y^k X(y;\alpha, \hbar)| \le C_k,\quad |\partial_y^k \partial_\alpha^\ell X(y;\alpha, \hbar)| \le C_{k,\ell}\, \hbar\, \alpha^{-2-\ell}
}
for all $k\ge0$, $\ell\ge1$. 
\end{corollary}
\begin{proof}
We set 
\EQ{\nn 
\log Y(x;\alpha,\hbar) &= \log x  - \frac{1}{1-2\hbar} \veps_1(x;\alpha,\hbar) \\
Y(x;\alpha,\hbar) &= x\exp\left(  -\frac{1}{1-2\hbar} \veps_1(x;\alpha,\hbar) \right). 
}
By Lemma~\ref{lem: Lemma 3.3 CDST}, $Y(x;\alpha,\hbar)$ is smooth as a function of~$x$ on a neighborhood of $0$ containing $J_0(\alpha,\hbar)$. Moreover, 
$$Y'(0;\alpha,\hbar)= \exp\left(  -\frac{1}{1-2\hbar} \veps_1(0;\alpha,\hbar) \right) $$
is bounded below uniformly in $\alpha\gtrsim 1$ and $\hbar\in(0,\frac13]$.  It follows that we can smoothly invert $y=Y(x;\alpha,\hbar)$  so that $x=X(y;\alpha,\hbar)$, uniformly in the parameters. Finally, we have the bounds
\begin{align}\label{higher deri Y}
|\partial_x^k Y(x;\alpha, \hbar)| \le C_k,\quad |\partial_x^k \partial_\alpha^\ell Y(x;\alpha, \hbar)| \le C_{k,\ell}\, \hbar\, \alpha^{-2-\ell}
\end{align}
for all $k\ge 0$ and $\ell\ge1$.  These bounds imply~\eqref{eq:diff bds}. Indeed, from $X(Y(x;\alpha,\hbar);\alpha,\hbar)=x$ we deduce
\begin{align}\label{higher deri X}
\begin{split}
\partial_y X \partial_x Y &=1 \\
- (\partial_x Y)^{-1} (\partial_y^2 X\partial_\alpha Y \partial_x Y  +  \partial_y X \partial_x \partial_\alpha Y )&=  \partial_y\partial_\alpha X  \\
-(\partial_x  Y)^{-1}\big(\partial_y^3 X (\partial_\alpha Y)^2 \partial_x Y + \partial_y^2\partial_\alpha X \partial_\alpha Y\partial_x Y + \partial_y^2 X \partial_\alpha^2 Y \partial_x Y + 
\partial_y^2 X \partial_\alpha Y \partial_x\partial_\alpha Y + & \\ 
\partial_y^2 \partial_\alpha X  \partial_\alpha  Y  \partial_x Y +  2\partial_y \partial_\alpha X \partial_x \partial_\alpha Y + \partial_y^2 X \partial_\alpha Y \partial_x\partial_\alpha Y+ \partial_y X \partial_x \partial_\alpha^2 Y \big)&=   \partial_y\partial_\alpha^2 X  
\end{split}
\end{align}
which ~\eqref{eq:diff bds}. The higher derivatives follow inductively. In fact, if we differentiate both sides of the first equation in \eqref{higher deri X} with respect to $x$, we obtain
	\begin{align*}
&		\partial^{2}_{y}X(\partial_{x}Y)^2+\partial_{y}X\partial_{x}^{2}Y=0,\\
&\partial^{3}_{y}X(\partial_{x}Y)^3+3\partial_{y}^{2}X\partial_{x}^{2}Y\partial_{x}Y+\partial_{y}X\partial_{x}^{3}Y=0,\\
&...
	\end{align*}
	Inductively, we have
	\begin{align}\label{higher deri X wrt y}
		\partial^{k}_{y}X(\partial_{x}Y)^k+P_{k}(\partial_{y}^{k-1}X,...,\partial_{y}X, \partial_xY, \partial^{2}_{x}Y,...,\partial_{x}^{k}Y)=0, \quad \textrm{for}\quad k\geq 2.
	\end{align}
	Here $P_{k}(\partial_{y}^{k-1}X,...,\partial_{y}X, \partial_xY, \partial^2_{x}Y,...,\partial_{x}^{k}Y)$ is a polynomial in the variables $\partial_{y}^{k-1}X,...,\partial_{y}X$ and $\partial_{x}Y,...,\partial_{x}^{k}Y$. 
	Then we use an induction argument and the lower bound on $\partial_{x}Y$ for small $x$ as well as the bounds \eqref{higher deri Y} to conclude 
	\begin{align}\label{y deri of X}
	|\partial_{y}^{k}X(y,\alpha;\hbar)|\leq C_{k}, \quad \textrm{for all}\quad k\geq 0.
	\end{align}
	For the $\alpha$-derivatives of $X$, we differentiate \eqref{higher deri X wrt y} in $\alpha$ on both sides to obtain
	\begin{align}\label{dy k dalpha X}
		\partial_{y}^{k}\partial_{\alpha}X(\partial_{x}Y)^k+k\partial_{y}^{k}X\partial_{x}\partial_{\alpha}Y(\partial_{x}Y)^{k-1}+\partial_{\alpha}\left(P_{k}(\partial_{y}^{k-1}X,...,\partial_{y}X, \partial_{x}Y,...,\partial_{x}^{k}Y)\right)=0,\quad \textrm{for}\quad k\geq 2.
	\end{align}
	Here $\partial_{\alpha}\left(P_{k}(\partial_{y}^{k-1}X,...,\partial_{y}X, \partial^{2}_{x}Y,...,\partial_{x}^{k}Y)\right)$ is a polynomial containing factors $\partial^{k-1}_{y}\partial_{\alpha}X,...,$ and $\partial_{y}\partial_{\alpha}X$, $\partial_{x}\partial_{\alpha}Y,...,\partial_{x}^{k}\partial_{\alpha}Y$, as well as $\partial_{y}^{k-1}X,...,\partial_{y}X$, $\partial_{x}Y,...,\partial_{x}^{k}Y$. Therefore the estimate on $\partial_{y}^{k}\partial_{\alpha}X$ follows from the estimates on the higher order $y$-derivatives of $X$, and the higher order derivatives of $Y$. Again, inductively we have
	\begin{align}\label{dy k dalpha ell X}
		\begin{split}
		\partial^{k}_{y}\partial^{\ell}_{\alpha}X\partial_{x}Y+Q_{k,\ell}=0.
		\end{split}
	\end{align}
	Here $Q_{k,\ell}$ is a polynomial with factors $\partial^{k'}_{y}\partial^{\ell'}_{\alpha}X$, as well as higher order derivatives of $Y$, where $k'\leq k, \ell'\leq \ell$ and at least one of $k'<k, \ell'<\ell$ holds.
Therefore the estimates on higher order derivatives $\partial_{y}^{k}\partial_{\alpha}^{\ell}X$ follow from an induction argument.
\end{proof}

Based on Lemma~\ref{lem: Lemma 3.3 CDST} and the previous corollary, we now analyze the behavior of $\tV$ for $0<x\leq\frac{1}{2}x_t(\alpha;\hbar)$ and $\alpha\gtrsim1$. 

\begin{proposition}
\label{prop: tV close to 0}
 For $0<x\leq\frac{1}{2}x_t(\alpha;\hbar)$ and $\alpha\gtrsim 1$, we have
 \begin{align*}
 |\partial_{\tau}^{k} \tV(\tau;\alpha, \hbar)|\leq C_{k}\, \tau^{-2-k},\qquad   |\partial_{\tau}^{k}\partial^{\ell}_{\alpha}\tV(\tau;\alpha,\hbar)|\leq C_{k,\ell}\hbar \, e^{\tau} \alpha^{-2-\ell}
 \end{align*}
 for all $k\ge0$ and $\ell\ge 1$. 
\end{proposition}
\begin{proof}
From Lemma~\ref{lem:tilVxtau} and \eqref{eq:Q2 form}
\EQ{\nn 
\tV(\tau;\alpha,\hbar) & = \frac{5}{16\tau^2} + \tau\, O_{\alpha,\hbar} \left(  X\left(\exp\left( - \frac{2}{3(1-2\hbar)}\left(-\tau\right)^{\frac{3}{2}}\right);\alpha, \hbar\right)^2\right) \\
&= \frac{5}{16\tau^2} + \tau\, O_{\alpha,\hbar} \left( \exp\left( - \frac{4}{3(1-2\hbar)}\left(-\tau\right)^{\frac{3}{2}}\right)  \right)
}
By inspection of the function $\mu$ in Lemma~\ref{lem:tilVxtau} and \eqref{eq:Q2 form} one sees that the  $O_{\alpha,\hbar}$-term here depends on $\alpha$ only through 
terms which are on the order of $\hbar$. Any derivative in $\alpha$ of the order $\ell$ will gain $\hbar \alpha^{-2-\ell}$ by the structure of $\mu$ and Corollary~\ref{cor:x tau}. 
\end{proof}

The exponential decay stated in the proposition is not optimal, it was chosen for convenience.

Next we turn to the case when $0<\alpha\ll1$.  An important difference arises here with respect to the potential $\tV(x;\alpha,\hbar)$. Indeed, one checks that
\[
x\mu'(x;\alpha,\hbar) - \mu(x;\alpha,\hbar)^2/4\Big|_{x=\alpha} = -\frac{\alpha^4 + \hbar^2 (1 - 2 \hbar + 5 \hbar^2) +
 2 \alpha^2 (2 - 9 \hbar + 7 \hbar^2)}{(1 - \alpha^2 - 2 \hbar)^2} <0
\]
for small $\hbar$ and $\alpha$.  This means that the term $\tau\varphi(x;\alpha,\hbar)$ in~\eqref{eq:tilVxtau} is large for small $\alpha$ when $x\simeq\alpha$, and dominates $5/(16\tau^2)$. 
We distinguish three parameter regimes: $0<x\ll\alpha\ll 1$, $0<\alpha\ll x\ll 1$, and $\alpha\simeq x\ll 1$.  We begin with the latter case. 

\begin{lemma}\label{lem:xsimal}
Let $0<k_0\ll 1\ll K_0$ be fixed constants. Then for all $x\in [k_0\alpha,K_0\alpha]$, 
\EQ{\label{eq:tau comp}
 \frac{2}{3}\left(-\tau(x;\alpha,\hbar)\right)^{\frac{3}{2}} = -\log x + \sigma(x;\alpha,\hbar) 
}
where 
\[
|\partial_x^k\partial_\alpha^\ell \sigma(x;\alpha,\hbar) |\le C_{k,\ell} \, \alpha^{-k-\ell}
\]
for all $k,\ell\ge0$. If $\ell\ge1$, then a factor of $\hbar$ is gained on the right-hand side. The constants depend on $k_0,K_0$, and $K_0\alpha\ll1$. 
\end{lemma}
\begin{proof}
We have
\EQ{\nn
 \frac{2}{3}\left(-\tau(x;\alpha,\hbar)\right)^{\frac{3}{2}}  & =  \int_x^{x_t} \left((1-2\hbar)^2 + \frac{4\hbar u^2}{u^2+\alpha^2} - \frac{4\hbar^2 u^2(u^2+3\alpha^2)}{(u^2+\alpha^2)^2}-u^2\right)^{\frac12}\, \frac{du}{u} \\
&=  \int_{x/\alpha}^{x_t/\alpha} \left((1-2\hbar)^2 + \frac{4\hbar v^2}{1+v^2} - \frac{4\hbar^2 v^2(v^2+3 )}{(1+v^2 )^2}-\alpha^2 v^2\right)^{\frac12}\, \frac{dv}{v} \\
&=  \int_{x/\alpha}^{x_t/\alpha} \left(1 - 4\hbar \frac{1-\hbar + (1+\hbar)v^2 }{(1+v^2 )^2}-\alpha^2 v^2\right)^{\frac12}\, \frac{dv}{v}
}
Since for $\hbar\in (0,\frac13]$ and $v>0$ 
\[
\frac{d}{dv} \frac{1-\hbar + (1+\hbar)v^2 }{(1+v^2 )^2} = - 2v(1+v^2)^{-3}(1-3\hbar+(1+\hbar)v^2)<0
\]
it follows that 
\[
4\hbar \frac{1-\hbar + (1+\hbar)v^2 }{(1+v^2 )^2} \le 4\hbar(1-\hbar)\le \frac89
\]
Thus, if $v\le \frac{1}{4\alpha}$, then 
\[
4\hbar \frac{1-\hbar + (1+\hbar)v^2 }{(1+v^2 )^2}+\alpha^2 v^2 \le \frac89+\frac{1}{16}<1
\]
and 
\EQ{\nn 
& \left(1 - 4\hbar \frac{1-\hbar + (1+\hbar)v^2 }{(1+v^2 )^2}-\alpha^2 v^2\right)^{\frac12}  = 1 - \sum_{k=1}^\infty (-1)^{k-1}\binom{\frac12}{k}\left( 4\hbar \frac{1-\hbar + (1+\hbar)v^2 }{(1+v^2 )^2}+\alpha^2 v^2\right)^k
}
converges uniformly on $0\le v\le \frac{1}{4\alpha}$, and $\hbar\in (0,\frac13]$. It is, however, not easy to integrate term-wise and estimate each term separately in the whole range of~$\hbar$. 
Instead, we proceed as follows:  if $4K_0 \alpha<1$, then
\EQ{\nn 
\frac{2}{3}\left(-\tau(x;\alpha,\hbar)\right)^{\frac{3}{2}} & = -\log(4x) + \omega_1(\alpha,\hbar) + \omega_2(x;\alpha,\hbar)\\ 
\omega_1(\alpha,\hbar) &:=  \int_{1/(4\alpha)}^{x_t/\alpha} \left(1 - 4\hbar \frac{1-\hbar + (1+\hbar)v^2 }{(1+v^2 )^2}-\alpha^2 v^2\right)^{\frac12}\, \frac{dv}{v}  \\
\omega_2(x;\alpha,\hbar)  &:= - \int_{x/\alpha}^{1/(4\alpha)} \frac{4\hbar \frac{1-\hbar + (1+\hbar)v^2 }{(1+v^2 )^2}+\alpha^2 v^2 }{1+\left(1 - 4\hbar \frac{1-\hbar + (1+\hbar)v^2 }{(1+v^2 )^2}-\alpha^2 v^2\right)^{\frac12}}\, \frac{dv}{v} 
}
Note that by Lemma~\ref{lem: monotonicity of root in alpha}, $x_t(\alpha;\hbar)>\frac13$. Lemma~\ref{lem: Lemma 3.2 CDST} applies to $\omega_1(\alpha,\hbar)$, i.e., 
\[
0\le \omega_1(\alpha,\hbar)\lesssim  1,  \qquad |\partial_\alpha^\ell \omega_1(\alpha,\hbar)|\le C_\ell\, \hbar \text{\ \ for all\ \ }\ell\ge1
\]
uniformly in the regime of parameters under consideration. Furthermore, 
\EQ{\nn
|\omega_2(x;\alpha,\hbar)|  &\le   \int_{x/\alpha}^{1/(4\alpha)} \left(4\hbar (1+\hbar)(1+v^2 )^{-1}+\alpha^2 v^2 \right)\, \frac{dv}{v}  \\
&= 2\hbar (1+\hbar) \log\frac{v^2}{1+v^2} + \frac12\alpha^2 v^2 \Big|^{1/(4\alpha)}_{x/\alpha} \\
&= 2\hbar (1+\hbar) \left( -\log(1+16\alpha^2) + \log(1+\alpha^2/x^2)\right) +\frac12\left(\frac{1}{16}-x^2\right) =O(1)
}
uniformly in $x\simeq \alpha\ll1$. Next, by inspection, for all $k\ge0$, 
\[
|\partial_x^k \omega_2(x;\alpha,\hbar)|\le C_k \alpha^{-k}
\]
Since $x\simeq \alpha$, we can replace $\alpha^{-k}$ with $x^{-k}$. For the derivatives relative to~$\alpha$ it is convenient to undo the scaling of $\omega_2$, to wit
\EQ{\nn 
\omega_2(x;\alpha,\hbar)  &:=  \int_{x}^{1/4} \frac{4\hbar \alpha^2 \frac{(1-\hbar)\alpha^2 + (1+\hbar)u^2 }{(\alpha^2+u^2 )^2}+u^2 }{1+\left(1 - 4\hbar\alpha^2 \frac{(1-\hbar)\alpha^2 + (1+\hbar)u^2 }{(\alpha^2+u^2 )^2}-  u^2\right)^{\frac12}}\, \frac{du}{u} 
}
From this expression one can derive that 
\[
|\partial_x^k \partial_\alpha^\ell \omega_2(x;\alpha,\hbar)|\le C_{k,\ell} \, \hbar \alpha^{-k-\ell}
\]
if $k\ge0$ and $\ell\ge1$. 
\end{proof}

Next, we describe $\tV$ in the regime of Lemma~\ref{lem:xsimal}. 

\begin{proposition}
\label{prop:tV xsimal}
Using the representation of Lemma~\ref{lem:tilVxtau}, we have 
\EQ{\label{eq:tVvarphi}
\tV(\tau;\alpha,\hbar) = \frac{5}{16\tau^2} - \tau \varphi(x;\alpha,\hbar)
}
where $|\varphi(x;\alpha,\hbar)|\le C(\hbar+\alpha^2)$ uniformly in $x\in [k_0\alpha,K_0\alpha]$, $0<\alpha\ll1$. Moreover,
\EQ{\label{eq:varphi der}
|\partial_\tau^k \partial_\alpha^\ell  \varphi(x;\alpha,\hbar)|\le C_{k,\ell}\, (\hbar+\alpha^2) (-\tau)^{\frac{k}{2}}\alpha^{-\ell}
}
for all $k,\ell\ge0$. If $\ell\ge1$, then a factor of $\hbar$ is gained in~\eqref{eq:varphi der}. The constants depend on $k_0,K_0$, and $\alpha_0>0$ where $K_0\alpha\le\alpha_0 \ll1$. 
\end{proposition}
\begin{proof}
By Lemma~\ref{lem:tilVxtau}, 
\[
\varphi(x;\alpha,\hbar) =   \frac{1}{4Q_2(x;\alpha,\hbar)} \Big( x\mu'(x;\alpha,\hbar)-  \frac14\mu^2(x;\alpha,\hbar)\Big)
\]
where $Q_2(x;\alpha,\hbar):=x^2 Q_0(x;\alpha,\hbar)$ and $\mu(x;\alpha,\hbar)= xQ_2'(x;\alpha,\hbar)/Q_2(x;\alpha,\hbar)$. 
Formulas~\eqref{eq:Q2 form} determine $\mu$. By inspection, $|\mu|\lesssim \hbar+\alpha^2$ and $|x\mu'|\lesssim \hbar+\alpha^2$ for all $x\simeq\alpha$. 
For the derivatives, recall that
\EQ{\label{eq:taudx1}
\frac{\partial}{\partial\tau}= q^{-\frac12} \frac{\partial}{\partial x}= (-\tau)^{\frac12} \frac{1}{\sqrt{Q_2(x;\alpha,\hbar)}} \, x\frac{\partial}{\partial x}
}
The $x\frac{\partial}{\partial x}$ operator does not change the bound on $\mu$ or $\varphi(x;\alpha,\hbar)$. The largest contribution to the higher derivatives comes from 
\EQ{\label{eq:taudx2}
\frac{\partial^k}{\partial\tau^k}= \left(q^{-\frac12} \frac{\partial}{\partial x}\right)^k= (-\tau)^{\frac{k}{2}} \left(\sqrt{Q_2(x;\alpha,\hbar)}\right)^{-k} \, \left(x\frac{\partial}{\partial x}\right)^k+ \text{\ lower order}
}
where ``lower order'' refers to terms involving fewer $x$-derivatives. 
The derivatives with respect to $\alpha$ bring out a (single) factor of $\hbar$, and lose a factor of $\alpha^{-1}$ each.  This yields~\eqref{eq:varphi der}. 
\end{proof}

Next, we analyze the case $0<x\ll\alpha\ll 1$. 

\begin{lemma}\label{lem: Lemma 3.3 harder}
 If $0<x\ll\alpha\ll 1$, then the function $\tau$ defined by \eqref{def tau} has the form
 \begin{align*}
  \frac{2}{3}\left(-\tau(x;\alpha,\hbar)\right)^{\frac{3}{2}}=  -\log x + 2\hbar\log(x/\alpha) + \psi(x/\alpha;\alpha,\hbar)  + \rho(\alpha;\hbar)
 \end{align*}
where $\psi(z;\alpha,\hbar)$ is analytic near $0$, and uniformly in $0<\alpha\ll1$, and $|z|\le r_0$ for some absolute constant $0<r_0\ll 1$, 
\EQ{\nn
|\partial_z^k \partial_\alpha^\ell \psi(z;\alpha,\hbar)| &\le C_{k,\ell}\,(\hbar+\alpha^2), \quad |\partial_\alpha^\ell\rho(\alpha;\hbar)|\le C_\ell\,\alpha^{-\ell}
}
for all $k,\ell\ge0$. 
\end{lemma}
\begin{proof}
 We have
\EQ{\label{eq:gc} 
 \frac{2}{3}\left(-\tau(x;\alpha,\hbar)\right)^{\frac{3}{2}} & =  \frac{2}{3}\left(-\tau(k_0\alpha;\alpha,\hbar)\right)^{\frac{3}{2}}  + \int_{x} ^{k_0\alpha} \left((1-2\hbar)^2 + \frac{4\hbar u^2}{u^2+\alpha^2} - \frac{4\hbar^2 u^2(u^2+3\alpha^2)}{(u^2+\alpha^2)^2}-u^2\right)^{\frac12}\, \frac{du}{u} \\
&=\frac{2}{3}\left(-\tau(k_0\alpha;\alpha,\hbar)\right)^{\frac{3}{2}}  +  \int_{x/\alpha}^{k_0} \left((1-2\hbar)^2 + \frac{4\hbar v^2}{1+v^2} - \frac{4\hbar^2 v^2(v^2+3 )}{(1+v^2 )^2}-\alpha^2 v^2\right)^{\frac12}\, \frac{dv}{v}  
}
We write, with analytic $f_1$ and $f_2$ on the unit disk, 
\EQ{\nn 
\frac{ v^2}{1+v^2}  &= v^2 f_1(v^2), \quad \frac{ v^2(v^2+3 )}{(1+v^2 )^2} = 3v^2 f_2(v^2),\quad f_1(0)=f_2(0)=1\\ 
}
whence 
\EQ{\nn
& \left((1-2\hbar)^2 + \frac{4\hbar v^2}{1+v^2} - \frac{4\hbar^2 v^2(v^2+3 )}{(1+v^2 )^2}-\alpha^2 v^2\right)^{\frac12} \\
&= (1-2\hbar) \left( 1+ 4\hbar (1-2\hbar)^{-2} v^2 f_1(v^2)  - 12\hbar^2(1-2\hbar)^{-2} v^2 f_2(v^2) -\alpha^2(1-2\hbar)^{-2}  v^2\right)^{\frac12} \\
&= (1-2\hbar) \left( 1 + v^2 f(v^2;\alpha,\hbar)\right)^{\frac12}   = (1-2\hbar)\left( 1 +  \sum_{n=1}^\infty  \binom{\frac12}{n} v^{2n} f(v^2;\alpha,\hbar)^n\right)
}
with 
\EQ{\nn
f(v^2;\alpha,\hbar) &= 4\hbar (1-2\hbar)^{-2} f_1(v^2)  - 12\hbar^2(1-2\hbar)^{-2}  f_2(v^2) -\alpha^2(1-2\hbar)^{-2}   \\
f(0;\alpha,\hbar) &= 4\hbar (1-2\hbar)^{-2} -12 \hbar^2(1-2\hbar)^{-2} -\alpha^2(1-2\hbar)^{-2}
}
as well as in the complex plane 
\EQ{\label{eq:fmax}
\max_{|z|\le \frac12} |f(z;\alpha,\hbar)| \le C(\hbar+ \alpha^2)
}
where $C$ is an absolute constant. With $n\ge1$, 
\EQ{\nn 
F_n(y;\alpha,\hbar) &:=  \int_{0}^{y} v^{2n-1} f(v^2;\alpha,\hbar)^n\, dv  \\
 \int_{x/\alpha}^{k_0} v^{2n-1} f(v^2;\alpha,\hbar)^n\, dv &= F_n(k_0;\alpha,\hbar) - F_n(x/\alpha;\alpha,\hbar) 
}
By \eqref{eq:fmax}, if $0\le y\le \frac12$, then 
\[
|F_n(y;\alpha,\hbar)|\le C^n(\hbar+ \alpha^2)^n y^{2n}
\]
Thus, by \eqref{eq:gc}, and Lemma~\ref{lem:xsimal}, and with $k_0$ some fixed small constant, 
\EQ{\nn
 \frac{2}{3}\left(-\tau(x;\alpha,\hbar)\right)^{\frac{3}{2}} & =  \frac{2}{3}\left(-\tau(k_0\alpha;\alpha,\hbar)\right)^{\frac{3}{2}}  + \int_{x/\alpha}^{k_0}  (1-2\hbar)\Big( 1 +  \sum_{n=1}^\infty  \binom{\frac12}{n} v^{2n} f(v^2;\alpha,\hbar)^n\Big) \, \frac{dv}{v}   \\
& = \frac{2}{3}\left(-\tau(k_{0}\alpha;\alpha,\hbar)\right)^{\frac{3}{2}}+(1-2\hbar)\log k_{0}-(1-2\hbar)\log\frac{x}{\alpha} \\
&\quad + (1-2\hbar) \sum_{n=1}^\infty  \binom{\frac12}{n} \big( F_n(k_0;\alpha,\hbar) - F_n(x/\alpha;\alpha,\hbar) \big) \\
& =  -\log x +  2\hbar\log(x/\alpha) + \psi(x/\alpha;\alpha,\hbar)  + \rho(\alpha;\hbar)
}
where $\psi$ and $\rho$ are given by

\begin{align*}
	\psi\left(z;\alpha,\hbar\right):=&-(1-2\hbar)\sum_{n=1}^{\infty} \binom{\frac12}{n}F_{n}\left(z;\alpha,\hbar\right),\\
	\rho(\alpha;\hbar):=&\frac{2}{3}\left(-\tau(k_{0}\alpha;\alpha,\hbar)\right)^{\frac{3}{2}}+(1-2\hbar)\log k_{0}+\log\alpha+(1-2\hbar)\sum_{n=1}^{\infty} \binom{\frac12}{n}F_{n}(k_{0};\alpha,\hbar),
\end{align*}
and have the properties stated in the lemma. In particular, the estimate on $\tau(k_{0}\alpha;\alpha,\hbar)$ follows from the result in Lemma \ref{lem:xsimal}.
\end{proof}

Based on Lemma~\ref{lem: Lemma 3.3 harder}, we now describe the potential $\tV(\tau;\alpha,\hbar)$ in the regime $0<x\ll\alpha\ll 1$. 

\begin{proposition}\label{prop: tV close to 0 harder}
In the parameter regime of the previous lemma,  and using the representation~\eqref{eq:tVvarphi}, 
\EQ{\nn 
|\partial_\tau^k \partial_\alpha^\ell\varphi(x;\alpha,\hbar)| 
&\le C_{k,\ell}\, (\hbar+\alpha^2)\,  \alpha^{-\ell} (  -\log x + 2\hbar\log(x/\alpha) )^{\frac{k}{3}}  (x/\alpha)^2 \\
&\le C_{k,\ell}\, (\hbar+\alpha^2)\,  \alpha^{-\ell} (  -\tau )^{\frac{k}{2}}  (x/\alpha)^2
}
for all $k,\ell\ge0$. 
\end{proposition}
\begin{proof}
By \eqref{eq:Q2 form},  with $\xi:=x/\alpha$, 
 \EQ{\nn 
 Q_2(x;\alpha,\hbar) &=  (1-2\hbar)^2 + \frac{4\hbar x^2}{x^2+\alpha^2} - \frac{4\hbar^2 x^2(x^2+3\alpha^2)}{(x^2+\alpha^2)^2} -x^2  \\
 & =  (1-2\hbar)^2 + \frac{4\hbar \xi^2}{1+\xi^2} - \frac{4\hbar^2 \xi^2(\xi^2+3)}{(1+\xi^2)^2} -\alpha^2 \xi^2 \\ 
xQ_2'(x;\alpha,\hbar) &= 8\hbar \alpha^2  \frac{x^2}{(x^2+\alpha^2)^2} + 8\hbar^2 \alpha^2 \frac{ x^2 (x^2-3\alpha^2)}{(x^2+\alpha^2)^3}-2x^2 \\
&= 8\hbar   \frac{\xi^2}{(1+\xi^2)^2} + 8\hbar^2  \frac{ \xi^2 (\xi^2-3)}{(1+\xi^2)^3}-2\alpha^2\xi^2 
}
Thus,  we have $|\mu|\le C(\hbar+\alpha^2)\xi^{2}$ and therefore 
\[
|\varphi(x;\alpha,\hbar)|\le C(\hbar+\alpha^2)\xi^2 ,\qquad |\partial_\alpha^\ell \varphi(x;\alpha,\hbar)|\le C(\hbar+\alpha^2)\alpha^{-\ell}\,\xi^2
\]
uniformly in the parameter regime of the proposition, for all $\ell\ge0$.  In view of \eqref{eq:taudx1}, \eqref{eq:taudx2}, 
$x\partial_x=\xi\partial_\xi$,  and Lemma~\ref{lem: Lemma 3.3 harder}, 
\EQ{\nn 
|\partial_\tau^k \varphi(x;\alpha,\hbar)| &\le C_k\, (\hbar+\alpha^2)\,  (\sqrt{-\tau})^{k} \xi^2\\
&\le C_k\, (\hbar+\alpha^2)\,  (  -\log x + 2\hbar\log(x/\alpha) )^{\frac{k}{3}}  \xi^2 
}
for all $k\ge0$, as well as 
\EQ{\nn 
|\partial_\tau^k \partial_\alpha^\ell\varphi(x;\alpha,\hbar)| 
&\le C_{k,\ell}\, (\hbar+\alpha^2)\,  \alpha^{-\ell} (  -\log x + 2\hbar\log(x/\alpha) )^{\frac{k}{3}}  \xi^2 
}
as claimed. 
\end{proof}

Finally we discuss the case when $0<\alpha\ll x \ll1$. 

\begin{lemma}\label{lem: Lemma 3.3 hardest}
The function $\tau$ defined by \eqref{def tau} has the form, with $\eta=\alpha/x$, 
\EQ{ \label{eq:taux final case}
\frac{2}{3}\left(-\tau(x;\alpha,\hbar)\right)^{\frac{3}{2}}=-(1-\alpha^2 f_3(\alpha;\hbar))\log x + C(\alpha;\hbar) + x^2 f_1(x,\eta,\alpha;\hbar) + \hbar\eta^2 f_2(x,\eta,\alpha;\hbar) 
}
for all $0<\alpha\ll x\ll 1$ and $0<\hbar\leq\frac{1}{3}$. Here $f_1$ and $f_2$ are analytic in small (complex) $x,\eta,\alpha$ and bounded on a small polydisc in $\C^3$ centered at $(0,0,0)$, uniformly in $0<\hbar\leq\frac{1}{3}$. The function $f_3$ is analytic and bounded on a small disk of small complex $\alpha$, uniformly in the same range of~$\hbar$. Finally, for all $\ell\ge0$ we have $|\partial_\alpha^\ell C(\alpha;\hbar)|\le C_\ell\, \alpha^{-\ell}$ uniformly in small $\alpha$ and $\hbar$ in the same range as before. 
\end{lemma}
\begin{proof}
With $\eta=\frac{\alpha}{u}$,  we have
\EQ{\label{Q0 small alpha}
 Q_{2}(u;\alpha, \hbar)&= 1 -\frac{4\hbar(1-\hbar)\alpha^{2}}{u^{2}+\alpha^{2}}-\frac{8\hbar^{2}\alpha^{2} u^2}{(u^{2}+\alpha^{2})^{2}}-u^2\\
 &= 1 -\frac{4\hbar(1-\hbar)\eta^{2}}{1+\eta^{2}}-\frac{8\hbar^{2}\eta^{2}}{(1+\eta^{2})^{2}}-u^2 \\
 &= 1- \hbar \eta^2 g(\eta^2;\hbar) - u^2
}
with 
\[
g(z;\hbar)= \frac{4(1-\hbar)}{1+z}+\frac{8\hbar}{(1+z)^{2}}, \quad g(0;\hbar) = 4(1+\hbar)
\]
Thus,  fixing some $0<x_0\ll 1$ so that $0\ll \alpha \ll x\le x_0$, we have 
\EQ{\label{eq:gc2} 
 \frac{2}{3}\left(-\tau(x;\alpha,\hbar)\right)^{\frac{3}{2}} & =  \frac{2}{3}\left(-\tau(x_0;\alpha,\hbar)\right)^{\frac{3}{2}} +\log x_0 - \log x \\
 &\qquad - \int_{x} ^{x_0}   \Big[\Big(1- \hbar \eta^2 g(\eta^2;\hbar) - u^2\Big)^{\frac12}-1\Big] \, \frac{du}{u}  \\
 &= \frac{2}{3}\left(-\tau(x_0;\alpha,\hbar)\right)^{\frac{3}{2}} +\log x_0 - \log x \\
 &\qquad + \sum_{n=1}^\infty (-1)^{n-1} \binom{\frac12}{n} \sum_{\ell=0}^n \binom{n}{\ell} \hbar^\ell  \alpha^{2\ell} \int_{x} ^{x_0}    g(\alpha^2 u^{-2};\hbar)^\ell u^{2(n-2\ell)-1} \,  du  
 }
 Plugging the Taylor expansion of $g$ into the integral leads to expressions of the form (because of smallness, uniform convergence holds and integrations and summations can be exchanged), with $j\ge0, 0\le \ell\le n$, 
 \EQ{\nn 
 \hbar^\ell  \alpha^{2(\ell+j)} \int_{x} ^{x_0} u^{2(n-2\ell-j)-1}\, du &=  \frac{\hbar^\ell  \alpha^{2(\ell+j)}}{2(n-2\ell-j)}(x_0^{2(n-2\ell-j)}- x^{2(n-2\ell-j)}) \\
 &= O(\hbar^\ell  \alpha^{2(\ell+j)}x_0^{2(n-2\ell-j)}) + O(\hbar^\ell x^{2(n-\ell)} \eta^{2(\ell+j)})
 }
 if $n-2\ell-j\ne0$. Observe that 
 \EQ{\label{eq:xoreta}
  O(\hbar^\ell x^{2(n-\ell)} \eta^{2(\ell+j)}) &= O(x^2) \text{\ \ if\ \ }\ell=0\\
 O(\hbar^\ell x^{2(n-\ell)} \eta^{2(\ell+j)}) &= \hbar \, O(\eta^2) \text{\ \ if\ \ }\ell>0
 }
 If  $n-2\ell-j=0$, then
  \EQ{\nn 
 \hbar^\ell  \alpha^{2(\ell+j)} \int_{x} ^{x_0} u^{2(n-2\ell-j)-1}\, du &=  \hbar^\ell  \alpha^{2(n-\ell)}  (\log x_0 -\log x) \\
 &= \hbar^\ell  \alpha^{2(\ell+j)}  (\log x_0 -\log x)
 } 
 We have $n-\ell=\ell+j\ge1$ (if $\ell=0$, then $j=n\ge1$), so the $-\log x$ is multiplied with the small factor $\alpha^2$. In summary, in view of~\eqref{eq:gc2}, \eqref{eq:xoreta}, and Lemma~\ref{lem: Lemma 3.2 CDST} we conclude that 
\EQ{ \nn 
\frac{2}{3}\left(-\tau(x;\alpha,\hbar)\right)^{\frac{3}{2}}=-\log x + C(\alpha;\hbar) + x^2 f_1(x,\eta,\alpha;\hbar) + \hbar\eta^2 f_2(x,\eta,\alpha;\hbar) +\alpha^2 f_3(\alpha;\hbar)\log x
}
where $\eta=\alpha/x$, and all functions have the properties stated in the lemma. In particular, in order to estimate $\tau(x_{0};\alpha,\hbar)$ contributing to $C(\alpha,\hbar)$, we use the results in Lemma \ref{lem: Lemma 3.2 CDST}, by choosing $x_{0}\in J_{1}$.
\end{proof}

Based on Lemma~\ref{lem: Lemma 3.3 hardest}, we are finally able to describe $\tV$ in the one remaining case. 

\begin{proposition}\label{prop: tV close to 0 hardest}
In the parameter regime of the previous lemma,  and using the representation~\eqref{eq:tVvarphi}, 
\EQ{\nn 
|\partial_\tau^k \partial_\alpha^\ell\varphi(x;\alpha,\hbar)| 
&\le C_{k,\ell}\, (\hbar\alpha^2/x^2+x^2)\,  \alpha^{-\ell} (\sqrt{-\tau})^{k} \\
&\le C_{k,\ell}\, (\hbar\alpha^2/x^2+x^2)\,  \alpha^{-\ell}(  -\log x  )^{\frac{k}{3}}
}
for all $k,\ell\ge0$. 
\end{proposition}
\begin{proof}
By \eqref{eq:Q2 form},  with $\eta:=\alpha/x$, 
 \EQ{\nn 
 Q_2(x;\alpha,\hbar) &=  (1-2\hbar)^2 + \frac{4\hbar x^2}{x^2+\alpha^2} - \frac{4\hbar^2 x^2(x^2+3\alpha^2)}{(x^2+\alpha^2)^2} -x^2  \\
 & =  (1-2\hbar)^2 + \frac{4\hbar}{1+\eta^2} - \frac{4\hbar^2 \eta^2(\eta^2+3)}{(1+\eta^2)^2} - x^2 \\ 
xQ_2'(x;\alpha,\hbar) &= 8\hbar \alpha^2  \frac{x^2}{(x^2+\alpha^2)^2} + 8\hbar^2 \alpha^2 \frac{ x^2 (x^2-3\alpha^2)}{(x^2+\alpha^2)^3}-2x^2 \\
&= 8\hbar   \frac{\eta^2}{(1+\eta^2)^2} + 8\hbar^2  \frac{ \eta^2 (1-3\eta^2)}{(1+\eta^2)^3}-2 x^2 
}
Thus, recalling Lemma~\ref{lem:tilVxtau} for the definition of $\mu$, we have $|\mu|\le C(\hbar\eta^2+x^2 )$ and therefore 
\[
|\varphi(x;\alpha,\hbar)|\le C(\hbar\eta^2+x^2 ) ,\qquad |\partial_\alpha^\ell \varphi(x;\alpha,\hbar)|\le C(\hbar\eta^2+x^2 )\alpha^{-\ell}
\]
uniformly in the parameter regime of the proposition, for all $\ell\ge0$.  In view of \eqref{eq:taudx1}, \eqref{eq:taudx2}, $x\partial_x=-\eta\partial_\eta$,  and Lemma~\ref{lem: Lemma 3.3 hardest}, 
\EQ{\nn 
|\partial_\tau^k \varphi(x;\alpha,\hbar)| &\le C_k\, (\hbar\eta^2+x^2)\,  (\sqrt{-\tau})^{k}\\
&\le C_k\, (\hbar\eta^2+x^2)\,  (-\log x)^{\frac{k}{3}}
}
for all $k\ge0$, as well as 
\EQ{\nn 
|\partial_\tau^k \partial_\alpha^\ell\varphi(x;\alpha,\hbar)| 
&\le C_{k,\ell}\, (\hbar\eta^2+x^2)\,  \alpha^{-\ell} (\sqrt{-\tau}   )^{k}  \\
&\le C_{k,\ell}\, (\hbar\eta^2+x^2)\,  \alpha^{-\ell}(-\log x)^{\frac{k}{3}}
}
as claimed. 
\end{proof}

Propositions \ref{prop: tV close xt}, \ref{prop: tV x large}, \ref{prop: tV close to 0}, \ref{prop: tV close to 0 harder} and \ref{prop: tV close to 0 hardest} provide the following complete description for the potential $\tV(\tau;\alpha,\hbar)$. We ignore a factor of $\hbar$ which might appear upon differentiation with respect to~$\alpha$. 

\begin{proposition}
	\label{prop:tV summary}
	There exists a  constant $\tau_*>0$ and a small constant $0<\alpha_*\ll1$ so that uniformly in $\hbar\in(0,\frac13]$, 
	\EQ{\label{eq:ptV}
		|\partial_\tau^k\partial_\alpha^\ell \tV(\tau;\alpha,\hbar)| \le C_{k,\ell}\, \langle  \alpha\rangle^{-\ell-1} \langle  \tau\rangle^{-2-k},\qquad \forall\; \tau\ge -\tau_*
	}
	for all $k,\ell\ge0$ and $\alpha>0$. Moreover, \eqref{eq:ptV} holds for $-\infty <\tau\le -\tau_*$, 
	all $k,\ell\ge0$ and $\alpha\ge \alpha_*$. Finally, if $0<\alpha\le \alpha_*$ and $-\infty <\tau\le -\tau_*$, then  
	\[
	\tV(\tau;\alpha,\hbar) = \frac{5}{16\tau^2} - \tau \varphi(x;\alpha,\hbar)
	\]
	where for all $k,\ell\ge0$
	\EQ{\label{eq:phi new}
		|\partial_\tau^k \partial_\alpha^\ell\varphi(x;\alpha,\hbar)| 
		&\le C_{k,\ell}\, \min\big(\hbar\, \alpha^2 x^{-2}+x^2, \hbar x^2/\alpha^2 + x^2\big)  \alpha^{-\ell} (  -\tau  )^{\frac{k}{2}} 
	}
	Here $x=x(\tau;\alpha,\hbar)$ is the inverse of the diffeomorphism $\tau=\tau(x;\alpha,\hbar)$ defined in~\eqref{def tau},  and satisfies
	\EQ{\label{eq:xtau}
		\frac{2}{3}\left(-\tau(x;\alpha,\hbar)\right)^{\frac{3}{2}}  = \left\{   \begin{array}{lll} & -(1-O(\alpha^2))\log x + O_{1}(1), & \text{\ \ if\ \ }  0< \alpha\le x\le x_*:=x(-\tau_*;\alpha,\hbar)  \\
			&  -\log x + 2\hbar\log(x/\alpha) + O_{2}(1),  &\text{\ \ if\ \ }    0< x\le \alpha\le \alpha_*
		\end{array} 
		\right. 
	}
	Here $O(\alpha^2)$ is analytic in complex $|\alpha|\le\alpha_*$, and bounded uniformly in  $\hbar\in(0,\frac13]$. Furthermore, the two terms $O_1(1)$, resp.~$O_2(1)$ refer to smooth functions of $\tau,\alpha$ (and thus also of $x$), uniformly bounded in $0<\alpha\le\alpha_*$, $-\infty< \tau\le-\tau_*$, and so that for all $k,\ell\ge0$ one has $\partial_x^k \partial_\alpha^\ell O_1(1)=O(x^{-k} \alpha^{-\ell})$ in the parameter regime of the first line of~\eqref{eq:xtau}, resp.\ $\partial_x^k \partial_\alpha^\ell O_2(1)=O(\alpha^{-k-\ell})$ in the regime of the second line of~\eqref{eq:xtau}.
\end{proposition}

\subsection{Fundamental system for the perturbed Airy equation}\label{sec:fund sys}
In this section we will analyze a fundamental system of solutions of the  perturbed Airy equation \eqref{perturbed Airy}, viz. 
\begin{align}\nn
\hbar^{2}\ddot{w}(\zeta;\hbar)=\zeta w(\zeta;\hbar)+\hbar^{2}\tV(-\zeta,\alpha;\hbar)w(\zeta;\hbar).
\end{align}
Here $\tV$ is as in Proposition~\ref{prop:tV summary}, but we switched to the independent variable $\zeta=-\tau$. 
With $\tau_*$ as in that proposition, we set $\zeta_*:=-(-\tau_*)=\tau_*$. 

Recall that the unperturbed semiclassical Airy equation 
\EQ{\label{eq:Airy0}
	\hbar^{2}\ddot{w}_0(\zeta;\hbar)=\zeta w_0(\zeta;\hbar)
}
has a fundamental system $\Ai(\hbar^{-\frac23} \zeta), \Bi(\hbar^{-\frac23} \zeta)$ which are positive for all $\zeta\ge0$. For $\zeta\le0$ we switch to the complex system 
\[
\Ai(\hbar^{-\frac23} \zeta)\pm i \Bi(\hbar^{-\frac23} \zeta)
\]
which does not vanish.  Throughout this section, we shall frequently use the following standard facts about Volterra integral equations
of the form
\EQ{\label{eq:Volt*}
 f(x)= g(x)+\int_{x}^{\infty} K(x,s) f(s) \, ds,
}
or
\EQ{\label{eq:Volt**}
 f(x)= g(x)+\int_{a}^{x} K(x,s) f(s) \, ds,
}
with some $g(x)\in L^{\infty}$ and $a\in \mathbb{R}$. 
These equations are solved by means of an iteration which crucially relies on the directedness of the variables.  This refers to $s>x$ in~\eqref{eq:Volt*} and  $s<x$ in~\eqref{eq:Volt**}. This ordering leads to a gain of $n!$ after $n$ iterations, as for the exponential series.

\begin{lemma}[Lemma 2.4 in~\cite{SSS1}]
Suppose $g(x)\in L^{\infty}([a,\infty))$ and 
\[\mu:= \int_{a}^\infty \sup_{a<x<s} |K(x,s)|\, ds<\infty\] 
There exists a unique solution to~\eqref{eq:Volt*} of the form 
\EQ{
\label{eq:volt_it} f(x) = g(x) + \sum_{n=1}^\infty \int_a^\infty
\ldots \int_a^\infty \prod_{i=1}^n \chi_{[x_{i-1}<x_i]}
K(x_{i-1},x_i) \; g(x_{n}) \, dx_{n}\ldots dx_1.
}
with $x_0:=x$.  It satisfies the bound
\[ \|f\|_{L^\infty(a,\infty)} \le e^\mu \|g\|_{L^\infty(a,\infty)}. \]
An analogous statement holds for \eqref{eq:Volt**}.
\end{lemma}

We refer the reader to~\cite{SSS1} (or elsewhere) for the elementary proof. 

\begin{lemma}
	\label{lem:solve pA}
	Let $w_0(\zeta;\hbar):=\Ai(\hbar^{-\frac23} \zeta)$ for $\zeta\ge0$ and $w_1(\zeta;\hbar):=\Ai(\hbar^{-\frac23} \zeta)+i \Bi(\hbar^{-\frac23} \zeta)$ for $\zeta\le0$. 
	Then the  Volterra integral equation  	
	\EQ{\label{eq:a0 system} 
		a_0(\zeta;\alpha,\hbar) & := \int_\zeta^\infty K_0(\zeta,s;\alpha,\hbar) (1+\hbar a_0(s;\alpha,\hbar))\, ds \\
		K_0(\zeta,s;\alpha,\hbar)  & = \hbar^{-1} \tV(-s;\alpha,\hbar) w_0^2(s; \hbar) \int_\zeta^s w_0^{-2}(t;\hbar)  \, dt 
	}
	has a unique bounded solution $a_0(\zeta;\alpha,\hbar)$ for all $\hbar\in(0,\frac13]$ and~$\alpha>0$, $\zeta\ge0$. One has $\lim_{\zeta\to\infty} a_0(\zeta;\alpha,\hbar)=0$ and 
	$w(\zeta;\alpha,\hbar):=w_0(\zeta;\hbar)(1+\hbar a_0(\zeta;\alpha,\hbar))$ is the unique solution of~\eqref{perturbed Airy} on $[0,\infty)$ with 
	$w(\zeta;\alpha,\hbar)\sim w_0(\zeta;\hbar)$ as $\zeta\to\infty$. 
	Analogously,  the  Volterra integral equation  	
	\EQ{\label{eq:a1 system} 
		a_1(\zeta;\alpha,\hbar) & := \int^\zeta_{-\infty} K_1(\zeta,s;\alpha,\hbar) (1+\hbar a_1(s;\alpha,\hbar))\, ds \\
		K_1(\zeta,s;\alpha,\hbar)  & = \hbar^{-1} \tV(-s;\alpha,\hbar) w_1^2(s;\hbar) \int^\zeta_s w_1^{-2}(t;\hbar)  \, dt 
	}
	has a unique bounded solution $a_1(\zeta;\alpha,\hbar)$ for all $\hbar\in(0,\frac13]$ and~$\alpha>0$, $\zeta\le0$. One has $\lim_{\zeta\to-\infty} a_1(\zeta;\alpha,\hbar)=0$ and 
	$w(\zeta;\alpha,\hbar):=w_1(\zeta;\hbar)(1+\hbar a_1(\zeta;\alpha,\hbar))$ is the unique solution of~\eqref{perturbed Airy} on $(-\infty,0]$ with $w(\zeta;\alpha,\hbar)\sim 
	w_1(\zeta;\hbar)$ as $\zeta\to-\infty$. 
\end{lemma}
\begin{proof}
	For simplicity, we suppress the parameters $\alpha,\hbar$ in the notation since they are fixed for the purposes of this lemma.  Suppose $w(\zeta):=w_0(\zeta)(1+\hbar a_0(\zeta))$ solves~\eqref{perturbed Airy}. Then 
	\EQ{\label{eq:a0ode}
		\ddot{w} &= \ddot{w}_0(1+\hbar a_{0}) + 2\hbar \dot{w}_0 \dot{a}_0 +\hbar w_0 \ddot{a}_0\\
		\hbar^2 \ddot{w} &= \zeta w(\zeta)  + \hbar^3( 2 \dot{w}_0 \dot{a}_0 + w_0 \ddot{a}_0 ) \\
		&= \zeta w(\zeta)  + \hbar^2 \tV(-\zeta) w_0(\zeta)(1+\hbar a_0(\zeta))
	} 
	whence 
	\EQ{\label{eq:w0a0}
		(w_0^2(\zeta) \dot{a}_0(\zeta))^{\dot{\ }} = \hbar^{-1} \tV(-\zeta) w_0^2(\zeta) (1+\hbar a_0(\zeta))
	}
	If a bounded solution to \eqref{eq:w0a0} exists, then it is given by 
	\EQ{\label{eq:Volt1}
		w_0^2(\zeta) \dot{a}_0(\zeta) & = -\hbar^{-1} \int_\zeta^\infty  \tV(-s) w_0^2(s) (1+\hbar a_0(s))\, ds\\
		a_0(\zeta) &= \hbar^{-1} \int_\zeta^\infty w_0^{-2}(t)  \int_t^\infty \tV(-s) w_0^2(s) (1+\hbar a_0(s))\, ds \, dt \\
		&= \hbar^{-1} \int_\zeta^\infty \int_\zeta^s w_0^{-2}(t) \, dt\,   \tV(-s) w_0^2(s) (1+\hbar a_0(s))\, ds   
	}
	Recall the well-known asymptotic behavior
	\begin{align}\label{eq:Ai}
	\Ai(x) =  (4\pi)^{-\frac12} x^{-\frac14} e^{-\frac23 x^{\frac32}}(1+\alpha(x)), \quad |\alpha^{(k)}(x)|\leq C_k x^{-\frac32 - k},\quad x\geq 1
	\end{align}
	Therefore, for large $s\ge \zeta\ge\hbar^{\frac{2}{3}}$, 
	\EQ{\label{eq:w0int}
		\int_\zeta^s  w_0^2(s) w_0^{-2}(t) \, dt &\le C \int_\zeta^s  e^{-\frac{4}{3\hbar}  (s^{\frac32}-t^{\frac32})} \, (t/s)^{\frac12}\, dt \\
			&\le Cs^{-\frac{1}{2}}\int_{\zeta}^{s}e^{-\frac{4}{3\hbar}  (s^{\frac32}-t^{\frac32})}dt^{\frac{3}{2}}\\
			&\le C\hbar s^{-\frac{1}{2}}.}
		
	By Proposition~\ref{prop:tV summary} we have $\tV(-s)=O(s^{-2})$ as $s\to\infty$ and thus
	\[
	\max_{\zeta\le s} \Big| \int_\zeta^s w_0^{-2}(t) \, dt\,   \tV(-s) w_0^2(s)\Big|  = O(\hbar s^{-\frac52}) \qquad s\to\infty
	\]
	By a standard Volterra iteration, \eqref{eq:Volt1} has a unique bounded solution, and in fact $a_0(\zeta)=O(\zeta^{-\frac32})$ as $\zeta\to\infty$.  
	The argument can be reversed, and defining $w(\zeta):=w_0(\zeta)(1+\hbar a_0(\zeta))$ gives a solution of~\eqref{perturbed Airy}  with $w(\zeta)\sim w_0(\zeta)$ as $\zeta\to\infty$. 
	Conversely, any solution of this type satisfies $w(\zeta)=w_0(\zeta)(1+b(\zeta))$ with $b=o(1)$ as $\zeta\to\infty$. Repeating the previous calculations they yields $b=\hbar a_0$, as claimed. 
	
	For negative $\zeta$ we proceed in an analogous fashion, arriving at the same equation~\eqref{eq:a0ode} with $a_1, w_1$ in place of $a_0, w_0$. Integrating it from 
	$-\infty$ to $\zeta\le0$ yields~\eqref{eq:a1 system}.  The well-known asymptotic behavior of $w_1$ is given by, see for example Corollary~C.4 in~\cite{CDST}, 
	\EQ{\label{eq:AiiBi}
		\Ai(-x)+i\Bi(-x) &= c x^{-\frac14} e^{-\frac{2i}{3}x^{\frac32}}[1+b(x)],\qquad x\ge1 \\
		|b^{(k)}(x)| &\le C\langle x\rangle^{-\frac32-k} \quad \forall\; k\ge0
	}
	whence for all $s\le\zeta\le -\hbar^{\frac{2}{3}}$ 
	\EQ{\nn
		\Big| w_1^2(s;\alpha, \hbar) \int^\zeta_s w_1^{-2}(t;\alpha,\hbar)  \, dt \Big| &\le \hbar |s|^{-\frac12} \Big|  \int_{-\hbar^{-\frac23}\zeta}^{-\hbar^{-\frac23} s}  t^{\frac12} e^{\frac{4i}{3}t^{\frac32}}[1+b(t)]^{-2} \, dt \Big| \\
		&\le  \hbar |s|^{-\frac12} \Big|  \int_{-\hbar^{-\frac23}\zeta}^{-\hbar^{-\frac23}s}  [1+b(t)]^{-2} \, d( e^{\frac{4i}{3}t^{\frac32}})  \Big| \le \hbar |s|^{-\frac12}
	}
	where the final inequality follows by integration by parts. 
	This is the same bound as~\eqref{eq:w0int}. Hence, the argument proceeds as before. 
\end{proof}

The previous lemma introduces the two fundamental systems of the perturbed Airy equation \eqref{perturbed Airy}. We have for $\zeta\le0$ the oscillatory solutions
\EQ{\label{eq:sys1}
	\big\{ w_1(\zeta;\hbar)(1+\hbar a_1(\zeta;\alpha,\hbar)),\; \overline{w_1(\zeta;\hbar)}(1+\hbar \overline{a_1(\zeta;\alpha,\hbar)}) \big\} 
}
as well as on $\zeta\ge0$ the real-valued solutions 
\EQ{\label{eq:sys0}
	\big\{ f_1(\zeta;\alpha,\hbar):=w_0(\zeta;\hbar)(1+\hbar a_0(\zeta;\alpha,\hbar)),\; f_2(\zeta;\alpha,\hbar) \big\} 
}
where $f_2$ is the growing solution obtained from the decaying one $f_1$ by the usual reduction ansatz: $f_2= gf_1$ where $f_1\ddot g + 2\dot g \dot{f}_1=0$. Thus,
\[
g(\zeta;\alpha,\hbar)= \int_{\zeta_0}^\zeta f_1(t;\alpha,\hbar)^{-2}\, dt, \quad f_2(\zeta;\alpha,\hbar) =  f_1(\zeta;\alpha,\hbar) \int_{\zeta_0}^\zeta f_1(t;\alpha,\hbar)^{-2}\, dt
\]
Here $\zeta_0\ge0$ is chosen such that $f_1(\zeta;\alpha,\hbar)>0$ for all $\zeta\ge\zeta_0$. The remainder of this section analyses the fundamental systems~\eqref{eq:sys0}, \eqref{eq:sys1} in more detail. For example, we need to show that $a_0, a_1$ remain uniformly bounded in all variables, obtain their decay in~$\zeta$, and we also need to  bound the derivatives of these functions.  The general treatment of Volterra equations with Airy kernels in the appendices in~\cite{CDST} does not  cover our problem due to the more delicate  behavior of~$\tV(-\zeta;\alpha,\hbar)$ in the regime $0<\alpha<1$ and $\zeta\ge0$.  We remark that the equations~\eqref{eq:a1 system}  hold not just on $\zeta\le0$ but also on, say, $\zeta\le\zeta_*$.     In the following lemma we treat the case $\zeta\le0 $ using the methods from the appendices in~\cite{CDST}. 

\begin{lemma}
	\label{lem:a1fine}
	The functions $a_1(\zeta;\alpha,\hbar)$ from Lemma~\ref{lem:solve pA} satisfy the bounds
	\EQ{
		\label{eq:a1 fine}
		|\partial_{\alpha}^{\ell}a_{1}(\zeta;\alpha,\hbar)| & \leq C_{\ell}\langle\alpha\rangle^{-\ell}
		\langle\zeta\rangle^{-\frac{3}{2}},\quad \zeta\leq 0, \\ 
		|\partial_\zeta^{k} \partial_\alpha^\ell a_1(\zeta;\alpha,\hbar)| & \le C_{k,\ell}\, \langle \alpha\rangle^{-\ell} \left\{ \begin{array}{cc} 
			|\zeta|^{-\frac{3}{2}-k} & -\infty<\zeta\le -1 \\
			|\zeta|^{\frac{1}{2}-k} & -1<\zeta\le -\hbar^{\frac23} \\
			\hbar^{\frac{1-2k}{3}}&  -\hbar^{\frac23}\le\zeta\le 0, 
		\end{array}\right. 
	}
	for all $\ell\ge0, k\geq 1$, $\alpha>0$, and $\hbar\in (0,\frac13]$. 
\end{lemma}
\begin{proof}
	In view of~\eqref{eq:a1 system} and \eqref{eq:AiiBi}, we have for all $s\le\zeta\le-\hbar^{\frac23}$, 
	\EQ{\nn 
		|K_1(\zeta,s;\alpha,\hbar)| &\le  C\hbar^{-1} \langle s\rangle^{-2} \langle\hbar^{-\frac23}s\rangle^{-\frac12}   (1+|b(-\hbar^{-\frac23}s)|)^2 
		\Big| \int_s^\zeta \langle\hbar^{-\frac23}t\rangle^{\frac12} e^{\frac{4i}{3\hbar}(-t)^{\frac32}} (1+b(-\hbar^{-\frac23}t))^{-2}\, dt \Big|  \\
		& = C \hbar^{-\frac13} \langle s\rangle^{-2} \langle\hbar^{-\frac23}s\rangle^{-\frac12}   (1+|b(-\hbar^{-\frac23}s)|)^2 
		\Big| \int_{\hbar^{-\frac23}s}^{\hbar^{-\frac23}\zeta} \langle t\rangle^{\frac12} e^{\frac{4i}{3}(-t)^{\frac32}} (1+b(-t))^{-2}\, dt \Big|
	}
	Next, 
	\EQ{\nn 
		& \int_{t_0}^{t_1} \langle t\rangle^{\frac12} e^{\frac{4i}{3}t^{\frac32}} (1+b(t))^{-2}\, dt = \frac{1}{2i}\int_{t_0}^{t_1} t^{-\frac12} \langle t\rangle^{\frac12}  (1+b(t))^{-2}\, d\Big(e^{\frac{4i}{3}t^{\frac32}}\Big)\\
		& = t^{-\frac12} \langle t\rangle^{\frac12}  (1+b(t))^{-2}e^{\frac{4i}{3}t^{\frac32}} \Big|_{t_0}^{t_1} -\frac{1}{2i}\int_{t_0}^{t_1}  e^{\frac{4i}{3}t^{\frac32}} \, d\big( t^{-\frac12} \langle t\rangle^{\frac12}  (1+b(t))^{-2}\big)  =O(1)
	}
	uniformly in $1\le t_0\le t_1$. Here we used that $\frac{d}{dt}\big( t^{-\frac12} \langle t\rangle^{\frac12}  (1+b(t))^{-2}\big)=O(t^{-1})$ as $t\rightarrow\infty$. It follows from the preceding that
	\[
	|K_1(\zeta,s;\alpha,\hbar)| \le  C\hbar^{-\frac13} \langle s\rangle^{-2} \langle\hbar^{-\frac23}s\rangle^{-\frac12}\le C\langle s\rangle^{-2}|s|^{-\frac12}
	\]
	for all $s\le\zeta\le-\hbar^{\frac23}$. If  $-\hbar^{\frac{2}{3}}\leq s\leq \zeta\leq 0$, then 
	\EQ{\label{eq:K1small} 
		|K_{1}(\zeta,s;\alpha,\hbar)|\leq C\hbar^{-1}\int_{s}^{\zeta}dt\leq C\hbar^{-\frac{1}{3}} 
	} 
	If $s\leq -\hbar^{\frac{2}{3}}\leq \zeta\leq 0$, we split the integral $\int_{s}^{\zeta}$ in the form 
	\begin{align*}
	\int_{s}^{\zeta}=\int_{s}^{-\hbar^{\frac{2}{3}}}+\int_{-\hbar^{\frac{2}{3}}}^{\zeta}.
	\end{align*}
	The contribution from the first integral is treated in the same way as the case $s\leq \zeta\leq -\hbar^{\frac{2}{3}}$, and the contribution from the second integral is bounded by~\eqref{eq:K1small}.
	By a standard Volterra iteration applied to the equation~\eqref{eq:a1 system} we obtain from this bound
	\[
	|a_1(\zeta;\alpha,\hbar)|\le C\langle \zeta\rangle^{-\frac32}\qquad\forall\; \zeta\le 0,
	\]
	uniformly in $\alpha>0$, and $\hbar\in (0,\frac13]$. The loss of $\hbar^{-\frac{1}{3}}$ in~\eqref{eq:K1small} is absorbed by an integration interval of length~$\hbar^{\frac{2}{3}}$. Combining with the contribution from the range $-\infty<s\leq -\hbar^{\frac{2}{3}}$, we obtain the above estimate. Taking derivatives in $\alpha$ yields
	\EQ{\label{eq:a1ell} 
		\partial_\alpha^\ell a_1(\zeta;\alpha,\hbar) & :=  \sum_{j=0}^\ell \binom{\ell}{j} \int^\zeta_{-\infty} \partial_\alpha^{j} K_1(\zeta,s;\alpha,\hbar) \partial_\alpha^{\ell-j} (1+\hbar a_1(s;\alpha,\hbar))\, ds 
	}
	We have 
	\[
	\partial_\alpha^{j} K_1(\zeta,s;\alpha,\hbar)   = \hbar^{-1} \partial_\alpha^{j}\tV(-s;\alpha,\hbar) w_1^2(s;\hbar) \int^\zeta_s w_1^{-2}(t;\hbar)  \, dt  
	\]
	and by Proposition~\ref{prop:tV summary}, 
	\[
	\Big| \partial_\alpha^{j} K_1(\zeta,s;\alpha,\hbar) \Big| \le C_j\, \langle s\rangle^{-2}|s|^{-\frac12}\langle \alpha\rangle^{-j} 
	\]
	for all $j\ge0$ and all $s\leq \zeta\leq 0$.     These give \eqref{eq:a1 fine} for $k=0$ by induction in $\ell$.  Indeed, \eqref{eq:a1ell} is a Volterra equation for $\partial_\alpha^\ell a_1(\zeta;\alpha,\hbar)$ and the lower order derivatives $\partial_\alpha^j a_1(\zeta;\alpha,\hbar)$ with $0\le j<\ell$ are estimated by means of the induction assumption. See~\cite[Proposition~B.1]{CDST} for more details. 
	
	The derivatives in~$\zeta$ are more delicate and we use the method from~\cite[Proposition~C.5]{CDST}. 
	Denote $S(x):=\Ai(x)+i\Bi(x)$ and $\tilde a_1(\zeta)=\tilde a_1(\zeta;\alpha,\hbar):= a_1(\hbar^{\frac23}\zeta)$. Then 
	\EQ{\nn 
		\tilde a_1(\zeta)  & := \hbar^{\frac23} \int^{\zeta}_{-\infty} K_1(\hbar^{\frac23} \zeta,\hbar^{\frac23} s;\alpha,\hbar) (1+\hbar \tilde a_1(s))\, ds   \\
		&= \hbar^{\frac13} \int^{ \zeta}_{-\infty}  \tV(-\hbar^{\frac23} s;\alpha,\hbar) w_1^2(\hbar^{\frac23} s;\hbar) \int^\zeta_s w_1^{-2}(\hbar^{\frac23} t;\hbar)  \, dt \; (1+\hbar \tilde a_1(s))\, ds \\
		&= \hbar^{\frac13} \int^{ \zeta}_{-\infty}  \tV(-\hbar^{\frac23} s;\alpha,\hbar) S(s)^2 \int^\zeta_s S(t)^{-2}  \, dt \; (1+\hbar \tilde a_1(s))\, ds 
	}
	If $\zeta\le-1$, then by \eqref{eq:AiiBi} we have 
	\EQ{\nn 
		\tilde a_1\left(\zeta\right) & = c \hbar^{\frac13} \int_{-\zeta}^{\infty}  \tV(\hbar^{\frac23} s;\alpha,\hbar) s^{-\frac12} e^{-\frac{4i}{3}s^{\frac32}}\big(1+b(s)\big)^2   \int_{-\zeta}^{s}  t^{\frac12} e^{\frac{4i}{3}t^{\frac32}}\big(1+b(t)\big)^{-2} \, dt \; (1+\hbar \tilde a_1(-s))\, ds 
	}
	Changing variables $s^{\frac32}=\sigma$, $t^{\frac32}=\tau$ and setting $a_2(u) = \tilde a_1(-u^{\frac23})=a_1(-(\hbar u)^{\frac23})$, $\beta(\sigma):=b(\sigma^{\frac23})$, and suppressing $\alpha,\hbar$ as arguments from the notation,  we obtain for all $u\ge1$
	\EQ{\nn
		& a_2(u) = c  \hbar^{\frac13} \int_{u}^{\infty}  \tV((\hbar \sigma)^{\frac23})  \big(1+\beta(\sigma) \big)^2   \int_{u}^{\sigma}  e^{\frac{4i}{3}(\tau-\sigma) }\big(1+\beta(\tau)\big)^{-2} \, d\tau \; (1+\hbar  a_2(\sigma))\, \sigma^{-\frac23}\, d\sigma \\
		&= c  \hbar^{\frac13} \int_{0}^{\infty}  \tV((\hbar (u+v))^{\frac23})  \big(1+\beta(u+v) \big)^2  \!\! \int_{0}^{v}  e^{\frac{4i}{3}(w-v) }\big(1+\beta(w+u)\big)^{-2}  dw \, (1+\hbar  a_2(u+v))\, (u+v)^{-\frac23}\, dv 
	}
	Note that the exterior variable $u$ does not appear in the phase of the complex exponential, as in~\cite[Proposition~B.1]{CDST}.  This is important as we differentiate in~$u$.  In fact, 
	\EQ{\nn
		& a_2'(u) = c_1  |\thbar| \int_{0}^{\infty}  \tV'((\hbar (u+v))^{\frac23})   \big(1+\beta(u+v) \big)^2  \!\! \int_{0}^{v}  e^{\frac{4i}{3}(w-v) }\big(1+\beta(w+u)\big)^{-2} \, dw \, (1+\hbar  a_2(u+v))\, (u+v)^{-1}\, dv \\
		& + c _2\hbar^{\frac13} \int_{0}^{\infty}  \tV((\hbar (u+v))^{\frac23})  \big(1+\beta(u+v) \big)\beta'(u+v)  \!\! \int_{0}^{v}  e^{\frac{4i}{3}(w-v) }\big(1+\beta(w+u)\big)^{-2} \,  dw \, (1+\hbar  a_2(u+v))\, (u+v)^{-\frac23}\, dv \\
		& + c _3\hbar^{\frac13} \int_{0}^{\infty}  \tV((\hbar (u+v))^{\frac23})  \big(1+\beta(u+v) \big)^{2}  \!\! \int_{0}^{v}  e^{\frac{4i}{3}(w-v) }\big(1+\beta(w+u)\big)^{-3} \beta'(u+w)\, dw \, (1+\hbar  a_2(u+v))\, (u+v)^{-\frac23}\, dv \\
		& + c_4  \hbar^{\frac13} \int_{0}^{\infty}  \tV((\hbar (u+v))^{\frac23})  \big(1+\beta(u+v) \big)^2  \!\! \int_{0}^{v}  e^{\frac{4i}{3}(w-v) }\big(1+\beta(w+u)\big)^{-2} \, dw \, (1+\hbar  a_2(u+v))\, (u+v)^{-\frac53}\, dv \\
		& + c_5  \hbar^{\frac43} \int_{0}^{\infty}  \tV((\hbar (u+v))^{\frac23})  \big(1+\beta(u+v) \big)^2  \!\! \int_{0}^{v}  e^{\frac{4i}{3}(w-v) }\big(1+\beta(w+u)\big)^{-2} \,  dw \,   a_2'(u+v)\, (u+v)^{-\frac23}\, dv 
	}
	By \eqref{eq:AiiBi} we have $|\beta^{(k)}(\sigma)|\le C_k\, \sigma^{-1-k}$ for all $k\ge0$ and $\sigma\ge1$.  Integrating by parts yields 
	\EQ{\nn
		\sup_{v\ge0} \Big| \int_{0}^{v}  e^{\frac{4i}{3}(w-v) }\big(1+\beta(w+u)\big)^{-2} \, dw \Big| &\lesssim 1 \\ 
		\sup_{v\ge0} \Big| \int_{0}^{v}  e^{\frac{4i}{3}(w-v) }\big(1+\beta(w+u)\big)^{-3} \beta'(u+w) \, dw \Big| & \lesssim u^{-2}
	}
	for all $u\ge1$.  It follows from these bounds, Proposition~\ref{prop:tV summary},  and the previous uniform bound on $a_1$,  and thus on $a_2$, that for all $u\ge1$, 
	\EQ{\nn
		|a_2'(u)|  &\lesssim   \hbar \int_{0}^{\infty}  \langle (\hbar (u+v))^{\frac23} \rangle^{-3}   \, (u+v)^{-1}\, dv +  \hbar^{\frac13} \int_{0}^{\infty}  \langle (\hbar (u+v))^{\frac23} \rangle^{-2}  \,  u^{-2}\, (u+v)^{-\frac23}\, dv   \\
		& \quad +      \hbar^{\frac13} \int_{0}^{\infty}  \langle (\hbar (u+v))^{\frac23} \rangle^{-2}  (u+v)^{-\frac53}\, dv   +   \hbar^{\frac43} \int_{0}^{\infty}  \langle (\hbar (u+v))^{\frac23} \rangle^{-2} \, (u+v)^{-\frac23}\, | a_2'(u+v)|\,  dv \\
		& \lesssim    \int_{u}^{\infty}  \big[ \hbar \langle (\hbar w)^{\frac23} \rangle^{-3}   \,w^{-1}  + \hbar^{\frac13} \langle (\hbar w)^{\frac23} \rangle^{-2}  \, ( u^{-2}\, w^{-\frac23}+ w^{-\frac53}) \big]\, dw + \hbar^{\frac43} \int_{u}^{\infty}  \langle (\hbar w)^{\frac23} \rangle^{-2} \, w^{-\frac23}\, | a_2'(w)|\,  dw
	}	
	If $u\gtrsim\hbar^{-1}$, we therefore have 
	\EQ{\nn
		|a_2'(u)|  
		& \lesssim   \hbar^{-1}  \int_{u}^{\infty}  \big(  w^{-3}  +  u^{-2}\, w^{-2}\big)\, dw +  \int_{u}^{\infty}   w^{-2}\, | a_2'(w)|\,  dw \\
		& \lesssim   \hbar^{-1}  u^{-2} +  \int_{u}^{\infty}   w^{-2}\, | a_2'(w)|\,  dw  
	}	
	which yields upon iteration that $ |a_2'(u)|\lesssim \hbar^{-1}  u^{-2}$ for all $u\ge\hbar^{-1}$. On the other hand, if $1\le u\ll \hbar^{-1}$, then 
	\EQ{\nn
		|a_2'(u)|  &\lesssim \int_{u}^{\hbar^{-1}}  \big[  \hbar w^{-1} + \hbar^{\frac13} (u^{-2} w^{-\frac23} +  w^{-\frac53})  \big] \, dw  + \hbar + \hbar^{\frac43} \int_u^{\hbar^{-1}} w^{-\frac23} | a_2'(w)|\,  dw  + \int_{\hbar^{-1}}^\infty w^{-2}  | a_2'(w)|\,  dw \\
		&\lesssim -  \hbar\log (\hbar u) + \hbar^{\frac13} u^{-\frac23} + \hbar^{\frac43} \int_u^{\hbar^{-1}} w^{-\frac23} | a_2'(w)|\,  dw  \\
		&\lesssim  \hbar^{\frac13} u^{-\frac23} + \hbar^{\frac43} \int_u^{\hbar^{-1}} w^{-\frac23} | a_2'(w)|\,  dw 
	} 
	By Volterra iteration, $|a_2'(u)|  \lesssim \hbar^{\frac13} u^{-\frac23} $ for all $1\le u\ll \hbar^{-1}$. 
	We have 
	\begin{align}\label{a1 a2 1st deri}
	a'_{2}(u)=\frac23 a'_{1}(-(\hbar u)^{\frac{2}{3}})(-\hbar^{\frac{2}{3}}) u^{-\frac13}\quad 
	\end{align}
	Thus, redefining $\zeta:=-\hbar^{\frac{2}{3}}u^{\frac23}$, we obtain
	\[
	|\dot a_1(\zeta)|\lesssim \left\{ \begin{array}{ll} |\zeta|^{-\frac12} &  \forall -1\le \zeta\le -\hbar^{\frac23} \\
	|\zeta|^{-\frac52} &  \forall \;  \zeta\le -1 
	\end{array}\right. 
	\]
	as claimed in \eqref{eq:a1 fine} for $\ell=0$.  Next, we need to discuss the case $-\hbar^{\frac23}\le\zeta\le 0$. For this we go back to the ODE for $\dot a_1$, i.e.,  the analogue of ~\eqref{eq:w0a0} for $a_{1}, w_{1}$.
	Integrating once we arrive at (again, not writing $\alpha,\hbar$ as arguments) 
	\begin{align}\label{a1 1st deri}
	\begin{split}
		\dot a_1(\zeta) &= w_1^{-2}(\zeta) w_1^2(-\hbar^{\frac23}) \dot a_1(-\hbar^{\frac23}) + \hbar^{-1}  w_1^{-2}(\zeta) \int_{-\hbar^{\frac23}}^\zeta \tilde V(-s) w_1^2(s)(1+\hbar a_1(s))\, ds \\
		&= w_1^{-2}(\zeta) w_1^2(-\hbar^{\frac23}) \dot a_1(-\hbar^{\frac23}) + \hbar^{-\frac13}  w_1^{-2}(\zeta) \int^{1}_{-\hbar^{-\frac23}\zeta} \tilde V(\hbar^{\frac23} s) S^2(-s)(1+\hbar a_1(-\hbar^{\frac23}s))\, ds \\
		& = O(\hbar^{-\frac13})
		\end{split}
	\end{align}
	by the already established bound for $ \dot a_1(-\hbar^{\frac23})$. 
	This concludes the proof of \eqref{eq:a1 fine} for $\ell=0$. 
	The previous method extends to all $\ell\ge1$, simply by combining the argument for the derivatives relative to~$\alpha$, see~\eqref{eq:a1ell},  with the previous derivation. Next we turn to the estimates for higher order derivatives in $\zeta$. The method is very similar to that for estimating $a'_{2}(u)$. Using Leibniz's rule, we have
	\begin{align*}
	a^{(k)}_{2}(u)=&\sum_{p+q+r+s+m= k}C_{p,q,r,s,m}\, \hbar^{\frac{1}{3}}\int_{0}^{\infty}\frac{d^{p}}{du^{p}}\left(\tV((\hbar (u+v))^{\frac23})\right)\frac{d^{q}}{du^{q}}\left(\big(1+\beta(u+v) \big)^2\right)\\
	&\cdot\int_{0}^{v}e^{\frac{4i}{3}(w-v) }\frac{d^{r}}{du^{r}}\left(\big(1+\beta(w+u)\big)^{-2}\right)dw\frac{d^{s}}{du^{s}}\left((1+\hbar  a_2(u+v)) \right)\frac{d^{m}}{du^{m}}\left((u+v)^{-\frac23}\right)dv.
\end{align*}
For each factor of the above formula, we have the following pointwise estimate:
\begin{align*}
	&\left|\frac{d^{m}}{du^{m}}\left((u+v)^{-\frac23}\right)\right|\lesssim (u+v)^{-\frac{2}{3}-m},\quad m\ge0, \\
	&\frac{d^{s}}{du^{s}}\left(1+\hbar  a_2(u+v) \right)=\hbar a^{(s)}_{2}(u+v),\quad s\geq 1,\\
	&\left|\frac{d^{q}}{du^{q}}\left(\big(1+\beta(u+v) \big)^2\right)\right|\lesssim (u+v)^{-1-q},\quad q\geq 1,\\
	&\left|\frac{d^{p}}{du^{p}}\left(\tV((\hbar (u+v))^{\frac23})\right)\right|\lesssim\left\langle\hbar^{\frac{4}{3}}(u+v)^{\frac{4}{3}}\right\rangle^{-1}(u+v)^{-p},\\
	&\left|\int_{0}^{v}e^{\frac{4i}{3}(w-v) }\frac{d^{r}}{du^{r}}\left(\big(1+\beta(w+u)\big)^{-2}\right)dw\right|\lesssim u^{-r-1},\quad r\geq 1.
\end{align*}
Therefore we obtain the following estimate:
\begin{align*}
	|a^{(k)}_{2}(u)|\lesssim& G^{(k)}(u)+\hbar^{\frac{4}{3}}\int_{0}^{\infty}\left\langle\hbar^{\frac{2}{3}}(u+v)^{\frac{2}{3}}\right\rangle^{-2}(u+v)^{-\frac{2}{3}}\big|a_{2}^{(k)}(u+v)\big|\,dv\\
	=&G^{(k)}(u)+\hbar^{\frac{4}{3}}\int_{u}^{\infty}\left\langle\hbar w\right\rangle^{-\frac{4}{3}}w^{-\frac{2}{3}}\big|a_{2}^{(k)}(w)\big|\,dw.
\end{align*}
Here $G^{(k)}(u)$ is given by
\begin{align*}
	G^{(k)}(u)=&G^{(k)}_{1}(u)+G^{(k)}_{2}(u)\\
	=:&\sum_{p+q+r+m= k}C_{p,q,r,m}\, \hbar^{\frac{1}{3}}\int_{0}^{\infty}\frac{d^{p}}{du^{p}}\left(\tV((\hbar (u+v))^{\frac23})\right)\frac{d^{q}}{du^{q}}\left(\big(1+\beta(u+v) \big)^2\right)\\
	&\cdot\int_{0}^{v}e^{\frac{4i}{3}(w-v) }\frac{d^{r}}{du^{r}}\left(\big(1+\beta(w+u)\big)^{-2}\right)dw \left(1+\hbar  a_2(u+v)\right)\frac{d^{m}}{du^{m}}\left((u+v)^{-\frac23}\right)\,dv\\
	&+\sum_{p+q+r+s+m= k, 1\leq s<k}C_{p,q,r,s,m}\hbar^{\frac{1}{3}}\int_{0}^{\infty}\frac{d^{p}}{du^{p}}\left(\tV((\hbar (u+v))^{\frac23})\right)\frac{d^{q}}{du^{q}}\left(\big(1+\beta(u+v) \big)^2\right)\\
	&\cdot\int_{0}^{v}e^{\frac{4i}{3}(w-v) }\frac{d^{r}}{du^{r}}\left(\big(1+\beta(w+u)\big)^{-2}\right)dw\frac{d^{s}}{du^{s}}\left(1+\hbar  a_2(u+v)\right)\frac{d^{m}}{du^{m}}\left((u+v)^{-\frac23}\right)\,dv.
\end{align*}
If $u\gtrsim \hbar^{-1}$, $G^{(k)}_{1}(u)$ is bounded by
\begin{align*}
	|G^{(k)}_{1}(u)|\lesssim& \hbar^{\frac{1}{3}}\int_{u}^{\infty}\left\langle\hbar^{\frac{2}{3}}w^{\frac{2}{3}}\right\rangle^{-2}w^{-\frac{2}{3}-k}dw+\sum_{1\leq r\leq k}\hbar^{\frac{1}{3}}u^{-r-1}\int_{u}^{\infty}\left\langle\hbar^{\frac{2}{3}}w^{\frac{2}{3}}\right\rangle^{-2}w^{-\frac{2}{3}-k+r}dw\\
	\lesssim&\hbar^{-1}\int_{u}^{\infty}w^{-2-k}dw+\sum_{1\leq r\leq k}\hbar^{-1}u^{-r-1}\int_{u}^{\infty}w^{-2-k+r}dw\\
	\lesssim&\hbar^{-1}u^{-1-k}.
\end{align*}
This computation shows that, for $u\gtrsim \hbar^{-1}$, we expect to obtain an estimate for $a^{(k)}_{2}(u)$ as $\left|a^{(k)}_{2}(u)\right|\lesssim \hbar u^{-k-1}$. Therefore, by induction, we can assume that $G^{(k)}_{2}(u)$, when $u\gtrsim \hbar^{-1}$, enjoys the same estimate as for $G^{(k)}_{1}(u)$. So we finally obtain
\begin{align}\label{a2 k 1}
	\left|a^{(k)}_{2}(u)\right|\lesssim \hbar u^{-k-1},\quad u\gtrsim\hbar^{-1}.
\end{align}
When $1\leq u\ll\hbar^{-1}$, we have
\begin{align*}
	|G^{(k)}_{1}(u)|\lesssim& \hbar^{\frac{1}{3}}\int_{u}^{\hbar^{-1}}w^{-\frac{2}{3}-k}dw+\sum_{1\leq r\leq k}C_{r}\hbar^{\frac{1}{3}}u^{-r-1}\int_{u}^{\hbar^{-1}}w^{-\frac{2}{3}-k+r}dw+\hbar^{k}\\
	\lesssim &\hbar^{\frac{1}{3}}u^{\frac{1}{3}-k}+\hbar^{k}\lesssim \hbar^{\frac{1}{3}}u^{\frac{1}{3}-k}.
\end{align*}
Again, this computation shows that, for $1\leq u\ll\hbar^{-1}$, we expect to obtain an estimate for $a^{(k)}_{2}(u)$ as $\left|a^{(k)}_{2}(u)\right|\lesssim \hbar^{\frac{1}{3}}u^{\frac{1}{3}-k}$. In fact, one verifies this as before by means of an induction argument. Based on the above discussion, we claim
\begin{align*}
	|\partial_{\zeta}^{k} a_1(\zeta)|\lesssim \left\{ \begin{array}{ll} |\zeta|^{-\frac12-k+1} &  \forall -1\le \zeta\le -\hbar^{\frac23} \\
	|\zeta|^{-\frac52-k+1} &  \forall \;  \zeta\le -1 
	\end{array}\right. 
\end{align*}
This can be proved using an induction argument. The case for $k=0$ is already proved. We assume that the estimate holds for $k'<k$, and prove the estimates for $k'=k$. In view of the relation \eqref{a1 a2 1st deri}, we have
\begin{align}\label{a1 a2 higher deri}
\begin{split}
a^{(k)}_{2}(u)=&c_{k}a^{(k)}_{1}\left(-(\hbar u)^{\frac{2}{3}}\right)\left(-\hbar^{\frac{2k}{3}}\right)u^{-\frac{k}{3}}+\sum_{1\leq k'<k}c_{k'}a^{(k')}_{1}\left(-(\hbar u)^{\frac{2}{3}}\right)\left(-\hbar^{\frac{2k'}{3}}\right)u^{-\frac{k'}{3}}\cdot u^{-k+k'}
\end{split}
\end{align}
Substituting the estimates for $a^{(k')}_{1}\left(-(\hbar u)^{\frac{2}{3}}\right)$ into the second term on the right hand side above, we obtain the desired result. 
Finally we consider the regime $-\hbar^{\frac{2}{3}}\leq \zeta\leq 0$. To this end, we simply differentiate the equation \eqref{a1 1st deri} with respect to $\zeta$:
\begin{align}\label{a1 2nd deri}
	\begin{split}
	\ddot{a}_{1}(\zeta)=&\left(w_{1}^{-2}(\zeta)\right)^{\cdot}w^{2}_{1}(-\hbar^{\frac{2}{3}})\dot{a}_{1}(-\hbar^{\frac{2}{3}})+\hbar^{-1}\tV(-\zeta)(1+\hbar a_{1}(\zeta))\\
	&+\hbar^{-1}\left(w_{1}^{-2}(\zeta)\right)^{\cdot}\int_{-\hbar^{\frac{2}{3}}}^{\zeta}\tV(-s)w^{2}_{1}(s)(1+\hbar a_{1}(s))ds\\
	=&O(\hbar^{-1}).
	\end{split}
\end{align}
The point here is that differentiating $w_{1}^{-2}(\zeta)$ once gives a factor of $\hbar^{-\frac{2}{3}}$, while integration over $[-\hbar^{\frac{2}{3}},\zeta]$ gains at least another factor $\hbar^{\frac23}$. Therefore using an induction argument, we have
\begin{align}\label{a1 higher deri}
	a^{(k)}_{1}(\zeta)=O(\hbar^{\frac{1-2k}{3}}).
\end{align}
This completes the proof of the lemma.
\end{proof}

The argument for  $\zeta\ge0 $ and $\alpha\ge1$ is very similar. 

\begin{lemma}
	\label{lem:a0fine1}
	The functions $a_0(\zeta;\alpha,\hbar)$ from Lemma~\ref{lem:solve pA} satisfy the bounds
	\EQ{
		\label{eq:a0 fine}
		|\partial_{\alpha}^{\ell}a_{0}(\zeta;\alpha,\hbar)| & \leq C_{\ell}\, \alpha^{-\ell}
		\langle\zeta\rangle^{-\frac{3}{2}},\quad \zeta\geq 0, \\ 
		|\partial_\zeta^{k} \partial_\alpha^\ell a_0(\zeta;\alpha,\hbar)| & \le C_{k,\ell}\,  \alpha^{-\ell} \left\{ \begin{array}{cc} 
			\zeta^{-\frac{3}{2}-k} &1<\zeta<\infty \\
			\zeta^{\frac{1}{2}-k} & \hbar^{\frac23} <\zeta\le 1 \\
			\hbar^{\frac{1-2k}{3}}&  0 \le\zeta\le \hbar^{\frac23} 
		\end{array}\right. 
	}
	for all $\ell\ge0,k\geq 1$, $\alpha\gtrsim 1 $, and $\hbar\in (0,\frac13]$. 
\end{lemma}
\begin{proof}
	In view of~\eqref{eq:a0 system} and \eqref{eq:Ai}, we have for all $s\ge\zeta\ge \hbar^{\frac23}$, 
	\EQ{\nn 
		|K_0 (\zeta,s;\alpha,\hbar)| &\le  C\hbar^{-1} \langle s\rangle^{-2} \langle\hbar^{-\frac23}s\rangle^{-\frac12}   (1+a(\hbar^{-\frac23}s))^2  
		e^{-\frac{4}{3\hbar}s^{\frac32}}
		\Big| \int^s_\zeta \langle\hbar^{-\frac23}t\rangle^{\frac12} e^{\frac{4}{3\hbar}t^{\frac32}} (1+a(\hbar^{-\frac23}t))^{-2}\, dt \Big|  \\
		& = C \hbar^{-\frac13} \langle s\rangle^{-2} \langle\hbar^{-\frac23}s\rangle^{-\frac12}   (1+a(\hbar^{-\frac23}s))^2  e^{-\frac{4}{3\hbar}s^{\frac32}}
		\int^{\hbar^{-\frac23}s}_{\hbar^{-\frac23}\zeta} \langle t\rangle^{\frac12} e^{\frac{4}{3} t^{\frac32}} (1+a(t))^{-2}\, dt 
	}
	Next, 
	\EQ{\nn 
		& e^{-\frac{4}{3}t_{1}^{\frac32}}\int_{t_0}^{t_1} \langle t\rangle^{\frac12} e^{\frac{4}{3}t^{\frac32}} (1+a(t))^{-2}\, dt = \frac{e^{-\frac{4}{3}t_{1}^{\frac32}}}{2}\int_{t_0}^{t_1} t^{-\frac12} \langle t\rangle^{\frac12}  (1+a(t))^{-2}\, d\Big(e^{\frac{4}{3}t^{\frac32}}\Big)\\
		& = \frac{e^{-\frac{4}{3}t_{1}^{\frac32}}}{2}t^{-\frac12} \langle t\rangle^{\frac12}  (1+a(t))^{-2}e^{\frac{4}{3}t^{\frac32}} \Big|_{t_0}^{t_1} -\frac{e^{-\frac{4}{3}t_{1}^{\frac32}}}{2}\int_{t_0}^{t_1}  e^{\frac{4}{3}t^{\frac32}} \, d\big( t^{-\frac12} \langle t\rangle^{\frac12}  (1+a(t))^{-2}\big)  =O(1)
	}
	uniformly in $1\le t_0\le t_1$. Here we used that $\frac{d}{dt}\big( t^{-\frac12} \langle t\rangle^{\frac12}  (1+a(t))^{-2}\big)=O(t^{-1})$ as $t\rightarrow\infty$. 
	It follows from the preceding that
	\begin{align}\label{eq:streich1}
	|K_0(\zeta,s;\alpha,\hbar)| \le  C\hbar^{-\frac13} \langle s\rangle^{-2} \langle\hbar^{-\frac23}s\rangle^{-\frac12}\le C\langle s\rangle^{-2}|s|^{-\frac12}
	\end{align}
	for all $s\ge\zeta\ge \hbar^{\frac23}$. If  $\hbar^{\frac{2}{3}}\geq s\geq \zeta\geq 0$, then 
	\EQ{\label{eq:streich2}
		|K_{0}(\zeta,s;\alpha,\hbar)|\leq C\hbar^{-1}\int^{s}_{\zeta}\, dt\leq C\hbar^{-\frac{1}{3}} 
	} 
	If $s\geq \hbar^{\frac{2}{3}}\geq \zeta\geq 0$, we split the integral $\int^{s}_{\zeta}$ in the form 
	\begin{align*}
	\int^{s}_{\zeta}=\int^{s}_{\hbar^{\frac{2}{3}}}+\int^{\hbar^{\frac{2}{3}}}_{\zeta}.
	\end{align*}
	The contribution from the first integral is treated in the same way as the case $s\geq \zeta\geq \hbar^{\frac{2}{3}}$, 
	and the contribution from the second integral is bounded by~\eqref{eq:streich2}.
	As in the oscillatory case,  a standard Volterra iteration applied to the equation~\eqref{eq:a0 system} implies 
	\[
	|a_0(\zeta;\alpha,\hbar)|\le C\langle \zeta\rangle^{-\frac32}\qquad\forall\; \zeta\ge 0,
	\]
	uniformly in $\alpha\gtrsim 1$, and $\hbar\in (0,\frac13]$.  Derivatives in $\alpha$ are handled as in the previous lemma.

	For the derivatives in~$\zeta$, denote  $\tilde a_0(\zeta)=\tilde a_0(\zeta;\alpha,\hbar):= a_0(\hbar^{\frac23}\zeta)$. Then 
	\EQ{\nn 
		\tilde a_0(\zeta)  & := \hbar^{\frac23} \int_{\zeta}^{\infty} K_0(\hbar^{\frac23} \zeta,\hbar^{\frac23} s;\alpha,\hbar) (1+\hbar \tilde a_0(s))\, ds   \\
		&= \hbar^{\frac13} \int_{ \zeta}^{\infty}  \tV(-\hbar^{\frac23} s;\alpha,\hbar) w_0^2(\hbar^{\frac23} s;\hbar) \int_\zeta^s w_0^{-2}(\hbar^{\frac23} t;\hbar)  \, dt \; (1+\hbar \tilde a_0(s))\, ds \\
		&= \hbar^{\frac13} \int_{ \zeta}^{\infty}  \tV(-\hbar^{\frac23} s;\alpha,\hbar) \Ai(s)^2 \int_\zeta^s \Ai(t)^{-2}  \, dt \; (1+\hbar \tilde a_0(s))\, ds 
	}
	If $\zeta\ge 1$, then by \eqref{eq:Ai} we have 
	\EQ{\nn 
		\tilde a_0\left(\zeta\right) & = c \hbar^{\frac13} \int_{\zeta}^{\infty}  \tV(-\hbar^{\frac23} s;\alpha,\hbar) s^{-\frac12} e^{-\frac{4}{3}s^{\frac32}}\big(1+a(s)\big)^2   \int_{\zeta}^{s}  t^{\frac12} e^{\frac{4}{3}t^{\frac32}}\big(1+a(t)\big)^{-2} \, dt \; (1+\hbar \tilde a_0(s))\, ds 
	}
	Changing variables $s^{\frac32}=\sigma$, $t^{\frac32}=\tau$ and setting $a_3(u) = \tilde a_0(u^{\frac23})=a_0((\hbar u)^{\frac23})$, $\gamma(\sigma):=a(\sigma^{\frac23})$, and suppressing $\alpha,\hbar$ as arguments  from the notation,  we obtain for all $u\ge1$
	\EQ{\nn
		& a_3(u) = c  \hbar^{\frac13} \int_{u}^{\infty}  \tV(-(\hbar \sigma)^{\frac23})  \big(1+\gamma(\sigma) \big)^2   \int_{u}^{\sigma}  e^{\frac{4}{3}(\tau-\sigma) }\big(1+\gamma(\tau)\big)^{-2} \, d\tau \; (1+\hbar  a_3(\sigma))\, \sigma^{-\frac23}\, d\sigma \\
		&= c  \hbar^{\frac13} \int_{0}^{\infty}  \tV(-(\hbar (u+v))^{\frac23})  \big(1+\gamma(u+v) \big)^2  \!\! \int_{0}^{v}  e^{\frac{4}{3}(w-v) }\big(1+\gamma(w+u)\big)^{-2}  \, dw \, (1+\hbar  a_3(u+v))\, (u+v)^{-\frac23}\, dv 
	}
	Note that the exterior variable $u$ does not appear in the phase of the complex exponential, as~\cite[Proposition~B.1]{CDST}.  This is important as we differentiate in~$u$.  In fact,  $a_3'(u)$ fulfills the following Volterra equation  
	\EQ{\nn
		a_3'(u) &= F(u)+ 
		c_5  |\thbar|^{\frac43} \int_{0}^{\infty}  \tV(-(\hbar (u+v))^{\frac23})  \big(1+\gamma(u+v) \big)^2  \!\! \int_{0}^{v}  e^{\frac{4}{3}(w-v) }\big(1+\gamma(w+u)\big)^{-2} \,  dw \,   a_3'(u+v)\, (u+v)^{-\frac23}\, dv 
	}
	where
	\EQ{\nn 
		&F(u) = c_1  \hbar \int_{0}^{\infty}  \tV'(-(\hbar (u+v))^{\frac23})   \big(1+\gamma(u+v) \big)^2  \!\! \int_{0}^{v}  e^{\frac{4}{3}(w-v) }\big(1+\gamma(w+u)\big)^{-2} \, dw \, (1+\hbar  a_3(u+v))\, (u+v)^{-1}\, dv \\
		& + c _2\hbar^{\frac13} \int_{0}^{\infty}  \tV(-(\hbar (u+v))^{\frac23})  \big(1+\gamma(u+v) \big)\gamma'(u+v)  \!\! \int_{0}^{v}  e^{\frac{4}{3}(w-v) }\big(1+\gamma(w+u)\big)^{-2} \,  dw \, (1+\hbar  a_3(u+v))\, (u+v)^{-\frac23}\, dv \\
		& + c _3\hbar^{\frac13} \int_{0}^{\infty}  \tV(-(\hbar (u+v))^{\frac23})  \big(1+\gamma(u+v) \big)^{2}  \!\! \int_{0}^{v}  e^{\frac{4}{3}(w-v) }\big(1+\gamma(w+u)\big)^{-3} \gamma'(u+w)\, dw \, (1+\hbar  a_3(u+v))\, (u+v)^{-\frac23}\, dv \\
		& + c_4  \hbar^{\frac13} \int_{0}^{\infty}  \tV(-(\hbar (u+v))^{\frac23})  \big(1+\gamma(u+v) \big)^2  \!\! \int_{0}^{v}  e^{\frac{4}{3}(w-v) }\big(1+\gamma(w+u)\big)^{-2} \, dw \, (1+\hbar  a_3(u+v))\, (u+v)^{-\frac53}\, dv \\
	}
	By \eqref{eq:Ai} we have $|\gamma^{(k)}(\sigma)|\le C_k\, \sigma^{-1-k}$ for all $k\ge0$ and $\sigma\ge1$.  Integrating by parts yields 
	\begin{align}\label{eq:2sup}
	\begin{split}
		\sup_{v\ge0} \Big| \int_{0}^{v}  e^{\frac{4}{3}(w-v) }\big(1+\gamma(w+u)\big)^{-2} \, dw \Big| &\lesssim 1 \\ 
		\sup_{v\ge0} \Big| \int_{0}^{v}  e^{\frac{4}{3}(w-v) }\big(1+\gamma(w+u)\big)^{-3} \gamma'(u+w) \, dw \Big| & \lesssim u^{-2}
		\end{split}
	\end{align}
	for all $u\ge1$. 
	
	For higher order derivatives in $\zeta$, we again use Leibniz's rule to obtain
\begin{align*}
	a^{(k)}_{3}(u)=&\sum_{p+q+r+s+m= k}C_{p,q,r,s,m}\,\hbar^{\frac{1}{3}}\int_{0}^{\infty}\frac{d^{p}}{du^{p}}\left(\tV(-(\hbar (u+v))^{\frac23})\right)\frac{d^{q}}{du^{q}}\left(\big(1+\gamma(u+v) \big)^2\right)\\
	&\cdot\int_{0}^{v}e^{\frac{4}{3}(w-v) }\frac{d^{r}}{du^{r}}\left(\big(1+\gamma(w+u)\big)^{-2}\right)dw\frac{d^{s}}{du^{s}}\left(1+\hbar  a_3(u+v)\right)\frac{d^{m}}{du^{m}}\left((u+v)^{-\frac23}\right)dv,
\end{align*}	
and we have the estimates
\EQ{\label{eq:qrsm}
	&\left|\frac{d^{m}}{du^{m}}\left((u+v)^{-\frac23}\right)\right|\lesssim (u+v)^{-\frac{2}{3}-m},\quad m\ge0\\
	&\frac{d^{s}}{du^{s}}\left(1+\hbar  a_3(u+v) \right)=\hbar a^{(s)}_{3}(u+v),\quad s\geq 1,\\
	&\left|\frac{d^{q}}{du^{q}}\left(\big(1+\gamma(u+v) \big)^2\right)\right|\lesssim (u+v)^{-1-q},\quad q\geq 1,\\
	&\left|\int_{0}^{v}e^{\frac{4}{3}(w-v) }\frac{d^{r}}{du^{r}}\left(\big(1+\gamma(w+u)\big)^{-2}\right)dw\right|\lesssim u^{-r-1},\quad r\geq 1.
}
as well as
\[
\left|\frac{d^{p}}{du^{p}}\left(\tV(-(\hbar (u+v))^{\frac23})\right)\right|\lesssim\left\langle\hbar^{\frac{2}{3}}(u+v)^{\frac{2}{3}}\right\rangle^{-2}(u+v)^{-p},\quad p\ge0.
\]
   The proof now concludes as in the previous lemma. 
\end{proof}

Finally, we consider the scenario which truly differs from the one in \cite{CDST}.  

\begin{lemma}
	\label{lem:a0fine2}
	The functions $a_0(\zeta;\alpha,\hbar)$ from Lemma~\ref{lem:solve pA} satisfy the bounds
	\EQ{
	\label{eq:a0 fine*}
	 |\partial_{\alpha}^{\ell}a_{0}(\zeta;\alpha,\hbar)| & \leq C_{\ell}\, \alpha^{-\ell}\big[\langle \zeta\rangle^{-\frac32}+ \min(1, x(-\zeta;\alpha,\hbar)^2/\alpha^2) \big]
		,\quad \zeta\geq 0, 		
	}
	for all $\ell\ge0$, $0<\alpha\ll 1 $, and $\hbar\in (0,\frac13]$. Here $x(\tau;\alpha,\hbar)$ is the diffeomorphism from Proposition~\ref{prop:tV summary}. 
	Furthermore,  in the same parameter regime, 
	\EQ{
	\label{eq:a0 fine**}
	 |\partial_{\alpha}^{\ell} \partial^{k}_\zeta\, a_{0}(\zeta;\alpha,\hbar)| & \leq  C_{k,\ell}\, \alpha^{-\ell} \zeta^{\frac{k}{2}}
	 \big[\langle \zeta\rangle^{-\frac32}+ \min(1, x(-\zeta;\alpha,\hbar)^2/\alpha^2) \big]
		,\quad \zeta\geq1, 		
	}
and 
\EQ{
	\label{eq:a0 fine***}
	 |\partial_{\alpha}^{\ell}\partial^{k}_\zeta\, a_{0}(\zeta;\alpha,\hbar)| & \le  C_{k,\ell}\,  \alpha^{-\ell} \left\{ \begin{array}{cc} 
	 	\zeta^{\frac{1}{2}-k} & \quad \hbar^{\frac23} <\zeta\le 1 \\
	 	\hbar^{\frac{1-2k}{3}}&  \quad 0 \le\zeta\le \hbar^{\frac23} 
	 \end{array}\right. 
}
\end{lemma}
\begin{proof}
By Proposition~\ref{prop:tV summary}, for all $-\tau\gtrsim1$ 
\EQ{\label{eq:tV516}
	\tV(\tau;\alpha,\hbar) = \frac{5}{16\tau^2} - \tau \varphi(x;\alpha,\hbar)
}
with \eqref{eq:phi new} and \eqref{eq:xtau} describing $\varphi(x;\alpha,\hbar)$ and the relation between $x$ and $\tau$.  In particular, from \eqref{eq:xtau}, 
\EQ{\label{eq:tau strich}
(-\tau(x))^{\frac12} \tau'(x) =  O(x^{-1}) 
}
uniformly in the parameters.  
By the preceding lemma, cf.~\eqref{eq:streich1}, \eqref{eq:streich2}, 
\[
k_0(s;\alpha,\hbar):=	\sup_{0\le \zeta\le s} |K_0(\zeta,s;\alpha,\hbar)| \le   C(\langle s\rangle^{-2}+ s |\varphi(x;\alpha,\hbar)|)  \min(|s|^{-\frac12},\hbar^{-\frac13})
\]
We adopt the convention that $\varphi(x;\alpha,\hbar)=0$ if $0\le\zeta\lesssim 1$ and we used that $\tV(\tau;\alpha,\hbar) $ is a smooth bounded function of $|\tau|\lesssim 1$ uniformly in the parameters, see~\eqref{eq:ptV}. 
The $s^{-2}$ here  contributes the exact same amount as Lemma~\ref{lem:a0fine1}, so we only need to consider~$s\varphi$. 
From
\[
|\varphi(x;\alpha,\hbar)| 
		\lesssim \, \min\big(\hbar\, \alpha^2 x^{-2}+x^2, (\hbar+\alpha^2)x^2/\alpha^2\big)  
\]
we conclude from \eqref{eq:tau strich} writing $ds=s'(x)\, dx$, 
\EQ{\label{eq:int K0}
& \int_\zeta^\infty |K_0(\zeta,s;\alpha,\hbar)|\, ds \le \int_\zeta^\infty k_0(s;\alpha,\hbar)\, ds \\
& \lesssim \int_\zeta^\infty (\langle s\rangle^{-2}+ s |\varphi(x;\alpha,\hbar)|)   \min(|s|^{-\frac12},\hbar^{-\frac13})\, ds  \\
&\lesssim \langle \zeta\rangle^{-\frac32} +  \int^{x(-\zeta)}_0 s \min\big(\hbar\, \alpha^2 x^{-2}+x^2, (\hbar+\alpha^2)x^2/\alpha^2\big)   \min(|s|^{-\frac12},\hbar^{-\frac13})\, s^{-\frac12} x^{-1}\, dx \\
&\lesssim \langle \zeta\rangle^{-\frac32} + \int^{x(-\zeta)}_0 \min\big(\hbar\, \alpha^2 x^{-3}+x, (\hbar+\alpha^2)x/\alpha^2\big) \, dx \\
& \lesssim \langle \zeta\rangle^{-\frac32} + \min(1, x(-\zeta)^2/\alpha^2). 
}
Recall that $x(-\zeta)\lesssim \exp(-c\,\zeta^{\frac32})$, thus $\langle \zeta\rangle^{-\frac32}$ decays more slowly than the final term as $\zeta\to\infty$.  This estimate controls the first term in the Volterra iteration computing the solution of~\eqref{eq:a0 system}. The full Volterra series now provides the bound
\EQ{\nn
|a_{0}(\zeta;\alpha,\hbar)| & \le \sum_{n=1}^\infty \frac{\hbar^{n-1} }{n!} \left(  \int_\zeta^\infty k_0(s;\alpha,\hbar)\, ds\right)^n  \\
&\le  \int_\zeta^\infty k_0(s;\alpha,\hbar)\, ds \exp\left(\hbar \int_\zeta^\infty k_0(s;\alpha,\hbar)\, ds \right) \\
&\lesssim \left\langle\zeta\right\rangle^{-\frac32}+  \min(1, x(-\zeta)^2/\alpha^2)
}
as claimed, see \eqref{eq:a0 fine*} with $\ell=0$. The derivatives in $\alpha$ are estimated in the same fashion, noting that each derivative in $\alpha$ loses a factor of~$\alpha$ by Proposition~\ref{prop:tV summary}. 
For the $\zeta$ derivatives, we use the same changes of variables as in the previous lemma. I.e.,  setting $a_3(u) = \tilde a_0(u^{\frac23})=a_0((\hbar u)^{\frac23})$ we obtain 
\EQ{\nn
		& a_3'(u) =  c_1  \hbar \int_{0}^{\infty}  \tV'(-(\hbar (u+v))^{\frac23})   \big(1+\gamma(u+v) \big)^2  \!\! \int_{0}^{v}  e^{\frac{4}{3}(w-v) }\big(1+\gamma(w+u)\big)^{-2} \, dw \, (1+\hbar  a_3(u+v))\, (u+v)^{-1}\, dv \\
		& + c _2\hbar^{\frac13} \int_{0}^{\infty}  \tV(-(\hbar (u+v))^{\frac23})  \big(1+\gamma(u+v) \big)\gamma'(u+v)  \!\! \int_{0}^{v}  e^{\frac{4}{3}(w-v) }\big(1+\gamma(w+u)\big)^{-2} \,  dw \, (1+\hbar  a_3(u+v))\, (u+v)^{-\frac23}\, dv \\
		& + c _3\hbar^{\frac13} \int_{0}^{\infty}  \tV(-(\hbar (u+v))^{\frac23})  \big(1+\gamma(u+v) \big)^{2}  \!\! \int_{0}^{v}  e^{\frac{4}{3}(w-v) }\big(1+\gamma(w+u)\big)^{-3} \gamma'(u+w)\, dw \, (1+\hbar  a_3(u+v))\, (u+v)^{-\frac23}\, dv \\
		& + c_4  \hbar^{\frac13} \int_{0}^{\infty}  \tV(-(\hbar (u+v))^{\frac23})  \big(1+\gamma(u+v) \big)^2  \!\! \int_{0}^{v}  e^{\frac{4}{3}(w-v) }\big(1+\gamma(w+u)\big)^{-2} \, dw \, (1+\hbar  a_3(u+v))\, (u+v)^{-\frac53}\, dv  \\ 
		& +   c_5  \hbar^{\frac43} \int_{0}^{\infty}  \tV(-(\hbar (u+v))^{\frac23})  \big(1+\gamma(u+v) \big)^2  \!\! \int_{0}^{v}  e^{\frac{4}{3}(w-v) }\big(1+\gamma(w+u)\big)^{-2} \,  dw \,   a_3'(u+v)\, (u+v)^{-\frac23}\, dv 
	}
	As before, \eqref{eq:Ai} implies that  $|\gamma^{(k)}(\sigma)|\le C_k\, \sigma^{-1-k}$ for all $k\ge0$ and $\sigma\ge1$.  Moreover, \eqref{eq:2sup} holds	for all $u\ge1$.  Hence, for all $u\ge1$, 
	\EQ{\nn
		 |a_3'(u)| &\lesssim   \hbar \int_{u}^{\infty} | \tV'(-(\hbar w)^{\frac23}) |  w^{-1}\, dw 
		 + \hbar^{\frac13} \int_{u}^{\infty}  |\tV(-(\hbar w)^{\frac23}) |  (w^{-1}+\hbar  |   a_3'(w)|  ) w^{-\frac23}\, dw 
	}
	Changing variables $u=\hbar^{-1} t^{\frac32}$, and $w=\hbar^{-1} s^{\frac32}$, and $|a_3'(u)|=:f(t)$ yields 
		\EQ{\label{eq:f Volterra}
		 f(t)  &\lesssim   \hbar \int_{{t}}^{\infty} | \tV'(-{s}) |  s^{-1}\, ds
		 + \hbar  \int_{{t}}^{\infty}  |\tV(-{s}) |  (  s^{{-\frac{3}{2}}}+  f(s)  ) s^{{-\frac12}}\, ds 
	}
	for all $t\ge {\hbar^{\frac{2}{3}}}$.  	By Proposition~\ref{prop:tV summary},  
\[
|\tV'(-\zeta)| \lesssim \langle \zeta\rangle^{-3} +  \langle \zeta\rangle^{\frac32}  \min\big(\hbar\, \alpha^2 x^{-2}+x^2, (\hbar+\alpha^2)x^2/\alpha^2\big)   
\]
for all $\zeta\ge0$. If $t\ge1$, we rewrite \eqref{eq:f Volterra} in the form using \eqref{eq:tau strich}, 
\EQ{\label{eq:f Volterra*}
		 f(t)   
		 &\lesssim  \hbar \int_{t}^{\infty} \big( \xi^{-4}  + \xi^{\frac12} \min\big(\hbar\, \alpha^2 x^{-2}+x^2, (\hbar+\alpha^2)x^2/\alpha^2\big)  \big)  \, d\xi \\
		 &
		\qquad + \hbar  \int_{t}^{\infty}  \big( \xi^{- {\frac{5}{2}}}  + {\xi^{\frac{1}{2}}}\min\big(\hbar\, \alpha^2 x^{-2}+x^2, (\hbar+\alpha^2)x^2/\alpha^2\big) \big)  f(\xi) \, d\xi \\
		&\lesssim  f_0(t) + \hbar  \int_t^\infty L(\xi;\alpha,\hbar)  f(\xi) \, d\xi 
	}
where
\[
f_0(t):=\hbar(t^{-3} + \min(1,x(t)^2/\alpha^2)),\quad L(\xi;\alpha,\hbar):= 
 \xi^{- {\frac{5}{2}}}  + {\xi^{\frac{1}{2}}}\min\big(\hbar\, \alpha^2 x^{-2}+x^2, (\hbar+\alpha^2)x^2/\alpha^2\big). 
 \]
 By Volterra iteration, \eqref{eq:f Volterra*} implies 
 \[
 f(t) \lesssim \hbar(t^{-\frac32} + \min(1,x(t)^2/\alpha^2))\qquad t\ge1
 \]
 If $\hbar^{\frac23}\le t\le 1$, then  \eqref{eq:f Volterra} becomes
 \[
 f(t)\lesssim \hbar\int_t^1 \xi^{-2}\, d\xi + \hbar \int_t^1 f(s)\, \frac{ds}{s^{\frac{1}{2}}}\lesssim \hbar t^{-1}
 \]
 The changes of variables above amount to 
 \EQ{
 \nn
 f(t) &= |a_3'(\hbar^{-1} t^{\frac32}) | = \frac23 \hbar^{\frac23} (\hbar^{-1} t^{\frac32})^{-\frac13} |a_0'(t)| = \frac23 \hbar t^{-\frac12}  |a_0'(t)|
 }
 whence 
 \[
 |a_0'(t)| \lesssim  t^{\frac12}\big[ t^{-\frac32}+\min(1,x(t)^2/\alpha^2) \big] \qquad\forall\; t\ge1
 \]
 and 
 \[
  |a_0'(t)| \lesssim t^{-\frac12} \qquad\forall\; \hbar^{\frac23}\le t\le 1
 \]
 Finally, if $0\le t\le  \hbar^{\frac23}$, then $|a_0'(t)|\lesssim \hbar^{-\frac13}$, as in Lemma~\ref{lem:a0fine1}. 
 
 Next we turn to the higher order $\zeta$-derivatives. We again use the formula
 \EQ{\label{eq:a3higher}
  a^{(k)}_{3}(u)=&\sum_{p+q+r+s+m= k}C_{p,q,r,s,m}\, \hbar^{\frac{1}{3}}\int_{0}^{\infty}\frac{d^{p}}{du^{p}}\left(\tV(-(\hbar (u+v))^{\frac23})\right)\frac{d^{q}}{du^{q}}\left(\big(1+\gamma(u+v) \big)^2\right)\\
	&\cdot\int_{0}^{v}e^{\frac{4}{3}(w-v) }\frac{d^{r}}{du^{r}}\left(\big(1+\gamma(w+u)\big)^{-2}\right)\, dw\frac{d^{s}}{du^{s}}\left(1+\hbar  a_3(u+v)\right)\frac{d^{m}}{du^{m}}\left((u+v)^{-\frac23}\right)\, dv,
 }
as well as the estimates~\eqref{eq:qrsm}. 
The true difference in this case are the bounds on higher order derivatives of $\tV$. Let us denote:
\begin{align*}
 \tV\left(-(\hbar(u+v))^{\frac23}\right)=&\tV_{1}\left(-(\hbar(u+v))^{\frac23}\right)+\tV_{2}\left(-(\hbar(u+v))^{\frac23}\right)\\
 :=&\frac{5}{16(\hbar(u+v))^{\frac{4}{3}}}-\left(\hbar(u+v)\right)^{\frac23}\vphi(x;\alpha,\hbar).
\end{align*}
By our convention, the contribution from $\tV_{2}$ appears only for $\hbar(u+v)\geq 1$. The contribution from $\tV_{1}$ is bounded the same way as in Lemma~\ref{lem:a0fine1}. For the contribution from $\tV_{2}$, we have, in view of Proposition~\ref{prop:tV summary},
\begin{align*}
	\left|\frac{d^{p}}{du^{p}}\tV_{2}\left(-\left(\hbar(u+v)\right)^{\frac{2}{3}}\right)\right|\lesssim \sum_{p'+p''=p}\hbar^{\frac{2}{3}}\cdot\hbar^{p'} (u+v)^{\frac{2}{3}-p''}\min\{\hbar\alpha^{2}x^{-2}+x^{2},(\hbar+\alpha^{2})x^{2}\alpha^{-2}\} ,\quad\textrm{for}\quad p\geq0.
\end{align*}
Based on these estimates, the contribution from $\tV_{2}$ is bounded by, for $\hbar u\geq 1$,
\begin{align*}
 \left|a^{(k)}_{3}(u)\right|\lesssim&\hbar\sum_{0\leq p'\leq k}\int_{u}^{\infty}\hbar^{p'}\left(w^{-k+p'}+\sum_{1\leq r\leq k-p'}u^{-r-1}w^{-k+p'+r}\right)\\
 &\quad \quad \quad \quad\cdot\min\{\hbar\alpha^{2}x^{-2}+x^{2},(\hbar+\alpha^{2})x^{2}\alpha^{-2}\} \, dw\\
 +&\hbar^2\sum_{0\leq p'\leq k'< k,k'\geq1}\int_{u}^{\infty}\hbar^{p'}\left(w^{-k'+p'}+\sum_{1\leq r\leq k'- p'}u^{-r-1}w^{-k'+p'+r}\right)\\
 &\quad \quad \quad \quad\cdot\min\{\hbar\alpha^{2}x^{-2}+x^{2},(\hbar+\alpha^{2})x^{2}\alpha^{-2}\}\left|a_{3}^{(k-k')}(w)\right| \, dw\\
 +&\hbar^{2}\int_{u}^{\infty}\min\{\hbar\alpha^{2}x^{-2}+x^{2},(\hbar+\alpha^{2})x^{2}\alpha^{-2}\}\left|a^{(k)}_{3}(w)\right| \, dw\\
 =:&I+II+III.
\end{align*}
The first term here is the contribution of no derivative falling on $ a_3(u+v)$ in \eqref{eq:a3higher}, the second one is if at least one falls on $a_3(u+v)$, but no more than~$k-1$, and the third, if all $k$ derivatives fall on~$a_3(u+v)$. 
Changing variables $u=\hbar^{-1}t^{\frac{3}{2}}, w=\hbar^{-1}s^{\frac{3}{2}}$, and $|a^{(k)}_{3}(u)|=:f^{(k)}(t)$ yields
\begin{align*}
I\lesssim &\hbar\sum_{0\leq p'\leq k}\int_{t}^{\infty}\left(\hbar^{k-1}s^{-\frac{3(k-p')-1}{2}}+\hbar^{k}\sum_{1\leq r\leq k- p'}t^{-\frac{3r+3}{2}}s^{-\frac{3k-3p'-3r-1}{2}}\right)\\
&\quad \quad \quad \quad\cdot \min\{\hbar\alpha^{2}x^{-2}+x^{2},(\hbar+\alpha^{2})x^{2}\alpha^{-2}\} \, ds,\\
II\lesssim &\hbar^2\sum_{0\leq p'\leq k'<k,k'\geq 1}\int_{t}^{\infty}\left(\hbar^{k'-1}s^{-\frac{3(k'-p')-1}{2}}+\hbar^{k'}\sum_{1\leq r\leq k'-p'}t^{-\frac{3r+3}{2}}s^{-\frac{3k'-3p'-3r-1}{2}}\right)\\
&\quad \quad \quad \quad\cdot \min\{\hbar\alpha^{2}x^{-2}+x^{2},(\hbar+\alpha^{2})x^{2}\alpha^{-2}\} f^{(k-k')}(s)\, ds,\\
III\lesssim&\hbar\int_{t}^{\infty}s^{\frac{1}{2}}\min\{\hbar\alpha^{2}x^{-2}+x^{2},(\hbar+\alpha^{2})x^{2}\alpha^{-2}\} f^{(k)}(s)\, ds.
\end{align*}
Recall the bound
\begin{align*}
	\int_{t}^{\infty}s^{\frac{1}{2}}\min\{\hbar\alpha^{2}x^{-2}+x^{2},(\hbar+\alpha^{2})x^{2}\alpha^{-2}\}ds\lesssim \min\left\{1,\frac{x(t)^{2}}{\alpha^{2}}\right\},
\end{align*}
Thus,  for $t\geq 1$,
\begin{align}\label{a3 k singular 1}
	I\lesssim \hbar^{k}\int_{t}^{\infty}s^{\frac{1}{2}}\min\{\hbar\alpha^{2}x^{-2}+x^{2},(\hbar+\alpha^{2})x^{2}\alpha^{-2}\}ds\lesssim \hbar^{k}\min\left\{1,\frac{x(t)^{2}}{\alpha^{2}}\right\}.
\end{align}
For the contribution from $\tV_{1}$, we have
\begin{align*}
	\left|\frac{d^{p}}{du^{p}}\tV_{1}\left(-\left(\hbar(u+v)\right)^{\frac{2}{3}}\right)\right|\lesssim \hbar^{-\frac{4}{3}}(u+v)^{-\frac{4}{3}-p}.
\end{align*}
Therefore this contribution is bounded in the same way as in Lemma \ref{lem:a1fine} and \ref{lem:a0fine1}. 

For the case when $\hbar^{\frac23}\zeta\leq 1$ and $0\leq \zeta\leq \hbar^{\frac23}$, by adopting the convention $\vphi=0, 0\leq \zeta\leq 1$, the corresponding argument in the proof of Lemma \ref{lem:a0fine1} applies here.
\end{proof}
\subsection{Distorted Fourier transform and the spectral measures for all $n\ge2$}
The machinery developed starting with Section~\ref{sec: ngeq2 linear} up until this point  allows us to determine the spectral measures associated to all operators $\calH_n$ with $n$ positive and large. In other words, we need $\hbar=\frac{1}{n+1}$ small, and  not just $0<\hbar\le\frac13$. The reason for this lies with the correction factors $1+\hbar a_j$, $j=0,1$ which are only useful if
$\hbar |a_j|\ll 1$.  The finitely many remaining $n$ will be treated below by the same analysis from~\cite{KST} which have already applied to $n=0,\pm1$. 

\subsubsection{$n\ge N_0$ with $N_0$ large and fixed} 

We wish to find the distorted Fourier transform as in Section~\ref{sec:DFT} associated with the spectral problem $-\calH_n^+ f = \xi f$, $\xi>0$, cf.~\eqref{eigen eq}. Recalling equations~\eqref{Schlag eq},  \eqref{Schlag eq alt}, \eqref{eq:qdef}, \eqref{def tau}, \eqref{perturbed Airy}, and Lemma~\ref{lem:solve pA}, the half-line problem $-\calH_n^+ f=E^2 f= \xi f$ in $L^2(dR)$, $R>0$ is converted into the perturbed Airy equation~\eqref{perturbed Airy} on the whole line 
\[
\hbar^2 \ddot w(\tau) = -\tau w(\tau) +\hbar^2 \tilde V(\tau;\alpha,\hbar) w(\tau)
\]
with $\alpha=\hbar E=\hbar \xi^{\frac12}$. Here $\tau$ and $x=\hbar ER=\hbar\xi^{\frac12}R$ are related by~\eqref{def tau}. See~\eqref{eq:xtau} for the asymptotic relation between $x$ and~$\tau$. Throughout, $\hbar=(n+1)^{-1}$. We also freely use the results and notations of the previous section, in particular Lemma~\ref{lem:solve pA} and the quantitative control on $a_0$ and $a_1$ obtained above.

\begin{proposition}
\label{prop:DFT nlarge}
The distorted Fourier transform associated with $-\calH_n^+$, $n\gg1$, takes the following form: 
\EQ{\label{eq:hatfn}
\hat{f}(\xi)  &=   \int_0^\infty \phi_n(R,\xi) f(R)\, dR,\qquad 
f(R)  = \int_0^\infty \phi_n(R,\xi) \hat{f}(\xi) \rho_n(d\xi)
}
with   
\EQ{\label{eq:phin} 
\phi_n(R;\xi) &= \hbar^{\frac13}\alpha^{-\frac12} q^{-\frac14}( \tau)\Ai(-\hbar^{-\frac23} \tau)(1+\hbar a_0(-\tau;\alpha,\hbar)) \\
&= \hbar^{-\frac16}\xi^{-\frac14}  q^{-\frac14}( \tau)\Ai(-\hbar^{-\frac23} \tau)(1+\hbar a_0(-\tau;\alpha,\hbar)) \\
 \tau & =\tau(x;\alpha,\hbar),  \quad x= \alpha R\le x_t,\quad \alpha=\hbar\xi^{\frac12}
}
Here $\tau=\tau(x;\alpha,\hbar)$ as in \eqref{def tau}, $q$ is defined in~\eqref{eq:qdef},  and $x_t=x_t(\alpha;\hbar)$ the turning point as in Lemma~\ref{lem: monotonicity of root in alpha}.  
One has $\phi_n(R;\xi)\sim \hbar^{n-\frac12}\xi^{\frac{n-1}{2}}R^{n-\frac12}$ as $ R\xi^{\frac12}\to0+$. To the right of the turning point we have the representation  
\begin{align}
\phi_n(R;\xi) &=  -c_1 \hbar^{-\frac16}\xi^{-\frac14} q^{-\frac14}( \tau) \Re\Big(\big(1+\hbar\Xi(\xi;\hbar)\big) \big(\Ai(-\hbar^{-\frac23}\tau)+i\Bi(-\hbar^{-\frac23}\tau)\big)(1+\hbar a_1(-\tau;\alpha,\hbar))\Big)  \nn \\
& \sim -c_2\xi^{-\frac14} \Re\Big(\big(1+\hbar\Xi(\xi;\hbar)\big) e^{i\frac{\pi}{4}}e^{i\xi^{\frac12} R}\Big) \text{\ \ as\ \ } R\xi^{\frac12}\to\infty,  \label{eq:phin Rlarge}
\end{align}
where $c_1, c_2>0$ are absolute constants and $|\partial_\xi^k \Xi(\xi;\hbar)|\le C_k\, \xi^{-k}$ for all $k\ge0$ uniformly in $\hbar$. The coefficients $a_0, a_1$ satisfy the bounds in Lemma~\ref{lem:a1fine}--Lemma~\ref{lem:a0fine2}. 
The spectral measure~$\rho_n$ is 
purely absolutely continuous with density satisfying 
 \EQ{\label{eq:rho n}
 \frac12 & \le \frac{d\rho_n(\xi)}{d\xi} \le 2,\qquad 
 \Big|\frac{d^\ell}{d\xi^\ell} \frac{d\rho_n(\xi)}{d\xi}\Big| \le C_\ell\,\xi^{-\ell},\qquad\forall\; \xi>0,\;\ell\ge0 
 }
 uniformly in $n$. For the higher order derivatives of $\phi_{n}(R,\xi)$ we have
 \begin{align*}
 	(R\partial_{R})^{k}\phi_{n}(R,\xi)\sim \hbar^{n-\frac12}\xi^{\frac{n-1}{2}}R^{n-\frac12},\quad (\xi\partial_{\xi})^{k}\phi_{n}(R,\xi)\sim \hbar^{n-\frac12}\xi^{\frac{n-1}{2}}R^{n-\frac12},\quad \textrm{as}\quad \hbar R\xi^{\frac12}\rightarrow0_{+},
 \end{align*}
 and 
 \begin{align*}
 	(R\partial_{R})^{k}\left(e^{-iR\xi^{\frac12}}\phi_{n}(R,\xi)\right)\sim c_{3} \xi^{-\frac14},\quad (\xi\partial_{\xi})^{k}\left(e^{-iR\xi^{\frac12}}\phi_{n}(R,\xi)\right)\sim c_{4} \xi^{-\frac14}\quad \textrm{as}\quad  R\xi^{\frac12}\rightarrow\infty.
 \end{align*}
\end{proposition} 
\begin{proof}
We fix $n\ge N_0$ large (the latter will be determined). By Section~\ref{sec:DFT}, specifically~\eqref{eq:gen m},
\[
m(z;\hbar) =\frac{W(\theta(\cdot,z),\psi(\cdot,z))}{W(\psi(\cdot,z),\phi(\cdot,z))}
\]
where $\psi(\cdot,z)$ is the exponentially decaying solution of 
\begin{align}\label{eq defining WT function}
-\calH_{n}^+ f =zf \quad \textrm{for} \quad\Im z>0,\quad \textrm{and}\quad |W(\psi,\overline{\psi})|\simeq 1.
\end{align}
Taking the limit $\xi=\xi+i0$ for $\xi>0$ we obtain
\EQ{\label{eq:m n}
m(\xi+i0;\hbar) = \frac{W(\theta(\cdot,\xi),\psi_+(\cdot,\xi))}{W(\psi_+(\cdot,\xi),\phi(\cdot,\xi))}
}
as density of the spectral measure, i.e., $\rho_n(d\xi)=m(\xi+i0)\,d\xi$. 
Let 
\EQ{\label{eq:w0w}
w(-\tau;\alpha,\hbar) = w_0(-\tau;\hbar)(1+\hbar a_0(-\tau;\alpha,\hbar))=\Ai(-\hbar^{-\frac23} \tau)(1+\hbar a_0(-\tau;\alpha,\hbar)) 
}
 be the solution from Lemma~\ref{lem:solve pA} (recall the change of sign $\zeta=-\tau$ which occurred at the beginning of Section~\ref{sec:fund sys}).  As $\tau\to-\infty$, which corresponds to $R\xi^{\frac12}\to0+$, this solution satisfies
 \EQ{\label{eq:wasymp}
 w(-\tau;\alpha,\hbar) \sim c \hbar^{\frac16} (-\tau)^{-\frac14} e^{-\frac{2}{3\hbar}(-\tau)^{\frac32}} \sim c'  \hbar^{\frac16} \left(-\log \left(\hbar R\xi^{\frac12}\right)\right)^{-\frac16} \left(\hbar R\xi^{\frac12}\right)^{\frac{1-2\hbar}{\hbar}}
 }
with some absolute constants~$c,c'$, cf.~\eqref{eq:xtau}. We note that $w(-\tau;\alpha,\hbar) >0$ for $\tau\le0$. 
Take $\hbar$ small enough so that $\hbar |a_0(-\tau;\alpha,\hbar)|\le\frac12$ on $(-\infty,0]$. Note that by the results of the previous section this then holds uniformly in $\tau\le0 $ and $ \alpha>0$. Define another solution of the perturbed Airy equation~\eqref{perturbed Airy}
\EQ{\label{eq:red w}
\tilde w(\zeta;\alpha,\hbar) := w(\zeta;\alpha,\hbar)\int_0^{\zeta} w(u;\alpha,\hbar)^{-2}\, du,\qquad \zeta\ge 0
}
Their Wronskian is 
$
W(\tilde w, w)=-1
$
as can be seen by evaluating it at $\tau=0$.   We pass between solutions of $-\calH_n f =\xi f$ on the one hand, and \eqref{perturbed Airy} viz.
\[
\hbar^2 \ddot w(\tau) = -\tau w(\tau) +\hbar^2 \tilde V(\tau;\alpha,\hbar) w(\tau)
\]
on the other hand, by means of the relation $f(R)=\tilde f(x) = q^{-\frac14} w(\tau)$,  see \eqref{Schlag eq}, \eqref{Schlag eq alt}, \eqref{perturbed Airy}, \eqref{eq:qdef}. In particular, $d\tau =\sqrt{q}\, dx$, $q=-Q_0/\tau$, and $x=\alpha R$. Therefore
\begin{align*}
	\frac{d}{dR} f(R)&=\alpha\widetilde{f}'(x)=\alpha\left(q^{-\frac{1}{4}}w\right)^{\prime}=\alpha\partial_{\tau}\left(q^{-\frac{1}{4}}w\right)\cdot \frac{d\tau}{dx}\\
	&=\alpha q^{\frac{1}{2}}\left(q^{-\frac{1}{4}}\dot{w}-\frac{1}{4}q^{-\frac{5}{4}}\dot{q}w\right).
\end{align*}
and $W(f_1,f_2)=\alpha W(w_1,w_2)$ for any pairs of related solutions of \eqref{Schlag eq}, respectively, \eqref{perturbed Airy}. 

We claim that, with $w,\tilde w$ as in \eqref{eq:w0w}, \eqref{eq:red w}
\EQ{\label{eq:phi theta n}
\phi(R;\xi) & =\phi_n(R,\xi) = \alpha^{-\frac12}\hbar^{\frac13} q^{-\frac14} w(-\tau;\alpha,\hbar) \\
\theta(R;\xi) & =\theta_n(R,\xi) = \alpha^{-\frac12}\hbar^{-\frac13} q^{-\frac14} \tilde w(-\tau;\alpha,\hbar)
}
are admissible choices for the pair $\phi,\theta$ in Section~\ref{sec:DFT}.   In fact, they solve $-\calH_n f =\xi f$, $\phi$ is $L^2$ near $R=0$,  they are real-valued for $\xi>0$, and their Wronskian is 
\[
W(\theta, \phi) = -W(\tilde w, w)=1
\]
From \eqref{def Q0} and \eqref{eq:wasymp}, up to multiplicative constants, 
\[
\phi_n(R,\xi)\sim (\hbar\xi^{\frac12})^{-\frac12} \hbar^{\frac13} (-\tau/Q_0)^{\frac14} \hbar^{\frac16} (-\tau)^{-\frac14} e^{-\frac{2}{3\hbar}(-\tau)^{\frac32}} \sim \hbar^{n-\frac{1}{2}}\xi^{\frac{n-1}{2}}R^{n-\frac12}
\]
as $ R\xi^{\frac12}\to0+$.  The Weyl-Titchmarsh solution $\psi_+(R,\xi;\hbar)$ defined in \eqref{eq defining WT function} is given by 
\begin{align}\label{eq:psi+ def}
\psi_+(R,\xi;\hbar) = \hbar^{\frac13} \alpha^{-\frac12} q^{-\frac14}\overline{w_+(-\tau;\alpha,\hbar)}
\end{align}
where with $\zeta=-\tau\le0$
\[
w_+(\zeta;\alpha,\hbar):= w_1(\zeta;\hbar)(1+\hbar a_1(\zeta;\alpha,\hbar)), \quad w_1(\zeta;\hbar):=\Ai(\hbar^{-\frac23} \zeta)+i \Bi(\hbar^{-\frac23} \zeta)
\]
is the oscillatory solution from Lemma~\ref{lem:solve pA}. Here we select $\hbar$ so small that $\hbar |a_1(\zeta;\alpha,\hbar)|\le \frac12$, uniformly in $\zeta\le0$ and $\alpha>0$. 
It follows from the asymptotic behavior  
\[
\Ai(-\tau)+i\Bi(-\tau)\sim \pi^{-\frac12}\tau^{-\frac14}e^{i\left(-\frac23\tau^{\frac32}+\frac{\pi}{4}\right)}\qquad \tau\to\infty
\]
and Lemma \ref{lem: Lemma 3.4 CDST} that 
 \[
 \psi_+(R;\xi,\hbar)\sim \pi^{-\frac12} e^{i\theta}\xi^{-\frac14} e^{i\xi^{\frac12} R} \qquad R\to\infty.
 \]
 for some $\theta=\theta(\alpha,\hbar)$.  In particular, evaluating at $R=\infty$ we find the Wronskian relation 
 \EQ{\label{eq:Wpsi n}
 W(\psi_+(\cdot;\xi,\hbar),\overline{\psi_+(\cdot;\xi,\hbar)}) = -2i\pi^{-1} \qquad \forall\; \xi>0
 }
 with a positive absolute constant $C_+$. 
Returning to the generalized $m$ function of~\eqref{eq:m n}, we conclude that (suppressing $\alpha,\hbar$ from $a_j$ for simplicity)
\begin{align*}
\begin{split}
m(\xi+i0;\hbar) &= \hbar^{-\frac23} \frac{w_+(0;\alpha,\hbar)}{w(0;\alpha,\hbar)(-w_+(0;\alpha,\hbar)w'(0;\alpha,\hbar)+w_+'(0;\alpha,\hbar)w(0;\alpha,\hbar))} \\
& = \frac{i(1+i\Bi(0)/\Ai(0))(1+\hbar a_0(0))^{-2}}{ W(\Ai,\Bi)+i\hbar^{\frac53}\Ai(0)(\Ai(0)+i\Bi(0))(a_0'(0)(1+\hbar a_0(0))^{-1}-a_1'(0)(1+\hbar a_1(0))^{-1}) } \\
& = \frac{\pi (-\Bi(0)/\Ai(0)+i)(1+\hbar a_0(0))^{-2}}{ 1+i\pi\hbar^{\frac53}\Ai(0)(\Ai(0)+i\Bi(0))(a_0'(0)(1+\hbar a_0(0))^{-1}-a_1'(0)(1+\hbar a_1(0))^{-1}) }
\end{split}
\end{align*}
By the previous section, $a_j(0;\alpha,\hbar)=O(1)$, $a_j'(0;\alpha,\hbar)=O(\hbar^{-\frac13})$ uniformly in $\xi>0$ and $0<\hbar\le\frac13$. 
Thus, 
\begin{align}\label{m n xi}
\begin{split}
\frac{1}{\pi} m(\xi+i0;\hbar) &= (i-\Bi(0)/\Ai(0)) (1+ i \hbar^{\frac43} \Omega(\alpha;\hbar))(1+\hbar a_0(0;\alpha,\hbar))^{-2} \\
&= \big[ i-\Bi(0)/\Ai(0)  -  \hbar^{\frac43} \Omega(\alpha;\hbar)  (1+i \Bi(0)/\Ai(0))  \big] (1+\hbar a_0(0;\alpha,\hbar))^{-2}   \\
|\partial_\alpha^\ell \Omega(\alpha;\hbar)| &\le C_\ell \,\alpha^{-\ell},\qquad   |\partial_\xi^\ell \Omega(\alpha;\hbar)| \le C_\ell \,\xi^{-\ell} \qquad\forall\ell\ge0
\end{split}
\end{align}
and with $\alpha=\hbar\xi^{\frac12}$
\EQ{\nn 
\frac{d\rho_n(\xi)}{d\xi} = \frac{1}{\pi} \Im m(\xi+i0;\hbar) = \big[ 1   -  \hbar^{\frac43} \Im\big( \Omega(\alpha;\hbar)  (1+i \Bi(0)/\Ai(0)) \big)  \big] (1+\hbar a_0(0;\alpha,\hbar))^{-2}  
}
 We conclude that  \eqref{eq:rho n} holds as claimed. From~\eqref{eq:phiapsi}, \eqref{eq:a rep rho},  and~\eqref{eq:Wpsi n} we conclude that
 \EQ{\label{eq:a size}
 \phi_n(R;\xi) &= 2\Re(a(\xi;\hbar)\psi_+(R;\xi,\hbar)) \\
|a(\xi;\hbar)|^{-2} &=  4 \frac{d\rho_n(\xi)}{d\xi} \simeq 1
 }
uniformly in $0<\hbar\ll 1$, $\alpha>0$ and $\xi>0$. In fact, 
\begin{align}\label{a xi n}
\begin{split}
a(\xi;\hbar) &= \frac{W(\phi,\overline{\psi_+})}{W(\psi_+,\overline{\psi_+})} = -iC_+^{-1} \hbar^{\frac23} \overline{W(w,w_+)} \\
& =  -iC_+^{-1} \hbar^{\frac23}(w(0)\overline{w_+'(0)}- w'(0)\overline{w_+(0)})  \\
& = -C_+^{-1} (1+\hbar a_0(0;\alpha,\hbar))(1+\hbar a_1(0;\alpha,\hbar)) W(\Ai,\Bi) + \hbar^{\frac43} \Xi_1(\xi;\hbar) \\
& = -(\pi C_+)^{-1} + \hbar \Xi_2(\xi;\hbar) 
\end{split}
\end{align}
where $|\partial_\xi^k \Xi_2(\xi;\hbar)|\le C_k\, \xi^{-k}$ for all $k\ge0$ uniformly in $\hbar$. This implies~\eqref{eq:phin Rlarge}. For the bounds on derivatives of $\phi_{n}(R,\xi)$, we observe that the behavior of the principal part follows by a direct calculation, and we only need to look at the contribution from $a_{0}(-\tau;\alpha,\hbar)$ and $a_{1}(\tau;\alpha,\hbar)$.  When $\hbar R\xi^{\frac12}\rightarrow0_{+}$, we have $\zeta=-\tau\rightarrow\infty$. Therefore we have, in view of  Proposition \ref{prop:tV summary},
\begin{align*}
	\frac{\partial a_{0}}{\partial R}=\frac{\partial a_{0}}{\partial \zeta}\cdot\frac{\partial \zeta}{\partial x}\cdot\frac{\partial x}{\partial R}\simeq\zeta^{-\frac12}R^{-1}\frac{\partial a_{0}}{\partial \zeta}. 
\end{align*}
For $\xi$-derivative, we have
\begin{align*}
	\frac{\partial a_{0}}{\partial \xi}=&\frac{\partial a_{0}}{\partial \zeta}\cdot\frac{\partial \zeta}{\partial x}\cdot\frac{\partial x}{\partial \xi}+\frac{\partial a_{0}}{\partial\alpha}\cdot \frac{\partial \alpha}{\partial \xi}
	\simeq\zeta^{-\frac12}\xi^{-1}\cdot \frac{\partial a_{0}}{\partial\zeta}+\hbar\xi^{-\frac12}\cdot \frac{\partial a_{0}}{\partial \alpha}.
\end{align*}
Then Lemma \ref{lem:a0fine1} and Lemma \ref{lem:a0fine2} imply
\begin{align*}
\left|R\frac{\partial a_{0}}{\partial R}\right|, \quad \left|\xi\frac{\partial a_{0}}{\partial\xi}\right|\lesssim \Gamma_{0},
\end{align*}
where $\Gamma_{0}$ denotes the bound on $|a_{0}|$. Therefore the bounds on $R\partial_{R}\phi_{n}(R,\xi)$ and $\xi\partial_{\xi}\phi_{n}(R,\xi)$ follow when $\hbar R\xi^{\frac12}\rightarrow0_{+}$. Now we turn to the case when $\hbar R\xi^{\frac12}\rightarrow\infty$. In this case $\zeta=-\tau\rightarrow-\infty$, and we have, in view of Lemma \ref{lem: Lemma 3.4 CDST}, 
\begin{align*}
	\frac{\partial a_{1}}{\partial R}=\frac{\partial a_{1}}{\partial\zeta}\cdot \frac{\partial \zeta}{\partial x}\cdot \frac{\partial x}{\partial R}\simeq -x^{-\frac13}\cdot \hbar\xi^{\frac12}\cdot \frac{\partial a_{1}}{\partial\zeta}.
\end{align*}
For $\xi$-derivative, we have
\begin{align*}
	\frac{\partial a_{1}}{\partial \xi}=\frac{\partial a_{1}}{\partial\alpha}\cdot\frac{\partial\alpha}{\partial\xi}+\frac{\partial a_{1}}{\partial\zeta}\cdot\frac{\partial\zeta}{\partial x}\cdot\frac{\partial x}{\partial\alpha}\simeq \hbar\xi^{-\frac12}\frac{\partial a_{1}}{\partial\alpha}-x^{-\frac13}\cdot \hbar\xi^{\frac12}\cdot\frac{\partial a_{1}}{\partial\zeta}.
\end{align*}
Then Lemma \ref{lem:a1fine} gives
\begin{align*}
	\left|R\frac{\partial a_{1}}{\partial R}\right|,\quad \left|\xi\frac{\partial a_{1}}{\partial \xi}\right|\lesssim \Gamma_{1},
\end{align*}
where $\Gamma_{1}$ denotes the bound on $|a_{1}|$. Therefore the bounds on $R\partial_{R}\phi_{n}(R,\xi)$ and $\xi\partial_{\xi}\phi_{n}(R,\xi)$ follow when $\hbar R\xi^{\frac12}\rightarrow\infty$. The estimates on higher order derivatives follow in a similar way.
\end{proof} 

\begin{remark}
	\label{rem:hbar32}
	The representations \eqref{eq:phin} and~\eqref{eq:phin Rlarge}, respectively,  
	both extend $\lesssim \hbar^{\frac32}$ across the turning point. 
\end{remark}

\begin{corollary}\label{cor:SM behavior large n}
For $\alpha\geq \alpha_{*}$ where $\alpha_{*}>0$ is given in Proposition \ref{prop:tV summary} and $\xi$ sufficiently large, we have the expansion
\begin{align*}
    \frac{d\rho_{n}(\xi)}{d\xi}=C_{0}+C_{1}\xi^{-1}+o(\xi^{-1})
\end{align*}
where $n$ is sufficiently large and the constants $C_{0}, C_{1}$ are independent of $n$. 
\end{corollary}
\begin{proof}
    According to Proposition \ref{prop:tV summary}, when $\alpha\geq \alpha_{*}$ we have 
    \begin{align*}
        \left|\partial_{\alpha}^{\ell}\widetilde{V}(\tau;\alpha,\hbar)\right|\leq C_{\ell}\langle\alpha\rangle^{-\ell-1},\quad \forall \tau\in \bbR.
    \end{align*}
    Then following the proof of Lemmas \ref{lem:solve pA}-\ref{lem:a0fine1} we have, for $\alpha\geq \alpha_{*}$,
    \begin{align*}
        \left|\partial_{\alpha}^{\ell}a_{j}(0;\alpha,\hbar)\right|\leq C_{\ell}\alpha^{-\ell-1},\quad \forall \ell\geq 1,\quad j=0,1.
    \end{align*}
    Therefore for $\Omega(\alpha,\hbar)$ in \eqref{m n xi}, we also have
    \begin{align*}
        \left|\partial_{\alpha}^{\ell}\Omega(\alpha;\hbar)\right|\leq C_{\ell}\alpha^{-1-\ell},\quad \left|\partial_{\xi}^{\ell}\Omega(\alpha;\hbar)\right|\leq C_{\ell}\xi^{-1-\ell},\quad \forall \ell\geq 1.
    \end{align*}
    This in turn implies that $\Xi_{2}(\xi;\hbar)$ in \eqref{a xi n} satisfies
    \begin{align}\label{a xi n remainder}
        \left|\partial_{\xi}^{k}\Xi_{2}(\xi;\hbar)\right|\leq C_{k}\xi^{-1-k},\quad \forall k\geq 1.
    \end{align}
    By the relation between $\frac{d\rho_{n}(\xi)}{d\xi}$ and $a(\xi;\hbar)$, we have
    \begin{align}\label{rho n xi remainder}
        \frac{d\rho_{n}(\xi)}{d\xi}=C_{0}+\widetilde{\Xi}_{2}(\xi;\hbar)
    \end{align}
    where $\widetilde{\Xi}_{2}(\xi;\hbar)$ satisfies the same estimates as $\Xi_{2}(\xi;\hbar)$. Then in order to obtain the desired result in the corollary, since $\xi$ is large, we can think of $\widetilde{\Xi}(\xi;\hbar)$ as a function of $\xi^{-1}$ and consider its Taylor expansion at ``$\xi^{-1}=0$" up to order $\xi^{-1}$.
\end{proof}

\subsubsection{The case of finitely many $2\le n\le N_0$, $N_0$ large and fixed}
Next we analyze the behavior of the spectral measure for $\hbar\simeq1$. We directly work with the operator appearing in \eqref{Schlag eq}:
\begin{align}\label{calH n operator}
 -\calH^{+}_{n}=-\partial_{R}^{2}+(n+1)^{2}V(R),\quad 2\leq n\leq N_{0}.
\end{align}
\begin{lemma}\label{lem:Fourier basis, large hbar} 
There exists a fundamental system $\phi_0, \theta_0$ of solutions for $\mathcal{H}_nf = 0$ with the asymptotic behavior 
\[
\phi_0(R)\sim R^{n-\frac{1}{2}},\quad\theta_0(R)\sim (2n-2)^{-1}R^{-n+\frac32},\quad\textrm{as}\quad R\rightarrow0+,
\]
and such that we have 
\begin{align}\label{Large R}
\phi_0(R)\sim c_{1,n}R^{n+\frac32},\quad\theta_0(R)\sim c_{2,n}R^{-n-\frac12},\quad\textrm{as}\quad R\rightarrow\infty
\end{align}
for some positive constants $c_{i,n},i=1,2$.
\end{lemma}
\begin{proof} It suffices to construct $\phi_0$ and then define $\theta_0(R)$ by means of $W(\theta_{0},\phi_{0})=1$.  We start by constructing a fundamental system $f_0, f_1$ near $R = +\infty$. Make the ansatz  $f_0 (R)= R^{n+\frac32} + \epsilon(R)$ where
	\begin{align}\label{f0 eq}
	\left(-\partial_R^2 + \frac{(n+1)^2 - \frac14}{R^2}\right)f_0 = \left(\frac{4n}{R^2(1+R^2)} + \frac{8}{(R^2+1)^2}\right)f_0. 
	\end{align}
	We solve for $\epsilon$ using the two-sided Green function  constructed from 
	 fundamental system $\{R^{n+\frac32},\,R^{-n-\frac12}\}$ for the operator on the left. Thus,   $\epsilon$ solves the Fredholm integral equation 
	\begin{align}\label{Fredholm inte eq}
	\begin{split}
	\epsilon(R) =& \frac{1}{2n+2} R^{n+\frac32}\int_R^\infty s^{-n-\frac12}\left(\frac{4n}{s^2(1+s^2)} + \frac{8}{(s^2+1)^2}\right)\left(s^{n+\frac32} + \epsilon(s)\right)\,ds\\
	& +  \frac{1}{2n+2} R^{-n-\frac12}\int_{R_{0}}^{R} s^{n+\frac32}\left(\frac{4n}{s^2(1+s^2)} + \frac{8}{(s^2+1)^2}\right)\left(s^{n+\frac32} + \epsilon(s)\right)\,ds
	\end{split}
	\end{align}
	where $R\ge R_{0}\ge 1$. Here $R_0$ is large enough to guarantee smallness of the integral operator in a suitable norm.  In fact,  \eqref{Fredholm inte eq} can be rewritten as follows:
\begin{align}\label{Fredholm abstract 1}
	\epsilon(R)=g_{n}(R)+\left(T_{n}\epsilon\right)(R),
\end{align}	
where
\begin{align}\label{Fredholm abstract 2}
	\begin{split}
g_{n}(R):=& \frac{1}{2n+2} R^{n+\frac32}\int_R^\infty s^{-n-\frac12}\left(\frac{4n}{s^2(1+s^2)} + \frac{8}{(s^2+1)^2}\right)s^{n+\frac32}\,ds\\
& +  \frac{1}{2n+2} R^{-n-\frac12}\int_{R_{0}}^{R} s^{n+\frac32}\left(\frac{4n}{s^2(1+s^2)} + \frac{8}{(s^2+1)^2}\right)s^{n+\frac32} \,ds,\\
\left(T_{n}\epsilon\right)(R):=& \frac{1}{2n+2} R^{n+\frac32}\int_R^\infty s^{-n-\frac12}\left(\frac{4n}{s^2(1+s^2)} + \frac{8}{(s^2+1)^2}\right)\epsilon(s)\,ds\\
& +  \frac{1}{2n+2} R^{-n-\frac12}\int_{R_{0}}^{R} s^{n+\frac32}\left(\frac{4n}{s^2(1+s^2)} + \frac{8}{(s^2+1)^2}\right)\epsilon(s) \,ds.
	\end{split}
\end{align}
We introduce a Banach space $X$ such that $f\in X$ if $R^{-n+\frac12}f(R)\in L^{\infty}([R_{0},\infty))$, and
\begin{align}\label{weighted Linfty}
	\|f\|_{X}:=\|R^{-n+\frac12}f(R)\|_{L^{\infty}([R_{0},\infty))}.
\end{align}
Then $\|g_n\|_X\le C$ uniformly in $n$. 
Writing $(T_n f)(R)=\int_{R_0}^\infty K_n(R,s)f(s)\, ds$ one has with some absolute constant $C$
\begin{align}\label{Fkernel}
	0&\le  R^{-n+\frac12}K_n(R,s)s^{n-\frac12} \le C(R^{-2n}s^{2n-3}\one_{[R_0<s<R]}+ R^2s^{-5}\one_{[s>R]})
\end{align}
This immediately implies that $\| T_n f\|_X\le CR_0^{-2}\|f\|_X$. Choosing  $R_0$ large enough $T_n$ is a contraction and there is a unique solution $\epsilon\in X$ to \eqref{Fredholm abstract 1}. Therefore we obtain a unique solution $f_{0}$ to \eqref{f0 eq} defined for all $R>0$ such that
\begin{align*}
	f_{0}(R)=R^{n+\frac32}+\epsilon(R)=R^{n+\frac32}+O(R^{n-\frac12}),\quad \textrm{for}\quad R\rightarrow\infty.
\end{align*}
We can then define $f_1$ for large $R$ by 
	\[
	f_1(R) =f_0(R)\cdot \int_R^\infty f_0^{-2}(s)\,ds,
	\]
	which behaves like $f_1(R)\sim c_{n} R^{-n-\frac12}$ for $R\rightarrow\infty$ and some positive constant $c_{n}$. This gives the asymptotic description \eqref{Large R}. We next construct the solution  $\phi_0(R)$ of 
	\begin{align*}
	\left(-\partial_R^2 + \frac{(n-1)^2 - \frac14}{R^2}\right)\phi_0 = \left(-\frac{4n}{(1+R^2)} + \frac{8}{(R^2+1)^2}\right)\phi_0,
	\end{align*}
	near the origin $R = 0$ by means of  the  ansatz $\phi_0(R) = R^{n-\frac12}+\gamma(R)$. Here $\gamma$  is a solution of the Volterra equation 
	\begin{align*}
	\gamma(R) =& -\frac{1}{2n-2} R^{n-\frac12}\int_0^R s^{-n+\frac32}\left(-\frac{4n}{(1+s^2)} + \frac{8}{(s^2+1)^2}\right)\left(s^{n-\frac12} + \gamma\right)\,ds\\
	& +  \frac{1}{2n-2} R^{-n+\frac32}\int_0^R s^{n-\frac12}\left(-\frac{4n}{(1+s^2)} + \frac{8}{(s^2+1)^2}\right)\left(s^{n-\frac12} + \gamma\right)\,ds\\
	 = &\gamma_0(R) + \int_0^R \tilde K_n(R,s) \gamma(s)\, ds
	\end{align*}
	One checks that $0\le \gamma_0(R)\le CR^{n+\frac32}$ and $$0\le R^{-n-\frac32}\tilde K_n(R,s)s^{n+\frac32}\le Cs$$ for all $0<s<R$. By a standard Volterra iteration,  $0\le \gamma(R)\lesssim R^{n+\frac32}$ for $R\ll 1$. 
	We claim that this $\phi_0$, when continued up to $R = +\infty$ has the desired asymptotics and is positive everywhere. For the latter assertion use that from the equation we have that 
	\[
	\partial_R^2\phi_0>0
	\]
	as long as $\phi_0(R)>0$. In fact, $\phi_{0},\partial_{R}\phi_{0}$ and $\partial_{R}^{2}\phi_{0}$ are all positive initially. Suppose that $\phi_{0}(R)$ becomes zero at $R=R_{0}>0$ for the first time. Then there is a $R_{1}\in (0,R_{0})$ such that $\partial_{R}\phi_{0}(R_{1})=0$. This means that $\phi_{0}$ loses convexity even before $R=R_{1}$. Therefore there is a $R_{2}\in (0,R_{1})$ such that $\partial_{R}^{2}\phi_{0}(R_{2})=0$. In view of the equation satisfied by $\phi_{0}$, we have $\phi(R_{2})=0$, which contradicts the fact that $R_{0}$ is the first point where $\phi_{0}$ vanishes. Therefore $\phi_{0}(R)>0$ for all $R>0$.

Since $\phi_0$ is a linear combination of $f_0, f_1$ for large $R$, we also immediately get the large $R$ asymptotics, due to the convexity just observed. 
\end{proof}
Next we construct a solution to $\left(-\calH^{+}_{n}-\xi\right)f=0$ for all $R,\xi>0$.
\begin{lemma}\label{lem: phin small R2xi small}
	There exists a smooth function $\phi(R, \xi)$ on $R,\xi>0$ satisfying $-\calH_{n}^{+} \phi(R, \xi)=\xi \phi(R, \xi)$ of the form 
	\[
	\phi(R, \xi) =  \phi_0(R)\big[1+\sum_{j\geq 1}\phi_j(R)(R^2\xi)^j\big],
	\]
	as an absolutely convergent series where $\phi_j$ are smooth functions of $R>0$ satisfying the bounds   $\big|(R\partial_R)^m \phi_j(R)\big|\leq \frac{C^j_{m}}{j!}$ for all $m\ge0$, $R>0$. 
\end{lemma}
\begin{proof} We find $\phi(R,\xi)$ by solving the Volterra equation for all $R>0$ (recall that $1 = W(\theta_0, \phi_0)$)
	\begin{align*}
	\phi(R, \xi) = \phi_0(R)-\xi\phi_0(R)\int_0^R\theta_0(s)\phi(s, \xi)\,ds
	 + \xi\theta_0(R)\int_0^R\phi_0(s)\phi(s, \xi)\,ds.
	\end{align*}
	Inserting the preceding ansatz, we obtain the identity 
	\begin{align*}
	\sum_{j\geq 1}\phi_j(R)(R^2\xi)^j = -\xi\int_0^R\theta_0(s) \phi_0(s)[1+\sum_{j\geq 1}\phi_j(s)(s^2\xi)^j]\,ds + \xi\frac{\theta_0(R)}{\phi_0(R)}\int_0^R\phi_0^2(s)[1+\sum_{j\geq 1}\phi_j(s)(s^2\xi)^j]\,ds
	\end{align*}
	This then defines the $\phi_j(R)$ inductively. To begin with, we have 
	\[
	\phi_1(R) =  - R^{-2}\int_0^R\theta_0(s) \phi_0(s)\,ds +  R^{-2}\frac{\theta_0(R)}{\phi_0(R)}\int_0^R\phi_0^2(s)\,ds.
	\]
	and the other coefficients are then determined by the following recursive formula:
	\begin{align*}
		\phi_{j}(R)=&-\frac{1}{R^{2j}}\int_{0}^{R}\theta_{0}(s)\phi_{0}(s)\phi_{j-1}(s)s^{2j-2}ds+\frac{\theta_{0}(R)}{R^{2j}\phi_{0}(R)}\int_{0}^{R}\phi_{0}^{2}(s)\phi_{j-1}(s)s^{2j-2}ds.
\end{align*}
By construction for all $s>0$ we have $0\le \theta_{0}(s)\phi_{0}(s)\le Cs$. For the second term on the right hand side above, we bound
\begin{align*}
	0< \frac{\theta_{0}(R)}{\phi_{0}(R)}\phi_{0}^{2}(s)\lesssim &\left\{ 
	\begin{array}{cc}  R^{-2n+2}\cdot s^{2n-1} &  0<s<R< 1 \\
	R^{-2n-2}\cdot s^{2n+3}&  1<s<R \\ 
	R^{-2n-2}s^{2n-1}&0<s<1<R
	\end{array}
	\right. \\
	\lesssim  &\,s,\quad \textrm{for all}\quad 0<s<R.
\end{align*}
Therefore the desired bound on $\phi_{j}(R)$ follows from an induction argument. 
\end{proof}

%
%
%
We next construct the Weyl-Titchmarsh solutions $\psi_{\pm}(R, \xi)$ for $R\xi^{\frac12}\gg 1$. 
\begin{lemma}\label{lem:psin R2xi large}
	There exists a pair of smooth functions $\psi_{\pm}(R,\xi)$ of $R,\xi>0$ and $R\xi^{\frac12}\gtrsim 1$ satisfying $-\calH_{n}^{+}\psi_{\pm}(R,\xi)=\xi\psi_{\pm}(R,\xi)$ of the form
	\begin{align*}
		\psi_{\pm}(R,\xi)=\frac{e^{\pm iR\xi^{\frac12}}}{\xi^{\frac14}}\left(1+g_{\pm}(R,\xi)\right).
	\end{align*}
	Here $g_{-}(R,\xi)=\overline{g_{+}(R,\xi)}$ and they satisfy the bounds (for some constants $c_{k}>0$)
	\begin{align*}
		|(R\partial_{R})^{k}g_{\pm}(R,\xi)|\leq c_{k}(R\xi^{\frac12})^{-1},\quad \textrm{and}\quad |(\xi\partial_{\xi})^{k}g_{\pm}(R,\xi)|\leq c_{k}(R\xi^{\frac12})^{-1}\quad \textrm{for}\quad R\xi^{\frac12}\gtrsim 1,
	\end{align*}
	and have the following asymptotic profile with some constant $C\in \bbR\setminus\{0\}$
	\begin{align*}
		g_{\pm}(R,\xi)\sim \frac{\psi_{1,\pm}(R)}{R\xi^{\frac12}}+O\left(\frac{1}{R^{2}\xi}\right),\quad \textrm{for}\quad R\xi^{\frac12}\rightarrow\infty,
	\end{align*}
	where
	\begin{align*}
		\psi_{1,\pm}(R)=\pm C\,i+O\left(\frac{1}{1+R^{2}}\right).
	\end{align*}
	In particular we have $W(\psi_+, \psi_-) = -2i$.
\end{lemma}
\begin{proof}
 We discuss $g_{+}(R,\xi)$ in detail and the argument for $g_{-}(R,\xi)$ is similar, therefore we omit it. Plugging the definition of  $g_{+}(R,\xi)$ into the equation satisfied by $\psi_{+}(R,\xi)$, we have
%
	\begin{align*}
	g_{+}(R, \xi) =& \int_R^\infty \frac{1-e^{2i(s-R)\xi^{\frac12}}}{2i\xi^{\frac12}}\left(1+g_{+}(s, \xi)\right)\cdot \left(\frac{(n+1)^2 - \frac14}{s^2}-\frac{4n}{s^{2}(1+s^{2})}-\frac{8}{(1+s^{2})^{2}}\right)\,ds\\
	:=&g_{+,0}(R,\xi)+\int_{R}^{\infty}K_{+,n}(R,s;\xi)g_{+}(s,\xi)\,ds
	\end{align*}
	where
	\begin{align*}
		g_{+,0}(R,\xi)=&\int_{R}^{\infty}\frac{1-e^{2i(s-R)\xi^{\frac12}}}{2i\xi^{\frac12}}\cdot \left(\frac{(n+1)^2 - \frac14}{s^2}-\frac{4n}{s^{2}(1+s^{2})}-\frac{8}{(1+s^{2})^{2}}\right)\,ds,\\
		K_{+,n}(R,s;\xi)=&\frac{1-e^{2i(s-R)\xi^{\frac12}}}{2i\xi^{\frac12}}\cdot \left(\frac{(n+1)^2 - \frac14}{s^2}-\frac{4n}{s^{2}(1+s^{2})}-\frac{8}{(1+s^{2})^{2}}\right).
	\end{align*}
	A direct calculation shows
	\begin{align*}
		|g_{+,0}(R,\xi)|\lesssim (R\xi^{\frac12})^{-1},\quad \sup_{R\leq s}\,(R\xi^{\frac12})\cdot |K_{+,n}(R,s;\xi)|(s\xi^{\frac12})^{-1}\lesssim \xi^{-\frac12}s^{-2}.
	\end{align*}
	A standard Volterra iteration gives the existence and uniqueness of such $g_{+}(R,\xi)$ and the estimates for $k=0$. To derive the bounds on derivatives, we introduce the notation $U_{n}(s):=\frac{(n+1)^2 - \frac14}{s^2}-\frac{4n}{s^{2}(1+s^{2})}-\frac{8}{(1+s^{2})^{2}}$, and the new variable $r:=s-R$. Therefore
	\begin{align*}
		&\left(R\partial_{R}g_{+,0}\right)(R,\xi)=\int_{0}^{\infty}\frac{1-e^{2ir\xi^{\frac12}}}{2i\xi^{\frac12}}\left(R\partial_{R}U_{n}\right)(r+R)dr,\\
		\Rightarrow\quad &|\left(R\partial_{R}g_{+,0}\right)(R,\xi)|\lesssim\int_{0}^{\infty}\xi^{-\frac12}(r+R)^{-2}dr\lesssim (R\xi^{\frac12})^{-1}.
	\end{align*}
	For the term involving $g_{+}(s,\xi)$, we have
	\begin{align*}
		R\partial_{R}\left(\int_{R}^{\infty}K_{+,n}(R,s;\xi)g_{+}(s,\xi)\,ds\right)=&\int_{0}^{\infty}\frac{1-e^{2ir\xi^{\frac12}}}{2i\xi^{\frac12}}R\partial_{R}\left(U_{n}(r+R)g_{+}(r+R,\xi)\right)dr\\
		=&\int_{0}^{\infty}\frac{1-e^{2ir\xi^{\frac12}}}{2i\xi^{\frac12}}R\partial_{R}\left(U_{n}(r+R)\right)g_{+}(r+R,\xi)dr\\&+\int_{0}^{\infty}\frac{1-e^{2ir\xi^{\frac12}}}{2i\xi^{\frac12}}U_{n}(r+R)\,R\partial_{R}\left(g_{+}(r+R,\xi)\right)dr\\
		:=&I_{R}+II_{R}.
	\end{align*}
	Therefore
	\begin{align*}
		|I_{R}|\lesssim &\int_{R}^{\infty}\xi^{-1}s^{-3}\,ds\lesssim (R^{2}\xi)^{-1},
	\end{align*}
	which is desired. For $II_{R}$ we have 
	\begin{align*}
		|II_{R}|\lesssim &\int_{R}^{\infty}\xi^{-\frac12}s^{-2}|s\partial_{s}g_{+}(s,\xi)|\,ds.
	\end{align*}
	Then a Volterra iteration argument gives the desired bound on 
	\begin{align*}
		R\partial_{R}g_{+}(R,\xi)=R\partial_{R}g_{+,0}(R,\xi)+R\partial_{R}\left(\int_{R}^{\infty}K_{+,n}(R,s;\xi)g_{+}(s,\xi)\,ds\right).
	\end{align*}
	
	Next we turn to $\xi\partial_{\xi}g_{+}(R,\xi)$. We have
	\begin{align*}
		\xi\partial_{\xi}g_{+,0}(R,\xi)=&\int_{0}^{\infty}\frac{i}{2}\xi\partial_{\xi}\left(e^{2ir\xi^{\frac{1}{2}}}\right)\xi^{-\frac{1}{2}}U_{n}(r+R)\,dr-\frac{i}{4}\int_{0}^{\infty}e^{2ir\xi^{\frac12}}\xi^{-\frac12}U_{n}(r+R)\,dr\\
		=&I_{\xi}+II_{\xi}.
	\end{align*}
	For $II_{\xi}$ we change back to $s$-variable:
	\begin{align*}
		II_{\xi}=-\frac{i}{4}\int_{R}^{\infty}e^{2i(s-R)\xi^{\frac12}}\xi^{-\frac12}U_{n}(s)\,ds,\quad \Rightarrow\quad |II_{\xi}|\lesssim (R\xi^{\frac12})^{-1}.
	\end{align*}
	For $I_{\xi}$, we rewrite it as
	\begin{align*}
		I_{\xi}=&\frac{i}{4}\int_{0}^{\infty}r\partial_{r}\left(e^{2ir\xi^{\frac12}}\right)\xi^{-\frac12}U_{n}(r+R)\,dr\\
		=&-\frac{i}{4}\int_{0}^{\infty}e^{2ir\xi^{\frac12}}\xi^{-\frac12}U_{n}(r+R)\,dr-\frac{i}{4}\int_{0}^{\infty}r\xi^{-\frac12}\partial_{r}U_{n}(r+R)\,dr
	\end{align*}
	A direct calculation implies $|I_{\xi}|\lesssim (R\xi^{\frac12})^{-1}$. For the term involving $g_{+}(s,\xi)$, we have
	\begin{align*}
		\xi\partial_{\xi}\left(\int_{R}^{\infty}K_{+,n}(R,s;\xi)g_{+}(s,\xi)ds\right)=:I'_{\xi}+II'_{\xi}+III'_{\xi},
		\end{align*}
		where
		\begin{align*}
			I'_{\xi}=&\frac{1}{2i}\int_{0}^{\infty}(1-e^{2ir\xi^{\frac12}})\xi\partial_{\xi}(\xi^{-\frac12})U_{n}(r+R)g_{+}(r+R,\xi)\,dr,\\
			II'_{\xi}=&\int_{0}^{\infty}\frac{1-e^{2ir\xi^{\frac12}}}{2i\xi^{\frac12}}U_{n}(r+R)\xi\partial_{\xi}\left(g_{+}(r+R,\xi)\right)\,dr\\
			III'_{\xi}=&\frac{i}{2\xi^{\frac12}}\int_{0}^{\infty}\xi\partial_{\xi}\left(e^{2ir\xi^{\frac12}}\right)U_{n}(r+R)g_{+}(r+R,\xi)\,dr.
		\end{align*}
		Based on the estimate for $g_{+}(s,\xi)$, the estimate for $I'_{\xi}$ is straightforward: $|I'_{\xi}|\lesssim (R^{2}\xi)^{-1}$. For $II'_{\xi}$ we use the Volterra iteration. For $III'_{\xi}$ we rewrite as
		\begin{align*}
			&\frac{i}{4\xi^{\frac12}}\int_{0}^{\infty}r\partial_{r}\left(e^{2ir\xi^{\frac12}}\right)U_{n}(r+R)g_{+}(r+R)\,dr.
		\end{align*}
		As in the estimate for $I_{\xi}$, we integrate by parts in $r$ and use the estimate for $(r+R)\partial_{r}g_{+}(r+R,\xi)$, which has already been proved, to obtain the desired bound on $|III'_{\xi}|\lesssim (R^{2}\xi)^{-1}$. The higher order derivative estimates follow in a similar vein.
		
		Next we turn to the asymptotic profile for $g_{+}(R,\xi)$. First, based on the estimate of the kernel and $g_{+}(s,\xi)$ we have
		\begin{align*}
			\left|\int_{R}^{\infty}K_{+,n}(R,s;\xi)g_{+}(s,\xi)\,ds\right|\lesssim (R\xi^{\frac12})^{-1}\cdot \xi^{-\frac12}\int_{R}^{\infty}s^{-2}\,ds\lesssim (R\xi^{\frac12})^{-2}.
		\end{align*}
		Therefore this contribution can be grouped into $O\left(\frac{1}{R^{2}\xi}\right)$. Next we turn to the contribution from $g_{+,0}(R,\xi)$. The contribution from $\int_{R}^{\infty}\frac{1}{2\,i\xi^{\frac12}}U_{n}(s)\, ds$ is already in the desired form of $\frac{\psi_{\pm,1}(R)}{R\xi^{\frac12}}$. For the contribution from $-\frac{1}{2i\xi^{\frac12}}\int_{R}^{\infty}e^{2i(s-R)\xi^{\frac12}}U_{n}(s)d\,s$, we write it as
		\begin{align*}
			-\frac{1}{2\,i\xi^{\frac12}}\int_{0}^{\infty}e^{2ir\xi^{\frac12}}U_{n}(r+R)d\,r=&-\frac{1}{2i\xi^{\frac12}}\int_{0}^{\infty}\frac{1}{2i\xi^{\frac12}}\partial_{r}\left(e^{2ir\xi^{\frac12}}\right)U_{n}(r+R)d\,r\\
			=&-\frac{1}{4\xi}U_{n}(R)-\frac{1}{4\xi}\int_{0}^{\infty}e^{2ir\xi^{\frac12}}\partial_{r}U_{n}(r+R)d\,r.
		\end{align*}
		The boundary term $-\frac{1}{4\xi}U_{n}(R)$ can be grouped into $O\left(\frac{1}{R^{2}\xi}\right)$. The integral term is bounded by
		\begin{align*}
			\left|\frac{1}{4\xi}\int_{0}^{\infty}e^{2ir\xi^{\frac12}}\partial_{r}U_{n}(r+R)d\,r\right|\lesssim \xi^{-1}\int_{0}^{\infty}(r+R)^{-3}d\,r\lesssim (R^{2}\xi)^{-1},
		\end{align*}
		and therefore is also grouped into $O\left(\left(R^2\xi\right)^{-1}\right)$.
	\end{proof}
 \begin{corollary}\label{cor:spectral measure expand small n}
     For $\xi$ sufficiently large and $1\leq n\leq N_{0}$, we have the expansion
     \begin{align*}
         \frac{d\rho_{n}(\xi)}{d\xi}=C_{0}+C_{1}\xi^{-1}+o(\xi^{-1})
     \end{align*}
     where $C_{0}, C_{1}$ are constants.
 \end{corollary}
 \begin{proof}
    In view of the behavior of $\phi(R,\xi)$ and $\Psi_{\pm}(R,\xi)$ constructed in Lemma \ref{lem: phin small R2xi small} and Lemma \ref{lem:psin R2xi large} respectively, the desired expansion follows by evaluating $\phi_{j}(R), j=0,1,2,...$ in Lemma \ref{lem: phin small R2xi small} and $\psi_{1,\pm}(R)$ in Lemma \ref{lem:psin R2xi large} at $R=\xi^{-\frac12}$ and then consider their Taylor expansions at ``$\xi^{-1}=0$".
 \end{proof}
	Next we compare the Fourier basis constructed in Proposition \ref{prop:DFT nlarge} and the solutions constructed in Lemma \ref{lem: phin small R2xi small} and Lemma \ref{lem:psin R2xi large}.
	\begin{proposition}\label{prop: FB match}
		Let $\phi(R,\xi)$ be as in Lemma \ref{lem: phin small R2xi small}. Then there exists a smooth complex nonvanishing function  $w_{n}(\xi)$ which has the following asymptotic behavior:
		\begin{align*}
			|w_{n}(\xi)|\sim \xi^{\frac{n+1}{2}},\quad \xi\rightarrow0+,\quad\textrm{and}\quad |w_{n}(\xi)|\sim \xi^{\frac{n-1}{2}}, \quad \xi\rightarrow\infty.
		\end{align*}
		such that for any fixed but large $n\leq N_{0}$, $w_{n}(\xi)\phi(R,\xi)$ and its higher order derivatives have the same asymptotic behavior as $\phi_{n}(R,\xi)$ constructed in Proposition \ref{prop:DFT nlarge} for $R\xi^{\frac12}\rightarrow0+$ and $ R\xi^{\frac12}\rightarrow\infty$. In particular, if  $\frac{d\rho_{n}(\xi)}{d\xi}$ is the spectral measure density associated to $w_{n}(\xi)\phi(R,\xi)$, then there exist absolute constants $C_{\ell}>0$ and $C>1$, such that
		\begin{align*}
			C^{-1}\leq \frac{d\rho_{n}(\xi)}{d\xi}\leq C,\quad \textrm{and}\quad \left|\frac{d^{\ell}}{d\xi^{\ell}} \frac{d\rho_{n}(\xi)}{d\xi}\right|\leq C_{\ell}\xi^{-\ell}\quad \textrm{for all}\quad \xi>0, \ell>0.
		\end{align*}
	\end{proposition}
\begin{proof}
	Since $\phi(R,\xi),\psi_{\pm}(R,\xi)$ are all solutions of $\left(\calH^{+}_{n}+\xi\right)f=0$, and $\phi(R,\xi)$ is real, there exists a complex function $a(\xi)$ such that
	\[
	\phi(R, \xi) = a(\xi)\psi_+(R, \xi) + \overline{a(\xi)}\psi_-(R, \xi),
	\]
	and thence 
	\[
	a(\xi) \simeq W(\phi(R, \xi), \psi_-(R, \xi)).
	\]
To obtain an estimate on $a(\xi)$, we choose $R, \xi$ such that $R\xi^{\frac12}\simeq 1$. Therefore if $\xi\ll1$, then $R\gg1$, and when $\xi\gtrsim 1$, we have $R\lesssim 1$. We start with the case when $\xi\ll1$. In order to get an upper bound on $|a(\xi)|$, we use that $|\phi(\xi^{-\frac12}, \xi)|\lesssim \xi^{-\frac{n}{2}-\frac34}$, $|(R\partial_R)\phi(\xi^{-\frac12}, \xi)|\lesssim \xi^{-\frac{n}{2}-\frac34}$, as well as 
	\[
	|\psi_\pm(R,\xi)| + |R\partial_R\psi_{\pm}(R,\xi)|\lesssim \xi^{-\frac14},
	\]
	whence 
	\begin{align}\label{a unnormalized 1}
|a(\xi)|\simeq	|W(\phi(\xi^{-\frac12}, \xi), \psi_-(\xi^{-\frac12}, \xi))|\lesssim \xi^{-\frac{n}{2}-\frac12}.
	\end{align} 
	In order to get the lower bound on $|a(\xi)|$, we evaluate our functions at $R \simeq \xi^{-\frac12}$. Then one finds a bound 
	\begin{align}\label{a unnormalized 2}
	|a(\xi)|\geq \frac{\big| \phi(\xi^{-\frac12}, \xi)\big|}{2\big|\psi_{\pm}(\xi^{-\frac12}, \xi)\big|}\gtrsim \frac{\xi^{-\frac{n}{2}-\frac34}}{\xi^{-\frac14}} = \xi^{-\frac{n+1}{2}}. 
	\end{align}
	This implies that the spectral measure of $w_{n}(\xi)\phi(R, \xi)$ behaves like 
	\begin{align}\label{spectral small n 1}
	\frac{d\rho_{n}(\xi)}{d\xi}\simeq \frac{1}{|w_{n}(\xi)a(\xi)|^2}\simeq 1,\quad \textrm{for}\quad \xi\ll1.
	\end{align}
	Next we turn to the case when $\xi\gtrsim 1$, we have $\phi(\xi^{-\frac12},\xi)\lesssim \xi^{-\frac{n}{2}+\frac14}$ and $(R\partial_{R})\phi(\xi^{-\frac12},\xi)\lesssim \xi^{-\frac{n}{2}+\frac14}$, and $|\psi_{\pm}(R,\xi)|\simeq \xi^{-\frac14}$. Therefore we have 
\begin{align}\label{a unnormalized 3}
	|a(\xi)|\simeq |W(\phi(\xi^{-\frac12},\xi),\psi_{-}(\xi^{-\frac12},\xi))|\lesssim \xi^{-\frac{n}{2}+\frac12},
\end{align}	
and 
\begin{align}\label{a unnormalized 4}
	|a(\xi)|\gtrsim \frac{|\phi(\xi^{-\frac12},\xi)|}{2|\psi_{\pm}(\xi^{-\frac12},\xi)|}\gtrsim\frac{\xi^{-\frac{n}{2}+\frac14}}{\xi^{-\frac14}}=\xi^{-\frac{n}{2}+\frac12}.
\end{align}
This again gives
\begin{align}\label{spectral small n 2}
\frac{d\rho_{n}(\xi)}{d\xi}\simeq \frac{1}{|w_{n}(\xi)a(\xi)|^2}\simeq 1,\quad \textrm{for}\quad \xi\gtrsim1.
\end{align}
Now \eqref{spectral small n 1} and \eqref{spectral small n 2} give the desired estimates on spectral measure. The bounds on the derivatives $\frac{d^{\ell}}{d\xi^{\ell}}\frac{d\rho_{n}(\xi)}{d\xi}$ follow directly from the construction of the function $a(\xi)$.

Next we analyze the asymptotic behavior of $w_{n}(\xi)\phi(R,\xi)$. Since now $\hbar$ has a positive lower bound, we equivalently consider the asymptotic regimes $R\xi^{\frac12}\rightarrow0+$ and $R\xi^{\frac12}\rightarrow\infty$. When $R\xi^{\frac12}\rightarrow0+$, Lemma \ref{lem:Fourier basis, large hbar} and Lemma \ref{lem: phin small R2xi small} imply, for $R\xi^{\frac12}\rightarrow0+$,
\begin{align}\label{FB small R2xi large hbar}
	\begin{split}
	w_{n}(\xi)\phi(R,\xi)\sim\begin{cases}
	&R^{n-\frac12}\xi^{\frac{n+1}{2}}\lesssim R^{n-\frac12}\xi^{\frac{n-1}{2}},\quad R\rightarrow0_{+},\quad \textrm{and}\quad \xi\rightarrow0+,\\
	&c_{1,n}R^{n+\frac32}\xi^{\frac{n+1}{2}}\lesssim R^{n-\frac12}\xi^{\frac{n-1}{2}},\quad R\rightarrow\infty,\quad\textrm{and}\quad  \xi\rightarrow0+,\\
	&R^{n-\frac12}\xi^{\frac{n-1}{2}},\quad R\rightarrow0_{+},\quad \textrm{and}\quad \xi\rightarrow\infty.
	\end{cases} 
	\end{split}
\end{align}
Therefore this agrees with the bound on $\phi_{n}(R,\xi)$ for $\hbar R\xi^{\frac12}\rightarrow0+$ in Proposition \ref{prop:DFT nlarge}. Now we look at the case when $R\xi^{\frac12}\rightarrow\infty$. We have, in view of Lemma \ref{lem:psin R2xi large}, and simply choosing $w_{n}(\xi)=\frac{1}{|a(\xi)|}$, 
\begin{align}\label{FB large R2xi large hbar}
	\begin{split}
	w_{n}(\xi)\phi(R,\xi)=&w_{n}(\xi)\left(a(\xi)\psi_{+}(R,\xi)+\overline{ a(\xi)}\psi_{-}(R,\xi)\right)\\
	=&\xi^{-\frac14}\Re\left(\frac{a(\xi)}{|a(\xi)|}e^{iR\xi^{\frac12}}(1+g_{+}(R,\xi))\right),
	\end{split}
\end{align}
	which again agrees with the asymptotic behavior of $\phi_{n}(R,\xi)$ in Proposition \ref{prop:DFT nlarge} as $\hbar R\xi^{\frac12}\rightarrow\infty$.
	The estimates on the higher order derivatives follow from the profile on $\phi(R,\xi)$ in Lemma \ref{lem: phin small R2xi small} and the derivative bounds on $g_{+}(R,\xi)$ given in Lemma \ref{lem:psin R2xi large}.
\end{proof}

Based on the analysis in Proposition \ref{prop:DFT nlarge} and Proposition \ref{prop: FB match}, we can justify the existence of the Fourier transform and its inverse transform appearing in Proposition \ref{prop:DFT nlarge}. 

\begin{proposition}\label{prop: rigor def FT ngeq2}
	The distorted Fourier transform associated with $\calH^{+}_{n},n\geq 2$ has the following property: for any $f\in C^{2}((0,\infty))$ with
	\begin{align*}
		\int_{0}^{\infty}\left(R^{-1}|f(R)|+|f'(R)|+R|f''(R)|\right)dR\leq M<\infty
	\end{align*}
	the Fourier transform
	\begin{align}\label{eq:hat f ngeq2 exist}
		\hat{f}(\xi)=\lim_{L\rightarrow\infty}\int_{0}^{L}\phi_{n}(R,\xi)f(R)dR
	\end{align}
	with $\phi_{n}(R,\xi)$ as in Proposition \ref{prop:DFT nlarge} for $n\geq N_{0}$ and $w_{n}(\xi)\phi(R,\xi)$ as in Proposition \ref{prop: FB match} for $2\leq n\leq N_{0}$, exists for all $\xi>0$ and 
	\begin{align}\label{inverse ngeq2 exist}
		\int_{0}^{\infty}|\hat{f}(\xi)||\phi_{n}(R,\xi)|d\xi\lesssim M.
	\end{align}
\end{proposition}
\begin{proof}
	Based on the analysis in Proposition \ref{prop:DFT nlarge} and Proposition \ref{prop: FB match}, we have
	\begin{align*}
		|\phi_{n}(R,\xi)|\lesssim \min\{R^{\frac12},\xi^{-\frac14}\},\quad \forall \xi>0, R>0.
	\end{align*}
	Thus, the left-hand side of \eqref{inverse ngeq2 exist} is bounded by $\int_{0}^{\infty}|\hat{f(\xi)}|\xi^{-\frac14}d\xi$ which we now proceed to estimate. Consider the partition of unity in \eqref{eq:POU} and define
	\begin{align*}
		A(\xi):=\int_{0}^{\infty}\chi_{0}(R^{2}\xi)\phi_{n}(R,\xi)f(R)dR,\quad B_{j}(\xi):=\int_{0}^{\infty}\chi(2^{-j}R^{2}\xi)\phi_{n}(R,\xi)f(R)dR.
	\end{align*}
	Therefore we have
	\begin{align*}
		&|A(\xi)|\lesssim \int_{0}^{\xi^{-\frac12}}R^{\frac12}|f(R)|dR,\\
		\Rightarrow\quad &\int_{0}^{\infty}\xi^{-\frac14}|A(\xi)|d\xi\lesssim \int_{0}^{\infty}\int_{0}^{R^{-2}}\xi^{-\frac14}d\xi|f(R)|R^{\frac12}dR\lesssim \int_{0}^{\infty}R^{-1}|f(R)|dR.
	\end{align*}
	To bound $B_{j}(\xi)$, we note that by Proposition \ref{prop:DFT nlarge} and Proposition \ref{prop: FB match}, for $R^{2}\xi\simeq 2^{j}$, one write $\phi_{n}(R,\xi)$ as
	\begin{align*}
		\phi_{n}(R,\xi)=\xi^{-\frac14}\Re\left(e^{iR\xi^{\frac14}}\sigma(R,\xi)\right),
	\end{align*}
	and $\sigma(R,\xi)$ satisfies
	\begin{align*}
		|(R\partial_{R})^{k}\sigma(R,\xi)|\lesssim 1,\quad |(\xi\partial_{\xi})^{k}\sigma(R,\xi)|\lesssim 1.
	\end{align*}
	Thus, we write
	\begin{align*}
		B_{j}(\xi)=\xi^{-\frac14}\int_{0}^{\infty}\chi(2^{-j}R^{2}\xi)e^{iR\xi^{\frac12}}\sigma(R,\xi)f(R)dR+\xi^{-\frac14}\int_{0}^{\infty}\chi(2^{-j}R^{2}\xi)e^{-iR\xi^{\frac12}}\sigma(R,\xi)f(R)dR.
	\end{align*}
	Using integration by parts, we have
	\begin{align*}
		|B_{j}(\xi)|\lesssim& \xi^{-\frac54}\int_{0}^{\infty}|\partial_{R}^{2}\left(\chi(2^{-j}R^{2}\xi)\sigma(R,\xi)f(R)\right)|dR\\\lesssim&\xi^{-\frac54}\int_{0}^{\infty}\one_{[R^2\xi\simeq 2^j]} \big(   |f''(R)| + R^{-1} |f'(R)|+R^{-2} |f(R)|\big)\, dR.
	\end{align*}
	Starting from this point we simply proceed as in the proof for Lemma \ref{lem:H1p} to obtain the finiteness of $\int_{0}^{\infty}\sum_{j=0}^{\infty}|B_{j}(\xi)|\xi^{-\frac14}d\xi$ and $\sum_{j=0}^{\infty}|B_{j}(\xi)|$. Therefore the proof is completed. 
\end{proof}

\subsection{The analysis for negative angular frequencies}

In this section we discuss the spectral theory for the operator $\calH_{-n}^{+}=\calH_{n}^{-}$ for $n\geq 2$. We start with the equation \eqref{eq:calHn}.
With the same notations as in \eqref{Schlag eq}, we consider the following eigenvalue problem associated to $\calH^{+}_{-n}$ for $n\ge2$:
\begin{align}\label{Schlag eq negative n}
\begin{split}
-\frac{1}{(n-1)^{2}}\partial_{R}^{2}f+V(R)f=&\frac{E^{2}}{(n-1)^{2}}f,\\
V(R):=&\frac{1}{R^{2}}-\frac{1}{4(n-1)^{2}R^{2}}+\frac{4n}{(n-1)^{2}R^{2}(R^{2}+1)}-\frac{8}{(n-1)^{2}(R^{2}+1)^{2}}.
\end{split}
\end{align}
By convention, we set $E<0$.
Introducing the notation $\thbar:=\frac{1}{1-n}$, we rewrite \eqref{Schlag eq negative n} as
\begin{align}\label{Schlag eq negative alt}
\begin{split}
-\thbar^{2}\partial_{R}^{2}f+V(R)f=&\thbar^{2}E^{2}f,\\
V(R)=&\frac{1}{R^{2}}\left(1+\frac{15\thbar^{2}}{4}-4\thbar\right)+\veps\left(R^{2},\thbar\right),
\end{split}
\end{align}
where recall 
\begin{align*}
\veps\left(R^{2},\thbar\right):=\frac{4\thbar}{R^{2}+1}-\frac{4\thbar^{2}(R^{2}+3)}{(R^{2}+1)^{2}}.
\end{align*}
The potential function looks identical to the one in \eqref{Schlag eq potential alter}, with $\hbar$ replaced by $\thbar$, and $\thbar =-1,-\frac12,-\frac13,...$. As for the positive angular frequencies, we distinguish between the regimes for $2\leq n\leq N_{0}$ and $n\geq N_{0}$, where $N_{0}$ is a large but fixed number. For $2\leq n\leq N_{0}$, we directly construct solutions to the eigenvalue problem associated to the operator 
\begin{align}\label{calH n operator minus}
-\calH^{+}_{-n}=-\partial_{R}^{2}+(n-1)^{2}V(R),\quad n\geq 2.
\end{align}
Following the same procedures as in the proofs of  Lemma~\ref{lem:Fourier basis, large hbar}, Lemma~\ref{lem: phin small R2xi small}, and Lemma~\ref{lem:psin R2xi large}, we obtain the following:
\begin{proposition}\label{prop: FB nega n large hbar}
	For $2\leq n\leq N_{0}$, there exists a solution $\phi_{-n}(R,\xi)$ to the eigenvalue problem $(\calH^{+}_{-n}+\xi)f=0$, which has the following asymptotic behavior for $R\xi^{\frac12}\rightarrow0_{+}$:
	\begin{align*}
		\phi_{-n}(R,\xi)\sim \begin{cases}
		&R^{n+\frac32}\xi^{\frac{n}{2}-\frac12}\leq R^{n-\frac12}\xi^{\frac{n}{2}-\frac12},\quad R\rightarrow0_{+}\quad \textrm{and}\quad \xi\rightarrow0_{+},\\
		&R^{n+\frac32}\xi^{\frac{n}{2}+\frac12}\leq R^{n-\frac12}\xi^{\frac{n}{2}-\frac12},\quad R\rightarrow0_{+}\quad \textrm{and}\quad \xi\rightarrow\infty,\\
		&c_{5}R^{n-\frac12}\xi^{\frac{n}{2}-\frac12},\quad R\rightarrow\infty\quad \textrm{and}\quad \xi\rightarrow0_{+}.
		\end{cases}
	\end{align*}
	For higher order derivatives, we have
	\begin{align*}
		|(R\partial_{R})^{k}\phi_{-n}(R,\xi)|\lesssim R^{n-\frac12}\xi^{\frac{n-1}{2}},\quad \textrm{and}\quad |(\xi\partial_{\xi})^{k}\phi_{-n}(R,\xi)|\lesssim R^{n-\frac12}\xi^{\frac{n-1}{2}},\quad \textrm{as}\quad R\xi^{\frac12}\rightarrow0_{+}.
	\end{align*}
	When $R\xi^{\frac12}\rightarrow\infty$, we have, for some complex function $a_{-}(\xi)$ which behaves $|a_{-}(\xi)|\simeq 1$,
	\begin{align*}
	\phi_{-n}(R,\xi)=2\xi^{-\frac14}\Re\left(a_{-}(\xi)e^{iR\xi^{\frac12}}(1+g^{-}(R,\xi))\right).
	\end{align*}
	Here$g^{-}(R,\xi)$ satisfies
	\begin{align*}
		g^{-}(R,\xi)= \frac{\psi_{2}(R)}{R\xi^{\frac12}}+O\left(\frac{1}{R^{2}\xi}\right),\quad \textrm{as}\quad R\xi^{\frac12}\rightarrow\infty,
	\end{align*}
	and the function $\psi_{2}(R)$ satisfies, for some constant $C\in \bbR\setminus\{0\}$,
	\begin{align*}
		\psi_{2}(R)=C\,i+O\left(\frac{1}{1+R^{2}}\right).
	\end{align*}
	For higher order derivatives, we have 
	\begin{align*}
		|(R\partial_{R})^{k}g^{-}(R,\xi)|\leq c_{k}(R\xi^{\frac12})^{-1},\quad \textrm{and}\quad |(\xi\partial_{\xi})^{k}g^{-}(R,\xi)|\leq c_{k}(R\xi^{\frac12})^{-1}, \quad \textrm{for}\quad R\xi^{\frac12}\gtrsim 1.
	\end{align*}
	This implies
	\begin{align*}
		(R\partial_{R})\left(e^{-iR\xi^{\frac12}}\phi_{-n}(R,\xi)\right)\sim c_{6}\xi^{-\frac14},\quad (\xi\partial_{\xi})^{k}\left(e^{-iR\xi^{\frac12}}\phi_{-n}(R,\xi)\right)\sim c_{7} \xi^{-\frac14},\quad \textrm{as}\quad R\xi^{\frac12}\rightarrow\infty.
	\end{align*}
\end{proposition}
For $n\geq N_{0}$, we introduce the notations
\begin{align*}
	\xbar:=\thbar ER,\quad \alpha:=\thbar E,\quad \tf(\xbar):=f(R).
\end{align*}
Then the equation \eqref{Schlag eq negative alt} transfers to 
\begin{align}\label{Schlag eq negative x}
\begin{split}
-\thbar^{2}\tf''(\xbar)+\tQ(\xbar)\tf(\xbar)=0,\quad \tQ(\xbar):=&\thbar^{-2}E^{-2}V^{-}\left(\frac{\xbar}{\thbar E}\right)-1\\
=&\xbar^{-2}\left(1+\frac{15\thbar^{2}}{4}-4\thbar\right)+\alpha^{-2}\veps_{-}\left(\frac{\xbar^{2}}{\alphab^{2}};\hbar\right)-1.
\end{split}
\end{align}
As in \eqref{def Q0}, we introduce
\begin{align}\label{def tQ0}
\tQ_{0}(\xbar;\alphab,\thbar):=\tQ(\xbar;\alphab,\thbar)+\frac{\thbar^{2}}{4\xbar^{2}}=\xbar^{-2}(1-2\thbar)^{2}+\alphab^{-2}\veps_{-}\left(\frac{\xbar^{2}}{\alphab^{2}};\thbar\right)-1.
\end{align}
As in \eqref{def tau} and \eqref{eq:qdef}, we introduce the variable
\begin{align}\label{def taub}
\taub(\xbar,\alpha;\thbar):=\sign(\xbar-\xbar_{t}(\alpha;\thbar))\left|\frac{3}{2}\int_{\xbar_{t}}^{\xbar}\sqrt{|\tQ_{0}(u,\alpha;\thbar)|}du\right|^{\frac23},
\end{align}
and the quantity
\begin{align}\label{def qb}
	\qb:=\frac{\tQ_{0}}{\taub},\quad \Rightarrow\quad \frac{d\taub}{d\xbar}=\taub'=\sqrt{\qb},\quad \frac{d}{d\taub}=\qb^{-\frac12}\frac{d}{d\xbar}.
\end{align}
We use ``$\prime$" to denote the derivative with respect to $\xbar$ and ``$\cdot$" to denote the derivative with respect to $\taub$. Then with the new variable $\wb:=\sqrt{\taub'}\tf$, \eqref{Schlag eq negative x} transfers to the following perturbed Airy equation:
\begin{align}\label{perturbed Airy negative}
	\thbar^{2}\ddot{\wb}(\taub)=\taub \wb(\taub)+\thbar^{2}\tVb(\taub;\alpha,\thbar) \wb(\taub)
\end{align}
where 
\begin{equation}\label{eq:master negative}
\tVb(\taub;\alpha,\thbar) = \ddot{\qb}^{\frac{1}{4}}\qb^{-\frac{1}{4}} - \frac{1}{4\xbar^{2}\qb}.
\end{equation}
We argue as in Section \ref{sec: ngeq2 linear}--\ref{sec:fund sys} to construct the Fourier basis associated to the operator $-\calH^{+}_{-n}$ for $n\geq N_{0}$. Since the argument is identical, we omit the details here.  Similar to 
 Lemma \ref{lem:solve pA}, we have the following (with $\zetab:=-\taub$):

\begin{lemma}
	\label{lem:solve pA negative}
	Let $w_0(\zetab;\thbar):=\Ai(|\thbar|^{-\frac23} \zetab)$ for $\zetab\ge0$ and $w_1(\zetab;\thbar):=\Ai(|\thbar|^{-\frac23} \zetab)+i \Bi(|\thbar|^{-\frac23} \zetab)$ for $\zetab\le0$. 
	Then the  Volterra integral equation  	
	\EQ{\label{eq:a0 system negative} 
		a_0(\zetab;\alpha,\thbar) & := \int_{\zetab}^\infty K_0(\zetab,s;\alpha,\thbar) (1+|\thbar| a_0(s;\alpha,|\thbar|))\, ds \\
		K_0(\zetab,s;\alpha,\thbar)  & = |\thbar|^{-1} \tVb(s;\alpha,\thbar) w_0^2(s; \thbar) \int_{\zetab}^{s} w_0^{-2}(t;\thbar)  \, dt 
	}
	has a unique bounded solution $a_0(\zetab;\alpha,\thbar)$ for all $\thbar=-\frac{1}{N_{0}-1},-\frac{1}{N_{0}},-\frac{1}{N_{0}+1},...$ and~$\alpha>0$, $\zetab\ge0$. One has $\lim_{\zetab\to\infty} a_0(\zetab;\alpha,\thbar)=0$ and 
	$w(\zetab;\alpha,\thbar):=w_0(\zetab;\thbar)(1+\thbar a_0(\zetab;\alpha,\thbar))$ is the unique solution of~\eqref{perturbed Airy negative} on $[0,\infty)$ with 
	$w(\zetab;\alpha,\thbar)\sim w_0(\zetab;\thbar)$ as $\zetab\to\infty$. 
	Analogously,  the Volterra integral equation  	
	\EQ{\label{eq:a1 system negative} 
		a_1(\zetab;\alpha,\thbar) & := \int^{\zetab}_{-\infty} K_1(\zetab,s;\alpha,\thbar) (1+\thbar a_1(s;\alpha,\thbar))\, ds \\
		K_1(\zetab,s;\alpha,\thbar)  & = |\thbar|^{-1} \tVb(s;\alpha,\thbar) w_1^2(s;\thbar) \int^{\zetab}_s w_1^{-2}(t;\thbar)  \, dt 
	}
	has a unique bounded solution $a_1(\zetab;\alpha,\thbar)$ for all $\thbar=-\frac{1}{N_{0}-1},-\frac{1}{N_{0}},-\frac{1}{N_{0}+1},...$ and~$\alpha>0$, $\zetab\le0$. One has $\lim_{\zetab\to-\infty} a_1(\zetab;\alpha,\thbar)=0$ and 
	$w(\zetab;\alpha,\thbar):=w_1(\zetab;\thbar)(1+\thbar a_1(\zetab;\alpha,\thbar))$ is the unique solution of~\eqref{perturbed Airy negative} on $(-\infty,0]$ with $w(\taub;\alpha,\thbar)\sim 
	w_1(\zetab;\thbar)$ as $\zetab\to-\infty$. 
\end{lemma}
We also have the analogue of  Proposition \ref{prop:tV summary}:
\begin{proposition}
	\label{prop:tVb summary}
	There exists a  constant $\taub_{*}>0$ and a small constant $0<\alpha_*\ll1$ so that uniformly in $\thbar=-\frac{1}{N_{0}-1},-\frac{1}{N_{0}},-\frac{1}{N_{0}+1},...$, 
	\EQ{\label{eq:ptVb}
		|\partial_{\taub}^k\partial_\alpha^\ell \tVb(\taub;\alpha,\thbar)| \le C_{k,\ell}\, \langle  \alpha\rangle^{-\ell-1} \langle  \taub\rangle^{-2-k},\qquad \forall\; \taub\ge -\taub_{*}
	}
	for all $k,\ell\ge0$ and $\alpha>0$. Moreover, \eqref{eq:ptVb} holds for $-\infty<\taub\leq-\taub_{*}$, 
	all $k,\ell\ge0$ and $\alpha\ge \alpha_{*}$. Finally, if $0<\alpha\leq \alpha_{*}$ and $-\infty\taub\leq -\tau_{*}$, then  
	\[
	\tVb(\taub;\alpha,\thbar) = \frac{5}{16\taub^{2}} -\taub \varphi(\xbar;\alpha,\thbar)
	\]
	where for all $k,\ell\ge0$
	\EQ{\label{eq:phi new negative}
		|\partial_{\taub}^k \partial_\alpha^\ell\varphi(\xbar;\alpha,\thbar)| 
		&\le C_{k,\ell}\, \min\big(|\thbar|\, \alpha^2 \xbar^{-2}+\xbar^{2}, |\thbar| \xbar^{2}/\alpha^2 + \xbar^{2}\big)  |\alpha|^{-\ell} (  -\taub  )^{\frac{k}{2}} 
	}
	Here $\xbar=\xbar(\taub;\alpha,\thbar)$ is the inverse of the diffeomorphism $\taub=\taub(\xbar;\alpha,\thbar)$ defined in~\eqref{def taub},  and satisfies
	\EQ{\label{eq:xbartaub}
		\frac{2}{3}\left(-\taub(\xbar;\alpha,\thbar)\right)^{\frac{3}{2}}  = \left\{   \begin{array}{lll} & -(1-O(\alpha^2))\log (\xbar) + O_{1}(1), & \text{\ \ if\ \ }  0< \alpha\le \xbar\le \xbar_*:=\xbar(-\tau_*;\alpha,\hbar) \\
			&  -\log (\xbar) + 2\thbar\log(\xbar/\alpha) + O_{2}(1),  &\text{\ \ if\ \ }    0< \xbar\le \alpha\le \alpha_*
		\end{array} 
		\right. 
	}
	Here $O(\alpha^2)$ is analytic in complex $|\alpha|\leq \alpha_{*}$, and bounded uniformly in  $\thbar=-\frac{1}{N_{0}-1},-\frac{1}{N_{0}},-\frac{1}{N_{0}+1}$. Furthermore, the two terms $O_1(1)$, resp.~$O_2(1)$ refer to smooth functions of $\taub,\alpha$ (and thus also of $\xbar$), uniformly bounded in $0<\alpha\le\alpha_{*}$, $-\infty<\tau\le-\tau_{*}$, and so that for all $k,\ell\ge0$ one has $\partial_{\xbar}^k \partial_\alpha^\ell O_1(1)=O(\xbar^{-k} \alpha^{-\ell})$ in the parameter regime of the first line of~\eqref{eq:xbartaub}, resp.\ $\partial_{\xbar}^k \partial_\alpha^\ell O_2(1)=O(\alpha^{-k-\ell})$ in the regime of the second line of~\eqref{eq:xbartaub}.
\end{proposition}

Based on Proposition \ref{prop:tVb summary}, we have the analogues of Lemma \ref{lem:a1fine}--Lemma \ref{lem:a0fine2}. 

\begin{lemma}
	\label{lem:a1fine negative}
	The functions $a_1(\zeta;\alpha,\thbar)$ from Lemma~\ref{lem:solve pA negative} satisfy the bounds
	\EQ{
		\label{eq:a1 fine negative}
		|\partial_{\alpha}^{\ell}a_{1}(\zetab;\alpha,\thbar)| & \leq C_{\ell}\langle\alpha\rangle^{-\ell}
		\langle\zetab\rangle^{-\frac{3}{2}},\quad \zetab\leq 0, \\ 
		|\partial_{\zetab}^{k} \partial_\alpha^\ell a_1(\zetab;\alpha,\thbar)| & \le C_{k,\ell}\, \langle \alpha\rangle^{-\ell} \left\{ \begin{array}{cc} 
			|\zetab|^{-\frac{3}{2}-k} & -\infty<\zetab\le -1 \\
			|\zetab|^{\frac{1}{2}-k} & -1<\zetab\le -|\thbar|^{\frac23} \\
			|\thbar|^{\frac{1-2k}{3}}&  -|\thbar|^{\frac23}\le\zetab\le 0, 
		\end{array}\right. 
	}
	for all $\ell\ge0, k\geq 1$, $\alpha>0$, and $\thbar=-\frac{1}{N_{0}-1},-\frac{1}{N_{0}},-\frac{1}{N_{0}+1},...$. 
\end{lemma}
\begin{lemma}
	\label{lem:a0fine1 negative}
	The functions $a_0(\zetab;\alpha,\thbar)$ from Lemma~\ref{lem:solve pA negative} satisfy the bounds
	\EQ{
		\label{eq:a0 fine negative}
		|\partial_{\alpha}^{\ell}a_{0}(\zetab;\alpha,\thbar)| & \leq C_{\ell}\, \alpha^{-\ell}
		\langle\zetab\rangle^{-\frac{3}{2}},\quad \zetab\geq 0, \\ 
		|\partial_{\zetab}^{k} \partial_\alpha^\ell a_0(\zetab;\alpha,\thbar)| & \le C_{k,\ell}\,  \alpha^{-\ell} \left\{ \begin{array}{cc} 
			\zetab^{-\frac{3}{2}-k} &1<\zetab<\infty \\
			\zetab^{\frac{1}{2}-k} & |\thbar|^{\frac23} <\zetab\le 1 \\
			|\thbar|^{\frac{1-2k}{3}}&  0 \le\zetab\le |\thbar|^{\frac23} 
		\end{array}\right. 
	}
	for all $\ell\ge0,k\geq 1$, $\alpha\gtrsim 1 $, and $\thbar=-\frac{1}{N_{0}-1},-\frac{1}{N_{0}},-\frac{1}{N_{0}+1},...$. 
\end{lemma}
and
\begin{lemma}
	\label{lem:a0fine2 negative}
	The functions $a_0(\zetab;\alpha,\thbar)$ from Lemma~\ref{lem:solve pA negative} satisfy the bounds
	\EQ{
		\label{eq:a0 fine* negative}
		|\partial_{\alpha}^{\ell}a_{0}(\zetab;\alpha,\thbar)| & \leq C_{\ell}\, \alpha^{-\ell}\big[\langle \zetab\rangle^{-\frac32}+ \min(1, \xbar(\zetab;\alpha,\thbar)^2/\alpha^2) \big]
		,\quad \zetab\geq 0, 		
	}
	for all $\ell\ge0$, $0<\alpha\ll 1 $, and $\thbar=-\frac{1}{N_{0}-1},-\frac{1}{N_{0}},-\frac{1}{N_{0}+1},...$. Here $\xbar(\zetab;\alpha,\thbar)$ is the diffeomorphism from Proposition~\ref{prop:tVb summary}. 
	Furthermore,  in the same parameter regime, 
	\EQ{
		\label{eq:a0 fine** negative}
		|\partial_{\alpha}^{\ell} \partial^{k}_{\zetab}\, a_{0}(\zetab;\alpha,\thbar)| & \leq  C_{k,\ell}\, \alpha^{-\ell} \zetab^{\frac{k}{2}}
		\big[\langle \zetab\rangle^{-\frac32}+ \min(1, \xbar(\zetab;\alpha,\thbar)^2/\alpha^2) \big]
		,\quad \zetab\geq1, 		
	}
	and 
	\EQ{
		\label{eq:a0 fine*** negative}
		|\partial_{\alpha}^{\ell}\partial^{k}_{\zetab}\, a_{0}(\zetab;\alpha,\thbar)| & \le  C_{k,\ell}\,  \alpha^{-\ell} \left\{ \begin{array}{cc} 
			\zetab^{\frac{1}{2}-k} & \quad |\thbar|^{\frac23} <\zetab\le 1 \\
			|\thbar|^{\frac{1-2k}{3}}&  \quad 0 \le\zetab\le |\thbar|^{\frac23} 
		\end{array}\right. 
	}
\end{lemma}
With the above preparations, we can now state the result on Fourier basis for $n\geq N_{0}$.
\begin{proposition}
	\label{prop:DFT nlarge negative}
	The distorted Fourier transform associated with $-\calH_{-n}^{+}$, $n\gg1$, takes the following form: 
	\EQ{\label{eq:hatfn negative}
		\hat{f}(\xi)  &=   \int_0^\infty \phi_{-n}(R,\xi) f(R)\, dR,\qquad 
		f(R)  = \int_0^\infty \phi_{-n}(R,\xi) \hat{f}(\xi) \rho_{-n}(d\xi)
	}
	with   
	\EQ{\label{eq:phin negative} 
		\phi_n(R;\xi) &= |\thbar|^{\frac13}\alpha^{-\frac12} \qb^{-\frac14}( \taub)\Ai(-|\thbar|^{-\frac23} \taub)(1+|\thbar| a_0(-\taub;\alpha,\thbar)) \\
		&= |\thbar|^{-\frac16}\xi^{-\frac14}  \qb^{-\frac14}( \taub)\Ai(-|\thbar|^{-\frac23} \taub)(1+|\thbar| a_0(-\taub;\alpha,\thbar)) \\
		\taub & =\taub(\xbar;\alpha,\thbar),  \quad \xbar= \alpha R\le 1,\quad \alpha=|\thbar|\xi^{\frac12}
	}
	Here $\taub=\taub(\xbar;\alpha,\thbar)$ as in \eqref{def taub}, $\qb$ is defined in~\eqref{def qb}. 
	One has $\phi_{-n}(R;\xi)\sim |\thbar|^{n-\frac12}\xi^{\frac{n-1}{2}}R^{n-\frac12}$ as $ R\xi^{\frac12}\to0+$. For $\xbar\ge1$
	\begin{align}
	\phi_{-n}(R;\xi) &=  -c_1 |\thbar|^{-\frac16}\xi^{-\frac14} \qb^{-\frac14}( \taub) \Re\Big(\big(1+|\thbar|\Xi(\xi;\thbar)\big) \big(\Ai(-|\thbar|^{-\frac23}\taub)+i\Bi(-|\thbar|^{-\frac23}\taub)\big)(1+|\thbar| a_1(\taub;\alpha,\thbar))\Big)  \nn \\
	& \sim -c_2\xi^{-\frac14} \Re\Big(\big(1+|\thbar|\Xi(\xi;\thbar)\big) e^{i\frac{\pi}{4}}e^{i\xi^{\frac12} R}\Big) \text{\ \ as\ \ } R\xi^{\frac12}\to\infty,  \label{eq:phin Rlarge negative}
	\end{align}
	where $c_1, c_2>0$ are absolute constants and $|\partial_\xi^k \Xi(\xi;\thbar)|\le C_k\, \xi^{-k}$ for all $k\ge0$ uniformly in $\thbar$. 
	The spectral measure~$\rho_{-n}$ is 
	purely absolutely continuous with density satisfying 
	\EQ{\label{eq:rho n negative}
		\frac12 & \le \frac{d\rho_{-n}(\xi)}{d\xi} \le 2,\qquad 
		\Big|\frac{d^\ell}{d\xi^\ell} \frac{d\rho_{-n}(\xi)}{d\xi}\Big| \le C_\ell\,\xi^{-\ell},\qquad\forall\; \xi>0,\;\ell\ge0 
	}
	uniformly in $n$. For the higher order derivatives of $\phi_{-n}(R,\xi)$ we have
	\begin{align*}
	(R\partial_{R})^{k}\phi_{-n}(R,\xi)\sim |\thbar|^{n-\frac12}\xi^{\frac{n-1}{2}}R^{n-\frac12},\quad (\xi\partial_{\xi})^{k}\phi_{-n}(R,\xi)\sim |\thbar|^{n-\frac12}\xi^{\frac{n-1}{2}}R^{n-\frac12},\quad \textrm{as}\quad  R\xi^{\frac12}\rightarrow0{+},
	\end{align*}
	and 
	\begin{align*}
	(R\partial_{R})^{k}\left(e^{-iR\xi^{\frac12}}\phi_{-n}(R,\xi)\right)\sim c_{8} \xi^{-\frac14},\quad (\xi\partial_{\xi})^{k}\left(e^{-iR\xi^{\frac12}}\phi_{-n}(R,\xi)\right)\sim c_{9} \xi^{-\frac14}\quad \textrm{as}\quad  R\xi^{\frac12}\rightarrow\infty.
	\end{align*}
\end{proposition} 
\begin{remark}\label{rmk: FB match negative}
	Proposition \ref{prop: FB nega n large hbar} and Proposition \ref{prop:DFT nlarge negative} show that the Fourier basis constructed for large and small $|\thbar|$ have the same asymptotic behavior as $R\xi^{\frac12}\rightarrow0_{+}$ and $R\xi^{\frac12}\rightarrow\infty$.
\end{remark}
\begin{remark}\label{rmk: FT exist negative}
	One can easily follow the same procedure in Proposition \ref{prop: rigor def FT ngeq2} to show that the Fourier transform and its inverse transform exist for $C^{2}$-functions $f$ satisfying
	\begin{align*}
		\int_{0}^{\infty}\left(R^{-1}|f(R)|+|f'(R)|+R|f''(R)|\right)dR<\infty.
	\end{align*}
\end{remark}
\begin{remark}\label{rmk: spectral measure negative}
	Based on Proposition \ref{prop: FB nega n large hbar} and Proposition \ref{prop:DFT nlarge negative}, one can follow similar procedures in the proof of these two propositions to show that the spectral measure density $\frac{d\rho_{-n}(\xi)}{d\xi}$ associated to $\phi_{-n}(R,\xi)$ satisfies
	\begin{align*}
		C^{-1}\leq\frac{d\rho_{-n }(\xi)}{d\xi}\leq C,\quad \left|\frac{d^{\ell}}{d\xi^{\ell}}\frac{d\rho_{-n}(\xi)}{d\xi}\right|\leq C_{\ell}\xi^{-\ell},\quad \forall \xi>0,\ell>0, n\ge2.
	\end{align*}
	Here the spectral measure density $\frac{d\rho_{-n}(\xi)}{d\xi}$ is given by
	\begin{align*}
		\frac{d\rho_{-n}(\xi)}{d\xi}=\frac{1}{4\pi|a_{-}(\xi)|^{2}}.
	\end{align*}
\end{remark}
\section{Transference identities}\label{sec:transference}

\subsection{The identities for angular momenta $n\geq2$ and $n\le -2$}
\subsubsection{$n\geq2$.}
Let us define the operator $\calK_\hbar$ as
\begin{align}\label{def calK}
\widehat{R\partial_{R}u}=-2\xi\partial_{\xi}\widehat{u}+\calK_\hbar\widehat{u}.
\end{align}
Here the Fourier transform is taken with respect to the Fourier basis $\phi(R,\xi;\hbar)=R^{-\frac12} \phi_*(R,\xi;\hbar)$ in $L^2(R\, dR)$ introduced in the previous section, with $\phi_{*}(R,\xi;\hbar)$ being the Fourier basis in $L^2(dR)$. Let $\rho_{n}(d\xi)$ be the  spectral 
measure associated with $\phi_{*}(R,\xi;\hbar)$, where $\hbar=\frac{1}{n+1}$. 
We have
\begin{align*}
\xb(\xi) &:=\calF_{*}\left(R^{\frac{1}{2}}f(R)\right)(\xi)=\langle R^{-\frac{1}{2}} \phi_{*}(R,\xi;\hbar),f(R)\rangle_{L^{2}_{RdR}},\\ f(R) &=\int_{0}^{\infty}\xb(\xi)R^{-\frac{1}{2}} \phi_{*}(R,\xi;\hbar)\rho_{n}(d\xi).
\end{align*}
Then one can easily check that the following Plancherel identity holds:
\begin{align*}
\|f\|^{2}_{L^{2}_{RdR}}=\|\langle f,R^{-\frac{1}{2}}\phi_{*}\rangle_{L^{2}_{RdR}}\|^{2}_{L^{2}_{\rho_{n}(d\xi)}}.
\end{align*}
Therefore $ \rho_{n}(d\xi)$ is also the spectral measure for $\phi(R,\xi;\hbar)$. The main result of this section is:

\begin{proposition}\label{prop: K operator}
For any $0<\hbar=\frac{1}{n+1}\le \frac13$ the operator $\calK_\hbar$ is given by
	\begin{align*}
	\calK_\hbar f(\xi)=-\left(2f(\eta)+ \frac{\eta\left(\frac{d\rho_{n}}{d\eta}(\eta)\right)^{\prime}}{\frac{d\rho_{n}(\eta)}{d\eta}} f(\eta)\right)\delta(\eta-\xi)+\left(\calK_\hbar^{(0)}f\right)(\xi)
	\end{align*}
	where the off-diagonal part $\calK^{(0)}_\hbar$ has a kernel $K_{0}(\xi,\eta;\hbar)$ given by
	\EQ{\label{eq:K0 kern}
		K_{0}(\xi,\eta;\hbar)= \frac{\frac{d\rho_{n}(\eta)}{d\eta}}{\xi-\eta}F(\xi,\eta;\hbar)
	}
	and the symmetric function $F(\xi,\eta;\hbar)$ satisfies (for any $0\leq k\leq k_{0}$ and sufficiently small $\hbar=\hbar(k_{0})$, where $k_{0}$ is arbitrary but fixed, and $\xi\leq \eta$)
	\begin{align}\label{F hbar bounds}
	\begin{split}
	|F(\xi,\eta;\hbar)|\lesssim& \left(\hbar\xi^{\frac12}\right)^{-1}\min\left\{1,\left(\hbar\xi^{\frac12}\right)^{3}\right\}\cdot G:=\Gamma_{\hbar}.
	\end{split}
	\end{align}
with 
\begin{align}\label{F hbar bounds auxi}
	G:=\begin{cases}
	&\min\left\{1,\left(\hbar\xi^{\frac12}\right)^{\frac14}\left|\eta^{\frac12}-\xi^{\frac12}\right|^{-\frac14}\right\},\quad \textrm{for}\quad \left|\frac{\eta^{\frac12}}{\xi^{\frac12}}-1\right|\lesssim1,\\
	&\hbar \left(\hbar\xi^{\frac12}\right)^{-1}\min\left\{1,\hbar\xi^{\frac12}\right\}\cdot\left(\frac{\xi^{\frac12}}{\eta^{\frac12}}\right)^{k},\quad \textrm{for}\quad \left|\frac{\eta^{\frac12}}{\xi^{\frac12}}-1\right|\gg1.
	\end{cases}
\end{align}
For $\hbar\gtrsim 1$, we have, for any $N>0$
	\begin{align}\label{F n bounds}
		\begin{split}
		|F(\xi,\eta;\hbar)|\lesssim \Gamma_{n}:= \begin{cases}
		&\eta|\log\eta|\cdot\frac{\xi^{\frac12}}{\eta^{\frac12}},\quad \textrm{for}\quad \eta\ll1,\\
		&\max\{\xi^{-\frac12}\eta^{-\frac12},\eta^{-\frac14}\xi^{-\frac14}\}\cdot\left(\frac{\xi^{\frac12}}{\eta^{\frac12}}\right)^{n}\cdot\left(\frac{\langle\xi\rangle}{\eta}\right)^{N},\quad \textrm{for}\quad \eta\gtrsim 1
		\end{cases}.
		\end{split}
	\end{align}

	For the derivatives of $F(\xi,\eta;\hbar)$, we have, for $k_{1}+k_{2}\leq 2$
		\begin{align}\label{deri F hbar and n bounds 3}
		\left|\partial_{\xi^{\frac12}}^{k_{1}}\partial_{\eta^{\frac12}}^{k_{2}}F(\xi,\eta;\hbar)\right|\lesssim \max\left\{1,\left(\hbar\xi^{\frac12}\right)^{-k_{1}}\right\}
		\cdot\max\left\{1,\left(\hbar\eta^{\frac12}\right)^{-k_{2}}\right\}\cdot \Gamma
		\end{align}
	Here $\Gamma$ is $\Gamma_{\hbar}$ if $\hbar\ll1$, and $\Gamma_{n}$ if $\hbar\gtrsim1$. The following estimate holds for trace-type derivatives for $\hbar\ll1$:
	\EQ{\label{eq:trace der}
		\left|\left(\partial_{\xi^{\frac12}}+\partial_{\eta^{\frac12}}\right)^{k}F(\xi,\eta;\hbar)\right| &\le C_k\, \xi^{-\frac{k}{2}} \cdot\Gamma,\quad \textrm{if}\quad \xi\simeq \eta,\; |\eta-\xi|\le \hbar^{\frac23}\xi
	}
for all $k\ge0$. 
		
\end{proposition}

\begin{proof}
	For simplicity,  we will often drop $\hbar$ from the notation. For example, we write $\phi(R,\xi)$ for $\phi(R,\xi;\hbar)$ etc. We also write $\rho_n(d\xi)=\rho(\xi)\, d\xi$. For $f\in C_{0}^{\infty}\left((0,\infty)\right)$ define 
	\[
	u(R) = \int_0^\infty \phi(R,\xi) f(\xi)\, \rho(\xi )\,d\xi 
	\]
	This function decays rapidly in $R$ due to the oscillations of $\phi(R,\xi) $, see Proposition~\ref{prop:DFT nlarge}. 
	Thus, 
	\begin{align}
	\widehat{R\partial_R u}(\xi) &= \left\langle \int_0^\infty R\partial_R  \phi(R,\eta) f(\eta) \rho(\eta)\, d\eta,  \phi(R,\xi) \right\rangle_{L^2(R\,dR)} \nn \\
	& = \left\langle \int_0^\infty [R\partial_R -2\eta\partial_\eta]  \phi(R,\eta) f(\eta) \rho(\eta)\, d\eta,    \phi(R,\xi) \right\rangle_{L^2(R\,dR)} \nn  \\
	&\qquad -2 \left\langle \int_0^\infty \phi(R,\eta)   {\partial_\eta \, (\eta f(\eta) \rho(\eta))}\, d\eta,    \phi(R,\xi) \right\rangle_{L^2(R\,dR)}    \nn     \\
	& = \left\langle \int_0^\infty [R\partial_R -2\eta\partial_\eta]  \phi(R,\eta) f(\eta) \rho(\eta)\, d\eta,    \phi(R,\xi) \right\rangle_{L^2(R\,dR)} \label{eq:main int}   \\
	&\qquad -2  f(\xi) -2 \frac{\xi\rho'(\xi)}{\rho(\xi)} f(\xi) - 2\xi f'(\xi) \nn 
	\end{align}
	We used that as limit of distributions 
	\begin{align}\label{delta kernel}
	\lim_{L\to\infty} \left\langle \chi(R/L)  \phi(R,\eta),  \phi(R,\xi) \right\rangle_{L^2(R\,dR)} = \rho(\xi)^{-1} \delta(\xi-\eta)
	\end{align}
	where $\chi$ is compactly supported smooth bump function which equals $1$ on $[0,1]$. 
	It follows that 
	\begin{align}\label{operator K pre}
	\begin{split}
	\calK f(\eta)=&\left\langle\int_{0}^{\infty}f(\xi)[R\partial_{R}-2\xi\partial_{\xi}]\phi(R,\xi)\rho(\xi)d\xi,\phi(R,\eta)\right\rangle_{L^{2}_{RdR}}\\
	&-2\left(1+\frac{\eta\rho'(\eta)}{\rho(\eta)}\right)f(\eta).
	\end{split}
	\end{align}
	To extract any  $\delta(\xi-\eta)$ appearing from the $RdR$ integral in the first line here we recall from Proposition~\ref{prop:DFT nlarge} and \eqref{eq:psi+ def} that 
	\EQ{ \label{eq:phi wieder}
		\phi(R,\xi) & = 2R^{-\frac12} \Re\left( a(\xi) \psi_+(R;\alpha)\right) \\
		&=    2 \hbar^{-\frac16} (R \xi^{\frac12})^{-\frac12} q^{-\frac14}( \tau) \Re\left( a(\xi) \left(\Ai(-\hbar^{-\frac23}\tau)-i\Bi(-\hbar^{-\frac23}\tau)\right)(1+\hbar \overline{a_1(-\tau;\alpha)})\right)   \\
		\psi_+(R,\xi;\hbar) & = \hbar^{\frac13} \alpha^{-\frac12} q^{-\frac14} (\Ai(-\hbar^{-\frac23} \tau)-i \Bi(-\hbar^{-\frac23} \tau)) (1+\hbar \overline{a_1(-\tau;\alpha,\hbar)} ) 
	}
	By Lemma~\ref{lem:a1fine} we have for all $\tau\ge0$
	\EQ{\label{eq:a1 absch}
		|\partial_{\alpha}^{\ell}\partial_\tau^k a_{1}(\tau;\alpha )|  \leq C_{\ell,k}\langle\alpha\rangle^{-\ell}
		\langle\tau\rangle^{-\frac{3}{2}-k} \qquad\forall\; k,\ell\ge0
	}
	and by \eqref{eq:AiiBi}, 
	\[
	\Ai(-\hbar^{-\frac23}\tau)-i\Bi(-\hbar^{-\frac23}\tau) = \pi^{-\frac12} e^{-i\frac{\pi}{4}} \hbar^{\frac16}\tau^{-\frac14} e^{i\frac{2}{3\hbar} \tau^{\frac{3}{2}}} \big(1 + \hbar O( \tau^{-\frac32})\big) \qquad\tau\to\infty 
	\]
	Inserting this into \eqref{eq:phi wieder} we obtain 
	\EQ{ \label{eq:phitauq}
		\phi(R,\xi) & =  2 \pi^{-\frac12}  (R \xi^{\frac12})^{-\frac12}  (\tau q(x;\alpha,\hbar))^{-\frac14}  \Re\Big(a(\xi) e^{i\frac{2}{3\hbar} \tau^{\frac{3}{2}}} e^{-i\frac{\pi}{4}}   (1+\hbar \tilde a_1(-\tau;\alpha))\Big) 
	}
	where $\tilde a_1$ satisfies the same bounds \eqref{eq:a1 absch} as $a_1$.  Here $x=\hbar \xi^{\frac12} R$ and $Q_0$ is defined in~\eqref{def Q0}. By~\eqref{eq:qtau4} and Lemma~\ref{lem: Lemma 3.4 CDST}
	\EQ{\nn 
		\tau^{-\frac14}q^{-\frac14} & = (1+\rho'(x;\alpha,\hbar))^{-\frac12} \\
		\frac{2}{3}\tau^{\frac32} &=  x-y(\alpha;\hbar)+\rho(x;\alpha,\hbar),
	}
	Writing $a_1$ in place of $\tilde a_1$ for simplicity (since both satisfy~\eqref{eq:a1 absch}) we conclude that 
	\EQ{ \label{eq:phi wieder2}
		\phi(R,\xi) & =  2 \pi^{-\frac12}  (R \xi^{\frac12})^{-\frac12}   (1+\rho'(\hbar \xi^{\frac12} R;\alpha,\hbar))^{-\frac12}   \Re \Big( a(\xi) e^{i\xi^{\frac12}R} e^{\frac{i}{\hbar}[-y(\alpha;\hbar)+\rho(\hbar \xi^{\frac12} R;\alpha,\hbar)]}    (1+\hbar   a_1(-\tau;\alpha)) \Big) 
	}
	For a function $g$ of the form $g=g(x;\alpha,\hbar)$, we have
	\begin{align}\label{rela R xi deri}
	R\partial_{R}g-2\xi\partial_{\xi}g=-\hbar\xi^{\frac12}\partial_{\alpha}g=-\alpha\partial_{\alpha}g.
	\end{align}
	Thus, applying this differential operator to~\eqref{eq:phi wieder2} yields 
	\EQ{ \nn 
		& ( R\partial_{R}-2\xi\partial_{\xi} )\phi(R,\xi)  =  -2 \pi^{-\frac12}  (R \xi^{\frac12})^{-\frac12}      \Re \Big( 2\xi a'(\xi) e^{i\xi^{\frac12}R} e^{\frac{i}{\hbar}[-y(\alpha;\hbar)+\rho(\hbar \xi^{\frac12} R;\alpha,\hbar)]}    (1+\hbar   a_1(-\tau;\alpha)) \Big) \\
		&  {+} 2 \hbar^{-1} \pi^{-\frac12} \alpha (R \xi^{\frac12})^{-\frac12}       \Re \Big(i \partial_\alpha y(\alpha;\hbar) a(\xi) e^{i\xi^{\frac12}R} e^{\frac{i}{\hbar}[-y(\alpha;\hbar)+\rho(\hbar \xi^{\frac12} R;\alpha,\hbar)]}    (1+\hbar   a_1(-\tau;\alpha)) \Big)  \\
		&  - 2 \hbar \pi^{-\frac12} \alpha (R \xi^{\frac12})^{-\frac12}     \Re \Big( a(\xi) e^{i\xi^{\frac12}R} e^{\frac{i}{\hbar}[-y(\alpha;\hbar)+\rho(\hbar \xi^{\frac12} R;\alpha,\hbar)]}    \partial_\alpha  a_1(-\tau;\alpha)  \Big)  + O(R^{-\frac52}) \\
		& = - 2 \pi^{-\frac12}  (R \xi^{\frac12})^{-\frac12}   \Re \Big\{ e^{i\hbar^{-1}\Psi(R;\xi,\hbar)} \Big( \big[  2\xi a'(\xi)  {-} i \xi^{\frac12} a(\xi) \partial_\alpha y(\alpha,\hbar)\big]  (1+\hbar   a_1(-\tau;\alpha)) + \hbar\alpha a(\xi)  \partial_\alpha  a_1(-\tau;\alpha) \Big) \Big\}  \\
		& \qquad + O(R^{-\frac52})
	}
	as $R\to\infty$  where
	\[
	\Psi(R;\xi,\hbar) = \hbar \xi^{\frac12}R -y(\alpha;\hbar)+\rho(\hbar \xi^{\frac12} R;\alpha,\hbar) 
	\]
	The final $O(R^{-\frac52})$ is a result of the derivatives $\alpha\partial_\alpha$ falling onto the $\rho$ terms. In addition, we used that $ (1+\rho'(\hbar \xi^{\frac12} R;\alpha,\hbar))^{-\frac12}  = 1+O(R^{-2})$.  It does not contribute to the diagonal $\xi=\eta$ and we therefore ignore it. 
	Using that $\Re z\Re w = \frac12 ( \Re(zw)+ \Re(z\overline{w}))$ we infer that the $\delta$ measure on the diagonal in the  integral, cf.~\eqref{eq:main int}
	\[
	\lim_{A\to \infty} \int_0^A [R\partial_R -2\xi\partial_\xi]  \phi(R,\xi)     \phi(R,\eta) R\,d R
	\]
	comes from the expression 
	\begin{align}
	- 2 \pi^{-1}  ( \xi\eta)^{-\frac14} \lim_{L\to\infty} \Re  \int_0^\infty  & e^{i\hbar^{-1}[\Psi(R;\xi,\hbar)-\Psi(R;\eta,\hbar)]} \Big( \big[  2\xi a'(\xi)  {-} i \xi^{\frac12} a(\xi) \partial_\alpha y(\alpha,\hbar)\big]  (1+\hbar   a_1(-\tau;\alpha))  \nn  \\
	&\quad   + \hbar\alpha a(\xi)  \partial_\alpha  a_1(-\tau;\alpha) \Big)    \overline{ a(\eta) }     (1+\hbar  \overline{  a_1(-\sigma;\beta) } ) \chi_1(R) \chi_2(R/L)   \, dR =: I(\xi,\eta;\hbar)
	\label{eq:Ixieta}
	\end{align}
	which is to be evaluated as a distributional limit as $L\to\infty$. Here $\chi_1$ is a s smooth cutoff which equals~$1$ near $\infty$ and which vanishes near~$0$, and $\chi_2=1-\chi_1$. It depends on the compact subset of~$(0,\infty)$ which holds $\xi,\eta$. Furthermore, 
	$\alpha$ is above, $\beta=\hbar \eta^{\frac12}$, and $\sigma$ is associated with $y=\hbar \eta^{\frac12} R$ the same way that $\tau$ is to $x$. We claim that we can discard all terms involving~$a_1$ from~$I(\xi,\eta;\hbar)$ as they do not contribute to $\delta$ on the diagonal. We note that by the relation between $\tau$ and $R$ we have $a_1(-\tau;\alpha) = O(R^{-1})$ as $R\to\infty$, with symbol behavior under differentiation in~$R$. Moreover this bound is uniform in $\xi$ in a compact set in $(0,\infty)$. 
	
	We first rewrite the phase as follows
	\EQ{ \label{eq:phase I} 
		\hbar^{-1}[\Psi(R;\xi,\hbar)-\Psi(R;\eta,\hbar)] &=  R(\xi^\frac12 - \eta^\frac12 ) [1+ \Omega(R;\xi,\eta,\hbar)  ]- \hbar^{-1}[ y(\hbar\xi^\frac12 ;\hbar) - y(\hbar\eta^\frac12 ;\hbar)] \\
		\Omega(R;\xi,\eta,\hbar) & = \hbar^{-1}R^{-1}(\xi^\frac12 - \eta^\frac12 )^{-1}[ \rho(\hbar \xi^{\frac12} R;\hbar\xi^\frac12,\hbar) - \rho(\hbar \eta^{\frac12} R;\hbar\eta^\frac12,\hbar) ] \\
		& = \int_0^1 [ \rho_x (R\gamma(s); \gamma(s),\hbar) + R^{-1} \rho_\alpha (R\gamma(s); \gamma(s),\hbar) ]\, ds\\
		\gamma(s) &= \hbar\eta^\frac12+ s\hbar (\xi^{\frac12} -\eta^{\frac12})
	}
	By Lemma~\ref{lem: Lemma 3.4 CDST}, one has $|\Omega(R;\xi,\eta,\hbar)|\lesssim R^{-2}$.   We introduce the new variable 
	\[
	\tilde R := R[1+ \Omega(R;\xi,\eta,\hbar)  ],\quad \frac{d\tilde R}{dR} = 1 + \partial_R (R\Omega(R;\xi,\eta,\hbar) ) = 1 + O(R^{-2})
	\]
	in \eqref{eq:Ixieta}. The $O(R^{-2})$ here as $R\to\infty$ is uniform in $\xi,\eta$ from a compact subset of $(0,\infty)$. Hence, this error term makes a bounded contribution to~\eqref{eq:Ixieta} and cannot contribute to any~$\delta$ measure. Consequently, we can ignore the $\Omega$ term in the phase of~$I(\xi,\eta;\hbar)$, see~\eqref{eq:phase I}. 
	
	Returning to the claim about $a_1$, note that an oscillatory integral of the form 
	\begin{align}\label{converge to L2}
	\int_0^\infty  e^{i R\zeta}   \omega(R) \chi_1(R) \chi_2(R/L)   \, dR  =: f_L(\zeta)
	\end{align}
	where $\omega(R)=O(R^{-1})$ has the property that $f_L$ converges strongly in $L^2(\R)$ to some $f\in L^2(\R)$. But this means in particular that any such expression cannot converge in the sense of distributions to an expression with a $\delta_0(\zeta)$ component. 
	
	In summary, any $\delta$ measure arising in~\eqref{eq:Ixieta} is already present in the simpler expression 
	\[
	- 2 \pi^{-1}  ( \xi\eta)^{-\frac14} \lim_{L\to\infty} \Re   \int_0^\infty   e^{i R(\xi^\frac12 - \eta^\frac12 )}    \big[  2\xi a'(\xi)  {-} i \xi^{\frac12} a(\xi) \partial_\alpha y(\alpha,\hbar)\big]\overline{ a(\eta) }    
	\chi_1(R) \chi_2(R/L)   \, dR =: I_0(\xi,\eta;\hbar)
	\]
	Note that we have also discarded the factor $e^{-i \hbar^{-1}[ y(\hbar\xi^\frac12 ;\hbar) - y(\hbar\eta^\frac12 ;\hbar)]}$ in the integral since it equals~$1$ on the diagonal. 
	By standard Fourier analysis,  since $y$ is real-valued 
	\EQ{\nn 
		I_0(\xi,\eta;\hbar) & = - 2 \pi^{-1} \pi   ( \xi\eta)^{-\frac14} \Re \big[  2\xi a'(\xi)  {-} i \xi^{\frac12} a(\xi) \partial_\alpha y(\alpha,\hbar)\big] \overline{ a(\eta) } \delta(\xi^{\frac12}-\eta^{\frac12}) +O(1)\\
		&= - 8 \xi \Re a'(\xi) \overline{ a(\xi) } \delta(\xi-\eta) +O(1)= -4 \xi \frac{d}{d\xi} \big( |a(\xi)|^2 \big) \delta(\xi-\eta) + O(1)
	} 
	as $\xi\to\eta$. 
	By \eqref{eq:a rep rho} and the asymptotics of $\psi_+$ in \eqref{eq:phi wieder}, 
	\EQ{\nn 
		W(\psi(\cdot,\xi), \overline{\psi(\cdot,\xi)}) &= -2i \pi^{-1}\\ 
		\rho(\xi) &= \frac{1 }{2i\pi |a(\xi)|^2  W(\psi(\cdot,\xi), \overline{\psi(\cdot,\xi)})}  = (4 |a(\xi)|^2)^{-1} 
	} 
	Recall that we are denoting the density of the spectral measure by~$\rho$.  Returning to \eqref{operator K pre} we finally obtain the exact $\delta$-measure contribution to the operator $\calK$.
	With an operator $\calK^{(0)}$ whose Schwartz kernel contains no $\delta$ measure on the diagonal, 
	\EQ{\label{operator K post}
		\calK f(\eta)= &- \rho(\eta) \eta\frac{d}{d\eta} \big( 4|a(\eta)|^2 \big) f(\eta) -2\left(1+\frac{\eta\rho'(\eta)}{\rho(\eta)}\right)f(\eta) + (\calK^{(0)} f)(\eta)  \\
		&= - \rho(\eta) \eta\frac{d}{d\eta} \big( \rho(\eta)^{-1} \big) f(\eta) -2\left(1+\frac{\eta\rho'(\eta)}{\rho(\eta)}\right)f(\eta) + (\calK^{(0)} f)(\eta) \\ 
		&=-\left(2+\frac{\eta\rho'(\eta)}{\rho(\eta)}\right)f(\eta) + (\calK^{(0)} f)(\eta)
	}
	in agreement with the statement of the proposition.

	To determine the off-diagonal part of the operator $\calK$, choose a  function $f\in C^{\infty}_{0}\left((0,\infty)\right)$ and define $u(R)$ to be the function
	\begin{align*}
	u(R):=\int_{0}^{\infty}f(\xi)\left(R\partial_{R}-2\xi\partial_{\xi}\right)\phi(R,\xi)\rho(\xi)d\xi,
	\end{align*}
	With $\phi(R,\xi)$ the Fourier basis of $H_n^+$, we have $-H_{n}^+\phi(R,\xi)=\xi\phi(R,\xi)$. Therefore
	\begin{align*}
		&\eta\left(\calK f\right)(\eta)-\calK\left(\xi f(\xi)\right)(\eta)\\
		=&\left\langle u,\eta\phi(R,\eta)\right\rangle_{L^{2}_{RdR}}-\left\langle\int_{0}^{\infty}\xi f(\xi)\left(R\partial_{R}-2\xi\partial_{\xi}\right)\phi(R,\xi)\rho(\xi)d\xi,\phi(R,\eta)\right\rangle_{L^{2}_{RdR}}\\
		=&-\left\langle H_{n}^{+}u,\phi(R,\eta)\right\rangle_{L^{2}_{RdR}}-\left\langle\int_{0}^{\infty}\xi f(\xi)\left(R\partial_{R}-2\xi\partial_{\xi}\right)\phi(R,\xi)\rho(\xi)d\xi,\phi(R,\eta)\right\rangle_{L^{2}_{RdR}}.
	\end{align*}
	On the other hand,
	\begin{align*}
	&\int_{0}^{\infty}\xi f(\xi)\left(R\partial_{R}-2\xi\partial_{\xi}\right)\phi(R,\xi)\rho(\xi)d\xi\\
	=&-\int_{0}^{\infty}f(\xi)\left(R\partial_{R}-2\xi\partial_{\xi}\right)\left(H_{n}^{+}\phi(R,\xi)\right)\rho(\xi)d\xi+2\int_{0}^{\infty}\xi f(\xi)\phi(R,\xi)\rho(\xi)d\xi\\
	=&\int_{0}^{\infty}f(\xi)\left[H_{n}^{+},R\partial_{R}\right]\phi(R,\xi)\rho(\xi)d\xi+2\int_{0}^{\infty}\xi f(\xi)\phi(R,\xi)\rho(\xi)d\xi-H_{n}^{+}u.
	\end{align*}
	The commutator $\left[H_{n}^{+},R\partial_{R}\right]$ is given by
	\begin{align*}
	[H^{+}_{n},R\partial_{R}]=2H_{n}^{+}+R(f_{n}- g_{n})'+2f_{n}-2g_{n}:=2H_{n}^{+}+W^{+}_{n}(R),
	\end{align*}
	and by the explicit expressions \eqref{fn gn expression} for $f_{n}(R)$ and $g_{n}(R)$, we have
	\begin{align}\label{W pre}
	\begin{split}
	&  R\partial_{R}f_{n}+2f_{n}=\frac{16}{(R^{2}+1)^{2}}-\frac{32}{(R^{2}+1)^{3}},\quad R\partial_{R}g_{n}+2g_{n}=-\frac{8n}{(R^{2}+1)^{2}},\\
	&W^{ {+}}_{n}(R)=\frac{16}{(R^{2}+1)^{2}}-\frac{32}{(R^{2}+1)^{3}}+\frac{8n}{(R^{2}+1)^{2}}.
	\end{split}
	\end{align}
	Therefore we conclude that
	\[
	\eta(\calK f)(\eta) - [\calK (\xi f(\xi))](\eta) = -\int_{0}^{\infty}f(\xi) \big\langle W^{+}_{n}(R) \phi(R,\xi),\phi(R,\eta)\rangle_{L^{2}_{RdR}} \, \rho(\xi)\, d\xi 
	\]
	The kernel function $F(\xi,\eta;\hbar)$ from \eqref{eq:K0 kern} is therefore given by
	\begin{align}\label{off diagonal F pre}
	F(\xi,\eta;\hbar):=&\int_{0}^{\infty}W^{ {+}}_{n}(R)\phi(R,\xi;\hbar)\phi(R,\eta;\hbar)R\, dR
	\end{align}
	Recall from~\eqref{eq:phitauq} and \eqref{eq:phi wieder2} the following representation of the oscillatory regime
	\EQ{ \label{eq:rechts von wp}
		\phi(R,\xi) & =  2 \hbar^{-\frac16} (R \xi^{\frac12})^{-\frac12} q^{-\frac14}( \tau) \Re\Big( a(\xi) \big(\Ai(-\hbar^{-\frac23}\tau)-i\Bi(-\hbar^{-\frac23}\tau)\big)(1+\hbar \overline{a_1(-\tau;\alpha)})\Big)  \\ 
		&= 2 \pi^{-\frac12}  (R \xi^{\frac12})^{-\frac12}  (\tau q(x;\alpha,\hbar))^{-\frac14}  \Re\Big(a(\xi) e^{i\frac{2}{3\hbar} \tau^{\frac{3}{2}}} e^{-i\frac{\pi}{4}}   (1+\hbar \tilde a_1(-\tau;\alpha))\Big) \\
		& =  2 \pi^{-\frac12}  (R \xi^{\frac12})^{-\frac12}   (1+\rho'(\hbar \xi^{\frac12} R;\alpha))^{-\frac12}   \Re \Big( a(\xi) e^{i\xi^{\frac12}R} e^{\frac{i}{\hbar}[-y(\alpha)+\rho(\hbar \xi^{\frac12} R;\alpha)]}    (1+\hbar  \tilde a_1(-\tau;\alpha)) \Big) 
	}
	where $\rho'(\hbar \xi^{\frac12} R;\alpha)$ is governed by~\eqref{eq:qtau4} and Lemma~\ref{lem: Lemma 3.4 CDST}, and $a(\xi)$ is as in~\eqref{eq:a size}. Here in view of \eqref{eq:AiiBi} and Lemma \ref{lem:a1fine}, $\tilde a_{1}(-\tau;\alpha)$ satisfies $|\tilde a_{1}(-\tau;\alpha)|\lesssim |\tau|^{-\frac32}$ as well as $\left|\partial_{\tau}^{\ell}\tilde a_{1}(-\tau;\alpha)\right|\lesssim |\tau|^{-\frac32-\ell}$ for all $\hbar^{\frac23}\lesssim \tau$. For ease of notation we suppress $\hbar$ from some of the notation (i.e., we write $\phi(R,\xi)$ in place of $\phi(R,\xi;\hbar)$, $a(\xi)$ instead of $a(\xi;\hbar)$ etc.). 
	
	On the other hand, to the left of the turning point, i.e., for $x= \hbar \xi^{\frac12} R\le x_t$, one has, cf.~\eqref{eq:phin} which we need to multiply by~$R^{-\frac12}$,  
	\EQ{\label{eq:links von wp}
\phi(R;\xi) &= \hbar^{\frac13}x^{-\frac12} q^{-\frac14}( \tau)\Ai(-\hbar^{-\frac23} \tau)(1+\hbar a_0(-\tau;\alpha,\hbar)) \\
&= \hbar^{\frac12} (x^2 Q_0(x,\alpha))^{-\frac14} e^{-\frac{2}{3\hbar}(-\tau)^{\frac32}} (1+\hbar \tilde a_0(-\tau;\alpha,\hbar)) 
}
Here $\tilde a_{0}(-\tau;\alpha,\hbar)$ satisfies $|\tilde a_{0}(-\tau;\alpha,\hbar)|\lesssim |\tau|^{-\frac32}$ and $\left|\partial_{\tau}^{\ell}\tilde a_{0}(-\tau;\alpha,\hbar)\right|\lesssim |\tau|^{-\frac32-\ell}$ for all $\tau\lesssim -\hbar^{\frac23}$.
The representation~\eqref{eq:rechts von wp} holds for $\tau\ge C\hbar^{\frac23}$, and \eqref{eq:links von wp} holds for $\tau\le -C\hbar^{\frac23}$ with some absolute constant~$C$. Here $\tau$ is the global variable defined by~\eqref{def tau}.   The third line of~\eqref{eq:rechts von wp} is only valid if $\tau\ge C$ since then also $x=\hbar \xi^{\frac12} R\gg 1$ and $\rho'(\hbar \xi^{\frac12} R;\alpha)=O(x^{-2})$, see Lemma~\ref{lem: Lemma 3.2 CDST}.  In the second line of~\eqref{eq:links von wp} we absorbed the error from the asymptotic expansion of $\Ai(-\hbar^{-\frac23} \tau)$ into the multiplicative correction $(1+\hbar \tilde a_0(-\tau;\alpha,\hbar))$. 
To pass to the second line we also used~\eqref{eq:qdef}, which means that we can use this second line only sufficiently far away from the unique zero of $Q_0$ which is the turning point~$x_t$. Thus, we need,  $0<x<\frac12 x_t$. On the other hand, if $\frac12 x_t\le x\le x_t$, then $x\simeq1$, $q\simeq 1$.

Without loss of generality, in the rest of the proof we assume $\xi\leq\eta$. We also first assume small $\hbar$.  To estimate \eqref{off diagonal F pre} we consider two regimes: ($\textbf{Case A}$) $\xi\simeq\eta$ and ($\textbf{Case B}$)  $\xi\ll\eta$.

	We start with {\bf Case A:}  {$\xi \simeq\eta$}. We break the integral in~\eqref{off diagonal F pre} into three separate regions, namely $x=\xi^{\frac12}\hbar R \ge C_1\gg1$ (case $A_3$), followed by $x\simeq 1$ (case $A_2$), and finally $0<x\ll 1$ (case $A_1$). 
	
	With $\chi(\xi^{\frac12}\hbar R)$ a smooth cutoff to the region in {\bf case~$\mathbf{A_3}$}, and with $\alpha=\hbar\xi^{\frac12}, \beta=\hbar\eta^{\frac12}$, and $\tau,\tilde\tau$ defined by \eqref{def tau} relative to $\xi$ and $\eta$, respectively, 
	\begin{align}
	&F_3(\xi,\eta;\hbar):=  \int_{0}^{\infty}W^{ {+}}_{n}(R)\phi(R,\xi;\hbar)\phi(R,\eta;\hbar)\chi(\xi^{\frac12}\hbar R) R\, dR \\
	& =4\pi^{-1}  \int_{0}^{\infty}W^{ {+}}_{n}(R) (R \xi^{\frac12})^{-\frac12}   (1+\rho'(\hbar \xi^{\frac12} R;\alpha))^{-\frac12}   \Re \Big( a(\xi) e^{i\xi^{\frac12}R} e^{\frac{i}{\hbar}[-y(\alpha)+\rho(\hbar \xi^{\frac12} R;\alpha)]}    (1+\hbar   \tilde a_1(-\tau;\alpha)) \Big)\nn  \\
	&  \qquad (R \eta^{\frac12})^{-\frac12}   (1+\rho'(\hbar \eta^{\frac12} R;\beta))^{-\frac12}   \Re \Big( a(\eta) e^{i\eta^{\frac12}R} e^{\frac{i}{\hbar}[-y(\beta)+\rho(\hbar \eta^{\frac12} R;\beta)]}    (1+\hbar   \tilde a_1(-\tilde\tau;\beta)) \Big) \chi(\xi^{\frac12}\hbar R) R\, dR\nn  \\
	& =4\pi^{-1}( \hbar\xi)^{-1} \zeta^{-\frac12} \int_{0}^{\infty}W^{ {+}}_{n}(x/\alpha)     (1+\rho'(x;\alpha))^{-\frac12}   \Re \Big( a(\xi)  e^{\frac{i}{\hbar}[x-y(\alpha)+\rho(x;\alpha)]}    (1+\hbar    O(x^{-1}) ) \Big) \nn \\
	&   \qquad \qquad\qquad(1+\rho'(x\zeta;\beta))^{-\frac12}   \Re \Big( a(\eta) e^{\frac{i}{\hbar}[x\zeta-y(\beta)+\rho(x\zeta;\beta)]}    (1+\hbar  O(x^{-1}) ) \Big) \chi(x) \, dx \label{eq:F3}
	\end{align}
	where $\zeta=(\eta/\xi)^{\frac12}\simeq1$. The third equality sign follows by a change of variables and the behavior of $a_1$ under differentiation in $\tau$, and thus in~$x$. The notation $ O(x^{-1}) $ refers to symbol-type behavior under differentiation. 
	Then, on the one hand, placing absolute values inside the integrals yields 
	\EQ{\label{eq:F3 0}
		|F_3(\xi,\eta;\hbar)| \lesssim \hbar^{-1}\xi^{-\frac12} \int_{\xi^{\frac12}\hbar R\gtrsim 1} \langle R\rangle^{-4}\, dR \lesssim (\hbar \xi^{\frac12})^{-1}  \min(1,  (\hbar  \xi^{\frac12} )^{3})  
	}
	On the other hand, integrating by parts in \eqref{eq:F3} yields
	\EQ{ \label{eq:F3 IBP}
		|F_3(\xi,\eta;\hbar)| \lesssim (\hbar \xi^{\frac12})^{-1}  |\xi^{\frac12}-\eta^{\frac12}|^{-1} \min (1,(\hbar \xi^{\frac12})^4)
	} 
	Indeed, the oscillatory terms in \eqref{eq:F3}  are of the form $e^{\pm \frac{i}{\hbar}\Phi_+(x,\alpha;\hbar,\zeta)}$ and 
	$e^{\pm  \frac{i}{\hbar}\Phi_-(x,\alpha;\hbar,\zeta)}$ 
	where
	\begin{align}\label{eq:F3 phases}
	\Phi_{\pm}(x,\alpha;\hbar,\zeta) = x(\zeta\pm 1) \pm \rho(x;\alpha) + \rho(x\zeta;\alpha\zeta)
	\end{align}
	We will only discuss the destructive interference given by $\Phi_{-}$ as $\Phi_+$ contributes less. One has
	\EQ{\label{eq:delx Phi-}
		\partial_x \Phi_{-}(x,\alpha;\hbar) &= (\zeta-1)\big\{1 +\int_0^1 \big[(\rho_x + y(s)\rho_{xx})(y(s),d(s)))
		+d(s)\rho_{x\alpha}(y(s),d(s))\big]\, ds \big\}\\
		&=  (\zeta-1)[1 + O(x^{-2})]
	}
	with $y(s)=x(1-s+s\zeta), d(s)=\alpha(1-s+s\zeta)$,
	where the final $O(x^{-2})$ is a consequence of  Lemma~\ref{lem: Lemma 3.4 CDST}.  Using that 
	\[(i\hbar^{-1} \partial_x \Phi_{-}(x,\alpha;\hbar,\zeta))^{-1}\partial_x e^{  \frac{i}{\hbar}\Phi_-(x,\alpha;\hbar,\zeta)} = e^{ \frac{i}{\hbar}\Phi_-(x,\alpha;\hbar,\zeta)}
	\]
	one integration by parts in \eqref{eq:F3} yields 
	\EQ{\nn
		|F_3(\xi,\eta;\hbar)| & \lesssim (\hbar \xi)^{-1} |1-\zeta|^{-1} \int_1^\infty \Big( \langle x/\alpha\rangle^{-6} x\alpha^{-2} + \langle x/\alpha\rangle^{-4} x^{-2}\Big)\, dx\\
		&\lesssim  ( \hbar\xi)^{-1} |1-\zeta|^{-1} \min (1,(\hbar \xi^{\frac12})^4),
	}
	which is the same as \eqref{eq:F3 IBP}.  After $\ell\ge1$ integrations by parts we obtain
	\EQ{\label{eq:F3 ell}
		|F_3(\xi,\eta;\hbar)| &\lesssim   ( \hbar\xi^{\frac12})^{-2} \hbar^\ell |1-\zeta|^{-\ell} \min (1,(\hbar \xi^{\frac12})^4) \\
		&\lesssim ( \hbar\xi^{\frac12})^{\ell-2} |\xi^{\frac12} -\eta^{\frac12} |^{-\ell} \min (1,(\hbar \xi^{\frac12})^4)
	}
	Combining \eqref{eq:F3 0} with \eqref{eq:F3 ell} (with $\ell=5$) one obtains with $\alpha=\hbar\xi^{\frac12}$, 
	\EQ{\label{eq:F3 summary}
		|F_3(\xi,\eta;\hbar)| &\lesssim \alpha^{-1}\min(1,\alpha^3)\min\big(1,\alpha^5 |\xi^{\frac12} -\eta^{\frac12} |^{-5} \big)
	}
	We now turn to {\bf case $\mathbf{A_2}$}, i.e., the second integral  for which $x\simeq 1$, i.e., $-C\le \tau\le C$ with some   arbitrary but fixed constant~$C$. Then, in view of \eqref{eq:links von wp} and \eqref{eq:rechts von wp}, with a smooth cutoff $\chi(x)$ to $x\simeq1$, 
	\begin{align}\label{eq:F2 pre}
	\begin{split}
	F_2(\xi,\eta;\hbar) &:=  \int_{0}^{\infty}W^{ {+}}_{n}(R)\phi(R,\xi;\hbar)\phi(R,\eta;\hbar)  \chi(x) R\, dR \\
	& = \int_{0}^{\infty}W^{ {+}}_{n}(R) \chi(\alpha R )R\,  \hbar^{\frac13}x^{-\frac12} q^{-\frac14}( \tau)\Big[  \Ai(-\hbar^{-\frac23} \tau)(1+\hbar a_0(-\tau;\alpha,\hbar)) \chi_{[\tau\lesssim  \hbar^{\frac23}]} + \\
	&\qquad + 2  \Re\big( a(\xi) \big(\Ai(-\hbar^{-\frac23}\tau)-i\Bi(-\hbar^{-\frac23}\tau)\big)(1+\hbar \overline{a_1(-\tau;\alpha)})  \big) \chi_{[\tau\gtrsim  \hbar^{\frac23}]}\Big] \\
	&\qquad \hbar^{\frac13}\tilde x^{-\frac12} \tilde q^{-\frac14}( \tilde \tau)\Big[  \Ai(-\hbar^{-\frac23}\tilde \tau)(1+\hbar a_0(-\tilde \tau;\beta,\hbar)) \chi_{[\tilde \tau\lesssim  \hbar^{\frac23}]} + \\
	&\qquad + 2  \Re\big( a(\eta) \big(\Ai(-\hbar^{-\frac23}\tilde\tau)-i\Bi(-\hbar^{-\frac23}\tilde\tau)\big)(1+\hbar \overline{a_1(-\tilde\tau;\beta)})  \big) \chi_{[\tilde\tau\gtrsim \hbar^{\frac23}]} \Big]\, dR
	\end{split}
	\end{align}
	Here $\chi_{[\tau\lesssim \hbar^{\frac23}]}$ and $\chi_{[\tau\gtrsim \hbar^{\frac23}]}$ are smooth cutoffs to the indicated regions and a partition of unity, i.e., $$\chi_{[\tau\lesssim \hbar^{\frac23}]}+\chi_{[\tau\gtrsim  \hbar^{\frac23}]}=1.$$
	In~\eqref{eq:F2 pre}, we are passing beyond the turning point in the non-oscillatory fundamental solution by an amount $\simeq \hbar^{\frac23}$, cf.~Remark~\ref{rem:hbar32}.  
	As usual, $\tilde\tau$ and $\beta=\zeta\alpha$ play the role of $\tau$ and $\alpha$, respectively, with $\eta$ in place of~$\xi$. 
	To handle case~$A_2$, 
	we analyze $\tilde\tau$ as a function of~$\tau$.  By Lemmas~\ref{lem: monotonicity of root in alpha}  and~\ref{lem: Lemma 3.2 CDST} we have 
	\EQ{\nn 
		\tau &= (x-x_t(\alpha))\Phi(x;\alpha),\qquad \tilde \tau = (\zeta x-x_t(\zeta\alpha))\Phi(\zeta x;\zeta \alpha) 
	}
	where $|\partial_\alpha \Phi(x,\alpha)|+|\partial_\alpha x_t(\alpha)|\lesssim \hbar \langle \alpha\rangle^{-3}$. In particular, $\tau=0$ corresponds to  (with $x_t=x_t(\alpha)$)
	\EQ{\label{eq:tau 0}
		\tilde\tau_0 &:= (\zeta x_t(\alpha)-x_t(\zeta\alpha))\Phi(\zeta{x_{t}};\zeta \alpha) \\
		& = (\zeta-1)\big[ x_{t} -\alpha\int_0^1 x_t'(s\zeta\alpha+(1-s)\alpha)\, ds \big]  {\Phi(\zeta x_{t};\zeta\alpha)}= (\zeta-1) \big[ 1 +O(\hbar) \big],
	}
	see \eqref{eq:xt ell}, and 
	\EQ{ \label{eq:tiltau tau} 
		\tilde \tau &= \zeta\tau \frac{\Phi(\zeta x;\zeta\alpha)}{\Phi(x;\alpha)} + \tilde\tau_0\cdot {\frac{\Phi(\zeta x;\zeta\alpha)}{\Phi(\zeta x_{t};\zeta\alpha)}}
	}
	Since $\zeta\ge 1$, we have $\tilde\tau_0\ge 0$.
	By Lemma \ref{lem: Lemma 3.2 CDST}, we have
	\begin{align*}
	\Phi(\zeta x;\zeta \alpha)\simeq 1,\quad \Phi(x;\alpha)\simeq 1
	\end{align*}
	and thus
	\EQ{
		\label{eq:tautil size}
		\tilde \tau\simeq \tau+\zeta-1
	}
	uniformly in the parameters. 
	It follows that, see~\eqref{eq:qdef}, 
	\EQ{\label{eq:tau p}
		\frac{\partial\tilde\tau}{\partial\tau} &= \zeta\frac{\Phi(\zeta x;\zeta\alpha)}{\Phi(x;\alpha)} +\zeta\tau \left( \frac{\Phi(\zeta x;\zeta\alpha)}{\Phi(x;\alpha)}\right)_x q^{-\frac12}+ {\tilde{\tau}_{0}\left(\frac{\Phi(\zeta x;\zeta\alpha)}{\Phi(\zeta x_{t};\zeta\alpha)}\right)_{x}q^{-\frac12} }\\
		\frac23\partial_\tau \left( \tilde\tau^{\frac32}-\tau^{\frac32} \right) &= \tau^{\frac12} \left( -1 + \frac{\partial\tilde\tau}{\partial\tau}\right) + \frac{\tilde\tau-\tau}{\tau^{\frac12}+\tilde\tau^{\frac12}}\frac{\partial\tilde\tau}{\partial\tau}
	}
	We claim that the functions $\calR, \calS$ defined by
	\begin{align}\label{def R S T}
		\begin{split}
		\zeta\frac{\Phi(\zeta x;\zeta\alpha)}{\Phi(x;\alpha)} - 1 & =: (\zeta-1) \calR(\zeta,x;\alpha,\hbar)  \\
		\left( \frac{\Phi(\zeta x;\zeta\alpha)}{\Phi(x;\alpha)}\right)_x &=: (\zeta-1) \calS(\zeta,x;\alpha,\hbar)\\
		{\left(\frac{\Phi(\zeta x;\zeta\alpha)}{\Phi(\zeta x_{t};\zeta\alpha)}\right)_{x}} & {=: \calT(\zeta, x;\alpha,\hbar)}
		\end{split}
	\end{align}
	satisfy uniformly in this regime of parameters 
	\EQ{\label{eq:RS claim}
		& |\calR(\zeta,x;\alpha,\hbar)| + |\calS(\zeta,x;\alpha,\hbar)| +  |\calT(\zeta,x;\alpha,\hbar) |\lesssim  1
	}
	To estimate $\calR$, we write
	\begin{align*}
	(\zeta-1)\Phi(\zeta x;\zeta\alpha)+\Phi(\zeta x;\zeta\alpha)-\Phi(x;\alpha)=(\zeta-1)\calR(\zeta,x;\alpha,\hbar)\Phi(x;\alpha).
	\end{align*}
	The difference on the left-hand side is 
	\EQ{\label{eq:Phi diff}
		\Phi(\zeta x;\zeta\alpha)-\Phi(x;\alpha)=&(\zeta-1)\int_{0}^{1}\left(x(\partial_{x}\Phi)+\alpha(\partial_{\alpha}\Phi)\right)(s\zeta x+(1-s)x;s\zeta\alpha+(1-s)\alpha)\, ds
	}
	By \eqref{eq:kell g0},  $\left|\alpha\partial_{\alpha}\Phi\right|\lesssim \hbar$ and $|\partial_{x}\Phi|\lesssim 1$, whence $|\calR(\zeta,x;\alpha,\hbar)|\lesssim 1$. 
	For $\calS(\zeta,x;\alpha,\hbar)$ we consider
	\begin{align*}
	\left(\frac{\Phi(\zeta x;\zeta\alpha)}{\Phi(x;\alpha)}\right)_{x}=\left(\frac{\Phi(\zeta x;\zeta\alpha)-\Phi(x;\alpha)}{\Phi(x;\alpha)}\right)_{x}
	=&-\frac{\Phi(\zeta x;\zeta\alpha)-\Phi(x;\alpha)}{\Phi(x;\alpha)^{2}}\cdot\partial_{x}\Phi(x;\alpha)+\frac{\left(\Phi(\zeta x;\zeta\alpha)-\Phi(x;\alpha)\right)_{x}}{\Phi(x;\alpha)}\\
	=:&(\zeta-1)(I+II).
	\end{align*}
	$|\partial_{x}\Phi(x;\alpha)|\lesssim 1$ and \eqref{eq:Phi diff} imply $|I|\lesssim 1$. Furthermore, 
	\begin{align}\label{Phizeta minus Phi dx}
	\begin{split}
	&	\partial_{x}(\Phi(\zeta x;\zeta\alpha)-\Phi(x;\alpha))\\
	=&(\zeta-1)\int_{0}^{1}\left(\left((\zeta-1)x^{2}(\partial^{2}_{x}\Phi)+(\zeta-1)x\alpha(\partial^{2}_{x\alpha}\Phi)+\partial_{x}\Phi\right)(s\zeta x+(1-s)x;s\zeta\alpha+(1-s)\alpha)\right)\, ds.
	\end{split}
	\end{align}
	Bounding the integrand by Lemma \ref{lem: Lemma 3.2 CDST} implies $|II|\lesssim 1$, whence  $|\calS(\zeta,x;\alpha,\hbar)| \lesssim 1$.
	The estimate on  $\calT(\zeta,x;\alpha,\hbar)$ follows from 
	\begin{align*}
	\calT(\zeta,x;\alpha,\hbar)=	\left(\frac{\Phi(\zeta x;\zeta\alpha)}{\Phi(\zeta x_{t};\zeta\alpha)}\right)_{x}=\frac{\zeta\partial_{x}\Phi(\zeta x;\zeta\alpha)}{\Phi(\zeta x_{t};\zeta\alpha)}.
	\end{align*}
	and \eqref{eq:RS claim} holds. 
	Returning to \eqref{eq:tau p}, we conclude that 
	\begin{align}\label{eq: tau minus tiltau dtiltau dtau}
	\begin{split}
	\frac{\partial\tilde\tau}{\partial\tau} -1 &= (\zeta-1) \big[ \calR(\zeta,x;\alpha,\hbar)  + \zeta\tau \calS(\zeta,x;\alpha,\hbar) q^{-\frac12} {+(1+O(\hbar))\calT(\zeta,x;\alpha,\hbar)q^{-\frac12}}\big]\\
	\tilde \tau -\tau &= \int_0^\tau \left [\frac{\partial\tilde\tau}{\partial\tau} -1\right ]\, dv + \tilde \tau_0 \\
	& = (\zeta-1) \big[ 1+ {O(\tau)}  +O(\hbar) \big], 
	\end{split}
	\end{align}
	where we invoked \eqref{eq:tau 0} to pass to the  third  line. 
	By~\eqref{eq:tau p}  and~\eqref{eq:tautil size}
	\EQ{\label{eq:32 diff}
		\partial_\tau \big( \tilde\tau^{\frac32}-\tau^{\frac32} \big)\simeq (\zeta-1) (\tau+\tilde\tau_0)^{-\frac12}\simeq(\zeta-1) (\tau+\zeta-1)^{-\frac12}
	}
	uniformly in $0\le \tau \ll1$, $\alpha>0$ and small $\hbar$ (recall that we are in the range $1\le\zeta\lesssim 1$). 
	We will also require the following bound on the second derivatives 
	\EQ{\label{eq:32 diff 2}
		\left| \partial_\tau^2 \left( \tilde\tau^{\frac32}-\tau^{\frac32} \right) \right|\lesssim (\zeta-1)\tau^{-\frac12} (\tau+\zeta-1)^{-1}
	}
	for the same regime of parameters. 
	We write the left-hand side here in the form
	\begin{align}\label{diff 326}
	\begin{split}
	\frac43\partial_{\tau}^{2}\left(\tiltau^{\frac32}-\tau^{\frac32}\right)&= 2\tiltau^{\frac12}\frac{\partial^{2}\tiltau}{\partial\tau^{2}} +\tiltau^{-\frac12} - \tau^{-\frac12} + \tiltau^{-\frac12}\left(\left( \frac{\partial\tiltau}{\partial\tau}\right)^{2}-1\right)
	\end{split}
	\end{align}
	The second term on the right-hand side is 
	\EQ{\label{eq:zweiter T}
		\tiltau^{-\frac12} - \tau^{-\frac12} & = \frac{\tau-\tiltau }{(\tiltau\tau)^{\frac12}(\tiltau^{\frac12} + \tau^{\frac12})}  \\
		&= O\left((\zeta-1)\tau^{-\frac12}(\tau+\zeta-1)^{-1}\right)
	} 
	and the third
	\EQ{\label{eq:dritter T}
		\tiltau^{-\frac12}\left(\left( \frac{\partial\tiltau}{\partial\tau}\right)^{2}-1\right) &= O\left((\zeta-1)(\tau+\zeta-1)^{-\frac12}\right)
	} 
	Taking a derivative of \eqref{eq: tau minus tiltau dtiltau dtau} yields 
	\[
	\frac{\partial^2\tilde\tau}{\partial\tau^2} = (\zeta-1) \partial_\tau\left[ \calR(\zeta,x;\alpha,\hbar)  + \zeta\tau \calS(\zeta,x;\alpha,\hbar) q^{-\frac12} {+(1+O(\hbar))\calT(\zeta,x;\alpha,\hbar)q^{-\frac12}}\right] = O(\zeta-1)
	\]
	The final bound here follows from the expressions for $\calR$ etc.\ which we obtained above, see the calculations leading from~\eqref{eq:RS claim} to~\eqref{eq: tau minus tiltau dtiltau dtau}. In conclusion, \eqref{eq:32 diff 2} holds. 
	We now break up \eqref{eq:F2 pre} into two pieces, the first being  over $\tau\lesssim  \hbar^{\frac23}$ (here $\ell$ stands for {\em left}): 
	\begin{align}\label{eq:F2 pre left}
	\begin{split}
	F_{2\ell}(\xi,\eta;\hbar) 
	& = \int_{0}^{\infty}W^{ {+}}_{n}(R)  \hbar^{\frac13}x^{-\frac12} q^{-\frac14}( \tau)  \Ai(-\hbar^{-\frac23} \tau)(1+\hbar a_0(-\tau;\alpha,\hbar)) \chi_{[\tau\lesssim  \hbar^{\frac23}]} \\
	&\qquad \hbar^{\frac13}\tilde x^{-\frac12} \tilde q^{-\frac14}( \tilde \tau)\Big[  \Ai(-\hbar^{-\frac23}\tilde \tau)(1+\hbar a_0(-\tilde \tau;\beta,\hbar)) \chi_{[\tilde \tau\lesssim  \hbar^{\frac23}]} + \\
	&\qquad + 2  \Re\big( a(\eta) \big(\Ai(-\hbar^{-\frac23}\tilde\tau)-i\Bi(-\hbar^{-\frac23}\tilde\tau)\big)(1+\hbar \overline{a_1(-\tilde\tau;\beta)})  \big) \chi_{[\tilde\tau\gtrsim \hbar^{\frac23}]} \Big]\, \chi(\alpha R )R\, dR
	\end{split}
	\end{align}
	This is estimated using only the size of the integrand: 
	\begin{align}\label{eq:F2 pre left size}
	\begin{split}
	|F_{2\ell}(\xi,\eta;\hbar) |
	& \lesssim  \hbar^{-\frac13}  \alpha^{-2} \langle \alpha^{-1}\rangle^{-4} \int_{-\infty}^{\infty}   \chi_{[-C\le \tau\lesssim  \hbar^{\frac23}]}  \left\langle \hbar^{-\frac23} \tau \right\rangle^{-\frac14} \exp\left(-\frac23\hbar^{-1}|\tau|^{\frac32} \right)  \left\langle \hbar^{-\frac23}\tilde\tau\right\rangle^{-\frac14}  \, d\tau
	\end{split}
	\end{align}
	We do not need to use the rapid decay of  $\Ai(-\hbar^{-\frac23}\tilde \tau)$ for negative $\tilde\tau$ since it does not improve the upper bound.  If $-C\le \tau\le -c$ with some small constant $c>0$, then we obtain a gain of $O(\hbar^\infty)$ from the exponential. If $-c\le\tau\lesssim \hbar^{\frac23}$, then we use~\eqref{eq: tau minus tiltau dtiltau dtau} to write $\tilde\tau=\tau+(\zeta-1)(1+o(1))$. 
	The dominant contribution to~\eqref{eq:F2 pre left size} therefore comes from the  interval $|\tau|\simeq \hbzwd$ whence 
	\EQ{
		\label{eq:F2ell bd}
		|F_{2\ell}(\xi,\eta;\hbar) |
		& \lesssim \alpha^{2}\langle\alpha\rangle^{-4} \min\big(\hbar^{\frac13}, \hbar^{\frac14} \Lambda^{\frac14}\big), 
	}
	where $\Lambda:=\hbar(\zeta-1)^{-1}=\frac{\hbar\xi^{\frac12}}{\eta^{\frac12}-\xi^{\frac12}}$.  By~\eqref{eq: tau minus tiltau dtiltau dtau},  $\tau\gtrsim \hbzwd$ implies the same for $\tilde\tau$.  Hence we  can write $F_{2r}:=F_2-F_{2\ell}$ in the form
	\EQ{
		\nn
		F_{2r}(\xi,\eta;\hbar)  
		& = 4  \hbar^{\frac23} \int_{0}^{\infty}W^{ {+}}_{n}(R)  x^{-\frac12} q^{-\frac14}( \tau) \tilde x^{-\frac12} \tilde q^{-\frac14}( \tilde \tau)  \Re\big( a(\xi) \big(\Ai(-\hbar^{-\frac23}\tau)-i\Bi(-\hbar^{-\frac23}\tau)\big)(1+\hbar \overline{a_1(-\tau;\alpha)})  \big)   \\
		&\qquad   \Re\big( a(\eta) \big(\Ai(-\hbar^{-\frac23}\tilde\tau)-i\Bi(-\hbar^{-\frac23}\tilde\tau)\big)(1+\hbar \overline{a_1(-\tilde\tau;\beta)})  \big) \chi_{[\tau\gtrsim \hbar^{\frac23}]} \, \chi(x)R\,dR
	}
	Integrating by parts using \eqref{eq:32 diff}   and~\eqref{eq:32 diff 2} requires that $\hbzwd \lesssim \tau\ll1$. Thus, we split $F_{2r}$ further by means of a partition of unity:
	\[
	F_{2r} = F_{2r1} + F_{2r2}.
	\]
	Here $F_{2r1}$ covers the integration over $\hbzwd \lesssim \tau\ll 1$, and  $F_{2r2}$ deals with $\tau\simeq1$.  Changing variables via  $dR=\alpha^{-1}q^{-\frac12}\, d\tau$ yields
	\EQ{\nn
		F_{2r1} (\xi,\eta;\hbar) &:=  
		4 \hbar^{\frac23}\alpha^{-2} \int_{-\infty}^{\infty}W^{ {+}}_{n}(R) \chi(x )  \Re\big( a(\xi) \big(\Ai(-\hbar^{-\frac23}\tau)-i\Bi(-\hbar^{-\frac23}\tau)\big)(1+\hbar \overline{a_1(-\tau;\alpha)})  \big)  \\
		&\qquad x^{\frac12} q^{-\frac34}( \tau)  \tilde x^{-\frac12} \tilde q^{-\frac14}( \tilde \tau)  \Re\big( a(\eta) \big(\Ai(-\hbar^{-\frac23}\tilde\tau)-i\Bi(-\hbar^{-\frac23}\tilde\tau)\big)(1+\hbar \overline{a_1(-\tilde\tau;\beta)})  \big) \chi_{[1\gg \tau\gtrsim \hbar^{\frac23}]}\,  d\tau
	}
	The dominant contributions are made by the resonant terms due to phases exhibiting opposite signs. Without loss of generality it suffices to bound 
	\EQ{\label{eq:Fs osc+-}
		F_{2r1}^{+-}(\xi,\eta;\hbar) &:=  
		{\hbar}\alpha^{-2} \int_{-\infty}^{\infty}W^{ {+}}_{n}(R)   a(\xi)\overline{a(\eta)}   e^{\frac{2i}{3\hbar} (\tilde \tau^{\frac32} -\tau^{\frac32})} (1+\hbar \overline{\tilde a_1(-\tau;\alpha)})  \big)  \\
		&\qquad x^{\frac12} q^{-\frac34}( \tau)  \tilde x^{-\frac12} \tilde q^{-\frac14}( \tilde \tau) {\tau^{-\frac14}\tiltau^{-\frac14}}(1+\hbar {\tilde a_1(-\tilde\tau;\beta)})  \big) \chi_{[1\gg \tau\gtrsim \hbar^{\frac23}]}\,   d\tau\\
		& = \frac32 i\hbar^{2}\alpha^{-2}  a(\xi)\overline{a(\eta)}\int_{-\infty}^{\infty}  e^{\frac{2i}{3\hbar} (\tilde \tau^{\frac32} -\tau^{\frac32})}\partial_\tau\Big\{ [ \partial_\tau \big( \tilde\tau^{\frac32}-\tau^{\frac32} \big)]^{-1} \Big[ W^{ {+}}_{n}(R)    (1+\hbar \overline{\tilde a_1(-\tau;\alpha)})   \\
		&\qquad x^{\frac12} q^{-\frac34}( \tau)  \tilde x^{-\frac12} \tilde q^{-\frac14}( \tilde \tau) {\tau^{-\frac14}\tiltau^{-\frac14}}(1+\hbar {\tilde a_1(-\tilde\tau;\beta)}) \chi_{[1\gg \tau\gtrsim \hbar^{\frac23}]}\Big]\Big\} \,   d\tau
	}
	The second expression, which was obtained by integration by parts, is only useful if $\xi$ and $\eta$ are not too close. We estimate the first term on the right-hand side of~\eqref{eq:Fs osc+-} by passing absolute values inside yielding 
	\begin{align}\label{eq:F2 osc 1}
	|F_{2r1}^{+-}(\xi,\eta;\hbar)|  &\lesssim \alpha^2\langle\alpha\rangle^{-4} 
	\end{align}
	To obtain the refined bounds in terms of $\Lambda$ provided $\Lambda\ll1$, we further split the interval $\hbar^{\frac23}\ll\tau\ll1$ according to whether $\tau^{\frac12}\lesssim \Lambda$ or $\tau^{\frac12}\gtrsim\Lambda$.  In other words, using a smooth partition of unity as before, we split one more time
	\[
	F_{2r1}^{+-}(\xi,\eta;\hbar) = F_{2r11}^{+-}(\xi,\eta;\hbar) + F_{2r12}^{+-}(\xi,\eta;\hbar)
	\]
	The first term $F_{2r11}^{+-}(\xi,\eta;\hbar)$ corresponds to $\hbar^{\frac23}\ll\tau\lesssim \Lambda^{2}$, and we bound it 
	by placing absolute values inside and using that $\tilde\tau^{-\frac14}\le\tau^{-\frac14}$, which yields
	\begin{align}
	\label{eq:F2r1 e}
	\left|F_{2r11}^{+-}(\xi,\eta;\hbar)\right|\lesssim \alpha^{2}\langle\alpha\rangle^{-4}\Lambda,\quad \textrm{for}\quad \Lambda\ll1.
	\end{align}
	The second term  $F_{2r12}^{+-}$ lives on $\Lambda^2\lesssim \tau \ll1$, and we integrate by parts as in the second line of \eqref{eq:Fs osc+-} but with $\chi_{[1\gg \tau\gg \hbar^{\frac23}]}$ replaced by~$\chi_{[1\gg \tau\gtrsim \Lambda^2]}$. 
	By \eqref{eq:32 diff}   and~\eqref{eq:32 diff 2} we have
	\[
	\frac{|\partial_{\tau}^{2}(\tiltau^{\frac32}-\tau^{\frac32})|}{|\partial_{\tau}(\tiltau^{\frac32}-\tau^{\frac32})|^2}  \lesssim (\zeta-1)^{-1} \tau^{-\frac12}
	\]
	Placing absolute values inside the integral in the second expression in~\eqref{eq:Fs osc+-} produces the bound
	\begin{align*}
	\left| F_{2r12}^{+-}(\xi,\eta;\hbar)\right|  & \lesssim \Lambda   \alpha^2\langle\alpha\rangle^{-4}  \Big\{ \int_{-\infty}^\infty    \tau^{-\frac34} (\tau+\zeta-1)^{-\frac14} \chi_{[1\gg \tau\gg \Lambda^{2}]} \,   d\tau \\
	& \qquad + \int_{-\infty}^\infty (\tau+\zeta-1)^{\frac14} \big (\hbar \tau^{-\frac34}  + \tau^{-\frac54} \big)  \chi_{[1\gg \tau\gg \Lambda^{2}]} \,   d\tau \Big\}.
	\end{align*}
	On the one hand,  $$\int \tau^{-\frac34} (\tau+\zeta-1)^{-\frac14} \chi_{[1\gg \tau\gg \Lambda^{2}]}\,d\tau \lesssim \int \tau^{-1}\chi_{[1\gg \tau\gg \Lambda^{2}]}\, d\tau \lesssim |\log\Lambda|. $$ 
	Note that without the condition $\tau\gg \Lambda^{2}$ we would have obtained $|\log(\zeta-1)|$ resulting in a loss of $\log\hbar$ which is inadmissible. 
	On the other hand,  $$\int (\tau+\zeta-1)^{\frac14} \tau^{-\frac54}\chi_{[1\gg \tau\gg \Lambda^{2}]}\,d\tau  \lesssim \int (\tau^{-1}+(\zeta-1)^{\frac14}\tau^{-\frac54})\chi_{[1\gg \tau\gg \Lambda^{2}]}
	\, d\tau \lesssim |\log\Lambda|+\hbar^{\frac14}\Lambda^{-\frac34},$$ 
	and 
	\[
	\int_{-\infty}^\infty (\tau+\zeta-1)^{\frac14} \hbar\tau^{-\frac34} \chi_{[1\gg\tau\gg\Lambda^{2}]} \, d\tau\lesssim \hbar 
	\]
	In combination with \eqref{eq:F2 osc 1} these estimates imply 
	\EQ{\nn 
		\left| F_{2r12}^{+-}(\xi,\eta;\hbar)\right| & \lesssim   \alpha^2\langle\alpha\rangle^{-4}  \Lambda \left[ \hbar + |\log\Lambda| + \hbar^{\frac14} \Lambda^{-\frac34} \right]  \lesssim   \alpha^2\langle\alpha\rangle^{-4}  \Lambda^{\frac14}
	}
	and in summary, 
	\EQ{\label{eq:F2r1 h}
		|F_{2r1} (\xi,\eta;\hbar)|&\lesssim   \alpha^2\langle\alpha\rangle^{-4} \min (1, \Lambda^{\frac14} ).
	}

Next we turn to $F_{2r2}$ which covers the integration over $\tau\simeq1$. We will proceed as in Case~$A_3$  above. However, since $x$ is no longer large, we need to add the following observation concerning the lower bound on the derivative of the phase in~\eqref{eq:delx Phi-}. 
By \eqref{eq:Q0darst} and \eqref{eq:rho def}, 
\[
y\rho_y = - \frac{1}{y+\sqrt{y^2-1}} + O(\hbar), \quad (y\rho_y)_y = \frac{1}{(y+\sqrt{y^2-1})\sqrt{y^2-1}} + O(\hbar)
\]
Note that the derivative on the right-hand side is large if $y$ is close to $1$, but it is positive. On the other hand, by Lemma~\ref{lem: Lemma 3.4 CDST}, $\alpha\rho_{\alpha y}=O(\hbar)$. This shows that \eqref{eq:delx Phi-} continues to hold in the form $\partial_x \Phi_-\simeq \zeta-1$. The same analysis we used above to bound~\eqref{eq:F3} thus still goes through whence
\EQ{\label{eq:F2r2 bd}
	|F_{2r2}(\xi,\eta;\hbar)| &\lesssim \alpha^{-1}\min(1,\alpha^3)\min\big(1,\Lambda^5 \big)
}
This completes the discussion for {\bf case $\mathbf{A_2}$}.
 
 \medskip

It remains to deal with {\bf case $\mathbf{A_1}$}, namely, $0<x\ll1$. Lemmas~\ref{lem: Lemma 3.3 CDST}, \ref{lem:xsimal}, \ref{lem: Lemma 3.3 harder} and \ref{lem: Lemma 3.3 hardest} imply that, for sufficiently small $\hbar$,
\begin{align*}
\left(-\tau\right)^{\frac32}\simeq -\log x
\end{align*}
uniformly in the parameters. 
By the definition of $q$ and \eqref{eq:Q2 form}, for $\tau\ll-1$ (i.e., $0<x \leq \frac12 x_{t}$), we have
\begin{align*}
q=-\frac{Q_0}{\tau}\simeq \frac{1}{x^{2}\tau}.
\end{align*}
This together with the asymptotic profile for the Airy function gives
\begin{align}\label{eq:F2r2 integrand}
\left|\hbar^{\frac13}x^{-\frac12}q^{-\frac14}(\tau)\Ai\left(-\hbar^{-\frac23}\tau\right)\left(1+\hbar a_{0}(-\tau,\alpha;\hbar)\right)\right|\lesssim  \hbar^{\frac12} e^{-\frac23\hbar^{-1}(-\tau)^{\frac32}}\lesssim \hbar^{\frac12}\cdot x^{c\hbar^{-1}} 
\end{align}
Since $\tilde x=x\zeta\simeq x$,  we have $\tilde{\tau}\ll-1$ and 
\begin{align*}
\left|\hbar^{\frac13}\tilde{x}^{-\frac12}q^{-\frac14}(\tiltau)\Ai\left(-\hbar^{-\frac23}\tiltau\right)\left(1+\hbar a_{0}(-\tiltau,\beta;\hbar)\right)\right|\lesssim  
\hbar^{\frac12}\cdot x^{c\hbar^{-1}}
\end{align*}
Therefore, with $p=c\hbar^{-1}$ and a fixed small constant $c>0$, 
\begin{align}\label{eq:F2r2}
\begin{split}
\Big|&\int_{0}^{\infty}\chi_{[0<x\ll 1]}\cdot W_{n}^{+}(R)\hbar^{\frac13}x^{-\frac12}q^{-\frac14}(\tau)\Ai(-\hbar^{-\frac23}\tau)(1+\hbar a_{0}(-\tau,\alpha;\hbar))\\
&\cdot\hbar^{\frac13}\tilde{x}^{-\frac12}q^{-\frac14}(\tilde{\tau})\Ai(-\hbar^{-\frac23}\tiltau)(1+\hbar a_{0}(-\tilde{\tau},\beta;\hbar))R\,dR\Big|\\
&\lesssim \alpha^{-2}\int_0^{\frac12} \langle x/\alpha\rangle^{-4} x^{p}\, dx = \alpha^{-1+p}\int_0^{(2\alpha)^{-1}} \langle y\rangle^{-4} y^p\, dy\\
&\lesssim \hbar \alpha^{2}\langle\alpha\rangle^{-4} \lesssim \alpha^{2}\langle\alpha\rangle^{-4}\min\left(1,\Lambda\right),
\end{split}
\end{align}
since $\Lambda\gtrsim \hbar$. Adding this estimate to the sum of the prior cases, i.e., \eqref{eq:F3 summary}, \eqref{eq:F2ell bd}, \eqref{eq:F2r1 h},  \eqref{eq:F2r2 bd}, establishes the bound stated in~\eqref{F hbar bounds}. The bound in the last line of \eqref{eq:F2r2} is good enough for our purposes, but crude. In fact, we gain an exponential factor of the form~$2^{-p}=e^{-c\hbar^{-1}}$.

\bigskip

Now we turn to {\bf Case $\mathbf{B}$} where ${\xi\ll\eta}$.   As before, we break the integral in~\eqref{off diagonal F pre} into three separate regions, namely $x\gg1$ (case $B_3$), followed by $x\simeq 1$ (case $B_2$), and finally $0<x\ll 1$ (case $B_1$).  Beginning with $\mathbf{B_3}$, we need to estimate the expression $F_3(\xi,\eta;\hbar)$ defined by~\eqref{eq:F3}. 
The bound~\eqref{eq:F3 0} is replaced by 
\[
|F_3(\xi,\eta;\hbar)|\lesssim \zeta^{-\frac12} \alpha^{-1}\min (1,\alpha^3)
\]
In Case~B, we have $\zeta-1\simeq\zeta= (\eta/\xi)^{\frac12}$ which implies that  the phases in \eqref{eq:F3 phases} satisfy
\[
|\partial_x \Phi_{\pm}(x,\alpha;\hbar,\zeta)|\simeq \zeta
\]
uniformly in the parameters.  Integrating by parts as in the calculation leading up to~\eqref{eq:F3 ell} yields 
\EQ{\label{eq:F3 B}
	|F_3(\xi,\eta;\hbar)|  \le C_\ell\, \zeta^{-\frac12} (\hbar/ \zeta)^{\ell}  \alpha^2\langle \alpha\rangle^{-4}
}
for any $\ell\ge1$. 
In {\bf Case $\mathbf{B_2}$}, we have $x\simeq 1$, i.e., $|\tau|\lesssim 1$ as well as $\tilde x\gg1$ and $\tilde \tau \gg1$. Thus, \eqref{eq:F2 pre} now takes the form 
	\begin{align}\label{eq:F2 pre B}
	\begin{split}
	F_2(\xi,\eta;\hbar) &:=  \int_{0}^{\infty}W^{ {+}}_{n}(R)\phi(R,\xi;\hbar)\phi(R,\eta;\hbar)  \chi(x) R\, dR \\
	& = 2\alpha^{-2}\int_{0}^{\infty}W^{ {+}}_{n}(R) \chi(\alpha R )  \hbar^{\frac13}x^{\frac12} q^{-\frac14}( \tau)\Big[  \Ai(-\hbar^{-\frac23} \tau)(1+\hbar a_0(-\tau;\alpha,\hbar)) \chi_{[-1\lesssim \tau\lesssim  \hbar^{\frac23}]} + \\
	&\quad + 2  \Re\big( a(\xi) \big(\Ai(-\hbar^{-\frac23}\tau)-i\Bi(-\hbar^{-\frac23}\tau)\big)(1+\hbar \overline{a_1(-\tau;\alpha)})  \big) \chi_{[1\gtrsim \tau\gtrsim  \hbar^{\frac23}]}\Big] \\
	&\qquad\qquad \hbar^{\frac13}\tilde x^{-\frac12} \tilde q^{-\frac14}( \tilde \tau)  \Re\big( a(\eta) \big(\Ai(-\hbar^{-\frac23}\tilde\tau)-i\Bi(-\hbar^{-\frac23}\tilde\tau)\big)(1+\hbar \overline{a_1(-\tilde\tau;\beta)})\big)  \, dx \\
	&:= F_{2\ell}(\xi,\eta;\hbar) + F_{2r}(\xi,\eta;\hbar)
	\end{split}
	\end{align}
where the summands in the last line correspond to the respective summands inside the brackets. Below we will simplify~\eqref{eq:F2 pre B} using  that $x^{\frac12} \tilde x^{-\frac12} = \zeta^{-\frac12}$, see~\eqref{eq:F2ell0 B}, \eqref{eq:F2ell1 B}.  The factors $a(\xi), a(\eta)$ are uniformly bounded and do not affect the bounds. 
By Lemma~\ref{lem: Lemma 3.4 CDST} we have
\EQ{\label{eq:Psi xalp}
	\Psi=\Psi(x;\alpha,\zeta,\hbar):=\frac23\tiltau^{\frac32}=\tilde{x}-y(\beta;\hbar)+\rho(\tilde{x},\beta;\hbar) = \zeta x-y(\beta;\hbar)+\rho(\zeta x,\beta;\hbar) \simeq \zeta x
}
and the oscillatory Airy functions are as follows, see~\eqref{eq:AiiBi}: 
\begin{align*}
\Ai(-\hbar^{-\frac23}\tiltau)\mp \,i\Bi(-\hbar^{-\frac23}\tiltau)=c_{\pm}(\hbar^{-\frac23}\tiltau)^{-\frac14}\cdot e^{\pm\,\frac{i}{\hbar} \Psi}\cdot(1+b(\hbar^{-\frac23}\tiltau)).
\end{align*}
By \eqref{eq:Ai}, $b(\hbar^{-\frac23}\tiltau) = \hbar \tilde b(\tilde x,\beta;\hbar)$ with $\partial_{\tilde x}^k \tilde b(\tilde x,\beta;\hbar) = O(\tilde x^{-k-1})$. We can thus absorb the factor $(1+b(\hbar^{-\frac23}\tiltau))$ into the factor $(1+\hbar \overline{a_1(\tilde\tau;\beta)})$ above. We shall carry out this step without further mention, also for the  non-oscillatory Airy function. 
Integration by parts in \eqref{eq:F2 pre B} is carried out by means of the identity 
\EQ{\label{eq:phase large tiltau}
	\calL\;  e^{\pm\,\frac{i}{\hbar} \Psi} &= \pm e^{\pm\,\frac{i}{\hbar} \Psi} \\
	\calL &:= \hbar \big[ i\,\partial_{x}\Psi(x;\alpha,\zeta,\hbar)\big]^{-1}\partial_x
}
Since $\tilde{x}\gg1$, the derivative of the phase satisfies
\EQ{\label{eq:phase deri large tiltau}
	\partial_{x}\Psi &=\zeta(1+\rho_{\tilde{x}}(\tilde{x},\beta;\hbar))=\zeta (1+O(\tilde{x}^{-2}))\simeq \zeta \\
	\partial_{x}^k\Psi &= \zeta^k O(\tilde x^{-1-k})\simeq_k \zeta^{-1}, \quad k\ge2
}
where the final $\simeq$ holds due to $\tilde x\simeq \zeta$. 
From \eqref{eq:qdef}, 
\begin{align*}
\hbar^{\frac13}\tilde{x}^{-\frac12}q^{-\frac14}(\tiltau)(\hbar^{-\frac23}\tilde{\tau})^{-\frac14}=\hbar^{\frac12}\tilde{x}^{-\frac12}(-Q_{0}(\tilde x,\beta;\hbar))^{-\frac14}, 
\end{align*}
with Corollary \ref{cor:Q0diff} providing the bounds 
\EQ{\label{eq:Q0 derxk}
	-Q_{0}(\tilde{x},\alpha;\hbar)\simeq 1,\quad \left|\partial_{\tilde x}^{k}Q_{0}(\tilde{x},\alpha;\hbar)\right|\lesssim \tilde{x}^{-k-2}\simeq\zeta^{-k-2},\; k\ge1.
}
We break up $F_{2\ell}=F_{2\ell0}+F_{2\ell1}$ further, by means of a smooth partition of unity adapted to $|x-x_t(\alpha,\hbar)|\lesssim \hbar^{\frac23}$, respectively $-1\lesssim x-x_t\lesssim - \hbar^{\frac23}$.  Note that by Lemma~\ref{lem: Lemma 3.2 CDST} we may interchange $\tau$ with $x-x_t(\alpha)$ in these conditions. 
Thus, up to a constant multiplicative factor which we ignore, 
\EQ{
		\label{eq:F2ell0 B}
		F_{2\ell0} (\xi,\eta;\hbar)& = \hbar^{\frac56} \alpha^{-2}a(\eta)\zeta^{-\frac12} \int_{0}^{\infty}  e^{\frac{i}{\hbar}\Psi}\; (\calL^*)^m\Big[ W^{+}_{n}(x/\alpha)      q^{-\frac14}( \tau)  \Ai(-\hbar^{-\frac23} \tau)(1+\hbar a_0(-\tau;\alpha,\hbar))  \\
		&\qquad\qquad\qquad   (-Q_{0}(\tilde x,\beta;\hbar))^{-\frac14} (1+\hbar \overline{\tilde a_1(-\tilde\tau;\beta)}) \chi_{[ |x-x_t(\alpha)|\lesssim  \hbar^{\frac23}]}  \Big] \, dx + \text{cc}
	}
where ``cc" stands for {\em complex conjugate}.  Similarly, 
\EQ{
		\label{eq:F2ell1 B}
		F_{2\ell1} (\xi,\eta;\hbar)& = \hbar \alpha^{-2}a(\eta)\zeta^{-\frac12}\int_{0}^{\infty}  e^{\frac{i}{\hbar}\Psi}\; (\calL^*)^m\Big[ W^{+}_{n}(x/\alpha)     q^{-\frac14}( \tau) (-\tau)^{-\frac14}  e^{-\frac{2}{3\hbar}(-\tau)^{\frac32}}(1+\hbar \tilde a_0(-\tau;\alpha,\hbar))   \\
		&\qquad\qquad\qquad    (-Q_{0}(\tilde x,\beta;\hbar))^{-\frac14} (1+\hbar \overline{\tilde a_1(-\tilde\tau;\beta)}) \chi_{[-1\lesssim x -x_t(\alpha)\lesssim  -\hbar^{\frac23}]} \Big] \, dx + \text{cc.}
	} 
We now analyze the  contribution of each factor as it arises by integrating by parts in $F_{2\ell 0}$, see~\eqref{eq:F2 pre B}. First, it follows by induction that 
\EQ{\label{eq:SOT} 
	(\calL^*)^m  &:= -\hbar^m \Big( \partial_x \big[ i\,\partial_{x}\Psi(x;\alpha,\zeta,\hbar)\big]^{-1}\Big)^m \\
	& =  \hbar^m \sum \text{coeff.} \frac{\partial_x^{j_1+1} \Psi \partial_x^{j_2+1} \Psi \cdot\ldots \cdot \partial_x^{j_\ell+1} \Psi }{(\partial_x \Psi )^{m+\ell}} \partial_x^k 
}
where the sum runs over those terms which obey $j_1+j_2+\ldots+j_\ell + k =m$, $0\le\ell\le m$, $\min(j_1,j_2,\ldots, j_\ell)\ge1$. The coefficients are some absolute constants. In view of~\eqref{eq:phase deri large tiltau}, the factors in the sum are of size
\EQ{\label{eq:of size}
	\frac{\partial_x^{j_1+1} \Psi \partial_x^{j_2+1} \Psi \cdot\ldots \cdot \partial_x^{j_\ell+1} \Psi }{(\partial_x \Psi )^{m+\ell}} &\simeq \frac{\zeta^{-(j_1+\ldots+j_\ell)-\ell}}{\zeta^{m+\ell}} 
	\simeq \zeta^{k-2(m+\ell)} \lesssim \zeta^{-m-3\ell}\lesssim \zeta^{-m} 
}
since $\ell+k\le m$. Thus, we gain $(\hbar/\zeta)^m$ in~\eqref{eq:F2ell0 B} and~\eqref{eq:F2ell0 B}, and it remains to estimate the effect of $\partial_x^k$, $0\le k\le m$, in these integrals. We now analyze $\partial_x^k$ derivatives of each of the factors in~\eqref{eq:F2ell0 B} for all $k\ge1$, and uniformly in the parameters such as $x\simeq1$:  
\begin{itemize}
		\item $\alpha^{-2}|\partial_x^k (W^{+}_{n}(x/\alpha))|\le C_k\hbar^{-1} \alpha^2 \langle \alpha\rangle^{-4}$, as one checks by differentiation of~\eqref{W pre}. 
		\item By Lemma~\ref{lem: Lemma 3.4 CDST}, $|\partial_x^k q^{-\frac14}( \tau)|\le C_k $.
		\item $|\partial_x^k \big( \Ai(-\hbar^{-\frac23} \tau) \chi_{[ |x-x_t(\alpha)|\lesssim  \hbar^{\frac23}]} \big) |\le C_k \hbar^{-\frac{2k}{3}}$, again using Lemma~\ref{lem: Lemma 3.4 CDST}. 
		\item By Lemmas~\ref{lem:a0fine1}, \ref{lem:a0fine2}, $|\partial_x^k (1+\hbar a_0(-\tau;\alpha,\hbar))|\le C_k \,\hbar^{\frac{4-2k}{3}}$ for $|\tau|\lesssim\hbar^{\frac23}$.
		\item By \eqref{eq:Ai}, $|\partial_x^k (1+\hbar \tilde a_0(-\tau;\alpha,\hbar))|\le C_k \,\hbar|\tau|^{-\frac32-k}$ for $\hbar^{\frac23}\ll|\tau|\lesssim 1$.
		\item By \eqref{eq:Q0 derxk} and $\frac{\partial\tilde x}{\partial x}=\zeta$, one has $|\partial_x^k  (-Q_{0}(\tilde x,\beta;\hbar))^{-\frac14}|\le C_k\, \zeta^{-2}$.
		\item By Lemma~\ref{lem:a1fine}, see~\eqref{eq:a1 fine}, and $\frac{\partial\tilde \tau}{\partial \tilde x}=(-Q_0(\tilde x)/{\tilde \tau})^{\frac12}$, one checks that $|\partial_x^k (1+\hbar \tilde a_1(-\tilde\tau;\beta,\hbar))|\le C_k \,\hbar\zeta^{-1}$. 
	\end{itemize}
We remark that $\zeta=-\tau$ in Lemmas~\ref{lem:a1fine},  \ref{lem:a0fine1}, \ref{lem:a0fine2}, and elsewhere in that section has nothing to do with $\zeta=(\eta/\xi)^{\frac12}$ as it appears here.
In summary, these bounds imply that for any $m\ge0$, 
\EQ{\label{eq:F2ell0 B bd}
	|F_{2\ell0} (\xi,\eta;\hbar)|&\le C_m\, \alpha^2 \langle \alpha\rangle^{-4} \, (\hbar/\zeta)^{\frac12} (\hbar^{\frac13}/\zeta)^m
}
The factor $\hbar^\frac12$ arises as $\hbar^{\frac56}\hbar^{-1} \hbar^{\frac23}$, the latter being the length of the integration interval. 
This same estimate also holds for $F_{2\ell1} (\xi,\eta;\hbar)$ as can be seen by a dyadic decomposition $-\tau\simeq 2^j\hbar^{\frac23}$. Indeed, one then has
\[
\hbar (-\tau)^{-\frac14} \exp\left(-\frac{2}{3\hbar}(-\tau)^{\frac32}\right)\lesssim \hbar^{\frac56} 2^{-\frac{j}{4}} e^{-\frac23 2^{3j/2}}
\]
which is rapidly decaying in $j$. We leave the remaining details to the reader. Thus, $F_{2\ell}$ also satisfies the bound~\eqref{eq:F2ell0 B bd}. 

The analysis of $F_{2r}(\xi,\eta;\hbar)$ is similar. Indeed, up to uniformly bounded factors, this term is the sum of the following four integrals:
\EQ{\nn
	F^{\pm}_{\pm}:=\alpha^{-2}\zeta^{-\frac12} \hbar \int_0^\infty &W_n^+(x/\alpha)(\tau q(\tau))^{-\frac14} (-Q_0(\tilde x;\beta,\hbar))^{-\frac14} \exp\left( \frac{2i}{3\hbar} (\pm\tau^{\frac32} \pm \tilde\tau^{\frac32})\right) \chi_{[\hbar^{\frac23}\lesssim x-x_t(\alpha)\lesssim 1]}\\
	&\cdot \left(1+\hbar \tilde{a}_{1}^{\pm}(-\tau;\alpha,\hbar)\right)\left(1+\hbar \tilde a_{1}^{\pm}(-\tiltau;\beta,\hbar)\right)\, dx
}
We write the complex exponential as $e^{\frac{i}{\hbar}\Phi^{\pm}_{\pm}}$ with phases 
\[
\Phi^{\pm}_\pm(x;\alpha,\zeta,\hbar) = \frac23\left( \pm \tau^{\frac32} \pm \tilde \tau^{\frac32} \right) 
\]
and $a_{1}^{+}:=a_{1},\quad a_{1}^{-}:=\overline{ a_{1}},\quad  \tilde a_{1}^{+}:=\tilde a_{1},\quad \tilde a_{1}^{-}:=\overline{ \tilde a_{1}}$.

It suffices to consider one choice of signs here due to the fact that $\tilde\tau\gg1$ and $0<\tau\lesssim1$, and we write $\Phi= \Phi^+_+$.   Then 
\EQ{\label{eq:Phi deri B}
	\partial_x \Phi &=  \tau^{\frac12} q^{\frac12} + \partial_x \Psi\simeq \zeta \\
	\partial_x^2 \Phi & =  \frac12 \tau^{-\frac12} q + \frac12\tau^{\frac12}  \partial_\tau q+ \partial_x^2 \Psi \simeq\tau^{-\frac12}\\
	|\partial_x^k \Phi| & \simeq_k \tau^{\frac32-k},\quad k\ge2
}
see~\eqref{eq:phase deri large tiltau}.  We redefine $\calL$ as in \eqref{eq:phase large tiltau} but with $\Phi$ in place of~$\Psi$. In analogy with~\eqref{eq:F2ell1 B} we now have
\begin{align*}
F^+_+ = \hbar\alpha^{-2}\zeta^{-\frac12}  \int_0^\infty e^{\frac{i}{\hbar}\Psi} (\calL^*)^m&\big[ W_n^+(x/\alpha)(\tau q(\tau))^{-\frac14} (-Q_0(\tilde x;\beta,\hbar))^{-\frac14}  \chi_{[\hbar^{\frac23}\lesssim x-x_t(\alpha)\lesssim 1]}\\
&\cdot\left(1+\hbar \tilde a_{1}(-\tau;\alpha,\hbar)\right)\left(1+\hbar \tilde a_{1}(-\tiltau;\beta,\hbar)\right)\big]\, dx
\end{align*}
where $(\calL^*)^m$ is given by \eqref{eq:SOT} but with $\Phi$ in place of $\Psi$.  In view of the preceding bullet list, and with $m\ge2$ and $k,\ell$ and the sum as in~\eqref{eq:SOT}, we have 
\EQ{\nn
	|F_{2r}(\xi,\eta;\hbar)| &\lesssim \hbar^{m+1} \alpha^{-2} \zeta^{-\frac12}  \int_0^\infty \sum \frac{\prod_{i=1}^\ell\tau^{ \frac12-j_i}}{\zeta^{m+\ell}} \Big|\partial_x^k \big[ W_n^+(x/\alpha)(-\tau q(\tau)Q_0(\tilde x;\beta))^{-\frac14}  \chi_{[\hbar^{\frac23}\lesssim x-x_t(\alpha)\lesssim 1]}\\
	&\qquad\qquad\qquad\qquad \left(1+\hbar \tilde a_{1}(-\tau;\alpha,\hbar)\right)\left(1+\hbar \tilde a_{1}(-\tiltau;\beta,\hbar)\right)\big]\Big| \, dx \\
	&\lesssim  \hbar^{m} \alpha^2\langle \alpha\rangle^{-4}  \sum \zeta^{-\frac12-m-\ell} \int_{\hbar^{\frac23}}^1 \tau^{\frac{\ell}{2}+k-m} \left(\tau^{-\frac14 -k}+\hbar\tau^{-\frac32-k}\right)\, d\tau \\
	&\lesssim \alpha^2\langle \alpha\rangle^{-4}  \sum  \hbar^{m} \zeta^{-\frac12-m-\ell} (\hbar^{\frac23})^{\frac{\ell}{2}-m+\frac34} \lesssim \alpha^2\langle \alpha\rangle^{-4} (\hbar/\zeta)^{\frac12} \sum_{\ell\ge0}  (\hbar^{\frac13}/\zeta)^{m+\ell}\\
	&\lesssim \alpha^2\langle \alpha\rangle^{-4} (\hbar/\zeta)^{\frac12}  (\hbar^{\frac13}/\zeta)^{m}
}
To pass to the third line we used that $\frac{\ell}{2}-m-\frac14<-1$ due to $\ell\le m$ and $m\ge2$. The final estimate exactly matches~\eqref{eq:F2ell0 B bd} whence 
\EQ{\label{eq:F2 B}
	|F_{2} (\xi,\eta;\hbar)|&\le C_m\, \alpha^2 \langle \alpha\rangle^{-4} \, (\hbar/\zeta)^{\frac12} (\hbar^{\frac13}/\zeta)^m,\quad \forall\; m\ge0
}

It remains to analyze the contribution of Case~$B_1$. Let $\chi(x)$ be a smooth cutoff which equals~$1$ on $0<x\leq \frac14$ and is supported on $0<x\leq \frac12$.  
	We split the integral 
	\begin{align}
	F_1(\xi,\eta;\hbar) &:=  \int_{0}^{\infty}W^{ {+}}_{n}(R)\phi(R,\xi;\hbar)\phi(R,\eta;\hbar)  \chi(x) R\, dR=\alpha^{-2}\int_{0}^{\infty}W^{ {+}}_{n}(x/\alpha)\phi(R,\xi;\hbar)\phi(R,\eta;\hbar)  \chi(x) x\, dx \nn \\
	&= \sum_{j=1}^3 \alpha^{-2}\int_{0}^{\infty} W^{ {+}}_{n}(x/\alpha) \hbar^{\frac13} x^{\frac12} q(\tau;\hbar)^{-\frac14} \Ai(-\hbar^{-\frac23}\tau)(1+\hbar a_0(-\tau;\alpha,\hbar))\phi(R,\eta;\hbar)  \chi_j(\tilde \tau;\hbar) \, dx \nn \\
	& = (F_{11} + F_{12} + F_{13})(\xi,\eta;\hbar) \label{eq:F1 pre B}
	\end{align}
	where $\sum_{j=1}^3 \chi_j(\tilde \tau;\hbar)=\chi(x)$ is a smooth partition of unity corresponding to $\tilde\tau\lesssim -1$, $|\tilde\tau|\lesssim 1$, $\tilde\tau\gtrsim1$, respectively. Then, on the one hand, by \eqref{eq:xtau} and with $N=(2\hbar)^{-1}$, 
	\EQ{\label{eq:F11 bd B}
		|F_{11}(\xi,\eta;\hbar)| &\lesssim  \alpha^{-2}\zeta^{-\frac12} \int_{0}^{\infty} \langle x/\alpha\rangle^{-4}  \left(Q_0(x;\alpha,\hbar)Q_0(\tilde x;\beta,\hbar)\right)^{-\frac14} e^{ - \frac{2}{3\hbar} ((-\tau)^{\frac32}+(-\tilde\tau)^{\frac32})}\\
		&\qquad \qquad\qquad |(1+\hbar \tilde a_0(-\tau;\alpha,\hbar))(1+\hbar \tilde a_0(-\tilde\tau;\beta,\hbar))|  \chi_1(\tilde \tau;\hbar) \, dx \\
		&\lesssim \alpha^{-2}\zeta^{-\frac12} \int_{0}^{\frac{1}{2\zeta}} \langle x/\alpha\rangle^{-4}  (x\tilde x)^{\frac12} (x\tilde x)^{N} \, dx 
		\lesssim \alpha^{-2}\zeta^{N} \int_{0}^{\frac{1}{2\zeta}} \langle x/\alpha\rangle^{-4}  x^{2N+1}   \, dx \\
		&\lesssim \hbar \alpha^{-2}\zeta^{-N-2} \min(1,(\alpha\zeta)^4)\\
		&\lesssim \hbar\alpha^{2}\langle\alpha\rangle^{-4}\zeta^{-N+2}
	}
	and, on the other hand, 
	\EQ{\label{eq:F13 bd B}
		F_{13}(\xi,\eta;\hbar) &= 2\alpha^{-2}\int_{0}^{\infty} W^{ {+}}_{n}(x/\alpha) \hbar^{\frac13} x^{\frac12} q(\tau;\hbar)^{-\frac14} \Ai(-\hbar^{-\frac23}\tau)(1+\hbar a_0(-\tau;\alpha,\hbar))  \\
		&\qquad \hbar^{\frac13} \tilde x^{-\frac12} \tilde q(\tilde \tau;\hbar)^{-\frac14} \Re\Big(a(\eta) (\Ai(-\hbar^{-\frac23}\tilde \tau)-i \Bi(-\hbar^{-\frac23}\tilde \tau)) (1+\hbar \overline{a_1}(-\tilde\tau;\alpha,\hbar)) \Big) \chi_3(\tilde \tau;\hbar) \, dx  \\
		& = 2\hbar \, \alpha^{-2}\zeta^{-\frac12} \int_{0}^{\infty} W^{ {+}}_{n}(x/\alpha)  \left(Q_0(x;\alpha,\hbar)\right)^{-\frac14}\left(-Q_0(\tilde x;\beta,\hbar)\right)^{-\frac14}  e^{ -\frac{2}{3\hbar} (-\tau)^{\frac32}}  (1+\hbar \tilde a_0(-\tau;\alpha,\hbar))  \\
		&\qquad \Re\Big(a(\eta) e^{\frac{i}{\hbar}\Psi}   (1+\hbar\overline{\tilde a_1}(-\tilde\tau;\beta,\hbar)) \Big)  \chi_{[\zeta^{-1}\ll x\ll 1]} \, dx 
	}
	where we wrote the smooth cutoff function $\chi_3(\tilde \tau;\hbar)=\chi_{[\zeta^{-1}\ll x\ll 1]}$, and $\Psi$ is defined in \eqref{eq:Psi xalp}.  Integrating by parts by means of $\calL$ defined in~\eqref{eq:phase large tiltau} we obtain
	\EQ{
		\label{eq:F13 B}
		F_{13}(\xi,\eta;\hbar)& = \hbar \alpha^{-2}a(\eta)\zeta^{-\frac12}\int_{0}^{\infty}  e^{\frac{i}{\hbar}\Psi}\; (\calL^*)^m\Big[ W^{+}_{n}(x/\alpha)   \left(Q_0(x;\alpha,\hbar)\right)^{-\frac14}\left(-Q_0(\tilde x;\beta,\hbar)\right)^{-\frac14} e^{ -\frac{2}{3\hbar} (-\tau)^{\frac32}}\\
		&\qquad  \qquad\qquad  (1+\hbar \tilde a_0(-\tau;\alpha,\hbar))     (1+\hbar\overline{\tilde a_1}(-\tilde\tau;\beta,\hbar))    \chi_{[\zeta^{-1}\ll x\ll 1]} \Big]\, dx + \text{cc.},
	}
	cf.~\eqref{eq:F2ell1 B}. In view of \eqref{eq:SOT} and~\eqref{eq:of size}, 
	\EQ{ \label{eq:F13 B bd}
		|F_{13}(\xi,\eta;\hbar)| & \lesssim \hbar \alpha^{-2} \zeta^{-\frac12}\sum_{\ell+k\le m} \zeta^{k-2(m+\ell)} \int_{0}^{\infty} \Big| \partial_x^k\Big\{  W^{+}_{n}(x/\alpha)   \left(Q_0(x;\alpha,\hbar)\right)^{-\frac14}\left(-Q_0(\tilde x;\beta,\hbar)\right)^{-\frac14} e^{ -\frac{2}{3\hbar} (-\tau)^{\frac32}}    \\
		&\qquad \left(\chi_{[x\lesssim\alpha]}+\chi_{[x\gtrsim\alpha]}\right)  (1+\hbar \tilde a_0(-\tau;\alpha,\hbar))  (1+\hbar\overline{\tilde a_1}(-\tilde\tau;\beta,\hbar))    \chi_{[\zeta^{-1}\ll x\ll 1]} \Big\}\Big| \, dx
	}
Here $\chi_{[x\lesssim\alpha]}$ and $\chi_{[x\gtrsim \alpha]}$ are smooth cutoffs.
 
	We now analyze $\partial_x^k$ derivatives of each of the factors in~\eqref{eq:F13 B bd} for all $k\ge1$, and uniformly in the parameters such as $x\ll1$:  
	\begin{itemize}
		\item By differentiation of \eqref{W pre}, $\alpha^{-2}\left|\partial_{x}^{k}\left(W_{n}^{+}(x/\alpha)\right)\right|\leq C_{k}\hbar^{-1}\alpha^{-2-k}$ if $x\lesssim \alpha$ and $\alpha^{-2}\left|\partial_{x}^{k}\left(W_{n}^{+}(x/\alpha)\right)\right|\leq C_{k}\hbar^{-1}\alpha^{2}x^{-4-k}$ if $x\gg\alpha$.
		\item $\left|\partial_{x}^{k}\left|Q_0(x;\alpha,\hbar)\right|^{-\frac14}\right|\leq C_{k}x^{\frac12-k}$ by the definition of $Q_{0}(x;\alpha,\hbar)$ and the fact that $0<x\ll1$
		\item $\left|\partial_{x}^{k}\left(-Q_0(\tilde{x};\beta,\hbar)\right)^{-\frac14}\right|\leq C_{k}\zeta^{-2}$ by \eqref{eq:Q0 derxk}
		\item $\left|\partial_{x}^{k}\left(\chi_{x\leq\alpha}e^{-\frac{2}{3\hbar}(-\tau)^{\frac32}}\right)\right|\leq C_{k}x^{N-k}$ by Lemma \ref{lem: Lemma 3.3 CDST} and the second equation in \eqref{eq:xtau}
		\item $\left|\partial_{x}^{k}\left(\chi_{x\geq\alpha}e^{-\frac{2}{3\hbar}(-\tau)^{\frac32}}\right)\right|\leq C_{k}x^{N-k}$ by the first equation in \eqref{eq:xtau}
		\item $\left|\partial_{x}^{k}(1+\hbar \tilde a_{1}(-\tilde{\tau};\beta,\hbar))\right|\leq C_{k}\hbar \zeta^{-1}$ by Lemma \ref{lem:a1fine}.
		\item By Lemma \ref{lem:a0fine1} $\left|\partial_{x}^{k}(1+\hbar \tilde a_{0}(-\tau;\alpha,\hbar))\right|\leq C_{k}\hbar (-\tau)^{-3}x^{-k}$ for $\alpha\gtrsim 1$  
		\item By Lemma \ref{lem:a0fine2} $\left|\partial_{x}^{k}(1+\hbar \tilde a_{0}(-\tau;\alpha,\hbar))\right|\leq C_{k}\hbar x^{-k}$ for $0<\alpha\ll1$.
	\end{itemize}
To obtain the last two estimates in the above list, besides Lemma \ref{lem:a0fine1} and Lemma \ref{lem:a0fine2}, we also use the estimate $\left|\frac{\partial(-\tau)}{\partial x}\right|\lesssim (-\tau)^{-\frac12}$.

This bullet list implies that for a constant $c\in(0,1)$,
\begin{align}\label{eq:F13 B bd final}
	\begin{split}
	\left|F_{13}(\xi,\eta;\hbar)\right| \lesssim&\hbar\alpha^{2}\langle\alpha\rangle^{-4} \hbar^{-1}\zeta^{-\frac12}\zeta^{-m}\int_{\zeta^{-1}}^{1}x^{c\hbar^{-1}}\,dx\lesssim \hbar\alpha^{2}\langle\alpha\rangle^{-4}\zeta^{-\frac12-m}.
	\end{split}
\end{align}
Finally we consider $F_{12}$. Similar to \eqref{eq:F13 bd B}, we have
\begin{align}\label{eq:F12 bd B}
	\begin{split}
	F_{12}(\xi,\eta;\hbar) &= 2\alpha^{-2}\int_{0}^{\infty} W^{ {+}}_{n}(x/\alpha) \hbar^{\frac13} x^{\frac12} q(\tau;\hbar)^{-\frac14} \Ai(-\hbar^{-\frac23}\tau)(1+\hbar a_0(-\tau;\alpha,\hbar))  \\
	&\qquad \hbar^{\frac13} \tilde x^{-\frac12} \tilde q(\tilde \tau;\hbar)^{-\frac14} \Re\Big(a(\eta) (\Ai(-\hbar^{-\frac23}\tilde \tau)-i \Bi(-\hbar^{-\frac23}\tilde \tau)) (1+\hbar \overline{a_1}(-\tilde\tau;\alpha,\hbar)) \Big) \chi_2(\tilde \tau;\hbar) \, dx 
	\end{split}
\end{align}
We simply bound the oscillatory Airy function in absolute value  by
\begin{align*}
	\left|\tilde x^{-\frac12} \tilde q(\tilde \tau;\hbar)^{-\frac14} \Re\Big(a(\eta) (\Ai(-\hbar^{-\frac23}\tilde \tau)-i \Bi(-\hbar^{-\frac23}\tilde \tau)) (1+\hbar \overline{a_1}(-\tilde\tau;\alpha,\hbar)) \Big)\right|\lesssim 1.
\end{align*}
Then similar as for $F_{11}$, we have
\begin{align}\label{eq:F12 bd B final}
	\begin{split}
	\left|F_{12}(\xi,\eta;\hbar)\right|\lesssim &\alpha^{-2}\int_{0}^{\frac{1}{\zeta}}\langle x/\alpha\rangle^{-4}\left|Q_{0}(x;\alpha,\hbar)\right|^{-\frac14}e^{-\frac{2}{3\hbar}(-\tau)^{\frac32}}\left|1+\hbar \tilde a_{0}(-\tau;\alpha,\hbar)\right|\chi_{2}(\tilde{\tau};\hbar)\,dx\\
	\lesssim &\alpha^{-2}\int_{0}^{\frac{1}{\zeta}}\langle x/\alpha\rangle^{-4}x^{N+\frac12}\,dx\\
	\lesssim &\hbar\alpha^{2}\langle\alpha\rangle^{-4}\zeta^{-N+1}.
	\end{split}
\end{align}
So the estimates \eqref{eq:F11 bd B},\eqref{eq:F13 B bd final} and \eqref{eq:F12 bd B final} finally conclude the proof of \eqref{F hbar bounds} and \eqref{F hbar bounds auxi}.
\\

We now turn to bounding the derivatives in the variables $\xi^{\frac12}$, resp., $\eta^{\frac12}$, and we consider $k_1=1, k_2=0$ in~\eqref{deri F hbar and n bounds 3}. We need to verify that differentiating once in $\xi^{\frac12}$ results in the same bound that we obtained above for the undifferentiated case, but with a loss of~$\alpha^{-1}$ if $0<\alpha<1$. We will instead differentiate in $\alpha$, resp., $\beta$ using that $\partial_{\xi^{\frac12}}=\hbar \partial_\alpha$, $\partial_{\eta^{\frac12}}=\hbar \partial_\beta$. 
Beginning with Case~$A_3$,  rewrite  \eqref{eq:F3} in the form 
\begin{align}
	&F_3(\xi,\eta;\hbar)  =4\pi^{-1} \hbar \alpha^{-2} \zeta^{-\frac12} \int_{0}^{\infty}W^{ {+}}_{n}(x/\alpha)     (1+\rho'(x;\alpha))^{-\frac12}   \Re \left( a(\xi)  e^{\frac{i}{\hbar}[x-y(\alpha)+\rho(x;\alpha)]}    (1+\hbar    O(x^{-1}) ) \right) \nn \\
	&   \qquad \qquad\qquad(1+\rho'(x\zeta;\beta))^{-\frac12}   \Re \left( a(\eta) e^{\frac{i}{\hbar}[x\zeta-y(\beta)+\rho(x\zeta;\beta)]}    (1+\hbar  O(x^{-1}) ) \right)\, \chi(x) \, dx \label{eq:F3*} \\
	& = 4\pi^{-1} \hbar \alpha^{-2} \zeta^{-\frac12} a(\xi)  \overline{a(\eta)}  \int_{0}^{\infty} e^{\frac{i}{\hbar} x(1-\zeta)}  W^{ {+}}_{n}(x/\alpha)     (1+\rho'(x;\alpha))^{-\frac12}  (1+\rho'(x\zeta;\beta))^{-\frac12} \nn \\
	& \qquad \qquad\qquad  \qquad\qquad \qquad\qquad e^{\frac{i}{\hbar}[y(\beta)-y(\alpha)+\rho(x;\alpha)-\rho(x\zeta;\beta)]} (1+\hbar  O(x^{-1}) )  \, \chi(x) \, dx  + \text{cc} \label{eq:F3 onebar} \\
	& \quad + 4\pi^{-1} \hbar \alpha^{-2} \zeta^{-\frac12} a(\xi)   {a(\eta)} \int_{0}^{\infty} e^{\frac{i}{\hbar} x(1+\zeta)}  W^{ {+}}_{n}(x/\alpha)     (1+\rho'(x;\alpha))^{-\frac12}  (1+\rho'(x\zeta;\beta))^{-\frac12}  \nn \\
	& \qquad \qquad\qquad  \qquad\qquad \qquad\qquad e^{\frac{i}{\hbar}[-y(\beta)-y(\alpha)+\rho(x;\alpha)+\rho(x\zeta;\beta)]} (1+\hbar  O(x^{-1}) )  \, \chi(x) \, dx  + \text{cc}\label{eq:F3 nobar}
\end{align}
where ``cc" stands for {\em complex conjugate}. As noted above, the $O(x^{-1})$ has symbol behavior under differentiation in $x$, and each derivative in $\alpha$ brings out a factor of $\alpha^{-1}$. 
The latter follows from~\eqref{eq:a1 fine}. 
We only treat the term \eqref{eq:F3 onebar}, since \eqref{eq:F3 nobar} satisfies better bounds.  By \eqref{eq:rho n} and~\eqref{eq:a size}, we have $\partial_\alpha a(\xi)= O(\alpha^{-1})$ uniformly in all parameters (while $a(\xi)= O(1)$). 
Writing $\zeta=\beta/\alpha$ implies $\zeta_\alpha=-\zeta/\alpha$ and $\zeta_\beta=\zeta/\beta$. Therefore, in view of Lemma~\ref{lem: Lemma 3.4 CDST}, differentiating 
\[
\alpha^{-2} \zeta^{-\frac12} a(\xi)\overline{a(\eta)}  (1+\rho'(x;\alpha))^{-\frac12}  (1+\rho'(x\zeta;\beta))^{-\frac12}  e^{\frac{i}{\hbar}[y(\beta)-y(\alpha)+\rho(x;\alpha)-\rho(x\zeta;\beta)]} (1+\hbar  O(x^{-1}) )
\]
in $\alpha$  produces a sum of terms of similar form, but multiplied with extra factors which are bounded by $\alpha^{-1}$ or even by $\hbar (1+\alpha)^{-1}$ (or $\hbar (1+\alpha)^{-3}$ as for $y'(\alpha)$). On the other hand, 
\[
\partial_\alpha \left( W^{ {+}}_{n}(x/\alpha) \right) = -\alpha^{-1} (W^{ {+}}_{n})'(x/\alpha) x/\alpha = \hbar^{-1} \alpha^{-1} O(\langle x/\alpha\rangle^{-4}) 
\]
which differs from $W^{ {+}}_{n}(x/\alpha)=  \hbar^{-1}  O(\langle x/\alpha\rangle^{-4}) $ only by the $\alpha^{-1}$ factor. Note that all these terms are better by a factor of $\hbar$ than what we need since we are actually taking a derivative $\hbar \partial_\alpha $. The most important contribution results from 
\[
\hbar \partial_\alpha \, e^{\frac{i}{\hbar} x(1-\zeta)}  =  i\zeta \frac{x}{\alpha} e^{\frac{i}{\hbar} x(1-\zeta)}  
\]
Since we are in the regime $\zeta\simeq 1$, the only  effect here is to multiply the potential $W^{ {+}}_{n}(x/\alpha) $ by one factor of~$\frac{x}{\alpha} $.  Consider first $\alpha>1$. Due to the fact that in our estimation of $F_3$ above the dominant contribution in the integrals $\int_1^\infty\ldots dx$ came from $x\simeq \alpha$, we see that we arrive at the same estimate as above by analogous arguments. The only essential feature here is that $W^{ {+}}_{n}(x/\alpha) x/\alpha$ remains integrable at $x=\infty$ since it decays like $x^{-3}$ (this changes if we were to take three derivatives). On the other hand, if $0<\alpha\lesssim 1$, then we lose a factor of $\alpha^{-1}$ since the dominant contribution to the integral is derived from the region $x\simeq 1$. 
In summary, the upshot from this discussion is that Case~$A_3$ satisfies the claimed estimate. 

\smallskip

Turning to Case~$A_2$, we rewrite \eqref{eq:F2 pre} in the form 
\begin{align}\label{eq:F2 pre*}
	\begin{split}
		F_2(\xi,\eta;\hbar) 
		& = \hbar^{\frac23} \alpha^{-2}\zeta^{-\frac12} \int_{0}^{\infty}W^{ {+}}_{n}(x/\alpha)  q^{-\frac14}( \tau)\tilde q^{-\frac14}( \tilde \tau)\Big[  \Ai(-\hbar^{-\frac23} \tau)(1+\hbar a_0(-\tau;\alpha,\hbar)) \chi_{[\tau\lesssim  \hbar^{\frac23}]} + \\
		&\qquad + 2  \Re\big( a(\xi) \big(\Ai(-\hbar^{-\frac23}\tau)-i\Bi(-\hbar^{-\frac23}\tau)\big)(1+\hbar \overline{a_1(-\tau;\alpha)})  \big) \chi_{[\tau\gtrsim  \hbar^{\frac23}]}\Big] \\
		&\qquad  \Big[  \Ai(-\hbar^{-\frac23}\tilde \tau)(1+\hbar a_0(-\tilde \tau;\beta,\hbar)) \chi_{[\tilde \tau\lesssim  \hbar^{\frac23}]} + \\
		&\qquad + 2  \Re\big( a(\eta) \big(\Ai(-\hbar^{-\frac23}\tilde\tau)-i\Bi(-\hbar^{-\frac23}\tilde\tau)\big)(1+\hbar \overline{a_1(-\tilde\tau;\beta)})  \big) \chi_{[\tilde\tau\gtrsim \hbar^{\frac23}]} \Big]\, \chi(x )\, dx
	\end{split}
\end{align}
where the integration runs over $x\simeq 1$. We differentiate this expression by $\partial_{\xi^{\frac12}}=\hbar \partial_\alpha$ using Lemma~\ref{lem: Lemma 3.2 CDST} in the regime $x\simeq1$. Viewed as functions of $x$ and $\alpha$, the derivative $\partial_\alpha \tau$  gains a factor of  $\hbar (1+\alpha)^{-3}$, uniformly in this regime, while 
\[
\partial_\alpha \tilde\tau = \partial_\alpha \tau(x\zeta,\beta)= - (\partial_1 \tau) (x\zeta,\beta) \zeta x/\alpha = O(\alpha^{-1})
\]
By Lemma~\ref{lem: Lemma 3.2 CDST} and the chain rule, both~$\partial_{\alpha}q(\tau(x,\alpha),\alpha)$ and $\partial_{\alpha}\tilde q(\tilde \tau)=\partial_{\alpha}q(\tau(x\zeta,\beta),\beta)$ are bounded by $O(\alpha^{-1})$.  
The potential term is treated the same as in Case~$A_3$ above, as are $ \alpha^{-2}\zeta^{-\frac12}$ and $a(\xi), a(\eta)$. We now consider the case when the $\partial_\alpha$ derivative falls on $a_j, j=0,1$. Note that both $\partial_\alpha$ and $a_j$ come with a factor of $\hbar$, so in total we gain a factor of~$\hbar^2$. First (suppressing $\hbar$ from the notation), by the preceding bound on $\tau_\alpha$ and Lemmas~\ref{lem:a0fine1} and~\ref{lem:a0fine2} (with $\zeta=-\tau$ in those lemmas having nothing to do with the $\zeta$ used in this proof),
\EQ{
	\nn
	\partial_\alpha a_0(-\tau(x,\alpha);\alpha) =  -\tau_\alpha \,\partial_1 a_0 (-\tau,\alpha) + \partial_\alpha a_0 (-\tau,\alpha) = O(\alpha^{-1}), \quad -1\lesssim \tau\lesssim \hbar^{\frac23}
}
Here we used that $\partial_1 a_0 (-\tau,\alpha) = O(\hbar^{-\frac13})$, see \eqref{eq:a0 fine} and~\eqref{eq:a0 fine***}, which is harmless due to the gains of factors of~$\hbar$. We remark that the main argument above for Case~$A_2$ only needs the {\em size of the integrand} in the regime of $\tau$ in which  $a_0$ is relevant, so these bounds suffice (without any further differentiation, in contrast to the oscillatory $a_1$ regime we integrate by parts once in $\tau$ as independent variable). Next,
\EQ{
	\nn
	\partial_\alpha a_0(-\tilde\tau;\beta)= \partial_\alpha \, a_0(-\tau(x\zeta,\beta);\beta) =  x\zeta\alpha^{-1}\partial_1\tau(x\zeta,\beta) \,\partial_1 a_0 (-\tilde\tau,\beta)  = O(\alpha^{-1}\hbar^{-\frac13}), \quad \frac12 \le x \le x_t +  \hbar^{\frac23}
}
By \eqref{eq:a1 fine} 
\[
\partial_\alpha a_1(-\tau(x,\alpha),\alpha) = -\tau_\alpha\, \partial_1 a_1(-\tau,\alpha) + \partial_2 a_1(-\tau,\alpha)  = O(\langle \alpha\rangle^{-1}) 
\]
Differentiating this once in $x$ (which is the same as differentiating in $\tau$ up to a factor of $q^{\frac12}\simeq 1$), leads to a loss of a factor of $\hbar^{-\frac13}$ which is harmless given the aforementioned~$\hbar^2$ whence $\partial_x\partial_\alpha a_1(-\tau(x,\alpha),\alpha) = O(\hbar^{-\frac13}\langle \alpha\rangle^{-1})$. 
A more delicate term is $\hbar a_1(-\tilde\tau,\beta)=\hbar a_1(-\tau(\zeta x,\beta),\beta)$ since it is involved in one integration by parts relative to~$\tau$, see~\eqref{eq:Fs osc+-}. On the one hand, 
\EQ{\label{eq:a1 delic}
	\hbar\partial_\alpha \hbar a_1(-\tilde\tau,\beta) =  \hbar^2 (\partial_1 a_1)(-\tilde\tau,\beta)(\partial_1 \tau) (x\zeta,\beta) \zeta x/\alpha = O\left(\hbar^2 \tilde\tau^{-\frac12}\alpha^{-1}\right)=O\left(\hbar^2 \hbar^{-\frac13}\alpha^{-1}\right),\quad \tilde\tau\gtrsim \hbar^{\frac23}
}
and on the other, by~\eqref{eq:a1 fine}, another derivative in $x$ (or $\tau$) loses a further factor of $\hbar^{-\frac23}$. 
In conclusion, those terms in which the $\partial_\alpha$ derivative {\em does not fall on the Airy functions} are treated in the exact same fashion as in Case~$A_2$ above, with the net effect of a factor of   $\hbar\alpha^{-1}$ for all $\alpha>0$ which multiplies the original estimate.

From \eqref{eq:AiiBi}, for all $x\ge1$, 
\EQ{\label{eq:AiiBi'}
	\left(\Ai(-x)+i\Bi(-x)\right)' &= c\, x^{\frac12} (1+b_1(x)) \left(\Ai(-x)+i\Bi(-x)\right) \\
	b_1^{(k)} (x) & = O(x^{-\frac32-k}), \quad k\ge0
}
whence
\EQ{ \label{eq:AiiBi tau'}
	\hbar\partial_\alpha \left(\Ai\left(-\hbar^{-\frac23} \tau\right)-i\Bi\left(-\hbar^{-\frac23} \tau\right)\right) &=  c\, \tau^{\frac12} \tau_\alpha \left(1+\hbar O\left(\tau^{-\frac32}\right) \right) \left(\Ai\left(-\hbar^{-\frac23} \tau\right)-i\Bi\left(-\hbar^{-\frac23} \tau\right)\right)
}
for $\tau\geq \hbar^{\frac23}$, with $|\tau_\alpha|\lesssim \hbar (1+\alpha)^{-3}$. Similarly, 
\EQ{\label{eq:Ai'}
	\Ai'(x) &= c\, x^{\frac12} (1+b_2(x)) \Ai(x) \\
	b_2^{(k)} (x) & = O(x^{-\frac32-k}), \quad k\ge0
}
leading to the same type of gain of the $\hbar (1+\alpha)^{-3}$ factor. Using that  
$
\hbar\partial_\alpha\tilde\tau  = O(\hbar \alpha^{-1}) 
$
from above, 
we conclude that $\hbar\partial_\alpha$ acting on the second bracket in~\eqref{eq:F2 pre*} yields  the oscillatory term 
\EQ{\label{eq:AiiBi tildetau'}
	\hbar\partial_\alpha \left(\Ai\left(-\hbar^{-\frac23} \tilde\tau\right)-i\Bi\left(-\hbar^{-\frac23} \tilde\tau\right)\right) &=  c\, \tilde \tau^{\frac12} O\left(\alpha^{-1}\right) \left(1+\hbar O\left(\tilde{\tau}^{-\frac32}\right) \right) \left(\Ai\left(-\hbar^{-\frac23} \tilde\tau\right)-i\Bi\left(-\hbar^{-\frac23}\tilde \tau\right)\right)
}
and analogously for the decaying Airy function. This shows that by repeating the exact same arguments as in the undifferentiated Case~$A_2$, see in particular~\eqref{eq:Fs osc+-}, we obtain the desired bound. 

\smallskip

It remains to consider Case~$A_1$.  Using that $q=-\frac{Q_0}{\tau}$, we rewrite the relevant integral in this case, see~\eqref{eq:F2r2}, in the form  
\begin{align}
	\begin{split}\nn 
		& \int_{0}^{x_0/\alpha} W_{n}^{+}(R)\hbar^{\frac13}x^{-\frac12}q^{-\frac14}(\tau)\Ai(-\hbar^{-\frac23}\tau)(1+\hbar a_{0}(-\tau,\alpha;\hbar))
		\hbar^{\frac13}\tilde{x}^{-\frac12}q^{-\frac14}(\tilde{\tau})\Ai(-\hbar^{-\frac23}\tiltau)(1+\hbar a_{0}(-\tilde{\tau},\beta;\hbar))R\,dR \\
		& =  \alpha^{-2} \hbar \zeta^{-\frac12} \!\!\int_0^{x_0}  W_{n}^{+}(x/\alpha) (Q_0(x,\alpha)Q_0(\zeta x,\beta))^{-\frac14}  e^{-\frac{2}{3\hbar}[(-\tau)^{\frac32} + (-\tilde\tau)^{\frac32} ]}
		(1+\hbar a_{0}(-\tau,\alpha;\hbar)) (1+\hbar a_{0}(-\tilde{\tau},\beta;\hbar))\,dx
	\end{split}
\end{align}
with some $0<x_0\ll 1$. The terms arising if $\partial_\alpha$ hits $\alpha^{-2}\zeta^{-\frac12}W_{n}^{+}(x/\alpha)$ are the same as above. Next, 
\EQ{\nn 
	\partial_\alpha (Q_0(x,\alpha)Q_0(\zeta x,\beta))^{-\frac14}  = -\frac14 (Q_0(x,\alpha)Q_0(\zeta x,\beta))^{-\frac14}  \left[ \frac{\partial_\alpha Q_0(x,\alpha)}{Q_0(x,\alpha)} -  \frac{\partial_1 Q_0(\zeta x,\beta)}{Q_0(\zeta x,\beta)} x \zeta/\alpha  \right]
} 
By \eqref{eq:Q0kell} the term in brackets is $O(\alpha^{-1})$. Here we also used the fact that if $0<x\ll1, 0<\zeta x\ll1$, then $\frac{1}{Q_{0}(x,\alpha)}\simeq x^{2}, \frac{1}{Q_{0}(\zeta x,\beta)}\simeq \zeta^{2}x^{2}$. The dominant contribution is (recall $\tau=\tau(x,\alpha;\hbar)$ and $\tilde \tau=\tau(x\zeta,\beta;\hbar)$)
\EQ{
	\hbar \partial_\alpha e^{-\frac{2}{3\hbar}[(-\tau)^{\frac32} + (-\tilde\tau)^{\frac32} ]}  & = e^{-\frac{2}{3\hbar}[(-\tau)^{\frac32} + (-\tilde\tau)^{\frac32} ]}  \left((-\tau)^{\frac12}\tau_\alpha -(-\tilde\tau)^{\frac12} \partial_1\tau(x\zeta,\beta)x \zeta/\alpha \right)
}
By \eqref{eq:xtau}, $\tau_\alpha=O((-\tau)^{-\frac12}\alpha^{-1})$ and $\partial_1\tau(x\zeta,\beta)= O(x^{-1}(-\tilde\tau)^{-\frac12})$ whence
\EQ{
	\hbar \partial_\alpha e^{-\frac{2}{3\hbar}[(-\tau)^{\frac32} + (-\tilde\tau)^{\frac32} ]}  & = O(\alpha^{-1}) e^{-\frac{2}{3\hbar}[(-\tau)^{\frac32} + (-\tilde\tau)^{\frac32} ]}  
}
as desired. Finally, by \eqref{eq:a0 fine} and \eqref{eq:a0 fine***}, 
\EQ{
	\nn
	\partial_\alpha a_0(-\tau(x,\alpha);\alpha) &=  -\tau_\alpha \,\partial_1 a_0 (-\tau,\alpha) + \partial_\alpha a_0 (-\tau,\alpha) \\
	&=  O((-\tau)^{-\frac12}\alpha^{-1}) (-\tau)^{\frac12} +  O(\alpha^{-1}) = O(\alpha^{-1}),  
}
and 
\EQ{
	\nn
	\partial_\alpha a_0(-\tilde\tau;\beta) &= \partial_\alpha \, a_0(-\tau(x\zeta,\beta);\beta) =  x\zeta\alpha^{-1}\partial_1\tau(x\zeta,\beta) \,\partial_1 a_0 (-\tilde\tau,\beta) \\
	& =   x\zeta\alpha^{-1} O((-\tilde\tau)^{-\frac12} x^{-1}) O((-\tilde\tau)^{\frac12}) = O(\alpha^{-1})
}
uniformly in the regime $0<x\ll1$ (which corresponds to $\tau\lesssim -1$). 

Case~$B$, which means $\eta\gg \xi$ or $\zeta\gg1$, now proceeds entirely analogously. For example, \eqref{eq:F3 B} needs to be replaced by 
\EQ{\label{eq:F3 B*}
	|\partial_{\xi^{\frac12}} F_3(\xi,\eta;\hbar)|  \le C_\ell\, \zeta^{-\frac12} (\hbar/ \zeta)^{\ell}  \alpha\langle \alpha\rangle^{-3}
}
for any $\ell\ge1$, reflecting a loss of a factor of $\max(1,\alpha^{-1})$.  This can be seen by applying $\hbar\partial_\alpha$ to \eqref{eq:F3 onebar} and~\eqref{eq:F3 nobar}, followed by the exact same integration by parts which yield~\eqref{eq:F3}. The analysis of $B_2$ and $B_1$ proceed in the same way and we skip the details. 
The other derivatives $\partial_{\xi^{\frac12}}^{k_1} \partial_{\eta^{\frac12}}^{k_1}$ which can appear in the regime $k_1+k_2\le2$ are also treated in the same fashion, leading to the estimates claimed in the proposition. 

We now turn to the expressions 
\[
\left(\partial_{\xi^{\frac12}} + \partial_{\eta^{\frac12}}\right)^k   F(\xi,\eta;\hbar) = \hbar^k \left(\partial_{\alpha} + \partial_{\beta}\right)^k   F(\xi,\eta;\hbar)
\]
assuming $\xi\simeq\eta$ or $\zeta\simeq1$. Thus it suffices to apply the operator $\hbar \left(\partial_{\alpha} + \partial_{\beta}\right)$ and its powers to Case~$A$ above. 
We begin with Case~$A_3$, see \eqref{eq:F3 onebar} and~\eqref{eq:F3 nobar}.   The key identities for \eqref{eq:F3 onebar} are
\EQ{\label{eq:spur1}
	-i\hbar\partial_x  e^{\frac{i}{\hbar}x(1-\zeta)} &= (1-\zeta)   e^{\frac{i}{\hbar}x(1-\zeta)}\\
	\hbar (\partial_{\alpha} + \partial_{\beta}) e^{\frac{i}{\hbar}x(1-\zeta)} &= -i\frac{x}{\alpha}(1-\zeta) e^{\frac{i}{\hbar}x(1-\zeta)}= -\hbar\frac{x}{\alpha}\partial_x  e^{\frac{i}{\hbar}x(1-\zeta)}
}
and for~\eqref{eq:F3 nobar} they take the form
\EQ{\label{eq:spur2}
	-i\hbar\partial_x  e^{\frac{i}{\hbar}x(1+\zeta)} &= (1+\zeta)   e^{\frac{i}{\hbar}x(1+\zeta)}\\
	\hbar (\partial_{\alpha} + \partial_{\beta}) e^{\frac{i}{\hbar}x(1+\zeta)} &= i\frac{x}{\alpha}(1-\zeta) e^{\frac{i}{\hbar}x(1+\zeta)}= \hbar\frac{x}{\alpha}\frac{1-\zeta}{1+\zeta}\partial_x  e^{\frac{i}{\hbar}x(1+\zeta)}. 
}
Thus, applying \eqref{eq:spur1} to \eqref{eq:F3 onebar}, respectively \eqref{eq:spur2} to~\eqref{eq:F3 nobar}, yields, upon integration by parts in $x$,
\begin{align}
	& \hbar \left(\partial_{\alpha} + \partial_{\beta}\right) F_3(\xi,\eta;\hbar)  = 4\pi^{-1} \hbar \alpha^{-2} \zeta^{-\frac12} a(\xi)  \overline{a(\eta)}  \int_{0}^{\infty} e^{\frac{i}{\hbar} x(1-\zeta)} \hbar \partial_x\left[  W^{ {+}}_{n}(x/\alpha) x/\alpha     \left(1+\rho'(x;\alpha)\right)^{-\frac12}  \left(1+\rho'(x\zeta;\beta)\right)^{-\frac12}\right. \nn \\
	& \qquad \qquad\qquad \left.e^{\frac{i}{\hbar}\left[y(\beta)-y(\alpha)+\rho(x;\alpha)-\rho(x\zeta;\beta)\right]} \left(1+\hbar  O(x^{-1}) \right)  \, \chi(x)\right] \, dx  + \text{cc}  \label{eq:F3 onebar*} \\
	& \quad \qquad + 4\pi^{-1} \hbar \alpha^{-2} \zeta^{-\frac12} a(\xi)   {a(\eta)} \frac{\zeta-1}{1+\zeta} \int_{0}^{\infty} e^{\frac{i}{\hbar} x(1+\zeta)}  \hbar \partial_x\left[   W^{ {+}}_{n}(x/\alpha) x/\alpha    (1+\rho'(x;\alpha))^{-\frac12}  (1+\rho'(x\zeta;\beta))^{-\frac12}\right.  \nn \\
	& \qquad \qquad\qquad  \left.  e^{\frac{i}{\hbar}\left[-y(\beta)-y(\alpha)+\rho(x;\alpha)+\rho(x\zeta;\beta)\right]} \left(1+\hbar  O(x^{-1}) \right)  \, \chi(x)\right] \, dx  + \text{cc} + O\left(\hbar \alpha^{-1} \Gamma\right) \label{eq:F3 nobar*}
\end{align}
The final $O$-term here is a result of those expressions in which the  $\hbar (\partial_{\alpha} + \partial_{\beta})$ derivatives fall on the non-oscillatory terms, and they are treated as above. It is essential here that we obtain $O(\hbar \alpha^{-1} \Gamma)$ rather than $O(\hbar \max(1,\alpha^{-1} )\Gamma)$ as in the case of $\hbar\partial_\alpha$ and~$\hbar\partial_\beta$, see~\eqref{deri F hbar and n bounds 3}. As explained in the paragraph preceding~\eqref{eq:F2 pre*}, the reason for the absence of a gain 
of~$\alpha^{-1}$  for $\alpha>1$ lies only with the extra $x$-factor resulting from $\partial_\alpha$ hitting the complex exponential. This term is absent in the $O(\cdot)$ in~\eqref{eq:F3 nobar*}.  Placing absolute values in the integrals in~\eqref{eq:F3 onebar*} and~\eqref{eq:F3 nobar*}, and using Lemma~\ref{lem: Lemma 3.4 CDST} as before,  yields a bound of~$O(\hbar \alpha^{-2})=O(\xi^{-\frac12}\Gamma)$. 
Iterating this process leads to~\eqref{eq:trace der}.

For Case~$A_2$, we return to the integral \eqref{eq:F2 pre*}. As we proved in the paragraph preceding \eqref{eq:AiiBi'}, if either $\hbar\partial_\alpha$ or, by symmetry, $\hbar\partial_\beta$, fall on any term other than the Airy functions in the integral in~\eqref{eq:F2 pre*}, then we obtain a net factor of $\hbar\alpha^{-1}$ as required for~\eqref{eq:trace der} (recall $\alpha\simeq\beta$).  Furthermore, if $\left(\hbar \partial_\alpha\right)^k$, respectively, $\left(\hbar \partial_\beta\right)^k$ with $k\ge2$ hits any term in~\eqref{eq:F2 pre*} other than the Airy functions and $a_0, a_1$, then by Lemmas~\ref{lem: monotonicity of root in alpha}, \ref{lem: Lemma 3.2 CDST}, \ref{lem: Lemma 3.4 CDST},  we obtain the desired $(\hbar\alpha^{-1})^k$ factor. For the case of $a_j(-\tau,\alpha)=a_j(-\tau(x,\alpha),\alpha)$, with $j=0,1$, we have from~\eqref{eq:a1 fine}, \eqref{eq:a0 fine}, and~\eqref{eq:a0 fine***}, 
\EQ{\nn 
	\hbar \partial_\alpha \,\hbar a_j(-\tau,\alpha) &= -\hbar^2\partial_\alpha \tau(x,\alpha) (\partial_1 a_j)(-\tau,\alpha) + \hbar^2(\partial_2 a_j)(-\tau,\alpha) \\
	& = O\left(\hbar^3\langle \alpha\rangle^{-3} |\tau|^{-\frac12}\right) + O\left(\hbar^2\alpha^{-1}\right) = O\left(\hbar^2\alpha^{-1}\right)
}
for $\hbar^{\frac23}\lesssim |\tau|\lesssim 1$ (and as the reader can check also for $|\tau|\le \hbar^{\frac23}$). The higher order derivatives $(\hbar \partial_\alpha)^{\ell} \,\hbar a_j(-\tau,\alpha)$ are linear combinations of terms of the form
\EQ{\label{eq:hell p 1}
	\hbar^{\ell+1}\, \partial^{k_1}_\alpha \tau\cdot\ldots\cdot \partial^{k_n}_\alpha \tau \; \partial_1^{n} \partial_2^{m} a_j(-\tau,\alpha),\qquad k_1+\ldots+k_n+m=\ell,\;k_i\ge1,\; n\ge0
}
By Lemmas~\ref{lem: monotonicity of root in alpha} and \ref{lem: Lemma 3.2 CDST}, $\partial^{k_i}_\alpha \tau= O(\hbar (1+\alpha)^{-k_i-2})$, whereas by the aforementioned bounds on $a_j$ and their derivatives, $\partial_1^{n} \partial_2^{m} a_j(-\tau,\alpha)= O\left(|\tau|^{\frac12-n}\alpha^{-m}\right)$.  In conclusion, the term in~\eqref{eq:hell p 1} is bounded by, for $\hbar^{\frac23}\le|\tau|\lesssim 1$,  
\[
\lesssim \hbar^{2\ell+1-m} (1+\alpha)^{-2n-(\ell-m)} |\tau|^{\frac12-n}\alpha^{-m}\lesssim \hbar^{\ell+1}\cdot\hbar^{\ell-m}(1+\alpha)^{-\ell+m}|\tau|^{\frac12-n}\alpha^{-m} \lesssim  \hbar^{\ell+1} \alpha^{-\ell}
\]
using that $m+n\le\ell$, and the same bounds also applies to $|\tilde\tau|\lesssim  \hbar^{\frac23}$. While we have so far not made use of the sum of the derivatives, 
higher derivatives of the more delicate $a_j(-\tilde\tau,\beta)=a_j\left(-\tau(x\zeta,\beta),\beta\right)$ do require the sum $\partial_\alpha+\partial_\beta$, cf.~\eqref{eq:a1 delic}. The point being that $(\hbar \partial_\alpha)^{\ell} \,\hbar a_j(-\tilde\tau,\beta)$  generates expressions of the form 
\[
\hbar^{\ell+1} \partial_1^\ell a_j (-\tilde\tau,\beta) (\partial_\alpha \tilde\tau)^{\ell} = O\left(\hbar^{\ell+1}  |\tilde\tau|^{\frac12-\ell} \alpha^{-\ell}\right)
\]
which violate our desired bounds if $\tilde\tau\simeq \hbar^{\frac23}$ and $\ell$ becomes large. Note that there are other contributions from the expression for $(\hbar \partial_\alpha)^{\ell} \,\hbar a_j(-\tilde\tau,\beta)$, which, however, contributes less powers in $|\tilde{\tau}|^{-1}$. Now 
\EQ{\label{eq:aj spur}
	\hbar(\partial_\alpha+\partial_\beta)\hbar a_j(-\tilde\tau,\beta) & =  \hbar^2 x\zeta \beta^{-1}(\zeta-1) (\partial_1 a_j)(-\tilde\tau,\beta)  \partial_1\tau(x\zeta,\beta)
	\\
	&\qquad -\hbar^2 (\partial_1 a_j)(-\tilde\tau,\beta) (\partial_\beta\tau)(x\zeta,\beta)+\hbar^2 \partial_2a_j(-\tilde\tau,\beta)
}
The final term is $O\left(\hbar^2\alpha^{-1}\right)$, which is better by a factor of~$\hbar$ than what we need. 
The derivative $\partial_1 a_j$ loses $\left(|\tilde\tau|+\hbar^{\frac23}\right)^{-\frac12}$. In the second term this loss is compensated for by 
$(\partial_\beta\tau)(x\zeta,\beta)=O\left(\hbar (1+\alpha)^{-3}\right)$. For the integration by parts in $\tau$, we need the following estimates (by Lemma \ref{lem: Lemma 3.2 CDST}, Lemma \ref{lem:a1fine} and \eqref{eq:dritter T}):
\begin{align*}
	\frac{\partial\tilde{\tau}}{\partial\tau}\simeq 1,\quad \frac{\partial x}{\partial\tau}\simeq 1,\quad \left|\partial_{1}^{2}a_{1}(-\tilde{\tau},\beta)\right|=O\left(\hbar^{-1}\right),\quad \left|\partial_{1}\partial_{2}a_{1}(-\tilde{\tau},\beta)\right|=O\left(\hbar^{-\frac13}\beta^{-1}\right).
\end{align*}
The powers in $\hbar^{-1}$ above are compensated by the extra power in $\hbar$.
On the other hand, by \eqref{eq:tautil size} we have $\tilde\tau\ge \zeta-1$ if $\tau\ge0$. Thus, in the first term of~\eqref{eq:aj spur} with $j=1$ one has $(\zeta-1) (\partial_1 a_1)(-\tilde\tau,\beta) = O\left(\tilde\tau^{\frac12}\right)$ (note that $a_1$ in~\eqref{eq:F2 pre*} only involves $\tau\ge0$). In the subsequent integration by parts in $\tau$, the operator $\partial_{\tau}$ hits the first term on the RHS of \eqref{eq:aj spur}. Again by Lemma \ref{lem: Lemma 3.2 CDST}, Lemma \ref{lem:a1fine} and \eqref{eq:dritter T}, the most troublesome contribution is $$\hbar^{2}x\zeta\beta^{-1}(\zeta-1)(\partial_{1}^{2}a_{1})(-\tilde{\tau},\beta)\partial_{1}\tau(x\zeta,\beta)\cdot\frac{\partial\tilde{\tau}}{\partial\tau}$$, which is bounded in absolute value by
\begin{align*}
	\lesssim \hbar^{2}\beta^{-1}(\zeta-1)|\tilde{\tau}|^{\frac12-2}\lesssim \hbar^{\frac53}\beta^{-1},
\end{align*}
which is better than what we need. This argument for $j=1$ can be iterated: the term $\tilde\tau^{\frac12-k}$ arising in $(\hbar(\partial_\alpha+\partial_\beta))^m \hbar a_j(-\tilde\tau,\beta)$ from~\eqref{eq:a1 fine} are compensated for by $(\zeta-1)^k$. The details are as in~\eqref{eq:hell p 1} above, and we skip them. The end result is that 
this particular expression derived from~\eqref{eq:F2 pre*} gives $O(\hbar^m \alpha^{-m})$ (the integration by parts in~$\tau$, see~\eqref{eq:Fs osc+-}, generates a loss of at most $\hbar^{-\frac23}$ which is absorbed by the extra $\hbar$ factor at our disposal). Here we also used the estimate on the higher order derivatives: $\left|\frac{\partial^{m}\tilde{\tau}}{\partial\tau^{m}}\right|\lesssim 1$, which is seen by Lemma \ref{lem: Lemma 3.2 CDST}, \eqref{def R S T}.

For $j=0$, i.e., $a_0(-\tilde\tau,\beta)=a_0(-\tau(x\zeta,\beta),\beta)$, the mechanism is slightly different. It is only the first term of~\eqref{eq:aj spur} 
that requires a different argument. Using \eqref{eq:tautil size} once again, we distinguish three cases: (i) if $\tau\ge -c(\zeta-1)$ with some absolute constant $0<c\ll1$, then $\hbar^{\frac23} \gtrsim \tilde\tau\gtrsim \zeta-1$ and the same argument as for $a_1$ applies. (ii) if $\tau\le -C(\zeta-1)$ with some absolute constant $C\gg1$, then $\tilde\tau \simeq \tau \le -C(\zeta-1)$ whence $|\tilde\tau|\gtrsim \zeta-1$ and we can argue again as for $a_1$. (iii) if $-\tau\simeq \zeta-1$, then we use the factor $\Ai\left(-\hbar^{-\frac23} \tau\right)$ from~\eqref{eq:F2 pre*} to bound 
\EQ{\nn
	& \left| \Ai\left(-\hbar^{-\frac23} \tau\right) (\zeta-1) \left(\partial_1 a_0\right)(-\tilde\tau,\beta)  \partial_1\tau(x\zeta,\beta)\right|
	\lesssim \left(\zeta-1\right) \hbar^{-\frac13}  e^{-\frac{c}{\hbar}|\zeta-1|^{\frac32}} \lesssim 1
}
uniformly in the parameters, which is sufficient (we can ignore the $-\frac14$ power preceding the exponential bound of $\Ai$). For the higher derivatives we encounter terms of the form 
\EQ{\nn
	\left|  \Ai\left(-\hbar^{-\frac23} \tau\right) \hbar^{m+1} \beta^{-m}  (\zeta-1)^m    \left(\partial_1^m a_j\right)(-\tilde\tau,\beta) \right|  
	&\lesssim \hbar^{m+1}\alpha^{-m} (\zeta-1)^m  (|\tilde\tau|+\hbar^{\frac23})^{\frac12-m}    e^{-\frac{c}{\hbar}|\zeta-1|^{\frac32}} \\
	& \lesssim \hbar^{m+1}\alpha^{-m} (\zeta-1)^m \hbar^{-\frac23m}  e^{-\frac{c}{\hbar}|\zeta-1|^{\frac32}} \lesssim \hbar^{m+1}\alpha^{-m}
}
It remains to consider the derivatives of the Airy functions, which are governed by \eqref{eq:AiiBi tau'} and~\eqref{eq:AiiBi tildetau'}. 
In view of these equations,  we lose a factor of $\hbar$ by differentiating $\Ai\left(-\hbar^{\frac23}\tau\right)$ and the other Airy function. However, if the argument is~$\tau$ and not $\tilde\tau$, then we regain this factor in $\tau_\alpha$, resp., $\tau_\beta$ which are each bounded by $\hbar (1+\alpha)^{-3}$. This process can be repeated and we gain the desired $(\hbar/\alpha)^m$ {\em without} having to exploit the trace type derivative $\partial_\alpha+\partial_\beta$. In contrast, we do not gain $\hbar$ in~\eqref{eq:AiiBi tildetau'} and instead invoke the following: 
\EQ{ \label{eq:AiiBi tildetau' spur}
	\hbar\left(\partial_\alpha+\partial_\beta\right) \left(\Ai\left(-\hbar^{-\frac23} \tilde\tau\right)-i\Bi\left(-\hbar^{-\frac23} \tilde\tau\right)\right) &=  c\, \tilde \tau^{\frac12} \left(\tilde\tau_\alpha+ \tilde\tau_\beta\right)\left(1+\hbar O(\tilde\tau^{-\frac32}) \right) \left(\Ai\left(-\hbar^{-\frac23} \tilde\tau\right)-i\Bi\left(-\hbar^{-\frac23}\tilde \tau\right)\right)\\
	\tilde\tau_\alpha+ \tilde\tau_\beta	&=  (\partial_1\tau)(x\zeta,\beta)x\zeta\beta^{-1}(1-\zeta)+(\partial_2\tau)(x\zeta,\beta)
}
The second term is $\hbar/\alpha$ and again regains the essential $\hbar$ factor, while the first one is of size $(\zeta-1)/\alpha$. By \eqref{eq:32 diff} and~\eqref{eq:tautil size}, and assuming $\tau\ge -c(\zeta-1)$ for simplicity as in case (i) in the preceding paragraph, 
\[
\left|\partial_\tau \left( \tilde\tau^{\frac32}-\tau^{\frac32} \right)\right|\simeq (\zeta-1) \tilde\tau^{-\frac12}
\]
Integrating by parts as in \eqref{eq:Fs osc+-} in~$\tau$ therefore gains a factor of $\hbar \tilde\tau/\alpha$ (note the $\tilde \tau^{\frac12}$ in~\eqref{eq:AiiBi tildetau' spur}). As we iterate this process, we accumulate a factor of $(\hbar/\alpha)^m  \tilde\tau^m$, albeit at the expense of $\partial_\tau^m$ hitting all terms in the integrand other than the complex exponential, cf.~\eqref{eq:Fs osc+-}. If the operator $\hbar(\partial_{\alpha}+\partial_{\beta})$ hits the coefficient of the oscillatory Airy function in \eqref{eq:AiiBi tildetau' spur}, we again obtain a bound like $\hbar\tilde{\tau}/\alpha$. In fact we have
\begin{align*}
	\hbar(\partial_{\alpha}+\partial_{\beta})\left(\tilde{\tau}^{\frac12}\right)\cdot\left(\tilde{\tau}_{\alpha}+\tilde{\tau}_{\beta}\right)\simeq \hbar\tilde{\tau}^{-\frac12}\left((\partial_{1}\tau)(x\zeta,\beta)\right)^{2}x^{2}\zeta^{2}\beta^{-2}(1-\zeta)^{2}.
\end{align*}
Then upon integration by parts twice, we obtain a bound as (in view of $\hbar^{\frac23}\lesssim \tilde{\tau}$, since we are in the oscillatory regime)
\begin{align*}
	\hbar^{3}\left((\partial_{1}\tau)(x\zeta,\beta)\right)^{2}x^{2}\zeta^{2}\beta^{-2}\tilde{\tau}\lesssim \hbar^{2}\tilde{\tau}^{2}/\beta^{-2}.
\end{align*}
The case when the derivative $\hbar(\partial_{\alpha}+\partial_{\beta})$ hits the factor $\tilde{\tau}_{\alpha}+\tilde{\tau}_{\beta}$ is handled similarly, and it is in fact lower order.

Note that in order to obtain the factor involving $\Lambda=\hbar (\zeta-1)^{-1}$ which is necessary for the~$\Gamma$ bound~\eqref{F hbar bounds}, we then need to carry out one more integration by parts exactly as in~\eqref{eq:Fs osc+-}. If the $\partial_\tau^m$ hits $\tilde\tau^{-\frac14} \tilde a_1(-\tilde\tau,\beta)$ in~\eqref{eq:Fs osc+-}, then we lose $\tilde\tau^{-m}$ which is exactly compensated for by the $\tilde\tau^m$ that we gained. However, if $\partial_\tau^m$ hits $\tau^{-\frac14} \tilde a_1(-\tau,\alpha)$, then we lose $\tau^{-m}$. Recall our conditions  $|\xi-\eta|\le\hbar^{\frac23}\xi$ (which we have not used before) and $\xi\simeq\eta$ for~\eqref{eq:trace der}. They imply  that $0\le \zeta-1\lesssim \hbar^{\frac23}$, which  ensures that $\tau\simeq\tilde\tau$, since we can assume that $\tau\ge\hbar^{\frac23}$ (If $|\tau|\lesssim \hbar^{\frac23}$, then we don't even need the integration by parts argument).

Finally, Case~$A_1$ is easy due to the exponential gain inherent in the regime $0<x\ll1$, see the comment following~\eqref{eq:F2r2}. We can therefore pass 
any number of $\hbar\partial_\alpha$, resp.~$\hbar\partial_\beta$ derivatives into the integral in~\eqref{eq:F2r2} without any concern about losing powers of~$\hbar$. The are swallowed by the gain of $\exp(-c\hbar^{-1})$ which we have at our disposal.  

Now we turn to the case when $2\leq n\leq N_{0}$ for some fixed $N_{0}$. We use $F_{n}(\xi,\eta)$ to denote the off-diagonal kernel \eqref{off diagonal F pre}:
\begin{align}\label{off-diagonal F large hbar}
	F_{n}(\xi,\eta):=\int_{0}^{\infty}W^{+}_{n}(R)\phi_{n}(R,\xi)\phi_{n}(R,\eta) \,dR.
\end{align}
Here $\phi_{n}(R,\cdot)=w_{n}(\cdot)\phi(R,\cdot)$ is as in Proposition \ref{prop: FB match} and was constructed with respect to the measure $dR$, which appears in \eqref{off-diagonal F large hbar}. We start with the bound on $F_{n}(\xi,\eta)$ itself. When $0\leq R\xi^{\frac12}\leq R\eta^{\frac12}\lesssim 1$, we have, by \eqref{FB small R2xi large hbar}, 
\begin{align}\label{Fn 1}
	\begin{split}
	\left|\int_{0}^{\infty}W^{+}_{n}(R)\chi_{R\xi^{\frac12}\lesssim\frac12}\chi_{R\eta^{\frac12}\lesssim\frac12}\phi_{n}(R,\xi)\phi_{n}(R,\eta)\,dR\right|
	\lesssim &\xi^{-\frac14}\eta^{-\frac14}\cdot\xi^{\frac{n}{2}-\frac14}\eta^{\frac{n}{2}-\frac14}\int_{0}^{\eta^{-\frac12}}\frac{R^{2n-1}}{(1+R^{2})^{2}}dR\\
	\lesssim
	&\begin{cases}
&\eta^{-1}\cdot\left(\frac{\xi^{\frac12}}{\eta^{\frac12}}\right)^{n-1},\quad \textrm{if}\quad \eta\gtrsim1,\\
&\eta\cdot|\log\eta|\cdot\frac{\xi^{\frac12}}{\eta^{\frac12}},\quad \textrm{if}\quad \eta\ll1.
	\end{cases}
	\end{split}
\end{align}
Next we turn to the case $0<R\xi^{\frac12}\lesssim1\ll R\eta^{\frac12}$. We have, by \eqref{FB small R2xi large hbar} and \eqref{FB large R2xi large hbar} (Here we only look at $e^{iR\xi^{\frac12}}(1+g_{+}(R,\xi))$, and its conjugate can be handled similarly),
\begin{align}\label{Fn 2 pre}
	\begin{split}
	&\int_{0}^{\infty}W^{+}_{n}(R)\chi_{R\xi^{\frac12}\lesssim1}\chi_{R\eta^{\frac12}\gg 1}\phi_{n}(R,\xi)\phi_{n}(R,\eta)\,dR\\
	=&\int_{0}^{\infty}W^{+}_{n}(R)\chi_{R\xi^{\frac12}\lesssim1}\phi_{n}(R,\xi)\cdot\chi_{R\eta^{\frac12}\gg1} \eta^{-\frac14}\frac{a(\eta)}{|a(\eta)|}e^{iR\eta^{\frac12}}(1+g_{+}(R,\eta))\,dR\\
	=&\int_{0}^{\infty}W^{+}_{n}(R)\chi_{R\xi^{\frac12}\lesssim1}\phi_{n}(R,\xi)\cdot\chi_{R\eta^{\frac12}\gg1} \eta^{-\frac14}\frac{a(\eta)}{|a(\eta)|}\cdot \frac{1}{i\eta^{\frac12}}\partial_{R}\left(e^{iR\eta^{\frac12}}\right)(1+g_{+}(R,\eta))\,dR
	\end{split}
\end{align}
By Lemma \ref{lem:psin R2xi large} and \eqref{FB small R2xi large hbar}, upon a differentiation in $R$ on the factor $W^{+}_{n}(R)\phi_{n}(R,\xi)(1+g_{+}(R,\eta))$, we obtain a factor of $\xi^{\frac12}$, compared to the undifferentiated one. Therefore after integration by parts, we obtain
\begin{align}\label{Fn 2}
	\begin{split}
	\left|\int_{0}^{\infty}W^{+}_{n}(R)\chi_{R\xi^{\frac12}\lesssim1}\chi_{R\eta^{\frac12}\gg 1}\phi_{n}(R,\xi)\phi_{n}(R,\eta)\,dR\right|\lesssim \frac{\min\left\{1,\eta^{\frac32}\right\}}{\xi^{\frac14}\eta^{\frac14}}\cdot\left(\frac{\xi^{\frac12}}{\eta^{\frac12}}\right)^{N},\quad \textrm{for}\quad 0<N\leq n.
	\end{split}
\end{align}
Now we consider the case when $R\xi^{\frac12}\gg1$ and $R\eta^{\frac12}\gg1$. This is similar to the previous case, because when $\partial_{R}$ hits $W^{+}_{n}(R)\phi_{n}(R,\xi)(1+g_{+}(R,\eta))$, we again obtain an extra factor of $\xi^{\frac12}$. Therefore we have
\begin{align}\label{Fn 3}
	\begin{split}
		\left|\int_{0}^{\infty}W^{+}_{n}(R)\chi_{R\xi^{\frac12}\gg1}\chi_{R\eta^{\frac12}\gg1}\phi_{n}(R,\xi)\phi_{n}(R,\eta)\,dR\right|\lesssim \frac{\min\{1,\xi^{\frac32}\}}{\xi^{\frac14}\eta^{\frac14}}\cdot\left(\frac{\xi^{\frac12}}{\eta^{\frac12}}\right)^{N},\quad \textrm{for any}\quad N>0.
	\end{split}
\end{align}
Now we refine the estimates when $R\xi^{\frac12}\ll1$ and $\eta\geq 1$ to obtain the rapid decay factor $\left(\frac{\langle\xi\rangle}{\eta}\right)^{N}$. We again start with the expression
\begin{align*}
	F_{n}(\xi,\eta):=\int_{0}^{\infty}W_{n}^{+}(R)\phi_{n}(R,\xi)\phi_{n}(R,\eta)R\,dR.
\end{align*}
Here $\phi_{n}(R,\xi)$ is the Fourier basis with respect to the measure $R\,dR$, and $W^{+}_{n}(R)$ is a linear combination of $\frac{1}{(1+R^{2})^{k}}$ where $k\geq 2$.  Since now $2\leq n\leq N_{0}$, the $n$-dependence of $W^{+}_{n}(R)$ is not a concern here. Therefore without loss of generality, we consider
\begin{align}\label{Fn refine pre}
\begin{split}
	\int_{0}^{\infty}\frac{1}{(1+R^{2})^{k}}\phi_{n}(R,\xi)\phi_{n}(R,\eta)R\,dR.
	\end{split}
\end{align}
A direct calculation gives 
\begin{align}\label{Fn refine 1}
\begin{split}
	\eta\int_{0}^{\infty}\frac{1}{(1+R^{2})^{k}}\phi_{n}(R,
	\xi)\phi_{n}(R,\eta)R\,dR=&-\int_{0}^{\infty}\frac{1}{(1+R^{2})^{k}}\phi_{n}(R,\xi)H_{n}^{+}\left(\phi_{n}(R,\eta)\right)R\,dR\\
	=&\xi\int_{0}^{\infty}\frac{1}{(1+R^{2})^{k}}\phi_{n}(R,\xi)\phi_{n}(R,\eta)R\,dR\\
	&+2\int_{0}^{\infty}\partial_{R}\left(\frac{1}{(1+R^{2})^{k}}\right)\partial_{R}\left(\phi_{n}(R,\xi)\right)\phi_{n}(R,\eta)R\,dR\\
	&+\int_{0}^{\infty}R^{-1}\partial_{R}\left(\frac{1}{(1+R^{2})^{k}}\right)\phi_{n}(R,\xi)\phi_{n}(R,\eta)R\,dR\\
	&-\int_{0}^{\infty}\partial_{R}^{2}\left(\frac{1}{(1+R^{2})^{k}}\right)\phi_{n}(R,\xi)\phi_{n}(R,\eta)R\,dR.
	\end{split}
\end{align}
Since $\partial_{R}^{2}\left(\frac{1}{(1+R^{2})^{k}}\right)$ and $R^{-1}\partial_{R}\left(\frac{1}{(1+R^{2})^{k}}\right)$ are again linear combinations of the functions of the form $\frac{1}{(1+R^{2})^{k}}$, the third and the fourth terms on the RHS of \eqref{Fn refine 1} can be the same way as \eqref{Fn refine pre}. So we focus on the second term on the RHS of \eqref{Fn refine 1}. With $\phi_{0}(R), \phi_{j}(R), j\geq 1$ given by Lemma \ref{lem: phin small R2xi small} and $w_{n}(\xi)$ given by Proposition \ref{prop: FB match}, we have, for $R\xi^{\frac12}\ll1$,
\begin{align*}
	\phi_{n}(R,\xi)=\phi_{0}(R)w_{n}(\xi)\left[1+\sum_{j\geq1}\phi_{j}(R)\left(R^{2}\xi\right)\right].
\end{align*}
It follows that, for a constant $C_{n}$ depending on $n$,
\begin{align*}
	\chi_{R\xi^{\frac12}\ll1}\partial_{R}\left(\phi_{n}(R,\xi)\right)=&C_{n}R^{-1}\phi_{n}(R,\xi)+R\phi_{n}(R,\xi)\cdot\left[\sum_{j\geq 0}\psi_{j}(R)\left(R^{2}\xi\right)^{j}\right]
\end{align*}
for certain coefficient functions $\psi_{j}(R)$ satisfying analogous bounds as the $\phi_{j}(R)$. In fact, if $R\xi^{\frac12}\ll1$ is sufficiently small, the reciprocal of $\left[1+\sum_{j\geq1}\phi_{j}(R)(R^{2}\xi)^{j}\right]$ can be also expanded into a power series in $R^{2}\xi$. So we conclude that, for some constant $C_{k,n}$ depending on $k,n$, and some constant $D_{k}$ depending on $k$,
\begin{align}\label{Fn refine 2}
\begin{split}
	&\int_{0}^{\infty}\chi_{R\xi^{\frac12}\ll1}\partial_{R}\left(\frac{1}{(1+R^{2})^{k}}\right)\partial_{R}\left(\phi_{n}(R,\xi)\right)\phi_{n}(R,\eta)R\,dR\\
	=&\int_{0}^{\infty}\chi_{R\xi^{\frac12}\ll1}C_{k,n}\frac{1}{(1+R^{2})^{k+1}}\phi_{n}(R,\xi)\phi_{n}(R,\eta)R\,dR\\
	&+\int_{0}^{\infty}\chi_{R\xi^{\frac12}\ll1}D_{k}\frac{1}{(1+R^{2})^{k}}\phi_{n}(R,\xi)\left[\sum_{j\geq0}\psi_{j}(R)\left(R^{2}\xi\right)^{j}\right]\phi_{n}(R,\eta)R\,dR.
	\end{split}
\end{align}
Using the fact that $\psi_{j}(R)$ has rapid decay in $j$ and $\psi_{j}(R)$ is a function of $R^{2}$, the terms on the RHS of \eqref{Fn refine 2} can be handled the same way as \eqref{Fn refine pre}. 

The refined derivative estimates follow in the same way, using the relations:
\begin{align}\label{Fn refine deri}
	\begin{split}
	H_{n}^{+}\left(\partial_{\xi}^{k}\phi_{n}(R,\xi)\right)=&-\xi\partial_{\xi}^{k}\phi_{n}(R,\xi)-k\partial_{\xi}^{k-1}\phi_{n}(R,\xi),\quad \textrm{for}\quad k\geq1.
	\end{split}
\end{align}
\\
Now we turn to the derivatives of $F_{n}(\xi,\eta)$, if $R\xi^{\frac12}\lesssim1$ and $R\eta^{\frac12}\lesssim1$, then by \eqref{FB small R2xi large hbar}, upon a differentiation in $\xi^{\frac12}$ (or $\eta^{\frac12}$), we obtain an extra factor of $\xi^{-\frac12}$ (or $\eta^{-\frac12}$). This gives all the derivative estimates in this regime. If $R\xi^{\frac12}\lesssim1$ and $R\eta^{\frac12}\gg1$, we use the fact 
\begin{align}\label{deri transfer FB}
R\partial_{R}\left(e^{iR\eta^{\frac12}}\right)=\eta^{\frac12}\partial_{\eta^{\frac12}}\left(e^{iR\eta^{\frac12}}\right),\quad \Rightarrow\quad \partial_{\eta^{\frac12}}\left(e^{iR\eta^{\frac12}}\right)=R\eta^{-\frac12}\partial_{R}\left(e^{iR\eta^{\frac12}}\right),
\end{align}
and integration by parts to obtain the derivative estimates in this regime. Finally we look at the more delicate case $1\ll R\xi^{\frac12}\leq R\eta^{\frac12}$.  For the pure $\partial_{\eta^{\frac12}}$-derivatives of $F_{n}(\xi,\eta)$, we again use \eqref{deri transfer FB} and integration by parts to obtain
\begin{align}\label{deri Fn eta hard 1}
	\begin{split}
	\left|\partial_{\eta^{\frac12}}^{k}F_{n}(\xi,\eta)\right|\lesssim \max\{\xi^{-\frac{k}{2}},1\}\cdot\Gamma_{n} \cdot\frac{\xi^{\frac{k}{2}}}{\eta^{\frac{k}{2}}},\quad \textrm{for}\quad k=1,2\quad \textrm{and}\quad 1\ll R\xi^{\frac12}\leq R\eta^{\frac12}.
	\end{split}
\end{align}
Here $\Gamma_{n}$ is the bound for the undifferentiated $F_{n}(\xi,\eta)$ in the same regime. For the derivatives in $\xi^{\frac12}$. We distinguish the discussion according to the order of the derivative in $\xi^{\frac12}$. We first consider the bound for $\partial^{k}_{\xi^{\frac12}}F_{n}(\xi,\eta),k=1,2$. It is straightforward to see that for $k=1,2$,
\begin{align}\label{deri Fn xi hard 1}
	\left|\partial^{k}_{\xi^{\frac12}}F_{n}(\xi,\eta)\right|\lesssim \xi^{-\frac{k}{2}}\cdot \Gamma_{n}\quad\textrm{for}\quad \xi\ll1,\quad \textrm{and}\quad \left|\partial^{k}_{\xi^{\frac12}}F_{n}(\xi,\eta)\right| \lesssim \Gamma_{n}\quad \textrm{for}\quad \xi\gtrsim 1.
\end{align}
$\partial_{\xi^{\frac12}}\partial_{\eta^{\frac12}}F_{n}(\x,\eta)$ is even more delicate. Without loss of generality, we consider the integral:
\begin{align}\label{deri Fn xi eta pre}
\begin{split}
\int_{0}^{\infty}W^{+}_{n}(R)\chi_{R\xi^{\frac12}\gg1}\chi_{R\eta^{\frac12}\gg1}\xi^{-\frac14}\eta^{-\frac14}\partial_{\xi^{\frac12}}\left(e^{iR\xi^{\frac12}}(1+g_{+}(R,\xi))\right)\cdot\partial_{\eta^{\frac12}}\left(e^{iR\eta^{\frac12}}(1+g_{+}(R,\eta))\right)\,dR,
\end{split}
\end{align}
and we only focus on the contribution when the derivatives $\partial_{\xi^{\frac12}}, \partial_{\eta^{\frac12}}$ hit the oscillatory factors. Therefore the integral in consideration is given by (omitting the constant coefficients)
\begin{align}\label{deri Fn xi eta}
	\begin{split}
	&\left|\partial_{\xi^{\frac12}}\partial_{\eta^{\frac12}}F_{n}(\xi,\eta)\right|\\
\simeq &	\xi^{-\frac14}\eta^{-\frac14}\left|\int_{0}^{\infty}W^{+}_{n}(R)\chi_{R\xi^{\frac12}\gg1}\chi_{R\eta^{\frac12}\gg1}R\cdot e^{iR\xi^{\frac12}}(1+g_{+}(R,\xi))\cdot\partial_{\eta^{\frac12}}\left(e^{iR\eta^{\frac12}}\right)\cdot(1+g_{+}(R,\eta))\,dR\right|\\
=&\xi^{-\frac14}\eta^{-\frac14}\left|\int_{0}^{\infty}W^{+}_{n}(R)\chi_{R\xi^{\frac12}\gg1}\chi_{R\eta^{\frac12}\gg1}R\cdot e^{iR\xi^{\frac12}}(1+g_{+}(R,\xi))\cdot \frac{R}{\eta^{\frac12}}\partial_{R}\left(e^{iR\eta^{\frac12}}\right)\cdot(1+g_{+}(R,\eta))\,dR\right|\\\lesssim &\max\{\xi^{-1},1\}\cdot\Gamma_{n}\cdot\frac{\xi^{\frac12}}{\eta^{\frac12}}.
	\end{split}
\end{align}
using integration by parts. At the end we consider the third order derivatives in the regime $1\ll R\xi^{\frac12}\leq R\eta^{\frac12}$. In this case we will encounter a logarithmic divergence any way, so we bound the contributions from all the possible third order derivatives by the absolute value of (Again here we only consider the phase $e^{iR(\xi^{\frac12}-\eta^{\frac12})}$, and the contribution from $e^{iR(\xi^{\frac12}+\eta^{\frac12})}$ can be treated similarly)
\begin{align}\label{deri Fn 3rd}
	\begin{split}
	&\xi^{-\frac14}\eta^{-\frac14}\left(\int_{\xi^{-\frac12}}^{|\xi^{\frac12}-\eta^{\frac12}|^{-1}}R^{-1}\,dR+\int_{|\xi^{\frac12}-\eta^{\frac12}|^{-1}}^{\infty}R^{-1}e^{iR(\eta^{\frac12}-\xi^{\frac12})}\,dR\right)\\
	\lesssim &\xi^{-\frac14}\eta^{-\frac14}\left(1+\left|\log\left(\xi^{\frac12}\left|\xi^{\frac12}-\eta^{\frac12}\right|^{-1}\right)\right|\right)
	\end{split}
\end{align}
\end{proof}
\subsubsection{Negative frequency}
Next we turn to the case when $n\leq -2$. Let us define the operator $\calK_{\thbar}$ as
\begin{align}\label{def calK thbar}
	\widehat{R\partial_{R}u}=-2\xi\partial_{\xi}\widehat{u}+\calK_{\thbar}\widehat{u}.
\end{align}
Here the Fourier transform is taken with respect to the Fourier basis $\phi(R,\xi;\thbar)=R^{-\frac12}\phi_{*}(R,\xi;\thbar)$ in $L^{2}(R\,dR)$ with $\phi_{*}(R,\xi;\thbar)$ being the Fourier basis in $L^{2}(dR)$ (See Propositions \ref{prop: FB nega n large hbar} and \ref{prop:DFT nlarge negative}). Let $\rho_{-n}(d\xi)$ be the spectral measure associated with $\phi_{*}(R,\xi;\thbar)$, where $\thbar=\frac{1}{1-n}$. As for the positive frequencies, $\rho_{-n}(d\xi)$ is also the spectral measure of $\phi(R,\xi;\thbar)$. 

Similar as Proposition \ref{prop: K operator}, we have
\begin{proposition}\label{prop: K operator negative}
	For any $-1\leq \thbar<0$ the operator $\calK_{\thbar}$ is given by
	\begin{align*}
		\calK_{\thbar}f(\xi)=-\left(2f(\eta)+\frac{\eta\left(\frac{d\rho_{-n}}{d\eta}(\eta)\right)^{\prime}}{\frac{d\rho_{-n}(\eta)}{d\eta}}f(\eta)\right)\delta(\eta-\xi)+\left(\calK_{\thbar}^{(0)}f\right)(\xi)
	\end{align*}
	where the off-diagonal part $\calK^{(0)}_{\thbar}$ has a kernel $K_{0}(\xi,\eta;\thbar)$ given by
	\EQ{\label{eq:K0 kern negative}
		K_{0}(\xi,\eta;\thbar)= \frac{\frac{d\rho_{-n}(\eta)}{d\eta}}{\xi-\eta}F(\xi,\eta;\thbar)
	}
	and the symmetric function $F(\xi,\eta;\thbar)$ satisfies (for any $0\leq k\leq k_{0}$ and sufficiently small $\hbar=\hbar(k_{0})$, where $k_{0}$ is arbitrary but fixed, and $\xi\leq \eta$)
	\begin{align}\label{F thbar bounds}
		\begin{split}
			|F(\xi,\eta;\thbar)|\lesssim& \left(\thbar\xi^{\frac12}\right)^{-1}\min\left\{1,\left(\thbar\xi^{\frac12}\right)^{3}\right\}\cdot G:=\Gamma_{\thbar}.
		\end{split}
	\end{align}
	with 
	\begin{align}\label{F thbar bounds auxi}
		G:=\begin{cases}
			&\min\left\{1,\left(\thbar\xi^{\frac12}\right)^{\frac14}\left|\eta^{\frac12}-\xi^{\frac12}\right|^{-\frac14}\right\},\quad \textrm{for}\quad \left|\frac{\eta^{\frac12}}{\xi^{\frac12}}-1\right|\lesssim1,\\
			&\hbar \left(\thbar\xi^{\frac12}\right)^{-1}\min\left\{1,\thbar\xi^{\frac12}\right\}\cdot\left(\frac{\xi^{\frac12}}{\eta^{\frac12}}\right)^{k},\quad \textrm{for}\quad \left|\frac{\eta^{\frac12}}{\xi^{\frac12}}-1\right|\gg1.
		\end{cases}
	\end{align}
	For $\thbar\gtrsim 1$, we have, for any $N>0$
	\begin{align}\label{F -n bounds}
		\begin{split}
			|F(\xi,\eta;\thbar)|\lesssim \Gamma_{n}:= \begin{cases}
				&\eta|\log\eta|\cdot\frac{\xi^{\frac12}}{\eta^{\frac12}},\quad \textrm{for}\quad \eta\ll1,\\
				&\max\{\xi^{-\frac12}\eta^{-\frac12},\eta^{-\frac14}\xi^{-\frac14}\}\cdot\left(\frac{\xi^{\frac12}}{\eta^{\frac12}}\right)^{n}\cdot\left(\frac{\langle\xi\rangle}{\eta}\right)^{N},\quad \textrm{for}\quad \eta\gtrsim 1
			\end{cases}.
		\end{split}
	\end{align}

	For the derivatives of $F(\xi,\eta;\thbar)$, we have, for $k_{1}+k_{2}\leq 2$
	\begin{align}\label{deri F thbar and n bounds 3}
		\left|\partial_{\xi^{\frac12}}^{k_{1}}\partial_{\eta^{\frac12}}^{k_{2}}F(\xi,\eta;\thbar)\right|\lesssim \max\left\{1,\left(\thbar\xi^{\frac12}\right)^{-k_{1}}\right\}
		\cdot\max\left\{1,\left(\thbar\eta^{\frac12}\right)^{-k_{2}}\right\}\cdot \Gamma
	\end{align}
	Here $\Gamma$ is $\Gamma_{\thbar}$ if $\hbar\ll1$, and $\Gamma_{n}$ if $\thbar\gtrsim1$. The following estimate holds for trace-type derivatives for $\thbar\ll1$:
	\EQ{\label{eq:trace der negative}
		\left|\left(\partial_{\xi^{\frac12}}+\partial_{\eta^{\frac12}}\right)^{k}F(\xi,\eta;\thbar)\right| &\le C_k\, \xi^{-\frac{k}{2}} \cdot\Gamma,\quad \textrm{if}\quad \xi\simeq \eta,\; |\eta-\xi|\le \thbar^{\frac23}\xi
	}
	for all $k\ge0$. 
\end{proposition}
\begin{proof}
	Again we will often drop $\thbar$ from the notation. The proof is identical to that of Proposition \ref{prop: K operator} except that for $2\leq n\leq N_{0}$ we use Proposition \ref{prop: FB nega n large hbar}, since the Airy functions are only used to construct $\phi_{-n}(R,\xi)$ for $n\geq N_{0}$. A similar calculation as in Proposition \ref{prop: K operator} gives
	\begin{align*}
		\calK f(\eta)=&\left\langle\int_{0}^{\infty}f(\xi)[R\partial_{R}-2\xi\partial_{\xi}]\phi(R,\xi)\rho(\xi)\,d\xi,\phi(R,\eta)\right\rangle_{L^{2}_{R\,dR}}\\
		&-2\left(1+\frac{\eta\rho'(\eta)}{\rho(\eta)}\right)f(\eta).
	\end{align*}
	To extract any $\delta(\xi-\eta)$ from the first line above we recall from Proposition \ref{prop: FB nega n large hbar},
	\begin{align*}
		\phi(R,\xi)=&2\xi^{-\frac14}R^{-\frac12}\Re\left(a_{-}(\xi)e^{iR\xi^{\frac12}}\left(1+g^{-}(R,\xi)\right)\right)
	\end{align*}
	with
	\begin{align*}
		&g^{-}(R,\xi)=\frac{1}{R\xi^{\frac12}}\left(C\,i+O\left(\frac{1}{1+R^{2}}\right)\right)+O\left(R^{-2}\xi^{-1}\right),\\
		&  \left|(R\partial_{R})^{k}g^{-}(R,\xi)\right|\leq c_{k}(R\xi^{\frac12})^{-1},\quad \textrm{and}\quad \left|(\xi\partial_{\xi})^{k}g^{-}(R,\xi)\right|\leq c_{k}(R\xi^{\frac12})^{-1}.
	\end{align*}
	Therefore we have
	\begin{align*}
		\left(R\partial_{R}-2\xi\partial_{\xi}\right)\phi(R,\xi)=&-4\xi^{-\frac14}R^{-\frac12}\Re\left(\xi a^{\prime}_{-}(\xi)e^{iR\xi^{\frac12}}\left(1+g^{-}(R,\xi)\right)\right)\\
		&+2\xi^{-\frac14}R^{-\frac12}\Re\left(a_{-}(\xi)e^{iR\xi^{\frac12}}\left(R\partial_{R}-2\xi\partial_{\xi}\right)g^{-}(R,\xi)\right).
	\end{align*}
	According to the profile of $g^{-}(R,\xi)$, we have
	\begin{align*}
		\left|\left(R\partial_{R}-2\xi\partial_{\xi}\right)g^{-}(R,\xi)\right|\lesssim&\frac{1}{R\xi^{\frac12}}O\left(\frac{R^{2}}{(1+R^{2})^{2}}\right) .
	\end{align*}
	Using that $\Re z\Re w=\frac12\left(\Re(zw)+\Re(z\overline{w})\right)$ we infer that the $\delta$ measure on the diagonal in the integral
	\begin{align*}
		\lim_{A\rightarrow\infty}\int_{0}^{A}\left(R\partial_{R}-2\xi\partial_{\xi}\right)\phi(R,\xi)\phi(R,\eta)R\,dR
	\end{align*}
	comes from the expression
	\begin{align*}
		&(\xi\eta)^{-\frac14}\lim_{L\rightarrow\infty}\Re\int_{0}^{\infty}-4\left(\xi a^{\prime}_{-}(\xi)\overline{a_{-}(\eta)}e^{iR(\xi^{\frac12}-\eta^{\frac12})}\left(1+g^{-}(R,\xi)\right)\left(\overline{1+g^{-}(R,\eta)}\right)\right)\chi_{1}(R)\chi_{2}(R/L)dR\\
		&+(\xi\eta)^{-\frac14}\lim_{L\rightarrow\infty}\Re\int_{0}^{\infty}2e^{iR(\xi^{\frac12}-\eta^{\frac12})}a_{-}(\xi)\overline{a_{-}(\eta)}\left(R\partial_{R}-2\xi\partial_{\xi}\right)g^{-}(R,\xi)\cdot\overline{\left(1+g^{-}(R,\eta)\right)}\cdot \chi_{1}(R)\chi_{2}(R/L)\,dR
	\end{align*}
	Here $\chi_{1}$ is a smooth cutoff which equals $1$ near $\infty$ and which vanishes near $0$, and $\chi_{2}=1-\chi_{1}$. According to the estimate on $\left(R\partial_{R}-2\xi\partial_{\xi}\right)g^{-}(R,\xi)$ and the profile of $g^{-}$, the $\delta$ measure comes from
	\begin{align*}
		&(\xi\eta)^{-\frac14}\lim_{L\rightarrow\infty}\Re\int_{0}^{\infty}-4\xi a^{\prime}_{-}(\xi)\overline{a_{-}(\eta)}e^{iR(\xi^{\frac12}-\eta^{\frac12})}\chi_{1}(R)\chi_{2}(R/L)dR\\
		=&-4\pi(\xi\eta)^{-\frac14}\Re\left(\xi a^{\prime}_{-}(\xi)\overline{  a_{-}(\eta)}\right)\delta(\xi^{\frac12}-\eta^{\frac12})\\
		=&-8\pi\xi\Re\left(a^{\prime}_{-}(\xi)\overline{  a_{-}(\xi)}\right)\delta(\xi-\eta)=-4\pi\xi\frac{d}{d\xi}\left(|a_{-}(\xi)|^{2}\right)=-\xi\frac{d}{d\xi}\left(\rho(\xi)^{-1}\right)\delta(\xi-\eta).
	\end{align*}
	This together with the calculation in \eqref{operator K post} gives the desired result.
\end{proof}

\subsection{Transference operator at angular momentum $n=1$}
The operator $\calK^{1}$ is defined as
\begin{align}\label{def trans n1}
\widehat{R\partial_{R}u}=-2\xi\partial_{\xi}\widehat{u}+\calK^{1}\widehat{u}.
\end{align}
In this section, the Fourier transform is taken with respect to the Fourier basis $\tilde{\phi}(R,\xi)=R^{-\frac12}\phi(R,\xi)$ in $L^{2}(RdR)$ where $\phi(R,\xi)$ is introduced in \eqref{modified funda sys}. Let $\rho_{1}(d\xi)$ be the spectral measure associated with $\phi(R,\xi)$. Using an argument similar to the $n\geq2$ case, we obtain the Plancherel identity
\begin{align*}
	\|f\|^{2}_{L^{2}_{RdR}}=\left\|\langle f,R^{-\frac12}\phi\rangle_{L^{2}_{RdR}}\right\|^{2}_{L^{2}_{\rho_{1}(d\xi)}},
\end{align*}
and it follows that $\rho_{1}(d\xi)$ is also the spectral measure for $\tilde{\phi}(R,\xi)$. The main result of this section is:
\begin{proposition}\label{prop: K operator n1}
 The operator $\calK^{1}$ is given by
 \begin{align*}
 \left( \calK^{1} f\right)(\xi)=-\left(2f(\eta)+\frac{\eta\rho'_{1}(\eta)}{\rho_{1}(\eta)}f(\eta)\right)\delta(\eta-\xi).
 \end{align*}
In other words, the off-diagonal part of $\calK^{1}$ vanishes. Here for simplicity we write the spectral measure density $\frac{d\rho_{1}(\xi)}{d\xi}$ as $\rho_{1}(\xi)$, and recall from \eqref{eq:an1} that $\rho_{1}(\xi)=\frac{2\xi}{\pi^{2}}$ for $\xi>0$.
\end{proposition}
\begin{proof}
Similar to the proof of Proposition \ref{prop: K operator}, for a function $f\in C_{0}^{\infty}\left((0,\infty)\right)$, we define $u(R)=\int_{0}^{\infty}\tilde{\phi}(R,\xi)f(\xi)\rho_{1}(\xi)d\xi$. A computation similar to \eqref{eq:main int} gives
\begin{align}\label{eq:main int 1}
\begin{split}
\widehat{R\partial_{R}u}(\xi)=&\left\langle\int_{0}^{\infty}[R\partial_{R}-2\eta\partial_{\eta}]\tilde{\phi}(R,\eta)f(\eta)\rho_{1}(\eta)d\eta,\tilde{\phi}(R,\xi)\right\rangle_{L^{2}(RdR)}\\
	&-2f(\xi)-2\frac{\xi\rho'_{1}(\xi)}{\rho_{1}(\xi)}f(\xi)-2\xi f'(\xi).
\end{split}
\end{align}
 Here we again used the representation \eqref{delta kernel} for Dirac measure. It follows that
 \begin{align}\label{operator K n1 pre}
 \begin{split}
  \left(\calK^{1} f\right)(\eta)=&\left\langle\int_{0}^{\infty}f(\xi)[R\partial_{R}-2\xi\partial_{\xi}]\tilde{\phi}(R,\xi)\rho_{1}(\xi)d\xi,\tilde{\phi}(R,\eta)\right\rangle_{L^{2}_{RdR}}\\
 &-2\left(1+\frac{\eta\rho'_{1}(\eta)}{\rho_{1}(\eta)}\right)f(\eta).
 \end{split}
 \end{align}
 To extract any $\delta(\xi-\eta)$-contribution from the $RdR$-integral above, we recall from \eqref{eq:Hosc}, \eqref{eq:phian1} and \eqref{eq:an1} (Here we adopt the notations from these equations):
 \begin{align}\label{eq:phi wieder n1}
 	\begin{split}
 	\tilde{\phi}(R,\xi)=&2R^{-\frac12}\Re\left(a(\xi)\psi(R,\xi)\right)\\
 	=&2\Re\left(-e^{-\frac{\pi i}{4}}\frac{\pi}{4i\xi^{\frac34}}R^{-\frac12}e^{iR\xi^{\frac12}}\left(\sqrt{\frac{2}{\pi}}+\frac{a_{1}}{R\xi^{\frac12}}+O((R^{2}\xi)^{-1})\right)\right).
 	\end{split}
 \end{align}
 Here $O((R^{2}\xi)^{-1})$ obeys the ``symbol-type" behavior upon differentiation. In view of \eqref{rela R xi deri}, applying the differential operator $R\partial_{R}-2\xi\partial_{\xi}$ to \eqref{eq:phi wieder n1} yields
 \begin{align}\label{eq:phi wieder n1 deri}
 \begin{split}
 	(R\partial_{R}-2\xi\partial_{\xi})\tilde{\phi}(R,\xi)=&-2\Re\left(e^{-\frac{\pi i}{4}}\frac{3\pi}{8i\xi^{\frac34}}R^{-\frac12}e^{iR\xi^{\frac12}}\left(\sqrt{\frac{2}{\pi}}+\frac{a_{1}}{R\xi^{\frac12}}+O((R^{2}\xi)^{-1})\right)\right)\\
 	&+\Re\left(e^{-\frac{\pi i}{4}}\frac{\pi}{4i\xi^{\frac34}}R^{-\frac12}e^{iR\xi^{\frac12}}\left(\sqrt{\frac{2}{\pi}}+\frac{a_{1}}{R\xi^{\frac12}}+O((R^{2}\xi)^{-1})\right)\right)\\
 	=&-\Re\left(e^{-\frac{\pi i}{4}}\frac{\pi}{2i\xi^{\frac34}}R^{-\frac12}e^{iR\xi^{\frac12}}\left(\sqrt{\frac{2}{\pi}}+\frac{a_{1}}{R\xi^{\frac12}}+O((R^{2}\xi)^{-1})\right)\right)
 	\end{split}
 \end{align}
 Now we again use the fact $\Re z\Re w=\frac12\left(\Re(zw)+\Re(z\overline{w})\right)$ to find the $\delta$ measure on the diagonal in the integral
 \begin{align*}
 	\lim_{A\rightarrow\infty}\int_{0}^{A}[R\partial_{R}-2\xi\partial_{\xi}]\tilde{\phi}(R,\xi)\tilde{\phi}(R,\eta)RdR.
 \end{align*}
 The contribution from the $\delta$ measure is given by
 \begin{align}\label{eq:Ixieta n1}
 	\begin{split}
 	\lim_{L\rightarrow\infty}\frac{\pi^{2}}{8}\cdot\frac{2}{\pi}\xi^{-\frac34}\eta^{-\frac34}\Re\int_{0}^{\infty}e^{iR(\xi^{\frac12}-\eta^{\frac12})}\chi_{1}(R)\chi_{2}(R/L)dR:=I(\xi,\eta;1).
 	\end{split}
 \end{align}
 Here the cutoff functions $\chi_{1},\chi_{2}$ are the same as in \eqref{eq:Ixieta}. The contribution from the profile $O((R^{2}\xi)^{-1})$ in \eqref{eq:phi wieder n1}-\eqref{eq:phi wieder n1 deri} is integrable in $R$, therefore is not part of the $\delta$ measure. The contribution from $\frac{a_{1}}{R\xi^{\frac12}}$ in \eqref{eq:phi wieder n1}-\eqref{eq:phi wieder n1 deri} gives an integral of the form \eqref{converge to L2}, therefore is not part of the $\delta$ measure either. Then a standard calculation gives
 \begin{align*}
 I(\xi,\eta;1)=\frac{\pi^{2}}{4}\xi^{-\frac34}\eta^{-\frac34}\delta(\xi^{\frac12}-\eta^{\frac12})=\frac{\pi^{2}}{2}\xi^{-\frac14}\eta^{-\frac34}\delta(\xi-\eta),
 \end{align*}
 which in term gives the diagonal contribution
 \begin{align*}
&\left\langle\int_{0}^{\infty}f(\xi)[R\partial_{R}-2\xi\partial_{\xi}]\tilde{\phi}(R,\xi)\rho_{1}(\xi)d\xi,\tilde{\phi}(R,\eta)\right\rangle_{L^{2}_{RdR}}\\
=& \left\langle\frac{\pi^{2}}{2}\xi^{-\frac14}\eta^{-\frac34}\delta(\xi-\eta), f(\xi)\rho_{1}(\xi)\right\rangle_{d\xi}=f(\eta)=\frac{\eta\rho'_{1}(\eta)}{\rho_{1}(\eta)}f(\eta).
 \end{align*}
 This gives the desired representation for the diagonal part of $\calK^{1}$. For the off-diagonal part, we proceed similarly as in Proposition \ref{prop: K operator} to obtain that the kernel of the off-diagonal part of $\calK^{1}$ is given by
 \begin{align*}
 	F(\xi,\eta;1)=\int_{0}^{\infty}W_{1}(R)\phi(R,\xi)\phi(R,\eta)RdR
 \end{align*}
 where 
 \begin{align*}
 	W_{1}(R):=[H_{1}^{+},R\partial_{R}]-2H_{1}^{+}.
 \end{align*}
 But a direct computation shows that $[H_{1}^{+},R\partial_{R}]-2H_{1}^{+}=0$. Therefore the off-diagonal part vanishes and the proof is completed.
\end{proof}
\subsection{Transference operator at angular momentum $n=-1$}\label{subsec:n=-1transference}
The operator $\calK_{-1}$ is defined as
\begin{align}\label{def trans n-1}
	\widehat{R\partial_{R}u}=-2\xi\partial_{\xi}\widehat{u}+\calK_{-1}\widehat{u}.
\end{align}
In this section, the Fourier transform is taken with respect to the Fourier basis $\phi_{-1}(R,\xi)=\xi^{-1}\calD_{-}\left(R^{-\frac12}\Phi^{-1}(R,\xi)\right)$ where  $\Phi^{-1}(R,\xi)$ is constructed in Proposition \ref{prop: basis n-1 small}. Let $\rho_{-1}(d\xi)$ be the spectral measure associated with $\Phi^{-1}(R,\xi)$. Let $f\in L^{2}(RdR)$ as in Lemma \ref{lem:H-1conv} and recall from \eqref{eq:D_f}, 
\begin{align*}
	\calD_{-}f(R)=\int_{0}^{\infty}x_{-1}(\xi)\phi_{-1}(R,\xi)\trho_{-1}(d\xi),\quad x_{-1}(\xi)=\left\langle \calD_{-}f,\phi_{-1}(R,\xi)\right\rangle_{L^{2}_{RdR}},
\end{align*}
where $\trho_{-1}(d\xi)=\xi\rho_{-1}(d\xi)$. Moreover, Proposition \ref{prop: spectral measure -1 diff} gives the Plancherel identity:
\begin{align*}
	\|\calD_{-}f\|^{2}_{L^{2}(RdR)}=\left\|\left\langle\calD_{-}f,\phi_{-1}(R,\xi)\right\rangle\right\|_{L^{2}(\trho_{-1})}^{2}.
\end{align*}
The main result of this section is 
\begin{proposition}\label{prop: K operator n-1}
	The operator $\calK_{-1}$ is given by 
	\begin{align*}
		\calK_{-1}f(\xi)=-\left(2f(\eta)+\frac{\eta\trho'_{-1}(\eta)}{\trho_{-1}(\eta)}f(\eta)\right)\delta(\eta-\xi)+\calK_{-1}^{0}f(\xi)
	\end{align*}
	where the off-diagonal part $\calK^{0}_{-1}$ has a kernel $K_{0}(\xi,\eta;-1)$ given by
	\begin{align}\label{eq:K0-1 kern}
		K_{0}(\xi,\eta;-1)=\frac{\trho_{-1}(\eta)}{\xi-\eta}F_{-1}(\xi,\eta)
	\end{align}
	where the symmetric function $F_{-1}(\xi,\eta)$ satisfies (for $\xi\leq \eta$, and any $N\in\bbN$)
	\begin{align}\label{F-1 bounds}
	\begin{split}
		|F_{-1}(\xi,\eta)|\lesssim 
		&\left(\frac{\langle\xi\rangle}{\langle\eta\rangle}\right)^{N}\left(\langle\eta\rangle^{-1}\cdot\frac{\min\{\xi^{\frac12},\xi^{-2}\}}{\eta^{\frac12}}+\langle\eta\rangle^{-1}\min\{1,\eta^{-\frac34}\}\langle\xi\rangle^{-\frac74}\right)\\&+\langle\xi\rangle^{-1}\langle\eta\rangle^{-1}\eta^{-\frac34}\cdot\min\{\xi^{\frac34},\xi^{-\frac34}\}\left(\frac{\xi^{\frac12}}{\eta^{\frac12}}\right)^{N}:=\Gamma_{-1}.
		\end{split}
	\end{align}
	For the derivatives of $F_{-1}(\xi,\eta)$, we have
	\begin{align}\label{F-1 deri bounds}
		\begin{split}
	&	\left|\partial_{\xi^{\frac12}}F_{-1}(\xi,\eta)\right|\lesssim \left(1+\xi^{-\frac12}\right)\cdot \Gamma_{-1},\quad \left|\partial_{\eta^{\frac12}}F_{-1}(\xi,\eta)\right|\lesssim \frac{\Gamma_{-1}}{\eta^{\frac12}}\left(1+\xi^{\frac12}\right),\\
	&\left|\partial^{2}_{\xi^{\frac12}}F_{-1}(\xi,\eta)\right|\lesssim \left(1+\xi^{-1}\right)\cdot\Gamma_{-1},\quad \left|\partial^{2}_{\eta^{\frac12}}F_{-1}(\xi,\eta)\right|\lesssim \frac{\Gamma_{-1}}{\eta}\cdot\left(1+\xi\right),\\
	&\left|\partial_{\xi^{\frac12}}\partial_{\eta^{\frac12}}F_{-1}(\xi,\eta)\right|\lesssim \frac{\Gamma_{-1}}{\eta^{\frac12}}\left(\xi^{\frac12}+\xi^{-\frac12}\right).
		\end{split}
	\end{align}
\end{proposition}

\begin{proof}
	Similar to the proof of Proposition \ref{prop: K operator} and Proposition \ref{prop: K operator n1}, for a function $f\in C^{\infty}_{0}((0,\infty))$, we define $u(R)=\int_{0}^{\infty}\phi_{-1}(R,\xi)f(\xi)\rho_{-1}(\xi)d\xi$. A computation similar to \eqref{eq:main int} and \eqref{operator K pre} gives
	\begin{align}\label{operator K n-1 pre}
		\begin{split}
		(\calK_{-1}f)(\eta)=&\left\langle\int_{0}^{\infty}f(\xi)[R\partial_{R}-2\xi\partial_{\xi}]\phi_{-1}(R,\xi)\trho_{-1}(\xi)d\xi,\phi_{-1}(R,\eta)\right\rangle_{L^{2}_{RdR}}\\&-2\left(1+\frac{\eta\trho'_{-1}(\eta)}{\trho_{-1}(\eta)}\right)f(\eta).
		\end{split}
	\end{align}
	To extract the $\delta$ measure, we recall from Proposition \ref{prop: basis n-1 large}, Proposition \ref{prop: spectral measure n-1} and \eqref{D minus and star}:
	\begin{align}\label{eq:phi wieder n-1}
		\begin{split}
		\phi_{-1}(R,\xi)=&\xi^{-1}\calD_{-}\left(2R^{-\frac12}\Re\left(a_{-1}(\xi)\Psi^{+}_{-1}(R,\xi)\right)\right)\\=&2\xi^{-1}\calD_{-}\left(\Re\left(R^{-\frac12}a_{-1}(\xi)\xi^{-\frac14}e^{iR\xi^{\frac12}}\left(1-\frac{i}{8R\xi^{\frac12}}+\frac{1}{R\xi^{\frac12}}O\left(\frac{1}{1+R^{2}}\right)+O\left(\frac{1}{R^{2}\xi}\right)\right)\right)\right)\\
		=&2\Re\left(R^{-\frac12}i\xi^{-\frac34}a_{-1}(\xi)e^{iR\xi^{\frac12}}\left(1+O\left(\frac{1}{R\xi^{\frac12}}\right)+O\left(\frac{1}{R^{2}}\right)\right)\right)
		\end{split}
	\end{align}
	Again the profiles $O\left(\frac{1}{R\xi^{\frac12}}\right)$ and $O\left(\frac{1}{R^{2}}\right)$ obey ``symbol-type" behavior upon differentiation. Applying the differential operator $R\partial_{R}-2\xi\partial_{\xi}$ to \eqref{eq:phi wieder n-1} yields
	\begin{align}\label{eq:phi wieder n-1 deri}
		\begin{split}
		&(R\partial_{R}-2\xi\partial_{\xi})\phi_{-1}(R,\xi)\\
		=&-R^{-\frac12}\Re\left(a_{-1}(\xi)i\xi^{-\frac34}e^{iR\xi^{\frac12}}\left(1+O\left(\frac{1}{R\xi^{\frac12}}\right)+O\left(\frac{1}{R^{2}}\right)\right)\right)\\&+\Re\left(R^{-\frac12}\left(3i\xi^{-\frac34}a_{-1}(\xi)-4i\xi^{\frac14}a'_{-1}(\xi)\right)e^{iR\xi^{\frac12}}\left(1+O\left(\frac{1}{R\xi^{\frac12}}\right)+O\left(\frac{1}{R^{2}}\right)\right)\right)\\
		=&\Re\left(R^{-\frac12}\left(2i\xi^{-\frac34}a_{-1}(\xi)-4i\xi^{\frac14}a'_{-1}(\xi)\right)e^{iR\xi^{\frac12}}\left(1+O\left(\frac{1}{R\xi^{\frac12}}\right)+O\left(\frac{1}{R^{2}}\right)\right)\right)
		\end{split}
	\end{align}
	We again use the fact $\Re z\Re w=\frac12(\Re(zw)+\Re(z\overline{w})$ to find the $\delta$ measure on the diagonal in the integral
	\begin{align*}
		\lim_{A\rightarrow\infty}\int_{0}^{A}[R\partial_{R}-2\xi\partial_{\xi}]\phi_{-1}(R,\xi)\phi_{-1}(R,\eta)RdR.
	\end{align*}
	As in the proof for Proposition \ref{prop: K operator} and Proposition \ref{prop: K operator n-1}, the profiles $O\left(\frac{1}{R\xi^{\frac12}}\right)$ and $O\left(\frac{1}{R^{2}}\right)$ in the parenthesis $\Re(...)$ in \eqref{eq:phi wieder n-1} and \eqref{eq:phi wieder n-1 deri} do not contribute to the $\delta$ measure. Therefore the contribution in consideration is given by
	\begin{align}\label{eq:Ixieta n-1}
		\begin{split}
		\lim_{L\rightarrow\infty}\int_{0}^{\infty}\Re\left(\left(\left(2\xi^{-\frac34}a_{-1}(\xi)-4\xi^{\frac14}a'_{-1}(\xi)\right)\eta^{-\frac34}\overline{  a_{-1}(\eta)}\right)e^{iR(\xi^{\frac12}-\eta^{\frac12})}\right)\chi_{1}(R)\chi_{2}(R/L)dR:=I(\xi,\eta;-1)
		\end{split}
	\end{align}
	where $\chi_{1},\chi_{2}$ is again the same as in \eqref{eq:Ixieta}. Now a routine calculation gives
	\begin{align*}
		I(\xi,\eta;-1)=&\pi\Re\left(\left(2\xi^{-\frac34}a_{-1}(\xi)-4\xi^{\frac14}a'_{-1}(\xi)\right)\eta^{-\frac34}\overline{  a_{-1}(\eta)}\right)\delta(\xi^{\frac12}-\eta^{\frac12})\\=&2\xi^{\frac12}\pi\Re\left(\left(2\xi^{-\frac34}a_{-1}(\xi)-4\xi^{\frac14}a'_{-1}(\xi)\right)\eta^{-\frac34}\overline{  a_{-1}(\eta)}\right)\delta(\xi-\eta),
	\end{align*}
	which in turn gives the diagonal contribution
	\begin{align*}
&		\left\langle\int_{0}^{\infty}f(\xi)[R\partial_{R}-2\xi\partial_{\xi}]\phi_{-1}(R,\xi)\trho_{-1}(\xi)d\xi,\phi_{-1}(R,\eta)\right\rangle_{L^{2}_{RdR}}\\
=&\left\langle2\xi^{\frac12}\pi\Re\left(\left(2\xi^{-\frac34}a_{-1}(\xi)-4\xi^{\frac14}a'_{-1}(\xi)\right)\eta^{-\frac34}\overline{  a_{-1}(\eta)}\right)\delta(\xi-\eta),f(\xi)\trho_{-1}(\xi)\right\rangle_{d\xi}\\
=&2\pi\eta^{\frac12}\Re\left(2\eta^{-\frac32}|a_{-1}(\eta)|^{2}-4\eta^{-\frac12}a'_{-1}(\eta)\overline{a_{-1}(\eta)}\right)f(\eta)\trho_{-1}(\eta)\\
=&\left(4\pi\eta^{-1}|a_{-1}(\eta)|^{2}-4\pi\left(|a_{-1}(\eta)|^{2}\right)^{\prime}\right)f(\eta)\trho_{-1}(\eta)\\
=&\left(\frac{1}{\eta\rho_{-1}(\eta)}+\frac{\rho'_{-1}(\eta)}{(\rho_{-1}(\eta))^{2}}\right)f(\eta)\trho_{-1}(\eta)=\frac{\eta\trho'_{-1}(\eta)}{\trho_{-1}(\eta)}f(\eta).
	\end{align*}
	This gives the desired result on the diagonal part of $\calK_{-1}$. Here we used the relation $\rho_{-1}(\xi)=\frac{1}{4\pi}|a_{-1}(\xi)|^{-2}$ given in Proposition \ref{prop: spectral measure n-1}.

For the off-diagonal part $\calK^{0}_{-1}$, a routing calculation as in the proof of Proposition \ref{prop: K operator} gives the expression for $F_{-1}(\xi,\eta)$:
\begin{align*}
	F_{-1}(\xi,\eta)=\int_{0}^{\infty}W_{-1}(R)\phi_{-1}(R,\xi),\phi_{-1}(R,\eta)R\,dR,\quad
	W_{-1}(R)=\left[\tilde H^{-}_{1},R\partial_{R}\right]-2\tilde H^{-}_{1}=\frac{16}{(1+R^{2})^{2}}.
\end{align*}
By Proposition \ref{prop: basis n-1 small} and $\calD_{-}\left(R^{-\frac12}\Phi^{-1}_{0}(R)\right)=0$, we have, for $R^{2}\xi\lesssim 1$
\begin{align}\label{phi -1 small R2xi}
	|\phi_{-1}(R,\xi)|\lesssim R^{3}\langle R^{2}\rangle^{-1}\lesssim \xi^{-\frac12}\langle\xi\rangle^{-1},\quad |\partial_{R}\phi_{-1}(R,\xi)|\lesssim R^{2}\langle R^{2}\rangle^{-1}\lesssim \langle\xi\rangle^{-1}.
\end{align}
For $R^{2}\xi\gg 1$, we have, by Propositions \ref{prop: basis n-1 large}, \ref{prop: spectral measure n-1} and \eqref{basis n-1 diff}, 
\begin{align}\label{phi-1 large R2xi}
\begin{split}
	\phi_{-1}(R,\xi)=\,&\xi^{-1}\calD_{-}\left(2\Re\left(a_{-1}(\xi)R^{-\frac12}\xi^{-\frac14}e^{iR\xi^{\frac12}}\sigma_{-1}(R\xi^{\frac12},R)\right)\right)\\\simeq \,&\xi^{-\frac54}\langle\xi\rangle^{-1}\calD_{-}\left(\Re\left(R^{-\frac12}e^{iR\xi^{\frac12}}\left(1+O\left(R^{-1}\xi^{-\frac12}\right)\right)\right)\right)\\
	\simeq \,&\xi^{-\frac34}\langle\xi\rangle^{-1}\Re\left(R^{-\frac12}e^{iR\xi^{\frac12}}\left(1+O\left(R^{-1}\xi^{-\frac12}\right)\right)\right)
	\end{split}
\end{align}
Without loss of generality, we again assume $\xi\leq \eta$ and set $\eta^{\frac12}=\zeta\xi^{\frac12}$. We discuss the following three regimes: $R\eta^{\frac12}\geq R\xi^{\frac12}\gg1$, $R\eta^{\frac12}\gg1\gtrsim R\xi^{\frac12}$ and $1\gtrsim R\eta^{\frac12}\geq R\xi^{\frac12}$. We start with first case where both $\phi_{-1}(R,\xi)$ and $\phi_{-1}(R,\eta)$ are in the oscillatory regime. In this case we consider the integral
\begin{align}\label{F-1 a pre}
	\begin{split}
	F_{-1}^{a}(\xi,\eta):=&\int_{0}^{\infty}\chi_{R\eta^{\frac12}\gg1}\cdot \chi_{R\xi^{\frac12}\gg1}\cdot W_{-1}(R)\phi_{-1}(R,\xi)\phi_{-1}(R,\eta)R\,dR\\
	=&\int_{0}^{\infty}\chi_{R\eta^{\frac12}\gg1} \chi_{R\xi^{\frac12}\gg1} W_{-1}(R)\xi^{-\frac34}\langle\xi\rangle^{-1}\eta^{-\frac34}\langle\eta\rangle^{-1}\\
	&\cdot\Re\left(e^{iR\xi^{\frac12}}\left(1+O\left(R^{-1}\xi^{-\frac12}\right)\right)\right)\Re\left(e^{iR\eta^{\frac12}}\left(1+O\left(R^{-1}\eta^{-\frac12}\right)\right)\right) \,dR.
	\end{split}
\end{align}
If we simply bound the oscillatory factors $e^{iR\xi^{\frac12}}$ and $e^{iR\eta^{\frac12}}$ in absolute value, then we have (taking into account the integral of $W_{-}(R)$)
\begin{align}\label{F-1 a 1}
	\begin{split}
	\left|F_{-1}^{a}(\xi,\eta)\right|\lesssim &\langle\xi\rangle^{-1}\langle\eta\rangle^{-1}\eta^{-\frac34}\cdot\min\{\xi^{\frac34},\xi^{-\frac34}\}.
	\end{split}
\end{align}
To gain any factor of $\frac{\xi^{\frac12}}{\eta^{\frac12}}$, we need to use integration by parts. 
We consider an integral of the following form. The other cases can be handled similarly:
\begin{align}\label{F-1 a 2 pre}
	\begin{split}
	\int_{0}^{\infty}&\chi_{R^{2}\eta\gg1}\chi_{R^{2}\xi\gg1}W_{-1}(R)\xi^{-\frac34}\langle\xi\rangle^{-1}\eta^{-\frac34}\langle\eta\rangle^{-1}\cdot e^{iR\xi^{\frac12}}\left(1+O\left(R^{-1}\xi^{-\frac12}\right)\right)\\&\cdot \frac{\partial_{R}(e^{iR\eta^{\frac12}})}{i\eta^{\frac12}}\left(1+O\left(R^{-1}\eta^{-\frac12}\right)\right)\,dR.
	\end{split}
\end{align}
Upon integration by parts, no matter $\partial_{R}$ hits $W_{-}(R)$ or $O(R^{-1}\xi^{-\frac12})$, we gain a factor of $R^{-1}\lesssim \xi^{\frac12}$. Therefore $F_{-1}^{a}(\xi,\eta)$ can be also bounded by
\begin{align}\label{F-1 a 2}
	|F_{-1}^{a}(\xi,\eta)|\lesssim \langle\xi\rangle^{-1}\langle\eta\rangle^{-1}\eta^{-\frac34}\cdot\min\{\xi^{\frac34},\xi^{-\frac34}\}\left(\frac{\xi^{\frac12}}{\eta^{\frac12}}\right)^{N},\quad \textrm{for any}\quad N>0.
\end{align}
Next we turn to the regime $R\eta^{\frac12}\gg1\gtrsim R\xi^{\frac12}$, in which we consider the integral
\begin{align}\label{F-1 b pre}
\begin{split}
	F_{-1}^{b}(\xi,\eta):=&\int_{0}^{\infty}\chi_{R^{2}\eta\gg1}\chi_{R^{2}\xi\lesssim1}W_{-1}(R)\phi_{-1}(R,\xi)\phi_{-1}(R,\eta)R\,dR\\=&\int_{0}^{\infty}\chi_{R^{2}\eta\gg1}\chi_{R^{2}\xi\lesssim1}W_{-1}(R)\eta^{-\frac34}\langle\eta\rangle^{-1}R^{\frac12}\\
	&\cdot \phi_{-1}(R,\xi)\Re\left(e^{iR\eta^{\frac12}}\left(1+O\left(R^{-1}\eta^{-\frac12}\right)\right)\right)\,dR.
	\end{split}
\end{align}
If $\xi\leq 1$, then we use the bound $|\phi_{-1}(R,\xi)|\lesssim \langle R\rangle$ in \eqref{phi -1 small R2xi} to bound this integral as
\begin{align}\label{F-1b 1}
	\begin{split}
	|F_{-1}^{b}(\xi,\eta)|\lesssim &\langle\eta\rangle^{-1}\min\{1,\eta^{-\frac34}\}.
	\end{split}
\end{align}
If $\xi\geq1$, then we use the bound $|\phi_{-1}(R,\xi)|\lesssim \xi^{-\frac32}$ in \eqref{phi -1 small R2xi} to bound the integral as
\begin{align}\label{F-1b 1 prime}
	|F_{-1}^{b}(\xi,\eta)|\lesssim \eta^{-\frac74}\xi^{-\frac74},
\end{align}
where we have bounded the factor $R^{\frac12}$ in \eqref{F-1 b pre} by $\xi^{-\frac14}$.
Finally we consider the regime $0<R\xi^{\frac12}\leq R\eta^{\frac12}\lesssim1$. The integral in consideration now is given by
\begin{align}\label{F-1 c pre}
	F^{c}_{-1}(\xi,\eta):=\int_{0}^{\infty}\chi_{R^{2}\eta\lesssim 1}\chi_{R^{2}\xi\lesssim 1}W_{-1}(R)\phi_{-1}(R,\xi)\phi_{-1}(R,\eta)R\,dR.
\end{align}
In view of \eqref{phi -1 small R2xi}, this can be bounded by
\begin{align}\label{F-1 c 1}
\begin{split}
	\left|F_{-1}^{c}(\xi,\eta)\right|\lesssim \left|\int_{0}^{\infty}\chi_{R^{2}\eta\lesssim 1}\chi_{R^{2}\xi\lesssim 1}W_{-1}(R)R^{3}\langle R^{2}\rangle^{-1}\eta^{-\frac12}\langle\eta\rangle^{-1}R\,dR\right|
	\lesssim\langle\eta\rangle^{-1}\cdot\frac{\min\{\xi^{\frac12},\xi^{-2}\}}{\eta^{\frac12}}.
	\end{split}
\end{align}
To obtain the rapid decay factor $\left(\frac{\langle\xi\rangle}{\langle\eta\rangle}\right)^{N}$ for any $N>0$, we consider
\begin{align*}
	\eta\int_{0}^{\infty}\chi_{R\xi^{\frac12}\ll1}W_{-1}(R)\phi_{-1}(R,\xi)\phi_{-1}(R,\eta)R\,dR=&-\int_{0}^{\infty}\chi_{R\xi^{\frac12}\ll1}W_{-1}(R)\phi_{-1}(R,\xi)H_{-1}^{+}\left(\phi_{-1}(R,\eta)\right)R\,dR,
\end{align*}
and use an argument similar to the proof of Proposition \ref{prop: K operator}.\\

Next we turn to the derivatives of $F_{-1}(\xi,\eta)$. By Proposition \ref{prop: basis n-1 small} and \eqref{basis n-1 diff}, if $R\xi^{\frac12}\lesssim 1$, then upon a differentiation in $\xi^{\frac12}$ for $\phi_{-1}(R,\xi)$, we gain an extra factor of $\xi^{-\frac12}$. This of course also holds for $\phi_{-1}(R,\eta)$ when $R\eta^{\frac12}\lesssim 1$. Therefore the desired estimates hold for $R\xi^{\frac12}\lesssim1$ and $R\eta^{\frac12}\lesssim 1$. If $R\eta^{\frac12}\gg1$, we use Propositions \ref{prop: basis n-1 large} and \ref{prop: spectral measure n-1}, as well as \eqref{basis n-1 diff}, to obtain that compared to the undifferentiated basis, we gain a factor of $R$ upon a differentiation in $\eta^{\frac12}$. In the estimates for the derivatives of $F_{-1}(\xi,\eta)$, this extra factor of $R$ is absorbed into the function $W_{-1}(R)$. Therefore if $\eta\ll1$, we simply obtain a factor of $\eta^{-\frac12}$. If $\eta\gtrsim1$ and $R\xi^{\frac12}\lesssim 1$, we use \eqref{deri transfer FB}and integration by parts in $R$ to obtain the desired estimates. The case when $1\ll R\xi^{\frac12}\leq R\eta^{\frac12}$ is more delicate. First, for $\partial_{\eta^{\frac12}}F_{-1}(\xi,\eta)$ and $\partial_{\eta^{\frac12}}^{2}F_{-1}(\xi,\eta)$, we again use \eqref{deri transfer FB} and the above integration by parts argument to obtain:
\begin{align}\label{deri F-1 eta hard 1}
	\begin{split}
\left|\partial^{k}_{\eta^{\frac12}}F_{-1}(\xi,\eta)\right|\lesssim \max\{1,\xi^{-k}\}\cdot\Gamma_{-1}\cdot \left(\frac{\xi^{\frac12}}{\eta^{\frac12}}\right)^{k},\quad \textrm{for}\quad k=1,2,\quad\textrm{and}\quad 1\ll R\xi^{\frac12}\leq R\eta^{\frac12}.
	\end{split}
\end{align}
Here $\Gamma_{-1}$ is the bound for the undifferentiated $F_{-1}(\xi,\eta)$ in the same regime. For the derivatives in $\xi^{\frac12}$, we again distinguish the cases when $\xi\ll1$ and $\xi\gtrsim1$.  If $\xi\ll1$, we simply obtain a factor of $\xi^{-\frac12}$ upon a differentiation in $\xi^{\frac12}$. If $\xi\gtrsim1$, then we obtain the same bound as the undifferentiated one. Therefore we have
\begin{align}\label{deri F-1 eta hard 2}
	\begin{split}
&	\left|\partial_{\xi^{\frac12}}F_{-1}(\xi,\eta)\right|\lesssim \Gamma_{-1},\quad \left|\partial^{2}_{\xi^{\frac12}}F_{-1}(\xi,\eta)\right|\lesssim \Gamma_{-1},\\
&\left|\partial_{\xi^{\frac12}}\partial_{\eta^{\frac12}}F_{-1}(\xi,\eta)\right|\lesssim\max\{1,\xi^{-1}\}\cdot \Gamma_{-1}\cdot\frac{\xi^{\frac12}}{\eta^{\frac12}}\quad \textrm{for}\quad 1\ll R\xi^{\frac12}\leq R\eta^{\frac12}.
	\end{split}
\end{align}
This completes the proof.
\end{proof}

\subsection{Transference operator at angular momentum $n=0$}\label{subsec:linear end}

The operator $\calK$ is defined as
\begin{align}\label{def trans n0}
\widehat{R\partial_{R}u}=-2\xi\partial_{\xi}\widehat{u}+\calK\widehat{u}.
\end{align}
In this section, the Fourier transform is taken with respect to the Fourier basis (See (4.1) in \cite{KMiao})
\begin{align}\label{FB n0}
\begin{split}
\phi_{0}(R,\xi):=\xi^{-1}\calD_{0}\left(R^{-\frac12}\phi_{KST}(R,\xi)\right)
\end{split}
\end{align}
Here $\phi_{KST}(R,\xi)$ (See \cite{KST}) is the Fourier basis associated to the operator
\begin{align*}
\calL_{KST}:=-\partial_{R}^{2}+\frac{3}{4R^{2}}-\frac{8}{(1+R^{2})^{2}}
\end{align*}
with respect to the measure $dR$. The operator $\calD_{0}$ (See \cite{KMiao}) is given by
\begin{align*}
\calD_{0}:=\partial_{R}+\frac{1}{R}-\frac{2}{R(R^{2}+1)}.
\end{align*}
As in \cite{KMiao}, we have
\begin{proposition}\label{prop: K n0}
	The operator $\calK$ is given by, with $\trho_{0}(\xi)$ being the spectral measure density associated to the Fourier basis $\phi_{0}(R,\xi)$,
	\begin{align*}
	\calK f(\xi)=-2f(\xi)-\frac{\trho^{\prime}_{0}(\xi)\xi}{\trho_{0}(\xi)}\cdot f(\xi)+\calK_{0}f(\xi)
	\end{align*}
	where the off-diagonal operator $\calK_{0}$ is an integral operator with kernel
	\begin{align*}
	\frac{\trho_{0}(\eta)F_{0}(\xi,\eta)}{\xi-\eta}.
	\end{align*}
	The symmetric $C^{2}$ function $F_{0}(\xi,\eta)$ satisfies the bounds
	\begin{align}\label{F0 bound}
	\begin{split}
	\left|F_{0}(\xi,\eta)\right|\lesssim& \min\{\xi^{\frac32},1\}\xi^{-\frac14}\log(1+\xi^{-1})\eta^{-\frac14}\log(1+\eta^{-1})\cdot\left(\frac{\xi^{\frac12}}{\eta^{\frac12}}\right)^{N}\\
	&+\frac{\eta\left(\log(1+\eta^{-1})\right)^{2}}{\langle\eta\rangle^{N+\frac14}}+\langle\eta\rangle^{-3}\left(\frac{\langle\xi\rangle}{\langle\eta\rangle}\right)^{N}:=\Gamma_{0}.
	\end{split}
	\end{align}
	For the derivatives of $F_{0}(\xi,\eta)$, we have
	\begin{align}\label{F0 deri bound}
		\begin{split}
		&	\left|\partial_{\xi^{\frac12}}F_{0}(\xi,\eta)\right|\lesssim \left(1+\xi^{-\frac12}\right)\cdot \Gamma_{0},\quad \left|\partial_{\eta^{\frac12}}F_{0}(\xi,\eta)\right|\lesssim \frac{\Gamma_{0}}{\eta^{\frac12}}\left(1+\xi^{\frac12}\right),\\
		&\left|\partial^{2}_{\xi^{\frac12}}F_{0}(\xi,\eta)\right|\lesssim \left(1+\xi^{-1}\right)\cdot\Gamma_{0},\quad \left|\partial^{2}_{\eta^{\frac12}}F_{0}(\xi,\eta)\right|\lesssim \frac{\Gamma_{0}}{\eta}\cdot\left(1+\xi\right),\\
		&\left|\partial_{\xi^{\frac12}}\partial_{\eta^{\frac12}}F_{0}(\xi,\eta)\right|\lesssim \frac{\Gamma_{0}}{\eta^{\frac12}}\left(\xi^{\frac12}+\xi^{-\frac12}\right).
		\end{split}
	\end{align}
\end{proposition}
\begin{proof}
	The diagonal part of the operator $\calK$ was calculated in \cite{KMiao}. Here we strive to derive a more precise bound on the off-diagonal kernel $F_{0}(\xi,\eta)$. A routing calculation gives (See (4.17 ) in \cite{KMiao})
	\begin{align}\label{def F0}
	F_{0}(\xi,\eta):=\int_{0}^{\infty}W_{0}(R)\phi_{0}(R,\xi)\phi_{0}(R,\eta)R\,dR,\quad W_{0}(R):=\frac{16}{(1+R^{2})^{2}}-\frac{32}{(1+R^{2})^{3}}.
	\end{align}
	Again in view of the bounds derived in \cite{KMiao}, we have
	\begin{align}\label{phi0 bounds 1}
	\left|\phi_{0}(R,\xi)\right|\lesssim \log(1+R^{2})\quad \textrm{for}\quad R^{2}\xi\lesssim1.
	\end{align}
	and
	\begin{align}\label{phi0 bounds 2}
	\begin{split}
	\phi_{0}(R,\xi)=&R^{-\frac12}\xi^{-\frac14}\log\left(1+\xi^{-1}\right)\Re\left(e^{iR\xi^{\frac12}}\left(1+O\left(R^{-1}\xi^{-\frac12}\right)\right)\right)\quad \textrm{for}\quad R\xi^{\frac12}\gg1.
	\end{split}
	\end{align}
	Again we split into the following three regimes: $R\eta^{\frac12}\geq R\xi^{\frac12}\gg 1, R\eta^{\frac12}\gg 1\gtrsim R\xi^{\frac12}$ and $1\gtrsim R\eta^{\frac12}\geq R\xi^{\frac12}$. We start with the case when both $\phi_{0}(R,\xi)$ and $\phi_{0}(R,\eta)$ are in the oscillatory regime:
	\begin{align}\label{F0 a pre}
	\begin{split}
	F_{0}^{a}(\xi,\eta):=&\int_{0}^{\infty}\chi_{R\eta^{\frac12}\gg1}\chi_{R\xi^{\frac12}\gg1} W_{0}(R)\phi_{0}(R,\xi)\phi_{0}(R,\eta)R\,dR\\
	=&\int_{0}^{\infty}\chi_{R\eta^{\frac12}\gg1}\chi_{R\xi^{\frac12}\gg1} W_{0}(R)\xi^{-\frac14}\log\left(1+\xi^{-1}\right)\eta^{-\frac14}\log\left(1+\eta^{-1}\right)\\
	&\cdot\Re\left(e^{iR\xi^{\frac12}}\left(1+O\left(R^{-1}\xi^{-\frac12}\right)\right)\right)\Re\left(e^{iR\eta^{\frac12}}\left(1+O\left(R^{-1}\eta^{-\frac12}\right)\right)\right)\,dR.
	\end{split}
	\end{align}
	If we simply bound the oscillatory factors in absolute value, then we have
	\begin{align*}
	\left|F_{0}^{a}(\xi,\eta)\right|\lesssim\min\{\xi^{\frac32},1\}\cdot\xi^{-\frac14}\log\left(1+\xi^{-1}\right)\eta^{-\frac14}\log\left(1+\eta^{-1}\right)
	\end{align*}
	Writing $e^{iR\eta^{\frac12}}=\frac{1}{i\eta^{\frac12}}\partial_{R}\left(e^{iR\eta^{\frac12}}\right)$, we can do integration by parts. No matter $\partial_{R}$ hits $W_{0}(R)$ or $O(R^{-1}\xi^{-\frac12})$, we gain a factor of $R^{-1}\lesssim \xi^{\frac12}$. Therefore $F^{a}_{0}(\xi,\eta)$ can be bounded by
	\begin{align}\label{F0 a}
	\begin{split}
	\left|F_{0}^{a}(\xi,\eta)\right|\lesssim \min\{\xi^{\frac32},1\}\cdot\xi^{-\frac14}\log\left(1+\xi^{-1}\right)\eta^{-\frac14}\log\left(1+\eta^{-1}\right)\cdot\left(\frac{\xi^{\frac12}}{\eta^{\frac12}}\right)^{N},\quad \textrm{for any}\quad N>0.
	\end{split}
	\end{align}
	Next we turn to the regime $R\eta^{\frac12}\gg1\gtrsim R\xi^{\frac12}$:
	\begin{align}\label{F0 b pre}
	\begin{split}
	F_{0}^{b}(\xi,\eta):=&\int_{0}^{\infty}\chi_{R\eta^{\frac12}\gg1}\chi_{R\xi^{\frac12}\lesssim 1}W_{0}(R)\phi_{0}(R,\xi)\phi_{0}(R,\eta)R\,dR\\=&\int_{0}^{\infty}\chi_{R\eta^{\frac12}\gg1}\chi_{R\xi^{\frac12}\lesssim 1}W_{0}(R)R^{\frac12}\log(1+R^{2})\eta^{-\frac14}\log(1+\eta^{-1})\Re\left(e^{iR\eta^{\frac12}}\left(1+O\left(R^{-1}\eta^{-\frac12}\right)\right)\right)\,dR,
	\end{split}
	\end{align}
	which is bounded by, upon integration by parts
	\begin{align}\label{F0 b}
	\left|F_{0}^{b}(\xi,\eta)\right|\lesssim \frac{\eta\left(\log\left(1+\eta^{-1}\right)\right)^{2}}{\langle\eta\rangle^{\frac14}}\langle\eta\rangle^{-N},\quad \textrm{for any}\quad N>0.
	\end{align}
	Finally we consider the regime $0<R\xi^{\frac12}\leq R\eta^{\frac12}\lesssim 1$, and we have
	\begin{align}\label{F0 c}
	\begin{split}
	\left|F^{c}_{0}(\xi,\eta)\right|\lesssim& \left|\int_{0}^{\infty}\chi_{R^{2}\eta\lesssim 1}\chi_{R^{2}\xi\lesssim 1}W_{0}(R)\phi_{0}(R,\xi)\phi_{0}(R,\eta)R\,dR\right|
	\lesssim \langle\eta\rangle^{-3}.
	\end{split}
	\end{align}
	To gain the rapid decay $\left(\frac{\langle\xi\rangle}{\langle\eta\rangle}\right)^{N}$ for any $N>0$, we consider
	\begin{align*}
	\eta\int_{0}^{\infty}\chi_{R\xi^{\frac12}\ll1}W_{0}(R)\phi_{0}(R,\xi)\phi_{0}(R,\eta)R\,dR=-\int_{0}^{\infty}\chi_{R\xi^{\frac12}\ll1}W_{0}(R)\phi_{0}(R,\xi)H_{0}\left(\phi_{0}(R,\eta)\right)R\,dR,
	\end{align*}
	and then a routing argument gives the desired estimate.\\
	
	The estimates on the derivatives of $F_{0}(\xi,\eta)$ are derived similarly as those for $F_{-1}(\xi,\eta)$.
\end{proof}

\section{Estimates for smooth sources for $n\geq2$}\label{sec:bilinlargensmooth}
\subsection{Precise equations-Revisit}
In this section, we derive the precise equations for each angular mode and their counterpart on Fourier side with respect to the $R$-variable. We start by taking the Fourier transform for the equation \eqref{coe eq tau R} with respect to the angular variable. To this end, we (formally) expand $\frakN(\vphi_{\ell}), \ell=1,2$ (recall the notation in \eqref{nonlinearity tau R}) into Fourier series
\begin{align}\label{Fourier expansion NL}
	\begin{split}
	&\frakN(\vphi_{\ell})(t,R,\theta)=\sum_{n}\hfrakN(\vphi_{\ell})(n,t,R)e^{in\theta},\quad \hfrakN(\vphi_{\ell})(n,t,R):=\int_{0}^{2\pi}\frakN(\vphi_{\ell})(t,R,\theta)e^{-in\theta}\frac{d\theta}{2\pi},\quad n\in\bbZ,\quad \ell=1,2
	\end{split}
\end{align}
Then we have the Fourier counterpart of the equation \eqref{coe eq tau R}, with $\veps_{\pm}(\tau,R,\theta)=\varphi_{1}(\tau,R,\theta)\mp i\varphi_{2}(\tau,R,\theta)$:
\begin{equation}\label{eq:RegularFinestructure1}\begin{split}
&-\left(\left(\partial_{\tau} + \frac{\lambda_{\tau}}{\lambda}R\partial_R\right)^2 + \frac{\lambda_{\tau}}{\lambda}\left(\partial_{\tau} + \frac{\lambda_{\tau}}{\lambda}R\partial_R\right)\right)\varepsilon_{\pm}(n)+ H_n^{\pm}\varepsilon_{\pm}(n) = F_{\pm}(n),
\end{split}\end{equation}
where the source term admits a decomposition into the following pieces: 
\begin{equation}\label{eq:RegularFinestructure2}\begin{split}
F_{\pm}(n) &= \lambda^{-2}N_{\pm}(n)  - 2\frac{\sin\left[2Q + \eps\right]\sin\eps}{R^2}\varepsilon_{\pm}(n) - 4\frac{\sin\left[Q+\frac{\eps}{2}\right]\sin\left[\frac{\eps}{2}\right]}{R^2}i\varepsilon_{\pm,\theta}(n)\\
&\pm i\left(\frac{2\partial_R\eps}{1+R^2} + \left(\partial_R\epsilon\right)^2 - \frac{\lambda_{\tau}}{\lambda}\frac{8R}{1+R^2}\left(\partial_{\tau}\eps + \frac{\lambda_{\tau}}{\lambda}R\partial_R\eps\right) - \left(\partial_{\tau}\eps + \frac{\lambda_{\tau}}{\lambda}R\partial_R\eps\right)^2\right)\varphi_2(n)\\
& \mp \frac{(1+\nu)^2}{\nu^2\tau^2}\cdot\frac{4R^2}{1+2R^2 + R^4}i\varphi_2(n),
\end{split}\end{equation}
and where the first component of $F_{\pm}(n)$ is the angular momentum $n$ projection (which we shall denote by $\Pi_n$) of $\lambda^{-2}N(\varphi_1) \mp i\lambda^{-2}N(\varphi_2)$: 
\begin{equation}\label{eq:RegularFinestructure3}\begin{split}
\lambda^{-2}N_{\pm}(n) &= \Pi_n\left(\mathcal{P}\varepsilon_{\pm}\right)\\
& + \Pi_n\left(\frac{2}{\sqrt{1-\left|\Pi_{\Phi^{\perp}}\varphi\right|^2}}\left[U_R\sum_{j=1}^2\varphi_j\varphi_{j,R} - \left(\partial_{\tau}+\frac{\lambda_{\tau}}{\lambda}R\partial_R\right)U\sum_{j=1}^2\varphi_j\left(\partial_{\tau}+\frac{\lambda_{\tau}}{\lambda}R\partial_R\right)\varphi_j\right]\right)\\
& \mp  \Pi_n\left(\frac{2\sin U}{R^2\sqrt{1-\left|\Pi_{\Phi^{\perp}}\varphi\right|^2}}\sum_{j=1}^2\varphi_j\varphi_{j,\theta}\right)\\
&\mathcal{P} = \big|\Psi_t\big|^2 - \big|\Phi_t\big|^2 -\big|\nabla\Psi\big|^2 + \big|\nabla\Phi\big|^2, 
\end{split}\end{equation}
and where the fine structure of the last term $\mathcal{P} $ is given in \eqref{non linear 1 2}. 

Now we derive the equations in $(\tau,\xi)$-variable for $|n|\geq 2$. We first introduce the following notations:
\begin{align}\label{notations xi variab}
\begin{split}
	\xh(\tau,\xi):=&\Fh(\veps_{+}(\tau,R,n)):=\int_{0}^{\infty}\phi(R,\xi;\hbar)\veps_{+}(\tau,R,n)R\,dR,\\ \xhb(\tau,\xi):=&\Fhb(\veps_{-}(\tau,R,-n)):=\int_{0}^{\infty}\phi(R,\xi;\thbar)\veps_{-}(\tau,R,-n)R\,dR.
	\end{split}
\end{align}
Here $\phi(R,\xi;\hbar):=R^{-\frac12}\phi_{n}(R,\xi)$ where $\phi_{n}(R,\xi)$ are introduced in Proposition \ref{prop:DFT nlarge} for $n\gg1$, and $\phi(R,\xi;\hbar):R^{-\frac12}w_{n}(\xi)\phi(R,\xi)$ where $w_{n}(\xi)\phi(R,\xi)$ are introduced in Proposition \ref{prop: FB match} for $2\leq n\leq N_{0}$. Similarly, $\phi(R,\xi;\thbar):=R^{-\frac12}\phi_{-n}(R,\xi)$ where $\phi_{-n}(R,\xi)$ are introduced in Proposition \ref{prop:DFT nlarge negative} for $n\gg1$, and $\phi(R,\xi;\thbar):=R^{-\frac12}\phi_{-n}(R,\xi)$ where $\phi_{-n}(R,\xi)$ are introduced in Proposition \ref{prop: FB nega n large hbar} for $2\leq n\leq N_{0}$.

Applying the Fourier transform $\Fh$ to both sides of the ``$+$"-component of \eqref{eq diag phys}, we obtain
\begin{align}\label{eq tau xi hbar pre}
	\begin{split}
	-\left(\partial_{\tau}-2\dftau\xi\partial_{\xi}+\dftau\calK_{\hbar}\right)^{2}\xh-\dftau\left(\partial_{\tau}-2\dftau\xi\partial_{\xi}+\dftau\calK_{\hbar}\right)\xh-\xi\xh=\Fh\left(F_{+}\right).
	\end{split}
\end{align}
By the definition of $\calK^{(0)}_{\hbar}$, the above equation becomes
\begin{align}\label{eq tau xi hbar D}
	\begin{split}
	&-\left(\calD_{\tau}^{2}+\dftau\calD_{\tau}+\xi\right)\xh\\
	=&\Fh\left(\frakN_{+}(\veps)(n)\right)\\
	&+2\dftau\calK_{\hbar}^{(0)}\calD_{\tau}\xh+\left(\dftau\right)^{\prime}\calK_{\hbar}^{(0)}\xh+\dftau\left[\calD_{\tau},\calK_{\hbar}^{(0)}\right]\xh+\left(\dftau\right)^{2}\left(\left(\calK_{\hbar}^{(0)}\right)^{2}+\calK_{\hbar}^{(0)}\right)\xh\\
	=:&\Fh\left(\frakN_{+}(\veps)(n)\right)+\calR_{\hbar}(\xh,\calD_{\tau}\xh).
	\end{split}
\end{align}
Here the operator $\calD_{\tau}$ is defined as
\begin{align}\label{def calD}
	\calD_{\tau}:=\partial_{\tau}-2\dftau\xi\partial_{\xi}-\dftau\frac{\rho'_{n}(\xi)\xi}{\rho_{n}(\xi)}-2\dftau.
\end{align}
For the second order linear operator appearing in the equation \eqref{eq tau xi hbar D}, we have, by direct calculation (see \cite{KS_FullRange}),
\begin{lemma}\label{lem: parametrix hbar}
	Let $\calD_{\tau}$ be as in \eqref{def calD}. If $\xb(\tau,\xi)$ satisfies
	\begin{align}\label{parametrix hbar homo eq}
		\left(\calD_{\tau}^{2}+\dftau\calD_{\tau}+\xi\right)\xb(\tau,\xi)=0;\quad \xb(\tau_{0},\xi)=\xb_{0}(\xi),\quad \calD_{\tau}\xb(\tau_{0},\xi)=\xb_{1}(\xi),
	\end{align}
	then $\xb(\tau,\xi)$ is given by
	\begin{align}\label{parametrix hbar homo}
		\begin{split}
		\xb(\tau,\xi)=&\frac{\lambda(\tau)^{2}}{\lambda(\tau_{0})^{2}}\frac{\rho^{\frac{1}{2}}_{n}\left(\frac{\lambda(\tau)^{2}}{\lambda(\tau_{0})^{2}}\xi\right)}{\rho^{\frac{1}{2}}_{n}(\xi)}\cos\left(\lambda(\tau)\xi^{\frac{1}{2}}\int_{\tau_{0}}^{\tau}\lambda(u)^{-1}du\right)\xb_{0}\left(\frac{\lambda(\tau)^{2}}{\lambda(\tau_{0})^{2}}\xi\right)\\
		&+\xi^{-\frac{1}{2}}\frac{\lambda(\tau)}{\lambda(\tau_{0})}\frac{\rho^{\frac{1}{2}}_{n}\left(\frac{\lambda(\tau)^{2}}{\lambda(\tau_{0})^{2}}\xi\right)}{\rho^{\frac{1}{2}}_{n}(\xi)}\sin\left(\lambda(\tau)\xi^{\frac{1}{2}}\int_{\tau_{0}}^{\tau}\lambda(u)^{-1}du\right)\xb_{1}\left(\frac{\lambda(\tau)^{2}}{\lambda(\tau_{0})^{2}}\xi\right).
		\end{split}
	\end{align}
	If $\xb(\tau,\xi)$ satisfies
	\begin{align}\label{parametrix hbar inhomo eq}
	\left(\calD_{\tau}^{2}+\dftau\calD_{\tau}+\xi\right)\xb(\tau,\xi)=f(\tau,\xi);\quad \xb(\tau_{0},\xi)=\calD_{\tau}\xb(\tau_{0},\xi)=0,
	\end{align}
	then 
	\begin{align}\label{parametrix hbar inhomo}
		\begin{split}
		\xb(\tau,\xi)=\xi^{-\frac{1}{2}}\int_{\tau_{0}}^{\tau}\frac{\lambda(\tau)}{\lambda(\sigma)}\frac{\rho^{\frac{1}{2}}_{n}\left(\frac{\lambda(\tau)^{2}}{\lambda(\tau_{0})^{2}}\xi\right)}{\rho^{\frac{1}{2}}_{n}(\xi)}\sin\left(\lambda(\tau)\xi^{\frac{1}{2}}\int_{\sigma}^{\tau}\lambda(u)^{-1}du\right)f\left(\sigma,\frac{\lambda(\tau)^{2}}{\lambda(\sigma)^{2}}\xi\right)d\sigma
		\end{split}
	\end{align}
\end{lemma}
\subsection{Spaces for smooth sources}
In this section we strive to prove the estimates on the ``regular" part of $\xh(\tau,\xi)$. Here ``regular" means in $H^{5+}(\bbR^{2})$. We first introduce the following space:
\begin{align}\label{def S0 hbar}
\|f(\xi)\|_{\Sh_{0}}:=\|(\hbar^{2}\xi)^{1-\frac{\delta}{2}}\langle\hbar^{2}\xi\rangle^{\frac32+\delta}f(\xi)\|_{L^{2}_{d\xi}}.
\end{align}
Here $\delta>0$ is a sufficiently small constant. 


We have the following decay estimate for \eqref{parametrix hbar homo} in $\|\cdot\|_{\Sh_{0}}$:
\begin{proposition}\label{prop: homo para decay hbar}
	Let $\xb(\tau,\xi)$ be as in \eqref{parametrix hbar homo} with $\xb_{0}(\xi)$ and $\xb_{1}(\xi)$ satisfying
	\begin{align*}
		\|\xb_{0}\|_{\Sh_{0}}<\infty,\quad \|\xb_{1}\|_{\Sh_{1}}<\infty,\quad \textrm{where}\quad \|\cdot\|_{\Sh_{1}}:=\|\xi^{-\frac12}\cdot\|_{\Sh_{0}}.
	\end{align*}
	Then we have, for some constant $c>0$,
	\begin{align*}
		\|\xb(\tau,\cdot)\|_{\Sh_{0}}\lesssim \left(\frac{\tau}{\tau_{0}}\right)^{-3+c\nu}\left(\|\xb_{0}\|_{\Sh_{0}}+\|\xb_{1}\|_{\Sh_{1}}\right),
	\end{align*}
	where the implicit constant is independent of $\hbar$.
\end{proposition}
\begin{proof}
	Without loss of generality, we may assume that $\xb_{1}=0$. A direct calculation gives
	\begin{align*}
		\|(\hbar^{2}\xi)^{1-\frac{\delta}{2}}\langle\hbar^{2}\xi\rangle^{\frac32+\delta}\xb(\tau,\xi)\|_{L^{2}_{d\xi}}
		\lesssim& \left(\frac{\lambda(\tau)}{\lambda(\tau_{0})}\right)^{-2+\delta}\left\|\left(\hbar^{2}\frac{\lambda(\tau)^{2}}{\lambda(\tau_{0})^{2}}\xi\right)^{1-\frac{\delta}{2}}\left\langle\hbar^{2}\frac{\lambda(\tau)^{2}}{\lambda(\tau_{0})^{2}}\xi\right\rangle^{\frac32+\delta}\xb(\tau,\xi)\right\|_{L^{2}_{d\xi}}\\
		\lesssim &\left(\frac{\lambda(\tau)}{\lambda(\tau_{0})}\right)^{-1+\delta}\|\xb_{0}\|_{\Sh_{0}}.
	\end{align*}
	For any $\nu\in\left(0,\frac12\right)$, we simply choose $\delta>0$ sufficiently small such that $\delta\lesssim \nu^{2}$ to obtain the desired result.
\end{proof}
We also have the following estimate for \eqref{parametrix hbar inhomo} in $\|\cdot\|_{\Sh_{0}}$:
\begin{proposition}\label{prop: inhomo para decay hbar}
	Let $\xb(\tau,\xi)$ be as in \eqref{parametrix hbar inhomo} with $f$ satisfying, for some constant $c>0$,
	\begin{align*}
		\left(\frac{\sigma}{\tau_{0}}\right)^{3-c\nu}\cdot \sigma\left\|f(\sigma,\cdot)\right\|_{\Sh_{1}}<\infty,\quad \textrm{for}\quad \sigma\in[\tau_{0},\tau].
	\end{align*}
	Then for the same constant $c>0$, we have
	\begin{align*}
	\left(\frac{\tau}{\tau_{0}}\right)^{3-c\nu}	\cdot\left\|\xb(\tau,\cdot)\right\|_{\Sh_{0}}\lesssim \sup_{\tau_{0}\leq \sigma\leq \tau}\left(\frac{\sigma}{\tau_{0}}\right)^{3-c\nu}\cdot\sigma\left\|f(\sigma,\cdot)\right\|_{\Sh_{1}}.
	\end{align*}
	Here the implicit constant depends only on $\nu$ and $c$.
\end{proposition}
\begin{proof}
	A direct calculation gives, for $\delta\ll \nu^{2}$ and some $c'\ll1$,
	\begin{align*}
\left(\frac{\tau}{\tau_{0}}\right)^{3-c\nu}\cdot\left\|\xb(\tau,\cdot)\right\|_{\Sh_{0}}\lesssim& \left(\frac{\tau}{\tau_{0}}\right)^{3-c\nu}\cdot\int_{\tau_{0}}^{\tau}\frac{\lambda(\tau)}{\lambda(\sigma)}\left\|\frac{(\hbar^{2}\xi)^{1-\frac{\delta}{2}}}{\xi^{\frac12}}\langle\hbar^{2}\xi\rangle^{\frac32+\delta}\cdot f\left(\sigma,\frac{\lambda(\tau)^{2}}{\lambda(\sigma)^{2}}\xi\right)\right\|_{L^{2}_{d\xi}}\,d\sigma\\
\lesssim&\left(\frac{\tau}{\tau_{0}}\right)^{3-c\nu}\cdot\int_{\tau_{0}}^{\tau}\left(\frac{\lambda(\tau)}{\lambda(\sigma)}\right)^{-1+\delta}\|f(\sigma,\cdot)\|_{\Sh_{1}}\,d\sigma\\
\lesssim &\left(\frac{\tau}{\tau_{0}}\right)^{3-c\nu}\cdot\int_{\tau_{0}}^{\tau}\left(\frac{\sigma}{\tau}\right)^{3-c\nu}\cdot \left(\frac{\sigma}{\tau}\right)^{c'\nu}\sigma^{-1}\cdot\sigma\|f(\sigma,\cdot)\|_{\Sh_{1}}\,d\sigma\\
\lesssim & \sup_{\tau_{0}\leq \sigma\leq \tau}\left(\frac{\sigma}{\tau_{0}}\right)^{3-c\nu}\cdot\sigma\left\|f(\sigma,\cdot)\right\|_{\Sh_{1}},
	\end{align*}
	as claimed.
\end{proof}
Based on the proof of the above proposition, we have the following:
\begin{corollary}\label{cor: inhomo para decay hbar}
		Let $\xb(\tau,\xi)$ be as in \eqref{parametrix hbar inhomo} with $f$ satisfying, for some constant $c>0$,
	\begin{align*}
	\left(\frac{\sigma}{\tau_{0}}\right)^{3-c\nu}\cdot \sigma^{2}\left\|f(\sigma,\cdot)\right\|_{\Sh_{1}}<\infty,\quad \textrm{for}\quad \sigma\in[\tau_{0},\tau].
	\end{align*}
	Then for the same constant $c>0$, we have
	\begin{align*}
	\left(\frac{\tau}{\tau_{0}}\right)^{3-c\nu}	\cdot\left\|\xb(\tau,\cdot)\right\|_{\Sh_{0}}\ll_{\tau_{0}} \sup_{\tau_{0}\leq \sigma\leq \tau}\left(\frac{\sigma}{\tau_{0}}\right)^{3-c\nu}\cdot\sigma^{2}\left\|f(\sigma,\cdot)\right\|_{\Sh_{1}}.
	\end{align*}
\end{corollary}

We also have the following estimate on the operator $\calK^{(0)}_{\hbar}$ acting on the space $\Sh_{0}$ and $\Sh_{1}$:
\begin{proposition}\label{prop: F hbar boundedness}
	Let $\calK^{(0)}_{\hbar}$ be as in Proposition \ref{prop: K operator}. We have the bounds:
	\begin{align*}
		\left\|\calK^{(0)}_{\hbar}f\right\|_{\Sh_{0}}\lesssim \left\|f\right\|_{\Sh_{0}},\quad \left\|\calK^{(0)}_{\hbar}f\right\|_{\Sh_{1}}\lesssim \left\|f\right\|_{\Sh_{1}},\quad \left\|\calK_{\hbar}^{(0)}f\right\|_{\Sh_{1}}\lesssim \left\|f\right\|_{\Sh_{0}}
	\end{align*}
	uniformly in $\hbar$.
\end{proposition}
\begin{proof}
The idea is to reduce the operator to a simple Hilbert operator, by exploiting the differentiability properties of $F(\xi,\eta;\hbar)$. For simplicity, we denote $F(\xi,\eta;\hbar)$ by $F(\xi,\eta)$. We first prove the bounds for $\|\cdot\|_{\Sh_{0}}\rightarrow\|\cdot\|_{\Sh_{0}}$ and $\|\cdot\|_{\Sh_{1}}\rightarrow\|\cdot\|_{\Sh_{1}}$:
\begin{itemize}
	\item $\underline{\xi\simeq\eta}$. We reformulate the norms in terms of the variable $\xi^{\frac12}$, using the fact $\left\|f\right\|_{L^{2}_{d\xi}}\simeq \left\|\xi^{\frac14}f\right\|_{L^{2}_{d\xi^{\frac12}}}$. So the integral can be rewritten as
	\begin{align*}
		\int_{0}^{\infty}\chi_{\xi\simeq\eta}\frac{F(\xi,\eta)\rho_{n}(\eta)}{\xi-\eta}f(\eta)\,d\eta=2\int_{0}^{\infty}\chi_{\xi^{\frac12}\simeq\eta^{\frac12}}\frac{\eta^{\frac12}}{\xi^{\frac12}+\eta^{\frac12}}\cdot\frac{F(\xi,\eta)\rho_{n}(\eta)}{\xi^{\frac12}-\eta^{\frac12}}f(\eta)\,d\eta^{\frac12}.
	\end{align*}
If $|\xi^{\frac12}-\eta^{\frac12}|\geq 1$, we have, by Proposition \ref{prop: K operator},
\begin{align*}
\left|\chi_{\xi\simeq\eta}\frac{(\hbar^{2}\xi)^{1-\frac{\delta}{2}}\langle\hbar^{2}\xi\rangle^{\frac32+\delta}}{(\hbar^{2}\eta)^{1-\frac{\delta}{2}}\langle\hbar^{2}\eta\rangle^{\frac32+\delta}}\cdot\frac{\eta^{\frac12}}{\xi^{\frac12}+\eta^{\frac12}}\cdot\frac{F(\xi,\eta)}{\xi^{\frac12}-\eta^{\frac12}}\right|\lesssim \left|\xi^{\frac12}-\eta^{\frac12}\right|^{-\frac54},
\end{align*}
which is in $L^{1}_{d\eta^{\frac12}}$. This gives
\begin{align*}
	\left\|\int_{0}^{\infty}\chi_{\xi\simeq\eta,|\xi^{\frac12}-\eta^{\frac12}|\geq 1}\frac{\eta^{\frac12}}{\xi^{\frac12}+\eta^{\frac12}}\cdot\frac{F(\xi,\eta)\rho_{n}(\eta)}{\xi^{\frac12}-\eta^{\frac12}}f(\eta)d\eta^{\frac12}\right\|_{\Sh_{0}}\lesssim \left\|f\right\|_{\Sh_{0}}.
\end{align*}
When $|\xi^{\frac12}-\eta^{\frac12}|<1$, we write
\begin{align*}
	F(\xi,\eta)=F(\xi,\xi)+\partial_{\eta^{\frac12}}F\cdot O\left(|\xi^{\frac12}-\eta^{\frac12}|\right).
\end{align*}
Since Proposition \ref{prop: K operator} gives that $\left|\partial_{\eta^{\frac12}}F\right|\lesssim 1$ uniformly in $\hbar$, the contribution from the error $\partial_{\eta^{\frac12}}F\cdot O\left(|\xi^{\frac12}-\eta^{\frac12}|\right)$ again leads to a $L^{1}_{d\eta^{\frac12}}$-kernel, which can be handled as for the case $|\xi^{\frac12}-\eta^{\frac12}|\geq 1$. The contribution from the main term $F(\xi,\xi)$ is handled via the boundedness of $F(\xi,\eta)$ and the estimate
\begin{align*}
	\left\|\int_{0}^{\infty}\chi_{\xi^{\frac12}\simeq\eta^{\frac12},|\xi^{\frac12}-\eta^{\frac12}|<1}\frac{\eta^{\frac12}}{\xi^{\frac12}+\eta^{\frac12}}\cdot\frac{f(\eta)}{\xi^{\frac12}-\eta^{\frac12}}d\eta^{\frac12}\right\|_{L^{2}_{d\xi^{\frac12}}}\lesssim \left\|f\right\|_{L^{2}_{d\eta^{\frac12}}}.
\end{align*}
The bound for $\|\cdot\|_{\Sh_{1}}\rightarrow\|\cdot\|_{\Sh_{1}}$ is handled similarly using the fact $\xi\simeq\eta$.

\item $\underline{\xi\gg\eta}$. We split the integral into the cases $\hbar^{2}\xi\geq 1$ and $\hbar^{2}\xi\leq 1$. For the latter we have
\begin{align*}
	&\left\|\chi_{\hbar^{2}\xi\leq 1}(\hbar^{2}\xi)^{1-\frac{\delta}{2}}\int_{0}^{\infty}\frac{F(\xi,\eta)\rho_{n}(\eta)}{(\xi-\eta)(\hbar^{2}\eta)^{1-\frac{\delta}{2}}}(\hbar^{2}\eta)^{1-\frac{\delta}{2}}f(\eta)d\eta\right\|_{L^{2}_{d\xi}}\\
	\lesssim &\hbar\left\|\chi_{\hbar^{2}\xi\leq 1}\frac{(\hbar^{2}\xi)^{1-\frac{\delta}{2}}}{\hbar^{2}\xi}\right\|_{L^{2}_{\hbar^{2}d\xi}}\cdot\int_{0}^{\infty}\chi_{\eta\ll\xi}(\hbar^{2}\eta)^{-\frac12+\frac{\delta}{2}}\cdot(\hbar^{2}\eta)^{1-\frac{\delta}{2}}f(\eta)d\eta.
\end{align*}
Here we have used the bound $|F(\xi,\eta)|\lesssim \hbar\eta^{\frac12}$ for $\hbar\eta^{\frac12}\leq 1$. The desired estimate follows from the bound
\begin{align*}
\left|	\int_{0}^{\infty}\chi_{\eta\ll\xi}(\hbar^{2}\eta)^{-\frac12+\frac{\delta}{2}}\cdot(\hbar^{2}\eta)^{1-\frac{\delta}{2}}f(\eta)d\eta\right|
\lesssim \left\|\chi_{\eta\ll\xi}(\hbar^{2}\eta)^{-\frac12+\frac{\delta}{2}}\right\|_{L^{2}_{d\eta}}\cdot\left\|f\right\|_{\Sh_{0}}\lesssim\hbar^{-1}\|f\|_{\Sh_{0}}.
\end{align*}
The estimate for $\|\cdot\|_{\Sh_{1}}\rightarrow\|\cdot\|_{\Sh_{1}}$ follows in the same way using the fact $\xi\geq \eta$. \\

For the case $\hbar^{2}\xi\geq1$, we have, by Proposition \ref{prop: K operator}, that for $\hbar^{2}\eta\leq 1$
\begin{align*}
	&\left\|\chi_{\hbar^{2}\xi\geq 1}(\hbar^{2}\xi)^{\frac52+\frac{\delta}{2}}\int_{0}^{\infty}\frac{F(\xi,\eta)\rho_{n}(\eta)}{(\xi-\eta)(\hbar^{2}\eta)^{1-\frac{\delta}{2}}}(\hbar^{2}\eta)^{1-\frac{\delta}{2}}f(\eta)d\eta\right\|_{L^{2}_{d\xi}}\\
	\lesssim&\hbar\left\|\chi_{\hbar^{2}\xi\geq 1}\frac{(\hbar^{2}\xi)^{\frac52+\frac{\delta}{2}}}{(\hbar^{2}\xi)^{\frac72}}\right\|_{L^{2}_{\hbar^{2}d\xi}}\cdot\int_{0}^{\infty}\chi_{\eta\ll\xi,\hbar^{2}\eta\leq 1}(\hbar^{2}\eta)^{-\frac12+\frac{\delta}{2}}\cdot(\hbar^{2}\eta)^{1-\frac{\delta}{2}}f(\eta)d\eta,
\end{align*}
which is bounded similarly as for the case $\hbar^{2}\xi\leq 1$.\\

 If $\hbar^{2}\eta\geq 1$, Proposition \ref{prop: K operator} implies
\begin{align*}
	\left|\frac{F(\xi,\eta)}{\xi-\eta}\right|\lesssim \hbar^{2}(\hbar^{2}\eta)^{\frac34}(\hbar^{2}\xi)^{-\frac94}\cdot (\hbar^{2}\eta)^{\frac54}(\hbar^{2}\xi)^{-\frac54}\leq\hbar^{2}(\hbar^{2}\eta)^{2}(\hbar^{2}\xi)^{-\frac72}
\end{align*}
Therefore we have the bound
\begin{align*}
	&\left\|\chi_{\hbar^{2}\xi\geq 1}(\hbar^{2}\xi)^{\frac52+\frac{\delta}{2}}\int_{0}^{\infty}\frac{F(\xi,\eta)\rho_{n}(\eta)}{(\xi-\eta)(\hbar^{2}\eta)^{\frac52+\frac{\delta}{2}}}(\hbar^{2}\eta)^{\frac52+\frac{\delta}{2}}f(\eta)d\eta\right\|_{L^{2}_{d\xi}}\\
	\lesssim &\hbar\left\|\chi_{\hbar^{2}\xi\geq 1}\frac{(\hbar^{2}\xi)^{\frac52+\frac{\delta}{2}}}{(\hbar^{2}\xi)^{\frac72}}\right\|_{L^{2}_{\hbar^{2}d\xi}}\cdot
\int_{0}^{\infty}\chi_{\eta\ll\xi,\hbar^{2}\eta\geq 1}(\hbar^{2}\eta)^{-\frac12-\frac{\delta}{2}}\cdot(\hbar^{2}\eta)^{\frac52+\frac{\delta}{2}}f(\eta)d\eta,
\end{align*}
which again is bounded in the similar way. The estimate for $\|\cdot\|_{\Sh_{1}}\rightarrow\|\cdot\|_{\Sh_{1}}$ is handled similarly using the fact $\xi\geq \eta$.

\item $\underline{\xi\ll\eta}$. We start with the case when $\hbar^{2}\xi\geq 1$. By Proposition \ref{prop: K operator}
\begin{align*}
\left|\frac{F(\xi,\eta)}{\eta-\xi}\right|\lesssim \hbar^{2}(\hbar^{2}\xi)^{\frac34}(\hbar^{2}\eta)^{-\frac94}\cdot (\hbar^{2}\xi)^{\frac54}(\hbar^{2}\eta)^{-\frac54}\leq \hbar^{2}(\hbar^{2}\xi)^{-\frac72}(\hbar^{2}\eta)^{2},
\end{align*}
which means that this case can be handled in the same way as the case when $\hbar^{2}\xi\gg\hbar^{2}\eta\geq 1$. If $\hbar^{2}\xi\leq 1$ and $\hbar^{2}\eta\geq1$, we have
\begin{align*}
	\left|\frac{F(\xi,\eta)}{\eta-\xi}\right|\lesssim \hbar^{2}(\hbar^{2}\xi)^{\frac12}(\hbar^{2}\eta)^{-\frac94},
\end{align*}
which is surely enough to bound the integral
\begin{align*}
	\left\|\chi_{\hbar^{2}\xi\leq 1,\hbar^{2}\eta\geq 1}(\hbar^{2}\xi)^{1-\frac{\delta}{2}}\int_{0}^{\infty}\frac{F(\xi,\eta)\rho_{n}(\eta)}{(\xi-\eta)(\hbar^{2}\eta)^{\frac52+\frac{\delta}{2}}}(\hbar^{2}\eta)^{\frac52+\frac{\delta}{2}}f(\eta)d\eta\right\|_{L^{2}_{d\xi}}
\end{align*}
Finally for the case $\hbar^{2}\xi\ll\hbar^{2}\eta\leq 1$, we have
\begin{align*}
	\left|\frac{F(\xi,\eta)}{\xi-\eta}\right|\lesssim \hbar^{2}(\hbar^{2}\xi)^{\frac12}(\hbar^{2}\eta)^{-\frac12}|\log\hbar^{2}\eta|\lesssim \hbar^{2}(\hbar^{2}\xi)^{-1}(\hbar^{2}\eta)|\log \hbar^{2}\eta|,
\end{align*}
which is the desired estimate. \\

For the bound for $\|\cdot\|_{\Sh_{1}}\rightarrow\|\cdot\|_{\Sh_{1}}$, we proceed as follows. For $\hbar\ll1$, since we can have as many powers of $\frac{\xi^{\frac12}}{\eta^{\frac12}}$ as we wish, the argument is the same as that for $\|\cdot\|_{\Sh_{0}}\rightarrow\|\cdot\|_{\Sh_{0}}$. For $\hbar\simeq1$, if $1\leq \hbar^{2}\xi\ll\hbar^{2}\eta$, then the argument is the same as for $\hbar\ll1$. For the remaining two cases, we use the bound:
\begin{align*}
		\left|\frac{F(\xi,\eta)}{\xi-\eta}\right|\lesssim &\hbar^{2}(\hbar^{2}\xi)^{-\frac12}(\hbar^{2}\eta)^{-\frac54},\quad\textrm{for}\quad \hbar^{2}\xi\leq 1\leq \hbar^{2}\eta,\\
		\left|\frac{F(\xi,\eta)}{\xi-\eta}\right|\lesssim &\hbar^{2}(\hbar^{2}\xi)^{-\frac12}(\hbar^{2}\eta)^{\frac12}|\log\hbar^{2}\eta|,\quad \textrm{for}\quad \hbar^{2}\xi\ll\hbar^{2}\eta\leq 1.
\end{align*}
\end{itemize}
The bound for $\|\cdot\|_{\Sh_{0}}\rightarrow\|\cdot\|_{\Sh_{1}}$ can be handled as follows: First, if $\hbar^{2}\xi\geq 1$, the desired bound follows from the bound for $\|\cdot\|_{\Sh_{0}}\rightarrow\|\cdot\|_{\Sh_{0}}$. If $\hbar^{2}\xi\leq 1$, the argument is similar to that of deriving the bound for $\|\cdot\|_{\Sh_{0}}\rightarrow\|\cdot\|_{\Sh_{0}}$, and here we only give an outline:
\begin{itemize}
	\item When $\xi\simeq\eta$ and $|\xi^{\frac12}-\eta^{\frac12}|<1$, for $\hbar\gtrsim1$ we have $\partial_{\eta^{\frac12}}F(\xi,\eta)\simeq \xi^{\frac12}|\log\eta|$. Here $O(|\xi^{\frac12}-\eta^{\frac12}|)$ cancels with singular denominator and $L^{1}_{d\eta^{\frac12}}$-integrability of $|\log\eta|$ gives the desired bound.
	
	\item When $\hbar^{2}\xi\leq 1$ and $\hbar^{2}\eta\geq 1$, the same bound
	\begin{align*}
		\left|\frac{F(\xi,\eta)}{\eta-\xi}\right|\lesssim \hbar^{2}(\hbar^{2}\xi)^{\frac12}(\hbar^{2}\eta)^{-\frac94}
	\end{align*}
	suffices.
	
	\item When $\hbar^{2}\xi\ll\hbar^{2}\eta\leq 1$, we have
	\begin{align*}
		\left|\frac{F(\xi,\eta)}{\xi-\eta}\right|\lesssim \hbar^{2}(\hbar^{2}\xi)^{\frac12}(\hbar^{2}\eta)^{-\frac12}|\log\hbar^{2}\eta|\lesssim \hbar^{2}(\hbar^{2}\xi)^{\frac12}(\hbar^{2}\xi)^{-\frac32+\delta}(\hbar^{2}\eta)^{1-\delta}|\log\hbar^{2}\eta|,
	\end{align*}
	which suffices.
\end{itemize}
\end{proof}
\subsection{Linear smooth sources}
Now we start to estimate the contribution of the smooth part of $\xh(\tau,\xi)$ to the RHS of the equation \eqref{eq tau xi hbar D}. We start with the contribution to $\calR_{\hbar}(\xh,\calD_{\tau}\xh)$:
\begin{proposition}\label{prop: xh good para linear}
	Let $S$ be any of the expressions in $\calR(\xh,\calD_{\tau}\xh)$, and we define
	\begin{align*}
		\txh(\tau,\xi):=\xi^{-\frac{1}{2}}\int_{\tau_{0}}^{\tau}\frac{\lambda(\tau)}{\lambda(\sigma)}\frac{\rho^{\frac{1}{2}}_{n}\left(\frac{\lambda(\tau)^{2}}{\lambda(\tau_{0})^{2}}\xi\right)}{\rho^{\frac{1}{2}}_{n}(\xi)}\sin\left(\lambda(\tau)\xi^{\frac{1}{2}}\int_{\sigma}^{\tau}\lambda(u)^{-1}du\right)S\left(\sigma,\frac{\lambda(\tau)^{2}}{\lambda(\sigma)^{2}}\xi\right)d\sigma.
	\end{align*}
	Then we have, for $S=2\dftau\calK^{(0)}_{\hbar}\xh$,
	\begin{align*}
		\left\|S\right\|_{\Sh_{0}}\lesssim \left\|\txh(\tau,\cdot)\right\|_{\Sh_{0}}+\left\|\calD_{\tau}\txh(\tau,\cdot)\right\|_{\Sh_{1}}.
	\end{align*}
If $S$ is one of the terms $\left(\dftau\right)^{\prime}\calK^{(0)}_{\hbar}\xh, \left(\dftau\right)^{2}\calK^{(0)}_{\hbar}\xh, \left(\dftau\right)^{\prime}\left(\calK^{(0)}_{\hbar}\right)^{2}\xh$ or $\dftau\left[\calD_{\tau},\calK^{(0)}_{\hbar}\right]\xh$, then we have
	\begin{align*}
\left\|S\right\|_{\Sh_{0}}\ll\left\|\txh(\tau,\cdot)\right\|_{\Sh_{0}}+\left\|\calD_{\tau}\txh(\tau,\cdot)\right\|_{\Sh_{1}}.
\end{align*}
\end{proposition}
\begin{proof}
	Except the commutator $\dftau\left[\calD_{\tau},\calK^{(0)}_{\hbar}\right]\xh$, the estimates for the terms in $\calR_{\hbar}(\xh,\calD_{\tau}\xh)$ follow directly from Proposition \ref{prop: inhomo para decay hbar} and Proposition \ref{prop: F hbar boundedness}. Now we consider the commutator term
	\begin{align*}
		\dftau\left[\calD_{\tau},\calK^{(0)}_{\hbar}\right]\xh=-2\left(\dftau\right)^{2}\left[\xi\partial_{\xi},\calK^{(0)}_{\hbar}\right]\xh-\left(\dftau\right)^{2}\left[\frac{\rho^{\prime}_{n}(\xi)\xi}{\rho_{n}(\xi)},\calK^{(0)}_{\hbar}\right]\xh.
	\end{align*}
	Since $\left|\frac{\rho^{\prime}_{n}(\xi)\xi}{\rho_{n}(\xi)}\right|\simeq 1$, the second commutator on the RHS above can be handled similarly as the other terms in $\calR_{\hbar}(\xh,\calD_{\tau}\xh)$. We focus on the first commutator. In the following calculation, we omit the factor $\left(\dftau\right)^{2}$ and write $2\xi\partial_{\xi}=\xi^{\frac12}\partial_{\xi^{\frac12}}$. For a function $f\in C^{\infty}_{c}(\bbR_{+})$, we consider
	\begin{align*}
	\left[\xi^{\frac12}\partial_{\xi^{\frac12}},\calK^{(0)}_{\hbar}\right]f(\xi)=&\int_{0}^{\infty}\xi^{\frac12}\partial_{\xi^{\frac12}}\left(\frac{F(\xi,\eta)}{\xi-\eta}\right)\rho_{n}(\eta)f(\eta)\,d\eta-\int_{0}^{\infty}\frac{F(\xi,\eta)}{\xi-\eta}\rho_{n}(\eta)\eta^{\frac12}\partial_{\eta^{\frac12}}f(\eta)d\eta:=I+II.
	\end{align*}
	For $II$ we use integration by parts to obtain
	\begin{align*}
		II=&2\int_{0}^{\infty}\partial_{\eta^{\frac12}}\left(\frac{F(\xi,\eta)\eta}{\xi-\eta}\rho_{n}(\eta)\right)d\eta^{\frac12}\\
		=&-2\int_{0}^{\infty}\partial_{\eta^{\frac12}}\left(F(\xi,\eta)\rho_{n}(\eta)\right)f(\eta)d\eta^{\frac12}+2\xi\int_{0}^{\infty}\partial_{\eta^{\frac12}}\left(\frac{F(\xi,\eta)}{\xi-\eta}\rho_{n}(\eta)\right)f(\eta)d\eta^{\frac12}\\
		=&-\int_{0}^{\infty}\frac{1}{\eta^{\frac12}}\partial_{\eta^{\frac12}}\left(\frac{F(\xi,\eta)\eta^{\frac12}}{\xi^{\frac12}+\eta^{\frac12}}\rho_{n}(\eta)\right)f(\eta)d\eta-\int_{0}^{\infty}\frac{\xi^{\frac12}}{\eta^{\frac12}}\partial_{\eta^{\frac12}}\left(\frac{F(\xi,\eta)}{\xi^{\frac12}+\eta^{\frac12}}\rho_{n}(\eta)\right)f(\eta)\,d\eta\\
		&+\int_{0}^{\infty}\frac{\xi}{\eta^{\frac12}}\partial_{\eta^{\frac12}}\left(\frac{F(\xi,\eta)\rho_{n}(\eta)}{(\xi^{\frac12}+\eta^{\frac12})(\xi^{\frac12}-\eta^{\frac12})}\right)d\eta:=II_{1}+II_{2}+II_{3}.
	\end{align*}
	We write $II_{2}$ as
	\begin{align*}
		II_{2}=-\int_{0}^{\infty}\frac{\xi}{\eta^{\frac12}}\partial_{\eta^{\frac12}}\left(\frac{F(\xi,\eta)\rho_{n}(\eta)}{(\xi^{\frac12}+\eta^{\frac12})(\xi^{\frac12}-\eta^{\frac12})}\right)f(\eta)d\eta+\int_{0}^{\infty}\frac{\xi^{\frac12}}{\eta^{\frac12}}\partial_{\eta^{\frac12}}\left(\frac{F(\xi,\eta)\eta^{\frac12}\rho_{n}(\eta)}{(\xi^{\frac12}+\eta^{\frac12})(\xi^{\frac12}-\eta^{\frac12})}\right)f(\eta)d\eta.
	\end{align*}
	Therefore
	\begin{align*}
		I+II_{2}+II_{3}=&\int_{0}^{\infty}\frac{\xi^{\frac12}}{\eta^{\frac12}}\left(\partial_{\xi^{\frac12}}+\partial_{\eta^{\frac12}}\right)\left(\frac{F(\xi,\eta)\eta^{\frac12}\rho_{n}(\eta)}{(\xi^{\frac12}+\eta^{\frac12})(\xi^{\frac12}-\eta^{\frac12})}\right)f(\eta)d\eta\\
		=&\int_{0}^{\infty}\frac{1}{\xi-\eta}\cdot\frac{\xi^{\frac12}(\xi^{\frac12}+\eta^{\frac12})}{\eta^{\frac12}}\left(\partial_{\xi^{\frac12}}+\partial_{\eta^{\frac12}}\right)\left(\frac{F(\xi,\eta)\eta^{\frac12}\rho_{n}(\eta)}{\xi^{\frac12}+\eta^{\frac12}}\right)f(\eta)d\eta.
	\end{align*}
	So we have
	\begin{align*}
		\left[\xi^{\frac12}\partial_{\xi^{\frac12}},\calK^{(0)}_{\hbar}\right]f(\xi)=\int_{0}^{\infty}\frac{H_{1}(\xi,\eta)}{\xi-\eta}f(\eta)d\eta+\int_{0}^{\infty}H_{2}(\xi,\eta)f(\eta)d\eta,
	\end{align*}
	where
	\begin{align}\label{H1 H2}
		H_{1}(\xi,\eta):=\frac{\xi^{\frac12}(\xi^{\frac12}+\eta^{\frac12})}{\eta^{\frac12}}\left(\partial_{\xi^{\frac12}}+\partial_{\eta^{\frac12}}\right)\left(\frac{F(\xi,\eta)\eta^{\frac12}\rho_{n}(\eta)}{\xi^{\frac12}+\eta^{\frac12}}\right),\quad H_{2}(\xi,\eta):=-\frac{1}{\eta^{\frac12}}\partial_{\eta^{\frac12}}\left(\frac{F(\xi,\eta)\eta^{\frac12}}{\xi^{\frac12}+\eta^{\frac12}}\rho_{n}(\eta)\right).
	\end{align}
	Let us first look at the contribution from $H_{2}(\xi,\eta)$, whose principal contribution is given by $-\frac{\partial_{\eta^{\frac12}}F(\xi,\eta)}{\xi^{\frac12}+\eta^{\frac12}}\rho_{n}(\eta)$. Here we try to mimic the proof of Proposition \ref{prop: F hbar boundedness}:
	\begin{itemize}
		\item $\underline{\xi\simeq \eta}$. After change of the variable $\eta\rightarrow\eta^{\frac12}$, we consider
		\begin{align}\label{kernel H2}
			\left|\chi_{\xi\simeq\eta}\hbar\frac{(\hbar^{2}\xi)^{\frac12-\frac{\delta}{2}}\langle\hbar^{2}\xi\rangle^{\frac32+\delta}}{(\hbar^{2}\eta)^{1-\frac{\delta}{2}}\langle\hbar^{2}\eta\rangle^{\frac32+\delta}}\cdot\frac{\eta^{\frac12}\partial_{\eta^{\frac12}}F(\xi,\eta)}{\xi^{\frac12}+\eta^{\frac12}}\right|.
		\end{align}
		If $\hbar^{2}\eta\gtrsim1$, Proposition \ref{prop: K operator} shows that the above absolute value is bounded by $\hbar\cdot(\hbar\eta^{\frac12})^{-2}$, which is in $L^{1}_{d\eta^{\frac12}}$. If $\hbar^{2}\eta\ll1$, Proposition \ref{prop: K operator} shows that the above absolute value is bounded by $\hbar\log(\hbar^{2}\eta)$, which is again in $L^{1}_{d\eta^{\frac12}}$.
		
		\item $\underline{\xi\ll\eta}$. This can be treated in the same way as the case $\xi\simeq\eta$.
		
		\item $\underline{\xi\gg\eta}$. If $\hbar^{2}\eta\geq 1$, Proposition \ref{prop: K operator} shows that \eqref{kernel H2}is bounded by
		\begin{align*}
			\hbar\cdot(\hbar^{2}\xi)^{2+\frac{\delta}{2}}\cdot(\hbar^{2}\xi)^{-\frac54}\cdot(\hbar^{2}\eta)^{-\frac52-\frac{\delta}{2}}\cdot(\hbar^{2}\eta)^{\frac34}\cdot(\hbar^{2}\xi)^{-\frac12}\cdot\left(\frac{\hbar^{2}\eta}{\hbar^{2}\xi}\right)^{N}\leq \hbar(\hbar^{2}\eta)^{-\frac32-\frac{\delta}{2}},
		\end{align*}
		which is in $L^{1}_{d\eta^{\frac12}}$. If $\hbar^{2}\eta\leq 1$ and $\hbar^{2}\xi\geq 1$, then \eqref{kernel H2} is bounded by
		\begin{align*}
			\hbar\cdot(\hbar^{2}\xi)^{\frac52+\frac{\delta}{2}}\cdot(\hbar^{2}\eta)^{-1+\frac{\delta}{2}}\cdot(\hbar^{2}\eta)^{-\frac12}\cdot (\hbar^{2}\xi)^{-\frac14}\frac{\hbar^{2}\eta}{\hbar^{2}\xi}\cdot (\hbar^{2}\xi)^{-N}\leq (\hbar^{2}\eta)^{-\frac12+\frac{\delta}{2}},
		\end{align*}
		which is again in $L^{1}_{d\eta}$. Finally if $\hbar^{2}\eta\ll\hbar^{2}\xi\leq 1$, Proposition 3.1 shows that \eqref{kernel H2} is bounded by
		\begin{align*}
			\hbar\cdot(\hbar^{2}\xi)^{\frac12-\frac{\delta}{2}}\cdot(\hbar^{2}\eta)^{-1+\frac{\delta}{2}}\cdot(\hbar^{2}\eta)^{\frac12}|\log\hbar^{2}\eta|\leq \hbar|\log\hbar^{2}\eta|,
		\end{align*}
		which is again in $L^{1}_{d\eta^{\frac12}}$. This completes the estimate for $H_{2}(\xi,\eta)$.
	\end{itemize}
Next we turn to the estimate for $\frac{H_{1}(\xi,\eta)}{\xi-\eta}$. Systematically $H_{1}(\xi,\eta)$ is a linear combination of the following terms:
\begin{align}\label{kernel H1 profile}
	\xi^{\frac12}\partial_{\xi^{\frac12}}F(\xi^{\frac12},\eta^{\frac12})\rho_{n}(\eta), \quad \xi^{\frac12}\partial_{\eta^{\frac12}}F(\xi^{\frac12},\eta^{\frac12})\rho_{n}(\eta),\quad \frac{\xi^{\frac12}}{\xi^{\frac12}+\eta^{\frac12}}F(\xi,\eta)\rho_{n}(\eta),\quad \frac{\hbar\xi^{\frac12}}{\hbar\eta^{\frac12}}F(\xi,
	\eta)\rho_{n}(\eta)
\end{align}
The contribution from the third term in \eqref{kernel H1 profile} can be treated exactly the same as $\calK^{(0)}_{\hbar}$ itself. For the fourth term in \eqref{kernel H1 profile}, $\hbar\xi^{\frac12}$ in numerator cancels with the $\xi^{-\frac12}$ in the definition of $\Sh_{1}$, giving an extra $\hbar$. In the denominator, if $\hbar\eta^{\frac12}\leq 1$, then we use the small factors of $\hbar\eta^{\frac12}$ in $F(\xi,\eta)$ to absorb it. Therefore the argument is similar as the one treating the $\Sh_{0}\rightarrow\Sh_{1}$ bound for $\calK^{(0)}_{\hbar}$. The second term in \eqref{kernel H1 profile} is treated in the same way, except that $\xi^{\frac12}$ cancels exactly with the $\xi^{-\frac12}$ in the definition of $\Sh_{1}$, without an extra $\hbar$. This approach also applies to the first term in \eqref{kernel H1 profile}. Therefore we obtain the following estimate on $\left[\calD_{\tau},\calK^{(0)}_{\hbar}\right]$:
\begin{align*}
	\left\|\left[\calD_{\tau},\calK^{(0)}_{\hbar}\right]f\right\|_{\Sh_{1}}\lesssim \left\|f\right\|_{\Sh_{0}},
\end{align*}
which completes the proof.
\end{proof}
\subsection{Nonlinear smooth sources}
Observing the expressions \eqref{coe eq precise 1} and \eqref{coe eq precise 2} for $N(\vphi_{1})$ and $N(\vphi_{2})$, there are two types of nonlinearities in $N(\vphi_{1})$ and $N(\vphi_{2})$: One is a quadratic expression of $\vphi_{i}, i=1,2$ and their derivatives multiplied by a spatial weight depending on the background solution, the other is $\vphi_{i}, i=1,2$ multiplied by a quadratic expression of the derivatives of $\vphi_{i},i=1.2$. Since in this section we only handle smooth sources and the background solution has limited smoothness, now we focus on the second type of nonlinearity. More precisely, in this section we derive multilinear estimates for products of regular functions, i.e., those whose Fourier transforms are in $\Sh_0$, where the $\hbar$ may differ amongst the factors. The key point here is that in these estimates, we do not lose in the smallest $\hbar$, i.e., the estimates are {\it{uniform with respect to the smallest $\hbar$ present}}. At the most basic level we need to bound the expressions
\[
\partial_R\phi_1\cdot \partial_R\phi_2,\quad \left(\partial_\tau + \frac{\lambda'}{\lambda}R\partial_R\right)\phi_1\cdot \left(\partial_\tau + \frac{\lambda'}{\lambda}R\partial_R\right)\phi_2, 
\]
We start with the following basic
\begin{proposition}\label{prop:derivative1} Let
	\[
	f(R): = \int_0^\infty \phi_{n}(R;\xi)\xb(\xi)\rho_{n}(\xi)\,d\xi,
	\]
	we have the bound, for $n\gg k$,
	\begin{align*}
	\left\|\partial_R^k f(R)\right\|_{L^2_{R\,dR}}\lesssim_k \left\|\xi^{\frac{k}{2}}\xb(\xi)\right\|_{L^2_{d\xi}}, 
	\end{align*}
	the implied constant being uniform in $\hbar$. 
\end{proposition}
\begin{proof}
	 We cover the case $n\gg 1$, the remaining case $n\lesssim 1$ being simpler since no attention needs to be paid to the $\hbar$-dependence, and can be handled similarly by Proposition \ref{prop: FB match}. Moreover, the cases $k\geq 2$ being similar to the one of $k = 1$, we treat this latter one in detail. We do this by explicit use of the asymptotics of $\phi(R;\xi,\hbar)$, decomposing  the frequency regime into a small frequency, turning point, and oscillatory regime:
	 \\
	 
	 {\it{Contribution of the small frequency regime}}. Here we consider the norm of the expression 
	 \[
	 \int_0^\infty\chi_{R\xi^{\frac12}\hbar<\frac{x_t}{2}}\partial_{R}\phi_{n}(R;\xi)x(\xi)\rho_{n}(\xi)\,d\xi,
	 \]
	 where the cutoff smoothly localizes to the indicated region. Then using the representation from Proposition \ref{prop:DFT nlarge}, we easily infer the bound 
	 \begin{align*}
	 	\left|\chi_{R\xi^{\frac12}\hbar<\frac{x_{t}}{2}}\partial_{R}\phi_{n}(R,\xi)\right|\lesssim R^{-1}\cdot\hbar^{n-\frac12}\xi^{\frac{n-1}{2}}R^{n-1}\leq \hbar^{\frac12}\cdot\hbar\xi^{\frac12}\cdot\left(\hbar\xi^{\frac12}R\right)^{n-2}
	 \end{align*}
	 Hence subdividing this frequency region into dyadic sub-intervals, we infer 
	 \begin{align*}
	 	\begin{split}
	 &	\left\|\int_0^\infty\chi_{R\xi^{\frac12}\hbar<\frac{x_t}{2}}\partial_{R}\phi_{n}(R;\xi)\xb(\xi)\rho_{n}(\xi)\,d\xi\right\|_{L^{2}_{R\,dR}}\\
	 \lesssim &\sum_{\mu<\frac{x_{t}}{2}}\mu^{n-2}\left\|\int_{0}^{\infty}\chi_{R\xi^{\frac12}\hbar\simeq \mu}\hbar^{\frac32}\xi^{\frac12}\cdot\xb(\xi)\rho_{n}(\xi)d\xi\right\|_{L^{2}_{R\,dR}}\\
	 \lesssim &\sum_{\mu<\frac{x_{t}}{2}}\mu^{n-2}\left(\sum_{\lambda}\left\|\int_{0}^{\infty}\chi_{R\xi^{\frac12}\hbar\simeq \mu}\hbar^{\frac32}\xi^{\frac12}\cdot\chi_{\xi^{\frac12}\simeq \lambda}\xb(\xi)\rho_{n}(\xi)d\xi\right\|_{L^{2}_{R\,dR}}\right)^{\frac12}\\
	 \lesssim&\sum_{\mu<\frac{x_{t}}{2}}\mu^{n-2}\left(\sum_{\lambda}\left(\hbar^{\frac32}\frac{\mu}{\hbar\lambda}\right)^{2}\cdot\lambda^{2}\left\|\chi_{\xi^{\frac12}\simeq\lambda}\xi^{\frac12}\xb(\xi)\right\|^{2}_{L^{2}_{d\xi}}\right)^{\frac12}\\
	 \lesssim &\hbar^{\frac12}\left\|\xi^{\frac12}\xb(\xi)\right\|_{L^{2}_{d\xi}}\sum_{\mu<\frac{x_{t}}{2}}\mu^{n-1}\lesssim\hbar^{\frac12}\left\|\xi^{\frac12}\xb(\xi)\right\|_{L^{2}_{d\xi}}.
	 	\end{split}
	 \end{align*}
	 We have used H\"older's inequality to obtain the second line from the bottom.
	 \\
	 
	 {\it{Contribution of the turning point}}. Here we include the smooth cutoff $\chi_{R\xi^{\frac12}\hbar\in\left[\frac{x_t}{2},2x_t\right]}$. By means of a sharp cutoff, we split this into two portions 
	 \begin{align*}
	 \int_0^\infty\chi_{R\xi^{\frac12}\hbar\in\left[\frac{x_t}{2},2x_t\right]}\partial_{R}\phi_{n}(R;\xi)\xb(\xi)\rho_{n}(\xi)\,d\xi &=  \int_0^\infty\chi_{R\xi^{\frac12}\hbar\in\left[\frac{x_t}{2},x_t\right]}\partial_{R}\phi_{n}(R;\xi)\xb(\xi)\rho_{n}(\xi)\,d\xi\\
	 &+ \int_0^\infty\chi_{R\xi^{\frac12}\hbar\in\left[x_t,2x_t\right]}\partial_{R}\phi_{n}(R;\xi)\xb(\xi)\rho_{n}(\xi)\,d\xi
	 \end{align*}
	 Correspondingly, we get a contribution from the non-oscillatory regime to the left of the turning point, and a contribution of the oscillatory regime to the right. 
	 \\
	 {\it{Non-oscillatory regime $R\xi^{\frac12}\hbar\in\left[\frac{x_t}{2},x_t\right]$.}}
	 Invoking again Proposition \ref{prop:DFT nlarge} and the estimates for $\tau$ in the regime $R\xi^{\frac12}\hbar\in [\frac{x_t}{2},x_t]$ which come from Lemma \ref{lem: Lemma 3.2 CDST}, we infer in this regime the bound 
	 \[
	 \left|\chi_{R\xi^{\frac12}\hbar\in[\frac{x_t}{2},x_t]}\partial_R\phi_{n}(R;\xi)\right|\lesssim \hbar^{\frac23}\xi^{\frac12}. 
	 \]
	 This time we have to be more careful to avoid any $\hbar$-losses since we no longer get the very rapid exponential decay in $\hbar^{-1}$ for free, compared to the preceding case. Thus for $\hbar\ll 1$ we have to use finer asymptotics of Airy functions. Setting $x_{t} - x = \hbar^{\frac23}\mu$ with $\lambda\geq 1,\,x = R\xi^{\frac12}\hbar$, we have 
	 \[
	 \left|\chi_{R\xi^{\frac12}\hbar\in[\frac{x_t}{2},x_t]}\partial_R\phi_{n}(R;\xi)\right|\lesssim \hbar^{\frac23}\xi^{\frac12}\mu^{-\frac14}e^{-\frac23\mu^{\frac32}}. 
	 \]
	 Now we divide the $\xi$ integral into one where $x_t-x\in \left[0,\hbar^{\frac23}\right]$ and its complement. For the contribution of the first of these, we observe that for fixed $R$, and using the fact that
	 \[
	 \chi_{R\xi^{\frac12}\hbar\in\left[\frac{x_t}{2},2x_t\right]}\partial_{\xi}\left(x - x_t\right)\simeq \chi_{R\xi^{\frac12}\hbar\in\left[\frac{x_t}{2},2x_t\right]}\cdot\xi^{-1}
	 \]
	 thanks to Lemma \ref{lem: monotonicity of root in alpha}, we have that $x_t-x\in \left[0,\hbar^{\frac23}\right]$ implies that $\xi$ is confined to an interval $I_R$ of length $\simeq \frac{\hbar^{\frac23}}{(R\hbar)^2}$. So if $R\simeq \lambda$, then 
	 \begin{align*}
	 \left|\int_0^\infty\chi_{x_t - x\in [0,\hbar^{\frac23}]}\partial_R\phi_{n}(R;\xi)\xb(\xi)\rho_{n}(\xi)\,d\xi\right|&\lesssim \hbar^{\frac23}\cdot \left|I_R\right|^{\frac12}\cdot\left\|\xi^{\frac12}\xb\right\|_{L^2(\xi\simeq \lambda^{-2}\hbar^{-2})}
	 \lesssim \frac{1}{\lambda}\cdot \left\|\xi^{\frac12}\xb\right\|_{L^2(\xi\simeq \lambda^{-2}\hbar^{-2})}. 
	 \end{align*}
	 Here we have used H\"older's inequality for the first inequality. We conclude, using orthogonality, that 
	 \begin{align*}
	 &\left\|\int_0^\infty\chi_{x_t - x\in \left[0,\hbar^{\frac23}\right]}\partial_{R}\phi_{n}(R;\xi)\xb(\xi)\rho_{n}(\xi)\,d\xi\right\|_{L^2_{R\,dR}}\\
	 &\lesssim \left(\sum_{\lambda}\left\|\int_0^\infty\chi_{x_t - x\in [0,\hbar^{\frac23}]}\partial_{R}\phi_{n}(R;\xi)\xb(\xi)\rho_{n}(\xi)\,d\xi\right\|_{L^2_{R\,dR}(R\simeq\lambda)}^2\right)^{\frac12}\\
	 &\lesssim  \left(\sum_{\lambda}\left\|\xi^{\frac12}\xb\right\|_{L^2(\xi\sim \lambda^{-2}\hbar^{-2})}^2\right)^{\frac12}\lesssim \left\|\xi^{\frac12}\xb(\xi)\right\|_{L^2_{d\xi}}.
	 \end{align*}
	 In the regime $x_t-x\in\left[\hbar^{\frac23}, \frac{x_t}{2}\right]$, we write $x_t-x\simeq \mu\hbar^{\frac23}$ for dyadic $\mu\geq 1$. Then fixing such $\mu$ the corresponding interval $I_R$ for $\xi$ is of length $\simeq   \frac{\mu\hbar^{\frac23}}{(R\hbar)^2}$, and we gain $\mu^{-\frac14}e^{-\frac32\mu^{\frac32}}$ from the decay of the Airy function. Proceeding as before one gets an analogous bound with (super)exponential decay in $\mu$ which can be summed up over dyadic $\mu\geq 1$. We omit the similar details.  
	 \\
	 {\it{Oscillatory regime $R\xi^{\frac12}\hbar\in\left[x_t, 2x_t\right]$.}} We start by observing that the bound in the transition regime $x-x_t\in \left[0,\hbar^{\frac23}\right]$ is handled exactly like the one for $x_t-x\in \left[0,\hbar^{\frac23}\right]$. There is, however, a difference in the regime $x_t-x\in \left[\hbar^{\frac23}, x_t\right]$, since we no longer get the rapid decay in terms of the re-scaled variable $\hbar^{-\frac23}(x-x_t)$, and so we have to argue more carefully, exploiting oscillation. Defining the kernel function
	 \begin{align*}
	 G(\xi, \eta): = \int_0^\infty \chi_{x\in \left[x_t+ \hbar^{\frac23}, 2x_t\right]}\partial_R\phi_{n}(R;\xi)\cdot\phi_{n}(R;\eta)\,R\,dR,
	 \end{align*}
	 we strive for an estimate like
	 \begin{align*}
	 	\left\|\int_{0}^{\infty}G(\xi,\eta)\xb(\eta)\,d\eta\right\|_{L^{2}_{d\xi}}\lesssim \left\|\eta^{\frac12}\xb(\eta)\right\|_{L^{2}_{d\eta}}.
	 \end{align*}
	 We start with the most difficult situation when the second factor $\phi_{n}(R;\eta)$ is also in the weakly oscillatory regime to the right of the turning point, i.e., the expression 
	 \begin{align*}
	 G_1(\xi, \eta): = \int_0^\infty \chi_{x\in \left[x_t+ \hbar^{\frac23}, 2x_t\right]} \chi_{\tilde{x}\in\left[x_t+ \hbar^{\frac23}, 2x_t\right]}\partial_R\phi_{n}(R;\xi)\cdot\phi_{n}(R;\eta)\,R\,dR,
	 \end{align*}
	 where we recall the notation $x = R\xi^{\frac12}\hbar, \tilde{x} = R\eta^{\frac12}\hbar$. Then by Proposition \ref{prop:DFT nlarge} we have
	 \begin{align*}
	 &\phi_{n}(R;\xi) = c\hbar^{\frac13}R^{-\frac12}\alpha^{-\frac12}q^{-\frac14}(\tau)\cdot\Re\left[a(\xi)\cdot\left(\Ai(\hbar^{-\frac23}\tau) -i\Bi(\hbar^{-\frac23}\tau)\right)\left(1+\hbar\overline{a_1(\tau;\alpha,\hbar)}\right)\right]\\
	 &\phi_{n}(R;\eta) = c\hbar^{\frac13}R^{-\frac12}\beta^{-\frac12}q^{-\frac14}(\tilde{\tau})\cdot\Re\left[a(\eta)\cdot\left(\Ai(\hbar^{-\frac23}\tilde{\tau}) -i\Bi(\hbar^{-\frac23}\tilde{\tau})\right)\left(1+\hbar\overline{a_1(\tilde{\tau};\beta,\hbar)}\right)\right].
	 \end{align*}
	 Using the standard asymptotics of Airy functions in the oscillatory regime, and only keeping track of the worst term with destructive resonance of the phases and $\partial_R$ falling on the oscillatory phase, we obtain the oscillatory integrals (omitting unnecessary constants)
	 \begin{equation}\label{eq:delicateoscillatorynearxt}
	 G_{1,\pm}(\xi,\eta): = \xi^{\frac12}\int_0^\infty \chi_{x\in \left[x_t+ \hbar^{\frac23}, 2x_t\right]} \chi_{\tilde{x}\in \left[x_t+ \hbar^{\frac23}, 2x_t\right]}(\xi\eta)^{-\frac14}a(\xi)\overline{a(\eta)}\tau^{\frac14}\frac{e^{\pm i\frac23\hbar^{-1}\left[\tau^{\frac32} - \tilde{\tau}^{\frac32}\right]}}{\tilde{\tau}^{\frac14}}\cdot(1+ b(\tau,\tilde{\tau},\alpha,\beta,\hbar))\,dR,
	 \end{equation}
	 where the function $b(\ldots)$ has symbol behavior with respect to all its variables. In order to obtain \eqref{eq:delicateoscillatorynearxt} we have used Lemma \ref{lem: monotonicity of root in alpha} and Lemma \ref{lem: Lemma 3.2 CDST}. We split the discussion according to $\ttau\gtrsim\tau$ and $\ttau\ll\tau$. By the presence of the factor $\frac{\tau^{\frac14}}{\ttau^{\frac14}}$, the second case is more difficult, which we discuss first.  Let us introduce $\tilde{\tau}: = \hbar^{\frac23}y^{\frac23}$, whence in our domain $y\gtrsim 1$, and recalling $\xi\simeq \eta$ in our region of integration, we find (by Lemma \ref{lem: Lemma 3.2 CDST})
	 \begin{align*}
	 \frac{\partial\tau}{\partial y} = \frac{\partial\tau}{\partial R}\cdot \frac{\partial R}{\partial\tilde{\tau}}\cdot \frac{\partial\tilde{\tau}}{\partial y}\sim y^{-\frac13}\hbar^{\frac23}, 
	 \end{align*}
	 and so (note that $\tau^{\frac12}\hbar^{-\frac13}y^{-\frac13}\gg1$) 
	 \begin{align*}
	 \chi_{\tau\gg\tilde{\tau}}\partial_y\left[e^{\pm i\frac23\hbar^{-1}\left[\tau^{\frac32} - \tilde{\tau}^{\frac32}\right]}\right]\simeq \tau^{\frac12}\hbar^{-\frac13}y^{-\frac13}\cdot e^{\pm i\frac23\hbar^{-1}\left[\tau^{\frac32} - \tilde{\tau}^{\frac32}\right]}.
	 \end{align*}
	 Thus we may replace $\chi_{\tau\gg\tilde{\tau}}e^{\pm i\frac23\hbar^{-1}\left[\tau^{\frac32} - \tilde{\tau}^{\frac32}\right]}$ by $\chi_{\tau\gg\tilde{\tau}}\frac{\tilde{\tau}^{\frac12}}{\tau^{\frac12}}\partial_y\left[e^{\pm i\frac23\hbar^{-1}\left[\tau^{\frac32} - \tilde{\tau}^{\frac32}\right]}\right]$, and performing integration by parts after passing to the variable $y$ (which does not cost anything since it produces additional factors $y^{-1}\lesssim 1$) and then reverting to the variable $R$, we have effectively replaced $\chi_{\tau\gg\tilde{\tau}}\frac{\tau^{\frac14}}{\tilde{\tau}^{\frac14}}$ by $\chi_{\tau\gg\tilde{\tau}}\left[\frac{\tau^{\frac14}}{\tilde{\tau}^{\frac14}}\right]^{-1}$. 
	 \\
	 Returning to \eqref{eq:delicateoscillatorynearxt}, we may henceforth assume that $\tilde{\tau}\gtrsim \tau$. Now introduce the new variable 
	 \[
	 \omega: = \frac23\hbar^{-1}\left[\tau^{\frac32} - \tilde{\tau}^{\frac32}\right]. 
	 \]
	 Observe that
	 \begin{align*}
	 \frac{\partial\omega}{\partial R} = \xi^{\frac12}\frac{\partial\tau}{\partial x}\tau^{\frac12} - \eta^{\frac12}\frac{\partial\tilde{\tau}}{\partial \tilde{x}}\tilde{\tau}^{\frac12}.
	 \end{align*}
	 By Lemma \ref{lem: Lemma 3.2 CDST}
	 \[
	 \frac{\partial\tau}{\partial x} = \Phi(x;\alpha, \tau) + \tau\cdot \frac{\partial_x \Phi(x;\alpha, \tau)}{\Phi(x;\alpha, \tau)}\simeq 1
	 \]
	 for $\tau$ sufficiently small. Therefore we have
	 \[
	 \tau^{\frac12} - \tilde{\tau}^{\frac12} = \frac{\tau - \tilde{\tau}}{\tau^{\frac12} + \tilde{\tau}^{\frac12}}\simeq \frac{\left[\xi^{\frac12} - \eta^{\frac12}\right]\cdot R\hbar}{\tau^{\frac12} + \tilde{\tau}^{\frac12}}\simeq  \frac{\left[\xi^{\frac12} - \eta^{\frac12}\right]}{\xi^{\frac12}(\tau^{\frac12} + \tilde{\tau}^{\frac12})}.
	 \]
	 On the other hand, we have
	 \begin{align*}
	 	\left|\frac{\partial}{\partial\xi^{\frac12}}\left(\xi^{\frac12}\frac{\partial\tau}{\partial x}\right)\right|=\left|\frac{\partial\tau}{\partial x}+\xi^{\frac12}\frac{\partial^{2}\tau}{\partial x^{2}}\cdot\hbar R\right|\lesssim 1.
	 \end{align*}
	 Therefore we have
	 \begin{align*}
	 \frac{\partial\omega}{\partial R} = \left(\xi^{\frac12}\frac{\partial\tau}{\partial x} - \eta^{\frac12}\frac{\partial\tilde{\tau}}{\partial \tilde{x}}\right)\cdot \tau^{\frac12} + \left(\tau^{\frac12} - \tilde{\tau}^{\frac12}\right)\cdot \eta^{\frac12}\frac{\partial\tilde{\tau}}{\partial \tilde{x}}\simeq \frac{\xi^{\frac12} - \eta^{\frac12}}{\tau^{\frac12} + \tilde{\tau}^{\frac12}},
	 \end{align*}
	 since 
	 \[
	 \left|\left(\xi^{\frac12}\frac{\partial\tau}{\partial x} - \eta^{\frac12}\frac{\partial\tilde{\tau}}{\partial \tilde{x}}\right)\cdot \tau^{\frac12}\right|\lesssim \left|\xi^{\frac12} - \eta^{\frac12}\right|\ll  \left|\frac{\xi^{\frac12} - \eta^{\frac12}}{\tau^{\frac12} + \tilde{\tau}^{\frac12}}\right|
	 \]
	 at last if $\tau + \tilde{\tau}\ll 1$, which we may as well assume(one can replace the bound $2x_t$ of the $x$-interval by any bound of the form $(1+c)x_t$ for our purposes). The integral under consideration can then be re-formulated as 
	 \begin{align*}
	 &\xi^{\frac12}\cdot(\xi\eta)^{-\frac14}\cdot \int_0^\infty  \chi_{x\in \left[x_t+ \hbar^{\frac23}, 2x_t\right]} \chi_{\tilde{x}\in \left[x_t+ \hbar^{\frac23}, 2x_t\right]}\chi_{\tau\lesssim \tilde{\tau}}\cdot e^{i\omega}\Psi(\omega;\alpha,\beta,\hbar)\frac{dR}{d\omega}\,d\omega\\
	 & = \frac{\xi^{\frac12}\cdot(\xi\eta)^{-\frac14}}{\xi^{\frac12} - \eta^{\frac12}}\cdot \int_0^\infty  \chi_{x\in \left[x_t+ \hbar^{\frac23}, 2x_t\right]} \chi_{\tilde{x}\in \left[x_t+ \hbar^{\frac23}, 2x_t\right]}\chi_{\tau\lesssim \tilde{\tau}}\cdot e^{i\omega}\tilde{\Psi}(\omega;\alpha,\beta,\hbar)\,d\omega, 
	 \end{align*}
	 where $\tilde{\Psi}(\omega;\alpha,\beta,\hbar)$ has symbol behavior with respect to all derivatives, and is uniformly bounded. In fact, observe that 
	 \begin{align*}
	 \frac{\partial\tilde{\Psi}(\omega;\alpha,\beta,\hbar)}{\partial\omega} = \frac{\partial\tilde{\Psi}}{\partial R}\cdot \frac{\partial R}{\partial \omega}
	 \end{align*}
	 and we have 
	 \begin{align*}
	 \omega \simeq \hbar^{-1}\left(\tau^{\frac12}-\tilde{\tau}^{\frac12}\right)\cdot (\tau+\tilde{\tau})\simeq \hbar^{-1}\frac{\xi^{\frac12}-\eta^{\frac12}}{\xi^{\frac12}}\cdot\left(\tau^{\frac12} + \tilde{\tau}^{\frac12}\right)\simeq R\cdot \left(\xi^{\frac12}-\eta^{\frac12}\right)\cdot\left(\tau^{\frac12} + \tilde{\tau}^{\frac12}\right),
	 \end{align*}
	 which implies
	 \begin{align*}
	 \left|\frac{\partial\tilde{\Psi}(\omega;\alpha,\beta,\hbar)}{\partial\omega} \right|\lesssim R^{-1}\cdot \frac{\tau^{\frac12}+\tilde{\tau}^{\frac12}}{\left|\xi^{\frac12} - \eta^{\frac12}\right|}\lesssim \left|\omega^{-1}\right|. 
	 \end{align*}
	 We have used here the easily verified symbol behavior of $\tilde{\Psi}$ with respect to $R$. The higher order derivatives are treated analogously. Finally we write
	 \begin{align*}
	 \int_0^\infty  \chi_{x\in \left[x_t+ \hbar^{\frac23}, 2x_t\right]} \chi_{\tilde{x}\in \left[x_t+ \hbar^{\frac23}, 2x_t\right]}\chi_{\tau\lesssim \tilde{\tau}}\cdot e^{i\omega}\tilde{\Psi}(\omega;\alpha,\beta,\hbar)\,d\omega =:H(\xi, \eta),
	 \end{align*}
	 where we have $\left|H(\xi, \eta)\right|\lesssim 1$ and symbol behavior with respect to $\xi, \eta$. But (keeping in mind that $\xi\simeq\eta$ on the support of $H$) it is easy to then verify the bound 
	 \begin{align*}
	 \left\|\int_0^\infty\frac{\xi^{\frac12} \chi_{\xi\simeq\eta}H(\xi,\eta)\cdot(\xi\eta)^{-\frac14}}{\xi^{\frac12} - \eta^{\frac12}}\cdot x(\eta)\,d\eta\right\|_{L^2_{d\xi}}\lesssim \left\|\eta^{\frac12}x\right\|_{L^2_{d\eta}}, 
	 \end{align*}
	 as desired. 
	 \\ 
	 {\it{Oscillatory regime}}: $R\xi^{\frac12}\hbar>2x_t$. Here, in analogy to the analysis of the transference operator, we have to understand the distribution valued kernel 
	 \begin{align*}
	 \lim_{A\rightarrow\infty}\int_0^A \partial_R\phi_{n}(R;\xi)\phi_{n}(R;\eta,\hbar)R\,dR =:F(\xi,\eta),
	 \end{align*}
	 and show that it induces an operator obeying the asserted bound. For this we again use the oscillatory expansion 
	 \begin{align*}
	 \phi_{n}(R;\xi) = 2\pi^{-\frac12}\left(R\xi^{\frac12}\right)^{-\frac12}\left(\tau q(x;\alpha,\hbar)\right)^{-\frac14}\Re\left(a(\xi)e^{i\frac{2}{3\hbar}\tau^{\frac32}}(1+\hbar a_1(\tau,\alpha,\hbar))\right),
	 \end{align*}
	 where the description of the variable $\tau$ in terms of $x,\alpha,\hbar$ is furnished by Lemma \ref{lem: Lemma 3.4 CDST}. We only deal with the most delicate case when $\xi\simeq \eta$ and the resonant case where the phases cancel. The combined phases arising will then be of the form, by Lemma \ref{lem: Lemma 3.4 CDST}, 
	 \[
	 e^{\pm i\left[R\left(\xi^{\frac12} - \eta^{\frac12}\right)-\hbar^{-1}y(\alpha;\hbar) + \hbar^{-1}y(\beta;\hbar) - \hbar^{-1}\rho(x;\alpha,\hbar) + \hbar^{-1}\rho(\tilde{x};\beta,\hbar)\right]}.
	 \]
When the operator $\partial_R$ hits the non-phase factor $(1+\hbar a_1(\tau,\alpha,\hbar))$, using Lemma \ref{lem:a1fine} we obtain an extra factor 
\[
\lesssim \hbar x^{-\frac13}R^{-1}
\]
and hence the corresponding term can be integrated over the region $R\xi^{\frac12}\hbar\gtrsim 1$ to give a contribution to $\chi_{\xi\simeq \eta}F(\xi,\eta)$ of size $\lesssim \xi^{-\frac12}$. By Minkowski inequality and the fact $\xi\simeq \eta$, such a bound for $\chi_{\xi\simeq \eta}F(\xi,\eta)$ gives
\[
\left\|\int_0^\infty\chi_{\xi\simeq \eta}F(\xi,\eta)x(\eta)\,d\eta\right\|_{L^2_{d\xi}}\lesssim \left\|\eta^{\frac12}x(\eta)\right\|_{L^2_{d\eta}}. 
\]
We henceforth only consider the contribution arising when $\partial_R$ hits the phase associated with $\xi$. This then results in the term (coming purely from the phase, and omitting constant coefficients)
\begin{equation}\label{eq:tough1}\begin{split}
&\Re\left(a(\xi)\overline{a(\eta)}\right)\xi^{\frac12}\cdot\left[1-\rho_x(x;\alpha,\hbar)\right]\\
&\cdot \sin\left[R(\xi^{\frac12} - \eta^{\frac12})-\hbar^{-1}y(\alpha;\hbar) + \hbar^{-1}y(\beta;\hbar) - \hbar^{-1}\rho(x;\alpha,\hbar) + \hbar^{-1}\rho(\tilde{x};\beta,\hbar)\right]\\
& + \xi^{\frac12}\cdot\left[1-\rho_x(x;\alpha,\hbar)\right]\Im\left(a(\eta)\overline{a(\xi)}\right)\\
&\cdot \cos\left[R(\xi^{\frac12} - \eta^{\frac12})-\hbar^{-1}y(\alpha;\hbar) + \hbar^{-1}y(\beta;\hbar) - \hbar^{-1}\rho(x;\alpha,\hbar) + \hbar^{-1}\rho(\tilde{x};\beta,\hbar)\right].
\end{split}
\end{equation}
Using the bounds on $a(\cdot)$ and its derivatives, we have
\begin{align*}
	\chi_{\xi\simeq\eta}\left(a(\eta)\overline{a(\xi)}-a(\xi)\overline{a(\eta)}\right)=&\chi_{\xi\simeq\eta}\left(\left(a(\eta)-a(\xi)\right)\overline{a(\xi)}+a(\xi)\left(\overline{a(\xi)}-\overline{a(\eta)}\right)\right)=\frac{(\xi^{\frac12} - \eta^{\frac12})}{\xi^{\frac12}}\cdot b(\xi,\eta),
\end{align*}
where $b$ itself obeys symbol behavior. 
Writing the phase as $\left(\xi^{\frac12} - \eta^{\frac12}\right)\cdot \tilde{R}$, $\tilde{R} = \tilde{R}(R;\alpha,\beta,\hbar)$, we have in our regime $\hbar\xi^{\frac12}R\simeq \hbar\eta^{\frac12}R>(1+c)x_t$, $c>0$ that 
\[
\frac{d\tilde{R}}{dR}\simeq 1
\]
uniformly in all parameters. The contribution of the second term in \eqref{eq:tough1} can then be handled by integration by parts: omitting constant terms and proceeding schematically, we reduce this to 
\begin{align*}
\frac{b(\xi,\eta)}{\xi^{\frac12}}\int_0^\infty \left[1-\rho_x(x;\alpha,\hbar)\right]\cdot (1+\hbar a_1(\tau,\alpha,\hbar))(1+\hbar a_1(\tilde{\tau},\beta,\hbar))\cdot \partial_{\tilde{R}}\left(\sin\left(\xi^{\frac12}-\eta^{\frac12}\right)\right)\cdot \left(\frac{d\tilde{R}}{dR}\right)^{-1}\,d\tilde{R}. 
\end{align*}
Carrying out the integration by parts again leads to a contribution to $\chi_{\xi\sim\eta}F(\xi,\eta)$ of size $\lesssim \xi^{-\frac12}$. 
\\
We can thus reduce things to the first, principal term in \eqref{eq:tough1}, which is the following contribution to $F(\xi, \eta)$:
\begin{equation}\label{eq:tough1main}
\begin{split}
&\frac{\Re\left(a(\xi)\overline{a(\eta)}\right)}{(\xi\eta)^{\frac14}}\xi^{\frac12}\\
&\cdot\lim_{A\rightarrow\infty}\int_0^A\chi_{R\xi^{\frac12}\hbar>2x_t}\left[1-\rho_x(x;\alpha,\hbar)\right]\cdot\frac{(1+\hbar a_1(\tau,\alpha,\hbar))}{(q(x;\alpha,\hbar)\tau)^{\frac14}}\frac{(1+\hbar a_1(\tilde{\tau},\beta,\hbar))}{(q(\tilde{x};\beta,\hbar)\tilde{\tau})^{\frac14}}\cdot \sin\left[\tilde{R}\left(\xi^{\frac12}-\eta^{\frac12}\right)\right]\cdot \left(\frac{d\tilde{R}}{dR}\right)^{-1}\,d\tilde{R}
\end{split}
\end{equation}
Writing
\begin{align*}
\chi_{R\xi^{\frac12}\hbar>2x_t}[1-\rho_x(x;\alpha,\hbar)]\cdot\frac{(1+\hbar \tilde{a}_1(\tau,\alpha,\hbar))}{(q(x;\alpha,\hbar)\tau)^{\frac14}}\frac{(1+\hbar \tilde{a}_1(\tilde{\tau},\beta,\hbar))}{(q(\tilde{x};\beta,\hbar)\tilde{\tau})^{\frac14}}\cdot \left(\frac{d\tilde{R}}{dR}\right)^{-1} = \Phi(\tilde{R};\alpha,\beta,\hbar)
\end{align*}
and
\begin{align*}
\Phi(\tilde{R};\alpha,\beta,\hbar) = 1 - \chi_{R\xi^{\frac12}\hbar<2x_t} +  \Psi(\tilde{R};\alpha,\beta,\hbar),
\end{align*}
with $\left| \Psi(\tilde{R};\alpha,\beta,\hbar)\right|\lesssim x^{-1}$ and symbol behavior with respect to its arguments. Then we have 
\begin{align*}
\chi_{\xi\simeq \eta}\frac{a(\xi)\overline{a(\eta)}}{(\xi\eta)^{\frac14}}\xi^{\frac12}\cdot\lim_{A\rightarrow\infty}\int_0^A\sin\left[\tilde{R}(\xi^{\frac12}-\eta^{\frac12})\right]\,d\tilde{R} = c\chi_{\xi\simeq \eta}\frac{a(\xi)\overline{a(\eta)}}{(\xi\eta)^{\frac14}}\xi^{\frac12}\cdot \left(\frac{1}{\xi^{\frac12} - \eta^{\frac12}}\right)_{P.V.}
\end{align*}
and acts like a Hilbert transform like operator which is easily seen to satisfy the desired bound. Furthermore we can write 
\begin{align*}
\lim_{A\rightarrow\infty}\int_0^A\chi_{R\xi^{\frac12}\hbar<2x_t}\sin\left[\tilde{R}(\xi^{\frac12}-\eta^{\frac12})\right]\,d\tilde{R} = c\left(\left(\frac{1}{\cdot}\right)_{P.V.}* m\right)(\xi^{\frac12} - \eta^{\frac12})
\end{align*}
for a function $m$ of bounded $L^1$-norm (independent of $\hbar$), and hence this contribution can be bounded analogously to the preceding one. 
\\
Finally, for the operator 
\begin{align*}
&\chi_{\xi\simeq \eta}\frac{a(\xi)\overline{a(\eta)}}{(\xi\eta)^{\frac14}}\xi^{\frac12}\cdot\lim_{A\rightarrow\infty}\int_0^A  \Psi(\tilde{R};\alpha,\beta,\hbar)\sin\left[\tilde{R}(\xi^{\frac12}-\eta^{\frac12})\right]\,d\tilde{R}\\
& = \chi_{\xi\simeq \eta}\frac{a(\xi)\overline{a(\eta)}}{(\xi\eta)^{\frac14}}\xi^{\frac12}\cdot\sum_{\lambda\gtrsim 2x_t}\lim_{A\rightarrow\infty}\int_0^A \chi_{x\simeq\lambda}\Psi(\tilde{R};\alpha,\beta,\hbar)\sin\left[\tilde{R}(\xi^{\frac12}-\eta^{\frac12})\right]\,d\tilde{R},\\
\end{align*}
where the sum is over dyadic $\lambda$, we use that 
\begin{align*}
\lim_{A\rightarrow\infty}\int_0^A \chi_{x\simeq\lambda}\Psi(\tilde{R};\alpha,\beta,\hbar)\sin\left[\tilde{R}(\xi^{\frac12}-\eta^{\frac12})\right]\,d\tilde{R} = c\left(\left(\frac{1}{\cdot}\right)_{P.V.}* m_{\lambda}\right)(\xi^{\frac12} - \eta^{\frac12})
\end{align*}
where $\left\|m_{\lambda}\right\|_{L^1}\lesssim \lambda^{-1}$, and we can bound this contribution as before and sum over all dyadic $\lambda\gtrsim x_t$. 
\end{proof}
In a similar vein we have the following: 
\begin{proposition}\label{prop:singularmultiplier} 
	Assuming, with $\hbar = \frac{1}{n+1}$, $|n|\geq 2$,
	\[
	f(R) = \int_0^\infty \phi_{n}(R;\xi)\xb(\xi)\rho_{n}(\xi)\,d\xi,
	\]
	we have for $k+\ell\ll n$ the bound 
	\begin{align*}
	\left\|\left(\frac{n}{R}\right)^k \partial_R^{\ell}f(R)\right\|_{L^2_{R\,dR}}\lesssim_{k,\ell} \left\|\xi^{\frac{k+\ell}{2}}\xb(\xi)\right\|_{L^2_{d\xi}}, 
	\end{align*}
	the implied constant being uniform in $\hbar$. 
\end{proposition}
\begin{proof}
	This is similar to the preceding proposition. One exploits that for $x = R\xi^{\frac12}\hbar\gtrsim 1$, we have 
	\[
	\left(\frac{n}{R}\right)^k\sim \left(\frac{1}{R\hbar}\right)^k\lesssim \xi^{\frac{k}{2}},
	\]
	while if $x = R\xi^{\frac12}\hbar<\frac{x_t}{2}$, we have 
	\[
	\left| \left(\frac{n}{R}\right)^k\phi_{n}(R;\xi)\right|\lesssim \hbar^{\frac12}\left(R\xi^{\frac12}\hbar\right)^{n-1-k}\cdot \xi^{\frac{k}{2}}.
	\]
	We omit the remaining similar details.
\end{proof}
The following lemma gives the $L^{\infty}$-bound in terms of the $\Sh_{0}$-norm:
\begin{lemma}\label{lem:derLinfty} Let
	\[
	f(R) = \int_0^\infty \phi_{n}(R;\xi)\xb(\xi)\rho_{n}(\xi)\,d\xi. 
	\]
	Then for $k\geq 1$ we can estimate 
	\begin{align*}
	\left\|\partial_R^k f\right\|_{L^\infty_{dR}}\lesssim \hbar^{-\frac32-\delta}\left\|(\xi\hbar^2)^{1-\frac{\delta}{2}}\langle\xi\hbar^2\rangle^{\delta}\xi^{\frac{k-1}{2}}\xb(\xi)\right\|_{L^2_{d\xi}}.
	\end{align*}
	In particular we have (with $\hbar = \frac{1}{n+1}$)
	\begin{align*}
	\left\|\partial_Rf\right\|_{L^\infty_{dR}}\lesssim \hbar^{-\frac32-\delta}\left\|\xb\right\|_{\Sh_{0}}.
	\end{align*}
	Finally, we also have, for $k+\ell\geq 1$,
	\begin{align*}
	\left\|\left(\frac{n}{R}\right)^{\ell}\partial_R^k f\right\|_{L^\infty_{dR}}\lesssim \hbar^{-\frac32-\delta}\left\|(\xi\hbar^2)^{1-\frac{\delta}{2}}\langle\xi\hbar^2\rangle^{\delta}\xi^{\frac{k+\ell-1}{2}}
	\xb(\xi)\right\|_{L^2_{d\xi}}.
	\end{align*}
	In particular, we have the estimate 
	\begin{align*}
	\left\|\frac{f}{R}\right\|_{L^\infty_{dR}}\lesssim \hbar^{-\frac12-\delta}\left\|(\xi\hbar^2)^{1-\frac{\delta}{2}}\langle\xi\hbar^2\rangle^{\delta}
	\xb(\xi)\right\|_{L^2_{d\xi}}.
	\end{align*}
\end{lemma}
\begin{proof}
	We treat the case $k=1$. The higher order derivatives can be handled similarly. From the proof of Proposition \ref{prop:derivative1}, we have
	\begin{align*}
		\left|\partial_{R}\phi_{n}(R,\xi)\right|\lesssim \hbar^{\frac12}\xi^{\frac12}.
	\end{align*}
	Then we have
	\begin{align*}
	\left\|\partial_Rf\right\|_{L^\infty_{dR}}&\lesssim \hbar^{\frac12}\int_0^\infty \xi^{\frac12}\left|\xb(\xi)\right|\,d\xi\\
	&\lesssim  \hbar^{\frac12}\left(\int_0^1 \left|\xi^{1-\frac{\delta}{2}}\xb(\xi)\right|^2\,d\xi\right)^{\frac12}\cdot \left(\int_0^1\xi^{-1+\delta}\,d\xi\right)^{\frac12}\\
	&+  \hbar^{\frac12}\left(\int_1^\infty \left|\xi^{1+\frac{\delta}{2}}\xb(\xi)\right|^2\,d\xi\right)^{\frac12}\cdot \left(\int_1^\infty\xi^{-1-\delta}\,d\xi\right)^{\frac12}\\
	&\lesssim \hbar^{-\frac32-\delta}\left\|(\xi\hbar^2)^{1-\frac{\delta}{2}}\langle\xi\hbar^2\rangle^{\delta}\xb(\xi)\right\|_{L^2_{d\xi}}\lesssim \hbar^{-\frac32-\delta}\left\|\xb\right\|_{\Sh_0}.
	\end{align*}
\end{proof}

We now state the first main proposition on multilinear estimates. The assumption here is that we have the following Fourier representations 
\begin{equation}\label{eq:Fourier1}
\phi_j(R) = \int_0^\infty \phi_{n_{j}}(R;\xi)\xb_j(\xi)\rho_{n_{j}}(\xi)\,d\xi,\quad j = 1,\,2,\quad \xb_j(\xi)\in \Shj_0.
\end{equation}

\begin{proposition}\label{prop:bilin1} 
	Assume that $\hbar_1\lesssim \hbar_2\ll 1$, and either $\hbar_3\simeq \hbar_1$ or $\hbar_3\gg \hbar_1, |\hbar_1|\simeq|\hbar_2|$. Then we have the bound 
	\begin{align*}
	\left\|\left\langle \phi_{n_{3}}(R;\xi),\,\partial_R\phi_1\cdot\partial_R\phi_2\right\rangle_{L^2_{R\,dR}}\right\|_{\Shth_1}\lesssim \hbar_2^{-2}\prod_{j=1,2}\left\|\xb_j\right\|_{\Shj_0}.
	\end{align*}
	We similarly have the bound
	\begin{align*}
	\left\|\left\langle \phi_{n_{3}}(R;\xi),\,\chi_{R\lesssim\tau}\phi_1\cdot \phi_2\right\rangle_{L^2_{R\,dR}}\right\|_{\Shth_0}\lesssim \tau\cdot \hbar_2^{-2}\prod_{j=1,2}\left\|\xb_j\right\|_{\Shj_0}.
	\end{align*}
\end{proposition}
\begin{proof} The second estimate is proved similarly to the first, which we do here. Note that the loss of $\tau$ in the second estimate comes from application of the $L^\infty$-bound fo $\frac{f}{R}$ in the preceding lemma.
	We treat the case $\hbar_1\simeq \hbar_3$ in detail. The second case is similar and simpler, in the sense that this case happens when $\hbar_{1}$ and $\hbar_{2}$ have opposite signs and $n_{1}, n_{2}$ almost cancel each other. Therefore there is extra decay in $|\hbar_{1}|$ and $|\hbar_{2}|$ in the product $\prod_{j=1,2}\left\|\xb_j\right\|_{\Shj_0}$. We distinguish between the case when the output frequency $\xi$ is less than the maximal frequency of the factors, and the case when it is larger. Call the output frequency $\xi$ and those of the factors $\phi_j$ $\xi_j$, $j = 1,2$. 
	\\
	
	{\it{(1): $\xi<\max\{\xi_1, \xi_2\}$.}} Due to the asymmetry of the situation, we split this further into two sub-cases.
	\\
	
	{\it{(1.a): $\xi<\xi_1$.}} We write this contribution as 
	\[
	\sum_{\mu<\lambda}\chi_{\xi\simeq\mu}\left\langle \phi_{n_{3}}(R;\xi),\,\partial_R\phi_{1,\lambda}\cdot\partial_R\phi_2\right\rangle_{L^2_{R\,dR}}
	\]
	where $\mu,\lambda$ range over dyadic numbers, and we use the notation 
	\[
	\phi_{1,\lambda} = \int_0^\infty \chi_{\xi_1\simeq\lambda}\phi_{n_{1}}(R;\xi_1)\xb_j(\xi_1)\rho_{n_{1}}(\xi_1)\,d\xi_1
	\]
	By Plancherel's theorem for the distorted Fourier transform we have (recall $\hbar_3\simeq \hbar_1$)
	\begin{align*}
	&\left\|\chi_{\xi\simeq\mu}\left\langle \phi_{n_{3}}(R;\xi),\,\partial_R\phi_{1,\lambda}\cdot\partial_R\phi_2\right\rangle_{L^2_{R\,dR}}\right\|_{\Shth_1}\\
	&\lesssim \hbar_1(\mu\hbar_1^2)^{\frac12-\frac{\delta}{2}}\langle \mu\hbar_1^2\rangle^{\delta + \frac32}\left\|\left\langle \phi_{n_{3}}(R;\xi),\,\partial_R\phi_{1,\lambda}\cdot\partial_R\phi_2\right\rangle_{L^2_{R\,dR}}\right\|_{L^2_{d\xi}}\\
	&\simeq \hbar_1(\mu\hbar_1^2)^{\frac12-\frac{\delta}{2}}\langle \mu\hbar_1^2\rangle^{\delta + \frac32}\left\|\left\langle \phi_{n_{3}}(R;\xi),\,\partial_R\phi_{1,\lambda}\cdot\partial_R\phi_2\right\rangle_{L^2_{R\,dR}}\right\|_{L^2_{\rho_{n_{3}}(\xi)d\xi}}\\
	&\lesssim \hbar_1(\mu\hbar_1^2)^{\frac12-\frac{\delta}{2}}\langle \mu\hbar_1^2\rangle^{\delta + \frac32}\left\|\partial_R\phi_{1,\lambda}\cdot\partial_R\phi_2\right\|_{L^2_{R\,dR}}.
	\end{align*}
	The last $L^{2}_{R\,dR}$-norm of the product can be estimated using Proposition \ref{prop:singularmultiplier} and Lemma \ref{lem:derLinfty}:
	\begin{align*}
	\begin{split}
	\left\|\partial_R\phi_{1,\lambda}\cdot\partial_R\phi_2\right\|_{L^2_{R\,dR}}&\lesssim \left\|\partial_R\phi_{1,\lambda}\right\|_{L^2_{R\,dR}}\cdot\left\|\partial_R\phi_2\right\|_{L^\infty_{R\,dR}}\\
	&\lesssim \lambda^{\frac12}\left\|\xb_1\right\|_{L^2_{d\xi_1}(\xi_1\simeq\lambda)}\cdot \hbar_2^{-\frac32-\delta}\cdot \left\|\xb_2\right\|_{\Sht_0}.
	\end{split}
	\end{align*}
	Combining the two preceding bounds, we infer that 
	\begin{align*}
	&\left\|\chi_{\xi\simeq\mu}\left\langle \phi(R;\xi,\hbar_3),\,\partial_R\phi_{1,\lambda}\cdot\partial_R\phi_2\right\rangle_{L^2_{R\,dR}}\right\|_{\Shth_1}\\
	&\lesssim  \hbar_2^{-\frac32-\delta}\lambda^{\frac12} \hbar_1(\mu\hbar_1^2)^{\frac12-\frac{\delta}{2}}\langle \mu\hbar_1^2\rangle^{\delta + \frac32}\left\|\xb_1\right\|_{L^2_{d\xi_1}(\xi_1\simeq\lambda)}\cdot  \left\|\xb_2\right\|_{\Sht_0}\\
	&\lesssim \hbar_2^{-\frac32-\delta}\left(\frac{\mu}{\lambda}\right)^{\frac12-\frac{\delta}{2}}\cdot \left\|\xb_1\right\|_{\Sho_0(\xi_1\simeq\lambda)}\cdot \left\|\xb_2\right\|_{\Sht_0}.
	\end{align*}
	Finally exploiting orthogonality as well as the Cauchy-Schwarz inequality, we infer that 
	\begin{align*}
	&\left\|\sum_{\mu<\lambda}\chi_{\xi\simeq\mu}\left\langle \phi_{n_{3}}(R;\xi),\,\partial_R\phi_{1,\lambda}\cdot\partial_R\phi_2\right\rangle_{L^2_{R\,dR}}\right\|_{\Shth_1}^2\\
	&\lesssim \sum_\mu\left\|\sum_{\lambda>\mu}\chi_{\xi\simeq\mu}\left\langle \phi_{n_{3}}(R;\xi),\,\partial_R\phi_{1,\lambda}\cdot\partial_R\phi_2\right\rangle_{L^2_{R\,dR}}\right\|_{\Shth_1}^2\\
	&\lesssim \hbar_2^{-3-2\delta}\sum_{\mu<\lambda}\left(\frac{\mu}{\lambda}\right)^{1-\delta}\cdot \left\|\xb_1\right\|_{\Sho_0(\xi_1\simeq\lambda)}^2\cdot \left\|\xb_2\right\|_{\Sht_0}^2\\
	&\lesssim \hbar_2^{-3-2\delta}\left\|\xb_1\right\|_{\Sho_0}^2\cdot \left\|\xb_2\right\|_{\Sht_0}^2,
	\end{align*}
	which is desired. 
	\\
	
	{\it{(1.b): $\xi_1\leq \xi<\xi_2$.}} Write this term as 
	\[
	\sum_{\mu<\lambda}\chi_{\xi\simeq\mu}\langle \phi_{n_{3}}(R;\xi),\,\partial_R\phi_{1,\leq \mu}\cdot\partial_R\phi_{2,\lambda}\rangle_{L^2_{R\,dR}}.
	\]
	This time we use Lemma~\ref{lem:derLinfty} to bound the factor $\partial_{R}\phi_{2,\lambda}$: as in the preceding case, it suffices to bound 
	\begin{align*}
	&\hbar_1(\mu\hbar_1^2)^{\frac12-\frac{\delta}{2}}\langle \mu\hbar_1^2\rangle^{\delta + \frac32}\cdot\left\|\partial_R\phi_{1,\leq \mu}\cdot\partial_R\phi_{2,\lambda}\right\|_{L^2_{R\,dR}}\\
	&\lesssim \hbar_1(\mu\hbar_1^2)^{\frac12-\frac{\delta}{2}}\langle \mu\hbar_1^2\rangle^{\delta + \frac32}\left\|\partial_R\phi_{1,\leq \mu}\right\|_{L^\infty_{R\,dR}}\cdot \left\|\partial_R\phi_{2,\lambda}\right\|_{L^2_{R\,dR}}\\
	&\lesssim \hbar_2^{-2+\delta}\left(\frac{\mu}{\lambda}\right)^{\frac12-\frac{\delta}{2}}\cdot\hbar_1^{\frac12-2\delta}\cdot (\lambda\hbar_2^2)^{1-\frac{\delta}{2}}\langle \lambda\hbar_2^2\rangle^{\delta + \frac32}\cdot \left\|\xb_1\right\|_{\Sho_0(\xi_1\leq \mu)}\cdot \left\|\xb_2\right\|_{L^2_{d\xi_{2}}(\xi_2\simeq\lambda)}\\
	&\simeq  \hbar_2^{-2+\delta}\left(\frac{\mu}{\lambda}\right)^{\frac12-\frac{\delta}{2}}\cdot\hbar_1^{\frac12-2\delta}\cdot \left\|\xb_1\right\|_{\Sho_0(\xi_1\leq \mu)}\cdot \left\|\xb_2\right\|_{\Sht_0(\xi_2\simeq\lambda)}.\\
	\end{align*}
	From this point the estimate can be completed just as in the preceding case. 
	\\
	
	We next deal with the case where the output frequency $\xi$ dominates the input frequencies:
	\\
	
	{\it{(2): $\xi\geq \max\{\xi_1, \xi_2\}$.}} Here we perform integration by parts in order to gain in terms of $\frac{ \max\{\xi_1, \xi_2\}}{\xi}$. We can do this rather carelessly in the present situation since we are assuming that all angular parameters $\hbar_j$, $j = 1,2$, are very small, which ensures sufficient vanishing at the origin $R = 0$.  We have to exploit that $ \phi_{n_{3}}(R;\xi)$ is a generalized eigenfunction of $H_{n_{3}}^{\pm}$, $\hbar_3 = \frac{1}{n_{3}+1}$ which we denote as $H_{n_{3}}$ for simplicity. Then we have 
	\[
	-H_{n_{3}} \phi_{n_{3}}(R;\xi) = \xi \phi_{n_{3}}(R;\xi)
	\]
	We apply $H_{n_{3}}$ three times and integrate by parts, resulting schematically in 
	\begin{align*}
	&\xi^3\left\langle \phi_{n_{3}}(R;\xi),\,\partial_R\phi_{1,<\xi}\cdot\partial_R\phi_{2,<\xi}\right\rangle_{L^2_{R\,dR}}\\& = 
	\sum_{i+j+k = 6}\left\langle \phi_{n_{3}}(R;\xi),\,\left(\frac{n_{3}}{R}\right)^i\partial_R^{1+j}\phi_{1,<\xi}\cdot\partial_R^{1+k}\phi_{2,<\xi}\right\rangle_{L^2_{R\,dR}},
	\end{align*}
	where we have only included the most singular terms, the remaining ones being much simpler to treat. 
	Importantly note that we combine the factors $\frac{n_{3}}{R}$ with $\phi_{1,<\xi}$, in order not to lose inverse powers of $\hbar_1\simeq \hbar_3$. In order to proceed we dyadically localize the output frequency $\xi\simeq \lambda$ 
	and the inner factors to dyadic frequencies $\mu_1,\mu_2$ respectively. Invoking orthogonality, we infer 
	\begin{align*}
	&\left\|\left\langle \phi_{n_{3}}(R;\xi),\,\partial_R\phi_{1,<\xi}\cdot\partial_R\phi_{2,<\xi}\right\rangle_{L^2_{R\,dR}}\right\|_{\Shth_1}^2\\
	& \simeq \sum_{\lambda}\left\|\left\langle \phi_{n_{3}}(R;\xi),\,\partial_R\phi_{1,<\xi}\cdot\partial_R\phi_{2,<\xi}\right\rangle_{L^2_{R\,dR}}\right\|_{\Shth_1(\xi\simeq\lambda)}^2\\
	&\lesssim \sum_{i+j+k = 6}\sum_{\lambda}\lambda^{-6}\left\|\left\langle \phi_{n_{3}}(R;\xi),\,\left(\frac{n_{3}}{R}\right)^i\partial_R^{1+j}\phi_{1,<\xi}\cdot\partial_R^{1+k}\phi_{2,<\xi}\right\rangle_{L^2_{R\,dR}}\right\|_{\Shth_1(\xi\simeq\lambda)}^2\\
	&\lesssim  \sum_{i+j+k = 6}\sum_{\lambda}\lambda^{-6}\left\|\sum_{\mu_1,\mu_2<\lambda}\left\langle \phi_{n_{3}}(R;\xi),\,\left(\frac{n_{3}}{R}\right)^i\partial_R^{1+j}\phi_{1,\mu_1}\cdot\partial_R^{1+k}\phi_{2,\mu_2}\right\rangle_{L^2_{R\,dR}}\right\|_{\Shth_1(\xi\simeq\lambda)}^2.
	\end{align*}
	Using Plancherel's theorem for the distorted Fourier transform we obtain, for $\mu_{2}\leq \mu_{1}$,
	\begin{align}\label{mu1 mu2 distinct}
	\begin{split}
	&\lambda^{-3}\left\|\sum_{\mu_1,\mu_2<\lambda}\left\langle \phi_{n_{3}}(R;\xi),\,\left(\frac{n_{3}}{R}\right)^i\partial_R^{1+j}\phi_{1,\mu_1}\cdot\partial_R^{1+k}\phi_{2,\mu_2}\right\rangle_{L^2_{R\,dR}}\right\|_{\Shth_1}\\
	&\lesssim \frac{(\hbar_3^2\lambda)^{1-\frac{\delta}{2}}\langle\hbar_3^2\lambda\rangle^{\delta+\frac32}}{\lambda^{\frac{7}{2}}}\sum_{\mu_{1,2}<\lambda}\left\|\left(\frac{n_{3}}{R}\right)^i\partial_R^{1+j}\phi_{1,\mu_1}\right\|_{L^2_{R\,dR}}\cdot\left\|\partial_R^{1+k}\phi_{2,\mu_2}\right\|_{L^\infty_{R\,dR}}.
	\end{split}
	\end{align}
	By Proposition~\ref{prop:singularmultiplier} as well as Lemma~\ref{lem:derLinfty}, we bound the preceding by 
	\begin{align*}
	\begin{split}
	&\hbar_{2}^{-\frac32-\delta}\frac{(\hbar_3^2\lambda)^{1-\frac{\delta}{2}}\langle\hbar_3^2\lambda\rangle^{\delta+\frac32}}{\lambda^\frac{7}{2}}\sum_{\mu_{1,2}<\lambda}\mu_1^{\frac{1+i+j}{2}}\left\|\xb_1\right\|_{L^2_{d\xi}(\xi_{1}\simeq\mu_1)}\cdot\mu_2^{\frac{k}{2}}\left\|\xb_{2}\right\|_{\Sht_{0}(\xi_{2}\simeq\mu_{2})}
	\\
	&\lesssim \hbar_2^{-\frac32-\delta}\sum_{\mu_{1,2}<\lambda}\frac{(\hbar_3^2\lambda)^{1-\frac{\delta}{2}}\langle\hbar_3^2\lambda\rangle^{\delta+\frac32}}{(\hbar_1^2\mu_1)^{1-\frac{\delta}{2}}\langle\hbar_1^2\mu_1\rangle^{\delta+\frac32}}\cdot \left(\frac{\mu_1}{\lambda}\right)^{\frac{1+i+j}{2}}\left\|\xb_1\right\|_{\Sho_0(\xi_{1}\simeq\mu_1)}\cdot \left(\frac{\mu_2}{\lambda}\right)^{\frac{k}{2}}\left\|\xb_2\right\|_{\Sht_0(\xi\simeq\mu_2)}.
	\end{split}
	\end{align*}
	Finally, observe that 
	\begin{align*}
	\frac{(\hbar_3^2\lambda)^{1-\frac{\delta}{2}}\langle\hbar_3^2\lambda\rangle^{\delta+\frac32}}{(\hbar_1^2\mu_1)^{1-\frac{\delta}{2}}\langle\hbar_1^2\mu_1\rangle^{\delta+\frac32}}\cdot \left(\frac{\mu_1}{\lambda}\right)^{\frac{1+i+j}{2}}\cdot \left(\frac{\mu_2}{\lambda}\right)^{\frac{k}{2}}\lesssim \prod_{j=1,2}\left(\frac{\mu_j}{\lambda}\right)^{\frac14}.
	\end{align*}
	We conclude by using the Cauchy-Schwarz inequality to infer that 
	\begin{align*}
	&\hbar_2^{-\frac32-\delta}\left(\sum_{\lambda}\left[\sum_{\mu_{1,2}<\lambda}\prod_{j=1,2}\left(\frac{\mu_j}{\lambda}\right)^{\frac14}\left\|\xb_1\right\|_{\Sho_0(\xi\simeq\mu_1)}\left\|\xb_2\right\|_{\Sht_0(\xi\simeq\mu_2)}\right]^2\right)^{\frac12}\\
	&\lesssim \hbar_2^{-\frac32-\delta}\left\|\xb_1\right\|_{\Sho_0}\cdot\left\|\xb_2\right\|_{\Sht_0},
	\end{align*}
	as desired. If $\mu_{1}\leq \mu_{2}$, then we bound \eqref{mu1 mu2 distinct} as
	\begin{align*}
		\frac{(\hbar_3^2\lambda)^{1-\frac{\delta}{2}}\langle\hbar_3^2\lambda\rangle^{\delta+\frac32}}{\lambda^{\frac{7}{2}}}\sum_{\mu_{1,2}<\lambda}\left\|\left(\frac{n_{3}}{R}\right)^i\partial_R^{1+j}\phi_{1,\mu_1}\right\|_{L^\infty_{R\,dR}}\cdot\left\|\partial_R^{1+k}\phi_{2,\mu_2}\right\|_{L^2_{R\,dR}}.
	\end{align*}
	By Proposition \ref{prop:derivative1} and Lemma \ref{lem:derLinfty} we bound this as
	\begin{align*}
		&\hbar_{1}^{-\frac32-\delta}\frac{(\hbar_{3}^{2}\lambda)^{1-\frac{\delta}{2}}\langle\hbar_{3}^{2}\lambda\rangle^{\delta+\frac32}}{\lambda^{\frac72}}\sum_{\mu_{1,2}<\lambda}\mu_{1}^{\frac{i+j}{2}}\left\|\xb_{1}\right\|_{\Sho_{1}\left(\xi_{1}\simeq\mu_{1}\right)}\cdot\mu_{2}^{\frac{k+1}{2}}\left\|\xb_{2}\right\|_{L^{2}_{d\xi}(\xi_{2}\simeq\mu_{2})}\\
		&\lesssim \hbar_{3}^{\frac12}\hbar_{2}^{-2+\delta}\sum_{\mu_{1,2}<\lambda}\frac{\lambda^{1-\frac{\delta}{2}}\langle\hbar_{3}^{2}\lambda\rangle^{\delta+\frac32}}{\mu_{2}^{1-\frac{\delta}{2}}\langle\hbar_{2}^{2}\mu_{2}\rangle^{\delta+\frac32}}\cdot \left(\frac{\mu_1}{\lambda}\right)^{\frac{i+j}{2}}\left\|\xb_1\right\|_{\Sho_0(\xi_{1}\simeq\mu_1)}\cdot \left(\frac{\mu_2}{\lambda}\right)^{\frac{1+k}{2}}\left\|\xb_2\right\|_{\Sht_0(\xi\simeq\mu_2)}.
	\end{align*}
	The rest argument is similar to the case when $\mu_{2}\leq \mu_{1}$.
\end{proof}
We next aim to derive an analogue of the preceding proposition but with all angular frequencies $\hbar_j\gtrsim 1$, again assuming the compatibility conditions on the $\hbar_j$ as in the proposition. Thus in that situation we no longer worry about losses in the $\hbar_j^{-1}$, but instead we need to worry about the last step in the proof, where we transferred the operator $H_n$ from the left side of the inner product to the right. 
\subsection{The good spaces for the exceptional angular momenta $n=0,\pm1$}\label{subsec:nexcspaces}
The next estimate will be formulated for arbitrary inputs, also those for angular momenta $n = 0,\pm 1$. Thus we need to introduce a norm for these as well. 
\\

{\it{$n = -1$}}. Recall from \eqref{eq:D_f}
\[
\mathcal{D}_-f(R) = \int_0^\infty y(\xi)\phi_{-}(R,\xi)\tilde{\rho}_{-1}(\xi)d\xi,\quad y(\xi): = \langle \mathcal{D}_-f(R),\,\phi_{-}(R,\xi)\rangle_{L^2_{R\,dR}}.
\]
To pass back from this formula to the underlying function $f$, write it as 
\begin{equation}\label{eq:levelup-rep}
f(R) = c_{-}\cdot \phi_{-1}(R) + \phi_{-1}(R)\cdot \int_0^R \left(\phi_{-1}(s)\right)^{-1}\mathcal{D}_-f(s)\,ds,\quad \phi_{-1}(R): = \frac{R^2}{1+R^2}. 
\end{equation}
Recall the asymptotics of the spectral measure 
\[
\tilde{\rho}_{-1}(\xi)\sim \xi\langle\xi\rangle^2. 
\]
In the sequel, it will also be important to recall the asymptotic bounds for the Fourier basis $\phi_{-1}(R,\xi)$, which gives $\big|\phi_{-1}(R,\xi)\big|\lesssim R^3/\langle R\rangle^2$, whence linear growth toward $R = +\infty$. 
Then we introduce the norm 
\begin{align}\label{eq:S0-1}
	\left\|y(\xi)\right\|_{S^{-}_{0}}:=\left\|\xi^{1-\frac{\delta}{2}}\langle\xi\rangle^{\frac52+\delta}y(\xi)\right\|_{L^{2}_{d\xi}}.
\end{align}
Then we shall describe the function $f(R)$ in terms of the pair $\big(c_{-}, y\big)$. For later reference we also introduce the space $S_1^{-}$ such that $\|\cdot\|_{S_{1}^{-}}:=\|\xi^{-\frac12}\cdot\|_{S^{-}_{0}}$. 

{\it{$n = 0$}}. We recall from \cite{KMiao} the representation 
\[
\mathcal{D}f(R) = \int_0^\infty y(\xi)\phi_{0}(R,\xi)\tilde{\rho}_{0}(\xi)d\xi,\quad y(\xi): = \langle \mathcal{D}f(R),\,\phi_{0}(R,\xi)\rangle_{L^2_{R\,dR}}.
\]
Here we use the notation $\phi_0(R,\xi) := \phi(R, \xi)$ where the latter Fourier basis is described in \cite{KMiao}, and the spectral measure $\tilde{\rho}_{0}$ corresponds to $\tilde{\rho}(\xi)$ in loc. cit.. Moreover, denoting by $\phi_0(R): = \frac{R}{1+R^2}$ the resonance at zero frequency, we have the representation 
\begin{equation}\label{eq:levelup0rep}
f(R) = c_{0}\cdot\phi_{0}(R) + \phi_{0}(R)\cdot \int_0^R \phi_{0}^{-1}(s)\mathcal{D}f(s)\,ds,\quad \phi_{0}(R): = \frac{R}{1+R^2}. 
\end{equation}
Here we recall the spectral asymptotics
\begin{align*}
&\tilde{\rho}_{0}(\xi)\sim \langle\log\xi\rangle^{-2},\quad \xi\lesssim 1,\\
&\tilde{\rho}_{0}(\xi)\sim \xi^2,\quad \xi\gg 1. 
\end{align*}
Then we introduce the norm 
\begin{equation}\label{eq:S00}
\left\|y(\xi)\right\|_{S_0^{0}}: = \left\|\left(\frac{\langle\log\langle\xi\rangle\rangle}{\langle\log\xi\rangle}\right)^{1+\delta}\xi^{\frac12}\langle\xi\rangle^{\frac52+\frac{\delta}{2}}y(\xi)\right\|_{L^2_{d\xi}}. 
\end{equation}
\begin{remark}\label{rmk: norm in KM}
	 This norm is similar to the one use in \cite{KMiao}, except it reflects $H^{4+}$ regularity (for the derivative $\mathcal{D}f$) rather than $H^{3+}$-regularity as in \cite{KMiao}; moreover, the logarithmic weight is slightly altered here in order for this norm to exactly match the regularity for the other $n$ at high frequencies. 
\end{remark}
{\it{$n = 1$}}. Recall from Lemma \ref{lem:H1p} the representation 
\[
\mathcal{D}_+f(R) = \frac{2}{\pi^{2}}\int_0^\infty y(\xi)\phi_{1}(R,\xi)\xi d\xi,\quad y(\xi): = \langle \mathcal{D}_+f(R),\,\phi_{1}(R,\xi)\rangle_{L^2_{R\,dR}}.
\]
where we set $\phi_{1}(R,\xi)$ is identical to $\phi(R, z)$ in \eqref{modified funda sys} with $z = \xi$. In the following we will denote by $\trho_{1}(\xi):=\frac{2}{\pi^{2}}\xi$ the spectral measure associated to $\phi_{1}(R,\xi)$. Then we have the representation
\begin{equation}\label{eq:levelup1rep}
f(R) = c_{1}\cdot\phi_{1}(R) + \phi_{1}(R)\cdot \int_0^R \phi_{1}(s)^{-1}\mathcal{D}_+f(s)\,ds,\quad\phi_{1}(R): = \frac{1}{1+R^2}. 
\end{equation}
Then we introduce the norm 
\begin{equation}\label{eq:S01}
\left\|y(\xi)\right\|_{S_0^{+}}: = \left\|\xi^{1-\frac{\delta}{2}}\langle\xi\rangle^{\frac32+\delta}y(\xi)\right\|_{L^2_{d\xi}}. 
\end{equation}
We then have the following basic analogous of Proposition~\ref{prop:derivative1} and Lemma~\ref{lem:derLinfty}:
\begin{proposition}\label{prop:exceptionalnderivative1}
	Let 
	\begin{align*}
	&f_{-1}(R) = c_{-1}\cdot \phi_{-1}(R) + \phi_{-1}(R)\cdot \int_0^R \phi_{-1}(s)^{-1}\mathcal{D}_-f_{-1}(s)\,ds\\
	&f_0(R) = c_{0}\cdot\phi_{0}(R) + \phi_{0}(R)\cdot \int_0^R \phi_{0}(s)^{-1}\mathcal{D}f_0(s)\,ds\\
	&f_1(R) = c_{1}\cdot\phi_{1}(R) + \phi_{1}(R)\cdot \int_0^R \phi_{1}(s)^{-1}\mathcal{D}_+f_1(s)\,ds
	\end{align*}
	Moreover, set 
	\begin{align*}
	y_{\mp 1}(\xi)= \langle \mathcal{D}_{\mp}f_{\mp 1}(R),\,\phi_{\mp 1}(R,\xi)\rangle_{L^2_{R\,dR}},\quad y_0(\xi): = \langle \mathcal{D}f_0(R),\,\phi_{0}(R,\xi)\rangle_{L^2_{R\,dR}}
	\end{align*}
	Then for $k\geq 1$ we have the derivative bounds (for $\tau\gg 1$)
	\begin{equation}\label{eq:exceptionalnderL2}\begin{split}
	&\left\|\partial_R^k f_{-1}(R)\right\|_{L^2_{R\,dR}}\lesssim \left|c_{-1}\right| + \left\|\xi^{\frac12-\frac{\delta}{2}}\langle\xi\rangle^{\frac12+\delta}\left(\tilde{\rho}_{-1}(\xi)\right)^{\frac12}y_{-1}(\xi)\right\|_{L^2_{d\xi}}\\
	&\hspace{3.2cm} +\left\|\xi^{\frac{k-1}{2}}\left(\tilde{\rho}_{-1}(\xi)\right)^{\frac12}y_{-1}(\xi)\right\|_{L^2_{d\xi}} ,\\
	&\left\|\partial_R^k f_{+1}(R)\right\|_{L^2_{R\,dR}(R\lesssim\tau)}\lesssim \left|c_{1}\right| + \tau^{\max\{2-k,0\}(1-\delta)}\left\|\xi^{\frac12-\frac{\delta}{2}}\langle\xi\rangle^\delta\left(\tilde{\rho}_{1}(\xi)\right)^{\frac12}y_{1}(\xi)\right\|_{L^2_{d\xi}}\\
	&\hspace{3.2cm} +\left\|\xi^{\frac{k-1}{2}}\left(\tilde{\rho}_{1}(\xi)\right)^{\frac12}y_{1}(\xi)\right\|_{L^2_{d\xi}} ,\\
	&\left\|\partial_R^k f_0(R)\right\|_{L^2_{R\,dR}}\lesssim \left|c_0\right| +\tau^{(1-\delta)\max\{2-k,0\}} \left\|\xi^{\frac12-\frac{\delta}{2}}\langle\xi\rangle^{\delta}\left(\tilde{\rho}_{0}(\xi)\right)^{\frac12}y_0(\xi)\right\|_{L^2_{d\xi}}\\
	&\hspace{3.2cm}+\left\|\xi^{\frac{k-1}{2}}\left(\tilde{\rho}_{0}(\xi)\right)^{\frac12}y_{0}(\xi)\right\|_{L^2_{d\xi}}
	\end{split}\end{equation}
	and in terms of the norms $S_0^{\pm}, S^{0}_{0}$, we have the bounds 
	\begin{align}\label{L2 inside cone}
	\begin{split}
	&\left\|\partial_R f_{\pm1}(R)\right\|_{L^2_{R\,dR}(R\lesssim \tau)}\lesssim \left|c_{\pm1}\right| + \tau^{1-\delta}\left\|y_{\pm1}\right\|_{S_0^{\pm}},\\
	&\left\|\partial_R f_0(R)\right\|_{L^2_{R\,dR}(R\lesssim \tau)}\lesssim \left|c_0\right| + \tau\left\|y_{0}\right\|_{S_0^0},\\
	&\left\|\partial_R^2 f_0(R)\right\|_{L^2_{R\,dR}(R\lesssim \tau)}\lesssim \left|c_0\right| + \sqrt{|\log\tau|}\left\|y_{0}\right\|_{S_0^0}.
	\end{split}
	\end{align}
	
	Moreover, for $k\geq 1$ we have the $L^\infty$-bounds
	\begin{align*}
	&\left\|\partial_R^kf_{-1}(R)\right\|_{L^\infty_{R\,dR}}\lesssim \left|c_{-1}\right| + \left\|\xi^{\frac{k}{2}-\frac{\delta}{2}}\langle\xi\rangle^{\delta}\tilde{\rho}_{-1}(\xi)y_j(\xi)\right\|_{L^2_{d\xi}}
	\end{align*}
	and in particular we have 
	\begin{align*}
	\left\|\partial_Rf_{-1}(R)\right\|_{L^\infty_{R\,dR}}\lesssim  \left|c_{-1}\right| +\left\|y_{-1}\right\|_{S_0^{-}}. 
	\end{align*}
	Similarly, we have 
	\begin{align*}
	&\left\|\partial_R^kf_{1}(R)\right\|_{L^\infty_{R\,dR}}\lesssim  \left|c_{1}\right| +\left\|\xi^{\frac{k}{2}-\frac{\delta}{2}}\langle\xi\rangle^{\delta}\tilde{\rho}_{1}(\xi)y_1(\xi)\right\|_{L^2_{d\xi}}
	\end{align*}
	and in particular we have 
	\begin{align*}
	\left\|\partial_Rf_{1}(R)\right\|_{L^\infty_{R\,dR}}\lesssim \left|c_{1}\right| +\left\|y_{1}\right\|_{S_0^{+}}. 
	\end{align*}
	Finally we have 
	\begin{align*}
	&\left\|\partial_R^kf_{0}(R)\right\|_{L^\infty_{R\,dR}}\lesssim \left|c_{0}\right| +\left\|\xi^{\frac{k}{2}}\langle\xi\rangle^{\frac{\delta}{2}}\left(\frac{\langle\log\xi\rangle}{\langle\log\langle\xi\rangle\rangle}\right)^{1-\delta}\tilde{\rho}_{0}(\xi)y_0(\xi)\right\|_{L^2_{d\xi}},
	\end{align*}
	and in particular we have 
	\begin{align*}
	\left\|\partial_Rf_{0}(R)\right\|_{L^\infty_{R\,dR}}\lesssim \left|c_{0}\right| +\left\|y_{0}\right\|_{S_0^{0}}. 
	\end{align*}
\end{proposition}

\begin{proof}
	 We start by proving \eqref{eq:exceptionalnderL2}. 
	\\
	
	{\it{First estimate of \eqref{eq:exceptionalnderL2}}}. To begin with, we observe that 
	\begin{align*}
		\left\|c_{-1}\partial_R^k\phi_{-1}(R)\right\|_{L^2_{R\,dR}}\lesssim \left|c_{-1}\right|,\quad k\geq 1. 
	\end{align*}
	Next, write 
	\begin{align*}
		\partial_R^k\left( \phi_{-1}(R)\cdot \int_0^R \phi_{-1}(s)^{-1}\mathcal{D}_-f_{-1}(s)\,ds\right) &= \sum_{i+j =k,\,j\geq 1}C_{i,j}\partial_R^i\left(\phi_{-1}(R)\right)\partial_R^{j-1}\left(\phi_{-1}(R)^{-1}\mathcal{D}_-f_{-1}(R)\right)\\
		& + \partial_R^k\left( \phi_{-1}(R)\right)\cdot \int_0^R \phi_{-1}(s)^{-1}\mathcal{D}_-f_{-1}(s)\,ds
	\end{align*}
	Here we can quickly dispose of the second term on the right by observing that 
	\[
	\left\|\langle R\rangle\cdot \partial_R^k\left( \phi_{-}(R)\right)\right\|_{L^2_{R\,dR}}<\infty,\,k\geq 1,
	\]
	as well as the bounds (valid for all positive $s$)
	\[
	\left| \phi_{-1}(s)^{-1}\phi_{-1}(s,\xi)\right|\lesssim \min\left\{s, s^{-2}(1+s^2)\langle\xi\rangle^{-1}\xi^{-\frac12}(s\xi^{\frac12})^{-\frac12}\right\}
	\]
	Then write 
	\begin{align*}
		&\int_0^R \phi_{-1}(s)^{-1}\mathcal{D}_-f_{-1}(s)\,ds\\
		& = \int_0^{\min\{1,R\}}\phi_{-1}(s)^{-1}\int_0^1 \phi_{-1}(s,\xi) y_{-1}(\xi)\tilde{\rho}_{-1}(\xi)\,d\xi\,ds\\
		& +  \int_0^{\min\{1,R\}}\phi_{-1}(s)^{-1}\int_1^{s^{-2}} \phi_{-1}(s,\xi) y_{-1}(\xi)\tilde{\rho}_{-1}(\xi)\,d\xi\,ds\\
		& + \int_0^{\min\{1,R\}}\phi_{-1}(s)^{-1}\int_{s^{-2}}^\infty \phi_{-1}(s,\xi) y_{-1}(\xi)\tilde{\rho}_{-1}(\xi)\,d\xi\,ds\\
		& + \int_{\min\{1,R\}}^R\phi_{-1}(s)^{-1}\mathcal{D}_-f_{-1}(s)\,ds\\
		& = \sum_{j=1}^4 S_j.
	\end{align*}
	We easily obtain 
	\begin{align*}
		\left|S_1\right|\lesssim \int_0^1\int_0^1 s\xi \left|y_{-1}(\xi)\right|\,d\xi\, ds\lesssim \left\|\xi^{\frac32-\delta}y_{-1}(\xi)\right\|_{L^2_{d\xi}(\xi<1)},
	\end{align*}
	which is better than the bound needed for \eqref{eq:exceptionalnderL2}. For $S_{2,3}$, we change the order of integration:
	\begin{align*}
		\left|S_2\right|&\lesssim \int_1^\infty \int_0^{\min\{R, \xi^{-\frac12}\}}s\big|y_{-1}(\xi)\big|\tilde{\rho}_{-1}(\xi)\,ds\,d\xi\lesssim \int_1^\infty \xi^2\left|y_{-1}(\xi)\right|\,d\xi\lesssim \left\|\xi^{\frac52+\frac{\delta}{2}}y(\xi)\right\|_{L^2_{d\xi}(\xi>1)}
	\end{align*}
	in accordance with the bound needed for  \eqref{eq:exceptionalnderL2}, and similarly 
	\begin{align*}
		\left|S_3\right|&\lesssim  \int_1^\infty \left(\int_{\xi^{-\frac12}}^R s^{-2}(1+s^2)\langle\xi\rangle^{-1}\xi^{-\frac12}(s\xi^{\frac12})^{-\frac12}\,ds\right)\left|y_{-1}(\xi)\right|\tilde{\rho}_{-1}(\xi)\,d\xi\\
		&\lesssim \langle R\rangle\cdot \int_1^\infty \xi^2\left|y_{-1}(\xi)\right|\,d\xi\lesssim \langle R\rangle\cdot \left\|\xi^{\frac52+\frac{\delta}{2}}y_{-1}(\xi)\right\|_{L^2_{d\xi}(\xi>1)}, 
	\end{align*}
	and so the contribution of the this term to the full expression can be estimated by 
	\begin{align*}
		\left\|\partial_R^k\left( \phi_{-1}(R)\right)\cdot S_3\right\|_{L^2_{R\,dR}}&\lesssim \left\|\langle R\rangle\cdot \partial_R^k\left( \phi_{-}(R)\right)\right\|_{L^2_{R\,dR}}\cdot \left\|\langle R\rangle^{-1}S_3\right\|_{L^\infty_{R\,dR}}\lesssim \left\|\xi^{\frac52+\frac{\delta}{2}}y_{-1}(\xi)\right\|_{L^2_{d\xi}(\xi>1)}, 
	\end{align*}
	which is consistent with \eqref{eq:exceptionalnderL2}. Finally, for the last term, we may assume $R>1$, of course, and then 
	\begin{align*}
		\left|S_4\right|&\lesssim \int_1^R \int_0^1\chi_{s\xi^{\frac12}\lesssim1}s\xi^{\frac12}\cdot \xi^{\frac12}|y_{-1}(\xi)\big|\,d\xi\,ds\\
		&+\int_1^R \int_0^1\chi_{s\xi^{\frac12}\gtrsim1} \xi^{-\frac12}(s\xi^{\frac12})^{-\frac12}\big|y_{-1}(\xi)\big|\,\xi d\xi\,ds\\
		& + \int_1^R \int_1^\infty \langle\xi\rangle^{-1}\xi^{-\frac12}(s\xi^{\frac12})^{-\frac12}\big|y_{-1}(\xi)\big|\xi^3\,d\xi\,ds,
	\end{align*}
	and all of these are seen to be bounded by 
	\[
	\lesssim \langle R\rangle \left\|\xi^{\frac12-\frac{\delta}{2}}\langle\xi\rangle^{\frac12+\delta}\left(\tilde{\rho}_{-1}(\xi)\right)^{\frac12}y_{-1}(\xi)\right\|_{L^2_{d\xi}},
	\]
	which then leads to the desired bound as see before. As far as the remaining term 
	\[
	\sum_{i+j =k,\,j\geq 1}C_{i,j}\partial_R^i\left(\phi_{-1}(R)\right)\partial_R^{j-1}\left(\phi_{-1}(R)^{-1}\mathcal{D}_-f_{-1}(R)\right),\quad k\geq 1,
	\]
	is concerned, using arguments analogous to the ones for the proof of Proposition~\ref{prop:derivative1}, we have for $i+j\leq k$ that 
	\begin{align*}
		\left\|\partial_R^i\left(\phi_{-1}(R)\right)\partial_R^{j-1}\left(\phi_{-1}(R)^{-1}\mathcal{D}_-f_{-1}(R)\right)\right\|_{L^2_{R\,dR}}\lesssim \left\|\xi^{\frac{k-1}{2}}\left(\tilde{\rho}_{-1}(\xi)\right)^{\frac12}y_{-1}(\xi)\right\|_{L^2_{d\xi}}.
	\end{align*}
	In fact, the term $\phi_{-1}(R)^{-1}$ being singular at the origin costs two powers of $\xi^{\frac12}$ for $\xi\gg 1$, but the factor $\partial_R^i\big(\phi_{-1}(R)\big)$ contributes at least $\max\{2-i,0\}$ powers of $R$, whence we have to pay 
	\[
	j+1 - \max\{2-i,0\} = k-i+1- \max\{2-i,0\}\leq k-1
	\]
	many powers of $\xi^{\frac12}$ in the high-frequency regime. In the low frequency regime $\xi\lesssim 1$ derivatives translate into gains of $R^{-1}$ which in turn translate into additional factors $\xi^{\frac12}$. When $R\xi^{\frac12}\gtrsim1$, this result is  straightforward. When $R\xi^{\frac12}\ll1$ and $k=1$, the result is also straightforward using Plancherel Theorem. Now we give a proof for  the case $R\xi^{\frac12}\ll1$ and $k=2$. We essentially need to prove the estimate:
		\begin{align*}
		\left\|\int_{0}^{\infty}\chi_{R\xi^{\frac12}\ll1}\partial_{R}\phi_{-1}(R,\xi)\,y_{-1}(\xi)\rho_{-1}(\xi)\,d\xi\right\|_{L^{2}_{R\,dR}}\lesssim \left\|\xi^{\frac12}y_{-1}(\xi)\left(\rho_{-1}(\xi)\right)^{\frac12}\right\|_{L^{2}_{d\xi}}.
		\end{align*}
	To prove this, we consider the dyadic localization $R\xi^{\frac12}\simeq \mu\ll1, \xi^{\frac12}\simeq \lambda$, which implies $R\simeq \frac{\mu}{\lambda}$. By Proposition \ref{prop: spectral measure -1 diff}, we have
	\begin{align*}
		\left|\chi_{R\xi^{\frac12}\ll1}\partial_{R}\phi_{-1}(R,\xi)\right|\lesssim R^{2}\left\langle R^{2}\right\rangle^{-1}.
	\end{align*}
Therefore we have
\begin{align*}
	&\left\|\int_{0}^{\infty}\chi_{R\xi^{\frac12}\ll1}\partial_{R}\phi_{-1}(R,\xi)\,y_{-1}(\xi)\rho_{-1}(\xi)\,d\xi\right\|_{L^{2}_{R\,dR}}\\
	\lesssim &\left(\sum_{\mu\ll1, \lambda}\left(\frac{\mu}{\lambda}\right)^{4}\left\langle\frac{\mu^{2}}{\lambda^{2}}\right\rangle^{-2}\left\|\int_{0}^{\infty}\chi_{R\xi^{\frac12}\simeq\mu}\chi_{\xi^{\frac12}\simeq\lambda}y_{-1}(\xi)\rho_{-1}(\xi)\,d\xi\right\|^{2}_{L^{2}_{R\,dR}}\right)^{\frac12}\\
	\lesssim &\left(\sum_{\mu\ll1,\lambda}\left(\frac{\mu}{\lambda}\right)^{4}\left\langle\frac{\mu^{2}}{\lambda^{2}}\right\rangle^{-2}\cdot\left(\frac{\mu}{\lambda}\right)^{2}\cdot\lambda^{2}\left\|\chi_{\xi^{\frac12}\simeq\lambda}\xi^{\frac12}y_{-1}(\xi)\left(\rho_{-1}(\xi)\right)^{\frac12}\right\|^{2}_{L^{2}_{d\xi}}\right)^{2}\\
	\lesssim &\left(\sum_{\mu\ll1,\lambda}\mu^{2}\left\|\chi_{\xi^{\frac12}\simeq\lambda}\xi^{\frac12}y_{-1}(\xi)\left(\rho_{-1}(\xi)\right)^{\frac12}\right\|^{2}_{L^{2}_{d\xi}}\right)^{\frac12},
\end{align*}
which is desired.
	This completes the proof of the first inequality in \eqref{eq:exceptionalnderL2}.
	\\
	
	{\it{Second estimate of \eqref{eq:exceptionalnderL2}}}. As in the preceding case the contribution of the root mode $c_{1}\phi_{1}(R)$ is trivial, and so we now treat the contribution of the integral expression. As before write
	\begin{align*}
		\partial_R^k\left( \phi_{1}(R)\cdot \int_0^R \phi_{1}(s)^{-1}\mathcal{D}_{+}f_{1}(s)\,ds\right) &= \sum_{i+j =k,\,j\geq 1}C_{i,j}\partial_R^i\left(\phi_{1}(R)\right)\partial_R^{j-1}\left(\phi_{+1}(R)^{-1}\mathcal{D}_{+}f_{1}(R)\right)\\
		& + \partial_R^k\left( \phi_{1}(R)\right)\cdot \int_0^R \phi_{1}(s)^{-1}\mathcal{D}_{+}f_{1}(s)\,ds
	\end{align*}
	The loss of a factor $\tau$ for $k = 1$ in the estimate is due to the case $k = 1$ when the integral does not get hit by $\partial_R$. Considering this case first, we get the bound 
	\begin{align*}
		&\left|\partial_R\left( \phi_{1}(R)\right)\cdot \int_0^R \phi_{1}(s)^{-1}\mathcal{D}_{+}f_{1}(s)\,ds\right|\\
		&\lesssim \left|\partial_R\left( \phi_{1}(R)\right)\right|\cdot \int_0^1\int_0^{\min\{R, \xi^{-\frac12}\}}(1+s^2)s\,ds\left|y_{1}(\xi)\right|\,\xi d\xi\\
		& +  \left|\partial_R\left( \phi_{1}(R)\right)\right|\cdot \int_0^1\int_{\xi^{-\frac12}}^R (1+s^2)\cdot \xi^{-\frac12}(s\xi^{\frac12})^{-\frac12}\,ds\left|y_{1}(\xi)\right|\,\xi d\xi\\
		& +  \left|\partial_R\left( \phi_{1}(R)\right)\right|\cdot \int_1^\infty \int_0^R (1+s^2)\min\{s,\,\xi^{-\frac12}(s\xi^{\frac12})^{-\frac12}\}\,ds\cdot \left|y_{1}(\xi)\right|\,\xi d\xi\\
	\end{align*}
	Each of these expressions can be easily bounded by 
	\[
	\lesssim\langle R\rangle^{-\delta}\left\|\xi^{\frac12-\frac{\delta}{2}}\langle\xi\rangle^\delta\left(\tilde{\rho}_{1}(\xi)\right)^{\frac12}y_{1}(\xi)\right\|_{L^2_{d\xi}}.
	\]
	Correspondingly the $L^2_{R\,dR}$-norm over the region $R\lesssim \tau$ is bounded by the second expression on the right of \eqref{eq:exceptionalnderL2} in the case $k = 1$. Observe that expression we obtain in the $k = 1$ case when $\partial_R$ falls on the integral is simply $\mathcal{D}_+f_{1}(R)$, and its $L^2_{R\,dR}$ norm is simply bounded by $\left\|\left(\tilde{\rho}_{1}(\xi)\right)^{\frac12}y_{1}(\xi)\right\|_{L^2_{d\xi}}$ by Plancherel for the distorted Fourier transform, again compatible with the second bound in  \eqref{eq:exceptionalnderL2}. 
	\\
	Next consider the case $k = 2$, again assuming that both derivatives fall on the first factor $\phi_{1}(R)$. For $R\gg1$, we have
	\begin{align*}
		\left|\partial^{2}_{R}\left(\phi_{1}(R)\right)\cdot\int_{0}^{R}\phi_{1}(s)^{-1}\calD_{+}f_{1}(s)\,ds\right|\lesssim \langle R\rangle^{-1-\delta}\left\|\xi^{\frac12-\frac{\delta}{2}}\langle\xi\rangle^\delta\left(\tilde{\rho}_{1}(\xi)\right)^{\frac12}y_{1}(\xi)\right\|_{L^2_{d\xi}},
	\end{align*}
	whose $L^{2}_{RdR}$-norm is bounded by the second expression on the right of \eqref{eq:exceptionalnderL2} in the case $k = 2$.
	The remaining terms arising in the $k = 2$ case are 
	\[
	\partial_R\left( \phi_{1}(R)\right)\cdot [\phi_{1}(R)]^{-1}\mathcal{D}_{+}f_{1}(R),\quad \phi_{1}(R)\cdot \partial_R\left([\phi_{1}(R)]^{-1}\mathcal{D}_{+}f_{1}(R)\right).
	\]
	In order to estimate them, we use the Plancherel's theorem for the distorted Fourier transform, as well as the bounds
	\begin{align*}
		&\left| \partial_R\left( \phi_{1}(R)\right) [\phi_{1}(R)]^{-1}\right| + \left|\phi_{1}(R)\cdot \partial_R\left([\phi_{1}(R)]^{-1}\right)\right|\lesssim 1,\\
		&\left\|\partial_R\mathcal{D}_{+}f_{1}(R)\right\|_{L^2_{R\,dR}}\lesssim \left\|\xi^{\frac12}\left(\tilde{\rho}_{1}(\xi)\right)^{\frac12}y_{1}(\xi)\right\|_{L^2_{d\xi}}, 
	\end{align*}
	which imply the desired estimate in accordance with the second inequality in \eqref{eq:exceptionalnderL2}. The case of derivatives of degree $k\geq 3$ is straightforward due to the fact that now 
	\[
	\partial_R^k\left( \phi_{1}(R)\right)\cdot\langle R\rangle\cdot  [\phi_{1}(R)]^{-1}\in L^2_{R\,dR}.
	\]
	We omit the straightforward details. 
	\\
When $n=0$, the argument is similar to that of $n=1$, and most delicate step is to estimate the contribution from
\begin{align*}
&	\left|\partial_{R}\left(\phi_{0}(R)\right)\cdot\int_{0}^{R}\phi_{0}(s)^{-1}\calD f_{0}(s)\,ds\right|\\
\lesssim &\left|\partial_R\left( \phi_{0}(R)\right)\right|\cdot \int_0^1\int_0^{\min\{1,R, \xi^{-\frac12}\}}\left(s+s^{-1}\right)s^{2}\,ds\left|y_{0}(\xi)\right|\,\langle\log\xi\rangle^{-2} d\xi\\
&+\left|\partial_R\left( \phi_{0}(R)\right)\right|\cdot \int_0^1\int_1^{\min\{R, \xi^{-\frac12}\}}\left(s+s^{-1}\right)\log s\,ds\left|y_{0}(\xi)\right|\,\langle\log\xi\rangle^{-2} d\xi\\
& +  \left|\partial_R\left( \phi_{0}(R)\right)\right|\cdot \int_0^1\int_{\xi^{-\frac12}}^R \left(s+s^{-1}\right)\cdot (s\xi^{\frac12})^{-\frac12}\langle\log\xi\rangle\,ds\left|y_{0}(\xi)\right|\,\langle\log\xi\rangle^{-2}d\xi\\
& +  \left|\partial_R\left( \phi_{0}(R)\right)\right|\cdot \int_1^\infty \int_0^R \left(s+s^{-1}\right)\min\{s^{2}, (s\xi^{\frac12})^{-\frac12}\xi^{-1}\}\,ds\cdot \left|y_{0}(\xi)\right|\,\langle\xi\rangle^{2}d\xi.
\end{align*}
Similar as the $n=1$ case, each of these expressions is bounded by
\begin{align*}
	\lesssim \langle R\rangle^{-\delta}\left\|\xi^{\frac12-\frac{\delta}{2}}\langle\xi\rangle^{\delta}\left(\trho_{0}(\xi)\right)^{\frac12}y_{0}(\xi)\right\|_{L^{2}_{d\xi}}.
\end{align*}
The rest argument is identical to the $n=1$ case, and we omit the details.
\\

Now we turn to the estimates \eqref{L2 inside cone}. Based on the estimates in \eqref{eq:exceptionalnderL2}, the proof for \eqref{L2 inside cone} is more straightforward and here we only give an outline of its proof. For instance we look at $\partial_{R}f_{1}$:
\begin{align*}
	\partial_{R}f_{1}(R)=c_{1}\partial_{R}\phi_{1}(R)+\calD_{+}f_{1}(R)+\left(\partial_{R}\phi_{1}(R)\right)\cdot\int_{0}^{R}\left(\phi_{1}(s)\right)^{-1}\calD_{+}f_{1}(s)\,ds
\end{align*}
The first and the last term on the RHS above can be estimated exactly the same way as for the proof of \eqref{eq:exceptionalnderL2}, in view of the fact that the $\|\cdot\|_{S_{0}^{+}}$-norm controls the norm $\left\|\xi^{\frac12-\frac{\delta}{2}}\langle\xi\rangle^{\delta}\left(\trho_{1}(\xi)\right)^{\frac12}\right\|_{L^{2}_{d\xi}}$. For the second term, we first estimate the pointwise bound for $\calD_{+}f_{1}(R)$, which decays for $R\gg1$, and take the $L^{2}_{RdR}$-norm over the regime $R\lesssim \tau$. More precisely, we have (assuming $R\gg1$, and omitting the constant coefficients)
\begin{align*}
	\calD_{+}f_{1}(R)=\int_{0}^{\infty}\phi_{1}(R,\xi)y_{1}(\xi)\xi\,d\xi=\left(\int_{0}^{R^{-2}}+\int_{R^{-2}}^{1}+\int_{1}^{\infty}\right)\phi_{1}(R,\xi)y_{1}(\xi)\xi\,d\xi.
\end{align*}
For the integral $\int_{0}^{R^{-2}}...$ we have
\begin{align*}
	\left|\int_{0}^{R^{-2}}\phi_{1}(R,\xi)y_{1}(\xi)\xi\,d\xi\right|\lesssim \int_{0}^{R^{-2}}R\xi|y_{1}(\xi)|\,d\xi\lesssim \left(\int_{0}^{R^{-2}}\xi^{-1+\delta}\,d\xi\right)^{\frac12}\left\|y_{1}\right\|_{S_{0}^{+}}\lesssim R^{-\delta}\left\|y_{1}\right\|_{S_{0}^{+}}.
\end{align*}
For the integral $\int_{R^{-2}}^{1}...$ we have
\begin{align*}
	\left|\int_{R^{-2}}^{1}\phi_{1}(R,\xi)y_{1}(\xi)\xi\,d\xi\right|\lesssim R^{-\frac12}\int_{R^{-2}}^{1}\xi^{\frac14}|y_{1}(\xi)|\,d\xi\lesssim R^{-\frac12}\left(\int_{R^{-2}}^{1}\xi^{-\frac32+\delta}\,d\xi\right)^{\frac12}\left\|y_{1}\right\|_{S_{0}^{+}}\lesssim R^{-\delta}\left\|y_{1}\right\|_{S_{0}^{+}}.
\end{align*}
Finally for the integral $\int_{1}^{\infty}...$ we have
\begin{align*}
	\left|\int_{1}^{\infty}\phi_{1}(R,\xi)y_{1}(\xi)\xi\,d\xi\right|\lesssim R^{-\frac12}\int_{1}^{\infty}\xi^{\frac14}|y_{1}(\xi)|\,d\xi\lesssim R^{-\frac12}\|y_{1}\|_{S_{0}^{+}}.
\end{align*}
Therefore the estimate on $\|\partial_{R}f_{1}\|_{L^{2}_{RdR}(R\lesssim \tau)}$ follows. The other estimates in \eqref{L2 inside cone} follow in the same way. 
\\

Finally the $L^{\infty}$-estimates on $\partial_{R}^{k}f_{j}$ can be proved directly using the Fourier representations for $\calD_{-}f_{-1},\calD f_{0}$ and $\calD_{+}f_{1}$. We omit the details here.
\end{proof}

In the sequel, we shall need certain weighted versions of inequalities in the preceding proposition which allow to mostly avoid the losses in $\tau$, which may be understood as a low-frequency issue:
\begin{lemma}\label{lem:exceptionalnderivativeparad1} We have the following weighted and paradifferentiated derivative bounds, where $y_j$ denotes the Fourier transform in analogy to the preceding proposition, and where we assume $c_j = 0$: 
	\begin{align*}
	&\left(\sum_{\lambda\lesssim 1}\lambda^{1-\delta}\left\|\partial_Rf_{j,[\lambda,1]}\right\|_{L^2_{R\,dR}}^2\right)^{\frac12}\lesssim \left\|y_j\right\|_{S_0^{j}},\quad j = \pm 1, \\
	&\left(\sum_{\lambda\lesssim 1}\lambda^{1-\delta}\left\|\partial_Rf_{0,[\lambda,1]}\right\|_{L^2_{R\,dR}(R\lesssim\tau)}^2\right)^{\frac12}\lesssim \tau^{\delta}\left\|y_j\right\|_{S_0^{0}},\quad \tau\gtrsim 1. 
	\end{align*}
	Here $\lambda$ ranges over dyadic frequencies. The same bound obtains if $[\lambda,1]$ is replaced by $[\lambda, a)$ where $a\in [1,\infty]$. Here we use the notations $S_{0}^{1}:=S_{0}^{+}, S_{0}^{-1}:=S_{0}^{-}$.
\end{lemma}
\begin{proof}  We observe that the difference between the cases $j = \pm 1$ and $j = 0$ arises due to the different weights in the small frequency region in the norms $\|\cdot\|_{S_0^{j}}$. We treat the cases $j =\mp 1$, the case $j = 0$ being similar. \\
	{\it{$j = -1$.}} Write as before 
	\begin{align*}
	\partial_R f_{-1,[\lambda,1]} = \partial_R\phi_{-1}(R)\cdot\int_0^R[\phi_{-1}(s)]^{-1}\cdot \mathcal{D}_{-}f_{-1,[\lambda,1]}(s)\,ds+ \mathcal{D}_{-}f_{-1,[\lambda,1]}(R)
	\end{align*}
	Using Plancherel's theorem for the distorted Fourier transform we infer for dyadic $\mu\in [\lambda, 1]$
	\begin{align*}
	\lambda^{\frac12-\frac{\delta}{2}}\left\|\mathcal{D}_{-}f_{-1, \mu}(R)\right\|_{L^2_{R\,dR}}&\lesssim \left(\frac{\lambda}{\mu}\right)^{\frac12-\frac{\delta}{2}}\cdot \mu^{\frac12-\frac{\delta}{2}}\left\|(\tilde{\rho}_{-1}(\xi))^{\frac12}y_{-1}(\xi)\right\|_{L^2_{d\xi}(\xi\simeq\mu)}\\
	&\lesssim \left(\frac{\lambda}{\mu}\right)^{\frac12-\frac{\delta}{2}}\cdot \left\|y_{-1}(\xi)\right\|_{S_0^{-1}(\xi\simeq\mu)}.
	\end{align*}
	It follows that 
	\begin{align*}
	\left(\sum_{\lambda\lesssim 1}\lambda^{1-\delta} \left\|\mathcal{D}_{-}f_{-1,[\lambda,1]}\right\|_{L^2_{R\,dR}}^2\right)^{\frac12}&\lesssim  \left(\sum_{\lambda\lesssim 1}\sum_{\lambda\lesssim\mu\lesssim1}\left(\frac{\lambda}{\mu}\right)^{1-\delta}\cdot \big\|y_{-1}(\xi)\big\|_{S_0^{-1}(\xi\simeq\mu)}^2\right)^{\frac12}\\
	&\lesssim \left\|y_{-1}(\xi)\right\|_{S_0^{-1}}.
	\end{align*}
	Next consider the integral term above.
	\[
	\partial_R\phi_{-1}(R)\cdot\int_0^R\left[\phi_{-1}(s)\right]^{-1}\cdot \mathcal{D}_{-1}f_{-1,[\lambda,1]}(s)\,ds.
	\]
	Following from the proof of the preceding proposition, we immediately obtain
	\[
	\left\|\partial_{R}\phi_{-1}(R)\cdot\int_0^R\left[\phi_{-1}(s)\right]^{-1}\cdot \mathcal{D}_{-}f_{-1,[\lambda,1]}(s)\,ds\right\|_{L^{2}_{RdR}}\lesssim \left\|y_{-1}\right\|_{S_0^{-1}},
	\]
	 and the desired weighted bound follows.
	\\
	{\it{$j = +1$.}} Write 
	\begin{align*}
	\partial_R f_{1,[\lambda,1]} &= \partial_R\phi_{1}(R)\cdot\int_0^R\left[\phi_{1}(s)\right]^{-1}\cdot \mathcal{D}_{+}f_{1,[\lambda,1]}(s)\,ds+ \mathcal{D}_{+}f_{1,[\lambda,1]}(R)
	\end{align*}
	The contribution of the second term is treated exactly as in the preceding case. However, the contribution of the first term is more complicated. Our main goal here is to control the growth in $R$ of its $L^{2}_{RdR}$-norm. Therefore without loss of generality, we assume that $s\in[1,R]$. Fix dyadic $\mu\in[\lambda, 1]$ and consider the localized term (where $y_{1,\mu}(\xi) = \chi_{\xi\simeq\mu}y_{1}(\xi)$)
	\begin{align*}
	&\partial_R\phi_{1}(R)\cdot\int_0^R\left[\phi_{1}(s)\right]^{-1}\cdot \mathcal{D}_{+}f_{1,\mu}(s)\,ds\\
	 =& \partial_R\phi_{1}(R)\cdot\int_0^R\left[\phi_{1}(s)\right]^{-1}\cdot\int_0^\infty \chi_{s\xi^{\frac12}\lesssim1}\phi_{1}(s,\xi)y_{1,\mu}(\xi)\tilde{\rho}_{1}(\xi)\,d\xi\,ds\\
	+& \partial_R\phi_{1}(R)\cdot\int_0^R\left[\phi_{1}(s)\right]^{-1}\cdot\int_0^\infty \chi_{s\xi^{\frac12}\gtrsim1}\phi_{1}(s,\xi)y_{1,\mu}(\xi)\tilde{\rho}_{1}(\xi)\,d\xi\,ds\\
	\end{align*}
	Then we bound the contribution of the first term on the right as follows, where we first localize $R$ to dyadic scale $R\simeq\rho$ (Recall that since $s\geq1$, we have $\xi\leq 1$ in this case.):
	\begin{align*}
	&\lambda^{\frac12-\frac{\delta}{2}}\left\|\partial_R\phi_{1}(R)\cdot\int_0^R[\phi_{1}(s)]^{-1}\cdot\int_0^\infty \chi_{s\xi^{\frac12}\lesssim1}\phi_{1}(s,\xi)y_{1,\mu}(\xi)\tilde{\rho}_{1}(\xi)\,d\xi\,ds\right\|_{L^2_{R\,dR}(R\simeq\rho)}\\
	&\lesssim \left(\frac{\lambda}{\mu}\right)^{\frac12-\frac{\delta}{2}}\cdot \left(\sum_{\kappa\lesssim \min\{\mu^{-\frac12}, \rho\}}\left(\frac{\kappa}{\rho}\right)^2\cdot \left(\mu^{\frac12}\kappa\right)^2\right)\cdot\left\|y_{1,\mu}\right\|_{S_0^{1}(\xi\simeq\mu)}
	\end{align*}
	where $\kappa$ runs over dyadic scales for $s\in[1,R]$, and we have 
	\begin{align*}
	\sum_{\kappa\lesssim \min\{\mu^{-\frac12}, \rho\}}\left(\frac{\kappa}{\rho}\right)^2\cdot \left(\mu^{\frac12}\kappa\right)^2\lesssim \frac{\mu^{\frac12}\rho}{\langle\mu^{\frac12}\rho\rangle^2}.
	\end{align*}
	We conclude that 
	\begin{align*}
	&\left(\sum_{\lambda\lesssim 1}\lambda^{1-\delta}\left\| \partial_R\phi_{1}(R)\cdot\int_0^R\left[\phi_{1}(s)\right]^{-1}\cdot\int_0^\infty \chi_{s\xi^{\frac12}\lesssim1}\phi_{1}(s,\xi)y_{1,[\lambda,1]}(\xi)\tilde{\rho}_{1}(\xi)\,d\xi\,ds\right\|_{L^2_{R\,dR}}^2\right)^{\frac12}\\
	&\lesssim \left(\sum_{\lambda\lesssim 1}\sum_{\rho}\left[\sum_{\mu\in [\lambda,1]}\left(\frac{\lambda}{\mu}\right)^{\frac12-\frac{\delta}{2}}\cdot \frac{\mu^{\frac12}\rho}{\langle\mu^{\frac12}\rho\rangle^2}\cdot \left\|y_{1,\mu}\right\|_{S_0^{+}(\xi\simeq\mu)}
	\right]^2\right)^{\frac12}\\
	&\lesssim  \left(\sum_{\lambda\lesssim 1}\sum_{\rho}\sum_{\mu\in [\lambda,1]}\left(\frac{\lambda}{\mu}\right)^{1-\delta}\cdot \frac{\mu\rho^{2}}{\langle\mu^{\frac12}\rho\rangle^4}\cdot \left\|y_{1,\mu}\right\|_{S_0^{+}(\xi\simeq\mu)}^2
	\right)^{\frac12}\\
	&\lesssim \left\|y_{1}\right\|_{S_0^{+}},
	\end{align*}
	which corresponds to the desired bound. As for the oscillatory term above where the integration is restricted to $s\xi^{\frac12}\gtrsim 1$, we have to perform integration by parts (for instance, three times) with respect to $s$, which leads to the analogous bound 
	\begin{align*}
	&\lambda^{\frac12-\frac{\delta}{2}}\left\|\partial_R\phi_{1}(R)\cdot\int_0^R\left[\phi_{1}(s)\right]^{-1}\cdot\int_0^\infty \chi_{s\xi^{\frac12}\gtrsim1}\phi_{1}(s,\xi)y_{1,\mu}(\xi)\tilde{\rho}_{1}(\xi)\,d\xi\,ds\right\|_{L^2_{R\,dR}(R\simeq\rho)}\\
	&\lesssim \left(\frac{\lambda}{\mu}\right)^{\frac12-\frac{\delta}{2}}\cdot \left(\sum_{\mu^{-\frac12}\lesssim\kappa\lesssim\rho}\left(\frac{\kappa}{\rho}\right)^2\cdot \left(\mu^{\frac12}\kappa\right)^{-2}\right)\cdot\left\|y_{1,\mu}\right\|_{S_0^{+}(\xi\simeq\mu)},
	\end{align*}
	and the estimate is concluded from here in the same way as before. 
	\\
	We omit the simple modifications of the preceding proof to also handle the generalised bounds involving $\partial_R f_{j,[\lambda,a)}$.

\end{proof}
\subsection{More general bilinear estimates}

We can now formulate the following proposition, in which we still restrict the output to `non-exceptional angular momenta', i.e., $|n_3|\geq 2$. 
\begin{proposition}\label{prop:bilin2} 
	Let $|n_1|\gg 1$, $|n_2|\leq |n_1|$, where $n_2$ is allowed to take any integer value. Also, assume that $|n_3|\geq 2$ and that either $n_3\simeq n_1$ or else $|n_3|\ll |n_1|$ as well as $n_1\simeq -n_2$. For $|n|\geq 2$, we say that $\phi(R)$ is an angular momentum $n$ function provided 
	\[
	\phi(R) = \int_0^\infty \phi_{n}(R;\xi)\xb(\xi)\tilde{\rho}_{n}(\xi)\,d\xi,
	\]
	in which case we set 
	\[
	\left\|\phi\right\|_{\tilde{S}^{(n)}}: = \left\|\xb(\xi)\right\|_{\Sh_{0}}.  
	\]
	For $n = 0,\pm 1$ we say that $\phi(R)$ is an angular momentum $n$ function provided 
	\[
	\phi(R) = c_n\phi_n(R) + \phi_n(R)\cdot\int_0^R\left[\phi_n(s)\right]^{-1}\tilde{\mathcal{D}}\phi(s)\,ds,\quad \tilde{\mathcal{D}}\phi(R) = \int_0^\infty x(\xi)\phi_{n}(R,\xi)\tilde{\rho}_{n}(\xi)\,d\xi,
	\]
	where $\tilde{\calD} = \mathcal{D}_-,\mathcal{D}, \mathcal{D}_+$ according to $n$ as in the preceding. Then set 
	\[
	\left\|\phi\right\|_{\tilde{S}^{(n)}}: = \left|c_n\right| + \left\| \xb(\xi)\right\|_{S_0^{(n)}},\quad \textrm{where}\quad S_{0}^{(0)}:=S^{0}_{0},\quad S_{0}^{(\pm1)}=S_{0}^{\pm}.
	\]
	With the conventions and under the assumptions on the $n_j$ stated at the beginning, and assuming $\tau\gg 1$, we have the following bound 
	\begin{equation}\label{eq:productlargen3alln2}
	\left\|\langle\phi_{n_{3}}(R;\xi),\,\chi_{R\lesssim\tau}\partial_R\phi_1\cdot\partial_R\phi_2\rangle_{L^2_{R\,dR}}\right\|_{\Shth_{1}}\lesssim \tau^{\frac{\delta}{2}}\langle n_2\rangle^{2}\prod_{j=1}^2\left\|\phi_j\right\|_{\tilde{S}^{(n_j)}}
	\end{equation}
	where we assume that the factors $\phi_j$ are angular momentum $n_j$ functions, $j = 1, 2$. 
	\\
	Next, assume that $|n_1|\lesssim 1$, while the other assumptions on $n_j$ stated above are still valid. Then if $\phi_{1,2}$ are angular momentum $n_j$ functions with finite $\|\cdot\|_{\tilde{S}^{(n_j)}}$-norm, the product 
	\[
	\partial_R\phi_1\cdot\partial_R\phi_2
	\]
	admits a third order Taylor expansion at $R = 0$ of the form
	\[
	P_3(R): = \sum_{\ell=0}^3\gamma_{\ell} R^\ell,
	\]
	where we have the bound
	\begin{equation}\label{eq:Taylorbound1}
	\sum_{\ell=0}^3\left|\gamma_\ell\right|\lesssim \prod_{j=1,2}\left\|\phi_j\right\|_{\tilde{S}^{(n_j)}}.
	\end{equation}
	Furthermore, we have the bound 
	\begin{equation}\label{eq:productsmalln1geq2alln2}
	\left\|\langle\phi_{n_{3}}(R;\xi),\,\chi_{R\lesssim\tau}\partial_R\phi_1\cdot\partial_R\phi_2 - \chi_{R\lesssim 1}P_3(R)\rangle_{L^2_{R\,dR}}\right\|_{S_1^{(\hbar_3)}}\lesssim \tau^{\frac{\delta}{2}}\prod_{j=1}^2\left\|\phi_j\right\|_{\tilde{S}^{(n_j)}}
	\end{equation}
	The last inequality remains correct if we subtract from $P_3$ those terms $\gamma_l R^l$ with $l\geq |n_3 - 1|, l\equiv n_3 - 1\text{(mod 2)}$, provided $|n_3 -1|\leq 3$. 
\end{proposition}

\begin{proof}
	{\it{Proof of \eqref{eq:productlargen3alln2}}}:\\
	In order to prove \eqref{eq:productlargen3alln2}, in light of Proposition~\ref{prop:bilin1}, it suffices to deal with the case when $n_2 = 0,\pm 1$, and so $n_3\simeq n_1$. 
	Then we repeat the proof of Proposition~\ref{prop:bilin1}. We split the proof into the same cases, and use the same terminology. Frequencies for the second factor $\phi_2$ correspond to the frequencies of the Fourier integral.
	\\
	
	{\it{(1.a): $\xi<\xi_1$. }} This proceeds exactly as in the earlier proof by taking advantage of the $L^\infty$-bounds coming from Proposition~\ref{prop:exceptionalnderivative1}. In particular, the contribution of the resonant/root part in $\phi_2$ leads to an admissible contribution, and we can omit it from now on. 
	\\
	
	{\it{(1.b): $\xi_2\geq \xi\geq \xi_1$. }} This we can write as in the earlier proof as 
	\begin{align*}
	&\sum_{\mu}\chi_{\xi\simeq\mu}\left\langle \phi_{n_{3}}(R;\xi),\,\chi_{R\lesssim\tau}\partial_R\phi_{1,\leq \mu}\cdot\partial_R\phi_{2,[\mu,\infty]}\right\rangle_{L^2_{R\,dR}}\\
	& = \sum_{\mu<\lambda}\chi_{\xi\simeq\mu}\left\langle \phi_{n_{3}}(R;\xi),\,\chi_{R\lesssim\tau}\partial_R\phi_{1,\leq \mu}\cdot[c_n\partial_R\phi_n(R) + \partial_R\phi_{2,\lambda}]\right\rangle_{L^2_{R\,dR}}.
	\end{align*}
	where we set 
	\begin{align*}
	&\phi_{2,\lambda}(R) = \phi_j(R)\cdot \int_0^R \left[\phi_j(s)\right]^{-1}\,\tilde{\mathcal{D}}\phi_{2,\lambda}(s)\,ds,\,\tilde{\mathcal{D}}\phi_{2,\lambda}(R) = \int_0^\infty \chi_{\xi\simeq\lambda}\xb_2(\xi)\phi_j(R,\xi)\tilde{\rho}_{j}(\xi)\,d\xi,\quad j = 0,\pm 1. 
	\end{align*}
	We distinguish between small output frequencies $\mu\lesssim 1$ and large ones. 
	\\
	
	{\it{Small output frequency $\mu\lesssim 1$.}} The contribution of the root/resonant part is straightforward to handle by means of Plancherel's Theorem and Lemma~\ref{lem:derLinfty}:
	\begin{align*}
	&\left\|\sum_{\mu \lesssim 1}\chi_{\xi\simeq\mu}\langle \phi_{n_{3}}(R;\xi),\,\chi_{R\lesssim\tau}\partial_R\phi_{1,\leq \mu}\cdot c_n\partial_R\phi_n(R)\rangle_{L^2_{R\,dR}}\right\|_{\Shth_{1}}\\
	&\lesssim \left|c_n\right|\cdot \sum_{\mu \lesssim 1}\hbar_3^{2-\delta}\mu^{\frac12-\frac{\delta}{2}}\left\|\partial_R\phi_{1,\leq \mu}\right\|_{L^\infty_{R\,dR}}\cdot \left\|\partial_R\phi_n(R)\right\|_{L^2_{R\,dR}}\\
	&\lesssim  \left|c_n\right|\cdot\left\|\phi_1\right\|_{\tilde{S}_0^{(n_1)}}\lesssim \prod_{j=1,2}\left\|\phi_j\right\|_{\tilde{S}^{(n_j)}}. 
	\end{align*}
	
	Next, denoting $\tilde{\phi}_2 = \phi_2 - c_n\phi_n(R)$ and taking advantage of Lemma~\ref{lem:exceptionalnderivativeparad1} as well as Lemma~\ref{lem:derLinfty} as well as orthogonality, we have 
	\begin{align*}
	&\left\|\sum_{\mu\lesssim 1}\chi_{\xi\simeq\mu}\left\langle \phi_{n_{3}}(R;\xi),\,\chi_{R\lesssim\tau}\partial_R\phi_{1,\leq \mu}\cdot\partial_R\tilde{\phi}_{2,[\mu,\infty]}\right\rangle_{L^2_{R\,dR}}\right\|_{\Shth_1}\\
	&\lesssim \left(\sum_{\mu\lesssim 1}\hbar_3^{4-\delta}\mu^{1-\delta}\left\|\chi_{R\lesssim\tau}\partial_R\phi_{1,\leq \mu}\cdot\partial_R\tilde{\phi}_{2,[\mu,\infty]}\right\|_{L^2_{R\,dR}}^2\right)^{\frac12}\\
	&\lesssim  \left(\sum_{\mu\lesssim 1}\hbar_3^{4-\delta}\mu^{1-\delta}\left\|\partial_R\phi_{1,\leq \mu}\right\|_{L^\infty_{R\,dR}}^2\cdot\left\|\partial_R\tilde{\phi}_{2,[\mu,\infty]}\right\|_{L^2_{R\,dR}(R\lesssim\tau)}^2\right)^{\frac12}\\
	&\lesssim \hbar_1^{\frac12}\tau^{\delta}\cdot \left\|\phi_1\right\|_{\tilde{S}^{(n_1)}}\cdot \left\|\phi_2\right\|_{\tilde{S}^{(n_2)}}
	\end{align*}
	
	{\it{Large output frequency $\mu\gg 1$.}} Here in principle we intend to absorb the frequency weight $\mu^{\frac12-\frac{\delta}{2}}\langle \mu\hbar_1^2\rangle^{\delta + \frac32}$ into the high frequency factor, which requires some care, due to the particular structure of $\partial_R\phi_2$. We quickly dispose of the contribution of the resonant/root part via integration by parts: assuming $\phi_2(R) = c_n\cdot\phi_{n}(R)$, $n = 0,\pm 1$, 
	\begin{align*}
	&c_n\chi_{\xi\simeq\mu}\left\langle \phi_{n_{3}}(R;\xi),\,\chi_{R\lesssim\tau}\partial_R\phi_{1,\leq \mu}\cdot\partial_R\phi_{n}(R)\right\rangle_{L^2_{R\,dR}}\\
	& = c_n\mu^{-3}\cdot  \chi_{\xi\simeq\mu}\left\langle \phi_{n_{3}}(R;\xi),\,H_{n_3}^3\left(\chi_{R\lesssim\tau}\partial_R\phi_{1,\leq \mu}\cdot\partial_R\phi_{n}(R)\right)\right\rangle_{L^2_{R\,dR}},
	\end{align*}
	where boundary terms at $R = 0$ play no role because $|n_1|\gg 1$ by assumption. Localize the term $\phi_{1,\leq \mu}$ to dyadic frequency $\simeq \lambda\lesssim\mu$ and write 
	\begin{align*}
	H_{n_3}^3\left(\chi_{R\lesssim\tau}\partial_R\phi_{1,\lambda}\cdot\partial_R\phi_{n}(R)\right) = \sum_{i+j+k = 6}\left(\frac{n_3}{R}\right)^i\partial_R^{j+1}\big(\phi_{1,\lambda}\big)\cdot \partial_R^k\big(\chi_{R\lesssim\tau}\partial_R\phi_{n}(R)\big).
	\end{align*}
	Then use Proposition~\ref{prop:singularmultiplier}, which gives (letting $\xb_1$ be the Fourier variable for $\phi_1$)
	\begin{align*}
	\left\|\left(\frac{n_3}{R}\right)^i\partial_R^{j+1}\left(\phi_{1,\lambda}\right)\cdot \partial_R^k\left(\chi_{R\lesssim\tau}\partial_R\phi_{n}(R)\right)\right\|_{L_{R\,dR}^2}\lesssim \left(\lambda^{\frac12}+\lambda^{\frac72}\right)\left\|\xb_1\right\|_{L^2_{d\xi}},\quad i+j\leq 6, 
	\end{align*}
	and so we infer 
	\begin{align*}
	&\left|c_n\right|\left\|\chi_{\xi\simeq\mu}\langle \phi_{n_{3}}(R;\xi),\,\chi_{R\lesssim\tau}\partial_R\phi_{1,\leq \mu}\cdot\partial_R\phi_{n}(R)\rangle_{L^2_{R\,dR}}\right\|_{\Shth_1}\\
	&\lesssim \left|c_n\right|\hbar_{3}^{5}\,\mu^{2+\frac{\delta}{2}}\cdot\mu^{-3}\cdot \left(\lambda^{\frac12}+\lambda^{\frac72}\right)\left\|\xb_1\right\|_{L^2_{d\xi}(\xi_1\simeq\lambda)}\\
	&\lesssim \left|c_n\right|\cdot\left(\frac{\lambda}{\mu}\right)^{1-\frac{\delta}{2}}\cdot \left\|\xb_1\right\|_{\Sho_0(\xi_1\simeq\lambda)}. 
	\end{align*}
	as long as we restrict $\lambda \gtrsim 1$. It follows that 
	\begin{align*}
	&\left|c_n\right|\left\|\sum_{\mu\gtrsim 1}\chi_{\xi\simeq\mu}\left\langle \phi_{n_{3}}(R;\xi),\,\chi_{R\lesssim\tau}\partial_R\phi_{1,1\leq\cdot\leq \mu}\cdot\partial_R\phi_{n}(R)\right\rangle_{L^2_{R\,dR}}\right\|_{\Shth_{1}}\\
	&\lesssim \left|c_n\right|\left(\sum_{\mu\gtrsim 1}\left\|\chi_{\xi\simeq\mu}\left\langle \phi_{n_{3}}(R;\xi),\,\chi_{R\lesssim\tau}\partial_R\phi_{1,1\leq\cdot\leq \mu}\cdot\partial_R\phi_{n}(R)\right\rangle_{L^2_{R\,dR}}\right\|^2_{\Shth_{1}}\right)^{\frac12}\\
	&\lesssim \left|c_n\right|\left(\sum_{\mu\gtrsim 1}\left[\sum_{1\lesssim\lambda\lesssim\mu}\left(\frac{\lambda}{\mu}\right)^{1-\frac{\delta}{2}}\cdot \left\|\xb_1\right\|_{\Sho_0(\xi_1\simeq\lambda)}\right]^2\right)^{\frac12}\\
	&\lesssim \left|c_n\right|\cdot \left\|\xb_1\right\|_{\Sho_{0}}\leq \left\|\phi_1\right\|_{\tilde{S}_0^{(n_1)}}\cdot \left\|\phi_2\right\|_{\tilde{S}_0^{(n_2)}}.
	\end{align*}
	The case when $\phi_1$ is at frequency $\lesssim 1$ is easily handled by placing this factor into $L^\infty_{R\,dR}$, we omit the details. 
	\\
	
	This reduces things to the contribution of the integral term, where we shall again have to perform some integration by parts in the inner product for very large output frequencies $\gtrsim \hbar_3^{-2}$. Assume henceforth that $\phi_2$ is given by the expression stated at the beginning of case {\it{(1.b)}}. 
	\\
	
	{\it{(i): intermediate output frequencies $1\lesssim \mu\lesssim \hbar_3^{-2}$.}} Here we have 
	\begin{align*}
	&\left\|\chi_{\xi\simeq\mu}\langle \phi_{n_{3}}(R;\xi),\,\chi_{R\lesssim\tau}\partial_R\phi_{1,\leq \mu}\cdot\partial_R\phi_{2,[\mu,\infty]}\rangle_{L^2_{R\,dR}}\right\|_{\Shth_{1}}\\
	&\lesssim \hbar_3^{2-\delta}\mu^{\frac12 - \frac{\delta}{2}}\sum_{\lambda\geq \mu}\left\|\partial_R\phi_{1,\leq \mu}\right\|_{L^\infty_{R\,dR}}\cdot \left\|\partial_R\phi_{2,\lambda}\right\|_{L^2_{R\,dR}}.
	\end{align*}
	where we can expand 
	\begin{align*}
	\partial_R\phi_{2,\lambda} &= \partial_R\left(\phi_n(R)\right)\cdot\int_0^R\left[\phi_n(s)\right]^{-1}\calD_n\phi_{2,\lambda}(s)\,ds\\
	& + \calD_n\phi_{2,\lambda}(R),\quad n = 0,\pm 1. 
	\end{align*}
	Using Plancherel's theorem for the distorted Fourier transform and the fact that $\lambda\geq \mu\gtrsim 1$, we obtain (with $\xb_2$ denoting the Fourier transform of $D_n\phi_2$)
	\begin{align*}
	\mu^{\frac12 - \frac{\delta}{2}}\left\| \calD_n\phi_{2,\lambda}(R)\right\|_{L^2_{R\,dR}}\lesssim \left(\frac{\mu}{\lambda}\right)^{\frac12-\frac{\delta}{2}}\left\|\xb_2\right\|_{\Sht_0(\xi\simeq\lambda)}.
	\end{align*}
	To bound the contribution of the integral term to $\partial_R\phi_{2,\lambda}$, we observe that in all cases $n = 0,\pm 1$ we have (under the hypothesis $\lambda\gtrsim\mu\gtrsim 1$)
	\begin{align*}
	&\mu^{\frac12 - \frac{\delta}{2}}\left\| \partial_R\left(\phi_n(R)\right)\cdot\int_0^R\left[\phi_n(s)\right]^{-1}\calD_n\phi_{2,\lambda}(s)\,ds\right\|_{L^2_{R\,dR}}\\
	&\lesssim \left(\frac{\mu}{\lambda}\right)^{\frac12-\frac{\delta}{2}}\left\|\xb_2\right\|_{\Sht_0(\xi\simeq\lambda)}.
	\end{align*}
In fact here we can restrict to the regime where $s\gg1$, since the other regime $s\lesssim1$ is much more straightforward. So in this case we have $s\lambda^{\frac12}\gg1$ and we are in the oscillatory regime. As an example, for $n=0$, we first localize $R$ to dyadic scale $R\simeq\rho$ and then perform integration by parts to obtain
\begin{align*}
	&\mu^{\frac12-\frac{\delta}{2}}\left\|\partial_{R}\phi_{0}(R)\cdot\int_{1}^{R}\left[\phi_{0}(s)\right]^{-1}\int_{0}^{\infty}\chi_{s\xi^{\frac12}\gtrsim1}\chi_{\xi\simeq\lambda}\phi_{0}(s,\xi)\xb_{2}(\xi)\trho_{0}(\xi)d\xi ds\right\|_{L^{2}_{RdR}(R\simeq\rho)}\\
	\lesssim &\left(\frac{\mu}{\lambda}\right)^{\frac12-\frac{\delta}{2}}
	\left(\sum_{\lambda^{-\frac12}\lesssim 1\lesssim \kappa\lesssim \rho}\left(\frac{\kappa}{\rho}\right)\cdot\left(\lambda^{\frac12}\kappa\right)^{-1}\right)\left\|\xb_{2}\right\|_{\Sht_0(\xi\simeq\lambda)}\\
	\lesssim &\left(\frac{\mu}{\lambda}\right)^{\frac12-\frac{\delta}{2}}\left\|\xb_{2}\right\|_{\Sht_0(\xi\simeq\lambda)}.
\end{align*}
If we combine these bounds with the usual $L^\infty$-bound 
\begin{align*}
\hbar_3^{2-\delta}\big\|\partial_R\phi_{1,\leq \mu}\big\|_{L^\infty_{R\,dR}}\lesssim \big\|\phi_1\big\|_{\tilde{S}_0^{(n_1)}}
\end{align*}
and to exploit orthogonality and the Cauchy-Schwarz inequality as usual to deduce the desired bound. 
\\

{\it{(ii): large  output frequencies $\mu\gg \hbar_3^{-2}$.}} Here the weight used for the norm $\|\cdot\|_{\Shth_1}$ becomes 
\[
\hbar_3\left(\mu\hbar_3^2\right)^{\frac12-\frac{\delta}{2}}\left\langle \mu\hbar_3^2\right\rangle^{\delta+\frac32}\simeq \hbar_3\left(\mu\hbar_3^2\right)^{2+\frac{\delta}{2}}
\]
We may again assume that 
\begin{align*}
\partial_R\phi_{2,[\mu,\infty]} = \partial_R\left(\phi_n(R)\right)\cdot\int_0^R\left[\phi_n(s)\right]^{-1}\calD_n\phi_{2,[\mu,\infty]}(s)\,ds
+ \calD_n\phi_{2,[\mu,\infty]}(R),\quad n = 0,\pm 1.
\end{align*}
Substituting the second term on the right for $\partial_R\phi_{2,[\mu,\infty]}$ leads to term that can be bounded directly via the Plancherel's theorem and Lemma~\ref{lem:derLinfty}:
\begin{align*}
&\left\|\chi_{\xi\simeq\mu}\left\langle \phi_{n_{3}}(R;\xi),\,\chi_{R\lesssim\tau}\partial_R\phi_{1,\leq \mu}\cdot  \calD_n\phi_{2,[\mu,\infty]}(R)\right\rangle_{L^2_{R\,dR}}\right\|_{\Shth_{1}}\\
&\lesssim \hbar_3(\mu\hbar_3^2)^{2+\frac{\delta}{2}}\cdot \left\|\chi_{R\lesssim\tau}\partial_R\phi_{1,\leq \mu}\right\|_{L^\infty_{R\,dR}}\cdot \left\|\calD_n\phi_{2,[\mu,\infty]}(R)\right\|_{L^2_{R\,dR}}\\
&\lesssim \sum_{\lambda\geq \mu}\left(\frac{\mu}{\lambda}\right)^{2+\frac{\delta}{2}}\left\|\phi_1\right\|_{\tilde{S}_0^{(n_1)}}\cdot \left\|\phi_2\right\|_{\tilde{S}_0^{(n_2)}(\xi\simeq\lambda)}.
\end{align*}
As for the contribution of the integral term to $\partial_R\phi_{2,[\mu,\infty]}$ we perform integration by parts as needed: Setting now 
\[
\partial_R\phi_{2,[\mu,\infty]} = \partial_R\left(\phi_n(R)\right)\cdot\int_0^R\left[\phi_n(s)\right]^{-1}\calD_n\phi_{2,[\mu,\infty]}(s)\,ds,
\]
we write 
\begin{equation}\label{eq:highfreqneedpartialint1}\begin{split}
&\chi_{\xi\simeq\mu}\left\langle \phi_{n_{3}}(R;\xi),\,\chi_{R\lesssim\tau}\partial_R\phi_{1,\leq \mu}\cdot\partial_R\phi_{2,[\mu,\infty]}\right\rangle_{L^2_{R\,dR}}\\
& = \mu^{-1}\chi_{\xi\simeq\mu}\left\langle \phi_{n_{3}}(R;\xi),\,H_{n_3}\left[\chi_{R\lesssim\tau}\partial_R\phi_{1,\leq \mu}\cdot\partial_R\phi_{2,[\mu,\infty]}\right]\right\rangle_{L^2_{R\,dR}}\\
& =  \mu^{-1}\chi_{\xi\simeq\mu}\left\langle \phi_{n_{3}}(R;\xi),\,H_{n_3}\left[\chi_{R\lesssim\tau}\partial_R\phi_{1,\leq \mu}\cdot \partial_R\left(\phi_n(R)\right)\right]\cdot\left(\int_0^R\left[\phi_n(s)\right]^{-1}\calD_n\phi_{2,[\mu,\infty]}(s)\,ds\right)\right\rangle_{L^2_{R\,dR}}\\
& +  \mu^{-1}\chi_{\xi\simeq\mu}\left\langle \phi_{n_{3}}(R;\xi),\,\partial_R\left(\chi_{R\lesssim\tau}\partial_R\phi_{1,\leq \mu}\cdot \partial_R\phi_n(R)\left[\phi_n(R)\right]^{-1}\cdot \calD_n\phi_{2,[\mu,\infty]}(R)\right)\right\rangle_{L^2_{R\,dR}}
\end{split}\end{equation}
The second term on the right can be bounded directly by using Plancherel's Theorem for the cases $n=0,\pm1$, as well as Leibinz' rule to expand things out more. Specifically, we have 
\begin{align*}
&\left\|\partial_R\left(\chi_{R\lesssim\tau}\partial_R\phi_{1,\leq \mu}\cdot \partial_R\left(\phi_n(R)\right)\left[\phi_n(R)\right]^{-1}\right)
\cdot \calD_n\phi_{2,[\mu,\infty]}(R)\right\|_{L^2_{R\,dR}}\\
&\lesssim \left\|\partial_R\left(\chi_{R\lesssim\tau}\partial_R\phi_{1,\leq \mu}\cdot \partial_R\left(\phi_n(R)\right)\left[\phi_n(R)\right]^{-1}\right)\right\|_{L^\infty_{R\,dR}}\cdot \left\| \calD_n\phi_{2,[\mu,\infty]}(R)\right\|_{L^2_{R\,dR}}\\
&\lesssim \hbar_1^{-4-\delta}\left\|\phi_1\right\|_{\tilde{S}_0^{(n_1)}}\cdot \left(\sum_{\lambda\geq \mu}\lambda^{-2-\frac{\delta}{2}}\left\|\phi_{2,\lambda}\right\|_{\tilde{S}_0^{(n_2)}(\xi\simeq\lambda)}\right),
\end{align*}
where we have also taken advantage of Lemma~\ref{lem:derLinfty}. Further, we have 
\begin{align*}
&\left\|\chi_{R\lesssim\tau}\partial_R\phi_{1,\leq \mu}\cdot \partial_R\left(\phi_n(R)\right)\left[\phi_n(R)\right]^{-1}
\cdot \partial_R\left(\calD_n\phi_{2,[\mu,\infty]}(R)\right)\right\|_{L^2_{R\,dR}}\\
&\lesssim \left\|\chi_{R\lesssim\tau}\partial_R\phi_{1,\leq \mu}\cdot \partial_R\left(\phi_n(R)\right)\left[\phi_n(R)\right]^{-1}\right\|_{L^\infty_{R\,dR}}\cdot \left\| \partial_R\left(\calD_n\phi_{2,[\mu,\infty]}(R)\right)\right\|_{L^2_{R\,dR}}\\
&\lesssim \hbar_1^{-4-\delta}\left\|\phi_1\right\|_{\tilde{S}_0^{(n_1)}}\cdot \left(\sum_{\lambda\geq \mu}\lambda^{-\frac32-\frac{\delta}{2}}\left\|\phi_{2,\lambda}\right\|_{\tilde{S}_0^{(n_2)}(\xi\simeq\lambda)}\right)
\end{align*}
Keeping in mind orthogonality, we can the  estimate the contribution of the second term on the right in \eqref{eq:highfreqneedpartialint1} to the $\|\cdot\|_{\Shth_{1}}$-norm by 
\begin{align*}
&\left( \sum_{\mu\gtrsim \hbar_3^{-2}}\left[\hbar_3(\mu\hbar_3^2)^{2+\frac{\delta}{2}}\cdot\mu^{-1}\hbar_1^{-4-\delta}\left\|\phi_1\right\|_{\tilde{S}_0^{(n_1)}}\cdot \left(\sum_{\lambda\geq \mu}\lambda^{-\frac32-\frac{\delta}{2}}\left\|\phi_{2,\lambda}\right\|_{\tilde{S}_0^{(n_2)}(\xi\sim\lambda)}\right)\right]^2\right)^{\frac12}\\
&\lesssim \hbar_3^{\frac12}\left\|\phi_1\right\|_{\tilde{S}_0^{(n_1)}}\cdot \left(\sum_{\mu\gtrsim \hbar_3^{-2}}\sum_{\lambda\geq\mu}\left(\frac{\mu}{\lambda}\right)^{\frac12}\left\|\phi_{2,\lambda}\right\|_{\tilde{S}_0^{(n_2)}(\xi\simeq\lambda)}^2\right)^{\frac12}\\
&\lesssim  \hbar_3^{\frac12}\left\|\phi_1\right\|_{\tilde{S}_0^{(n_1)}}\cdot \left\|\phi_{2}\right\|_{\tilde{S}_0^{(n_2)}},
\end{align*}
which is better than the bound we are striving to establish. \\

As for the contribution of the first term in \eqref{eq:highfreqneedpartialint1}, this requires another application of integration by parts. In fact, we can write 
\begin{align*}
&\mu^{-1}\chi_{\xi\simeq\mu}\left\langle \phi_{n_{3}}(R;\xi),\,H_{n_3}\left[\chi_{R\lesssim\tau}\partial_R\phi_{1,\leq \mu}\cdot \partial_R\left(\phi_n(R)\right)\right]\cdot\left(\int_0^R\left[\phi_n(s)\right]^{-1}\calD_n\phi_{2,[\mu,\infty]}(s)\,ds\right)\right\rangle_{L^2_{R\,dR}}\\
& = \mu^{-2}\chi_{\xi\simeq\mu}\left\langle \phi_{n_{3}}(R;\xi),\,H_{n_3}^2\left[\chi_{R\lesssim\tau}\partial_R\phi_{1,\leq \mu}\cdot \partial_R\left(\phi_n(R)\right)\right]\cdot\left(\int_0^R\left[\phi_n(s)\right]^{-1}\calD_n\phi_{2,[\mu,\infty]}(s)\,ds\right)\right\rangle_{L^2_{R\,dR}}\\
& +  \mu^{-2}\chi_{\xi\simeq\mu}\left\langle \phi_{n_{3}}(R;\xi),\,\partial_R\left(H_{n_3}\left[\chi_{R\lesssim\tau}\partial_R\phi_{1,\leq \mu}\cdot \partial_R\left(\phi_n(R)\right)\right]\left[\phi_n(R)\right]^{-1}\cdot \calD_n\phi_{2,[\mu,\infty]}(R)\right)\right\rangle_{L^2_{R\,dR}}.\\
\end{align*}
Here both terms on the right can be bounded directly. The second term is analogous (with the same exponent $\hbar_{1}^{-4-\delta}$) to the second term in \eqref{eq:highfreqneedpartialint1} and hence omitted, while for the first term we use (with $\alpha_{-1} = 1,\alpha_0 = 2,\alpha_1 = 3$).
\begin{align*}
&\left\|\langle R\rangle^{\alpha_n}H_{n_3}^2\left[\chi_{R\lesssim\tau}\partial_R\phi_{1,\leq \mu}\cdot \partial_R\left(\phi_n(R)\right)\right]\right\|_{L^2_{R\,dR}}\lesssim \hbar_1^{-5-\delta}\left\|\phi_1\right\|_{\tilde{S}_0^{(n_1)}},\\
&\left\|\langle R\rangle^{-\alpha_n}\int_0^R\left[\phi_n(s)\right]^{-1}\calD_n\phi_{2,[\mu,\infty]}(s)\,ds\right\|_{L^\infty_{R\,dR}}\lesssim \mu^{-1}\left\|\phi_2\right\|_{\tilde{S}_0^{(n_2)}}, 
\end{align*}
where for the first bound we again take advantage of Proposition~\ref{prop:derivative1}  and Lemma~\ref{lem:derLinfty}. These bounds imply 
\begin{align*}
&\left\| \mu^{-2}\chi_{\xi\simeq\mu}\left\langle \phi_{n_{3}}(R;\xi),\,H_{n_3}^2\left[\chi_{R\lesssim\tau}\partial_R\phi_{1,\leq \mu}\cdot \partial_R\left(\phi_n(R)\right)\right]\cdot\left(\int_0^R\left[\phi_n(s)\right]^{-1}\calD_n\phi_{2,[\mu,\infty]}(s)\,ds\right)\right\rangle_{L^2_{R\,dR}}\right\|_{\Shth_{1}}\\
&\lesssim \hbar_3^{5+\delta}\mu^{2+\frac{\delta}{2}}\\
&\cdot \left\| \mu^{-2}\chi_{\xi\simeq\mu}\left\langle \phi_{n_{3}}(R;\xi),\,H_{n_3}^2\left[\chi_{R\lesssim\tau}\partial_R\phi_{1,\leq \mu}\cdot \partial_R\left(\phi_n(R)\right)\right] \cdot\left(\int_0^R
\left[\phi_n(s)\right]^{-1}\calD_n\phi_{2,[\mu,\infty]}(s)\,ds\right)\right\rangle_{L^2_{R\,dR}}\right\|_{L^2_{d\xi}}\\
&\lesssim \hbar_3^{\delta}\mu^{-(1-\frac{\delta}{2})}\cdot  \big\|\phi_1\big\|_{\tilde{S}_0^{(n_1)}}\cdot  \big\|\phi_2\big\|_{\tilde{S}_0^{(n_2)}},
\end{align*}
and this can be summed over dyadic $\mu\gtrsim \hbar_3^{-2}$. This finally concludes the estimate for case {\it{(1.b)}}.
\\

{\it{(2): $\xi\geq \max\{\xi_1,\xi_2\}$. }} This is again completely analogous to {\it{(2)}} in the proof of Proposition~\ref{prop:bilin1}, taking advantage of Proposition~\ref{prop:exceptionalnderivative1} as well as repeated integrations by parts as in the preceding case, and hence omitted here. 
\\

{\it{Proof of \eqref{eq:productsmalln1geq2alln2}}}. We start with the case when the output frequency $\xi$ dominates the two input frequencies $\xi_{1,2}$, which corresponded to the case {\it{(2)}} in the preceding argument. 
\\
In order to define $P_3$, we formally expand the product $\chi_{R\lesssim\tau}\partial_R\phi_1\cdot\partial_R\phi_2$ in a Taylor series around $R = 0$ and stop at order three terms. For this we expand (observe that we do not include a cutoff here)
\[
\phi_1(R) = c\int_0^\infty \left(R\xi^{\frac12}\right)^{|n_1|\mp 1}\left(1+ \sum_{j\geq 1}\phi_j(R^2)\cdot \left(R^2\xi\right)^j\right)\xb_1(\xi)\tilde{\rho}_{n_{1}}(\xi)\,d\xi,
\]
provided $\phi_1$ is an angular momentum $n_1$ function with $|n_1|\geq 2$, while if $|n_1| = 1,\,0$, we use an analogous formula at the level of $\mathcal{D}_j\phi_1(R)$, $j = 0,\,\pm 1$, which gets then inserted into the formula giving  $\phi_1(R)$ in terms of $\left(c_j,\,\mathcal{D}_j\phi_1(R)\right)$. Keeping track only of terms up to order $R^3$, we easily infer the bound \eqref{eq:Taylorbound1}. Now write the bad term corresponding to the case {\it{(2)}} from before as 
\begin{equation}\label{eq:Taylorcorrectioneffect1}\begin{split}
&\sum_{\lambda}\chi_{\xi\simeq\lambda}\left\langle \phi_{n_{3}}(R;\xi),\,\chi_{R\lesssim\tau}\partial_R\phi_{1,\ll\lambda}\cdot\partial_R\phi_{2,\ll\lambda} - \chi_{R\lesssim 1}P_3(R)\right\rangle_{L^2_{R\,dR}}\\
& =  \sum_{\lambda}\chi_{\xi\simeq\lambda}\left\langle\chi_{R\xi^{\frac12}\gtrsim 1}\phi_{n_{3}}(R;\xi),\,\chi_{R\lesssim\tau}\partial_R\phi_{1,\ll \lambda}\cdot\partial_R\phi_{2,\ll\lambda} - \chi_{R\lesssim 1}P_3(R)\right\rangle_{L^2_{R\,dR}}\\
& +  \sum_{\lambda}\chi_{\xi\simeq\lambda}\left\langle\chi_{R\xi^{\frac12}\lesssim 1}\phi_{n_{3}}(R;\xi),\,\chi_{R\lesssim\tau}\partial_R\phi_{1,\ll \lambda}\cdot\partial_R\phi_{2,\ll\lambda} - \chi_{R\lesssim 1}P_3(R)\right\rangle_{L^2_{R\,dR}}\\
\end{split}\end{equation}
Observe that we no longer keep careful track of the $\hbar_3$-dependence, since $n_3 = O(1)$ now, and similarly for all other angular momenta. 
\\
To deal with the first term on the right of \eqref{eq:Taylorcorrectioneffect1}, observe that we can replace the cutoff $\chi_{R\lesssim\tau}$ by $\chi_{R\lesssim 1}$, since including a cutoff $\chi_{1\lesssim R\lesssim\tau}$ the resulting term can be handled like {\it{(2)}} in the proof of Proposition~\ref{prop:bilin1} via integration by parts without incurring problematic boundary terms. Assuming for simplicity that both factors $\phi_{1,2}$ are angular momentum $n_j$ functions with $|n_j|\geq 2$, $j = 1,2$, we can then schematically write 
\begin{align*}
&\chi_{R\lesssim 1}\partial_R\phi_{1,\ll\lambda}\cdot\partial_R\phi_{2,\ll\lambda} - \chi_{R\lesssim 1}P_3(R) = \sum_{j=1,2}X_j -\chi_{R\lesssim 1} \tilde{P}_3(R)\\
&X_1 = \chi_{R\lesssim 1}\sum_{i+j\geq 4}\left(\int_0^\infty \chi_{R\xi_1^{\frac12}\lesssim 1}(R\xi_1^{\frac12})^i\left(1+\sum_{k\geq 1}\phi_k(R)(R\xi_1^{\frac12})^k\right)\xi_1^{\frac12}\chi_{\xi_1\ll\lambda}\xb_1(\xi_1)\tilde{\rho}_{n_{1}}(\xi_1)\,d\xi_1\right)\\
&\hspace{5cm}\cdot\left(\int_0^\infty\chi_{R\xi_2^{\frac12}\lesssim 1} (R\xi_2^{\frac12})^j\left(1+\sum_{l\geq 1}\phi_l(R)(R\xi_2^{\frac12})^l\right)\xi_2^{\frac12}\chi_{\xi_2\ll\lambda}\xb_2(\xi_2)\tilde{\rho}_{n_{2}}(\xi_2)\,d\xi_2\right)\\
&X_2 =  \chi_{R\lesssim 1}\left(\int_0^\infty\chi_{R\xi_1^{\frac12}\gtrsim 1}\partial_R\phi_{n_{1}}(R;\xi_1)\chi_{\xi_1\ll\lambda}\xb_1(\xi_1)\tilde{\rho}_{n_1}(\xi_1)\,d\xi_1\right)\\&\hspace{7cm}\cdot \left(\int_0^\infty\partial_R\phi_{n_{2}}(R;\xi_2)\chi_{\xi_2\ll\lambda}\xb_2(\xi_2)\tilde{\rho}_{n_2}(\xi_2)\,d\xi_2\right)\\
& +  \chi_{R\lesssim 1}\left(\int_0^\infty\chi_{R\xi_1^{\frac12}\lesssim 1}\partial_R\phi_{n_{1}}(R;\xi_1)\chi_{\xi_1\ll\lambda}\xb_1(\xi_1)\tilde{\rho}_{n_{1}}(\xi_1)\,d\xi_1\right)\\&\hspace{7cm}\cdot \left(\int_0^\infty\chi_{R\xi_2^{\frac12}\gtrsim 1}\partial_R\phi_{n_{2}}(R;\xi_2)\chi_{\xi_2\ll\lambda}\xb_2(\xi_2)\tilde{\rho}_{n_{2}}(\xi_2)\,d\xi_2\right)\\
\end{align*}
where the functions $\phi_{l,k}(R)$ in the definition of $X_1$ satisfy bounds stated in Lemma 2.34. Furthermore, we can write schematically 
\begin{align*}
\tilde{P}_3(R) &= \sum_{i+j\leq 3}C_{i,j}\left(\int_0^\infty \chi_{R\xi_1^{\frac12}\gtrsim 1}\left(R\xi_1^{\frac12}\right)^i\xi_1^{\frac12}\xb_1(\xi_1)\tilde{\rho}_{n_{1}}(\xi_1)\,d\xi_1\right)\cdot \left(\int_0^\infty \left(R\xi_2^{\frac12}\right)^j\xi_2^{\frac12}\xb_2(\xi_2)\tilde{\rho}_{n_{2}}(\xi_2)\,d\xi_2\right)\\
&+ \sum_{i+j\leq 3}C_{i,j}\left(\int_0^\infty \chi_{R\xi_1^{\frac12}\lesssim 1}\left(R\xi_1^{\frac12}\right)^i\xi_1^{\frac12}\xb_1(\xi_1)\tilde{\rho}_{n_{1}}(\xi_1)\,d\xi_1\right)\cdot \left(\int_0^\infty\chi_{R\xi_2^{\frac12}\gtrsim 1} \left(R\xi_2^{\frac12}\right)^j\xi_2^{\frac12}\xb_2(\xi_2)\tilde{\rho}_{n_{2}}(\xi_2)\,d\xi_2\right)\\
&+ \ldots,
\end{align*}
and where the terms denoted by $\ldots$ refer to similar products for $R\xi_{1}^{\frac12}\lesssim 1, R\xi^{\frac12}_{2}\lesssim 1$ where at least one cutoff $\chi_{\xi_{1,2}\gtrsim\lambda}$ is included into one of the factors.
\\
Then we handle their contribution to the first term on the right of \eqref{eq:Taylorcorrectioneffect1} as follows: for the contribution of $X_1$, i.e., the term 
\begin{align*}
\left\|\sum_{\lambda}\chi_{\xi\simeq\lambda}\left\langle \chi_{R\xi^{\frac12}\gtrsim 1}\chi_{R\simeq\kappa}\phi_{n_{3}}(R;\xi),\,X_1\right\rangle_{L^2_{R\,dR}}\right\|_{\Shth_{1}}, 
\end{align*}
localize $R$ to dyadic size $\kappa\gtrsim \lambda^{-\frac12}$. Exploiting the oscillatory nature of $\chi_{R\xi^{\frac12}\gtrsim 1}\phi_{n_{3}}(R;\xi)$ from Lemma \ref{lem:psin R2xi large}, Proposition \ref{prop: FB match} and their analogues for negative $n$, and performing integration by parts with respect to $R$, we gain $\left(\frac{1}{\lambda^{\frac12}\kappa}\right)^N$. Then localizing the frequencies $\xi_{1,2}$ to dyadic values $\mu_{1,2}\ll\lambda$, and calling $X_{1,\mu_{1,2}}$ the corresponding contribution to $X_1$, we have  
\begin{align*}
&\left\|\chi_{\xi\simeq\lambda}\left\langle \chi_{R\xi^{\frac12}\gtrsim 1}\chi_{R\simeq\kappa}\phi_{n_{3}}(R;\xi),\,X_{1,\mu_{1,2}}\right\rangle_{L^2_{R\,dR}}\right\|_{\Shth_{1}}\\
&\lesssim \lambda^{\frac12-\frac{\delta}{2}}\langle\lambda\rangle^{\frac32+\delta}\left\|\chi_{\xi\simeq\lambda}\left\langle \chi_{R\xi^{\frac12}\gtrsim 1}\chi_{R\simeq\kappa}\phi_{n_{3}}(R;\xi),\,X_{1,\mu_{1,2}}\right\rangle_{L^2_{R\,dR}}\right\|_{L^2_{d\xi}}\\
&\lesssim \lambda^{1-\frac{\delta}{2}}\langle\lambda\rangle^{\frac32+\delta}\left\|\chi_{\xi\simeq\lambda}\left\langle \chi_{R\xi^{\frac12}\gtrsim 1}\chi_{R\simeq\kappa}\phi_{n_{3}}(R;\xi),\,X_{1,\mu_{1,2}}\right\rangle_{L^2_{R\,dR}}\right\|_{L^\infty_{d\xi}},\\
\end{align*}
and further for fixed $\xi\simeq\lambda$ we have, by Lemma \ref{lem:psin R2xi large} and the profile of $X_{1}$,
\begin{align*}
&\left|\left\langle \chi_{R\xi^{\frac12}\gtrsim 1}\chi_{R\simeq\kappa}\phi_{n_{3}}(R;\xi),\,X_{1,\mu_{1,2}}\right\rangle_{L^2_{R\,dR}}\right|\\
&\lesssim \frac{\kappa^6}{\left(\lambda^{\frac12}\kappa\right)^N}\cdot\left(\mu_1^{\frac12-\frac{\delta}{2}}\frac{\mu_2^{\frac{\delta}{2}}}{\langle\mu_2\rangle^{\frac32+\delta}} + \mu_2^{\frac12-\frac{\delta}{2}}\frac{\mu_1^{\frac{\delta}{2}}}{\langle\mu_1\rangle^{\frac32+\delta}}\right)\prod_{j=1,2}\left\|\xb_j\right\|_{\Shj_0(\xi_j\simeq\mu_j)}
\end{align*}
Summing over $\kappa\gtrsim\lambda^{-\frac12}$ we obtain 
\begin{align*}
&\sum_{\kappa\gtrsim\lambda^{-\frac12}} \frac{\kappa^6}{\left(\lambda^{\frac12}\kappa\right)^N}\left(\mu_1^{\frac12-\frac{\delta}{2}}\frac{\mu_2^{\frac{\delta}{2}}}{\langle\mu_2\rangle^{\frac32+\delta}} + \mu_2^{\frac12-\frac{\delta}{2}}\frac{\mu_1^{\frac{\delta}{2}}}{\langle\mu_1\rangle^{\frac32+\delta}}\right)\cdot\prod_{j=1,2}\left\|x_j\right\|_{\Shj_0(\xi_j\simeq\mu_j)}\\&\lesssim 
\frac{\mu_1^{\frac12-\frac{\delta}{2}}\frac{\mu_2^{\frac{\delta}{2}}}{\langle\mu_2\rangle^{\frac32+\delta}} + \mu_2^{\frac12-\frac{\delta}{2}}\frac{\mu_1^{\frac{\delta}{2}}}{\langle\mu_1\rangle^{\frac32+\delta}}}{\lambda^3}\prod_{j=1,2}\left\|\xb_j\right\|_{\Shj_0(\xi_j\simeq\mu_j)}.
\end{align*}
Combining with the inequality further above we obtain 
\begin{align*}
\left\|\chi_{\xi\simeq\lambda}\left\langle \chi_{R\xi^{\frac12}\gtrsim 1}\phi_{n_{3}}(R;\xi),\,X_{1,\mu_{1,2}}\right\rangle_{L^2_{R\,dR}}\right\|_{\Shth_1}
\lesssim\frac{\mu_1^{\frac12-\frac{\delta}{2}}\frac{\mu_2^{\frac{\delta}{2}}}{\langle\mu_2\rangle^{\delta}} + \mu_2^{\frac12-\frac{\delta}{2}}\frac{\mu_1^{\frac{\delta}{2}}}{\langle\mu_1\rangle^{\delta}}}{\lambda^{\frac12-\frac{\delta}{2}}}\prod_{j=1,2}\left\|x_j\right\|_{\Shj_0(\xi_j\simeq\mu_j)}
\end{align*}
provided $\lambda\gtrsim 1$. 
Finally, exploiting orthogonality and Cauchy-Schwarz inequality, we infer that 
\begin{align*}
\left\|\left\langle \chi_{R\xi^{\frac12}\gtrsim 1}\phi_{n_{3}}(R;\xi),\,X_{1,\mu_{1,2}}\right\rangle_{L^2_{R\,dR}}\right\|_{\Shth_1(\xi\gtrsim 1)}\lesssim \left(\sum_{\lambda}\left\|\chi_{\xi\simeq\lambda}\langle \chi_{R\xi^{\frac12}\gtrsim 1}\phi_{n_{3}}(R;\xi),\,X_{1,\mu_{1,2}}\rangle_{L^2_{R\,dR}}\right\|_{\Shth_1}^2\right)^{\frac12}\lesssim \prod_{j=1}^2\left\|\phi_j\right\|_{\tilde{S}^{(n_j)}}.
\end{align*}
The estimate for low output frequencies $\xi\lesssim 1$ is much easier due to the restriction on $R$ in the definition of $X_1$ and omitted here. 
\\
The contribution of the term $X_2$ to the first term on the right hand side in \eqref{eq:Taylorcorrectioneffect1} is handled similarly, except that now one has to use the oscillatory expansion for $\chi_{R\xi_1^{\frac12}\gtrsim 1}\phi_{n_{1}}(R;\xi_1)$ and the other oscillatory terms and combine the phases with the one of $\chi_{R\xi^{\frac12}\gtrsim 1}\phi_{n_{3}}(R;\xi)$ before performing integration by parts. We omit the similar details.  
\\
In order to complete the estimate for the first term on the right hand side of \eqref{eq:Taylorcorrectioneffect1}, it thus suffices to deal with the contribution of $\chi_{R\lesssim 1} \tilde{P}_3(R)$. We observe right away that the restrictions $R\xi^{\frac12}\gtrsim 1,\,R\lesssim 1$ imply $\xi\gtrsim 1$. Consider the first term in the definition of $\tilde{P}_3(R)$, and call this $\tilde{P}_{3}^{(1)}(R)$. Further localise the frequencies of the two factors to dyadic size $\mu_{1,2}$, respectively, resulting in $\tilde{P}_{3,\mu_{1,2}}^{(1)}(R)$. Then we bound for $\kappa\gtrsim \mu_1^{-\frac12}$
\begin{align*}
&\left|\chi_{\xi\simeq\lambda}\left\langle  \chi_{R\xi^{\frac12}\gtrsim 1}\phi_{n_{3}}(R;\xi),\,\chi_{R\simeq\kappa}\chi_{R\lesssim 1} \tilde{P}_{3,\mu_{1,2}}^{(1)}(R)\right\rangle_{L^2_{R\,dR}}\right|\\
&\lesssim \sum_{i+j\leq 3}\frac{\kappa^{2+i+j}}{(\lambda^{\frac12}\kappa)^N}\cdot\left|\int_0^\infty \chi_{\xi_1\simeq\mu_1}\xi_1^{\frac{1+i}{2}}\xb_1(\xi_1)\tilde{\rho}_{n_{1}}(\xi_1)\,d\xi_1\right|\cdot \left|\int_0^\infty \chi_{\xi_2\simeq\mu_2}\xi_2^{\frac{1+j}{2}}\xb_2(\xi_2)\tilde{\rho}_{n_{2}}(\xi_2)\,d\xi_2\right|.
\end{align*}
The integrals are easily bounded by 
\begin{align*}
&\left|\int_0^\infty \chi_{\xi_1\simeq\mu_1}\xi_1^{\frac{1+i}{2}}\xb_1(\xi_1)\tilde{\rho}_{n_{1}}(\xi_1)\,d\xi_1\right|\lesssim \frac{\mu_1^{\frac{\delta}{2}}}{\langle\mu_1\rangle^{\delta}}\cdot\langle\mu_1\rangle^{-\frac32+\frac{i}{2}}\cdot\left\|\xb_1\right\|_{\Sho_{0}(\xi_1\simeq\mu_1)},\\
&\left|\int_0^\infty \chi_{\xi_2\simeq\mu_1}\xi_2^{\frac{1+j}{2}}\xb_2(\xi_2)\tilde{\rho}_{n_{2}}(\xi_2)\,d\xi_2\right|\lesssim \frac{\mu_2^{\frac{\delta}{2}}}{\langle\mu_2\rangle^{\delta}}\cdot\langle\mu_2\rangle^{-\frac32+\frac{j}{2}}\cdot\left\|\xb_2\right\|_{\Sho_{0}(\xi_2\simeq\mu_2)}.
\end{align*}
We conclude that 
\begin{align*}
&\sum_{\kappa\gtrsim \lambda^{-\frac12}}\left|\chi_{\xi\simeq\lambda}\left\langle  \chi_{R\xi^{\frac12}\gtrsim 1}\phi_{n_{3}}(R;\xi),\,\chi_{R\simeq\kappa}\chi_{R\lesssim 1} \tilde{P}_{3,\mu_{1,2}}^{(1)}(R)\right\rangle_{L^2_{R\,dR}}\right|\\
&\lesssim \lambda^{-1-\frac{3}{2}}\cdot\left(\frac{\mu_1}{\lambda}\right)^{N^{\prime}}\cdot  \frac{\mu_1^{\frac{\delta}{2}}}{\langle\mu_1\rangle^{\delta}}\cdot\langle\mu_1\rangle^{-\frac32+\frac{i}{2}}\frac{\mu_2^{\frac{\delta}{2}}}{\langle\mu_2\rangle^{\delta}}\cdot\langle\mu_2\rangle^{-\frac32+\frac{j}{2}}\cdot \prod_{k=1,2}\left\|\xb_k\right\|_{\Shk_{0}(\xi_k\simeq\mu_k)}.
\end{align*}
Note that the factor $\left(\frac{\mu_{1}}{\lambda}\right)^{N'}$ is also obtained using integration by parts.
Finally, we infer 
\begin{align*}
&\left\|\chi_{\xi\simeq\lambda}\left\langle  \chi_{R\xi^{\frac12}\gtrsim 1}\phi_{n_{3}}(R;\xi),\,\chi_{R\lesssim 1} \tilde{P}_{3,\mu_{1,2}}^{(1)}(R)\right\rangle_{L^2_{R\,dR}}\right\|_{\Shth_{1}}\\
&\lesssim \lambda^{\frac12-\frac{\delta}{2}}\langle\lambda\rangle^{\frac32 + \delta}\left\|\chi_{\xi\simeq\lambda}\left\langle  \chi_{R\xi^{\frac12}\gtrsim 1}\phi_{n_{3}}(R;\xi),\,\chi_{R\lesssim 1} \tilde{P}_{3,\mu_{1,2}}^{(1)}(R)\right\rangle_{L^2_{R\,dR}}\right\|_{L^2_{d\xi}}\\
&\lesssim  \lambda^{1-\frac{\delta}{2}}\langle\lambda\rangle^{\frac32 + \delta}\left\|\chi_{\xi\simeq\lambda}\left\langle  \chi_{R\xi^{\frac12}\gtrsim 1}\phi_{n_{3}}(R;\xi),\,\chi_{R\lesssim 1} \tilde{P}_{3,\mu_{1,2}}^{(1)}(R)\right\rangle_{L^2_{R\,dR}}\right\|_{L^\infty_{d\xi}}\\
\end{align*}
whence using the preceding bound we get 
\begin{align*}
&\left\|\chi_{\xi\simeq\lambda}\langle  \chi_{R\xi^{\frac12}\gtrsim 1}\phi_{n_{3}}(R;\xi),\,\chi_{R\lesssim 1} \tilde{P}_{3,\mu_{1,2}}^{(1)}(R)\rangle_{L^2_{R\,dR}}\right\|_{\Shth_{1}}\\
&\lesssim \left(\frac{\mu_1}{\lambda}\right)^{N^{\prime}}\cdot \frac{(\lambda\mu_1)^{\frac{\delta}{2}}}{\langle\mu_1\rangle^{\delta}}\cdot \frac{\mu_2^{\frac{\delta}{2}}}{\langle\mu_2\rangle^{\delta}}\cdot  \prod_{k=1,2}\left\|\xb_k\right\|_{\Shk_{0}(\xi_k\simeq\mu_k)}.
\end{align*}
But then exploiting orthogonality and the Cauchy-Schwarz inequality as usual we infer (here $\tilde{P}_{3,\ll\lambda,\ll\lambda}^{(1)}(R)$ is the term arising from $\tilde{P}_3^{(1)}$ after restricting both factors to frequency $\ll \lambda$) 
\begin{align*}
&\left\|\sum_{\lambda}\left\langle \chi_{\xi\simeq\lambda}\chi_{R\xi^{\frac12}\gtrsim 1}\phi_{n_{3}}(R;\xi),\,\chi_{R\lesssim 1} \tilde{P}_{3,\ll\lambda,\ll\lambda}^{(1)}(R)\right\rangle_{L^2_{R\,dR}}\right\|_{S_1^{(\hbar_3)}}\\
&\lesssim \left(\sum_{\lambda}\left[\sum_{\mu_{1,2}\ll\lambda}\left\|\left\langle  \chi_{R\xi^{\frac12}\gtrsim 1}\phi_{n_{3}}(R;\xi),\,\chi_{R\lesssim 1} \tilde{P}_{3,\mu_{1,2}}^{(1)}(R)\right\rangle_{L^2_{R\,dR}}\right\|_{\Shth_{1}(\xi\simeq\lambda)}\right]^2\right)^{\frac12}\\
&\lesssim  \left(\sum_{\lambda}\sum_{\mu_{1,2}\ll\lambda} \left(\frac{\mu_1}{\lambda}\right)^{N^{\prime}}\cdot \frac{(\lambda\mu_1)^{\frac{\delta}{2}}}{\langle\mu_1\rangle^{\delta}}\cdot \frac{\mu_2^{\frac{\delta}{2}}}{\langle\mu_2\rangle^{\delta}}\cdot  \prod_{k=1,2}\left\|\xb_k\right\|_{\Shk_{0}(\xi_k\simeq\mu_k)}^2\right)^{\frac12}\\
&\lesssim  \prod_{k=1,2}\left\|\xb_k\right\|_{\Shk_{0}}
\end{align*}
This still leaves the contribution of 
\[
\chi_{R\lesssim 1}[\tilde{P}_3^{(1)}(R) -  \tilde{P}_{3,\ll\lambda,\ll\lambda}^{(1)}(R)]
\]
to be bounded, which will be done like below, when treating the contribution of the final term constituting $\tilde{P}_3(R)$. The second term in the decomposition of $\tilde{P}_3(R)$ is handled analogously to the first, and so we only need to bound the contribution of the last term, $\tilde{P}_3^{(3)}(R)$ to finish the bound for the first term on the right of \eqref{eq:Taylorcorrectioneffect1}. 
This contribution is then given by\footnote{In addition there is a similar term with the roles of $x_1, x_2$ interchanged.} 
\begin{align*}
&\sum_{i+j\leq 3}C_{i,j}\sum_{\lambda}\Big\langle \chi_{\xi\simeq\lambda}\chi_{R\xi^{\frac12}\gtrsim 1}\phi_{n_{3}}(R;\xi),\,\chi_{R\lesssim 1}\left(\int_0^\infty\left(R\xi_1^{\frac12}\right)^i\chi_{\xi_1\gtrsim\lambda}\xi_1^{\frac12}\xb_1(\xi_1)\tilde{\rho}_{n_{1}}(\xi_1)\,d\xi_1\right)\\&\hspace{9cm}\cdot \left(\int_0^\infty \left(R\xi_2^{\frac12}\right)^j\xi_2^{\frac12}\xb_2(\xi_2)\tilde{\rho}_{n_{2}}(\xi_2)\,d\xi_2\right)\Big\rangle_{L^2_{R\,dR}}.
\end{align*}
Then we use the integral bounds (for $i+j\leq 3$)
\begin{align*}
&\left|\int_0^\infty\xi_1^{\frac{i+1}{2}}\chi_{\xi_1\gtrsim\lambda}\xb_1(\xi_1)\tilde{\rho}_{n_{1}}(\xi_1)\,d\xi_1\right|\lesssim \sum_{\mu\gtrsim\lambda}\mu^{-\frac32-\frac{\delta}{2}}\mu^{\frac{i}{2}}\cdot\left\|\xb_1\right\|_{\Sho_{0}(\xi\simeq\mu)},\\
&\left|\int_0^\infty\xi_2^{\frac{j+1}{2}}\xb_2(\xi_2)\tilde{\rho}_{n_{1}}(\xi_2)\,d\xi_2\right|\lesssim \left\|\xb_2\right\|_{\Sht_{0}},
\end{align*}
and so localizing $R$ to $\kappa\lesssim 1$, $\kappa\gtrsim \lambda^{-\frac12}$, we infer after repeated integrations by parts (denoting by $A_{\lambda}^{i,j}(R)$ the product of the two integrals in the long expression above)
\begin{align*}
&\left|\langle \chi_{\xi\simeq\lambda}\chi_{R\xi^{\frac12}\gtrsim 1}\phi_{n_{3}}(R;\xi),\,\chi_{R\simeq\kappa}A_{\lambda}^{i,j}(R)\rangle_{L^2_{R\,dR}}\right|\\
&\lesssim \frac{\kappa^{2+i+j}}{(\lambda^{\frac12}\kappa)^N}\cdot\left(\sum_{\mu\gtrsim\lambda}\mu^{-\frac32-\frac{\delta}{2}}\mu^{\frac{i}{2}}\cdot\left\|\xb_1\right\|_{\Sho_{0}(\xi\simeq\mu)}\right)\cdot \left\|\xb_2\right\|_{\Sht_{0}},
\end{align*}
whence summing over $\kappa$
\begin{align*}
&\sum_{\lambda^{-\frac12}\lesssim\kappa\lesssim1}\left|\left\langle \chi_{\xi\simeq\lambda}\chi_{R\xi^{\frac12}\gtrsim 1}\phi_{n_{3}}(R;\xi),\,\chi_{R\simeq\kappa}A_{\lambda}^{i,j}(R)\right\rangle_{L^2_{R\,dR}}\right|\\
&\lesssim \lambda^{-1-\frac{3}{2}}\cdot\left(\sum_{\mu\gtrsim\lambda}\mu^{-\frac32-\frac{\delta}{2}}\mu^{\frac{i}{2}}\cdot\left\|\xb_1\right\|_{\Sho_{0}(\xi\simeq\mu)}\right)\cdot \left\|\xb_2\right\|_{\Sht_0}.
\end{align*}
Finally, we infer 
\begin{align*}
&\left\|\sum_{i+j\leq 3}C_{i,j}\sum_{\lambda}\langle \chi_{\xi\simeq\lambda}\chi_{R\xi^{\frac12}\gtrsim 1}\phi_{n_{3}}(R;\xi),\,\chi_{R\lesssim 1}A_{\lambda}^{i,j}(R)\rangle_{L^2_{R\,dR}}\right\|_{\Shth_1}\\
&\lesssim \sum_{i+j\leq 3}\left(\sum_{\lambda\gtrsim 1}\left\|\langle \chi_{\xi\simeq\lambda}\chi_{R\xi^{\frac12}\gtrsim 1}\phi_{n_{3}}(R;\xi),\,\chi_{R\lesssim 1}A_{\lambda}^{i,j}(R)\rangle_{L^2_{R\,dR}}\right\|_{\Shth_1}^2\right)^{\frac12}\\
&\lesssim \sum_{i+j\leq 3}\left(\sum_{\lambda\gtrsim 1}\left[\lambda^{1-\frac{\delta}{2}}\langle\lambda\rangle^{\frac32+\delta}\left\|\left\langle \chi_{\xi\simeq\lambda}\chi_{R\xi^{\frac12}\gtrsim 1}\phi_{n_{3}}(R;\xi),\,\chi_{R\lesssim 1}A_{\lambda}^{i,j}(R)\right\rangle_{L^2_{R\,dR}}\right\|_{L^\infty_{d\xi}}\right]^2\right)^{\frac12},
\end{align*}
and the preceding expression can be bounded by 
\begin{align*}
&\sum_{i+j\leq 3}\left(\sum_{\lambda\gtrsim 1}\left[\lambda^{1-\frac{\delta}{2}}\langle\lambda\rangle^{\frac32+\delta}\left\|\left\langle \chi_{\xi\simeq\lambda}\chi_{R\xi^{\frac12}\gtrsim 1}\phi_{n_{3}}(R;\xi),\,\chi_{R\lesssim 1}A_{\lambda}^{i,j}(R)\right\rangle_{L^2_{R\,dR}}\right\|_{L^\infty_{d\xi}}\right]^2\right)^{\frac12}\\
&\lesssim \left(\sum_{\lambda\gtrsim 1}\left[\sum_{\mu\gtrsim\lambda}\left(\frac{\lambda}{\mu}\right)^{\delta}\cdot\left\|\xb_1\right\|_{\Sho_{0}(\xi\simeq\mu)}\right]^2\right)^{\frac12}\left\|\xb_2\right\|_{\Sht_0}\lesssim \prod_{j=1,2}\left\|\xb_j\right\|_{\Shj_0}.
\end{align*}
This finally concludes the bound for the first term on the right hand side of \eqref{eq:Taylorcorrectioneffect1}, provided that both factors $\phi_{1,2}$ are angular momentum $n_j$-functions with $|n_j|\geq 2$. The case $|n_j|<2$ for at least one $j$ is handled similarly and omitted.
\\

As for the second term on the right hand side of \eqref{eq:Taylorcorrectioneffect1}, here we have to take advantage of the high degree of vanishing of the expression in the inner product at $R = 0$, as well as the shortness of the interval of integration. To begin with we can reduce $\chi_{R\lesssim\tau}\partial_R\phi_{1,\ll \lambda}\cdot\partial_R\phi_{2,\ll\lambda}$ to $\chi_{R\lesssim 1}\partial_R\phi_{1,\ll \lambda}\cdot\partial_R\phi_{2,\ll\lambda}$. In light of the cutoff $\chi_{R\xi^{\frac12}\lesssim 1}$ in the second term in \eqref{eq:Taylorcorrectioneffect1}, this further localization is automatic if $\xi\gtrsim 1$. Thus consider now 
\[
\sum_{\lambda\lesssim 1}\chi_{\xi\simeq\lambda}\left\langle\chi_{R\xi^{\frac12}\lesssim 1}\phi_{n_{3}}(R;\xi),\,\chi_{1\lesssim R\lesssim\tau}\partial_R\phi_{1,\ll \lambda}\cdot\partial_R\phi_{2,\ll\lambda}\right\rangle_{L^2_{R\,dR}}
\]
We shall again assume for simplicity that both factors $\phi_{1,2}$ are angular momentum $n_j$ functions with $|n_j|\geq 2$, the remaining cases being dealt with similarly. A straightforward sharpening of Lemma~\ref{lem:derLinfty} furnishes the bounds
\begin{align*}
\left|\partial_R\phi_{1,\mu_1}\right|\lesssim \mu_1^{\frac{\delta}{2}}\left\|\xb_1\right\|_{\Sho_{0}(\xi_1\simeq\mu_1)},\quad\left|\partial_R\phi_{2,\ll\lambda}\right|\lesssim \left\|\xb_2\right\|_{\Sht_0}, 
\end{align*}
and so we infer 
\begin{align*}
&\left\|\chi_{\xi\simeq\lambda}\left\langle\chi_{R\xi^{\frac12}\lesssim 1}\phi_{n_{3}}(R;\xi),\,\chi_{1\lesssim R\lesssim\tau}\partial_R\phi_{1,\ll \lambda}\cdot\partial_R\phi_{2,\ll\lambda}\right\rangle_{L^2_{R\,dR}}\right\|_{\Shth_1}\\
&\lesssim \sum_{\mu_1\ll\lambda}\lambda^{1-\frac{\delta}{2}}\cdot\lambda^{-1}\cdot\left\|\partial_R\phi_{1,\mu_1}\cdot\partial_R\phi_{2,\ll\lambda}\right\|_{L^\infty_{R\,dR}}\\
&\lesssim \sum_{\mu_1\ll\lambda}\left(\frac{\mu_{1}}{\lambda}\right)^{\frac{\delta}{2}}\left\|\xb_1\right\|_{\Sho_{0}(\xi_1\simeq\mu_1)}\cdot \left\|\xb_2\right\|_{\Sht_{0}}.
\end{align*}
Square-summing over $\lambda$ yields the desired bound upon applying Cauchy-Schwarz and orthogonality. 
\\
It follows that its suffices to bound 
\begin{align*}
\left\|\sum_{\lambda}\chi_{\xi\simeq\lambda}\left\langle\chi_{R\xi^{\frac12}\lesssim 1}\phi_{n_{3}}(R;\xi),\,\chi_{R\lesssim 1}\left[\partial_R\phi_{1,\ll \lambda}\cdot\partial_R\phi_{2,\ll\lambda} - P_3(R)\right]\right\rangle_{L^2_{R\,dR}}\right\|_{\Shth_1}. 
\end{align*}
Observe that in the inner product the factors $\partial_R\phi_{j,\ll \lambda}$ are automatically in the non-oscillatory regime due to the restrictions on their frequencies and the cutoff $\chi_{R\xi^{\frac12}\lesssim 1}$. Thus we can write 
\begin{align*}
\chi_{\xi\simeq\lambda}\chi_{R\xi^{\frac12}\lesssim 1}\chi_{R\lesssim 1}\left[\partial_R\phi_{1,\ll \lambda}\cdot\partial_R\phi_{2,\ll\lambda} - P_3(R)\right] = X_3 + \tilde{P}_4,
\end{align*}
where we can write schematically
\begin{align*}
X_3 &=  \chi_{R\lesssim 1}\chi_{\xi\simeq\lambda}\chi_{R\xi^{\frac12}\lesssim 1}\sum_{i+j\geq 4}\left(\int_0^\infty (R\xi_1^{\frac12})^i\left(1+\sum_{k\geq 1}\phi_k(R)(R\xi_1^{\frac12})^k\right)\xi_1^{\frac12}\chi_{\xi_1\ll\lambda}\xb_1(\xi_1)\tilde{\rho}_{n_{1}}(\xi_1)\,d\xi_1\right)\\
&\hspace{5cm}\cdot\left(\int_0^\infty (R\xi_2^{\frac12})^j\left(1+\sum_{l\geq 1}\phi_l(R)(R\xi_2^{\frac12})^l\right)\xi_2^{\frac12}\chi_{\xi_2\ll\lambda}\xb_2(\xi_2)\tilde{\rho}_{n_{2}}(\xi_2)\,d\xi_2\right)\\
\end{align*}
as well as 
\begin{align*}
\tilde{P}_4& = \sum_{i+j\leq 3}C_{i,j}\left(\int_0^\infty \chi_{\xi_1\gtrsim \lambda}\left(R\xi_1^{\frac12}\right)^i\xi_1^{\frac12}\xb_1(\xi_1)\tilde{\rho}_{n_{1}}(\xi_1)\,d\xi_1\right)\cdot \left(\int_0^\infty \left(R\xi_2^{\frac12}\right)^j\xi_2^{\frac12}\xb_2(\xi_2)\tilde{\rho}_{n_{2}}(\xi_2)\,d\xi_2\right)\\
&+ \sum_{i+j\leq 3}C_{i,j}\left(\int_0^\infty \chi_{\xi_1\lesssim \lambda}\left(R\xi_1^{\frac12}\right)^i\xi_1^{\frac12}\xb_1(\xi_1)\tilde{\rho}_{n_{1}}(\xi_1)\,d\xi_1\right)\cdot \left(\int_0^\infty\chi_{\xi_2\gtrsim \lambda} \left(R\xi_2^{\frac12}\right)^j\xi_2^{\frac12}\xb_2(\xi_2)\tilde{\rho}_{n_{2}}(\xi_2)\,d\xi_2\right).
\end{align*}
We bound the contributions of $X_3, \tilde{P}_4$ as follows: 
\\
{\it{Contribution of $X_3$}}: First, the low frequencies $\lambda\lesssim 1$ are handled as follows. Localizing the variable $\xb_1(\xi_1)$ to dyadic frequency $\mu_1\ll \lambda$, and calling the resulting expression $X_{3,\mu_1}$, we get 
\begin{align*}
\left|\chi_{\xi\simeq\lambda}\left\langle\chi_{R\xi^{\frac12}\lesssim 1}\phi_{n_{3}}(R;\xi),\,\chi_{R\lesssim 1}X_{3,\mu_1}\right\rangle_{L^2_{R\,dR}}\right|&\lesssim \lambda^{-1}\cdot \left\|\chi_{R\lesssim 1}X_{3,\mu_1}\right\|_{L^\infty_{R\,dR}}\\
&\lesssim \lambda^{-1}\mu_1^{\frac{\delta}{2}}\left\|\xb_1\right\|_{\Sho_{0}(\xi_1\simeq\mu_1)}\cdot \left\|\xb_2\right\|_{\Sht_0}.
\end{align*}
From here we conclude that 
\begin{align*}
&\left\|\chi_{\xi\simeq\lambda}\left\langle\chi_{R\xi^{\frac12}\lesssim 1}\phi_{n_{3}}(R;\xi),\,\chi_{R\lesssim 1}X_{3,\mu_1}\right\rangle_{L^2_{R\,dR}}\right\|_{\Shth_1}\\
&\lesssim \lambda^{1-\frac{\delta}{2}}\left\|\chi_{\xi\simeq\lambda}\left\langle\chi_{R\xi^{\frac12}\lesssim 1}\phi_{n_{3}}(R;\xi),\,\chi_{R\lesssim 1}X_{3,\mu_1}\right\rangle_{L^2_{R\,dR}}\right\|_{L^\infty_{d\xi}}\\
&\lesssim \left(\frac{\mu_1}{\lambda}\right)^{\frac{\delta}{2}}\cdot \left\|\xb_1\right\|_{\Sho_0(\xi_1\simeq\mu_1)}\cdot \left\|\xb_2\right\|_{\Sht_0}.
\end{align*}
Square-summing over $\lambda$ and invoking the Cauchy-Schwarz inequality as usual gives the desired bound.  For high frequencies $\lambda\gtrsim 1$, we take advantage of the high degree of vanishing at $R = 0$ of $X_3$: specifically, localising the frequencies of the factors to $\mu_{1,2}\ll\lambda$ as well as fixing $i, j$ for simplicity and calling the resulting expression $X_{3,\mu_{1,2}}^{(i,j)}$, we get the poitnwise bound 
\begin{align*}
&\left|\chi_{\xi\simeq\lambda}\left\langle\chi_{R\xi^{\frac12}\lesssim 1}\phi_{n_{3}}(R;\xi),\,\chi_{R\lesssim 1}X_{3,\mu_{1,2}}^{(i,j)}\right\rangle_{L^2_{R\,dR}}\right|\\
&\lesssim \lambda^{-1-\frac{i+j}{2}}\cdot \mu_1^{\max\{\frac{i-3}{2}-\frac{\delta}{2},0\}}\cdot \mu_2^{\max\{\frac{j-3}{2}-\frac{\delta}{2},0\}}\cdot\prod_{k=1,2}\left\|\xb_k\right\|_{\Shk_{0}(\xi_k\simeq\mu_k)}
\end{align*}
whence 
\begin{align*}
&\left\|\chi_{\xi\simeq\lambda}\left\langle\chi_{R\xi^{\frac12}\lesssim 1}\phi_{n_{3}}(R;\xi),\,\chi_{R\lesssim 1}X_{3,\mu_{1,2}}^{(i,j)}\right\rangle_{L^2_{R\,dR}}\right\|_{\Shth_1}\\
&\lesssim \lambda^{\frac32 + \frac{\delta}{2}-\frac{i+j}{2}}\cdot \mu_1^{\max\{\frac{i-3}{2}-\frac{\delta}{2},0\}}\cdot \mu_2^{\max\{\frac{j-3}{2}-\frac{\delta}{2},0\}}\cdot\prod_{k=1,2}\left\|\xb_k\right\|_{\Shk_{0}(\xi_k\simeq\mu_k)},
\end{align*}
where we keep the condition $i+j\geq 4$ in mind. Square-summing over $\lambda\gtrsim 1$ and exploiting simple orthogonal arguments, the desired bound easily follows. 
\\
As for the contribution of $\tilde{P}_4$, where we only treat the first line, the second being more of the same, we get smallness from the large frequencies inside at least one of the factors in this expression. Using that $i+j\leq 3$, and fixing $i, j$ as well as the frequency $\xi_1\simeq\mu_1$ in the first factor, which results in $\tilde{P}_{4,\mu_1}^{(i,j)}$, we can bound, for large frequencies $\lambda\gtrsim 1$,
\begin{align*}
\left|\chi_{\xi\simeq\lambda}\left\langle\chi_{R\xi^{\frac12}\lesssim 1}\phi_{n_{3}}(R;\xi), \tilde{P}_{4,\mu_1}^{(i,j)}\right\rangle_{L^2_{R\,dR}}\right|\lesssim \lambda^{-1-\frac{i}{2}}\cdot\mu_1^{-(\frac32-\frac{i}{2}+\frac{\delta}{2})}\cdot \left\|\xb_1\right\|_{\Sho_{0}(\xi_1\simeq\mu_1)}\left\|\xb_2\right\|_{\Sht_0}. 
\end{align*}
This in turn implies that 
\begin{align*}
\left\|\chi_{\xi\simeq\lambda}\left\langle\chi_{R\xi^{\frac12}\lesssim 1}\phi_{n_{3}}(R;\xi), \tilde{P}_{4,\mu_1}^{(i,j)}\right\rangle_{L^2_{R\,dR}}\right\|_{\Shth_1}
\lesssim \left(\frac{\lambda}{\mu_1}\right)^{(\frac32-\frac{i}{2}+\frac{\delta}{2})}\cdot \left\|\xb_1\right\|_{\Sho_{0}(\xi_1\simeq\mu_1)}\left\|\xb_2\right\|_{\Sht_0}.
\end{align*}
But from here, keeping in mind that $i\leq 3$, we infer that 
\begin{align*}
\left(\sum_{\lambda\gtrsim 1}\left\|\sum_{\mu_1\gtrsim\lambda}\chi_{\xi\simeq\lambda}\left\langle\chi_{R\xi^{\frac12}\lesssim 1}\phi_{n_{3}}(R;\xi), \tilde{P}_{4,\mu_1}^{(i,j)}\right\rangle_{L^2_{R\,dR}}\right\|_{\Shth_1}^2\right)^{\frac12}\lesssim \prod_{k=1,2}\left\|\xb_k\right\|_{\Shk_{0}}. 
\end{align*}
For small frequencies $\lambda\lesssim 1$, we distinguish two cases $\mu_{1}\lesssim 1$ and $\mu_{1}\gtrsim1$. For the first case, we have
\begin{align*}
\left|\chi_{\xi\simeq\lambda}\left\langle\chi_{R\xi^{\frac12}\lesssim 1}\phi_{n_{3}}(R;\xi), \tilde{P}_{4,\mu_1}^{(i,j)}\right\rangle_{L^2_{R\,dR}}\right|\lesssim \mu_1^{\frac{\delta}{2}}\cdot \left\|\xb_1\right\|_{\Sho_{0}(\xi_1\simeq\mu_1)}\left\|\xb_2\right\|_{\Sht_0},
\end{align*}
which implies
\begin{align*}
\left\|\chi_{\xi\simeq\lambda}\left\langle\chi_{R\xi^{\frac12}\lesssim 1}\phi_{n_{3}}(R;\xi), \tilde{P}_{4,\mu_1}^{(i,j)}\right\rangle_{L^2_{R\,dR}}\right\|_{\Shth_1}
\lesssim \lambda^{1-\frac{\delta}{2}}\mu_{1}^{\frac{\delta}{2}}\cdot \left\|\xb_1\right\|_{\Sho_{0}(\xi_1\simeq\mu_1)}\left\|\xb_2\right\|_{\Sht_0}.
\end{align*}
Then a routing argument gives the desired result. For the second case where $\mu_{1}\gtrsim1$, we have
\begin{align*}
\left|\chi_{\xi\simeq\lambda}\left\langle\chi_{R\xi^{\frac12}\lesssim 1}\phi_{n_{3}}(R;\xi), \tilde{P}_{4,\mu_1}^{(i,j)}\right\rangle_{L^2_{R\,dR}}\right|\lesssim \mu_1^{-\left(\frac32-\frac{i}{2}+\frac{\delta}{2}\right)}\cdot \left\|\xb_1\right\|_{\Sho_{0}(\xi_1\simeq\mu_1)}\left\|\xb_2\right\|_{\Sht_0},
\end{align*}
which implies
\begin{align*}
\left\|\chi_{\xi\simeq\lambda}\left\langle\chi_{R\xi^{\frac12}\lesssim 1}\phi_{n_{3}}(R;\xi), \tilde{P}_{4,\mu_1}^{(i,j)}\right\rangle_{L^2_{R\,dR}}\right\|_{\Shth_1}
\lesssim \lambda^{1-\frac{\delta}{2}}\mu_1^{-\left(\frac32-\frac{i}{2}+\frac{\delta}{2}\right)}\cdot \left\|\xb_1\right\|_{\Sho_{0}(\xi_1\simeq\mu_1)}\left\|\xb_2\right\|_{\Sht_0}.
\end{align*}
Again a routing argument gives the desired result. 
\\
Now we turn to the case when $\xi_{1}\leq \xi\leq \xi_{2}$. Reviewing the argument handling the case $|n_{3}|\gg1$, there are two places where we used the integration by parts regarding the operator $H_{n_{3}}$. The first place is when we handle the contribution from
\begin{align*}
	\chi_{\xi\simeq\mu}\left\langle\phi_{n_{3}}(R,\xi),\,\chi_{R\lesssim\tau}\partial_{R}\phi_{1,\leq\mu}\cdot \partial_{R}\phi_{n}(R)\right\rangle_{L^{2}_{RdR}}.
\end{align*}
Note that if the vanishing order at $R=0$ of $\partial_{R}\phi_{1,\leq\mu}$ is greater or equal to $4$, then we can still perform the integration by parts argument. If the vanishing order is less or equal to $3$, then we can follow the argument estimating the contribution from 
\begin{align*}
	\chi_{\xi\simeq\mu}\left\langle\phi_{n_{3}}(R,\xi),\,\chi_{R\lesssim1}\tilde{P}_{3}(R)\right\rangle_{L^{2}_{RdR}}.
\end{align*}
The other place where we used the integration by parts argument is when estimate the contribution from 
\begin{align*}
	\chi_{\xi\simeq\mu}\left\langle\phi_{n_{3}}(R,\xi),\,\chi_{R\lesssim\tau}\partial_{R}\phi_{1,\leq\mu}\cdot\partial_{R}\phi_{2,[\mu,\infty]}\right\rangle_{L^{2}_{RdR}},
\end{align*}
where, with $n=0,\pm1$,
\begin{align*}
	\partial_{R}\phi_{2,[\mu,\infty]}=\partial_{R}\left(\phi_{n}(R)\right)\cdot\int_{0}^{R}\left[\phi_{n}(s)\right]^{-1}\calD_{n}\phi_{2,[\mu,\infty]}(s)\,ds.
\end{align*}
Note  that the vanishing order at $R=0$ in the expression $\chi_{\xi\simeq\mu}\left\langle\phi_{n_{3}}(R,\xi),\,\chi_{R\lesssim1}\tilde{P}_{3}(R)\right\rangle_{L^{2}_{RdR}}$ is at least of $O(R^{\frac92})$, which is enough for integration by parts regarding $H_{n_{3}}$ twice, as for the $|n_{3}|\gg1$ case. This completes the proof of the proposition. 
\end{proof}

In order to bound the quadratic ``null-forms" arising in the nonlinearity, we also need to deal with the terms involving temporal derivatives, which in terms of the $(\tau, R)$-coordinates involve the operator $\partial_{\tau} + \frac{\lambda_{\tau}}{\lambda}R\partial_R$. In order to understand the effect of this operator on the exceptional low angular momentum terms, it is useful to determine the effect of this operator on the expressions in \eqref{eq:levelup-rep}, \eqref{eq:levelup0rep}, \eqref{eq:levelup1rep}. Letting $n = 0,\pm1$, we get (now the function $\phi$ also depends on $\tau$)
\begin{equation}\label{eq:Dtaueffecton0pm1term}\begin{split}
&\left(\partial_{\tau} + \frac{\lambda_{\tau}}{\lambda}R\partial_R\right)\left(c_n(\tau)\phi_n(R) + \phi_n(R)\cdot\int_0^R\left[\phi_n(s)\right]^{-1}\mathcal{D}_n\phi(\tau, s)\,ds\right)\\
& = c_n'(\tau)\phi_n(R) +  c_n(\tau)\frac{\lambda_{\tau}}{\lambda}R\partial_R\phi_n(R) + \frac{\lambda_{\tau}}{\lambda}\left(R\partial_R\phi_n(R)\right)\cdot\int_0^R\left[\phi_n(s)\right]^{-1}\mathcal{D}_n\phi(\tau,s)\,ds\\
& +  \frac{\lambda_{\tau}}{\lambda}R\mathcal{D}_n\phi(\tau, R)  +  \phi_n(R)\cdot\int_0^R\left[\phi_n(s)\right]^{-1}\partial_{\tau}\left(\mathcal{D}_n\phi(\tau, s)\right)\,ds.
\end{split}\end{equation}
Here in light of the explicit algebraic nature of $\phi_n(R)$, the first three terms at the end are essentially like the terms in the original formula for $\phi(\tau, R)$, except that they come with an extra factor $\tau^{-1}\sim\frac{\lambda_{\tau}}{\lambda}$. As for the last term, write it as 
\begin{equation}\label{eq:Dtaueffecton0pm1worstterm}\begin{split}
\phi_n(R)\cdot\int_0^R\left[\phi_n(s)\right]^{-1}\partial_{\tau}\left(\mathcal{D}_n\phi(\tau, s)\right)\,ds &=  \phi_n(R)\cdot\int_0^R\left[\phi_n(s)\right]^{-1}\left(\partial_{\tau}+\frac{\lambda_{\tau}}{\lambda}s\partial_s\right)\left(\mathcal{D}_n\phi(\tau, s)\right)\,ds\\
& +  \frac{\lambda_{\tau}}{\lambda}\phi_n(R)\cdot\int_0^R\partial_s\left(s\left[\phi_n(s)\right]^{-1}\right)\left(\mathcal{D}_n\phi(\tau, s)\right)\,ds\\
& - \frac{\lambda_{\tau}}{\lambda}R\mathcal{D}_n\phi(\tau, R),
\end{split}\end{equation}
where the last term cancels against the fourth term in the earlier identity. The middle term is completely analogous to the term 
\[
\phi_n(R)\cdot\int_0^R\left[\phi_n(s)\right]^{-1}\mathcal{D}_n\phi(\tau, s)\,ds,
\]
The first term on the right is more delicate, and will be handled by formulating it on the Fourier side. Specifically, writing 
\begin{align*}
\mathcal{D}_n\phi(\tau, R) = \int_0^\infty \phi_n(R,\xi)\xb(\tau, \xi)\tilde{\rho}_{n}(\xi)\,d\xi,
\end{align*}
we express the effect of the operator $\partial_{\tau}+\frac{\lambda_{\tau}}{\lambda}R\partial_R$ by 
\begin{equation}\label{eq:dilationoperatorFourier1}\begin{split}
\left(\partial_{\tau}+\frac{\lambda_{\tau}}{\lambda}R\partial_R\right)\mathcal{D}_n\phi(\tau, R) = \int_0^\infty \phi_n(R,\xi)\mathcal{D}_{\tau}^{(n)}\xb(\tau, \xi)\tilde{\rho}_{n}(\xi)\,d\xi +  \int_0^\infty \phi_n(R,\xi)\mathcal{K}_n^{(0)}\xb(\tau, \xi)\tilde{\rho}_{n}(\xi)\,d\xi,
\end{split}\end{equation}
where $\mathcal{K}_n^{(0)}$ is the non-diagonal part of the transference operator, which vanishes when $n = 1$. Moreover, the dilation type operator $\mathcal{D}_{\tau}^{(n)}$ is given by the formula 
\begin{align*}
\mathcal{D}_{\tau}^{(n)} = \partial_{\tau} - 2\frac{\lambda_{\tau}}{\lambda}\xi\partial_{\xi} - \frac{\lambda_{\tau}}{\lambda}\frac{\left(\tilde{\rho}_{n}(\xi)\right)'\xi}{\tilde{\rho}_{n}(\xi)} - \frac{\lambda_{\tau}}{\lambda}. 
\end{align*}
Next, we consider the effect of $\partial_{\tau} + \frac{\lambda_{\tau}}{\lambda}R\partial_R$ on large angular momentum functions, where things are more direct. In fact, writing 
\[
\phi(\tau, R) = \int_0^\infty \phi_{n}(R;\xi)\xb(\tau, \xi)\tilde{\rho}_{n}(\xi)\,d\xi, 
\]
we obtain 
\begin{equation}\label{eq:dilationoperatorFourier2}\begin{split}
\left(\partial_{\tau}+\frac{\lambda_{\tau}}{\lambda}R\partial_R\right)\phi(\tau, R) &=  \int_0^\infty \phi_{n}(R;\xi)\mathcal{D}_{\tau}^{(\hbar)}\xb(\tau, \xi)\tilde{\rho}_{n}(\xi)\,d\xi\\
& +  \int_0^\infty \phi_{n}(R;\xi)\mathcal{K}_{\hbar}^{(0)}\xb(\tau, \xi)\tilde{\rho}_{n}(\xi)\,d\xi,
\end{split}\end{equation}
where the dilation type operator is given by the formula 
\begin{align*}
\mathcal{D}_{\tau}^{(\hbar)} = \partial_{\tau} - 2\frac{\lambda_{\tau}}{\lambda}\xi\partial_{\xi} - \frac{\lambda_{\tau}}{\lambda}\frac{\left(\tilde{\rho}_{n}(\xi)\right)'\xi}{\tilde{\rho}_{n}(\xi)} - 2\frac{\lambda_{\tau}}{\lambda}.
\end{align*}
We shall then use bounds on $\xb(\tau, \xi),\,\mathcal{D}_{\tau}^{(\hbar)}\xb(\tau, \xi)$ or else $\mathcal{D}_{\tau}^{(n)}\xb(\tau, \xi)$ in order to control the bilinear terms involving time derivatives. We shall require an analogue of the $L^\infty$-bounds in Proposition~\ref{prop:exceptionalnderivative1}  without derivatives and Fourier coefficients that are better behaved in the low frequency regime:
\begin{lemma}\label{lem:exceptionalnderivative2} Let the functions $f_j,\,j = \pm 1,0$ be as in the statement of  Proposition~\ref{prop:exceptionalnderivative1}. Then we have the bounds
	\begin{align*}
	\left\|f_j\right\|_{L^\infty_{R\,dR}}\lesssim \left|c_j\right| + \left\|y_j\right\|_{S_1^{(j)}},\quad S_1^{(j)} = \xi^{\frac12}S_0^{(j)}. 
	\end{align*}
	Here $S^{(1)}_{0}:=S^{+}_{0}, S^{(0)}_{0}:=S_{0}^{0}, S^{(-j)}_{0}:=S_{0}^{-}$.
\end{lemma}
\begin{proof} We consider the integral expression contributing to $f_j$. Thus set now 
	\begin{align*}
	f_j(R) &= \phi_j(R)\cdot\int_0^R\left[\phi_j(s)\right]^{-1}\cdot \left(\int_0^\infty \phi_j(s,\xi)y_j(\xi)\tilde{\rho}_{j}(\xi)\,d\xi\right)\,ds\\
	& = \phi_j(R)\cdot\int_0^R\left[\phi_j(s)\right]^{-1}\cdot \left(\int_0^\infty\chi_{s\xi^{\frac12}\lesssim 1} \phi_j(s,\xi)y_j(\xi)\tilde{\rho}_{j}(\xi)\,d\xi\right)\,ds\\
	& + \phi_j(R)\cdot\int_0^R\left[\phi_j(s)\right]^{-1}\cdot \left(\int_0^\infty\chi_{s\xi^{\frac12}\gtrsim 1} \phi_j(s,\xi)y_j(\xi)\tilde{\rho}_{j}(\xi)\,d\xi\right)\,ds.
	\end{align*}
	Consider for example the case $j = 1$, the others being similar. Then recall the bound 
	\[
	\left|\phi_1(s,\xi)\right|\lesssim \min\{s, \xi^{-\frac12}(s\xi^{\frac12})^{-\frac12}\}, 
	\]
	and so we get the bound 
	\begin{align*}
	\left|\chi_{s\xi^{\frac12}\lesssim 1}\phi_1(s,\xi)\tilde{\rho}_{1}(\xi)\right|\lesssim s^{-1}\cdot (s\xi^{\frac12})
	\end{align*}
	We conclude that 
	\begin{align*}
	&\left|\phi_1(R)\cdot\int_0^R\left[\phi_1(s)\right]^{-1}\cdot \left(\int_0^\infty\chi_{s\xi^{\frac12}\lesssim 1} \phi_1(s,\xi)y_1(\xi)\tilde{\rho}_{1}(\xi)\,d\xi\right)\,ds\right|\\
	&\lesssim \int_0^\infty \langle\xi\rangle^{\frac12}\left|y_1(\xi)\right|\,d\xi\lesssim \left\|y_1\right\|_{S_1^{(1)}}. 
	\end{align*}
	The second integral expression is bounded similarly, by using the oscillatory nature of $\chi_{s\xi^{\frac12}\gtrsim 1} \phi_j(s,\xi)$ and performing integration by parts. 
	\end{proof}
We also have the following estimates for large frequencies:
\begin{lemma}\label{lem:exceptionalnderivatie2 large}
	Let $\lambda\geq 1$ and $f_{j,\geq \lambda},j=0,\pm1$ be defined as
	\begin{align}\label{def: fj ge lambda}
		f_{j,\geq \lambda}(R):=\phi_j(R)\cdot\int_0^R\left[\phi_j(s)\right]^{-1}\cdot \left(\int_{\lambda}^{\infty} \phi_j(s,\xi)y_j(\xi)\tilde{\rho}_{j}(\xi)\,d\xi\right)\,ds.
	\end{align}
	Then we have the following $L^{\infty}$-estimate
	\begin{align*}
		\left\|f_{j,\geq\lambda}\right\|_{L^{\infty}_{RdR}}\lesssim \lambda^{-\frac12}\left\|y_{j}\right\|_{S_{1}^{(j)}}, \quad \left\|\partial_{R}f_{j,\geq \lambda}\right\|_{L^{\infty}_{RdR}}\lesssim \sum_{\mu\geq \lambda}\mu^{-1-\frac{\delta}{2}}\left\|y_{j}(\xi)\right\|_{S_{1}^{(j)}(\xi\simeq\mu)}.
	\end{align*}
	and the following $L^{2}$-estimates
	\begin{align*}
		\left\|\partial_{R}f_{j,\geq \lambda}\right\|_{L^{2}_{RdR}}\lesssim \sum_{\mu\geq \lambda}\mu^{-1-\frac{\delta}{2}}\left\|y_{j}(\xi)\right\|_{S_{1}^{(j)}(\xi\simeq\mu)},\quad \textrm{for}\quad j=0,1,
	\end{align*}
	\begin{align*}
		\left\|\partial_{R}f_{-1,\geq \lambda}\right\|_{L^{2}_{RdR}(R\lesssim \lambda^{-\frac12})}\lesssim \lambda^{-1-\frac{\delta}{2}}\left\|y_{-1}(\xi)\right\|_{S_{1}^{(-1)}},\quad \left\|\partial_{R}f_{-1,\geq \lambda}\right\|_{L^{2}_{RdR}(R\gtrsim \lambda^{-\frac12})}\lesssim \lambda^{-\frac12-\frac{\delta}{2}}\left\|y_{-1}(\xi)\right\|_{S_{1}^{(-1)}}.
	\end{align*}
\end{lemma}
\begin{proof}
The proof is similar to that of Lemma \ref{lem:exceptionalnderivative2} and we omit the details.
\end{proof}

\begin{proposition}\label{prop:bilin3} Let $\left|n_1\right|\gg 1$, $\left|n_2\right|\leq \left|n_1\right|$, where $n_2$ is allowed to take any integer value. Also, assume that $\left|n_3\right|\geq 2$ and that either $n_3\simeq n_1$ or else $\left|n_3\right|\ll \left|n_1\right|$ as well as $n_1\simeq -n_2$. For $\left|n\right|\geq 2$, assuming that $\phi(R)$ is an angular momentum $n$ function with
	\[
	\phi(\tau, R) = \int_0^\infty \phi_{n}(R;\xi)\xb(\tau,\xi)\tilde{\rho}_{n}(\xi)\,d\xi,\quad\hbar = \frac{1}{n+1}, 
	\]
	we set (recalling $\Sh_{1}= \xi^{\frac12}\Sh_{0}$)
	\[
	\left\|\phi(\tau, \cdot)\right\|_{\tilde{S}_1^{(n)}}: = \left\|\mathcal{D}_{\tau}\xb(\tau,\xi)\right\|_{\Sh_{1}}.  
	\]
	For $n = 0,\pm 1$ assuming that $\phi(\tau,R)$ is an angular momentum $n$ function provided 
	\[
	\phi(\tau,R) = c_n(\tau)\phi_n(R) + \phi_n(R)\cdot\int_0^R\left[\phi_n(s)\right]^{-1}\tilde{\mathcal{D}}\phi(\tau, s)\,ds,\quad\tilde{\mathcal{D}}\phi(\tau,R) = \int_0^\infty \xb(\tau,\xi)\phi_{n}(R,\xi)\tilde{\rho}_{n}(\xi)\,d\xi,
	\]
	where $\tilde{\calD} = \mathcal{D}_-,\mathcal{D}, \mathcal{D}_+$ according to $n$ as in the preceding. Then set 
	\[
	\left\|\phi(\tau, \cdot)\right\|_{\tilde{S}_1^{(n)}}: = \left|c_n'(\tau)\right| + \left\| \mathcal{D}_{\tau}\xb(\tau,\xi)\right\|_{S_1^{(n)}}. 
	\]
	With the conventions and under the assumptions on the $n_j$ stated at the beginning, and assuming $\tau\gg 1$, we have the following bound 
	\begin{equation}\label{eq:productlargen3alln20}
	\left\|\left\langle\phi_{n_{3}}(R;\xi),\,\chi_{R\lesssim\tau}\left(\partial_\tau + \frac{\lambda_{\tau}}{\lambda}R\partial_R\right)\phi_1\cdot\left(\partial_\tau + \frac{\lambda_{\tau}}{\lambda}R\partial_R\right)\phi_2\right\rangle_{L^2_{R\,dR}}\right\|_{\Shth_1}\lesssim \tau^{\delta}\left\langle n_2\right\rangle^{2}\prod_{j=1}^2\left\|\phi_j\right\|_{\tilde{S}_0^{(n_j)}\cap\tilde{S}_1^{(n_j)}}
	\end{equation}
	where we assume that the factors $\phi_j$ are angular momentum $n_j$ functions, $j = 1, 2$. 
	\\
	Next, assume that $|n_1|\lesssim 1$, while the other assumptions on $n_j$ stated above are still valid. Then if $\phi_{1,2}$ are angular momentum $n_j$ functions with finite $\|\cdot\|_{\cap_{k=0,1}\tilde{S}_{k}^{(n_j)}}$-norm, the product 
	\[
	\left(\partial_{\tau}+\frac{\lambda_{\tau}}{\lambda}R\partial_R\right)\phi_1\cdot\left(\partial_{\tau}+\frac{\lambda_{\tau}}{\lambda}R\partial_R\right)\phi_2
	\]
	admits a third order Taylor expansion at $R = 0$ of the form
	\[
	P_3(R): = \sum_{l=0}^3\gamma_l R^l,
	\]
	where we have the bound
	\begin{equation}\label{eq:Taylorbound10}
	\sum_{l=0}^3\left|\gamma_l\right|\lesssim \prod_{j=1,2}\left\|\phi_j\right\|_{\tilde{S}_0^{(n_j)}\cap \tilde{S}_1^{(n_j)}}.
	\end{equation}
	Furthermore, we have the bound 
	\begin{equation}\label{eq:productsmalln1geq2alln20}
	\left\|\left\langle\phi_{n_{3}}(R;\xi),\,\chi_{R\lesssim\tau}\left(\partial_{\tau}+\frac{\lambda_{\tau}}{\lambda}R\partial_R\right)\phi_1\cdot\left(\partial_{\tau}+\frac{\lambda_{\tau}}{\lambda}R\partial_R\right)\phi_2 - \chi_{R\lesssim 1}P_3(R)\right\rangle_{L^2_{R\,dR}}\right\|_{\Shth_1}\lesssim \tau^{\delta}\prod_{j=1}^2\left\|\phi_j\right\|_{\tilde{S}_0^{(n_j)}\cap \tilde{S}_1^{(n_j)}}
	\end{equation}
	The last inequality remains correct if we subtract from $P_3$ those terms $\gamma_l R^l$ with $l\geq |n_3 - 1|, l\equiv n_3 - 1\text{(mod 2)}$, provided $|n_3 -1|\leq 3$. 
\end{proposition}
\begin{proof} {\it{First inequality}}. Start with the case when $|n_2|\geq 2$, i.e., we can use the representation \eqref{eq:dilationoperatorFourier2}. Observe that since $\mathcal{K}_{\hbar}^{0}$ maps $\Sh_{0}$ into $\Sh_{1}$, and we assume that $\mathcal{D}_{\tau}\xb\in \Sh_{1}$, we can treat this case exactly as in the proof of Proposition~\ref{prop:bilin1}, observing that there we always deal with the product $\xi^{\frac12}\xb(\tau, \xi)\in \Sh_{1}$ for the Fourier transform. Thus to conclude the proof of \eqref{eq:productlargen3alln20}, it suffices to deal with the case $|n_2|\leq 1$, where we have to take advantage of \eqref{eq:Dtaueffecton0pm1term}, \eqref{eq:Dtaueffecton0pm1worstterm} as well as \eqref{eq:dilationoperatorFourier2}. In the following, we omit the contributions coming directly from the resonance/root part, i.e., the first, second term on the right hand side in \eqref{eq:Dtaueffecton0pm1term}, since their contribution is straightforward to bound. Throughout we have $n_3\simeq n_1$ due to the assumptions on the angular momenta for the first inequality.   
	\\
	
	{\it{Contribution of the third term on the right hand side of  \eqref{eq:Dtaueffecton0pm1term} and the second term on the right hand side in \eqref{eq:Dtaueffecton0pm1worstterm}.}} For all intents and purposes, we may assume here that 
	\[
	\left(\partial_{\tau}+\frac{\lambda_{\tau}}{\lambda}R\partial_R\right)\phi_2(R) = \frac{\lambda_{\tau}}{\lambda}\phi_{n}(R)\cdot \int_0^R \left[\phi_n(s)\right]^{-1}\mathcal{D}_n\phi_2(s)\,ds,\quad n = 0,\pm 1, 
	\]
	where $\mathcal{D}_n\phi_2$ admits the representation
	\begin{align*}
	\mathcal{D}_n\phi_2(R) = \int_0^\infty \phi_n(R,\xi) \xb(\xi)\tilde{\rho}_{n}(\xi)\,d\xi. 
	\end{align*}
	On the other hand, $\phi_1(\tau, R)$ is assumed to admit a representation as at the beginning of the proposition, corresponding to large angular momentum, and we then have \eqref{eq:dilationoperatorFourier2}. We can then essentially repeat the proof of \eqref{eq:productlargen3alln2} in Prop.~\ref{prop:bilin2} to deal with this case. 
	\\
	
	This reduces things to the case when we substitute the two terms on the right of \eqref{eq:dilationoperatorFourier1} for $\left(\partial_{\tau}+\frac{\lambda_{\tau}}{\lambda}s\partial_s\right)\mathcal{D}_n\phi_2$, which in turn gets substituted into \eqref{eq:Dtaueffecton0pm1worstterm}. On the other hand, we still use the representation \eqref{eq:dilationoperatorFourier2} for $\left(\partial_{\tau}+\frac{\lambda_{\tau}}{\lambda}R\partial_R\right)\phi_1$. We follow the same steps as for the proof of \eqref{eq:productlargen3alln2}. As before we call $\xi$ the output frequency, and $\xi_{1,2}$ the frequencies of the factors (more precisely, in case of $\left(\partial_{\tau}+\frac{\lambda_{\tau}}{\lambda}R\partial_R\right)\phi_2$, the frequency of $\mathcal{D}_n\phi_2$). 
	\\
	
	{\it{(1.a) $\xi\leq \xi_1$}}. We take advantage of Lemma~\ref{lem:exceptionalnderivative2}, which under our current hypotheses on $\left(\partial_{\tau}+\frac{\lambda_{\tau}}{\lambda}R\partial_R\right)\phi_2$ give
	\begin{align*}
	\left\|\left(\partial_{\tau}+\frac{\lambda_{\tau}}{\lambda}R\partial_R\right)\phi_2\right\|_{L^\infty_{R\,dR}}&\lesssim \left\|\mathcal{D}_{\tau}^{(n)}\xb(\tau,\xi)\right\|_{S_1^{(n)}} + \left\|\mathcal{K}_n^{(0)}\xb\right\|_{S_1^{(n)}}\\
	&\lesssim  \left\|\mathcal{D}_{\tau}^{(n)}\xb(\tau,\xi)\right\|_{S_1^{(n)}} + \left\|\xb\right\|_{S_0^{(n)}}\\
	&\lesssim \left\|\phi_2\right\|_{\tilde{S}_0^{(n)}\cap \tilde{S}_1^{(n)}}. 
	\end{align*}
	We need to bound the expression
	\begin{align*}
	&\left\|\sum_{\lambda}\chi_{\xi\simeq\lambda}\left\langle \phi_{n_{3}}(R;\xi),\,\left(\partial_{\tau}+\frac{\lambda_{\tau}}{\lambda}R\partial_R\right)\phi_{1\,\geq\lambda}\cdot \left(\partial_{\tau}+\frac{\lambda_{\tau}}{\lambda}R\partial_R\right)\phi_2\right\rangle_{L^2_{R\,dR}}\right\|_{\Shth_1}\\
	&\lesssim \left(\sum_{\lambda}\left\|\left\langle \phi_{n_{3}}(R;\xi),\,\left(\partial_{\tau}+\frac{\lambda_{\tau}}{\lambda}R\partial_R\right)\phi_{1\,\geq\lambda}\cdot \left(\partial_{\tau}+\frac{\lambda_{\tau}}{\lambda}R\partial_R\right)\phi_2\right\rangle_{L^2_{R\,dR}}\right\|_{S_1^{(\hbar_3)}(\xi\simeq\lambda)}^2\right)^{\frac12}\\
	&\lesssim \hbar_3^{2-\delta}\left(\sum_{\lambda}\lambda^{1-\delta}\langle\lambda\hbar_3^2\rangle^{3+2\delta}\left\|\left(\partial_{\tau}+\frac{\lambda_{\tau}}{\lambda}R\partial_R\right)\phi_{1\,\geq\lambda}\cdot \left(\partial_{\tau}+\frac{\lambda_{\tau}}{\lambda}R\partial_R\right)\phi_2\right\|_{L^2_{R\,dR}}^2\right)^{\frac12},
	\end{align*}
	where we have exploited orthogonality and the Plancherel's theorem for the distorted Fourier transform. Recalling \eqref{eq:dilationoperatorFourier2} we get (letting $\xb_1$ the Fourier transform at angular momentum $n_1$ of $\phi_1$) 
	\begin{align*}
	\left\|\left(\partial_{\tau}+\frac{\lambda_{\tau}}{\lambda}R\partial_R\right)\phi_{1\,\geq\lambda}\right\|_{L^2_{R\,dR}}&\lesssim \sum_{\mu\geq \lambda}\left[\left\|\mathcal{D}_{\tau}^{(\hbar_1)}\xb_1\right\|_{L^2_{d\xi}} + \left\|\mathcal{K}_{\hbar_1}^{(0)}\xb_1\right\|_{L^2_{d\xi}}\right]\\
	&\lesssim \sum_{\mu\geq \lambda}\hbar_1^{-1}(\mu\hbar_1^2)^{-\frac12+\frac{\delta}{2}}\langle\mu\hbar_1^2\rangle^{-\frac32-\delta}\cdot\left[\left\|\xb_1\right\|_{\Sho_{0}(\xi\simeq\mu)}+\left\|\mathcal{D}_{\tau}^{(\hbar_1)}\xb_1\right\|_{\Sho_{1}(\xi\simeq\mu)}\right].
	\end{align*}
	Recalling our hypothesis that $n_1\simeq n_3$ and using Holder's inequality and also recalling the above $L^\infty$ bound for the second factor, we then have 
	\begin{align*}
	&\hbar_3^{2-\delta}\lambda^{\frac12-\frac{\delta}{2}}\langle\lambda\hbar_3^2\rangle^{\frac32+\delta}\cdot \left\|\left(\partial_{\tau}+\frac{\lambda_{\tau}}{\lambda}R\partial_R\right)\phi_{1\,\geq\lambda}\cdot \left(\partial_{\tau}+\frac{\lambda_{\tau}}{\lambda}R\partial_R\right)\phi_2\right\|_{L^2_{R\,dR}}\\
	&\lesssim   \left\|\phi_2\right\|_{\tilde{S}_0^{(n)}\cap \tilde{S}_1^{(n)}}\sum_{\mu\geq \lambda}\left(\frac{\lambda}{\mu}\right)^{\frac12-\frac{\delta}{2}}\cdot \left[\left\|\xb_1\right\|_{\Sho_{0}(\xi\simeq\mu)}+\left\|\mathcal{D}_{\tau}^{(\hbar_1)}\xb_1\right\|_{\Sho_{1}(\xi\simeq\mu)}\right].
	\end{align*}
	Using the Cauchy-Schwarz inequality as well as orthogonality, it follows that 
	\begin{align*}
	&\hbar_3^{2-\delta}\left(\sum_{\lambda}\lambda^{1-\delta}\langle\lambda\hbar_3^2\rangle^{3+2\delta}\left\|\left(\partial_{\tau}+\frac{\lambda_{\tau}}{\lambda}R\partial_R\right)\phi_{1\,\geq\lambda}\cdot \left(\partial_{\tau}+\frac{\lambda_{\tau}}{\lambda}R\partial_R\right)\phi_2\right\|_{L^2_{R\,dR}}^2\right)^{\frac12}\\
	&\lesssim  \left\|\phi_2\right\|_{\tilde{S}_0^{(n)}\cap \tilde{S}_1^{(n)}}\left(\sum_{\lambda}\left[\sum_{\mu\geq \lambda}\left(\frac{\lambda}{\mu}\right)^{\frac12-\frac{\delta}{2}}\cdot \left[\left\|\xb_1\right\|_{\Sho_{0}(\xi\simeq\mu)}+\left\|\mathcal{D}_{\tau}^{(\hbar_1)}\xb_1\right\|_{\Sho_{1}(\xi\simeq\mu)}\right]\right]^2\right)^{\frac12}\\
	&\lesssim  \left\|\phi_2\right\|_{\tilde{S}_0^{(n)}\cap \tilde{S}_1^{(n)}}\left(\sum_{\lambda}\sum_{\mu\geq \lambda}\left(\frac{\lambda}{\mu}\right)^{\frac12-\frac{\delta}{2}}\cdot \left[\left\|\xb_1\right\|_{\Sho_{1}(\xi\simeq\mu)}+\left\|\mathcal{D}_{\tau}^{(\hbar_1)}\xb_1\right\|_{\Sho_{1}(\xi\simeq\mu)}\right]^2\right)^{\frac12}\\
	&\lesssim \prod_{j=1,2}\left\|\phi_j\right\|_{\tilde{S}_0^{(n_j)}\cap \tilde{S}_1^{(n_j)}}.
	\end{align*}
	
	{\it{(1.b) $\xi_1<\xi\leq\xi_2$}}. We distinguish between different frequency ranges for the output frequency $\xi$. 
	\\
	
	{\it{Small output frequency $\xi\lesssim 1$.}} Localizing to dyadic $\xi\simeq\lambda\lesssim 1$, consider the expression 
	\[
	\chi_{\xi\simeq\lambda}\left\langle\phi_{n_{3}}(R;\xi),\,\chi_{R\lesssim\tau}\left(\partial_{\tau}+\frac{\lambda_{\tau}}{\lambda}R\partial_R\right)\phi_{1\,<\lambda}\cdot \left(\partial_{\tau}+\frac{\lambda_{\tau}}{\lambda}R\partial_R\right)\phi_{2,\geq \lambda}\right\rangle_{L^2_{R\,dR}}, 
	\]
	The subscript in $\phi_{2,\geq \lambda}$ refers to the frequency variable occurring in the representation \eqref{eq:dilationoperatorFourier1}, which in turn gets substituted into \eqref{eq:Dtaueffecton0pm1worstterm}, \eqref{eq:Dtaueffecton0pm1term}. For this small frequency regime, this localization does not play an important role, however. We decompose 
	\begin{equation}\label{eq:timederivprodlowfreq1}\begin{split}
	&\chi_{\xi\simeq\lambda}\left\langle\phi_{n_{3}}(R;\xi),\,\chi_{R\lesssim\tau}\left(\partial_{\tau}+\frac{\lambda_{\tau}}{\lambda}R\partial_R\right)\phi_{1\,<\lambda}\cdot \left(\partial_{\tau}+\frac{\lambda_{\tau}}{\lambda}R\partial_R\right)\phi_{2,\geq \lambda}\right\rangle_{L^2_{R\,dR}}\\
	&=\chi_{\xi\simeq\lambda}\left\langle\chi_{\hbar_{3}R\lambda^{\frac12}\leq \frac12}\phi_{n_{3}}(R;\xi),\,\chi_{R\lesssim\tau}\left(\partial_{\tau}+\frac{\lambda_{\tau}}{\lambda}R\partial_R\right)\phi_{1\,<\lambda}\cdot \left(\partial_{\tau}+\frac{\lambda_{\tau}}{\lambda}R\partial_R\right)\phi_{2,\geq \lambda}\right\rangle_{L^2_{R\,dR}}\\
	& + \chi_{\xi\simeq\lambda}\left\langle\chi_{\hbar_{3}R\lambda^{\frac12}\geq \frac12}\phi_{n_{3}}(R;\xi),\,\chi_{R\lesssim\tau}\left(\partial_{\tau}+\frac{\lambda_{\tau}}{\lambda}R\partial_R\right)\phi_{1\,<\lambda}\cdot \left(\partial_{\tau}+\frac{\lambda_{\tau}}{\lambda}R\partial_R\right)\phi_{2,\geq \lambda}\right\rangle_{L^2_{R\,dR}}
	\end{split}\end{equation}
	By Proposition \ref{prop:DFT nlarge}, the Fourier basis $\phi_{n_{3}}(R,\xi)$ satisfies the bound $\left|\chi_{\hbar_{3}R\xi^{\frac12}\leq \frac12}\phi_{n_{3}}(R,\xi)\right|\lesssim \left(\frac12\right)^{c\hbar_{3}^{-1}}$ for some absolute constant $c>0$. This together with a slightly sharpened version of Lemma~\ref{lem:derLinfty} gives the following estimate for $\lambda\lesssim 1$ and an absolute constant $c>0$ 
	\begin{align*}
	\left\|\int_0^\infty \chi_{\hbar_{3} R\lambda^{\frac12}\leq \frac12}\phi_{n_{3}}(R;\xi)\chi_{\xi\lesssim \lambda}\xb(\xi)\tilde{\rho}_{n}(\xi)\,d\xi\right\|_{L^\infty_{R\,dR}}\lesssim \hbar^{-2+\delta}_{3}\lambda^{\frac{\delta}{2}}\cdot \left(\frac12\right)^{c\hbar^{-1}_{3}}\cdot \left\|\xb\right\|_{S^{\hbar_{3}}_{1}}.  
	\end{align*}
Since $\hbar_{1}\simeq \hbar_{3}$, we conclude, using Plancherel's theorem for the distorted Fourier transform as well as Holder's inequality that 
	\begin{align*}
	&\left\|\chi_{\xi\simeq\lambda}\left\langle\chi_{\hbar_{3}R\lambda^{\frac12}\leq \frac12}\phi_{n_{3}}(R;\xi),\,\chi_{R\lesssim\tau}\left(\partial_{\tau}+\frac{\lambda_{\tau}}{\lambda}R\partial_R\right)\phi_{1\,<\lambda}\cdot \left(\partial_{\tau}+\frac{\lambda_{\tau}}{\lambda}R\partial_R\right)\phi_{2,\geq \lambda}\right\rangle_{L^2_{R\,dR}}\right\|_{\Shth_1}\\
	&\lesssim \hbar_3^{2-\delta}\lambda^{\frac12-\frac{\delta}{2}}\cdot \left\|\chi_{R\lesssim\min\{\hbar_{3}^{-1}\lambda^{-\frac12},\tau\}}\left(\partial_{\tau}+\frac{\lambda_{\tau}}{\lambda}R\partial_R\right)\phi_{1\,<\lambda}\cdot \left(\partial_{\tau}+\frac{\lambda_{\tau}}{\lambda}R\partial_R\right)\phi_{2,\geq \lambda}\right\|_{L^2_{R\,dR}}\\
	&\lesssim \hbar_3^{2-\delta}\lambda^{\frac12-\frac{\delta}{2}}\cdot\hbar_{3}^{-1+\delta} \lambda^{-\frac12+\frac{\delta}{2}}\tau^{\delta}\cdot \left\|\left(\partial_{\tau}+\frac{\lambda_{\tau}}{\lambda}R\partial_R\right)\phi_{1\,<\lambda}\right\|_{L^\infty_{R\,dR}}\cdot \left\|\left(\partial_{\tau}+\frac{\lambda_{\tau}}{\lambda}R\partial_R\right)\phi_{2,\geq \lambda}\right\|_{L^\infty_{R\,dR}}. 
	\end{align*}
	Taking advantage of the preceding sharpened $L^\infty$-bound, the fact that $\hbar_1\simeq \hbar_3$, as well as Lemma~\ref{lem:exceptionalnderivative2}, and also keeping in mind the representation \eqref{eq:dilationoperatorFourier2}, we can bound the preceding by 
	\begin{align*}
	& \hbar_3^{2-\delta}\lambda^{\frac12-\frac{\delta}{2}}\cdot \lambda^{-\frac12+\frac{\delta}{2}}\hbar_{3}^{-1+\delta}\tau^{\delta}\cdot \left\|\left(\partial_{\tau}+\frac{\lambda_{\tau}}{\lambda}R\partial_R\right)\phi_{1\,<\lambda}\right\|_{L^\infty_{R\,dR}}\cdot \left\|\left(\partial_{\tau}+\frac{\lambda_{\tau}}{\lambda}R\partial_R\right)\phi_{2,\geq \lambda}\right\|_{L^\infty_{R\,dR}}\\
	&\lesssim \tau^{\delta}\cdot \lambda^{\frac{\delta}{2}}\cdot\left\|\phi_1\right\|_{\tilde{S}_1^{(n_1)}\cap\tilde{S}_0^{(n_1)}}\cdot \left\|\phi_2\right\|_{\tilde{S}_1^{(n_2)}\cap\tilde{S}_0^{(n_2)}}.
	\end{align*}
	This can then be summed over dyadic $\lambda\lesssim 1$ to yield the desired bound. 
	\\
	As for the second term on the right in \eqref{eq:timederivprodlowfreq1}, we perform integration by parts in the inner product, replacing it schematically by 
	\begin{align*}
	\lambda^{-\frac12}\chi_{\xi\simeq\lambda}\left\langle\chi_{\hbar_{3}R\lambda^{\frac12}\geq \frac12}\phi_{n_{3}}(R;\xi),\,\partial_R\left[\chi_{R\lesssim\tau}\left(\partial_{\tau}+\frac{\lambda_{\tau}}{\lambda}R\partial_R\right)\phi_{1\,<\lambda}\cdot \left(\partial_{\tau}+\frac{\lambda_{\tau}}{\lambda}R\partial_R\right)\phi_{2,\geq \lambda}\right]\right\rangle_{L^2_{R\,dR}}
	\end{align*}
	To bound this, we recall our assumptions on the factors $\left(\partial_{\tau}+\frac{\lambda_{\tau}}{\lambda}R\partial_R\right)\phi_j,\,j = 1,2$, and use the following $L^2$-bounds 
	\begin{align*}
	&\left\|\partial_R\left(\partial_{\tau}+\frac{\lambda_{\tau}}{\lambda}R\partial_R\right)\phi_{1\,<\lambda}\right\|_{L^2_{R\,dR}}\lesssim \left\|\xi^{\frac12}\mathcal{D}_{\tau}^{(\hbar_1)}\xb_1\right\|_{L^2_{d\xi}(\xi<\lambda)} + \left\|\xi^{\frac12}\mathcal{K}_{\hbar_1}^{(0)}\xb_1\right\|_{L^2_{d\xi}(\xi<\lambda)}\\
	&\hspace{5cm}\lesssim \hbar_1^{-2+\delta}\lambda^{\frac{\delta}{2}}\left\|\phi_1\right\|_{\tilde{S}_1^{(\hbar_1)}\cap \tilde{S}_0^{(\hbar_1)}}\\
	&\left\|\partial_R\left(\partial_{\tau}+\frac{\lambda_{\tau}}{\lambda}R\partial_R\right)\phi_{2,\geq \lambda}\right\|_{L^2_{R\,dR}(R\lesssim\tau)}\lesssim \left\|\phi_2\right\|_{\tilde{S}_1^{(n_2)}\cap \tilde{S}_0^{(n_2)}}.
	\end{align*}
	In fact, to get the last bound for $n=0,\pm1$, write (recall the discussion between \eqref{eq:Dtaueffecton0pm1worstterm} and \eqref{eq:dilationoperatorFourier1})
	\begin{align*}
	\partial_R\left(\partial_{\tau}+\frac{\lambda_{\tau}}{\lambda}R\partial_R\right)\phi_{2,\geq \lambda} &= \partial_R\phi_n(R)\cdot\int_0^R\left[\phi_n(s)\right]^{-1}\left(\partial_\tau + \frac{\lambda_{\tau}}{\lambda}s\partial_s\right)\mathcal{D}_n\phi_{2,\geq\lambda}(\tau,s)\,ds\\
	& + \left(\partial_\tau + \frac{\lambda_{\tau}}{\lambda}R\partial_R\right)\mathcal{D}_n\phi_{2,\geq\lambda}(\tau,R)
	\end{align*}
	and also recall the Fourier representation \eqref{eq:dilationoperatorFourier1}, which in light of the Plancherel's theorem for the distorted Fourier transform implies the desired bound for the second term on the right. As for the first term on the right, one can argue similarly to the proof of Lemma~\ref{lem:exceptionalnderivative2} to get the desired bound. If we add the easily verified bounds (see Lemma~\ref{lem:exceptionalnderivative2})
	\begin{align*}
	&\left\|\left(\partial_{\tau}+\frac{\lambda_{\tau}}{\lambda}R\partial_R\right)\phi_{2,\geq \lambda}\right\|_{L^\infty_{R\,dR}}\lesssim \left\|\phi_2\right\|_{\tilde{S}_1^{(n_2)}\cap \tilde{S}_0^{(n_2)}},\\
	&\left\|\left(\partial_{\tau}+\frac{\lambda_{\tau}}{\lambda}R\partial_R\right)\phi_{1,<\lambda}\right\|_{L^\infty_{R\,dR}}\lesssim \hbar_1^{-2+\delta}\lambda^{\frac{\delta}{2}}\left\|\phi_1\right\|_{\tilde{S}_1^{(n_1)}\cap\tilde{S}_0^{(n_1)}}
	\end{align*}
	and using Holder's inequality and the Leibniz product rule, we infer that (recall the restriction $\lambda\lesssim 1$)
	\begin{align*}
	&\left\| \chi_{\xi\simeq\lambda}\left\langle\chi_{\hbar_{3}R\lambda^{\frac12}\geq \frac12}\phi_{n_{3}}(R;\xi),\,\chi_{R\lesssim\tau}\left(\partial_{\tau}+\frac{\lambda_{\tau}}{\lambda}R\partial_R\right)\phi_{1\,<\lambda}\cdot \left(\partial_{\tau}+\frac{\lambda_{\tau}}{\lambda}R\partial_R\right)\phi_{2,\geq \lambda}\right\rangle_{L^2_{R\,dR}}\right\|_{\Shth_{1}}\\
	&\lesssim \hbar_3^{2-\delta}\lambda^{\frac12-\frac{\delta}{2}}\cdot \left\| \chi_{\xi\simeq\lambda}\left\langle\chi_{\hbar_{3}R\lambda^{\frac12}\geq \frac12}\phi_{n_{3}}(R;\xi),\,\chi_{R\lesssim\tau}\left(\partial_{\tau}+\frac{\lambda_{\tau}}{\lambda}R\partial_R\right)\phi_{1\,<\lambda}\cdot \left(\partial_{\tau}+\frac{\lambda_{\tau}}{\lambda}R\partial_R\right)\phi_{2,\geq \lambda}\right\rangle_{L^2_{R\,dR}}\right\|_{L^2_{d\xi}}\\
	&\lesssim \hbar_3^{2-\delta}\tau^{\delta}\cdot \left\|\partial_R\left[\chi_{R\lesssim\tau}\left(\partial_{\tau}+\frac{\lambda_{\tau}}{\lambda}R\partial_R\right)\phi_{1\,<\lambda}\cdot \left(\partial_{\tau}+\frac{\lambda_{\tau}}{\lambda}R\partial_R\right)\phi_{2,\geq \lambda}\right]\right\|_{L^2_{R\,dR}}\\
	&\lesssim \tau^{\delta}\lambda^{\frac{\delta}{2}}\cdot \prod_{j=1,2}\left\|\phi_j\right\|_{\tilde{S}_0^{(n_j)}\cap\tilde{S}_1^{(n_j)}},
	\end{align*}
	which can also be summed over dyadic $\lambda\lesssim 1$, giving the desired bound. This concludes the case {\it{(1.b)}} in the small output regime. 
	\\
	
	{\it{Intermediate output frequency $\hbar_3^{-2}\gtrsim \xi\gtrsim 1$.}} This case can be handled similarly as the preceding case since the weight in the norm $\|\cdot\|_{\Shth_1}$ is the same as in the small frequency regime. However, since now the frequency $\lambda\geq 1$, we must proceed differently when we sum over $\lambda$. In the non-oscillatory regime $\hbar_{3}R\lambda^{\frac12}\leq\frac12$ we use the bound
\begin{align}\label{eq: L infty rapid decay non osc}
	\left\|\chi_{\hbar_{1}R\lambda^{\frac12}\leq\frac12}\left(\partial_{\tau}+\frac{\lambda_{\tau}}{\lambda}R\partial_{R}\right)\phi_{1,\leq \lambda}\right\|_{L^{\infty}_{RdR}}\lesssim \hbar_{1}^{-2+\delta}\lambda^{\frac{\delta}{2}}\cdot\left(\frac12\right)^{c\hbar_{1}^{-1}}\left\|\phi_{1}\right\|_{\tilde{S}_{1}^{(n_{1})}\cap\tilde{S}_{0}^{(n_{1})}}.
\end{align}	
The rapid decaying factor $\left(\frac12\right)^{c\hbar_{1}^{-1}}$ absorbs the growth from summing over $1\leq\lambda\lesssim \hbar_{3}^{-2}$, using the fact $\hbar_{1}\simeq\hbar_{3}$.  For the oscillatory regime $\hbar_{3}R\lambda^{\frac12}\geq \frac12$, in addition to the estimate \eqref{eq: L infty rapid decay non osc} (without the rapid decaying factor $\left(\frac12\right)^{c\hbar_{1}^{-1}}$), we also need to use the following refined estimates for $\left(\partial_{\tau}+\frac{\lambda_{\tau}}{\lambda}R\partial_{R}\right)\phi_{2,\geq\lambda}$:
\begin{align*}
	\left\|\left(\partial_{\tau}+\frac{\lambda_{\tau}}{\lambda}R\partial_{R}\right)\phi_{2,\geq\lambda}\right\|_{L^{\infty}_{RdR}}\lesssim \lambda^{-\frac12}\left\|\phi_{2}\right\|_{\tilde{S}_{1}^{(n_{2})}\cap\tilde{S}_{0}^{(n_{2})}},\quad \left\|\partial_{R}\left(\partial_{\tau}+\frac{\lambda_{\tau}}{\lambda}R\partial_{R}\right)\phi_{2,\geq\lambda}\right\|_{L^{2}_{RdR}}\lesssim\lambda^{-\frac12-\frac{\delta}{2}}\left\|\phi_{2}\right\|_{\tilde{S}_{1}^{(n_{2})}\cap\tilde{S}_{0}^{(n_{2})}},
\end{align*}
which are obtained by Lemma \ref{lem:exceptionalnderivatie2 large}. The decaying factors $\lambda^{-\frac12}$ absorb the growth from summing over $1\leq \lambda\lesssim \hbar_{3}^{-2}$.
	\\
	
	{\it{Large output frequency $\xi\gtrsim \hbar_3^{-2}$.}} This is accomplished by integration by parts: schematically we have 
	\begin{equation}\label{eq:timederivprodHLH}\begin{split}
	& \chi_{\xi\simeq\lambda}\left\langle \phi_{n_{3}}(R;\xi),\,\chi_{R\lesssim\tau}\left(\partial_{\tau}+\frac{\lambda_{\tau}}{\lambda}R\partial_R\right)\phi_{1\,<\lambda}\cdot \left(\partial_{\tau}+\frac{\lambda_{\tau}}{\lambda}R\partial_R\right)\phi_{2,\geq \lambda}\right\rangle_{L^2_{R\,dR}}\\
	& = \xi^{-1} \chi_{\xi\simeq\lambda}\left\langle \phi_{n_{3}}(R;\xi),\,H_{n_3}\left[\chi_{R\lesssim\tau}\left(\partial_{\tau}+\frac{\lambda_{\tau}}{\lambda}R\partial_R\right)\phi_{1\,<\lambda}\cdot \left(\partial_{\tau}+\frac{\lambda_{\tau}}{\lambda}R\partial_R\right)\phi_{2,\geq \lambda}\right]\right\rangle_{L^2_{R\,dR}}\\
	& = \xi^{-1} \chi_{\xi\simeq\lambda}\left\langle \phi_{n_{3}}(R;\xi),\,H_{n_3}\left[\chi_{R\lesssim\tau}\left(\partial_{\tau}+\frac{\lambda_{\tau}}{\lambda}R\partial_R\right)\phi_{1\,<\lambda}\right]\cdot \left(\partial_{\tau}+\frac{\lambda_{\tau}}{\lambda}R\partial_R\right)\phi_{2,\geq \lambda}\right\rangle_{L^2_{R\,dR}}\\
	& +  \sum_{i+j=2,i\leq 1}\xi^{-1} \chi_{\xi\simeq\lambda}\left\langle \phi_{n_{3}}(R;\xi),\,\left(2\partial_R+\frac1R\right)^i\left[\chi_{R\lesssim\tau}\left(\partial_{\tau}+\frac{\lambda_{\tau}}{\lambda}R\partial_R\right)\phi_{1\,<\lambda}\right]\cdot \partial_R^j\left(\partial_{\tau}+\frac{\lambda_{\tau}}{\lambda}R\partial_R\right)\phi_{2,\geq \lambda}\right\rangle_{L^2_{R\,dR}}
	\end{split}\end{equation}
	Taking advantage of Proposition~\ref{prop:singularmultiplier} as well as Lemma~\ref{lem:exceptionalnderivatie2 large}, we can bound the first term on the right: 
	\begin{align*}
	&\left\|\xi^{-1} \chi_{\xi\simeq\lambda}\left\langle \phi_{n_{3}}(R;\xi),\,H_{n_3}\left[\chi_{R\lesssim\tau}\left(\partial_{\tau}+\frac{\lambda_{\tau}}{\lambda}R\partial_R\right)\phi_{1\,<\lambda}\right]\cdot \left(\partial_{\tau}+\frac{\lambda_{\tau}}{\lambda}R\partial_R\right)\phi_{2,\geq \lambda}\right\rangle_{L^2_{R\,dR}}\right\|_{\Shth_1}\\
	&\lesssim \hbar_3^{2-\delta}\lambda^{-\frac12 - \frac{\delta}{2}}\left\langle\lambda\hbar_3^2\right\rangle^{\delta+\frac32}\cdot \left\|H_{n_3}\left[\chi_{R\lesssim\tau}\left(\partial_{\tau}+\frac{\lambda_{\tau}}{\lambda}R\partial_R\right)\phi_{1\,<\lambda}\right]\right\|_{L^2_{R\,dR}}\cdot \left\|\left(\partial_{\tau}+\frac{\lambda_{\tau}}{\lambda}R\partial_R\right)\phi_{2,\geq \lambda}\right\|_{L^\infty_{R\,dR}}\\
	&\lesssim \lambda^{-\frac12}\prod_{j=1,2}\left\|\tilde{\phi}_j\right\|_{\tilde{S}_0^{(n_j)}\cap\tilde{S}_1^{(n_j)}},
	\end{align*}
	where we have taken advantage of the bounds (with $\xb_j$ denoting the Fourier transforms as explained in Prop.~\ref{prop:bilin2})
	\begin{align*}
	&\hbar_3^{2-\delta}\lambda^{-\frac12 - \frac{\delta}{2}}\left\langle\lambda\hbar_3^2\right\rangle^{\delta+\frac32}\cdot \left\|H_{n_3}\left[\chi_{R\lesssim\tau}\left(\partial_{\tau}+\frac{\lambda_{\tau}}{\lambda}R\partial_R\right)\phi_{1\,<\lambda}\right]\right\|_{L^2_{R\,dR}}\lesssim \left[\left\|\mathcal{D}_{\tau}^{(\hbar_1)}\xb_1\right\|_{\Sho_{1}} + \left\|\mathcal{K}_{\hbar_1}^{(0)}\xb_1\right\|_{\Sho_{1}}\right]
	\\&\hspace{10cm}\lesssim \left\|\tilde{\phi}_1\right\|_{\tilde{S}_0^{(n_1)}\cap\tilde{S}_1^{(n_1)}}\\
	&\left\|\left(\partial_{\tau}+\frac{\lambda_{\tau}}{\lambda}R\partial_R\right)\phi_{2,\geq \lambda}\right\|_{L^\infty_{R\,dR}}\lesssim \lambda^{-\frac12}\cdot \left[\left\|\mathcal{D}_{\tau}^{(\hbar_1)}\xb_2\right\|_{S_1^{(n_1)}} + \left\|\mathcal{K}_{(n_1)}^{(0)}\xb_1\right\|_{S_1^{(n_1)}}\right]\lesssim  \lambda^{-\frac12}\cdot \left\|\tilde{\phi}_2\right\|_{\tilde{S}_0^{(n_2}\cap\tilde{S}_1^{(n_2)}}.
	\end{align*}
	In particular, summing over dyadic $\lambda\gtrsim \hbar_3^{-2}$ furnishes a bound for the first of the last two terms in \eqref{eq:timederivprodHLH} of the desired form. 
	\\
	As for the final term in \eqref{eq:timederivprodHLH}, we get (under our current hypothesis on $\left(\partial_{\tau}+\frac{\lambda_{\tau}}{\lambda}R\partial_R\right)\phi_{2}$) the schematic (only keeping the term $\calD_{\tau}^{(n_{2})}\xb_{2}$) decompositions
		\begin{align*}
		\partial_R\left(\partial_{\tau}+\frac{\lambda_{\tau}}{\lambda}R\partial_R\right)\phi_{2,\geq \lambda} &= \partial_R\phi_{n_{2}}(R)\cdot \left(\int_0^R\left[\phi_{n_{2}}(s)\right]^{-1}\int_0^\infty\chi_{\xi\geq \lambda}\phi_{n_2}(s,\xi)\mathcal{D}_{\tau}^{(n_2)}\xb_2(\tau,\xi)\tilde{\rho}_{n_2}(\xi)\,d\xi\right)\,ds\\
		& + \int_0^\infty \chi_{\xi\geq \lambda}\phi_{n_2}(R,\xi)\mathcal{D}^{(n_{2})}_{\tau}\xb_2(\tau,\xi)\tilde{\rho}_{n_2}(\xi)\,d\xi,
		\end{align*}
		\begin{align*}
		\partial_R^2\left(\partial_{\tau}+\frac{\lambda_{\tau}}{\lambda}R\partial_R\right)\phi_{2,\geq \lambda} &= \partial_R^2\phi_{n_2}(R)\cdot \left(\int_0^R\left[\phi_{n_2}(s)\right]^{-1}\int_0^\infty\chi_{\xi\geq \lambda}\phi_{n_2}(s,\xi)\mathcal{D}_{\tau}^{(n_2)}\xb_2(\tau,\xi)\tilde{\rho}_{n_2}(\xi)\,d\xi\right)\,ds\\
		& + \partial_R\phi_{n_{2}}(R)\cdot\left[\phi_{n_{2}}(R)\right]^{-1}\cdot \left(\int_0^\infty\chi_{\xi\geq \lambda}\phi_{n_{2}}(R,\xi)\mathcal{D}_{\tau}^{(n_2)}\xb_2(\tau,\xi)\tilde{\rho}_{n_2}(\xi)\,d\xi\right)\\
		& + \partial_R\left(\int_0^\infty\chi_{\xi\geq \lambda}\phi_{n_{2}}(R,\xi)\mathcal{D}_{\tau}^{(n_2)}\xb_2(\tau,\xi)\tilde{\rho}_{n_2}(\xi)\,d\xi\right).
		\end{align*}
	
	For $j=0,1$ we directly use Lemma \ref{lem:exceptionalnderivatie2 large} to obtain for $\lambda>1$
	\begin{align*}
	\left\|\partial_R\left(\partial_{\tau}+\frac{\lambda_{\tau}}{\lambda}R\partial_R\right)\phi_{2,\geq \lambda}\right\|_{L^2_{R\,dR}}\lesssim \sum_{\mu\geq \lambda}\mu^{-1-\frac{\delta}{2}}\left[\left\|\mathcal{D}_{\tau}^{(n_2)}\xb_2(\tau,\xi)\right\|_{S_1^{(n_2)}(\xi\simeq\mu)} + \left\|\mathcal{K}_{n_2}^{(0)}\xb_2\right\|_{S_1^{(n_2)}(\xi\simeq\mu)}\right].
	\end{align*}
	For $j=-1$, Lemma \ref{lem:exceptionalnderivatie2 large} gives the desired estimate for $R\lesssim \lambda^{-\frac12}$. For $R\lambda^{\frac12}\gtrsim 1$ we will treat differently at the end. 
	
	Furthermore, for the high-angular momentum term, in light of Lemma~\ref{lem:derLinfty} we have the $L^\infty$-bound (for $i\leq 1$)
	\begin{align*}
	\left\|\left(\partial_R+\frac1R\right)^i\left[\chi_{R\lesssim\tau}\left(\partial_{\tau}+\frac{\lambda_{\tau}}{\lambda}R\partial_R\right)\phi_{1\,<\lambda}\right]\right\|_{L^\infty_{R\,dR}}
	\lesssim \hbar_1^{-\frac52-\delta}\cdot \left[\left\|\calD_\tau^{(\hbar_1)} \xb_1\right\|_{\Sho_{1}} + \left\|\mathcal{K}_{\hbar_1}^{(0)}\xb_1\right\|_{\Sho_{1}}\right]
	\end{align*}
Combining these bounds leads to the estimate (except for $j=-1$ with $R\lambda^{\frac12}\gtrsim1$)
\begin{align*}
&\left\|\left(\text{last term of \eqref{eq:timederivprodHLH}}\right)\right\|_{\Shth_1(\xi\gtrsim\hbar_3^{-2})}\\
&\lesssim \hbar_3^{5+\frac{\delta}{2}}\cdot \left(\sum_{\lambda\gtrsim \hbar_3^{-2}}\lambda^{-2}\lambda^{2+\delta}\left\|(\text{last term of \eqref{eq:timederivprodHLH}})\right\|_{L^{2}_{d\xi}}^2\right)^{\frac12}\\
&\lesssim \hbar_1\cdot \left[\left\|\calD_\tau^{(\hbar_1)} \xb_1\right\|_{S_1^{(\hbar_1)}} + \left\|\mathcal{K}_{\hbar_1}^{(0)}\xb_1\right\|_{\Sho_{1}}\right]\\&\hspace{4cm}\cdot \left(\sum_{\lambda\gtrsim\hbar_{3}^{-2}}\lambda^{-2}\left(\sum_{\mu\geq \lambda}\left(\frac{\lambda}{\mu}\right)^{1+\frac{\delta}{2}}\left[\left\|\mathcal{D}_{\tau}^{(n_2)}\xb_2(\tau,\xi)\right\|_{S_1^{(n_2)}(\xi\simeq\mu)} + \left\|\mathcal{K}_{n_2}^{(0)}\xb_2\right\|_{S_1^{(n_2)}(\xi\simeq\mu)}\right]\right)^2\right)^{\frac12}
\end{align*}
For $j=-1$ and $R\lambda^{\frac12}\gtrsim1$, we simply write $\lambda^{-2}=\lambda^{-1}\cdot\lambda^{-1}$ and use one of the two $\lambda^{-1}$ to make up the discrepancy of the decay in $\mu$. Applying the Cauchy-Schwarz inequality and exploiting orthogonality allows to bound the preceding by 
\[
\lesssim \hbar_1\cdot\prod_{j=1,2}\left\|\phi_j\right\|_{\tilde{S}_1^{(n_j)}\cap \tilde{S}_0^{(n_j)}}
\]
This concludes the case of large output frequencies for the case {\it{(1.b)}}. 
\\

{\it{(2): $\xi\geq \max\{\xi_1,\xi_2\}$. Output frequency dominates both input frequencies.}} We proceed in analogy to case {\it{(2)}} in the proof Proposition~\ref{prop:bilin1}. The case of output frequency $\xi\lesssim 1$ here is handled exactly like the small output frequency case in {\it{(1.b)}} before. We shall henceforth restrict to output frequency $\xi\gtrsim 1$. We shall again exploit multi-fold integration by parts in order to shift derivatives around. Precisely, we write, always keeping the assumed underlying fine structure of $\left(\partial_{\tau}+\frac{\lambda_{\tau}}{\lambda}R\partial_R\right)\phi_{2,<\lambda}$ in mind
\begin{equation}\label{eq:bilintimederintbyparts}\begin{split}
&\chi_{\xi\simeq\lambda}\left\langle \phi_{n_{3}}(R;\xi),\,\chi_{R\lesssim\tau}\left(\partial_{\tau}+\frac{\lambda_{\tau}}{\lambda}R\partial_R\right)\phi_{1,<\lambda}\cdot \left(\partial_{\tau}+\frac{\lambda_{\tau}}{\lambda}R\partial_R\right)\phi_{2,<\lambda}\right\rangle_{L^2_{R\,dR}}\\
& = \xi^{-3}\chi_{\xi\sim\lambda}\left\langle \phi(R;\xi,\hbar_3),\,H_{n_3}^3\left[\chi_{R\lesssim\tau}\left(\partial_{\tau}+\frac{\lambda_{\tau}}{\lambda}R\partial_R\right)\phi_{1,<\lambda}\cdot \left(\partial_{\tau}+\frac{\lambda_{\tau}}{\lambda}R\partial_R\right)\phi_{2,<\lambda}\right]\right\rangle_{L^2_{R\,dR}}\\
& = \sum_{i+j = 6}C_{i,j}\xi^{-3}\chi_{\xi\simeq\lambda}\left\langle \phi_{n_{3}}(R;\xi),\,\partial_R^i \left(\frac{n_3}{R}\right)^j\left[\chi_{R\lesssim\tau}\left(\partial_{\tau}+\frac{\lambda_{\tau}}{\lambda}R\partial_R\right)\phi_{1,<\lambda}\right]\cdot \left(\partial_{\tau}+\frac{\lambda_{\tau}}{\lambda}R\partial_R\right)\phi_{2,<\lambda}\right\rangle_{L^2_{R\,dR}}\\
& + \sum_{\substack{l+i+j = 6\\l\geq 1}}\xi^{-3}\chi_{\xi\simeq\lambda}\left\langle \phi_{n_{3}}(R;\xi),\,\partial_R^i \left(\frac{n_3}{R}\right)^j\left[\chi_{R\lesssim\tau}\left(\partial_{\tau}+\frac{\lambda_{\tau}}{\lambda}R\partial_R\right)\phi_{1,<\lambda}\right]\cdot\partial_R^l\left(\phi_n(R)\right)\cdot \int_0^R\left[\phi_n(s)\right]^{-1}f_{<\lambda}(s)\,ds\right\rangle_{L^2_{R\,dR}}\\
& + \sum_{\substack{l+i+j+k = 6\\l\geq 1}}\xi^{-3}\chi_{\xi\simeq\lambda}\left\langle \phi_{n_{3}}(R;\xi),\,\partial_R^i \left(\frac{n_3}{R}\right)^j\left[\chi_{R\lesssim\tau}\left(\partial_{\tau}+\frac{\lambda_{\tau}}{\lambda}R\partial_R\right)\phi_{1,<\lambda}\right]\cdot\partial_R^{l-1}\left[\partial_R^k\left(\phi_n(R)\right)\left[\phi_n(R)\right]^{-1}\cdot f_{<\lambda}(R)\right]\right\rangle_{L^2_{R\,dR}},
\end{split}\end{equation}
where we have introduced the quantity 
\[
f_{<\lambda}(R): = \int_0^\infty\chi_{\xi<\lambda}\,\phi_n(R,\xi)\left(\mathcal{D}_{\tau}^{(n_2)}\xb_2(\tau,\xi)+\calK_{\hbar_{2}}^{(0)}\xb_2(\tau,\xi)\right)\tilde{\rho}_{n_2}(\xi)\,d\xi. 
\]

Then we estimate the last three terms as follows: for the first term at the end in \eqref{eq:bilintimederintbyparts}, we use Lemma~\ref{lem:exceptionalnderivative2} as well Proposition~\ref{prop:singularmultiplier} to bound the first, respectively the second term below
\[
\left\|\left(\partial_{\tau}+\frac{\lambda_{\tau}}{\lambda}R\partial_R\right)\phi_{2,<\lambda}\right\|_{L^\infty_{R\,dR}}, \quad\partial_R^i \left(\frac{n_3}{R}\right)^j\left[\chi_{R\lesssim\tau}\left(\partial_{\tau}+\frac{\lambda_{\tau}}{\lambda}R\partial_R\right)\phi_{1,<\lambda}\right]. 
\]
Invoking Plancherel's theorem for the distorted Fourier transform, this leads to the bound 
\begin{align*}
&\left\|\xi^{-3}\chi_{\xi\simeq\lambda}\left\langle \phi_{n_{3}}(R;\xi),\,\partial_R^i \left(\frac{n_3}{R}\right)^j\left[\chi_{R\lesssim\tau}\left(\partial_{\tau}+\frac{\lambda_{\tau}}{\lambda}R\partial_R\right)\phi_{1,<\lambda}\right]\cdot \left(\partial_{\tau}+\frac{\lambda_{\tau}}{\lambda}R\partial_R\right)\phi_{2,<\lambda}\right\rangle_{L^2_{R\,dR}}\right\|_{\Shth_1}\\
&\lesssim  \frac{\hbar_3^{2-\delta}\lambda^{\frac12-\frac{\delta}{2}}\langle\lambda\hbar_3^2\rangle^{\delta + \frac32}}{\lambda^3}\cdot \left\|\partial_R^i \left(\frac{n_3}{R}\right)^j\left[\chi_{R\lesssim\tau}\left(\partial_{\tau}+\frac{\lambda_{\tau}}{\lambda}R\partial_R\right)\phi_{1,<\lambda}\right]\right\|_{L^2_{R\,dR}}\cdot \left\|\left(\partial_{\tau}+\frac{\lambda_{\tau}}{\lambda}R\partial_R\right)\phi_{2,<\lambda}\right\|_{L^\infty_{R\,dR}}\\
&\lesssim \frac{\hbar_3^{2-\delta}\lambda^{\frac12-\frac{\delta}{2}}\left\langle\lambda\hbar_3^2\right\rangle^{\delta + \frac32}}{\lambda^3}\cdot \left(\sum_{\mu<\lambda} \left\|\xi^3\mathcal{D}_{\tau}^{(\hbar_1)}\xb_1\right\|_{L^2_{d\xi}(\xi\simeq\mu)} + \left\|\xi^3\mathcal{K}_{\hbar_1}^{(0)}\xb_1\right\|_{L^2_{d\xi}(\xi\simeq\mu)}\right)\cdot\left\|\phi_{2}\right\|_{\tilde{S}_{0}^{(n_{2})}\cap\tilde{S}_1^{(n_2)}}.
\end{align*}
Here we can bound the product of the first two expressions by (recalling the hypothesis $\hbar_1\simeq\hbar_3$)
\begin{equation}\label{eq:neededlater1}\begin{split}
&\frac{\hbar_3^{2-\delta}\lambda^{\frac12-\frac{\delta}{2}}\langle\lambda\hbar_3^2\rangle^{\delta + \frac32}}{\lambda^3}\cdot \left(\sum_{\mu<\lambda} \left\|\xi^3\mathcal{D}_{\tau}^{(\hbar_1)}\xb_1\right\|_{L^2_{d\xi}(\xi\simeq\mu)} + \left\|\xi^3\mathcal{K}_{\hbar_1}^{(0)}\xb_1\right\|_{L^2_{d\xi}(\xi\simeq\mu)}\right)\\
&\lesssim \sum_{\mu<\lambda} \left(\frac{\mu}{\lambda}\right)^{\frac12-\frac{\delta}{2}}\cdot\left[\left\|\mathcal{D}_{\tau}^{(\hbar_1)}\xb_1\right\|_{\Sho_{1}(\xi\simeq \mu)} + \left\|\xb_1\right\|_{\Sho_{0}(\xi\simeq \mu)}\right].
\end{split}\end{equation}
Substituting this bound into the preceding one, we easily infer that 
\begin{align*}
&\left(\sum_{\lambda\gtrsim 1}\left\|(\text{First term of \eqref{eq:bilintimederintbyparts}})\right\|_{\Shth_{1}}^2\right)^{\frac12}\\
&\lesssim \left\|\phi_{2}\right\|_{\tilde{S}_{0}^{(n_{2})}\cap\tilde{S}_1^{(n_2)}}\cdot \left(\sum_{\lambda\gtrsim 1}\left[ \sum_{\mu<\lambda} \left(\frac{\mu}{\lambda}\right)^{\frac12-\frac{\delta}{2}}\cdot\left[\left\|\mathcal{D}_{\tau}^{(\hbar_1)}\xb_1\right\|_{\Sho_{1}(\xi\simeq \mu)} + \left\|\xb_1\right\|_{\Sho_{0}(\xi\simeq \mu)}\right]\right]^2\right)^{\frac12}\\
&\lesssim \prod_{j=1,2}\left\|\phi_j\right\|_{\tilde{S}_1^{(n_j)}\cap \tilde{S}_0^{(n_j)}},
\end{align*}
which is as desired, and concludes the bound for the first term in \eqref{eq:bilintimederintbyparts}. For the second term, Lemma~\ref{lem:exceptionalnderivative2} implies 
\[
\left\|\partial_R^l\left(\phi_{n_{2}}(R)\right)\cdot \int_0^R\left[\phi_{n_{2}}(s)\right]^{-1}f_{<\lambda}(s)\,ds\right\|_{L^{\infty}_{R\,dR}}\lesssim \left\|\mathcal{D}^{(n_{2})}_{\tau}\xb_1\right\|_{S_1^{(n_2)}}+\left\|\mathcal{K}^{(0)}_{\hbar_{2}}\xb_1\right\|_{S_1^{(n_2)}},\quad l\geq 0, 
\]
As in the previous case, with $i+j\leq 6$, we have
\[
\frac{\hbar_3^{2-\delta}\lambda^{\frac12-\frac{\delta}{2}}\langle\lambda\hbar_3^2\rangle^{\delta + \frac32}}{\lambda^3}\left\|\partial_R^i \left(\frac{n_3}{R}\right)^j\left[\chi_{R\lesssim\tau}(\partial_{\tau}+\frac{\lambda_{\tau}}{\lambda}R\partial_R)\phi_{1,< \lambda}\right]\right\|_{L^2_{R\,dR}}
\]
can be bounded by \eqref{eq:neededlater1}. Applying Plancherel's theorem as well as Holder's inequality suitably we then infer 
\begin{align*}
&\left(\sum_{\lambda\gtrsim 1}\left\|(\text{Second term of \eqref{eq:bilintimederintbyparts}})\right\|_{\Shth_1}^2\right)^{\frac12}\\
&\lesssim \left\|\phi_{2}\right\|_{\tilde{S}_{0}^{(n_{2})}\cap\tilde{S}_1^{(n_2)}}\cdot \left(\sum_{\lambda\gtrsim 1}\left[ \sum_{\mu<\lambda} \left(\frac{\mu}{\lambda}\right)^{\frac12-\frac{\delta}{2}}\cdot\left[\left\|\mathcal{D}_{\tau}^{(\hbar_1)}\xb_1\right\|_{\Sho_{1}(\xi\simeq \mu)} + \left\|\xb_1\right\|_{\Sho_{0}(\xi\simeq \mu)}\right]\right]^2\right)^{\frac12},
\end{align*}
which implies the desired bound. 
\\
The final term in \eqref{eq:bilintimederintbyparts} is bounded by applying Leibniz' rule to the last term. Considering the case $n = -1$, the cases $n = +1,0$ being handled similarly, we write 
\begin{align*}
&\partial_R^{l-1}\left[\partial_R^k\left(\phi_n(R)\right)\left[\phi_n(R)\right]^{-1}\cdot f_{<\lambda}(R)\right]\\& = \sum_{l_1+l_2+l_3=l-1}C_{l_{1,2,3}}\partial_R^{k+l_1}\left(\phi_n(R)\right)\cdot \partial_R^{l_2}\left([\phi_n(R)]^{-1}\right)\partial_R^{l_3}\left( f_{<\lambda}(R)\right)\\
& = \sum_{l_1+l_2+l_3=l-1}g_{k,l_{1,2}}(R)\cdot \partial_R^{l_3}\left( f_{<\lambda}(R)\right), 
\end{align*}
where the function $g_{k,l_{1,2}}(R)$ satisfies the bound
\begin{align*}
\left|g_{k,l_{1,2}}(R)\right|\lesssim \max\{R^{-2-l_2}\cdot R^{\max\{2-k-l_1, 0\}}, \langle R\rangle^{-(k+l_1+l_2)}\},
\end{align*}
and so we get the bound 
\begin{align*}
&\left|\partial_R^i \left(\frac{n_3}{R}\right)^j\left[\chi_{R\lesssim\tau}\left(\partial_{\tau}+\frac{\lambda_{\tau}}{\lambda}R\partial_R\right)\phi_{1,<\lambda}\right]\cdot\partial_R^{l-1}\left[\partial_R^k\left(\phi_n(R)\right)[\phi_n(R)]^{-1}\cdot f_{<\lambda}(R)\right]\right|\\
&\lesssim \sum_{l\leq 5-i-j-k}\left|R^{-(5-i-j-l)}\partial_R^i \left(\frac{n_3}{R}\right)^j\left[\chi_{R\lesssim\tau}\left(\partial_{\tau}+\frac{\lambda_{\tau}}{\lambda}R\partial_R\right)\phi_{1,<\lambda}\right]\right|\cdot\left|\partial_R^{l}\left( f_{<\lambda}(R)\right)\right|.
\end{align*}
Taking advantage of Prop.~\ref{prop:singularmultiplier}, Lemma~\ref{lem:derLinfty} as well as Prop.~\ref{prop:exceptionalnderivative1} and Lemma~\ref{lem:derLinfty} and placing the factors 
\[
R^{-(5-i-j-l)}\partial_R^i \left(\frac{n_3}{R}\right)^j\left[\ldots\right],\quad \partial_R^{l}\left( f_{<\lambda}(R)\right)
\]
into $L_{R\,dR}^\infty, L_{R\,dR}^{2}$ or the other way round according to whether $l\geq 3$, $l<3$, we infer the bound 
\begin{align*}
&\left\|\xi^{-3}\chi_{\xi\simeq\lambda}\left\langle \phi_{n_{3}}(R;\xi),\,\partial_R^i \left(\frac{n_3}{R}\right)^j\left[\chi_{R\lesssim\tau}\left(\partial_{\tau}+\frac{\lambda_{\tau}}{\lambda}R\partial_R\right)\phi_{1,<\lambda}\right]\cdot\partial_R^{l-1}\left[\partial_R^k\left(\phi_n(R)\right)[\phi_n(R)]^{-1}\cdot f_{<\lambda}(R)\right]\right\rangle_{L^2_{R\,dR}}\right\|_{\Shth_1}\\
&\lesssim \sum_{\mu_{1,2}<\lambda}\prod_{j=1,2}\min\{\mu_j,1\}^{\frac{\delta}{2}}\cdot\left(\frac{\max\{\mu_{1,2}\}}{\lambda}\right)\cdot \left[\left\|\mathcal{D}_{\tau}^{(\hbar_1)}\xb_1\right\|_{\Sho_{1}(\xi\simeq\mu_1)} + \left\|\xb_1\right\|_{\Sho_{0}(\xi\simeq\mu_1)}\right]\\
&\cdot \left(\left\|\mathcal{D}_{\tau}^{(n_2)}\xb_2\right\|_{S_1^{(n_2)}(\xi\simeq\mu_2)}+\left\|\mathcal{K}_{\hbar_2}^{(0)}\xb_2\right\|_{S_1^{(n_2)}(\xi\simeq\mu_2)}\right)
\end{align*}
The desired bound results from here in the usual fashion by square summing over $\lambda$ and exploiting Cauchy-Schwarz as well as orthogonality. This concludes the proof of \eqref{eq:productlargen3alln20}. 
\\
The proof of the remaining inequalities \eqref{eq:Taylorbound10}, \eqref{eq:productsmalln1geq2alln20} proceeds in analogy to the proof of \eqref{eq:Taylorbound1},
\eqref{eq:productsmalln1geq2alln2}. This concludes the proof of Proposition~\ref{prop:bilin3}. 
	\end{proof}
In order to deal with the higher order nonlinear source terms we need to pass from the above basic estimates to estimates for higher order terms. For this we have the following 
\begin{proposition}\label{prop:bilin4} Assume that $F(R)$ is a function on $[0,\infty)$ admitting a third order Taylor development $\sum_{j=0}^3\gamma_j R^j$ at $R = 0$ and such that 
	\[
	\left\|F(R) - \chi_{R\lesssim1}\sum_{j=0}^3\gamma_j R^j\right\|_{\tilde{S}_1^{(n_1)}} + \sum_j\left|\gamma_j\right| =: \Lambda_1<\infty
	\]
	for some $|n_1|\geq 2$, and where $\|\cdot\|_{\tilde{S}_1^{(n_1)}}$ is defined like $\left\|\cdot\right\|_{\tilde{S}_0^{(n_1)}}$ in Proposition~\ref{prop:bilin2}, except that $\left\|\cdot\right\|_{\Sh_{0}}$ there is replaced by $\left\|\cdot\right\|_{\Sh_{1}}$. Also, assume that all $\gamma_j = 0$ provided $|n_1|\geq K$. 
	Next assume that $\phi_2$ is an angular momentum $n_2$ function (for arbitrary integral $n_2$) with  
	\[
	\left\|\phi_2\right\|_{\tilde{S}_0^{(n_2)}} = :\Lambda_2<\infty. 
	\]
	Then the function $F\phi_2$ admits a Taylor expansion $P_3 = \sum_{j=0}^3\tilde{\gamma}_jR^j$ of order three at $R = 0$ with 
	\[
	\sum_j\left|\tilde{\gamma}_j\right|\lesssim \Lambda_1\cdot\Lambda_2,
	\]
	and such that, with $|n_3|\geq 2$ and either (i) $n_1\simeq n_3$ and $|n_2|\lesssim |n_1|$, or (ii) $|n_3|\ll |n_1|$ and $n_1\simeq -n_2$, or (iii) $|n_3|\gg|n_1|$ and $n_3\simeq n_2$, we have  
	\begin{equation}\label{eq:somewhatpoorbound}
	\left\|\chi_{R\lesssim\tau}\left[F(R)\phi_2(R) -  \chi_{R\lesssim1}\sum_{j=0}^3\tilde{\gamma}_jR^j\right]\right\|_{\Shth_1}\lesssim \langle\chi_{|n_2|<K}|n_2|\rangle^{3+\delta}\langle \min\{n_1,n_2\}\rangle^2\tau^{1+\delta}\prod_{k=1,2}\Lambda_k,\,\tau\gg 1. 
	\end{equation}
	The last inequality remains correct if we subtract from $P_3$ those terms $\tilde{\gamma}_l R^l$ with $l\geq |n_3 - 1|, l\equiv n_3 - 1\text{(mod 2)}$, provided $|n_3 -1|\leq 3$. 
	Note the particular choice $n_3 = n_3(n_1,n_2)$ where $n_3(n_1,n_2) = n_1 + n_2$ if $|n_1+n_2|\geq 2$ and $n_3(n_1,n_2) = 2$ otherwise, satisfies (i) - (iii).
	\end{proposition}
\begin{proof}
	This is in fact completely analogous to the one of Prop.~\ref{prop:bilin2}, observing that we implicitly exploited there that $\partial_R\phi_j\in \tilde{S}_1^{(n_j)}$.  Then if $|n_2|\gg 1$, we shall set $\tilde{\gamma}_j = 0$ for all $j$ and we decompose 
	\begin{align*}
	\chi_{R\lesssim\tau}F(R)\phi_2(R) &= \left(\chi_{R\lesssim\tau}F(R) - \chi_{R\lesssim 1}\sum_{j\leq 3}\gamma_jR^j\right)\phi_2(R)\\
	& + \left( \chi_{R\lesssim 1}\sum_{j\leq 3}\gamma_jR^j\right)\phi_2(R).\\ 
	\end{align*}
	Here the first term on the right is bounded precisely like in the preceding argument, taking advantage of the last bound of Lemma~\ref{lem:derLinfty}(which causes the loss of a factor $\tau$ in the estimate). For the second term, one expands out 
	\[
	\phi_2(R) = \int_0^\infty \phi_{n_{2}}(R;\xi_2)\xb_2(\xi_2)\tilde{\rho}_{n_{2}}(\xi_2)\,d\xi_2
	\]
	and divides into two cases depending on the relation of the output frequency $\xi$ to $\xi_2$. Thus write
	\begin{equation}\label{eq:technical12}\begin{split}
	&\left\langle\phi_{n_{3}}(R;\xi),\, \left(\chi_{R\lesssim 1}\sum_{j\leq 3}\gamma_jR^j\right)\phi_2(R)\right\rangle_{L^2_{R\,dR}} \\
	&= \sum_{\lambda}\chi_{\xi\simeq\lambda}\left\langle\phi_{n_{3}}(R;\xi),\, \left(\chi_{R\lesssim 1}\sum_{j\leq 3}\gamma_jR^j\right)\phi_{2,\ll\lambda}(R)\right\rangle_{L^2_{R\,dR}}\\
	& +   \sum_{\lambda}\chi_{\xi\simeq\lambda}\left\langle\phi_{n_{3}}(R;\xi),\, \left(\chi_{R\lesssim 1}\sum_{j\leq 3}\gamma_jR^j\right)\phi_{2,\gtrsim\lambda}(R)\right\rangle_{L^2_{R\,dR}}
	\end{split}\end{equation}
	To deal with the second term on the right, we further localize $\phi_2$ to frequency $\mu_2\gtrsim\lambda$, which then gives for fixed $\lambda$ 
	\begin{align*}
	&\left\|\chi_{\xi\simeq\lambda}\left\langle\phi_{n_{3}}(R;\xi),\, \left(\chi_{R\lesssim 1}\sum_{j\leq 3}\gamma_jR^j\right)\phi_{2,\mu_2}(R)\right\rangle_{L^2_{R\,dR}}\right\|_{\Shth_1}\\
	&\lesssim \lambda^{-\frac12}(\lambda\hbar_3^2)^{1-\frac{\delta}{2}}\langle\lambda\hbar_3^2\rangle^{\delta + \frac32}\cdot \left\|\left(\chi_{R\lesssim 1}\sum_{j\leq 3}\gamma_jR^j\right)R^{-1}\phi_{2,\mu_2}(R)\right\|_{L^2_{R\,dR}},
	\end{align*}
	where we have inserted the extra factor $R^{-1}$ on purpose, since Proposition~\ref{prop:singularmultiplier} gives 
	\begin{align*}
	\left\|\left(\chi_{R\lesssim 1}\sum_{j\leq 3}\gamma_jR^j\right)R^{-1}\phi_{2,\mu_2}(R)\right\|_{L^2_{R\,dR}}\lesssim \left\|\xi_2^{\frac12}\xb_2(\xi_2)\right\|_{L^2(\xi_2\simeq\mu_2)},
	\end{align*}
	and so we infer 
	\begin{align*}
	&\left\|\chi_{\xi\simeq\lambda}\left\langle\phi_{n_{3}}(R;\xi),\, \left(\chi_{R\lesssim 1}\sum_{j\leq 3}\gamma_jR^j\right)\phi_{2,\mu_2}(R)\right\rangle_{L^2_{R\,dR}}\right\|_{\Shth_{1}}\\
	&\lesssim \hbar_3\mu^{\frac12}_{2}\cdot\left(\lambda\hbar_3^2\right)^{\frac12-\frac{\delta}{2}}\left\langle\lambda\hbar_3^2\right\rangle^{\delta+\frac32}\cdot \left\|\xb_2(\xi_2)\right\|_{L^2(\xi_2\simeq\mu_2)}\\
	&\lesssim\left[\chi_{|n_2|<K}|n_2|\right]^{3+\delta}\min\{n_1, n_2\}^2\left(\frac{\lambda}{\mu_2}\right)^{\frac12 - \frac{\delta}{2}}\cdot \left\|\xb_2\right\|_{\Sht_{0}(\xi_2\simeq \mu_2)},
	\end{align*}
	where we exploited our assumption that if $|n_2|\gg |n_3|$, then $K>|n_1|\simeq |n_2|$, since if $|n_1|\geq K$ the coefficients $\gamma_j$ all vanish . Square summing over $\lambda$ and exploiting Cauchy-Schwarz and orthogonality as usual leads to the desired bound. 
	Bounding the first term in \eqref{eq:technical12} involves integrating by parts in a manner similar to the one in the proof of the previous proposition, we omit the details. 
	\\
	There remains the case when all involved angular momenta $n_j$ are of small size, $|n_j|\lesssim 1$, $j = 1,2,3$. There we use a slightly refined decomposition
	\begin{align*}
	\chi_{R\lesssim\tau}F(R)\phi_2(R) &= \left(\chi_{R\lesssim\tau}F(R) - \chi_{R\lesssim 1}\sum_{j\leq 3}\gamma_jR^j\right)\phi_2(0)\\
	& +  \left(\chi_{R\lesssim\tau}F(R) - \chi_{R\lesssim 1}\sum_{j\leq 3}\gamma_jR^j\right)[\phi_2(R) - \phi_{2}(0)]\\
	& + \left( \chi_{R\lesssim 1}\sum_{j\leq 3}\gamma_jR^j\right)\phi_2(R).\\ 
	\end{align*}
	Then proceeding as in the proof of the preceding proposition one shows that the first two terms on the right satisfy the desired bound
	\begin{align*}
	&\left\| \left(\chi_{R\lesssim\tau}F(R) - \chi_{R\lesssim 1}\sum_{j\leq 3}\gamma_jR^j\right)\phi_2(0)\right\|_{\Shth_1}\lesssim \cdot\prod_{k=1,2}\Lambda_k\\
	&\left\| \left(\chi_{R\lesssim\tau}F(R) - \chi_{R\lesssim 1}\sum_{j\leq 3}\gamma_jR^j\right)\left[\phi_2(R) - \phi_{2}(0)\right]\right\|_{\Shth_1}\lesssim \tau^{1+\frac{\delta}{2}}\cdot\prod_{k=1,2}\Lambda_k.\\
	\end{align*}
	Finally, as for the remaining term, it is straightforward to check that it admits a degree Taylor expansion around $R = 0$:
	\[
	\sum_{j=1}^3\tilde{\gamma}_jR^j,\quad \sum_{j}|\tilde{\gamma}_j|\lesssim \prod_{k=1,2}\Lambda_k, 
	\]
	and such that 
	\begin{align*}
	\left\|\left( \chi_{R\lesssim 1}\sum_{j\leq 3}\gamma_jR^j\right)\phi_2(R) - \chi_{R\lesssim 1}\sum_{j=1}^3\tilde{\gamma}_jR^j\right\|_{\Shth_{1}}\lesssim \prod_{k=1,2}\Lambda_k. 
	\end{align*}
\end{proof}
A consequence of the preceding proposition is the following multilinear estimate. 
\begin{proposition}\label{prop:multilingen1} For an integer $n$ define $n_{|\cdot|\geq 2}: = n$ if $|n|\geq 2$ and $n_{|\cdot|\geq 2}: =2$ otherwise. Assume that $\phi_1,\phi_2,\ldots, \phi_k$ are angular momentum $m_j$ functions, $j = 1,2,\ldots,k$, $m_j\in \Z$, where we make the same structural assumption on each of them as in Proposition~\ref{prop:bilin4}. Further, assume that $F$ is a angular momentum $n_0$ function with $|n_0|\geq 2$,  admitting a third order Taylor development $\sum_{j=0}^3\gamma_j R^j$ at $R = 0$ and such that 
	\[
	\left\|F(R) - \chi_{R\lesssim1}\sum_{j=0}^3\gamma_j R^j\right\|_{\tilde{S}_1^{(n_0)}} + \sum_j\left|\gamma_j\right| =: \Lambda_0<\infty.
	\]
	Assume that all $\gamma_j = 0$ if $|n_0|\geq C_0$. 
	Then the function 
	\[
	F(R)\cdot\prod_{j=1}^k\phi_j(R)
	\]
	admits a third order Taylor expansion $\sum_{j=0}^3\tilde{\gamma}_j R^j$ at $R = 0$, where we have the coefficient bound 
	\[
	\sum_{j=0}^3\left|\tilde{\gamma}_j\right|\leq C_1^k\Lambda_0\cdot\prod_{j=1}^k\left\|\phi_j\right\|_{\tilde{S}_0^{(m_j)}}. 
	\]
	for a constant $C_1(C_0)$. Further, setting 
	\[
	n: = n_0 + \sum_{j=1}^k m_j,
	\]
	we have the bound 
	\begin{align*}
	&n_{|\cdot|\geq 2}^4\left\|\chi_{R\lesssim\tau}F(R)\cdot\prod_{j=1}^n\phi_j(R) - \chi_{R\lesssim1}\sum_{j=0}^3\tilde{\gamma}_j R^j\right\|_{\tilde{S}_1^{(n_{|\cdot|\geq 2})}}\\&\leq C_2^kn_0^4\Lambda_0\cdot\prod_{j=1}^k\left(\tau^{1+\delta}\langle m_j\rangle^4\left\|\phi_j\right\|_{\tilde{S}_0^{(m_j)}}\right). 
	\end{align*}
	Similar bounds obtain when the exponent $4$ is replaced by $p\geq 4$. 
\end{proposition}
\begin{proof} By reordering the factors we may assume that $|m_1|\leq|m_2|\leq \ldots\leq |m_k|$, and either (i) $|n_0|\gtrsim \frac{\left|n_{|\cdot|\geq 2}\right|}{k}$ or else (ii) $|m_{k}|\gtrsim \frac{\left|n_{|\cdot|\geq 2}\right|}{k}$. Then we apply Prop.~\ref{prop:bilin4} consecutively to the sequence of products 
	\[
	F(R)\phi_1,\,F(R)\phi_1\cdot\phi_2,\ldots, F(R)\cdot\prod_{j=1}^k\phi_j, 
	\]
	where the function $F(R)\phi_1\cdot\ldots\phi_l$ will be interpreted as angular momentum $\left(n_0+\sum_{j=1}^l m_j\right)_{|\cdot|\geq 2}$ function. Here we have to be careful not to lose on account of the factor 
	\[
	\chi_{|n_2|<K}|n_2|^{3+\delta}
	\]
	in \eqref{eq:somewhatpoorbound}, where we note that the angular momenta of these products (and thus $K$) can in principle grow. However, if one of the factors $\phi_l$ is an angular momentum $n$ function with $|n|\geq 5$,  all products $F\cdot\prod_{j=1}^{l'}\phi_j, l'\geq l$ have trivial third order Taylor polynomial around $R =0$, and so the term in the proof of Prop.~\ref{prop:bilin4} responsible for the preceding factor is not present. It follows that there is {\it{at most one case}} where one loses a factor $\chi_{|n_2|<K}|n_2|^{3+\delta}$, namely the first instance where $\phi_l$ is of absolute angular momentum $\geq 5$, and we lose here at most a factor $Cl^{3+\delta}$. Such a factor can be absorbed into the factors $C^l$ which automatically occur. 
	In case (i) we inductively get the bounds 
	\begin{equation}\label{eq:inductivemultilin1}\begin{split}
	\left\|\chi_{R\lesssim\tau}F(R)\cdot\prod_{j=1}^{l}\phi_j - \chi_{R\lesssim 1}\sum_{j=1}^3 \tilde{\gamma}^{(l)}R^j\right\|_{\tilde{S}_1^{\left((n_0+\sum_{j=1}^lm_j)_{|\cdot|\geq 2}\right)}}\leq
	C^l\Lambda_0\cdot\prod_{j=1}^l\left(\tau^{1+\delta}\langle m_j\rangle^2\left\|\phi_j\right\|_{\tilde{S}_0^{(m_j)}}\right)
	\end{split}\end{equation}
	where $\sum_{j=1}^3 \tilde{\gamma}^{(l)}R^j$ is the Taylor polynomial of $\chi_{R\lesssim\tau}F(R)\cdot\prod_{j=1}^{l}\phi_j$ at $R = 0$, and the desired conclusion follows for $l = k$ by using assumption (i) and $k^4C^k\leq C_2^k$ for suitable $C_2 = C_2(C)$. 
	\\
	In case (ii), use 
	\[
	\left|n_0 + \sum_{j=1}^{k-1}m_j\right|\leq \left|n_0\right| + (k-1)\left|m_{k-1}\right|, 
	\]
	and distinguish between $|n_0|<(k-1)\left|m_{k-1}\right|, |n_0|\geq (k-1)\left|m_{k-1}\right|$. in the former case, use \eqref{eq:inductivemultilin1} up to $l = k-1$ and further 
	\begin{align*}
	&\left\|\chi_{R\lesssim\tau}F(R)\cdot\prod_{j=1}^{k}\phi_j - \chi_{R\lesssim 1}\sum_{j=1}^3 \tilde{\gamma}^{(k)}R^j\right\|_{\tilde{S}_1^{\left((n_0+\sum_{j=1}^k m_j)_{|\cdot|\geq 2}\right)}}\\&\leq
	\langle 2(k-1)\left|m_{k-1}\right|\rangle^2C^k\Lambda_0\cdot\prod_{j=1}^{k-1}\left(\tau^{1+\delta}\langle m_j\rangle^2\left\|\phi_j\right\|_{\tilde{S}_0^{(m_j)}}\right)\cdot \left(\tau^{1+\delta}\left\|\phi_k\right\|_{\tilde{S}_0^{(m_k)}}\right)\\
	&\leq C_1^k\Lambda_0\cdot\prod_{j=1}^{k-1}\left(\tau^{1+\delta}\langle m_j\rangle^4\left\|\phi_j\right\|_{\tilde{S}_0^{(m_j)}}\right)\cdot\left(\tau^{1+\delta}\left\|\phi_k\right\|_{\tilde{S}_0^{(m_k)}}\right)\\
	\end{align*}
	from which the desired bound follows due to assumption (ii). In the case $|n_0|\geq (k-1)\left|m_{k-1}\right|$ it suffices to replace the factor $\langle 2(k-1)\left|m_{k-1}\right|\rangle$ by $2|n_0|$. 
\end{proof}

In a similar vein, we can bound products where each factor is a sum of different angular momentum components: 
\begin{corollary}\label{cor:multilingen2} Assume that each $\phi_j$ can be written as a sum of angular momentum $l$ functions $\phi_j^{(l)}$ satisfying the bound 
	\[
	\tau^{1+\delta}\sum_{l\in \Z}\langle l\rangle^4\left\|\phi_j^{(l)}\right\|_{\tilde{S}_0^{(l)}}: = \Lambda_j<\infty,\,j = 1,2,\ldots, k. 
	\]
	Further, assume that $F$ is a sum of angular momentum $l_{|\cdot|\geq 2}$-functions $F = \sum_{l\in\Z}F^{(l)}$ each of which has a third order Taylor polynomial $\sum_{j=1}^3\gamma_j^{(l)}R^j$ at $R = 0$ and such that 
	\[
	\sum_{l\in\Z}\langle l\rangle^4\left[\left\|\chi_{R\lesssim\tau}F^{(l)} - \chi_{R\lesssim 1}\sum_{j=1}^3\gamma_j^{(l)}R^j\right\|_{\tilde{S}_1^{(l_{|\cdot|\geq 2})}} + \sum_{j=1}^3\left|\gamma_j^{(l)}\right|\right]: = \Lambda_0<\infty. 
	\]
	If we then set 
	\begin{align*}
	\Psi^{(n)}_k: = \sum_{l+\sum_{j=1}^k m_j = n}F^{(l)}\cdot \prod_{j=1}^k \phi_j^{(m_j)},
	\end{align*}
	then each $\Psi^{(n)}_k$ admits a third order Taylor expansion $\sum_{j=0}^3\tilde{\gamma}^{(n)}_jR^j$ around $R = 0$ and we have, with $C_3$ a universal constant, the bounds 
	\begin{align*}
	&\sum_{n\in \Z}\langle n\rangle^4\left\|\chi_{R\lesssim\tau}\Psi^{(n)}_k - \chi_{R\lesssim1}\sum_{j=0}^3\tilde{\gamma}^{(n)}_jR^j\right\|_{\tilde{S}_1^{(n_{|\cdot|\geq 2})}}\leq C_3^k\Lambda_0\prod_{j=1}^k\left(\tau^{1+\delta}\Lambda_j\right).
	\end{align*}
	The same bound obtains if we subtract from the $\sum_{j=0}^3\tilde{\gamma}^{(n)}_jR^j$ those terms $\tilde{\gamma}^{(n)}_lR^l$ with $l\geq |n_{|\cdot|\geq 2}-1|, l\equiv (n_{|\cdot|\geq 2}-1)\text{mod}(2)$, provided $|n_{|\cdot|\geq 2}-1|\leq 3$. 
	The exponent $4$ may be replaced by $p\geq 4$ throughout.  
\end{corollary}
\subsection{Control over all source terms for angular momenta $|n|\geq 2$ and away from the light cone}\label{subsec:allsourcetermsngeq2}
We use here the results of the preceding section to bound all the source terms arising in the equations for the angular momenta $|n|\geq 2$ in the interior of the light cone, away from the shock region. Dealing with the latter will require a different set of estimates exploiting the fine structure of the shock on the light cone. We will be using the equations \eqref{eq:RegularFinestructure1}-\eqref{eq:RegularFinestructure3}. We observe that all terms on the right hand side of \eqref{eq:RegularFinestructure2} with the exception of the first term are linear in $\varepsilon_{1,2}$, and in fact only depend on $\varepsilon_{1,2}(n)$. The fine structure of the correction $\epsilon = \epsilon(\tau, R)$ in turn is adopted from the works \cite{GaoK,KST}, and we adopt the following function space essentially from the latter reference:
\begin{definition}\label{def:KST1} Let $b_1: = \frac{\left(\log(1+R^2)\right)^2}{(t\lambda)^2},\,b_2 = \frac{1}{(t\lambda)^2}$. Then we denote by 
	\[
	S^m\left(R^k(\log R)^l\right)
	\]
	the class of analytic functions $v: [0,\infty)\times [0, b_0]^2\longrightarrow \R$, $b_0>0$ a fixed small positive number, such that 
	\begin{itemize}
		\item $v(R, b_1, b_2)$ vanishes to order $m$ at $R =0$, where we interpret $b_1$ as a function of $R$ and $t$.
		\item $v$ admits a convergent expansion at $R = \infty$
		\[
		v(R, b_1, b_2) = \sum_{\substack{0\leq j\leq l+i\\ i\geq 0}}c_{ij}(b_1, b_2)R^{k-i}(\log R)^j
		\]
		and the functions $c_{ij}$ are analytic on $[0,b_0]^2$. 
	\end{itemize}
\end{definition}
We observe that this class of functions captures precisely the corrections $\epsilon$ which define $Q$ in the region away from the light cone, say $r<\frac{t}{2}$, provided $m = 3, k = 1$. Specifically, we have
\\

\begin{theorem}\label{thm:interiorepsilonstructure}(\cite{GaoK,KST}) The correction $\epsilon$ can be chosen in the form 
	\begin{align*}
	\chi_{r<\frac{t}{2}}\epsilon &= \epsilon_1 + \epsilon_2\\
	&\epsilon_1 = \sum_{i=1}^N v_i,\quad v_{2k-1}\in \frac{1}{(t\lambda)^{2k}}S^3\left(R(\log R)^{2k-1}\right),\quad v_{2k}\in \frac{1}{(t\lambda)^{2k+2}}S^3\left(R^3(\log R)^{2k-1}\right)
	\end{align*}
	where $N$ can be chosen arbitrarily large, and 
	such that $\epsilon_2$ is a $C^\infty$ function admitting a Taylor expansion around $R = 0$ in terms of odd powers of $R$, as well as the bounds 
	\begin{align*}
	\left\|\nabla_{\tau, R}^l\epsilon_2\right\|_{L^\infty_{dR}}\lesssim_l \frac{1}{(t\lambda)^{N_1}},
	\end{align*}
	for any $l\geq 0$, where $N_1 = N_1(N)$ grows linearly in $N$. 
\end{theorem}
Using this structural ingredient, we can immediately bound all the linear terms in \eqref{eq:RegularFinestructure2}. 
\begin{proposition}\label{prop:smoothlinearsource} Let $|n|\geq 2$ and assume that $\varepsilon$ is an angular momentum $n$ function. Then denoting $F_j$, $j = 1,\ldots, 4$ the final four terms in \eqref{eq:RegularFinestructure2}, we have the bounds (with $\hbar = \frac{1}{|n|+1}$)
	\begin{align*}
	\left\|\left\langle \phi(R;\xi,\hbar), \chi_{r<\frac{t}{2}}F_j\right\rangle_{L^2_{R\,dR}}\right\|_{\Sh_{1}}\lesssim \frac{1}{\tau^2}\sum_{\pm}\left\|\varepsilon_{\pm}(n)\right\|_{\tilde{S}_0^{(n)}},
	\end{align*}
	where we recall that $\phi_2(n) = i\cdot\left(\varepsilon_+(n) - \varepsilon_-(n)\right)$, $\tau = \int_t^\infty \lambda(s)\,ds$. 
\end{proposition}
\begin{proof}
	We give details for the estimate of the second term $F_2$, as the others can be handled similarly. To begin with, Theorem~\ref{thm:interiorepsilonstructure} implies that we have the bound 
	\begin{align*}
	\left\|4\frac{\sin\left[Q+\frac{\eps}{2}\right]\sin\left[\frac{\eps}{2}\right]}{R^2}\right\|_{L^2_{R\,dR}\cap L^\infty_{R\,dR}}\lesssim \frac{1}{(t\lambda)^2},\quad\left\|4\langle R\rangle\frac{\sin\left[Q+\frac{\eps}{2}\right]\sin\left[\frac{\eps}{2}\right]}{R^2}\right\|_{L^\infty_{R\,dR}}\lesssim \frac{1}{(t\lambda)^2}
	\end{align*}
	Also, we can obviously write $i\varepsilon_{\pm, \theta}(n) = -n\varepsilon_{\pm}(n)$. Then we localize the output frequency to dyadic size $\xi\simeq\lambda$ and split the source term into the sum of two pieces: 
	\begin{equation}\label{eq:proofofpropsmoothlinearsource1}
	\chi_{r<\frac{t}{2}}F_2 = \chi_{r<\frac{t}{2}}4n\frac{\sin\left[Q+\frac{\eps}{2}\right]\sin\left[\frac{\eps}{2}\right]}{R^2}\varepsilon_{\pm,\geq \lambda}(n) + \chi_{r<\frac{t}{2}}4n\frac{\sin\left[Q+\frac{\eps}{2}\right]\sin\left[\frac{\eps}{2}\right]}{R^2}\varepsilon_{\pm,< \lambda}(n)
	\end{equation}
	Correspondingly we have to bound two contributions: 
	\\
	
	{\it{(i) Contribution of first term in \eqref{eq:proofofpropsmoothlinearsource1}}}. This is the expression 
	\begin{align*}
	&\chi_{\xi\simeq\lambda}\left\langle\phi_{n}(R;\xi),\quad\chi_{r<\frac{t}{2}}4n\frac{\sin\left[Q+\frac{\eps}{2}\right]\sin\left[\frac{\eps}{2}\right]}{R^2}\varepsilon_{\pm,\geq \lambda}(n)\right\rangle_{L^2_{R\,dR}}\\
	& = \sum_{\mu\geq\lambda} \chi_{\xi\simeq\lambda}\left\langle\phi_{n}(R;\xi),\quad\chi_{r<\frac{t}{2}}4n\frac{\sin\left[Q+\frac{\eps}{2}\right]\sin\left[\frac{\eps}{2}\right]}{R^2}\varepsilon_{\pm,\mu}(n)\right\rangle_{L^2_{R\,dR}}
	\end{align*}
	Using Plancherel's theorem for the distorted Fourier transform, the triangle inequality as well as the point wise bound above, we infer (recalling that $\lambda(t)\cdot t\sim \tau$)
	\begin{align*}
	&\left\|\chi_{\xi\simeq\lambda}\left\langle\phi_{n}(R;\xi),\,\chi_{r<\frac{t}{2}}4n\frac{\sin\left[Q+\frac{\eps}{2}\right]\sin\left[\frac{\eps}{2}\right]}{R^2}\varepsilon_{\pm,\geq \lambda}(n)\right\rangle_{L^2_{R\,dR}}\right\|_{\Sh_{1}}\\
	&\lesssim \frac{(\lambda\hbar^2)^{1-\frac{\delta}{2}}\cdot\left\langle\lambda\hbar^2\right\rangle^{\delta+\frac32}}{\lambda^{\frac12}}\sum_{\mu\geq\lambda}\left\|\chi_{r<\frac{t}{2}}4n\frac{\sin\left[Q+\frac{\eps}{2}\right]\sin\left[\frac{\eps}{2}\right]}{R^2}\varepsilon_{\pm,\mu}(n)\right\|_{L^2_{R\,dR}}\\
	&\lesssim \frac{1}{\tau^2}(\lambda\hbar^2)^{\frac12-\frac{\delta}{2}}\cdot\left\langle\lambda\hbar^2\right\rangle^{\delta+\frac32}\sum_{\mu\geq\lambda}\left\|\langle R\rangle^{-1}\varepsilon_{\pm,\mu}(n)\right\|_{L^2_{R\,dR}}
	\end{align*}
	The last term can be bounded by taking advantage of Proposition~\ref{prop:singularmultiplier}: letting $x(\xi)$ denote the distorted Fourier transform of $\varepsilon_{\pm}$ as angular momentum-$n$ function, we have 
	\begin{align*}
	&\frac{1}{\tau^2}(\lambda\hbar^2)^{\frac12-\frac{\delta}{2}}\cdot\langle\lambda\hbar^2\rangle^{\delta+\frac32}\sum_{\mu\geq\lambda}\left\|\langle R\rangle^{-1}\varepsilon_{\pm,\mu}(n)\right\|_{L^2_{R\,dR}}\\
	&\lesssim \frac{1}{\tau^2}(\lambda\hbar^2)^{\frac12-\frac{\delta}{2}}\cdot\langle\lambda\hbar^2\rangle^{\delta+\frac32}\sum_{\mu\geq\lambda}|n|^{-1}\left\|\xi^{\frac12}\xb\right\|_{L^2_{d\xi}(\xi\simeq\mu)}\\
	&\lesssim \frac{1}{\tau^2}\sum_{\mu\geq\lambda}\left(\frac{\lambda}{\mu}\right)^{\frac12-\frac{\delta}{2}}\cdot \left\|\xb\right\|_{\Sh_{0}(\xi\simeq\mu)}. 
	\end{align*}
	Finally, exploiting orthogonality, we infer the bound 
	\begin{align*}
	&\left\|\left\langle\phi_{n}(R;\xi),\,\chi_{r<\frac{t}{2}}4n\frac{\sin\left[Q+\frac{\eps}{2}\right]\sin\left[\frac{\eps}{2}\right]}{R^2}\varepsilon_{\pm,\geq \lambda}(n)\right\rangle_{L^2_{R\,dR}}\right\|_{\Sh_{1}}\\
	&\lesssim \left(\sum_{\lambda}\left\|\chi_{\xi\simeq\lambda}\left\langle\phi_{n}(R;\xi),\,\chi_{r<\frac{t}{2}}4n\frac{\sin\left[Q+\frac{\eps}{2}\right]\sin\left[\frac{\eps}{2}\right]}{R^2}\varepsilon_{\pm,\geq \lambda}(n)\right\rangle_{L^2_{R\,dR}}\right\|_{\Sh_{1}}^2\right)^{\frac12}\\
	&\lesssim  \tau^{-2}\left(\sum_{\lambda}\left[\sum_{\mu\geq\lambda}\left(\frac{\lambda}{\mu}\right)^{\frac12-\frac{\delta}{2}}\cdot \left\|\xb\right\|_{\Sh_{0}(\xi\simeq\mu)}\right]^2\right)^{\frac12}\\
	&\lesssim \tau^{-2}\left(\sum_{\lambda}\sum_{\mu\geq\lambda}\left(\frac{\lambda}{\mu}\right)^{\frac12-\frac{\delta}{2}}\cdot \left\|\xb\right\|_{\Sh_{0}(\xi\simeq\mu)}^2\right)^{\frac12}\\
	&\lesssim \tau^{-2}\left\|\xb\right\|_{\Sh_{0}}, 
	\end{align*}
	which confirms the estimate of the proposition for this contribution. 
	\\
	
	{\it{(ii) Contribution of second term in \eqref{eq:proofofpropsmoothlinearsource1}}}. This case is more delicate, and we have to distinguish between different output frequency regimes. 
	\\
	
	{\it{(ii.a) $\hbar^2\lambda<1$.}} Here know that $\varepsilon_{\pm,< \lambda}(n)$ is in a low frequency regime, and we have to avoid losing factors $\hbar^{-1}$ when bounding this term. For this we take again advantage of Proposition~\ref{prop:singularmultiplier} but use a multiplier $\langle R\rangle^{-2+}$ this time, exploiting the fact that 
	\[
	\left|\frac{\sin\left[Q+\frac{\eps}{2}\right]\sin\left[\frac{\eps}{2}\right]}{R^2}\right|\lesssim \tau^{-2}\langle R\rangle^{-2}. 
	\]
	Then we obtain 
	\begin{align*}
	&\left\|\left\langle\chi_{\xi\simeq\lambda}\phi_{n}(R;\xi),\,\chi_{r<\frac{t}{2}}4n\frac{\sin\left[Q+\frac{\eps}{2}\right]\sin\left[\frac{\eps}{2}\right]}{R^2}\varepsilon_{\pm,<\lambda}(n)\right\rangle_{L^2_{R\,dR}}\right\|_{\Sh_{1}}\\
	&\lesssim \tau^{-2}(\lambda\hbar^2)^{\frac12-\frac{\delta}{2}}\cdot \left\|\langle R\rangle^{-2}\varepsilon_{\pm,<\lambda}(n)\right\|_{L^2_{R\,dR}}\\
	&\lesssim \tau^{-2}(\lambda\hbar^2)^{\frac12-\frac{\delta}{2}}\cdot\left\|\xb\right\|_{\Sh_{0}(\xi<\lambda)}.
	\end{align*}
	This can be $l^1$-summed over dyadic scales $\lambda<\hbar^{-2}$, giving the desired bound. 
	\\
	
	{\it{(ii.b) $\hbar^2\lambda\geq 1$.}} Here we perform integration by parts in order to absorb the outer weight into the expression. Importantly, observe that there will be no issues with boundary values at $R = 0$ since we use the same angular momentum for the output as well as for the factor $\varepsilon_{\pm,<\lambda}(n)$. Write 
	\begin{align*}
	&\chi_{\xi\simeq\lambda}\left\langle\phi_{n}(R;\xi),\,\chi_{r<\frac{t}{2}}4n\frac{\sin\left[Q+\frac{\eps}{2}\right]\sin\left[\frac{\eps}{2}\right]}{R^2}\varepsilon_{\pm,<\lambda}(n)\right\rangle_{L^2_{R\,dR}}\\
	& = \xi^{-3}\chi_{\xi\simeq\lambda}\left\langle\phi_{n}(R;\xi),\,\left(H_n^{\pm}\right)^3\left[\chi_{r<\frac{t}{2}}4n\frac{\sin\left[Q+\frac{\eps}{2}\right]\sin\left[\frac{\eps}{2}\right]}{R^2}\varepsilon_{\pm,<\lambda}(n)\right]\right\rangle_{L^2_{R\,dR}}\\
	& = \sum_{i+j+k = 6}C_{i,j,k}\xi^{-3}\chi_{\xi\simeq\lambda}\left\langle\phi_{n}(R;\xi),\,\partial_R^i\left[\chi_{r<\frac{t}{2}}4n\frac{\sin\left[Q+\frac{\eps}{2}\right]\sin\left[\frac{\eps}{2}\right]}{R^2}\right]\left(\frac{n}{R}\right)^j\partial_R^k\varepsilon_{\pm,<\lambda}(n)\right\rangle_{L^2_{R\,dR}}\\
	\end{align*}
	To bound the expression on the right in the inner product, we use that 
	\begin{align*}
	\langle R\rangle^{2-\delta}\left|\left(R\partial_R\right)^i\left[\chi_{r<\frac{t}{2}}4n\frac{\sin[Q+\frac{\eps}{2}]\sin[\frac{\eps}{2}]}{R^2}\right]\right|\lesssim_{i,\delta} \frac{n}{\tau^2}
	\end{align*}
	for any $i\geq 0$, and so we have the bound 
	\begin{align*}
	&\left\|\partial_R^i\left[\chi_{r<\frac{t}{2}}4n\frac{\sin\left[Q+\frac{\eps}{2}\right]\sin\left[\frac{\eps}{2}\right]}{R^2}\right]\left(\frac{n}{R}\right)^j\partial_R^k\varepsilon_{\pm,<\lambda}(n)\right\|_{L^2_{R\,dR}}\\&
	\lesssim \tau^{-2}n\left\|R^{-i}\left\langle R\right\rangle^{-2+\delta}\left(\frac{n}{R}\right)^j\partial_R^k\varepsilon_{\pm,<\lambda}(n)\right\|_{L^2_{R\,dR}}
	\end{align*}
	It follows that 
	\begin{align*}
	&\left\|\chi_{\xi\simeq\lambda}\left\langle\phi_{n}(R;\xi),\,\chi_{r<\frac{t}{2}}4n\frac{\sin\left[Q+\frac{\eps}{2}
		\right]\sin\left[\frac{\eps}{2}\right]}{R^2}\varepsilon_{\pm,<\lambda}(n)\right\rangle_{L^2_{R\,dR}}\right\|_{\Sh_{1}}\\
	&\lesssim \tau^{-2}\frac{(\lambda\hbar^2)^{\frac12-\frac{\delta}{2}}\langle\lambda\hbar^2\rangle^{\delta+\frac32}}{\lambda^3}\cdot \sum_{j+k\leq 6}\left\|R^{-i}\langle R\rangle^{-2+\delta}\left(\frac{n}{R}\right)^j\partial_R^k\varepsilon_{\pm,<\lambda}(n)\right\|_{L^2_{R\,dR}}\\
	&\lesssim  \tau^{-2}\frac{(\lambda\hbar^2)^{\frac12-\frac{\delta}{2}}\langle\lambda\hbar^2\rangle^{\delta+\frac32}}{\lambda^3}\cdot \sum_{j+k\leq 6}\left\|\xi^{3}x(\xi)\right\|_{L^2_{d\xi}(\xi<\lambda)}\\
	&\lesssim  \tau^{-2}\sum_{\mu<\lambda}\left(\frac{\mu}{\lambda}\right)^{\frac12-\frac{\delta}{2}}\cdot \left\|\xb(\xi)\right\|_{\Sh_{0}(\xi\simeq\mu)}.
	\end{align*}
	The desired estimate follows from this by square summing over dyadic scales $\lambda\geq \hbar^{-2}$ and applying the Cauchy-Schwarz inequality and orthogonality in the usual manner. 
\end{proof}
The preceding proposition reveals that in order to control the right hand side in \eqref{eq:RegularFinestructure2}, it suffices to deal with the more difficult first term $\lambda^{-2}N_{\pm}(n)$, which contains all the nonlinear interactions. We commence by estimating the second and third term on the right hand side in \eqref{eq:RegularFinestructure3}. Recall that $U = Q+\epsilon$. 
\begin{proposition}\label{prop:bilinwithUregular1} For $n_1, n_2$ arbitrary integers, let $\phi$ be an angular momentum $n_1$ function and $\psi$ and angular momentum $n_2$ function in the sense of Prop.~\ref{prop:bilin2}. Further, 
	let $n_3, |n_3|\geq 2$, be an integer satisfying either  (i) $n_1\simeq n_3$ and $|n_2|\lesssim |n_1|$, or (ii) $|n_3|\ll |n_1|$ and $n_1\simeq -n_2$, or (iii) $|n_3|\gg|n_1|$ and $n_3\simeq n_2$. Then if $\min\{|n_1|, |n_2|\}\gg1$, 
	the functions
	\[
	F_1: = \chi_{R\ll\tau}U_R\cdot\phi\cdot\psi_{R},\quad F_2 = \chi_{R\ll\tau}\left(\partial_{\tau}+\frac{\lambda_{\tau}}{\lambda}R\partial_R\right)U\cdot\phi\cdot\left(\partial_{\tau}+\frac{\lambda_{\tau}}{\lambda}R\partial_R\right)\psi,\quad F_3 = \chi_{R\ll\tau}\frac{2\sin U}{R^2}(n_1-n_2)\phi\cdot\overline{\psi}
	\]
	satisfy the bound 
	\begin{align*}
	\left\|F_l\right\|_{\tilde{S}_1^{(n_3)}}\lesssim \langle\min\{|n_1|, |n_2|\}\rangle^6\left\|\phi\right\|_{\tilde{S}_0^{(n_1)}}\cdot\left\|\psi\right\|_{\tilde{S}_0^{(n_2)}\cap\tilde{S}_1^{(n_2)}}.
	\end{align*}
	If $|n_j|\lesssim 1$ for all $j$, then the $F_l$ admit third order Taylor expansions around $R = 0$ of the form $P_3^{(l)} = \sum_{j=0}^3\gamma_j^{(l)}R^j$, $l = 1,2,3$, with 
	\[
	\sum_j\left|\gamma_j^{(l)}\right|\lesssim \left\|\phi\right\|_{\tilde{S}_0^{(n_1)}}\cdot\left\|\psi\right\|_{\tilde{S}_0^{(n_2)}\cap\tilde{S}_1^{(n_2)}},
	\]
	and such that we have the bounds 
	\begin{align*}
	\left\|F_l - \chi_{R\lesssim 1} \sum_{j=0}^3\gamma_j^{(l)}R^j\right\|_{\tilde{S}_1^{(n_3)}}\lesssim \left\|\phi\right\|_{\tilde{S}_0^{(n_1)}}\cdot\left\|\psi\right\|_{\tilde{S}_0^{(n_2)}\cap\tilde{S}_1^{(n_2)}}. 
	\end{align*}
\end{proposition}
\begin{remark}\label{rem:singularityatzero1} The reason for the form of $F_3$ is the term 
	\[
	\frac{2\sin U}{R^2}\sum_{j=1}^2\varphi_j\varphi_{j,\theta} = \sum_{n,m}\frac{\sin U}{R^2}\left(\varepsilon_+(n)\varepsilon_-(m)e^{i(n+m)\theta}\right)_{\theta},
	\]
	which occurs up to a small perturbative factor in \eqref{eq:RegularFinestructure3} and where we recall the relation $\varepsilon_+ = \overline{\varepsilon_-}$. Expanding out $\varepsilon_{\pm}$ in angular Fourier modes, we observe that there is a delicate cancellation in this term at $R = 0$, since the product $e^{i(n+m)\theta)}\varepsilon_+(n)\varepsilon_-(m)$ does not vanish at $R = 0$ precisely in the case $n = -m = 1$, in which case the $\theta$-derivative vanishes. 
\end{remark}
\begin{proof}
We give details for the expressions $F_{1,3}$, the remaining $F_2$ being handled similarly. 
\\

{\it{(1): Bounding $F_1$.}} We distinguish between different situations involving the angular momentum and regular frequencies. 
\\

{\it{(1.1): $\max\{|n_j|\}\gg 1$. Trivial Taylor polynomial for output.}} Localizing the output frequency $\xi\simeq\lambda$ for dyadic $\lambda$, we need to estimate 
\begin{align*}
\left\|\chi_{\xi\simeq\lambda}\left\langle \phi_{n_{3}}(R;\xi),\,\chi_{R\ll\tau}U_R\cdot\phi\cdot\psi_{R}\right\rangle_{L^2_{R\,dR}}\right\|_{\Shth_1},
\end{align*}
where the output angular momentum $n_3$ is constrained by conditions (i) - (iii) in the statement of the proposition. Then decompose the expression as follows:
\begin{equation}\label{eq:propbilinwithUregular11}\begin{split}
\chi_{\xi\simeq\lambda}\left\langle \phi_{n_{3}}(R;\xi),\,\chi_{R\lesssim\tau}U_R\phi\psi_{R}\right\rangle_{L^2_{R\,dR}} &= \chi_{\xi\simeq\lambda}\left\langle \phi_{n_{3}}(R;\xi),\,\chi_{R\ll\tau}U_R\phi\psi_{R,\geq \lambda}\right\rangle_{L^2_{R\,dR}}\\
& + \chi_{\xi\simeq\lambda}\left\langle \phi_{n_{3}}(R;\xi),\,\chi_{R\ll\tau}U_R\phi_{\geq\lambda}\psi_{R,<\lambda}\right\rangle_{L^2_{R\,dR}}\\
& + \chi_{\xi\simeq\lambda}\left\langle \phi_{n_{3}}(R;\xi),\,\chi_{R\ll\tau}U_R\phi_{<\lambda}\psi_{R,<\lambda}\right\rangle_{L^2_{R\,dR}}
\end{split}\end{equation}
We use the simple bound $\left|U_R\right|\lesssim \langle R\rangle^{-2}$. Then the first term on the right is bounded by 
\begin{align*}
&\left\|\chi_{\xi\simeq\lambda}\left\langle \phi_{n_{3}}(R;\xi),\,\chi_{R\ll\tau}U_R\phi\psi_{R,\geq \lambda}\right\rangle_{L^2_{R\,dR}}\right\|_{\Shth_1}\\
&\lesssim \hbar_3(\lambda\hbar_3^2)^{\frac12-\frac{\delta}{2}}\cdot\langle\lambda\hbar_3^2\rangle^{\delta+\frac32}\cdot \left\|\chi_{R\ll\tau}U_R\phi\psi_{R,\geq \lambda}\right\|_{L^2_{R\,dR}}. 
\end{align*}
Then in situations (iii) we have $|n_3|\gtrsim |n_2|\geq |n_1|$, and taking advantage of Lemma~\ref{lem:derLinfty} (the last estimate) and the estimate on $U_{R}$, 
we infer in this case the bound 
\begin{align*}
&\hbar_3(\lambda\hbar_3^2)^{\frac12-\frac{\delta}{2}}\cdot\langle\lambda\hbar_3^2\rangle^{\delta+\frac32}\cdot \left\|\chi_{R\ll\tau}U_R\phi\psi_{R,\geq \lambda}\right\|_{L^2_{R\,dR}}\\&
\lesssim \left\|U_R\phi\right\|_{L^\infty_{R\,dR}}\cdot \hbar_3(\lambda\hbar_3^2)^{\frac12-\frac{\delta}{2}}\cdot\langle\lambda\hbar_3^2\rangle^{\delta+\frac32}\left\|\psi_{R,\geq \lambda}\right\|_{L^2_{R\,dR}}\\
&\lesssim |n_1|^{\frac12+\delta}\left\|\phi_1\right\|_{\tilde{S}_0^{(n_1)}}\cdot \sum_{\mu\geq\lambda}\left(\frac{\lambda}{\mu}\right)^{\frac12-\frac{\delta}{2}}\cdot\left\|\xb_2\right\|_{\Sht_0(\xi\simeq\mu)},
\end{align*}
where $\xb_2$ stands for the Fourier transform of $\psi$ interpreted as angular momentum $n_2$-function, while in case (i), exploiting the bound 
\[
\hbar_3\left\|U_R\phi\right\|_{L^\infty_{R\,dR}}\lesssim \left\|\xb_1\right\|_{\Sho_{0}}, 
\]
we get 
\begin{align*}
&\hbar_3(\lambda\hbar_3^2)^{\frac12-\frac{\delta}{2}}\cdot\langle\lambda\hbar_3^2\rangle^{\delta+\frac32}\cdot \left\|\chi_{R\ll\tau}U_R\phi\psi_{R,\geq \lambda}\right\|_{L^2_{R\,dR}}\\&
\lesssim |n_2|\left\|\phi_1\right\|_{\tilde{S}_0^{(n_1)}}\cdot \sum_{\mu\geq\lambda}\left(\frac{\lambda}{\mu}\right)^{\frac12-\frac{\delta}{2}}\cdot\left\|\xb_2\right\|_{\Sht_0(\xi\simeq\mu)}.
\end{align*}
In the remaining case (ii), we have $|n_3|\ll |n_1|\simeq |n_2|$, and so we can bound 
\begin{align*}
&\hbar_3(\lambda\hbar_3^2)^{\frac12-\frac{\delta}{2}}\cdot\langle\lambda\hbar_3^2\rangle^{\delta+\frac32}\cdot \left\|\chi_{R\ll\tau}U_R\phi\psi_{R,\geq \lambda}\right\|_{L^2_{R\,dR}}\\&
\lesssim \left\|U_R\phi\right\|_{L^\infty_{R\,dR}}\cdot \hbar_3(\lambda\hbar_3^2)^{\frac12-\frac{\delta}{2}}\cdot\langle\lambda\hbar_3^2\rangle^{\delta+\frac32}\left\|\psi_{R,\geq \lambda}\right\|_{L^2_{R\,dR}}\\
&\lesssim |n_1|^{\frac12+\delta}\left\|\phi_1\right\|_{\tilde{S}_0^{(n_1)}}\cdot |n_2|^{5+\delta}\sum_{\mu\geq\lambda}\left(\frac{\lambda}{\mu}\right)^{\frac12-\frac{\delta}{2}}\cdot\left\|\xb_2\right\|_{\Sht_{0}(\xi\simeq\mu)}\\
&\lesssim |n_1|^{6}\left\|\phi_1\right\|_{\tilde{S}_0^{(n_1)}}\cdot\sum_{\mu\geq\lambda}\left(\frac{\lambda}{\mu}\right)^{\frac12-\frac{\delta}{2}}\cdot\left\|\xb_2\right\|_{\Sht_{0}(\xi\simeq\mu)}. 
\end{align*}
The desired bound follows for this case by square-summing over dyadic $\lambda$ and exploiting Cauchy-Schwarz and orthogonality. 
\\
In order to bound the second term in \eqref{eq:propbilinwithUregular11}, use that according to Prop.~\ref{prop:singularmultiplier} we have 
\[
\left\|\langle R\rangle^{-1}\phi_{\geq\lambda}\right\|_{L^2_{R\,dR}}\lesssim |n_1|^{-1}\left\|\xi^{\frac12}\xb_1\right\|_{L^2_{d\xi}},\,
\]
and so in situations (i) we have the bound
\begin{align*}
\hbar_3(\lambda\hbar_3^2)^{\frac12-\frac{\delta}{2}}\cdot\langle\lambda\hbar_3^2\rangle^{\delta+\frac32}\cdot \left\|U_R\phi_{\geq\lambda}\right\|_{L^2_{R\,dR}}\lesssim \hbar_3\sum_{\mu\geq\lambda}\left(\frac{\lambda}{\mu}\right)^{\frac12-\frac{\delta}{2}}\cdot \left\|\xb_1\right\|_{\Sho_{0}(\xi\simeq\mu)},
\end{align*}
and so using $\hbar_3\left\|\psi_{R,<\lambda}\right\|_{L^\infty_{R\,dR}}\lesssim \left|n_2\right|^{\frac12+\delta}\left\|\xb_2\right\|_{\Sht_{0}}$ according to Lemma~\ref{lem:derLinfty}, we get the bound 
\begin{align*}
&\left\|\chi_{\xi\simeq\lambda}\left\langle \phi_{n_{3}}(R;\xi),\,\chi_{R\ll\tau}U_R\phi_{\geq\lambda}\psi_{R,<\lambda}\right\rangle_{L^2_{R\,dR}}\right\|_{\Shth_1}\\
&\lesssim \left|n_2\right|^{\frac12+\delta}\cdot \left\|\xb_2\right\|_{\Sht_{0}}\cdot \sum_{\mu\geq\lambda}\left(\frac{\lambda}{\mu}\right)^{\frac12-\frac{\delta}{2}}\cdot \left\|\xb_1\right\|_{\Sho_{0}(\xi\simeq\mu)}.
\end{align*}
In situation (iii), we use (Lemma~\ref{lem:derLinfty})
\begin{align*}
\hbar_3^{2-\delta}\left\|\psi_{R,<\lambda}\right\|_{L^\infty_{R\,dR}}\lesssim \left\|\xb_2\right\|_{\Sht_{0}}
\end{align*}
as well as (Prop.~\ref{prop:singularmultiplier})
\begin{align*}
\lambda^{\frac12-\frac{\delta}{2}}\cdot\langle\lambda\hbar_3^2\rangle^{\delta+\frac32}\cdot \left\|U_R\phi_{\geq\lambda}\right\|_{L^2_{R\,dR}}\lesssim\left|n_1\right|^{1-\delta}\cdot \sum_{\mu\geq\lambda}\left(\frac{\lambda}{\mu}\right)^{\frac12-\frac{\delta}{2}}\cdot \left\|\xb_1\right\|_{\Sho_{0}(\xi\simeq\mu)}, 
\end{align*}
in order to infer the bound 
\begin{align*}
&\left\|\chi_{\xi\simeq\lambda}\left\langle \phi_{n_{3}}(R;\xi),\,\chi_{R\ll\tau}U_R\phi_{\geq\lambda}\psi_{R,<\lambda}\right\rangle_{L^2_{R\,dR}}\right\|_{\Sht_{1}}\\
&\lesssim \left|n_1\right|^{1-\delta}\cdot\left\|\xb_2\right\|_{\Sht_{0}}\cdot \sum_{\mu\geq\lambda}\left(\frac{\lambda}{\mu}\right)^{\frac12-\frac{\delta}{2}}\cdot \left\|\xb_1\right\|_{\Sho_{0}(\xi\simeq\mu)},
\end{align*}
while case (ii) is treated as for the first term in \eqref{eq:propbilinwithUregular11}. 
The desired bound again follows for this case by square-summing over dyadic $\lambda$ and exploiting Cauchy-Schwarz and orthogonality. 
\\
The last term in \eqref{eq:propbilinwithUregular11} is treated in the customary fashion by integration by parts (the fact that we assume high vanishing at the origin in case (1.1) preventing problems with boundary terms at $R = 0$), and we omit the details. 
\\

{\it{(1.2): $\max\{|n_j|\}\lesssim 1$.}} Here we have to take advantage of subtracting off the Taylor polynomial around $R = 0$. We shall treat the exceptional case when $n_j\in\{0,\pm1\}$, $j = 1,2$. This means we assume that $\phi, \psi$ admit representations in terms of the root/resonant modes $\phi_n$ as detailed in the statement of Proposition~\ref{prop:bilin2}:
\begin{align*}
\phi = c_1\cdot\phi_{n_1}(R) + \phi_{n_1}(R)\cdot \int_0^R\left[\phi_{n_1}(s)\right]^{-1}\mathcal{D}\phi(s)\,ds,\quad \psi = c_2\cdot\phi_{n_2}(R) + \phi_{n_2}(R)\cdot \int_0^R\left[\phi_{n_2}(s)\right]^{-1}\tilde{\mathcal{D}}\phi(s)\,ds
\end{align*}

Also, observe that by assumption we set $n_3 = O(1)$. 
To begin with, the case when $\phi = c_1\cdot\phi_{n_1}(R)$, $\psi = c_2\cdot\phi_{n_2}(R)$, is easy to handle. In this case the third order Taylor polynomial $\sum_{j=0}^3\gamma_j R^j$ around $R = 0$ of $ \chi_{R\ll\tau}U_R\phi\psi_{R}$ satisfies 
\[
\sum_j\left|\gamma_j\right|\lesssim \prod_{k=1,2}\left|c_k\right|.
\]
We then need to verify the bound 
\begin{align*}
\left\| \chi_{R\ll\tau}U_R\phi\psi_{R} - \chi_{R\lesssim 1}\sum_{j=0}^3\gamma_j R^j\right\|_{\Shth_1}\lesssim \prod_{k=1,2}\left|c_k\right|, 
\end{align*}
which in light of the smoothness properties of the resonances/root modes reduces to verifying sufficient decay in $R$. In fact, we have 
\begin{align*}
&\left\|\left\langle \phi_{n_{3}}(R;\xi),\,\chi_{R\ll\tau}U_R\phi\psi_{R} - \chi_{R\lesssim 1}\sum_{j=0}^3\gamma_j R^j\right\rangle_{L^2_{R\,dR}}\right\|_{\Shth_{1}}\\&\leq  \left\|\left\langle \phi_{n_{3}}(R;\xi),\,\chi_{R\ll\tau}U_R\phi\psi_{R}- \chi_{R\lesssim 1}\sum_{j=0}^3\gamma_j R^j\right\rangle_{L^2_{R\,dR}}\right\|_{\Shth_1(\xi<\hbar_3^{-2})}\\
& +  \left\|\left\langle \phi_{n_{3}}(R;\xi),\,\chi_{R\ll\tau}U_R\phi\psi_{R} - \chi_{R\lesssim 1}\sum_{j=0}^3\gamma_j R^j\right\rangle_{L^2_{R\,dR}}\right\|_{\Shth_1(\xi\geq \hbar_3^{-2})},
\end{align*}
and the first term on the right can be bounded by neglecting the weight defining $\left\|\cdot\right\|_{\Shth_1}$ and using Plancherel's theorem for the distorted Fourier transform: 
\begin{align*}
&\left\|\left\langle \phi_{n_{3}}(R;\xi),\, \chi_{R\ll\tau}U_R\phi\psi_{R}- \chi_{R\lesssim 1}\sum_{j=0}^3\gamma_j R^j \right\rangle\right\|_{\Shth_{1}(\xi<\hbar_3^{-2})}\\
&\lesssim  \prod_{k=1,2}\left|c_k\right|\cdot \left\| \chi_{R\ll\tau}U_R\phi_{n_1}\phi_{n_2, R}- \chi_{R\lesssim 1}\sum_{j=0}^3\gamma_j R^j\right\|_{L^2_{R\,dR}}\lesssim  \prod_{k=1,2}\left|c_k\right|.
\end{align*}
For the second, large frequency contribution above, we use integration by parts: the fourth order vanishing of the term on the right in the inner product implies that 
\begin{align*}
&\left\langle \phi_{n_{3}}(R;\xi),\,\chi_{R\ll\tau}U_R\phi\psi_{R} - \chi_{R\lesssim 1}\sum_{j=0}^3\gamma_j R^j\right\rangle_{L^2_{R\,dR}}\\& = \xi^{-2}\left\langle \phi_{n_{3}}(R;\xi),\,\left(H_{n_3}^{\pm}\right)^{2}\left[\chi_{R\ll\tau}U_R\phi\psi_{R} - \chi_{R\lesssim 1}\sum_{j=0}^3\gamma_j R^j\right]\right\rangle_{L^2_{R\,dR}}
\end{align*}
Then we have the simple bound 
\begin{align*}
\left|\left(H_{n_3}^{\pm}\right)^{2}\left[\chi_{R\ll\tau}U_R\phi\psi_{R} - \chi_{R\lesssim 1}\sum_{j=0}^3\gamma_j R^j\right]\right|\lesssim \frac{1}{\langle R\rangle^8}\cdot \prod_{k=1,2}\left|c_k\right|, 
\end{align*}
and using the asymptotics of $ \phi_{n_{3}}(R;\xi)$ both in the oscillatory and non-oscillatory regime as well as simple integration by parts, we infer the bound 
\begin{align*}
\left|\xi^{-2}\left\langle \phi_{n_{3}}(R;\xi),\,\left(H_{n_3}^{\pm}\right)^{2}\left[\chi_{R\ll\tau}U_R\phi\psi_{R} - \chi_{R\lesssim 1}\sum_{j=0}^3\gamma_j R^j\right]\right\rangle_{L^2_{R\,dR}}\right|\lesssim \xi^{-3}\langle\log\xi\rangle\cdot n_3^4\prod_{k=1,2}\left|c_k\right|. 
\end{align*}
We conclude that 
\begin{align*}
&\left\|\left\langle \phi_{n_{3}}(R;\xi),\,\chi_{R\ll\tau}U_R\phi\psi_{R} - \chi_{R\lesssim 1}\sum_{j=0}^3\gamma_j R^j\right\rangle_{L^2_{R\,dR}}\right\|_{\Shth_{1}(\xi\geq \hbar_3^{-2})}\\&
\lesssim  \prod_{k=1,2}\left|c_k\right|\cdot \left\|\hbar_3(\xi\hbar_3^2)^{\frac12-\frac{\delta}{2}}\langle\xi\hbar_3^2\rangle^{\delta + \frac32}\cdot  \xi^{-3}\langle\log\xi\rangle\right\|_{L^2_{d\xi}(\xi>\hbar_3^{-2})}\\
&\lesssim \prod_{k=1,2}\left|c_k\right|. 
\end{align*}
Next we consider the mixed case, where, say, 
\[
\phi = \phi_{n_1}(R)\cdot \int_0^R[\phi_{n_1}(s)]^{-1}\mathcal{D}\phi(s)\,ds,
\]
while $\psi = c_2\phi_{n_2}(R)$ is a multiple of the root/resonant mode. By Proposition \ref{prop:bilin2} (See \eqref{eq:Taylorbound1}), the third order Taylor polynomial $\sum_{j=0}^3\gamma_jR^j$ of $\chi_{R\ll\tau}U_R\phi\psi_{R}$ around $R = 0$ satisfies the bound (here $n\in \{0, \pm 1\}$, and $\xb_1$ refers to the distorted Fourier coefficient at the level of $\calD\phi$)
\[
\sum_j\left|\gamma_j\right|\lesssim \left|c_2\right|\cdot \left\|\xb_1\right\|_{\Sho_{0}},
\]
The low frequency regime of the expression is again straightforward to estimate, since 
\begin{align*}
\left\|\left\langle \phi_{n_{3}}(R;\xi),\,\chi_{R\ll\tau}U_R\phi\psi_{R} - \chi_{R\lesssim1}\sum_{j=0}^3\gamma_jR^j\right\rangle_{L^2_{R\,dR}}\right\|_{\Shth_{1}(\xi<\hbar_3^{-2})}
\lesssim \left\|\xb_1\right\|_{\Sho_{0}}\cdot \left|c_2\right|, 
\end{align*}
where we have used the fact that 
\begin{align*}
\left\|U_R\phi\right\|_{L^2_{R\,dR}}\lesssim \left\|\langle R\rangle^{-(1-\delta)}\phi\right\|_{L^2_{R\,dR}}\lesssim  \left\|\xb_1\right\|_{\Sho_{0}}, 
\end{align*}
recalling the proof of Prop.~\ref{prop:exceptionalnderivative1}, as well as the bound $\left\|\chi_{R\lesssim1}\sum_{j=0}^3\gamma_jR^j\right\|_{L^2_{R\,dR}}\lesssim \left\|\xb_1\right\|_{\Sho_0}\cdot \left|c_2\right|$. For the large frequency contribution, the interplay of the output frequency and the frequency of $\phi$ needs to be analyzed. Restricting the output frequency to dyadic scale $\lambda\geq\hbar_3^{-2}$, we have 
\begin{equation}\label{eq:propbilinwithUregular12}\begin{split}
&\chi_{\xi\simeq\lambda}\left\langle \phi_{n_{3}}(R;\xi),\,\chi_{R\ll\tau}U_R\phi\psi_{R} - \chi_{R\lesssim1}\sum_{j=0}^3\gamma_jR^j\right\rangle_{L^2_{R\,dR}}\\
& = \chi_{\xi\simeq\lambda}\left\langle \phi_{n_{3}}(R;\xi),\,\chi_{R\ll\tau}U_R\phi_{\geq\lambda}\psi_{R} \right\rangle_{L^2_{R\,dR}}\\
& + \chi_{\xi\simeq\lambda}\left\langle \phi_{n_{3}}(R;\xi),\,\chi_{R\ll\tau}U_R\phi_{<\lambda}\psi_{R} - \chi_{R\lesssim1}\sum_{j=0}^3\gamma_jR^j\right\rangle_{L^2_{R\,dR}}
\end{split}\end{equation}
To bound the first term on the right, we integrate by parts twice, using the fact that our choice for $\phi$ means that $\phi_{\geq\lambda}$ vanishes to order at least two at the origin. This gives 
\begin{align*}
&\chi_{\xi\simeq\lambda}\left\langle \phi_{n_{3}}(R;\xi),\,\chi_{R\ll\tau}U_R\phi_{\geq\lambda}\psi_{R} \right\rangle_{L^2_{R\,dR}}\\
&=\lambda^{-1} \chi_{\xi\simeq\lambda}\left\langle \phi_{n_{3}}(R;\xi),\,H_{n_3}^{\pm}\left[\chi_{R\ll\tau}U_R\phi_{\geq\lambda}\psi_{R}\right]\right\rangle_{L^2_{R\,dR}}\\
& = \lambda^{-1} \chi_{\xi\simeq\lambda}\left\langle \phi_{n_{3}}(R;\xi),\,H_{n_3}^{\pm}\left[\chi_{R\ll\tau}U_R\psi_{R}\right]\phi_{\geq\lambda}\right\rangle_{L^2_{R\,dR}}\\
& + \sum_{\substack{i+j =2\\j\geq 1}}C_{i,j}\lambda^{-1} \chi_{\xi\simeq\lambda}\left\langle \phi_{n_{3}}(R;\xi),\,\partial_R^i\left(\chi_{R\ll\tau}U_R\psi_{R}\right)\left(\partial_R^j\phi_n(R)\right)\left(\int_0^R[\phi_n(s)]^{-1}\mathcal{D}\phi_{\geq\lambda}(s)\,ds\right)\right\rangle_{L^2_{R\,dR}}\\
& +  \lambda^{-1} \chi_{\xi\simeq\lambda}\left\langle \phi_{n_{3}}(R;\xi),\,\partial_R\left[\chi_{R\ll\tau}U_R\psi_{R}\mathcal{D}\phi_{\geq\lambda}(R)\right]\right\rangle_{L^2_{R\,dR}}\\
& =: \sum_{l=1}^3 A_l. 
\end{align*}

To bound $A_1$, use a simple variation on Prop.~\ref{prop:exceptionalnderivative1} which is 
\[
\left\|\langle R\rangle^{-2}\phi_{\geq\lambda}\right\|_{L^\infty_{R\,dR}}\lesssim \sum_{\mu\geq\lambda}\left\|\xi\left[\tilde{\rho}_{n_1}(\xi)\right]^{\frac12}\xb_1(\xi)\right\|_{L^2_{d\xi}(\xi\simeq\mu)},
\]
and so we infer the bound (recall that the output frequency is restricted to $\xi\geq \hbar_3^{-2}$)
\begin{align*}
\left\|A_1\right\|_{\Shth_{1}}&\lesssim \lambda^{1+\frac{\delta}{2}}\cdot \left\|\langle R\rangle^2H_{n_3}^{\pm}\left[\chi_{R\ll\tau}U_R\psi_{R}\right]\right\|_{L^2_{R\,dR}}\cdot \left\|\langle R\rangle^{-2}\phi_{\geq\lambda}\right\|_{L^\infty_{R\,dR}}\\
&\lesssim \left|c_2\right|\sum_{\mu\geq\lambda}\left(\frac{\lambda}{\mu}\right)^{1+\frac{\delta}{2}}\cdot \left\|\xb_1\right\|_{S_0^{(n_1)}(\xi\simeq\mu)}
\end{align*}
The term $A_2$ leads to a similar bound, we omit the details. As for the last term $A_3$, expand it as follows
\begin{align*}
A_3 &=  \lambda^{-1} \chi_{\xi\simeq\lambda}\left\langle \phi_{n_{3}}(R;\xi),\,\partial_R\left[\chi_{R\ll\tau}U_R\psi_{R}\right]\mathcal{D}\phi_{\geq\lambda}(R)\right\rangle_{L^2_{R\,dR}}\\
&+  \lambda^{-1} \chi_{\xi\simeq\lambda}\left\langle \phi_{n_{3}}(R;\xi),\,\chi_{R\ll\tau}U_R\psi_{R}\partial_R\left[\mathcal{D}\phi_{\geq\lambda}(R)\right]\right\rangle_{L^2_{R\,dR}}\\
\end{align*}
Then use Prop.~\ref{prop:derivative1} to conclude that 
\begin{align*}
\left\|\mathcal{D}\phi_{\geq\lambda}(R)\right\|_{L^2_{R\,dR}}\lesssim \left\|\xb_1\right\|_{L^2_{d\xi}(\xi\geq\lambda)},\quad\left\|\partial_R\left[\mathcal{D}\phi_{\geq\lambda}(R)\right]\right\|_{L^2_{R\,dR}}\lesssim \left\|\xi^{\frac12}\xb_1\right\|_{L^2_{d\xi}(\xi\geq\lambda)}, 
\end{align*}
and so we can crudely bound 
\begin{align*}
\left\|A_3\right\|_{\Shth_1}&\lesssim \left\|\langle\partial_R\rangle\left[\chi_{R\ll\tau}U_R\psi_{R}\right]\right\|_{L^\infty_{R\,dR}}\cdot\lambda^{1+\frac{\delta}{2}}\cdot\left\|\xi^{\frac12}\xb_1\right\|_{L^2_{d\xi}(\xi\geq\lambda)}\\
&\lesssim \left|c_2\right|\cdot \sum_{\mu\geq\lambda}\left(\frac{\lambda}{\mu}\right)^{1+\frac{\delta}{2}}\left\|\xb_1\right\|_{\Sho_{0}(\xi\simeq\mu)} 
\end{align*}
Combining the preceding bounds for $A_j,\,j=1,2,3$, and square-summing over $\lambda$ as well as exploiting Cauchy-Schwarz and orthogonality gives the desired bound for the contribution of the first term on the right in \eqref{eq:propbilinwithUregular12}.
\\

Consider then the second term in \eqref{eq:propbilinwithUregular12}, where we have to take advantage of a partial cancellation between the principal term and the truncated Taylor polynomial. To begin with, observe that in the case $n = -1$, the third order Taylor polynomial at $R = 0$ is trivial, since 
\[
\left|\phi_{-1}(R)\cdot\int_0^R\left[\phi_{-1}(s)\right]^{-1}\mathcal{D}_{-1}\phi(s)\,ds\right|\lesssim R^{4+\delta}\left\|\xb_1\right\|_{S_0^{(n_1)}}
\]
Thus in this case a direct integration by parts argument works: 
\begin{align*}
&\chi_{\xi\simeq\lambda}\left\langle \phi_{n_{3}}(R;\xi),\,\chi_{R\ll\tau}U_R\phi_{<\lambda}\psi_{R} - \chi_{R\lesssim1}\sum_{j=0}^3\gamma_jR^j\right\rangle_{L^2_{R\,dR}}\\
& =  \chi_{\xi\simeq\lambda}\left\langle \phi_{n_3}(R;\xi),\,\chi_{R\ll\tau}U_R\phi_{<\lambda}\psi_{R} \right\rangle_{L^2_{R\,dR}}\\
& = \lambda^{-2}\chi_{\xi\simeq\lambda}\left\langle \phi_{n_{3}}(R;\xi),\,\left(H_{n_3}^{\pm}\right)^2\left[\chi_{R\ll\tau}U_R\phi_{<\lambda}\psi_{R}\right]\right\rangle_{L^2_{R\,dR}}
\end{align*}
Expand
\begin{align*}
&\left(H_{n_3}^{\pm}\right)^2\left[\chi_{R\ll\tau}U_R\phi_{<\lambda}\psi_{R}\right]\\& = \sum_{i+j = 4}C_{i,j}\left(\frac{n_3}{R}\right)^i\partial_R^j\left(\chi_{R\ll\tau}U_R\psi_R\phi_{-1}(R)\right)\cdot\int_0^R\left[\phi_{-1}(s)\right]^{-1}\mathcal{D}_{-1}\phi_{<\lambda}(s)\,ds\\
& + \sum_{\substack{i+j+k = 4\\k\geq 1}}C_{i,j,k}\left(\frac{n_3}{R}\right)^i\partial_R^j\left(\chi_{R\ll\tau}U_R\psi_R\phi_{-1}(R)\right)\cdot\partial_R^{k-1}\left(\left[\phi_{-1}(R)\right]^{-1}\mathcal{D}_{-1}\phi_{<\lambda}(R)\right),
\end{align*}
and we have the crude bound (compare Prop.~\ref{prop:exceptionalnderivative1} and its proof)
\begin{align*}
&\left|\langle R\rangle^{-1}\int_0^R\left[\phi_{-1}(s)\right]^{-1}\mathcal{D}_{-1}\phi_{<\lambda}(s)\,ds\right|\lesssim \left\|\xb_1\right\|_{S_0^{(-1)}(\xi<\lambda)},
\end{align*}
while for the second term we have the more sophisticated bound (under the constraint $i+j+k = 4$)
\begin{align*}
&\left|\left(\frac{n_3}{R}\right)^i\partial_R^j\left(\chi_{R\ll\tau}U_R\psi_{R}\phi_{-1}(R)\right)\partial_R^{k-1}\left(\left[\phi_{-1}(R)\right]^{-1}\mathcal{D}_{-1}\phi_{<\lambda}(R)\right)\right|\lesssim |c_{2}|\cdot\left\|\xb_1\right\|_{S_0^{(-1)}(\xi<\lambda)}.
\end{align*}
Back to the inner product $\langle\cdot\rangle_{L^2_{R\,dR}}$, we split this via a smooth cutoff into the regions $R\lambda^{\frac12}<1,\,R\lambda^{\frac12}\geq 1$, and perform further integration by parts as needed in the latter region (of course without generating boundary terms). Note that the vanishing order of $\phi_{-1}(R)\cdot\int_{0}^{R}\left[\phi_{-1}(s)\right]^{-1}\,\calD_{-1}\phi(s)\,ds$ at $R=0$ together with the measure $R\,dR$ allows us to perform one more integration by parts. This finally leads to the bound 
\begin{align*}
&\left|\lambda^{-2}\chi_{\xi\simeq\lambda}\left\langle \phi_{n_{3}}(R;\xi),\,\left(H_{n_3}^{\pm}\right)^2\left[\chi_{R\ll\tau}U_R\phi_{<\lambda}\psi_{R}\right]\right\rangle_{L^2_{R\,dR}}\right|
\lesssim \lambda^{-3}\cdot \left\|\xb_1\right\|_{S_0^{(-1)}(\xi<\lambda)}\cdot\left|c_2\right|. 
\end{align*}
If we now apply the $\|\cdot\|_{\Shth_{1}}$-norm to this term and neglect the weight $\hbar_3\simeq 1$, we obtain 
\begin{align*}
&\left\|\lambda^{-2}\chi_{\xi\simeq\lambda}\left\langle \phi_{n_{3}}(R;\xi),\,\left(H_{n_3}^{\pm}\right)^2\left[\chi_{R\ll\tau}U_R\phi_{<\lambda}\psi_{R}\right]\right\rangle_{L^2_{R\,dR}}\right\|_{\Shth_{1}}\\
&\lesssim \lambda^{\frac52+\frac{\delta}{2}}\cdot\lambda^{-3}\cdot\left\|\xb_1\right\|_{S_0^{(-1)}(\xi<\lambda)} \cdot\left|c_2\right|
\end{align*}
This can then even be $l^1$-summed over dyadic scales $\lambda\geq\hbar_3^{-2}$ to result in the desired bound. 
\\

Let us next consider the case $n_1 = 1$, where the Taylor polynomial of third order is not necessarily trivial. Then we write 
\begin{align*}
\chi_{R\ll\tau}U_R\phi_{<\lambda}\psi_{R} - \chi_{R\lesssim1}\sum_{j=0}^3\gamma_jR^j
&= \chi_{1\lesssim R\ll\tau}U_R\phi_{<\lambda}\psi_{R}\\
&+  \chi_{R\lesssim1}U_R\psi_{R}\phi_{1}(R)\cdot\int_0^R\left[\phi_1(s)\right]^{-1}\mathcal{D}_1\phi_{<\lambda}(s)\,ds -  \chi_{R\lesssim1}\sum_{j=0}^3\gamma_jR^j,
\end{align*}
where throughout all cutoffs $ \chi_{1\lesssim R\ll\tau}$ are chosen smoothly. Then the contribution of the first term on the right, which is 
\[
\chi_{\xi\simeq\lambda}\left\langle \phi_{n_{3}}(R;\xi),\, \chi_{1\lesssim R\ll\tau}U_R\phi_{<\lambda}\psi_{R}\right\rangle_{L^2_{R\,dR}},
\]
can be handled in the customary fashion by shifting three copies of $H_{n_3}^{\pm}$ from left to right, since the term vanishes in a neighborhood of $R  =0$, and we omit the details. For the second term on the right, we note that here in fact $\sum_{j=0}^3\gamma_jR^j = \gamma_2R^2$ or $\gamma_3R^3$(depending on $n_2$),  and we have 
\begin{equation}\label{eq:propbilinwithUregular13}\begin{split}
&\chi_{R\lesssim1}U_R\psi_{R}\phi_{1}(R)\cdot\int_0^R\left[\phi_1(s)\right]^{-1}\mathcal{D}_1\phi_{<\lambda}(s)\,ds -  \chi_{R\lesssim1}\sum_{j=0}^3\gamma_jR^j\\
&= \chi_{R\lesssim1}U_R\psi_{R}\phi_{1}(R)\cdot\int_0^R\left[\phi_1(s)\right]^{-1}\widetilde{\mathcal{D}_1\phi_{<\lambda}}(s)\,ds+ E, 
\end{split}\end{equation}
and where we define the expression $\widetilde{\mathcal{D}_1\phi_{<\lambda}}(s)$ as 
\begin{align*}
\widetilde{\mathcal{D}_1\phi_{<\lambda}}(s) =\int_0^\infty \chi_{\xi<\lambda}\left[\phi_1(s,\xi) - \frac{\pi}{4}s\right]\xb_1(\xi)\tilde{\rho}_{1}(\xi)\,d\xi, 
\end{align*}
and we recall that $\frac{\pi}{4}s$ is the linear part of $\phi(s,\xi)$ near $s = 0$. Then we claim that both terms in \eqref{eq:propbilinwithUregular13} lead to good contributions. In fact, for the first term, splitting the inner product $\langle\cdot\rangle_{L^2_{R\,dR}}$ smoothly into the ranges $R\lambda^{\frac12}<1, R\lambda^{\frac12}>1$ and performing integration by parts in the second regime as needed, and using the bound 
\begin{align*}
\left|\chi_{R\lesssim1}U_R\psi_{R}\phi_{1}(R)\cdot\int_0^R\left[\phi_1(s)\right]^{-1}\widetilde{\mathcal{D}_1\phi_{<\lambda}}(s)\,ds\right|\lesssim R^4\cdot\left\|\xb_1\right\|_{\Sho_{0}}
\end{align*}
as well as simple analogues for its derivatives, we infer 
\begin{align*}
&\left|\chi_{\xi\simeq\lambda}\left\langle \phi_{n_{3}}(R;\xi),\,\chi_{R\lesssim1}U_R\psi_{R}\phi_{1}(R)\cdot\int_0^R\left[\phi_1(s)\right]^{-1}\widetilde{\mathcal{D}_1\phi_{<\lambda}}(s)\,ds\right\rangle_{L^2_{R\,dR}}\right|\\
&\lesssim \lambda^{-\frac52}\cdot |c_2|\cdot \big\|x_1\big\|_{\Sho_{0}}, 
\end{align*}
which in combination with H\"older's inequality gives 
\begin{align*}
&\left\|\chi_{\xi\simeq\lambda}\left\langle \phi_{n_{3}}(R;\xi),\,\chi_{R\lesssim1}U_R\psi_{R}\phi_{1}(R)\cdot\int_0^R\left[\phi_1(s)\right]^{-1}\widetilde{\mathcal{D}_1\phi_{<\lambda}}(s)\,ds\right\rangle_{L^2_{R\,dR}}\right\|_{\Shth_1}\\
&\lesssim \lambda^{-\frac12+\frac{\delta}{2}}\cdot |c_2|\cdot \left\|\xb_1\right\|_{\Sho_{0}}.
\end{align*}
This bound can again be $l^1$-summed over dyadic scales $\lambda\geq \hbar_3^{-2}$. It remains to deal with the term $E$ in \eqref{eq:propbilinwithUregular13}. Assuming $n_2 = \pm 1$ to be concrete(the case $n_2 = 0$ being handled analogously), this term can be written as 
\begin{align*}
E &= \chi_{R\lesssim1}U_R\psi_{R}\phi_{1}(R)\cdot\int_0^R\left[\phi_1(s)\right]^{-1}\cdot\frac{\pi}{4}s\,ds\cdot \int_0^\infty \chi_{\xi<\lambda}\xb_1(\xi)\tilde{\rho}_{1}(\xi)\,d\xi\\
& - \chi_{R\lesssim1}U_R(0)R\psi_{RR}(0)\phi_{1}(0)\cdot\int_0^R\frac{\pi}{4}s\,ds\cdot \int_0^\infty \chi_{\xi<\lambda}\xb_1(\xi)\tilde{\rho}_{1}(\xi)\,d\xi\\
& -  \chi_{R\lesssim1}U_R(0)R\psi_{RR}(0)\phi_{1}(0)\cdot\int_0^R\frac{\pi}{4}s\,ds\cdot \int_0^\infty \chi_{\xi\geq \lambda}\xb_1(\xi)\tilde{\rho}_{1}(\xi)\,d\xi
\end{align*}
The difference of the first two terms is a $C^\infty$ function on $[0,\infty)$ of size $O_{\left|c_{2}\right|\left\|\xb_1\right\|_{\Sho_{0}}}\left(R^4\right)$ near the origin, and its contribution is handled analogously to the one of the first term on the right in \eqref{eq:propbilinwithUregular13}. There remains the contribution of the last term, which is a smooth function of size 
\[
O_{\left|c_2\right|\cdot\left\|\xb_1\right\|_{\Sho_{0}}}\left(\lambda^{-1-\frac{\delta}{2}}R^3\right).
\]
Calling this last term $E_3$, we then conclude that 
\begin{align*}
\left|\chi_{\xi\simeq\lambda}\left\langle \phi_{n_{3}}(R;\xi),\,E_3\right\rangle_{L^2_{R\,dR}}\right|\lesssim |c_2|\cdot\left\|\xb_1\right\|_{\Sho_{0}}\cdot \lambda^{-3-\frac{\delta}{2}}, 
\end{align*}
from which the desired $\|\cdot\|_{\Shth_{1}}$ follows as usual via H\"older and summation over dyadic scales $\lambda\geq \hbar_3^{-2}$.\\
The case when $\phi = c_1\phi_{n_1}(R)$ but $\psi = \phi_{n_2}(R)\cdot \int_0^R\left[\phi_{n_2}(s)\right]^{-1}\tilde{\mathcal{D}}\phi(s)\,ds$ is handled analogously, as is the case when both $\phi, \psi$ are of the latter form. 
This concludes case {\it{(1.2)}} and thereby case {\it{(1)}}, i.e., bounding the term $F_1$. 
\\

{\it{(2): Bounding $F_3$.}} This is largely analogous to the bound for $F_1$, since $\left|\frac{\sin U}{R}\right|$ has the same asymptotics as $R\rightarrow 0$ or $R\rightarrow\infty$ as $\left|U_R\right|$, and $\psi_R$ gets replaced by $(n_1 - n_2)\frac{\psi}{R}$ or $(n_1-n_2)\frac{\phi}{R}$. As remarked earlier after the statement of the proposition, the precise choice of the coefficient $n_1 - n_2$ ensures that the function is never singular at $R = 0$. 
\end{proof}
\begin{corollary}\label{cor:bilinwithUregular1} Assume that $\phi = \sum_{n\in\Z}\phi^{(n)}$, where $\phi^{(n)}$ is an angular momentum $n$ function, and similarly $\psi = \sum_{m\in\Z}\psi^{(m)}$, and we assume the bound 
	\[
	\sum_n\langle n\rangle^{12}\left\|\phi^{(n)}\right\|_{\tilde{S}_0^{(n)}\cap \tilde{S}_1^{(n)}} = \Lambda_1<\infty,\quad \sum_m\langle m\rangle^{12}\left\|\psi^{(m)}\right\|_{\tilde{S}_0^{(n)}\cap \tilde{S}_1^{(n)}} = \Lambda_2<\infty.
	\]
	Then if 
	\[
	F_l = F_l(\phi, \psi) = \sum_{n,m}F_l(\phi^{(n)}, \psi^{(m)}) =\sum_k\sum_{n+m = k}F_l(\phi^{(n)}, \psi^{(m)}) =:\sum_kF_l^{(k)}
	\]
	where $F_l$ is one of the bilinear expressions in Proposition~\ref{prop:bilinwithUregular1}, with $n, m$ replacing $n_1, n_2$ if $l = 3$, then for each $k$ the function $F_l^{(k)}$ admits a third order Taylor expansion 
	\[
	\sum_{j=0}^3\gamma^{(l,k)}_jR^j
	\]
	around $0$ and we have 
	\begin{align*}
	\sum_{k\in \Z}\langle k\rangle^{12}\left[\left\|\chi_{R\lesssim\tau}F_l^{(k)} - \chi_{R\lesssim 1}\sum_{j=0}^3\gamma^{(l,k)}_jR^j\right\|_{\tilde{S}_1^{k_{|\cdot|\geq 2}}} + \sum_j\left|\gamma^{(l,k)}_j\right|\right]\lesssim \prod_{j=1,2}\Lambda_j. 
	\end{align*}
\end{corollary}
We can now control the second and third term on the right in \eqref{eq:RegularFinestructure3}. For this we have to keep in mind that the angular decompositions and estimates will be at the level of the functions $\varepsilon_{\pm}$, which in turn determine the functions $\varphi_{1,2}$ by means of  
\begin{align*}
\varphi_1 = \frac12\left[\varepsilon_{+} + \varepsilon_{-}\right],\quad\varphi_2 = \frac{1}{2i}\left[\varepsilon_{-} -\varepsilon_{+}\right].
\end{align*}
Then we decompose 
\begin{equation}\label{eq:varepsilonpmangulardecomp}
\varepsilon_{+} = \sum_{n\in\Z}\varepsilon_{+}(n)e^{in\theta},\quad\varepsilon_{-} = \sum_{n\in\Z}\varepsilon_{-}(m)e^{im\theta},
\end{equation}
where $\varepsilon_{+}(n)$ is an angular momentum $n$ function in reference to $H_n^+$, while $\varepsilon_{-}(m)$ is an angular momentum $m$ function in reference to $H_m^{-}$. We shall assume the bounds 
\begin{equation}\label{eq:varepsilonpmangulardecomp1}
\sum_{n\in\Z}\langle n\rangle^{12}\left\|\varepsilon_{+}(n)\right\|_{\tilde{S}_0^{(n)}\cap\tilde{S}_1^{(n)}} +  \sum_{m\in\Z}\langle m\rangle^{12}\left\|\varepsilon_{-}(m)\right\|_{\tilde{S}_0^{(m)}\cap\tilde{S}_1^{(m)}} =:\Lambda<\infty. 
\end{equation}
Furthermore, we shall impose the smallness condition 
\begin{equation}\label{eq:varepsilonpmangulardecomp2}
\tau^{1+\delta}\Lambda\ll 1.
\end{equation}
\begin{proposition}\label{prop:bilinwithUregular2} Denote by $F_{1,2}$ the functions 
	\begin{align*}
	&F_ 1 = \chi_{R\ll\tau}\frac{2}{\sqrt{1-\left|\Pi_{\Phi^{\perp}}\varphi\right|^2}}\left[U_R\sum_{j=1}^2\vphi_j\vphi_{j,R} - \left(\partial_{\tau}+\frac{\lambda_{\tau}}{\lambda}R\partial_R\right)U\sum_{j=1}^2\vphi_j\left(\partial_{\tau}+\frac{\lambda_{\tau}}{\lambda}R\partial_R\right)\vphi_j\right]\\
	&F_2 = \frac{2\sin U}{R^2\sqrt{1-\left|\Pi_{\Phi^{\perp}}\vphi\right|^2}}\sum_{j=1}^2\vphi_j\vphi_{j,\theta}
	\end{align*}
	Then assuming \eqref{eq:varepsilonpmangulardecomp1}, \eqref{eq:varepsilonpmangulardecomp2}, and setting 
	\[
	F_j = \sum_{n}F_j^{(n)}e^{in\theta},\quad j = 1,2. 
	\]
	we have the bound 
	\begin{align*}
	\sum_{n\in\Z}\langle n\rangle^{12}\left[\left\|F_j^{(n)} - \chi_{R\lesssim 1}\sum_{l=0}^3 \gamma_l^{(n,j)}R^l\right\|_{\tilde{S}_1^{(n_{|\cdot|\geq 2})}} + \sum_{l=0}^3\left|\gamma_l^{(n,j)}\right|\right]\lesssim \Lambda^2. 
	\end{align*}
	where $\sum_{l=0}^3 \gamma_l^{(n,j)}R^l$ is the third order Taylor development of $F_j^{(n)}$ at $R = 0$. The same bound obtains if we subtract from the $\sum_{l=0}^3\gamma_l^{(n,j)}R^l$ those terms $\gamma_l^{(n,j)}R^l$ with $l\geq |n_{|\cdot|\geq 2}-1|, l\equiv (n_{|\cdot|\geq 2}-1)\, \text{mod} \,2$, provided $|n_{|\cdot|\geq 2}-1|\leq 3$. 
\end{proposition}
\begin{proof} This is a consequence of Corollary~\ref{cor:bilinwithUregular1} and Corollary~\ref{cor:multilingen2} after expanding $\frac{2}{\sqrt{1-\left|\Pi_{\Phi^{\perp}}\vphi\right|^2}}$ in a power series. 
\end{proof}
The problem of controlling the source terms \eqref{eq:RegularFinestructure2}, \eqref{eq:RegularFinestructure3} in the regular regime away from the light cone is then reduced to bounding the term $\mathcal{P}\varepsilon_{\pm}$. 
\\

We now consider the last remaining  term, for which we have to rely on the identities \eqref{non linear 2b}-\eqref{non linear 1 2}. Here we encounter a number of terms which appear singular at the origin $R = 0$, but these are of course spurious singularities that disappear when taking the algebraic fine structure into account. To formulate the necessary estimates allowing us to deal with this term, we first have to render these cancellations explicit. We shall sort these according to their degree. 
\\

{\it{(i) Linear terms singular at $R = 0$.}} These are given by 
\[
-2\left(1+a\left(\Pi_{\Phi^{\perp}}\vphi\right)\right)\frac{\sin U \cos U}{R^2}\vphi_1 - 2\frac{\sin U}{R^2}\left(1+a\left(\Pi_{\Phi^{\perp}}\vphi\right)\right)\vphi_{2,\theta}. 
\]
Since $1-\cos U =O(R^2)$, we can replace this expression by the sum of a  simpler one and a term regular at the origin: 
\begin{align*}
&-2\left(1+a\left(\Pi_{\Phi^{\perp}}\vphi\right)\right)\frac{\sin U \cos U}{R^2}\vphi_1 - 2\frac{\sin U}{R^2}\left(1+a\left(\Pi_{\Phi^{\perp}}\vphi\right)\right)\vphi_{2,\theta}\\
& = -\frac{2\sin U}{R^2}\left(1+a\left(\Pi_{\Phi^{\perp}}\vphi\right)\right)\left[\vphi_1 + \vphi_{2,\theta}\right] + \frac{2\sin U}{R^2}(1-\cos U)\left(1+a\left(\Pi_{\Phi^{\perp}}\vphi\right)\right)\vphi_1
\end{align*}
Here the first, singular term can be written in terms of $\varepsilon_{\pm}$ as follows:
\begin{align*}
\frac{2\sin U}{R^2}\left(1+a\left(\Pi_{\Phi^{\perp}}\vphi\right)\right)\left[\vphi_1 + \vphi_{2,\theta}\right] &= \frac{\sin U}{R^2}\left(1+a\left(\Pi_{\Phi^{\perp}}\vphi\right)\right)\left[\varepsilon_+ + i\varepsilon_{+,\theta}\right]\\
& +  \frac{\sin U}{R^2}\left(1+a\left(\Pi_{\Phi^{\perp}}\vphi\right)\right)[\varepsilon_- - i\varepsilon_{-,\theta}]
\end{align*}
Observe that the expression $[\varepsilon_+ + i\varepsilon_{+,\theta}]$ vanishes at angular momentum $n = 1$, which is precisely the value for which $\varepsilon_+$ does not vanish at the origin $R = 0$. 
Similarly $[\varepsilon_- - i\varepsilon_{-,\theta}]$ vanishes at angular momentum $n = -1$, which is precisely the value at which $\varepsilon_-$ does not vanish at $R = 0$. It follows that both expressions on the right are in fact non-singular at $R = 0$. 
\\

{\it{(ii) Quadratic terms singular at $R = 0$.}} These are 
\[
- \frac{1}{R^2}\left[\vphi_{1,\theta}^2 + \vphi_{2,\theta}^2\right]-\frac{1}{R^2}\left[\vphi_2^2 + \cos^2U\vphi_1^2\right] + \frac{2\cos U}{R^2}\left(\vphi_{1,\theta}\vphi_2 - \vphi_1\vphi_{2,\theta}\right)
\]
Again we can split this into a 'principal term' and two additional terms which are clearly non-singular at $R = 0$: 
\begin{align*}
& - \frac{1}{R^2}\left[\vphi_{1,\theta}^2 + \vphi_{2,\theta}^2\right]-\frac{1}{R^2}\left[\vphi_2^2 + \vphi_1^2\right] + \frac{2}{R^2}\left(\vphi_{1,\theta}\vphi_2 - \vphi_1\vphi_{2,\theta}\right)\\
&+\frac{1}{R^2}(1- \cos^2U)\vphi_1^2 +  \frac{2(\cos U-1)}{R^2}\left(\vphi_{1,\theta}\vphi_2 - \vphi_1\vphi_{2,\theta}\right)
\end{align*}
Here the first term can be written as 
\begin{align*}
& - \frac{1}{R^2}\left[\vphi_{1,\theta}^2 + \vphi_{2,\theta}^2\right]-\frac{1}{R^2}\left[\vphi_2^2 + \vphi_1^2\right] + \frac{2}{R^2}\left(\vphi_{1,\theta}\vphi_2 - \vphi_1\vphi_{2,\theta}\right)\\
& = -\frac{1}{R^2}\left(\vphi_1 + \vphi_{2,\theta}\right)^2 - \frac{1}{R^2}\left(\vphi_2 - \vphi_{1,\theta}\right)^2 
\end{align*}
and since 
\[
\vphi_1 + \vphi_{2,\theta} = \frac{\varepsilon_+ +  i\varepsilon_{+,\theta}}{2} +  \frac{\varepsilon_- -  i\varepsilon_{-,\theta}}{2},\quad\vphi_2 - \vphi_{1,\theta} = i\frac{\varepsilon_+ +  i\varepsilon_{+,\theta}}{2} -  i\frac{\varepsilon_- -  i\varepsilon_{-,\theta}}{2},
\]
we again verify that the preceding quadratic expression is non-singular at $R = 0$. 
\\

{\it{(iii) Higher order singular terms at $R = 0$}}: 
\[
-\frac{1}{R^2}\left[\left(a\left(\Pi_{\Phi^{\perp}}\vphi\right)\right)_{\theta}\right]^2 + \frac{2\sin U}{R^2}\left(a\left(\Pi_{\Phi^{\perp}}\vphi\right)\right)_{\theta}\vphi_2
\]
Recall that $a\left(\Pi_{\Phi^{\perp}}\vphi\right) = \sqrt{1-\left|\Pi_{\Phi^{\perp}}\vphi\right|^2} - 1$, and so 
\begin{align*}
\left(a\left(\Pi_{\Phi^{\perp}}\vphi\right)\right)_{\theta} = -\frac{\left(\left|\Pi_{\Phi^{\perp}}\vphi\right|^2\right)_{\theta}}{\sqrt{1-\left|\Pi_{\Phi^{\perp}}\vphi\right|^2}}. 
\end{align*}
Recalling Remark~\ref{rem:singularityatzero1} we see that this expression vanishes at $R= 0$, and hence the above terms are all non-singular at $R =0$. 
\\
Finally, control over the term $\mathcal{P}\varepsilon_{\pm}$ in the region away from the light cone will follow from 
\begin{proposition}\label{prop:bilinwithUregular3} Assume \eqref{eq:varepsilonpmangulardecomp}, \eqref{eq:varepsilonpmangulardecomp1}, \eqref{eq:varepsilonpmangulardecomp2}. Define the functions (\eqref{non linear 2b}-\eqref{non linear 1 2})
	\begin{align*}
	F_1 =\chi_{R\ll\tau}&\Bigg[-\left(2a\left(\Pi_{\Phi^{\perp}}\vphi\right) + \left(a\left(\Pi_{\Phi^{\perp}}\vphi\right)\right)^2\right)\left(U_R^2 - \left[\left(\partial_{\tau}+\frac{\lambda_{\tau}}{\lambda}R\partial_R
	\right)U\right]^2 + \frac{\sin^2U}{R^2}\right)\\
	&+ \left[\left(\partial_{\tau}+\frac{\lambda_{\tau}}{\lambda}R\partial_R\right)\left(a\left(\Pi_{\Phi^{\perp}}\vphi\right)\right)\right]^2 -  \left[\partial_R\left(a\left(\Pi_{\Phi^{\perp}}\vphi\right)\right)\right]^2\\
	&+ \sum_{j=1,2}\left[\left(\partial_{\tau}+\frac{\lambda_{\tau}}{\lambda}R\partial_R\right)\vphi_j\right]^2 -  \sum_{j=1,2}\left[\partial_R\vphi_j\right]^2\\
	&+2\phi_1\left(\left(a\left(\Pi_{\Phi^{\perp}}\vphi\right)\right)_RU_R - \left(\partial_{\tau}+\frac{\lambda_{\tau}}{\lambda}R\partial_R\right)\left(a\left(\Pi_{\Phi^{\perp}}\vphi\right)\right)\left(\partial_{\tau}+\frac{\lambda_{\tau}}{\lambda}R\partial_R\right)U\right)\\
	& + \left(U_R^2 - \left[\left(\partial_{\tau}+\frac{\lambda_{\tau}}{\lambda}R\partial_R\right)U\right]^2\right)\vphi_1^2\\
	& + 2\left(1+a\left(\Pi_{\Phi^{\perp}}\vphi\right)\right)\left( \left(\partial_{\tau}+\frac{\lambda_{\tau}}{\lambda}R\partial_R\right)U\cdot\left(\partial_{\tau}+\frac{\lambda_{\tau}}{\lambda}R\partial_R\right)\vphi_1 - U_R\cdot\vphi_{1,R}\right)\Bigg]\cdot\varepsilon_{\pm},
	\end{align*}
	\begin{align*}
	F_2 = \chi_{R\ll\tau}&\Bigg[\frac{2\sin U}{R^2}\left(1-\cos U\right)\left(1+a\left(\Pi_{\Phi^{\perp}}\vphi\right)\right)\vphi_1\\
	&+\frac{1}{R^2}\left(1- \cos^2U\right)\vphi_1^2 +  \frac{2\left(\cos U-1\right)}{R^2}\left(\vphi_{1,\theta}\vphi_2 - \vphi_1\vphi_{2,\theta}\right)\Bigg]\cdot\varepsilon_{\pm}, 
	\end{align*}
	\begin{align*}
	F_3 = \chi_{R\ll\tau}&\Bigg[\frac{\sin U}{R^2}\left(1+a\left(\Pi_{\Phi^{\perp}}\vphi\right)\right)\left[\varepsilon_+ + i\varepsilon_{+,\theta}\right]\\
	& +  \frac{\sin U}{R^2}\left(1+a\left(\Pi_{\Phi^{\perp}}\vphi\right)\right)\left[\varepsilon_- - i\varepsilon_{-,\theta}\right]\Bigg]\cdot\varepsilon_{\pm}, 
	\end{align*}
	\begin{align*}
	F_4 = \chi_{R\ll\tau}&\left[ -\frac{1}{R^2}\left(\vphi_1 + \vphi_{2,\theta}\right)^2 - \frac{1}{R^2}\left(\vphi_2 - \vphi_{1,\theta}\right)^2\right]\cdot\varepsilon_{\pm},  
	\end{align*}
	\begin{align*}
	F_5 = \chi_{R\ll\tau}&\left[\left(a\left(\Pi_{\Phi^{\perp}}\vphi\right)\right)_{\theta}\right]\cdot\varepsilon_{\pm} = -\frac{\left(\left|\Pi_{\Phi^{\perp}}\vphi\right|^2\right)_{\theta}}{\sqrt{1-\left|\Pi_{\Phi^{\perp}}\vphi\right|^2}}\cdot\varepsilon_{\pm}. 
	\end{align*}
	Then setting 
	\[
	F_j = \sum_{n}F_j^{(n)}e^{in\theta},\quad j = 1,\ldots,5, 
	\]
	we have the bound 
	\begin{align*}
	\sum_{n\in\Z}\langle n\rangle^{12}\left[\left\|F_j^{(n)} - \chi_{R\lesssim 1}\sum_{l=0}^3 \gamma_l^{(n,j)}R^l\right\|_{\tilde{S}_1^{(n_{|\cdot|\geq 2})}} + \sum_{l=0}^3\left|\gamma_l^{(n,j)}\right|\right]\lesssim \Lambda^2. 
	\end{align*}
	where $\sum_{l=0}^3 \gamma_l^{(n,j)}R^l$ is the third order Taylor development of $F_j^{(n)}$ at $R = 0$. The same bound obtains if we subtract from the $\sum_{l=0}^3\gamma_l^{(n,j)}R^l$ those terms $\gamma_l^{(n,j)}R^l$ with $l\geq |n_{|\cdot|\geq 2}-1|, l\equiv (n_{|\cdot|\geq 2}-1)\text{mod}(2)$, provided $|n_{|\cdot|\geq 2}-1|\leq 3$. 
	\end{proposition}
\begin{proof}
		The proof is similar to that of Proposition \ref{prop:bilinwithUregular1} and we omit the details here.
\end{proof}
We observe that Proposition~\ref{prop:smoothlinearsource}, Proposition~\ref{prop:bilinwithUregular2} and Proposition~\ref{prop:bilinwithUregular3} complete control over the source terms in the non-singular region for the $|n|\geq 2$ angular momentum  modes. 
\section{Estimates for singular sources}\label{sec:multilinestimatesnearcone}
\subsection{Description of the shock}\label{subsec:introshock}
Here we briefly recall the function spaces from \cite{KST,GaoK} which in analogy to Definition~\ref{def:KST1} describe the correction $\epsilon$ used to build $U$ close to the light cone: 
\begin{defi}\label{defi:Q_n}For $i\in \mathbb{N}$, let $j(i) = i$ if $\nu$ is irrational, respectively $j(i) = 2i^2$ if $\nu$ is rational. 
	Then
	\begin{itemize}
		\item $\mathcal{Q}$ is the algebra of continuous functions $q: [0,1]\rightarrow \R$ with the following properties: 
		\\
		(i) $q$ is analytic in $[0,1)$ with even expansion around $a = 0$.\\
		(ii) near $a = 1$ we have an absolutely convergent expansion of the form
		\begin{align*}
		q(a) = &q_0(a) + \sum_{i=1}^\infty (1-a)^{\beta(i)+\frac{1}{2}}\sum_{j=0}^{j(i)}q_{i,j}(a)\big(\log(1-a)\big)^j\\
		&+ \sum_{i=1}^\infty (1-a)^{\tilde{\beta}(i)+\frac{1}{2}}\sum_{j=0}^{j(i)}\tilde{q}_{i,j}(a)\big(\log(1-a)\big)^j
		\end{align*}
		with analytic coefficients $q_0, q_{i,j}$, and $\beta(i) = i\nu$, $\tilde{\beta}(i) = \nu i+\frac{1}{2}$.  
		\item $\mathcal{Q}_n$ is the algebra which is defined similarly, but also requiring $q_{i,j}(1) = 0$ if $i\geq 2n+1$. 
		\item $\mathcal{Q}_n^{(N)}$ is the vector space of functions in $\mathcal{Q}_n$ for which $q_{i,j}(a) = 0$ for $i\geq N$. 
	\end{itemize}
\end{defi}
For later purposes, we also define the space of functions obtained by differentiating $\mathcal{Q}_n$: 
\begin{defi}\label{defi:Q'_n}
	Define $\mathcal{Q}'$ as in the preceding definition but replacing $\beta(i)$ by $\beta'(i): = \beta(i) - 1$, and similarly for $\mathcal{Q}'_n$.  
\end{defi}

\begin{defi}\label{defi:S^mQ_n}  Pick $t$ sufficiently small such that both $b_1, b_2$ (recall their definitions in Definition \ref{def:KST1}), when restricted to the light cone $r\leq t$ are of size at most $b_0$. 
	\begin{itemize}
		\item $S^m(R^k(\log R)^l, \mathcal{Q}_n)$ is the class of analytic functions $v: [0,\infty)\times [0,1)\times [0,b_0]^2\rightarrow\R$ so that\\
		(i) $v$ is analytic as a function of $R, b_1, b_2$, 
		\[
		v: [0,\infty)\times [0, b_0]^2\rightarrow {\mathcal{Q}}_n
		\]
		(ii) $v$ vanishes to order $m$ at $R = 0$.\\
		(iii) $v$ admits a convergent expansion at $R = \infty$, 
		\[
		v(R,\cdot,b_1, b_2) = \sum_{\substack{0\leq j\leq l+i\\ i\geq 0}}c_{ij}(\cdot, b_1, b_2)R^{k-i}(\log R)^j  
		\]
		where the coefficients $c_{ij}: [0, b_0]^2\rightarrow \mathcal{Q}_n$ are analytic with respect to $b_{1,2}$. 
		\item $IS^m(R^k(\log R)^l, \mathcal{Q}_n)$ is the class of analytic functions $w$ inside the cone $r<t$ which can be represented as 
		\[
		w(t, r) = v(R, a, b_1, b_2),\quad v\in S^m(R^k(\log R)^l, \mathcal{Q}_n)
		\]
		and $t>0$ sufficiently small.  
		\item Define $IS^m(R^k(\log R)^l, \mathcal{Q}_n^{(N)})$ analogously by replacing $\mathcal{Q}_n$ by $ \mathcal{Q}_n^{(N)}$. 
	\end{itemize}
\end{defi}

Then the core content of \cite{KST,GaoK}, customized for the application we have in mind, can be encapsulated by the following theorem:
\begin{theorem}\label{thm:KSTGaoKBase} Let $\nu>0$. Then, given $N>0$, there exist $t_0 = t_0(\nu, N)>0$, $M = M(\nu,N)$, $P = P(\nu)$ and a co-rotational blow up solution $U = Q(\lambda(t)r) + \epsilon(t, r)$,\,$\lambda(t) = t^{-1-\nu}$, on $(0, t_0]\times \R^2$, and such that 
	$\epsilon(t, r)$ admits the following fine structure:
	\begin{align*}
	\epsilon = \epsilon_1 + \epsilon_2, 
	\end{align*}
	where we have\footnote{The condition on $\epsilon_2$ is chosen to conform with our requirements; one could impose any $H^s_{R\,dR}$ condition here.}  (with $R = \lambda(t)r,\,\tau = \int_t^\infty\lambda(s)\,ds$)
	\[
	\left\|\epsilon_2(t, \cdot)\right\|_{\ H^{5+\delta}_{R\,dR}}\lesssim \frac{1}{\tau^N},\quad \tau = \tau(t)\in [\tau_0,\infty),\quad \tau_0 = \int_{t_0}^\infty\lambda(s)\,ds,
	\]
	for the 'small but unstructured term' $\epsilon_2$, while we have 
	\[
	\epsilon_1\in \sum_{k=1}^M\frac{1}{(t\lambda)^{2k}}IS^3\left(R(\log R)^{2k-1}, \mathcal{Q}_{k-1}^{(P)}\right) + \sum_{k=1}^M\frac{1}{(t\lambda)^{2k+2}}IS^3\left(R^3(\log R)^{2k-1}, \mathcal{Q}_{k}^{(P)}\right).
	\]
\end{theorem}
The additional restriction here replacing $ \mathcal{Q}_{k}$ by $\mathcal{Q}_{k}^{(P)}$ here is of course simply a reflection of the fact that those terms the expansion of a general function $q(a)\in \mathcal{Q}_{k}$ with sufficiently large power of $(1-a)$ are automatically in $ H^{5+\delta}_{R\,dR}$. 
\\

In order to fit the functions described in the preceding into our calculus set up on the Fourier side, we have to implement a concise translation of the algebraic structure on the physical side to the (distorted) Fourier side, for each angular momentum $n$. This is what we do in the next subsection. It is also here where we start to keep track of the {\it{temporal decay}} (as measured in terms of the variable $\tau$), which will start to play an important role for the final iterative scheme, of course. 
\subsection{Description of the shock on the distorted Fourier side I: prototypical expansions}

In order to handle the source terms near the light cone, we shall have to introduce a suitable algebra of functions which contains all the possible singular terms, is formulated on the distorted Fourier side, and flexible enough that it is preserved under certain para-differential operations, as well as the transference operators and the solution operator for the wave equation on the Fourier side. Due to the somewhat complicated nature of the function spaces concerned, we do this in a two-step fashion. First, we introduce a prototype which handles certain key aspects (it gives the correct translation from the physical singularity as detailed in the preceding subsection to the Fourier side at fixed times), but is still not adequate due to its incompatibility with the solution operator for the wave operator on the Fourier side. Nonetheless it will be useful to describe the nonlinear source terms on the distorted Fourier side. Once we have studied some of the key properties in this simplified setting, we can finally introduce the spaces that will overcome all hurdles. We already introduce aspects of the somewhat complicated temporal bounds for these function spaces, although strictly speaking this is not important yet at this stage. It will play a very prominent role later for the multilinear estimates, of course. 

\begin{definition}\label{defi:xsingulartermsngeq2proto} We call a function $\xb(\tau,\xi)$, $\xi\in [0,\infty)$, an prototype singular part at angular momentum $n, |n|\geq 2$, provided it allows an expansion (here all cutoffs are smooth localizers to the indicated regions)
\begin{align*}
&\xb(\tau, \xi)\\& = \sum_{\pm}\sum_{k=1}^{N}\sum_{j=0}^{N_1}\chi_{\xi\geq \hbar^{-2}}\hbar^{-1}\frac{e^{\pm i\nu\tau\xi^{\frac12}}}{\xi^{1+k\frac{\nu}{2}}}\left(\log\xi\right)^j\cdot a_{k,j}^{\pm}(\tau) + \sum_{\pm}\sum_{k=1}^{N}\sum_{j=0}^{N_1}\chi_{\xi\geq \hbar^{-2}}\hbar^{-1}\frac{e^{\pm i(\nu\tau\xi^{\frac12} + \hbar^{-1}\rho(x_{\tau};\alpha,\hbar))}}{\xi^{1+k\frac{\nu}{2}}}\left(\log\xi\right)^j\cdot \tilde{a}_{k,j}^{\pm}(\tau)\\&  + \sum_{\pm}\sum_{l =1}^7\sum_{k=1}^{N}\sum_{j=0}^{N_1}\chi_{\xi\geq \hbar^{-2}}\hbar^{-1}\langle\hbar^2\xi\rangle^{-\frac{l}{4}}\frac{e^{\pm i\nu\tau\xi^{\frac12}}}{\xi^{1+k\frac{\nu}{2}}}\left(\log\xi\right)^j\cdot F_{l,k,j}^{\pm}(\tau, \xi)\\
&\hspace{4.5cm} + \sum_{\pm}\sum_{l =1}^7\sum_{k=1}^{N}\sum_{j=0}^{N_1}\chi_{\xi\geq \hbar^{-2}}\hbar^{-1}\langle\hbar^2\xi\rangle^{-\frac{l}{4}}\frac{e^{\pm i(\nu\tau\xi^{\frac12} + \hbar^{-1}\rho(x_{\tau};\alpha,\hbar))}}{\xi^{1+k\frac{\nu}{2}}}\left(\log\xi\right)^j\cdot \tilde{F}_{l,k,j}^{\pm}(\tau, \xi)
\end{align*}
	Here $N_1 = N_1(N,\nu)$ is sufficiently large, $N = \lfloor\frac{5}{\nu}\rfloor + 1$,  $x_{\tau} = \xi^{\frac12}\hbar\cdot\nu\tau$ and $\rho$ is as in Lemma \ref{lem: Lemma 3.4 CDST}, and the functions $a_{k,j}^{\pm}(\tau), \tilde{a}_{k,j}^{\pm}(\tau), F_{k,j}^{\pm}(\tau, \xi), F_{l,k,j}^{\pm}(\tau, \xi), \tilde{F}_{l,k,j}^{\pm}(\tau, \xi)$ have the following properties:
	\begin{itemize}
		\item  Pointwise bounds with (weak) temporal decay $\left| a_{k,j}^{\pm}(\tau)\right| + \left|\tilde{a}_{k,j}^{\pm}(\tau)\right| \lesssim \left(\log\tau\right)^{N_1-j}\cdot \tau^{-1-\nu}$,  as well as 
		$\left|F_{l,k,j}^{\pm}(\tau, \xi)\right| + \left|\tilde{F}_{l,k,j}^{\pm}(\tau, \xi)\right|\lesssim \left(\log\tau\right)^{N_1-j}\tau^{-1-\nu}$. 
		\item Symbol type behavior: with $k_1\in \{0,1\}$, $k_2\in \{0,1,\ldots,20 - l\}$, $l\in \{1,2,\ldots,7\}$, and $\delta, \delta_1$ small absolute positive constants,  
		\begin{align*}
		&\left|\partial_{\tau}^{k_1}a_{k,j}^{\pm}(\tau)\right| + \left|\partial_{\tau}^{k_1}\tilde{a}_{k,j}^{\pm}(\tau)\right| \lesssim \left(\log\tau\right)^{N_1-j}\cdot \tau^{-1-\nu-k_1}\\
		&\left|\partial_{\tau}^{k_1}\partial_{\xi}^{k_2}F_{l,k,j}^{\pm}(\tau, \xi)\right| + \left|\partial_{\tau}^{k_1}\partial_{\xi}^{k_2}\tilde{F}_{l,k,j}^{\pm}(\tau, \xi)\right|\lesssim \left(\log\tau\right)^{N_1-j}\cdot \hbar^{-k_1}\tau^{-1-\nu-k_1}\cdot \hbar^{-\delta_1\cdot k_2}\xi^{-k_2},
		\end{align*}
		as well as the `closure bounds'
		\[
		\left\|\xi^{20 - l +\delta }\partial_{\tau}^{k_1}\partial_{\xi}^{20 - l}F_{l,k,j}^{\pm}(\tau, \xi)\right\|_{\dot{C}^{\delta}_{\xi}} + \left\|\xi^{20 - l +\delta-k_1}\partial_{\tau}^{k_1}\partial_{\xi}^{20 - l-k_1}\tilde{F}_{l,k,j}^{\pm}(\tau, \xi)\right\|_{\dot{C}^{\delta}_{\xi}}\lesssim  \left(\log\tau\right)^{N_1-j}\cdot \hbar^{-k_*}\cdot \tau^{-1-\nu - k_1},
		\]
		where $k_*= k_*(\delta_1,\delta,l,k_1): = (20-l)\delta_1 + \delta + k_1$, for some $\delta\in (0,1)$ we  set 
		\[
		\left\|g\right\|_{\dot{C}^{\delta}}: = \sup_{x\neq y}\frac{|g(x) - g(y)|}{|x-y|^{\delta}}
		\]
		We shall refer to the first two sums 
		\[
		\sum_{\pm}\sum_{k=1}^{N}\sum_{j=0}^{N_1}\chi_{\xi\geq \hbar^{-2}}\hbar^{-1}\frac{e^{\pm i\nu\tau\xi^{\frac12}}}{\xi^{1+k\frac{\nu}{2}}}\left(\log\xi\right)^j\cdot a_{k,j}^{\pm}(\tau),\quad \sum_{\pm}\sum_{k=1}^{N}\sum_{j=0}^{N_1}\chi_{\xi\geq \hbar^{-2}}\hbar^{-1}\frac{e^{\pm i(\nu\tau\xi^{\frac12} + \hbar^{-1}\rho(x_{\tau};\alpha,\hbar))}}{\xi^{1+k\frac{\nu}{2}}}\left(\log\xi\right)^j\cdot \tilde{a}_{k,j}^{\pm}(\tau) 
		\]
		as the principal singular part, and the second sum as the connecting singular part. This is because the second sum involves terms which are still too rough to be included into the space $\Sh_{0}$, but still much less singular than the principal singular part. Here `singular' refers of course to the physical side, and is measured in terms of decay of the Fourier transform with respect to $\xi$ for large frequencies. 
	\end{itemize}
	We say that a function $f(\tau, R)$ is a function at angular momentum $n$, $|n|\geq 2$, and prototypical singular part, provided it can be written as 
	\[
	f(\tau, R) = \int_0^\infty \phi_{n}(R;\xi)\cdot \xb(\tau, \xi)\cdot\rho_{n}(\xi)\,d\xi,
	\]
	and $\xb = \xb_1 + \xb_2$ where $\xb_1\in \Sh_{0}$ for each $\tau$, while $\xb_2$ is a prototypical singular part at angular momentum $n$. 
	\\
	Finally, we shall say that the principal singular part is of restricted type, provided we have the more specific structure 
	\[
	a_{k,i}^{\pm}(\tau) = \sum_{l=i}^{N_1}c_{l,k,i}^{\pm}\cdot \tau^{-k(1+\nu)}\cdot \left(\log\tau\right)^{l-i} + b_{k,i}^{\pm}(\tau),\quad k\in\{1,2,3\}
	\]
	for certain constants $c_{l,k,i}^{\pm}$, where $\left| b_{k,i}^{\pm}(\tau)\right| + \left|\tau\cdot  (b_{k,i}^{\pm})'(\tau)\right|\lesssim \tau^{-3-\nu}\left(\log\tau\right)^{N_1}$, and similarly for $\tilde{a}_{k,i}^{\pm}(\tau)$, and for $k\geq 4$ we have 
	\[
	\left|a_{k,i}^{\pm}(\tau)\right| + \left|\tau\cdot (a_{k,i}^{\pm})'(\tau)\right|\lesssim \tau^{-3-\nu}\left(\log\tau\right)^{N_1}. 
	\]
	Finally, we say that $\overline{y}(\tau,\xi)$ is a prototype singular source term, provided $\xi^{-\frac12}\cdot \overline{y}(\tau,\xi)$ is prototype singular in the previous sense, except we replace $\tau^{-1-\nu}$ by $\tau^{-2-\nu}$, and we set $k_1 = 0$ throughout. We also define the concept of a prototype singular source term with principal part of restricted type, in analogy to the preceding. 
\end{definition}
\begin{remark}\label{rem:xsingulartermsngeq2} It is important that the functions $a_{k,i}(\tau), \tilde{a}_{k,i}(\tau)$ do not depend on $\xi$. A similar property for the final admissible function space will play a key role when applying the modulation techniques for the exceptional modes $n = 0,\,\pm 1$. 
\end{remark}
\begin{remark}\label{rem:xsingulartermsngeq2proto2} The phase $e^{\pm i(\nu\tau\xi^{\frac12} + \hbar^{-1}\rho(x_{\tau};\alpha,\hbar))}$ appears quite natural in light of the asymptotic properties of the angular momentum $n$ Fourier basis in the oscillatory regime and near $R = \nu\tau$. However, the properties of the Duhamel parametrix will force a `deformation' of the function $ \hbar^{-1}\rho(x_{\tau};\alpha,\hbar)$ toward smaller and smaller values, effectively making the phase converge toward $e^{\pm i \nu\tau\xi^{\frac12}}$. Thus the `true phase' that will appear in the admissible function space defined below will in fact be an interpolate between these two phases, involving a continuously varying family of phases. The fact that we can work with such in a certain sense `unnatural phases' comes from an important monotonicity property of the function $\rho$ with respect to the variable $x$, as well as the fact that we always restriction attention to the interior of the light cone $0\leq R<\nu\tau$. 
\end{remark}

The preceding concept of function comes equipped with a natural concept of `norm': 
\begin{definition}\label{defy: protofunctionnorm} Let $\xb(\tau, \xi)$ be a prototype singular part at angular momentum $|n|\geq 2$. Then define 
	\begin{align*}
		\left\|\xb\right\|_{proto}: &= \sum_{k_1\in\{0,1\}}\sum_{k=1}^N\sum_{i=0}^{N_1}\left\|\tau^{1+\nu+k_1}\left(\log\tau\right)^{-N_1+i}\partial_\tau^{k_1} a_{k,i}^{\pm}(\tau)\right\|_{L^\infty([\tau_0,\infty))}\\
		& +  \sum_{k_1\in\{0,1\}}\sum_{k=1}^N\sum_{i=0}^{N_1}\left\|\tau^{1+\nu+k_1}\left(\log\tau\right)^{-N_1+i}\partial_\tau^{k_1} \tilde{a}_{k,i}^{\pm}(\tau)\right\|_{L^\infty([\tau_0,\infty))}\\
		& + \sum_{k_1\in\{0,1\}}\sum_{l=1}^7\sum_{k_2\in \{0,\ldots,20 - l\}}\hbar^{k_2\cdot\delta_1 + k_1}\left\| \tau^{1+\nu+k_1}\left(\log\tau\right)^{-N_1+i}\xi^k_2 \partial_\tau^{k_1}\partial_{\xi}^{k_2}F_{l,k,i}^{\pm}(\tau, \xi)\right\|_{L_{\tau,\xi}^{\infty}([\tau_0,\infty)\times [0,\infty))}\\
		& + \sum_{k_1\in\{0,1\}}\sum_{l=1}^7\sum_{k=1}^N\sum_{i=0}^{N_1}\sum_{k_2\in \{0,\ldots,20 - l\}}\hbar^{k_2\cdot\delta_1 + k_1}\left\| \tau^{1+\nu+k_1}\left(\log\tau\right)^{-N_1+i}\xi^k_2 \partial_\tau^{k_1}\partial_{\xi}^{k_2}\tilde{F}_{l,k,i}^{\pm}(\tau, \xi)\right\|_{L_{\tau,\xi}^{\infty}([\tau_0,\infty)\times [0,\infty))}\\
		& +  \sum_{k_1\in\{0,1\}}\sum_{l=1}^7\sum_{k=1}^N\sum_{i=0}^{N_1}\left\|\sup_{\lambda>0}\hbar^{k_*}\left\| \tau^{1+\nu+k_1}\left(\log\tau\right)^{-N_1+i}\xi^{20-l+\delta} \partial_\tau^{k_1}\partial_{\xi}^{20 - l}F_{l,k,i}^{\pm}(\tau, \xi)\right\|_{\dot{C}^{\delta_l}(\xi\simeq\lambda)}\right\|_{L^\infty_{\tau}([\tau_0,\infty))}\\
		& + \sum_{k_1\in\{0,1\}}\sum_{l=1}^7\sum_{k=1}^N\sum_{i=0}^{N_1}\hbar^{k_*}\left\|\sup_{\lambda>0}\left\| \tau^{1+\nu+k_1}\left(\log\tau\right)^{-N_1+i}\xi^{20-l}\partial_\tau^{k_1}\partial_{\xi}^{20-l}\tilde{F}_{l,k,i}^{\pm}(\tau, \xi)\right\|_{\dot{C}^{\delta_l}(\xi\sim\lambda)}\right\|_{L^\infty_{\tau}([\tau_0,\infty))}
	\end{align*}
	In case the function $\xb$ has principal part of restricted type, we also use the norm $\left\|\xb\right\|_{proto(r)}$, where (keeping in mind the preceding definition) we replace the first expression on the right by 
	\[
	\sum_{k_1\in\{0,1\}}\sum_{k=1}^N\sum_{i=0}^{N_1}\left\|\tau^{3+\nu+k_1}\left(\log\tau\right)^{-N_1+i}\partial_\tau^{k_1} b_{k,i}^{\pm}(\tau)\right\|_{L^\infty([\tau_0,\infty))} + \sum_{k=1}^N\sum_{i=0}^{N_1}\sum_{l=1}^7\left|c_{l,k,i}^{\pm}\right|. 
	\]
	We define an analogous norm for a prototype singular source term, resulting in $\big\|\overline{y}\big\|_{sourceproto}, \big\|\overline{y}\big\|_{sourceproto(r)}$. 
\end{definition}

The first order of the day shall be to translate the information of the preceding definition to the physical side, i.e., identify a vector space of functions which correspond to the above Fourier description: 
\begin{lemma}\label{lem:singFouriertiphysicalngeq2} Assume that $\xb$ is a prototypical singular part at angular momentum $n, |n|\geq 2$. Then the associated function 
	\[
	f(\tau, R): = \int_0^\infty \phi_{n}(R;\xi)\cdot \xb(\tau, \xi)\cdot\rho_{n}(\xi)\,d\xi
	\]
	restricted to the light cone $R\leq \nu\tau$, can be decomposed as 
	\[
	f = f_1 + f_2+f_3
	\]
	where $f_1 = f_1(\tau, R)$ is a $C^\infty$-function supported in $\nu\tau - R\gtrsim 1$ and satisfying\footnote{The conclusion for $f_1$ is not sharp, it is stated in such a way that it implies in particular the function is in $\tilde{S}_0^{\hbar}$ at fixed time $\tau$.} 
	\[
	\big|\nabla_R^{k_2} \partial_{\tau}^{k_1}f_1(\tau, R)\big|\lesssim \hbar^{-\frac12+}\left(\log\tau\right)^{N_1}\cdot \tau^{-\frac32-\nu}\left|\nu\tau - R\right|^{-5}\quad k_1\in \{0,1\},\quad k_1 + k_2\leq 5, 
	\]
	while $f_2 =\sum_{l=1}^8 f_{2l}$ where we have the explicit form
	\begin{align*}
	&f_{2l}(\tau, R) = \chi_{|\nu\tau-R|\lesssim \hbar}\sum_{k=1}^N\sum_{i=0}^{N_1}\frac{G_{k,l,i}(\tau,\nu\tau - R)}{\tau^{\frac12}}\hbar^{-\frac{l+1}{2}}\left[\nu\tau - R\right]^{\frac{l}{2}+k\nu}\left(\log(\nu\tau - R)\right)^i
	\end{align*}
	Here the function $G_{k,i}(\tau, x)$ has symbol type behavior with respect to $R, x$, as follows:
	\begin{align*}
	\hbar^{\delta_1\cdot \min\{k_2, 20-l\}}\left|\partial_{\tau}^{k_1}\partial_x^{k_2}G_{k,l,i}(\tau, x)\right| \lesssim  \left(\log\tau\right)^{N_1-i}\cdot\hbar^{-k_1}\tau^{-1-\nu-k_1}x^{-k_2},\quad k_1\in \{0,1\},\quad k_2\in\{0,\ldots,20 - l + \lfloor \frac{l}{2}+k\nu\rfloor\},
	\end{align*}
	and we have the bound 
	\begin{align*}
	\hbar^{k_*}\cdot \left\|x^{20-l_*+\delta_3}\partial_{\tau}^{k_1}\partial_x^{20 - l}G_{k,l,i}(\tau, x)\right\|_{\dot{C}^{\delta_3}} \lesssim \left(\log\tau\right)^{N_1-i}\cdot\tau^{-1-\nu-k_1},\quad k_1\in \{0,1\},\,&\delta_3 = \delta_2 + \frac{l}{2}+k\nu - \lfloor \frac{l}{2}+k\nu\rfloor\\
	&l_*: =  l - \lfloor \frac{l}{2}+k\nu\rfloor,
	\end{align*}
	and $\delta_2>\delta_1$. Also, as before we set $k_* = (20-l)\delta_1+ \delta + k_1$. 
		
	Finally, the remaining function $f_3$ is $C^{\infty}$ and supported in $|\nu\tau - R|\lesssim 1$ and satisfies
	 \[
	\left\|f_3\right\|_{\tilde{S}_0^{(n)}}\lesssim \left(\log\tau\right)^{N_1}\cdot\tau^{-1-\nu},\quad \left\|\partial_{\tau}f_3\right\|_{\tilde{S}_1^{(n)}}\lesssim  \left(\log\tau\right)^{N_1}\cdot \tau^{-1-\nu}.
	\]
	Moreover, we have the bounds 
	\begin{align*}
	&\left|\partial_R^{k_1} f_3\right|\lesssim \left(\log\tau\right)^{N_1}\tau^{-\frac32-\nu}\cdot\hbar^{-\frac12+}\min\{(\nu\tau - R)^{-k_1},\hbar^{-k_1}\},\\& \left|\partial_R^{k_2} \partial_{\tau}f_3\right|\lesssim \left(\log\tau\right)^{N_1}\tau^{-\frac32-\nu}\cdot\hbar^{-\frac12+}\min\{(\nu\tau - R)^{-k_2-1},\hbar^{-k_2-1}\},\\&
	0\leq k_1\leq 5,\,0\leq k_2\leq 4. 
	\end{align*}
\end{lemma}
\begin{proof}
Due to linearity, we can separately consider the contributions of the different parts constituting $\xb$. 
\\
We deal with the contribution of\footnote{Here we omit the factor $\hbar^{-1}$ for simplicity.} $ \sum_{\pm}\sum_{k=1}^{N}\sum_{j=0}^{N_1}\chi_{\xi\geq \hbar^{-2}}\frac{e^{\pm i\nu\tau\xi^{\frac12}}}{\xi^{1+k\frac{\nu}{2}}}\left(\log\xi\right)^j\cdot a_{k,j}^{\pm}(\tau)$, and later explain what
changes are necessary to handle the terms in the second to the fourth sum constituting $\xb(\tau,\xi)$, according to the definition. We split the Fourier integral into a number of regions reflecting the oscillatory/non-oscillatory nature of $\phi_{n}(R;\xi)$. Also we may as well omit the factors $a_{k,j}^{\pm}(\tau)$, modifying the required bounds accordingly. We split the integral for $f(\tau,R)$ into a number of regions.
\\

{\it{(1) The rapidly decaying region: $R\xi^{\frac12}\hbar<\frac{x_t}{2}$. }} Let $\psi(x)$ be a smooth function which is identically $1$ for $x\leq \frac12$ and vanishes beyond $x = 1$. Then consider the integral 
\begin{align*}
\int_0^\infty\psi\left(\frac{2x}{x_t}\right)\phi_{n}(R;\xi)\cdot \chi_{\xi\geq \hbar^{-2}}\frac{e^{\pm i\nu\tau\xi^{\frac12}}}{\xi^{1+k\frac{\nu}{2}}}\left(\log\xi\right)^j\cdot\rho_{n}(\xi)\,d\xi,\quad x = R\xi^{\frac12}\hbar.
\end{align*}
Taking advantage of the asymptotics for $\phi_{n}(R;\xi)$ in the non-oscillatory regime away from the turning point $x_t$, as displayed in Proposition~\ref{prop:DFT nlarge} but with an additional factor $R^{-\frac12}$ to transform things back into the context of $\mathbb{R}^2$,  this becomes 
\begin{align*}
\int_0^\infty\psi\left(\frac{2x}{x_t}\right)\hbar^{\frac13}x^{-\frac12}q^{-\frac14}(\zeta)\cdot\Ai\left(\hbar^{-\frac23}\zeta\right)\cdot\left(1+\hbar a_0(-\zeta;\alpha,\hbar)\right) \chi_{\xi\geq \hbar^{-2}}\frac{e^{\pm i\nu\tau\xi^{\frac12}}}{\xi^{1+k\frac{\nu}{2}}}\left(\log\xi\right)^j\cdot\rho_{n}(\xi)\,d\xi
\end{align*}
We claim that the preceding function is smooth with respect to $R$ and $\tau$, and that the function as well as its derivatives can be bounded by $\lesssim_N \tau^{-N}\lesssim \left|R-\nu\tau\right|^{-N}$ for any $N$. To see this, distinguish between the cases $R\ll 1$ and $R\gtrsim 1$, and use that $\alpha = \xi^{\frac12}\hbar\gtrsim 1$ on the support of the integrand . Then use Lemma \ref{lem: Lemma 3.3 CDST} for the description of $\zeta$, and perform integration by parts with respect to $\xi^{\frac12}$, where we also have to take advantage of Lemma \ref{lem:a0fine1}. To initiate the integration by parts, write 
\[
e^{\pm i\nu\tau\xi^{\frac12}} = \frac{1}{\pm i\nu\tau}\cdot \partial_{\xi^{\frac12}}\big(e^{\pm i\nu\tau\xi^{\frac12}} \big). 
\]
 Also, use Proposition \ref{prop:DFT nlarge} for the symbol behavior of the spectral measure. It follows from these sources that the worst case occurs when $\partial_{\xi^{\frac12}}$ hits $\Ai\left(\hbar^{-\frac23}\zeta\right)$, which in light of the Airy function asymptotics leads to a factor 
\[
\hbar^{-1}\cdot x^{-1}\cdot R\hbar = \frac{R}{x} = \frac{1}{\xi^{\frac12}\hbar}. 
\]
Then performing $N$ times integration by parts with respect to $\xi^{\frac12}$, we get the bound 
\begin{align*}
&\left|\int_0^\infty\psi\left(\frac{2x}{x_t}\right)\hbar^{\frac13}x^{-\frac12}q^{-\frac14}(\zeta)\cdot\Ai\left(\hbar^{-\frac23}\zeta\right)\cdot\left(1+\hbar a_0(-\zeta;\alpha,\hbar)\right) \chi_{\xi\geq \hbar^{-2}}\frac{e^{\pm i\nu\tau\xi^{\frac12}}}{\xi^{1+k\frac{\nu}{2}}}\left(\log\xi\right)^j\cdot\rho_{n}(\xi)\,d\xi\right|\\
&\lesssim_{N}\tau^{-N}\int_0^\infty (\xi^{\frac12}\hbar)^{-N}\cdot x^{\hbar^{-1} - 3}\cdot  \frac{\chi_{\xi\geq 1}}{\xi^{1+k\frac{\nu}{2}}}\left(\log\xi\right)^{j}\,d\xi\lesssim \tau^{-N}
\end{align*}
since we have $\sup_{0<\hbar\leq \frac13}\hbar^{-N} x^{\hbar^{-1} - 3}\lesssim_N 1$.  Derivatives with respect to $R,\tau$ are handled similarly, as the additional factors $\xi^{\frac12}$ can be easily absorbed. 
\\

{\it{(2) Region near the turning point, still non-oscillatory: $x_t + \hbar^{\frac23}>R\xi^{\frac12}\hbar\geq \frac{x_t}{2}$. }} Let $\vphi(x)$ be non-vanishing on $[-1, 1]$ and zero on $[-2,2]^c$, chosen in such a way that 
\[
\sum_{\lambda}\vphi\left(\frac{x - \left[x_t - \lambda \hbar^{\frac23}\right]}{\frac12\langle\lambda\rangle \hbar^{\frac23}}\right) = 1,\quad x\in (0, x_{t}+\hbar^{\frac23}]. 
\]
where $\lambda$ ranges over $0\cup\{2^j\}_{j\geq 0}$. Then, in light of the fact that the above asymptotics for $\phi_{n}(R;\xi)$ are still valid in this regime, but with a different law for $\tau$, we reduce to bounding 
\begin{align*}
\sum_{\lambda}\int_0^\infty\left[1-\psi\left(\frac{2x}{x_t}\right)\right]\cdot \vphi\left(\frac{x - \left[x_t - \lambda \hbar^{\frac23}\right]}{\frac12\langle\lambda\rangle \hbar^{\frac23}}\right)\cdot \Xi(R;\xi,\hbar) \chi_{\xi\geq \hbar^{-2}}\frac{e^{\pm i\nu\tau\xi^{\frac12}}}{\xi^{1+k\frac{\nu}{2}}}\left(\log\xi\right)^j\rho_{n}(\xi)\,d\xi,\quad x = R\xi^{\frac12}\hbar.  
\end{align*}
where we set 
\[
\Xi(R;\xi,\hbar) = \hbar^{\frac13}x^{-\frac12}q^{-\frac14}(\zeta)\cdot\Ai\left(\hbar^{-\frac23}\zeta\right)\cdot\left(1+\hbar a_0(-\zeta;\alpha,\hbar)\right). 
\]
Also, $\zeta$ is described by Lemma \ref{lem: Lemma 3.2 CDST}. Then since 
\[
\hbar^{-\frac23}\partial_{\xi^{\frac12}}\zeta\simeq \hbar^{\frac13}R,\quad R\xi^{\frac12}\hbar\simeq 1
\]
in the support of the integrand, we easily infer using integration by parts and taking into account the asymptotics of Airy functions, that 
\begin{align*}
&\left|\int_0^\infty\left[1-\psi\left(\frac{2x}{x_t}\right)\right]\cdot \vphi\left(\frac{x - \left[x_t - \lambda \hbar^{\frac23}\right]}{\frac12\langle\lambda\rangle \hbar^{\frac23}}\right)\cdot \Xi(R;\xi,\hbar) \chi_{\xi\geq \hbar^{-2}}\frac{e^{\pm i\nu\tau\xi^{\frac12}}}{\xi^{1+k\frac{\nu}{2}}}\left(\log\xi\right)^j\rho_{n}(\xi)\,d\xi\right|\\
&\lesssim_N \left(\hbar^{\frac13}R\right)^N\cdot \tau^{-N}\cdot e^{-c\lambda^{\frac32}}\lesssim \left(\frac{R}{\tau}\right)^N\cdot \left(\frac{1}{R\xi^{\frac12}}\right)^{\frac{N}{3}}\cdot e^{-c\lambda^{\frac32}}\lesssim \tau^{-\frac{N}{3}}\cdot\xi^{-\frac{N}{6}}\cdot e^{-c\lambda^{\frac32}}.
\end{align*}
Again derivatives with respect to $R$ and $\tau$ are easily absorbed by the gain in $\xi^{-\frac12}$, and the rapid exponential decay in $\lambda$ allows for summation over dyadic $\lambda$. Here we also used the fact that inside the cone $r\leq t$ we have $R\leq \nu\tau< \tau$. The contributions from the derivatives on $a_{0}(-\zeta;\alpha,\hbar)$ are handled similarly by Lemma \ref{lem:a0fine1} and Lemma \ref{lem:a0fine2}.
\\

{\it{(3) Region near the turning point, oscillatory regime: $x_t + \hbar^{\frac23}\leq R\xi^{\frac12}\hbar <(1+\gamma)x_t, 0<\gamma\ll1$. }} Here the asymptotics of the $\phi_{n}(R;\xi)$ are of the oscillatory kind, with a weaker decay. Proceeding as in the preceding case and letting $\psi$ a smooth function which equals $1$ on $[0,1]$ but vanishes beyond $2$, we reduce to bounding (with $\lambda$ again taking dyadic values)
\begin{align*}
\sum_{\lambda\geq 1}\int_0^\infty\psi\left(\frac{x}{(1+\gamma)x_t}\right)\cdot \vphi\left(\frac{x - [x_t +\lambda \hbar^{\frac23}]}{\lambda \hbar^{\frac23}}\right)\cdot \Xi(R;\xi,\hbar) \chi_{\xi\geq \hbar^{-2}}\frac{e^{\pm i\nu\tau\xi^{\frac12}}}{\xi^{1+k\frac{\nu}{2}}}\left(\log\xi\right)^j\rho_{n}(\xi)\,d\xi,\quad x = R\xi^{\frac12}\hbar, 
\end{align*}
where we set this time 
\begin{align*}
\Xi(R;\xi,\hbar) = \hbar^{\frac13}x^{-\frac12}q^{-\frac14}(\zeta)\Re\left[2a(\xi)\left(\Ai(-\hbar^{-\frac23}\zeta) - i\Bi(-\hbar^{-\frac23}\zeta)\right)(1+\hbar \overline{a_1(\zeta;\alpha)})\right]
\end{align*}
recalling Proposition \ref{prop:DFT nlarge}. Then using the oscillatory asymptotics of the complex Airy function, we get a phase $e^{\pm\frac23 i\hbar^{-1}\zeta^{\frac32}}$, and we have (recall Lemma \ref{lem: Lemma 3.2 CDST})
\begin{align*}
\partial_{\xi^{\frac12}}\left[\hbar^{-1}\zeta^{\frac32}\right] &= \frac32\hbar^{-1}\zeta^{\frac12}\cdot\left[R\hbar\cdot \Phi(x;\alpha,\hbar) + (1+O(x-x_t))\cdot O(\hbar(1+\alpha)^{-3})+O(x-x_{t})\cdot R\hbar\right]\\
& \simeq \zeta^{\frac12}\cdot R, 
\end{align*}
where we have used that $\left|O(\hbar(1+\alpha)^{-3}) \right|\lesssim \hbar \cdot \left(\hbar\xi^{\frac12}\right)^{-1} = O(\hbar R)$. In particular, for the combined phase $e^{\pm\frac23 i\hbar^{-1}\zeta^{\frac32}\pm i\nu\tau\xi^{\frac12}}$, we have the identity 
\begin{align*}
e^{\pm\frac23 i\hbar^{-1}\zeta^{\frac32}\pm i\nu\tau\xi^{\frac12}} &= \frac{1}{\partial_{\xi^{\frac12}}\left[\pm\frac23 i\hbar^{-1}\zeta^{\frac32}\pm i\nu\tau\xi^{\frac12}\right]}\cdot \partial_{\xi^{\frac12}}\left(e^{\pm\frac23 i\hbar^{-1}\zeta^{\frac32}\pm i\nu\tau\xi^{\frac12}} \right)\\
& \simeq \frac{1}{\pm\nu\tau\pm \zeta^{\frac12}R}\cdot \partial_{\xi^{\frac12}}\left(e^{\pm\frac23 i\hbar^{-1}\zeta^{\frac32}\pm i\nu\tau\xi^{\frac12}} \right),
\end{align*}
and since $\left|\zeta\right|\ll 1$ on the support of the integrand, we have $ \frac{1}{\pm\nu\tau\pm \zeta^{\frac12}R}\simeq \frac{1}{\tau}$ for $R\leq\nu\tau$. If we then again take advantage of the bound $\hbar^{-\frac23}\partial_{\xi^{\frac12}}\zeta\simeq \hbar^{\frac13}R$ and perform repeated integration by parts (after combining  the oscillatory phases), we obtain the gain (similar as in the previous regime)
\begin{align*}
&\left|\int_0^\infty\psi\left(\frac{x}{(1+\gamma)x_t}\right)\cdot \vphi\left(\frac{x - \left[x_t +\lambda \hbar^{\frac23}\right]}{\lambda \hbar^{\frac23}}\right)\cdot \Xi(R;\xi,\hbar) \chi_{\xi\geq 1}\frac{e^{\pm i\nu\tau\xi^{\frac12}}}{\xi^{1+k\frac{\nu}{2}}}\left(\log\xi\right)^j\rho_{n}(\xi)\,d\xi\right|\\
&\lesssim \lambda^{-\frac14}\cdot \tau^{-\frac{N}{3}}\cdot\xi^{-\frac{N}{6}}, 
\end{align*}
and this can again be summed over dyadic $\lambda\geq 1$. Derivatives are again handled similarly. Observe that up to this point all terms have been of the rapidly decaying type. 
\\

{\it{(4) Region away from the turning point and in oscillatory regime: $R\xi^{\frac12}\hbar \geq \left(1+\gamma\right)x_t$. }} Here, recalling Lemma \ref{lem: Lemma 3.4 CDST}, we have to consider the integral 
\begin{align*}
&\int_0^\infty \left[1- \psi\left(\frac{x}{(1+\gamma)x_t}\right)\right]\cdot\hbar^{\frac13}x^{-\frac12}q^{-\frac14}(\zeta)\cdot \left(\hbar^{-\frac23}\zeta\right)^{-\frac14}\cdot e^{\pm i\hbar^{-1}\left[x - y(\alpha;\hbar) + \rho(x;\alpha,\hbar)\right]}\\&\hspace{6cm}\cdot\left(1+b(x;\alpha,\hbar)\right) \chi_{\xi\geq \hbar^{-2}}\frac{e^{\pm i\nu\tau\xi^{\frac12}}}{\xi^{1+k\frac{\nu}{2}}}\left(\log\xi\right)^j\rho_{n}(\xi)\,d\xi,
\end{align*}
where $1+b(x;\alpha,\hbar)$ accounts for the factor $1+\hbar \overline{a_1(\zeta;\alpha)}$ as well as the corrections entailed by the Airy function asymptotics. Here we can combine the phases into 
\[
e^{i(\pm\nu\tau\pm R)\xi^{\frac12}}\cdot e^{\pm i\hbar^{-1}[\rho(x;\alpha,\hbar) - y(\alpha;\hbar)]}
\]
In the case of non-resonance, i.e., when the phases $e^{\pm i(\nu\tau + R)\xi^{\frac12}}$ occur, we can prove rapid decay for the contribution as well as all its derivatives with respect to $\tau$, whence also $|\nu\tau - R|$, within a dilate of the light cone. For this we also need to take into account the phase function $\hbar^{-1}\rho(x;\alpha,\hbar)$, for which we have, interpreting $x = R\hbar\xi^{\frac12}$ as a function of $R$ and $\xi$,
\[
\partial_{\xi^{\frac12}}\left(\hbar^{-1}\rho(x;\alpha,\hbar)-\hbar^{-1}y(\alpha,\hbar)\right) = -\frac{Rx^{-2}(1-2\hbar)^2}{1+\sqrt{-Q_0\big(x;\alpha,\hbar)}} + O\left(\hbar\right),
\]
where we note that the integration takes place over the region $x\geq (1+\gamma)x_t>1$. We conclude that 
\begin{align*}
\partial_{\xi^{\frac12}}\left[(\nu\tau + R)\xi^{\frac12}  - \hbar^{-1}y(\alpha,\hbar)+  \hbar^{-1}\rho(x;\alpha,\hbar)\right]\gtrsim \tau,
\end{align*}
and the claim follows after repeated integration by parts with respect to $\xi^{\frac12}$. 
\\
We may hence assume that there is destructive resonance, meaning the phases $e^{\pm i(\nu\tau - R)\xi^{\frac12}}$ occur. Denoting for now 
\begin{equation}\label{eq:Psixalphahbardef}
\left[1- \psi\left(\frac{x}{(1+\gamma)x_t}\right)\right]\cdot\hbar^{\frac13}x^{-\frac12}q^{-\frac14}(\zeta)\cdot \left(\hbar^{-\frac23}\zeta\right)^{-\frac14}\cdot\left(1+b(x;\alpha,\hbar)\right)\cdot\chi_{\xi\geq \hbar^{-2}}\frac{\left(\log\xi\right)^j\rho_{n}(\xi)}{\xi^{1+k\frac{\nu}{2}}} =:\Psi(x;\alpha,\hbar), 
\end{equation}
we reduce to bounding 
\begin{align*}
\int_0^\infty e^{\pm i(\nu\tau - R)\xi^{\frac12}}\cdot e^{\mp i\hbar^{-1}\left[\rho(x;\alpha,\hbar) - y(\alpha;\hbar)\right]}\cdot \Psi(x;\alpha,\hbar)\,d\xi
\end{align*}
We shall decompose this into a number of pieces with different structure. First, observe that if $\nu\tau - R\gg 1$, then the expression is a $C^\infty$ function, which decays, in addition to all its derivatives, like $|\nu\tau-R|^{-N}\tau^{-\frac12}$. 
This follows as before by considering the phase $(\nu\tau - R)\xi^{\frac12} + \hbar^{-1}y(\alpha,\hbar)-  \hbar^{-1}\rho(x;\alpha,\hbar)$ and using the bound 
\begin{equation}\label{eq:phaseboundcombinedoscillatory}
\left|\partial_{\xi^{\frac12}}\left[(\nu\tau - R)\xi^{\frac12} + \hbar^{-1}y(\alpha,\hbar)-  \hbar^{-1}\rho(x;\alpha,\hbar)\right]\right|\gtrsim \nu\tau - R, 
\end{equation}
provided $\nu\tau-R\gg 1$, and in the regime $x>(1+\gamma)x_t>1$. Note that to derive \eqref{eq:phaseboundcombinedoscillatory}, we can use the fact that for $\hbar\ll1$, $\partial_{x}\rho(x;\alpha,\hbar)\leq0$\\
We next peel off a function in $\tilde{S}_0^{(n)}$, by considering the intermediate region $\hbar\ll \nu\tau - R\lesssim 1$. Observe that 
\begin{align*}
&\left|\partial_R\left[(\nu\tau - R)\xi^{\frac12} -  \hbar^{-1}\rho(x;\alpha,\hbar)\right]\right|\lesssim \xi^{\frac12} + \hbar^{-1}x^{-2}\cdot \hbar \xi^{\frac12}\lesssim \xi^{\frac12},\\
&\left|\partial_R \Psi(x;\alpha,\hbar)\right|\lesssim R^{-\frac32}\lesssim \tau^{-\frac32},
\end{align*}
under our assumptions on $R, \tau$, and similarly for the $j$-th derivatives with the right hand sides replaced by the $j$-th powers.  Furthermore, we have the estimate \eqref{eq:phaseboundcombinedoscillatory}, where the right hand side can be replaced by $\hbar$ under the current assumption. Moreover, we have the now somewhat more delicate relation 
\begin{align*}
\partial_{\xi^{\frac12}}\left[(\nu\tau - R)\xi^{\frac12} + \hbar^{-1}y(\alpha,\hbar)-  \hbar^{-1}\rho(x;\alpha,\hbar)\right] &= (\nu\tau - R) +\hbar^{-1}\frac{\partial\alpha}{\partial\xi^{\frac12}}\cdot\partial_{\alpha}y(\alpha,\hbar) - R\cdot\rho_x(x;\alpha,\hbar)\\
&\hspace{5.5cm} - \rho_{\alpha}(x;\alpha,\hbar)\\
&\geq (\nu\tau - R) +\hbar^{-1}\frac{\partial\alpha}{\partial\xi^{\frac12}}\cdot\partial_{\alpha}y(\alpha,\hbar) -  \rho_{\alpha}(x;\alpha,\hbar)\\
&\gtrsim \nu\tau - R
\end{align*}
due to Lemma \ref{lem: Lemma 3.4 CDST} and the assumption $\hbar\ll \nu\tau - R\lesssim 1$.

Again using integration by parts, we conclude that  
\begin{align*}
&\left\|\chi_{1\gtrsim\nu\tau - R\gg\hbar}\int_0^\infty \chi_{\xi\geq\hbar^{-2}}\cdot e^{\pm i(\nu\tau - R)\xi^{\frac12}}\cdot e^{\mp i\hbar^{-1}[\rho(x;\alpha,\hbar) - y(\alpha;\hbar)]}\cdot \Psi(x;\alpha,\hbar)\,d\xi\right\|_{\tilde{S}_0^{(n)}}\lesssim 1.
\end{align*}
Here the large bound $(\nu\tau-R)^{-k}\ll\hbar^{-k}$ is compensated by the power of $\hbar$ in the definition of the norm $\Sh_{0}$. It follows that this contribution can be absorbed into $f_3$ (in fact, it is easily seen that we also have the point wise bounds required in the lemma), and so we reduce to controlling the remaining part. \\

Next we need to distinguish between different regions in terms of the functions $(\nu\tau - R)\xi^{\frac12}, \hbar^{-1}\rho(x;\alpha,\hbar)$. 
\\

{\it{(4.i): $\hbar^{-1}\rho(x;\alpha,\hbar)\lesssim 1$.}} We enforce this condition by inclusion of a smooth cutoff $\chi_{\lesssim 1}\left(\hbar^{-1}\rho(x;\alpha,\hbar)\right)$, and we shall expand 
\[
\chi_{\lesssim 1}\left(\hbar^{-1}\rho(x;\alpha,\hbar)\right)e^{\mp i\hbar^{-1}\rho(x;\alpha,\hbar)} = \chi_{\lesssim 1}\left(\hbar^{-1}\rho(x;\alpha,\hbar)\right)\sum_{l=0}^\infty\frac{1}{l!}\left(\mp i\hbar^{-1}\rho(x;\alpha,\hbar)\right)^l. 
\]
Consider then a term (for $l\geq 0$)
\begin{align*}
&\chi_{\nu\tau - R\lesssim\hbar}\int_0^\infty\chi_{\lesssim 1}\left(\hbar^{-1}\rho(x;\alpha,\hbar)\right) e^{\pm i(\nu\tau - R)\xi^{\frac12}}\cdot \left(\mp i\hbar^{-1}\rho(x;\alpha,\hbar)\right)^l e^{\pm i\hbar^{-1}y(\alpha;\hbar)}\cdot \Psi(x;\alpha,\hbar)\,d\xi\\
& = \chi_{\nu\tau - R\lesssim\hbar}\int_0^\infty\chi_{\lesssim 1}\left((\nu\tau - R)\xi^{\frac12}\right)\chi_{\lesssim 1}\left(\hbar^{-1}\rho(x;\alpha,\hbar)\right) e^{\pm i(\nu\tau - R)\xi^{\frac12}}\cdot \left(\mp i\hbar^{-1}\rho(x;\alpha,\hbar)\right)^l e^{\pm i\hbar^{-1}y(\alpha;\hbar)}\cdot \Psi(x;\alpha,\hbar)\,d\xi\\
& + \chi_{\nu\tau - R\lesssim\hbar}\int_0^\infty\chi_{\gtrsim1}\left((\nu\tau - R)\xi^{\frac12}\right)\chi_{\lesssim 1}\left(\hbar^{-1}\rho(x;\alpha,\hbar)\right) e^{\pm i(\nu\tau - R)\xi^{\frac12}}\cdot \left(\mp i\hbar^{-1}\rho(x;\alpha,\hbar)\right)^l e^{\pm i\hbar^{-1}y(\alpha;\hbar)}\cdot \Psi(x;\alpha,\hbar)\,d\xi.
\end{align*}
For the first term here on the right containing the cutoff $\chi_{\lesssim 1}\left((\nu\tau - R)\xi^{\frac12}\right)$, we also expand the exponential $e^{\pm i(\nu\tau - R)\xi^{\frac12}}$ into a Taylor series, which leads us to consider the expressions (for $l\geq 0, k\geq 0$)
\begin{align*}
\chi_{\nu\tau - R\lesssim\hbar}\int_0^\infty\chi_{\lesssim 1}\left((\nu\tau - R)\xi^{\frac12}\right)\chi_{\lesssim 1}\left(\hbar^{-1}\rho(x;\alpha,\hbar)\right) \left[i(\nu\tau - R)\xi^{\frac12}\right]^k\cdot \left(\mp i\hbar^{-1}\rho(x;\alpha,\hbar)\right)^l e^{\pm i\hbar^{-1}y(\alpha;\hbar)}\cdot \Psi(x;\alpha,\hbar)\,d\xi.
\end{align*}
If $k = 0$, writing 
\begin{align}\label{smooth plus singular}
	\begin{split}
&\chi_{\nu\tau - R\lesssim\hbar}\int_0^\infty\chi_{\lesssim 1}\left((\nu\tau - R)\xi^{\frac12}\right)\chi_{\lesssim 1}\left(\hbar^{-1}\rho(x;\alpha,\hbar)\right)\left(\mp i\hbar^{-1}\rho(x;\alpha,\hbar)\right)^l e^{\pm i\hbar^{-1}y(\alpha;\hbar)}\cdot \Psi(x;\alpha,\hbar)\,d\xi\\
& = \chi_{\nu\tau - R\lesssim\hbar}\int_0^\infty\chi_{\lesssim 1}\left(\hbar^{-1}\rho(x;\alpha,\hbar)\right)\left(\mp i\hbar^{-1}\rho(x;\alpha,\hbar)\right)^l e^{\pm i\hbar^{-1}y(\alpha;\hbar)}\cdot \Psi(x;\alpha,\hbar)\,d\xi\\
& - \chi_{\nu\tau - R\lesssim\hbar}\int_0^\infty\chi_{\gtrsim 1}\left((\nu\tau - R)\xi^{\frac12}\right)\chi_{\lesssim 1}\left(\hbar^{-1}\rho(x;\alpha,\hbar)\right)\left(\mp i\hbar^{-1}\rho(x;\alpha,\hbar)\right)^l e^{\pm i\hbar^{-1}y(\alpha;\hbar)}\cdot \Psi(x;\alpha,\hbar)\,d\xi.
\end{split}
\end{align}
and using the simple bounds 
\begin{align*}
\left|\partial_R^k\left[\chi_{\lesssim 1}\left(\hbar^{-1}\rho(x;\alpha,\hbar)\right)\left(\mp i\hbar^{-1}\rho(x;\alpha,\hbar)\right)^l\right]\right|&\lesssim_k C^l\cdot l^k\cdot \tilde{\chi}_{\lesssim 1}\left(\hbar^{-1}\rho(x;\alpha,\hbar)\right)\frac{\left(\hbar\xi^{\frac12}\right)^{k}}{x^{k}}
\lesssim l^k\cdot \tau^{-k}
\end{align*}
\begin{align*}
\left|\partial_R^k\left( \chi_{\nu\tau - R\lesssim\hbar}\right)\right|\lesssim \hbar^{-k}, 
\end{align*}
we can easily place the first term on the RHS of \eqref{smooth plus singular} into $\tilde{S}_0^{(n)}$, with norm bound $\lesssim \log \tau)^{N_1 - j}\cdot\tau^{-1-\nu}$, and the term is as $f_3$ in the lemma. For the second term in \eqref{smooth plus singular} we can write it as 
\begin{align*}
& \chi_{\nu\tau - R\lesssim\hbar}\int_0^\infty\chi_{\gtrsim 1}\left((\nu\tau - R)\xi^{\frac12}\right)\chi_{\lesssim 1}\left(\hbar^{-1}\rho(x;\alpha,\hbar)\right)\left(\mp i\hbar^{-1}\rho(x;\alpha,\hbar)\right)^l e^{\pm i\hbar^{-1}y(\alpha;\hbar)}\cdot \Psi(x;\alpha,\hbar)\,d\xi\\
& = \left(\nu\tau - R\right)^{\frac12+k\nu}\left(\log\left(\nu\tau-R\right)\right)^j\cdot \frac{H(\tau, \nu\tau-R)}{\tau^{\frac12}},
\end{align*}
where the term $H(\tau, x)$ has the symbol behavior asserted in the lemma. In fact, note that the restricting $\hbar^{-1}\rho(x;\alpha,\hbar)\lesssim 1$ implies $x\gtrsim \hbar^{-1}$ in light of Lemma~\ref{lem: Lemma 3.4 CDST}, and then we infer 
\begin{align*}
\big|\partial_R\left(\mp i\hbar^{-1}\rho(x;\alpha,\hbar)\right)^l\big|\lesssim l\cdot \big|\hbar^{-1}\rho_x\big|\cdot \xi^{\frac12}\hbar\lesssim \big|x\cdot\hbar^{-1}\rho_x\big|\cdot R^{-1}\lesssim R^{-1}, 
\end{align*}
again using Lemma~\ref{lem: Lemma 3.4 CDST}, and similarly for higher order derivatives. Hence this expression has better than symbol behavior with respect to differentiation with respect to $\nu\tau - R$. This is also the case for $ \Psi(x;\alpha,\hbar)$ in light of \eqref{eq:Psixalphahbardef}, while the function $\chi_{\gtrsim 1}\left((\nu\tau - R)\xi^{\frac12}\right)$ has symbol type bounds vis-a-vis $\nu\tau - R$.\\

Next for the term with oscillatory phase $e^{\pm i(\nu\tau - R)\xi^{\frac12}}$ in the regime $(\nu\tau - R)\xi^{\frac12}\gtrsim 1$, pass to the new integration variable $\tilde{\xi} = (\nu\tau - R)\xi^{\frac12}$, which also gives the structure 
\begin{align*}
	& \chi_{\nu\tau - R\lesssim\hbar}\int_0^\infty\chi_{\gtrsim1}\left((\nu\tau - R)\xi^{\frac12}\right)\chi_{\lesssim 1}\left(\hbar^{-1}\rho(x;\alpha,\hbar)\right) e^{\pm i(\nu\tau - R)\xi^{\frac12}}\cdot \left(\mp i\hbar^{-1}\rho(x;\alpha,\hbar)\right)^l e^{\pm i\hbar^{-1}y(\alpha;\hbar)}\cdot \Psi(x;\alpha,\hbar)\,d\xi\\
	& =  \chi_{\nu\tau - R\lesssim\hbar}\cdot(\nu\tau - R)^{\frac12 + k\nu}\cdot \frac{H(\tau, \nu\tau - R)}{\tau^{\frac12}}.
\end{align*}
 This also applies to the terms with $k\geq 1$ (In fact if $k\geq1$ is sufficiently large, we do not write $\chi_{(\nu\tau-R)\xi^{\frac12}\lesssim 1}$ as $1-\chi_{(\nu\tau-R)\xi^{\frac12}\gtrsim1}$. Instead, we directly integrate the expression in the variable $\tilde{\xi}$, because now the power in $\tilde{\xi}$ from the $k\geq 1$ expansion kills some powers in the denominator of the expression. Therefore we are able to integrate over the range $\tilde{\xi}\lesssim1$.) 
This concludes the case {\it{(4.i)}}. 
\\

{\it{(4.ii): $\hbar^{-1}\rho(x;\alpha,\hbar)\gg 1$, $(\nu\tau - R)\xi^{\frac12}\ll \hbar^{-1}\rho(x;\alpha,\hbar)$.}} We effect the additional restriction by means of the following smooth cutoff
\[
\chi_{\gg 1}\left(\frac{\hbar^{-1}\rho(x;\alpha,\hbar)}{(\nu\tau - R)\xi^{\frac12}}\right)\cdot \chi_{\gg 1}\left(\hbar^{-1}\rho(x;\alpha,\hbar)\right) =: \Lambda(x;\alpha,\hbar). 
\]
Observe that then thinking of $\frac{x}{\xi^{\frac12}\hbar} = \nu\tau - (\nu\tau - R)$ as a function of $\tau, \nu\tau - R$, 
\begin{align*}
\left|\frac{\partial \Lambda}{\partial\xi^{\frac12}}\right|\lesssim \xi^{-\frac12},\quad \left|\frac{\partial \Lambda}{\partial (\nu\tau - R)}\right|\lesssim (\nu\tau - R)^{-1},\quad \left|\frac{\partial \Lambda}{\partial\tau}\right|\lesssim\tau^{-1}. 
\end{align*}
We decompose 
\begin{align*}
&\chi_{\nu\tau - R\lesssim\hbar}\int_0^\infty  \Lambda(x;\alpha,\hbar)e^{\pm i(\nu\tau - R)\xi^{\frac12}}\cdot e^{\mp i\hbar^{-1}[\rho(x;\alpha,\hbar) - y(\alpha;\hbar)]}\cdot \Psi(x;\alpha,\hbar)\,d\xi\\
& = \chi_{\nu\tau - R\lesssim\hbar}\int_0^\infty\chi_{\lesssim 1}((\nu\tau - R)\xi^{\frac12})\Lambda(x;\alpha,\hbar)e^{\pm i(\nu\tau - R)\xi^{\frac12}}\cdot e^{\mp i\hbar^{-1}[\rho(x;\alpha,\hbar) - y(\alpha;\hbar)]}\cdot \Psi(x;\alpha,\hbar)\,d\xi\\
& +  \chi_{\nu\tau - R\lesssim\hbar}\int_0^\infty\chi_{\gtrsim 1}((\nu\tau - R)\xi^{\frac12})\Lambda(x;\alpha,\hbar)e^{\pm i(\nu\tau - R)\xi^{\frac12}}\cdot e^{\mp i\hbar^{-1}[\rho(x;\alpha,\hbar) - y(\alpha;\hbar)]}\cdot \Psi(x;\alpha,\hbar)\,d\xi
\end{align*}
For the first term on the right, we further decompose it into 
\begin{align*}
& \chi_{\nu\tau - R\lesssim\hbar}\int_0^\infty\chi_{\lesssim 1}((\nu\tau - R)\xi^{\frac12})\Lambda(x;\alpha,\hbar)e^{\pm i(\nu\tau - R)\xi^{\frac12}}\cdot e^{\mp i\hbar^{-1}[\rho(x;\alpha,\hbar) - y(\alpha;\hbar)]}\cdot \Psi(x;\alpha,\hbar)\,d\xi\\
& =  \chi_{\nu\tau - R\lesssim\hbar}\int_0^\infty\chi_{\lesssim 1}((\nu\tau - R)\xi^{\frac12})\Lambda(x;\alpha,\hbar)\cdot e^{\mp i\hbar^{-1}[\rho(x;\alpha,\hbar) - y(\alpha;\hbar)]}\cdot \Psi(x;\alpha,\hbar)\,d\xi\\
& +  \chi_{\nu\tau - R\lesssim\hbar}\int_0^\infty\chi_{\lesssim 1}((\nu\tau - R)\xi^{\frac12})\Lambda(x;\alpha,\hbar)[e^{\pm i(\nu\tau - R)\xi^{\frac12}}-1]\cdot e^{\mp i\hbar^{-1}[\rho(x;\alpha,\hbar) - y(\alpha;\hbar)]}\cdot \Psi(x;\alpha,\hbar)\,d\xi\\
& =  \chi_{\nu\tau - R\lesssim\hbar}\int_0^\infty\Lambda(x;\alpha,\hbar)\cdot e^{\mp i\hbar^{-1}[\rho(x;\alpha,\hbar) - y(\alpha;\hbar)]}\cdot \Psi(x;\alpha,\hbar)\,d\xi\\
& -  \chi_{\nu\tau - R\lesssim\hbar}\int_0^\infty\chi_{\gtrsim 1}((\nu\tau - R)\xi^{\frac12})\Lambda(x;\alpha,\hbar)\cdot e^{\mp i\hbar^{-1}[\rho(x;\alpha,\hbar) - y(\alpha;\hbar)]}\cdot \Psi(x;\alpha,\hbar)\,d\xi\\
& +  \chi_{\nu\tau - R\lesssim\hbar}\int_0^\infty\chi_{\lesssim 1}((\nu\tau - R)\xi^{\frac12})\Lambda(x;\alpha,\hbar)[e^{\pm i(\nu\tau - R)\xi^{\frac12}}-1]\cdot e^{\mp i\hbar^{-1}[\rho(x;\alpha,\hbar) - y(\alpha;\hbar)]}\cdot \Psi(x;\alpha,\hbar)\,d\xi
\end{align*}
Here the first of the last three terms can be reformulated as 
\begin{align*}
& \chi_{\nu\tau - R\lesssim\hbar}\int_0^\infty \chi_{\gg 1}\left(\hbar^{-1}\rho(x;\alpha,\hbar)\right)\cdot e^{\mp i\hbar^{-1}[\rho(x;\alpha,\hbar) - y(\alpha;\hbar)]}\cdot \Psi(x;\alpha,\hbar)\,d\xi\\
& - \chi_{\nu\tau - R\lesssim\hbar}\int_0^\infty \chi_{\lesssim 1}\left(\frac{\hbar^{-1}\rho(x;\alpha,\hbar)}{(\nu\tau - R)\xi^{\frac12}}\right)\chi_{\gg 1}\left(\hbar^{-1}\rho(x;\alpha,\hbar)\right)\cdot e^{\mp i\hbar^{-1}[\rho(x;\alpha,\hbar) - y(\alpha;\hbar)]}\cdot \Psi(x;\alpha,\hbar)\,d\xi,
\end{align*}
of which the first term is easily seen to be in $\tilde{S}_0^{(n)}$ and to also satisfy the bounds required for $f_3$ in the lemma. For the second term, since differentiation in $R$ would generate a factor of $(\nu\tau-R)^{-1}$, this term only admits the representation
\[
\left(\nu\tau - R\right)^{\frac12+k\nu}\log\left(\nu\tau-R\right)^j\cdot \frac{H(\tau, \nu\tau-R)}{\tau^{\frac12}}
\]
where $H$ has the desired symbol behavior and satisfies the bounds required of $G_{k,0,j}$. Here we note that interpreting $x = x(\tau,\nu\tau - R)\sim \hbar\tau\xi^{\frac12}$, we have the estimate 
\begin{align*}
\big|\partial_{\nu\tau - R}\big(e^{\mp i\hbar^{-1}[\rho(x;\alpha,\hbar)}\big)\big|\lesssim \hbar^{-1}\cdot \hbar\xi^{\frac12}\cdot \big|\rho_x\big|\sim \frac{\xi^{\frac12}}{(\hbar\tau\xi^{\frac12})^2}\lesssim (\nu\tau - R)^{-1}\cdot \frac{\nu\tau - R}{\hbar}\cdot \tau^{-1}\ll (\nu\tau-R)^{-1}\cdot\tau^{-1},
\end{align*}
and similar estimates apply to higher order derivatives. 
\\
The latter conclusion also applies to the terms 
\begin{align*}
&\chi_{\nu\tau - R\lesssim\hbar}\int_0^\infty\chi_{\gtrsim 1}((\nu\tau - R)\xi^{\frac12})\Lambda(x;\alpha,\hbar)\cdot e^{\mp i\hbar^{-1}[\rho(x;\alpha,\hbar) - y(\alpha;\hbar)]}\cdot \Psi(x;\alpha,\hbar)\,d\xi\\
&\chi_{\nu\tau - R\lesssim\hbar}\int_0^\infty\chi_{\lesssim 1}((\nu\tau - R)\xi^{\frac12})\Lambda(x;\alpha,\hbar)[e^{\pm i(\nu\tau - R)\xi^{\frac12}}-1]\cdot e^{\mp i\hbar^{-1}[\rho(x;\alpha,\hbar) - y(\alpha;\hbar)]}\cdot \Psi(x;\alpha,\hbar)\,d\xi
\end{align*}
In fact here we simply change the variable to $\tilde{\xi}:=(\nu\tau-R)^{2}\xi$.
\\

It remains to deal with the term 
\[
\chi_{\nu\tau - R\lesssim\hbar}\int_0^\infty\chi_{\gtrsim 1}\left((\nu\tau - R)\xi^{\frac12}\right)\Lambda(x;\alpha,\hbar)e^{\pm i(\nu\tau - R)\xi^{\frac12}}\cdot e^{\mp i\hbar^{-1}[\rho(x;\alpha,\hbar) - y(\alpha;\hbar)]}\cdot \Psi(x;\alpha,\hbar)\,d\xi.
\]
Introduce the phase function 
\[
\tilde{\xi}: = -(\nu\tau - R)\xi^{\frac12} + \hbar^{-1}\rho(x;\alpha,\hbar),
\]
whence 
\begin{align*}
\frac{\partial\tilde{\xi}}{\partial\xi^{\frac12}} = -(\nu\tau - R) + \hbar^{-1}\rho_x\cdot x\xi^{-\frac12} + \hbar^{-1}\rho_{\alpha}\cdot\hbar,
\end{align*}
and so since $\rho_x\cdot x\simeq -\rho$  on the support of the integrand and further $\hbar^{-1}\rho\gg 1$ there, as well as $\left|\rho_{\alpha}\right|\lesssim \hbar\alpha^{-1}\lesssim \xi^{-\frac12}$, we infer, using the fact $\left|\nu\tau-R\right|\ll \hbar^{-1}\rho\xi^{-\frac12}$,
\[
\left|\frac{\partial\tilde{\xi}}{\partial\xi^{\frac12}}\right|\simeq \hbar^{-1}\rho\xi^{-\frac12}
\]
on the support of the integrand. Then we write the preceding integral as 
\begin{align*}
\chi_{\nu\tau - R\lesssim\hbar}\int_0^\infty\chi_{\gtrsim 1}\left((\nu\tau - R)\xi^{\frac12}\right)\Lambda(x;\alpha,\hbar)e^{\mp i\tilde{\xi}}\cdot e^{\pm i\hbar^{-1}y(\alpha;\hbar)}\cdot \tilde{\Psi}(x;\alpha,\hbar)\,d\tilde{\xi}
\end{align*}
where we set 
\[
\tilde{\Psi}(x;\alpha,\hbar) = 2\xi^{\frac12}\frac{\partial\xi^{\frac12}}{\partial\tilde{\xi}}\cdot \Psi(x;\alpha,\hbar).
\]
Then we observe that 
\begin{align*}
&\left|\partial_{\tilde{\xi}}^l\tilde{\Psi}(x;\alpha,\hbar)\right|\lesssim \tau^{-\frac12}\cdot \xi^{-\frac14-\frac{k\nu}{2}}(\log\xi)^j\cdot \left(\hbar^{-1}\rho\right)^{-l-1},\quad \left|\partial_{\tilde{\xi}}^l\chi_{\gtrsim 1}\left((\nu\tau - R)\xi^{\frac12}\right)\right|\lesssim \left(\hbar^{-1}\rho\right)^{-l},\\&\left|\partial_{\tilde{\xi}}^l\Lambda(x;\alpha,\hbar)\right|\lesssim \left(\hbar^{-1}\rho\right)^{-l},\quad
\left|\partial_{\tilde{\xi}}^l\left(e^{\pm i\hbar^{-1}y(\alpha;\hbar)}\right)\right|\lesssim \left(\hbar^{-1}\rho\right)^{-l}.
\end{align*}
Note that applying a derivative $\partial_{\nu\tau - R}$ to the original expression (formulated as integral with respect to $\xi$, and interpreted as function of $\nu\tau - R$ and $\tau$) results at worst in a loss of 
\[
\lesssim \xi^{\frac12} + \hbar^{-1}\rho\cdot\tau^{-1}\lesssim \left(\hbar^{-1}\rho\right)\cdot \left((\nu\tau - R)^{-1} + \tau^{-1}\right),
\]
while applying $\partial_{\tau}$ results in at worst the loss $\hbar^{-1}\rho\cdot\tau^{-1}$. 
We conclude that after performing integration by parts in the variable $\tilde{\xi}$ and then reverting to the original integration variable $\xi$, we have 
\begin{align*}
&\chi_{\nu\tau - R\lesssim\hbar}\int_0^\infty\chi_{\gtrsim 1}\left((\nu\tau - R)\xi^{\frac12}\right)\Lambda(x;\alpha,\hbar)e^{\pm i(\nu\tau - R)\xi^{\frac12}}\cdot e^{\mp i\hbar^{-1}[\rho(x;\alpha,\hbar) - y(\alpha;\hbar)]}\cdot \Psi(x;\alpha,\hbar)\,d\xi\\
&= (\nu\tau - R)^{\frac12 + k\nu}\left(\log(\nu\tau-R)\right)^{j}\frac{H(\tau, \nu\tau - R)}{\tau^{\frac12}},
\end{align*}
where $H$ has the symbol behavior asserted in the lemma. (Here the factor $(\nu\tau-R)^{\frac12+k\nu}\cdot \left(\log(\nu\tau-R)\right)^{j}$ comes from the factor $\xi^{-\frac14-\frac{k\nu}{2}}(\log\xi)^{j}$ and the cutoff. The rapid decay in $\hbar^{-1}\rho$ produced by integration by parts gives the integrability.)
\\

{\it{(4.iii): $\hbar^{-1}\rho(x;\alpha,\hbar)\gg 1$, $(\nu\tau - R)\xi^{\frac12}\gtrsim\hbar^{-1}\rho(x;\alpha,\hbar)$.}}
Start by considering the term 
\begin{align*}
\chi_{\nu\tau - R\lesssim\hbar}\int_0^\infty \tilde{\Lambda}(x;\alpha,\hbar)\cdot e^{\pm i\hbar^{-1}y(\alpha;\hbar)}\cdot \Psi(x;\alpha,\hbar)\,d\xi,
\end{align*}
where we set 
\begin{align*}
\tilde{\Lambda}(x;\alpha,\hbar): = \chi_{\lesssim 1}\left(\frac{\hbar^{-1}\rho(x;\alpha,\hbar)}{(\nu\tau - R)\xi^{\frac12}}\right)\cdot \chi_{\gg 1}\left(\hbar^{-1}\rho(x;\alpha,\hbar)\right).
\end{align*}
This term can be seen to be of the explicit form
\begin{align*}
	\left(\nu\tau - R\right)^{\frac12+k\nu}\left(\log\left(\nu\tau-R\right)\right)^j\cdot \frac{H(\tau, \nu\tau-R)}{\tau^{\frac12}}
\end{align*}
 by arguing as for the term 
\[
\chi_{\nu\tau - R\lesssim\hbar}\int_0^\infty\Lambda(x;\alpha,\hbar)\cdot e^{\mp i\hbar^{-1}\left[\rho(x;\alpha,\hbar) - y(\alpha;\hbar)\right]}\cdot \Psi(x;\alpha,\hbar)\,d\xi.
\]
This leaves us with the term 
\[
\chi_{\nu\tau - R\lesssim\hbar}\int_0^\infty\tilde{\Lambda}(x;\alpha,\hbar) \left[e^{\pm i(\nu\tau - R)\xi^{\frac12}}\cdot e^{\mp i\hbar^{-1}\rho(x;\alpha,\hbar)}-1\right]e^{\pm i\hbar^{-1}y(\alpha;\hbar)}\cdot \Psi(x;\alpha,\hbar)\,d\xi.
\]
Introduce the variable 
\[
\tilde{\xi}: = (\nu\tau - R)\xi^{\frac12} - \hbar^{-1}\rho(x;\alpha,\hbar),
\]
and so (for suitable $c>0$)
\begin{align*}
\frac{\partial\tilde{\xi}}{\partial\xi^{\frac12}} \geq \nu\tau - R + c\hbar^{-1}x^{-1}\xi^{-\frac12}\geq \nu\tau - R,\quad \xi^{\frac12} = \left[\tilde{\xi} + \hbar^{-1}\rho(x;\alpha,\hbar)\right](\nu\tau - R)^{-1}.
\end{align*}
Also, the cutoff $\tilde{\Lambda}(x;\alpha,\hbar)$ ensures that $(\nu\tau - R)\xi^{\frac12} = \tilde{\xi} + \hbar^{-1}\rho(x;\alpha,\hbar)\gtrsim \left|\tilde{\xi}\right|$ on the support of the integrand. Then, assuming $k\nu<\frac12$, we can write (with $c(\tau, \nu\tau - R) = \hbar^{-1}(\nu\tau - R) - \hbar^{-1}\rho(R; 1,\hbar)$) 
\begin{align*}
&\chi_{\nu\tau - R\lesssim\hbar}\int_0^\infty \tilde{\Lambda}(x;\alpha,\hbar)\left[e^{\pm i(\nu\tau - R)\xi^{\frac12}}\cdot e^{\mp i\hbar^{-1}\rho(x;\alpha,\hbar)}-1\right]e^{\pm i\hbar^{-1}y(\alpha;\hbar)}\cdot \Psi(x;\alpha,\hbar)\,d\xi\\
& = \sum_{j'+j''= j}C_{j',j''}\chi_{\nu\tau - R\lesssim\hbar}\cdot(\nu\tau - R)^{\frac32+k\nu}\left[\log(\nu\tau - R)\right]^{j'}\cdot\\
&\int_{c(\tau, \nu\tau - R)}^\infty \frac{e^{\pm i\tilde{\xi}}-1}{[\tilde{\xi}+ \hbar^{-1}\rho(x;\alpha,\hbar)]^{\frac32+k\nu}}\left(\log\left(\tilde{\xi}+\hbar^{-1}\rho\right)\right)^{j''}\tilde{\Psi}(x;\alpha,\hbar)\,d\tilde{\xi}, 
\end{align*}
where we put 
\begin{align*}
\tilde{\Psi}(x;\alpha,\hbar) = \tilde{\Lambda}(x;\alpha,\hbar)e^{\pm i\hbar^{-1}y(\alpha;\hbar)}\cdot \Psi(x;\alpha,\hbar)\cdot \xi^{\frac54+\frac{k\nu}{2}}(\log\xi)^{-j}\cdot \frac{\partial\xi^{\frac12}}{\partial\tilde{\xi}}.
\end{align*}
Then from the definitions and in particular \eqref{eq:Psixalphahbardef} we directly infer the bound 
\begin{align*}
\left|\tilde{\Psi}(x;\alpha,\hbar)\right|\lesssim \tau^{-\frac12}\cdot(\nu\tau-R)^{-1}
\end{align*}
provided $\hbar\gtrsim\nu\tau-R>0$, and putting 
\[
G(\tau, \nu\tau-R) = \tau^{\frac12}(\nu\tau-R)\cdot\int_{c(\tau, \nu\tau - R)}^\infty \frac{e^{\pm i\tilde{\xi}}-1}{\left[\tilde{\xi}+ \hbar^{-1}\rho(x;\alpha,\hbar)\right]^{\frac32+k\nu}}\left(\log\left(\tilde{\xi}+\hbar^{-1}\rho\right)\right)^{j''}\tilde{\Psi}(x;\alpha,\hbar)\,d\tilde{\xi},
\]
we easily infer the desired symbol type behavior asserted in the lemma. In fact, note that upon writing $x = \hbar\xi^{\frac12}\left[\nu\tau - (\nu\tau - R)\right]$, 
\begin{align*}
\partial_{(\nu\tau - R)}\left(\frac{1}{\left[\tilde{\xi}+ \hbar^{-1}\rho(x;\alpha,\hbar)\right]^{\frac32+k\nu}}\right) = -\left(\frac32+k\nu\right)\cdot \frac{\hbar^{-1}\rho_x\cdot \hbar\xi^{\frac12}}{\left[\tilde{\xi}+ \hbar^{-1}\rho(x;\alpha,\hbar)\right]^{\frac52+k\nu}},
\end{align*}
and so (again using $\tau\simeq R$)
\begin{align*}
\left|\partial_{(\nu\tau - R)}\left(\frac{1}{\left[\tilde{\xi}+ \hbar^{-1}\rho(x;\alpha,\hbar)\right]^{\frac32+k\nu}}\right)\right|&\lesssim  \frac{\hbar^{-1}\rho_x\cdot x\tau^{-1}}{\left[\tilde{\xi}+ \hbar^{-1}\rho(x;\alpha,\hbar)\right]^{\frac52+k\nu}}\\
&\lesssim \frac{\tau^{-1}}{\left[\tilde{\xi}+ \hbar^{-1}\rho(x;\alpha,\hbar)\right]^{\frac32+k\nu}}
\end{align*}
on the support of the integrand, which leads to better than symbol behavior in the case the derivative $\partial_{(\nu\tau - R)}$ falls on this term in the integral. However, this is not the case for the contribution when this derivative falls on the factor 
\[
\rho_{n}(\xi) = \rho_{n}(\xi(\tilde{\xi},\tau, R;\hbar)),
\]
for which we have to use Proposition \ref{prop:DFT nlarge}. The case $k\nu\geq \frac12$ is handled by subtracting further terms from $e^{\pm i\tilde{\xi}}$ to handle the term $\left[\tilde{\xi}+ \hbar^{-1}\rho(x;\alpha,\hbar)\right]^{-(\frac32+k\nu)}$ when $\tilde{\xi}\ll 1$. 
\\

It remains to explain how to handle the contributions of the second to fourth sum constituting $\xb(\tau, \xi)$ in Definition~\ref{defi:xsingulartermsngeq2proto}. For the third sum of terms, one needs to invoke Lemma~\ref{lem:FtoGimprovement} in the oscillatory regime $(\nu\tau - R)\xi^{\frac12}\gtrsim 1$, with $x$ in the lemma representing $\nu\tau - R$. This results in the gain of regularity of the functions $G_{k,l,j}(\tau,\nu\tau - R)$ with respect to the second variable, compared to the functions $F_{l,k,j}^{\pm}$. 
Due to the restriction of the expressions to the high-frequency region $\xi\gtrsim \hbar^{-2}$, one may replace the factors $\langle\hbar^2\xi\rangle^{-\frac{l}{4}}$ by $(\hbar^2\xi)^{-\frac{l}{4}}$ for all intents and purposes, which upon following the steps in the preceding is seen to lead to additional factors of the form $\left(\hbar^{-1}[\nu\tau - R]\right)^{\frac{l-1}{2}}$ for the terms of explicit form (with symbol behavior with respect to $\nu\tau - R)$). Moreover, the extra factors $F_{l,k,j}(\tau, \xi)$ are seen to not lead to any additional complications due to their symbol type bounds. 
\\
It remains to deal with the contribution of terms in the second and fourth sum, which differ more significantly due to the presence of the phase function $e^{\pm i\hbar^{-1}\rho(x_{\tau};\alpha,\hbar)}$. This requires revisiting all terms where integration by parts was required. To do so we briefly revisit the cases {\it{(1)}} - {\it{(4)}} in the preceding, explaining the differences. For simplicity we shall set $l = j = 0$, the more general case being similar, and so we can unambiguously refer to these cases from before. 
\\
In cases {\it{(1), (2)}}, one replaces the term $e^{\pm i\nu\tau\xi^{\frac12}}$ by $e^{\pm i(\nu\tau\xi^{\frac12} + \hbar^{-1}\rho(x_{\tau};\alpha,\hbar))}$. Keeping in mind that $\xi\hbar^2\gtrsim 1$ for these terms, one has
\begin{align*}
\partial_{\xi^{\frac12}}\left(\nu\tau\xi^{\frac12} + \hbar^{-1}\rho(x_{\tau};\alpha,\hbar)\right) &= \nu\tau + \hbar^{-1}\rho_x(x_{\tau};\alpha,\hbar)\cdot \hbar\nu\tau + \hbar^{-1}\rho_{\alpha}(x_{\tau};\alpha,\hbar)\cdot\hbar\\
& = \nu\tau + O(1),
\end{align*}
since we have $x_{\tau}\gtrsim \tau$ on the support of the expression, and so $\hbar^{-1}\rho_x(x_{\tau};\alpha,\hbar)\cdot \hbar\nu\tau\lesssim \tau^{-1}$, while $\left|\rho_{\alpha}\right|\lesssim\hbar$. Thus the integration by parts with respect to $\xi^{\frac12}$ works in the same way here. 
\\
Case {\it{(3)}} also basically doesn't change, since here we encounter the phase function 
\[
\pm \frac23\hbar^{-1}\zeta^{\frac32} \pm (\nu\tau\xi^{\frac12} + \hbar^{-1}\rho(x_{\tau};\alpha,\hbar)), 
\]
Here we have the same bound 
\begin{align*}
\left|\partial_{\xi^{\frac12}}\left(\pm \frac23\hbar^{-1}\zeta^{\frac32} \pm \left(\nu\tau\xi^{\frac12} + \hbar^{-1}\rho(x_{\tau};\alpha,\hbar)\right)\right)\right|\gtrsim \tau, 
\end{align*}
on account of the above bound for $\partial_{\xi^{\frac12}}\left( \hbar^{-1}\rho(x_{\tau};\alpha,\hbar)\right) =  \hbar^{-1}\rho_x(x_{\tau};\alpha,\hbar)\cdot \hbar\nu\tau+\hbar^{-1}\rho_{\alpha}(x_{\tau};\alpha,\hbar)\cdot\hbar$, and then the argument proceeds as before. \\
In case {\it{(4)}}, we reduce as before to the situation $\nu\tau - R\lesssim 1$. Now to deal with the intermediate case $\hbar\ll\nu\tau - R\lesssim 1$, we encounter the phase 
\begin{align*}
\left[(\nu\tau - R)\xi^{\frac12} + \hbar^{-1}y(\alpha,\hbar)-  \hbar^{-1}\rho(x;\alpha,\hbar) + \hbar^{-1}\rho(x_{\tau};\alpha,\hbar)\right], 
\end{align*}
for which we have the relation 
\begin{align*}
&\partial_{\xi^{\frac12}}\left[(\nu\tau - R)\xi^{\frac12} + \hbar^{-1}y(\alpha,\hbar) -  \hbar^{-1}\rho(x;\alpha,\hbar) + \hbar^{-1}\rho(x_{\tau};\alpha,\hbar)\right]\\
& = (\nu\tau - R) + y_{\alpha}(\alpha,\hbar) - \hbar^{-1}\partial_{\xi^{\frac12}}\left(\int_0^1\rho_x( sx+(1-s)x_{\tau};\alpha,\hbar)\cdot(x-x_{\tau})\,ds\right),
\end{align*}
and we have the relation (for $\xi\gtrsim \hbar^{-2}$)
\begin{align*}
&\hbar^{-1}\partial_{\xi^{\frac12}}\left(\int_0^1\rho_x( sx+(1-s)x_{\tau};\alpha,\hbar)\cdot(x-x_{\tau})\,ds\right)\\
& = (R-\nu\tau)\cdot\xi^{\frac12}\int_0^1\left[\rho_{xx}( sx+(1-s)x_{\tau};\alpha,\hbar)\cdot (s\hbar R + (1-s)\hbar\nu\tau) + \hbar\rho_{x\alpha}( sx+(1-s)x_{\tau};\alpha,\hbar)\right]\,ds\\
&\hspace{2cm} + (R-\nu\tau)\cdot\int_0^1\rho_x( sx+(1-s)x_{\tau};\alpha,\hbar)\,ds\\
& = (R-\nu\tau)\cdot O\left(\tau^{-2}\right). 
\end{align*}
In light of $ y_{\alpha}(\alpha,\hbar) = O(\hbar)$ because of Lemma \ref{lem: Lemma 3.4 CDST}, we conclude that 
\begin{align*}
\left|\partial_{\xi^{\frac12}}\left[(\nu\tau - R)\xi^{\frac12} + \hbar^{-1}y(\alpha,\hbar) -  \hbar^{-1}\rho(x;\alpha,\hbar) + \hbar^{-1}\rho(x_{\tau};\alpha,\hbar)\right]\right|\gtrsim \nu\tau - R, 
\end{align*}
and one can again argue as in the earlier situation to handle the regime $\hbar\ll\nu\tau - R\lesssim 1$. 
\\
Finally, consider the expression 
\begin{align*}
\chi_{\nu\tau - R\lesssim \hbar}\int_0^\infty e^{\pm i(\nu\tau - R)\xi^{\frac12}}\cdot e^{\mp i\hbar^{-1}[\rho(x;\alpha,\hbar) - \rho(x_{\tau};\alpha,\hbar) - y(\alpha;\hbar)]}\cdot \chi_{\xi\hbar^2\gtrsim1}\Psi(x;\alpha,\hbar)\,d\xi
\end{align*}
where $\Psi(x;\alpha,\hbar)$ is again defined as in the earlier case {\it{(4)}}. Dealing with this is now in fact a bit simpler than in the earlier situation. Write the preceding as 
\begin{align*}
&\chi_{\nu\tau - R\lesssim \hbar}\int_0^\infty e^{\pm i(\nu\tau - R)\xi^{\frac12}}\cdot e^{\mp i\hbar^{-1}\left[\rho(x;\alpha,\hbar) - \rho(x_{\tau};\alpha,\hbar) - y(\alpha;\hbar)\right]}\cdot  \chi_{\xi\hbar^2\gtrsim1}\Psi(x;\alpha,\hbar)\,d\xi\\
& = \chi_{\nu\tau - R\lesssim \hbar}\int_0^\infty\chi_{(\nu\tau-R)\xi^{\frac12}\lesssim 1} e^{\pm i(\nu\tau - R)\xi^{\frac12}}\cdot e^{\mp i\hbar^{-1}\left[\rho(x;\alpha,\hbar) - \rho(x_{\tau};\alpha,\hbar)\right]} e^{\pm i\hbar^{-1}y(\alpha;\hbar)}\cdot  \chi_{\xi\hbar^2\gtrsim1}\Psi(x;\alpha,\hbar)\,d\xi\\
& +  \chi_{\nu\tau - R\lesssim \hbar}\int_0^\infty\chi_{(\nu\tau-R)\xi^{\frac12}\gtrsim 1} e^{\pm i(\nu\tau - R)\xi^{\frac12}}\cdot e^{\mp i\hbar^{-1}\left[\rho(x;\alpha,\hbar) - \rho(x_{\tau};\alpha,\hbar)\right]} e^{\pm i\hbar^{-1}y(\alpha;\hbar)}\cdot  \chi_{\xi\hbar^2\gtrsim1}\Psi(x;\alpha,\hbar)\,d\xi\\
& =: X_1 + X_2.
\end{align*}
In order to treat the first term, write it as (assuming $\frac12 + k\nu <1$, the case $\frac12 + k\nu \geq 1$ being handled analogously by subtracting more terms of the Taylor expansion of $e^z,\,z = \pm i(\nu\tau - R)\xi^{\frac12}\mp i\hbar^{-1}\left[\rho(x;\alpha,\hbar) - \rho(x_{\tau};\alpha,\hbar)\right]$)
\begin{align*}
 X_1 = &\chi_{\nu\tau - R\lesssim \hbar}\int_0^\infty\chi_{(\nu\tau-R)\xi^{\frac12}\lesssim 1} \left[e^{\pm i(\nu\tau - R)\xi^{\frac12}}\cdot e^{\mp i\hbar^{-1}[\rho(x;\alpha,\hbar) - \rho(x_{\tau};\alpha,\hbar)]} - 1\right]\tilde{\Psi}(x;\alpha,\hbar)\,d\xi\\
& + \chi_{\nu\tau - R\lesssim \hbar}\int_0^\infty\chi_{(\nu\tau-R)\xi^{\frac12}\lesssim 1}\tilde{\Psi}(x;\alpha,\hbar)\,d\xi,
\end{align*}
and we set $\tilde{\Psi}(x;\alpha,\hbar) = e^{\pm i\hbar^{-1}y(\alpha;\hbar)}\cdot  \chi_{\xi\hbar^2\gtrsim1}\Psi(x;\alpha,\hbar)$. The second term on the right, call it $X_{12}$, is easily seen to be handled the same way as for \eqref{smooth plus singular} 
The first term, call it $X_{11}$, can be written as 
\begin{align*}
X_{11} =& \chi_{\nu\tau - R\lesssim \hbar}(\nu\tau - R)\int_0^\infty\chi_{(\nu\tau-R)\xi^{\frac12}\lesssim 1} \frac{\left[e^{\pm i(\nu\tau - R)\xi^{\frac12}}\cdot e^{\mp i\hbar^{-1}\left[\rho(x;\alpha,\hbar) - \rho(x_{\tau};\alpha,\hbar)\right]} - 1\right]}{\left[\xi^{\frac12}(\nu\tau - R)\right]}\xi^{\frac12}\tilde{\Psi}(x;\alpha,\hbar)\,d\xi\\
 =:& \chi_{\nu\tau - R\lesssim \hbar}(\nu\tau - R)^{\frac12 + k\nu}\cdot\frac{H(\tau, \nu\tau - R)}{\tau^{\frac12}}.
\end{align*}
Here we integrate with respect to the variable $\tilde{\xi}:=\left(\nu\tau-R\right)^{2}\xi$.
We claim that the function $H(\tau, x)$ satisfies the symbol bounds $\left|\partial_{\tau}^{k_1}\partial_x^{k_2}H(\tau, x)\right|\lesssim \tau^{-k_1}x^{-k_2}$. To begin with, note that the function 
\[
\frac{\left[e^{\pm i(\nu\tau - R)\xi^{\frac12}}\cdot e^{\mp i\hbar^{-1}[\rho(x;\alpha,\hbar) - \rho(x_{\tau};\alpha,\hbar)]} - 1\right]}{\left[\xi^{\frac12}(\nu\tau - R)\right]}
\]
is bounded on the support of the integrand, since 
\[
\left|\pm i(\nu\tau - R)\xi^{\frac12}\mp i\hbar^{-1}\left[\rho(x;\alpha,\hbar) - \rho(x_{\tau};\alpha,\hbar)\right]\right|\lesssim (\nu\tau - R)\xi^{\frac12}
\]
there. As for the derivatives, write 
\begin{align*}
\hbar^{-1}\left[\rho(x;\alpha,\hbar) - \rho(x_{\tau};\alpha,\hbar)\right] = \int_0^1\rho_x\left(sx + (1-s)x_{\tau};\alpha,\hbar\right)\,ds\cdot (R-\nu\tau)\xi^{\frac12}, 
\end{align*}
whence (as usual interpreting $x$ as function of $\tau$ and $\nu\tau - R$), we infer the relation
\begin{align*}
&\left|\partial_{\tau}\left[ \frac{\left[e^{\pm i(\nu\tau - R)\xi^{\frac12}}\cdot e^{\mp i\hbar^{-1}[\rho(x;\alpha,\hbar) - \rho(x_{\tau};\alpha,\hbar)]} - 1\right]}{[\xi^{\frac12}(\nu\tau - R)]}\right]\right|\\
& \leq \left|1+\int_0^1\rho_x\left(sx + (1-s)x_{\tau};\alpha,\hbar\right)\,ds\right|\cdot \left|Z'\left((\nu\tau - R)\xi^{\frac12}\cdot(1+\int_0^1\rho_x\big(sx + (1-s)x_{\tau};\alpha,\hbar\big)\,ds)\right)\right|\\
&\hspace{7cm}\cdot \left|(\nu\tau - R)\xi^{\frac12}\cdot \partial_{\tau}\left(\int_0^1\rho_x\left(sx + (1-s)x_{\tau};\alpha,\hbar\right)\,ds\right)\right|\\
& + \left|\partial_{\tau}\left(\int_0^1\rho_x\left(sx + (1-s)x_{\tau};\alpha,\hbar\right)\,ds\right)\right|\cdot \left|Z\left((\nu\tau - R)\xi^{\frac12}\cdot(1+\int_0^1\rho_x\left(sx + (1-s)x_{\tau};\alpha,\hbar\right)\,ds)\right)\right|,
\end{align*}
where we set $Z(x): = \frac{e^{\pm ix} - 1}{x}$. Using the crude bound 
\[
\left|\partial_{\tau}\left(\int_0^1\rho_x\left(sx + (1-s)x_{\tau};\alpha,\hbar\right)\,ds\right)\right|\lesssim \tau^{-1}, 
\]
we easily infer, also using the definition of $\tilde{\Psi}(x;\alpha,\hbar)$, that on the support of the integrand for the integral in the definition of $X_{11}$, whence in particular with $(\nu\tau - R)\xi^{\frac12}\lesssim 1$, $x\gtrsim 1$, we have 
\[
\left|\partial_{\tau}H(\tau, \nu\tau - R)\right|\lesssim \tau^{-1}, 
\]
and analogously for higher derivatives with respect to $\tau$. Derivatives with respect to $\nu\tau - R$ are handled analogously.
\\
It remains to deal with the term $X_2$, which is 
\[
X_2 =  \chi_{\nu\tau - R\lesssim \hbar}\int_0^\infty\chi_{(\nu\tau-R)\xi^{\frac12}\gtrsim 1} e^{\pm i(\nu\tau - R)\xi^{\frac12}}\cdot e^{\mp i\hbar^{-1}\left[\rho(x;\alpha,\hbar) - \rho(x_{\tau};\alpha,\hbar)\right]} e^{\pm i\hbar^{-1}y(\alpha;\hbar)}\cdot  \chi_{\xi\hbar^2\gtrsim1}\Psi(x;\alpha,\hbar)\,d\xi.
\]
We claim that this is also of the form $\chi_{\nu\tau - R\lesssim \hbar}(\nu\tau - R)^{\frac12 + k\nu}\cdot\frac{H(\tau, \nu\tau - R)}{\tau^{\frac12}}$, with $H$ having the appropriate symbol behavior. For this, introduce the integration variable 
\[
\tilde{\xi}: = (\nu\tau - R)\xi^{\frac12} - \hbar^{-1}[\rho(x;\alpha,\hbar) - \rho(x_{\tau};\alpha,\hbar)].
\]
Then the preceding formula for the difference $\hbar^{-1}[\rho(x;\alpha,\hbar) - \rho(x_{\tau};\alpha,\hbar)]$ and the fact that $sx + (1-s)x_{\tau} \simeq \tau$ on the support of the integrand for $X_2$ imply that 
\[
\tilde{\xi}\simeq (\nu\tau - R)\xi^{\frac12},\quad \frac{\partial\tilde{\xi}}{\partial \xi^{\frac12}}\simeq \nu\tau - R
\]
on the support of the integrand. We can then write 
\begin{align*}
X_2 =  \chi_{\nu\tau - R\lesssim \hbar}(\nu\tau - R)^{\frac32+k\nu}\int_0^\infty\chi_{(\nu\tau-R)\xi^{\frac12}\gtrsim 1} \frac{e^{\mp i\tilde{\xi}}}{\left(\tilde{\xi}+\hbar^{-1}\left[\rho(x;\alpha,\hbar) - \rho(x_{\tau};\alpha,\hbar)\right]\right)^{\frac32+k\nu}}\cdot \frac{\partial\xi^{\frac12}}{\partial\tilde{\xi}}\tilde{\tilde{\Psi}}\,d\tilde{\xi}, 
\end{align*}
where we have set $\tilde{\tilde{\Psi}}: = 2\xi^{\frac54+\frac{k\nu}{2}}e^{\pm i\hbar^{-1}y(\alpha;\hbar)}\cdot  \chi_{\xi\hbar^2\gtrsim1}\Psi(x;\alpha,\hbar)$, and we recall that for simplicity of exposition we assume $l = j = 0$ (recall the definition of $\Psi$). Then since $\left|\tilde{\tilde{\Psi}}(x;\alpha,\hbar)\right|\lesssim \tau^{-\frac12}$ on the support of the integrand, writing $X_2 = \chi_{\nu\tau - R\lesssim \hbar}(\nu\tau - R)^{\frac12 + k\nu}\cdot\frac{H(\tau, \nu\tau - R)}{\tau^{\frac12}}$ we infer $\left|H(\tau, \nu\tau - R)\right|\lesssim 1$. Further using 
\begin{align*}
\frac{\partial\tilde{\xi}}{\partial\xi^{\frac12}} = &(\nu\tau - R)\cdot \left(1 + \int_0^1\rho_x(sx+(1-s)x_{\tau};\alpha,\hbar)\,ds + \int_0^1\rho_{xx}(sx+(1-s)x_{\tau};\alpha,\hbar)\cdot (sx+(1-s)x_{\tau})\,ds\right)\\
& + (\nu\tau - R)\cdot\hbar\int_0^1\rho_{x\alpha}(sx+(1-s)x_{\tau};\alpha,\hbar)\,ds, 
\end{align*}
taking advantage of Lemma \ref{lem: Lemma 3.4 CDST}, we find 
\[
\left|\partial_{\tau}^{k_1}\partial_{\nu\tau - R}^{k_2}\left(\frac{\partial\xi^{\frac12}}{\partial\tilde{\xi}}\right)\right|\lesssim \tau^{-k_1}(\nu\tau - R)^{-k_2-1}, 
\]
and in conjunction with the definition of $\tilde{\tilde{\Psi}}$, the desired symbol behavior of $H$ follows. 
\end{proof}
The preceding lemma admits a counterpart translating a physical singular expansion into a Fourier expansion. Note that this is not a simple inverse, and adapted to the specific needs later on:
\begin{lemma}\label{lem:singPhysicaltoFourierngeq2} Assume that $f = f(\tau, R)$ is a function supported on $\frac{\nu\tau}{2}<R\leq\nu\tau$ which admits an expansion 
	\[
	f = \sum_{j=1}^3 f_j
	\]
	where the $f_1$ is as in the statement of the preceding lemma, $f_2$ can be written as $f_2 = \sum_{l=0}^7f_{2l}$ while the $f_{2l}$ vanish on the set $\nu\tau<R$ and are defined on $R<\nu\tau$ as in the preceding lemma\footnote{In particular, we assume that the coefficient functions $G_{k,l,j}(\tau, x)$ admit derivatives up to order $20 - l + \lfloor \frac{l}{2}+k\nu\rfloor$ with respect to $x$.} except that $\hbar^{-\frac{l+1}{2}}$ is replaced by $\hbar^{-\frac{l}{2}-1}$, and $f_3$ is $C^\infty$ with $\|f_3\|_{\tilde{S}_1^{(\hbar)}}\lesssim (\log\tau)^{N_1}\cdot \tau^{-1-\nu}$ and satisfies the pointwise bounds 
	\[
	\left|\partial_R^kf_3(\tau, R)\right|\lesssim (\log\tau)^{N_1}\cdot\tau^{-\frac32-\nu}\cdot \hbar^{-1}\min\{\left|\nu\tau - R\right|^{-k},\hbar^{-k}\},\quad 0\leq k\leq 5. 
	\]
	Then setting 
	\[
	\xb(\tau, \xi): = \int_0^\infty \phi_{n}(R;\xi)\cdot f(\tau, R)\cdot \rho_{n}(\xi)\,d\xi, 
	\]
	we have $\xb = \xb_1 + \xb_2$, where 
	\[
	\left\|\xb_1\right\|_{\Sh_{1}}\lesssim (\log\tau)^{N_1}\cdot \tau^{-1-\nu},
	\]
	and furthermore we can write 
	\[
	\xb_1(\tau, \xi) = \sum_{\pm} e^{\pm i\nu\tau\xi^{\frac12}}\cdot g_{\pm}(\tau, \xi)
	\]
	where $g_{\pm}(\tau, \xi)$ is $C^\infty$ with respect to the second variable, and satisfies the bounds 
	\[
	\left|\partial_{\xi}^{k_2}g_{\pm}(\tau, \xi)\right|\lesssim(\log\tau)^{N_1}\cdot \tau^{-1-\nu}\cdot\hbar^{-1}\xi^{-\frac12}\langle\hbar\xi^{\frac12}\rangle^{-4},\quad 0\leq k_2\leq 4. 
	\]
	Furthermore, $\xb_2$ admits the symbol expansion (with $\rho$ as in Lemma \ref{lem: Lemma 3.4 CDST} and $x_{\tau} = \hbar\xi^{\frac12}\cdot\nu\tau$)
	\[
	\xb_2(\tau, \xi) = \sum_{\pm}\sum_{l =0}^7\sum_{k=1}^{N}\sum_{i=0}^{N_1}\chi_{\xi\geq \hbar^{-2}}\hbar^{-1}\langle\hbar^2\xi\rangle^{-\frac{l}{4}}\frac{e^{\pm i\left(\nu\tau\xi^{\frac12}+\hbar^{-1}\rho(x_{\tau};\alpha,\hbar)\right)}}{\xi^{\frac34+k\frac{\nu}{2}}}\left(\log\xi\right)^j\cdot F_{l,k,j}(\tau, \xi)
	\]
	and we have the symbol type bounds (with $k_1\in \{0,1\}, k_2\in\{0,\ldots, 20 - l\}$)
	\[
	\hbar^{\delta_1\cdot k_2}\cdot \left|\partial_{\tau}^{k_1}\partial_{\xi}^{k_2}F_{l,k,j}(\tau, \xi)\right|\lesssim (\log\tau)^{N_1 - j}\cdot\hbar^{-k_1}\tau^{-1-\nu-k_1}\cdot \xi^{-k_2}, 
	\]
	as well as the closure bound (with $k_*$ defined as in the statement of the preceding lemma)
	\[
	\hbar^{k_*}\left\|\xi^{20 - l+\delta}\partial_{\tau}^{k_1}\partial_{\xi}^{20 - l}F_{l,k,j}(\tau, \xi)\right\|_{\dot{C}^{\delta}}\lesssim (\log\tau)^{N_1 - j}\cdot\tau^{-1-\nu-k_1}.
	\]
\end{lemma}
\begin{proof}
	\[
	\xb(\tau,\xi) = \sum_{j=1}^3 \yb_j(\tau, \xi), 
	\]
	where we set 
	\[
	\yb_j(\tau, \xi) = \int_0^\infty \phi_{n}(R;\xi)f_j(\tau,R)\,R\,dR. 
	\]
	Then it is immediate that $\left\|\yb_3\right\|_{\Sh_{1}}\lesssim (\log\tau)^{N_1}\tau^{-1-\nu}$, whence $\yb_3$ can be incorporated into $\xb_1$. The same holds for $f_{1}$ by using a similar argument and the cutoff $\nu\tau-R\gtrsim 1$.
	\\
	It thus remains to deal with $f_2 = \sum_{l=0}^7 f_{2l}$, where it suffices to set (with $l,k,i$ as in the specified ranges above)
	\[
	g(\tau,R): = \chi_{|\nu\tau-R|\lesssim \hbar}\frac{G_{k,l,i}(\tau,\nu\tau - R)}{\tau^{\frac12}}\hbar^{-\frac{l}{2}-1}[\nu\tau - R]^{\frac{l}{2}+k\nu}(\log(\nu\tau - R))^i
	\]
	We split the distorted Fourier transform of this term into several pieces. To begin with, write
	\begin{align*}
	&\left\langle \phi_{n}(R;\xi),\,\chi_{|\nu\tau-R|\lesssim \hbar}\frac{G_{k,l,j}(\tau,\nu\tau - R)}{\tau^{\frac12}}\hbar^{-\frac{l}{2}-1}\left[\nu\tau - R\right]^{\frac{l}{2}+k\nu}\left(\log(\nu\tau - R)\right)^j\right\rangle_{L^2_{R\,dR}}\\
	& = \chi_{\xi\hbar^2\lesssim 1}\left\langle\ldots\right\rangle_{L^2_{R\,dR}} +  \chi_{\xi\hbar^2\gtrsim 1}\left\langle\ldots\right\rangle_{L^2_{R\,dR}}.
	\end{align*}
	Then the low-frequency term can again be included into $\Sh_{1}$: 
	\begin{align*}
	\left\|\chi_{\xi\hbar^2\lesssim 1}\langle\ldots\rangle_{L^2_{R\,dR}}\right\|_{\Sh_{1}}\lesssim \hbar\cdot \left\|g\right\|_{L^2_{R\,dR}}\lesssim (\log\tau)^{N_1}\cdot\tau^{-1-\nu},
	\end{align*}
	where we take advantage of the fact that 
	\[
	\left|\chi_{|\nu\tau-R|\lesssim \hbar}\hbar^{-\frac{l}{2}}\cdot |\nu\tau - R|^{\frac{l}{2}}\right|\lesssim 1. 
	\]
	It is also straightforward to check that $\chi_{\xi\hbar^2\lesssim 1}\langle\ldots\rangle_{L^2_{R\,dR}}$ can be written as $\sum_{\pm}e^{\pm i\nu\tau\xi^{\frac12}}\cdot g_{\pm}(\tau,\xi)$ where the functions $g_{\pm}$ satisfy the assertions of the lemma. 
	We can thus reduce to considering 
	\begin{align*}
	\xb_2(\tau,\xi): = \chi_{\xi\hbar^2\gtrsim 1}\int_0^\infty  \phi_{n}(R;\xi)\chi_{|\nu\tau-R|\lesssim \hbar}\frac{G_{k,l,j}(\tau,\nu\tau - R)}{\tau^{\frac12}}\hbar^{-\frac{l}{2}-1}\left[\nu\tau - R\right]^{\frac{l}{2}+k\nu}\left(\log(\nu\tau - R)\right)^j R\,dR
	\end{align*}
	Observe that since $\xi^{\frac12}\hbar R\simeq \xi^{\frac12}\hbar\tau\gg 1$ on the support of the above integrand, the function $ \phi_{n}(R;\xi)$ is in the oscillatory regime, and we may replace this function by 
	\[
	\hbar^{\frac13}x^{-\frac12}\left(\hbar^{-\frac23}\zeta\right)^{-\frac14}q^{-\frac14}(\zeta)\Re\left[a(\xi)e^{\pm i\left(R\xi^{\frac12} - y(\alpha;\hbar) + \hbar^{-1}\rho(x;\alpha,\hbar)\right)}\cdot\left(1+\hbar \overline{\tilde{a}_1(\zeta;\alpha)}\right)\right],
	\]
	where as usual we have incorporated the asymptotic correction terms for the oscillatory Airy functions into the term $1+\hbar \overline{\tilde{a}_1(\zeta;\alpha)}$, and we also recall Proposition~\ref{prop:DFT nlarge}. Write 
	\begin{align*}
	\phi_{n}(R;\xi) &= \sum_{\pm} e^{\pm i\nu\tau\xi^{\frac12}}\cdot \frac{e^{\mp i(\nu\tau - R)\xi^{\frac12}\pm i\hbar^{-1}\rho(x;\alpha,\hbar)}}{(R\xi^{\frac12})^{\frac12}}\cdot \psi(x;\alpha,\hbar)\cdot(1+\hbar \overline{\tilde{a}_1(\zeta;\alpha)})\\
	& =: \sum_{\pm} e^{\pm i\nu\tau\xi^{\frac12}}\cdot \phi_{\pm,\tau}(R;\xi,\hbar), 
	\end{align*}
	where $\psi$ has symbol behavior with respect to $x,\alpha$. Then we can write 
	\begin{align*}
	&\xb_2(\tau,\xi)\\& = \sum_{\pm}  \chi_{\xi\hbar^2\gtrsim 1}e^{\pm i\nu\tau\xi^{\frac12}}\cdot \int_0^\infty  \phi_{\pm,\tau}(R;\xi,\hbar)\chi_{|\nu\tau-R|\lesssim \hbar}\frac{G_{k,l,j}(\tau,\nu\tau - R)}{\tau^{\frac12}}\hbar^{-\frac{l}{2}-1}\left[\nu\tau - R\right]^{\frac{l}{2}+k\nu}\left(\log\left(\nu\tau - R\right)\right)^j R\,dR\\
	& = \sum_{\pm}  \chi_{\xi\hbar^2\gtrsim 1}e^{\pm i\nu\tau\xi^{\frac12}}\cdot \int_0^\infty  \phi_{\pm,\tau}(R;\xi,\hbar)\chi_{|\nu\tau-R|\lesssim \xi^{-\frac12}}\frac{G_{k,l,j}(\tau,\nu\tau - R)}{\tau^{\frac12}}\hbar^{-\frac{l}{2}-1}\left[\nu\tau - R\right]^{\frac{l}{2}+k\nu}\left(\log\left(\nu\tau - R\right)\right)^j R\,dR\\
	& +  \sum_{\pm}  \chi_{\xi\hbar^2\gtrsim 1}e^{\pm i\nu\tau\xi^{\frac12}}\cdot \int_0^\infty  \phi_{\pm,\tau}(R;\xi,\hbar)\chi_{ \xi^{-\frac12}\lesssim|\nu\tau-R|\lesssim\hbar}\frac{G_{k,l,j}(\tau,\nu\tau - R)}{\tau^{\frac12}}\hbar^{-\frac{l}{2}-1}\left[\nu\tau - R\right]^{\frac{l}{2}+k\nu}\left(\log\left(\nu\tau - R\right)\right)^j R\,dR\\
	&=:\xb_{21} + \xb_{22}. 
	\end{align*}
	To deal with $\xb_{21}$, we simply write it as 
	\[
	\xb_{21} =  \sum_{\pm}  \chi_{\xi\hbar^2\gtrsim 1}\hbar^{-1}\left\langle \xi\hbar^2\right\rangle^{-\frac{l}{4}}\frac{e^{\pm i\nu\tau\xi^{\frac12}}}{\xi^{\frac34+\frac{k\nu}{2}}}\left(\log\xi\right)^j\cdot F^{(1),\pm}_{l,k,j}(\tau,\xi), 
	\]
	where we define 
	\begin{align*}
	&F^{(1),\pm}_{l,k,j}(\tau,\xi)\\ =& \tilde{\chi}_{\xi\hbar^2\gtrsim 1}\frac{\left\langle \xi\hbar^2\right\rangle^{\frac{l}{4}}}{\left(\xi\hbar^2\right)^{\frac{l}{4}}}\cdot \frac{\xi^{\frac34}}{\left(\log\xi\right)^j}\\
	&\cdot\int_0^\infty  \phi_{\pm,\tau}(R;\xi,\hbar)\chi_{|\nu\tau-R|\lesssim \xi^{-\frac12}}\frac{G_{k,l,j}(\tau,\nu\tau - R)}{\tau^{\frac12}}\left[\xi^{\frac12}(\nu\tau - R)\right]^{\frac{l}{2}+k\nu}\left(\log\left(\nu\tau - R\right)\right)^j R\,dR
	\end{align*}
	 Then it is straightforward to check that 
	\[
	\left|F^{(1),\pm}_{l,k,j}(\tau,\xi)\right|\lesssim \big(\log\tau\big)^{N_1-j}\cdot \tau^{-1-\nu},
	\]
	due to the fact that $\left| \phi_{\pm,\tau}(R;\xi,\hbar)\right|\lesssim \tau^{-\frac12}\xi^{-\frac14}$ on the support of the integrand, as well as the fact that the $R$-integration interval is of length $\lesssim \xi^{-\frac12}$. We claim that after multiplication with $e^{\mp i\hbar^{-1}\rho(x_{\tau};\alpha,\hbar)}$, the function $F^{(1),\pm}_{l,k,i}(\tau,\xi)$ also has the kind of symbol behavior asserted in the lemma. The main difficulty here comes from dealing with the phase $e^{\pm i\hbar^{-1}\rho(x;\alpha,\hbar)}$. The idea is to exploit the localization near the light cone expressed by the condition $(\nu\tau - R)\xi^{\frac12}\lesssim 1$ to expand this phase into Taylor expansion around $R = \nu\tau$. In fact, setting $x_{\tau}: = \xi^{\frac12}\hbar\nu\tau$, we can write\footnote{In fact, this holds for $\tau\gg 1$ sufficiently large under the condition $|x-x_\tau|\lesssim 1$, thanks to the analyticity properties of $\rho$ implied by Lemma~\ref{lem: Lemma 3.4 CDST} and its proof.}
	\begin{equation}\label{eq:goodTaylorBadrho}
	\hbar^{-1}\rho(x;\alpha,\hbar) = \hbar^{-1}\rho(x_{\tau};\alpha,\hbar) + \hbar^{-1}\sum_{j\geq 1}\frac{1}{j!}(\partial_x^j\rho)(x_{\tau};\alpha,\hbar)\cdot (x-x_\tau)^j
	\end{equation}
	 Write 
	\begin{align*}
	\xb_{21} =  \sum_{\pm}  \chi_{\xi\hbar^2\gtrsim 1}\hbar^{-1}\left\langle \xi\hbar^2\right\rangle^{-\frac{l}{4}}\frac{e^{\pm i\left[\nu\tau\xi^{\frac12}+\hbar^{-1}\rho(x_{\tau};\alpha,\hbar)\right]}}{\xi^{\frac34+\frac{k\nu}{2}}}\left(\log\xi\right)^j\cdot \tilde{F}^{(1),\pm}_{l,k,j}(\tau,\xi), 
	\end{align*}
	where we have 
	\begin{align*}
	&\tilde{F}^{(1),\pm}_{l,k,j}(\tau,\xi)\\& = \tilde{\chi}_{\xi\hbar^2\gtrsim 1}\frac{\left\langle \xi\hbar^2\right\rangle^{\frac{l}{4}}}{\left(\xi\hbar^2\right)^{\frac{l}{4}}}\cdot \frac{\xi^{\frac34}}{\left(\log\xi\right)^j}\int_0^\infty  \tilde{\phi}_{\pm,\tau}(R;\xi,\hbar)\chi_{\left|\nu\tau-R\right|\lesssim \xi^{-\frac12}}\frac{G_{k,l,j}(\tau,\nu\tau - R)}{\tau^{\frac12}}\left[\xi^{\frac12}(\nu\tau - R)\right]^{\frac{l}{2}+k\nu}\left(\log\left(\nu\tau - R\right)\right)^j R\,dR,\\
	&\tilde{\phi}_{\pm,\tau}(R;\xi,\hbar): = \frac{e^{\mp i(\nu\tau - R)\xi^{\frac12}\pm i\hbar^{-1}\left[\rho(x;\alpha,\hbar) - \rho(x_{\tau};\alpha,\hbar)\right]}}{\left(R\xi^{\frac12}\right)^{\frac12}}\cdot \psi(x;\alpha,\hbar)\cdot\left(1+\hbar \overline{\tilde{a}_1(\zeta;\alpha)}\right)
	\end{align*}
	But then taking advantage of the localization of the integral as well as \eqref{eq:goodTaylorBadrho} it is straightforward to check that $\tilde{F}^{(1),\pm}_{l,k,i}(\tau,\xi)$ has the desired symbol type behavior. 
	\\
	
	It remains to deal with the term $\xb_{22}(\tau, \xi)$. Proceeding in analogy to $\xb_{21}$, write this as 
	\begin{align*}
	\xb_{22} =  \sum_{\pm}  \chi_{\xi\hbar^2\gtrsim 1}\hbar^{-1}\left\langle \xi\hbar^2\right\rangle^{-\frac{l}{4}}\frac{e^{\pm i\left[\nu\tau\xi^{\frac12}+\hbar^{-1}\rho(x_{\tau};\alpha,\hbar)\right]}}{\xi^{\frac34+\frac{k\nu}{2}}}\left(\log\xi\right)^j\cdot \tilde{F}^{(2),\pm}_{l,k,j}(\tau,\xi), 
	\end{align*}
	where we have 
	\begin{align*}
	&\tilde{F}^{(2),\pm}_{l,k,j}(\tau,\xi)\\& = \tilde{\chi}_{\xi\hbar^2\gtrsim 1}\frac{\left\langle \xi\hbar^2\right\rangle^{\frac{l}{4}}}{\left(\xi\hbar^2\right)^{\frac{l}{4}}}\cdot \frac{\xi^{\frac34}}{\left(\log\xi\right)^j}\int_0^\infty \chi_{ \xi^{-\frac12}\lesssim|\nu\tau-R|\lesssim\hbar}\tilde{\phi}_{\pm,\tau}(R;\xi,\hbar)\frac{G_{k,l,j}(\tau,\nu\tau - R)}{\tau^{\frac12}}\left[\xi^{\frac12}(\nu\tau - R)\right]^{\frac{l}{2}+k\nu}\left(\log\left(\nu\tau - R\right)\right)^j R\,dR,\\
	&\tilde{\phi}_{\pm,\tau}(R;\xi,\hbar): = \frac{e^{\mp i(\nu\tau - R)\xi^{\frac12}\pm i\hbar^{-1}\left[\rho(x;\alpha,\hbar) - \rho(x_{\tau};\alpha,\hbar)\right]}}{\left(R\xi^{\frac12}\right)^{\frac12}}\cdot \psi(x;\alpha,\hbar)\cdot\left(1+\hbar \overline{\tilde{a}_1(\zeta;\alpha)}\right)
	\end{align*}
	Then introduce the new variable 
	\begin{align}\label{def tildeR}
	\tilde{R}:= (\nu\tau - R)\xi^{\frac12}- \hbar^{-1}\left[\rho(x;\alpha,\hbar) - \rho(x_{\tau};\alpha,\hbar)\right]. 
	\end{align}
	Then observe that on account of $x\gg 1$ on the support of the integrand, we have 
	\[
	\frac{\partial\tilde{R}}{\partial R} = -\xi^{\frac12} - \hbar^{-1}\rho_x\cdot\hbar\xi^{\frac12}\simeq -\xi^{\frac12}.
	\]
	Moreover, due to 
	\begin{align*}
	\hbar^{-1}\left[\rho(x;\alpha,\hbar) - \rho(x_{\tau};\alpha,\hbar)\right] &= \hbar^{-1}\left(\int_0^1\rho_x(sx+(1-s)x_\tau;\alpha,\hbar)\,ds\right)\cdot(x-x_\tau)\\
	& = \left(\int_0^1\rho_x(sx+(1-s)x_\tau;\alpha,\hbar)\,ds\right)\cdot(R-\nu\tau)\xi^{\frac12}
	\end{align*}
	and $\left|\rho_x(sx+(1-s)x_\tau;\alpha,\hbar)\right|\lesssim \tau^{-2}$ on the region $|R-\nu\tau|<1,\,x\gtrsim 1$, we see that 
	\[
	\tilde{R}\simeq (\nu\tau - R)\xi^{\frac12}
	\]
	on the support of the integrand. Furthermore, interpreting $\nu\tau - R = (\nu\tau - R)(\tilde{R}, \xi,\tau)$, we get 
	\begin{align}\label{nutau minus R deri}
		\begin{split}
	\left[(\nu\tau - R)\xi^{\frac12}\right]_{,\xi^{\frac12}} &= \hbar^{-1}\left[\rho(x;\alpha,\hbar) - \rho(x_{\tau};\alpha,\hbar)\right]_{,\xi^{\frac12}}\\
	& = \hbar^{-1}\cdot\left(\int_0^1\rho_x(sx+(1-s)x_\tau;\alpha,\hbar)\,ds \cdot(x-x_{\tau})\right)_{,\xi^{\frac12}}\\
	& = \xi^{\frac12}\int_0^1\rho_x(sx+(1-s)x_\tau;\alpha,\hbar)\,ds \cdot(R-\nu\tau)_{,\xi^{\frac12}}\\
	& + \left(\xi^{\frac12}\int_0^1\rho_x(sx+(1-s)x_\tau;\alpha,\hbar)\,ds\right)_{,\xi^{\frac12}}\cdot (R-\nu\tau). 
	\end{split}
	\end{align}
	This relation easily implies that 
	\begin{align*}
	\left|(R-\nu\tau)_{,\xi^{\frac12}}\right|\lesssim \frac{\nu\tau - R}{\xi^{\frac12}}
	\end{align*}
	on the support of the integrand defining $\tilde{F}^{(2),\pm}_{l,k,j}(\tau,\xi)$. In fact the identity \eqref{nutau minus R deri} implies
\begin{align*}
	\left(\nu\tau-R\right)_{,\xi^{\frac12}}\cdot\xi^{\frac12}+\left(\nu\tau-R\right)\simeq \xi^{\frac12}\cdot\frac{1}{\hbar^{2}\tau^{2}\xi}\cdot\left(R-\nu\tau\right)_{,\xi^{\frac12}}+\frac{1}{\hbar^{2}\tau^{2}\xi}\cdot\left(R-\nu\tau\right).
\end{align*}	
The desired estimate follows in view of the fact $\frac{1}{\hbar^{2}\tau^{2}\xi}\ll1$. On the other hand, using the fact $\tilde{R}\simeq (\nu\tau-R)\xi^{\frac12}$ we have
\begin{align*}
	\left|(\nu\tau-R)_{,\tilde{R}}\right|\lesssim \xi^{-\frac12}\simeq \frac{\nu\tau-R}{\tilde{R}}.
\end{align*}	
In fact by definition \eqref{def tildeR}, we have
\begin{align*}
	\left(\nu\tau-R\right)\xi^{\frac12}=\tilde{R}+\left(\int_{0}^{1}\rho_{x}\left(sx+(1-s)x_{\tau};\alpha,\hbar\right)\,ds\right)\cdot \left(R-\nu\tau\right)\xi^{\frac12},
\end{align*}
which implies
\begin{align*}
	\left(\nu\tau-R\right)_{,\tilde{R}}\cdot\xi^{\frac12}=&1+\left(\int_{0}^{1}\rho_{x}\left(sx+(1-s)x_{\tau};\alpha,\hbar\right)\,ds\right)\cdot \left(R-\nu\tau\right)_{,\tilde{R}}\cdot\xi^{\frac12}\\
	&+\left(\int_{0}^{1}\rho_{xx}\left(sx+(1-s)x_{\tau};\alpha,\hbar\right)\cdot s\hbar\xi^{\frac12}\,ds\right)\cdot (R-\nu\tau)_{,\tilde{R}}\cdot\xi^{\frac12}\cdot \left(R-\nu\tau\right).
\end{align*}
Note that
\begin{align*}
	\left|\int_{0}^{1}\rho_{x}\left(sx+(1-s)x_{\tau};\alpha,\hbar\right)\,ds\right|\ll1,\quad \left|\int_{0}^{1}\rho_{xx}\left(sx+(1-s)x_{\tau};\alpha,\hbar\right)\cdot s\hbar\xi^{\frac12}\left(R-\nu\tau\right)\,ds\right|\ll1,
\end{align*}
where the second estimate above follows from the fact $\hbar\left(R-\nu\tau\right)\xi^{\frac12}\ll x$. Therefore the desired estimate $\left(\nu\tau-R\right)_{,\tilde{R}}\simeq  \xi^{-\frac12}$ follows.
Then expressing the integral for $\tilde{F}^{(2),\pm}_{l,k,j}(\tau,\xi)$ in terms of $\tilde{R}$ and performing sufficiently many integrations by parts with respect to $\tilde{R}$, as well as exploiting the symbol behavior of $G_{k,l,j}(\tau, \nu\tau - R)$ and the preceding inequality, the desired symbol behavior for $\tilde{F}^{(2),\pm}_{l,k,j}(\tau,\xi)$ follows easily. In fact, this is a consequence of a slight variation of Lemma~\ref{lem:fromGtoF}. 
\end{proof}
\subsection{Description of the shock on the distorted Fourier side II: admissible expansions}
When formulating the precise asymptotic expansions we shall use for the description of the singular part on the Fourier side, we have to take into account the action of the solution of the wave equation at angular momentum $n$, $|n|\geq 2$, and formulated on the (distorted) Fourier side, given by the Duhamel formula (here we recall Lemma~\ref{lem: parametrix hbar})
\[
\int_{\tau_0}^{\tau}\frac{\lambda(\tau)}{\lambda(\sigma)}\cdot \left(\frac{\rho_{n}\left(\frac{\lambda^2(\tau)}{\lambda^2(\sigma)}\xi\right)}{\rho_{n}(\xi)}\right)^{\frac12}\cdot \frac{\sin\left[\lambda(\tau)\xi^{\frac12}\int_{\sigma}^{\tau}\lambda^{-1}(u)\,du\right]}{\xi^{\frac12}}\cdot F\left(\sigma, \frac{\lambda^2(\tau)}{\lambda^2(\sigma)}\xi\right)\,d\sigma
\]
If $F(\sigma, \xi)$ is a singular source term admitting an expansion as in the preceding subsection, then the re-scaled term $F\left(\sigma, \frac{\lambda^2(\tau)}{\lambda^2(\sigma)}\xi\right)$ will come in general with a phase 
\[
e^{\pm i\left(\nu\sigma\cdot\frac{\lambda(\tau)}{\lambda(\sigma)}\xi^{\frac12} + \hbar^{-1}\rho\left(x_{\sigma\cdot\frac{\lambda(\tau)}{\lambda(\sigma)}};\alpha\cdot\frac{\lambda(\tau)}{\lambda(\sigma)},\hbar\right)\right)},
\]
while the other oscillatory term $\sin\left[\lambda(\tau)\xi^{\frac12}\int_{\sigma}^{\tau}\lambda^{-1}(u)\,du\right]$ will contribute 
\[
e^{\pm i\nu\left(\tau - \frac{\lambda(\tau)}{\lambda(\sigma)}\sigma\right)\cdot\xi^{\frac12}}.
\]
These phases will either enter `destructive resonance', resulting in a singularity propagating precisely along the light cone characteristic $R = \nu\tau$, or else enter into a `constructive resonance' resulting in a singularity propagating outside of the light cone. The latter situation, when translated to the interior of the light cone, will result in smoother terms, in principle analogous to the connecting singular terms from the preceding section, but with a more complicated algebraic structure, which we will have to keep track of. Precisely, in essence we will encounter integrals of the form (where as usual $\alpha = \hbar\xi^{\frac12}$)
\begin{align*}
&e^{\pm i\nu\tau\xi^{\frac12}}\int_{\tau_0}^{\tau}e^{\pm i\hbar^{-1}\rho\left(x_{\sigma\cdot\frac{\lambda(\tau)}{\lambda(\sigma)}};\alpha\cdot\frac{\lambda(\tau)}{\lambda(\sigma)},\hbar\right)}\cdot G\left(\sigma, \frac{\lambda^2(\tau)}{\lambda^2(\sigma)}\xi\right)\,d\sigma,\\
&\int_{\tau_0}^{\tau} e^{\pm i\left(\left(\nu\tau - 2\frac{\lambda(\tau)}{\lambda(\sigma)}\nu\sigma\right)\xi^{\frac12} - \hbar^{-1}\rho\left(x_{\sigma\cdot\frac{\lambda(\tau)}{\lambda(\sigma)}};\alpha\cdot\frac{\lambda(\tau)}{\lambda(\sigma)},\hbar\right)\right)}\cdot G\left(\sigma, \frac{\lambda^2(\tau)}{\lambda^2(\sigma)}\xi\right)\,d\sigma,
\end{align*}
where the two $\pm$-signs in the first integral are synchronized, and where the factors $G\left(\sigma, \frac{\lambda^2(\tau)}{\lambda^2(\sigma)}\xi\right)$ have suitable symbol behaviour. The key point concerning the second integral represents a function whose singularity (on the physical side) is located at $R = 2\frac{\lambda(\tau)}{\lambda(\sigma)}\nu\sigma - \nu\tau\geq \nu\tau$ due to $\sigma\leq \tau$, and the integration over $\sigma$ ensures this function, upon restriction to the interior of the light cone, is of regularity $H^{2+}$. 
\\
Of course, the preceding integration expressions need to be re-inserted into suitable source terms and again subjected to the wave parametrix, resulting in more complex expressions in principle. However, a manageable formalism is possible to capture all the terms that arise in this fashion, which is realized in the following definition. We emphasize that the bounds of the coefficient functions closely mimic the ones of Definition~\ref{defi:xsingulartermsngeq2proto}:
\begin{definition}\label{defi:xsingulartermsngeq2adm} We call a function $\xb(\tau,\xi)$, $\xi\in [0,\infty)$, an admissible singular part at angular momentum $n, |n|\geq 2$, provided it allows a representation 
	\[
	\xb(\tau, \xi) = \xb_{\text{in}}(\tau, \xi) + \xb_{\text{out}}(\tau, \xi)
	\]
	where $\xb_{\text{in}}$ represents the {\it{ingoing part}} of the singularity, which constitutes the dominant part, while $\xb_{\text{out}}$ represents the {\it{outgoing part}}, and the functions on the right admit the following expansions: 
	\\
	We can write (in the following all cutoffs $\chi$ are smooth)
	\begin{equation}\label{eq:xsingin}\begin{split}
	\xb_{\text{in}}(\tau, \xi) =& \sum_{\pm}\sum_{k=1}^{N}\sum_{j=0}^{N_1}\chi_{\xi\geq \hbar^{-2}}\hbar^{-1}\frac{e^{\pm i\nu\tau\xi^{\frac12}}}{\xi^{1+k\frac{\nu}{2}}}\left(\log\xi\right)^j\cdot \int_{\tau_0}^{\tau}e^{\pm i\hbar^{-1}\rho\left(x_{\sigma\cdot\frac{\lambda(\tau)}{\lambda(\sigma)}};\alpha\cdot\frac{\lambda(\tau)}{\lambda(\sigma)},\hbar\right)}\cdot a_{k,j}^{(\pm)}(\tau,\sigma)\,d\sigma\\
	+&  \sum_{\pm}\sum_{l =1}^7\sum_{k=1}^{N}\sum_{j=0}^{N_1}\chi^{(l)}_{\xi\geq \hbar^{-2}}\hbar^{-1}\left\langle\hbar^2\xi\right\rangle^{-\frac{l}{4}}\frac{e^{\pm i\nu\tau\xi^{\frac12}}}{\xi^{1+k\frac{\nu}{2}}}\left(\log\xi\right)^j\\
	&\hspace{4.5cm}\cdot \int_{\tau_0}^{\tau}e^{\pm i\hbar^{-1}\rho\left(x_{\sigma\cdot\frac{\lambda(\tau)}{\lambda(\sigma)}};\alpha\cdot\frac{\lambda(\tau)}{\lambda(\sigma)},\hbar\right)}\cdot  F_{l,k,j}^{(\pm)}\left(\tau, \sigma, \frac{\lambda^2(\tau)}{\lambda^2(\sigma)}\xi\right)\,d\sigma,\\
	\end{split}\end{equation}
	where the $\pm$-signs in each expression on the right are synchronized, and we have the following bounds, where $\delta_1, \delta$ are small positive numbers : 
	\begin{align*}
	&\left|a_{k,j}^{(\pm)}(\tau,\sigma)\right| + \tau^{k_1}\left|\partial_{\tau}^{k_1}a_{k,j}^{(\pm)}(\tau,\sigma)\right|\lesssim \left(\log\tau\right)^{N_1-j}\tau^{-1-\nu}\cdot \sigma^{-3},\,k_1\in\{0, 1\}\\
	&\hbar^{\delta_1\cdot k_2}\cdot \left|\xi^{k_2}\partial_{\tau}^{\iota}\partial_{\xi}^{k_2} F_{l,k,j}^{(\pm)}(\tau, \sigma, \xi)\right|\lesssim_{k_1}  \left(\log\tau\right)^{N_1-j}\hbar^{-\iota}\tau^{-1-\nu-\iota}\cdot\sigma^{-1}\cdot \left[\sigma^{-2} + \kappa\left(\hbar\xi^{\frac12}\right)\right],\quad  0\leq k_2\leq 20-l,\quad \iota\in\{0,1\}\\
	&\hbar^{k_*}\cdot\left\|\xi^{20-l+\delta}\partial_{\tau}^{\iota}\partial_{\xi}^{20 -l} F_{l,k,i}^{(\pm)}(\tau, \sigma, \xi)\right\|_{\dot{C}^{\delta}_{\xi}(\xi\simeq\mu)}\lesssim_{k_1}  \left(\log\tau\right)^{N_1-j}\tau^{-1-\nu-\iota}\cdot\sigma^{-1}\cdot \left[\sigma^{-2} + \kappa\left(\hbar\mu^{\frac12}\right)\right],\quad \iota\in\{0,1\}
	\end{align*}
	where $k_* = k_*(\delta_1,\delta,l,k_1)$ is defined as in Definition~\ref{defi:xsingulartermsngeq2proto}, and we set $\kappa(x): = \frac{x}{1+x^2}$, while $\mu$ ranges over all positive numbers. We shall again call the first sum on the right in \eqref{eq:xsingin} the principal part of the ingoing singularity, while we call the second sum the connecting part of the ingoing singularity.
	\\
	For the outgoing part, we can split it into 
	\[
	\xb_{\text{out}} = \xb_{\text{out}, 1} + \xb_{\text{out}, 2}
	\]
	where the first term on the right represents the `outgoing seed singularity' given by 
	\begin{equation}\label{eq:xsingout1}\begin{split}
	\xb_{\text{out}, 1}(\tau, \xi) =  &\sum_{\pm}\sum_{l =0}^7\sum_{k=1}^{N}\sum_{j=0}^{N_1}\chi_{\xi\geq \hbar^{-2}}\hbar^{-1}\left\langle\hbar^2\xi\right\rangle^{-\frac{l}{4}}\frac{\left(\log\xi\right)^j}{\xi^{1+\frac{k\nu}{2}}}\\&\hspace{2cm}\cdot\int_{\tau_0}^{\tau}e^{\pm i\left[\left(\nu\tau - 2\frac{\lambda(\tau)}{\lambda(\sigma)}\nu\sigma\right)\xi^{\frac12} - \hbar^{-1}\rho\left(x_{\sigma\cdot\frac{\lambda(\tau)}{\lambda(\sigma)}};\alpha\cdot\frac{\lambda(\tau)}{\lambda(\sigma)},\hbar\right)\right]}\cdot \tilde F^{\pm}_{l,k,j}\left(\tau,\sigma, \frac{\lambda^2(\tau)}{\lambda^2(\sigma)}\xi\right)\,d\sigma, 
	\end{split}\end{equation}
	where $\tilde{F}^{\pm}_{0,k,i}\left(\tau,\sigma, \frac{\lambda^2(\tau)}{\lambda^2(\sigma)}\xi\right) = b^{\pm}_{k,i}(\tau, \sigma)$, and we have the bounds 
	\begin{align*}
	&\hbar^{\delta_1\cdot r}\cdot\left|\xi^{r}\partial_{\tau}^{\iota}\partial_{\xi}^{r} \tilde{F}_{l,k,j}^{(\pm)}(\tau, \sigma, \xi)\right|\lesssim \left(\log\tau\right)^{N_1-j}\hbar^{-\iota}\tau^{-1-\nu-\iota}\cdot\sigma^{-1}\cdot \left[\sigma^{-2} + \kappa\left(\hbar\xi^{\frac12}\right)\right], \quad 0\leq r\leq 20-l,\quad \iota\in\{0,1\}\\
	&\hbar^{k_*}\cdot\left\|\xi^{20-l+\delta}\partial_{\tau}^{\iota}\partial_{\xi}^{20-l} \tilde{F}_{l,k,j}^{(\pm)}(\tau,\sigma, \xi)\right\|_{\dot{C}^{\delta}_{\xi}(\xi\simeq\mu)}\lesssim \left(\log\tau\right)^{N_1-j}\tau^{-1-\nu-\iota}\cdot\sigma^{-1}\cdot \left[\sigma^{-2} + \kappa\left(\hbar\mu^{\frac12}\right)\right],\quad \iota\in\{0,1\}.
	\end{align*}
	where at the end $k_* = k_*(\delta_1,\delta,l,\iota)$. 
	The second term $\xb_{\text{out}, 2}$, which is the `outgoing perpetuated singularity' admits the description
	\begin{equation}\label{eq:xsingout2}\begin{split}
	\xb_{\text{out}, 2}(\tau, \xi) =  &\sum_{\pm}\sum_{l =0}^7\sum_{k=1}^{N}\sum_{j=0}^{N_1}\chi_{\xi\geq \hbar^{-2}}\hbar^{-1}\left\langle\hbar^2\xi\right\rangle^{-\frac{l}{4}}\frac{\left(\log\xi\right)^j}{\xi^{1+k\frac{\nu}{2}}}\\&\cdot\int_0^\infty \int_{\tau_0}^{\tau}e^{\pm i\left[\nu\left(\frac{\lambda(\tau)}{\lambda(\sigma)}x + \tau\right)\xi^{\frac12} + \hbar^{-1}\rho\left(x_{\sigma\cdot\frac{\lambda(\tau)}{\lambda(\sigma)}};\alpha\cdot\frac{\lambda(\tau)}{\lambda(\sigma)},\hbar\right)\right]}\cdot G^{\pm}_{l,k,j}\left(\tau, \sigma,x,\frac{\lambda^2(\tau)}{\lambda^2(\sigma)}\xi\right)\,d\sigma dx,
	\end{split}\end{equation}
	where $G^{\pm}_{0,k,i}\left(\tau, \sigma,x,\frac{\lambda^2(\tau)}{\lambda^2(\sigma)}\xi\right) = c^{\pm}_{k,i}(\tau,\sigma,x)$, and we have the bounds 
	\begin{align*}
	&\hbar^{\delta_1\cdot k_2}\left\|\xi^{k_2}\partial_{\xi}^{k_2} \partial_{\tau}^{\iota}G_{l,k,i}^{(\pm)}(\tau, \sigma, x, \xi)\right\|_{L_x^1}\lesssim  \left(\log\tau\right)^{N_1-i}\hbar^{-\iota}\tau^{-1-\nu-\iota}\cdot\sigma^{-1}\cdot \left[\sigma^{-2} + \kappa\left(\hbar\xi^{\frac12}\right)\right],\quad 0\leq k_2\leq 20-l,\,\iota\in \{0,1\}\\
	&\hbar^{k_*}\cdot\left\|\left\|\xi^{20-l+\delta_l}\partial_{\xi}^{20-l} \partial_{\tau}G_{l,k,i}^{(\pm)}(\tau, \sigma, x, \xi)\right\|_{\dot{C}^{\delta_l}_{\xi}(\xi\simeq\lambda)}\right\|_{L_x^1}\lesssim  \left(\log\tau\right)^{N_1-i}\tau^{-1-\nu-\iota}\cdot\sigma^{-1}\cdot \left[\sigma^{-2} + \kappa\left(\hbar\lambda^{\frac12}\right)\right],\,\iota\in \{0,1\}
	\end{align*}
	We say that the principal singular part is of restricted type, provided for $k\in\{1,2,3\}$ we have 
	\begin{align*}
	&a_{k,i}(\tau, \sigma) = \sum_{l=i}^{N_1}c_{l,k,i}^{(\pm)}(\sigma)\cdot \tau^{-k(1+\nu)}\cdot \left(\log\tau\right)^{l-i} + b_{k,i}^{(\pm)}(\tau, \sigma),\\
	&\left|c_{l,k,i}^{(\pm)}(\sigma)\right|\lesssim \sigma^{-3},\quad \left|b_{k,i}(\tau, \sigma)\right| + \tau\cdot\left|b_{k,i}'(\tau, \sigma)\right|\lesssim \tau^{-3-\nu}\left(\log\tau\right)^{N_1}\cdot\sigma^{-3},
	\end{align*}
	while for $k\geq 4$ we have 
	\begin{align*}
	\left|a_{k,i}^{(\pm)}(\tau,\sigma)\right| + \tau\left|\partial_{\tau}a_{k,i}^{(\pm)}(\tau,\sigma)\right|\lesssim \left(\log\tau\right)^{N_1-i}\tau^{-3-\nu}\cdot \sigma^{-3}.
	\end{align*}
	Finally, we say that a function $\overline{y}(\tau,\xi)$ is source admissibly singular, provided we can write $\overline{y} = \overline{y}_1 + \overline{y}_2$, and where $\xi^{-\frac12}\cdot \overline{y}_1(\tau,\xi)$ is a admissibly singular except that we replace $\tau^{-1-\nu}$ by $\tau^{-2-\nu}$ and we also set $k_1 = 0$ throughout; furthermore, we require $\overline{y}_2$ to be a prototype singular source term. We analogously define source admissibly singular functions of restricted type. 
\end{definition}
\begin{remark}\label{rem:admsingdef} The reason for including the weight $\kappa\left(\hbar\xi^{\frac12}\right)$ in the estimates comes from the action of the transference operator on admissibly singular terms, and subsequent application of the Duhamel parametrix. These compositions will be required to solve the angular momentum $n$ wave equations when translated to the distorted Fourier side, as becomes apparent upon inspecting \eqref{eq:ngeq2Fourier2}. The function $\kappa$ reflects the decay properties of the kernel $F_n(\xi, \eta)$ as described in Proposition~\ref{prop: K operator}. The product $\sigma^{-1}\cdot \kappa\left(\hbar\frac{\lambda(\tau)}{\lambda(\sigma)}\xi^{\frac12}\right)$ is absolutely integrable with respect $\sigma $ on $[\tau_0,\infty)$, which is important in controlling the size of the $\sigma$-integral. We also note that the small power of $\hbar$ weights $\hbar^{\delta_1\cdot k_2}$, are due to the action of the transference operator. 
\end{remark}
The preceding space of functions comes equipped with a natural norm, given in the following 
\begin{definition}\label{defi:xsingulartermsngeq2admnorm}  Assume that $\xb(\tau, \xi)$ is an admissible singular part at angular momentum $n, |n|\geq 2$. Then we set 
	\begin{align*}
	&\lVert x\rVert_{\text{adm}} : = \\
	&\sum_{\pm}\sum_{k=1}^{N}\sum_{i=0}^{N_1}\left\|\left[\left|a_{k,i}^{(\pm)}(\tau,\sigma)\right| + \tau\left|\partial_{\tau}a_{k,i}^{(\pm)}(\tau,\sigma)\right|\right]\cdot\left(\log\tau\right)^{-N_1+i}\tau^{1+\nu}\cdot \sigma^{3}\right\|_{L^\infty\left([\tau_0,\infty)\right)L^\infty\left([\tau_0,\tau]\right)}\\
	&+\sum_{\pm}\sum_{l =0}^7\sum_{k=1}^{N}\sum_{i=0}^{N_1}\sum_{\substack{0\leq k_2\leq 20-l\\ \iota\in\{0,1\}}}\hbar^{\iota+\delta_1\cdot k_2}\left\|\left(\log\tau\right)^{-N_1+i}\tau^{1+\nu+\iota}\sigma\cdot \left\|\frac{\xi^{k_2}\partial_{\xi}^{k_2}\partial_{\tau}^{\iota}F^{(\pm)}_{l,k,i}(\tau,\sigma,\xi)}{\left(\sigma^{-2} + \kappa(\hbar\xi^{\frac12})\right)}\right\|_{L_{\xi}^\infty\left([0,\infty)\right)}\right\|_{L_{\tau}^\infty\left([\tau_0,\infty)\right)L_{\sigma}^\infty\left([\tau_0,\tau]\right)}\\
	& + \sum_{\pm}\sum_{l =0}^7\sum_{k=1}^{N}\sum_{i=0}^{N_1}\sum_{\substack{0\leq k_2\leq 20-l\\ \iota\in\{0,1\}}}\hbar^{\iota+\delta_1\cdot k_2}\left\|\left(\log\tau\right)^{-N_1+i}\tau^{1+\nu+\iota}\sigma\cdot \left\|\frac{\xi^{k_2}\partial_{\xi}^{k_2}\partial_{\tau}^{\iota}\tilde{F}^{(\pm)}_{l,k,i}(\tau,\sigma,\xi)}{\left(\sigma^{-2} + \kappa(\hbar\xi^{\frac12})\right)}\right\|_{L_{\xi}^\infty\left([0,\infty)\right)}\right\|_{L_{\tau}^\infty([\tau_0,\infty)L_{\sigma}^\infty([\tau_0,\tau]}\\
	& + \sum_{\pm}\sum_{l =0}^7\sum_{k=1}^{N}\sum_{i=0}^{N_1}\sum_{\iota}\hbar^{k_*}\left\|\left(\log\tau\right)^{-N_1+i}\tau^{1+\nu+\iota}\sigma\cdot \left\|\frac{\xi^{20-l+\delta}\partial_{\xi}^{20-l}\partial_{\tau}^{\iota}F^{(\pm)}_{l,k,i}(\tau,\sigma,\xi)}{\left(\sigma^{-2} + \kappa(\hbar\xi^{\frac12})\right)}\right\|_{\dot{C}_{\xi}^{\delta}}\right\|_{L_{\tau}^\infty\left([\tau_0,\infty)\right)L_{\sigma}^\infty\left([\tau_0,\tau]\right)}\\
	& +  \sum_{\pm}\sum_{l =0}^7\sum_{k=1}^{N}\sum_{i=0}^{N_1}\sum_{\iota}\hbar^{k_*}\left\|\left(\log\tau\right)^{-N_1+i}\tau^{1+\nu+\iota}\sigma\cdot \left\|\frac{\xi^{20-l+\delta}\partial_{\xi}^{20-l}\partial_{\tau}^{\iota}\tilde{F}^{(\pm)}_{l,k,i}(\tau,\sigma,\xi)}{\left(\sigma^{-2} + \kappa(\hbar\xi^{\frac12})\right)}\right\|_{\dot{C}_{\xi}^{\delta}}\right\|_{L_{\tau}^\infty\left([\tau_0,\infty)\right)L_{\sigma}^\infty\left([\tau_0,\tau]\right)}\\
	& + \sum_{\pm}\sum_{l =0}^7\sum_{k=1}^{N}\sum_{i=0}^{N_1}\sum_{0\leq k_2\leq 20-l}\sum_{\iota}\hbar^{\iota+\delta_1\cdot k_2}\left\|\left(\log\tau\right)^{-N_1+i}\tau^{1+\nu+\iota}\sigma\cdot \left\|\frac{\xi^{k_2}\partial_{\xi}^{k_2}\partial_{\tau}^{\iota}G^{(\pm)}_{l,k,i}(\tau,\sigma,x,\xi)}{\left(\sigma^{-2} + \kappa(\hbar\xi^{\frac12})\right)}\right\|_{L_{\xi}^\infty([0,\infty)}\right\|_{L_x^1\left([0,\infty)\right)L_{\tau}^\infty\left(([\tau_0,\infty)\right)L_{\sigma}^\infty\left(([\tau_0,\tau]\right)}\\
	& +\sum_{\pm}\sum_{l =0}^7\sum_{k=1}^{N}\sum_{i=0}^{N_1}\sum_{\iota} \hbar^{k_*}\left\|\left(\log\tau\right)^{-N_1+i}\tau^{1+\nu+\iota}\sigma\cdot \left\|\frac{\xi^{20-l}\partial_{\xi}^{20-l}\partial_{\tau}^{\iota}G^{(\pm)}_{l,k,i}(\tau,\sigma,x,\xi)}{\left(\sigma^{-2} + \kappa(\hbar\xi^{\frac12})\right)}\right\|_{\dot{C}_{\xi}^{\delta}}\right\|_{L_x^1\left([0,\infty)\right)L_{\tau}^\infty\left([\tau_0,\infty)\right)L_{\sigma}^\infty\left([\tau_0,\tau]\right)}.\\
	\end{align*}
	If the principal singular part is of restricted type, we shall replace the first term 
	\[
	\sum_{\pm}\sum_{k=1}^{N}\sum_{i=0}^{N_1}\left\|\left[\left|a_{k,i}^{(\pm)}(\tau,\sigma)\right| + \tau\left|\partial_{\tau}a_{k,i}^{(\pm)}(\tau,\sigma)\right|\right]\cdot\left(\log\tau\right)^{-N_1+i}\tau^{1+\nu}\cdot \sigma^{3}\right\|_{L^\infty\left([\tau_0,\infty)\right)L^\infty\left([\tau_0,\tau]\right)}
	\]
	by 
	\begin{align*}
	&\sum_{\pm}\sum_{l =0}^7\sum_{k=1,2,3}\left\|\sigma^{3}\cdot c_{l,k,i}^{(\pm)}(\sigma)\right\|_{L^\infty\left([\tau_0,\infty)\right)} + \sum_{\pm}\sum_{k=1,2}\sum_{i=0}^{N_1}\left\|b_{k,i}^{(\pm)}(\tau, \sigma)\cdot\left(\log\tau\right)^{-N_1+i}\tau^{3+\nu}\cdot \sigma^{3}\right\|_{L^\infty\left([\tau_0,\infty)\right)L^\infty\left([\tau_0,\tau]\right)}\\
	& + \sum_{\pm}\sum_{k=4}^{N}\sum_{i=0}^{N_1}\left\|\left[\left|a_{k,i}^{(\pm)}(\tau,\sigma)\right| + \tau\left|\partial_{\tau}a_{k,i}^{(\pm)}(\tau,\sigma)\right|\right]\cdot\left(\log\tau\right)^{-N_1+i}\tau^{3+\nu}\cdot \sigma^{3}\right\|_{L^\infty\left([\tau_0,\infty)\right)L^\infty\left([\tau_0,\tau]\right)}\\
	\end{align*}
	and call the resulting norm $\lVert\cdot\rVert_{\text{adm}(r)}$. If $\overline{y} = \overline{y}_1 + \overline{y}_2$ is source admissibly singular, we define $\Big\|\overline{y}_1\Big\|$ in analogy to the preceding, and then set 
	\[
	\big\|\overline{y}\big\|_{sourceadm} = \Big\|\overline{y}_1\Big\| + \big\|\overline{y}_2\big\|_{sourceproto},
	\]
	with the latter norm as in Def.~\ref{defy: protofunctionnorm}.
\end{definition}

In the rest of this subsection as well as the next one, we verify that this choice of function space is compatible with the key operations that will arise during the iterative scheme. The first order of the day shall be to recover an analogue of Lemma~\ref{lem:singFouriertiphysicalngeq2}. The key point is that the next lemma has an almost identical conclusion.

 \begin{lemma}\label{lem:singFouriertiphysicalngeq2adm} Assume that $\xb$ is an admissible singular part at angular momentum $n, |n|\geq 2$. Then the associated function 
	\[
	f(\tau, R): = \int_0^\infty \phi_{n}(R;\xi)\cdot \xb(\tau, \xi)\cdot\rho_{n}(\xi)\,d\xi,
	\]
	restricted to the light cone $R<\nu\tau$, can be decomposed as 
	\[
	f = f_1 + f_2+f_3
	\]
	where $f_1 = f_1(\tau, R)$ is a $C^5$-function supported in $\nu\tau - R\gtrsim 1$ and satisfying 
	\[
	\nabla_R^{k_2}\partial_{\tau}^{k_1}f_1(\tau, R)\lesssim \hbar^{-\frac12+}\left(\log\tau\right)^{N_1}\cdot \tau^{-\frac32-\nu}\left|\nu\tau - R\right|^{-5},\quad 0\leq k_1 + k_2\leq 5, 
	\]
	while $f_2 =\sum_{l=1}^8 f_{2l}$ where we have the explicit form
	\begin{align*}
	f_{2l}(\tau, R) = \chi_{|\nu\tau-R|\lesssim \hbar}\sum_{k=1}^N\sum_{j=0}^{N_1}\frac{G_{k,l,j}(\tau,\nu\tau - R)}{\tau^{\frac12}}\hbar^{-\frac{l+1}{2}}\left[\nu\tau - R\right]^{\frac{l}{2}+k\nu}\left(\log(\nu\tau - R)\right)^j
	\end{align*}
	Here the function $G_{k,j}(\tau, x)$ has symbol type behavior with respect to $x$, as follows:
	\begin{align*}
	\hbar^{\delta_1\cdot (20-l)}\cdot \left|\partial_x^{k_2}\partial_{\tau}^{\iota}G_{k,l,j}(\tau, x)\right| \lesssim  \left(\log\tau\right)^{N_1-j}\cdot\tau^{-1-\nu-\iota}x^{-k_2},\quad 0\leq k_2\leq 20 - l + \lfloor \frac{l}{2}+k\nu\rfloor,\,\iota\in \{0, 1\},
	\end{align*}
	and we have the bound 
	\begin{align*}
	\hbar^{k_*}\cdot\left\|x^{20-l_*+\delta_3}\partial_x^{20-l_*}\partial_{\tau}^{\iota}G_{k,l,j}(\tau, x)\right\|_{\dot{C}^{\delta_3}} \lesssim \left(\log\tau\right)^{N_1-j}\cdot \tau^{-1-\nu-\iota},\,\iota\in \{0, 1\},
	\end{align*}
	where $k_*, l_*, \delta_3$ are defined in analogy to the statement of Lemma~\ref{lem:singFouriertiphysicalngeq2} .
	Finally, the remaining function $f_3$ is also $C^5$ and supported in $|\nu\tau - R|\lesssim 1$ and satisfies
	\[
	\left\|f_3\right\|_{\tilde{S}_0^{(\hbar)}}\lesssim \left(\log\tau\right)^{N_1}\cdot\tau^{-1-\nu},\quad  \left\|\partial_{\tau}f_3\right\|_{\tilde{S}_1^{(\hbar)}}\lesssim  \left(\log\tau\right)^{N_1}\cdot \tau^{-1-\nu}.
	\]
	Moreover, we have the bounds 
	\begin{align*}
	&\left|\partial_R^k f_3\right|\lesssim \left(\log\tau\right)^{N_1}\tau^{-\frac32-\nu}\cdot\hbar^{-\frac12+}\min\{(\nu\tau - R)^{-k},\hbar^{-k}\},\quad
	0\leq k\leq 5,\\
	&\left|\partial_R^k\partial_{\tau} f_3\right|\lesssim \left(\log\tau\right)^{N_1}\tau^{-\frac32-\nu}\cdot\hbar^{-\frac12+}\min\{(\nu\tau - R)^{-k-1},\hbar^{-k-1}\},\quad
	0\leq k\leq 4.
	\end{align*}
	Furthermore, the terms $\xb_{\text{out}}$ enjoy a smoothness gain visible on the physical side upon restriction to the interior of the light cone: defining 
	\[
	g(\tau, R): = \int_0^\infty \phi_{n}(R;\xi)\cdot \xb_{\text{out}}(\tau, \xi)\cdot\rho_{n}(\xi)\,d\xi,
	\]
	we can write 
	\[
	g = g_1 + g_2+g_3, 
	\]
	where $g_{1,3}$ have the same properties as $f_{1,3}$, while $g_2$ admits the representation 
	\[
	g_2 = \chi_{|\nu\tau-R|\lesssim \hbar}\sum_{l=3}^{8}\sum_{k=1}^N\sum_{j=0}^{N_1}\frac{H_{k,l,j}(\tau,\nu\tau - R)}{\tau^{\frac12}}\hbar^{-\frac{l+1}{2}}\left[\nu\tau - R\right]^{\frac{l}{2}+k\nu}\left(\log(\nu\tau - R)\right)^j
	\]
	and we have the symbol type bounds 
	\begin{align*}
	\hbar^{\delta_1\cdot(k_2+3)}\left|\partial_x^{k_2}\partial_{\tau}^{\iota}H_{k,l,j}(\tau, x)\right| \lesssim  \left(\log\tau\right)^{N_1-j}\cdot\hbar^{-1}\tau^{-1-\nu-\iota}x^{-k_2},\quad k_2\in\{0,\ldots,15-l + \lfloor \frac{l}{2}+k\nu\rfloor\},
	\end{align*}
	as well as
	\begin{align*}
	\hbar^{k_*}\cdot\left\|x^{15-l_*+\delta_3}\partial_x^{15-l_*+\delta_3}\partial_{\tau}^{\iota}H_{k,l,j}(\tau, x)\right\|_{\dot{C}^{\delta_3}} \lesssim \left(\log\tau\right)^{N_1-j}\cdot \tau^{-1-\nu-\iota}.
	\end{align*}
\end{lemma}
\begin{proof}
	We shall consider the first term in \eqref{eq:xsingin}, and explain later how to get the improved regularity for the outgoing terms. The contribution of the second term in \eqref{eq:xsingin} requires use of Lemma~\ref{lem:FtoGimprovement}, but is otherwise treated in the same manner, as are the subsequent terms.\\
	  We can then work at fixed time $\sigma$ and deduce the final result by integration over $\sigma$. The proof follows exactly the one of  Lemma~\ref{lem:singFouriertiphysicalngeq2}, and we explain the minor differences arising.
	\\
	{\it{The case of ingoing singularity, i.e., the contribution of \eqref{eq:xsingin}.}} The argument is the same for the two integral expressions up to minor details. Fixing $\sigma\in [\tau_0,\tau]$, we encounter the oscillatory expression 
	\[
	e^{\pm i\left(\nu\tau\xi^{\frac12} + \hbar^{-1}\rho(x_{\sigma'};\alpha',\hbar)\right)},\quad \sigma' = \sigma\cdot \frac{\lambda(\tau)}{\lambda(\sigma)}\geq\tau,\quad \alpha' = \alpha\cdot\frac{\lambda(\tau)}{\lambda(\sigma)}. 
	\]
	Then going through the case {\it{(1)}} - {\it{(4)}} in the proof of Lemma~\ref{lem:singFouriertiphysicalngeq2}, in cases  {\it{(1)}},  {\it{(2)}} we replace $e^{\pm i\nu\tau\xi^{\frac12}}$ by $e^{\pm i\left(\nu\tau + \hbar^{-1}\rho(x_{\sigma'};\alpha',\hbar)\right)\xi^{\frac12}}$, then we use the bound 
	\begin{align*}
	\partial_{\xi^{\frac12}}\left(\nu\tau\xi^{\frac12} + \hbar^{-1}\rho(x_{\sigma'};\alpha',\hbar)\right) &= \nu\tau + \hbar^{-1}\rho_x(x_{\sigma'};\alpha',\hbar)\cdot \hbar\nu\sigma' + \hbar^{-1}\rho_{\alpha}(x_{\sigma'};\alpha',\hbar)\cdot \hbar\\
	& = \nu\tau + O\left(\sigma'^{-1}+\hbar\right)\gtrsim \tau
	\end{align*}
	on the support of the expression which we recall is $\xi\gtrsim \hbar^{-2}$. 
	\\
	In case {\it{(3)}}, repeating the notation from the end of the proof of Lemma~\ref{lem:singFouriertiphysicalngeq2}, we encounter the phase 
	\[
	\pm \frac23\hbar^{-1}\zeta^{\frac32} \pm \left(\nu\tau\xi^{\frac12} + \hbar^{-1}\rho(x_{\sigma'};\alpha,\hbar)\right), 
	\]
	and thanks to $\sigma'\geq \tau$ and the condition $\xi\gtrsim \hbar^{-2}$ on the support, we again have  
	\begin{align*}
	\left|\partial_{\xi^{\frac12}}\left(\pm \frac23\hbar^{-1}\zeta^{\frac32} \pm \left(\nu\tau\xi^{\frac12} + \hbar^{-1}\rho(x_{\sigma'};\alpha,\hbar)\right)\right)\right|\gtrsim \tau.
	\end{align*}
	Case {\it{(4)}}, which we recall means $x\geq \left(1+\gamma\right)x_t$, is again the most complicated, and we argue as follows. Here in the regime $\nu\tau - R\gg 1$ and the most delicate case of resonant phases we encounter the phase function 
	\begin{align*}
	&\pm(\nu\tau - R)\xi^{\frac12}\mp \hbar^{-1}\left[\rho(x;\alpha,\hbar) - \rho(x_{\sigma'};\alpha',\hbar) - y(\alpha;\hbar)\right]\\
	& = \pm(\nu\tau - R)\xi^{\frac12}\mp \hbar^{-1}\left[\rho(x;\alpha,\hbar) - \rho(x_{\tau};\alpha,\hbar) - y(\alpha;\hbar)\right] \pm \hbar^{-1}\left[ \rho(x_{\sigma'};\alpha',\hbar)  -  \rho(x_{\tau};\alpha,\hbar)\right]\\
	&=:\pm\Omega(\tau,R,\alpha,\hbar)
	\end{align*}
	and the following lower bound obtains: 
	\begin{align*}
	\partial_{\xi^{\frac12}}\Omega&\gtrsim \nu\tau - R + \hbar^{-1}\rho_x(x_{\sigma'};\alpha',\hbar)\cdot \hbar\nu\sigma' - \hbar^{-1}\rho_x(x_{\tau};\alpha,\hbar)\cdot \hbar\nu\tau + O(\hbar)\\
	&\gtrsim \nu\tau - R
	\end{align*}
	provided $\nu\tau - R\gg \hbar$, since we have the important positivity property 
	\begin{align}\label{monotonicity rho deri}
	\hbar^{-1}\rho_x(x_{\sigma'};\alpha',\hbar)\cdot \hbar\nu\sigma' - \hbar^{-1}\rho_x(x_{\tau};\alpha,\hbar)\cdot \hbar\nu\tau\geq 0
	\end{align}
	on account of $\sigma'\geq\tau,\,\alpha'\geq \alpha$. This implies that the term (we re-use the notation from case {\it{(4)}} in the proof of Lemma~\ref{lem:singFouriertiphysicalngeq2}) and where the temporally decaying factor $a_{k,j}^{\pm}(\tau, \sigma)$ according to Definition~\ref{defi:xsingulartermsngeq2adm} has been included into $ \Psi(x;\alpha,\hbar)$ (recall \eqref{eq:Psixalphahbardef}), whence we denote this now by $\Psi_{\tau,\sigma}(x;\alpha,\hbar)$, 
	\[
	\chi_{\nu\tau-R\gg 1}\int_0^\infty e^{\pm i(\nu\tau - R)\xi^{\frac12}}\cdot e^{\mp i\hbar^{-1}\left[\rho(x;\alpha,\hbar) - y(\alpha;\hbar) - \rho(x_{\sigma'};\alpha',\hbar)\right]}\cdot \Psi_{\tau,\sigma}(x;\alpha,\hbar)\,d\xi
	\]
	can be included into $f_1(\tau, R)$, and similarly the term 
	\[
	\chi_{1\gtrsim\nu\tau-R\gg \hbar}\int_0^\infty e^{\pm i(\nu\tau - R)\xi^{\frac12}}\cdot e^{\mp i\hbar^{-1}\left[\rho(x;\alpha,\hbar) - y(\alpha;\hbar) - \rho(x_{\sigma'};\alpha',\hbar)\right]}\cdot \Psi_{\tau,\sigma}(x;\alpha,\hbar)\,d\xi
	\]
	can be incorporated into $f_3(\tau, R)$. 
	\\
	This reduces things to controlling the integral 
	\[
	\chi_{\nu\tau-R\lesssim\hbar}\int_0^\infty e^{\pm i(\nu\tau - R)\xi^{\frac12}}\cdot e^{\mp i\hbar^{-1}\left[\rho(x;\alpha,\hbar) - y(\alpha;\hbar) - \rho(x_{\sigma'};\alpha',\hbar)\right]}\cdot \Psi_{\tau,\sigma}(x;\alpha,\hbar)\,d\xi
	\]
	Write this as 
	\begin{equation}\label{eq:maincontribxin}\begin{split}
	&\chi_{\nu\tau-R\lesssim\hbar}\int_0^\infty e^{\pm i(\nu\tau - R)\xi^{\frac12}}\cdot e^{\mp i\hbar^{-1}\left[\rho(x;\alpha,\hbar) - y(\alpha;\hbar) - \rho(x_{\sigma'};\alpha',\hbar)\right]}\cdot \Psi_{\tau,\sigma}(x;\alpha,\hbar)\,d\xi\\&
	=\chi_{\nu\tau-R\lesssim\hbar}\int_0^\infty\chi_{(\nu\tau - R)\xi^{\frac12}\lesssim1} e^{\pm i(\nu\tau - R)\xi^{\frac12}}\cdot e^{\mp i\hbar^{-1}\left[\rho(x;\alpha,\hbar) - y(\alpha;\hbar) - \rho(x_{\sigma'};\alpha',\hbar)\right]}\cdot \Psi_{\tau,\sigma}(x;\alpha,\hbar)\,d\xi\\
	& + \chi_{\nu\tau-R\lesssim\hbar}\int_0^\infty\chi_{(\nu\tau - R)\xi^{\frac12}\gtrsim 1} e^{\pm i(\nu\tau - R)\xi^{\frac12}}\cdot e^{\mp i\hbar^{-1}\left[\rho(x;\alpha,\hbar) - y(\alpha;\hbar) - \rho(x_{\sigma'};\alpha',\hbar)\right]}\cdot \Psi_{\tau,\sigma}(x;\alpha,\hbar)\,d\xi
	\end{split}\end{equation}
	Split the first term on the right into the following 
	\begin{align*}
	&\chi_{\nu\tau-R\lesssim\hbar}\int_0^\infty\chi_{(\nu\tau - R)\xi^{\frac12}\lesssim1} e^{\pm i(\nu\tau - R)\xi^{\frac12}}\cdot e^{\mp i\hbar^{-1}\left[\rho(x;\alpha,\hbar) - y(\alpha;\hbar) - \rho(x_{\sigma'};\alpha',\hbar)\right]}\cdot \Psi_{\tau,\sigma}(x;\alpha,\hbar)\,d\xi\\
	& = \chi_{\nu\tau-R\lesssim\hbar}\int_0^\infty\chi_{(\nu\tau - R)\xi^{\frac12}\lesssim1}\cdot e^{\mp i\hbar^{-1}\left[\rho(x_{\tau};\alpha,\hbar) - y(\alpha;\hbar) - \rho(x_{\sigma'};\alpha',\hbar)\right]}\cdot \Psi_{\tau,\sigma}(x;\alpha,\hbar)\,d\xi\\
	& + \chi_{\nu\tau-R\lesssim\hbar}\int_0^\infty\chi_{(\nu\tau - R)\xi^{\frac12}\lesssim1} \left(e^{\pm i\left[(\nu\tau - R)\xi^{\frac12} + \hbar^{-1}(\rho(x_{\tau};\alpha,\hbar) - \rho(x;\alpha,\hbar))\right]}-1\right)\\&\hspace{5cm}\cdot e^{\mp i\hbar^{-1}\left[\rho(x_{\tau};\alpha,\hbar) - y(\alpha;\hbar) - \rho(x_{\sigma'};\alpha',\hbar)\right]}\cdot \Psi_{\tau,\sigma}(x;\alpha,\hbar)\,d\xi\\
	& =:A_1+A_2.
	\end{align*}
	Then the first term $A_1$ can be split into a regular term and a explicit term via 
	\begin{align*}
	& \chi_{\nu\tau-R\lesssim\hbar}\int_0^\infty\chi_{(\nu\tau - R)\xi^{\frac12}\lesssim1}\cdot e^{\mp i\hbar^{-1}\left[\rho(x_{\tau};\alpha,\hbar) - y(\alpha;\hbar) - \rho(x_{\sigma'};\alpha',\hbar)\right]}\cdot \Psi_{\tau,\sigma}(x;\alpha,\hbar)\,d\xi\\
	& =  \chi_{\nu\tau-R\lesssim\hbar}\int_0^\infty e^{\mp i\hbar^{-1}\left[\rho(x_{\tau};\alpha,\hbar) - y(\alpha;\hbar) - \rho(x_{\sigma'};\alpha',\hbar)\right]}\cdot \Psi_{\tau,\sigma}(x;\alpha,\hbar)\,d\xi\\
	& -  \chi_{\nu\tau-R\lesssim\hbar}\int_0^\infty\chi_{(\nu\tau - R)\xi^{\frac12}\gtrsim 1}\cdot e^{\mp i\hbar^{-1}\left[\rho(x_{\tau};\alpha,\hbar) - y(\alpha;\hbar) - \rho(x_{\sigma'};\alpha',\hbar)\right]}\cdot \Psi_{\tau,\sigma}(x;\alpha,\hbar)\,d\xi.
	\end{align*}
	Here the first of the last two terms is of the required form for $f_3$, as the oscillatory phase is now independent of $R$, while the second term admits the representation required of terms of the form $f_2$ in the lemma. In fact, one can write 
	\begin{align*}
	\chi_{(\nu\tau - R)\xi^{\frac12}\gtrsim 1}\cdot \Psi_{\tau,\sigma}(x;\alpha,\hbar) = \big(\nu\tau - R\big)^{\frac12+k\nu}\cdot \big(\log(\nu\tau - R)\big)^j\cdot \frac{H(\tau,\sigma,\nu\tau - R,\xi;\hbar)}{\tau^{\frac12}},
	\end{align*}
	where the function $H$ satisfies 
	\begin{align*}
	&\Big\|\partial_{\nu\tau - R}^l\big(\int_0^\infty e^{\mp i\hbar^{-1}\left[\rho(x_{\tau};\alpha,\hbar) - y(\alpha;\hbar) - \rho(x_{\sigma'};\alpha',\hbar)\right]}\cdot H(\tau,\sigma,\nu\tau - R,\xi;\hbar)\,d\xi\big)\Big\|_{L^1_{d\sigma}}\lesssim_l \hbar^{-1}\frac{\log\tau)^{N_1-j}}{\tau^{1+\nu}}\cdot(\nu\tau - R)^{-l},\,l\geq 0,\\
	&\Big\|\partial_{\nu\tau - R}^l\partial_{\tau}\big(\int_0^\infty e^{\mp i\hbar^{-1}\left[\rho(x_{\tau};\alpha,\hbar) - y(\alpha;\hbar) - \rho(x_{\sigma'};\alpha',\hbar)\right]}\cdot H(\tau,\sigma,\nu\tau - R,\xi;\hbar)\,d\xi\big)\Big\|_{L^1_{d\sigma}}\lesssim_l  \hbar^{-1}\frac{\log\tau)^{N_1-j}}{\tau^{2+\nu}}\cdot(\nu\tau - R)^{-l},\,l\geq 0.
	\end{align*}
	The term $A_2$ is also of the explicit form, i. e. can be written like the preceding term, by expressing it as 
	\begin{align*} 
	A_2 =  &(\nu\tau - R)\chi_{\nu\tau-R\lesssim\hbar}\int_0^\infty\chi_{(\nu\tau - R)\xi^{\frac12}\lesssim1} \frac{\left(e^{\pm i\left[(\nu\tau - R)\xi^{\frac12} + \hbar^{-1}\left(\rho(x_{\tau};\alpha,\hbar) - \rho(x;\alpha,\hbar)\right)\right]}-1\right)}{(\nu\tau - R)\xi^{\frac12}}\\&\hspace{7cm}\cdot e^{\mp i\hbar^{-1}\left[\rho(x_{\tau};\alpha,\hbar) - y(\alpha;\hbar) - \rho(x_{\sigma'};\alpha',\hbar)\right]}\cdot \xi^{\frac12}\Psi(x;\alpha,\hbar)\,d\xi,
	\end{align*}
	and arguing as in the proof of Lemma~\ref{lem:singFouriertiphysicalngeq2}.
	\\ 
	It remains to deal with the integral 
	\[
	\chi_{\nu\tau-R\lesssim\hbar}\int_0^\infty\chi_{(\nu\tau - R)\xi^{\frac12}\gtrsim 1} e^{\pm i(\nu\tau - R)\xi^{\frac12}}\cdot e^{\mp i\hbar^{-1}\left[\rho(x;\alpha,\hbar) - y(\alpha;\hbar) - \rho(x_{\sigma'};\alpha',\hbar)\right]}\cdot \Psi(x;\alpha,\hbar)\,d\xi,
	\]
	which we claim is of the explicit form, i.e., can be included into $f_2$. This can easily be seen by combining the term $e^{\pm i\hbar^{-1}y(\alpha;\hbar)}$ with $\Psi(x;\alpha,\hbar)$, and taking advantage of the inequality, in view of the positivity property \eqref{monotonicity rho deri},
	\begin{align*}
	\partial_{\xi^{\frac12}}\left[(\nu\tau - R)\xi^{\frac12} - \hbar^{-1}\left[\rho(x;\alpha,\hbar) - \rho(x_{\sigma'};\alpha',\hbar)\right]\right]\gtrsim \nu\tau - R,
	\end{align*}
	and invoking integration by parts with respect to $\xi^{\frac12}$. This concludes the argument for the contribution of the terms \eqref{eq:xsingin}. 
	\\
	
	{\it{The case of outgoing singularity, i.e., the contribution of \eqref{eq:xsingout1}, \eqref{eq:xsingout2}.}} Here the argument is very similar except that we have to replace the phases 
	\[
	e^{\pm i\left[\left(\nu\tau - R\right)\xi^{\frac12} + \hbar^{-1} \rho\left(x_{\sigma'};\alpha',\hbar\right)\right]}
	\]
	by 
	\[
	e^{\pm i\left[\left(\nu\tau - 2\frac{\lambda(\tau)}{\lambda(\sigma)}\nu\sigma + R\right)\xi^{\frac12} - \hbar^{-1} \rho\left(x_{\sigma'};\alpha',\hbar\right)\right]},\quad e^{\pm i\left[\nu\left(\frac{\lambda(\tau)}{\lambda(\sigma)}x + \tau\right)\xi^{\frac12} - R\xi^{\frac12} + \hbar^{-1}\rho\left(x_{\sigma\cdot\frac{\lambda(\tau)}{\lambda(\sigma)}};\alpha\cdot\frac{\lambda(\tau)}{\lambda(\sigma)},\hbar\right)\right]}
	\]
	for the contributions of  \eqref{eq:xsingout1}, \eqref{eq:xsingout2}, respectively, and on account of (recall $\sigma\leq \tau, x\geq 0$)
	\[
	\nu\tau - 2\frac{\lambda(\tau)}{\lambda(\sigma)}\nu\sigma\leq -\nu\tau,\quad \nu\left(\frac{\lambda(\tau)}{\lambda(\sigma)}x + \tau\right) - R\geq \nu\tau - R,
	\]
	we can repeat the arguments from before since all functions will be restricted to the interior of the light cone, and hence we encounter functions of the type 
	\[
	\left(\nu\tau - R + y\right)^{\frac{l}{2} + k\nu}H(\tau, \nu\tau - R+y),\quad y\geq 0,
	\]
	where $H$ has symbol behavior with respect to the second variable. In fact we claim that the resulting explicit terms are smoother than the ones obtained for $\xb_{\text{in}}$, as explained in the second part of the lemma. We briefly outline the arguments for the explicit singular terms generated by \eqref{eq:xsingout2}. In fact, consider a schematically written term (where we omit several factors which do not affect its smoothness but only its decay with respect to time)
	\[
	\chi_{\nu\tau-R\lesssim\hbar} \int_0^\infty\chi_{\lesssim 1}\left(\frac{\lambda(\tau)}{\lambda(\sigma)}x\right) \left[\nu\tau - R + \frac{\lambda(\tau)}{\lambda(\sigma)}x\right]^{\frac12+k\nu}\left(\log\left(\nu\tau - R + \frac{\lambda(\tau)}{\lambda(\sigma)}x\right)\right)^i H\left(\nu\tau - R + \frac{\lambda(\tau)}{\lambda(\sigma)}x\right)g(x)\,dx,
	\]
	where $g$ is bounded as well as in $L^1([0,\infty))$, $H(\cdot)$ is bounded and has symbol type behavior, and we have included an extra cutoff $\chi_{\lesssim 1}\left(\frac{\lambda(\tau)}{\lambda(\sigma)}x\right)$ since the contribution of \eqref{eq:xsingout2} where $\frac{\lambda(\tau)}{\lambda(\sigma)}x\gtrsim 1$ is seen to lead to terms of type $f_3$. Then splitting 
	\begin{align*}
	& \chi_{\nu\tau-R\lesssim\hbar} \int_0^\infty\chi_{\lesssim 1}\left(\frac{\lambda(\tau)}{\lambda(\sigma)}x\right) \left[\nu\tau - R + \frac{\lambda(\tau)}{\lambda(\sigma)}x\right]^{\frac12+k\nu}\left(\log\left(\nu\tau - R + \frac{\lambda(\tau)}{\lambda(\sigma)}x\right)\right)^i H\left(\nu\tau - R + \frac{\lambda(\tau)}{\lambda(\sigma)}x\right)g(x)\,dx\\
	& =  \chi_{\nu\tau-R\lesssim\hbar} \int_0^\infty\chi_{\lesssim 1}\left(\frac{\lambda(\tau)}{\lambda(\sigma)}x\right) \left[\frac{\lambda(\tau)}{\lambda(\sigma)}x\right]^{\frac12+k\nu}\left(\log\left(\frac{\lambda(\tau)}{\lambda(\sigma)}x\right)\right)^i H\left(\frac{\lambda(\tau)}{\lambda(\sigma)}x\right)g(x)\,dx\\
	& + \chi_{\nu\tau-R\lesssim\hbar} \int_0^\infty\chi_{\lesssim 1}\left(\frac{\lambda(\tau)}{\lambda(\sigma)}x\right) \left[\nu\tau - R + \frac{\lambda(\tau)}{\lambda(\sigma)}x\right]^{\frac12+k\nu}\left(\log\left(\nu\tau - R + \frac{\lambda(\tau)}{\lambda(\sigma)}x\right)\right)^i H\left(\nu\tau - R + \frac{\lambda(\tau)}{\lambda(\sigma)}x\right)g(x)\,dx\\
	& - \chi_{\nu\tau-R\lesssim\hbar} \int_0^\infty\chi_{\lesssim 1}\left(\frac{\lambda(\tau)}{\lambda(\sigma)}x\right) \left[\frac{\lambda(\tau)}{\lambda(\sigma)}x\right]^{\frac12+k\nu}\left(\log\left(\frac{\lambda(\tau)}{\lambda(\sigma)}x\right)\right)^i H\left(\frac{\lambda(\tau)}{\lambda(\sigma)}x\right)g(x)\,dx,
	\end{align*}
	the first term on the right is again seen to lead to a contribution which can be absorbed into $f_3$ (upon also taking into account the currently omitted factors ensuring the right temporal decay). The second and the third (difference) term on the right can be written as 
	\begin{align*}
	\chi_{\nu\tau-R\lesssim\hbar}(\nu\tau - R)\cdot\int_0^\infty \int_0^1 Z'\left(s(\nu\tau - R) + x\right)\cdot \chi_{\lesssim 1}\left(\frac{\lambda(\tau)}{\lambda(\sigma)}x\right)g(x)\,ds dx,
	\end{align*}
	which of course involves a derivative falling onto the term $H(\cdot)$, which will then be responsible for the derivative loss which  the `upgrade of smoothness' for the contribution of $\xb_{\text{out}}$ entails. Write the preceding integral as 
	\begin{align*}
	&\chi_{\nu\tau-R\lesssim\hbar}(\nu\tau - R)\cdot\int_0^\infty \int_0^1 Z'\left(s(\nu\tau - R) + x\right)\cdot \chi_{\lesssim 1}\left(\frac{\lambda(\tau)}{\lambda(\sigma)}x\right)g(x)\,ds dx\\
	& = \chi_{\nu\tau-R\lesssim\hbar}\cdot(\nu\tau - R)\cdot\int_0^{\nu\tau - R} \int_0^1 Z'\left(s(\nu\tau - R) + x\right)\cdot \chi_{\lesssim 1}\left(\frac{\lambda(\tau)}{\lambda(\sigma)}x\right)g(x)\,ds dx\\
	& + \chi_{\nu\tau-R\lesssim\hbar}\cdot(\nu\tau - R)\cdot\int_{\nu\tau - R}^\infty\int_0^1 Z'\left(s(\nu\tau - R) + x\right)\cdot \chi_{\lesssim 1}\left(\frac{\lambda(\tau)}{\lambda(\sigma)}x\right)g(x)\,ds dx\\
	\end{align*}
	Observe that 
	\[
	Z'(z) = z^{k\nu-\frac12}(\log z)^i\cdot \tilde{H}(z), 
	\]
	where $\tilde{H}$ has properties like $H$ in the region $z>0$. It easily follows that the first integral on the right admits an explicit representation 
	\begin{align*}
	&\chi_{\nu\tau-R\lesssim\hbar}\cdot(\nu\tau - R)\cdot\int_0^{\nu\tau - R} \int_0^1 Z'\left(s(\nu\tau - R) + x\right)\cdot \chi_{\lesssim 1}\left(\frac{\lambda(\tau)}{\lambda(\sigma)}x\right)g(x)\,ds dx\\& = \left(\nu\tau - R\right)^{\frac32+k\nu}\log\left(\nu\tau -R\right)^i \cdot H_1\left(\tau, \sigma, \nu\tau - R\right), 
	\end{align*}
	where $H_1$ is bounded and has symbol type bound with respect to all its variables. For the second integral above over the range $x\in [\nu\tau - R,\infty)$, we use one more splitting of a similar kind: 
	\begin{align*}
	Z'\left(s(\nu\tau - R) + x\right) = Z'(x) + s(\nu\tau - R)\cdot \int_0^1Z''\left(s_1s(\nu\tau - R) + x\right)\,ds_1,
	\end{align*}
	which then yields
	\begin{align*}
	&\chi_{\nu\tau-R\lesssim\hbar}\cdot(\nu\tau - R)\cdot\int_{\nu\tau - R}^\infty \int_0^1 Z'\left(s(\nu\tau - R) + x\right)\cdot \chi_{\lesssim 1}\left(\frac{\lambda(\tau)}{\lambda(\sigma)}x\right)g(x)\,ds dx\\&
	=  \chi_{\nu\tau-R\lesssim\hbar}\cdot(\nu\tau - R)\cdot\left[\int_0^\infty \chi_{\lesssim 1}\left(\frac{\lambda(\tau)}{\lambda(\sigma)}x\right)Z'(x)g(x)\,dx - \int_0^{\nu\tau - R} \chi_{\lesssim 1}\left(\frac{\lambda(\tau)}{\lambda(\sigma)}x\right)Z'(x)g(x)\,dx\right]\\
	& + \chi_{\nu\tau-R\lesssim\hbar}\cdot(\nu\tau - R)^2\cdot\int_{\nu\tau - R}^\infty\int_0^1\int_0^1s Z''\left(s_1s(\nu\tau - R) + x\right)g(x)\,ds_1 ds dx
	\end{align*}
	Here, the first term on the right is the difference of a function which can be included into $f_3$ and a function of the explicit type, while for the second term on the right, due to the inequality five derivatives of regularity for the function $H_{k,l,j}$ in the statement of the lemma, compared to the functions $G_{k,l,j}$. 
	\[
	\left|Z''(z)\right|\lesssim z^{k\nu-\frac32}\left|\log z\right|^i,
	\]
	it will be of the explicit form 
	\begin{align*}
	&\chi_{\nu\tau-R\lesssim\hbar}\cdot(\nu\tau - R)^2\cdot\int_{\nu\tau - R}^\infty\int_0^1\int_0^1s Z''\left(s_1s(\nu\tau - R) + x\right)g(x)\,ds_1 ds dx \\= &\left(\nu\tau - R\right)^{\frac32+k\nu}\log\left(\nu\tau -R\right)^i \cdot H_2(\tau, \sigma, \nu\tau - R)
	\end{align*}
	as long as $k\nu<\frac12$. If $k\nu\geq \frac12$, then we repeat the above process at most four more times (since $k\nu\leq 5$), which loses at most five derivatives. 
	\\
	The contribution of the `seed outgoing singular terms' \eqref{eq:xsingout1} is handled similarly by replacing the role of $x$ by $\sigma$. 
\end{proof}

\begin{remark}\label{rem:outgoing smoothing} 
	We will only take advantage of the additional smoothness of the physical incarnation of the outgoing singular terms for the contribution of the seed outgoing part with $l = 0$, where the coefficients $F_{l,k,i}$ do not depend on the frequency variable, and hence no smoothness loss will be incurred. 
\end{remark}

The adequacy of the physical space expansions of the singular part of admissible functions detailed in Lemma~\ref{lem:singFouriertiphysicalngeq2adm}  in the context of our null-form nonlinearities follows from the next lemma, which makes the connection to  Lemma~\ref{lem:singPhysicaltoFourierngeq2} :
\begin{lemma}\label{lem:basicnullformforsingularinputs} Assume that the functions $f_{2l}, \tilde{f}_{2l'}$ are as in Lemma~\ref{lem:singFouriertiphysicalngeq2adm}  where the coefficient functions 
\begin{align*}
G_{k,l,j}, \tilde{G}_{k',l',j'}
\end{align*}
satisfy (with $\delta>0$)
\begin{align*}
&\big|x^{k_2}\cdot \partial_x^{k_2}G_{k,l,j}\big|\lesssim (\log\tau)^{N_1-j+1}\cdot \tau^{-1-\nu},\,k_2\leq 20 - l - \lfloor 2\nu k\rfloor + \lfloor\frac{l}{2}+\nu k\rfloor,\\
&\big\|x^{k_3+\delta_1}\partial_{x}^{k_3}G_{k,l,j}\big\|_{\dot{C}^{\delta_1}}\lesssim (\log\tau)^{N_1-j+1}\cdot \tau^{-1-\nu},\,k_3 = 20 - l - \lfloor 2\nu k\rfloor + \lfloor\frac{l}{2}+\nu k\rfloor,\,\delta_1 = \frac{l}{2}+\nu k - \lfloor\frac{l}{2}+\nu k\rfloor + \delta,
\end{align*}
and analogously for $\tilde{G}_{k',l',j'}$. Then letting $(l_*, k_*,j_*)\in \{(l, k,j),\,(l',k',j')\}$ be such that $l_* + \lfloor 2\nu k_*\rfloor = \max\{l+\lfloor 2\nu k\rfloor ,\,l'+\lfloor 2\nu k'\rfloor\}$, and setting 
\begin{align*}
&\phi: = \chi_{|\nu\tau - R|\lesssim \hbar}\cdot \frac{G_{k,l,j}(\tau, \nu\tau - R)}{\tau^{\frac12}}\cdot \hbar^{-\frac{l+1}{2}}[\nu\tau - R]^{\frac{l}{2}+k\nu}\big(\log(\nu\tau - R)\big)^j\\
&\tilde{\phi}: = \chi_{|\nu\tau - R|\lesssim \hbar}\cdot \frac{G_{k,l,j}(\tau, \nu\tau - R)}{\tau^{\frac12}}\cdot \hbar^{-\frac{l'+1}{2}}[\nu\tau - R]^{\frac{l'}{2}+k'\nu}\big(\log(\nu\tau - R)\big)^{j'},\\
\end{align*}
we have that 
\begin{align*}
H: = \big(\partial_{\tau} + \frac{\lambda_{\tau}}{\lambda}R\partial_R - \partial_R\big)\phi\cdot \big(\partial_{\tau} + \frac{\lambda_{\tau}}{\lambda}R\partial_R + \partial_R\big)\phi'
\end{align*}
can be written in the form $f_{2\tilde{l}}$, $\tilde{l}\geq 0$, as in the statement of Lemma~\ref{lem:singPhysicaltoFourierngeq2} , i. e. 
\begin{align*}
H =  \tilde{\chi}_{|\nu\tau - R|\lesssim \hbar}\cdot \frac{G_{*k_*,l_*, j_*}(\tau, \nu\tau - R)}{\tau^{\frac12}}\cdot \hbar^{-\frac{l_*+1}{2}}[\nu\tau - R]^{\frac{l_*}{2}+k_*\nu}\big(\log(\nu\tau - R)\big)^{j_*}
\end{align*}
and we have the bounds 
\begin{align*}
&\big|x^{k_2}\cdot \partial_x^{k_2}G_{k_*,l_*,j_*}\big|\lesssim \hbar^{-1}(\log\tau)^{2N_1-j_*+2}\cdot \tau^{-3-2\nu},\,k_2\leq 19 - l_* - \lfloor 2\nu k_*\rfloor + \lfloor\frac{l_*}{2}+\nu k_*\rfloor,\\
&\big\|x^{k_3+\delta_1}\partial_{x}^{k_3}G_{k_*,l_*,j_*}\big\|_{\dot{C}^{\delta_1}}\lesssim \hbar^{-1}(\log\tau)^{2N_1-j_*+2}\cdot \tau^{-3-2\nu},\,k_3 = 19 - l_* - \lfloor 2\nu k_*\rfloor + \lfloor\frac{l_*}{2}+\nu k_*\rfloor,\\&\hspace{4cm}\delta_1 = \frac{l_*}{2}+\nu k_* - \lfloor\frac{l_*}{2}+\nu k_*\rfloor + \delta.
\end{align*}
\end{lemma}
\begin{proof} It suffices to observe that 
\begin{align*}
[\frac{\nu\tau - R}{\hbar}]^{\frac{l}{2}+k\nu}\cdot [\frac{\nu\tau - R}{\hbar}]^{\frac{l'}{2}+k'\nu - 1} = [\frac{\nu\tau - R}{\hbar}]^{\frac{l_*-1}{2}+k_*\nu}\cdot [\frac{\nu\tau - R}{\hbar}]^{\frac{l+l'-l_*-1}{2}+(k+k' - k_*)\nu},
\end{align*}
and we have $l+l'-l_*-1\geq 0$, $k+k' - k_*>0$. Furthermore, the 'good' derivative 
\[
\partial_{\tau} +  \frac{\lambda_{\tau}}{\lambda}R\partial_R - \partial_R = \partial_{\tau} + \nu\partial_R + (1+\nu^{-1})\cdot\frac{R - \nu\tau}{\tau}\partial_R
\]
does not alter the order of vanishing at $R = \nu\tau$ of either $\phi$ or $\phi'$, but does result in applying a scaling operator $(\nu\tau - R)\partial_{\nu\tau - R}$ to $G(\nu\tau - R)$, reducing its smoothness, which is reflected by the final derivative bounds. We note that the added factor of $\hbar^{-1}$ in the final derivative bounds arises since an operator $\partial_{\tau}$ may apply to the first variable $\tau$ in $G(\tau, \nu\tau - R)$. We note that this factor does not cause difficulties since the null-form singularity is of non-principal source admissibly singular type, and hence we can absorb it, see Lemma~\ref{lem:singPhysicaltoFourierngeq2}. 
\end{proof}

\subsection{Operations on functions with admissible singular part}

Our definition of admissible singular part is chosen to be flexible enough that important operations, such as frequency localizations as well as the transference operator, essentially preserve such functions. Furthermore, restrictions of admissibly singular functions away from the light cone lead to functions in $S_0^{\hbar}$ with good temporal decay.

\subsubsection{Restriction away from the light cone}

\begin{lemma}\label{lem:admsingawayfromshock} Let $\overline{y} = \overline{y}(\tau,\xi)$ be an admissibly singular function at angular momentum $|n|\geq 2$, and let 
\[
f(\tau, R) = \int_0^\infty \phi_n(R;\xi)\cdot \overline{y}(\tau,\xi)\cdot\rho_n(\xi)\,d\xi
\]
its distorted inverse Fourier transform. Then we have the bound 
\begin{align*}
\tau^4\cdot\big\|\chi_{R\ll\nu\tau}\cdot f(\tau, R)\big\|_{\tilde{S}_0^{\hbar}}\lesssim \big\|\overline{y}\big\|_{adm}.
\end{align*}
\end{lemma}
\begin{proof} We rely on Lemma~\ref{lem:singFouriertiphysicalngeq2adm}. Expand $f(\tau, R)$ as in that lemma, and notice that only the function $f_1(\tau, R)$ contributes in the region $R\ll\nu\tau$. Then we have 
\begin{align*}
&\big\|\chi_{R\ll\nu\tau}\cdot f(\tau, R)\big\|_{\tilde{S}_0^{\hbar}} = \big\|\langle \phi_n(R;\xi),\,\chi_{R\ll\nu\tau}\cdot f_1(\tau, R)\rangle_{L^2_{R\,dR}}\big\|_{S_0^{\hbar}}\\
&\lesssim \hbar\cdot\big\|\partial_R\langle \partial_R\rangle^4\big(\chi_{R\ll\nu\tau}\cdot f_1(\tau, R)\big)\big\|_{L^2_{R\,dR}}. 
\end{align*}
The desired bound is then a direct consequence of the pointwise bound for $f_1$ in  Lemma~\ref{lem:singFouriertiphysicalngeq2adm} .

\end{proof}

\subsubsection{The effect of the transference operator}

Recall from Proposition \ref{prop: K operator} the precise form of the transference operator at angular momentum $n, |n|\geq 2$, and consisting of a diagonal and an off-diagonal part. Two key structural properties ensure that this preserves an admissible singular part: On the one hand, the spectral density $\rho_n(\xi)$ admits an asymptotic expansion toward $\xi = +\infty$ of Hankel type. On the other hand the off-diagonal part $\mathcal{K}_{\hbar}^{(0)}$ improves decay toward $\xi = +\infty$ and its kernel has further remarkable regularity properties uniformly in $n$. The latter property means that the off-diagonal transference operator will send the principal part of the ingoing singular term into its connecting part, which is less rigid structurally but decays better at large frequencies. 
The following lemma gives the basic bounds for the Hilbert transform concerning the weighted H\"older type norms we shall work with: 
\begin{lemma}\label{lem:weightedHolder1} Introduce the norm (for $\delta\in (0,1)$)
	\[
	\left\|f\right\|_{wC^{\delta}}: = \left\|x f(x)\right\|_{L^\infty([0,\infty))} + \sup_{\lambda>0}\lambda^{1+\delta}\left\|\chi_{x\simeq\lambda}f\right\|_{\dot{C}^{\delta}}
	\]
	Assume that $f$ is supported at $\eta\simeq\lambda\gtrsim \hbar^{-2}$, and consider the function\footnote{The factor $\langle\xi^{\frac12} - \eta^{\frac12}\rangle^\gamma$ is motivated by the off-diagonal decay of the kernel $F_n(\xi,\eta)$ in Prop.~\ref{prop: K operator}.} 
	\begin{align*}
	\left(\mathcal{H}_{\text{twisted}}^{\pm}f\right)(\xi;\tau,\sigma'): = \int_0^\infty \frac{e^{\pm i\left[\nu\tau\eta^{\frac12} + \hbar^{-1}\rho\left(x_{\sigma'}^{(\eta)};\beta',\hbar\right)\right]}}{(\xi - \eta)\langle\xi^{\frac12} - \eta^{\frac12}\rangle^\gamma}\,f(\eta)\,d\eta,\quad\gamma>0, 
	\end{align*}
	where we use the notation $\beta' = \hbar\eta^{\frac12}\cdot\frac{\lambda(\tau)}{\lambda(\sigma)},\quad x_{\sigma'}^{(\eta)} = \hbar\eta^{\frac12}\cdot\nu\sigma\cdot\frac{\lambda(\tau)}{\lambda(\sigma)}$. Then for any $\delta'<\delta$ we have the relation 
	\begin{align*}
	\left(\mathcal{H}_{\text{twisted}}^{\pm}f\right)(\xi;\tau,\sigma') = e^{\pm i\left[\nu\tau\xi^{\frac12} + \hbar^{-1}\rho\left(x_{\sigma'}^{(\xi)};\alpha',\hbar\right)\right]}\cdot g_{\pm}(\xi;\tau, \sigma'),\quad \left\|\chi_{\xi\simeq\lambda}g(\xi;\tau, \sigma')\right\|_{wC^{\delta'}_{\xi}}\lesssim_{\delta',\delta,\gamma} \left\|f\right\|_{wC^{\delta}},
	\end{align*}
	where $\left\|\cdot\right\|_{wC^{\delta'}_{\xi}}$ indicates that the norm is with respect to the variable $\xi$, we use the notations $\alpha'=\hbar\xi^{\frac12}\cdot\frac{\lambda(\tau)}{\lambda(\sigma)}$ and $x_{\sigma'}^{(\xi)}= \hbar\xi^{\frac12}\cdot\nu\sigma\cdot\frac{\lambda(\tau)}{\lambda(\sigma)}$, and the bound is uniform in $\hbar, \tau, \sigma'$. 
\end{lemma}
\begin{proof}
	By definition we have 
	\begin{align*}
	g_{\pm}(\xi;\tau, \sigma') = \int_0^\infty \frac{e^{i\left[\nu\tau\left(\eta^{\frac12} - \xi^{\frac12}\right) + \hbar^{-1}\left[\rho(x_{\sigma'}^{(\eta)};\beta',\hbar) - \rho(x_{\sigma'}^{(\xi)};\alpha',\hbar)\right]\right]}}{(\xi - \eta)\langle\xi^{\frac12} - \eta^{\frac12}\rangle^\gamma}f(\eta)\,d\eta
	\end{align*}
	Introduce the variable $\tilde{\eta} = \nu\tau(\eta^{\frac12} - \xi^{\frac12})$, and write the phase as (we suppress the $\hbar$-dependence in the notation since $\hbar\ll 1$ is fixed)
	\begin{equation}\label{eq:Psixietatausigmaprime}\begin{split}
	&\Psi(\xi,\eta,\tau, \sigma'): = \nu\tau(\eta^{\frac12} - \xi^{\frac12}) + \hbar^{-1}\left[\rho(x_{\sigma'}^{(\eta)};\beta',\hbar) - \rho(x_{\sigma'}^{(\xi)};\alpha',\hbar)\right]\\&
	= \tilde{\eta} + \hbar^{-1}\int_0^1 \left[\rho_x\left(sx_{\sigma'}^{(\eta)} + (1-s)x_{\sigma'}^{(\xi)};s\beta' + (1-s)\alpha',\hbar\right)\cdot\hbar\frac{\sigma'}{\tau}\tilde{\eta}\right.
	\\&\hspace{5cm}+ \left.\rho_\alpha\left(sx_{\sigma'}^{(\eta)} + (1-s)x_{\sigma'}^{(\xi)};s\beta' + (1-s)\alpha',\hbar\right)\cdot\hbar\frac{\tilde{\eta}}{\nu\tau}\cdot\frac{\lambda(\tau)}{\lambda(\sigma)}\right]\,ds
	\end{split}\end{equation}
	By the bounds on $\rho_{x}$ and $\rho_{\alpha}$ in Lemma \ref{lem: Lemma 3.4 CDST}, if we pick $\xi\simeq\eta$, the phase is in effect bounded by $\left|\Psi(\xi,\eta,\tau, \sigma')\right|\lesssim \tilde{\eta}$ under those conditions. Write $f(x^2) = g(x)$, restricted to the positive real axis, and switch variables in the integral  
	\begin{align*}
	g_{\pm}(\xi;\tau, \sigma') = \int_{-\infty}^{\infty} \frac{e^{i\Psi(\xi,\eta,\tau, \sigma')}}{\tilde{\eta}\left\langle\frac{\tilde{\eta}}{\nu\tau}\right\rangle^{\gamma}}\cdot \frac{2\eta^{\frac12}}{\xi^{\frac12}+\eta^{\frac12}}\cdot g\left(\frac{\tilde{\eta}}{\nu\tau} + \xi^{\frac12}\right)\,d\tilde{\eta},
	\end{align*}
	where of course $\eta^{\frac12} = \frac{\tilde{\eta}}{\nu\tau} + \xi^{\frac12}$. In fact, restricting $\xi\simeq \lambda$, the integration limits may be set to be $\pm \lambda^{\frac12}\tau$. We now prove the two bounds required for $g$: 
\\

	{\it{The weighted $L^\infty$-bound.}} Denote by $P_{>\kappa},\,\kappa\in (0,\infty)$ a Littlewood-Paley type frequency cutoff to frequencies $\gtrsim \kappa$, and which can be realized by convolution with the function $\kappa\hat{\chi}(\eta\kappa)$, where $\chi$ is a smooth function supported on $(-\infty, -\frac12]\cup[\frac12,\infty)$ and identically $1$ on $[1,\infty)$. Note that the Fourier transform of $\chi$ contains a delta measure. Therefore instead of $P_{>\kappa}$ we consider a localization $P_{\simeq 2^{j}\kappa}$ with $j\in \bbZ^{+}$, by assuming that $\chi$ is supported on the interval $[1, 2)$. Then we have 
	\begin{equation}\label{eq:freqcutoff}
	P_{\simeq 2^{j}\kappa}g(\eta) = \int_{-\infty}^\infty 2^{j}\kappa\hat{\chi}\left((\eta-\zeta)2^{j}\kappa\right)g(\zeta)\,d\zeta =  \int_{-\infty}^\infty 2^{j}\kappa\hat{\chi}\left((\eta-\zeta)2^{j}\kappa\right)\left[g(\zeta) - g(\eta)\right]\,d\zeta,
	\end{equation}
	Here we used the fact that $\chi(\cdot)$ is supported away from the origin while the delta measure is supported at the origin. Setting $\kappa = \left(\frac{\tau^{2}}{\lambda}\right)^{\frac{1}{4}}$ and taking advantage of the rapid decay of the function $\hat{\chi}$ as well as the bound 
	\[
	\left|g(x_1) - g(x_2)\right| = \left|f(x_1^2) - f(x_2^2)\right|\lesssim \lambda^{-1-\delta}\left\|f\right\|_{wC^{\delta}}\cdot \left|x_1^2 - x_2^2\right|^{\delta}\lesssim  \lambda^{-1-\frac{\delta}{2}}\left\|f\right\|_{wC^{\delta}}\cdot \left|x_1 - x_2\right|^{\delta}
	\]
	provided $x_{1,2}\simeq\lambda^{\frac12}$, we infer 
	\begin{align*}
		\left|P_{\simeq2^{j}\left(\frac{\tau^{2}}{\lambda}\right)^{\frac{1}{4}}}g(\eta) \right|\lesssim 2^{-j\delta} \tau^{-\frac{\delta}{2}}\lambda^{-1-\frac{\delta}{4}}\cdot \left\|f\right\|_{wC^{\delta}},
	\end{align*}
which further implies
	\begin{align*}
	\left|P_{>\left(\frac{\tau^{2}}{\lambda}\right)^{\frac{1}{4}}}g(\eta) \right|\lesssim  \tau^{-\frac{\delta}{2}}\lambda^{-1-\frac{\delta}{4}}\cdot \left\|f\right\|_{wC^{\delta}},
	\end{align*}
	and so we infer the bound 
	\begin{align*}
	\left| \int_{-\infty}^{\infty}\psi_{|\tilde{\eta}|\gtrsim 1}\psi_{\eta\simeq\lambda} \frac{e^{i\Psi(\xi,\eta,\tau, \sigma')}}{\tilde{\eta}\langle\frac{\tilde{\eta}}{\nu\tau}\rangle^{\gamma}}\cdot \frac{2\eta^{\frac12}}{\xi^{\frac12}+\eta^{\frac12}}\cdot P_{>\left(\frac{\tau^{2}}{\lambda}\right)^{\frac{1}{4}}}g\left(\frac{\tilde{\eta}}{\nu\tau} + \xi^{\frac12}\right)\,d\tilde{\eta}\right|\lesssim\lambda^{-1} \frac{\log(\lambda^{\frac12}\tau)}{ \tau^{\frac{\delta}{2}}\lambda^{\frac{\delta}{4}}}\cdot  \left\|f\right\|_{wC^{\delta}},
	\end{align*}
	where the extra localizer $\psi_{\eta\simeq\lambda}$ ensures the correct support of the integral after application of the frequency cutoff. 
	\\
	In order to handle the remaining integral with the additional restriction $|\tilde{\eta}|\gtrsim 1$ and involving $P_{\leq\left(\frac{\tau^{2}}{\lambda}\right)^{\frac{1}{4}}}g$, we perform integration by parts with respect to $\tilde{\eta}$, using 
	\[
	\left|\partial_{\tilde{\eta}}\Psi(\xi,\eta,\tau, \sigma')\right|\gtrsim 1,
	\]
	as well as straightforward higher derivative bounds with respect to $\tilde{\eta}$. On the other hand, using the crude $L^{\infty}$-bound, we have
	\begin{align*}
		\left|P_{\leq\left(\frac{\tau^{2}}{\lambda}\right)^{\frac{1}{4}}}g(\eta)\right|\lesssim \lambda^{-1}\|f\|_{wC^{\delta}}.
	\end{align*}
	Then if the derivative falls on $P_{\leq\left(\frac{\tau^{2}}{\lambda}\right)^{\frac{1}{4}}}g\left(\frac{\tilde{\eta}}{\nu\tau} + \xi^{\frac12}\right)$, we gain a factor 
	\[
	\left(\frac{\tau^{2}}{\lambda}\right)^{\frac{1}{4}}\cdot \tau^{-1},
	\]
	which, together with the bound $\lambda^{-1}\|f\|_{wC^{\delta}}$, gives a bound
	\begin{align*}
		\left(\tau^{\frac12}\lambda^{\frac14}\right)^{-1}\cdot\lambda^{-1}\|f\|_{wC^{\delta}}
	\end{align*}
	 This again beats the loss of $\log(\lambda^{\frac12}\tau)$. While if the derivative falls on $\tilde{\eta}^{-1}\left(\partial_{\tilde{\eta}}\Psi(\xi,\eta,\tau, \sigma')\right)^{-1}\cdot \left(\frac{\tilde{\eta}}{\nu\tau} + 2\xi^{\frac12}\right)$, we gain $\tilde{\eta}^{-1}$, whence forcing integrability of the expression. 
	 \\
	 This reduces things to bounding the integral 
	 \begin{align*}
	 &\int_{-\infty}^{\infty}\psi_{|\tilde{\eta}|\lesssim1} \frac{e^{i\Psi\left(\xi,\eta,\tau, \sigma'\right)}}{\tilde{\eta}\left\langle\frac{\tilde{\eta}}{\tau}\right\rangle^{\gamma}}\cdot \frac{2\eta^{\frac12}}{\xi^{\frac12}+\eta^{\frac12}}\cdot g\left(\frac{\tilde{\eta}}{\nu\tau} + \xi^{\frac12}\right)\,d\tilde{\eta}\\
	 & =\sum_{\pm}(-1)^{\pm} \int_{0}^{\infty}\psi_{|\tilde{\eta}|\lesssim1} \frac{e^{i\Psi\left(\xi,\eta(\pm\tilde{\eta},\xi,\tau),\tau, \sigma'\right)}}{\tilde{\eta}\left\langle\frac{\tilde{\eta}}{\tau}\right\rangle^{\gamma}}\cdot \frac{2\eta^{\frac12}\left(\pm\tilde{\eta},\xi,\tau\right)}{\xi^{\frac12}+\eta^{\frac12}\left(\pm\tilde{\eta},\xi,\tau\right)}\cdot g\left(\frac{\pm\tilde{\eta}}{\nu\tau} + \xi^{\frac12}\right)\,d\tilde{\eta},
	 \end{align*}
	 where we have made explicit the fact that $\eta$ is a function of $\tilde{\eta},\xi, \tau$. But here differencing easily leads to gains of factors $\lesssim \tilde{\eta}^{\delta}$, which counteracts the singular term $\tilde{\eta}^{-1}$, which easily concludes the weighted $L^\infty$-bound. Here we also used the fact that the kernel $\tilde{\eta}$ is odd and the other factors in the integral except $g\left(\frac{\pm\tilde{\eta}}{\nu\tau} + \xi^{\frac12}\right)$ are smooth in $\tilde{\eta}$.
	 \\
	 
	 {\it{H\"older differencing bound}}. Fix two values $\xi_1, \xi_2\in \R_+, \xi_1\simeq \xi_2\simeq \lambda$, and set $l: = \left|\xi_1^{\frac12} - \xi_2^{\frac12}\right|$. As before set $\tilde{\eta} = \nu\tau\left(\eta^{\frac12} - \xi^{\frac12}\right)$, and decompose the integrals giving $\left(\mathcal{H}_{\text{twisted}}^{\pm}f\right)(\xi_j;\tau,\sigma'),\,j = 1,2$, into (where for now $\xi = \xi_{1,2}$)
	 \begin{equation}\label{eq:smalllargetildeetasplit}\begin{split}
	 &\int_{-\infty}^{\infty}\psi_{|\tilde{\eta}|\lesssim1} \frac{e^{i\Psi(\xi,\eta,\tau, \sigma')}}{\tilde{\eta}\left\langle\frac{|\tilde{\eta}|}{\tau}\right\rangle^{\gamma}}\cdot \frac{2\eta^{\frac12}}{\xi^{\frac12}+\eta^{\frac12}}\cdot g\left(\frac{\tilde{\eta}}{\nu\tau} + \xi^{\frac12}\right)\,d\tilde{\eta}\\
	 & + \int_{-\infty}^{\infty}\psi_{|\tilde{\eta}|\gtrsim 1} \frac{e^{i\Psi(\xi,\eta,\tau, \sigma')}}{\tilde{\eta}\left\langle\frac{|\tilde{\eta}|}{\tau}\right\rangle^{\gamma}}\cdot \frac{2\eta^{\frac12}}{\xi^{\frac12}+\eta^{\frac12}}\cdot g\left(\frac{\tilde{\eta}}{\nu\tau} + \xi^{\frac12}\right)\,d\tilde{\eta}\\
	 \end{split}\end{equation}
	 {\it{First integral above, small $\tilde{\eta}$.}} Due the  bound $\left|\Psi(\xi,\eta,\tau, \sigma')\right|\lesssim |\tilde{\eta}|\lesssim 1$, we can split the exponential into 
	 \[
	 e^{i\Psi(\xi,\eta,\tau, \sigma')} = 1 + (e^{i\Psi(\xi,\eta,\tau, \sigma')} - 1),\quad \left|(e^{i\Psi(\xi,\eta,\tau, \sigma')} - 1)\right|\lesssim |\tilde{\eta}|, 
	 \]
	 and it is also easy to check that under our support conditions $\xi_1\simeq \xi_2\simeq\eta\simeq\lambda\gtrsim \hbar^{-2}$ we have 
	 \[
	 \left|\partial_{\xi}e^{i\Psi(\xi,\eta,\tau, \sigma')}\right|\lesssim \xi^{-1}|\tilde{\eta}|,
	 \]
	 which furnishes the desired bound for the modified contribution with the exponential replaced by $e^{i\Psi(\xi,\eta,\tau, \sigma')} - 1$: 
	 \begin{align*}
	 \left|\sum_{j = 1,2}(-1)^j\int_{-\infty}^{\infty}\psi_{|\tilde{\eta}|\lesssim1} \frac{e^{i\Psi(\xi_j,\eta,\tau, \sigma')}-1}{\tilde{\eta}\left\langle\frac{|\tilde{\eta}|}{\tau}\right\rangle^{\gamma}}\cdot \frac{2\eta^{\frac12}}{\xi_j^{\frac12}+\eta^{\frac12}}\cdot g\left(\frac{\tilde{\eta}}{\nu\tau} + \xi_j^{\frac12}\right)\,d\tilde{\eta}\right|\lesssim &\lambda^{-1-\frac{\delta}{2}}l^{\delta}\|f\|_{wC^{\delta}}\\
	 \lesssim& \lambda^{-1-\delta}|\xi_1 - \xi_2|^{\delta}\|f\|_{wC^{\delta}}. 
	 \end{align*}
	 This allows us to reduce the case of small $\tilde{\eta}$ to the following integral, where we may re-arrange things if we pick the cutoff $\psi$ symmetrically around the origin: 
	 \begin{align*}
	 \sum_{j = 1,2}(-1)^j\int_{0}^{\infty}\psi_{|\tilde{\eta}|\lesssim1} \frac{1}{\left\langle\frac{|\tilde{\eta}|}{\tau}\right\rangle^{\gamma}}\cdot \frac{2\eta^{\frac12}}{\xi_j^{\frac12}+\eta^{\frac12}}\cdot \frac{g\left(\frac{\tilde{\eta}}{\nu\tau} + \xi_j^{\frac12}\right) - g\left(\frac{-\tilde{\eta}}{\nu\tau} + \xi_j^{\frac12}\right)}{\tilde{\eta}}\,d\tilde{\eta}.
	 \end{align*}
	 We split the above integrals further by including smooth cutoffs $\chi_{|\tilde{\eta}|\lesssim \tau l},  \chi_{|\tilde{\eta}|\gtrsim \tau l}$. Including the former, we take advantage of the bound 
	 \[
	 \left|g\left(\frac{\tilde{\eta}}{\nu\tau} + \xi_j^{\frac12}\right) - g\left(\frac{-\tilde{\eta}}{\nu\tau} + \xi_j^{\frac12}\right)\right|\lesssim \left|\frac{\tilde{\eta}}{\nu\tau}\right|^{\delta}\cdot \lambda^{-1-\frac{\delta}{2}}\cdot \left\|f\right\|_{wC^{\delta}},
	 \]
	 which leads to the desired bound upon integration over $0\leq \tilde{\eta}\lesssim \tau l$: 
	 \begin{align*}
	 &\left|\int_{0}^{\infty}\psi_{ |\tilde{\eta}|\lesssim \tau l} \frac{1}{\left\langle\frac{|\tilde{\eta}|}{\tau}\right\rangle^{\gamma}}\cdot \frac{2\eta^{\frac12}}{\xi_j^{\frac12}+\eta^{\frac12}}\cdot \frac{g\left(\frac{\tilde{\eta}}{\nu\tau} + \xi_j^{\frac12}\right) - g\left(\frac{-\tilde{\eta}}{\nu\tau} + \xi_j^{\frac12}\right)}{\tilde{\eta}}\,d\tilde{\eta}\right|\\
	 &\lesssim \left(\frac{\tau l}{\tau}\right)^{\delta}\cdot\lambda^{-1-\frac{\delta}{2}}\cdot \left\|f\right\|_{wC^{\delta}} = l^{\delta}\cdot\lambda^{-1-\frac{\delta}{2}}\cdot \left\|f\right\|_{wC^{\delta}}.
	 \end{align*}
	 On the other hand, in the regime $|\tilde{\eta}|\gtrsim \tau l$, the difference structure due to the simmation over $j$ will become important. Split 
	 \[
	 g\left(\frac{\tilde{\eta}}{\nu\tau} + \xi_j^{\frac12}\right) = P_{< l^{-1}}g\left(\frac{\tilde{\eta}}{\nu\tau} + \xi_j^{\frac12}\right) + P_{\geq l^{-1}}g\left(\frac{\tilde{\eta}}{\nu\tau} + \xi_j^{\frac12}\right)
	 \]
	 where the subscripts denote frequency cutoffs defined in analogy to \eqref{eq:freqcutoff}. Then arguing as earlier, we get the bound 
	 \[
	 \left|P_{\geq l^{-1}}g\left(\frac{\tilde{\eta}}{\nu\tau} + \xi_j^{\frac12}\right)\right|\lesssim l^\delta\cdot \lambda^{-1-\frac{\delta}{2}}\cdot \left\|f\right\|_{wC^{\delta}},.
	 \]
	 Moreover, using orthogonality (or Plancherel identity) and the fact that $P_{\geq l^{-1}}g\left(\frac{\tilde{\eta}}{\nu\tau} + \xi_j^{\frac12}\right) = P^{(\tilde{\eta})}_{\geq l^{-1}\left(\nu\tau\right)^{-1}}\left(g\left(\frac{\tilde{\eta}}{\nu\tau} + \xi_j^{\frac12}\right)\right)$ where $P^{(\tilde{\eta})}_{\geq l^{-1}\left(\nu\tau\right)^{-1}}$ means the frequency cutoff applied to the following expression interpreted as a function of $\tilde{\eta}$, we infer 
	 \begin{align*}
	 &\int_{-\infty}^{\infty}\psi_{\tau l\lesssim |\tilde{\eta}|\lesssim1} \frac{1}{\tilde{\eta}\left\langle\frac{|\tilde{\eta}|}{\tau}\right\rangle^{\gamma}}\cdot \frac{2\eta^{\frac12}}{\xi_j^{\frac12}+\eta^{\frac12}}\cdot P_{\geq l^{-1}}g\left(\frac{\tilde{\eta}}{\nu\tau} + \xi_j^{\frac12}\right)\,d\tilde{\eta}\\
	 & = \int_{-\infty}^{\infty}P^{(\tilde{\eta})}_{\geq l^{-1}\left(\nu\tau\right)^{-1}}\left(\psi_{\tau l\lesssim |\tilde{\eta}|\lesssim1} \frac{1}{\tilde{\eta}\left\langle\frac{|\tilde{\eta}|}{\tau}\right\rangle^{\gamma}}\cdot \frac{2\eta^{\frac12}}{\xi_j^{\frac12}+\eta^{\frac12}}\right)\cdot P_{\geq l^{-1}}g\left(\frac{\tilde{\eta}}{\nu\tau} + \xi_j^{\frac12}\right)\,d\tilde{\eta}
	 \end{align*}
	 Then note that for any $k\geq 0$ we have 
	 \begin{align}\label{dyadic L1 bound}
	 \left\|P^{(\tilde{\eta})}_{\geq l^{-1}(\nu\tau)^{-1}}\left(\psi_{2^k\tau l\simeq |\tilde{\eta}|} \frac{1}{\tilde{\eta}\left\langle\frac{|\tilde{\eta}|}{\tau}\right\rangle^{\gamma}}\cdot \frac{2\eta^{\frac12}}{\xi_j^{\frac12}+\eta^{\frac12}}\right)\right\|_{L^1_{d\tilde{\eta}}}\lesssim 2^{-k}.
	 \end{align}
	 In fact we denote by
\begin{align*}
	h(\zeta):=\psi_{2^{k}\tau l\simeq|\zeta|}\frac{1}{\zeta\left\langle\frac{|\zeta|}{\tau}\right\rangle^{\gamma}}\cdot \frac{2\eta^{\frac12}(\zeta)}{\xi_j^{\frac12}+\eta^{\frac12}(\zeta)}.
\end{align*}	 
By definition of the operator $P^{(\tilde{\eta})}_{\geq l^{-1}(\nu\tau)^{-1}}$, we have
\begin{align*}
	\left|\left(P^{(\tilde{\eta})}_{\geq l^{-1}(\nu\tau)^{-1}}h\right)(\tilde{\eta})\right|\lesssim \nu\tau l\cdot\left|\partial_{\tilde{\eta}}h(\tilde{\eta})\right|.
\end{align*}
On the other hand, we have, by the fact $\left|\partial_{\tilde{\eta}}h(\tilde{\eta})\right|\lesssim\tilde{\eta}^{-2}$,
\begin{align*}
	\left\|\partial_{\tilde{\eta}}h(\tilde{\eta})\right\|_{L^{1}_{d\tilde{\eta}}}\lesssim 2^{-k}\tau^{-1}l^{-1}.
\end{align*}
Therefore the estimate \eqref{dyadic L1 bound} follows. \eqref{dyadic L1 bound} in turn implies 
\begin{align}\label{dyadic l1}
	\begin{split}
&\left\|P^{(\tilde{\eta})}_{\geq l^{-1}(\nu\tau)^{-1}}\left(\psi_{\tau l\lesssim |\tilde{\eta}|\lesssim1} \frac{1}{\tilde{\eta}\left\langle\frac{|\tilde{\eta}|}{\tau}\right\rangle^{\gamma}}\cdot \frac{2\eta^{\frac12}}{\xi_j^{\frac12}+\eta^{\frac12}}\right)\right\|_{L^1_{d\tilde{\eta}}}\\
&\lesssim \sum_{k\geq 0}\left\|P^{(\tilde{\eta})}_{\geq l^{-1}(\nu\tau)^{-1}}\left(\psi_{2^k\tau l\simeq |\tilde{\eta}|} \frac{1}{\tilde{\eta}\langle\frac{|\tilde{\eta}|}{\tau}\rangle^{\gamma}}\cdot \frac{2\eta^{\frac12}}{\xi_j^{\frac12}+\eta^{\frac12}}\right)\right\|_{L^1_{d\tilde{\eta}}}\lesssim 1. 
\end{split}
\end{align}
It follows that 
\begin{align*}
&\left| \int_{-\infty}^{\infty}P^{(\tilde{\eta})}_{\geq l^{-1}(\nu\tau)^{-1}}\left(\psi_{\tau l\lesssim |\tilde{\eta}|\lesssim1} \frac{1}{\tilde{\eta}\left\langle\frac{|\tilde{\eta}|}{\tau}\right\rangle^{\gamma}}\cdot \frac{2\eta^{\frac12}}{\xi_j^{\frac12}+\eta^{\frac12}}\right)\cdot P_{\geq l^{-1}}g\left(\frac{\tilde{\eta}}{\nu\tau} + \xi_j^{\frac12}\right)\,d\tilde{\eta}\right|\\
&\lesssim \left\|P^{(\tilde{\eta})}_{\geq l^{-1}(\nu\tau)^{-1}}\left(\psi_{\tau l\lesssim |\tilde{\eta}|\lesssim1} \frac{1}{\tilde{\eta}\left\langle\frac{|\tilde{\eta}|}{\tau}\right\rangle^{\gamma}}\cdot \frac{2\eta^{\frac12}}{\xi_j^{\frac12}+\eta^{\frac12}}\right)\right\|_{L^1_{d\tilde{\eta}}}\cdot \left\|P_{\geq l^{-1}}g\left(\frac{\tilde{\eta}}{\nu\tau} + \xi_j^{\frac12}\right)\right\|_{L^\infty_{d\tilde{\eta}}}\\
&\lesssim  l^\delta\cdot \lambda^{-1-\frac{\delta}{2}}\cdot \left\|f\right\|_{wC^{\delta}},
\end{align*}
as desired.  
\\

 Next, to deal with the contribution of the low frequency part of $g$, namely $P_{< l^{-1}}g\left(\frac{\tilde{\eta}}{\nu\tau} + \xi_j^{\frac12}\right)$, we localize this further to frequency $2^j l^{-1},\,j<0$, whence we get 
\begin{align*}
\left|P_{2^j  l^{-1}}g\left(\frac{\tilde{\eta}}{\nu\tau} + \xi_1^{\frac12}\right) - P_{2^j  l^{-1}}g\left(\frac{\tilde{\eta}}{\nu\tau} + \xi_2^{\frac12}\right)\right|\lesssim 2^j\cdot 2^{-\delta j}l^\delta\cdot\lambda^{-1-\frac{\delta}{2}}\cdot \left\|f\right\|_{wC^{\delta}}.
\end{align*}
Here we used the following fact about H\"older functions: Let $P_{j}$ be the standard dyadic Littlewood-Paley projection and $f$ be a H\"older function. Then we have
\begin{align*}
	\left|P_{j}f(x+y)-P_{j}f(x)\right|\lesssim |y|\left\|\partial P_{j}f\right\|_{L^{\infty}}\lesssim  |y|[f]_{\alpha}2^{j(1-\alpha)}.
\end{align*}
Here $[\cdot]_{\alpha}$ is the H\"older difference norm for a H\"older function in $C^{\alpha}$. For our use, recall that $[g]_{\alpha}\lesssim \lambda^{-1-\frac{\delta}{2}}\left\|f\right\|_{wC^{\delta}}$.

 On the other hand, it is easily seen that (recalling that $\eta\sim\lambda$ throughout)
\begin{align*}
\left\|P^{(\tilde{\eta})}_{2^j l^{-1}(\nu\tau)^{-1}}\left(\psi_{\tau l\lesssim |\tilde{\eta}|\lesssim1} \frac{1}{\tilde{\eta}\left\langle\frac{|\tilde{\eta}|}{\tau}\right\rangle^{\gamma}}\cdot \frac{2\eta^{\frac12}}{\xi_r^{\frac12}+\eta^{\frac12}}\right)\right\|_{L^1_{d\tilde{\eta}}}\lesssim& |j|,\\
\left\|\sum_{r=1,2}(-1)^rP^{(\tilde{\eta})}_{<l^{-1}(\nu\tau)^{-1}}\left(\psi_{\tau l\lesssim |\tilde{\eta}|\lesssim1} \frac{1}{\tilde{\eta}\left\langle\frac{|\tilde{\eta}|}{\tau}\right\rangle^{\gamma}}\cdot \frac{2\eta^{\frac12}}{\xi_r^{\frac12}+\eta^{\frac12}}\right)\right\|_{L^1_{d\tilde{\eta}}}
\lesssim& l|\log \tau l|\cdot\lambda^{-\frac12}\lesssim l|\log  l|\cdot\lambda^{-\frac12}.
\end{align*}
For the first estimate above, we simply write the cutoff $\psi_{\tau l\lesssim|\tilde{\eta}|\lesssim 1}=\psi_{\tau l\lesssim|\tilde{\eta}|\lesssim 2^{-j}\tau l}+\psi_{2^{-j}\tau l\lesssim |\tilde{\eta}|\lesssim 1}$. The contribution from $\psi_{\tau l\lesssim|\tilde{\eta}|\lesssim 2^{-j}\tau l}$ is bounded by $|j|$, upon integrating $|\tilde{\eta}|^{-1}$. While the contribution from $\psi_{2^{-j}\tau l\lesssim |\tilde{\eta}|\lesssim 1}$ is bounded by $1$, using the argument deriving \eqref{dyadic l1}. For the second estimate above, we use the crude bound by ignoring the cutoff $P^{(\tilde{\eta})}_{<l^{-1}(\nu\tau)^{-1}}$, and using the Lagrange Mean Value Theorem for the function inside the parenthesis.
\\

In total, the remaining contribution to the case $|\tilde{\eta}|\lesssim 1$ is then estimated as follows: 
\begin{align*}
&\left|\sum_{r=1,2}(-1)^r \int_{-\infty}^{\infty}\psi_{\tau l\lesssim |\tilde{\eta}|\lesssim1} \frac{1}{\tilde{\eta}\left\langle\frac{|\tilde{\eta}|}{\tau}\right\rangle^{\gamma}}\cdot \frac{2\eta^{\frac12}}{\xi_r^{\frac12}+\eta^{\frac12}}\cdot P_{< l^{-1}}g\left(\frac{\tilde{\eta}}{\nu\tau} + \xi_r^{\frac12}\right)\,d\tilde{\eta}\right|\\
&\lesssim \sum_{j<0}\left| \int_{-\infty}^{\infty}P^{(\tilde{\eta})}_{2^jl^{-1}(\nu\tau)^{-1}}\left(\psi_{\tau l\lesssim |\tilde{\eta}|\lesssim1} \frac{1}{\tilde{\eta}\left\langle\frac{|\tilde{\eta}|}{\tau}\right\rangle^{\gamma}}\cdot \frac{2\eta^{\frac12}}{\xi_1^{\frac12}+\eta^{\frac12}}\right)\cdot \sum_{r=1,2}(-1)^rP_{2^j l^{-1}}g\left(\frac{\tilde{\eta}}{\nu\tau} + \xi_r^{\frac12}\right)\,d\tilde{\eta}\right|\\
& + \left|\int_{-\infty}^{\infty} \sum_{r=1,2}(-1)^rP^{(\tilde{\eta})}_{<l^{-1}(\nu\tau)^{-1}}\left(\psi_{\tau l\lesssim |\tilde{\eta}|\lesssim1} \frac{1}{\tilde{\eta}\left\langle\frac{|\tilde{\eta}|}{\tau}\right\rangle^{\gamma}}\cdot \frac{2\eta^{\frac12}}{\xi_r^{\frac12}+\eta^{\frac12}}\right)\cdot P_{< l^{-1}}g\left(\frac{\tilde{\eta}}{\nu\tau} + \xi_2^{\frac12}\right)\,d\tilde{\eta}\right|\\
&\lesssim  \sum_{j<0}|j| 2^{(1-\delta)j}\cdot l^\delta\cdot\lambda^{-1-\frac{\delta}{2}}\cdot \left\|f\right\|_{wC^{\delta}} +  \chi_{l\ll 1}l|\log l|\cdot\lambda^{-\frac32}\cdot \left\|f\right\|_{wC^{\delta}}\\
&\lesssim  l^\delta\cdot\lambda^{-1-\frac{\delta}{2}}\cdot \left\|f\right\|_{wC^{\delta}}, 
\end{align*}
as desired. 
\\

{\it{Second integral on the right in \eqref{eq:smalllargetildeetasplit}, the contribution of large $\tilde{\eta}$.}} Here we can no longer neglect the oscillatory term $e^{i\Psi(\xi,\eta,\tau, \sigma')}$, and in fact we shall exploit that 
\[
\Psi(\xi,\eta,\tau, \sigma') = \left(1 + O\left(\frac{1}{\tau}\right)\right)\cdot\tilde{\eta}.
\]
By directly evaluating the Fourier transform, we infer for $r\geq 0$ that 
\begin{align*}
\left\|P^{(\tilde{\eta})}_{\ll 1}\left(\psi_{|\tilde{\eta}|\simeq 2^r} \frac{e^{i\Psi(\xi,\eta,\tau, \sigma')}}{\tilde{\eta}\left\langle\frac{|\tilde{\eta}|}{\tau}\right\rangle^{\gamma}}\cdot \frac{2\eta^{\frac12}}{\xi_j^{\frac12}+\eta^{\frac12}}\right)\right\|_{L^1_{d\tilde{\eta}}} = O_N\left(\frac{1}{2^{Nr}}\right)
\end{align*}
for any $N>0$, and in particular by summing over $r\geq 0$ we find that 
\begin{equation}\label{eq:lowfreqgoodL1bound}
\left\|P^{(\tilde{\eta})}_{\ll 1}\left(\psi_{|\tilde{\eta}|\gtrsim 1} \frac{e^{i\Psi(\xi,\eta,\tau, \sigma')}}{\tilde{\eta}\left\langle\frac{|\tilde{\eta}|}{\tau}\right\rangle^{\gamma}}\cdot \frac{2\eta^{\frac12}}{\xi_j^{\frac12}+\eta^{\frac12}}\right)\right\|_{L^1_{d\tilde{\eta}}}\lesssim 1.
\end{equation}
\\
Furthermore, we have the difference bound 
\begin{align*}
\left\|\sum_{r = 1,2}(-1)^rP^{(\tilde{\eta})}_{\ll 1}\left(\psi_{|\tilde{\eta}|\gtrsim 1} \frac{e^{i\Psi(\xi_r,\eta,\tau, \sigma')}}{\tilde{\eta}\left\langle\frac{|\tilde{\eta}|}{\tau}\right\rangle^{\gamma}}\cdot \frac{2\eta^{\frac12}}{\xi_r^{\frac12}+\eta^{\frac12}}\right)\right\|_{L^1_{d\tilde{\eta}}}\lesssim \lambda^{-\frac12}\cdot l. 
\end{align*}
Now we break the second integral on the right in \eqref{eq:smalllargetildeetasplit} into two contributions: restricting the function $g$ to large frequency $>C^{-1}\tau$, we get the term
\begin{align*}
\sum_{r = 1,2}(-1)^r\int_{-\infty}^{\infty}\psi_{|\tilde{\eta}|\gtrsim 1} \frac{e^{i\Psi(\xi_r,\eta,\tau, \sigma')}}{\tilde{\eta}\left\langle\frac{|\tilde{\eta}|}{\tau}\right\rangle^{\gamma}}\cdot \frac{2\eta^{\frac12}}{\xi_r^{\frac12}+\eta^{\frac12}}\cdot P_{>C^{-1}\tau}g\left(\frac{\tilde{\eta}}{\nu\tau} + \xi_r^{\frac12}\right)\,d\tilde{\eta}
\end{align*}
and using \eqref{eq:freqcutoff} we find 
\begin{align*}
&\left|\sum_{r = 1,2}(-1)^r P_{>C^{-1}\tau}g\left(\frac{\tilde{\eta}}{\nu\tau} + \xi_r^{\frac12}\right)\right|\lesssim \min\left\{l^\delta\cdot\lambda^{-1-\frac{\delta}{2}},\,\lambda^{-1-\frac{\delta}{2}}\cdot\tau^{-\delta}\right\}\cdot \left\|f\right\|_{wC^{\delta}}\\
&\sum_{r = 1,2}\left|P_{>C^{-1}\tau}g\left(\frac{\tilde{\eta}}{\nu\tau} + \xi_r^{\frac12}\right)\right|\lesssim \lambda^{-1-\frac{\delta}{2}}\cdot\tau^{-\delta}\cdot \left\|f\right\|_{wC^{\delta}}, 
\end{align*}
and in particular for any $\delta'<\delta$ we have 
\[
\left|\sum_{r = 1,2}(-1)^r P_{>C^{-1}\tau}g\left(\frac{\tilde{\eta}}{\nu\tau} + \xi_r^{\frac12}\right)\right|\lesssim_{\delta',\delta} l^{\delta'}\cdot\lambda^{-1-\frac{\delta'}{2}}\cdot \left(\lambda^{\frac12}\tau\right)^{-\tilde{\delta}}\cdot \left\|f\right\|_{wC^{\delta}}
\]
for suitable $\tilde{\delta} = \tilde{\delta}(\delta',\delta)>0$. Since we also have 
\[
\left|\sum_{r = 1,2}(-1)^r e^{i\Psi(\xi_r,\eta,\tau, \sigma')} \frac{2\eta^{\frac12}}{\xi_r^{\frac12}+\eta^{\frac12}}\right|\lesssim \lambda^{-\frac12}\cdot l\cdot\langle\tilde{\eta}\rangle\lesssim l\cdot\tau
\]
in the domain $\xi_r\simeq \lambda, r = 1,2$, we easily conclude that 
\begin{align*}
&\left|\sum_{r = 1,2}(-1)^r\int_{-\infty}^{\infty}\psi_{|\tilde{\eta}|\gtrsim 1} \frac{e^{i\Psi(\xi_r,\eta,\tau, \sigma')}}{\tilde{\eta}\left\langle\frac{|\tilde{\eta}|}{\tau}\right\rangle^{\gamma}}\cdot \frac{2\eta^{\frac12}}{\xi_r^{\frac12}+\eta^{\frac12}}\cdot P_{>C^{-1}\tau}g\left(\frac{\tilde{\eta}}{\nu\tau} + \xi_r^{\frac12}\right)\,d\tilde{\eta}\right|\\
&\lesssim_{\delta',\delta}  l^{\delta'}\cdot\lambda^{-1-\frac{\delta'}{2}}\cdot \left\|f\right\|_{wC^{\delta}}
\end{align*}
It remains to deal with the contribution of the term arising when $g$ is replaced by $P_{\leq C^{-1}\tau}g$. Using Plancherel's theorem, we can write this as 
\begin{align*}
&\sum_{r = 1,2}(-1)^r\int_{-\infty}^{\infty}\psi_{|\tilde{\eta}|\gtrsim 1} \frac{e^{i\Psi(\xi_r,\eta,\tau, \sigma')}}{\tilde{\eta}\left\langle\frac{|\tilde{\eta}|}{\tau}\right\rangle^{\gamma}}\cdot \frac{2\eta^{\frac12}}{\xi_r^{\frac12}+\eta^{\frac12}}\cdot P_{\leq C^{-1}\tau}g\left(\frac{\tilde{\eta}}{\nu\tau} + \xi_r^{\frac12}\right)\,d\tilde{\eta}\\
& = \sum_{r = 1,2}(-1)^r\int_{-\infty}^{\infty}P^{(\tilde{\eta})}_{\lesssim C^{-1}}\left(\psi_{|\tilde{\eta}|\gtrsim 1} \frac{e^{i\Psi(\xi_r,\eta,\tau, \sigma')}}{\tilde{\eta}\left\langle\frac{|\tilde{\eta}|}{\tau}\right\rangle^{\gamma}}\cdot \frac{2\eta^{\frac12}}{\xi_r^{\frac12}+\eta^{\frac12}}\right)\cdot P_{\leq C^{-1}\tau}g\left(\frac{\tilde{\eta}}{\nu\tau} + \xi_r^{\frac12}\right)\,d\tilde{\eta}.
\end{align*}
Here we can take advantage of \eqref{eq:lowfreqgoodL1bound} and the bound following it provided $C\gg 1$, which gives 
\begin{align*}
&\left| \sum_{r = 1,2}(-1)^r\int_{-\infty}^{\infty}P^{(\tilde{\eta})}_{\lesssim C^{-1}}\left(\psi_{|\tilde{\eta}|\gtrsim 1} \frac{e^{i\Psi(\xi_r,\eta,\tau, \sigma')}}{\tilde{\eta}\left\langle\frac{|\tilde{\eta}|}{\tau}\right\rangle^{\gamma}}\cdot \frac{2\eta^{\frac12}}{\xi_r^{\frac12}+\eta^{\frac12}}\right)\cdot P_{\leq C^{-1}\tau}g\left(\frac{\tilde{\eta}}{\nu\tau} + \xi_r^{\frac12}\right)\,d\tilde{\eta}\right|\\
&\lesssim \left(l^\delta\cdot\lambda^{-1-\frac{\delta}{2}} +  \chi_{l\lesssim\lambda^{\frac12}}\lambda^{-\frac32}\cdot l\right)\cdot \left\|f\right\|_{wC^{\delta}}\\
&\lesssim l^\delta\cdot\lambda^{-1-\frac{\delta}{2}} \cdot \left\|f\right\|_{wC^{\delta}}.
\end{align*}
This concludes the proof. 
\end{proof}

\begin{remark}\label{rem:lemweightedHolder1} The loss of Holder regularity is only due to the contribution of the large frequency part $P_{>C^{-1}\tau}g$ of $g$. In particular, if $g$ is very smooth, then this contribution will be very small since $\tau\gg 1$. 
\end{remark}

In a similar vein, we can improve the conclusion if the function $f(\eta)$ has added differentiability:

\begin{lemma}\label{lem:weightedHolder2} Assume that $f$ is supported at $\eta\simeq\lambda\gtrsim \hbar^{-2}$, and consider the function 
	\begin{align*}
	\left(\mathcal{H}_{\text{twisted}}^{\pm}f\right)(\xi;\tau,\sigma'): = \int_0^\infty \frac{e^{\pm i\left[\nu\tau\eta^{\frac12} + \hbar^{-1}\rho\left(x_{\sigma'}^{(\eta)};\beta',\hbar\right)\right]}}{\left(\xi - \eta\right)\left\langle\xi^{\frac12} - \eta^{\frac12}\right\rangle^\gamma}\,f(\eta)\,d\eta,\quad \gamma>0, 
	\end{align*}
	where we use the notation $\beta' = \hbar\eta^{\frac12}\cdot\frac{\lambda(\tau)}{\lambda(\sigma)},\quad x_{\sigma'}^{(\eta)} = \hbar\eta^{\frac12}\cdot\nu\sigma\cdot\frac{\lambda(\tau)}{\lambda(\sigma)}$. Then for any $\delta'<\delta$ we have the relation 
	\begin{align*}
	&\left(\mathcal{H}_{\text{twisted}}^{\pm}f\right)(\xi;\tau,\sigma') = e^{\pm i\left[\nu\tau\xi^{\frac12} + \hbar^{-1}\rho\left(x_{\sigma'}^{(\xi)};\alpha',\hbar\right)\right]}\cdot g_{\pm}(\xi;\tau, \sigma'),\\
	&\sum_{0\leq r\leq N}\sup_{\xi>0}\left|\chi_{\xi\simeq\lambda}\xi^r\partial_{\xi}^rg(\xi;\tau, \sigma')\right| + \left\|\chi_{\xi\simeq\lambda}\xi^{N+\delta'}\partial_{\xi}^Ng(\xi;\tau, \sigma')\right\|_{wC^{\delta'}_{\xi}}\\&\lesssim_{\delta',\delta,\gamma}
	\sum_{0\leq r\leq N}\sup_{\xi>0}\left|\chi_{\xi\simeq\lambda}\xi^r\partial_{\xi}^rf(\xi;\tau, \sigma')\right| + \left\|\chi_{\xi\simeq\lambda}\xi^{N+\delta}\partial_{\xi}^Nf(\xi;\tau, \sigma')\right\|_{wC^{\delta}}.
	\end{align*}
\end{lemma}
The proof is analogous to the one of the preceding lemma. 
\\

Using the preceding lemmas, we can now analyze the effect of the transference operator at angular momentum $n, |n|\geq 2$ on admissible singular terms. For the following proposition, recall the terminology from Proposition \ref{prop: K operator}. 

\begin{proposition}\label{prop:transferenceonsingularngeq2} Assume that $\xb(\tau, \xi)$ is an admissible singular part at angular momentum $n, |n|\geq 2$, as in Definition~\ref{defi:xsingulartermsngeq2adm}. Then we have 
	\[
	\chi_{\hbar^2\xi\gtrsim 1}\cdot \left(\mathcal{K}_{\hbar}^{(0)}\xb\right)(\tau, \xi) = \xb_1(\tau, \xi) + \xb_2(\tau, \xi) + \xb_3(\tau, \xi)
	\]
	where $\xi^{\frac12}\cdot \tau^{-1}\cdot\xb_1$ is source admissibly singular with vanishing principal ingoing part(recall Def.~\ref{defi:xsingulartermsngeq2adm}), while $\xb_2$ satisfies 
	\[
	\left\|\xb_2\right\|_{\Sh_{0}}\lesssim \tau^{-5}. 
	\]
	Furthermore $\xb_3 =  \sum_{\pm} e^{\pm i\nu\tau\xi^{\frac12}}\cdot g_{\pm}(\tau, \xi)$
	where $g_{\pm}(\tau, \xi)$ is $C^\infty$ with respect to the second variable, and satisfies the bounds 
	\[
	\hbar^{10\delta_1}\cdot \left|\xi^{k_2}\cdot\partial_{\xi}^{k_2}\partial_{\tau}^{\iota}g_{\pm}(\tau, \xi)\right|\lesssim\log\tau\cdot \hbar^{-\iota}\tau^{-1-\nu-\iota}\cdot\hbar^{-1}(\log \xi)^{N_1}\cdot \xi^{-1-\nu}\left\langle\hbar\xi^{\frac12}\right\rangle^{-4},\quad 0\leq k_2\leq 10. 
	\]
	As a corollary, we see that 
	\[
	\chi_{\hbar^2\xi\gtrsim 1}\cdot \left(\mathcal{K}_{\hbar}\xb\right)(\tau, \xi) = \xb_1(\tau, \xi) + \xb_2(\tau, \xi) +  \xb_3(\tau, \xi)
	\]
	with $\xb_1$ admissibly singular and $\xb_2, \xb_3$ as before. 
	Furthermore, we can write 
	\[
	\chi_{\hbar^2\xi\gtrsim 1}\cdot \left(\mathcal{K}_{\hbar}^{(0)}\mathcal{D}_{\tau}\xb\right)(\tau, \xi) = \yb_1(\tau, \xi) + \yb_2(\tau, \xi) + \yb_3(\tau, \xi)
	\]
	where $\tau^{-1}\cdot \yb_1$ is source admissibly singular, $ \yb_2$ satisfies a bound like $ \xb_2$ but with $\Sh_{0}$ replaced by $\Sh_{1}$ and $\xi^{-\frac12}\cdot \yb_3$ is like $\xb_3$. \\
	The same conclusion applies if we replace $\mathcal{K}_{\hbar}^{(0)}$ by $\big[\mathcal{D}_{\tau}, \mathcal{K}_{\hbar}^{(0)}\big]$, $\big(\mathcal{K}_{\hbar}^{(0)}\big)^2$. 
\end{proposition}
\begin{remark}\label{rem:transferenceonsingularngeq20} Notice that source admissibility of $\xi^{\frac12}\cdot \tau^{-1}\cdot\xb_1$ implies that we needn't apply the temporal derivatives $\partial_{\tau}^{\iota}, \iota\in \{0, 1\}$. 
\end{remark}
\begin{remark}\label{rem:transferenceonsingularngeq21} The term $ \xb_3(\tau, \xi)$ arises when applying the transference operator to a connecting singular term with $l = 7$ in Definition~\ref{defi:xsingulartermsngeq2adm}: the added frequency decay due to the transference operator kernel ensures that the 'output' is then in $\Sh_{0}$, with weak temporal decay but an added structure. 
\end{remark}
\begin{remark}\label{rem:proptransferenceonsingularngeq2} The reason for the presence of the term $\xb_2$ is the fact that application of the transference operator always causes a small loss of differentiability of the coefficients $F_{l,k,i}^{(\pm)}(\tau, \sigma,\xi)$ with respect to the last variable. However, the fact that the preceding lemma and Remark~\ref{rem:lemweightedHolder1} identify the precise source for this loss of differentiability allow us to gain an additional temporal smallness for it and place it into $\xb_2$. This added temporal decay will be crucial, since otherwise the fine structure would be required to obtain acceptable bounds on the physical side. 
	Similarly, the term $\xb_3$ is also in the good space $\Sh_{0}$ but has poor temporal decay (its presence is forced by very technical reasons, related to the cutoff in the definition of admissibly singular function), and so we have to retain enough fine structure to be able to handle its contributions, when dealing with the modulation theory for the exceptional modes. 
\end{remark}

\begin{proof}
	We verify the conclusion  for both $\xb_{in}$ and $\xb_{out}$ according to Definition~\ref{defi:xsingulartermsngeq2adm}. 
	\\
	
	{\it{Contribution of the principal ingoing singular part, i.e., the first term on the right in \eqref{eq:xsingin}.}} We fix a sign $\pm$ as well as $k, i$ and the time $\sigma$, as we can move the temporal integral to the beginning of the resulting expression after applying the transference operator. The transference operator shall then be given by the kernel 
	\[
	\frac{F(\xi,\eta;\hbar)\cdot \rho_n(\eta)}{\xi - \eta}.
	\]
	In fact, a technical complication caused by the precise statement of the trace-type derivatives in \eqref{eq:trace der} of Proposition \ref{prop: K operator} forces us to replace $F(\xi,\eta;\hbar)$ by $\chi_{\lesssim 1}\left(\frac{\left|\xi^{\frac12} - \eta^{\frac12}\right|}{\xi^{\delta}}\right)\cdot F(\xi,\eta;\hbar)$ for some small $\delta>0$. That we may do so while only generating an error of type $\xb_2$ follows from Lemma~\ref{lem:transferenceonadmsingoffdiag}, and we shall henceforth tacitly assume that the kernel $F(\xi,\eta;\hbar)$ satisfies the trace-type bounds without restriction on $\xi, \eta$, as we may after inclusion of the preceding cutoff. We conclude that the integral is then restricted to the region $\xi\simeq \eta$, and we can take advantage of the preceding lemma, in conjunction to refined bounds on the kernel function $F$ and the spectral density $\rho_n$.  
	\\
	Since $\tau, \sigma$ are fixed, we may as well move the function $a_{k,i}^{(\pm)}(\tau, \sigma)$ to the outside of the $\eta$-integral. Then we need to show that the expression (where we use the same notation as in the proof of the preceding lemma)
	\begin{align*}
	\Xi_1^{(\pm)}(\tau, \xi): = \chi_{\xi\geq \hbar^{-2}}\int_{\tau_0}^{\tau}a_{k,i}^{(\pm)}(\tau, \sigma)\cdot \int_0^\infty\chi_{\xi\sim\eta} \frac{F(\xi,\eta;\hbar)\cdot \rho_n(\eta)}{\xi - \eta}\cdot \chi_{\eta\geq \hbar^{-2}}\hbar^{-1}\frac{e^{\pm i\left[\nu\tau\eta^{\frac12} + \hbar^{-1}\rho\left(x_{\sigma'}^{(\eta)};\beta',\hbar\right)\right]}}{\eta^{1+\frac{k\nu}{2}}}\left(\log\eta\right)^i\,d\eta\,d\sigma
	\end{align*}
	is as asserted in the proposition. Write the preceding expression as
	\begin{align*}
	\Xi_1^{(\pm)}(\tau, \xi) &=  \chi_{\xi\geq \hbar^{-2}}\int_{\tau_0}^{\tau}a_{k,i}^{(\pm)}(\tau, \sigma)e^{\pm i\left[\nu\tau\xi^{\frac12} + \hbar^{-1}\rho(x_{\sigma'}^{(\xi)};\alpha',\hbar)\right]}\cdot \tilde{\Xi}_1^{(\pm)}(\tau,\sigma,\xi)\,d\sigma
        \end{align*}
	where we have (recall \eqref{eq:Psixietatausigmaprime} for the definition of $\Psi(\xi,\eta,\tau,\sigma')$)
	\[
	\tilde{\Xi}_1^{(\pm)}(\tau,\sigma,\xi) = \int_0^\infty\chi_{\xi\sim\eta} \frac{F(\xi,\eta;\hbar)\cdot \rho_n(\eta)}{\xi - \eta}\cdot \chi_{\eta\geq \hbar^{-2}}\hbar^{-1}\frac{e^{\pm i\Psi(\xi,\eta,\tau,\sigma')}}{\eta^{1+\frac{k\nu}{2}}}\left(\log\eta\right)^i\,d\eta
	\]
	We claim that we may write 
	\[
	a_{k,i}^{(\pm)}(\tau, \sigma)\cdot \tilde{\Xi}_1^{(\pm)}(\tau,\sigma,\xi) = \hbar^{-1}\left\langle\hbar^2\xi\right\rangle^{-\frac14}\cdot \xi^{-1-\frac{k\nu}{2}}\left(\log\xi\right)^i\cdot F_{1,k,i}^{\pm}\left(\tau,\sigma,\frac{\lambda^2(\tau)}{\lambda^2(\sigma)}\xi\right),
	\]
	where the function $F_{1,k,i}^{\pm}$ satisfies the bounds required in Definition~\ref{defi:xsingulartermsngeq2adm} and with $l = 1$. Due to the scaling invariance of these bounds, it suffices to prove them for the function 
	\[
	Z_1^{\pm}(\tau,\sigma,\xi): = \tilde{\chi}_{\xi\geq \hbar^{-2}}a_{k,i}^{(\pm)}(\tau, \sigma)\cdot \tilde{\Xi}_1^{(\pm)}(\tau,\sigma,\xi)\cdot \hbar\left\langle\hbar^2\xi\right\rangle^{\frac14}\cdot \xi^{1+\frac{k\nu}{2}}\left(\log\xi\right)^{-i}
	\]
	where $\tilde{\chi}_{\xi\geq \hbar^{-2}}\chi_{\xi\geq \hbar^{-2}} = \chi_{\xi\geq \hbar^{-2}}$.	In order to get the derivative bounds with respect to $\xi$, introduce the variable 
	\[
	\tilde{\tilde{\eta}}: = \Psi(\xi,\eta,\tau,\sigma'),
	\]
	whence in terms of the preceding variable $\tilde{\eta} = \nu\tau\left(\eta^{\frac12} - \xi^{\frac12}\right)$, we have in light of the proof of the preceding lemma as well as Lemma \ref{lem: Lemma 3.4 CDST}
	\[
	\tilde{\tilde{\eta}} = \tilde{\eta}\cdot\left(1 + \kappa\left(\tilde{\tilde{\eta}},\xi,\tau, \sigma'\right)\right),
	\]
	where on the support of the integrand (i.e., $\eta\geq \hbar^{-2}$ and also $\xi\geq \hbar^{-2}$) we have 
	\[
	\left|\kappa\left(\tilde{\tilde{\eta}},\xi,\tau, \sigma'\right)\right| = O\left(\left[\sigma'\tau\right]^{-1}\right)\ll 1,\quad  \left|\partial_{\xi}^{k_1}\partial_{\tau}^{k_2}\kappa\left(\tilde{\tilde{\eta}},\xi,\tau, \sigma'\right)\right|\lesssim \tau^{-1-k_2}\sigma'^{-1}\xi^{-k_1},
	\]
	and furthermore we get 
	\[
	\frac{\partial\tilde{\tilde{\eta}}}{\partial\tilde{\eta}}\simeq 1,\quad \left|\partial_{\xi}^{k_1}\partial_{\tau}^{k_2}\left(\left[\frac{\partial\tilde{\tilde{\eta}}}{\partial\tilde{\eta}}\right]^{\iota}\right)\right|\lesssim \tau^{-1-k_2}\sigma'^{-1}\xi^{-k_1},\quad \iota = \pm 1,\,k_{1,2}\geq 0. 
	\]
	Then write 
	\[
	F(\xi,\eta;\hbar) = \tilde{F}\left(\x^{\frac12},\eta^{\frac12};\hbar\right) =  \tilde{F}\left(\x^{\frac12},\xi^{\frac12} + \frac{\tilde{\tilde{\eta}}}{\nu\tau(1+\kappa)};\hbar\right)
	\]
	and interpret $\eta = \left(\xi^{\frac12} + \frac{\tilde{\tilde{\eta}}}{\nu\tau(1+\kappa)}\right)^2$ as function of $\xi,\tilde{\tilde{\eta}},\tau, \sigma'$. Thus we now have 
	\begin{align*}
	\tilde{\Xi}_1^{(\pm)}(\tau,\sigma,\xi) = \int_0^\infty\chi_{\hbar^{-2}\lesssim\xi\ll\eta} \frac{\tilde{F}\left(\x^{\frac12},\xi^{\frac12} + \frac{\tilde{\tilde{\eta}}}{\nu\tau(1+\kappa)};\hbar\right)}{\tilde{\tilde{\eta}}}\cdot e^{\pm i\tilde{\tilde{\eta}}}\cdot\zeta(\tilde{\tilde{\eta}},\xi,\tau, \sigma')\,d\tilde{\tilde{\eta}},
	\end{align*}
	where we have set 
	\[
	\zeta(\tilde{\tilde{\eta}},\xi,\tau, \sigma') =  \chi_{\eta\geq \hbar^{-2}}\hbar^{-1}\frac{(\log\eta)^i}{\eta^{1+\frac{k\nu}{2}}}\rho_n(\eta)\cdot\frac{2\eta^{\frac12}}{\eta^{\frac12}+\xi^{\frac12}}\cdot\frac{\partial\tilde{\eta}}{\partial \tilde{\tilde{\eta}}}.
	\]
	Using the bounds from above, it is then easy to verify that 
	\begin{align*}
	\left|\partial_{\xi}^{k_1}\left(\zeta(\tilde{\tilde{\eta}},\xi,\tau, \sigma')\cdot \hbar\left\langle\hbar^2\xi\right\rangle^{\frac14}\cdot \xi^{1+\frac{k\nu}{2}}\left(\log\xi\right)^{-i}\right)\right|\lesssim \left\langle\hbar^2\xi\right\rangle^{\frac14}\cdot \tau^{-1}\xi^{-k_1},\,k_1>0
	\end{align*}
	on the support of the integrand. Furthermore, observe the relation 
	\begin{align*}
	\partial_{\xi^{\frac12}}\left(\tilde{F}\left(\x^{\frac12},\xi^{\frac12} + \frac{\tilde{\tilde{\eta}}}{\nu\tau(1+\kappa)};\hbar\right)\right) &= \left(\left(\partial_{\xi^{\frac12}} + \partial_{\eta^{\frac12}}\right)\tilde{F}\right)\left(\xi^{\frac12},\xi^{\frac12} + \frac{\tilde{\tilde{\eta}}}{\nu\tau(1+\kappa)};\hbar\right)\\
	& - \frac{\tilde{\tilde{\eta}}\cdot\partial_{\xi^{\frac12}}\kappa}{\nu\tau(1+\kappa)^2}\cdot \left(\partial_{\eta^{\frac12}}\tilde{F}\right)\left(\xi^{\frac12},\xi^{\frac12} + \frac{\tilde{\tilde{\eta}}}{\nu\tau(1+\kappa)};\hbar\right),
	\end{align*}
	and we can bound each term on the right due to the refined derivative bounds for the kernel of the transference operator(recalling Prop.~\ref{prop: K operator}): 
	\begin{align*}
	&\left|\left(\left(\partial_{\xi^{\frac12}} + \partial_{\eta^{\frac12}}\right)\tilde{F}\right)\left(\x^{\frac12},\xi^{\frac12} + \frac{\tilde{\tilde{\eta}}}{\nu\tau(1+\kappa)};\hbar\right)\right|\lesssim \xi^{-\frac12}\cdot\langle\hbar^2\xi\rangle^{-\frac12},\\
	&\left|\frac{\tilde{\tilde{\eta}}\cdot\partial_{\xi^{\frac12}}\kappa}{\nu\tau(1+\kappa)^2}\cdot \left(\partial_{\eta^{\frac12}}\tilde{F}\right)\left(\xi^{\frac12},\xi^{\frac12} + \frac{\tilde{\tilde{\eta}}}{\nu\tau(1+\kappa)};\hbar\right)\right|
	\lesssim  \left|\partial_{\xi^{\frac12}}\kappa\right|\cdot \left|\left(\left(\eta^{\frac12}-\xi^{\frac12}\right)\partial_{\eta^{\frac12}}\tilde{F}\right)\left(\xi^{\frac12},\xi^{\frac12} + \frac{\tilde{\tilde{\eta}}}{\nu\tau(1+\kappa)};\hbar\right)\right|\\
	&\hspace{7cm}\lesssim \left(\tau\sigma'\right)^{-1}\xi^{-\frac12+\delta}\cdot\left\langle\hbar^2\xi\right\rangle^{-\frac12}.
	\end{align*}

Analogous bounds are obtained for higher derivatives, although here we need to pay careful attention to not let too many 'non-diagonal' derivatives act on the kernel $F(\xi, \eta;\hbar)$. This can be achieved via integration by parts: 
\begin{align*}
&\partial_{\xi^{\frac12}}^l\tilde{\Xi}_1^{(\pm)}(\tau,\sigma,\xi)\\
& = \sum_{\sum_{j=1}^4 l_j = l}C_{l_{1,2,3,4}} \int_0^\infty \partial_{\xi^{\frac12}}^{l_1}\big(\chi_{\hbar^{-2}\lesssim\xi\ll\eta}\big) \frac{T_{l_{2,3}}\tilde{F}\left(\x^{\frac12},\xi^{\frac12} + \frac{\tilde{\tilde{\eta}}}{\nu\tau(1+\kappa)};\hbar\right)}{\tilde{\tilde{\eta}}}\cdot e^{\pm i\tilde{\tilde{\eta}}}\cdot \partial_{\xi^{\frac12}}^{l_4}\zeta(\tilde{\tilde{\eta}},\xi,\tau, \sigma')\,d\tilde{\tilde{\eta}},
\end{align*}
where we set $T_{l_{2,3}} = \sum_{l_2'\leq l_2}D_{l_2, l_2'}(\xi)\cdot \big(\big(\partial_{\xi^{\frac12}}\kappa\cdot \frac{\tilde{\tilde{\eta}}}{\nu\tau(1+\kappa)}\big)\partial_{\eta^{\frac12}}\big)^{l_3}\big(\partial_{\xi^{\frac12}} + \partial_{\eta}^{\frac12}\big)^{l_2'}$ , where we have the estimate 
\begin{align*}
\big|D_{l_2, l_2'}(\xi)\big|\lesssim \xi^{-\frac{l_2-l_2'}{2}}
\end{align*}
and further $D_{l_2, l_2'}$ is smooth and obeys symbol type bounds. Then we note that 
\begin{align*}
\Big(\big(\partial_{\xi^{\frac12}}\kappa\cdot \frac{\tilde{\tilde{\eta}}}{\nu\tau(1+\kappa)}\big)\partial_{\eta^{\frac12}}\tilde{F}\Big)\left(\x^{\frac12},\xi^{\frac12} + \frac{\tilde{\tilde{\eta}}}{\nu\tau(1+\kappa)};\hbar\right) = \partial_{\xi^{\frac12}}\kappa\cdot\big(\tilde{\tilde{\eta}}\partial_{\tilde{\tilde{\eta}}}\big)\Big(\tilde{F}\left(\x^{\frac12},\xi^{\frac12} + \frac{\tilde{\tilde{\eta}}}{\nu\tau(1+\kappa)};\hbar\right)\Big). 
\end{align*}
We can then perform integration by parts in the $\tilde{\tilde{\eta}}$-integral for $\partial_{\xi^{\frac12}}^l\tilde{\Xi}_1^{(\pm)}(\tau,\sigma,\xi)$. Note that amongst the other factors only the phase $e^{\pm i\tilde{\tilde{\eta}}}$ does not have symbol behavior with respect to $\tilde{\tilde{\eta}}$. However, the tacit inclusion of the cutoff $\chi_{\lesssim 1}\big(\frac{|\xi^{\frac12} - \eta^{\frac12}|}{\xi^{\delta}}\big)$ at the beginning implies $\big|\tilde{\tilde{\eta}}\big|\lesssim \xi^{\delta}\cdot \tau$, and so 
\begin{align*}
\Big|\partial_{\xi^{\frac12}}\kappa\cdot\big(\tilde{\tilde{\eta}}\partial_{\tilde{\tilde{\eta}}}\big)e^{\pm i\tilde{\tilde{\eta}}}\Big|\lesssim \xi^{\delta - \frac12}
\end{align*}
on the support of the integrand. Performing $l_3$ integrations by parts thus results in a loss of $\xi^{\delta\cdot l_3}$; the smallness gain $[\hbar\xi^{\frac12}]^{-1}$ for $F(\xi,\eta;\hbar)$ and its diagonal derivatives implied by Prop.~\ref{prop: K operator} allows us to replace the loss of $\xi^{\delta l_3}$ by $\hbar^{-2\delta l_3}$, which is an acceptable loss in light of the required derivative bounds. 
\\

In light of the preceding, it suffices to prove the $L^\infty$ and Holder type bound for the un-differentiated expression $\tilde{\Xi}_1^{(\pm)}(\tau,\sigma,\xi)$, modulo replacing the kernel $F$ by a trace-type derivative, which does not affect things due to our implicit restriction $|\xi^{\frac12} - \eta^{\frac12}|\lesssim \xi^{\delta}$, as well as Prop.~\ref{prop: K operator}.
To proceed, we decompose $\tilde{\Xi}_1^{(\pm)}(\tau,\sigma,\xi)$ further into a `weakly diagonal' and a `strongly diagonal' part via inclusion of smooth cutoffs $\chi_{\left|\xi^{\frac12}-\eta^{\frac12}\right|\geq 1}$, $\chi_{\left|\xi^{\frac12}-\eta^{\frac12}\right|< 1}$.
\\

{\it{Weakly diagonal part}}. Set 
\begin{align*}
\tilde{\Xi}_{11}^{(\pm)}(\tau,\sigma,\xi): =  \int_0^\infty\chi_{\hbar^{-2}\lesssim\xi\simeq\eta}\cdot \chi_{\left|\xi^{\frac12}-\eta^{\frac12}\right|\geq 1}\frac{F(\xi,\eta;\hbar)\cdot \rho_n(\eta)}{\xi - \eta}\cdot \chi_{\eta\geq \hbar^{-2}}\hbar^{-1}\frac{e^{\pm i\Psi(\xi,\eta,\tau,\sigma')}}{\eta^{1+\frac{k\nu}{2}}}\left(\log\eta\right)^i\,d\eta
\end{align*}
Here we proceed as above and, using the same notation we get 
\begin{align*}
\tilde{\Xi}_{11}^{(\pm)}(\tau,\sigma,\xi) =  \int_{-\infty}^\infty\chi_{\xi\simeq\eta} \chi_{\left|\xi^{\frac12}-\eta^{\frac12}\right|\geq 1} \frac{\tilde{F}\left(\x^{\frac12},\xi^{\frac12} + \frac{\tilde{\tilde{\eta}}}{\nu\tau(1+\kappa)};\hbar\right)}{\tilde{\tilde{\eta}}}\cdot e^{\pm i\tilde{\tilde{\eta}}}\cdot\zeta(\tilde{\tilde{\eta}},\xi,\tau, \sigma')\,d\tilde{\tilde{\eta}},
\end{align*}
where $\eta$ is thought of as function of $\xi,\tilde{\tilde{\eta}}, \tau, \sigma'$. Note that the restriction $\left|\xi^{\frac12}-\eta^{\frac12}\right|\geq 1$ implies $\left|\tilde{\tilde{\eta}}\right|\gtrsim \tau$, and that we can write 
\[
\chi_{\left|\xi^{\frac12}-\eta^{\frac12}\right|\geq 1} =  \chi_{\left|\frac{\tilde{\tilde{\eta}}}{\nu\tau(1+\kappa)}\right|\geq 1},
\]
and we have 
\[
\left|\partial_{\xi}\left(\chi_{\left|\frac{\tilde{\tilde{\eta}}}{\nu\tau(1+\kappa)}\right|\geq 1}\right)\right|\lesssim \xi^{-1}. 
\]
Taking advantage of the refined asymptotics of $F(\xi,\eta;\hbar)$ and its derivatives, as detailed in Proposition \ref{prop: K operator}, we infer the bound 
\begin{align*}
\left|\partial_{\xi}^{k_1}\left( \tilde{\Xi}_{11}^{(\pm)}(\tau,\sigma,\xi)\cdot \hbar\left\langle\hbar^2\xi\right\rangle^{\frac14}\cdot \xi^{1+\frac{k\nu}{2}}\left(\log\xi\right)^{-i}\right)\right|\lesssim\xi^{-k_1},\,k_1\in \{0, 1\}.  
\end{align*}
In fact, integrability over $\left|\tilde{\tilde{\eta}}\right|\gtrsim \tau$ follows from the refined off-diagonal decay
\begin{align*}
\big|\partial_{\xi^{\frac12}}^{k_1}\tilde{F}\left(\x^{\frac12},\xi^{\frac12} + \frac{\tilde{\tilde{\eta}}}{\nu\tau(1+\kappa)};\hbar\right)\big|\lesssim \xi^{-\frac{k_1}{2}}\cdot\big| \frac{\tilde{\tilde{\eta}}}{\tau}\big|^{-\gamma},\,\gamma>0,\,k_1\in \{0, 1\}. 
\end{align*}
The required weighted Holder bound follows from simple interpolation.
\\

 {\it{Strongly diagonal part}}. This is the term 
\begin{align*}
\tilde{\Xi}_{12}^{(\pm)}(\tau,\sigma,\xi): =  \int_0^\infty\chi_{\hbar^{-2}\lesssim \xi\simeq\eta}\cdot \chi_{\left|\xi^{\frac12}-\eta^{\frac12}\right|<1}\frac{F(\xi,\eta;\hbar)\cdot \rho_n(\eta)}{\xi - \eta}\cdot \chi_{\eta\geq \hbar^{-2}}\hbar^{-1}\frac{e^{\pm i\Psi(\xi,\eta,\tau,\sigma')}}{\eta^{1+\frac{k\nu}{2}}}\left(\log\eta\right)^i\,d\eta
\end{align*}
We split this into a term which can be handled directly by the preceding lemma, as well as a better error term. For this split $F(\xi,\eta;\hbar) = \tilde{F}\left(\xi^{\frac12},\eta^{\frac12};\hbar\right)$ as follows:
\[
\tilde{F}\left(\xi^{\frac12},\eta^{\frac12};\hbar\right) =  \tilde{F}\left(\xi^{\frac12},\xi^{\frac12};\hbar\right) + \int_0^1 \left(\xi^{\frac12} - \eta^{\frac12}\right)\cdot\partial_{\eta^{\frac12}}\tilde{F}\left(\xi^{\frac12}, s\xi^{\frac12} + (1-s)\eta^{\frac12};\hbar\right)\,ds,
\]
and so we infer 
\begin{align*}
\tilde{\Xi}_{12}^{(\pm)}(\tau,\sigma,\xi) = \tilde{\Xi}_{12}^{(\pm),a}(\tau,\sigma,\xi) +  \tilde{\Xi}_{12}^{(\pm),b}(\tau,\sigma,\xi),
\end{align*}
where we set 
\begin{align*}
&\tilde{\Xi}_{12}^{(\pm),a}(\tau,\sigma,\xi)  =  \tilde{F}\left(\xi^{\frac12},\xi^{\frac12};\hbar\right)\cdot\int_0^\infty\chi_{\hbar^{-2}\lesssim\xi\simeq\eta}\cdot \chi_{\left|\xi^{\frac12}-\eta^{\frac12}\right|<1}\frac{\rho_n(\eta)}{\xi - \eta}\cdot \chi_{\eta\geq \hbar^{-2}}\hbar^{-1}\frac{e^{\pm i\Psi(\xi,\eta,\tau,\sigma')}}{\eta^{1+\frac{k\nu}{2}}}\left(\log\eta\right)^i\,d\eta\\
&\tilde{\Xi}_{12}^{(\pm),b}(\tau,\sigma,\xi)  =  \int_0^\infty\chi_{\hbar^{-2}\lesssim\xi\simeq\eta}\cdot \chi_{\left|\xi^{\frac12}-\eta^{\frac12}\right|<1} \tilde{\tilde{F}}\left(\xi^{\frac12},\eta^{\frac12};\hbar\right)\frac{\rho_n(\eta)}{\xi^{\frac12}+\eta^{\frac12}}\cdot \chi_{\eta\geq \hbar^{-2}}\hbar^{-1}\frac{e^{\pm i\Psi(\xi,\eta,\tau,\sigma')}}{\eta^{1+\frac{k\nu}{2}}}\left(\log\eta\right)^i\,d\eta,\\
\end{align*}
where in the second integral we use the notation 
\[
\tilde{\tilde{F}}\left(\xi^{\frac12},\eta^{\frac12};\hbar\right) =  \int_0^1 \partial_{\eta^{\frac12}}\tilde{F}\left(\xi^{\frac12}, s\xi^{\frac12} + (1-s)\eta^{\frac12};\hbar\right)\,ds
\]
Using the special derivative bounds for the transference kernel, as well as the preceding lemma, we then easily infer the better than required bounds
\begin{align*}
\left|\xi^{k_2}\partial_{\xi}^{k_2}\left(a^{(\pm)}_{k,i}(\tau,\sigma)\tilde{\Xi}_{12}^{(\pm),a}(\tau,\sigma,\xi) \cdot  \hbar\left\langle\hbar^2\xi\right\rangle^{\frac14}\cdot \xi^{1+\frac{k\nu}{2}}\left(\log\xi\right)^{-i}\right)\right|\lesssim \left(\log\tau\right)^{N_1-i}\tau^{-1-\nu}\cdot\sigma^{-3},\,k_2\in\{0, 1\},
\end{align*}
which implies the desired Holder bound via interpolation. 
The same bounds obtain also for the second term $\tilde{\Xi}_{12}^{(\pm),b}(\tau,\sigma,\xi)$, where one can simply repeat the change of coordinates from the weakly diagonal part and follow the same arguments as in that case.  
\\

{\it{Contribution of the connecting ingoing singular part, i.e., the second term on the right in \eqref{eq:xsingin}.}} Here one can essentially follow the same outline as in the preceding case. 
Consider the expression
\begin{align*}
\int_{\tau_0}^{\tau}\int_0^\infty\chi_{\hbar^{-2}\lesssim\xi\sim\eta} \frac{F(\xi,\eta;\hbar)\cdot \rho_n(\eta)}{\xi - \eta}\cdot \chi_{\eta\geq \hbar^{-2}}\hbar^{-1}\frac{e^{\pm i\left[\nu\tau\eta^{\frac12} + \hbar^{-1}\rho\left(x_{\sigma'}^{(\eta)};\beta',\hbar\right)\right]}}{\eta^{1+\frac{k\nu}{2}}}\left(\log\eta\right)^i\cdot F_{l,k,i}^{(\pm)}\left(\tau,\sigma,\frac{\lambda^2(\tau)}{\lambda^2(\sigma)}\eta\right)\,d\eta d\sigma.
\end{align*}
Using the same change of variables as before, i.e., $\tilde{\tilde{\eta}}:=\Psi(\xi,\eta,\tau,\sigma')$, we can write the preceding term as (abusing notation in order to emphasize the analogy) 
\begin{align*}
\chi_{\xi\geq \hbar^{-2}} \int_{\tau_0}^{\tau} e^{\pm i\left[\nu\tau\xi^{\frac12} + \hbar^{-1}\rho\left(x_{\sigma'}^{(\xi)};\alpha',\hbar\right)\right]}\cdot \tilde{\Xi}_1^{(\pm)}(\tau,\sigma,\xi)\,d\sigma,
\end{align*}
where this time we set 
\begin{align*}
\tilde{\Xi}_1^{(\pm)}(\tau,\sigma,\xi) = \int_0^\infty\chi_{\hbar^{-2}\lesssim\xi\sim\eta} \frac{\tilde{F}\left(\xi^{\frac12},\xi^{\frac12} + \frac{\tilde{\tilde{\eta}}}{\nu\tau(1+\kappa)};\hbar\right)}{\tilde{\tilde{\eta}}}\cdot e^{\pm i\tilde{\tilde{\eta}}}\cdot\zeta\left(\tilde{\tilde{\eta}},\xi,\tau, \sigma'\right)\,d\tilde{\tilde{\eta}}
\end{align*}
and we define 
\begin{align*}
\zeta\left(\tilde{\tilde{\eta}},\xi,\tau, \sigma'\right) =  \chi_{\eta\geq \hbar^{-2}}\hbar^{-1}\frac{\left(\log\eta\right)^i}{\eta^{1+\frac{k\nu}{2}}}\rho_n(\eta)\cdot\frac{2\eta^{\frac12}}{\eta^{\frac12}+\xi^{\frac12}}\cdot\frac{\partial\tilde{\eta}}{\partial \tilde{\tilde{\eta}}}\cdot F_{l,k,i}^{(\pm)}\left(\tau, \sigma, \frac{\lambda^{2}(\tau)}{\lambda^{2}(\sigma)}\left(\xi^{\frac12}+\frac{\tilde{\tilde{\eta}}}{\nu\tau\left(1+\kappa\right)}\right)^{2}\right).
\end{align*}
In the following estimate $\tilde{\kappa}$ corresponds to $\kappa$ in Definition~\ref{defi:xsingulartermsngeq2adm}. We then obtain the following bounds for $\zeta$ in the region $\hbar^{-2}\lesssim\xi\sim\eta$: 
\begin{align*}
&\hbar^{\delta_1\cdot k_2}\cdot\left|\partial_{\xi}^{k_1}\left(\zeta\left(\tilde{\tilde{\eta}},\xi,\tau, \sigma'\right)\cdot \hbar\langle\hbar^2\xi\rangle^{\frac14}\cdot \xi^{1+\frac{k\nu}{2}}\left(\log\xi\right)^{-i}\right)\right|\\&\lesssim \langle\hbar^2\xi\rangle^{\frac14}\cdot \tau^{-1-\nu}\left(\log\tau\right)^{N_1 - i}\cdot \sigma^{-1}\left[\sigma^{-2} + \tilde{\kappa}\left(\hbar\frac{\lambda(\tau)}{\lambda(\sigma)}\eta^{\frac12}\right)\right]\cdot \xi^{-k_1},\quad 0\leq k_1\leq 20-l, 
\end{align*}
as well as the `closing H\"older bounds' 
\begin{align*}
&\hbar^{\delta_1\cdot(2-0l)}\left\|\partial_{\xi}^{20-l}\left(\zeta\left(\tilde{\tilde{\eta}},\xi,\tau, \sigma'\right)\cdot \hbar\left\langle\hbar^2\xi\right\rangle^{\frac14}\cdot \xi^{1+\frac{k\nu}{2}}\left(\log\xi\right)^{-i}\right)\right\|_{\dot{C}^{\delta}_{\xi}\left(\xi\simeq\mu\right)}\\&\lesssim \mu^{-(20-l+\delta)} \left\langle\hbar^2\mu\right\rangle^{\frac14}\cdot \tau^{-1-\nu}\left(\log\tau\right)^{N_1 - i}\cdot \sigma^{-1}\left[\sigma^{-2} + \tilde{\kappa}\left(\hbar\frac{\lambda(\tau)}{\lambda(\sigma)}\mu^{\frac12}\right)\right]
\end{align*}
The desired bounds for the transference operator can then be inferred as for the first term on the right in \eqref{eq:xsingin}.
\\

{\it{Contribution of the outgoing singular part, i.e., the terms \eqref{eq:xsingout1}, \eqref{eq:xsingout2}.}} This can again be handled by small modifications of the preceding cases. Consider for example the contribution of an `outgoing perpetuated singularity' \eqref{eq:xsingout2}. To handle it, we need to use a slightly modified version of Lemma~\ref{lem:weightedHolder1}, Lemma~\ref{lem:weightedHolder2}, where the oscillatory factor 
\[
e^{\pm i\nu\tau\eta^{\frac12}}
\]
gets replaced by 
\[
e^{\pm i\nu\left(\frac{\lambda(\tau)}{\lambda(\sigma)}x + \tau\right)\eta^{\frac12}},\quad x\geq 0. 
\]
This forces utilization of a different variable $\tilde{\eta}$ in the proof of these lemmas, namely 
\[
\tilde{\eta} = \nu\left(\frac{\lambda(\tau)}{\lambda(\sigma)}x + \tau\right)\left(\eta^{\frac12} - \xi^{\frac12}\right), 
\]
while the remaining steps are again identical to the ones in the original proofs. Observe in particular that in Remark~\ref{rem:lemweightedHolder1} the function $P_{>C^{-1}\tau}g$ is then replaced by 
\[
P_{>C^{-1}\cdot\left(\frac{\lambda(\tau)}{\lambda(\sigma)}x + \tau\right)}g.
\]
To obtain the last conclusion of the proposition, one notices that $\mathcal{D}_{\tau}$ acts trivially on the phases 
\[
e^{\pm i\hbar^{-1}\rho\big(x_{\sigma\cdot\frac{\lambda(\tau)}{\lambda(\sigma)}};\alpha\cdot\frac{\lambda(\tau)}{\lambda(\sigma)},\hbar\big)}, 
\]
and we have 
\begin{align*}
\mathcal{D}_{\tau}\Big(F_{l,k,j}^{(\pm)}\big(\tau,\sigma, \frac{\lambda^2(\tau)}{\lambda^2(\sigma)}\xi\big)\Big) = \big(\partial_{\tau}F_{l,k,j}^{(\pm)}\big)\big(\tau,\sigma, \frac{\lambda^2(\tau)}{\lambda^2(\sigma)}\xi\big).
\end{align*}
We carefully note that bounding $\partial_{\tau}F_{l,k,j}^{(\pm)}$ 'costs' $\hbar^{-1}$ (in light of the estimates in Def.~\ref{defi:xsingulartermsngeq2adm}), but it also gains a factor $\xi^{-\frac12}$ decay at large frequencies, and hence we can combine these two factors to obtain an admissibly singular souce term. Finally, when $\mathcal{D}_{\tau}$ lands on $e^{\pm i\nu\tau\xi^{\frac12}}$, we simply lose a factor $\xi^{\frac12}$. 
\\
The final conclusion from the proposition follows by using Prop.~\ref{prop: K operator} to show that the off-diagonal and diagonal parts of $\big[\mathcal{D}_{\tau}, \mathcal{K}^{(\hbar)}\big],\, \big(\mathcal{K}^{(\hbar)}\big)^2$ are similar to the ones of $ \mathcal{K}^{(\hbar)}$ and hence the preceding argument applies as well. See also the proof of Prop.~\ref{prop: xh good para linear}.
\end{proof}
\subsubsection{The effect of derivatives on functions with admissibly singular distorted Fourier transform at angular momentum $n, |n|\geq 2$.} Recall that in the coordinates $(\tau, R)$, the derivatives arising in the non-linear source terms are of type $\partial_\tau + \frac{\lambda_{\tau}}{\lambda}R\partial_R \pm \partial_R$. In fact, in the original $(t, r)$-coordinates, upon inspection of \eqref{non linear 1 2}, \eqref{coe eq precise 1},  \eqref{coe eq precise 2}, we encounter the schematically written quadratic null-forms
\[
g_t\cdot f_t - g_r\cdot f_r = \frac12(g_t - g_r)\cdot (f_t + f_r) + \frac12(g_t + g_r)\cdot (f_t - f_r),
\]
and we have 
\[
\lambda^{-1}\cdot (\partial_t \pm \partial_r) =  -\partial_\tau - \frac{\lambda_{\tau}}{\lambda}R\partial_R \pm \partial_R
\]
Here we establish how these operators act on admissibly singular functions, taking advantage of the preceding subsection. This gets expressed in the form of a lemma analogous to Lemma~\ref{lem:singFouriertiphysicalngeq2adm}:
\begin{lemma}\label{lem:singFouriertiphysicalngeq2admDeriv} Assume that $f(\tau, R)$ is an angular momentum $n, |n|\geq 2$ function, represented by 
	\[
	f(\tau, R) = \int_0^\infty \phi_{n}(R,\xi)\cdot \xb(\tau, \xi)\cdot\rho_{n}(\xi)\,d\xi, 
	\]
	with $\xb(\tau, \xi)$ admissibly singular (at angular momentum $n$). Then for the action of the `good derivative' $\partial_{\tau} + \frac{\lambda_\tau}{\lambda}R\partial_R - \partial_R$ we have the following representation on the physical side: restricting to $R<\nu\tau$ we can write 
	\begin{align*}
	\left(\partial_{\tau} + \frac{\lambda_\tau}{\lambda}R\partial_R - \partial_R\right)f =  f_1 + f_2+f_3
	\end{align*}
	where $f_1 = f_1(\tau, R)$ is a $C^5$-function supported in $\nu\tau - R\gtrsim 1$ and satisfying 
	\[
	\nabla_R^{k_2}f_1(\tau, R)\lesssim \hbar^{-\frac32+}\left(\log\tau\right)^{N_1}\cdot \tau^{-\frac52-\nu}\left|\nu\tau - R\right|^{-5},\quad 0\leq k_2\leq 5, 
	\]
	while $f_2 =\sum_{l=1}^8 f_{2l}$ where we have the explicit form
	\begin{align*}
	&f_{2l}(\tau, R) = \chi_{|\nu\tau-R|\lesssim \hbar}\sum_{k=1}^N\sum_{i=0}^{N_1}\frac{G_{k,l,i}(\tau,\nu\tau - R)}{\tau^{\frac12}}\hbar^{-\frac{l+3}{2}}\left[\nu\tau - R\right]^{\frac{l}{2}+k\nu}\left(\log\left(\nu\tau - R\right)\right)^i\\
	\end{align*}
	Here the function $G_{k,l,i}(\tau, x)$ has symbol type behavior with respect to $x$, as follows:
	\begin{align*}
	\hbar^{\delta_1\cdot (k_2+1)}\left|\partial_x^{k_2}G_{k,l,i}(\tau, x)\right| \lesssim  \left(\log\tau\right)^{N_1-i}\cdot\tau^{-2-\nu}x^{-k_2},\quad 0\leq k_2\leq 19 - l
	\end{align*}
	and we have the Holder type bound 
	\begin{align*}
	\hbar^{(20-l)\delta_1 + \delta}\cdot\left\|x^{19-l+\delta}\partial_x^{19-l}G_{k,l,i}(\tau, x)\right\|_{\dot{C}^{\delta}} \lesssim \left(\log\tau\right)^{N_1-i}\cdot \tau^{-2-\nu}.
	\end{align*}
	Finally, the remaining function $f_3$ is also $C^5$ and supported in $|\nu\tau - R|\lesssim 1$ and satisfies
	\[
	\left\|f_3\right\|_{\Sh_{0}}\lesssim \left(\log\tau\right)^{N_1}\cdot\tau^{-2-\nu}.
	\]
	Moreover, we have the bounds 
	\begin{align*}
	&\left|\partial_R^k f_3\right|\lesssim \left(\log\tau\right)^{N_1}\tau^{-\frac52-\nu}\cdot\hbar^{-\frac12+}\min\{\left(\nu\tau - R\right)^{-k},\hbar^{-k}\},\quad 
	0\leq k\leq 5. 
	\end{align*}
	For the action of the `bad derivative' $\partial_{\tau} + \frac{\lambda_\tau}{\lambda}R\partial_R + \partial_R$, we have analogous conclusion except that the expression $\hbar^{-\frac{l+3}{2}}\cdot\left[\nu\tau - R\right]^{\frac{l}{2}+k\nu}$ is replaced by 
	\[
	\hbar^{-\frac{l+1}{2}}\cdot\left[\nu\tau - R\right]^{\frac{l}{2}+k\nu-1},
	\]
	while $19 - l$ gets replaced by $20-l$. Furthermore, all estimates lose one power of decay with respect to $\tau$. 
\end{lemma}
\begin{remark}\label{rem:singFouriertiphysicalngeq2admDeriv} This lemma shows that applying the 'good' derivative $\partial_{\tau} + \frac{\lambda_\tau}{\lambda}R\partial_R - \partial_R$ in some sense leaves admissibly singular terms unchanged, although the differentiated functions needs to be interpreted as 'source admissibly singular'. This is then consistent with the slightly changed differentiability properties of the coefficient functions $G_{l,k,i}$. 
\end{remark}
\begin{proof}  The idea is to decompose the derivative operator into 
	\begin{align*}
	\partial_{\tau} + \frac{\lambda_\tau}{\lambda}R\partial_R \mp \partial_R &= \partial_{\tau} + \left[\left(1+\nu\right)\mp 1\right]\partial_R + \left(1+\nu^{-1}\right)\frac{R-\nu\tau}{\tau}\partial_R
	\end{align*}
	Fixing the $-$-sign first, this becomes 
	\[
	\partial_{\tau} +\nu\partial_R +  \left(1+\nu^{-1}\right)\frac{R-\nu\tau}{\tau}\partial_R. 
	\]
	The operator $\partial_{\tau} +\nu\partial_R$ acts trivially on functions of the form $f(\nu\tau - R)$. Thus if we recall Lemma~\ref{lem:singFouriertiphysicalngeq2adm}, the conclusion of the first part of the present lemma follows easily for the contribution arising by applying  $\partial_{\tau} +\nu\partial_R$, while the contribution arising upon applying $(1+\nu^{-1})\frac{R-\nu\tau}{\tau}\partial_R$ is handled by invoking the symbol behavior of the coefficients 
	\[
	G_{k,l,i}(\tau, \nu\tau - R)
	\]
	with respect to the second variable. The case of the $+$-sign is handled similarly. 
\end{proof}
Due to the technical difficulty of precisely translating the physical properties of functions to the Fourier properties, we give a more precise statement for the action of $\partial_{\tau} + \frac{\lambda_\tau}{\lambda}R\partial_R + \partial_R$ on the physical realisation of the principal ingoing part:
\begin{lemma}\label{lem:DerivonPrincSing} Assume that $f(\tau, R)$ is an angular momentum $n, |n|\geq 2$ function, represented by 
	\[
	f(\tau, R) = \int_0^\infty \phi_{n}(R,\xi)\cdot \xb(\tau, \xi)\cdot\rho_{n}(\xi)\,d\xi, 
	\]
	with $\xb(\tau, \xi)$ given by the first term on the right in \eqref{eq:xsingin}. Then we have 
	\[
	\left(\partial_{\tau} + \frac{\lambda_\tau}{\lambda}R\partial_R +\partial_R\right) f(\tau, R)|_{R<\nu\tau} =  \big(\int_0^\infty \phi_{n}(R,\xi)\cdot \yb(\tau, \xi)\cdot\rho_{n}(\xi)\,d\xi\big)|_{R<\nu\tau} 
	\]
	where $\yb = c\cdot\xi^{\frac12}\cdot \xb + \sum_{j=1}^2\yb_j$, with $c = c(\nu)$ a suitable constant, $\yb_1$ a source admissibly singular term with vanishing principal part, while $\yb_2$ is in $\Sh_{0}$, satisfying the bound 
	\begin{align*}
	\big\|\yb_1\big\|_{sourceadm}\lesssim \big\| \xb\big\|_{adm},\,\big\|\yb_2\big\|_{\Sh_{0}}\lesssim \tau^{-2-\nu}\cdot(\log\tau)^{N_1}\cdot  \big\| \xb\big\|_{adm}.
	\end{align*}
\end{lemma} 
\begin{proof} Write 
	\begin{align*}
	\left(\partial_{\tau} + \frac{\lambda_\tau}{\lambda}R\partial_R +\partial_R\right) f &= \left(\partial_{\tau} + \frac{\lambda_\tau}{\lambda}R\partial_R -\partial_R\right) f + 2\partial_R f\\
	& =  \left(\partial_{\tau} + \frac{\lambda_\tau}{\lambda}R\partial_R -\partial_R\right) f + \frac{2}{\nu}\big(\partial_{\tau}+\nu\partial_R\big) f -  \frac{2}{\nu}\cdot \partial_{\tau}f. 
	\end{align*}
	The first two terms at the end fall under the purview of the preceding lemma, and so it suffices to consider the last term $  \frac{2}{\nu}\cdot \partial_{\tau}f$, which as far as 'incoming admissibly singular functions' are concerned we do on the distorted Fourier side, keeping in mind Def.~\ref{defi:xsingulartermsngeq2adm}. It suffices to observe that (here $\sigma' = \sigma\cdot\frac{\lambda(\tau)}{\lambda(\sigma)}$)
	\begin{align*}
	\partial_{\tau}\Big(e^{\pm i\hbar^{-1}\rho\left(x_{\sigma\cdot\frac{\lambda(\tau)}{\lambda(\sigma)}};\alpha\cdot\frac{\lambda(\tau)}{\lambda(\sigma)},\hbar\right)}\Big) = c_1\cdot \frac{\sigma'}{\tau}\cdot \xi^{\frac12}\cdot\rho_x\left(x_{\sigma\cdot\frac{\lambda(\tau)}{\lambda(\sigma)}};\alpha\cdot\frac{\lambda(\tau)}{\lambda(\sigma)},\hbar\right) + c_2\tau^{-1}\cdot \xi^{\frac12}\cdot \rho_{\alpha}\left(x_{\sigma\cdot\frac{\lambda(\tau)}{\lambda(\sigma)}};\alpha\cdot\frac{\lambda(\tau)}{\lambda(\sigma)},\hbar\right)
	\end{align*}
	Using Lemma~\ref{lem: Lemma 3.4 CDST}, we can write
	\begin{align*}
	\chi_{\hbar^2\xi\gtrsim 1}\cdot \frac{\sigma'}{\tau}\cdot \xi^{\frac12}\cdot\rho_x\left(x_{\sigma\cdot\frac{\lambda(\tau)}{\lambda(\sigma)}};\alpha\cdot\frac{\lambda(\tau)}{\lambda(\sigma)},\hbar\right) = \tau^{-1}\cdot \xi^{\frac12}\cdot \langle \hbar\xi^{\frac12}\rangle^{-2}\cdot H\big(\tau, \sigma, \frac{\lambda^2(\tau)}{\lambda^2(\sigma)}\xi\big), 
	\end{align*}
	where $H$ satisfies the same deivative bounds as the function $F$ in  Def.~\ref{defi:xsingulartermsngeq2adm}. The same applies to the expression 
	\[
	\tau^{-1}\cdot \xi^{\frac12}\cdot \rho_{\alpha}\left(x_{\sigma\cdot\frac{\lambda(\tau)}{\lambda(\sigma)}};\alpha\cdot\frac{\lambda(\tau)}{\lambda(\sigma)},\hbar\right).
	\]
	Further, we note (with $F_{k,l,j}$ again as in Def.~\ref{defi:xsingulartermsngeq2adm}) we have 
	\begin{align*}
	\partial_{\tau}\Big(F_{k,l,j}\big(\tau, \sigma, \frac{\lambda^2(\tau)}{\lambda^2(\sigma)}\xi\big)\Big) = \xi^{\frac12}\cdot (\hbar\xi^{\frac12})^{-1}\cdot \big(\hbar\partial_{\tau}F_{k,l,j}\big)\big(\tau, \sigma, \frac{\lambda^2(\tau)}{\lambda^2(\sigma)}\xi\big) + c_3\cdot \tau^{-1}\big((\xi\partial_{\xi})F_{k,l,j}\big)\big(\tau, \sigma, \frac{\lambda^2(\tau)}{\lambda^2(\sigma)}\xi\big),
	\end{align*}
	and the functions $ \big(\hbar\partial_{\tau}F_{k,l,j}\big)\big(\tau, \sigma, \frac{\lambda^2(\tau)}{\lambda^2(\sigma)}\xi\big)$, $\big((\xi\partial_{\xi})F_{k,l,j}\big)\big(\tau, \sigma, \frac{\lambda^2(\tau)}{\lambda^2(\sigma)}\xi\big)$ satisfy the coefficient bounds as required for source admissibly singular terms. We conclude that if $\overline{x}$ is admissibly singular of incoming type, when $\partial_{\tau}$ falls either on the oscillatory phase $(e^{\pm i\hbar^{-1}\rho\left(x_{\sigma\cdot\frac{\lambda(\tau)}{\lambda(\sigma)}};\alpha\cdot\frac{\lambda(\tau)}{\lambda(\sigma)},\hbar\right)}$ or on a coefficient function $F_{k,l,j}\big(\tau, \sigma, \frac{\lambda^2(\tau)}{\lambda^2(\sigma)}\xi\big)$, or also on a coefficient $a_{k,j}^{(\pm)}(\tau,\sigma)$, we arrive at a source admissibly singular term of non-principal type. On the other hand, if $\partial_{\tau}$ falls on $e^{\pm i\nu\tau\xi^{\frac12}|}$, we clearly obtain the Fourier coefficient $c\cdot\xi^{\frac12}\cdot\overline{x}$. 
	\\
	As far as outgoing admissibly singular $\overline{x}$ are concerned, we use the 'physical representation' of such functions via Lemma~\ref{lem:singFouriertiphysicalngeq2adm}, Remark~\ref{rem:outgoing smoothing}, as well as Lemma~\ref{lem:singPhysicaltoFourierngeq2} to translate the physcial realization back to the distorted Fourier side. It easily follows that these functions are all source admissibly singular of non-principal type, up to errors of type $\overline{y}_2$.  	
\end{proof}
\begin{remark}\label{rem:lem:DerivonPrincSing} The physical realizations 
\[
\big(\int_0^\infty \phi_n(;\xi)\cdot \overline{y}_2(\tau,\xi)\cdot \rho_n(\xi)\,d\xi\big)|_{R<\nu\tau}
\]
admit pointwise bounds like the function $f_3$ in Lemma~\ref{lem:singFouriertiphysicalngeq2admDeriv}. 
\end{remark}
\subsection{Multilinear estimates near the light cone with singular inputs at angular momentum $|n|\geq 2$}
In this subsection, we finally control the source terms near the light cone with singular inputs, which is particularly delicate for the null-form source terms. In a first stage, we shall assume that all factors (inputs) are angular momentum $|n|\geq 2$ functions, as the general case will be a rather straightforward extension of this case. However, at this stage we only consider the source terms at angular momentum $|n|\geq 2$, as the exceptional angular momenta $n= 0,\pm 1$ will be treated in a separate section at the end. 
\\
We note that the null-form estimates are delicate, since a priori the regularity we are dealing with is only $H^{1+}$, whence we are at the limit of the strong local well-posedness regime. While the general theory requires the use of $H^{s,\delta}$ spaces in this setting, see \cite{KlMach}, we can and have to take advantage of the very particular structure of our solutions involving the shock on the light cone, which turns out to be very naturally adapted to the null-form structure of the most singular source terms. 
\subsubsection{Basic product estimates for angular momentum $|n|\geq 2$ functions with admissibly singular distorted Fourier transform}
Here we show that our concept of admissible singularity leads to good product estimates, and that these concepts are also compatible with forming paraproducts. These will arise naturally when proving the basic null-form estimates needed to handle the source terms.

\begin{proposition}\label{prop:admsinbproduct1} Let $n_j, j = 1,2,3$ obey the same conditions as in Prop.~\ref{prop:bilin2} . Assume that the functions $f_j(\tau, R), j = 1,2$ are angular momentum $n_j, |n_j|\geq 2$ functions admitting representations
	\[
	f_j(\tau, R) = \int_0^\infty \phi_{n_{j}}(R,\xi)\cdot \xb_j(\tau, \xi)\cdot \rho_{n_{j}}(\xi)\,d\xi,\quad j = 1,2, 
	\]
	where the distorted Fourier transforms $\xb_j,\,j = 1,2$ each can be written as $\xb_j = \yb_j + \zb_j$ with $\yb_j$ admissibly singular (at angular momentum $n_j$) and $\zb_j\in S_0^{\hbar_j}$. Then we have\footnote{We have to restrict things away from the spatial origin to avoid complications involving the vanishing order to ensure the function is regular as an angular momentum $n_3$ function. This is not important since near the origin, admissibly singular functions are smooth and decay rapidly with respect to time.}
	\[
	\chi_{R\gtrsim 1}\prod_{j=1,2}f_j|_{R<\nu\tau} = g(\tau, R),
	\]
	where, $g$ admits an angular momentum $n_3$ representation, with $\big|n_3 - \sum_{j=1,2} n_j\big|\leq O(1)$, $|n_3|\geq 2$, 
	\[
	g(\tau, R) = \int_0^\infty \phi_{n_{3}}(R,\xi)\xb_3(\tau, \xi)\rho_{n_{3}}(\xi)\,d\xi,
	\]
	where $\xb_3 = \yb_3  + \tilde{\yb}_3 + \zb_3$, and where $\yb_3$ is admissibly singular of principal ingoing type, $\tilde{\yb}_3$ is of prototypical singular type with vanishing principal part, and $\zb_3\in S_{0}^{\hbar_{3}}$. If we further assume that\footnote{We fomulate the decay for the regular part in this way since this will be the decay rate we will work with in the sequel. For the purposes of this proposition, weaker decay rates would work as well.}
	\[
	\big\|\zb_j(\tau,\cdot)\big\|_{S_{0}^{\hbar_{j}}} + \big\|\mathcal{D}_{\tau}^{(\hbar_j)}\zb_j(\tau,\cdot)\big\|_{S_{1}^{\hbar_{j}}} \lesssim \tau^{-3},\,j = 1, 2, 
	\]
	then we can make the quantitative bounds (recall Def.~\ref{defy: protofunctionnorm} )
	\[
	\big\| \yb_3\big\|_{adm} + \big\| \tilde{\yb}_3\big\|_{proto}\lesssim \min_{j = 1,2}\{|n_j|^C\}\cdot \prod_{j=1}^2\big(\big\|\yb_j\big\|_{adm} + \sup_{\tau\geq \tau_0}\tau^3\cdot \big\|\zb_j(\tau,\cdot)\big\|_{S_{0}^{\hbar_{j}}}\big)
	\]
	Furthermore, we can decompose $ \zb_3 =  \zb_{31} +  \zb_{32}$ where 
	\[
	\sup_{\tau\geq \tau_0}\tau^3\big(\big\|\zb_{31}\big\|_{S_{0}^{\hbar_{3}}} + \big\|\mathcal{D}_{\tau}\zb_{31}\big\|_{S_{1}^{\hbar_{3}}}\big)\lesssim \min_{j = 1,2}\{|n_j|^C\}\cdot \prod_{j=1}^2\big(\big\|\yb_j\big\|_{adm} + \sup_{\tau\geq \tau_0}\tau^3\cdot \big(\big\|\zb_j(\tau,\cdot)\big\|_{S_{0}^{\hbar_{j}}} +  \big\|\mathcal{D}_{\tau}^{(\hbar_j)}\zb_j(\tau,\cdot)\big\|_{S_{1}^{\hbar_{j}}}\big)\big),
	\]
	while the physical realisation of $ \zb_{32}$, when restricted to $R<\nu\tau$, satisfies the same bound as $f_3$ in Lemma~\ref{lem:singFouriertiphysicalngeq2adm}, with an extra factor 
	\[
	\min_{j = 1,2}\{|n_j|^C\}\cdot \prod_{j=1}^2\big(\big\|\yb_j\big\|_{adm} + \sup_{\tau\geq \tau_0}\tau^3\cdot \big(\big\|\zb_j(\tau,\cdot)\big\|_{S_{0}^{\hbar_{j}}} +  \big\|\mathcal{D}_{\tau}^{(\hbar_j)}\zb_j(\tau,\cdot)\big\|_{S_{1}^{\hbar_{j}}}\big)\big)
	\]
	on the right. 
\end{proposition}
\begin{remark}\label{rem:prop:admsinbproduct1} The reason for the distinction between $ \yb_3,\tilde{\yb}_3$ is that we do not truncate the physical realization of the function corresponding to the principal singular part, as we need to maintain its precise structure, while we use Lemma~\ref{lem:singFouriertiphysicalngeq2adm} and Lemma~\ref{lem:singPhysicaltoFourierngeq2} to deal with products of admissibly singular functions at least one of which is of non-principal type. 
\end{remark}
\begin{proof}
	Assume that $\hbar_2\gg\hbar_1$, say, whence the product is at angular momentum $n_3$ with $n_3\simeq n_1$, the remaining case being handled similarly. In addition, first, assume that both $\yb_{1,2}$ have vanishing principal ingoing part.  
	Write 
	\[
	f_j(\tau, R)  = \tilde{f}_j(\tau, R) + \tilde{\tilde{f}}_j(\tau, R),\,j = 1,2, 
	\]
	where we define 
	\[
	\tilde{f}_j(\tau, R) =  \int_0^\infty \phi_{n_{j}}(R,\xi)\cdot \yb_j(\tau, \xi)\cdot \rho_{n_{j}}(\xi)\,d\xi,
	\]
	 i.e., corresponding to the singular part. Using Lemma~\ref{lem:singFouriertiphysicalngeq2adm}, as well as Lemma~\ref{lem:singPhysicaltoFourierngeq2}, the conclusion follows readily for the product of the singular parts\footnote{This means the products of two factors of type $f_2$ in Lemma~\ref{lem:singFouriertiphysicalngeq2adm}.} of the $\tilde{f}_j$, and well as and a simple version of the basic product estimates such as Prop.~\ref{prop:bilin1} gives the conclusion for the product of the $\tilde{\tilde{f}}_j$, or for the products of the regular parts\footnote{This means the products of two factors of type $f_{1,3}$ in Lemma~\ref{lem:singFouriertiphysicalngeq2adm}.} of the $\tilde{f}_j$. It remains to consider the mixed case, i.e., the products
	\[
	\tilde{\tilde{f}}_1\cdot \tilde{f}_2,\,\tilde{f}_1\cdot\tilde{\tilde{f}}_2,
	\]
	or the products of a regular term and a singular term in the decomposition in Lemma~\ref{lem:singFouriertiphysicalngeq2adm}. Both situations are handled analogously. \\
{\it{The product $\tilde{\tilde{f}}_1\cdot \tilde{f}_2$}}. To begin with, we can reduce $ \tilde{f}_2$ to $\chi_{\nu\tau-R< \hbar_1} \tilde{f}_2$, by means of the following technical 
\begin{lemma}\label{lem:inproofofpropadmsinbproduct1} We have 
	\begin{align*}
	\tilde{\tilde{f}}_1\cdot \chi_{\nu\tau-R\geq \hbar_1}\tilde{f}_2\in \tilde{S}_0^{(\hbar_3)},\quad \left(\partial_{\tau}+\frac{\lambda'}{\lambda}R\partial_{R}\right)\left( \tilde{\tilde{f}}_1\cdot \chi_{\nu\tau-R\geq \hbar_1}\tilde{f}_2\right)\in \tilde{S}_1^{(\hbar_3)}
	\end{align*}
	Moreover, we have the bounds 
	\begin{align*}
	\big\|\big(\tilde{\tilde{f}}_1\cdot \chi_{\nu\tau-R\geq \hbar_1}\tilde{f}_2\big)(\tau,\cdot)\big\|_{\tilde{S}_0^{(\hbar_3)}}\lesssim \tau^{-3}\cdot |n_2|^C\cdot X,
	\end{align*}
	whee $X$ denotes the product of norms at the end of the statement of Prop.~\ref{prop:admsinbproduct1}, and similarly for the second product. 
\end{lemma}
\begin{proof}
	(lemma) We sketch the argument, which is very similar to the one for Prop.~\ref{prop:bilin1}: for the undifferentiated term, we need to bound 
	\begin{align*}
	\left\|\left\langle \phi_{n_{3}}(R,\xi),\, \tilde{\tilde{f}}_1\cdot \chi_{\nu\tau-R\geq \hbar_1}\tilde{f}_2\right\rangle_{L^2_{R\,dR}}\right\|_{S_0^{\hbar_3}}
	\end{align*}
	Labeling $\xi_1$ the frequency in the angular momentum $n_1$ Fourier representation of $\tilde{\tilde{f}}_1$, due to $\hbar_3\simeq \hbar_1$ the case $\xi_1\geq \xi$ is easily handled by invoking Lemma~\ref{lem:singFouriertiphysicalngeq2adm} to bound the $L^\infty$-norm of $\chi_{\nu\tau-R\geq \hbar_1}\tilde{f}_2$ (for instance, the method to handle the case \emph{(1)} in the proof of Proposition \ref{prop:bilin1} can be directly applied here). It remains to deal with the case $\xi_1<\xi$, where we have to perform integration by parts. For this, in the low frequency regime $\hbar_3^2\xi<1$, write (recall \eqref{diag operator redefine})
	\begin{align*}
	\left\langle \phi_{n_{3}}(R,\xi),\, \tilde{\tilde{f}}_1\cdot \chi_{\nu\tau-R\geq \hbar_1}\tilde{f}_2\right\rangle_{L^2_{R\,dR}} = \frac{1}{\xi}\left\langle H_{n_3}\phi_{n_{3}}(R,\xi),\, \tilde{\tilde{f}}_1\cdot \chi_{\nu\tau-R\geq \hbar_1}\tilde{f}_2\right\rangle_{L^2_{R\,dR}}
	\end{align*}
	Then the terms\footnote{The situation where one derivative hits each of $\tilde{\tilde{f}}_1, \tilde{f}_2$ is handled similarly.}
	\begin{align*}
	\frac{1}{\xi}\left\langle\phi_{n_{3}}(R,\xi),\, \partial_R^2\left(\tilde{\tilde{f}}_1\right)\cdot \chi_{\nu\tau-R\geq \hbar_1}\tilde{f}_2\right\rangle_{L^2_{R\,dR}},\quad \frac{1}{\xi}\left\langle\phi_{n_{3}}(R,\xi),\, \frac{n_3^2}{R^2}\left(\tilde{\tilde{f}}_1\right)\cdot \chi_{\nu\tau-R\geq \hbar_1}\tilde{f}_2\right\rangle_{L^2_{R\,dR}}
	\end{align*}
	are bounded by means of Prop.~\ref{prop:singularmultiplier} (by again placing $\frac{n_3^2}{R^2}\tilde{\tilde{f}}_1, \partial_R^2\tilde{\tilde{f}}_1$ into $L^2$ and $ \chi_{\nu\tau-R\geq \hbar_1}\tilde{f}_2$ into $L^\infty$). In fact, we have\footnote{We denote by $P_{<\xi_1}$ the Fourier localizer restricting the frequency (here for an angular momentum $n_1$ function) to size $ <\xi$.}
	\begin{align*}
	&\Big\|\frac{1}{\xi}\left\langle\phi_{n_{3}}(R,\xi),\, \partial_R^2\left(P_{<\xi}\tilde{\tilde{f}}_1\right)\cdot \chi_{\nu\tau-R\geq \hbar_1}\tilde{f}_2\right\rangle_{L^2_{R\,dR}}\Big\|_{S_0^{(\hbar_3)}(\xi\hbar_3^2\lesssim1)}\\
	&\lesssim \hbar_3^{2-\delta}\cdot\Big(\sum_{\lambda<\hbar_3^{-2}}\big[\sum_{\lambda_1<\lambda_1}\frac{\lambda_1}{\lambda}\cdot\big\|P_{\lambda_1}\tilde{\tilde{f}}_1\big\|_{L^2_{R\,dR}}\big]^2]\Big)^{\frac12}\cdot \big\|\chi_{\nu\tau-R\geq \hbar_1}\tilde{f}_2\big\|_{L^\infty_{dR}}\\
	&\lesssim \hbar_2^{-1}\cdot \big\|\tilde{\tilde{f}}_1\big\|_{S_0^{(\hbar_1)}}\cdot \big\|\overline{y}_2\big\|_{adm}, 
	\end{align*}
	where we have taken advantage of the assumption $\hbar_3\sim \hbar_1$ and the Cauchy-Schwarz inequality as well as orthogonality to bound the square sum by means of $\big\|\tilde{\tilde{f}}_1\big\|_{S_0^{(\hbar_1)}}$. 
	\\
	The term 
	\begin{align*}
	\frac{1}{\xi}\left\langle\phi_{n_{3}}(R,\xi),\, \tilde{\tilde{f}}_1\cdot \partial_R^2\left(\chi_{\nu\tau-R\geq \hbar_1}\tilde{f}_2\right)\right\rangle_{L^2_{R\,dR}}
	\end{align*}
	is estimated by using the bound (see Proposition \ref{prop:derivative1})
	\begin{align*}
	\left|\tilde{\tilde{f}}_1(\tau, R)\right|\lesssim \tau\hbar_1^{-\frac12-\delta}\cdot\left\|\tilde{\tilde{f}}_1\right\|_{\tilde{S}_0^{(\hbar_1)}},\quad \text{for}\quad  R\lesssim \tau, 
	\end{align*}
	in conjunction with (see Lemma \ref{lem:singFouriertiphysicalngeq2adm} and again Proposition \ref{prop:derivative1})
	\begin{align*}
	\left\|\partial_R^2\left(\chi_{\nu\tau-R\geq \hbar_1}\tilde{f}_2\right)\right\|_{L^2_{R\,dR}}\lesssim \hbar_1^{-1}\cdot\hbar_{2}^{-1+\delta}\cdot\tau^{-1-\nu}\left(\log\tau\right)^{N_1}. 
	\end{align*}
Note that the negative powers in $\hbar_{1}$ are canceled by positive powers in $\hbar_{3}$ in considering the $\left\|\cdot\right\|_{S_{0}^{\hbar_{3}}}$-norm. To deal with the high frequency regime $\xi>\hbar_3^{-2}$, further integration by parts are required, which can be handled analogously to the preceding. 	
\end{proof}
It remains to deal with $\tilde{\tilde{f}}_1\cdot \chi_{\nu\tau-R<\hbar_1}\tilde{f}_2$. Here it suffices to split 
\[
\tilde{\tilde{f}}_1(\tau, R) = \sum_{j=0}^3 \left(\nu\tau - R\right)^j\cdot  \frac{\tilde{\tilde{f}}_1^{(j)}(\tau, \nu\tau)}{j!} +  \tilde{\tilde{g}}_1(\tau, R) = :P_3 \tilde{\tilde{f}}_1(\tau, R) + \tilde{\tilde{g}}_1(\tau, R)
\]
Also, we may assume that in accordance with Lemma~\ref{lem:singFouriertiphysicalngeq2adm} we have (with $l\geq 1$)
\[
\chi_{\nu\tau-R<\hbar_1}\tilde{f}_2(\tau, R) = \chi_{|\nu\tau-R|\lesssim \hbar_1}\sum_{k=1}^N\sum_{i=0}^{N_1}\frac{G_{k,l,i}(\tau,\nu\tau - R)}{\tau^{\frac12}}\hbar_2^{-\frac{l+1}{2}}\left[\nu\tau - R\right]^{\frac{l}{2}+k\nu}\left(\log\left(\nu\tau - R\right)\right)^i
\]
with bounds as stated there for the coefficients $G_{k,l,i}(\tau,\nu\tau - R)$. We have the bounds 
\begin{align*}
\left|\tilde{\tilde{f}}_1^{(j)}(\tau, \nu\tau)\right| &= \left|\partial_R^j\int_0^\infty \phi_{n_{1}}(R,\xi)\cdot \zb_1(\tau,\xi)\rho_{n_{1}}(\xi)\,d\xi\right|_{R = \nu\tau}\\
& \lesssim \tau\cdot \hbar_1^{-\frac12-j}\cdot\left\|\zb_1(\tau,\cdot)\right\|_{\Sho_{0}},\quad 0\leq j\leq 3,
\end{align*}
\begin{align*}
\left|\left(\partial_{\tau}+\frac{\lambda'}{\lambda}R\partial_{R}\right)\tilde{\tilde{f}}_1^{(j)}(\tau, \nu\tau)\right|\lesssim \hbar_1^{-\frac32-j}\cdot\left(\left\|\zb_{1}(\tau,\cdot)\right\|_{\Sho_{0}}+\left\|\calD_{\tau}\zb_1(\tau,\cdot)\right\|_{\Sho_{1}}\right),\quad 0\leq j\leq 3.
\end{align*}

Then the product of the pure polynomial of the singular part is given by 
\begin{align*}
&\left(P_3 \tilde{\tilde{f}}_1(\tau, R)\right)\cdot   \chi_{\nu\tau-R<\hbar_1}\tilde{f}_2(\tau, R)\\
&=  \chi_{|\nu\tau-R|\lesssim \hbar_1}\sum_{j=0}^3\sum_{k=1}^N\sum_{i=0}^{N_1}\frac{\tilde{G}_{k,l,j,i}(\tau,\nu\tau - R)}{\tau^{\frac12}}\hbar_1^{-\frac{l+2j+1}{2}}\left[\nu\tau - R\right]^{\frac{l+2j}{2}+k\nu}\left(\log(\nu\tau - R)\right)^i ,
\end{align*}
where we set 
\[
\tilde{G}_{k,l,j,i}(\tau,\nu\tau - R) = \hbar_1^{\frac{2j+l+1}{2}}\cdot \hbar_2^{-\frac{1+l}{2}}\cdot \frac{1}{j!}\tilde{\tilde{f}}_1^{(j)}(\tau, \nu\tau)\cdot G_{k,l,i}(\tau,\nu\tau - R)
\]
and which is seen to satisfy the same bounds as the coefficients in $f_{2(l+2j)}$ in Lemma~\ref{lem:singFouriertiphysicalngeq2adm} up to a factor $\hbar_2^{-C}$. In fact, the gain of $\hbar^{j+\frac{l+1}{2}}$ more than compensates for the loss of $\hbar^{-\frac12-j}$ occurring when estimating $\tilde{\tilde{f}}_1^{(j)}(\tau, \nu\tau)$. 
Next we consider the product 
	\[
	E(\tau, R): = \tilde{\tilde{g}}_1\cdot  \chi_{\nu\tau-R<\hbar_1}\tilde{f}_2
	\]
	We claim that this function can be placed into $\tilde{S}_0^{(\hbar_1)}$, with time derivative in $\tilde{S}_1^{(\hbar_1)}$. This can be seen by using Taylor's theorem:
	\[
	\tilde{\tilde{g}}_1(\tau, R) = \int_0^1 \tilde{\tilde{f}}^{(4)}_{1}(\tau, \nu\tau + s(R - \nu\tau))\cdot(\nu\tau - R)^4\,ds = \tilde{\tilde{f}}_{1}(\tau, R) - P_3\tilde{\tilde{f}}_{1}(\tau, R).
	\]
	We also have 
	\begin{align*}
	\partial_R^k \tilde{\tilde{g}}_1(\tau, R) = \partial_R^k \tilde{\tilde{f}}_{1}(\tau, R) - P_{3-k}\big(\partial_R^k\tilde{\tilde{f}}_{1}\big)(\tau, R),\,0\leq k\leq 4,
	\end{align*}
	where we set $P_{-1}f = 0$, and the integral relation 
	\begin{align*}
	\partial_R^k \tilde{\tilde{g}}_1(\tau, R) = \int_0^1\big(\partial_R^4\tilde{\tilde{f}}_{1}\big)(\tau, \nu\tau + s(R - \nu\tau)))\cdot(\nu\tau - R)^{4-k}\,ds,\,0\leq k\leq 3. 
	\end{align*}
	It follows that 
	\begin{align*}
	\partial_R^5 E(\tau, R) = \sum_{p = 0}^5 C_p\cdot \partial_R^p \tilde{\tilde{g}}_1(\tau, R)\cdot \partial_R^{5-p}\big(\chi_{\nu\tau-R<\hbar_1}\tilde{f}_2\big),
	\end{align*}
	and we have 
	\begin{align*}
	&\hbar_3^{5+}\cdot \Big\| \partial_R^5 \tilde{\tilde{g}}_1(\tau, R)\cdot\chi_{\nu\tau-R<\hbar_1}\tilde{f}_2\Big\|_{H_{R\,dR}^{0+}}\lesssim \big\|\tilde{\tilde{f}}_{1}\big\|_{\tilde{S}_0^{(\hbar_1)}}\cdot \big\|\chi_{\nu\tau-R<\hbar_1}\tilde{f}_2\big\|_{L^\infty\cap W^{\infty,0+}},\\
	&\hbar_3^{5+}\cdot \Big\| \partial_R^4 \tilde{\tilde{g}}_1(\tau, R)\cdot\partial_R\big(\chi_{\nu\tau-R<\hbar_1}\tilde{f}_2\big)\Big\|_{H_{R\,dR}^{0+}}\lesssim \hbar_3^{5+}\cdot\big\| \partial_R^4 \tilde{\tilde{g}}_1(\tau,\cdot)\big\|_{L^\infty\cap W^{\infty,0+}}\cdot \big\|\partial_R\big(\chi_{\nu\tau-R<\hbar_1}\tilde{f}_2\big)\big\|_{H_{R\,dR}^{0+}}\\
	&\hspace{6cm}\lesssim \hbar_2^{-1}\cdot  \big\|\tilde{\tilde{f}}_{1}\big\|_{\tilde{S}_0^{(\hbar_1)}}.
	\end{align*}
	For $0\leq p\leq 3$ we get 
	\begin{align*}
	\partial_R^p \tilde{\tilde{g}}_1(\tau, R)\cdot \partial_R^{5-p}\big(\chi_{\nu\tau-R<\hbar_1}\tilde{f}_2\big) =  \int_0^1\big(\partial_R^4\tilde{\tilde{f}}_{1}\big)(\tau, \nu\tau + s(R - \nu\tau)))\,ds\cdot(\nu\tau - R)^{4-p}\cdot \partial_R^{5-p}\big(\chi_{\nu\tau-R<\hbar_1}\tilde{f}_2\big), 
	\end{align*}
	and we have the estimates 
\begin{align*}
 &\hbar_1^{5+}\cdot \Big\|\int_0^1\big(\partial_R^4\tilde{\tilde{f}}_{1}\big)(\tau, \nu\tau + s(R - \nu\tau)))\,ds\Big\|_{L^\infty_{RdR}\cap W^{\infty, 0+}_{R\,dR}}\lesssim  \big\|\tilde{\tilde{f}}_{1}\big\|_{\tilde{S}_0^{(\hbar_1)}}, \\
 &\Big\|(\nu\tau - R)^{4-p}\cdot \partial_R^{5-p}\big(\chi_{\nu\tau-R<\hbar_1}\tilde{f}_2\big)\Big\|_{H_{R\,dR}^{0+}}\lesssim \hbar_2^{-1}\cdot\tau^{-1-\nu}\cdot\big(\log\tau\big)^{N_1}.
\end{align*}
This again yields via the fractional Leibniz rule that 
\begin{align*}
\hbar_3^{5+}\cdot \Big\| \sum_{p = 0}^3 C_p\cdot \partial_R^p \tilde{\tilde{g}}_1(\tau, R)\cdot \partial_R^{5-p}\big(\chi_{\nu\tau-R<\hbar_1}\tilde{f}_2\big)\Big\|_{H^{0+}_{R\,dR}}\lesssim  \hbar_2^{-1}\cdot  \big\|\tilde{\tilde{f}}_{1}\big\|_{\tilde{S}_0^{(\hbar_1)}}.
\end{align*}


{\it{The product $\tilde{f}_1\cdot \tilde{\tilde{f}}_2$}}. Here the large angular momentum term $\tilde{f}_1$ is of singular type, while the small angular momentum term is of the smooth type. The argument then proceeds analogously to the preceding case by splitting $\tilde{\tilde{f}}_2$ into its third order Taylor polynomial centered at $R = \nu\tau$ and a remainder term. 
\\

Now we assume that one of $\yb_j$, say $\yb_1$, is of principal ingoing part. Then if $\yb_2$ is also of principal ingoing type, write 
\[
\tilde{f}_j(\tau, R) = \chi_{\nu\tau - R\geq \hbar_j}\tilde{f}_j(\tau, R)  +  \chi_{\nu\tau - R<\hbar_j}\left(\tilde{f}_j(\tau, R) - c_j(\tau)\right) + \chi_{\nu\tau - R<\hbar_j} c_j(\tau),\quad c_j(\tau) = \tilde{f}_j(\tau, \nu\tau),\quad j = 1,2,
\]
where we let 
\[
\tilde{f}_j(\tau, R) = \int_0^\infty \phi_{n_{j}}(R,\xi)\cdot \yb_j(\tau, \xi)\cdot\rho_{n_{j}}(\xi)\,d\xi. 
\]
Here $\chi_{\nu\tau - R\geq \hbar_j}\tilde{f}_j(\tau, R)\in \tilde{S}_0^{(\hbar_j)}$, and furthermore the preceding argument implies, using the fact that $\hbar_{1}\simeq\hbar_{3}$ and the product is in the regime $\nu\tau-R\gtrsim \hbar_{3}$,
\[
\chi_{\nu\tau - R\geq \hbar_j}\tilde{f}_j(\tau, R)\cdot  \chi_{\nu\tau - R<\hbar_k}\tilde{f}_k(\tau, R) \in \tilde{S}_0^{(\hbar_3)},\quad \{j,k\} = \{1,2\}. 
\]
More precisely, using Lemma~\ref{lem:singFouriertiphysicalngeq2adm}, we see that we can write the preceding product in the form of the function $f_3$ of that lemma and satisfying the bounds stated there with $\hbar$ replaced by $\hbar_3$, and an extra factor 
\[
\hbar_2^{-\frac12+}\cdot \prod_{j=1,2}\big\|\overline{y}_j\big\|_{adm}
\]
on the right. 
 It follows that to complete analysis of the product $\prod_{j=1,2}\tilde{f}_j(\tau, R)$, it suffices to consider 
\begin{align*}
&\prod_{j=1,2}\left[ \chi_{\nu\tau - R<\hbar_j}\left(\tilde{f}_j(\tau, R) - c_j(\tau)\right) + \chi_{\nu\tau - R<\hbar_j} c_j(\tau)\right]\\
& = \prod_{j=1,2}\chi_{\nu\tau - R<\hbar_j}\left(\tilde{f}_j(\tau, R) - c_j(\tau)\right)\\
& + \sum_{\{j,k\} = \{1,2\}}\chi_{\nu\tau - R<\hbar_j}\left(\tilde{f}_j(\tau, R) - c_j(\tau)\right)\cdot \chi_{\nu\tau - R<\hbar_k} c_k(\tau)\\
& + \prod_{j=1,2}\chi_{\nu\tau - R<\hbar_j} c_j(\tau)
\end{align*}
Using Lemma~\ref{lem:singFouriertiphysicalngeq2adm} and Lemma~\ref{lem:singPhysicaltoFourierngeq2}, the first term on the right is seen to be of prototype connecting singular type or smoother. More precisely, using the argument of Lemma~\ref{lem:inproofofpropadmsinbproduct1} to handle smooth/non-smooth products, we can write 
\begin{align*}
 \prod_{j=1,2}\chi_{\nu\tau - R<\hbar_j}\left(\tilde{f}_j(\tau, R) - c_j(\tau)\right)\big|_{<\nu\tau} = \big(g_2 + g_3\big)\big|_{<\nu\tau},
\end{align*}
 where $g_2$ has the same structure as the term $f_2$ in Lemma~\ref{lem:singFouriertiphysicalngeq2adm} with the added restriction $l\geq 2$ and $\hbar = \hbar_3$, while $g_3$ is like the term $f_3$, and furthermore, the same bounds obtain as in that lemma, up to an extra factor
 \[
\hbar_2^{-\frac12+}\cdot \prod_{j=1,2}\big\|\overline{y}_j\big\|_{adm}
\]
on the right. \\
For the second term above on the right, considering the case $j = 1, k=2$, say, we have 
\begin{align*}
&\chi_{\nu\tau - R<\hbar_1}\left(\tilde{f}_1(\tau, R) - c_1(\tau)\right)\cdot \chi_{\nu\tau - R<\hbar_2} c_2(\tau)\\
& = \chi_{\nu\tau - R<\hbar_1}\tilde{f}_1(\tau, R)\cdot \chi_{\nu\tau - R<\hbar_2} c_2(\tau) -  \chi_{\nu\tau - R<\hbar_1}c_1(\tau)\cdot \chi_{\nu\tau - R<\hbar_2} c_2(\tau)\\
& = \tilde{f}_1(\tau, R)\cdot c_2(\tau)\\
& - \chi_{\nu\tau - R\geq \hbar_1}\tilde{f}_1(\tau, R)\cdot c_2(\tau)\\
& - \chi_{\nu\tau - R<\hbar_1}\tilde{f}_1(\tau, R)\cdot \chi_{\nu\tau - R\geq\hbar_2} c_2(\tau)\\
&-  \chi_{\nu\tau - R<\hbar_1}c_1(\tau)\cdot \chi_{\nu\tau - R<\hbar_2} c_2(\tau).\\ 
\end{align*}
Here the last three terms, call them $X_j$, $j = 1, 2, 3$, are easily seen to be in $\tilde{S}_0^{(\hbar_3)}$ (Note that the third term on the RHS above actually vanishes, since $\hbar_{2}\gg\hbar_{1}$), satisfying the bounds
\begin{align*}
\big\|X_j\big\|_{\tilde{S}_0^{(\hbar_3)}}\lesssim \hbar_2^{-\frac12+}\cdot \prod_{j=1,2}\big\|\overline{y}_j\big\|_{adm},
\end{align*}
 while for the first of the last four terms, this is clearly a function which is of admissibly singular principal ingoing type {\it{when interpreted as angular momentum $n_1$ function}}. However, this function needs to be interpreted as an angular momentum $n_3$ function, for which we need a `translation device'. To begin with, using Lemma~\ref{lem:singFouriertiphysicalngeq2adm}, we easily conclude that 
\[
\chi_{R<\frac{\nu\tau}{2}}\tilde{f}_1(\tau, R)\cdot c_2(\tau)\in \tilde{S}_0^{(\hbar_3)}. 
\]
It thus suffices to understand 
\[
\left\langle \phi_{n_{3}}(R,\xi),\,\chi_{R\geq\frac{\nu\tau}{2}}\tilde{f}_1(\tau, R)\cdot c_2(\tau)\right\rangle_{L^2_{R\,dR}}
\]
For this we use the following lemma:
\begin{lemma}\label{lem:n1ton2translate} Let $f(R)$ be an angular momentum $n_1$ function, $|n_1|\geq 2$, let $|n_2|\geq 2$, and setting 
	\[
	f(R) = \int_0^\infty\phi_{n_{1}}(R,\xi) \xb(\xi)\rho_{n_{1}}(\xi)\,d\xi, 
	\]
	let the function $\mathcal{K}^{n_1,n_2}_{\tau}\xb$ be defined by the relation 
	\[
	\left(\mathcal{K}^{n_1,n_2}_{\tau}\xb\right)(\eta): = \left\langle \phi_{n_{2}}(R,\eta),\,\chi_{R\geq \frac{\nu\tau}{2}}f(R)\right\rangle_{L^2_{R\,dR}}. 
	\]
	Then we have in analogy to the transference operator the distributional identity 
	\[
	\mathcal{K}^{n_1,n_2}_{\tau}(\xi,\eta) = \frac{\overline{ a_{n_{2}}(\eta) }}{\overline{ a_{n_{1}}(\xi) }}\delta(\xi - \eta) + \tilde{\mathcal{K}}^{n_1,n_2}_{\tau}(\xi,\eta),
	\]
	where the operator $\tilde{\mathcal{K}}^{n_1,n_2}_{\tau}$ acts via integration against a kernel $\frac{F(\xi,\eta;\tau,n_1,n_2)\rho_{n_{1}}(\eta)}{\xi-\eta}$
	\[
	\left(\tilde{\mathcal{K}}^{n_1,n_2}_{\tau}f\right)(\eta) = \int_0^\infty \frac{F(\xi,\eta;\tau,n_1,n_2)\rho_{n_{1}}(\xi)}{\xi-\eta}f(\xi)\,d\xi. 
	\]
	The kernel function $F(\xi,\eta;\tau,n_1,n_2)$ can be decomposed as
	\begin{align*}
		F(\xi,\eta;\tau,n_1,n_2)=F_{P}(\xi,\eta;\tau,n_{1},n_{2})+F_{N}(\xi,\eta;\tau,n_{1},n_{2})
	\end{align*}
	where $F_{P}(\xi,\eta;\tau,n_{1},n_{2})$ satisfies (assuming $\xi\leq\eta$) for $k\leq k_{0}(n_{1})$,
	\begin{align}\label{FP n1 n2 bound 1}
		\begin{split}
		\left|F_{P}(\xi,\eta;\tau,n_{1},n_{2})\right|\lesssim P(|n_{1}-n_{2}|)\left(\hbar_{1}\xi^{\frac12}\right)^{-1}\min\left\{1,\left(\hbar_{1}\xi^{\frac12}\right)\right\}\cdot G,
		\end{split}
	\end{align}
	with
	\begin{align}\label{FP n1 n2 bound 1 auxi}
		G:=\begin{cases}
		&\min\left\{1,\left(\hbar_{1}\xi^{\frac12}\right)^{\frac14}\left|\eta^{\frac12}-\xi^{\frac12}\right|^{-\frac14}\right\},\quad \text{for}\quad \left|\frac{\eta^{\frac12}}{\xi^{\frac12}}-1\right|\lesssim 1,\\
		&\hbar_{1}\left(\hbar_{1}\xi^{\frac12}\right)^{-1}\min\left\{1,\hbar_{1}\xi^{\frac12}\right\}\cdot \left(\frac{\xi^{\frac12}}{\eta^{\frac12}}\right)^{k}\cdot\left(\frac{\langle\xi\rangle}{\eta}\right)^{N},\quad \text{for}\quad \left|\frac{\eta^{\frac12}}{\xi^{\frac12}}-1\right|\gg1,
		\end{cases}
	\end{align}
	for arbitrary $N>0$ and $\eta\geq 1$.
	Here $P(\cdot)$ is a polynomial. For the derivatives of $F_{P}(\xi,\eta;\tau,n_{1},n_{2})$, we have
	\begin{align}\label{FP deri 1}
	\begin{split}
		&\left|\partial_{\xi^{\frac12}}F_{P}(\xi,\eta;\tau,n_{1},n_{2})\right|,\quad \left|\partial_{\eta^{\frac12}}F_{P}(\xi,\eta;\tau,n_{1},n_{2})\right|\\
		\lesssim &\frac{P(|n_{1}-n_{2}|)}{\left(\hbar_{1}\xi^{\frac12}\right)^{\frac12}\left(\hbar_{2}\eta^{\frac12}\right)^{\frac12}}\left(1+\left|\log\left(\hbar_{1}\xi^{\frac12}\left|\xi^{\frac12}-\eta^{\frac12}\right|^{-1}\right)\right|\right)\quad \text{if}\quad \xi\simeq\eta,
		\end{split}
	\end{align}
	The operator corresponding to $F_{N}(\xi,\eta;n_{1},n_{2},\tau)$ has a stronger smoothing effect, which maps an admissibly singular function to a smooth function.
\end{lemma}
\begin{proof}
	The proof is similar to that of Proposition \ref{prop: K operator}, and here we only give an outline. According to its definition, the operator $\left(\mathcal{K}^{n_1,n_2}_{\tau}\xb\right)(\eta)$ is given by
	\begin{align*}
	\left(\mathcal{K}^{n_1,n_2}_{\tau}\xb\right)(\eta)=\int_{0}^{\infty}\left\langle\chi_{R\geq\frac{\nu\tau}{2}}\phi_{n_{1}}(R,\xi),\phi_{n_{2}}(R,\eta)\right\rangle_{L^{2}_{R\,dR}}\cdot\rho_{n_{1}}(\xi)f(\xi)\,d\xi.
	\end{align*}
	In view of \eqref{eq:phi wieder2}, we infer that the $\delta$ measure on the diagonal in the integral
	\begin{align*}
		\lim_{A\rightarrow\infty}\int_{0}^{A}\chi_{R\geq\frac{\nu\tau}{2}}\phi_{n_{1}}(R,\xi)\phi_{n_{2}}(R,\eta)R\,dR
	\end{align*}
	comes from the expression
	\begin{align*}
		&2\pi^{-1}\left(\xi\eta\right)^{-\frac14}\lim_{L\rightarrow\infty}\Re\int_{0}^{\infty}e^{i\left(\hbar_{1}^{-1}\Psi(R,\xi,\hbar_{1})-\hbar_{2}^{-1}\Psi(R,\eta,\hbar_{2})\right)}a_{n_{1}}(\xi)\left(1+\hbar a_{1}(-\tau,\alpha)\right)\\
		&\cdot\overline{ a_{n_{2}}(\eta) }\left(1+\hbar\overline{  a_1(-\sigma;\beta) }\right)\chi_{R\geq\frac{\nu\tau}{2}}\chi_{2}(R/L)\,dR.
	\end{align*}
	Here 
	\begin{align*}
		&\alpha:=\hbar_{1}\xi^{\frac12},\quad  \Psi(R,\xi,\hbar_{1}):=\hbar_{1}\xi^{\frac12}R-y(\alpha,\hbar_{1})+\rho(\hbar_{1}\xi^{\frac12}R,\alpha,\hbar_{1})\\
		&\beta:=\hbar_{2}\eta^{\frac12},\quad \Psi(R,\eta,\hbar_{2}):=\hbar_{2}\eta^{\frac12}R-y(\beta,\hbar_{2})+\rho(\hbar_{2}\eta^{\frac12}R,\alpha,\hbar_{2}),
	\end{align*}
	and $\chi_{2}$ is defined the same as in \eqref{eq:Ixieta}. The result for the diagonal part follows in a similar way as in the proof for Proposition \ref{prop: K operator}.
	\\
	Now we turn to the off-diagonal kernel $\tilde{\calK}^{n_{1},n_{2}}_{\tau}(\xi,\eta)$. By a routing calculation similar to the one in Proposition \ref{prop: K operator}, we have
	\begin{align*}
		\eta\left(\calK^{n_{1},n_{2}}_{\tau}f\right)(\eta)-\calK^{n_{1},n_{2}}_{\tau}\left(\xi f(\xi)\right)(\eta)=&-\int_{0}^{\infty}\left\langle\chi_{R\geq\frac{\nu\tau}{2}}\phi_{n_{1}}(R,\xi),H_{n_{2}}^{+}\phi_{n_{2}}(R,\eta)\right\rangle_{L^{2}_{R\,dR}}\rho_{n_{1}}(\xi)f(\xi)\,d\xi\\
		&+\int_{0}^{\infty}\left\langle \chi_{R\geq\frac{\nu\tau}{2}}\phi_{n_{1}}(R,\xi),H_{n_{1}}^{+}\phi_{n_{2}}(R,\eta)\right\rangle_{L^{2}_{R\,dR}}\rho_{n_{1}}(\xi)f(\xi)\,d\xi\\
		&+\int_{0}^{\infty}\left\langle\left(\partial_{R}^{2}+R^{-1}\partial_{R}\right)\chi_{R\geq\frac{\nu\tau}{2}}\cdot\phi_{n_{2}}(R,\eta),\phi_{n_{1}}(R,\xi)\right\rangle_{L^{2}_{R\,dR}}\rho_{n_{1}}(\xi) f(\xi)\,d\xi	\\
		&+\int_{0}^{\infty}\left\langle\partial_{R}\chi_{R\geq\frac{\nu\tau}{2}}\partial_{R}\phi_{n_{2}}(R,\eta),\phi_{n_{1}}(R,\xi)\right\rangle_{L^{2}_{R\,dR}}\rho_{n_{1}}(\xi)f(\xi)\,d\xi.
	\end{align*}
	Therefore the function $F(\xi,\eta;\tau, n_{1}, n_{2})$ is given by
\begin{align}\label{F n1 n2}
\begin{split}
	F(\xi,\eta;\tau,n_{1},n_{2})=&\left\langle\left(H^{+}_{n_{1}}-H^{+}_{n_{2}}\right)\phi_{n_{2}}(R,\eta),\chi_{R\geq\frac{\nu\tau}{2}}\phi_{n_{1}}(R,\xi)\right\rangle_{L^{2}_{R\,dR}}\\&+\left\langle\left(\partial_{R}^{2}+R^{-1}\partial_{R}\right)\chi_{R\geq\frac{\nu\tau}{2}}\cdot\phi_{n_{2}}(R,\eta),\phi_{n_{1}}(R,\xi)\right\rangle_{L^{2}_{R\,dR}}\\
	&+2\left\langle\partial_{R}\chi_{R\geq\frac{\nu\tau}{2}}\cdot \partial_{R}\phi_{n_{2}}(R,\eta),\phi_{n_{1}}(R,\xi)\right\rangle_{L^{2}_{R\,dR}}.
	\end{split}
\end{align}
We denote the first term on the RHS of \eqref{F n1 n2} by $F_{P}(\xi,\eta; n_{1}, n_{2},\tau)$ and the rest three terms by $F_{N}(\xi,\eta; n_{1},n_{2},\tau)$.
Note that $F_{P}(\xi,\eta; n_{1}, n_{2},\tau)$ has a similar structure as the function $F(\xi,\eta;\hbar)$ in Proposition \ref{prop: K operator}, except that 1) the potential $H_{n_{1}}^{+}-H_{n_{2}}^{+}$ decays as $\langle R\rangle^{-2}$ only, 2) the coefficient of the potential is $n_{1}^{2}-n_{2}^{2}=(n_{1}-n_{2})(n_{1}+n_{2})$, and 3) the Fourier basis $\phi_{n_{1}}(R,\xi)$ and $\phi_{n_{2}}(R,\eta)$ are for different angular modes. If $\hbar_{1}\leq \hbar_{2}$ (the vice-versa being similar), then 
\begin{align*}
	\hbar_{2}=\hbar_{1}\cdot\frac{\hbar_{2}}{\hbar_{1}}=\hbar_{1}\cdot\frac{n_{1}+1}{n_{2}+1}=\hbar_{1}\cdot \left(1+\frac{n_{1}-n_{2}}{n_{2}+1}\right).
\end{align*}
Therefore the estimates on $F_{P}(\xi,\eta;n_{1},n_{2},\tau)$ follows in a similar way as in Proposition \ref{prop: K operator}.
\end{proof}
Using the preceding lemma, we infer that 
\begin{align*}
\left\langle \phi_{n_{3}}(R,\xi),\,\chi_{R\geq\frac{\nu\tau}{2}}\tilde{f}_1(\tau, R)\cdot c_2(\tau)\right\rangle_{L^2_{R\,dR}} = c_2(\tau)\cdot \yb_1(\tau) + \tilde{\yb}_1(\tau,\xi), 
\end{align*}
where $ \tilde{\yb}_1(\tau,\xi)$ is of connecting singular type, and hence more regular. This completes the outline of the proof of the proposition. 
\end{proof}

\begin{remark}\label{rem:prop:admsinbproduct1} The same conclusion applies if we replace $\overline{y}_j$ in the statement of the preceding proposition by a linear combination of an admissibly singular term and a prototype singular term of connecting type (i. e. non-principal type). For this is suffices to use Lemma~\ref{lem:singFouriertiphysicalngeq2} in addition to Lemma~\ref{lem:singFouriertiphysicalngeq2adm}. We can also weaken our assumptions on the regular parts $\overline{z}_j$, and admit functions $\tilde{\overline{z}}_j\in S_0^{(\hbar_j)}$ whose 'physical realization' is like the function $f_3$ in Lemma~\ref{lem:singFouriertiphysicalngeq2adm}, at angular momentum $n_j$. This is a consequence of the preceding proof, where such functions arise when representing the singular parts $\overline{y}_j$ on the physical side. 
We also importantly note that the preceding proof implies that if the $\overline{y}_j$ are of restricted admissibly singular type, then so will be the principal part of the product. 
\end{remark}
Applying the preceding proposition inductively, we can then also handle more general products: 
\begin{corollary}\label{cor:admsinbproduct2} Let $\phi_j,\,j = 1,2,\ldots,k$, $k\geq 2$, be angular momentum $n_j,\,|n_j|\geq 2$ functions whose distorted Fourier transform $\overline{x}_j$ (in the angular momentum $n_j$-sense) admits a description as in the preceding proposition. Let the angular momenta $m_l,\,l = 1,\ldots, k-1$, be determined such that $|m_l|\geq 2$ and  $\{m_2,n_1,n_2\}$ satisfy the conditions in of  Prop.~\ref{prop:bilin2} (instead of  satisfying the conditions $\{n_3,n_1,n_2\}$ in that proposition), and similarly for the triples $\{m_r, n_r,m_{r-1}\}$. Then the product 
	\[
	\left(\prod_{j=1}^k\phi_j\right)|_{<\nu\tau} = g|_{R<\nu\tau}
	\]
	where $g$ may be represented as an angular momentum $m = m_{k-1}$-function whose Fourier transform is a linear combination of an admissible function, a prototype singular function of connecting type, and a function in $S_0^{\tilde{\hbar}_{k-1}}$, where $\tilde{\hbar  }_l: = \frac{1}{m_{l}+1}$. If all $\phi_j$ have restricted admissibly singular part (in terms of their distorted Fourier representation), then so does $g$.  More quantitatively, we can write 
	\begin{align*}
	g = \sum_{j=1}^3 g_j,\,g_j(\tau, R) = \int_0^\infty \phi_{m}(R;\xi)\cdot\overline{p}_j(\tau,\xi)\cdot\rho_m(\xi)\,d\xi,
	\end{align*} 
	such that (with $\hbar = |m|^{-1}$, and $C, D$ suitable absolute constants at most depending on $\nu$)
	\begin{align*}
	\left\|\overline{p}_1\big\|_{adm} + \tau^3\cdot \big\|\overline{p}_2(\tau,\cdot)\right\|_{S_0^{\hbar}} + \tau^{1+\nu-}\cdot \left\|\overline{p}_3(\tau,\cdot)\right\|_{S_0^{\hbar}}\leq D^k\cdot \left(\max_{j=1,\ldots k}\{|n_j|^C\}\right)^{-1}\cdot \prod_{j=1}^k |n_j|^C\cdot\left(\left\|\overline{y}_j\right\|_{adm} + \sup_{\tau\geq \tau_0}\tau^3\left\|\overline{z}_j\right\|_{S_0^{\hbar_j}} \right).
	\end{align*}
	Furthermore, the function $g_3$ admits a description like the function $f_3$ in Lemma~\ref{lem:singFouriertiphysicalngeq2adm}, but where the right hand product in the immediately preceding inequality occurs as additional factor on the right. 
	\end{corollary}
	\begin{remark}\label{rem:cor:admsinbproduct2} A less precise statement is also possible, by replacing $\sup_{\tau\geq \tau_0}\tau^3\big\|\overline{z}_j\big\|$ by $\sup_{\tau\geq \tau_0}\tau^{1+}\left\|\overline{z}_j\right\|_{S_0^{\hbar_j}}$ on the right and $p_2+p_3$ by $p_4$ and the sum of the last two terms on the left in the preceding inequality by 
	\[
	 \tau^{1+}\cdot \left\|\overline{p}_4(\tau,\cdot)\right\|_{S_0^{\hbar}}.
	 \]
	The argument is essentially the same. 
	\end{remark}

\subsubsection{Fourier localization on admissibly singular functions}\label{subsubsec:Fourierlocal}

One nice feature of our concept of admissibly singular functions is that except for the principal ingoing  singular part, they are compatible with the application of Fourier localization operators 
\[
\chi_{\xi^{\frac12}>\mu},\quad\chi_{\xi^{\frac12}<\mu},\quad  \chi_{\xi^{\frac12}\simeq\mu}. 
\]
In fact, if $F^{(\pm)}_{l,k,i}\left(\tau, \sigma,\frac{\lambda^2(\tau)}{\lambda^2(\sigma)}\xi\right)$ is as in Definition~\ref{defi:xsingulartermsngeq2adm}, then the function 
\begin{align*}
\chi_{\xi^{\frac12}>\mu}F^{(\pm)}_{l,k,i}\left(\tau, \sigma,\frac{\lambda^2(\tau)}{\lambda^2(\sigma)}\xi\right) &= \chi_{\frac{\lambda(\sigma)}{\lambda(\tau)}\cdot\left(\frac{\lambda(\tau)}{\lambda(\sigma)}\xi^{\frac12}\right)>\mu}F^{(\pm)}_{l,k,i}\left(\tau, \sigma,\frac{\lambda^2(\tau)}{\lambda^2(\sigma)}\xi\right)\\
&=:\tilde{F}^{(\pm)}_{l,k,i}\left(\tau, \sigma,\frac{\lambda^2(\tau)}{\lambda^2(\sigma)}\xi\right)
\end{align*}
is easily seen to satisfy the same estimates as the function $F^{(\pm)}_{l,k,i}\left(\tau, \sigma,\frac{\lambda^2(\tau)}{\lambda^2(\sigma)}\xi\right)$, and similarly for the other Fourier cutoffs. Moreover, these cutoffs act boundedly on $\Sh_{0}, \Sh_{1}$ for trivial reasons. 
\\
In the sequel, if $f$ is an angular momentum $n, |n|\geq 2$ function represented by 
\[
f(R) = \int_0^\infty\phi_{n}(R,\xi)\xb(\xi)\rho_{n}(\xi)\,d\xi,
\]
then we define the Fourier localization operators $P_{<\mu}, P_{\mu}, P_{>\mu}$ by means of 
\[
P_{>\mu}f(R) = \int_0^\infty\phi_{n}(R,\xi)\chi_{\xi^{\frac12}>\mu}\xb(\xi)\rho_{n}(\xi)\,d\xi,
\]
and similarly for the other operators. Thus we commit abuse of notion here, in that these operators tacitly depend on the angular momentum $n$. In the sequel, it will always be understood what the underlying angular momentum is. 
However, the principal ingoing part will be modified into inadmissible form by application of a Fourier cutoff; nonetheless, we shall be able to show that certain paraproducts of admissibly singular inputs remain admissibly singular (up to smoother errors). 

\subsubsection{Paraproduct estimates for angular momentum $|n|\geq 2$ functions with admissibly singular distorted Fourier transform}

Recall that in the expressions of the form $\sum_{\mu}P_{<\mu}f P_{\mu}g$ $\mu$ ranges over dyadic frequencies. Here we show that this concept is compatible with our concept of admissible singularity:
\begin{lemma}\label{lem:simpleparaproduct} Assume that $\phi_{1,2}$ are angular momentum $n_l$, $l = 1, 2$, functions, $|n_l|\geq 2$, and such that their $n_l$ distorted Fourier transforms $\overline{x_{n_l}}$ admit decompositions $\overline{x_{n_l}} = \overline{y_{n_l}} + \overline{z_{n_l}}$ into an admissibly singular piece $\overline{y_{n_l}}$ and a component in $S_0^{\hbar_l}$, $\hbar_l = \frac{1}{|n_l|}$. Then setting $n = \max\{2, |n_1 + n_2|\}$, the paraproduct 
\begin{align*}
\sum_{\lambda>1}P_{<\lambda}\phi_1\cdot P_{\lambda}\phi_2
\end{align*}
has angular momentum $n$-distorted Fourier transform $\overline{x}_n = \overline{y}_n + \overline{z}_n$, where $ \overline{y}_n $ is $n$-admissibly singular and $\overline{z}_n\in S_0^{\hbar}$, $\hbar = |n_2|^{-1}$. Moreover, we have the bounds 
\begin{align*}
\big\|\overline{y}_n\big\|_{adm} + \tau^{1+}\cdot\big\|\overline{z}_n\big\|_{S_0^{\hbar}}\lesssim \min\{|n_1|, |n_2|\}^C\cdot \prod_{l=1,2}\big(\big\|\overline{y}_{n_l}\big\|_{adm} + \sup_{\tau\geq \tau_0}\tau^{1+}\cdot \big\|\overline{z_{n_l}}\big\|_{S_0^{(\hbar_{n_l}})}\big). 
\end{align*}
If both functions $\phi_j$ have singular part of restricted type, so does the paraproduct. 
\end{lemma}
\begin{remark}\label{rem:lem:simpleparaproduct} A more refined statement analogous to Prop.~\ref{prop:admsinbproduct1} holds, with a similar proof. 
\end{remark}
\begin{proof} We distinguish between two cases:
\\

{\it{(1): $\phi_2$ of principal ingoing singular type.}} Here we reformulate the paraproduct as follows:
\begin{align*}
\sum_{\lambda>1}P_{<\lambda}\phi_1\cdot P_{\lambda}\phi_2 = \phi_1\cdot\phi_2 - \sum_{\lambda>1}P_{\geq\lambda}\phi_1\cdot P_{\lambda}\phi_2. 
\end{align*}
The first term on the right can be handled by means of Prop.~\ref{prop:admsinbproduct1}. For the second term on the right, we expand 
\begin{align*}
P_{\geq\lambda}\phi_1 = f_1^{(\geq\lambda)} + f_2^{(\geq\lambda)} + f_3^{(\geq\lambda)} 
\end{align*}
as in Lemma~\ref{lem:singFouriertiphysicalngeq2adm}, where we have (with $\hbar_1 = |n_1|^{-1}$)
\begin{align*}
\big\| f_1^{(\geq\lambda)} \big\|_{\tilde{S}_0^{\hbar_1}} + \big\| f_3^{(\geq\lambda)} \big\|_{\tilde{S}_0^{\hbar_1}}\lesssim \tau^{-1-\nu}\cdot\Big\|\overline{x_{n_1}}\Big\|,
\end{align*}
where for this proof  we use the notation $\Big\|\cdot\Big\|$ for the norm on the right in the statement of the lemma, while $ f_2^{(\geq\lambda)}$ admits an expansion 
\begin{align*}
f_2^{(\geq\lambda)} = \chi_{|\nu\tau - R|\lesssim \hbar_1}\sum_{l=1}^8\sum_{k=1}^{N}\sum_{j=0}^{N_1}\frac{G^{(\geq\lambda)}_{k,l,j}(\tau,\nu\tau - R)}{\tau^{\frac12}}\hbar_1^{-\frac{l+1}{2}}[\nu\tau - R]^{\frac{l}{2}+k\nu}\big(\log(\nu\tau - R)\big)^j
\end{align*}
with the coefficient functions $G^{(\geq\lambda)}_{k,l,j}(\tau,\nu\tau - R)$ satisfying the bounds stated in Lemma~\ref{lem:singFouriertiphysicalngeq2adm} with an extra factor $ \Big\|\overline{x_{n_1}}\Big\|$ on the right. Similarly, we can write 
\begin{align*}
P_{\lambda}\phi_1 = f_1^{(\lambda)} + f_2^{(\lambda)} + f_3^{(\lambda)},  
\end{align*}
where 
\begin{align*}
\big\| f_1^{(\lambda)} \big\|_{\tilde{S}_0^{\hbar_1}} + \big\| f_3^{(\lambda)} \big\|_{\tilde{S}_0^{\hbar_1}}\lesssim \tau^{-1-\nu}\cdot\Big\|\overline{x_{n_2}}\Big\|
\end{align*}
and further (recalling the assumption on $\phi_2$)
\begin{align*}
f_2^{(\lambda)} = \chi_{|\nu\tau - R|\lesssim \hbar_2}\sum_{k=1}^{N}\sum_{j=0}^{N_1}\frac{G^{(\lambda)}_{k,1,j}(\tau,\nu\tau - R)}{\tau^{\frac12}}\hbar_2^{-1}[\nu\tau - R]^{\frac{1}{2}+k\nu}\big(\log(\nu\tau - R)\big)^j.
\end{align*}
The coefficient functions $G^{(\lambda)}_{k,1,j}$ satisfy the bounds in Lemma~\ref{lem:singFouriertiphysicalngeq2adm},  Lemma~\ref{lem:lem:singFouriertiphysicalngeq2admrefined}. Then distinguish between the following cases:
\\

{\it{(1.a): Smooth-smooth interactions}}.
\\

We have the bound (using Prop.~\ref{prop:bilin1}, Lemma~\ref{lem:lem:singFouriertiphysicalngeq2admrefined})
\begin{align*}
&\Big\|\chi_{R\lesssim\tau}\big(f_1^{(\lambda)} + f_3^{(\lambda)}\big)\cdot \big( f_1^{(\geq\lambda)} + f_3^{(\geq\lambda)} \big)\Big\|_{\tilde{S}_0^{\hbar}}\\&\lesssim\max\{\hbar_1,\hbar_2\}^{-C}\cdot \tau\lambda^{-(0+)}\cdot \big(\big\| f_1^{(\lambda)} \big\|_{\tilde{S}_0^{\hbar_1}} + \big\| f_3^{(\lambda)} \big\|_{\tilde{S}_0^{\hbar_1}}\big)\cdot \big(\big\| f_1^{(\geq\lambda)} \big\|_{\tilde{S}_0^{\hbar_1}} + \big\| f_3^{(\geq\lambda)} \big\|_{\tilde{S}_0^{\hbar_1}}\big)\\
&\lesssim \max\{\hbar_1,\hbar_2\}^{-C}\cdot \tau^{-1-2\nu}\lambda^{-(0+)}\cdot \Big\|\overline{x_{n_1}}\Big\|\cdot \Big\|\overline{x_{n_2}}\Big\|.
\end{align*}
This can be summed over dyadic $\lambda>1$ to give the desired estimate.
\\

{\it{(1.b): Smooth-singular interactions}}. Consider the product 
\begin{align*}
\big(f_1^{(\geq\lambda)} + f_3^{(\geq\lambda)}\big)\cdot  f_2^{(\lambda)} . 
\end{align*}
We handle this via finite Taylor expansion (around $R = \nu\tau$) of the smooth factor $f_1^{(\geq\lambda)} + f_3^{(\geq\lambda)}$. Thus write 
\begin{align*}
f_1^{(\geq\lambda)} + f_3^{(\geq\lambda)} = f_1^{(\geq\lambda)} + f_3^{(\geq\lambda)} - P_4\big(f_1^{(\geq\lambda)} + f_3^{(\geq\lambda)}\big) + P_4\big(f_1^{(\geq\lambda)} + f_3^{(\geq\lambda)}\big). 
\end{align*}
Then we have the estimate
\begin{align*}
\Big\|\Big( f_1^{(\geq\lambda)} + f_3^{(\geq\lambda)} - P_4\big(f_1^{(\geq\lambda)} + f_3^{(\geq\lambda)}\big)\Big)\cdot f_2^{(\lambda)}\Big\|_{\tilde{S}_0^{\hbar}}\lesssim  \max\{\hbar_1,\hbar_2\}^{-C}\cdot \tau^{-1-2\nu}\lambda^{-(0+)}\cdot \big\|\overline{x_{n_1}}\big\|_{good}\cdot \big\|\overline{x_{n_2}}\big\|_{good}
\end{align*}
on account of Lemma~\ref{lem:lem:singFouriertiphysicalngeq2admrefined} and the proof of Prop.~\ref{prop:admsinbproduct1}.  Furthermore, we can write
\begin{align*}
 P_4\big(f_1^{(\geq\lambda)} + f_3^{(\geq\lambda)}\big)\cdot  f_2^{(\lambda)} &=  \chi_{|\nu\tau - R|\lesssim \hbar_2}\sum_{\alpha = 0}^4a_\alpha^{(\geq\lambda)}(\tau)\sum_{k=1}^{N}\sum_{j=0}^{N_1}\frac{G^{(\lambda)}_{k,1,j}(\tau,\nu\tau - R)}{\tau^{\frac12}}\hbar_2^{-1}[\nu\tau - R]^{\frac{1}{2}+\alpha + k\nu}\big(\log(\nu\tau - R)\big)^j\\
\end{align*}
where $a_{\alpha}{(\geq\lambda)}(\tau) = \partial_R^{\alpha}\big(f_1^{(\geq\lambda)} + f_3^{(\geq\lambda)}\big)(\tau,\nu\tau)$. Then if $\hbar_1\ll\hbar_2$, say, we can write this as 
\begin{align*}
\hbar_2^{-1}\chi_{|\nu\tau - R|\lesssim \hbar_2}\sum_{\alpha = 0}^4\hbar_1^{1+\alpha}a_\alpha{(\geq\lambda)}(\tau)\sum_{k=1}^{N}\sum_{j=0}^{N_1}\frac{G^{(\lambda)}_{k,1,j}(\tau,\nu\tau - R)}{\tau^{\frac12}}\hbar_1^{-1-\alpha}[\nu\tau - R]^{\frac{1}{2}+\alpha + k\nu}\big(\log(\nu\tau - R)\big)^j,
\end{align*}
and we have the estimate 
\begin{align*}
\lambda^{\frac12+}\cdot\Big(\big|\hbar_1^{1+\alpha}a_\alpha{(\geq\lambda)}(\tau)\big| + \big|\hbar_1^{2+\alpha}\partial_{\tau}a_\alpha{(\geq\lambda)}(\tau)\big|\Big)\lesssim \tau\cdot\big\|\overline{z}_{n_1}\big\|_{S_0^{\hbar_1}} + \tau^{-\frac32-}\cdot \big\|\overline{y}_{n_1}\big\|_{adm},\,\alpha\in \{0,1,\ldots,4\},
\end{align*}
taking advantage of Lemma~\ref{lem:lem:singFouriertiphysicalngeq2admrefined}. The preceding sum still contains the principal incoming singularities corresponding to $\alpha = 0$, but when replacing $P_{\lambda}\phi_2$ by $\lambda^{-\frac12}P_{\lambda}\phi_2$ we can effectively replace the term 
\[
\frac{G^{(\lambda)}_{k,1,j}(\tau,\nu\tau - R)}{\tau^{\frac12}}\hbar_2^{-1}[\nu\tau - R]^{\frac{1}{2}+k\nu}
\]
by 
\[
\frac{\tilde{G}^{(\lambda)}_{k,1,j}(\tau,\nu\tau - R)}{\tau^{\frac12}}\hbar_2^{-1}[\nu\tau - R]^{1+k\nu},
\]
which is no longer of principal incoming type. We note that the case $\hbar_2\lesssim\hbar_1$ is simpler since we don't have to carefully track the polynomial losses in $\hbar_1^{-1}$. 
\\

{\it{(1.c): Singular-smooth interactions}}.  Here we consider the product 
\[
f_2^{(\geq\lambda)}\cdot \big(f_1^{(\lambda)} + f_3^{(\lambda)}\big). 
\]
We note that $f_2^{(\geq\lambda)}$ is singular part of the physical realisation of the admissibly singular part of the distorted Fourier transform of $\phi_1$ (at angular momentum $n_1$), localized to frequencies $\xi\gtrsim \lambda$. On the other hand, $f_1^{(\lambda)} + f_3^{(\lambda)}$ is the smooth part of the physical realisation of $\phi_2$, localized to frequency $\xi\sim\lambda$. We then perform partial summation 
\begin{align*}
\sum_{\lambda>1}f_2^{(\geq\lambda)}\cdot \big(f_1^{(\lambda)} + f_3^{(\lambda)}\big) = f_2\cdot (f_1 + f_3) - \sum_{\lambda>1}f_2^{(<\lambda)}\cdot \big(f_1^{(\lambda)} + f_3^{(\lambda)}\big). 
\end{align*}
The desired assertion for the first term follows from Prop.~\ref{prop:admsinbproduct1}, while the second term on the right is treated as in case {\it{(1.b)}}. 
\\

{\it{(1.d): Singular-singular interactions.}} Finally we consider the product 
\[
\sum_{\lambda>1}f_2^{(\geq\lambda)}\cdot f_2^{(\lambda)}. 
\]
Inserting the expansions for $f_2^{(\geq\lambda)},  f_2^{(\lambda)}$ from before, we arrive at a sum involving the products 
\begin{align*}
\hbar_1^{-\frac{l+1}{2}}\cdot\hbar_2^{-\frac{l'+1}{2}}\cdot \frac{G^{(\lambda)}_{k,l,j}(\tau,\nu\tau - R)\cdot G^{(\geq \lambda)}_{k',l',j'}(\tau,\nu\tau - R)}{\tau}\cdot [\nu\tau - R]^{\frac{l+l'}{2}+(k+k')\nu}\cdot \log^{j+j'}(\nu\tau - R).
\end{align*}
Assuming, say, that $\hbar_1\lesssim\hbar_2$, we can write this alternatively in the form 
\begin{align*}
\hbar_1^{-\frac{l+l'+1}{2}}\cdot\hbar_2^{-\frac{1}{2}}\cdot \left(\frac{\hbar_1}{\hbar_2}\right)^{\frac{l'}{2}}\cdot\frac{G^{(\lambda)}_{k,l,j}(\tau,\nu\tau - R)\cdot G^{(\geq \lambda)}_{k',l',j'}(\tau,\nu\tau - R)}{\tau}\cdot [\nu\tau - R]^{\frac{l+l'}{2}+k\nu}\cdot[\nu\tau - R]^{k'\nu}\cdot \log^{j+j'}(\nu\tau - R).
\end{align*}
It is then straightforward to verify that identifying 
\begin{align*}
\tilde{G}_{k+k', l+l', j+j'}(\tau, \nu\tau - R)
\end{align*}
with the sum over all the 
\[
\hbar_2^{-\frac{1}{2}}\cdot \left(\frac{\hbar_1}{\hbar_2}\right)^{\frac{l'}{2}}\cdot\frac{G^{(\lambda)}_{k,l,j}(\tau,\nu\tau - R)\cdot G^{(\geq \lambda)}_{k',l',j'}(\tau,\nu\tau - R)}{\tau^{\frac12}}\cdot [\nu\tau - R]^{k'\nu}
\]
and given values of $k+k', l+l', j+j'$, we have the desired bound upon summing over $\lambda>1$. In fact, using Lemma~\ref{lem:lem:singFouriertiphysicalngeq2admrefined}, we can use the factor $[\nu\tau - R]^{k'\nu}$ to obtain exponential decay in $\lambda^{-1}$. 
\\

{\it{(2): $\phi_2$ not of principal ingoing singular type.}} Here one decomposes $P_{\lambda}\phi_2$ as well as $P_{<\lambda}\phi_1$ as in Lemma~\ref{lem:singFouriertiphysicalngeq2adm} and multiplies the terms directly, which is easily seen to lead to acceptable terms by following the procedure in the preceding case {\it{(1)}}.
\end{proof}

The preceding proof also reveals the following additional detail:
\begin{corollary}\label{cor:simpleparaproduct} Under the same hypotheses as for the preceding lemma, the singular part of 
\[
\sum_{\lambda>1}\lambda^{-\frac12}\cdot P_{<\lambda}\phi_1\cdot P_{\lambda}\phi_2
\]
is of level $l\geq 1$. Moreover, if the singular part of $\phi_2$ is of level $l\geq 1$ (i. e. not of principal incoming type), then so is the singular part of the preceding sum. 
\end{corollary}

\subsubsection{Iterated paraproducts}\label{subsubsec:paraproditerated} When controlling the quadratic null-forms arising in the source terms, we shall have to take advantage of an iterative process involving a kind of normal transforms, and resulting in iterated paraproducts which naturally generalize the paraproduct treated in the preceding subsection. Specifically, here we consider expressions of the form 
\begin{align*}
\Pi_{\lambda_k}(\phi_1,\phi_2,\ldots, \phi_{2k}) = \sum_{1<\lambda_1<\lambda_2<\ldots<\lambda_k}P_{<\lambda_1}\phi_1 P_{\lambda_1}\phi_2\ldots P_{<\lambda_k}\phi_{2k-1}\cdot P_{\lambda_k}\phi_{2k},
\end{align*}
where it is implied that the summation is over all tuples $\{\lambda_1, \lambda_2,\ldots, \lambda_{k-1}\}$ satisfying the indicated constraints, while $\lambda_k$ is held fixed. The $\lambda_j$ range over dyadic numbers $2^{\mathbb{Z}}$. 

\begin{lemma}\label{lem:longparaprod1} Assume that the functions $\phi_l$ are angular momentum $n_l$ functions, $|n_l|\geq 2$, and such that their $n_l$ distorted Fourier transforms $\overline{x}_{n_l}$ admit decompositions $\overline{x}_{n_l}= \overline{y}_{n_l} + \overline{z}_{n_l}$ into admissibly an singular piece $\overline{y_{n_l}}$ and a component in $S_0^{\hbar_l}$, $\hbar_l = \frac{1}{|n_l|}$. Then if $n = \sum_l n_l$, $|n|\geq 2$, the iterated paraproduct $\Pi(\phi_1,\phi_2,\ldots, \phi_{2k}): = \sum_{\lambda_k>1}\Pi_{\lambda_k}(\phi_1,\phi_2,\ldots, \phi_{2k})$ has distorted Fourier transform $\overline{x_n}$ admitting a decomposition into an $n$ admissibly singular part $\overline{y_n}$ and a part $\overline{z_n}\in S_0^{\hbar}$, $\hbar = \frac{1}{|n|}$. Moreover, we have the norm bound 
\begin{align*}
\Big\|\overline{x}_n\Big\|\leq D^k\cdot \max\{|n_l|\}^{-C}\cdot \prod_{l=1}^{2n}\big(|n_l|^C\cdot \Big\|\overline{x_{n_l}}\Big\|\big).
\end{align*}
Here we set as in the proof of the preceding lemma
\begin{align*}
\Big\|\overline{x_{n_l}}\Big\|: = \big\|\overline{y}_{n_l} \big\|_{adm} + \tau^{1+}\cdot \big(\big\|\overline{z}_{n_l}\big\|_{S_0^{\hbar_l}} + \big\|\mathcal{D}_{\tau}\overline{z}_{n_l}\big\|_{S_1^{\hbar_l}}\big),\,\hbar_l = |n_l|^{-1}. 
\end{align*}
\end{lemma}
\begin{proof} Let us introduce the notation $\Pi_{<\lambda}(\phi_1,\phi_2,\ldots, \phi_{2k}) = \sum_{\lambda_k<\lambda}\Pi_{\lambda_k}(\phi_1,\phi_2,\ldots, \phi_{2k})$ and analogously for $\Pi_{\geq \lambda}(\phi_1,\phi_2,\ldots, \phi_{2k})$. We need to show that 
\begin{equation}\label{eq:iteratedparaproduct}\begin{split}
&\Pi_{<\lambda_k}(\phi_1,\phi_2,\ldots, \phi_{2k-2})\cdot P_{<\lambda_k}\phi_{2k-1}\cdot P_{\lambda_k}\phi_{2k}\\
& = \Pi(\phi_1,\phi_2,\ldots, \phi_{2k-2})\cdot P_{<\lambda_k}\phi_{2k-1}\cdot P_{\lambda_k}\phi_{2k} - \Pi_{\geq \lambda_k}(\phi_1,\phi_2,\ldots, \phi_{2k-2})\cdot P_{<\lambda_k}\phi_{2k-1}\cdot P_{\lambda_k}\phi_{2k}\\
\end{split}\end{equation}
when summed over $\lambda_k>0$ is as asserted. We show inductively that the assertion of the lemma holds. According to the preceding lemma and corollary, the assertion holds when $k = 1$, so we now assume them for some $k-1\geq 1$. Let $\Pi_{\mu}: = \Pi_{\mu}\big(\phi_1, \phi_2,\ldots, \phi_{2k-2}\big),\,\Pi = \sum_{\mu>1}\Pi_{\mu}$. Then if $\phi_{2k}$ is {\it{not}} of principal incoming singular type, we can decompose both 
\[
\Pi_{<\lambda_k},\, P_{<\lambda_k}\phi_{2k-1}\cdot P_{\lambda_k}\phi_{2k},\,
\]
according to Lemma~\ref{lem:singFouriertiphysicalngeq2adm} and multiply out the individual terms, which easily gives the desired conclusion (after summation over $\lambda_k>1$) by exploiting exponential gains in $\lambda_k$ as in cases {\it{(1.a) - (1.d)}} before. Let us then assume that $\phi_{2k}$ is of principal incoming type. Using \eqref{eq:iteratedparaproduct}, the induction hypothesis, the preceding lemma as well as Prop.~\ref{prop:admsinbproduct1}, imply that 
\begin{equation}\label{eq:non-paraterm}
\sum_{\lambda_k>1}\Pi(\phi_1,\phi_2,\ldots, \phi_{2k-2})\cdot P_{<\lambda_k}\phi_{2k-1}\cdot P_{\lambda_k}\phi_{2k} = \Pi(\phi_1,\phi_2,\ldots, \phi_{2k-2})\cdot\sum_{\lambda_k>1}P_{<\lambda_k}\phi_{2k-1}\cdot P_{\lambda_k}\phi_{2k}.
\end{equation}
satisfies the conclusion of the lemma. For the remaining term, we replicate the argument for case {\it{(1.b)}} from the proof of Lemma~\ref{lem:simpleparaproduct}. The factor $ D^k\cdot \max\{|n_l|\}^{-C}\cdot \prod_{l=1}^{2n}|n_l|^C$ arises by simple case distinctions. For example in case 
\[
|n_{2k}| = \max_{l=1,\ldots,2k}\{|n_l|\}
\]
but $|n| = \big|\sum_{l<2k-1}n_l\big|>|n_{2k}|$, application of the bilinear estimate Prop.~\ref{prop:admsinbproduct1} and induction hypothesis as well as Lemma~\ref{lem:simpleparaproduct} to \eqref{eq:non-paraterm} results in a factor 
\begin{align*}
&C_1\cdot\max_{l=1,\ldots,2k-2}\{|n_l|\}^{-C}\cdot \prod_{l=1}^{2k-2}|n_l|^C\cdot  |n_{2k-1}|^C\cdot |n_{2k}|^C\\
&\leq C_1\cdot\max_{l=1,\ldots,2k-2}\{|n_l|\}^{-C}\cdot \prod_{l=1}^{2k-2}|n_l|^C\cdot  |n_{2k-1}|^C\cdot \big|\sum_{l=1}^{2k-2}n_l\big|^C\\
&\leq C_1\cdot (2k-2)^C\cdot  \max_{l=1,\ldots,2k}\{|n_l|\}^{-C}\prod_{l=1}^{2k}|n_l|^C,
\end{align*}
and of course we can bound $C_1\cdot (2k-2)^C\leq D^k$ for suitable $D$. 
\end{proof}

\subsection{Spaces for the source terms in the nonlinearity} Our next major task is to deduce the nonlinear estimates which will allow us to control all the source terms arising in \eqref{coe eq precise 1}, \eqref{coe eq precise 2}. 
In analogy to the `good spaces' $S_1^{\hbar}$ which are used to bound the smooth source terms, we shall  take advantage of the admissibly singular source terms already introduced at the end of Definition~\ref{defi:xsingulartermsngeq2adm}. We recall this here:
\begin{definition}\label{defi:xsingularsourcetermsngeq2adm} Let $F(\tau, R)$ be an angular momentum $n, |n|\geq 2$ function. Then we say that $F$ is an admissibly singular source term, provided its distorted Fourier transform $\yb(\tau, \xi)$ (at angular momentum $n$) has the property that 
	\[
	\yb(\tau, \xi) = \yb_1(\tau, \xi) + \yb_2(\tau, \xi) 
	\]
	and where $\xi^{-\frac12}\cdot \overline{y}_1(\tau,\xi)$ is a admissibly singular except that we replace $\tau^{-1-\nu}$ by $\tau^{-2-\nu}$ and we also set $k_1 = 0$ throughout(according to the details in Definition~\ref{defi:xsingulartermsngeq2adm}); furthermore, we require $\overline{y}_2$ to be a prototype singular source term.
	In particular, we have natural concepts of principally ingoing singular part and principal singular part of restricted type, for the admissibly singular source terms as well. If the principal ingoing parts of $\yb_j,\,j = 1, 2$ vanish and the sums $\sum_{l}\ldots$ occurring in the connecting incoming and outgoing parts are restricted to $l\geq l_1\geq 1$, we say that $\yb$ is an admissibly singular source term of level $l_1$.  \\ 
	We define $\big\|\yb\big\|_{sourceadm}$ as in Definition~\ref{defi:xsingulartermsngeq2adm}. 
	\end{definition}
	\begin{remark}\label{rem:defi:xsingularsourcetermsngeq2adm} The preceding definition implies that admissibly singular source terms involve both functions 'of admissibly singular type', meaning their definition involves integrals over an auxiliary temporal variable $\sigma$, as well prototypical singular type terms without such an integral. This is necessary since applying the operator $\mathcal{D}_{\tau}$ to an admissibly singular function according to Definition~\ref{defi:xsingulartermsngeq2adm} also results in functions of 'prototype'. 
	\end{remark}
	
We shall require estimates for products of source terms and admissibly singular functions. The following lemma treats both products as well as para-products:
\begin{lemma}\label{lem:sourceadmtimesgood} Assume that the angular momentum $n_1$ function $F$ has distorted Fourier transform which is the sum of a function $\overline{x}_1$ in $S_1^{\hbar_1}$ and a source admissibly singular function $\overline{y}_1$. Then if the distorted Fourier transform of the angular momentum $n_2$-function $\phi$ is the sum of a function $\overline{x}_2$ in $S_0^{\hbar_2}$ and an admissibly singular function $\overline{y}_2$, then the product 
\[
\chi_{1\lesssim R\lesssim \nu\tau}F\cdot \phi
\]
is source admissibly singular at angular momentum $m$,  provided we have $|m|\geq 2$ and $\{m, n_1, n_2\}$ is an admissible triple (meaning $\big|m - \sum_j n_j\big| = O(1)$) Specifically, its distorted Fourier transform can be written in the form 
\begin{align*}
\langle \phi_m(R;\xi), \chi_{1\lesssim R\lesssim \nu\tau}F\cdot \phi\rangle_{L^2_{R\,dR}} = \overline{y} + \overline{z}
\end{align*}
where we have the bound (with $\hbar = |m|^{-1}$)
\begin{align*}
\big\|\overline{y}\big\|_{sourceadm} + \tau^4\cdot\big\|\overline{z}\big\|_{S_1^{\hbar}}\lesssim \min\{|n_{1,2}|\}^C\cdot \big(\big\|\overline{y}_1\big\|_{sourceadm} + \tau^4\cdot\big\|\overline{x}_1\big\|_{S_1^{\hbar}}\big)\cdot \big(\tau^3\big\|\overline{x}_2\big\|_{S_0^{\hbar_2}} + \big\|\overline{y}_2\big\|_{adm}\big).
\end{align*}

The same conclusion applies to the para-products
\[
\sum_{\lambda>1}\chi_{1\lesssim R\lesssim \nu\tau}P_{<\lambda} F\cdot P_{\lambda}\phi,\,\sum_{\lambda>1}\chi_{1\lesssim R\lesssim \nu\tau}P_{\geq\lambda} F\cdot P_{\lambda}\phi.
\]
\end{lemma}
\begin{proof} {\it{First expression $ \chi_{1\lesssim R\lesssim \nu\tau}F\cdot \phi$.}} If $F$ has distorted Fourier transform in $S_1^{\hbar_1}$ and $\phi$ in $S_0^{\hbar_2}$, this is a consequence of Prop.~\ref{prop:bilin1}. If either $F$ or $\phi$ is admissibly (source) singular, while the other factor is in $S_1^{\hbar_1}$ (resp. $S_0^{\hbar_2}$), one uses Lemma~\ref{lem:singFouriertiphysicalngeq2adm} to represent the singular term as a sum of smooth terms plus an explicit non-smooth term (namely $f_2$). To handle the product of $f_2$ with the smooth factor, one develops the latter into a third, resp. fourth order Taylor expansion around $R= \nu\tau$ plus a Cauchy error. This leads to terms which are either directly seen to be in $S_1^{\hbar}$, or else of the explicit singular type for which application of Lemma~\ref{lem:singPhysicaltoFourierngeq2} leads to the desired conclusion. \\
If the distorted Fourier transform of $F$ is source admissibly while the one of $\phi$ is admissibly singular, this is also a consequence of taking advantage of Lemma~\ref{lem:singFouriertiphysicalngeq2adm} and multiplying the individual terms. 
\\

{\it{Para-product $\sum_{\lambda>1}\chi_{1\lesssim R\lesssim \nu\tau}P_{<\lambda} F\cdot P_{\lambda}\phi$.}} If $F\in S_1^{\hbar_1}$ while $\phi\in S_0^{\hbar_2}$, then this is a consequence of Prop.~\ref{prop:bilin1}. If $F$ is source admissibly singular while $\mathcal{F}(\phi)\in S_0^{\hbar_2}$, write 
\[
P_{\lambda}\phi = P_4\big(P_{\lambda}\phi\big)+ \big[P_{\lambda}\phi - P_4\big(P_{\lambda}\phi\big)\big],
\]
where $P_4\big(P_{\lambda}\phi_2\big)$ is the fourth order Taylor expansion (with respect to $R$) of $P_{\lambda}\phi_2$ around $R = \nu\tau$. Since $F_{<\lambda}$ is also source-admissibly singular, we have that 
\begin{align*}
 \big[P_{\lambda}\phi - P_4\big(P_{\lambda}\phi\big)\big]\cdot P_{<\lambda}F\in \tilde{S}_1^{\hbar}
\end{align*}
where $\hbar = \frac{1}{|m|}$. Moreover, frequency localizing the product as an angular momentum $m$-function and using integration by parts if needed (as in the proof of Prop.~\ref{prop:bilin1}), we have 
\begin{align*}
&\big\|\langle \phi_m(R;\xi),\,P_{\mu}\Big(\big[P_{\lambda}\phi - P_4\big(P_{\lambda}\phi\big)\big]\cdot P_{<\lambda}F\Big)\rangle\Big\|_{S_1^{\hbar}}\\&\lesssim \tau^{-4}\min\{\frac{\mu}{\lambda}, \frac{\lambda}{\mu}\}\cdot \min\{n_{1,2}\}^C\cdot \big\|\mathcal{F}(F)\big\|_{sourceadm}\cdot \big\|\mathcal{F}(P_{\lambda}\phi)\big\|_{S_0^{\hbar_2}}. 
\end{align*} 
The desired bound follows by square-summing over $\mu$ and taking advantage of the Cauchy-Schwarz inequality.
\\
On the other hand, the product 
\begin{align*}
P_4\big(P_{\lambda}\phi\big)\cdot P_{<\lambda} F
\end{align*}
is source admissibly singular (using Lemma~\ref{lem:singFouriertiphysicalngeq2adm}, Lemma~\ref{lem:singPhysicaltoFourierngeq2}) and we have the bound
\begin{align*}
\big\|\langle \phi_m(R;\xi),\,P_4\big(P_{\lambda}\phi\big)\cdot P_{<\lambda} F\rangle\big\|_{sourceadm}\lesssim \min\{|n_{1,2}|\}^C\cdot \lambda^{-(0+)}\cdot \big\|\mathcal{F}(F)\big\|_{sourceadm}\cdot \big\|\mathcal{F}(\phi)\big\|_{S_0^{\hbar_2}}. 
\end{align*}
The remaining cases when $\phi_2$ is admissibly singular while $F$ is either in $S_1^{\hbar_1}$ or admissibly singular are handled analogously. The assertion concerning the last para-product follows also from the preceding, as it is the difference of $\chi_{1\lesssim R\lesssim \nu\tau}F$ and the first para-product. 
\end{proof}


The preceding subsections furnish the functional framework to estimate the many source terms occurring in \eqref{coe eq tau R}. Due to the complexity of the many terms and the norms involved, we do this in several stages. To begin with, we restrict to expressions exclusively defined in terms of $|n|\geq 2$ angular momentum functions. It will not be difficult to generalise things to more general 'input' functions, once we have defined the suitable space analogues for the singular parts of the exceptional angular momenta. We shall also assume that each 'input' function has distorted Fourier transform which can be decomposed into a part in $S_0^{\hbar}$ with suitable temporal decay, as well as an admissibly singular part. This will have to be generalised a bit later on in order to arrive at a functional framework which is also compatible with the wave parametrix. 
The main difficulties will arise when treating the quadratic null-form terms, which we do in a separate subsection.

\subsection{Estimates for source terms not involving quadratic null-forms}

As a preliminary step to Prop.~\ref{prop:generalsourcetermsfirstboundnearlightcone}, which will give the main source term estimates at angular momentum $|n|\geq 2$, we need the following lemma, which handles the source terms other than the critical null-form, with all factors of angular momentum $|n|\geq 2$. Then we will need another analogous lemma for general inputs:
\begin{lemma}\label{lem:sourcetermsotherthannullformprop738prep} Let $F(\phi_1,\phi_2,\ldots, \phi_k)$ be any one\footnote{We expand the term $a\big(\Pi_{\Phi^{\perp}}(\varphi)\big)$ into a Taylor series which results in monomials of arbitrarily high degree in $\varphi$.} of the source terms in \eqref{coe eq tau R}, which are not of null-form type. This means we exclude the terms which involve a factor  
\begin{align*}
U_r^2 - U_t^2, \varphi_{j,t}^2 -\varphi_{j,r}^2, \big(a\big(\Pi_{\Phi^{\perp}}\varphi\big)\big)_{t}U_t - \big(a\big(\Pi_{\Phi^{\perp}}\varphi\big)\big)_{r}U_r, U_t\varphi_{j,t} - U_r\varphi_{j,r}, 
\end{align*}
in \eqref{non linear 1 2}, \eqref{coe eq precise 1}, \eqref{coe eq precise 2}, and further the terms involving a factor $(\partial_R\epsilon)^2 - \big(\partial_\tau\epsilon + \frac{\lambda_{\tau}}{\lambda}R\partial_R\epsilon\big)^2$ in \eqref{nonlinearity tau R}. Then assuming that all factors $\phi_j$ are angular momentum $n_j$ functions, $|n_j|\geq 2$, whose distorted Fourier transform admits a representation as sum $\overline{y}_j + \overline{z}_j$ with 
\begin{align*}
\overline{z}_j\in S_0^{\hbar_j},\,\hbar_j = |n_j|^{-1},
\end{align*}
and $\overline{y}_j $ admissibly singular (at angular momentum $n_j$), and assuming $|m|\geq 2$, $\big|\sum_j n_j - m\big|\lesssim 1$, $\hbar = |m|^{-1}$, we have 
\begin{align*}
\langle \phi_m(R;\xi),\,\chi_{1\lesssim R\lesssim\nu\tau}F(\phi_1,\phi_2,\ldots, \phi_k)\rangle_{L^2_{R\,dR}} = \overline{y} + \overline{z},
\end{align*}
where we have the bound 
\begin{align*}
\big\| \overline{y}\big\|_{sourceadm} + \tau^4\cdot\big\|\overline{z}\big\|_{S_1^{\hbar}}\lesssim D^k\cdot \max\{|n_j|\}^{-C}\cdot \prod_{j=1}^k |n_j|^C\Big(\big\|\overline{y}_j\big\|_{adm} + \tau^3\big(\big\|\overline{z}_j\big\|_{S_0^{\hbar_j}} + \big\|\mathcal{D}_{\tau}\overline{z}_j\big\|_{S_1^{\hbar_j}}\big)\Big)
\end{align*}
for a suitable constant $D = D(\nu)$. 
\end{lemma}
\begin{proof} Recalling \eqref{non linear 1 2}, which contributes to \eqref{coe eq precise 1}, \eqref{coe eq precise 2},, and omitting the terms involving a $\mathcal{Q}_0$-null form, we need to bound the schematically written terms (after multiplication by $\lambda^{-2}$ as in\eqref{nonlinearity tau R})
\begin{align*}
&a^{\kappa}\big(\Pi_{\Phi^{\perp}}(\varphi)\big)\cdot \frac{\sin^2U}{R^2}\cdot \varphi,\,\kappa\in \{1,2\},\,a\big(\Pi_{\Phi^{\perp}}(\varphi)\big)\cdot (\partial_{\tau} + \frac{\lambda_{\tau}}{\lambda}R\partial_R\big)a\big(\Pi_{\Phi^{\perp}}(\varphi)\big)\cdot \varphi,\\
&a\big(\Pi_{\Phi^{\perp}}(\varphi)\big)\cdot \partial_R a\big(\Pi_{\Phi^{\perp}}(\varphi)\big)\cdot \varphi,\,R^{-2}\big(\partial_\theta a\big(\Pi_{\Phi^{\perp}}(\varphi)\big)\big)^2\cdot \varphi,\,R^{-2}\varphi_{\theta}^2\cdot \varphi,\,R^{-2}\cos^{2\kappa} U\varphi^3,\kappa\in\{0,1\},\\
&\big(1 + a\big(\Pi_{\Phi^{\perp}}(\varphi)\big)\big)\cdot \frac{\sin U\cos U}{R^2}\varphi^2,\,\frac{\cos U}{R^2}\cdot \varphi_{\theta}\cdot\varphi^2,\,\frac{\sin U}{R^2}\cdot\varphi_{\theta}\cdot\varphi,\,\frac{\sin U}{R^2}\cdot \big(a\big(\Pi_{\Phi^{\perp}}(\varphi)\big)\big)_{\theta}\cdot \varphi^2,
\end{align*}
which we shall label $E_j,\,j = 1,\ldots, 10$. To this we add the term
\begin{align*}
E_{11}: = \frac{\sin U}{R^2\sqrt{1 - |\Pi_{\Phi^{\perp}}\varphi|^2}}\cdot \varphi_{\theta}\cdot\varphi
\end{align*}
from \eqref{coe eq precise 2}, as well as the terms 
\begin{align*}
\frac{\sin(2Q + \epsilon)\cdot\sin\epsilon}{R^2}\varphi,\,\frac{\sin(Q+\frac{\epsilon}{2})\cdot \sin(\frac{\epsilon}{2})}{R^2}\cdot\varphi_{\theta},\,\frac{\partial_R\epsilon}{1+R^2}\cdot\varphi,\,\frac{\lambda_{\tau}}{\lambda}\cdot \frac{R}{1+R^2}\cdot \big(\partial_\tau + \frac{\lambda_{\tau}}{\lambda}R\partial_R\big)\epsilon\cdot \varphi,
\end{align*}
which we label $E_j,\,j = 12,\ldots, 15$. 
\\
We shall now verify the assertion of the proposition for all of these terms. 
\\

{\it{$E_1$.}} Since the functions $\varphi$ are of small amplitude, we can expand $a\big(\Pi_{\Phi^{\perp}}(\varphi)\big)$ in terms of powers of $\varphi$. It follows that this term is schematically of the form
\[
 \frac{\sin^2U}{R^2}\cdot \varphi^k,\,k\geq 3,
 \]
 where we shall replace $\varphi^k$ by $\prod_{j=1}^k \phi_j$. Then using Lemma~\ref{lem:termswithU1} in conjunction with Lemma~\ref{lem:sourceadmtimesgood} (used iteratively) as well as Corollary~\ref{cor:multilingen2}, we infer that under our current assumptions on the $\phi_j$, and assuming that 
 \[
 |m|\geq 2,\,\big|m - \sum_{j=1}^k n_j\big|\lesssim 1
 \]
 we can write 
 \begin{align*}
\langle \phi_m(R;\xi),\,\frac{\sin^2U}{R^2}\cdot \prod_{j=1}^k \phi_j\rangle_{L^2_{R\,dR}} = \overline{y} + \overline{z}, 
 \end{align*}
 where we have the bound 
 \begin{align*}
 \big\|\overline{y}\big\|_{sourceadm} + \sup_{\tau\geq \tau_0}\tau^4\cdot \big\|\overline{z}\big\|_{S_1^{\hbar}}\leq D^k\cdot \max\{|n_j|\}^{-C}\cdot \prod_{j=1}^k |n_j|^C\big(\big\|\overline{y}_j\big\|_{adm} + \tau^3\big\|\overline{z}_j\big\|_{S_0^{\hbar_j}}\big)
 \end{align*}
 As we can assume the smallness of the quantities $|n_j|^C\big(\big\|\overline{y}_j\big\|_{adm} + \tau^3\big\|\overline{z}_j\big\|_{S_0^{\hbar_j}}\big)$, the right hand side will be $<\delta^k$ for some $0<\delta\ll 1$, and so the sum over $k$ will converge, given that $a(\phi)$ is analytic for $|\phi|\ll 1$.
\\

{\it{$E_2$.}} Again expanding $a\big(\Pi_{\Phi^{\perp}}(\varphi)$ into a power series in terms of $\varphi$, we reduce here to bounding terms of the form 
\begin{align*}
(\partial_{\tau} + \frac{\lambda_{\tau}}{\lambda}R\partial_R\big)\varphi\cdot \varphi^k,\,k\geq 4. 
\end{align*}
Replacing this by 
\begin{align*}
(\partial_{\tau} + \frac{\lambda_{\tau}}{\lambda}R\partial_R\big)\phi_1\cdot \prod_{j=2}^k\phi_j,\,k\geq 5,\,
\end{align*}
a conclusion analogous to the one for $E_1$ follows by taking advantage of Lemma~\ref{lem:singFouriertiphysicalngeq2admDeriv}, Lemma~\ref{lem:DerivonPrincSing}, together with Lemma~\ref{lem:sourceadmtimesgood}.  
\\
The term $E_3$ is similar. 
\\

For the terms $E_4, E_5$, we write them schematically as 
\begin{align*}
R^{-2}\big(\partial_\theta a\big(\Pi_{\Phi^{\perp}}(\varphi)\big)\big)^2\cdot \varphi = R^{-2}(\partial_{\theta}\varphi)^2\cdot \varphi\cdot b(\varphi)
\end{align*}
for a smooth function $b(\cdot)$ which is analytic near $\varphi = 0$. Hence as before we reduce to estimating 
\[
 R^{-2}\partial_{\theta}\phi_1\cdot\partial_{\theta}\phi_2\cdot \prod_{j=3}^k\phi_j,\,k\geq 3. 
 \]
Using Lemma~\ref{lem:basicproductsource},  the desired estimate follows from Lemma~\ref{lem:sourceadmtimesgood} used iteratively. 
\\

The remaining terms $E_j,\,j = 6,\ldots, 11$ can be handled analogously. 
\\

For the terms $E_{12} - E_{15}$, we notice that the factors 
\begin{align*}
\frac{\sin(2Q + \epsilon)\cdot\sin\epsilon}{R^2},\,\frac{\partial_R\epsilon}{1+R^2},\,\frac{\lambda_{\tau}}{\lambda}\cdot \frac{R}{1+R^2}\cdot \big(\partial_\tau + \frac{\lambda_{\tau}}{\lambda}R\partial_R\big)\epsilon
\end{align*}
are restricted to angular frequency $n = 0$. Thus the angular frequency of the expressions $E_j,\,j\in\{12,\ldots, 15\}$ is the same as the one for the factor $\varphi$. 
\\

For $E_{12}$, assuming that $\varphi = \phi_1$ where we can write 
\begin{align*}
\langle \phi_{n_1}(R;\xi),\,\phi_1\rangle_{L^2_{R\,dR}} = \overline{y} + \overline{z}, \big\|\overline{y}\big\|_{adm} + \sup_{\tau\geq\tau_0}\tau^3\cdot\big\|\overline{z}\big\|_{S_0^{\hbar}}<\infty, 
\end{align*}
we have the estimate (with $\hbar = n_1^{-1}$)
\begin{align*}
\sup_{\tau\geq\tau_0}\tau^4\big\|\chi_{R\ll\nu\tau}\frac{\sin(2Q + \epsilon)\cdot\sin\epsilon}{R^2}\cdot \phi_1\big\|_{\tilde{S}_1^{\hbar}}\lesssim \big(\sup_{\tau\geq\tau_0}\tau^3\big\|\phi_1\big\|_{\tilde{S}_0^{\hbar_1}}\big). 
\end{align*}
In fact, using Theorem \ref{thm:KSTGao precise}, we have 
\begin{align*}
\left\|\chi_{R\ll\nu\tau}\sin(2Q + \epsilon)\cdot\sin\epsilon\right\|_{\tilde{S}_0^{\frac13}}\lesssim \tau^{-2}. 
\end{align*}
The desired bound then follows from Lemma~\ref{lem:basicproductsource} and Remark~\ref{rem::basicproductsource}. \\
Near the light cone $R = \nu\tau$, we use that again from  Theorem \ref{thm:KSTGao precise}, we have 
\begin{align*}
\frac{\sin(2Q + \epsilon)\cdot\sin\epsilon}{R^2} = f_1 + f_2
\end{align*}
where 
\begin{align*}
\left\|f_1\right\|_{\tilde{S}_0^{\frac13}}\lesssim \tau^{-3}, 
\end{align*}
while the function $f_2$ admits an expansion $f_2 = \sum_{l=1}^8f_{2l}$ where the functions $f_{2l}$ in turn can be written like in Lemma~\ref{lem:singFouriertiphysicalngeq2adm}  but with $\chi_{|\nu\tau - R|\lesssim \hbar}$ replaced by $\chi_{|\nu\tau - R|\lesssim 1}$ and with the decay factor 
\[
\tau^{-1-\nu}
\]
in the bound for the coefficients $G_{k,l,j}(\tau, \nu\tau - R)$ replaced by the better factor $\tau^{-3-\nu}$. The desired conclusion for the product 
\[
\chi_{R\sim\nu\tau}\frac{\sin(2Q + \epsilon)\cdot\sin\epsilon}{R^2}\cdot \phi_1
\]
is  obtained in analogy to the proof of Lemma~\ref{lem:basicproductsource}. The other terms $E_j,\,j = 13, 14, 15$ are handled similarly. 
\end{proof}

The preceding lemma deals with source terms localized away from $R = 0$. Near the spatial origin all functions are smooth, but we require the proper vanishing rates near the origin: 

\begin{lemma}\label{lem:sourcetermsotherthannullformprop738prepnearorigin} Let $F(\phi_1,\phi_2,\ldots, \phi_k)$ be as in the preceding lemma. Then under the same hypotheses on the functions $\phi_j$, the parameters $n_j, m$, there is a polynomial 
\[
P_3(R) = \sum_{l=0}^3\gamma_lR^l,\,\gamma_l = \gamma_l(\tau), 
\]
such that we have the following bounds:
\begin{align*}
&\tau^4\cdot\big\|\chi_{R\lesssim 1}\cdot F(\phi_1,\phi_2,\ldots, \phi_k)\big\|_{\tilde{S}_1^{\hbar}}\lesssim D^k\cdot \max\{|n_j|\}^{-C}\cdot \prod_{j=1}^k |n_j|^C\big(\big\|\overline{y}_j\big\|_{adm} + \tau^3\big\|\overline{z}_j\big\|_{S_0^{\hbar_j}}\big),\\
& \tau^4\cdot \sum_{l=0}^3|\gamma_l|\lesssim D^k\cdot \max\{|n_j|\}^{-C}\cdot \prod_{j=1}^k |n_j|^C\big(\big\|\overline{y}_j\big\|_{adm} + \tau^3\big\|\overline{z}_j\big\|_{S_0^{\hbar_j}}\big).
\end{align*}
 \end{lemma}

The proof of this follows from the estimates in subsection~\ref{subsec:allsourcetermsngeq2} for the restrictions of the source terms away from the light cone, in conjunction with Lemma~\ref{lem:admsingawayfromshock}.

\subsection{Null-form estimates near the light cone; only angular momenta $|n|\geq 2$ involved.}\label{subset:hardnullformestimates}

We shall now deal with the most delicate type of source term due to the derivatives it contains, which we cast here as a trilinear expression 
\[
\mathcal{N}_0(\phi_1,\phi_2,\phi_3): = \phi_1\cdot\left[\left(\partial_{\tau}+ \frac{\lambda_{\tau}}{\lambda}R\partial_R\right)\phi_2\cdot \left(\partial_{\tau}+ \frac{\lambda_{\tau}}{\lambda}R\partial_R\right)\phi_3 - \partial_R\phi_2\cdot\partial_R\phi_3\right]
\]
As hinted at in the title of this subsection, we shall assume for now that all factors as well as the expression itself will be at angular momentum $|n|\geq 2$, i.e., we shall use Fourier representations with respect to such angular momenta. This will simplify the presentation a bit, and the extension to general factors will be rather routine. 
Our way to deal with this expression will be to take advantage of Fourier localizations, as set up in the preceding subsection. We shall face one particular technical difficulty, which is intimately tied to our functional setup, and which will require a bit of detour to handle. To understand it, let us write the above term in terms of the `good' derivative $\partial_{\tau}+ \frac{\lambda_{\tau}}{\lambda}R\partial_R - \partial_R$ as well as the `bad' derivative $\partial_{\tau}+ \frac{\lambda_{\tau}}{\lambda}R\partial_R+\partial_R$: 
\begin{equation}\label{eq:TheNullForm}\begin{split}
\mathcal{N}_0(\phi_1,\phi_2,\phi_3) &= \frac12\phi_1\cdot \left(\partial_{\tau}+ \frac{\lambda_{\tau}}{\lambda}R\partial_R - \partial_R\right)\phi_2\cdot  \left(\partial_{\tau}+ \frac{\lambda_{\tau}}{\lambda}R\partial_R + \partial_R\right)\phi_3\\
& +\frac12 \phi_1\cdot \left(\partial_{\tau}+ \frac{\lambda_{\tau}}{\lambda}R\partial_R + \partial_R\right)\phi_2\cdot  \left(\partial_{\tau}+ \frac{\lambda_{\tau}}{\lambda}R\partial_R - \partial_R\right)\phi_3
\end{split}\end{equation}
Then the strategy is to show that if the factors have Fourier transform which is either admissibly singular or in $\Sh_{0}$ (with $\hbar$ in accordance with their angular momentum), then the source term, upon application of the wave parametrix, leads again to such terms. This works for most of the situations, except in the case where a `bad' derivative hits a singular term (i.e., with Fourier transform admissibly singular), while the `good' derivative hits a term in $\Sh_{0}$. The problem is that these terms are still too rough to lead to a term in $\Sh_{1}$, but they also lack the precise fine structure that characterises our admissibly singular terms(or more precisely, their distorted Fourier transform). 
\\
The trick to deal with this problem is to use a type of {\it{'normal form transform'.}} Specifically, we `prepare' the wave equation we are trying to solve a bit, by modifying unknown variable (for us $\veps_{\pm}(n)$) by subtracting a suitable correction term (whose distorted Fourier transform happens to be a linear combination of an admissibly singular and a $S_0^{(\hbar)}$-term) from it which modifies the equation in such a way that the `troublesome' terms disappear, up to other troublesome terms of similar structure but which are smaller. Then the same method can be re-iterated to the remaining troublesome terms, until they eventually disappear. What makes this method work is the very special structure of the $\mathcal{N}_0$ null-structure. To begin with, we `micro-localize' the expression as follows: 
\begin{equation}\label{eq:keytrilinnullformdecomp}\begin{split}
\mathcal{N}_0(\phi_1,\phi_2,\phi_3) &= \sum_{\lambda_2,\lambda_3}\mathcal{N}_0\left(P_{<\min\{\lambda_2,\lambda_3\}}\phi_1,P_{\lambda_2}\phi_{2}, P_{\lambda_3}\phi_{3}\right)\\
& +  \sum_{\lambda_2,\lambda_3}\mathcal{N}_0\left(P_{\geq \min\{\lambda_2,\lambda_3\}}\phi_1,P_{\lambda_2}\phi_{2}, P_{\lambda_3}\phi_{3}\right)\\
& =: \mathcal{N}_{01}(\phi_1,\phi_2,\phi_3) + \mathcal{N}_{02}(\phi_1,\phi_2,\phi_3).
\end{split}\end{equation}
Of these it turns out that the second term $\mathcal{N}_{02}(\phi_1,\phi_2,\phi_3)$ is of the `good kind', as exemplified by the following 
\begin{proposition}\label{prop:goodcubicnullform1} Let $|n_j|\geq 2,\,j = 1,2,3$, and assume that $\{m_2, m_1\}$ satisfy $|m_j|\geq 2$ and the triples $\{m_2,m_1,n_1\}$, $\{m_1,n_2,n_3\}$ satisfy the conditions in of  Prop.~\ref{prop:bilin2} (instead of $\{n_3,n_1,n_2\}$ there). Assume that each $\phi_j$ admits a distorted Fourier transform $\overline{x}_j$ at angular momentum $n_j$ which can be written in the form $\overline{x}_j = \overline{y}_j + \overline{z}_j$ with $\overline{y}_j$ admissibly singular while $\overline{z}_j\in S_0^{\hbar_j}$, $\hbar_j = |n_j|^{-1}$. 
Then the function 
	\[
	\chi_{R\gtrsim \tau}\mathcal{N}_{02}(\phi_1,\phi_2,\phi_3)|_{R<\nu\tau} = g(\tau, R)|_{R<\nu\tau},
	\]
where $g(\tau, R)$ admits a angular momentum $m_2$ distorted Fourier representation 
	\[
	g(\tau, R) = \int_0^\infty \phi_{m_{2}}(R,\xi)\yb(\tau, \xi)\rho_{m_{2}}(\xi)\,d\xi,
	\]
	where $\yb =\tilde{ \yb}_1 +\tilde{\yb}_2$ and $\tilde{\yb}_1$ is an admissibly singular source term at angular momentum $m_2$, while $\tilde{\yb}_2\in \Sht_{1}$, $\hbar_2 = \frac{1}{m_2+1}$, enjoying the bound
	\begin{align*}
	\big\|\tilde{\yb}_1\big\|_{sourceadm} + \tau^4\cdot \big\|\tilde{\yb}_2\big\|_{ \Sht_{1}}\lesssim \max\{|n_j|\}^{-C}\cdot \prod_{j=1}^3|n_j|^C\Big(\big\| \overline{y}_j\big\|_{adm} + \sup_{\tau\geq\tau_0}\tau^3\cdot \big(\big\| \overline{z}_j\big\|_{S_0^{\hbar_j}} + \big\| \mathcal{D}_{\tau}\overline{z}_j\big\|_{S_1^{\hbar_j}}\big)\Big)
	\end{align*}
\end{proposition}
\begin{proof}
	For simplicity of notation, we shall indicate frequency localizations by subscripts, thus $P_{\lambda}\phi = \phi_{\lambda}$ etc. Expand the term out as 
	\begin{align*}
	\mathcal{N}_{02}(\phi_1,\phi_2,\phi_3) &= \sum_{\lambda_3<\lambda_2} \phi_{1,\geq\lambda_3}\cdot \left(\partial_{\tau}+ \frac{\lambda_{\tau}}{\lambda}R\partial_R - \partial_R\right)\phi_{2,\lambda_2}\cdot  \left(\partial_{\tau}+ \frac{\lambda_{\tau}}{\lambda}R\partial_R + \partial_R\right)\phi_{3,\lambda_3}\\
	& + \sum_{\lambda_3<\lambda_2} \phi_{1,\geq\lambda_3}\cdot \left(\partial_{\tau}+ \frac{\lambda_{\tau}}{\lambda}R\partial_R + \partial_R\right)\phi_{2,\lambda_2}\cdot  \left(\partial_{\tau}+ \frac{\lambda_{\tau}}{\lambda}R\partial_R - \partial_R\right)\phi_{3,\lambda_3}\\
	& + \sum_{\lambda_3\geq \lambda_2}\phi_{1,\geq\lambda_2}\cdot \left(\partial_{\tau}+ \frac{\lambda_{\tau}}{\lambda}R\partial_R - \partial_R\right)\phi_{2,\lambda_2}\cdot  \left(\partial_{\tau}+ \frac{\lambda_{\tau}}{\lambda}R\partial_R + \partial_R\right)\phi_{3,\lambda_3}\\
	& +  \sum_{\lambda_3\geq \lambda_2}\phi_{1,\geq\lambda_2}\cdot \left(\partial_{\tau}+ \frac{\lambda_{\tau}}{\lambda}R\partial_R + \partial_R\right)\phi_{2,\lambda_2}\cdot  \left(\partial_{\tau}+ \frac{\lambda_{\tau}}{\lambda}R\partial_R - \partial_R\right)\phi_{3,\lambda_3}.
	\end{align*}
	By symmetry it suffices to treat the first two terms on the right. Call them $\mathcal{N}_{02a}, \mathcal{N}_{02b}$ respectively. 
	\\
	
	{\it{The contribution of $\mathcal{N}_{02a}$.}} Observe that the `bad' derivative falls on the low frequency term $\phi_{3,\lambda_3}$. We shall derive the desired assertion by means of twofold application of suitable bilinear estimates. To begin with, we need a modification of the pure $\Sh_{j}$-based bilinear estimate Prop.~\ref{prop:bilin4}:
	\begin{lemma}\label{lem:propbilin4para1} Assume that $\phi_j,\,j = 1,2$ are angular momentum $n_j$, $|n_j|\geq 2$ functions, with $\phi_1\in \tilde{S}_0^{(\hbar_1)},\,\phi_2\in \tilde{S}_1^{(\hbar_2)}$. Let $\{m,n_1,n_2\}$, $|m|\geq 2$, be an admissible triple of angular momenta, i. e. $\big|m - \sum_j n_j\big| = O(1)$. Then 
		\[
		\chi_{1\lesssim R<\nu\tau}\sum_{\lambda>1}P_{\geq\lambda}\phi_1\cdot P_{\lambda}\phi_2\in \left\langle\xi\right\rangle^{-\frac12}S_1^{\hbar},\quad \hbar = \frac{1}{m+1}
		\]
		where the sum is over dyadic frequencies, and we have the bound 
		\begin{align*}
		\left\|\left\langle\phi_{m}(R,\xi),\,\chi_{R<\nu\tau}\sum_{\lambda>1}P_{\geq\lambda}\phi_1\cdot P_{\lambda}\phi_2\right\rangle_{L^2_{R\,dR}}\right\|_{\left\langle\xi\right\rangle^{-\frac12}S_1^{\hbar}}\lesssim \tau\cdot \max\{\hbar_1,\hbar_2\}^{-2}\cdot \left\|\phi_1\right\|_{\tilde{S}_0^{(\hbar_1)}}\cdot\left\|\phi_2\right\|_{\tilde{S}_0^{(\hbar_2)}}.
		\end{align*}
	\end{lemma}
	\begin{proof}.  Let $\hbar_j = \frac{1}{|n_j|}$, $\hbar = \frac{1}{|m|+1}$. We need to show that 
\begin{align*}
\langle \phi_m(R;\xi),\nabla_R\big(\chi_{R\lesssim \nu\tau}\sum_{\lambda>1}P_{\geq\lambda}\phi_1\cdot P_{\lambda}\phi_2\big)\rangle_{L^2_{R\,dR}}\in S_1^{\hbar}
\end{align*}
According to the Leibniz rule it suffices to show 
\begin{align*}
&\langle \phi_m(R;\xi),\nabla_R\big(\chi_{R\lesssim \nu\tau}\big)\sum_{\lambda>1}P_{\geq\lambda}\phi_1\cdot P_{\lambda}\phi_2\rangle_{L^2_{R\,dR}}\in S_1^{\hbar}\\
&\langle \phi_m(R;\xi), \chi_{R\lesssim \nu\tau}\sum_{\lambda>1}\nabla_R\big(P_{\geq\lambda}\phi_1\big)\cdot P_{\lambda}\phi_2\rangle_{L^2_{R\,dR}}\in S_1^{\hbar}\\
&\langle \phi_m(R;\xi), \chi_{R\lesssim \nu\tau}\sum_{\lambda>1}P_{\geq\lambda}\phi_1\cdot \nabla_R\big(P_{\lambda}\phi_2\big)\rangle_{L^2_{R\,dR}}\in S_1^{\hbar}\\
\end{align*}
For the first two of these expressions this is essentially a direct consequence of Prop.~\ref{prop:bilin1}. For the last term, we can trade a derivative on the low frequency term $P_{\lambda}\phi_2$ for a derivative on the high frequency term $P_{\geq\lambda}\phi_1$, after localizing the frequencies of both terms to dyadic size (in the sense of angular momentum $n_j$-functions). This results in a gain $\frac{\lambda}{\mu}$ where $\mu$ is the dyadic frequency of $P_{\geq\lambda}\phi_2$. One can then argue the same way as for Prop.~\ref{prop:bilin1}, and summation over all dyadic frequencies $\lambda,\mu$ reduces to square summation over these frequencies by exploiting the gain $\frac{\lambda}{\mu}$ as well as the Cauchy-Schwarz inequality. 
\end{proof}

	We now state the first main ingredient in the proof of the proposition:
	\begin{lemma}\label{lem:N02abilin1} Let $\phi_1, \phi_3$ be angular  momentum $n_j,\,j = 1,3$ functions (with $|n_j|\geq 2$) and with the property that their angular momentum $n_j$ distorted Fourier transforms $\yb_j(\tau, \xi)$ can be split into 
		\[
		\yb_j(\tau, \xi) = \tilde{\yb}_j(\tau, \xi) + \tilde{\tilde{\yb}}_j(\tau, \xi),
		\]
		where $ \tilde{\yb}_j(\tau, \xi)\in \Shj_0,\,\mathcal{D}_{\tau}\tilde{\yb}_j(\tau, \xi)\in \Shj_{1}$, and further $\tilde{\tilde{\yb}}_j$ is admissibly singular. Then we can write 
		\[
		\chi_{R\gtrsim \tau}\phi_{1,\geq\lambda_3}\cdot \left(\partial_{\tau}+ \frac{\lambda_{\tau}}{\lambda}R\partial_R + \partial_R\right)\phi_{3,\lambda_3}|_{R<\nu\tau} = g(\tau, R)|_{R<\nu\tau}
		\]
		where, if $\{m,n_1,n_3\}$ with $|m|\geq 2$ is an admissible angular momentum triple in the sense of Prop.~\ref{prop:bilin2}, we have an angular momentum $m$ Fourier representation of $g(\tau, R)$ 
		\[
		g(\tau, R) = \int_0^\infty \phi_{m}(R,\xi)\zb(\tau, \xi)\rho_{m}(\xi)\,d\xi,\quad \hbar = \frac{1}{m+1}
		\]
		where $\zb$ is the sum of a function in $\zb_1\in \left\langle\xi\right\rangle^{-\frac14-}S_1^{\hbar}$ and a function $\zb_2$ which is source admissibly singular at angular momentum $m$ and level $l_1\geq 1$. Furthermore, we have the estimate 
		\begin{align*}
		\tau^3\big\|\zb_1(\tau,\cdot)\big\|_{ \left\langle\xi\right\rangle^{-\frac14-}S_1^{\hbar}} + \big\|\zb_2\big\|_{sourceadm}\lesssim  \max_{j=1,3}\{|n_j|\}^{-C}\cdot \prod_{j=1,3}|n_j|^C\Big(\big\| \tilde{\tilde{\yb}}_j\big\|_{adm} + \sup_{\tau\geq\tau_0}\tau^3\cdot \big(\big\| \tilde{\yb}_j(\tau,\cdot)\big\|_{S_0^{\hbar_j}} + \big\| \mathcal{D}_{\tau}\tilde{\yb}_j(\tau,\cdot)\big\|_{S_1^{\hbar_j}}\big)\Big)
               \end{align*}
	\end{lemma}
\begin{proof}
	(lemma) We distinguish between different situations, depending on the nature of the factors. 
	\\
	
	{\it{(I): both factors of type $ \tilde{\yb}_j$.}} Here we write
	\begin{align*}
	\left(\partial_{\tau}+ \frac{\lambda_{\tau}}{\lambda}R\partial_R + \partial_R\right)f(\tau, R)=\left(\partial_{\tau}+ \frac{\lambda_{\tau}}{\lambda}R\partial_R \right)f+\partial_{R}f.
	\end{align*}
	Then the desired conclusion follows from Prop.~\ref{prop:bilin4}, where the complications concerning the spatial origin $R = 0$ are avoided here since we localize things away from the origin. The estimate stated at the end of the lemma is a direct consequence of the bound \eqref{eq:somewhatpoorbound} resulting in a better than needed decay rate $\tau^{-4-}$ for the product. 
	\\
	
	{\it{(II): first factor of type $\tilde{\tilde{\yb}}_1$, second factor of type $\tilde{\yb}_3$.}} We may reduce the first factor to $\chi_{\nu\tau-R<\min\{\hbar_1,\hbar_3\}}\phi_{1,\geq\lambda_3}$, $\hbar_1 = \frac{1}{|n_1|+1}$, since else the assertion follows from the preceding case, in light of Lemma~\ref{lem:singFouriertiphysicalngeq2adm}. Also, again using the preceding case and Lemma~\ref{lem:singFouriertiphysicalngeq2adm}, we may assume that $\phi_{1,\geq\lambda_3}$ is like the middle term $f_2$ in that lemma. We write $\partial_{\tau}+ \frac{\lambda_{\tau}}{\lambda}R\partial_R + \partial_R = D_\tau^+$, and decompose 
	\[
	D_\tau^+\phi_{3,\lambda_3} =\left[D_\tau^+\phi_{3,\lambda_3} - P_3\left(D_\tau^+\phi_{3,\lambda_3}\right)\right] + P_3\left(D_\tau^+\phi_{3,\lambda_3}\right),
	\]
	where $P_3\left(D_\tau^+\phi_{3,\lambda_3}\right)$ denotes the third order Taylor polynomial of $D_\tau^+\phi_{3,\lambda_3}$ centered at $R = \nu\tau$. Then we claim that 
	\[
	\chi_{\nu\tau-R<\min\{\hbar_1,\hbar_3\}}\phi_{1,\geq\lambda_3}\cdot P_3\left(D_\tau^+\phi_{3,\lambda_3}\right)|_{R<\nu\tau}
	\]
	coincides with the restriction there of an admissibly singular source term of level $l_1\geq 2$ up to a function in $\tilde{S}_1^{(\hbar)}$, $\hbar = \frac{1}{|m|+1}$. To see this, assume, say, that $\hbar_3\ll\hbar_1$. Then write 
	\[
	P_3\left(D_\tau^+\phi_{3,\lambda_3}\right) = \sum_{j=0}^3 \frac{1}{j!}\left(D_\tau^+\phi_{3,\lambda_3}\right)^{(j)}(\tau, \nu\tau)(\nu\tau - R)^j,
	\]
	where we have the coefficient bounds $\left|\left(D_\tau^+\phi_{3,\lambda_3}\right)^{(j)}(\tau, \nu\tau)\right|\lesssim \tau\hbar_3^{-\frac32-j}\cdot \left\|D_\tau^+\phi_{3,\lambda_3}\right\|_{\tilde{S}_1^{(\hbar_3)}}$.
	In fact we schematically\footnote{We omit the now routine details required to handle the turning point} estimate for $R\simeq \nu\tau$
	\begin{align*}
		\left|D_\tau^+\phi_{3,\lambda_3}(\tau, R)\right| &= \left|\int_0^\infty \mathcal{F}\left(D_\tau^+\phi_{3,\lambda_3}\right)(\xi)\cdot\phi_{n_{3}}(R,\xi)\cdot\rho_{n_{3}}(\xi)\,d\xi\right|\\
		&\lesssim \int_0^\infty \left|\mathcal{F}\left(D_\tau^+\phi_{3,\lambda_3}\right)(\xi)\right|\cdot\frac{1}{\xi^{\frac14}\tau^{\frac12}}\cdot\min\{\tau\xi^{\frac12}\hbar_3,1\}\,d\xi\\
		&\lesssim \hbar_3^{\frac12}\cdot \hbar_3^{-2}\left\|\mathcal{F}\left(D_\tau^+\phi_{3,\lambda_3}\right)\right\|_{\Shth_1},
	\end{align*}
Indeed, the second line on the RHS above can be written as (for $\hbar_{3}\xi^{\frac12}\ll1$)
	\begin{align*}
		&\hbar_{3}^{\frac12}\int_0^\infty \left|\mathcal{F}\left(D_\tau^+\phi_{3,\lambda_3}\right)(\xi)\right|\cdot\frac{1}{\hbar_{3}^{\frac12}\xi^{\frac14}\tau^{\frac12}}\cdot\min\{\tau\xi^{\frac12}\hbar_3,1\}\,d\xi\\
		\lesssim&\hbar_{3}^{\frac12}\cdot\hbar^{-2+\delta}_{3}\int_0^\infty\xi^{-1+\frac{\delta}{2}}\cdot\xi^{\frac12}\cdot \left(\hbar_{3}^{2}\xi\right)^{1-\frac{\delta}{2}}\cdot\xi^{-\frac12}\left|\mathcal{F}\left(D_\tau^+\phi_{3,\lambda_3}\right)(\xi)\right|\,d\xi.
	\end{align*}	
The case when $\hbar_{3}\xi^{\frac12}\gtrsim 1$ is handled similarly.
	Further, expand 
	\begin{align*}
	\chi_{\nu\tau-R<\min\{\hbar_1,\hbar_3\}}\phi_{1,\geq\lambda_3} = \chi_{|\nu\tau-R|\lesssim \hbar_3}\sum_{l=1}^8 \sum_{k=1}^N\sum_{i=0}^{N_1}\frac{G_{k,l,i}(\tau,\nu\tau - R)}{\tau^{\frac12}}\hbar_1^{-\frac{l+1}{2}}\left[\nu\tau - R\right]^{\frac{l}{2}+k\nu}\left(\log(\nu\tau - R)\right)^i
	\end{align*}
	Then the desired conclusion follows by multiplying the preceding expressions, reformulating things in the form of $f_2$ in  Lemma~\ref{lem:singFouriertiphysicalngeq2adm} with respect to angular momentum $n_3$, and using Lemma~\ref{lem:singPhysicaltoFourierngeq2}. Observe that the resulting function is actually admissibly regular source of level $l_1\geq 2$ (since it is one derivative more regular than a minimal regularity source), 
	and this information lets us gain another power $\hbar_3$ to counteract the loss of $\hbar_3^{-\frac32}$ to bound the coefficients of the Taylor polynomial. The situation $\hbar_3\gtrsim \hbar_1$ is handled similarly. 
	\\
	
	Next, consider the product 
	\[
	\chi_{\nu\tau-R<\min\{\hbar_1,\hbar_3\}}\phi_{1,\geq\lambda_3}\cdot\left[D_\tau^+\phi_{3,\lambda_3} - P_3\left(D_\tau^+\phi_{3,\lambda_3}\right)\right]. 
	\]
	We claim that this is in $\tilde{S}_1^{(\hbar)}$. For this let us again assume that $\hbar_3\ll\hbar_1$, say, which implies the `output angular momentum' $\hbar\simeq \hbar_3$. Then split the expression into a low and a high-frequency part:
	\begin{align*}
	&\chi_{\nu\tau-R<\min\{\hbar_1,\hbar_3\}}\phi_{1,\geq\lambda_3}\cdot\left[D_\tau^+\phi_{3,\lambda_3} - P_3\left(D_\tau^+\phi_{3,\lambda_3}\right)\right]\\
	& = P_{<\lambda_3}\left(\chi_{\nu\tau-R<\min\{\hbar_1,\hbar_3\}}\phi_{1,\geq\lambda_3}\cdot\left[D_\tau^+\phi_{3,\lambda_3} - P_3\left(D_\tau^+\phi_{3,\lambda_3}\right)\right]\right)\\
	& + P_{\geq\lambda_3}\left(\chi_{\nu\tau-R<\min\{\hbar_1,\hbar_3\}}\phi_{1,\geq\lambda_3}\cdot\left[D_\tau^+\phi_{3,\lambda_3} - P_3\left(D_\tau^+\phi_{3,\lambda_3}\right)\right]\right)\\
	&=:A + B. 
	\end{align*}
	To bound $A$, use 
	\begin{align*}
	&\left\|P_{<\lambda_3}\left(\chi_{\nu\tau-R<\min\{\hbar_1,\hbar_3\}}\phi_{1,\geq\lambda_3}\cdot\left[D_\tau^+\phi_{3,\lambda_3} - P_3\left(D_\tau^+\phi_{3,\lambda_3}\right)\right]\right)\right\|_{\tilde{S}_1^{(\hbar)}}\\
	& \lesssim \hbar\cdot\left(\sum_{\mu<\lambda_3}\left(\hbar^2\mu\right)^{1-\delta}\cdot \left\langle\hbar^2\mu\right\rangle^{3+2\delta}\cdot\left\|\chi_{\nu\tau-R<\min\{\hbar_1,\hbar_3\}}\phi_{1,\geq\lambda_3}\cdot D_\tau^+\phi_{3,\lambda_3} \right\|_{L^2_{R\,dR}}^2\right)^{\frac12}\\
	& + \hbar\cdot\left(\sum_{\mu<\lambda_3}\left(\hbar^2\mu\right)^{1-\delta}\cdot \left\langle\hbar^2\mu\right\rangle^{3+2\delta}\cdot\left\|\chi_{\nu\tau-R<\min\{\hbar_1,\hbar_3\}}\phi_{1,\geq\lambda_3}\cdot P_3\left(D_\tau^+\phi_{3,\lambda_3}\right)\right\|_{L^2_{R\,dR}}^2\right)^{\frac12}
	\end{align*}
	where $\mu$ ranges over dyadic scales. To bound the first term on the right, use the bound
	\begin{align*}
	\left\|D_\tau^+\phi_{3,\lambda_3}\right\|_{L^2_{R\,dR}}\lesssim \hbar_3^{-1}\left(\hbar_3^2\lambda_3\right)^{-\frac12+\frac{\delta}{2}}\left\langle \hbar_3^2\lambda_3\right\rangle^{-\frac32-\delta}\cdot \left\|D_\tau^+\phi_{3,\lambda_3}\right\|_{\tilde{S}_1^{(\hbar_3)}}.
	\end{align*}
	Further using the bound 
	\[
	\left\|\chi_{\nu\tau-R<\min\{\hbar_1,\hbar_3\}}\phi_{1,\geq\lambda_3}\right\|_{L^\infty_{R\,dR}}\lesssim_{\hbar_1} \lambda_3^{-\frac14}\cdot \tau^{-\frac32-\nu}(\log\tau)^C\cdot \big\|\tilde{\tilde{\yb}}_1\big\|_{adm},
	\]
	and taking advantage of H\"older's inequality, we bound the first sum by 
	\begin{align*}
	&\hbar\cdot\left(\sum_{\mu<\lambda_3}\left(\hbar^2\mu\right)^{1-\delta}\cdot \left\langle\hbar^2\mu\right\rangle^{3+2\delta}\cdot\left\|\chi_{\nu\tau-R<\min\{\hbar_1,\hbar_3\}}\phi_{1,\geq\lambda_3}\cdot D_\tau^+\phi_{3,\lambda_3} \right\|_{L^2_{R\,dR}}^2\right)^{\frac12}\\
	&\lesssim_{\hbar_1}  \lambda_3^{-\frac14}\cdot \tau^{-\frac32-\nu}(\log\tau)^C\cdot  \left\|D_\tau^+\phi_{3,\lambda_3}\right\|_{\tilde{S}_1^{(\hbar_3)}}\cdot \big\|\tilde{\tilde{\yb}}_1\big\|_{adm},
	\end{align*}
	where we have also taken advantage of the fact that $\hbar\simeq\hbar_3$. Here the decay factor $\lambda_{3}^{-\frac14}$ comes from the power $(\nu\tau-R)^{\frac12}$ in the singular profile of $\chi_{\nu\tau-R<\min\{\hbar_1,\hbar_3\}}\phi_{1,\geq\lambda_3}$. In fact if $(\nu\tau-R)\lesssim \lambda_{3}^{-\frac12}$, this is straightforward. If $(\nu\tau-R)\lambda_{3}^{\frac12}\gtrsim 1$, then we use integration by parts to gain negative powers in $(\nu\tau-R)\lambda_{3}^{\frac12}$ and use the power $(\nu\tau-R)^{\frac12}$ to gain the factor $\lambda_{3}^{-\frac14}$.
	\\
The other square sum above involving the Taylor polynomial can be treated  similarly, by bounding the Taylor polynomial in $L^{2}_{RdR}$. 

	We use the fact that the function $\phi_{1,\geq\lambda_3}$ decays rapidly beyond scale $\lambda_3^{-\frac12}$ in $\nu\tau-R$, 
	as follows by straightforward integration by parts. In fact, upon changing the integration variable, say $\tilde{\xi}:=(\nu\tau-R)\xi^{\frac12}$, we gain powers in $(\nu\tau-R)^{-1}\lambda^{-\frac12}_{3}$. For each term $(\nu\tau-R)^{j}$ in the Taylor expansion, we use $(\nu\tau-R)^{-j}\lambda_{3}^{-\frac{j}{2}}$ to gain a decay of $\lambda_{3}^{-\frac{j}{2}}$, which will finally cancel the growth $\hbar_{3}^{-1}\leq \lambda_{3}^{\frac12}$. Thus we may for all intents and purposes restrict to the region $\nu\tau - R\lesssim \lambda_3^{-\frac12}$, and there we have the bound (using the fact $\hbar^{-1}_{3}\leq \lambda_{3}^{\frac12}$)
	\begin{align*}
	\left\|\chi_{\nu\tau-R<\lambda_3^{-\frac12}}P_3\left(D_\tau^+\phi_{3,\lambda_3}\right)\right\|_{L^2_{R\,dR}}&\lesssim \sum_{j=0}^3\left\|\chi_{\nu\tau-R<\lambda_3^{-\frac12}}\left(D_\tau^+\phi_{3,\lambda_3}\right)^{(j)}(\tau, \nu\tau)\cdot (\nu\tau -R)^j\right\|_{L^2_{R\,dR}}\\
	&\lesssim \lambda_3^{-\frac14}\tau^{\frac12}\cdot \sum_{j=0}^3\left|\left(D_\tau^+\phi_{3,\lambda_3}\right)^{(j)}(\tau, \nu\tau)\right|\cdot\lambda_3^{-\frac{j}{2}}\\
	&\lesssim \lambda_{3}^{-\frac14}\left(\sum_{j=0}^3\lambda_3^{\frac{j}{2}}\cdot\lambda_3^{-\frac{j}{2}}\right)\cdot  \hbar_3^{-1}\left(\hbar_3^2\lambda_3\right)^{-\frac12+\frac{\delta}{2}}\left\langle \hbar_3^2\lambda_3\right\rangle^{-\frac32-\delta}\cdot \left\|D_\tau^+\phi_{3,\lambda_3}\right\|_{\tilde{S}_1^{(\hbar_3)}}.
	\end{align*}
	\\
	Combining with the $L^\infty$-bound for $\chi_{\nu\tau-R<\min\{\hbar_1,\hbar_3\}}\phi_{1,\geq\lambda_3}$ from before, we infer the bound for the second square sum for $A$
	\begin{align*}
	& \hbar\cdot\left(\sum_{\mu<\lambda_3}\left(\hbar^2\mu\right)^{1-\delta}\cdot \left\langle\hbar^2\mu\right\rangle^{3+2\delta}\cdot\left\|\chi_{\nu\tau-R<\min\{\hbar_1,\hbar_3\}}\phi_{1,\geq\lambda_3}\cdot P_3\left(D_\tau^+\phi_{3,\lambda_3}\right)\right\|_{L^2_{R\,dR}}^2\right)^{\frac12}\\
	&\lesssim_{\hbar_1} \lambda_3^{-\frac14}\cdot \tau^{-\frac32-\nu}(\log\tau)^C\cdot  \left\|D_\tau^+\phi_{3,\lambda_3}\right\|_{\tilde{S}_1^{(\hbar_3)}}\cdot \big\|\tilde{\tilde{\yb}}_1\big\|_{adm}.
	\end{align*}
	This concludes the bound for $A$.
	\\
	
	To get the desired bound for $B$, split it into 
	\begin{equation}\label{eq:lemmaN0abB1B2}\begin{split}
	&P_{\geq\lambda_3}\left(\chi_{\nu\tau-R<\min\{\hbar_1,\hbar_3\}}\phi_{1,\geq\lambda_3}\cdot\left[D_\tau^+\phi_{3,\lambda_3} - P_3\left(D_\tau^+\phi_{3,\lambda_3}\right)\right]\right)\\
	& = \sum_{\mu\geq\lambda_3}P_{\mu}\left(\chi_{\nu\tau-R<\min\{\hbar_1,\hbar_3\}}\phi_{1,[\lambda_3,\mu]}\cdot\left[D_\tau^+\phi_{3,\lambda_3} - P_3\left(D_\tau^+\phi_{3,\lambda_3}\right)\right]\right)\\
	& +  \sum_{\mu\geq\lambda_3}P_{\mu}\left(\chi_{\nu\tau-R<\min\{\hbar_1,\hbar_3\}}\phi_{1,\geq \mu}\cdot\left[D_\tau^+\phi_{3,\lambda_3} - P_3\left(D_\tau^+\phi_{3,\lambda_3}\right)\right]\right)\\
	& = :B_1 + B_2. 
	\end{split}\end{equation}
	To handle the term $B_1$ we perform integration by parts after a further subdivision to shift derivatives (inherent in the definition of the norm) from the outside to the inner lower frequency factors, as in the proof of Prop.~\ref{prop:bilin2}. Specifically, write 
	\begin{align*}
	&\mathcal{F}(B_1)(\xi)\\& = \sum_{\hbar^{-2}\geq \mu\geq\lambda_3}\chi_{\xi\simeq\mu}\left\langle\phi_{m}(R,\xi), \left(\chi_{\nu\tau-R<\min\{\hbar_1,\hbar_3\}}\phi_{1,[\lambda_3,\mu]}\cdot\left[D_\tau^+\phi_{3,\lambda_3} - P_3\left(D_\tau^+\phi_{3,\lambda_3}\right)\right]\right)\right\rangle_{L^2_{R\,dR}}\\
	& +  \sum_{\mu\geq\max\{\lambda_3,\hbar^{-2}\}}\xi^{-3}\chi_{\xi\simeq\mu}\left\langle\phi_{m}(R,\xi), H_{m}^3\left(\chi_{\nu\tau-R<\min\{\hbar_1,\hbar_3\}}\phi_{1,[\lambda_3,\mu]}\cdot\left[D_\tau^+\phi_{3,\lambda_3} - P_3\left(D_\tau^+\phi_{3,\lambda_3}\right)\right]\right)\right\rangle_{L^2_{R\,dR}}
	\end{align*}
	The first term on the right is estimated by using a point wise bound on the inner product and using Holder's inequality to bound the $S_1^{\hbar}$-norm of the output. In fact, note that we have 
	\begin{align*}
	&\left\|\sum_{\hbar^{-2}\geq \mu\geq\lambda_3}\chi_{\xi\simeq\mu}\left\langle\phi_{m}(R,\xi), \left(\chi_{\nu\tau-R<\min\{\hbar_1,\hbar_3\}}\phi_{1,[\lambda_3,\mu]}\cdot\left[D_\tau^+\phi_{3,\lambda_3} - P_3\left(D_\tau^+\phi_{3,\lambda_3}\right)\right]\right)\right\rangle_{L^2_{R\,dR}}\right\|_{S_1^{\hbar}}\\
	&\lesssim \sum_{\hbar^{-2}\geq \mu\geq\lambda_3}\hbar\cdot\left(\hbar^2\mu\right)^{\frac12-\frac{\delta}{2}}\cdot\mu^{\frac12}\cdot \left\|\left\langle\phi_{m}(R,\xi), \left(\chi_{\nu\tau-R<\min\{\hbar_1,\hbar_3\}}\phi_{1,[\lambda_3,\mu]}\cdot\left[D_\tau^+\phi_{3,\lambda_3} - P_3\left(D_\tau^+\phi_{3,\lambda_3}\right)\right]\right)\right\rangle_{L^2_{R\,dR}}\right\|_{L^\infty_{d\xi}}.
	\end{align*}
In view of the estimates (see Lemma~\ref{lem:extraLinftyrefinedbound})
\begin{align*}
	\left|\calD_{\tau}^{+}\phi_{3,\lambda_{3}}(\tau,R)\right|\lesssim \hbar_{3}^{-\frac32}\left\|\calF\left(\calD_{\tau}^{+}\phi_{3,\lambda_{3}}\right)\right\|_{\Shth_{1}},\quad \left|\chi_{\nu\tau-R<\min\{\hbar_1,\hbar_3\}}P_3\left(D_\tau^+\phi_{3,\lambda_3}\right)\right|\lesssim  \hbar_3^{-\frac32}\left\|\mathcal{F}\left(D_\tau^+\phi_{3,\lambda_3}\right)\right\|_{S_1^{\hbar_3}},
\end{align*}
as well as 
the fact that the factors $(\nu\tau - R)^j$ compensate for the frequency loss $\lambda_3^{\frac{j}{2}}$ from the coefficients in the Taylor polynomial due to the restrictions $\nu\tau-R<\min\{\hbar_1,\hbar_3\} = \hbar_3$ and $\hbar_3\lesssim \lambda_3^{-\frac12}$, we can finally estimate (for $R\simeq\tau$)
\begin{align*}
&\left\|\left\langle\phi_{m}(R,\xi), \left(\chi_{\nu\tau-R<\min\{\hbar_1,\hbar_3\}}\phi_{1,[\lambda_3,\mu]}\cdot\left[D_\tau^+\phi_{3,\lambda_3} - P_3\left(D_\tau^+\phi_{3,\lambda_3}\right)\right]\right)\right\rangle_{L^2_{R\,dR}}\right\|_{L^\infty_{d\xi}}\\
&\lesssim \tau\hbar_3\cdot\left\|\phi_{m}(R,\xi)\right\|_{L^\infty_{R\,dR}}\cdot \left\|\chi_{\nu\tau-R<\min\{\hbar_1,\hbar_3\}}\phi_{1,[\lambda_3,\mu]}\right\|_{L^\infty_{R\,dR}}\cdot\left\|D_\tau^+\phi_{3,\lambda_3} - P_3\left(D_\tau^+\phi_{3,\lambda_3}\right)\right\|_{L^\infty_{R\,dR}}, 
\end{align*}
which in the regime where $\nu R\xi^{\frac12}\hbar$ is away from the turning point (at output angular momentum $m$) can be bounded by 
\begin{align*}
\hbar_1^{-1}\cdot \tau\hbar_3\cdot \hbar^{\frac12}\cdot \tau^{-\frac32-\nu}(\log\tau)^C\cdot \hbar_3^{-\frac32}\left\|\mathcal{F}\left(D_\tau^+\phi_{3,\lambda_3}\right)\right\|_{S_1^{\hbar_3}}\cdot \big\|\tilde{\tilde{\yb}}_1\big\|_{adm}\lesssim_{\hbar_1} \tau^{-\frac12-\nu}(\log\tau)^C\cdot\left\|\mathcal{F}\left(D_\tau^+\phi_{3,\lambda_3}\right)\right\|_{S_1^{\hbar_3}}\cdot\big\|\tilde{\tilde{\yb}}_1\big\|_{adm},
\end{align*}
where we have taken advantage of the assumption that $\hbar\simeq\hbar_3$. If $\nu R\xi^{\frac12}\hbar$ is near the turning point, one a priori only gains $\hbar^{\frac13}$ for $\left|\phi_{m}(R,\xi)\right|$, and an extra power $\hbar^{\frac16}$ needs to be gained from the shortness of the $\xi$-integration interval, in analogy to arguments in the proof of Prop.~\ref{prop:derivative1}. In fact, recall from the proof for Proposition \ref{prop:derivative1}, $\xi$ is confined in an interval is length $\simeq\frac{\hbar^{\frac23}}{\left(R\hbar\right)^{2}}$. Then in using H\"older inequality to get the pointwise bound for $D_{\tau}^{+}\phi_{3,\lambda_3}$, we encounter the following ($I_{R}$ being the interval where $\xi$ lies in)
\begin{align*}
	\left(\int_{I_{R}}\xi^{-1+\delta}\,d\xi\right)^{\frac12}\lesssim \hbar_{3}^{\frac13}\cdot\hbar_{3}^{-\delta},
\end{align*}
which together with the factor $\hbar^{\frac13}_{3}\cdot\hbar_{3}^{-2+\delta}$ in front of the $\xi$-integral gives the desired extra factor $\hbar_{3}^{\frac13}$ (more than needed).
\\
To conclude the estimate for $B_1$, we still need to handle the case of large output frequencies $\mu\geq \max\{\hbar^{-2},\lambda_3\}$, i.e., the sum 
\begin{align*}
&\sum_{\mu\geq\max\{\lambda_3,\hbar^{-2}\}}\xi^{-3}\chi_{\xi\simeq\mu}\left\langle\phi_{m}(R,\xi), H_{m}^3\left(\chi_{\nu\tau-R<\min\{\hbar_1,\hbar_3\}}\phi_{1,[\lambda_3,\mu]}\cdot\left[D_\tau^+\phi_{3,\lambda_3} - P_3\left(D_\tau^+\phi_{3,\lambda_3}\right)\right]\right)\right\rangle_{L^2_{R\,dR}}
\end{align*}
Using the Leibniz rule, this can be written as linear combination of terms of the form 
\begin{align*}
&\sum_{\mu\geq\max\{\lambda_3,\hbar^{-2}\}}\xi^{-3}\chi_{\xi\simeq\mu}\left\langle\phi_{m}(R,\xi), \left(\frac{m^2}{R^2}\right)^l\left(\partial_R^i\chi_{\nu\tau-R<\min\{\hbar_1,\hbar_3\}}\partial_R^j\phi_{1,[\lambda_3,\mu]}\cdot\partial_R^k\left[D_\tau^+\phi_{3,\lambda_3} - P_3\left(D_\tau^+\phi_{3,\lambda_3}\right)\right]\right)\right\rangle_{L^2_{R\,dR}}
\end{align*}
where $2l+i+j+k = 6$. We consider the extreme cases $2l = 6, j = 6, k = 6$, the other situations being handled similarly. 
\\
{\it{$2l = 6$.}} Observe that 
\begin{align*}
\left\|\sum_{\mu\geq\max\{\lambda_3,\hbar^{-2}\}}\xi^{-3}\chi_{\xi\simeq\mu}\left\langle\ldots\right\rangle_{L^2_{R\,dR}}\right\|_{S_1^{\hbar}}\leq \sum_{\mu\geq\max\{\lambda_3,\hbar^{-2}\}}\hbar^7\cdot \left\|\left(\xi\hbar^2\right)^{-1+\frac{\delta}{2}}\left\langle\ldots\right\rangle_{L^2_{R\,dR}}\right\|_{L^2_{d\xi}(\xi\simeq\mu)}
\end{align*}
Using Planhcerel's theorem for the distorted Fourier transform at angular momentum $m$, we can then bound 
\begin{align*}
&\hbar^7\cdot \left\|\left(\xi\hbar^2\right)^{-1+\frac{\delta}{2}}\left\langle\ldots\right\rangle_{L^2_{R\,dR}}\right\|_{L^2_{d\xi}(\xi\simeq\mu)}\\&\lesssim \hbar^7\cdot \left(\mu\hbar^2\right)^{-1+\frac{\delta}{2}}\cdot \left\|\frac{m^6}{R^6}\chi_{\nu\tau-R<\min\{\hbar_1,\hbar_3\}}\phi_{1,[\lambda_3,\mu]}\cdot\left[D_\tau^+\phi_{3,\lambda_3} - P_3\left(D_\tau^+\phi_{3,\lambda_3}\right)\right]\right\|_{L^2_{R\,dR}},
\end{align*}
and we bound the last term by using Lemma~\ref{lem:extraLinftyrefinedbound} and Holder's inequality:
\begin{align*}
&\left\|\chi_{\nu\tau-R<\min\{\hbar_1,\hbar_3\}}D_\tau^+\phi_{3,\lambda_3}\right\|_{L^2_{R\,dR}}\lesssim  \hbar_3^{-1}\left\|D_\tau^+\phi_{3,\lambda_3}\right\|_{\tilde{S}_1^{(\hbar_3)}},\\
&\left\|\chi_{\nu\tau-R<\min\{\hbar_1,\hbar_3\}}P_3\left(D_\tau^+\phi_{3,\lambda_3}\right)\right\|_{L^2_{R\,dR}}\lesssim \hbar_3^{\frac12}\cdot \hbar_3^{\frac12}\cdot \hbar_3^{-2}\cdot\left\|D_\tau^+\phi_{3,\lambda_3}\right\|_{\tilde{S}_1^{(\hbar_3)}},
\end{align*}
where the second bound follows from Holder's inequality and the estimates 
\begin{align*}
\left|\left(D_\tau^+\phi_{3,\lambda_3}\right)^{(j)}(\tau, \nu\tau)\right|\lesssim \hbar_3^{\frac12}\cdot\hbar_3^{-2-j}\cdot \left\|D_\tau^+\phi_{3,\lambda_3}\right\|_{\tilde{S}_1^{(\hbar_3)}},\quad 0\leq j\leq 3. 
\end{align*}
Combining the preceding estimates with 
\[
\left\|\chi_{\nu\tau-R<\min\{\hbar_1,\hbar_3\}}\phi_{1,[\lambda_3,\mu]}\right\|_{L^\infty_{R\,dR}}\lesssim_{\hbar_1} \lambda_3^{-\frac14}\cdot\tau^{-\frac12}\cdot\tau^{-1-\nu}(\log\tau)^C\cdot \big\|\tilde{\tilde{\yb}}_1\big\|_{adm}, 
\]
and further taking advantage of $\hbar\simeq\hbar_3$ by our assumption as well as $m\simeq\hbar^{-1}$, we infer 
\begin{align*}
\hbar^7\cdot \left\|\left(\xi\hbar^2\right)^{-1+\frac{\delta}{2}}\left\langle\ldots\right\rangle_{L^2_{R\,dR}}\right\|_{L^2_{d\xi}(\xi\simeq\mu)}\lesssim_{\hbar_1} \tau^{-\frac{13}{2}-\nu}(\log\tau)^C\cdot\left\|\left(D_\tau^+\phi_{3,\lambda_3}\right)\right\|_{S_1^{\hbar_3}}\cdot \big\|\tilde{\tilde{\yb}}_1\big\|_{adm},
\end{align*}
where we exploit that $R\simeq \tau$ on the support of the expression. Note that in the case $i = 6$, which is formally similar we only gain a factor $\tau^{-\frac{1}{2}-\nu}(\log\tau)^C$. 
\\

{\it{$j = 6$.}} Here we exploit the fact that the factor $D_\tau^+\phi_{3,\lambda_3} - P_3\left(D_\tau^+\phi_{3,\lambda_3}\right)$ compensates for the singularity\footnote{Of course this function is actually $C^\infty$ due to the frequency cutoff but we need bounds which are uniform in $\mu$} of $\partial_R^6\phi_{1,[\lambda_3,\mu]}$. Precisely, write 
\begin{equation}\label{eq:Tayloridentity}
\left[D_\tau^+\phi_{3,\lambda_3} - P_3\left(D_\tau^+\phi_{3,\lambda_3}\right)\right](\tau, R) = C\left(\nu\tau - R\right)^4\cdot\int_0^1\left(D_\tau^+\phi_{3,\lambda_3}\right)^{(4)}(\tau, \nu\tau - s(\nu\tau - R))\,ds,
\end{equation}
where using Prop.~\ref{prop:derivative1} as well as Holder's inequality and a simple change of variables we have 
\begin{align*}
\left|\int_0^1\left(D_\tau^+\phi_{3,\lambda_3}\right)^{(4)}(\tau, \nu\tau - s(\nu\tau - R))\,ds\right|\lesssim \hbar_3^{-5}\cdot \left\|D_\tau^+\phi_{3,\lambda_3}\right\|_{\tilde{S}_1^{(\hbar_3)}}\cdot (\nu\tau - R)^{-\frac12}. 
\end{align*}
Furthermore using Lemma~\ref{lem:singFouriertiphysicalngeq2adm} we infer the bound (by re-arranging the powers on $(\nu\tau-R)$ and $\mu$)
\[
\left|\chi_{\nu\tau-R<\min\{\hbar_1,\hbar_3\}}\partial_R^6\phi_{1,[\lambda_3,\mu]}\right|\lesssim \mu^{1-}\tau^{-\frac32-\nu}(\log\tau)^C\cdot \chi_{\nu\tau-R<\min\{\hbar_1,\hbar_3\}}\cdot(\nu\tau - R)^{-\frac72+\nu-}\cdot \big\|\tilde{\tilde{\yb}}_1\big\|_{adm},
\]
whence in total we get the bound 
\begin{align*}
\left\|\chi_{\nu\tau-R<\min\{\hbar_1,\hbar_3\}}\partial_R^6\phi_{1,[\lambda_3,\mu]}\left[D_\tau^+\phi_{3,\lambda_3} - P_3\left(D_\tau^+\phi_{3,\lambda_3}\right)\right]\right\|_{L^2_{R\,dR}}\lesssim \hbar_3^{\frac12-5}\cdot \mu^{1-}\tau^{-1-\nu}(\log\tau)^C\cdot \left\|D_\tau^+\phi_{3,\lambda_3}\right\|_{\tilde{S}_1^{(\hbar_3)}}\cdot \big\|\tilde{\tilde{\yb}}_1\big\|_{adm}.
\end{align*}
Using Plancherel's theorem for the distorted Fourier transform, we infer 
\begin{align*}
&\left\|\sum_{\mu\geq\max\{\lambda_3,\hbar^{-2}\}}\xi^{-3}\chi_{\xi\simeq\mu}\left\langle\phi_{m}(R,\xi), \left(\chi_{\nu\tau-R<\min\{\hbar_1,\hbar_3\}}\partial_R^6\phi_{1,[\lambda_3,\mu]}\cdot\left[D_\tau^+\phi_{3,\lambda_3} - P_3\left(D_\tau^+\phi_{3,\lambda_3}\right)\right]\right)\right\rangle_{L^2_{R\,dR}}\right\|_{S_1^{\hbar}}\\
&\lesssim \sum_{\mu\geq\max\{\lambda_3,\hbar^{-2}\}}\hbar^{5+\frac{\delta}{2}}\mu^{-(1-)}\cdot\left\|\chi_{\nu\tau-R<\min\{\hbar_1,\hbar_3\}}\partial_R^6\phi_{1,[\lambda_3,\mu]}\cdot\left[D_\tau^+\phi_{3,\lambda_3} - P_3\left(D_\tau^+\phi_{3,\lambda_3}\right)\right]\right\|_{L^2_{R\,dR}}\\
&\lesssim \sum_{\mu\geq\max\{\lambda_3,\hbar^{-2}\}}\mu^{-(0+)}\cdot \tau^{-1-\nu}(\log\tau)^C\cdot \left\|D_\tau^+\phi_{3,\lambda_3}\right\|_{\tilde{S}_1^{(\hbar_3)}}\cdot \big\|\tilde{\tilde{\yb}}_1\big\|_{adm}.
\end{align*}

{\it{$k =6$}}. Here we use the estimate (note that $\hbar_{3}^{2}\lambda_{3}>1$) 
\begin{align*}
\left\|\partial_R^6\left[D_\tau^+\phi_{3,\lambda_3} - P_3\left(D_\tau^+\phi_{3,\lambda_3}\right)\right]\right\|_{L^2_{R\,dR}}\lesssim \hbar_3^{-5-\delta}\cdot \lambda_3^{1-\frac{\delta}{2}}\cdot \left\|D_\tau^+\phi_{3,\lambda_3}\right\|_{\tilde{S}_1^{(\hbar_3)}}. 
\end{align*}
Then using the simple point wise bound 
\[
\left\|\chi_{\nu\tau-R<\min\{\hbar_1,\hbar_3\}}\phi_{1,[\lambda_3,\mu]}\right\|_{L^\infty_{R\,dR}}\lesssim_{\hbar_1} \tau^{-\frac32-\nu}(\log\tau)^C\cdot \big\|\tilde{\tilde{\yb}}_1\big\|_{adm},
\]
we obtain by using the Plancherel's theorem for the distorted Fourier transform 
\begin{align*}
&\left\|\sum_{\mu\geq\max\{\lambda_3,\hbar^{-2}\}}\xi^{-3}\chi_{\xi\simeq\mu}\langle\phi_{m}(R,\xi), \left(\chi_{\nu\tau-R<\min\{\hbar_1,\hbar_3\}}\phi_{1,[\lambda_3,\mu]}\cdot\partial_R^6\left[D_\tau^+\phi_{3,\lambda_3} - P_3\left(D_\tau^+\phi_{3,\lambda_3}\right)\right]\right)\rangle_{L^2_{R\,dR}}\right\|_{S_1^{\hbar}}\\
&\lesssim \Big(\sum_{\mu\geq\max\{\lambda_3,\hbar^{-2}\}}\big[\hbar^{5+\delta}\mu^{-(1-\frac{\delta}{2})}\cdot \left\|\chi_{\nu\tau-R<\min\{\hbar_1,\hbar_3\}}\phi_{1,[\lambda_3,\mu]}\right\|_{L^\infty_{R\,dR}}\cdot \left\|\partial_R^6\left[D_\tau^+\phi_{3,\lambda_3} - P_3\left(D_\tau^+\phi_{3,\lambda_3}\right)\right]\right\|_{L^2_{R\,dR}}\big]^2\Big)^{\frac12}\\
&\lesssim_{\hbar_1} \tau^{-\frac32-\nu}(\log\tau)^C\cdot \left\|D_\tau^+\phi_{3,\lambda_3}\right\|_{\tilde{S}_1^{(\hbar_3)}}\cdot \big\|\tilde{\tilde{\yb}}_1\big\|_{adm},
\end{align*}
again exploiting that $\hbar\simeq \hbar_3$. 
\\
This concludes the bound for the term $B_1$, recalling \eqref{eq:lemmaN0abB1B2}. 
\\

We next turn to the term $B_2$ there, where we have to take advantage of the large frequency in the singular term. Precisely, we use that $\phi_{1,\geq\mu}$ decays rapidly beyond scale $\mu^{-\frac12}$ with respect to $\nu\tau - R$. If we then again invoke \eqref{eq:Tayloridentity}, we find the bound 
\begin{align*}
\left\|\chi_{\nu\tau-R<\min\{\hbar_1,\hbar_3\}}\phi_{1,\geq \mu}\cdot\left[D_\tau^+\phi_{3,\lambda_3} - P_3\left(D_\tau^+\phi_{3,\lambda_3}\right)\right]\right\|_{L^2_{R\,dR}}\lesssim_{\hbar_1} \hbar_3^{-5-\delta}\mu^{-2-\frac{\delta}{2}-\frac14}\cdot \tau^{-1-\nu}(\log\tau)^C\cdot \left\|D_\tau^+\phi_{3,\lambda_3}\right\|_{\tilde{S}_1^{(\hbar_3)}}\cdot \big\|\tilde{\tilde{\yb}}_1\big\|_{adm}.
\end{align*}
Here following the proof for Lemma \ref{lem:singFouriertiphysicalngeq2adm}, we integrate by parts to gain the decay in $\mu^{-1}$. In the meantime we also see powers in $(\nu\tau-R)^{-1}$. We can have up to $(\nu\tau-R)^{-4}$, to be compensated by the power $(\nu\tau-R)^{4}$ from the expression for $D_\tau^+\phi_{3,\lambda_3} - P_3\left(D_\tau^+\phi_{3,\lambda_3}\right)$. \\

Then using the Plancherel's theorem for the distorted Fourier transform, and the fact that by assumption $\hbar\simeq\hbar_3$, we get 
\begin{align*}
&\left\|\sum_{\mu\geq\lambda_3}P_{\mu}\left(\chi_{\nu\tau-R<\min\{\hbar_1,\hbar_3\}}\phi_{1,\geq \mu}\cdot\left[D_\tau^+\phi_{3,\lambda_3} - P_3\left(D_\tau^+\phi_{3,\lambda_3}\right)\right]\right)\right\|_{\tilde{S}_1^{(\hbar)}}\\
&\lesssim_{\hbar_1} \Big(\sum_{\mu}\big|\sum_{\lambda_3\leq \mu}\mu^{-\frac14}\cdot \tau^{-1-\nu}(\log\tau)^C\cdot \left\|D_\tau^+\phi_{3,\lambda_3}\right\|_{\tilde{S}_1^{(\hbar_3)}}\big|^2\Big)^{\frac12}\cdot \big\|\tilde{\tilde{\yb}}_1\big\|_{adm}\\
&\lesssim\Big(\sum_{\lambda_3}\big( \lambda_3^{-\frac14}\cdot \tau^{-1-\nu}(\log\tau)^C\cdot \left\|D_\tau^+\phi_{3,\lambda_3}\right\|_{\tilde{S}_1^{(\hbar_3)}}\big)^2\Big)^{\frac12}\cdot \big\|\tilde{\tilde{\yb}}_1\big\|_{adm}\\
&\lesssim \tau^{-1-\nu}(\log\tau)^C\cdot \left\|D_\tau^+\phi_{3}\right\|_{\langle\xi\rangle^{\frac14}\cdot \tilde{S}_1^{(\hbar_3)}}\cdot \big\|\tilde{\tilde{\yb}}_1\big\|_{adm}.
\end{align*}
This completes the estimate for $B$ and thereby finally for case (II) under the assumption $\hbar_3\ll\hbar_1$. The remaining cases $\hbar_3\simeq \hbar_1$, $\hbar_3\gg\hbar_1$ are handled analogously. 
\\

{\it{(III): first factor of type $\tilde{\yb}_1$, second factor of type $\tilde{\tilde{\yb}}_3$.}} It is in this case where the particular paraproduct structure of the expression in the lemma becomes important. Write 
\begin{align*}
&\chi_{R\gtrsim\tau}\phi_{1,\geq\lambda_3}\cdot \left(\partial_{\tau}+ \frac{\lambda_{\tau}}{\lambda}R\partial_R + \partial_R\right)\phi_{3,\lambda_3}\\
&= \chi_{R\gtrsim\tau}\phi_{1,\geq\lambda_3}(\tau, \nu\tau)\cdot \left(\partial_{\tau}+ \frac{\lambda_{\tau}}{\lambda}R\partial_R + \partial_R\right)\phi_{3,\lambda_3}\\
&+\chi_{R\gtrsim\tau}\left[\phi_{1,\geq\lambda_3}(\tau, R) - \phi_{1,\geq\lambda_3}(\tau, \nu\tau)\right]\cdot \left(\partial_{\tau}+ \frac{\lambda_{\tau}}{\lambda}R\partial_R + \partial_R\right)\phi_{3,\lambda_3},\\
\end{align*}
where in light of Lemma~\ref{lem:singFouriertiphysicalngeq2adm} we may assume that $\phi_{3,\lambda_3}$ admits the expansion of the middle term in that lemma. Then for the first term on the right, we argue exactly as in the proof of Lemma~\ref{lem:simpleparaproduct}, to conclude that it is an admissibly singular source term of level $l_1\geq 1$. This can be seen as follows: Suppose $\phi_{3,\lambda_{3}}$ contains an principal ingoing part. Then the Fourier transform of $\left(\partial_{\tau}+ \frac{\lambda_{\tau}}{\lambda}R\partial_R + \partial_R\right)\phi_{3,\lambda_{3}}$ contains a factor of $\lambda_{3}^{\frac12}$. On the other hand for $\chi_{R\gtrsim\tau}\phi_{1,\geq\lambda_3}$, using the pointwise decay for Fourier basis, we gain  a factor of $\lambda_{3}^{-\frac14}$. Therefore in total we have $\lambda_{3}^{\frac14}$. Since this is a source term, the power appearing in the admissible term is $\lambda_{3}^{-\frac14}$, which corresponds  to  $l_{1}=1$. The second term on the right is again in $\tilde{S}_1^{(\hbar)}$ by essentially the same argument as the one used at the end in the preceding case. 
\\

{\it{(IV): Both factors of type $\tilde{\tilde{y}}_j$, i.e., admissibly singular.}} This is handled by using Lemma~\ref{lem:singFouriertiphysicalngeq2admDeriv}, Lemma~\ref{lem:singFouriertiphysicalngeq2adm}, and translating things back to the Fourier side via Lemma~\ref{lem:singPhysicaltoFourierngeq2}. 
\\

This concludes the proof of the lemma. 
\end{proof}
In order to complete the proof of the proposition as far as the contribution of the term $\mathcal{N}_{02a}$ is concerned, we need one more bilinear lemma
\begin{lemma}\label{lem:N02abilin2} Let $F(\tau, R)$ be an angular momentum $n_{1}$ function and $\{m,n_1,n_2\}$ be an admissible triple of momenta\footnote{We recall that this means $|m - n_1 - n_2| = O(1)$.} all of absolute value $\geq 2$, and assume that $F = F_1+F_2$ is an angular momentum $n_1$ function which can be written as the sum of an angular momentum $n_1$ function $F_1\in \tilde{S}_1^{(\hbar_1)}$, while $F_2$ is an angular momentum $n_1$ admissibly singular source function of level $l_1\geq 1$, where we recall Definition~\ref{defi:xsingularsourcetermsngeq2adm}. 
	Then if $\phi$ is an angular momentum $n_2$ function as in Prop.~\ref{prop:goodcubicnullform1}, then 
	\[
	\left(\partial_{\tau}+ \frac{\lambda_{\tau}}{\lambda}R\partial_R - \partial_R\right)\phi\cdot F|_{R<\nu\tau} = g(\tau, R)|_{R<\nu\tau}, 
	\]
	where (with $\hbar = \frac{1}{1+m}$)
	\[
	g(\tau, R) = \int_0^\infty \phi_{m}(R,\xi)\cdot \zb(\tau, \xi)\rho_{m}(\xi)\,d\xi
	\]
	with $\zb = \zb_1+\zb_2$ and $\zb_1$ an admissibly singular source term at level $l_1\geq 1$, while $\zb_2\in S_1^{\hbar}$. Moreover, we have the bound 
	\begin{align*}
	&\big\| \zb_1\big\|_{sourceadm} + \tau^4\cdot \big\|\zb_2\big\|_{S_1^{\hbar}}\\&\lesssim \min\{|n_{1,2}|\}^{C}\cdot\Big(\tau^2\cdot\big\|F_1\big\|_{\tilde{S}_1^{(\hbar_1)}} + \big\|F_2\big\|_{sourceadm}\Big)\cdot \Big(\big\| \overline{y}\big\|_{adm} + \sup_{\tau\geq\tau_0}\tau^3\cdot \big(\big\| \overline{z}\big\|_{S_0^{\hbar_j}} + \big\| \mathcal{D}_{\tau}\overline{z}\big\|_{S_1^{\hbar_j}}\big)\Big),
	\end{align*}
where $\overline{x} = \overline{y} + \overline{z}$ is the distorted Fourier transform of $\phi$. 
\end{lemma}
\begin{proof}(sketch) This follows again by distinguishing between different cases, as in the proof of the preceding lemma. Note that the `good derivative' $\partial_{\tau}+ \frac{\lambda_{\tau}}{\lambda}R\partial_R - \partial_R$ sends admissibly singular functions into functions of almost the same regularity, albeit at the loss of regularity of some of the coefficient functions in the definition of admissibly singular functions. For this recall  Lemma~\ref{lem:singFouriertiphysicalngeq2admDeriv}. These functions are then good enough to be source admissible, see also Lemma~\ref{lem:basicnullformforsingularinputs} and its proof.  The details are rather similar to the ones in the preceding proof. Note that if $\phi$ is admissibly singular and $F\in \tilde{S}_1^{(\hbar_1)}$ and $\hbar_1\ll\hbar_2$, then since $\left(\partial_{\tau}+ \frac{\lambda_{\tau}}{\lambda}R\partial_R - \partial_R\right)\phi$ will be of essentially the same regularity(recalling Lemma~\ref{lem:singFouriertiphysicalngeq2admDeriv}), we can absorb the $\hbar_1^{-\frac32}$-loss coming from controlling $F$ via Lemma~\ref{lem:extraLinftyrefinedbound}  since we are allowed to lose a factor $\hbar_1^{-1}$ in light of the definition of source admissibly singular terms, and we can absorb another factor $\hbar_1^{-1}$ since the resulting expression will be of level $l_1\geq 2$. 
\end{proof}
Combining the preceding two lemmas and summing over $\lambda_3<\lambda_2$ completes the proof for the contribution of $\mathcal{N}_{02a}$. Observe that a conclusion of the preceding considerations is that $\mathcal{N}_{02a}$ is above minimal regularity (as expressed by the statement concerning the level $l_1\geq 1$). 
\\

{\it{The contribution of $\mathcal{N}_{02b}$}}. This is handled in the same manner, breaking things into a number of sub-steps. To begin with, we re-arrange the terms as follows:
\begin{equation}\label{eq:N02brearranged}\begin{split}
&\sum_{\lambda_3<\lambda_2} \phi_{1,\geq\lambda_3}\cdot \left(\partial_{\tau}+ \frac{\lambda_{\tau}}{\lambda}R\partial_R + \partial_R\right)\phi_{2,\lambda_2}\cdot  \left(\partial_{\tau}+ \frac{\lambda_{\tau}}{\lambda}R\partial_R - \partial_R\right)\phi_{3,\lambda_3}\\
& =  \left(\partial_{\tau}+ \frac{\lambda_{\tau}}{\lambda}R\partial_R + \partial_R\right)\phi_{2}\cdot \left(\sum_{\lambda_3} \phi_{1,\geq\lambda_3}\cdot  \left(\partial_{\tau}+ \frac{\lambda_{\tau}}{\lambda}R\partial_R - \partial_R\right)\phi_{3,\lambda_3}\right)\\
& - \sum_{\lambda_3\geq\lambda_2} \phi_{1,\geq\lambda_3}\cdot \left(\partial_{\tau}+ \frac{\lambda_{\tau}}{\lambda}R\partial_R + \partial_R\right)\phi_{2,\lambda_2}\cdot  \left(\partial_{\tau}+ \frac{\lambda_{\tau}}{\lambda}R\partial_R - \partial_R\right)\phi_{3,\lambda_3}
\end{split}\end{equation}
Here the second term on the right will again have better than minimum regularity (it will be seen to be at least admissibly singular source of level $l_1\geq 1$), while the first term on the right may be of minimal regularity. First, we need the following lemma analogous to Lemma~\ref{lem:N02abilin1}:
\begin{lemma}\label{lem:N02bbilin1} Let $\phi_1, \phi_3$ be angular  momentum $n_j,\,j = 1,3$ functions (with $|n_j|\geq 2$) and with the property that their angular momentum $n_j$ distorted Fourier transforms $\yb_j(\tau, \xi)$ can be split into 
	\[
	\yb_j(\tau, \xi) = \tilde{\yb}_j(\tau, \xi) + \tilde{\tilde{\yb}}_j(\tau, \xi),
	\]
	where $ \tilde{\yb}_j(\tau, \xi)\in S_0^{\hbar_j},\,\mathcal{D}_{\tau}\tilde{\yb}_j(\tau, \xi)\in S_1^{\hbar_j}$, and further $\tilde{\tilde{\yb}}_j$ is admissibly singular. Then we can write 
	\[
	\chi_{R\gtrsim \tau}\phi_{1,\geq\lambda_3}\cdot \left(\partial_{\tau}+ \frac{\lambda_{\tau}}{\lambda}R\partial_R - \partial_R\right)\phi_{3,\lambda_3}|_{R<\nu\tau} = g(\tau, R)|_{R<\nu\tau}
	\]
	where, if $\{m,n_1,n_3\}$ with $|m|\geq 2$ is an admissible angular momentum triple in the sense of Prop.~\ref{prop:bilin2}, we have an angular momentum $m$ Fourier representation of $g(\tau, R)$ 
	\[
	g(\tau, R) = \int_0^\infty \phi_{m}(R,\xi)\zb(\tau, \xi)\rho_{m}(\xi)\,d\xi,\quad \hbar = \frac{1}{m+1}
	\]
	where $\zb$ is the sum of a function in $\zb_1\in \left\langle\xi\right\rangle^{-\frac14-}S_1^{\hbar}$ and a function $\zb_2$ which is source admissibly singular at angular momentum $m$ at level $l_1\geq 2$. Furthermore, we have the bound 
	\begin{align*}
	\tau^3\cdot \big\|\zb_1\big\|_{ \left\langle\xi\right\rangle^{-\frac14-}S_1^{\hbar}} + \big\|\zb_2\big\|_{sourceadm} \lesssim  \max_{j=1,3}\{|n_j|\}^{-C}\cdot \prod_{j=1,3}|n_j|^C\Big(\big\| \tilde{\tilde{\yb}}_j\big\|_{adm} + \sup_{\tau\geq\tau_0}\tau^3\cdot \big(\big\| \tilde{\yb}_j(\tau,\cdot)\big\|_{S_0^{\hbar_j}} + \big\| \mathcal{D}_{\tau}\tilde{\yb}_j(\tau,\cdot)\big\|_{S_1^{\hbar_j}}\big)\Big).
	\end{align*}
\end{lemma}
The proof of this proceeds in analogy to the one of Lemma~\ref{lem:N02abilin1} and we omit it here. 
\\

To use the preceding lemma, the following refined $L^\infty$-type estimate shall be useful:
\begin{lemma}\label{lem:extraLinftyrefinedbound} Assume that $f(R)$ is an angular momentum $n, |n|\geq 2$ function represented by 
	\[
	f(R) = \int_0^\infty \phi_{n}(R,\xi)\cdot \xb(\xi)\cdot\rho_{n}(\xi)\,d\xi. 
	\]
	Then we have the point wise bound (with $\hbar=\frac{1}{n+1}$)
	\[
	\left|f(R)\right|\lesssim \hbar^{-1}\tau^{\frac12}\cdot\left\|\xb\right\|_{\langle\xi\rangle^{-\frac14-}S_1^{\hbar}},\quad R\lesssim \tau. 
	\]
	We also have the estimate 
	\begin{align*}
	\left|f(R)\right|\lesssim \hbar^{-\frac32}\cdot\left\|\xb\right\|_{S_1^{\hbar}},\,\left|\partial_R^jf(R)\right|\lesssim \hbar^{-\frac32-j}\cdot\left\|\xb\right\|_{S_1^{\hbar}},\,j = \{1, 2, 3\}. 
	\end{align*}
\end{lemma}
\begin{proof}
	(outline) We prove the first inequality, the second being similar. We first recall that away from the turning point we have the bound $\left|\phi_{n}(R,\xi)\right|\lesssim \hbar^{\frac12}$. Write 
	\begin{align*}
	\int_0^\infty \phi_{n}(R,\xi)\cdot \xb(\xi)\cdot\rho_{n}(\xi)\,d\xi &= \int_0^\infty \chi_{R\xi^{\frac12}\hbar\ll1}\phi_{n}(R,\xi)\cdot \xb(\xi)\cdot\rho_{n}(\xi)\,d\xi\\
	& + \int_0^\infty \chi_{R\xi^{\frac12}\hbar\gtrsim 1}\phi_{n}(R,\xi)\cdot \xb(\xi)\cdot\rho_{n}(\xi)\,d\xi.
	\end{align*}
	Then we can bound the first integral on the right by 
	\begin{align*}
	\left|\int_0^\infty \chi_{R\xi^{\frac12}\hbar\ll1}\phi_{n}(R,\xi)\cdot \xb(\xi)\cdot\rho_{n}(\xi)\,d\xi\right|&\lesssim \int_0^1  \chi_{R\xi^{\frac12}\hbar\ll1}\left|\phi_{n}(R,\xi)\right|\left|\xb(\xi)\right|\,d\xi\\
	& + \int_1^\infty  \chi_{R\xi^{\frac12}\hbar\ll1}\left|\phi_{n}(R,\xi)\right|\left|\xb(\xi)\right|\,d\xi\\
	&\lesssim c^{\hbar^{-1}}\left\|\xi^{\frac12-\frac{\delta}{2}}\xb(\xi)\right\|_{L^2_{d\xi}(\xi<1)}+ c^{\hbar^{-1}}\left\|\xi^{\frac12+\frac{\delta}{2}}\xb(\xi)\right\|_{L^2_{d\xi}(\xi\geq 1)}
	\end{align*}
	for some positive $c<1$, and the last two terms are easily seen to be bounded by $\ll \left\|\xb\right\|_{S_1^{\hbar}}$, which is better than what is needed. To control 
	\[
	\int_0^\infty \chi_{R\xi^{\frac12}\hbar\gtrsim 1}\phi_{n}(R,\xi)\cdot \xb(\xi)\cdot\rho_{n}(\xi)\,d\xi,
	\]
	use that $\tau\xi^{\frac12}\hbar\gtrsim R\xi^{\frac12}\hbar\gtrsim 1$, whence $\left(\xi^{\frac12}\hbar\right)^{-\frac12}\lesssim \tau^{\frac12}$, and so away from the turning point $R\xi^{\frac12}\hbar = x_t(\hbar)$, which we excise by the cutoff $\chi_1(\xi,R)$, we find by means of the Cauchy Schwarz inequality 
	\begin{align*}
	\left|\int_0^\infty \chi_1(\xi,R)\chi_{R\xi^{\frac12}\hbar\gtrsim 1}\phi_{n}(R,\xi)\cdot \xb(\xi)\cdot\rho_{n}(\xi)\,d\xi\right|
	&\lesssim \hbar^{\frac12}\left\|\chi_{\tau\xi^{\frac12}\hbar\gtrsim 1}\xi^{\frac12-\frac{\delta}{2}}\xb(\xi)\right\|_{L^2_{d\xi}(\xi<1)}\\
	& + \hbar^{\frac12-\delta}\left\|\chi_{\tau\xi^{\frac12}\hbar\gtrsim 1}\xi^{\frac12-\frac{\delta}{2}}\xb(\xi)\right\|_{L^2_{d\xi}\left(\hbar^{-2}\geq\xi\geq 1\right)}\\
	&+  \hbar^{\frac12+\delta}\left\|\chi_{\tau\xi^{\frac12}\hbar\gtrsim 1}\xi^{\frac12+\frac{\delta}{2}}\xb(\xi)\right\|_{L^2_{d\xi}\left(\xi\geq \hbar^{-2}\right)}\\
	\end{align*}
	and furthermore we have 
	\begin{align*}
	\hbar^{\frac12}\left\|\chi_{\tau\xi^{\frac12}\hbar\gtrsim 1}\xi^{\frac12-\frac{\delta}{2}}\xb(\xi)\right\|_{L^2_{d\xi}(\xi<1)}&\lesssim  \tau^{\frac12}\hbar\left\|\chi_{\tau\xi^{\frac12}\hbar\gtrsim 1}\xi^{\frac34-\frac{\delta}{2}}\xb(\xi)\right\|_{L^2_{d\xi}(\xi<1)}\\
	&\lesssim \tau^{\frac12}\hbar^{-1}\cdot\left\|\xb\right\|_{\langle\xi\rangle^{-\frac14-}S_1^{\hbar}},
	\end{align*}
	On the other hand for the large frequency contributions
	\begin{align*}
	\hbar^{\frac12-\delta}\left\|\chi_{\tau\xi^{\frac12}\hbar\gtrsim 1}\xi^{\frac12-\frac{\delta}{2}}\xb(\xi)\right\|_{L^2_{d\xi}\left(\hbar^{-2}\geq\xi\geq 1\right)}&\lesssim \tau^{\frac12}\hbar^{1-\delta}\left\|\chi_{\tau\xi^{\frac12}\hbar\gtrsim 1}\xi^{\frac34-\frac{\delta}{2}}\xb(\xi)\right\|_{L^2_{d\xi}(\hbar^{-2}\geq \xi\geq 1)}\\
	&\lesssim \tau^{\frac12}\hbar^{-1}\cdot \left\|\xb\right\|_{\left\langle\xi\right\rangle^{-\frac14}S_1^{\hbar}},
	\end{align*}
	and similarly 
	\begin{align*}
	\hbar^{\frac12+\delta}\left\|\chi_{\tau\xi^{\frac12}\hbar\gtrsim 1}\xi^{\frac12+\frac{\delta}{2}}\xb(\xi)\right\|_{L^2_{d\xi}\left(\xi\geq \hbar^{-2}\right)}&\lesssim \tau^{\frac12}\hbar^{1+\delta}\left\|\chi_{\tau\xi^{\frac12}\hbar\gtrsim 1}\xi^{\frac34+\frac{\delta}{2}}\xb(\xi)\right\|_{L^2_{d\xi}\left(\xi\geq \hbar^{-2}\right)}\\
	&\lesssim \tau^{\frac12}\hbar^{-1}\cdot \left\|\xb\right\|_{\left\langle\xi\right\rangle^{-\frac14}S_1^{\hbar}},
	\end{align*}
	which is as desired. It remains to deal with the turning point, which follows by means of the usual refined asymsptotics of $\phi_{n}(R,\xi)$ for $x = R\xi^{\frac12}\hbar$ near $x_t(\alpha;\hbar)$. This again can be achieved using the fact that $\xi$ lies in an interval of length $\simeq \frac{\hbar^{\frac23}}{(R\hbar)^{2}}$.
\end{proof}
To finish the argument for the first of the terms on the right of \eqref{eq:N02brearranged}, we need the following analogue of Lemma~\ref{lem:N02abilin2}:
\begin{lemma}\label{lem:N02bbilin2}
	Let $F(\tau, R)$ be an angular momentum $n_{1}$ function and $\{m,n_1,n_2\}$ be an admissible triple of momenta\footnote{Recall that this means $|m - n_1 - n_2| = O(1)$.} all of absolute value $\geq 2$, and assume that $F = F_1+F_2$ is an angular momentum $n_1$ function which can be written as the sum of an angular momentum $n_1$ function $F_1$ whose distorted Fourier transform satisfies $\overline{f}_1\in \left\langle\xi\right\rangle^{-\frac14}S_1^{\hbar_1}$, while $F_2$ is an angular momentum $n_1$ admissibly singular source function of level $l_1\geq 2$. 
	Then if $\phi_2$ is an angular momentum $n_2$ function as in Lemma~\ref{lem:N02abilin1} , then 
	\[
	\left(\partial_{\tau}+ \frac{\lambda_{\tau}}{\lambda}R\partial_R + \partial_R\right)\phi_{2}\cdot F|_{R<\nu\tau} = g(\tau, R)|_{R<\nu\tau}, 
	\]
	where (with $\hbar = \frac{1}{1+|m|}$)
	\[
	g(\tau, R) = \int_0^\infty \phi_{m}(R,\xi)\cdot \zb(\tau, \xi)\rho_{m}(\xi)\,d\xi
	\]
	with $\zb = \zb_1+\zb_2$ and $\zb_1$ an admissibly singular source term, while $\zb_2\in S_1^{\hbar}$. Furthermore, the following bound obtains: 
	\begin{align*}
	&\tau^4\cdot \big\|\zb_2\big\|_{ S_1^{\hbar}} + \big\| \zb_1\big\|_{sourceadm}\lesssim \min_{j=1,2}\{|n_j|\}^C\cdot \Big(\tau^3\cdot\big\|\overline{f}_1\big\|_{\left\langle\xi\right\rangle^{-\frac14}S_1^{\hbar_1}} + \big\|\overline{f}_2\big\|_{sourceadm}\Big)\\&\hspace{6cm}\cdot \Big(\big\| \tilde{\tilde{\yb}}_2\big\|_{adm} + \sup_{\tau\geq\tau_0}\tau^3\cdot \big(\big\| \tilde{\yb}_2(\tau,\cdot)\big\|_{S_0^{\hbar_2}} + \big\| \mathcal{D}_{\tau}\tilde{\yb}_2(\tau,\cdot)\big\|_{S_1^{\hbar_2}}\big)\Big)
	\end{align*}
	Here $ \tilde{\yb}_2 +  \tilde{\tilde{\yb}}_2$ is the distorted Fourier transform of $\phi_2$, and $\overline{f}_2$ is the distorted Fourier transform of $F_2$. 
\end{lemma}
\begin{proof}
	(sketch) This follows again by considering the various combinations of inputs allowed. The interaction of two $S_1^{\hbar_j}$-functions being routine by now in light of Prop~\ref{prop:bilin1}, we assume that $\phi_2$ is admissibly singular, and also $F_1\in \left\langle\xi\right\rangle^{-\frac14}S_1^{\hbar_1}$. Assuming $n_1\gg n_2$, say, whence $\hbar = \frac{1}{|m|+1}\simeq \hbar_1$, we split
	\begin{align*}
	\left(\partial_{\tau}+ \frac{\lambda_{\tau}}{\lambda}R\partial_R + \partial_R\right)\phi_{2}\cdot F &= \chi_{\nu\tau - R\gtrsim \hbar_1}\left(\partial_{\tau}+ \frac{\lambda_{\tau}}{\lambda}R\partial_R + \partial_R\right)\phi_{2}\cdot F\\
	& +  \chi_{\nu\tau - R\lesssim\hbar_1}\left(\partial_{\tau}+ \frac{\lambda_{\tau}}{\lambda}R\partial_R + \partial_R\right)\phi_{2}\cdot F
	\end{align*}
	The first term is easy to handle since it is the product of two $S_1^{\hbar_j}$-functions. For the second term, we may assume that $ \chi_{\nu\tau - R\lesssim\hbar_1}\left(\partial_{\tau}+ \frac{\lambda_{\tau}}{\lambda}R\partial_R + \partial_R\right)\phi_{2}$ is given in accordance with Lemma~\ref{lem:singFouriertiphysicalngeq2admDeriv} for the `connecting part' of the singularity, while we have to invoke Lemma~\ref{lem:DerivonPrincSing} to describe the contribution of the principal incoming singular part. Now if $\phi_2$ is of principal ingoing singular type, we decompose 
	\begin{align*}
	\chi_{\nu\tau - R\lesssim\hbar_1}\left(\partial_{\tau}+ \frac{\lambda_{\tau}}{\lambda}R\partial_R + \partial_R\right)\phi_{2}\cdot F &= \chi_{\nu\tau - R\lesssim\hbar_1}\left(\partial_{\tau}+ \frac{\lambda_{\tau}}{\lambda}R\partial_R + \partial_R\right)\phi_{2}\cdot F(\tau, \nu\tau)\\
	& +  \chi_{\nu\tau - R\lesssim\hbar_1}\left(\partial_{\tau}+ \frac{\lambda_{\tau}}{\lambda}R\partial_R + \partial_R\right)\phi_{2}\cdot\left[ F(\tau, R) - F(\tau, \nu\tau)\right]
	\end{align*}
	Here the second term on the right is at least source admissibly singular at level $l_1\ge 2$, and is handled as in the proof of Lemma~\ref{lem:N02abilin1}. 
	\\
	As for the first term on the right, in accordance with the conclusion of Lemma~\ref{lem:DerivonPrincSing} we may assume (i.e., this representation is valid on $R<\nu\tau$)
	\[
	\left(\partial_{\tau}+ \frac{\lambda_{\tau}}{\lambda}R\partial_R + \partial_R\right)\phi_{2}(\tau, R) = \int_0^\infty \phi_{n_{2}}(R,\xi)\cdot\xi^{\frac12}\cdot \xb(\tau,\xi)\rho_{n_{2}}(\xi)\,d\xi + \int_0^\infty  \phi_{n_{2}}(R,\xi)\cdot \yb(\tau,\xi)\rho_{n_{2}}(\xi)\,d\xi,
	\]
	where $\xb$ is principal singular in the sense of Definition~\ref{defi:xsingularsourcetermsngeq2adm}, while $\yb(\tau, \xi)$ is source admissibly singular  or in $S_1^{\hbar}$. It follows that 
	\[
	\chi_{\nu\tau - R\lesssim\hbar_1}\left(\partial_{\tau}+ \frac{\lambda_{\tau}}{\lambda}R\partial_R + \partial_R\right)\phi_{2}\cdot F(\tau, \nu\tau)
	\]
	agrees on the inner light cone $R<\nu\tau$ with the function 
	\[
	F(\tau, \nu\tau)\cdot\int_0^\infty \phi_{n_{2}}(R,\xi)\xi^{\frac12}\cdot \xb(\tau,\xi)\rho_{n_{2}}(\xi)\,d\xi 
	\]
	up to errors either of connection singular source type 
	or in $\tilde{S}_1^{\hbar}$. The size control of the term follows by invoking Lemma~\ref{lem:extraLinftyrefinedbound} with $R = \nu\tau$. 
	\\
	Next, if $\phi_2$ is as before but $F = F_2$ is source admissibly singular of level $l\geq 2$, we expand $F$ according to Lemma~\ref{lem:singFouriertiphysicalngeq2adm}. Referring to the notation there, the interaction of $\phi_2$ with $f_{1,3}$ is handled as before, while the product of $\left(\partial_{\tau}+ \frac{\lambda_{\tau}}{\lambda}R\partial_R + \partial_R\right)\phi_{2}$ and $f_2$ is handled by expanding $\left(\partial_{\tau}+ \frac{\lambda_{\tau}}{\lambda}R\partial_R + \partial_R\right)\phi_{2}$ according to Lemma~\ref{lem:singFouriertiphysicalngeq2adm} and multiplying out the terms. The situation with $\phi_2\in S_0^{\hbar}$ and $F$ admissibly source singular is more of the same. 
\end{proof}
Combining Lemma~\ref{lem:N02bbilin1}, Lemma~\ref{lem:N02bbilin2} easily implies the conclusion of the proposition for the first term on the right in \eqref{eq:N02brearranged}, while for the second term there, one easily checks that the microlocalization forces it to be of at connecting admissibly singular source type of level $l_1\geq 1$. This then completes the desired conclusion for $\mathcal{N}_{02b}$, and hence the proposition. 
\end{proof}
Returning to \eqref{eq:keytrilinnullformdecomp}, we have dealt with the second term on the right there, and now need a way to handle the first term, which we write as 
\begin{align*}
\sum_{\lambda_2,\lambda_3}\mathcal{N}_0\left(P_{<\min\{\lambda_2,\lambda_3\}}\phi_1,P_{\lambda_2}\phi_{2}, P_{\lambda_3}\phi_{3}\right) &=  \sum_{\lambda_2<\lambda_3}\mathcal{N}_0\left(P_{<\lambda_2}\phi_1,P_{\lambda_2}\phi_{2}, P_{\lambda_3}\phi_{3}\right)\\
& + \sum_{\lambda_2\geq\lambda_3}\mathcal{N}_0\left(P_{<\lambda_3}\phi_1,P_{\lambda_2}\phi_{2}, P_{\lambda_3}\phi_{3}\right)
\end{align*}
Introduce the wave type operator 
\[
\Box'_n: = \left(\partial_{\tau}+ \frac{\lambda_{\tau}}{\lambda}R\partial_R\right)^2 - \partial_R^2 + \frac{n^2}{R^2}. 
\]
Then we can write the preceding terms in the following manner, where we apply a cutoff $\chi_{R\gtrsim\tau}$ as we may in light of the results in subsection~\ref{subsec:allsourcetermsngeq2} as well as Lemma~\ref{lem:admsingawayfromshock}; in fact the latter shows that admissibly singular functions, when restricted far away from the light cone, are in $S_0^{\hbar}$ at fixed times, with fast temporal decay, and hence the results of subsection~\ref{subsec:allsourcetermsngeq2} also apply to source terms involving adfmissibly singular functions:
\begin{equation}\label{eq:nullformidentity1}\begin{split}
\chi_{R\gtrsim\tau}\sum_{\lambda_2<\lambda_3}\mathcal{N}_0\left(P_{<\lambda_2}\phi_1,P_{\lambda_2}\phi_{2}, P_{\lambda_3}\phi_{3}\right) &= \Box'_n\left(\chi_{R\gtrsim\tau}\left(\sum_{\lambda_2<\lambda_3}\phi_{1,<\lambda_2}\phi_{2,\lambda_2}\phi_{3,\lambda_3}\right)\right)\\
&+\left(\sum_{k=1}^3\frac{n_k^2}{R^2} - \frac{n^2}{R^2}\right)\chi_{R\gtrsim\tau}\left(\sum_{\lambda_2<\lambda_3}\phi_{1,<\lambda_2}\phi_{2,\lambda_2}\phi_{3,\lambda_3}\right)\\
& - \left[ \Box'_n,\chi_{R\gtrsim\tau}\right]\left(\sum_{\lambda_2<\lambda_3}\phi_{1,<\lambda_2}\phi_{2,\lambda_2}\phi_{3,\lambda_3}\right)\\
& - \chi_{R\gtrsim\tau}\sum_{\lambda_2<\lambda_3}\mathcal{N}_0\left(\phi_{3,\lambda_3}\phi_{1,<\lambda_2}\phi_{2,\lambda_2}\right)\\
& -  \chi_{R\gtrsim\tau}\sum_{\lambda_2<\lambda_3}\mathcal{N}_0\left(\phi_{2,\lambda_2}\phi_{1,<\lambda_2}\phi_{3,\lambda_3}\right)\\
& -  \chi_{R\gtrsim\tau}\sum_{\lambda_2<\lambda_3}\Box'_{n_1}\phi_{1,<\lambda_2}\phi_{2,\lambda_2}\phi_{3,\lambda_3}\\
& -  \chi_{R\gtrsim\tau}\sum_{\lambda_2<\lambda_3}\phi_{1,<\lambda_2}\Box'_{n_2}\phi_{2,\lambda_2}\phi_{3,\lambda_3}\\
& -  \chi_{R\gtrsim\tau}\sum_{\lambda_2<\lambda_3}\phi_{1,<\lambda_2}\phi_{2,\lambda_2}\Box'_{n_3}\phi_{3,\lambda_3}\\
& =:\sum_{j=1}^8 E_j
\end{split}\end{equation}
We shall then treat each of the terms $E_j$ separately. Some of these are good source terms and can be directly bounded, others can only be dealt with under the assumption that $\phi_j$ itself satisfy a wave equation (which is automatically fulfilled in an iterative scheme), one term will require further transformation ($E_8$), and one term will be used to modify the wave equation ($E_1$), thereby generating some additional but harmless source terms. 
\\

{\it{$E_1$}}: Recall that the wave equation at angular momentum $n$, $|n|\geq 2$, is given by 
\begin{equation}\label{eq:RegularFinestructure11}\begin{split}
&-\left(\left(\partial_{\tau} + \frac{\lambda_{\tau}}{\lambda}R\partial_R\right)^2 + \frac{\lambda_{\tau}}{\lambda}\left(\partial_{\tau} + \frac{\lambda_{\tau}}{\lambda}R\partial_R\right)\right)\varepsilon_{\pm}(n)+ H_n^{\pm}\varepsilon_{\pm}(n) = F_{\pm}(n),
\end{split}\end{equation}
where 
\[
H_n^{\pm} = \partial_R^2+\frac{1}{R}\partial_R - f_n(R)\pm g_n(R). 
\]
Call the wave operator on the left $-\Box_n$. Then we write the first term $E_1$ as 
\begin{align*}
E_1 &= \Box_n\left(\chi_{R\gtrsim\tau}\left(\sum_{\lambda_2<\lambda_3}\phi_{1,<\lambda_2}\phi_{2,\lambda_2}\phi_{3,\lambda_3}\right)\right) - \frac{\lambda_{\tau}}{\lambda}\left(\partial_{\tau} + \frac{\lambda_{\tau}}{\lambda}R\partial_R\right)\left(\chi_{R\gtrsim\tau}\left(\sum_{\lambda_2<\lambda_3}\phi_{1,<\lambda_2}\phi_{2,\lambda_2}\phi_{3,\lambda_3}\right)\right)\\
& + \left(\frac{1}{R}\partial_R - f_n(R)\pm g_n(R)+\frac{n^2}{R^2}\right)\left(\chi_{R\gtrsim\tau}\left(\sum_{\lambda_2<\lambda_3}\phi_{1,<\lambda_2}\phi_{2,\lambda_2}\phi_{3,\lambda_3}\right)\right)\\
&=: \Box_n\psi  + E_{11} + E_{12}. 
\end{align*}
Here the first term on the right $\Box_n\psi $ will be incorporated into the left-hand part of the wave equation \eqref{eq:RegularFinestructure11}, while the remaining two terms are good source terms: 
\begin{lemma}\label{lem:E1errorterms} Assume that $\phi_j$ are angular momentum $n_j$ functions, $|n_j|\geq 2$, $j = 1,2,3$, each of which admits a representation 
	\[
	\phi_j(\tau, R) = \int_0^\infty \phi_{n_{j}}(R,\xi)\cdot \xb_j(\tau, \xi)\cdot \rho_{n_{j}}(\xi)\,d\xi,\quad j = 1,2,3, 
	\]
	where the distorted Fourier transforms $\xb_j,\,j = 1,2$ each can be written as $\xb_j = \yb_j + \zb_j$ with $\yb_j$ admissibly singular (at angular momentum $n_j$) and $\zb_j\in S_0^{\hbar_j}, \mathcal{D}_{\tau}\zb_j\in S_1^{\hbar_j}$. 
	Then if $m$, $|m|\geq 2$, is an admissible angular momentum for the output of each of $E_{11}, E_{12}$, then each of $E_{11}, E_{12}$ satisfies 
	\[
	E_{1j}|_{R<\nu\tau} = g_j(\tau, R)|_{R<\nu\tau}, 
	\]
	where $g_j$ admits an angular momentum $m$ distorted Fourier representation 
	\[
	g_j(\tau, R) = \int_0^\infty \phi_{n}(R,\xi)\cdot r(\tau, \xi)\rho_{n}(\xi)\,d\xi,
	\]
	with $r = r_1 + r_2$ where $r_1$ is source admissibly singular, while $r_2\in S_1^{\hbar}$. Furthermore, we have the bound 
	\begin{align*}
	\tau^4\cdot\big\|r_2(\tau,\cdot)\big\|_{ S_1^{\hbar}} + \big\|r_1\big\|_{sourceadm}\lesssim \max_{j=1,2,3}\{|n_j|\}^{-C}\cdot \prod_{j=1}^3 |n_j|^C\Big(\big\|\overline{y}_j\big\|_{adm} + \tau^3\big(\big\|\overline{z}_j\big\|_{S_0^{\hbar_j}} + \big\|\mathcal{D}_{\tau}\overline{z}_j\big\|_{S_1^{\hbar_j}}\big)\Big)
        \end{align*}
	\end{lemma}
\begin{proof}
	Note that we avoid the key difficulty for the original null-form when the `good' derivative hits a $S_0^{\hbar}$-term while the `bad' derivative hits a singular term, because we have only one derivative to begin with here. The assertion for the term $E_{11}$ then follows by combining Lemma~\ref{lem:N02bbilin2} with Lemma~\ref{lem:longparaprod1}.  For the term $E_{12}$ we need in particular estimate the term (since $- f_n(R)\pm g_n(R)+\frac{n^2}{R^2} = O(\frac{n}{R^2})$) 
	\[
	\frac{n}{R^2}\chi_{R\gtrsim\tau}\cdot\left(\sum_{\lambda_2<\lambda_3}\phi_{1,<\lambda_2}\phi_{2,\lambda_2}\phi_{3,\lambda_3}\right),
	\]
	where $n$ refers to the angular momentum of the expression. Here we have to be careful not to lose in $n$, which we may assume in the most difficult case is the largest of all the angular momenta present, comparable to only one of the angular momenta $n_j$ in the $\phi_j$. In case all factors are in $\tilde{S}_0^{\hbar}$ this follows by the arguments used to prove Prop.~\ref{prop:bilin1}. In case that one of the factors is admissibly singular, and the expression is a singular source term, we can absorb the outer weight $n$ since the term is of level $l_1\geq 2$, the result following from Lemma~\ref{lem:longparaprod1}.
\end{proof}
{\it{$E_2$}}: Here we have to be more precise about the output angular momentum $n$, which we fix as $n = \sum_{j=1}^3n_j$, the value it takes in the iterative scheme. Then we have 
\[
\sum_{k=1}^3\frac{n_k^2}{R^2} - \frac{n^2}{R^2} = \frac{1}{R^2}\cdot O\left(\max\{|n_j|\}\cdot \max\{\{|n_j|\}\backslash\max\{|n_j|\}\}\right)
\]
and this shows that the contribution of $E_2$ can be handled analogously to the one of $E_1$. 
\\

{\it{$E_3$}}: This term is essentially 
\[
\frac{1}{R^2}\chi_{R\gtrsim\tau}\cdot\left(\sum_{\lambda_2<\lambda_3}\phi_{1,<\lambda_2}\phi_{2,\lambda_2}\phi_{3,\lambda_3}\right),
\]
and hence a simpler variant of the preceding terms. 
\\

{\it{$E_j,\,j\in \{6,7\}$}}: These terms can be handled if we recall the equation satisfied by $\phi_1,\phi_2$, these being linear combinations of the $\veps_{\pm}(n)$ (see \eqref{diag unknown}), which in turn satisfy \eqref{eq diag phys}. The right hand sides in these equations of course contain the same kinds of null-forms. However, the paraproduct structure inherent in the terms $E_{6, 7}$ makes the terms arising upon substituting the null-forms for the d'Alembertian factors better than the bare cubic null-form from before. To handle these terms, we (i) first have to move the operators $\Box_{n_j}'$ past the localizations $P_{<\lambda_2}, P_{\lambda_2}$, and then we (ii) have to invoke a paraproduct estimate to bound them. 
\\

{\it{(i): Commuting the wave operator at angular momentum $n,\,|n|\geq 2$, past localization operators.}} Here the following lemma mostly resolves the issue:
\begin{lemma}\label{lem:commutatorwithderivatives1} Let $P_{\mu}, P_{<\mu}$, $\lambda\geq 1$, be dyadic frequency localization operators acting on angular momentum $n,\,|n|\geq 2$ functions. Then if $f(\tau,R)$ is a function admitting the distorted Fourier representation 
	\[
	f(\tau, R) = \int_0^\infty \phi_{n}(R,\xi)\cdot \xb(\tau, \xi)\cdot \rho_{n}(\xi)\,d\xi, 
	\]
	with $\xb = \xb_1 + \xb_2$ where $\xb_1\in \Sh_{0}$, $\mathcal{D}_{\tau}\xb_1\in  S_1^{\hbar}$, and $\xb_2$ admissibly singular of level $l_1\geq 1$, then (recalling the notation from above)
	\[
	\mathcal{F}\left(\left[\Box_n,P_\mu\right]f\right) = \yb_1 + \yb_2
	\]
	where $\yb_1\in S_1^{\hbar}$, and $\yb_2$ is source admissibly singular. We have the bound
	\begin{align*}
	\tau^3\cdot\big\|\yb_1\big\|_{ S_1^{\hbar}} + \big\|\yb_2\big\|_{sourceadm}\lesssim \sup_{\tau\geq \tau_0}\tau^3\cdot\big(\big\|\xb_1(\tau,\cdot)\big\|_{\Sh_{0}} + \big\|\mathcal{D}_{\tau}\xb_1(\tau,\cdot)\big\|_{S_1^{\hbar}}\big) + \big\|\xb_2\big\|_{adm}.
	\end{align*}
	The same conclusion applies if $P_{\mu}$ is replaced by $P_{<\mu}$. 
\end{lemma}
\begin{proof}
	Since $H_n^{\pm}$ commutes with $P_{\mu}$, it suffices to consider the temporal part of $\Box_n$, which consists of the operators 
	\[
	\left(\partial_{\tau} + \frac{\lambda_{\tau}}{\lambda}R\partial_R\right)^2,\quad \frac{\lambda_{\tau}}{\lambda}\left(\partial_{\tau} + \frac{\lambda_{\tau}}{\lambda}R\partial_R\right).
	\]
	Write 
	\begin{align*}
		\left[\left(\partial_{\tau} + \frac{\lambda_{\tau}}{\lambda}R\partial_R\right)^2, P_{\mu}\right] &= \left(\partial_{\tau} + \frac{\lambda_{\tau}}{\lambda}R\partial_R\right)\left[\left(\partial_{\tau} + \frac{\lambda_{\tau}}{\lambda}R\partial_R\right), P_{\mu}\right] \\
		& +\left[\left(\partial_{\tau} + \frac{\lambda_{\tau}}{\lambda}R\partial_R\right), P_{\mu}\right] \left(\partial_{\tau} + \frac{\lambda_{\tau}}{\lambda}R\partial_R\right)\\
		& = \left(\partial_{\tau} + \frac{\lambda_{\tau}}{\lambda}R\partial_R\right)\left[\frac{\lambda_{\tau}}{\lambda}R\partial_R,P_{\mu}\right] +  \left[\frac{\lambda_{\tau}}{\lambda}R\partial_R,P_{\mu}\right]\left(\partial_{\tau} + \frac{\lambda_{\tau}}{\lambda}R\partial_R\right)\\
		& = \left[\frac{\lambda_{\tau}}{\lambda}R\partial_R,\left[\frac{\lambda_{\tau}}{\lambda}R\partial_R,P_{\mu}\right] \right] + 2\left[\frac{\lambda_{\tau}}{\lambda}R\partial_R,P_{\mu}\right]\left(\partial_{\tau} + \frac{\lambda_{\tau}}{\lambda}R\partial_R\right)\\
		&+\left(\frac{\lambda_{\tau}}{\lambda}\right)^{\prime}\left[R\partial_{R},P_{\mu}\right]
	\end{align*}
In order to proceed, we translate the commutator to the Fourier side:
\begin{align*}
	\mathcal{F}\circ \left[\frac{\lambda_{\tau}}{\lambda}R\partial_R,P_{\mu}\right] = \frac{\lambda_{\tau}}{\lambda}\left[-2\xi\partial_\xi +\mathcal{K}^{(0)}_{\hbar},\,\chi_{\xi\simeq\mu}\right]\circ \mathcal{F}
	= -2\frac{\lambda_{\tau}}{\lambda}\xi\partial_\xi(\chi_{\xi\simeq\mu})\mathcal{F} + \frac{\lambda_{\tau}}{\lambda}\left[\mathcal{K}^{(0)}_{\hbar},\,\chi_{\xi\simeq\mu}\right]\circ \mathcal{F}.
\end{align*}
Here the operator 
\[
\left[\mathcal{K}^{(0)}_{\hbar},\,\chi_{\xi\simeq\mu}\right]
\]
is seen to have similar properties to $\mathcal{K}^{(0)}_{\hbar}$. In fact, the kernel of this commutator is given by
\begin{align*}
	\frac{\left(\chi_{\eta\simeq\mu}-\chi_{\xi\simeq\mu}\right)F(\xi,\eta)\rho_{n}(\eta)}{\xi-\eta},
\end{align*}
which behaves at least as good as the kernel of $\calK_{\hbar}^{(0)}$. If we denote by 
\begin{align*}
	\calF\circ\left[\frac{\lambda_{\tau}}{\lambda}R\partial_{R},P_{\mu}\right]:=\calK_{C},
\end{align*}
then 
\begin{align*}
	\calF\circ\left[\frac{\lambda_{\tau}}{\lambda}R\partial_{R},\left[\frac{\lambda_{\tau}}{\lambda}R\partial_{R},P_{\mu}\right]\right]=\frac{\lambda_{\tau}}{\lambda}\left[-2\xi\partial_{\xi}+\calK^{(0)}_{\hbar},\calK_{C}\right]\circ\calF,
\end{align*}
which can be handled similarly as the commutator $\left[\calD_{\tau},\calK_{\hbar}^{(0)}\right]$ (see Proposition \ref{prop: xh good para linear}). 
It easily follows that if $\xb_1\in S_0^{\hbar}$ (or $\xb_{1}\in \Sh_{1}$) , then 
\[
\mathcal{F}\circ\left[\left(\partial_{\tau} + \frac{\lambda_{\tau}}{\lambda}R\partial_R\right)^2, P_{\mu}\right]\circ\mathcal{F}^{-1}\xb_1\in S_1^{\hbar}. 
\]
The conclusion for 
\[
\mathcal{F}\circ\left[\left(\partial_{\tau} + \frac{\lambda_{\tau}}{\lambda}R\partial_R\right)^2, P_{\mu}\right]\circ\mathcal{F}^{-1}\xb_2 
\]
with $\xb_2$ admissibly singular of level $l_1\geq 1$ follows by combining the observation in subsection~\ref{subsubsec:Fourierlocal} and Proposition~\ref{prop:transferenceonsingularngeq2}. 
\end{proof}

{\it{(ii): Using the paraproduct structure.}} The preceding lemma is not strong enough to handle the case when $\xb_2$ is of principal incoming singular type. This is where we have to take advantage of the paraproduct nature of the terms $E_6, E_7$: 
\begin{lemma}\label{lem:commutatorwithderivatives2} Let $f(R)$ be an angular momentum $n_1$ function as in the preceding lemma. 
	Let $g$, an angular momentum $n_2$ function with $|n_2|\geq 2$, have the same properties as $f$, and assume that $\{m,n_1,n_2\}$ with $|m|\geq 2$ is an admissible angular momentum triple.  
	Then 
	\[
	\sum_{\lambda>0}\left[\Box_n,P_\lambda\right]f P_{\geq \lambda}g|_{R<\nu\tau} = h(\tau, R)|_{R<\nu\tau}, 
	\]
	where the distorted Fourier transform of $h(\tau, R)$ can be decomposed as the sum of a function $\zb_1\in S_1^{\hbar}$ and a function $\zb_2$ which is source admissibly singular and of level $l_1\geq 1$. Moreover, assuming that the distorted Fourier transforms of $f, g$ are given by $\overline{x}_j,\,j = 1, 2$, such that we can write $\overline{x}_j = \overline{y}_j + \overline{z}_j$, with $\overline{y}_j$ admissibly singular and $\overline{z}_j\in S_0^{\hbar},\,\mathcal{D}_{\tau}\overline{z}_j\in S_1^{\hbar}$, then we have the bound
	\begin{align*}
	\tau^4\cdot \big\|\zb_1\big\|_{S_1^{\hbar}} + \big\|\zb_2\big\|_{sourceadm}\lesssim \min\{|n_{1,2}|\}^C\cdot \prod_{j=1,2}\Big(\big\|\overline{y}_j\big\|_{adm} + \sup_{\tau\geq\tau_0}\tau^3\cdot \big(\big\|\overline{z}_j\big\|_{S_0^{\hbar}} + \big\|\mathcal{D}_{\tau}\overline{z}_j\big\|_{S_1^{\hbar}}\big)\Big).
	\end{align*}
\end{lemma}
\begin{proof} We distinguish between different cases for $g$, whose distorted Fourier transform is given by an admissibly singular function term $\overline{y}_2$ and a function $\overline{z}_2\in S_0^{\hbar_2}$ . We also assume that the distorted Fourier transform of $f$ is the sum of $\overline{y}_1,\overline{z}_1$ with analogous properties. Also assume that $|n_2|\geq |n_1|$, where $n_1$ is the angular momentum of $f$ and $n_2$ the one of $g$, the case $|n_1|>|n_2|$ being similar. Assume that $\{m, n_1,n_2\}$, $|m|\geq 2$, is an admissibly triple of angular momenta and that $\hbar = \frac{1}{|m|}$. 
\\

{\it{(i): $g$ admissibly singular}}. Decompose $g$ in accordance with Lemma~\ref{lem:singFouriertiphysicalngeq2adm}. Calling this decomposition 
\[
P_{\geq \lambda}g = g_1 + g_2 + g_3,
\]
where on account of the frequency localization we have that $\lambda^{\frac12}\cdot g_j\in S_1^{\hbar_2}$, $j = 1, 3$, $\hbar_2  = \frac{1}{|n_2|}$. According to the proof of Lemma~\ref{lem:commutatorwithderivatives1}, the commutator expression 
\[
\lambda^{-\frac12}\cdot \big[\Box_{n_1}, P_{\lambda}\big]f
\]
is the sum of a function in $S_0^{\hbar_1}$ and an admissibly singular source term of level $l\geq 2$, which can be decomposed in analogy to  Lemma~\ref{lem:commutatorwithderivatives1}. Using a third order Taylor expansion for $\lambda^{\frac12}\cdot (g_1 + g_3)$ around $R = \nu\tau$, multiplying out the terms in the product and taking advantage of Lemma~\ref{lem:singPhysicaltoFourierngeq2}, we infer that there is a function $g$ with 
\begin{align*}
g|_{R<\nu\tau} = \lambda^{-\frac12}\cdot \big[\Box_{n_1}, P_{\lambda}\big]f\cdot \lambda^{\frac12}\cdot (g_1 + g_3)|_{R<\nu\tau},
\end{align*}
and such that 
\[
\langle \phi_m(R;\xi),\,g\rangle_{L^2_{R\,dR}}
\]
 is the sum of a source admissibly singular term $\overline{y}$ of level $l\geq 1$ as well as a term in $\overline{z}\in S_1^{\hbar}$ satisfying 
 \begin{align*}
 \big\|\overline{y}\big\|_{sourceadm} + \tau^3\cdot \big\|\overline{z}\big\|_{S_1^{\hbar}}\lesssim |n_1|^C\cdot \big(\big\|\overline{y}_1\big\|_{adm} + \tau^3\big\|\overline{z}_1\big\|_{S_0^{\hbar_1}}\big)\cdot \big\|\overline{y}_2\big\|_{adm}. 
 \end{align*}
 Indeed, letting $\mathbb{P}_3\big(\lambda^{\frac12}\cdot (g_1 + g_3)\big)$ the third order Taylor polynomial of $\lambda^{\frac12}\cdot (g_1 + g_3)$, the product 
 \[
 \lambda^{-\frac12}\cdot \big[\Box_{n_1}, P_{\lambda}\big]f\cdot \mathbb{P}_3\big(\lambda^{\frac12}\cdot (g_1 + g_3)\big)
 \]
 when restricted to $R<\nu\tau$, has a singular part with source admissibly singular distorted Fourier transform at level $l\geq 2$. If $\mathbb{E}_3\big(\lambda^{\frac12}\cdot (g_1 + g_3)\big)$ is the error of the third order Taylor expansion, then 
 \begin{align*}
  \lambda^{-\frac12}\cdot \big[\Box_{n_1}, P_{\lambda}\big]f\cdot \mathbb{E}_3\big(\lambda^{\frac12}\cdot (g_1 + g_3)\big)\in S_1^{\hbar}. 
 \end{align*}
 \\
 
 {\it{(ii) $g$ in $\tilde{S}_0^{\hbar_2}$.}} One argues as for the contribution of $g_1 + g_3$ in the previous case. This time we have the bound 
 \begin{align*}
 \big\|\lambda^{\frac12}\cdot P_{\geq \lambda}g\big\|_{\tilde{S}_1^{\hbar_2}}\lesssim \big\|\overline{z}_2\big\|_{S_0^{\hbar_2}}, 
 \end{align*}
 and writing 
 \begin{align*}
 \mathbb{P}_3\big(\lambda^{\frac12}\cdot P_{\geq \lambda}g\big) = \sum_{j=0}^3a_j(\tau)\cdot (R-\nu\tau)^j,
 \end{align*}
we have the bounds (Lemma~\ref{lem:extraLinftyrefinedbound})
\begin{align*}
\big|a_j(\tau)\big|\lesssim |n_2|^{\frac32+\delta+j}\cdot \big\|\overline{z}_2\big\|_{S_0^{\hbar_2}}. 
\end{align*}
The product of the singular part of $\lambda^{-\frac12}\cdot \big[\Box_{n_1}, P_{\lambda}\big]f$ and $ \mathbb{P}_3\big(\lambda^{\frac12}\cdot P_{\geq \lambda}g\big) $ can then be bounded as a source admissible term as asserted in the Lemma since it is of level $l\geq 2$ and we can absorb the loss in terms of $|n_2|$. 
\\
For the Taylor error $ \mathbb{E}_3\big(\lambda^{\frac12}\cdot P_{\geq \lambda}g\big) $, we can write 
\begin{align*}
\partial_R^k\mathbb{E}_3\big(\lambda^{\frac12}\cdot P_{\geq \lambda}g\big) = \big(R-\nu\tau\big)^{\frac72-k}\cdot \zeta_k(\tau, R), \,k = 0,\ldots, 3,\,\partial_R^4\mathbb{E}_3\big(\lambda^{\frac12}\cdot P_{\geq \lambda}g\big) =  \lambda^{\frac12}\cdot \zeta_4(\tau, R),
\end{align*}
 where we have the bound 
 \begin{align*}
 \big|\zeta_k\big|\lesssim |n_2|^{5+\delta}\cdot  \big\|\overline{z}_2\big\|_{S_0^{\hbar_2}},\,k = 0,\ldots, 4. 
 \end{align*}
 The estimate 
 \begin{align*}
 \tau^3\big\|\lambda^{-\frac12}\cdot \big[\Box_{n_1}, P_{\lambda}\big]f\cdot \mathbb{E}_3\big(\lambda^{\frac12}\cdot P_{\geq \lambda}g\big)\big\|_{\tilde{S}_1^{\hbar}}\lesssim |n_1|^C\cdot  \big(\big\|\overline{y}_1\big\|_{adm} + \tau^3\big\|\overline{z}_1\big\|_{S_0^{\hbar_1}}\big)\cdot \tau^3\big\|\overline{z}_2\big\|_{S_0^{\hbar_2}}
 \end{align*}
 then follows by applying the Leibniz rule to distribute the derivatives in the definition of $\tilde{S}_1^{\hbar}$, use the preceding bounds and observe that the $\hbar$ weights for the 'output norm' counteract the loss of $|n_2|^{5+\delta}$ provided $|n_2|\gg |n_1|$. 
\end{proof}

Taking advantage of the preceding lemma, as well as the fact that the error terms arising upon replacing $\Box_{n_j}'$ by $\Box_{n_j}$ are again see to be admissible source terms of at least level $l_1\geq 1$, we can now replace the terms $E_6, E_7$ by the following ones, respectively: 
\begin{align*}
	\chi_{R\gtrsim\tau}\sum_{\lambda_2<\lambda_3}P_{<\lambda_2}\left(\Box_{n_1}\phi_{1}\right)\phi_{2,\lambda_2}\phi_{3,\lambda_3},\quad \chi_{R\gtrsim\tau}\sum_{\lambda_2<\lambda_3}\phi_{1,<\lambda_2}P_{\lambda_2}\left(\Box_{n_2}\phi_{2}\right)\phi_{3,\lambda_3}
\end{align*}
Then we re-iterate application of the equations for $\phi_2,\phi_1$. The fact that these terms don't pose difficulties in spite of the fact that the same null-forms occur in $\Box_{n_1}\phi_{1}, \Box_{n_2}\phi_{2}$ are now consequences of the following 
\begin{lemma}\label{lem:E67paraprodterms}  Let $F$ be an angular momentum $n_1$ function, $|n_1|\geq 2$, whose distorted Fourier transform is the sum of a term in $\xi^{\frac14}S_1^{\hbar_1}$ and an admissibly singular source term. Also, assume that $\phi_2$ is an angular momentum $n_2, |n_2|\geq 2$ function, whose distorted Fourier transform is the sum of a function in $S_0^{\hbar_2}$ and an admissibly singular term. Finally, assume that $\{m,n_1,n_2\}$ is an admissible triple of angular momenta. Then 
	\[
	\chi_{R\gtrsim \tau}\sum_{0<\lambda}P_{<\lambda}F P_{\lambda}\phi_2|_{R<\nu\tau} = g(\tau, R)|_{R<\nu\tau}
	\]
	where $g$ is an angular momentum $m$ function whose distorted Fourier transform is the sum of a term in $S_1^{\hbar}$ and an admissibly singular source term of level $l_1\geq 1$. 
\end{lemma}
\begin{proof} We distinguish between the different possibilities for the factors. Write $\hbar_j = |n_j|^{-1},\,j = 1, 2$, where $n_1, n_2$ are the angular momenta of $F, \phi$, respectively. Also, assume that $\{m, n_1, n_2\}$ is an admissible triple, where we interpret the para-product as an angular momentum $m$ function. 
\\

{\it{(i): $\mathcal{F}_{n_1}(F)\in \langle\xi\rangle^{\frac14}S_1^{\hbar_1}$}}. If $\phi_2$ is admissibly singular (in the sense that its distorted Fourier transform is admissibly singular), then $\lambda^{\frac14}\cdot P_{\lambda}\phi_2$ is source admissibly singular of level $l\geq 1$. On the other hand, 
\[
\lambda^{-\frac14}P_{<\lambda}F\in \tilde{S}_1^{\hbar_1}.
\]
Using a third order Taylor approximation for $\lambda^{-\frac14}P_{<\lambda}F$ around $R = \nu\tau$, and denoting by $\mathbb{P}_3\big(\lambda^{-\frac14}P_{<\lambda}F\big)$ the third order Taylor polynomial while $\mathbb{E}_3\big(\lambda^{-\frac14}P_{<\lambda}F\big)$ is the error term, we have upon using Lemma~\ref{lem:singFouriertiphysicalngeq2adm} for $\lambda^{\frac14}\cdot P_{\lambda}\phi_2$ that the product 
\begin{align*}
\chi_{1\lesssim R}\cdot \mathbb{P}_3\big(\lambda^{-\frac14}P_{<\lambda}F\big)\cdot \lambda^{\frac14}\cdot P_{\lambda}\phi_2|_{R<\nu\tau}
\end{align*}
 agrees with a function $g|_{R<\nu\tau}$ whose distorted Fourier transform (at angular momentum $m$) is the sum of a source admissibly function $\overline{y}$ at level $l\geq 1$ as well as a function in $\overline{z}\in S_1^{\hbar}$ (with $\hbar = |m|^{-1}$), and we have the estimates (taking advantage of Lemma~\ref{lem:singPhysicaltoFourierngeq2} 
 \begin{align*}
 &\big\|\overline{y}\big\|_{sourceadm}\lesssim \lambda^{-(0+)}\cdot \min\{|n_1|, |n_2|\}^C\cdot \tau^3\cdot\big\|\langle\xi\rangle^{-\frac14}\mathcal{F}_{n_1}(F)\big\|_{S_1^{\hbar_1}}\cdot \big\|\mathcal{F}_{n_2}\big(\phi_2\big)\big\|_{adm}\\
 &\tau^4\cdot \big\|\overline{z}\big\|_{S_1^{\hbar}}\lesssim  \lambda^{-(0+)}\cdot\min\{|n_1|, |n_2|\}^C\cdot \tau^3\cdot\big\|\langle\xi\rangle^{-\frac14}\mathcal{F}_{n_1}(F)\big\|_{S_1^{\hbar_1}}\cdot \big\|\mathcal{F}_{n_2}\big(\phi_2\big)\big\|_{adm}.\\
 \end{align*}
 The contribution of $\chi_{1\lesssim R}\cdot \mathbb{E}_3\big(\lambda^{-\frac14}P_{<\lambda}F\big)$ satisfies 
 \begin{align*}
\tau^4\cdot  \big\|\chi_{1\lesssim R}\cdot \mathbb{E}_3\big(\lambda^{-\frac14}P_{<\lambda}F\big)\cdot \lambda^{\frac14}\cdot P_{\lambda}\phi_2|_{R<\nu\tau}\big\|_{\tilde{S}_1^{\hbar}}\lesssim  \lambda^{-(0+)}\cdot\min\{|n_1|, |n_2|\}^C\cdot \tau^3\cdot\big\|\langle\xi\rangle^{-\frac14}\mathcal{F}_{n_1}(F)\big\|_{S_1^{\hbar_1}}\cdot \big\|\mathcal{F}_{n_2}\big(\phi_2\big)\big\|_{adm}.\\
 \end{align*}
 {\it{(ii): $\mathcal{F}_{n_1}(F)$ source admissibly singular.}} Here we use that $\lambda^{-\frac14}P_{<\lambda}F$ is source admissibly singular of level $l\geq 1$, and distinguish between the case when $\phi_2$ is admissibly singular and when it is in $S_0^{\hbar_2}$. The main observation is that 
the product of two admissibly singular source terms of level $l\geq 1$ is again such a term, and localizing one of them to frequency $\lambda$ results in a small gain $\lambda^{-(0+)}$ for the corresponding norm. The details are similar to the preceding case.  
\end{proof}

To complete the estimates for the terms $E_{6,7}$, it suffices to combine the preceding lemma with Lemma~\ref{lem:sourcetermsotherthannullformprop738prep} for the non-null form terms in $\Box_{n_j}\phi_j$, while the next two lemmas deal with the contribution of the null-form terms, {\it{all under the hypothesis that only angular momentum $|n|\geq 2$ functions appear throughout}}. The general case can and will be discussed later, and follows the exact same pattern. 
\begin{lemma}\label{lem:badnullformpara} If $\phi_1, \phi_2, \phi_3$ are as in Prop.~\ref{prop:goodcubicnullform1}, as is $\phi_4$, then
\begin{align*}
\sum_{\lambda>1}P_{<\lambda}\big(\mathcal{N}_0\big(\phi_1, \phi_2, \phi_3\big)\big)\cdot P_{\lambda}\phi_4
\end{align*}
satisfies the conclusion of Prop.~\ref{prop:goodcubicnullform1}. In particular, we can write this expression (at angular momentum $m = \sum_{j=1}^4 n_j$), when restricted to $R<\nu\tau$, as 
\begin{align*}
\int_0^\infty \phi_m(R;\xi)\cdot \overline{z}(\tau,\xi)\rho_m(\xi)\,d\xi,
\end{align*}
where the distorted Fourier transform $\overline{z} = \overline{z}_1 + \overline{z}_2$ with $\overline{z}_1$ admissibly singular and $\overline{z}_2\in S_1^{\hbar},\,\hbar = \frac{1}{|m|+1}$ and we have the bound 
\begin{align*}
\big\|\overline{z}_1\big\|_{adm} + \tau^4\big\|\overline{z}_2\big\|_{S_1^{\hbar}}\lesssim \max_{j=1,\ldots, 4}\{|n_j|\}^{-C}\cdot \prod_{j=1}^4 |n_j|^C\big(\big\|\overline{y}_j\big\|_{adm} + \tau^3\big\|\overline{z}_j\big\|_{S_0^{\hbar_j}}\big)
\end{align*}

\end{lemma}
\begin{proof} The key point is that the only obstruction to placing $\mathcal{N}_0\big(\phi_1, \phi_2, \phi_3\big)$ directly into the admissible source space is the case where a principal incoming singular term is hit by $\partial_{\tau} + \frac{\lambda_{\tau}}{\lambda}R\partial_R +\partial_R$ while 
an unstructured term with Fourier transform in $S_0^{\hbar}$ is hit by $\partial_{\tau} + \frac{\lambda_{\tau}}{\lambda}R\partial_R - \partial_R$. For this recall the discussion at the beginning of subsection~\ref{subset:hardnullformestimates}. Nonetheless, by straightforward modification of the arguments for Prop.~\ref{prop:goodcubicnullform1}, we see that in this particular bad case we can still place $\mathcal{N}_0\big(\phi_1, \phi_2, \phi_3\big)$ into $\langle \xi^{\frac14}\rangle\cdot S_1^{\tilde{\hbar}}$ (with $\tilde{\hbar} = (|\sum_{j=1}^3 n_j| + 1)^{-1}$), and we have the bound 
\begin{align*}
\tau^3\big\|\mathcal{N}_0\big(\phi_1, \phi_2, \phi_3\big)\big\|_{\langle \xi^{\frac14}\rangle\cdot S_1^{\tilde{\hbar}}}\lesssim \max_{j=1,\ldots, 3}\{|n_j|\}^{-C}\cdot \prod_{j=1}^3 |n_j|^C\big(\big\|\overline{y}_j\big\|_{adm} + \tau^3\big\|\overline{z}_j\big\|_{S_0^{\hbar_j}}\big).
\end{align*}
The desired bound for this bad contribution is then a consequence of the preceding lemma in case that $\phi_4$ has distorted Fourier transform in $S_0^{\hbar_4}$, by writing 
\begin{align*}
\sum_{\lambda>1}P_{<\lambda}\big(\mathcal{N}_0\big(\phi_1, \phi_2, \phi_3\big)\big)\cdot P_{\lambda}\phi_4 = \sum_{\lambda>1}\lambda^{-\frac12}P_{<\lambda}\big(\mathcal{N}_0\big(\phi_1, \phi_2, \phi_3\big)\big)\cdot \lambda P_{\lambda}\phi_4
\end{align*}
If $\phi_4$  has admissibly singular distorted Fourier transform, we can invoke Lemma~\ref{lem:singFouriertiphysicalngeq2adm}; we still claim that if we replace\footnote{We may assume that $\mathcal{N}_0(\phi_1, \phi_2, \phi_3)$ is the first term on the right in \eqref{eq:TheNullForm}.} $(\partial_{\tau} + \frac{\lambda_{\tau}}{\lambda}R\partial_R - \partial_R)\phi_2$, with $\phi_2$ assumed to have distorted Fourier transform in $S_0^{\hbar_2}$, by the Cauchy error term in the third order Taylor approximation around $R = \nu\tau$, while $\phi_3$ is principal incoming singular (which is precisely the problematic case), then we have that 
\[
\langle \phi_m(R;\xi),\,\sum_{\lambda>1}P_{<\lambda}\big(\mathcal{N}_0\big(\phi_1, \phi_2, \phi_3\big)\big)\cdot P_{\lambda}\phi_4\rangle_{L^2_{R\,dR}}\in S_1^{\hbar}. 
\]
In fact, we have 
\begin{align*}
\big\|\sum_{\lambda>1}\partial_R^4 P_{<\lambda}\big(\mathcal{N}_0\big(\phi_1, \phi_2, \phi_3\big)\big)\cdot P_{\lambda}\phi_4\big\|_{H^{0+}}\lesssim \sum_{\lambda>1}\big\|\partial_R^3P_{<\lambda}\big(\mathcal{N}_0\big(\phi_1, \phi_2, \phi_3\big)\big)\big\|_{L^\infty\cap H^{0+}}\cdot \big\|\lambda P_{\lambda}\phi_4\big\|_{H^{0+}},
\end{align*}
which under our current assumptions can be bounded in terms of $ \prod_{j=1}^3 |n_j|^C\big(\big\|\overline{y}_j\big\|_{adm} + \tau^3\big\|\overline{z}_j\big\|_{S_0^{\hbar_j}}\big)$, as is easily verified(taking advantage of Lemma~\ref{lem:singFouriertiphysicalngeq2adm}) if $\phi_1$ is admissibly singular). 
On the other hand, for $k\in \{1,2,3,4\}$, we write 
\begin{equation}\label{eq:mixedderivativespara}\begin{split}
\sum_{\lambda>1}\partial_R^{4-k} P_{<\lambda}\big(\mathcal{N}_0\big(\phi_1, \phi_2, \phi_3\big)\big)\cdot \partial_R^kP_{\lambda}\phi_4 &= \partial_R^{4-k} \big(\mathcal{N}_0\big(\phi_1, \phi_2, \phi_3\big)\big)\cdot \partial_R^k \phi_4\\& - \sum_{\lambda>1}\partial_R^{4-k} P_{\geq\lambda}\big(\mathcal{N}_0\big(\phi_1, \phi_2, \phi_3\big)\big)\cdot \partial_R^kP_{\lambda}\phi_4.
\end{split}\end{equation}
Writing $\phi_4$ like the middle term $f_2$ in Lemma~\ref{lem:singFouriertiphysicalngeq2adm}, it is easily verified that the first term on the right is in fact in $H^{0+}$. For the second term on the right we can use that 
\begin{align*}
\mu^{k-1+\delta}\cdot \tau^3\big\| \partial_R^{4-k} P_{\mu}\big(\mathcal{N}_0\big(\phi_1, \phi_2, \phi_3\big)\big\|_{L^\infty},\,\delta>0,\,4\geq k\geq 1,
\end{align*}
is bounded in terms of $\prod_{j=1}^3 |n_j|^C\big(\big\|\overline{y}_j\big\|_{adm} + \tau^3\big\|\overline{z}_j\big\|_{S_0^{\hbar_j}}\big)$, while we also have 
\begin{align*}
\tau\big\|\partial_R^kP_{\lambda}\phi_4\big\|_{L^2_{R\,dR}}\lesssim \lambda^{k-1}\cdot \big\|\phi_4\big\|_{adm}.
\end{align*}
This shows that \eqref{eq:mixedderivativespara} can also be placed into $S_1^{\hbar}$. 
\end{proof}

The following lemma can be proved in similar fashion:
\begin{lemma}\label{lem:badnullformpara1} If $\phi_1, \phi_2, \phi_3$ are as in Prop.~\ref{prop:goodcubicnullform1}, as is $\phi_j, j = 4, 5$, then
\begin{align*}
\sum_{\lambda>1}P_{<\lambda}\big(\mathcal{N}_0\big(\phi_1, \phi_2, \phi_3\big)\big)\cdot P_{<\lambda}\phi_4\cdot P_{\lambda}\phi_5
\end{align*}
satisfies the conclusion of Prop.~\ref{prop:goodcubicnullform1}. In particular, we can write this expression (at angular momentum $m = \sum_{j=1}^5 n_j$), when restricted to $R<\nu\tau$, as 
\begin{align*}
\int_0^\infty \phi_m(R;\xi)\cdot \overline{z}(\tau,\xi)\rho_m(\xi)\,d\xi,
\end{align*}
where the distorted Fourier transform $\overline{z} = \overline{z}_1 + \overline{z}_2$ with $\overline{z}_1$ admissibly singular and $\overline{z}_2\in S_1^{\hbar},\,\hbar = \frac{1}{|m|+1}$ and we have the bound 
\begin{align*}
\big\|\overline{z}_1\big\|_{adm} + \tau^4\big\|\overline{z}_2\big\|_{S_1^{\hbar}}\lesssim \max_{j=1,\ldots, 4}\{|n_j|\}^{-C}\cdot \prod_{j=1}^5 |n_j|^C\big(\big\|\overline{y}_j\big\|_{adm} + \tau^3\big\|\overline{z}_j\big\|_{S_0^{\hbar_j}}\big)
\end{align*}
The same bound applies to the expression
\begin{align*}
\sum_{1<\lambda_1<\lambda_2}P_{\lambda_1}\big(\mathcal{N}_0\big(\phi_1, \phi_2, \phi_3\big)\big)\cdot P_{<\lambda_1}\phi_4\cdot P_{\lambda_2}\phi_5
\end{align*}
\end{lemma}
\begin{proof} We sketch how to deal with the last expression. We may assume that $\phi_2$ in $(\partial_{\tau} + \frac{\lambda_{\tau}}{\lambda}R\partial_R - \partial_R)\phi_2$ is in $S_0^{\hbar_2}$ and that the latter expression is replaced by its Cauchy error in the third order Taylor development around $R = \nu\tau$, while $\phi_3$ in $(\partial_{\tau} + \frac{\lambda_{\tau}}{\lambda}R\partial_R +\partial_R)\phi_3$ is of principal incoming singular type. Furthermore, we may assume $\phi_1$ is in $S_0^{\hbar_1}$. It is then easy to see that we may replace the term 
\[
P_{\lambda_1}\big(\mathcal{N}_0\big(\phi_1, \phi_2, \phi_3\big)\big)
\]
by the Cauchy error $\mathbb{E}_3\big(P_{\lambda_1}\big(\mathcal{N}_0\big(\phi_1, \phi_2, \phi_3\big)\big)\big)$ of its third order Taylor development around $R = \nu\tau$, by applying Lemma~\ref{lem:singFouriertiphysicalngeq2adm} to the terms $ P_{<\lambda_1}\phi_4, P_{\lambda_2}\phi_5$ if they are of admissibly singular type. This error being of size $O(|R-\nu\tau|^{3+})$ and having bounded derivatives up to third order, one can again use the representation Lemma~\ref{lem:singFouriertiphysicalngeq2adm} to the last two terms is needed to conclude that 
\begin{align*}
\sum_{1<\lambda_1<\lambda_2}\mathbb{E}_3\big(P_{\lambda_1}\big(\mathcal{N}_0\big(\phi_1, \phi_2, \phi_3\big)\big)\big)\cdot P_{<\lambda_1}\phi_4\cdot P_{\lambda_2}\phi_5\in S_1^{\hbar}.
\end{align*}
\end{proof}

{\it{$E_8$}}. Here we finally arrive at a term (up to easier error terms) where the same procedure is re-iterated, but with additional factors which gain smallness. Again we may replace $\Box'_{n_3}$ by $\Box_{n_3}$ up to good source error terms, and we may commute the localizer $P_{\lambda_3}$ and $\Box_{n_3}$, leading to the term 
\[
\chi_{R\gtrsim\tau}\sum_{\lambda_2<\lambda_3}\phi_{1,<\lambda_2}\phi_{2,\lambda_2}P_{\lambda_3}\left(\Box_{n_3}\phi_{3}\right)
\]
Writing 
\begin{align*}
	\chi_{R\gtrsim\tau}\sum_{\lambda_2<\lambda_3}\phi_{1,<\lambda_2}\phi_{2,\lambda_2}P_{\lambda_3}\left(\Box_{n_3}\phi_{3}\right) =&  -\chi_{R\gtrsim\tau}\sum_{\lambda_2\geq \lambda_3}\phi_{1,<\lambda_2}\phi_{2,\lambda_2}P_{\lambda_3}\left(\Box_{n_3}\phi_{3}\right)\\
	& +  \chi_{R\gtrsim\tau}\sum_{\lambda_2}\phi_{1,<\lambda_2}\phi_{2,\lambda_2}\left(\Box_{n_3}\phi_{3}\right),
\end{align*}
and taking advantage of Lemma~\ref{lem:E67paraprodterms} and the agument for the terms $E_{6,7}$, we can discard the first term on the right. Consider then the second term on the right, where we may assume that the term $\Box_{n_3}\phi_{3}$ is again given by a null-form of the type under consideration, the remaining terms constituting $\Box_{n_3}\phi_{3}$ leading to admissible source terms in light of Lemma~\ref{lem:sourceadmtimesgood}, Lemma~\ref{lem:simpleparaproduct} as well as Lemma~\ref{lem:sourcetermsotherthannullformprop738prep}. If we reduce this null-form to the paraproduct version which can again not be handled via Prop.~\ref{prop:goodcubicnullform1}, we arrive at the quintilinear expression
\begin{align*}
	&\chi_{R\gtrsim\tau}\left(\sum_{\lambda_2}\phi_{1,<\lambda_2}\phi_{2,\lambda_2}\right)\sum_{\lambda_5,\lambda_6}\mathcal{N}_0\left(P_{<\min\{\lambda_5,\lambda_6\}}\phi_4, \phi_{5,\lambda_5}, \phi_{6,\lambda_6}\right)\\
	& = \chi_{R\gtrsim\tau}\left(\sum_{\lambda_2}\phi_{1,<\lambda_2}\phi_{2,\lambda_2}\right)\sum_{\lambda_5<\lambda_6}\mathcal{N}_0\left(P_{<\lambda_5}\phi_4, \phi_{5,\lambda_5}, \phi_{6,\lambda_6}\right)\\& + \chi_{R\gtrsim\tau}\left(\sum_{\lambda_2}\phi_{1,<\lambda_2}\phi_{2,\lambda_2}\right)\sum_{\lambda_5\geq\lambda_6}\mathcal{N}_0\left(P_{<\lambda_6}\phi_4, \phi_{5,\lambda_5}, \phi_{6,\lambda_6}\right)
\end{align*}
At this point, again taking advantage of Lemma~\ref{lem:E67paraprodterms} , it is clear how to proceed: taking the first term on the right (the second being handled analogously), if $\lambda_2\geq \lambda_5$, this is a good source term (we can simply switch the positions between $\phi_{2}$ and $\phi_{4}$ in the quintilinear expression). Thus we reduce this term to 
\begin{align*}
	\chi_{R\gtrsim\tau}\sum_{\lambda_5<\lambda_6}\left(\sum_{\lambda_2<\lambda_5}\phi_{1,<\lambda_2}\phi_{2,\lambda_2}\right)\mathcal{N}_0\left(P_{<\lambda_5}\phi_4, \phi_{5,\lambda_5}, \phi_{6,\lambda_6}\right).
\end{align*}
Using the same argument as the one deriving \eqref{eq:nullformidentity1}, the problem reduces to bounding 
\[
\chi_{R\gtrsim\tau}\sum_{\lambda_5<\lambda_6}\left(\sum_{\lambda_2<\lambda_5}\phi_{1,<\lambda_2}\phi_{2,\lambda_2}\right)\phi_{4,<\lambda_5}\phi_{5,\lambda_5}P_{\lambda_6}\left(\Box_{n_6}\phi_{6}\right),
\]
which, up to good error terms, can be replaced by 
\[
\chi_{R\gtrsim\tau}\sum_{\lambda_5}\left(\sum_{\lambda_2<\lambda_5}\phi_{1,<\lambda_2}\phi_{2,\lambda_2}\right)\left(\phi_{4,<\lambda_5}\phi_{5,\lambda_5}\right)\left(\Box_{n_6}\phi_{6}\right).
\]
Now the process is repeated. The general step of this inductive procedure is encapsulated by the following 
\begin{proposition}\label{prop:nullformgeneralindictivestepwithnormaltransform} Assume that the functions $\phi_j$ are angular momentum $n_j$-functions, $|n_j|\geq 2$, $j = 1,\,2,\ldots,2k+1$, each of which has distorted Fourier transform (at angular momentum $n_j$) admitting a splitting 
\[
\langle \phi_{n_j}(R;\xi),\,\phi_j\rangle_{L^2_{R\,dR}} = \overline{y}_j + \overline{z}_j
\]
such that $ \overline{y}_j$ is admissibly singular while $\overline{z}_j\in S_0^{\hbar_j}$,\,$\hbar_j = \frac{1}{|n_j|}$. Furthermore assume that each angular momentum $n_j$-function $\phi_j$ satisfies a wave equation 
\begin{align*}
\Box_{n_j}\phi_j = F_j,
\end{align*}
where the angular momentum $n_j$-functions $F_j$ have the property that they are either source admissible\footnote{This means their distorted Fourier transform is the sum of a source admissibly singular function and a function in $S_1^{\hbar_j}$, $\hbar_j = |n_j|^{-1}$.}, or at least that 
\[
\sum_{\lambda>1}P_{\lambda}F_j\cdot P_{\geq \lambda}\psi
\]
is source admissible\footnote{By Lemma~\ref{lem:badnullformpara}, this is the case when $F_j$ is a null-form $\mathcal{N}_0$.}, provided $\psi$ is any of the $\phi_l$. Here we assume that the preceding expression is of size $O(1)$
by reference to the norms used to bound source admissible terms \\
Then if we define $\mathcal{N}_0\big(\phi_1, \phi_2, \phi_3\big)$ as at the beginning of subsection~\ref{subset:hardnullformestimates}, then assuming $|m - \sum_j n_j|\lesssim 1$, $|m|\geq 2$, there is an angular momentum $m$ function $\psi$ with the property that 
\begin{align*}
\langle \phi_m(R;\xi),\,\psi\rangle = \overline{y}_{\psi} + \overline{z}_{\psi}
\end{align*}
and such that (with $\hbar = \frac{1}{|m|}$)
\begin{align*}
&\big\| \overline{y}_{\psi}\big\|_{adm}\leq D^k\cdot \max_{j=1,\ldots,2k+1}\{|n_j|\}^{-C}\cdot \prod_{j=1}^{2k+1}|n_j|^C\big(\big\|\overline{y}_j \big\|_{adm} + \tau\big\| \overline{z}_j\big\|_{S_0^{\hbar_j}}\big)\\
&\tau^3\cdot \big\| \overline{z}_{\psi}\big\|_{S_0^{\hbar}}\lesssim D^k\cdot \max_{j=1,\ldots,2k+1}\{|n_j|\}^{-C}\cdot \prod_{j=1}^{2k+1}|n_j|^C\big(\big\|\overline{y}_j \big\|_{adm} + \tau\big\| \overline{z}_j\big\|_{S_0^{\hbar_j}}\big),\\
\end{align*}
and such that the following holds: there is a function $g(\tau, R)$ with the property that\footnote{Recall the notation from subsection~\ref{subsubsec:paraproditerated}.} 
\begin{align*}
\sum_{\lambda_k<\lambda_{k+1}}\Big(\chi_{R\gtrsim 1}\Pi_{<\lambda_k}\big(\phi_1,\ldots,\phi_{2k-2}\big)\cdot \mathcal{N}_0\big(P_{<\lambda_k}\phi_{2k-1},\,P_{\lambda_k}\phi_{2k}, P_{\lambda_{k+1}}\phi_{2k+1}\big) - \Box_m \psi\Big)|_{R<\nu\tau} = g|_{R<\nu\tau}, 
\end{align*}
and such that we can write 
\[
g = g_1 + g_2
\]
where the term $g_1$ is of a form analogous to the term $E_8$ in  \eqref{eq:nullformidentity1}, and given by 
\begin{align*}
g_1 = \sum_{\lambda_k}\Pi_{\lambda_k}\big(\phi_1,\ldots, \phi_{2k}\big)\cdot \Box_{n_{2k+1}} \big(\phi_{2k+1}\big),
\end{align*}
while $g_2$ has angular momentum $m$ distorted Fourier transform which can be decomposed into 
\[
\langle \phi_m(R;\xi),\,g_2\rangle_{L^2_{R\,dR}} = \overline{y}_g + \overline{z}_g
\]
with $\overline{y}_g$ source admissibly singular and $\overline{z}_g\in S_1^{\hbar}$, and we have the bounds
\begin{align*}
&\big\| \overline{y}_{g}\big\|_{sourceadm}\leq D^k\cdot \max_{j=1,\ldots,2k+1}\{|n_j|\}^{-C}\cdot \prod_{j=1}^{2k+1}|n_j|^C\big(\big\|\overline{y}_j \big\|_{adm} + \tau\big\| \overline{z}_j\big\|_{S_0^{\hbar_j}}\big)\\
&\tau^4\cdot \big\| \overline{z}_{g}\big\|_{S_1^{\hbar}}\lesssim D^k\cdot \max_{j=1,\ldots,2k+1}\{|n_j|\}^{-C}\cdot \prod_{j=1}^{2k+1}|n_j|^C\big(\big\|\overline{y}_j \big\|_{adm} + \tau\big\| \overline{z}_j\big\|_{S_0^{\hbar_j}}\big),\\
\end{align*}
\end{proposition}
\begin{proof} We use an analogue of the identity \eqref{eq:nullformidentity1}. Thus write 
\begin{equation}\label{eq:7.24analogue}\begin{split}
&\sum_{\lambda_k<\lambda_{k+1}}\chi_{\nu\tau\gtrsim R\gtrsim 1}\Pi_{<\lambda_k}\big(\phi_1,\ldots,\phi_{2k-2}\big)\cdot \mathcal{N}_0\big(P_{<\lambda_k}\phi_{2k-1},\,P_{\lambda_k}\phi_{2k}, P_{\lambda_{k+1}}\phi_{2k+1}\big)\\
& = \sum_{\lambda_{k+1}}\Box_m'\big(\chi_{\nu\tau\gtrsim R\gtrsim 1}\Pi_{\lambda_{k+1}}\big(\phi_1,\ldots,\phi_{2k+1}\big)\big)\\
& + \sum_{\lambda_{k+1}}\Big(\sum_{l=1}^{2k+1}\frac{n_l^2}{R^2} - \frac{m^2}{R^2}\Big)\cdot \chi_{\nu\tau\gtrsim R\gtrsim 1}\cdot \Pi_{\lambda_{k+1}}\big(\phi_1,\ldots,\phi_{2k+1}\big)\big)\\
& -\sum_{\lambda_{k+1}} \big[\Box_m', \chi_{\nu\tau\gtrsim R\gtrsim 1}\big]\cdot \Pi_{\lambda_{k+1}}\big(\phi_1,\ldots,\phi_{2k+1}\big)\\
& - \sum_* \sum'_{\lambda_{i_1}\leq\lambda_{i_2}\leq \lambda_{i_3}\leq \lambda_{i_4}\leq \lambda_{2k+1}}\mathcal{N}_0\big(\Pi_{*,\lambda_{i_1}}(\phi_1,\ldots,\phi_{i_1}), P_{*,\lambda_{i_2}}\phi_{i_2},\,P_{*,\lambda_{i_3}}\phi_{i_3}\big)\\&\hspace{7cm}\cdot \Pi_{*,\lambda_{i_2},\lambda_{i_3}}(\phi_{i_2+1},\ldots,\phi_{i_3-1})\cdot \Pi_{*,\lambda_{i_4},\lambda_{2k+1}}\big(\phi_{i_4},\ldots, \phi_{2k+1}\big)\\
& - \sum_*\sum'_{\lambda_{i_1}\leq \lambda_{i_2}\leq \lambda_{i_3}}\Pi_{*,\lambda_{i_1}}\big(\phi_1,\ldots, \phi_{i_1}\big)\cdot \Box'_{n_{i_2}}P_{*,\lambda_{i_2}}\phi_{i_2}\cdot \Pi_{*,\lambda_{i_3},\lambda_{i_{2k+1}}}\big(\phi_{i_3},\ldots,\phi_{2k+1}\big)\\
& - \sum_{\lambda_k<\lambda_{k+1}}\chi_{\nu\tau\gtrsim R\gtrsim 1}\Pi_{\lambda_{k}}\big(\phi_1,\ldots,\phi_{2k}\big)\cdot \Box_{n_{2k+1}}'P_{\lambda_{k+1}}\phi_{2k+1}. 
\end{split}\end{equation}
Let us label the sums on the right $E_j, 1\leq j\leq 6$. The operator $\Box_n'$ is defined as in \eqref{eq:nullformidentity1} The complicated sums $E_4, E_5$ are defined as follows: 
\\
{\it{$E_4$}}: in the sum $\sum'$, we sum over all quadruples $1\leq i_1<i_2<i_3<i_4\leq 2k+1$. The symbol $*$ in $P_{*,\lambda}$ indicates that this multiplier is either $P_{<\lambda}$ or $P_{\lambda}$, in the sense of the distorted Fourier transform at the angular momentum of the term it is applied to. If $\lambda_{i_2} = \lambda_{i_3}$, then we necessarily have $P_{*,\lambda_{i_2}}\phi_{i_2} = P_{<\lambda_{i_2}}\phi_{i_2},\,P_{*,\lambda_{i_3}}\phi_{i_3} = P_{\lambda_{i_3}}\phi_3$. The symbol $\Pi_{*,\lambda_{i_1}}(\phi_1,\ldots,\phi_{i_1})$ stands for either $\Pi_{<\lambda_{i_1}}(\phi_1,\ldots,\phi_{i_1})$ or for $\Pi_{\lambda_{i_1}}(\phi_1,\ldots,\phi_{i_1})$. If $\lambda_{i_1} = \lambda_{i_2}$, then necessarily we have $\Pi_{*,\lambda_{i_1}}(\phi_1,\ldots,\phi_{i_1}) = \Pi_{<\lambda_{i_1}}(\phi_1,\ldots,\phi_{i_1})$, as well as $P_{*,\lambda_{i_2}}\phi_{i_2} = P_{\lambda_{i_2}}\phi_{i_2}$. We also let $ \Pi_{*,\lambda_{i_2},\lambda_{i_3}}(\phi_{i_2+1},\ldots,\phi_{i_3-1})$ be the 'truncated' para-product $ \Pi_{*,\lambda_{i_3}}(\phi_{i_2+1},\ldots,\phi_{i_3-1})$ which is given by 
\[
\tilde{P}_{*,\lambda_{i_2}}\phi_{i_2+1}\ldots \tilde{P}_{*,\lambda_{i_3}}\phi_{i_3-1},
\]
where $\tilde{P}_{*,\lambda_{i_2}} = P_{\lambda_{i_2}}$ or $\tilde{P}_{*,\lambda_{i_2}} = P_{<\lambda_{i_2+1}}$, while $\tilde{P}_{*,\lambda_{i_3}} = P_{<\lambda_{i_3}}$ or $\tilde{P}_{*,\lambda_{i_3}} = P_{\lambda_{i_3-1}}$. The paraproduct $\Pi_{*,\lambda_{i_4},\lambda_{2k+1}}\big(\phi_{i_4},\ldots, \phi_{2k+1}\big)$ is defined analogously, and finally, the outermost summation $\sum_{*}$ indicates that we need to sum over all the different possibilities in terms of frequency localisations just indicated. 
\\
{\it{$E_5$}}: Here we sum over $1\leq i_1<i_2<i_3$ and if $\lambda_{i_1} = \lambda_{i_2}$ we have $\Pi_{*,\lambda_{i_1}}\big(\phi_1,\ldots, \phi_{i_1}\big) = \Pi_{<\lambda_{i_1}}\big(\phi_1,\ldots, \phi_{i_1}\big)$, $P_{*,\lambda_{i_2}}\phi_{i_2} = P_{\lambda_{i_2}}\phi_{i_2}$. The truncated para-product $\Pi_{*,\lambda_{i_3},\lambda_{i_{2k+1}}}\big(\phi_{i_3},\ldots,\phi_{2k+1}\big)$ is defined as for the term $E_4$. 
\\

We can now control each of the terms $E_j$, $j = 1,\ldots, 5$, in analogy to the terms in \eqref{eq:nullformidentity1}: 
\\

{\it{$E_1$}}. Set 
\[
\psi: =  \sum_{\lambda_{k+1}}\big(\chi_{\nu\tau\gtrsim R\gtrsim 1}\Pi_{\lambda_{k+1}}\big(\phi_1,\ldots,\phi_{2k+1}\big)\big),
\]
and then 
\begin{align*}
E_1 = \Box_m\psi + \big(\Box_m' - \Box_m\big)\psi. 
\end{align*}
It then suffices to apply the argument of Lemma~\ref{lem:E1errorterms} in conjunction with Lemma~\ref{lem:longparaprod1} to infer the desired decomposition and estimates for $\big(\Box_m' - \Box_m\big)\psi$. This latter term can then be incorporated into $g_2$.
\\

{\it{$E_2$}}. This also contributes to $g_2$. One argues as for the term $E_2$ in  \eqref{eq:nullformidentity1}. First, we note that (recalling $|m - \sum n_l|\lesssim 1$ by assumption)
\begin{align*}
\big|\sum_{l=1}^{2k+1}n_l^2 - m^2\big|\lesssim k^2\cdot\big(\max_{l=1}^{2k+1}|n_l|\big)\cdot \max\{\{|n_l|\}_{l=1}^{2k+1}\backslash \{\max_{l=1}^{2k+1}|n_l|\}\}
\end{align*}
Taking advantage of Lemma~\ref{lem:longparaprod1} write 
\begin{align*}
\langle \phi_m(R;\xi),\,\sum_{\lambda_{k+1}>1} \Pi_{\lambda_{k+1}}\big(\phi_1,\ldots,\phi_{2k+1}\big)\rangle = \overline{y} + \overline{z}, 
\end{align*}
where we have the bound (with $\hbar = \frac{1}{|m|}$)
\begin{align*}
\big\| \overline{y}\big\|_{adm} + \tau^4\cdot \big\|\overline{z}\big\|_{S_0^{\hbar}}\lesssim \max_{l=1}^{2k+1}\{|n_l|\}^{-C}\cdot \prod_{l=1}^{2k+1}\big(|n_l|^C\cdot \big\|\overline{x_{n_l}}\big\|_{good}\big).
\end{align*}
Now calling $X_1$ the contribution to the generalized paraproduct by $\overline{z}$, we have 
\begin{align*}
&\max\{|n_l|\}_{l=1}^{2k+1}\cdot \Big\|\langle \phi_m(R;\xi),\, \chi_{\nu\tau\gtrsim R\gtrsim 1}\cdot R^{-2}X_1\rangle_{L^2_{R\,dR}}\Big\|_{S_1^{\hbar}}\\&
\lesssim k\cdot  \max\{\{|n_l|\}_{l=1}^{2k+1}\backslash \{\max_{l=1}^{2k+1}|n_l|\}\}\cdot \big\|\overline{z}\big\|_{S_0^{\hbar}}, 
\end{align*}
where we have taken advantage of Prop.~\ref{prop:singularmultiplier} together with the elementary fact that 
\[
\max\{|n_l|\}_{l=1}^{2k+1}\lesssim k\cdot |m|\cdot  \max\{\{|n_l|\}_{l=1}^{2k+1}\backslash \{\max_{l=1}^{2k+1}|n_l|\}\}.
\]
This implies the desired estimate for the contribution of $X_1$. As for the contribution of $X_2$, first develop $R^{-2}$ into a Taylor series of fourth order around $R = \nu\tau$. Taking advantage of Lemma~\ref{lem:singFouriertiphysicalngeq2adm}, and placing the product of the Taylor polynomial with $X_2$ into the admissibly singular source part and the product of the Taylor error and $X_2$ into $S_1^{\hbar}$, the desired bound then follows from Lemma~\ref{lem:singPhysicaltoFourierngeq2} and an estimate as before to absorb the factor $\max\{|n_l|\}_{l=1}^{2k+1}$. 
\\

{\it{$E_3$}}. This is similar but simpler than the preceding term since $\big[\Box_m', \chi_{\nu\tau\gtrsim R\gtrsim 1}\big]\sim \tilde{\chi}\cdot R^{-2}$ for a suitable smooth cutoff $\tilde{\chi}$ supported at $R\sim 1$ and $R\sim \nu\tau$. 
\\

{\it{$E_4$}}. Since there is at least one factor $P_{\lambda}\phi_j$ whose frequency dominates the frequencies in the null-form $\mathcal{N}_0$, the conclusion for this term follows by combining Lemma~\ref{lem:longparaprod1} with Lemma~\ref{lem:badnullformpara}, leading to a contribution to $g_2$. 
\\

{\it{$E_5$}}. The estimate (as required for a contribution to $g_2$) for this term follows from the assumption on $\Box_{n_j}\phi_j$ as well as Lemma~\ref{lem:longparaprod1}, once we have replaced $\Box_{n_j}'$ by $\Box_{n_j}$ and moved $\Box_{n_j}$ past $P_{\lambda_j}$. Bounding the contributions of the error terms thereby generated follows from Lemma~\ref{lem:commutatorwithderivatives2}  and the argument for Lemma~\ref{lem:E1errorterms}, in conjunction with Lemma~\ref{lem:longparaprod1}.
\\

{\it{$E_6$}}. This is the term that contributes to $g_1$, as well as to $g_2$. Precisely, arguing as for the preceding term, we can replace $\Box_{n_{2k+1}}'$ by $\Box_{n_{2k+1}}$. Then we write 
\begin{align*}
&\sum_{\lambda_k<\lambda_{k+1}}\chi_{\nu\tau\gtrsim R\gtrsim 1}\Pi_{\lambda_{k}}\big(\phi_1,\ldots,\phi_{2k}\big)\cdot \Box_{n_{2k+1}}P_{\lambda_{k+1}}\phi_{2k+1}\\
& = \Big(\sum_{\lambda>1}\chi_{\nu\tau\gtrsim R\gtrsim 1}\Pi_{\lambda}\big(\phi_1,\ldots,\phi_{2k}\big)\Big)\cdot \Box_{n_{2k+1}}\phi_{2k+1}\\
& - \sum_{\lambda_k\geq \lambda_{k+1}}\chi_{\nu\tau\gtrsim R\gtrsim 1}\Pi_{\lambda_{k}}\big(\phi_1,\ldots,\phi_{2k}\big)\cdot \Box_{n_{2k+1}}P_{\lambda_{k+1}}\phi_{2k+1}\\
\end{align*}
For the second term on the right, taking advantage of Lemma~\ref{lem:longparaprod1} in conjunction with Lemma~\ref{lem:commutatorwithderivatives2} allows us to commute $ \Box_{n_{2k+1}}$ and $P_{\lambda_{k+1}}$ up to a term of type $g_2$, and the remaining term is also an acceptable contribution to $g_2$ by the hypothesis on $\phi_{2k+1}$. The first term on the right equals $g_1$. 
\end{proof}

We note that the expression 
\begin{align*}
\sum_{\lambda_k<\lambda_{k+1}}\Big(\chi_{R\gtrsim 1}\Pi_{<\lambda_k}\big(\phi_1,\ldots,\phi_{2k-2}\big)\cdot \mathcal{N}_0\big(P_{<\lambda_k}\phi_{2k-1},\,P_{\lambda_k}\phi_{2k}, P_{\lambda_{k+1}}\phi_{2k+1}\big)
\end{align*}
differs from the expression obtained upon replacing $\Box_{n_{2k+1}}(\phi_{2k+1})$ by a trilinear null-form $\mathcal{N}_0\big(\phi_{2k+1},\,\phi_{2k+2},\,\phi_{2k+3}\big)$, but one reduces to the preceding type of expression upon observing the following 
\begin{lemma}\label{lem:ullformgeneralindictivestepwithnormaltransformclosingstep} Assume that the $\phi_j,\,m,\,n_j$, are as in the preceding Prop.~\ref{prop:nullformgeneralindictivestepwithnormaltransform}. Then we can write\footnote{The restriction to $\lambda_{k+1}\leq \lambda_{k+2}$ can be achieved by simple symmetry considerations.} 
\begin{align*}
\sum_{\min\{\lambda_k,\,\lambda_{k+2}\}\geq \lambda_{k+1}} \chi_{R\gtrsim 1}\Pi_{\lambda_{k}}\big(\phi_1,\ldots,\phi_{2k}\big)\cdot\mathcal{N}_0\big(\phi_{2k+1},\,P_{\lambda_{k+1}}\phi_{2k+2},\,P_{\lambda_{k+2}}\phi_{2k+3}\big)\big|_{R<\nu\tau} = g\big|_{R<\nu\tau}, 
\end{align*}
and such that we can write 
\begin{align*}
\langle \phi_m(R;\xi),\,g\rangle = \overline{y}_g + \overline{z}_g, 
\end{align*}
with $\overline{y}_g$ source admissibly singular and $\overline{z}_g\in S_1^{\hbar}$, and we have the bounds
\begin{align*}
&\big\| \overline{y}_{g}\big\|_{sourceadm}\leq D^k\cdot \max_{j=1,\ldots,2k+3}\{|n_j|\}^{-C}\cdot \prod_{j=1}^{2k+3}|n_j|^C\big(\big\|\overline{y}_j \big\|_{adm} + \tau\big\| \overline{z}_j\big\|_{S_0^{\hbar_j}}\big)\\
&\tau^4\cdot \big\| \overline{z}_{g}\big\|_{S_1^{\hbar}}\lesssim D^k\cdot \max_{j=1,\ldots,2k+3}\{|n_j|\}^{-C}\cdot \prod_{j=1}^{2k+3}|n_j|^C\big(\big\|\overline{y}_j \big\|_{adm} + \tau\big\| \overline{z}_j\big\|_{S_0^{\hbar_j}}\big),\\
\end{align*}
The same conclusion applies to the expression 
\begin{align*}
\sum_{\lambda_{k+2}\geq \lambda_{k+1}} \chi_{R\gtrsim 1}\Pi_{\lambda_{k}}\big(\phi_1,\ldots,\phi_{2k}\big)\cdot\mathcal{N}_0\big(P_{\geq \lambda_{k+1}}\phi_{2k+1},\,P_{\lambda_{k+1}}\phi_{2k+2},\,P_{\lambda_{k+2}}\phi_{2k+3}\big)\big|_{R<\nu\tau} = g\big|_{R<\nu\tau}, 
\end{align*}
\end{lemma}
\begin{proof} This follows by combining Lemma~\ref{lem:longparaprod1} with Prop.~\ref{prop:goodcubicnullform1} in the region $R\gtrsim \tau$ and by using Prop.~\ref{prop:bilin1}  in the region $1\lesssim R\ll\tau$. For the latte region we also require Lemma~\ref{lem:admsingawayfromshock} to interpret admissibly singular functions as regular functions in $S_0^{\hbar}$ with enough temporal decay. 

\end{proof}

In order to complete the treatment of the null-form terms but still restricted to factors $\phi_j$ at angular momentum $|n_j|\geq 2$, we also have to control the source terms involving the singular functions $U, \epsilon$ (recall subsection~\ref{subsec:maintheorem} for the definition of $U, \epsilon$), which come from the unperturbed solution at angular momentum $n = 0$. The details for this are completely analogous to the proof of Proposition~\ref{prop:nullformgeneralindictivestepwithnormaltransform}, and so we simply formulate the analogous conclusion in the following 
\begin{proposition}\label{prop:nullformgeneralindictivestepwithnormaltransformandUepsiloningredients} Let us define the expressions 
\[
P_{<\lambda}\big(\chi_{R\lesssim\nu\tau}\cdot U\big),\,P_{<\lambda}\big(\chi_{R\lesssim\nu\tau}\cdot \epsilon\big)
\]
in analogy to subsection~\ref{subsubsec:Fourierlocal} via inclusion of cutoffs $\chi_{\xi<\lambda^2}$ in the distorted Fourier representation, at angular momentum $n =0$, as in \cite{KST}. Moreover, assume that the $\phi_j, F_j$, $j = 1,\ldots, 2k+1$ below are as in the statement of Proposition~\ref{prop:nullformgeneralindictivestepwithnormaltransform}. Then if we define $\mathcal{N}_0\big(\psi_1, \psi_2, \psi_3\big)$ as at the beginning of subsection~\ref{subset:hardnullformestimates}, but where now either $\psi_j = \phi_j$ or 
\[
\psi_j\in \{\chi_{R\lesssim\nu\tau}\cdot U,\,\chi_{R\lesssim\nu\tau}\cdot \epsilon\}, 
\]
then assuming $|m - \sum_j n_j|\lesssim 1$, $|m|\geq 2$, there is an angular momentum $m$ function $\psi$ with the property that 
\begin{align*}
\langle \phi_m(R;\xi),\,\psi\rangle = \overline{y}_{\psi} + \overline{z}_{\psi}
\end{align*}
and such that (with $\hbar = \frac{1}{|m|}$)
\begin{align*}
&\big\| \overline{y}_{\psi}\big\|_{adm}\leq D^k\cdot \max_{j=1,\ldots,2k'+1}\{|n_j|\}^{-C}\cdot \prod_{j=1}^{2k+1}|n_j|^C\big(\big\|\overline{y}_j \big\|_{adm} + \tau\big\| \overline{z}_j\big\|_{S_0^{\hbar_j}}\big)\\
&\tau^3\cdot \big\| \overline{z}_{\psi}\big\|_{S_0^{\hbar}}\lesssim D^k\cdot \max_{j=1,\ldots,2k'+1}\{|n_j|\}^{-C}\cdot \prod_{j=1}^{2k+1}|n_j|^C\big(\big\|\overline{y}_j \big\|_{adm} + \tau\big\| \overline{z}_j\big\|_{S_0^{\hbar_j}}\big),\\
\end{align*}
with $k'$ indicating the number of factors $\phi_j$ below, and such that the following holds: there is a function $g(\tau, R)$ with the property that 
\begin{align*}
\sum_{\lambda_k<\lambda_{k+1}}\Big(\chi_{R\gtrsim 1}\Pi_{<\lambda_k}\big(\phi_1,\ldots,\phi_{2k-2}\big)\cdot \mathcal{N}_0\big(P_{<\lambda_k}\psi_{2k-1},\,P_{\lambda_k}\psi_{2k}, P_{\lambda_{k+1}}\psi_{2k+1}\big) - \Box_m \psi\Big)|_{R<\nu\tau} = g|_{R<\nu\tau}, 
\end{align*}
and such that we can write 
\[
g = g_1 + g_2
\]
where the term $g_1$ is of a form analogous to the term $E_8$ in \eqref{eq:nullformidentity1}, and given by 
\begin{align*}
g_1 = \sum_{\lambda_k}\Pi_{\lambda_k}\big(\phi_1,\ldots, \phi_{2k}\big)\cdot \Box_{n_{2k+1}} \big(\phi_{2k+1}\big),
\end{align*}
while $g_2$ has angular momentum $m$ distorted Fourier transform which can be decomposed into 
\[
\langle \phi_m(R;\xi),\,g_2\rangle_{L^2_{R\,dR}} = \overline{y}_g + \overline{z}_g
\]
with $\overline{y}_g$ source admissibly singular and $\overline{z}_g\in S_1^{\hbar}$, and we have the bounds
\begin{align*}
&\big\| \overline{y}_{g}\big\|_{sourceadm}\leq D^k\cdot \max_{j=1,\ldots,2k+1}\{|n_j|\}^{-C}\cdot \prod_{j=1}^{2k+1}|n_j|^C\Big(\big\|\overline{y}_j \big\|_{adm} + \tau^3\big(\big\| \overline{z}_j\big\|_{S_0^{\hbar_j}} + \big\|\mathcal{D}_{\tau}\overline{z}_j\big\|_{S_1^{\hbar_j}}\big)\Big)\\
&\tau^4\cdot \big\| \overline{z}_{g}\big\|_{S_1^{\hbar}}\lesssim D^k\cdot \max_{j=1,\ldots,2k+1}\{|n_j|\}^{-C}\cdot \prod_{j=1}^{2k+1}|n_j|^C\big(\big\|\overline{y}_j \big\|_{adm} + \tau^3\big(\big\| \overline{z}_j\big\|_{S_0^{\hbar_j}} + \big\|\mathcal{D}_{\tau}\overline{z}_j\big\|_{S_1^{\hbar_j}}\big)\Big),\\
\end{align*}
We also have the analogue of Lemma~\ref{lem:ullformgeneralindictivestepwithnormaltransformclosingstep}. 
\end{proposition}

We note that the function $U$, which is the profile for the unperturbed blow up solution under the co-rotational ansatz, satisfies the wave equation 
\[
-\Big(\big(\partial_\tau + \frac{\lambda_{\tau}}{\lambda}R\partial_R\big)^2 + \frac{\lambda_{\tau}}{\lambda}\big(\partial_\tau + \frac{\lambda_{\tau}}{\lambda}R\partial_R\big)\Big)U + \big(\partial_R^2 + \frac{1}{R}\partial_R\big) U = \frac{\sin 2U}{2R^2},
\]
while $\epsilon = U - Q$. One uses this relation, instead of the equation $\Box_{n_j}\phi_j = F_j$ (as in the statement of  Proposition~\ref{prop:nullformgeneralindictivestepwithnormaltransform}), when expanding the cubic null-form. 
\\

The preceding proposition, together with Proposition~\ref{prop:nullformgeneralindictivestepwithnormaltransform} as well as Lemma~\ref{lem:sourcetermsotherthannullformprop738prep} allow us to inductively control all the source terms in \eqref{nonlinearity tau R}, which in turn relies on \eqref{non linear 1 2} - \eqref{coe eq precise 2}, {\it{provided all input functions $\phi_j$ are at angular momentum $|n_j|\geq 2$}}. We note here that all of these estimates still rely on 'input functions' $\phi_j$ whose distorted Fourier transform can be decomposed into an admissibly singular part and a regular part in $S_0^{\hbar}$ with good temporal decay properties. This framework is not quite compatible yet with the Duhamel wave propagator (at fixed angular momentum $|n|\geq 2$), and we need to implement a final refinement of the spaces for the 'input functions', characterised in terms of their distorted Fourier transform, in Definition~\ref{defi:goodfourierrepngeq2}. This will require a final round of minor adjustments of the nonlinear estimates at the very end. 
\\
Before implementing this, we explain how to also include factors depending on the exceptional angular momenta $n = 0,\,\pm 1$.

\subsection{Inclusion of the factors with exceptional angular momentum $n\in \{0,\pm1\}$. Definition of admissibly singular terms}\label{subsec:nexcepadmfns}
Recall from subsection~\ref{subsec:nexcspaces} the structure of functions with angular momentum $n = 0,\,\pm1$, namely 
\[
f(\tau, R) = c_j(\tau)\cdot\phi_j(R) + \phi_j(R)\cdot\int_0^R \left[\phi_j(s)\right]^{-1}\cdot\mathcal{D}_j f(s)\,ds,\quad j = 0,\pm 1, 
\]
and furthermore, we have the Fourier representation 
\begin{align*}
	\mathcal{D}_j f(\tau, R) = \int_0^\infty \phi_j(R,\xi)\cdot \xb(\tau, \xi)\cdot\tilde{\rho}_{j}(\xi)\,d\xi, 
\end{align*}
where the functions $\xb(\tau, \xi)$ got measured for the bilinear estimates away from the light cone via the norms \eqref{eq:S01}, \eqref{eq:S00}, \eqref{eq:S0-1}. Here we introduce the analogues of Definition~\ref{defi:xsingulartermsngeq2proto} and Definition~\ref{defi:xsingulartermsngeq2adm}, which allow to characterize the singularity across the light cone for the exceptional angular momentum functions on the distorted Fourier side. Things are simplified somewhat due to the fact that in the low angular momentum setting we can completely neglect the $n$-dependence, and in particular the somewhat cumbersome terms 
\[
e^{\pm i\hbar^{-1}\rho\left(x_{\sigma'};\alpha\cdot\frac{\lambda(\tau)}{\lambda(\sigma)},\hbar\right)}
\]
disappear from the presentation. On the other hand, we will have to treat each exceptional case on its own, due to the different spectral measures. It is to be kept in mind that the functions $\xb(\tau, \xi)$ in the following definition represent the derivative $\mathcal{D}_j f$. We observe that the improved temporal decay seen here compared to the angular modes $|n|\geq 2$ comes from the modulation theory developed later, and will play a crucial role in being able to control the evolution of the unstable modes, i.e., the functions $c_j(\tau)$, $j = 0,\pm 1$. 
\begin{definition}\label{defi:xsingulartermsnless2proto} For functions of angular momentum $n =1$, we say that $\xb(\tau, \xi)$ is a prototype singular function provided it admits the representation 
	\begin{align*}
		\xb(\tau, \xi) = \sum_{\pm}\sum_{k=1}^{N}\sum_{i=0}^{N_1}\chi_{\xi\geq 1}\frac{e^{\pm i\nu\tau\xi^{\frac12}}}{\xi^{1+k\frac{\nu}{2}}}\left(\log\xi\right)^i\cdot a_{k,i}(\tau) +  \sum_{\pm}\sum_{l =1}^7\sum_{k=1}^{N}\sum_{i=0}^{N_1}\chi_{\xi\geq 1}\frac{e^{\pm i\nu\tau\xi^{\frac12}}}{\xi^{1+k\frac{\nu}{2}+\frac{l}{4}}}\left(\log\xi\right)^i\cdot F_{l,k,i}(\tau, \xi),\\
	\end{align*}
	and furthermore we have the bounds with $k_1\in \{0,1\}$, $k_2\in \{0,1,\ldots,20 -l\}$ 
	\begin{align*}
		&\left|\partial_{\tau}^{k_1}a_{k,i}(\tau)\right|\lesssim \left(\log\tau\right)^{N_1-i}\cdot \tau^{-3-\nu-k_1}\\
		&\left|\partial_{\tau}^{k_1}\partial_{\xi}^{k_2}F_{l,k,i}(\tau, \xi)\right| \lesssim \left(\log\tau\right)^{N_1-i}\cdot \tau^{-3-\nu-k_1}\cdot \xi^{-k_2},
	\end{align*}
	as well as the `closure bounds'
	\[
	\left\|\xi^{20-l+\delta}\partial_{\tau}^{k_1}\partial_{\xi}^{20-l}F_{l,k,i}(\tau, \xi)\right\|_{\dot{C}^{\delta}_{\xi}} \lesssim  \left(\log\tau\right)^{N_1-i}\cdot \tau^{-3-\nu - k_1}.
	\]
	For functions of angular momentum $n = 0$, we say that $\xb(\tau, \xi)$ is a prototype singular function provided it admits the representation 
	\begin{align*}
		\xb(\tau, \xi) = \sum_{\pm}\sum_{k=1}^{N}\sum_{i=0}^{N_1}\chi_{\xi\geq 1}\frac{e^{\pm i\nu\tau\xi^{\frac12}}}{\xi^{\frac32+k\frac{\nu}{2}}}\left(\log\xi\right)^i\cdot a_{k,i}(\tau) +  \sum_{\pm}\sum_{l =1}^7\sum_{k=1}^{N}\sum_{i=0}^{N_1}\chi_{\xi\geq 1}\frac{e^{\pm i\nu\tau\xi^{\frac12}}}{\xi^{\frac32+k\frac{\nu}{2}+\frac{l}{4}}}\left(\log\xi\right)^i\cdot F_{l,k,i}(\tau, \xi).
	\end{align*}
	The coefficient functions are to satisfy the same bounds as in the case $n = 1$. 
	For functions of angular momentum $n = -1$, we say that $\xb(\tau, \xi)$ is a prototype singular function provided it admits the representation 
	\begin{align*}
		\xb(\tau, \xi) = \sum_{\pm}\sum_{k=1}^{N}\sum_{i=0}^{N_1}\chi_{\xi\geq 1}\frac{e^{\pm i\nu\tau\xi^{\frac12}}}{\xi^{2+k\frac{\nu}{2}}}\left(\log\xi\right)^i\cdot a_{k,i}(\tau) +  \sum_{\pm}\sum_{l =1}^7\sum_{k=1}^{N}\sum_{i=0}^{N_1}\chi_{\xi\geq 1}\frac{e^{\pm i\nu\tau\xi^{\frac12}}}{\xi^{2+k\frac{\nu}{2}+\frac{l}{4}}}\left(\log\xi\right)^i\cdot F_{l,k,i}(\tau, \xi)
	\end{align*}
	The same bounds as in the case $n = 1$ apply.  \\
	We define source prototype functions to have the same algebraic structure but to have $\xi^{\frac12}$ less decay towards $\xi = +\infty$, while we set $k_1 = 0$ throughout in the preceding.\\
	Similarly, we have the concept of admissible functions: 
	For functions of angular momentum $n = 1$, we say that $\xb(\tau, \xi)$ is an admissible singular function provided it admits the representation
	\[
	\xb =\xb_{\text{in}} + \xb_{\text{out}}
	\]
	where we have the following representations for the incoming and outgoing parts: 
	\begin{align*}
		\xb_{\text{in}}(\tau, \xi) = &\sum_{\pm}\sum_{k=1}^{N}\sum_{j=0}^{N_{1}}\chi_{\xi\geq 1}\frac{e^{\pm i\nu\tau\xi^{\frac12}}}{\xi^{1+k\frac{\nu}{2}}}\left(\log\xi\right)^{j}\cdot\int_{\tau_{0}}^{\tau}a_{k,j}^{(\pm)}(\tau,\sigma)\,d\sigma\\&+\sum_{\pm}\sum_{l =1}^7\sum_{k=1}^{N}\sum_{i=0}^{N_1}\chi_{\xi\geq 1}^{(l)}\left\langle\xi\right\rangle^{-\frac{l}{4}}\frac{e^{\pm i\nu\tau\xi^{\frac12}}}{\xi^{1+k\frac{\nu}{2}}}\left(\log\xi\right)^i\cdot \int_{\tau_0}^{\tau} F_{l,k,i}^{(\pm)}\left(\tau, \sigma, \frac{\lambda^2(\tau)}{\lambda^2(\sigma)}\xi\right)\,d\sigma,
	\end{align*}
	where the $\pm$-signs in each expression on the right are synchronized, and we have the following bounds, where the $\delta_l$ are small positive numbers decreasing in $l$: 
	\begin{align*}
		&\left|a^{(\pm)}_{k,j}(\tau,\sigma)\right|+\tau^{k_{1}}\left|\partial_{\tau}^{k_{1}}a_{k,j}^{(\pm)}(\tau,\sigma)\right|\lesssim_{k_{1}}\left(\log\tau\right)^{N_{1}-j}\tau^{-3-\nu}\cdot \sigma^{-3},\\
		&\left|\xi^{k_2}\partial_{\tau}^{\iota}\partial_{\xi}^{k_2} F_{l,k,i}^{(\pm)}(\tau, \sigma, \xi)\right|\lesssim  \left(\log\tau\right)^{N_1-i}\tau^{-3-\nu}\cdot\sigma^{-1}\cdot \left[\sigma^{-2} + \kappa\left(\xi^{\frac12}\right)\right],\quad 0\leq k_2\leq 20-l,\quad\iota\in \{0,1\}\\
		&\left\|\xi^{20-l+\delta_l}\partial_{\tau}^{\iota}\partial_{\xi}^{5} F_{l,k,i}^{(\pm)}(\tau, \sigma, \xi)\right\|_{\dot{C}^{\delta_l}_{\xi}(\xi\simeq\lambda)}\lesssim \left(\log\tau\right)^{N_1-i}\tau^{-3-\nu}\cdot\sigma^{-1}\cdot \left[\sigma^{-2} + \kappa\left(\lambda^{\frac12}\right)\right],\quad\iota\in \{0,1\},
	\end{align*}
	and further $\xb_{\text{out}} = \xb_{\text{out},1} + \xb_{\text{out},2}$, with
	\begin{equation}\label{eq:xsingout1n=1}\begin{split}
			\xb_{\text{out}, 1}(\tau, \xi) =  &\sum_{\pm}\sum_{l =0}^7\sum_{k=1}^{N}\sum_{i=0}^{N_1}\chi_{\xi\geq 1}\left\langle\xi\right\rangle^{-\frac{l}{4}}\frac{\left(\log\xi\right)^i}{\xi^{1+\frac{k\nu}{2}}}\\&\hspace{2cm}\cdot\int_{\tau_0}^{\tau}e^{\pm i\left[\left(\nu\tau - 2\frac{\lambda(\tau)}{\lambda(\sigma)}\nu\sigma\right)\xi^{\frac12} \right]}\cdot F^{\pm}_{l,k,i}\left(\tau,\sigma, \frac{\lambda^2(\tau)}{\lambda^2(\sigma)}\xi\right)\,d\sigma, 
	\end{split}\end{equation}
	and  the bounds 
	\begin{align*}
		&\left|\xi^{r}\partial_{\tau}^{\iota}\partial_{\xi}^{r} F_{l,k,i}^{(\pm)}(\tau, \sigma, \xi)\right|\lesssim \left(\log\tau\right)^{N_1-i}\tau^{-3-\nu}\cdot\sigma^{-1}\cdot \left[\sigma^{-2} + \kappa\left(\xi^{\frac12}\right)\right],\quad  0\leq r\leq 20-l,\quad\iota\in \{0,1\}\\
		&\left\|\xi^{20-l+\delta_l}\partial_{\tau}^{\iota}\partial_{\xi}^{20-l} F_{l,k,i}^{(\pm)}(\tau,\sigma, \xi)\right\|_{\dot{C}^{\delta_l}_{\xi}(\xi\simeq\lambda)}\lesssim  \left(\log\tau\right)^{N_1-i}\tau^{-3-\nu}\cdot\sigma^{-1}\cdot \left[\sigma^{-2} + \kappa\left(\lambda^{\frac12}\right)\right],\quad \iota\in \{0,1\}.
	\end{align*}
	The second term $\xb_{\text{out}, 2}$, which is the `outgoing perpetuated singularity' admits the description
	\begin{equation}\label{eq:xsingout2n=1}\begin{split}
			\xb_{\text{out}, 2}(\tau, \xi) =  &\sum_{\pm}\sum_{l =0}^7\sum_{k=1}^{N}\sum_{i=0}^{N_1}\chi_{\xi\geq 1}\left\langle\xi\right\rangle^{-\frac{l}{4}}\frac{\left(\log\xi\right)^i}{\xi^{1+k\frac{\nu}{2}}}\\&\cdot\int_0^\infty \int_{\tau_0}^{\tau}e^{\pm i\left[\nu\left(\frac{\lambda(\tau)}{\lambda(\sigma)}x + \tau\right)\xi^{\frac12} \right]}\cdot G^{\pm}_{l,k,i}\left(\tau, \sigma,x,\frac{\lambda^2(\tau)}{\lambda^2(\sigma)}\xi\right)\,d\sigma dx,
	\end{split}\end{equation}
	where $F^{\pm}_{0,k,i}\left(\tau,\sigma, \frac{\lambda^2(\tau)}{\lambda^2(\sigma)}\xi\right) = b^{\pm}_{k,i}(\tau, \sigma)$, $G^{\pm}_{0,k,i}\left(\tau, \sigma,x,\frac{\lambda^2(\tau)}{\lambda^2(\sigma)}\xi\right) = c^{\pm}_{k,i}(\tau,\sigma,x)$, and we have the bounds 
	\begin{align*}
		&\left\|\xi^{k_2}\partial_{\xi}^{k_2} G_{l,k,i}^{(\pm)}(\tau, \sigma, x, \xi)\right\|_{L_x^1}\lesssim  \left(\log\tau\right)^{N_1-i}\tau^{-3-\nu}\cdot\sigma^{-1}\cdot \left[\sigma^{-2} + \kappa\left(\xi^{\frac12}\right)\right],\quad 0\leq k_2\leq 20-l,\\
		&\left\|\left\|\xi^{20-l+\delta_l}\partial_{\xi}^{20-l} G_{l,k,i}^{(\pm)}(\tau, \sigma, x, \xi)\right\|_{\dot{C}^{\delta_l}_{\xi}(\xi\simeq\lambda)}\right\|_{L_x^1}\lesssim  \left(\log\tau\right)^{N_1-i}\tau^{-3-\nu}\cdot\sigma^{-1}\cdot \left[\sigma^{-2} + \kappa\left(\lambda^{\frac12}\right)\right],
	\end{align*}
	We similarly define admissibly singular functions at angular momentum $n = 0, n = -1$, except that the expression $\xi^{1+k\frac{\nu}{2}}$ gets replaced by $\xi^{\frac32+k\frac{\nu}{2}}, \xi^{2+k\frac{\nu}{2}}$, respectively. 
	\\
	We also define the norms $\lVert x\rVert_{\text{adm}}$ for each of these types of angular momentum $n$, $n\in \{0,\pm 1\}$ functions, in perfect analogy to Definition~\ref{defi:xsingulartermsngeq2admnorm} but with the obvious modifications taking into account the different decay properties. Finally, we define the restricted type of principal singular function in analogy to the case of angular momentum $|n|\geq 2$, and the corresponding norms $\lVert x\rVert_{\text{adm}(r)}$, and similarly for source admissible functions: We say that a function $\overline{y}(\tau,\xi)$ is source admissibly singular at angular momentum $n\in \{0,\,\pm 1\}$,  provided we can write $\overline{y} = \overline{y}_1 + \overline{y}_2$, and where $\xi^{-\frac12}\cdot \overline{y}_1(\tau,\xi)$ is a admissibly singular except that we replace $\tau^{-3-\nu}$ by $\tau^{-4-\nu}$ and we also set $k_1 = 0$ throughout; furthermore, we require $\overline{y}_2$ to be a prototype singular source term.
\end{definition}
In analogy to Lemma~\ref{lem:singFouriertiphysicalngeq2adm}, we have the following lemma for $n=0$. Note that since we work at the level of the derivative for the exceptional angular momenta, we lose one degree of smoothness near the light cone: 
\begin{lemma}\label{lem:singFouriertiphysicaln0adm}
    Assume that $\xb$ is an admissible singular part at angular momentum $n=0$. Then the associated function
    \begin{align*}
        f(\tau,R):=\int_{0}^{\infty}\phi_{0}(R,\xi)\cdot\xb(\tau,\xi)\cdot\trho_{0}(\xi)\,d\xi,
    \end{align*}
    restricted to the light cone $R<\nu\tau$, can be decomposed as
    \begin{align*}
        f=f_{1}+f_{2}
    \end{align*}
    where $f_{1}=f_{1}(\tau,R)$ is a $C^{5}$-function supported in $\nu\tau-R\gtrsim 1$ and satisfying
    \begin{align*}
        \nabla_{R}^{k_{2}}\partial_{\tau}^{k_1}f_{1}(\tau,R)\lesssim (\log\tau)^{N_{1}}\cdot \tau^{-\frac72-\nu}|\nu\tau-R|^{-5},\quad 0\leq k_1+k_{2}\leq 5,
    \end{align*}
    while $f_{2}$ have the explicit form
    \begin{align*}
        f_{2}(\tau,R)=\chi_{|\nu\tau-R|\lesssim 1}\sum_{k=1}^{N}\sum_{j=1}^{N_{1}}(\nu\tau-R)^{-\frac12+\frac{l}{2}+k\nu}\left(\log(\nu\tau-R)\right)^{j}\frac{G_{k,l, j}(\tau,\nu\tau-R)}{\tau^{\frac12}}.
    \end{align*}
    Here the function $G_{k,j}(\tau,x)$ has symbol type behavior with respect to $x$, as follows:
    \begin{align*}
        \left|\partial_{x}^{k_{2}}\partial_{\tau}^{\iota}G_{k,j}(\tau,x)\right|\lesssim \left(\log\tau\right)^{N_{1}-j}\cdot \tau^{-3-\nu}x^{-k_{2}},\quad 0\leq k_{2}\leq 20-l + \lfloor\frac{l}{2}+k\nu\rfloor,\,\iota\in \{0,\,1\},
    \end{align*}
    and we have the bound 
    \begin{align*}
        \left\|x^{20-l+\delta_3}\partial_{x}^{20-l}\partial_{\tau}^{\iota}G_{k,j}(\tau,x)\right\|_{\dot{C}^{\delta_3}}\lesssim \left(\log\tau\right)^{N_{1}-j}\cdot\tau^{-3-\nu},\,\,\iota\in \{0,\,1\}.
    \end{align*}
    If $\xb$ is of restricted admissibly singular type, then the leading part of the singularity with $l = 0$ admits a simpler description since we can set 
    \[
    G_{k,0, j}(\tau,\nu\tau-R) = g_{k,j}(\tau). 
    \]
      Similar conclusions apply to angular momentum $n = \pm 1$ admissibly singular functions. 
   \end{lemma}
       The proof is analogous to the one for Lemma \ref{lem:singFouriertiphysicalngeq2adm}. The last part of the lemma is a consequence of the fact that we do not have the complicated phase shifts $e^{\pm i\hbar^{-1}\rho(x_{\tau},\alpha;\hbar)}$ at low angular momenta. 
\subsection{Source term estimates near the light cone for outputs at angular momentum $|n|\geq 2$ but with arbitrary factors}

Using Lemma~\ref{lem:singFouriertiphysicaln0adm} instead of Lemma~\ref{lem:singFouriertiphysicalngeq2adm} and recalling from subsection~\ref{subsec:nexcspaces} that for the exceptional angular modes, the passage from the differentiated function $\mathcal{D}_nf$ to the original function is effected by means of the formal relation 
\begin{equation}\label{eq:exzcmomrecoveryformula}
f = c_n\cdot \phi_n(R) + \phi_n(R)\cdot \int_0^R\phi_n^{-1}(s)\cdot \mathcal{D}_nf(s)\,ds,
\end{equation}
one deduces in close analogy to the preceding the following
\begin{proposition}\label{prop:preliminarytrilinearnullformgeneralinputsngeq2} 

Assume that $F = F(\psi_1, \psi_2,\ldots, \psi_n)$ is as in Lemma~\ref{lem:sourcetermsotherthannullformprop738prep}, but where now the $\psi_j$ are angular momentum $n_j$ functions for arbitrary $n_j$, including the exceptional ones. Assume that for those $\psi_j$ with $|n_j|\geq 2$  the hypotheses of Lemma~\ref{lem:sourcetermsotherthannullformprop738prep} obtain, while for the exceptional $n_j\in \{0, \pm 1\}$, we $\psi_j$ is given by \eqref{eq:exzcmomrecoveryformula}, where we assume for now that $c_{n_j} = 0$. Then for the distorted Fourier transform$\overline{x}_j$ of $\mathcal{D}_{n_j}\psi_j$, we assume that $\overline{x}_j = \overline{y}_j + \overline{z}_j$ with $\overline{y}_j$ admissibly singular while $\overline{z}_j\in S_0^{n_j}$, $\mathcal{D}_{\tau}\overline{z}_j\in S_1^{n_j}$. In the formula below, we shall simply write $S_{0,1}^{\hbar_j}$ for these spaces also in the case of exceptional angular momenta. Then for $F$ at angular momentum $m,\,|m|\geq 2$, the same conclusion as in Lemma~\ref{lem:sourcetermsotherthannullformprop738prep} obtains, i. e. writing 
\begin{align*}
\langle \phi_m(R;\xi),\,\chi_{1\lesssim R\lesssim\nu\tau}F(\phi_1,\phi_2,\ldots, \phi_k)\rangle_{L^2_{R\,dR}} = \overline{y} + \overline{z},
\end{align*}
we have the bound 
\begin{align*}
\big\| \overline{y}\big\|_{sourceadm} + \tau^4\cdot\big\|\overline{z}\big\|_{S_1^{\hbar}}\lesssim D^k\cdot \max\{|n_j|\}^{-C}\prod_{j=1}^k |n_j|^C\Big(\big\|\overline{y}_j\big\|_{adm} + \tau^3\big(\big\|\overline{z}_j\big\|_{S_0^{\hbar_j}} + \big\|\mathcal{D}_{\tau}\overline{z}_j\big\|_{S_1^{\hbar_j}}\big)\Big)
\end{align*}
for a suitable constant $D = D(\nu)$. A similar analogue holds for Lemma~\ref{lem:sourcetermsotherthannullformprop738prepnearorigin}.\\
Furthermore, assuming that for exceptional angular momentum functions $\psi_j$ we have that $\Box\big(\mathcal{D}_{n_j}\psi_j\big) =:F_j$ is either source admissibly singular (at exceptional angular momentum $n_j$), or that 
\[
\sum_{\lambda>1}P_{\lambda}\big(\Box \psi_j\big)\cdot P_{\geq \lambda}\psi
\]
is source admissibly singular for $\psi$ any of the $\psi_l$, then the conclusions of Prop.~\ref{prop:nullformgeneralindictivestepwithnormaltransform}, Lemma~\ref{lem:ullformgeneralindictivestepwithnormaltransformclosingstep} , Prop.~\ref{prop:nullformgeneralindictivestepwithnormaltransformandUepsiloningredients} hold.

\end{proposition}
\subsection{Estimates for the remaining source terms near the light cone, and for output at angular momentum $|n|\geq 2$}
To complete the first version of the source term estimates near the light cone, we now consider all terms arising in Prop.~\ref{prop:smoothlinearsource}, Prop.~\ref{prop:bilinwithUregular1}, Prop.~\ref{prop:bilinwithUregular2}, Prop.~\ref{prop:bilinwithUregular3}, where we again assume that all inputs $\phi_j$ are the sum of an admissibly singular function and an (unstructured) smooth function. Specifically, we have 

\begin{proposition}\label{prop:generalsourcetermsfirstboundnearlightcone} Write $\varphi_j = \sum_{n\in Z}\varphi_j(n)e^{in\theta},\,j = 1,2$, as well as 
	\[
	\varepsilon_{\pm}(n) = \varphi_1(n)\mp i\varphi_2(n),\quad \varepsilon_-(n) = \overline{\varepsilon_+(-n)}. 
	\]
	Assume that the $\varepsilon_{\pm}(n) $ solve \eqref{eq diag phys} with the source terms $F_{\pm}$ defined in the discussion preceding \eqref{eq diag phys}.
	Furthermore, assume that for $|n|\geq 2$ we have (with $\hbar = \frac{1}{n+1}$)
	\[
	\varepsilon_{+}(n) = \int_0^\infty \phi_{n}(R,\xi)x_n(\tau,\xi)\rho_{n}(\xi)\,d\xi,
	\]
	where $\xb_n = \yb_n+ \zb_n$ with $\yb_n$ admissibly singular, and $\zb_n\in S_0^{\hbar}$, $\mathcal{D}_{\tau}\zb_n\in S_1^{\hbar}$ at each time $\tau$. For the exceptional angular momenta $n\in \{0,\pm 1\}$, we write 
	\[
	\mathcal{D}_n\varepsilon_{+}(n) = \int_0^\infty \phi_n(R,\xi)\xb_n(\tau,\xi)\tilde{\rho}_{n}(\xi)\,d\xi,
	\]
	with $\xb_n = \yb_n+ \zb_n$ and $\yb_n$ admissibly singular in the sense of Definition~\ref{defi:xsingulartermsnless2proto}, and $\zb_n\in S_0^{(n)}, \mathcal{D}_{\tau}\zb_n\in S_1^{(n)}$. Finally, assume that 
	\[
	\varepsilon_{+}(n) = \phi_n(R)\cdot\int_0^R\left[\phi_n(s)\right]^{-1}\mathcal{D}_n\varepsilon_{+}(n)(\tau,s)\,ds,
	\]
	i.e., the unstable part $c_n$ is assumed to vanish. Finally, we assume that for sufficiently large $C>1$
	\begin{align*}
		\sum_{|n|\geq 2}n^{C}\left(\left\|\yb_n\right\|_{adm} + \tau^3\left[\left\|\zb_n\right\|_{S_0^{\hbar}} + \left\|\mathcal{D}_{\tau}\zb_n\right\|_{S_1^{\hbar}}\right]\right) + \sum_{|n|<2}\left\|\yb_n\right\|_{adm} + \tau^3\left[\left\|\zb_n\right\|_{S_0^{(n)}}+ \left\|\mathcal{D}_{\tau}\zb_n\right\|_{S_1^{(n)}}\right] =:\Lambda\ll 1. 
	\end{align*}
	Then if $F_j$ is any one of the source terms in Prop.~\ref{prop:smoothlinearsource}, Prop.~\ref{prop:bilinwithUregular1}, Prop.~\ref{prop:bilinwithUregular2}, Prop.~\ref{prop:bilinwithUregular3}, writing 
	\[
	F_j = \sum_{n\in Z}F_j(n)e^{in\theta}, 
	\]
	then for each $n$ with $|n|\geq 2$, there is $\psi_n(\tau, R)$, $H_n(\tau, R)$, such that 
	\begin{align*}
		&\psi_n(\tau, R) = \int_0^\infty \phi_n(R,\xi)\tilde{\xb}_n(\tau,\xi)\rho_{n}(\xi)\,d\xi,\quad \tilde{\xb}_n = \tilde{\yb}_n + \tilde{\zb}_n,\\
		& H_n(\tau, R) = \int_0^\infty \phi_n(R,\xi)\tilde{\tilde{\xb}}_n(\tau,\xi)\rho_{n}(\xi)\,d\xi,\quad \tilde{\tilde{\xb}}_n = \tilde{\tilde{\yb}}_n + \tilde{\tilde{\zb}}_n,
	\end{align*}
	with 
	\begin{align*}
		&\sum_{|n|\geq 2}n^{C}\left(\left\|\tilde{\yb}_n\right\|_{adm} + \tau^4\left[\left\|\tilde{\zb}_n\right\|_{S_0^{\hbar}} + \left\|\mathcal{D}_{\tau}\tilde{\zb}_n\right\|_{S_1^{\hbar}}\right]\right)\leq \Lambda^2 + \tau_0^{-1}\Lambda, \\
		&\sum_{|n|\geq 2}n^{C}\left(\left\|\tilde{\tilde{\yb}}_n\right\|_{\text{source adm}} + \tau^4\left\|\tilde{\tilde{\zb}}_n\right\|_{S_1^{\hbar}}\right)\leq \Lambda^2 + \tau_0^{-1}\Lambda, 
	\end{align*}
	and such that 
	\[
	\left(F_j(n) - \Box_n\psi_n\right)|_{R<\nu\tau} = H_n(\tau, R)|_{R<\nu\tau}.
	\]
\end{proposition}

\begin{remark}\label{rem:generalsourcetermsfirstboundnearlightcone} We note that depending on the type of term $F_j$, the correction terms $\psi_n$ may all vanish. In fact, they are only required for terms with the characteristic $Q_0$ null-form structure. 
	
\end{remark}
\begin{proof} This is a consequence of Prop.~\ref{prop:preliminarytrilinearnullformgeneralinputsngeq2}. More precisely, for the terms involving the quadratic null-form, it is a consequence of iteratively applying Prop.~\ref{prop:nullformgeneralindictivestepwithnormaltransform}, Lemma~\ref{lem:ullformgeneralindictivestepwithnormaltransformclosingstep}. 

\end{proof}

\section{Closing the estimates for angular momenta $|n|\geq 2$: parametrix bounds for source admissibly singular source terms, and solving the wave equation by Fourier methods}\label{sec:ngeq2close}

\subsection{Parametrix bounds for source admissibly singular terms}\label{subsec:ngeqfourierparametrixsection}
In this subsection, we finally show that the functional framework introduced in Definition~\ref{defi:xsingulartermsngeq2adm}, Definition~\ref{defi:xsingularsourcetermsngeq2adm}, is compatible with the wave parametrix at angular momentum $n, |n|\geq 2$. Recall that this parametrix is explicitly given by 
\begin{align*}
	\int_{\tau_0}^\tau U^{(n)}(\tau, \sigma,\xi)\cdot f\left(\sigma,\frac{\lambda^2(\tau)}{\lambda^2(\sigma)}\xi\right)\,d\sigma,\quad U^{(n)}(\tau, \sigma,\xi) = \frac{\lambda(\tau)}{\lambda(\sigma)}\cdot \frac{\rho_n^{\frac12}\left(\frac{\lambda^2(\tau)}{\lambda^2(\sigma)}\xi\right)}{\rho_n^{\frac12}(\xi)}\cdot \frac{\sin\left[\lambda(\tau)\xi^\frac12\int_{\sigma}^{\tau}\lambda^{-1}(u)\,du\right]}{\xi^{\frac12}}
\end{align*}
Then we have the following 
\begin{proposition}\label{prop:parametrixonadmissiblesource} Let $\yb(\tau, \xi)$ be source admissibly singular at angular momentum $n,\,|n|\geq 2$. Then 
	\[
	\xb(\tau, \xi) := \int_{\tau_0}^\tau U^{(n)}(\tau, \sigma,\xi)\cdot \yb\left(\sigma,\frac{\lambda^2(\tau)}{\lambda^2(\sigma)}\xi\right)\,d\sigma
	\]
	can be written as $\xb = \xb_1 + \xb_2$ where $\xb_1$ is an  admissibly singular function and $\xb_2(\tau,]\cdot)\in S_0^{\hbar},\,\mathcal{D}_{\tau}\xb_2(\tau,\cdot)\in S_1^{\hbar}$ for each $\tau\in [\tau_0, \infty)$. Moreover, if $\yb$ is of restricted singular type, then so is $\xb_1$. We have the bound 
\begin{align*}
\big\|\xb_1\big\|_{adm}\lesssim \big\|\yb\big|_{sourceadm}. 
\end{align*}
\end{proposition}
\begin{remark}\label{rem:prop:parametrixonadmissiblesource} The regular part $\xb_2$ does not have the decay rate $\tau^{-3}$ which we have used for the input functions in most of the preceding source term estimates. In fact, we shall implement a final refinement of the spaces below, at which point we get a more complete statement for the Duhamel propagator. 
\end{remark}
\begin{proof}
	We verify this for the various parts involved in an admissibly singular source term. In the following, when we talk about the `contribution of a certain part' of $\xb(\tau, \xi)$ as displayed in Definition~\ref{defi:xsingulartermsngeq2adm}, we mean the parametrix applied to $f(\tau, \xi) = \tau^{-1}\xi^{\frac12}\cdot \xb(\tau, \xi)$. 
	\\
	
	{\it{The contribution of the first term in \eqref{eq:xsingin}.}} Observe that explicitly spelled out this is the function 
	\begin{align*}
		&\sum_{\pm}\sum_{k=1}^{N}\sum_{i=0}^{N_1}\int_{\tau_0}^\tau U^{(n)}(\tau, \sigma,\xi)\cdot \chi_{\frac{\lambda^2(\tau)}{\lambda^2(\sigma)}\xi\geq \hbar^{-2}}\hbar^{-1}\frac{e^{\pm i\nu\sigma\cdot\frac{\lambda(\tau)}{\lambda(\sigma)}\xi^{\frac12}}}{\left(\frac{\lambda^2(\tau)}{\lambda^2(\sigma)}\xi\right)^{\frac12+k\frac{\nu}{2}}}\left(\log\left(\frac{\lambda^2(\tau)}{\lambda^2(\sigma)}\xi\right)\right)^i\\&\hspace{4cm}\cdot \left(\int_{\tau_0}^{\sigma}e^{\pm i\hbar^{-1}\rho\left(x_{\sigma_1\cdot\frac{\lambda(\tau)}{\lambda(\sigma)}\cdot\frac{\lambda(\sigma)}{\lambda(\sigma_1)}};\alpha\cdot\frac{\lambda(\tau)}{\lambda(\sigma)}\cdot\frac{\lambda(\sigma)}{\lambda(\sigma_1)},\hbar\right)}\cdot \sigma^{-1}a_{k,i}^{(\pm)}(\sigma,\sigma_1)\,d\sigma_1\right)\,d\sigma
	\end{align*}
Here it is important to observe that the inner phase simplifies, which suggests switching the orders of integration:
\[
e^{\pm i\hbar^{-1}\rho\left(x_{\sigma_1\cdot\frac{\lambda(\tau)}{\lambda(\sigma)}\cdot\frac{\lambda(\sigma)}{\lambda(\sigma_1)}};\alpha\cdot\frac{\lambda(\tau)}{\lambda(\sigma)}\cdot\frac{\lambda(\sigma)}{\lambda(\sigma_1)},\hbar\right)} = e^{\pm i\hbar^{-1}\rho\left(x_{\sigma_1\cdot\frac{\lambda(\tau)}{\lambda(\sigma_1)}};\alpha\cdot\frac{\lambda(\tau)}{\lambda(\sigma_1)},\hbar\right)}
\]
The kernel $U^{(n)}(\tau, \sigma,\xi)$ of the parametrix has an oscillatory factor which can be written as 
\[
\frac{e^{-i\left(\nu\tau - \nu\sigma\frac{\lambda(\tau)}{\lambda(\sigma)}\right)\xi^{\frac12}} - e^{i\left(\nu\tau - \nu\sigma\frac{\lambda(\tau)}{\lambda(\sigma)}\right)\xi^{\frac12}}}{2i},
\]
whence combined with the oscillatory phase $e^{\pm i\nu\sigma\cdot\frac{\lambda(\tau)}{\lambda(\sigma)}\xi^{\frac12}}$ we either produce the phase $e^{\pm i\nu\tau\xi^{\frac12}}$ corresponding to an incoming singular term, or else the phase $e^{\pm i\left(\nu\tau - 2\sigma\cdot\frac{\lambda(\tau)}{\lambda(\sigma)}\right)\xi^{\frac12}}$, corresponding to an outgoing singularity. 
\\

{\it{Incoming singularity}}: Explicitly, this is the expression 
\begin{align*}
	&\sum_{\pm}\sum_{k=1}^{N}\sum_{i=0}^{N_1}\hbar^{-1}\frac{e^{\pm i\nu\tau\xi^{\frac12}}}{\xi^{\frac12}}\cdot \int_{\tau_0}^\tau e^{\pm i\hbar^{-1}\rho\left(x_{\sigma_1\cdot\frac{\lambda(\tau)}{\lambda(\sigma_1)}};\alpha\cdot\frac{\lambda(\tau)}{\lambda(\sigma_1)},\hbar\right)}\\&\hspace{3cm}\cdot\left(\int_{\sigma_1}^{\tau} \chi_{\frac{\lambda^2(\tau)}{\lambda^2(\sigma)}\xi\geq \hbar^{-2}}\frac{\lambda(\tau)}{\lambda(\sigma)}\cdot \frac{\rho_n^{\frac12}\left(\frac{\lambda^2(\tau)}{\lambda^2(\sigma)}\xi\right)}{\rho_n^{\frac12}(\xi)}\cdot \frac{\left(\log\left(\frac{\lambda^2(\tau)}{\lambda^2(\sigma)}\xi\right)\right)^i}{\left(\frac{\lambda^2(\tau)}{\lambda^2(\sigma)}\xi\right)^{\frac12+k\frac{\nu}{2}}}\cdot \sigma^{-1}a_{k,i}^{(\pm)}(\sigma,\sigma_1)\,d\sigma\right)\,d\sigma_1
\end{align*}
We claim that up to better terms this is again of incoming principal type. To show this, we need to first expand $\frac{\rho_n^{\frac12}\left(\frac{\lambda^2(\tau)}{\lambda^2(\sigma)}\xi\right)}{\rho_n^{\frac12}(\xi)}$ in a Hankel type expansion 
towards $\xi = +\infty$, i.e., write (see Cor.~\ref{cor:spectral measure expand small n}, Cor.~\ref{cor:SM behavior large n})
\[
\frac{\rho_n^{\frac12}\left(\frac{\lambda^2(\tau)}{\lambda^2(\sigma)}\xi\right)}{\rho_n^{\frac12}(\xi)} = 1 + \left\langle\hbar\xi^{\frac12}\right\rangle^{-1}\cdot g_n(\tau, \sigma, \xi)
\]
where $g_n$ has symbol type behavior and is bounded. The contribution if the term $\langle\hbar\xi^{\frac12}\rangle^{-1}\cdot g_n(\tau, \sigma, \xi)$ is then easily seen to lead to a connecting incoming singular term after multiplication by $\chi_{\hbar^2\xi\gtrsim 1}$ and a term in $S_0^{\hbar}$ after multiplication by $\chi_{\hbar^2\xi\lesssim 1}$. In fact, it suffices to put 
\begin{align*}
F_{l,k,j}^{(\pm)}\big(\tau, \sigma, \frac{\lambda^2(\tau)}{\lambda^2(\sigma)}\xi\big) =C_{i,j}\cdot \chi_{\hbar^2\xi\gtrsim 1}\cdot \int_{\sigma_1}^{\tau} \chi_{\frac{\lambda^2(\tau)}{\lambda^2(\sigma)}\xi\geq \hbar^{-2}}\frac{\lambda(\tau)}{\lambda(\sigma)}\cdot  g_n(\tau, \sigma, \xi)\cdot \frac{\left(\log\left(\frac{\lambda^2(\tau)}{\lambda^2(\sigma)}\right)\right)^{i-j}}{\left(\frac{\lambda^2(\tau)}{\lambda^2(\sigma)}\right)^{\frac12+k\frac{\nu}{2}}}\cdot \sigma^{-1}a_{k,i}^{(\pm)}(\sigma,\sigma_1)\,d\sigma,
\end{align*}
where the required bounds in accordance with Def.~\ref{defi:xsingulartermsngeq2adm} are straightforward to verify. 
Further, we have (where $\ldots$ are the remaining terms in the second integrand of the incoming singularity but with $ \frac{\rho_n^{\frac12}\left(\frac{\lambda^2(\tau)}{\lambda^2(\sigma)}\xi\right)}{\rho_n^{\frac12}(\xi)}$ replaced by $1$)
\begin{align*}
	&\sum_{\pm}\sum_{k=1}^{N}\sum_{i=0}^{N_1}\hbar^{-1}\frac{e^{\pm i\nu\tau\xi^{\frac12}}}{\xi^{\frac12}}\cdot \int_{\tau_0}^\tau e^{\pm i\hbar^{-1}\rho\left(x_{\sigma_1\cdot\frac{\lambda(\tau)}{\lambda(\sigma_1)}};\alpha\cdot\frac{\lambda(\tau)}{\lambda(\sigma_1)},\hbar\right)}\cdot\left(\int_{\sigma_1}^{\tau} \chi_{\frac{\lambda^2(\tau)}{\lambda^2(\sigma)}\xi\geq \hbar^{-2}}\ldots \right)d\sigma\\
	& = \sum_{\pm}\sum_{k=1}^{N}\sum_{i=0}^{N_1}\hbar^{-1} \chi_{\xi\geq \hbar^{-2}}\frac{e^{\pm i\nu\tau\xi^{\frac12}}}{\xi^{\frac12}}\cdot \int_{\tau_0}^\tau e^{\pm i\hbar^{-1}\rho\left(x_{\sigma_1\cdot\frac{\lambda(\tau)}{\lambda(\sigma_1)}};\alpha\cdot\frac{\lambda(\tau)}{\lambda(\sigma_1)},\hbar\right)}\cdot\left(\int_{\sigma_1}^{\tau} \chi_{\frac{\lambda^2(\tau)}{\lambda^2(\sigma)}\xi\geq \hbar^{-2}}\ldots \right)d\sigma\\
	& + \sum_{\pm}\sum_{k=1}^{N}\sum_{i=0}^{N_1}\hbar^{-1}\left(1-\chi_{\xi\geq \hbar^{-2}}\right)\frac{e^{\pm i\nu\tau\xi^{\frac12}}}{\xi^{\frac12}}\cdot \int_{\tau_0}^\tau e^{\pm i\hbar^{-1}\rho\left(x_{\sigma_1\cdot\frac{\lambda(\tau)}{\lambda(\sigma_1)}};\alpha\cdot\frac{\lambda(\tau)}{\lambda(\sigma_1)},\hbar\right)}\cdot\left(\int_{\sigma_1}^{\tau} \chi_{\frac{\lambda^2(\tau)}{\lambda^2(\sigma)}\xi\geq \hbar^{-2}}\ldots\right) d\sigma
\end{align*}
and here the last term is easily seen to be of type $\xb_2$ (as in the statement of the proposition), while the first term is decomposed as follows:
\begin{align*}
	&\sum_{\pm}\sum_{k=1}^{N}\sum_{i=0}^{N_1}\hbar^{-1} \chi_{\xi\geq \hbar^{-2}}\frac{e^{\pm i\nu\tau\xi^{\frac12}}}{\xi^{\frac12}}\cdot \int_{\tau_0}^\tau e^{\pm i\hbar^{-1}\rho\left(x_{\sigma_1\cdot\frac{\lambda(\tau)}{\lambda(\sigma_1)}};\alpha\cdot\frac{\lambda(\tau)}{\lambda(\sigma_1)},\hbar\right)}\cdot\left(\int_{\sigma_1}^{\tau} \chi_{\frac{\lambda^2(\tau)}{\lambda^2(\sigma)}\xi\geq \hbar^{-2}}\ldots\right) d\sigma\\
	& = \sum_{\pm}\sum_{k=1}^{N}\sum_{i=0}^{N_1}\hbar^{-1} \chi_{\xi\geq \hbar^{-2}}\frac{e^{\pm i\nu\tau\xi^{\frac12}}}{\xi^{\frac12}}\cdot \int_{\tau_0}^\tau e^{\pm i\hbar^{-1}\rho\left(x_{\sigma_1\cdot\frac{\lambda(\tau)}{\lambda(\sigma_1)}};\alpha\cdot\frac{\lambda(\tau)}{\lambda(\sigma_1)},\hbar\right)}\cdot\left(\int_{\sigma_1}^{\tau} \ldots\right) d\sigma\\
	& - \sum_{\pm}\sum_{k=1}^{N}\sum_{i=0}^{N_1}\hbar^{-1} \chi_{\xi\geq \hbar^{-2}}\frac{e^{\pm i\nu\tau\xi^{\frac12}}}{\xi^{\frac12}}\cdot \int_{\tau_0}^\tau e^{\pm i\hbar^{-1}\rho\left(x_{\sigma_1\cdot\frac{\lambda(\tau)}{\lambda(\sigma_1)}};\alpha\cdot\frac{\lambda(\tau)}{\lambda(\sigma_1)},\hbar\right)}\cdot\left(\int_{\sigma_1}^{\tau}\left(1-\chi_{\frac{\lambda^2(\tau)}{\lambda^2(\sigma)}\xi\geq \hbar^{-2}}\right)\ldots\right) d\sigma\\
\end{align*}
and here the second term on the right is again of type $\xb_2$ since the support of 
\[
\chi_{\xi\geq \hbar^{-2}}\cdot \left(1-\chi_{\frac{\lambda^2(\tau)}{\lambda^2(\sigma)}\xi\geq \hbar^{-2}}\right),\quad \sigma\leq \tau,
\]
is contained in $\xi\simeq \hbar^{-2}$. 
\\
We conclude that up to terms of type $\xb_2$ or of connecting incoming type, the incoming singularity is given by the function 
\begin{align*}
	&\sum_{\pm}\sum_{k=1}^{N}\sum_{i=0}^{N_1}\hbar^{-1}\chi_{\xi\geq \hbar^{-2}}\frac{e^{\pm i\nu\tau\xi^{\frac12}}}{\xi^{1+k\frac{\nu}{2}}}\cdot \int_{\tau_0}^\tau e^{\pm i\hbar^{-1}\rho\left(x_{\sigma_1\cdot\frac{\lambda(\tau)}{\lambda(\sigma_1)}};\alpha\cdot\frac{\lambda(\tau)}{\lambda(\sigma_1)},\hbar\right)}\\&\hspace{3cm}\cdot\left(\int_{\sigma_1}^{\tau}\frac{\lambda(\tau)}{\lambda(\sigma)}\cdot \frac{\left(\log\left(\frac{\lambda^2(\tau)}{\lambda^2(\sigma)}\xi\right)\right)^i}{\left(\frac{\lambda^2(\tau)}{\lambda^2(\sigma)}\right)^{\frac12+k\frac{\nu}{2}}}\cdot \sigma^{-1}a_{k,i}^{(\pm)}(\sigma,\sigma_1)\,d\sigma\right)\,d\sigma_1,
\end{align*}
where we can expand out the logarithm as 
\[
\left(\log\left(\frac{\lambda^2(\tau)}{\lambda^2(\sigma)}\xi\right)\right)^i = \sum_{k'+l = i}C_{k',l}\left(\log\left(\frac{\lambda^2(\tau)}{\lambda^2(\sigma)}\right)\right)^{k'}\cdot \left(\log\xi\right)^l
\]
and it then suffices to set for fixed $k, i$ and $p\leq i$
\begin{align*}
	\tilde{a}_{k,p}^{(\pm)}(\tau,\sigma_1): = \int_{\sigma_1}^{\tau}\frac{C_{k,l}\left(\log\left(\frac{\lambda^2(\tau)}{\lambda^2(\sigma)}\right)\right)^{i-p}}{\left(\frac{\lambda^2(\tau)}{\lambda^2(\sigma)}\right)^{k\frac{\nu}{2}}}\cdot \sigma^{-1}a_{k,i}^{(\pm)}(\sigma,\sigma_1)\,d\sigma,
\end{align*}
and to verify that in light of Definition~\ref{defi:xsingularsourcetermsngeq2adm} we get the bounds needed according to Definition~\ref{defi:xsingulartermsngeq2adm}
\begin{align*}
	\left|\tilde{a}_{k,p}^{(\pm)}(\tau,\sigma_1)\right| + \tau\left|\partial_{\tau}\tilde{a}_{k,p}^{(\pm)}(\tau,\sigma_1)\right|\lesssim \left(\log\tau\right)^{N_1-p}\tau^{-1-\nu}\cdot \sigma_1^{-3}.
\end{align*}

{\it{Outgoing singularity}}: Explicitly this is the expression 
\begin{align*}
	&\sum_{\pm}\sum_{k=1}^{N}\sum_{i=0}^{N_1}\hbar^{-1}\frac{e^{\pm i\nu\tau\xi^{\frac12}}}{\xi^{\frac12}}\cdot \int_{\tau_0}^\tau e^{\pm i\hbar^{-1}\rho\left(x_{\sigma_1\cdot\frac{\lambda(\tau)}{\lambda(\sigma_1)}};\alpha\cdot\frac{\lambda(\tau)}{\lambda(\sigma_1)},\hbar\right)}\\&\hspace{3cm}\cdot\left(\int_{\sigma_1}^{\tau}e^{\mp 2i\nu\sigma\frac{\lambda(\tau)}{\lambda(\sigma)}\xi^{\frac12}} \chi_{\frac{\lambda^2(\tau)}{\lambda^2(\sigma)}\xi\geq \hbar^{-2}}^{(l)}\frac{\lambda(\tau)}{\lambda(\sigma)}\cdot \frac{\rho_n^{\frac12}\left(\frac{\lambda^2(\tau)}{\lambda^2(\sigma)}\xi\right)}{\rho_n^{\frac12}(\xi)}\cdot \frac{\left(\log\left(\frac{\lambda^2(\tau)}{\lambda^2(\sigma)}\xi\right)\right)^i}{\left(\frac{\lambda^2(\tau)}{\lambda^2(\sigma)}\xi\right)^{\frac12+k\frac{\nu}{2}}}\cdot \sigma^{-1}a_{k,i}^{(\pm)}(\sigma,\sigma_1)\,d\sigma\right)\,d\sigma_1
\end{align*}
We claim that this can be interpreted in the form $\xb_{out,2}$, recalling Definition~\ref{defi:xsingulartermsngeq2adm}. In fact, write 
\begin{align*}
	e^{\pm i\nu\tau\xi^{\frac12}}\cdot e^{\mp 2i\nu\sigma\frac{\lambda(\tau)}{\lambda(\sigma)}\xi^{\frac12}}  = e^{\mp i\nu\left(\tau + 2\sigma\frac{\lambda(\tau)}{\lambda(\sigma)} - 2\tau\right)\xi^{\frac12}}, 
\end{align*}
and introduce the variable 
\[
x: = \frac{2\sigma\frac{\lambda(\tau)}{\lambda(\sigma)} - 2\tau}{\frac{\lambda(\tau)}{\lambda(\sigma_1)}},
\]
which takes non-negative values for $\sigma\leq \tau$. Then set, with $\sigma = \sigma(x,\tau, \sigma_1)$ 
\[
G_{k,0,p}(\tau, \sigma_1,x,\xi): = \chi_{\xi\geq \hbar^{-2}}\cdot\frac{\rho_n^{\frac12}(\xi)}{\rho_n^{\frac12}\left(\frac{\lambda^2(\sigma)}{\lambda^2(\tau)}\xi\right)}\frac{C_{i-p,l}\left(\log\left(\frac{\lambda^2(\tau)}{\lambda^2(\sigma)}\right)\right)^{i-p}}{\left(\frac{\lambda^2(\tau)}{\lambda^2(\sigma)}\right)^{k\frac{\nu}{2}}}\cdot \sigma^{-1}a_{k,i}^{(\pm)}(\sigma,\sigma_1)\cdot\frac{\partial\sigma}{\partial x}\cdot \chi^{(sharp)}_{\sigma\in [\sigma_1,\tau]}, 
\]
where $\chi^{(sharp)}_{\sigma\in [\sigma_1,\tau]}$ denotes the characteristic function corresponding to the interval $[\sigma_1,\tau]$ with respect to $\sigma$, which in turn is interpreted as a function of $x,\sigma_1, \tau$ via the above coordinate change. Here it may be objected that the function $G_{k,0,p}$ thus defined does not quite satisfy the requirements in Definition~\ref{defi:xsingulartermsngeq2adm} that it should not depend on $\xi$ (since $l = 0$), but this can be easily remedied by the arguments in the incoming case, by replacing 
\[
\frac{\rho_n^{\frac12}(\xi)}{\rho_n^{\frac12}\left(\frac{\lambda^2(\sigma)}{\lambda^2(\tau)}\xi\right)}
\]
by $1$ up to connecting singular terms, and abolishing the smooth cutoff $ \chi_{\xi\geq \hbar^{-2}}$ up to terms of more regular type $S_0^{\hbar}$.  It is now straightforward to check that the term can be written in the form \eqref{eq:xsingout2}(but with $\sigma_1$ replacing the variable $\sigma$ there), and that the required bounds for $G_{k,0,p}$ are satisfied.  
\\

The remaining terms in Definition~\ref{defi:xsingulartermsngeq2adm} are of course handled similarly, let us consider 
\\

{\it{The contribution of the term \eqref{eq:xsingout2}.}} We again distinguish between an incoming and an outgoing term depending on the interaction of the oscillatory phase $U^{(n)}(\tau, \sigma,\xi)$ and the phase $e^{\pm 2i\nu\sigma\frac{\lambda(\tau)}{\lambda(\sigma)}\xi^{\frac12}}$ from the source term. 
\\

{\it{Incoming singularity}}: In this case the combination of the phase from $U^{(n)}(\tau, \sigma,\xi)$ and the phase $e^{\pm 2i\nu\sigma\frac{\lambda(\tau)}{\lambda(\sigma)}\xi^{\frac12}}$ simplifies to $e^{\pm i\nu\tau\xi^{\frac12}}$, and so the phase 
\[
e^{\pm i\left[\nu\left(\frac{\lambda(\tau)}{\lambda(\sigma)}x + \tau\right)\xi^{\frac12} + \hbar^{-1}\rho\left(x_{\sigma\cdot\frac{\lambda(\tau)}{\lambda(\sigma)}};\alpha\cdot\frac{\lambda(\tau)}{\lambda(\sigma)},\hbar\right)\right]}
\]
after effecting the re-scaling $\xi\longrightarrow \frac{\lambda^2(\tau)}{\lambda^2(\sigma)}\xi$ gets transformed into 
\[
e^{\pm i\left[\nu\left(\frac{\lambda(\tau)}{\lambda(\sigma_1)}x + \tau\right)\xi^{\frac12} + \hbar^{-1}\rho\left(x_{\sigma_1\cdot\frac{\lambda(\tau)}{\lambda(\sigma_1)}};\alpha\cdot\frac{\lambda(\tau)}{\lambda(\sigma_1)},\hbar\right)\right]}
\]
Considering the case $l = 0$ for simplicity, we obtain the following expression (after splitting up the power of the logarithm $\left(\log\left(\frac{\lambda^2(\tau)}{\lambda^2(\sigma)}\xi\right)\right)^i$ as before)
\begin{align*}
	&\hbar^{-1}\chi_{\xi\geq \hbar^{-2}}\frac{e^{\pm i\nu\tau\xi^{\frac12}}}{\xi^{1+k\frac{\nu}{2}}}\left(\log\xi\right)^{i-p}\cdot\int_0^\infty \int_{\tau_0}^{\tau}e^{\pm i\left[\nu\left(\frac{\lambda(\tau)}{\lambda(\sigma_1)}x + \tau\right)\xi^{\frac12} + \hbar^{-1}\rho\left(x_{\sigma_1\cdot\frac{\lambda(\tau)}{\lambda(\sigma_1)}};\alpha\cdot\frac{\lambda(\tau)}{\lambda(\sigma_1)},\hbar\right)\right]}
	\\&\hspace{2cm}\cdot\left(\int_{\sigma_1}^{\tau} \chi_{\frac{\lambda^2(\tau)}{\lambda^2(\sigma)}\xi\geq \hbar^{-2}}\frac{\lambda(\tau)}{\lambda(\sigma)}\cdot \frac{\rho_n^{\frac12}\left(\frac{\lambda^2(\tau)}{\lambda^2(\sigma)}\xi\right)}{\rho_n^{\frac12}(\xi)}\cdot \frac{\left(\log(\frac{\lambda^2(\tau)}{\lambda^2(\sigma)})\right)^p}{\left(\frac{\lambda^2(\tau)}{\lambda^2(\sigma)}\right)^{\frac12+k\frac{\nu}{2}}}\cdot \sigma^{-1}G_{0,k,i}^{(\pm)}\left(\sigma,\sigma_1,x,\frac{\lambda^2(\tau)}{\lambda^2(\sigma_1)}\xi\right)\,d\sigma\right)\,d\sigma_1\,dx,
\end{align*}
where inclusion of the outer cutoff $\chi_{\xi\geq \hbar^{-2}}$ is again legitimate up to an error in $S_0^{\hbar}$. But then setting 
\begin{align*}
	\tilde{G}_{0,k,i}^{(\pm)}\left(\tau,\sigma_1,x,\frac{\lambda^2(\tau)}{\lambda^2(\sigma_1)}\xi\right): = \int_{\sigma_1}^{\tau} \chi_{\frac{\lambda^2(\tau)}{\lambda^2(\sigma)}\xi\geq \hbar^{-2}}\frac{\lambda(\tau)}{\lambda(\sigma)}\cdot \frac{\rho_n^{\frac12}\left(\frac{\lambda^2(\tau)}{\lambda^2(\sigma)}\xi\right)}{\rho_n^{\frac12}(\xi)}\cdot \frac{\left(\log\left(\frac{\lambda^2(\tau)}{\lambda^2(\sigma)}\right)\right)^p}{\left(\frac{\lambda^2(\tau)}{\lambda^2(\sigma)}\right)^{\frac12+k\frac{\nu}{2}}}\cdot \sigma^{-1}G_{0,k,i}^{(\pm)}\left(\sigma,\sigma_1,x,\frac{\lambda^2(\tau)}{\lambda^2(\sigma_1)}\xi\right)\,d\sigma,
\end{align*}
the desired conclusion follows upon checking the routine bounds according to Definition~\ref{defi:xsingulartermsngeq2adm}. Furthermore, applying $\partial_{\tau} - 2\frac{\lambda_{\tau}}{\lambda}\xi\partial_{\xi}$ is seen to led to an expression which is up to a factor $\xi^{\frac12}$ of the same form.
\\

{\it{Outgoing singularity}}: Here we encounter the more complex expressions 
\begin{align*}
	&\hbar^{-1}\chi_{\xi\geq \hbar^{-2}}\frac{e^{\pm i\nu\tau\xi^{\frac12}}}{\xi^{1+k\frac{\nu}{2}}}\left(\log\xi\right)^{i-p}\cdot\int_0^\infty \int_{\tau_0}^{\tau}\int_{\sigma_1}^{\tau}e^{\mp i\nu\cdot 2\sigma\frac{\lambda(\tau)}{\lambda(\sigma)}\xi^{\frac12} \mp i\nu\frac{\lambda(\tau)}{\lambda(\sigma_1)}x\xi^{\frac12}\mp i\hbar^{-1}\rho\left(x_{\sigma_1\cdot\frac{\lambda(\tau)}{\lambda(\sigma_1)}};\alpha\cdot\frac{\lambda(\tau)}{\lambda(\sigma_1)},\hbar\right)}\\
	&\hspace{2cm}\cdot \chi_{\frac{\lambda^2(\tau)}{\lambda^2(\sigma)}\xi\geq \hbar^{-2}}\frac{\lambda(\tau)}{\lambda(\sigma)}\cdot \frac{\rho_n^{\frac12}\left(\frac{\lambda^2(\tau)}{\lambda^2(\sigma)}\xi\right)}{\rho_n^{\frac12}(\xi)}\cdot \frac{\left(\log\left(\frac{\lambda^2(\tau)}{\lambda^2(\sigma)}\right)\right)^p}{\left(\frac{\lambda^2(\tau)}{\lambda^2(\sigma)}\right)^{\frac12+k\frac{\nu}{2}}}\cdot \sigma^{-1}G_{0,k,i}^{(\pm)}\left(\sigma,\sigma_1,x,\frac{\lambda^2(\tau)}{\lambda^2(\sigma_1)}\xi\right)\,d\sigma\,d\sigma_1\,dx,
\end{align*}
where all the signs are synchronized. In order to cast this in the mold required by Definition~\ref{defi:xsingulartermsngeq2adm}, introduce the new variable 
\[
\tilde{x}: = x+\frac{2\sigma\cdot\frac{\lambda(\tau)}{\lambda(\sigma)} -2\tau}{\frac{\lambda(\tau)}{\lambda(\sigma_1)}},
\]
which only takes non-negative values if $x\geq 0$ and $\sigma\leq \tau$. Then if 
\[
\tilde{x}\geq \frac{2\sigma_1\cdot\frac{\lambda(\tau)}{\lambda(\sigma_1)} -2\tau}{\frac{\lambda(\tau)}{\lambda(\sigma_1)}},
\]
the restriction on $\sigma$ remains $\tau\geq \sigma\geq \sigma_1$, while if 
\[
\tilde{x}<\frac{2\sigma_1\cdot\frac{\lambda(\tau)}{\lambda(\sigma_1)} -2\tau}{\frac{\lambda(\tau)}{\lambda(\sigma_1)}},
\]
the lower limit $\sigma_{min} = \sigma_{min}(\tau, \sigma_1,\tilde{x})$ for $\sigma$ is given by the condition 
\[
\tilde{x}=\frac{2\sigma_{min}\cdot\frac{\lambda(\tau)}{\lambda(\sigma_{min})} -2\tau}{\frac{\lambda(\tau)}{\lambda(\sigma_1)}}.
\]
Denoting as before by $\chi^{sharp}$ the sharp (characteristic function) cutoff, we can then set 
\begin{align*}
	&\tilde{G}_{0,k,i}^{(\pm)}\left(\tau,\sigma_1,\tilde{x},\frac{\lambda^2(\tau)}{\lambda^2(\sigma_1)}\xi\right): = \\&
	\int_{\sigma_1}^{\tau}\chi^{sharp}_{\tilde{x}\geq \frac{2\sigma_1\cdot\frac{\lambda(\tau)}{\lambda(\sigma_1)} -2\tau}{\frac{\lambda(\tau)}{\lambda(\sigma_1)}}} \chi_{\frac{\lambda^2(\tau)}{\lambda^2(\sigma)}\xi\geq \hbar^{-2}}\frac{\lambda(\tau)}{\lambda(\sigma)}\cdot \frac{\rho_n^{\frac12}\left(\frac{\lambda^2(\tau)}{\lambda^2(\sigma)}\xi\right)}{\rho_n^{\frac12}(\xi)}\cdot \frac{\left(\log\left(\frac{\lambda^2(\tau)}{\lambda^2(\sigma)}\right)\right)^p}{\left(\frac{\lambda^2(\tau)}{\lambda^2(\sigma)}\right)^{\frac12+k\frac{\nu}{2}}}\\&\hspace{4cm}\cdot \sigma^{-1}G_{0,k,i}^{(\pm)}\left(\sigma,\sigma_1,\tilde{x}-\frac{2\sigma\cdot\frac{\lambda(\tau)}{\lambda(\sigma)} -2\tau}{\frac{\lambda(\tau)}{\lambda(\sigma_1)}},\frac{\lambda^2(\tau)}{\lambda^2(\sigma_1)}\xi\right)\,d\sigma\\
	& + \int_{\sigma_{min}}^{\tau}\chi^{sharp}_{\tilde{x}< \frac{2\sigma_1\cdot\frac{\lambda(\tau)}{\lambda(\sigma_1)} -2\tau}{\frac{\lambda(\tau)}{\lambda(\sigma_1)}}}\ldots\,d\sigma, 
\end{align*}
which results in the desired representation 
\begin{align*}
	&\hbar^{-1}\chi_{\xi\geq \hbar^{-2}}\frac{1}{\xi^{1+k\frac{\nu}{2}}}\left(\log\xi\right)^{i-p}\cdot\int_0^\infty \int_{\tau_0}^{\tau}e^{\mp i\left[\nu\left(\frac{\lambda(\tau)}{\lambda(\sigma)}\tilde{x} + \tau\right)\xi^{\frac12} + \hbar^{-1}\rho\left(x_{\sigma_1\cdot\frac{\lambda(\tau)}{\lambda(\sigma_1)}};\alpha\cdot\frac{\lambda(\tau)}{\lambda(\sigma_1)},\hbar\right)\right]}\cdot \tilde{G}^{\pm}_{0,k,i}\left(\tau, \sigma_1,\tilde{x},\frac{\lambda^2(\tau)}{\lambda^2(\sigma_1)}\xi\right)\,d\sigma_1 d\tilde{x}.
\end{align*}
Verification of the required bounds for $\tilde{G}_{0,k,i}^{(\pm)}(\tau,\sigma_1,\tilde{x},\xi)$ is straightforward and omitted. 
\end{proof}
\subsection{Solution of the inhomogeneous wave equation on the distorted Fourier side}\label{subsec:inhomfouriersolution}

\subsubsection{The wave equation for angular momentum $|n|\geq 2$, formulation on the Fourier side}
Our point of departure is the fundamental equation \eqref{eq:RegularFinestructure1}, which we re-cast on the Fourier side. We emphasize at this point that for the technical reasons explained in subsection~\ref{subset:hardnullformestimates}, it will be necessary to modify this equation by changing the variable $\varepsilon_{\pm}(n)$ to eliminate certain bad source terms. This, however, has no bearing on the subsequent Fourier methods, and so we shall formally work with the `wrong equation'  \eqref{eq:RegularFinestructure1}. Direct translation of this action to the Fourier side by applying $\mathcal{F}^{\hbar}$ and setting ($\hbar=\frac{1}{n+1}$)
\[
\varepsilon_{\pm}(n) = \int_0^\infty \phi_{n}(R,\xi)\cdot \xb^{\hbar}(\tau, \xi)\rho_{n}(\xi)\,d\xi
\]
results in (see \eqref{eq tau xi hbar pre}, \eqref{eq tau xi hbar D} and \eqref{def calD})
\begin{equation}\label{eq:ngeq2Fourier1}
	-\left(\partial_{\tau} - 2\frac{\lambda'(\tau)}{\lambda(\tau)}\xi\partial_{\xi} + \frac{\lambda'(\tau)}{\lambda(\tau)}\mathcal{K}_{\hbar}\right)^2\overline{x}^{\hbar} - \frac{\lambda'(\tau)}{\lambda(\tau)}\left(\partial_{\tau} - 2\frac{\lambda'(\tau)}{\lambda(\tau)}\xi\partial_{\xi} + \frac{\lambda'(\tau)}{\lambda(\tau)}\mathcal{K}_{\hbar}\right)\overline{x}^{\hbar} - \xi\overline{x}^{\hbar} = \mathcal{F}^{\hbar}\left(F_\pm(n)\right)
\end{equation}
Introducing the important time derivative dilation type operator\footnote{We suppress the dependence on $n$ here to simplify the notation} 
\begin{equation}\label{eq:Dtaungeq2definition}
	\mathcal{D}_{\tau}: = \partial_{\tau} - 2\frac{\lambda'(\tau)}{\lambda(\tau)}\xi\partial_{\xi} - \frac{\lambda'(\tau)}{\lambda(\tau)}\frac{\rho_n'(\xi)\xi}{\rho_n(\xi)} - 2\frac{\lambda'(\tau)}{\lambda(\tau)}, 
\end{equation}
the preceding equation gets recast in the form 
\begin{equation}\label{eq:ngeq2Fourier2}\begin{split}
		&-\left(\mathcal{D}_{\tau}^2 + \frac{\lambda'(\tau)}{\lambda(\tau)}\mathcal{D}_{\tau} + \xi \right)\overline{x}^{\hbar}\\
		& =  \mathcal{F}^{\hbar}\left(F_\pm(n)\right) +  2\frac{\lambda'(\tau)}{\lambda(\tau)}\mathcal{K}_{\hbar}^{(0)}\mathcal{D}_{\tau}\overline{x}^{\hbar} + \left(\frac{\lambda'(\tau)}{\lambda(\tau)}\right)'\mathcal{K}_{\hbar}^{(0)}\overline{x}^{\hbar} + \frac{\lambda'(\tau)}{\lambda(\tau)}\left[\mathcal{D}_{\tau},\mathcal{K}_{\hbar}^{(0)}\right]\overline{x}^{\hbar} + \left(\frac{\lambda'(\tau)}{\lambda(\tau)}\right)^2\left(\left(\mathcal{K}_{\hbar}^{(0)}\right)^2 + \mathcal{K}_{\hbar}^{(0)}\right)\overline{x}^{\hbar} \\
		& =:\sum_{j=1}^5 f_j(\tau, \xi).  
\end{split}\end{equation}
The first order of the day to solve this equation will be to show that applying the Duhamel parametrix to the last four terms on the right hand side will send admissible functions into admissible functions. In fact, due to the smoothing property of the operator $\mathcal{K}_{\hbar}^{(0)}$, the resulting functions won't have a principal incoming part anymore. 
\begin{proposition}\label{prop:sourcetermswithtransferenceopngeq2} Assume that $\overline{x}^{\hbar}(\tau, \xi)$ is admissibly singular, and that $ \mathcal{F}^{\hbar}\left(F_\pm(n)\right)$ is an admissibly singular source term. Then 
	\[
	\int_{\tau_0}^\tau U^{(n)}(\tau, \sigma,\xi)\cdot f_j\left(\sigma,\frac{\lambda^2(\tau)}{\lambda^2(\sigma)}\xi\right)\,d\sigma,\quad j = 1,2,\ldots,5,
	\]
	is admissibly singular up to functions $\zb_j(\tau, \xi)$ in $S_0^{\hbar}$ and with $\mathcal{D}_{\tau}\zb_j\in S_1^{\hbar}$. Furthermore, calling $\overline{y}_j$, $j = 1, 2,\ldots,5$ the admissibly singular part for each contribution, then we have the bounds 
	\begin{align*}
	\big\|\overline{y}_1\big\|_{adm}\lesssim \big\|\mathcal{F}^{\hbar}\left(F_\pm(n)\right)\big\|_{sourceadm},\,\big\|\overline{y}_j\big\|_{adm}\lesssim \big\|\overline{x}^{\hbar}\big\|_{adm},\,j = 2,\ldots, 5. 
	\end{align*}
	If $\mathcal{F}^{\hbar}\left(F_\pm(n)\right)$ is of restricted type, then so is the principal part of $\overline{y}_1$, while the $\overline{y}_j$ are all of connecting type. 
\end{proposition}
\begin{proof} 
	According to the preceding proposition it suffices to check that the terms $f_j$, $j = 2,\ldots,5$ are source admissibly singular. It follows immediately from Definition~\ref{defi:xsingulartermsngeq2adm} that $\mathcal{D}_{\tau}\overline{x}^{\hbar}$ is of the form $\xi^{\frac12}\cdot \yb(\tau,\xi)$ with $\yb$ admissibly singular, up to better terms in $S_0^{\hbar}$. Using Proposition~\ref{prop:transferenceonsingularngeq2}, we conclude that all the terms $f_j,\,j\geq 2$, are source admissibly singular or better (in $S_1^{\hbar}$). Then the desired conclusion follows upon application of Proposition~\ref{prop:parametrixonadmissiblesource}. In fact, the last property follows by combining the preceding proposition and the fact that the off-diagonal part of the transference operator transforms the principal part into a connecting admissibly singular part. 
\end{proof}

\subsubsection{A final refinement of the spaces, with a view toward dealing with the exceptional angular modes}

As seen in many instances before, admissibly singular terms lead to error terms in the good space $S_0^{\hbar}$, for example when applying the transference operator, or when translating things from the Fourier side to the physical side as in Lemma~\ref{lem:singFouriertiphysicalngeq2adm}. Since admissibly singular terms have poor temporal decay properties (of the order $\tau^{-1-\nu}(\log\tau)^C$), this will mean the generation of poorly decaying terms with good regularity properties, but these will not be good enough to control the evolution of the instabilities of the exceptional modes. This problem forces us to refine the part of the distorted Fourier transform of $\varepsilon(n;\tau, R)$ which is in $S_0^{\hbar}$ by splitting it into a structured but poorly decaying part and an un-structured but well-decaying part, as follows: 
\begin{definition}\label{defi:goodfourierrepngeq2} We say that the function $\xb(\tau, \xi)$ is a good Fourier representation (or just a good function) at angular momentum $|n|\geq 2$, provided we can write 
	\begin{align*}
		\xb = \xb_{\text{smooth}} + \xb_{\text{singular}},
	\end{align*}
	where $\xb_{\text{singular}}$ is admissibly singular, and $\xb_{\text{smooth}}(\tau,\cdot)\in S_0^{\hbar}$, $\mathcal{D}_{\tau}\xb_{\text{smooth}}(\tau,\cdot)\in S_1^{\hbar}$, $\forall\tau\geq\tau_0$, and furthermore, we can write 
	\begin{align*}
		\xb_{\text{smooth}} = \xb_{1,\text{smooth}} +  \xb_{2,\text{smooth}}, 
	\end{align*}
	where the first function on the right is structured but temporally slowly decaying. The structure is similar to the one for the admissibly singular terms, but we require a minimum of smoothness for the coefficient functions:
	\begin{align*}
		\xb_{1,\text{smooth}}(\tau, \xi) &=\hbar^{-1}\cdot\langle\hbar^{2}\xi\rangle^{-4}\frac{e^{\pm i\nu\tau\xi^{\frac12}}}{\xi^{1+\delta}}\cdot \int_{\tau_0}^{\tau}\chi_{x_{\sigma\cdot\frac{\lambda(\tau)}{\lambda(\sigma)}}\gtrsim 1}\cdot e^{\pm i\hbar^{-1}\rho\left(x_{\sigma\cdot\frac{\lambda(\tau)}{\lambda(\sigma)}};\alpha\cdot\frac{\lambda(\tau)}{\lambda(\sigma)},\hbar\right)}\cdot  F^{(\pm)}\left(\tau, \sigma, \frac{\lambda^2(\tau)}{\lambda^2(\sigma)}\xi\right)\,d\sigma\\
		& + \hbar^{-1}\cdot\langle\hbar^{2}\xi\rangle^{-4}\frac{1}{\xi^{1+\delta}}\cdot\int_{\tau_0}^{\tau}\chi_{x_{\sigma\cdot\frac{\lambda(\tau)}{\lambda(\sigma)}}\gtrsim 1}\cdot e^{\pm i\left[\left(\nu\tau - 2\frac{\lambda(\tau)}{\lambda(\sigma)}\nu\sigma\right)\xi^{\frac12} - \hbar^{-1}\rho\left(x_{\sigma\cdot\frac{\lambda(\tau)}{\lambda(\sigma)}};\alpha\cdot\frac{\lambda(\tau)}{\lambda(\sigma)},\hbar\right)\right]}\cdot \tilde{F}^{(\pm)}\left(\tau,\sigma, \frac{\lambda^2(\tau)}{\lambda^2(\sigma)}\xi\right)\,d\sigma\\
		& + \hbar^{-1}\cdot\langle\hbar^{2}\xi\rangle^{-4}\frac{1}{\xi^{1+\delta}}\cdot\int_0^\infty \int_{\tau_0}^{\tau}\chi_{x_{\sigma\cdot\frac{\lambda(\tau)}{\lambda(\sigma)}}\gtrsim 1}\cdot e^{\pm i\left[\nu\left(\frac{\lambda(\tau)}{\lambda(\sigma)}x + \tau\right)\xi^{\frac12} + \hbar^{-1}\rho\left(x_{\sigma\cdot\frac{\lambda(\tau)}{\lambda(\sigma)}};\alpha\cdot\frac{\lambda(\tau)}{\lambda(\sigma)},\hbar\right)\right]}\cdot G^{\pm}\left(\tau, \sigma,x,\frac{\lambda^2(\tau)}{\lambda^2(\sigma)}\xi\right)\,d\sigma dx,
	\end{align*}
	and we have the bounds (with $\delta, \delta_1$ as in Def.~\ref{defi:xsingulartermsngeq2adm})
	\begin{align*}
		&\left|\xi^{k_2}\partial_{\tau}^{\iota}\partial_{\xi}^{k_2} F^{(\pm)}(\tau, \sigma, \xi)\right|\lesssim  \left(\log\tau\right)^{N_1}\hbar^{-\iota - \delta_1\cdot k_2}\tau^{-1-\iota - \nu}\cdot\sigma^{-1}\cdot \left[\sigma^{-2} + \kappa\left(\hbar\xi^{\frac12}\right)\right],\quad 0\leq k_2\leq 10,\quad\iota\in \{0,1\}\\
		&\left|\xi^{k_2}\partial_{\tau}^{\iota}\partial_{\xi}^{k_2} \tilde{F}^{(\pm)}(\tau, \sigma, \xi)\right|\lesssim  \left(\log\tau\right)^{N_1}\hbar^{-\iota - \delta_1\cdot k_2}\tau^{-1-\iota-\nu}\cdot\sigma^{-1}\cdot \left[\sigma^{-2} + \kappa\left(\hbar\xi^{\frac12}\right)\right],\quad 0\leq k_2\leq 10,\quad \iota\in \{0,1\}\\
		&\left\|\xi^{10+\delta}\partial_{\tau}^{\iota}\partial_{\xi}^{10} F^{(\pm)}(\tau, \sigma, \xi)\right\|_{\dot{C}^{\delta}_{\xi}(\xi\simeq\lambda)}\lesssim \left(\log\tau\right)^{N_1}\hbar^{-\iota - 10\delta_1}\tau^{-1-\iota-\nu}\cdot\sigma^{-1}\cdot \left[\sigma^{-2} + \kappa\left(\hbar\lambda^{\frac12}\right)\right],\quad \iota\in \{0,1\}\\
		&\left\|\xi^{10+\delta}\partial_{\tau}^{\iota}\partial_{\xi}^{10} \tilde{F}^{(\pm)}(\tau, \sigma, \xi)\right\|_{\dot{C}^{\delta}_{\xi}(\xi\simeq\lambda)}\lesssim \left(\log\tau\right)^{N_1}\hbar^{-\iota- 10\delta_1}\tau^{-1-\iota-\nu}\cdot\sigma^{-1}\cdot \left[\sigma^{-2} + \kappa\left(\hbar\lambda^{\frac12}\right)\right],\quad \iota\in \{0,1\}\\
		&\left\|\xi^{k_2}\partial_{\xi}^{k_2} \partial_{\tau}^{\iota}G^{(\pm)}(\tau, \sigma, x, \xi)\right\|_{L_x^1}\lesssim  \left(\log\tau\right)^{N_1}\hbar^{-\iota - \delta_1\cdot k_2}\tau^{-1-\iota-\nu}\cdot\sigma^{-1}\cdot \left[\sigma^{-2} + \kappa\left(\hbar\xi^{\frac12}\right)\right],\quad 0\leq k_2\leq 10,\\
		&\left\|\left\|\xi^{10+\delta}\partial_{\xi}^{10} \partial_{\tau}^{\iota}G^{(\pm)}(\tau, \sigma, x, \xi)\right\|_{\dot{C}^{\delta}_{\xi}(\xi\simeq\lambda)}\right\|_{L_x^1}\lesssim  \left(\log\tau\right)^{N_1}\hbar^{-\iota- 10\delta_1}\tau^{-1-\iota-\nu}\cdot\sigma^{-1}\cdot \left[\sigma^{-2} + \kappa\left(\hbar\lambda^{\frac12}\right)\right],
	\end{align*}
	On the other hand, for the `unstructured term' $\xb_{2,\text{smooth}}$ we have the bounds
	\[
	\left\|\xb_{2,\text{smooth}}(\tau, \cdot)\right\|_{S_0^{\hbar}} + \left\|\mathcal{D}_{\tau}\xb_{2,\text{smooth}}(\tau, \cdot)\right\|_{S_1^{\hbar}}\lesssim \tau^{-3}. 
	\]
	Similarly, we call a function $\xb(\tau, \xi)$ a good source term at angular momentum $|n|\geq 2$, provided it admits the representation
	\[
	\xb = \xb_{\text{smooth}} + \xb_{\text{singular}}
	\]
	where $\xb_{\text{singular}}$ is source admissibly singular, while the smooth part 
	\begin{align*}
		\xb_{\text{smooth}}(\tau, \xi) = \xb_{1,\text{smooth}} +  \xb_{2,\text{smooth}} +  \xb_{3,\text{smooth}},  
	\end{align*}
	where 
	\[
	\xb_{1,\text{smooth}}(\tau, \xi) =\hbar^{-1}\langle\hbar^{2}\xi\rangle^{-4}\chi_{x_{\tau\gtrsim 1}}\frac{e^{\pm i\left(\nu\tau\xi^{\frac12} + \hbar^{-1}\rho\left(x_{\tau};\alpha,\hbar\right)\right)}}{\xi^{\frac12+\delta}}\left(\log\xi\right)^i\cdot \tilde{F}(\tau, \xi)
	\]
	with the coefficient bounds 
	\begin{align*}
		&\left|\partial_{\xi}^{k_2}\tilde{F}^{(\pm)}(\tau, \xi)\right| \lesssim \hbar^{-\delta_1\cdot k_2}\cdot \left(\log\tau\right)^{N_1}\cdot \left[\tau^{-2-\nu-\delta}\cdot \xi^{-k_2} +  \tau^{-2-\nu}\cdot \xi^{-k_2}\cdot\kappa\left(\hbar\xi^{\frac12}\right)\right],\quad 0\leq k_2\leq 10,\\
		& \left\|\xi^{10+\delta}\partial_{\xi}^{10} F^{(\pm)}(\tau, \sigma, \xi)\right\|_{\dot{C}^{\delta}_{\xi}(\xi\simeq\lambda)}\lesssim \hbar^{-10\delta_1}\cdot\left[\tau^{-2-\nu-\delta} +  \tau^{-2-\nu}\cdot \kappa\left(\hbar\lambda^{\frac12}\right)\right]
	\end{align*}
	as well as 
	\begin{align*}
		\xb_{2,\text{smooth}} &= \hbar^{-1}\langle\hbar^{2}\xi\rangle^{-4}\frac{e^{\pm i\nu\tau\xi^{\frac12}}}{\xi^{\frac12+\delta}}\cdot \int_{\tau_0}^{\tau}e^{\pm i\hbar^{-1}\rho\left(x_{\sigma\cdot\frac{\lambda(\tau)}{\lambda(\sigma)}};\alpha\cdot\frac{\lambda(\tau)}{\lambda(\sigma)},\hbar\right)}\cdot  F^{(\pm)}\left(\tau, \sigma, \frac{\lambda^2(\tau)}{\lambda^2(\sigma)}\xi\right)\,d\sigma\\
		& +\hbar^{-1}\langle\hbar^{2}\xi\rangle^{-4}\frac{1}{\xi^{\frac12+\delta}}\cdot\int_{\tau_0}^{\tau}e^{\pm i\left[\left(\nu\tau - 2\frac{\lambda(\tau)}{\lambda(\sigma)}\nu\sigma\right)\xi^{\frac12} - \hbar^{-1}\rho\left(x_{\sigma\cdot\frac{\lambda(\tau)}{\lambda(\sigma)}};\alpha\cdot\frac{\lambda(\tau)}{\lambda(\sigma)},\hbar\right)\right]}\cdot \tilde{F}^{(\pm)}\left(\tau,\sigma, \frac{\lambda^2(\tau)}{\lambda^2(\sigma)}\xi\right)\,d\sigma\\
		& +\hbar^{-1}\langle\hbar^{2}\xi\rangle^{-4}\frac{1}{\xi^{\frac12+\delta}}\cdot\int_0^\infty \int_{\tau_0}^{\tau}e^{\pm i\left[\nu\left(\frac{\lambda(\tau)}{\lambda(\sigma)}x + \tau\right)\xi^{\frac12} + \hbar^{-1}\rho\left(x_{\sigma\cdot\frac{\lambda(\tau)}{\lambda(\sigma)}};\alpha\cdot\frac{\lambda(\tau)}{\lambda(\sigma)},\hbar\right)\right]}\cdot G^{\pm}\left(\tau, \sigma,x,\frac{\lambda^2(\tau)}{\lambda^2(\sigma)}\xi\right)\,d\sigma dx,
	\end{align*}
	and we have the bounds 
	\begin{align*}
		&\left|\xi^{k_2}\partial_{\xi}^{k_2} F^{(\pm)}(\tau, \sigma, \xi)\right|\lesssim  \hbar^{-k_2\delta_1}\cdot\left(\log\tau\right)^{N_1}\tau^{-1-\nu}\cdot\sigma^{-1}\cdot \left[\sigma^{-2} + \kappa\left(\hbar\xi^{\frac12}\right)\right]\cdot\left[\tau^{-1-\delta} +\tau^{-1} \kappa\left(\hbar\frac{\lambda(\sigma)}{\lambda(\tau)}\xi^{\frac12}\right)\right],\quad 0\leq k_2\leq 10,\\
		&\left|\xi^{k_2}\partial_{\xi}^{k_2} \tilde{F}^{(\pm)}(\tau, \sigma, \xi)\right|\lesssim \hbar^{-k_2\delta_1}\cdot\left(\log\tau\right)^{N_1}\tau^{-1-\nu}\cdot\sigma^{-1}\cdot \left[\sigma^{-2} + \kappa\left(\hbar\xi^{\frac12}\right)\right]\cdot\left[\tau^{-1-\delta} + \tau^{-1}\kappa\left(\hbar\frac{\lambda(\sigma)}{\lambda(\tau)}\xi^{\frac12}\right)\right],\quad 0\leq k_2\leq 10,\\
		&\left\|\xi^{10+\delta}\partial_{\xi}^{10} F^{(\pm)}(\tau, \sigma, \xi)\right\|_{\dot{C}^{\delta}_{\xi}(\xi\simeq\lambda)}\lesssim \hbar^{-10\delta_1}\cdot\left(\log\tau\right)^{N_1}\tau^{-1-\nu}\cdot\sigma^{-1}\cdot \left[\sigma^{-2} + \kappa\left(\hbar\lambda^{\frac12}\right)\right]\cdot\left[\tau^{-1-\delta} +\tau^{-1} \kappa\left(\hbar\frac{\lambda(\sigma)}{\lambda(\tau)}\lambda^{\frac12}\right)\right],\\
		&\left\|\xi^{10+\delta}\partial_{\xi}^{10} \tilde{F}^{(\pm)}(\tau, \sigma, \xi)\right\|_{\dot{C}^{\delta}_{\xi}(\xi\simeq\lambda)}\lesssim \hbar^{-10\delta_1}\cdot\left(\log\tau\right)^{N_1}\tau^{-1-\nu}\cdot\sigma^{-1}\cdot \left[\sigma^{-2} + \kappa\left(\hbar\lambda^{\frac12}\right)\right]\cdot\left[\tau^{-1-\delta} +\tau^{-1} \kappa\left(\hbar\frac{\lambda(\sigma)}{\lambda(\tau)}\lambda^{\frac12}\right)\right],\,\\
		&\left\|\xi^{k_2}\partial_{\xi}^{k_2} G^{(\pm)}(\tau, \sigma, x, \xi)\right\|_{L_x^1}\lesssim  \hbar^{-k_2\delta_1}\cdot\left(\log\tau\right)^{N_1}\tau^{-1-\nu}\cdot\sigma^{-1}\cdot \left[\sigma^{-2} + \kappa\left(\hbar\xi^{\frac12}\right)\right]\cdot\left[\tau^{-1-\delta} +\tau^{-1} \kappa\left(\hbar\frac{\lambda(\sigma)}{\lambda(\tau)}\xi^{\frac12}\right)\right],\quad 0\leq k_2\leq 10,\\
		&\left\|\left\|\xi^{10+\delta}\partial_{\xi}^{10} G^{(\pm)}(\tau, \sigma, x, \xi)\right\|_{\dot{C}^{\delta}_{\xi}(\xi\simeq\lambda)}\right\|_{L_x^1}\lesssim  \hbar^{-10\delta_1}\cdot\left(\log\tau\right)^{N_1}\tau^{-1-\nu}\cdot\sigma^{-1}\cdot \left[\sigma^{-2} + \kappa\left(\hbar\lambda^{\frac12}\right)\right]\cdot\left[\tau^{-1-\delta} +\tau^{-1} \kappa\left(\hbar\frac{\lambda(\sigma)}{\lambda(\tau)}\lambda^{\frac12}\right)\right].
	\end{align*}
	Finally
	\begin{align*}
		\left\| \xb_{3,\text{smooth}}(\tau, \cdot)\right\|_{S_1^{\hbar}}\lesssim \tau^{-4}. 
	\end{align*}
	We shall call $\xb_{\text{smooth}}$ the regular or smooth part of the good function $\xb(\tau, \xi)$, and similarly for  good source function. 
	\\
	Finally, we shall say that $\xb$ is a good function with restricted principal singular part, provided that $\xb_{\text{singular}}$ has restricted singular part. 
\end{definition}

\begin{remark}\label{rem:defi:goodfourierrepngeq2} We shall sometimes use the splitting 
\[
 \xb_{1,\text{smooth}} =  \chi_{\hbar^2\xi\lesssim1}\xb_{1,\text{smooth}} +  \chi_{\hbar^2\xi\gtrsim1}\xb_{1,\text{smooth}} =:  \xb_{1,\text{smooth}}^{(l)} +  \xb_{1,\text{smooth}}^{(h)}. 
 \]
Here the term $ \xb_{1,\text{smooth}}^{(h)}$ can be handled for purposes such as pointwise estimates of the corresponding 'physical' function like an admissibly singular function, while the low frequency part $\xb_{1,\text{smooth}}^{(l)}$ needs to be treated separately. 
\end{remark}

It is then natural to introduce a `norm' on good functions, as well as source functions, as follows:
\begin{definition}\label{defy:goodfunctionsnorm} Assume that $\xb(\tau, \xi)$ is a good function at angular momentum $n, |n|\geq 2$, with the representation 
	\[
	\xb = \xb_{\text{smooth}} + \xb_{\text{singular}}
	\]
	and the implied representation for $\xb_{\text{smooth}} = \xb_{1,\text{smooth}} +  \xb_{2,\text{smooth}}$ as in the preceding definition. Then we set 
	\begin{align*}
		\lVert \xb\rVert_{\text{good}}: = \lVert \xb_{\text{smooth}}\rVert_{\text{smooth}} + \lVert \xb_{\text{singular}}\rVert_{\text{adm}},
	\end{align*}
	and we define 
	\begin{align*}
		&\lVert \xb_{\text{smooth}}\rVert_{\text{smooth}} : = \sup_{\tau\geq\tau_0}\tau^3\cdot\big\|\xb_{2,\text{smooth}}(\tau,\cdot)\big\|_{S_0^{\hbar}} + \sup_{\tau\geq\tau_0}\tau^3\cdot \big\|\mathcal{D}_{\tau}\xb_{2,\text{smooth}}(\tau,\cdot)\big\|_{S_1^{\hbar}} \\
		&\sum_{\substack{0\leq k_2\leq 10\\ \iota\in\{0,1\}}}\hbar^{\iota + \delta_1\cdot k_2}\left\|\left(\log\tau\right)^{-N_1}\tau^{1+\iota+\nu}\sigma\cdot \left\|\frac{\xi^{k_2}\partial_{\xi}^{k_2}\partial_{\tau}^{\iota}F^{(\pm)}(\tau,\sigma,\xi)}{\left(\sigma^{-2} + \kappa\left(\hbar\xi^{\frac12}\right)\right)}\right\|_{L_{\xi}^\infty([0,\infty)}\right\|_{L_{\tau}^\infty([\tau_0,\infty)L_{\sigma}^\infty([\tau_0,\tau]}\\
		& + \sum_{\substack{0\leq k_2\leq 10\\ \iota\in\{0,1\}}}\hbar^{\iota + \delta_1\cdot k_2}\left\|\left(\log\tau\right)^{-N_1}\tau^{1+\iota+\nu}\sigma\cdot \left\|\frac{\xi^{k_2}\partial_{\xi}^{k_2}\partial_{\tau}^{\iota}\tilde{F}^{(\pm)}(\tau,\sigma,\xi)}{\left(\sigma^{-2} + \kappa\left(\hbar\xi^{\frac12}\right)\right)}\right\|_{L_{\xi}^\infty([0,\infty)}\right\|_{L_{\tau}^\infty([\tau_0,\infty)L_{\sigma}^\infty([\tau_0,\tau]}\\
		& + \sum_{\iota}\hbar^{\iota + 10\delta_1}\left\|\left(\log\tau\right)^{-N_1}\tau^{1+\iota+\nu}\sigma\cdot \left\|\frac{\xi^{10+\delta}\partial_{\xi}^{10}\partial_{\tau}^{\iota}F^{(\pm)}(\tau,\sigma,\xi)}{\left(\sigma^{-2} + \kappa\left(\hbar\xi^{\frac12}\right)\right)}\right\|_{\dot{C}_{\xi}^{\delta}}\right\|_{L_{\tau}^\infty([\tau_0,\infty)L_{\sigma}^\infty([\tau_0,\tau]}\\
		& +  \sum_{\iota}\hbar^{\iota + 10\delta_1}\left\|\left(\log\tau\right)^{-N_1}\tau^{1+\iota+\nu}\sigma\cdot \left\|\frac{\xi^{10+\delta}\partial_{\xi}^{10}\partial_{\tau}^{\iota}\tilde{F}^{(\pm)}(\tau,\sigma,\xi)}{\left(\sigma^{-2} + \kappa\left(\hbar\xi^{\frac12}\right)\right)}\right\|_{\dot{C}_{\xi}^{\delta}}\right\|_{L_{\tau}^\infty([\tau_0,\infty)L_{\sigma}^\infty([\tau_0,\tau]}\\
		& + \sum_{0\leq k_2\leq 10}\hbar^{\iota + \delta_1\cdot k_2}\left\|\left(\log\tau\right)^{-N_1}\tau^{1+\iota + \nu}\sigma\cdot \left\|\frac{\xi^{k_2}\partial_{\xi}^{k_2}\partial_{\tau}^{\iota}G^{(\pm)}(\tau,\sigma,x,\xi)}{\left(\sigma^{-2} + \kappa\left(\hbar\xi^{\frac12}\right)\right)}\right\|_{L_{\xi}^\infty([0,\infty)}\right\|_{L_x^1([0,\infty)L_{\tau}^\infty([\tau_0,\infty)L_{\sigma}^\infty([\tau_0,\tau]}\\
		& + \hbar^{\iota + 10\delta_1}\left\|\left(\log\tau\right)^{-N_1}\tau^{1+\iota + \nu}\sigma\cdot \left\|\frac{\xi^{10+\delta}\partial_{\xi}^{10}\partial_{\tau}^{\iota}G^{(\pm)}(\tau,\sigma,x,\xi)}{\left(\sigma^{-2} + \kappa\left(\hbar\xi^{\frac12}\right)\right)}\right\|_{\dot{C}_{\xi}^{\delta}}\right\|_{L_x^1L_{\tau}^\infty([\tau_0,\infty)L_{\sigma}^\infty([\tau_0,\tau]}.\\
	\end{align*}
	\\
	Using obvious modifications we similarly define $\lVert x\rVert_{\text{good source}}$. Finally, we define the corresponding norms 
	\[
	\lVert \xb\rVert_{\text{good}(r)},\quad \lVert \xb\rVert_{\text{good source}(r)}
	\]
	if the principal singular part is of restricted type, by replacing $\lVert\cdot\rVert_{\text{adm}}$ by $\lVert\cdot\rVert_{\text{adm}(r)}$, in accordance with Definition~\ref{defi:xsingulartermsngeq2admnorm}. 
\end{definition}

In perfect analogy to Proposition~\ref{prop:parametrixonadmissiblesource}, we have the following 

\begin{proposition}\label{prop:parametrixongoodsource}Let $\yb(\tau, \xi)$ be a good source at angular momentum $n,\,|n|\geq 2$. Then 
	\[
	\xb(\tau, \xi) := \int_{\tau_0}^\tau U^{(n)}(\tau, \sigma,\xi)\cdot \yb\left(\sigma,\frac{\lambda^2(\tau)}{\lambda^2(\sigma)}\xi\right)\,d\sigma
	\]
	is a good function. Moreover, we have the bound
	\begin{align*}
		\lVert \xb\rVert_{\text{good}}\lesssim \lVert \yb\rVert_{\text{goodsource}}
	\end{align*}
	where the implied constant is uniform in $n$ as well as $\tau_0\gg 1$.  
\end{proposition}
\begin{proof} In light of Proposition~\ref{prop:parametrixonadmissiblesource} and the earlier considerations on $S_1^{\hbar}$, it suffices to check this when the source term is of type $\yb_{\text{smooth}}$, and here the argument follows exactly the same lines as the one for Prop.~\ref{prop:sourcetermswithtransferenceopngeq2}.  As in Definition \ref{defi:goodfourierrepngeq2}, we decompose $\yb_{\text{smooth}}$ as
	\begin{align*}
		\yb_{\text{smooth}}=\sum_{j=1}^{3}\yb_{j,\text{smooth}}.
	\end{align*}
The contribution from $\yb_{1,\text{smooth}}$ to the inhomogeneous solution is
\begin{align*}
\xb_{1,\text{smooth}}(\tau,\xi):=&\int_{\tau_{0}}^{\tau}U^{(n)}(\tau,\sigma,\xi)\cdot \yb_{1,\text{smooth}}\left(\sigma,\frac{\lambda^{2}(\tau)}{\lambda^{2}(\sigma)}\xi\right)\,d\sigma\\
=&\hbar^{-1}\int_{\tau_{0}}^{\tau}\chi_{\hbar^{2}\frac{\lambda(\tau)^{2}}{\lambda(\sigma)^{2}}\xi\geq1}U^{(n)}(\tau,\sigma,\xi)\cdot\left(\hbar^{2}\frac{\lambda^{2}(\tau)}{\lambda^{2}(\sigma)}\xi\right)^{-4}\\&\hspace{2cm}\cdot\frac{e^{\pm i\left(\nu\sigma\frac{\lambda(\tau)}{\lambda(\sigma)}\xi^{\frac12}+\hbar^{-1}\rho\left(x_{\sigma\frac{\lambda(\tau)}{\lambda(\sigma)}};\alpha\frac{\lambda(\tau)}{\lambda(\sigma)},\hbar\right)\right)}}{\left(\frac{\lambda^{2}(\tau)}{\lambda^{2}(\sigma)}\xi\right)^{\frac12+\delta}}\left(\log\left(\frac{\lambda^{2}(\tau)}{\lambda^{2}(\sigma)}\xi\right)\right)^{i}\cdot F\left(\sigma,\frac{\lambda^{2}(\tau)}{\lambda^{2}(\sigma)}\xi\right)\,d\sigma.\\
\end{align*}
The kernel $U^{(n)}(\tau,\sigma,\xi)$ can be written as
\begin{align*}
	U^{(n)}(\tau,\sigma,\xi)=&\frac{\lambda(\tau)}{\lambda(\sigma)}\cdot\frac{\rho_{n}^{\frac12}\left(\frac{\lambda^{2}(\tau)}{\lambda^{2}(\sigma)}\xi\right)}{\rho^{\frac12}_{n}(\xi)}\cdot\frac{e^{-i\left(\nu\tau-\nu\sigma\frac{\lambda(\tau)}{\lambda(\sigma)}\right)\xi^{\frac12}}-e^{i\left(\nu\tau-\nu\sigma\frac{\lambda(\tau)}{\lambda(\sigma)}\right)\xi^{\frac12}}}{2i\xi^{\frac12}}.
\end{align*}
Therefore this contribution to $\xb_{1,\text{smooth}}(\tau,\xi)$ consists of the following expressions (up to a multiple of constants):
\begin{align}\label{xb 1 smooth 1}
	\begin{split}
	&\sum_{\pm}\frac{e^{\pm i\nu\tau\xi^{\frac12}}}{\hbar\xi^{\frac12}}\int_{\tau_{0}}^{\tau}\left(\hbar^{2}\frac{\lambda^{2}(\tau)}{\lambda^{2}(\sigma)}\xi\right)^{-4}e^{\pm i\hbar^{-1}\rho\left(x_{\sigma\frac{\lambda(\tau)}{\lambda(\sigma)}};\alpha\frac{\lambda(\tau)}{\lambda(\sigma)},\hbar\right)}\\&\hspace{3cm}\cdot \frac{\lambda(\tau)}{\lambda(\sigma)}\frac{\rho_{n}^{\frac12}\left(\frac{\lambda^{2}(\tau)}{\lambda^{2}(\sigma)}\xi\right)}{\rho^{\frac12}_{n}(\xi)}\cdot\frac{\left(\log\left(\frac{\lambda^{2}(\tau)}{\lambda^{2}(\sigma)}\xi\right)\right)^{i}}{\left(\frac{\lambda^{2}(\tau)}{\lambda^{2}(\sigma)}\xi\right)^{\frac12+\delta}}\cdot F\left(\sigma,\frac{\lambda^{2}(\tau)}{\lambda^{2}(\sigma)}\xi\right)\,d\sigma,
	\end{split}
\end{align}
and
\begin{align}\label{xb 1 smooth 2}
	\begin{split}
	&\sum_{\pm}\chi_{\hbar^{2}\xi\geq1}\frac{e^{\pm i\nu\tau\xi^{\frac12}}}{\hbar\xi^{\frac12}}\int_{\tau_{0}}^{\tau}\left(\hbar^{2}\frac{\lambda^{2}(\tau)}{\lambda^{2}(\sigma)}\xi\right)^{-4}e^{\pm i\hbar^{-1}\rho\left(x_{\sigma\frac{\lambda(\tau)}{\lambda(\sigma)}};\alpha\frac{\lambda(\tau)}{\lambda(\sigma)},\hbar\right)\mp 2i\nu\sigma\frac{\lambda(\tau)}{\lambda(\sigma)}\xi^{\frac12}}\\
	&\hspace{3cm}\cdot \frac{\lambda(\tau)}{\lambda(\sigma)}\frac{\rho_{n}^{\frac12}\left(\frac{\lambda^{2}(\tau)}{\lambda^{2}(\sigma)}\xi\right)}{\rho^{\frac12}_{n}(\xi)}\cdot\frac{\left(\log\left(\frac{\lambda^{2}(\tau)}{\lambda^{2}(\sigma)}\xi\right)\right)^{i}}{\left(\frac{\lambda^{2}(\tau)}{\lambda^{2}(\sigma)}\xi\right)^{\frac12+\delta}}\cdot F\left(\sigma,\frac{\lambda^{2}(\tau)}{\lambda^{2}(\sigma)}\xi\right)\,d\sigma
	\end{split}
\end{align}
Let us set
\begin{align*}
	F\left(\tau,\sigma,\frac{\lambda^{2}(\tau)}{\lambda^{2}(\sigma)}\xi\right):= \frac{\rho_{n}^{\frac12}\left(\frac{\lambda^{2}(\tau)}{\lambda^{2}(\sigma)}\xi\right)}{\rho^{\frac12}_{n}(\xi)}\cdot\left(\log\left(\frac{\lambda^{2}(\tau)}{\lambda^{2}(\sigma)}\xi\right)\right)^{i}\cdot F\left(\sigma,\frac{\lambda^{2}(\tau)}{\lambda^{2}(\sigma)}\xi\right).
\end{align*}
In view of the pointwise bounds for $\tilde{F}(\tau,\xi)$ and $F^{(\pm)}(\tau,\sigma,\xi), \tilde{F}^{(\pm)}(\tau,\sigma,\xi)$ in Definition \ref{defi:goodfourierrepngeq2}, as well as Prop.~\ref{prop:DFT nlarge}, Prop.~\ref{prop: FB match} and their analogues for negative angular momenta, we obtain the desired estimates. 

Next we consider the contribution of $\yb_{2,\text{smooth}}$ in the ``good source".  We start with the 1st term in its expression, which is:
\begin{align}\label{good source 2a}
	\begin{split}
	\yb^{a}_{2,\text{smooth}}\left(\sigma,\frac{\lambda^{2}(\tau)}{\lambda^{2}(\sigma)}\xi\right):=&\hbar^{-1}\chi_{\hbar^{2}\frac{\lambda(\tau)^{2}}{\lambda(\sigma)^{2}}\xi\geq1}\left(\hbar^{2}\frac{\lambda^{2}(\tau)}{\lambda^{2}(\sigma)}\xi\right)^{-4}\frac{e^{\pm i\nu\sigma\frac{\lambda(\tau)}{\lambda(\sigma)}\xi^{\frac12}}}{\left(\frac{\lambda^{2}(\tau)}{\lambda^{2}(\sigma)}\xi\right)^{\frac12+\delta}}\\
	&\cdot \int_{\tau_{0}}^{\sigma}e^{\pm i\hbar^{-1}\rho\left(x_{\sigma_{1}\cdot\frac{\lambda(\tau)}{\lambda(\sigma_{1})}};\alpha\cdot\frac{\lambda(\tau)}{\lambda(\sigma_{1})},\hbar\right)}\cdot F^{(\pm)}\left(\sigma,\sigma_{1},\frac{\lambda^{2}(\tau)}{\lambda^{2}(\sigma_{1})}\xi\right)\,d\sigma_{1}. 
	\end{split}
\end{align}
Upon applying the parametrix, the contribution to $\xb_{1,\text{smooth}}^{a}(\tau,\xi)$ consists of the following two expressions (up to a multiple of constants):
\begin{align}\label{good output 2a1}
	\begin{split}
		&\sum_{\pm}\frac{e^{\pm i\nu\tau\xi^{\frac12}}}{\hbar\xi^{\frac12}}\int_{\tau_{0}}^{\tau}\left(\hbar^{2}\frac{\lambda^{2}(\tau)}{\lambda^{2}(\sigma)}\xi\right)^{-4}\cdot \frac{\lambda(\tau)}{\lambda(\sigma)}\frac{\rho_{n}^{\frac12}\left(\frac{\lambda^{2}(\tau)}{\lambda^{2}(\sigma)}\xi\right)}{\rho^{\frac12}_{n}(\xi)}\cdot\frac{\left(\log\left(\frac{\lambda^{2}(\tau)}{\lambda^{2}(\sigma)}\xi\right)\right)^{i}}{\left(\frac{\lambda^{2}(\tau)}{\lambda^{2}(\sigma)}\xi\right)^{\frac12+\delta}}\\
		&\hspace{3cm}\cdot \int_{\tau_{0}}^{\sigma}e^{\pm i\hbar^{-1}\rho\left(x_{\sigma_{1}\cdot\frac{\lambda(\tau)}{\lambda(\sigma_{1})}};\alpha\cdot\frac{\lambda(\tau)}{\lambda(\sigma_{1})},\hbar\right)}\cdot F^{(\pm)}\left(\sigma,\sigma_{1},\frac{\lambda^{2}(\tau)}{\lambda^{2}(\sigma_{1})}\xi\right)\,d\sigma_{1}\,d\sigma
	\end{split}
\end{align}
and
\begin{align}\label{good output 2a2}
	\begin{split}
		&\sum_{\pm}\frac{e^{\pm i\nu\tau\xi^{\frac12}}}{\hbar\xi^{\frac12}}\int_{\tau_{0}}^{\tau}\left(\hbar^{2}\frac{\lambda^{2}(\tau)}{\lambda^{2}(\sigma)}\xi\right)^{-4}\cdot \frac{\lambda(\tau)}{\lambda(\sigma)}\frac{\rho_{n}^{\frac12}\left(\frac{\lambda^{2}(\tau)}{\lambda^{2}(\sigma)}\xi\right)}{\rho^{\frac12}_{n}(\xi)}\cdot\frac{\left(\log\left(\frac{\lambda^{2}(\tau)}{\lambda^{2}(\sigma)}\xi\right)\right)^{i}}{\left(\frac{\lambda^{2}(\tau)}{\lambda^{2}(\sigma)}\xi\right)^{\frac12+\delta}}\cdot e^{\mp 2i\nu\sigma\frac{\lambda(\tau)}{\lambda(\sigma)}\xi^{\frac12}}\\
		&\hspace{3cm}\cdot \int_{\tau_{0}}^{\sigma}e^{\pm i\hbar^{-1}\rho\left(x_{\sigma_{1}\cdot\frac{\lambda(\tau)}{\lambda(\sigma_{1})}};\alpha\cdot\frac{\lambda(\tau)}{\lambda(\sigma_{1})},\hbar\right)}\cdot F^{(\pm)}\left(\sigma,\sigma_{1},\frac{\lambda^{2}(\tau)}{\lambda^{2}(\sigma_{1})}\xi\right)\,d\sigma_{1}\,d\sigma
	\end{split}
\end{align}
Changing the order of integral variables, the above two terms become
\begin{align}\label{good output 2a1prime}
	\begin{split}
	&\sum_{\pm}\frac{e^{\pm i\nu\tau\xi^{\frac12}}}{\hbar\xi^{\frac12}}\int_{\tau_{0}}^{\tau}e^{\pm i\hbar^{-1}\rho\left(x_{\sigma_{1}\cdot\frac{\lambda(\tau)}{\lambda(\sigma_{1})}};\alpha\cdot\frac{\lambda(\tau)}{\lambda(\sigma_{1})},\hbar\right)}\int_{\sigma_{1}}^{\tau}\left(\hbar^{2}\frac{\lambda^{2}(\tau)}{\lambda^{2}(\sigma)}\xi\right)^{-4}\\
	&\hspace{2cm}\cdot \frac{\lambda(\tau)}{\lambda(\sigma)}\frac{\rho_{n}^{\frac12}\left(\frac{\lambda^{2}(\tau)}{\lambda^{2}(\sigma)}\xi\right)}{\rho^{\frac12}_{n}(\xi)}\cdot\frac{\left(\log\left(\frac{\lambda^{2}(\tau)}{\lambda^{2}(\sigma)}\xi\right)\right)^{i}}{\left(\frac{\lambda^{2}(\tau)}{\lambda^{2}(\sigma)}\xi\right)^{\frac12+\delta}}\cdot F^{(\pm)}\left(\sigma,\sigma_{1},\frac{\lambda^{2}(\tau)}{\lambda^{2}(\sigma_{1})}\xi\right)\,d\sigma\,d\sigma_{1}
	\end{split}
\end{align}
and 
\begin{align}\label{good output 2a2prime}
	\begin{split}
	&\sum_{\pm}\frac{e^{\pm i\nu\tau\xi^{\frac12}}}{\hbar\xi^{\frac12}}\int_{\tau_{0}}^{\tau}e^{\pm i\hbar^{-1}\rho\left(x_{\sigma_{1}\cdot\frac{\lambda(\tau)}{\lambda(\sigma_{1})}};\alpha\cdot\frac{\lambda(\tau)}{\lambda(\sigma_{1})},\hbar\right)}\int_{\sigma_{1}}^{\tau}\left(\hbar^{2}\frac{\lambda^{2}(\tau)}{\lambda^{2}(\sigma)}\xi\right)^{-4}\\
	&\hspace{1cm}\cdot \frac{\lambda(\tau)}{\lambda(\sigma)}\frac{\rho_{n}^{\frac12}\left(\frac{\lambda^{2}(\tau)}{\lambda^{2}(\sigma)}\xi\right)}{\rho^{\frac12}_{n}(\xi)}\cdot\frac{\left(\log\left(\frac{\lambda^{2}(\tau)}{\lambda^{2}(\sigma)}\xi\right)\right)^{i}}{\left(\frac{\lambda^{2}(\tau)}{\lambda^{2}(\sigma)}\xi\right)^{\frac12+\delta}}\cdot e^{\mp 2i\nu\sigma\frac{\lambda(\tau)}{\lambda(\sigma)}\xi^{\frac12}}\cdot F^{(\pm)}\left(\sigma,\sigma_{1},\frac{\lambda^{2}(\tau)}{\lambda^{2}(\sigma_{1})}\xi\right)\,d\sigma\,d\sigma_{1}.
	\end{split}
\end{align}
For the contribution in \eqref{good output 2a1prime}, we set
\begin{align*}
	\tilde{F}^{(\pm)}\left(\tau,\sigma_{1},\xi\right):=&\int_{\sigma_{1}}^{\tau}\left(\frac{\lambda(\tau)}{\lambda(\sigma)}\right)^{-\frac{15+\delta}{2}}\cdot \frac{\rho_{n}^{\frac12}\left(\frac{\lambda^{2}(\tau)}{\lambda^{2}(\sigma)}\xi\right)}{\rho^{\frac12}_{n}(\xi)}\cdot\left(\log\left(\frac{\lambda^{2}(\tau)}{\lambda^{2}(\sigma)}\xi\right)\right)^{i}\cdot F^{(\pm)}\left(\sigma,\sigma_{1},\frac{\lambda^{2}(\tau)}{\lambda^{2}(\sigma_{1})}\xi\right)\,d\sigma
\end{align*}
and for the contribution in \eqref{good output 2a2prime}, we set
\begin{align*}
	\tilde{F}^{(\pm)}\left(\tau,\sigma_{1},\xi\right):=&\int_{\sigma_{1}}^{\tau}\left(\frac{\lambda(\tau)}{\lambda(\sigma)}\right)^{-\frac{15+\delta}{2}}\cdot \frac{\rho_{n}^{\frac12}\left(\frac{\lambda^{2}(\tau)}{\lambda^{2}(\sigma)}\xi\right)}{\rho^{\frac12}_{n}(\xi)}\cdot\left(\log\left(\frac{\lambda^{2}(\tau)}{\lambda^{2}(\sigma)}\xi\right)\right)^{i}\cdot e^{\mp 2i\nu\sigma\frac{\lambda(\tau)}{\lambda(\sigma)}\xi^{\frac12}}\cdot F^{(\pm)}\left(\sigma,\sigma_{1},\frac{\lambda^{2}(\tau)}{\lambda^{2}(\sigma_{1})}\xi\right)\,d\sigma,
\end{align*}
then the desired result follows.

The contribution from the 2nd term in the ``good source function" becomes:
\begin{align}\label{good source 2b}
	\begin{split}
		\yb^{b}_{2,\text{smooth}}\left(\sigma,\frac{\lambda^{2}(\tau)}{\lambda^{2}(\sigma)}\xi\right):=&\chi_{\hbar^{2}\frac{\lambda(\tau)^{2}}{\lambda(\sigma)^{2}}\xi\geq1}\hbar^{-1}\left(\hbar^{2}\frac{\lambda^{2}(\tau)}{\lambda^{2}(\sigma)}\xi\right)^{-4}\frac{1}{\left(\frac{\lambda^{2}(\tau)}{\lambda^{2}(\sigma)}\xi\right)^{\frac12+\delta}}\\
		&\cdot \int_{\tau_{0}}^{\sigma}e^{\pm i\left[\left(\nu\sigma\frac{\lambda(\tau)}{\lambda(\sigma)}-2\nu\sigma_{1}\frac{\lambda(\tau)}{\lambda(\sigma_{1})}\right)\xi^{\frac12}-\hbar^{-1}\rho\left(x_{\sigma_{1}\cdot\frac{\lambda(\tau)}{\lambda(\sigma_{1})}};\alpha\cdot\frac{\lambda(\tau)}{\lambda(\sigma_{1})},\hbar\right)\right]}\cdot \tilde{F}^{(\pm)}\left(\sigma,\sigma_{1},\frac{\lambda^{2}(\tau)}{\lambda^{2}(\sigma_{1})}\xi\right)\,d\sigma_{1}.
	\end{split}
\end{align}
We use an identical argument as we handle \eqref{good source 2a} to obtain the desired result. Here the only difference is that now the phase function is $e^{\pm i\left(\nu\tau-2\nu\sigma_{1}\frac{\lambda(\tau)}{\lambda(\sigma_{1})}\right)\xi^{\frac12}}$ instead of $e^{\pm i\nu\tau\xi^{\frac12}}$.

Finally the contribution from the 3rd term in the ``good source function" becomes
\begin{align}\label{good source 2c}
	\begin{split}
		\yb_{2,\text{smooth}}^{c}\left(\sigma,\frac{\lambda^{2}(\tau)}{\lambda^{2}(\sigma)}\xi\right):=&\chi_{\hbar^{2}\frac{\lambda(\tau)^{2}}{\lambda(\sigma)^{2}}\xi\geq1}\hbar^{-1}\left(\hbar^{2}\frac{\lambda^{2}(\tau)}{\lambda^{2}(\sigma)}\xi\right)^{-4}\frac{1}{\left(\frac{\lambda^{2}(\tau)}{\lambda^{2}(\sigma)}\xi\right)^{\frac12+\delta}}\\&\cdot \int_{0}^{\infty}\int_{\tau_{0}}^{\sigma}e^{\pm i\left[\nu\left(\frac{\lambda(\tau)}{\lambda(\sigma_{1})}x+\sigma\frac{\lambda(\tau)}{\lambda(\sigma)}\right)\xi^{\frac12}+\hbar^{-1}\rho\left(x_{\sigma_{1}\cdot\frac{\lambda(\tau)}{\lambda(\sigma_{1})}};\alpha\cdot\frac{\lambda(\tau)}{\lambda(\sigma_{1})},\hbar\right)\right]}\cdot G^{\pm}\left(\sigma,\sigma_{1},x,\frac{\lambda^{2}(\tau)}{\lambda^{2}(\sigma_{1})}\xi\right)\,d\sigma_{1}\,dx,
	\end{split}
\end{align}
and its contribution to the inhomogeneous parametrix is handled in a similar vein as that for the contribution from \eqref{good source 2a} and \eqref{good source 2b}.
\end{proof}

We have the following analogue of Proposition~\ref{prop:sourcetermswithtransferenceopngeq2}: 
\begin{proposition}\label{prop:sourcetermswithtransferenceopngeq2specialregular} Assume that $\overline{x}^{\hbar}(\tau, \xi)$ is the regular part of a good source function , and that $ \mathcal{F}^{\hbar}\left(F_\pm(n)\right)$ is an a regular good source term. Then referring to \eqref{eq:ngeq2Fourier2}, 
	\[
	\int_{\tau_0}^\tau U^{(n)}(\tau, \sigma,\xi)\cdot f_j\left(\sigma,\frac{\lambda^2(\tau)}{\lambda^2(\sigma)}\xi\right)\,d\sigma,\quad j = 1,2,\ldots,5,
	\]
	is a regular good function, and we have the bound 
	\[
	\left\|\int_{\tau_0}^\tau U^{(n)}(\tau, \sigma,\xi)\cdot f_j\left(\sigma,\frac{\lambda^2(\tau)}{\lambda^2(\sigma)}\xi\right)\,d\sigma\right\|_{\text{smooth}}\lesssim \left\| \overline{x}^{\hbar}\right\|_{\text{smooth}},
	\]
	with uniform implied constant. 
	In conjunction with  Prop.~\ref{prop:sourcetermswithtransferenceopngeq2}, we conclude that if $\overline{x}^{\hbar}(\tau, \xi)$ is good and $ \mathcal{F}^{\hbar}\left(F_\pm(n)\right)$ is a good source, then 
	\[
	\int_{\tau_0}^\tau U^{(n)}(\tau, \sigma,\xi)\cdot f_j\left(\sigma,\frac{\lambda^2(\tau)}{\lambda^2(\sigma)}\xi\right)\,d\sigma,\quad j = 1,2,\ldots,5,
	\]  
	is a good function, and we have 
	\[
	\left\|\int_{\tau_0}^\tau U^{(n)}(\tau, \sigma,\xi)\cdot f_j\left(\sigma,\frac{\lambda^2(\tau)}{\lambda^2(\sigma)}\xi\right)\,d\sigma\right\|_{\text{good}}\lesssim \left\| \overline{x}^{\hbar}\right\|_{\text{good}}. 
	\]
\end{proposition}
\begin{proof}
In view of Propositions \ref{prop: xh good para linear} and \ref{prop:transferenceonsingularngeq2}, we only need to consider the case when $\xb^{\hbar}(\tau,\xi)$ is of the form $\xb_{1,\text{smooth}}$ in Definition \ref{defi:goodfourierrepngeq2}. More precisely, we need to show that the transference operator maps a function of the form $\xb_{1,\text{smooth}}$ into a ``good source function" defined in Definition \ref{defi:goodfourierrepngeq2}. 
For this we have the following analogue of Prop.~\ref{prop:transferenceonsingularngeq2}:
\begin{lemma}\label{lem:transferenceongoodngeq2} Assume that $\xb(\tau, \xi) = \xb_{1,smooth}(\tau, \xi) $ is a good function of smooth and structured type as in Definition~\ref{defi:goodfourierrepngeq2}, at angular momentum $n, |n|\geq 2$. Then the function
\begin{align*}
\langle\xi^{\frac12}\rangle\cdot \tau^{-1}\cdot \left(\mathcal{K}_{\hbar}^{(0)}\xb\right)
\end{align*}
is a good source function in the sense of Definition \ref{defi:goodfourierrepngeq2}, as are the functions obtained by replacing $\mathcal{K}_{\hbar}^{(0)}$ by $\big[\mathcal{D}_{\tau}, \mathcal{K}_{\hbar}^{(0)}\big]$, $\big(\mathcal{K}_{\hbar}^{(0)}\big)^2$. Furthermore, for $\xb(\tau, \xi) $, the function 
\[
\langle\xi^{\frac12}\rangle\cdot \tau^{-1}\cdot \chi_{\xi\hbar^2\lesssim 1}\left(\mathcal{K}_{\hbar}^{(0)}\xb\right)
\]
is a good source function, and similarly when replacing $\mathcal{K}_{\hbar}^{(0)}$ by $\big[\mathcal{D}_{\tau}, \mathcal{K}_{\hbar}^{(0)}\big]$, $\big(\mathcal{K}_{\hbar}^{(0)}\big)^2$.
\end{lemma}

This is proven in analogy to Proposition~\ref{prop:transferenceonsingularngeq2}, taking advantage of Lemma~\ref{lem:transferenceonadmsingoffdiaglowfreq} to reduce to a regime where the diagonal derivative bounds of Prop.~\ref{prop: K operator} apply, as well as Lemma~\ref{lem:transferenceonadmsingoffdiaglowfreqverylowoutputfreq}. \\
The proof of Prop.~\ref{prop:sourcetermswithtransferenceopngeq2specialregular} is then completed by applying Prop.~\ref{prop:parametrixongoodsource}.

\end{proof}
\begin{remark}\label{rem:prop:sourcetermswithtransferenceopngeq2specialregular} It is easily seen that in the estimate of the proposition, one actually gains a smallness factor $\tau_0^{-1}$ for the contributions of $f_j,\,j = 3, 4, 5$. 

\end{remark}

\subsubsection{Solution of the Fourier wave equation \eqref{eq:ngeq2Fourier1}}

Given the functional setup from the preceding, we can now detail the (re)-iterative scheme which leads to a solution of our Fourier wave equation, in analogy to \cite{BKmemo, Kmemo, KMiao}. As usual the main difficulty stems from the term 
\[
2\frac{\lambda'(\tau)}{\lambda(\tau)}\mathcal{K}_{\hbar}^{(0)}\mathcal{D}_{\tau}\overline{x}^{\hbar},
\]
since the temporal weight $\frac{\lambda'(\tau)}{\lambda(\tau)}\simeq \tau^{-1}$ only furnishes enough decay to recover any previous decay assumptions upon application of the wave parametrix, and in particular we cannot force a smallness gain by simply picking the initial time $\tau_0$ large enough. This issue does not occur for the terms $f_j,\,j = 3,4,5$ in \eqref{eq:ngeq2Fourier2}, as follows by the following simple sharpening of the preceding proposition: 
\begin{lemma}\label{lem:sourcetermswithtransferenceopngeq2sharpening1} Let $\overline{x}^{\hbar}(\tau, \xi)$ be a good function. Then 
	\begin{align*}
		\left\|\int_{\tau_0}^\tau U^{(n)}(\tau, \sigma,\xi)\cdot f_j\left(\sigma,\frac{\lambda^2(\tau)}{\lambda^2(\sigma)}\xi\right)\,d\sigma\right\|_{\text{good}}\lesssim \tau_0^{-1}\left\| \overline{x}^{\hbar}\right\|_{\text{good}},\quad j = 3,4,5.
	\end{align*} 
	In particular, given any $\gamma>0$, there is $\tau_0 = \tau_0(\gamma)$ large enough such that 
	\begin{align*}
		\left\|\int_{\tau_0}^\tau U^{(n)}(\tau, \sigma,\xi)\cdot f_j\left(\sigma,\frac{\lambda^2(\tau)}{\lambda^2(\sigma)}\xi\right)\,d\sigma\right\|_{\text{good}}\leq \gamma\cdot \left\| \overline{x}^{\hbar}\right\|_{\text{good}},\quad j = 3,4,5.
	\end{align*} 
\end{lemma}
The proof is a direct consequence of the one for Prop.~\ref{prop:sourcetermswithtransferenceopngeq2specialregular}, Prop.~\ref{prop:sourcetermswithtransferenceopngeq2}. 
\\

In order to cope with the term $\frac{\lambda'(\tau)}{\lambda(\tau)}\mathcal{K}_{\hbar}^{(0)}\mathcal{D}_{\tau}\overline{x}^{\hbar}$, we use the same method as in \cite{Kmemo, BKmemo,KMiao}, namely manifold iteration, introduce the auxiliary composite operator 
\begin{equation}\label{eq:badexplicitpropagator}
	\Phi(f): = \int_{\tau_0}^{\tau}\frac{\lambda(\tau)}{\lambda(\sigma)}\cdot\frac{\rho^{\frac12}_n\left(\frac{\lambda^2(\tau)}{\lambda^2(\sigma)}\xi\right)}{\rho^{\frac12}_n(\xi)}\cdot \cos\left[\lambda(\tau)\xi^{\frac12}\cdot\int_{\sigma}^{\tau}\lambda^{-1}(u)\,du\right]\cdot \beta_{\nu}(\sigma)\cdot \mathcal{K}_{\hbar}^{(0)}\left(\mathcal{D}_{\tau}f\right)\left(\sigma,\frac{\lambda^2(\tau)}{\lambda^2(\sigma)}\xi\right)\,d\sigma,
\end{equation}
where we set $\beta_{\nu}(\sigma) = \frac{\lambda'(\sigma)}{\lambda(\sigma)}\simeq \sigma^{-1}$. 
Observe that this operator arises upon applying the operator $\mathcal{D}_{\tau}$ to the expression arising upon applying the parametrix to $\beta_{\nu}(\sigma)\cdot \mathcal{K}_{\hbar}^{(0)}\mathcal{D}_{\tau}f$.  Then the following key proposition gives the desired smallness gain upon manifold application of $\Phi$: 
\begin{proposition}\label{prop:keyreiteratengeq2} Let $\epsilon>0$ sufficiently small. Then there is $\gamma_*>0$ such that for any $k\geq 1$, there is $\tau_{0}=\tau_0(\epsilon, k)$ large enough such that 
	\begin{align*}
		\left\|\beta_{\nu}(\tau)\mathcal{K}_{\hbar}^{(0)}\Phi^k f\right\|_{\text{goodsource}}\leq \epsilon^{\gamma k}\cdot e^{\epsilon^{-2}}\cdot \left\| f\right\|_{\text{good}}. 
	\end{align*}
\end{proposition}
\begin{proof}
	The idea is to decompose the transference operator $\mathcal{K}_{\hbar}^{(0)}$ into a number of pieces, and most pivotally into a `diagonal' and `off diagonal' piece 
	\[
	\mathcal{K}_{\hbar}^{(0)} = \mathcal{K}_{\hbar,d}^{(0)} + \mathcal{K}_{\hbar,nd}^{(0)}, 
	\]
	which arise as follows: fixing the $\epsilon>0$, $k$, we pick $\ell$ sufficiently large (specifically $\ell\simeq k$) and define $\mathcal{K}_{\hbar,d}^{(0)}$ by including a smooth cutoff $\chi_{|1-\frac{\xi}{\eta}|\leq \frac{1}{\ell}}$ into the kernel of 
	$\mathcal{K}_{\hbar}^{(0)}$, whence the off-diagonal part $\mathcal{K}_{\hbar,nd}^{(0)}$ is defined by including a cutoff $\chi_{|1-\frac{\xi}{\eta}|>\frac{1}{\ell}}$. Call $\Phi^{(d)}, \Phi^{(nd)}$ the operators arising by these changes.  Then the idea is to reduce things to the diagonal part, by observing that in any strings of $\Phi$s the presence of a single $ \mathcal{K}_{\hbar,nd}^{(0)}$ results in a smallness gain by picking $\tau_0$ large. This follows from the following 
	\begin{lemma}\label{lem:offdiagonaltransferencegainngeq2} We have the bound
		\begin{align*}
			\left\|\beta_{\nu}(\tau)\mathcal{K}_{\hbar}^{(0)} \Phi^{(nd)}\beta_\nu(\tau) \Phi f\right\|_{\text{goodsource}}\lesssim \tau_0^{-\gamma_1}\left\| f\right\|_{\text{good}}.
		\end{align*}
		for some $\gamma_1>0$, where the implied constant is independent of the angular momentum $n$ and the time $\tau_0$, but will depend on $\ell$. 
	\end{lemma} 
\begin{proof}
	(lemma) The idea is to perform integration by parts in the time variable for $\Phi^{(nd)}$. Specifically, write 
	\begin{align*}
		&\left( \Phi^{(nd)}\beta_\nu(\tau)\Phi\right)(g)\\&
		= c\sum_{\pm,\pm}\int_{\tau_0}^{\tau}\tilde{\rho}_n(\tau, \sigma,\xi)\cdot e^{\pm i\left(\nu\tau - \nu\sigma\frac{\lambda(\tau)}{\lambda(\sigma)}\right)\xi^{\frac12}}\cdot\int_0^\infty \frac{F_n^{(nd)}\left(\frac{\lambda^2(\tau)}{\lambda^2(\sigma)}\xi, \eta\right)\cdot\rho_n(\eta)}{\frac{\lambda^2(\tau)}{\lambda^2(\sigma)}\xi - \eta}\\&\hspace{4cm}\cdot \int_{\tau_0}^{\sigma}\tilde{\rho}_n(\sigma, \sigma_1,\eta)\cdot e^{\pm i\left(\nu\sigma - \nu\sigma_1\frac{\lambda(\sigma)}{\lambda(\sigma_1)}\right)\eta^{\frac12}}\cdot \tilde{g}\left(\sigma_1,\frac{\lambda^2(\sigma)}{\lambda^2(\sigma_1)}\eta\right)\,d\sigma_1\, d\eta \,d\sigma
	\end{align*}
where $c$ is a suitable constant, and we have set $\tilde{g} = \mathcal{K}^{(0)}_{\hbar}g$, $\tilde{\rho}_n(\tau, \sigma,\xi) = \beta_{\nu}(\sigma)\cdot \frac{\lambda(\tau)}{\lambda(\sigma)}\cdot \frac{\rho_n^{\frac12}\left(\frac{\lambda^2(\tau)}{\lambda^2(\sigma)}\xi\right)}{\rho_n^{\frac12}(\xi)}$. Also, $F^{(nd)}$ means inclusion of the cutoff $\chi_{\left|1-\frac{\xi}{\eta}\right|>\frac{1}{\ell}}$, which means here inclusion of $\chi_{\left|1-\frac{\frac{\lambda^2(\tau)}{\lambda^2(\sigma)}\xi}{\eta}\right|>\frac{1}{\ell}}$ due to the change of scale. To begin with, we note that if we include a further smooth cutoff $\chi_{\left|\sigma - \sigma_1\right|\lesssim \tau_0^{-\frac12}\sigma}$ in the second time integral, we can easily force smallness due to the restriction of the integration interval. This means we can include a smooth cutoff $\chi_{\left|\sigma - \sigma_1\right|\gtrsim\tau_0^{-\frac12}\sigma}$. Furthermore, due to the asymptotic bounds for 
$
F_n(\xi,\eta), 
$
we may also include a cutoff $\chi_{\frac{\lambda^2(\tau)}{\lambda^2(\sigma)}\xi\gtrsim 1}$, since else smallness  can again be forced by $\frac{\lambda^2(\tau)}{\lambda^2(\sigma)}\xi\ll 1,\Rightarrow\hbar^{2}\frac{\lambda^2(\tau)}{\lambda^2(\sigma)}\xi\ll 1$. The extra smallness gain also occurs if $\frac{\lambda(\tau)^{2}}{\lambda(\sigma)^{2}}\xi\gg\eta$ or $\frac{\lambda(\tau)^{2}}{\lambda(\sigma)^{2}}\xi\ll\eta$, in view of the bounds in Proposition \ref{prop: K operator}. Therefore without loss of generality, we assume $\frac{\lambda(\tau)^{2}}{\lambda(\sigma)^{2}}\xi\simeq\eta$.

Switching to the new variable $\tilde{\eta}$ defined via 
\[
\eta = \frac{\lambda^2(\tau)}{\lambda^2(\sigma)}\tilde{\eta}, 
\]
we replace the function $\tilde{g}\left(\sigma_1,\frac{\lambda^2(\sigma)}{\lambda^2(\sigma_1)}\eta\right)$ by $\tilde{g}\left(\sigma_1,\frac{\lambda^2(\tau)}{\lambda^2(\sigma_1)}\tilde{\eta}\right)$, and the integration kernel 
\[
\frac{F_n\left(\frac{\lambda^2(\tau)}{\lambda^2(\sigma)}\xi, \eta\right)\cdot\rho_n(\eta)}{\frac{\lambda^2(\tau)}{\lambda^2(\sigma)}\xi - \eta}
\]
is replaced by the more symmetrical 
\[
\frac{F_n\left(\frac{\lambda^2(\tau)}{\lambda^2(\sigma)}\xi, \frac{\lambda^2(\tau)}{\lambda^2(\sigma)}\tilde{\eta}\right)\cdot\rho_n\left(\frac{\lambda^2(\tau)}{\lambda^2(\sigma)}\tilde{\eta}\right)}{\xi - \tilde{\eta}}. 
\]
Finally, the oscillatory phase becomes 
\[
e^{\pm i\left(\nu\tau - \nu\sigma\frac{\lambda(\tau)}{\lambda(\sigma)}\right)\xi^{\frac12}}\cdot e^{\pm i\left(\nu\sigma\frac{\lambda(\tau)}{\lambda(\sigma)} - \nu\sigma_1\frac{\lambda(\tau)}{\lambda(\sigma_1)}\right)\tilde{\eta}^{\frac12}}
= e^{\pm i \nu\tau\xi^{\frac12}\mp i \nu\sigma_1\frac{\lambda(\tau)}{\lambda(\sigma_1)}\tilde{\eta}^{\frac12}}\cdot e^{ i\nu\sigma\frac{\lambda(\tau)}{\lambda(\sigma)}\left(\mp\xi^{\frac12}\pm \tilde{\eta}^{\frac12}\right)}
\]
For later reference, we observe that in case of anti-alignment of the pass, we get the $\sigma$-dependent phase $e^{\pm i\nu\sigma\frac{\lambda(\tau)}{\lambda(\sigma)}\cdot \left(\xi^{\frac12}+\tilde{\eta}^{\frac12}\right)}$, in which case we needn't even take advantage of the off-diagonal condition. In the worst case of destructive resonance, we get the $\sigma$-dependent phase $e^{\pm i\nu\sigma\frac{\lambda(\tau)}{\lambda(\sigma)}\left(\xi^{\frac12} - \tilde{\eta}^{\frac12}\right)}$, where our assumption implies that $\left|\frac{\lambda(\tau)}{\lambda(\sigma)}\xi^{\frac12} - \frac{\lambda(\tau)}{\lambda(\sigma)}\tilde{\eta}^{\frac12}\right|> \frac{1}{\ell}\cdot\frac{\lambda(\tau)}{\lambda(\sigma)}\xi^{\frac12}$. Then we change the order of temporal integration in the above triple integral, so that the $\sigma$-integral comes first, and perform integration by parts with respect to $\sigma$. We have
\begin{align*}
	\partial_{\sigma}\left(e^{\pm i\nu\sigma\frac{\lambda(\tau)}{\lambda(\sigma)}\left(\xi^{\frac12}-\tilde{\eta}^{\frac12}\right)}\right)=\mp i\frac{\lambda(\tau)}{\lambda(\sigma)}\left(\xi^{\frac12}-\tilde{\eta}^{\frac12}\right)\cdot e^{\pm i\nu\sigma\frac{\lambda(\tau)}{\lambda(\sigma)}\left(\xi^{\frac12}-\tilde{\eta}^{\frac12}\right)}
\end{align*}
This produces an additional factor 
\[
\sim \ell\cdot \frac{1}{\sigma}\cdot \frac{1}{\frac{\lambda(\tau)}{\lambda(\sigma)}\xi^{\frac12}}\lesssim \frac{1}{\sigma}, 
\]
where we have taken into account the additional localization $\frac{\lambda^2(\tau)}{\lambda^2(\sigma)}\xi\gtrsim 1$, the derivative bounds of the transference kernel $F_n(\xi,\eta)$, as well as the fact that the additional time cutoff above prevents boundary terms. In fact in the regime $|\sigma-\sigma_{1}|\gtrsim\tau_{0}^{-\frac12}\sigma$, we have $\frac{\sigma_{1}}{1-\tau_{0}^{-\frac12}}\lesssim \sigma\leq \tau$. Therefore the boundary term at $\sigma=\sigma_{1}$ vanishes. For the boundary term at $\sigma=\tau$, we get something like (up to multiple of a constant)
\begin{align*}
\int_{\tau_{0}}^{\tau}\beta_{\nu}(\tau)\beta_{\nu}(\sigma_{1})\cdot\frac{\lambda(\tau)}{\lambda(\sigma_{1})}\int_{0}^{\infty}\frac{\rho_{n}^{\frac12}\left(\frac{\lambda^{2}(\tau)}{\lambda^{2}(\sigma_{1})}\eta\right)}{\rho_{n}^{\frac12}(\eta)}\frac{F_{n}^{(nd)}(\xi,\eta)\cdot \rho_{n}(\eta)}{\xi-\eta}\cdot e^{\pm i\nu\left(\tau-\sigma_{1}\frac{\lambda(\tau)}{\lambda(\sigma_{1})}\right)\eta}\cdot \tilde{g}\left(\sigma_{1},\frac{\lambda^{2}(\tau)}{\lambda^{2}(\sigma_{1})}\eta\right)\,d\eta\,d\sigma_{1}.
\end{align*}
Here the growth $\frac{\lambda(\tau)}{\lambda(\sigma_{1})}$ is absorbed by the decay of $\tilde{g}(\sigma_{1},\cdot)$ and the extra smallness comes from the factor $\beta_{\nu}(\tau)\beta_{\nu}(\sigma_{1})$.
For the integral in $\sigma$, the smallness gain follows since $\sigma\geq \tau_0$, and the fact that one can place the resulting function in the good source space of functions by re-arranging the phases as before and invoking Proposition~\ref{prop:sourcetermswithtransferenceopngeq2specialregular}. 
\end{proof}
Note that by the preceding proof, we may assume for the proof of Proposition~\ref{prop:keyreiteratengeq2} that the oscillatory phases (due to the factors $\cos\left[\lambda(\tau)\xi^{\frac12}\int_{\sigma}^{\tau}\lambda^{-1}(u)\,du\right]$, written as a sum of exponentials), all have the same sign, since otherwise, two adjacent opposite signs allow us to utilize the gain from the preceding lemma to gain smallness. This implies (following the proof of Proposition~\ref{prop:parametrixonadmissiblesource}) that starting with one of the three expressions in Definition~\ref{defi:goodfourierrepngeq2}, or alternatively one of the expressions in $\xb_{in}, \xb_{out}$ in Definition~\ref{defi:xsingulartermsngeq2adm}, we always reproduce the same kind of expression, with coefficient functions which are given by suitable iterated integrals. Moreover, the preceding lemma allows us to replace $\Phi$ by $\Phi^{(d)}$ throughout.
One then follows the proof of Proposition 11.2 in \cite{KMiao} to reduce the problem to integration over simplices which forces smallness via a combinatorial argument. Here we show how to gain smallness from the contribution of $\Phi^{(d)}$ independent of $\hbar$. Let us denote by $\calK_{\hbar}^{(0),d}$ the transference operator with the restriction $\left|\frac{\xi}{\eta}-1\right|\leq \frac{1}{\ell}$, and decompose $\calK_{\hbar}^{(0),d}$ as
\begin{align*}
	\calK_{\hbar}^{(0),d}=\calK_{\hbar,\epsilon}^{(0),d,1}+\calK_{\hbar,\epsilon}^{(0),d,2}+\calK_{\hbar,\epsilon}^{(0),d,3}
\end{align*}
We use $K_{\hbar}^{(0),d}, K_{\hbar,\epsilon}^{(0),d,1}, K_{\hbar,\epsilon}^{(0),d,2}, K_{\hbar,\epsilon}^{(0),d,3}$ to denote the kernels of $\calK_{\hbar}^{(0),d}, \calK_{\hbar,\epsilon}^{(0),d,1}, \calK_{\hbar,\epsilon}^{(0),d,2}, \calK_{\hbar,\epsilon}^{(0),d,3}$ respectively, where
\begin{align*}
	K_{\hbar,\epsilon}^{(0),d,1}:=\chi_{\hbar^{2}\xi<\epsilon}\cdot K_{\hbar}^{(0),d},\quad K_{\hbar,\epsilon}^{(0),d,3}:=\chi_{\hbar^{2}\xi>\epsilon^{-1}}\cdot K_{\hbar}^{(0),d},\quad K_{\hbar,\epsilon}^{(0),d,2}:=\cdot K_{\hbar}^{(0),d}-
\end{align*}
According to the kernel bounds in Proposition \ref{prop: F hbar boundedness}, there is a smallness gain from the operators $\calK_{\hbar,\epsilon}^{(0),d,1}, \calK_{\hbar,\epsilon}^{(0),d,3}$. On the other hand, define an operator $\calD_{\tau}U$ as
\begin{align}\label{def DtauU}
	\left(\calD_{\tau}U h\right)(\tau,\xi):=\int_{\tau_{0}}^{\tau}\frac{\lambda(\tau)}{\lambda(\sigma)}\frac{\rho_{n}^{\frac12}\left(\frac{\lambda^{2}(\tau)}{\lambda^{2}(\sigma)}\xi\right)}{\rho_{n}^{\frac12}(\xi)}\cos\left[\lambda(\tau)\xi^{\frac12}\int_{\tau}^{\sigma}\lambda^{-1}(u)\,du\right]h\left(\sigma,\frac{\lambda^{2}(\tau)}{\lambda^{2}(\sigma)}\xi\right)\,d\sigma.
\end{align}
Then we have the following vanishing relations:
\begin{align}\label{vanishing Kd}
	\begin{split}
		\calK_{\hbar,\epsilon}^{(0),d,3}\calD_{\tau}U\calK_{\hbar,\left(1+\frac{1}{\ell}\right)\epsilon}^{(0),d,2}=0,\quad \calK_{\hbar,\left(1+\frac{1}{\ell}\right)\epsilon}^{(0),d,2}\calD_{\tau}U\calK_{\hbar,\epsilon}^{(0),d,1}=0.
	\end{split}
\end{align}
This can be seen as follows. Let us consider the first relation in \eqref{vanishing Kd} for example. Regarding to the operator $\calK_{\hbar,\epsilon}^{(0),d,3}$, the output frequency $\xi$ and the input frequency $\eta$ satisfy
\begin{align*}
	\hbar^{2}\xi>\epsilon^{-1},\quad \left(1+\frac{1}{\ell}\right)\epsilon\leq\frac{\lambda^{2}(\tau)}{\lambda^{2}(\sigma)}\hbar^{2}\eta\leq\left(1+\frac{1}{\ell}\right)^{-1}\epsilon^{-1}.
\end{align*}
On the other hand, the diagonal restriction implies
\begin{align*}
	\left|\frac{\hbar^{2}\xi}{\frac{\lambda^{2}(\tau)}{\lambda^{2}(\sigma)}\hbar^{2}\eta}-1\right|\leq \frac{1}{\ell},\quad \Rightarrow\quad \hbar^{2}\xi\leq\left(1+\frac{1}{\ell}\right)\frac{\lambda^{2}(\tau)}{\lambda^{2}(\sigma)}\hbar^{2}\eta\leq \epsilon^{-1},
\end{align*}
which violates the restriction $\hbar^{2}\xi>\epsilon^{-1}$. Therefore the first relation in \eqref{vanishing Kd} holds. The second one can be analyzed similarly. Using this argument we also infer relation
\begin{align}\label{vanishing Kd alt}
	\begin{split}
		\calK_{\hbar,\epsilon}^{d,3}\calD_{\tau}U\calK_{\hbar,\left(1+\frac{1}{\ell}\right)\epsilon}^{d,1}=0.
	\end{split}
\end{align}
Following the argument in \cite{Kmemo,KMiao}, we expand
\begin{align*}
	\left(\beta_{\nu}(\tau)\calK_{\hbar}^{(0),d}\calD_{\tau}U\right)^{n}=&\left(\beta_{\nu}(\tau)\calK_{\hbar,\epsilon}^{(0),d,2}\calD_{\tau}U\right)^{n}\\
	&+\sum_{k=1}^{n-1}\left(\beta_{\nu}(\tau)\calK_{\hbar,\epsilon}^{(0),d,2}\calD_{\tau}U\right)^{k}\left(\beta_{\nu}(\tau)\calK_{\hbar,\epsilon}^{(0),d,1}\calD_{\tau}U\right)\left(\beta_{\nu}(\tau)\calK_{\hbar}^{(0),d}\calD_{\tau}U\right)^{n-k-1}\\
	&+\sum_{k=1}^{n-1}\left(\beta_{\nu}(\tau)\calK_{\hbar,\epsilon}^{(0),d,2}\calD_{\tau}U\right)^{k}\left(\beta_{\nu}(\tau)\calK_{\hbar,\epsilon}^{(0),d,3}\calD_{\tau}U\right)\left(\beta_{\nu}(\tau)\calK_{\hbar}^{(0),d}\calD_{\tau}U\right)^{n-k-1}\\
	&+\left(\beta_{\nu}(\tau)\calK_{\hbar,\epsilon}^{(0),d,1}\calD_{\tau}U\right)\left(\beta_{\nu}(\tau)\calK_{\hbar}^{(0),d}\calD_{\tau}U\right)^{n-1}\\
	&+\left(\beta_{\nu}(\tau)\calK_{\hbar,\epsilon}^{(0),d,3}\calD_{\tau}U\right)\left(\beta_{\nu}(\tau)\calK_{\hbar}^{(0),d}\calD_{\tau}U\right)^{n-1}\\
	=:&A+B+C+D+E.
\end{align*}
For $C$, we have, by \eqref{vanishing Kd} and \eqref{vanishing Kd alt},
\begin{align*}
	\left(\beta_{\nu}(\tau)\calK_{\hbar,\epsilon}^{(0),d,3}\calD_{\tau}U\right)\left(\beta_{\nu}(\tau)\calK_{\hbar}^{(0),d}\calD_{\tau}U\right)^{n-k-1}
	=\left(\beta_{\nu}(\tau)\calK_{\hbar,\epsilon}^{(0),d,3}\calD_{\tau}U\right)\left(\beta_{\nu}(\tau)\calK_{\hbar,4\epsilon}^{(0),d,3}\calD_{\tau}U\right)^{n-k-1},
\end{align*}
which implies
\begin{align*}
	C=\sum_{k=1}^{n-1}\left(\beta_{\nu}(\tau)\calK_{\hbar,\epsilon}^{d,2}\calD_{\tau}U\right)^{k}\left(\beta_{\nu}(\tau)\calK_{\hbar,\epsilon}^{(0),d,3}\calD_{\tau}U\right)\left(\beta_{\nu}(\tau)\calK_{\hbar,4\epsilon}^{(0),d,3}\calD_{\tau}U\right)^{n-k-1}.
\end{align*}
For $B$, we have
\begin{align*}
	B=&\sum_{k=1}^{n-1}\left(\beta_{\nu}(\tau)\calK_{\hbar,\epsilon}^{(0),d,2}\calD_{\tau}U\right)^{k}\left(\beta_{\nu}(\tau)\calK_{\hbar,\epsilon}^{(0),d,1}\calD_{\tau}U\right)\left(\beta_{\nu}(\tau)\calK_{\hbar}^{(0),d}\calD_{\tau}U\right)^{n-k-1}\\
	=&\sum_{k=1}^{n-1}\left(\beta_{\nu}(\tau)\calK_{\hbar,\epsilon}^{(0),d,2}\calD_{\tau}U\right)^{k}\left(\beta_{\nu}(\tau)\calK_{\hbar,\epsilon}^{(0),d,1}\calD_{\tau}U\right)\left(\beta_{\nu}(\tau)\calK_{\hbar,\frac{\epsilon}{4}}^{(0),d,2}\calD_{\tau}U\right)^{n-k-1}\\
	&+\sum_{1\leq k\leq n-1,j\leq n-k-2}\left(\beta_{\nu}(\tau)\calK_{\hbar,\epsilon}^{(0),d,2}\calD_{\tau}U\right)^{k}\left(\beta_{\nu}(\tau)\calK_{\hbar,\epsilon}^{(0),d,1}\calD_{\tau}U\right)\left(\beta_{\nu}(\tau)\calK_{\hbar,\frac{\epsilon}{4}}^{(0),d,2}\calD_{\tau}U\right)^{j}\\
	&\cdot \left(\beta_{\nu}(\tau)\calK_{\hbar,\frac{\epsilon}{4}}^{(0),d,3}\calD_{\tau}U\right)\left(\beta_{\nu}(\tau)\calK_{\hbar,\epsilon}^{(0),d,3}\calD_{\tau}U\right)^{n-k-1}.
\end{align*}
Here we have used the first relation in \eqref{vanishing Kd} and the vanishing property
\begin{align}\label{vanishing Kd alt 1}
	\begin{split}
		\left(\beta_{\nu}(\tau)\calK_{\hbar,\epsilon}^{(0),d,2}\calD_{\tau}U\right)\left(\beta_{\nu}(\tau)\calK_{\hbar,\epsilon}^{(0),d,1}\calD_{\tau}U\right)\left(\beta_{\nu}(\tau)\calK_{\hbar,\frac{\epsilon}{4}}^{(0),d,1}\calD_{\tau}U\right)=0.
	\end{split}
\end{align}
\eqref{vanishing Kd alt 1} can be seen as follows. Let $(\tau,\xi), (\sigma,\eta), (\sigma_{1},\eta_{1})$ be the output variables. Then we have
\begin{align*}
	\hbar^{2}\eta_{1}\leq \frac{\lambda^{2}(\sigma)}{\lambda^{2}(\sigma_{1})}\hbar^{2}\eta_{1}<\frac{\epsilon}{4}.
\end{align*}
The diagonal restriction for the operator $\calK_{\hbar,\epsilon}^{(0),d,1}$ implies
\begin{align*}
	\hbar^{2}\eta\leq \frac{\lambda^{2}(\tau)}{\lambda^{2}(\sigma)}\hbar^{2}\eta\leq \left(1+\frac{1}{\ell}\right)\hbar^{2}\eta_{1}<\left(1+\frac{1}{\ell}\right)\cdot\frac{\epsilon}{4}.
\end{align*}
Similarly the diagonal restriction for the operator $\calK_{\hbar,\epsilon}^{(0),d,2}$ implies
\begin{align*}
	\hbar^{2}\xi\leq \left(1+\frac{1}{\ell}\right)\hbar^{2}\eta<\left(1+\frac{1}{\ell}\right)^{2}\frac{\epsilon}{4}.
\end{align*}
But this violates the restriction $\hbar^{2}\xi\geq \epsilon$, which proves \eqref{vanishing Kd alt 1}.

For $D$ we write it as
\begin{align*}
	D=&\left(\beta_{\nu}(\tau)\calK_{\hbar,\epsilon}^{(0),d,1}\calD_{\tau}U\right)\left(\beta_{\nu}(\tau)\calK_{\hbar}^{(0),d}\calD_{\tau}U\right)^{n-1}\\
	=&\sum_{j=1}^{n}\left(\beta_{\nu}(\tau)\calK_{\hbar,\epsilon}^{(0),d,1}\calD_{\tau}U\right)^{j}\left[A^{(n-j)}+B^{(n-j)}+C^{(n-j)}+E^{(n-j)}\right],
\end{align*}
where the superscript indicates that these terms are defined just as in $A,B,C$ and $E$ but with $n$ replaced by $n-j$.

Since the operators $\left(\beta_{\nu}(\tau)\calK_{\hbar,\epsilon}^{(0),d,1}\calD_{\tau}U\right)$ and $\left(\beta_{\nu}(\tau)\calK_{\hbar,\epsilon}^{(0),d,3}\calD_{\tau}U\right)$ produce extra smallness, we have reduced the problem of bounding $(\beta_{\nu}(\tau)\calK_{\hbar}^{(0)}\calD_{\tau}U)^{n}$ to the problem of bounding $(\beta_{\nu}(\tau)\calK_{\hbar,\epsilon}^{(0),d,2}\calD_{\tau}U)^{n}$. Given a function $f(\tau,\xi)$, we consider
$\left(\beta_{\nu}(\tau)\calK_{\hbar,\epsilon}^{(0),d,2}\calD_{\tau}U\right)^{j}\left(\beta_{\nu}(\tau)\calK_{\hbar,\epsilon}^{(0),d,2}f\right)$:
\begin{align}\label{reiteration main}
	\begin{split}
	&\left(\beta_{\nu}(\tau)\calK_{\hbar,\epsilon}^{(0),d,2}\calD_{\tau}U\right)^{j}\left(\beta_{\nu}(\tau)\calK_{\hbar,\epsilon}^{(0),d,2}f\right)\\
	=&\beta_{\nu}(\tau)\calK_{\hbar,\epsilon}^{(0),d,2}\int_{\tau_{0}}^{\tau}\int_{0}^{\infty}\,d\sigma_{1}\,d\eta_{1}\,\frac{\lambda(\tau)}{\lambda(\sigma_{1})}\frac{\rho_{n}^{\frac12}\left(\frac{\lambda^{2}(\tau)}{\lambda^{2}(\sigma_{1})}\xi\right)}{\rho_{n}^{\frac12}(\xi)}\beta_{\nu}(\sigma_{1})\cos\left[\lambda(\tau)\xi^{\frac12}\int_{\tau}^{\sigma_{1}}\lambda^{-1}(u)\,du\right]\\
	&\cdot K_{\hbar,\epsilon}^{(0),d,2}\left(\frac{\lambda^{2}(\tau)}{\lambda^{2}(\sigma_{1})}\xi,\eta_{1}\right)\int_{\tau_{0}}^{\sigma_{1}}\int_{0}^{\infty}\,d\sigma_{2}\,d\eta_{2}\,\frac{\lambda(\sigma_{1})}{\lambda(\sigma_{2})}\frac{\rho_{n}^{\frac12}\left(\frac{\lambda^{2}(\sigma_{1})}{\lambda^{2}(\sigma_{2})}\eta_{1}\right)}{\rho_{n}^{\frac12}(\eta_{1})}\beta_{\nu}(\sigma_{2})\cos\left[\lambda(\sigma_{1})\eta_{1}^{\frac12}\int_{\sigma_{1}}^{\sigma_{2}}\lambda^{-1}(u)\,du\right]\\
	&...\\
	&\cdot K_{\hbar,\epsilon}^{(0),d,2}\left(\frac{\lambda^{2}(\sigma_{j-2})}{\lambda^{2}(\sigma_{j-1})}\eta_{j-2},\eta_{j-1}\right)\int_{\tau_{0}}^{\sigma_{j-1}}\,d\sigma_{j}\frac{\lambda(\sigma_{j-1})}{\lambda(\sigma_{j})}\frac{\rho_{n}^{\frac12}\left(\frac{\lambda^{2}(\sigma_{j-1})}{\lambda^{2}(\sigma_{j})}\eta_{j-1}\right)}{\rho_{n}^{\frac12}(\eta_{j-1})}\beta_{\nu}(\sigma_{j})\\
	&\cdot\cos\left[\lambda(\sigma_{j-1})\eta_{j-1}^{\frac12}\int_{\sigma_{j-1}}^{\sigma_{j}}\lambda^{-1}(u)\,du\right]\left(\beta_{\nu}(\sigma_{j})\calK_{\hbar,\epsilon}^{(0),d,2}f\right)\left(\sigma_{j},\frac{\lambda^{2}(\sigma_{j-1})}{\lambda^{2}(\sigma_{j})}\eta_{j-1}\right).
	\end{split}
\end{align}
The diagonal restriction gives
\begin{align*}
	\left|\frac{\lambda^{2}(\sigma_{k})\eta_{k}}{\lambda^{2}(\sigma_{k+1})\eta_{k+1}}-1\right|\leq \frac{1}{\ell},\quad 1\leq k\leq j-1<n,
\end{align*}
which implies
\begin{align*}
	\frac{\lambda^{2}(\tau)\xi}{\lambda^{2}(\sigma_{1})}\leq \eta_{1}\left(1+\frac{1}{\ell}\right),\quad \frac{\lambda^{2}(\sigma_{1})\eta_{1}}{\lambda^{2}(\sigma_{2})}\leq \eta_{2}\left(1+\frac{1}{\ell}\right),\quad ...\, ,\quad \frac{\lambda^{2}(\sigma_{k-1})\eta_{k-1}}{\lambda^{2}(\sigma_{k})}\leq \eta_{k}\left(1+\frac{1}{\ell}\right),\quad 1\leq k\leq j-1.
\end{align*}
Multiplying all above inequalities, we have
\begin{align*}
	\frac{\lambda^{2}(\tau)}{\lambda^{2}(\sigma_{k})}\xi\eta_{1}...\eta_{k-1}\leq \eta_{1}...\eta_{k-1}\eta_{k}\left(1+\frac{1}{\ell}\right)^{k},\quad 1\leq k\leq j -1,
\end{align*} 
which finally implies
\begin{align}\label{compare xi etak}
	\frac{1}{\left(1+\frac{1}{\ell}\right)^{k}}\lambda^{2}(\tau)\xi\leq \lambda^{2}(\sigma_{k})\eta_{k},\quad 1\leq k\leq j-1.
\end{align}
Now we claim that on the support of the full expression \eqref{reiteration main}, we have 
\begin{align}\label{xi etak epsilon}
\hbar^{2}\xi\gtrsim\epsilon,\quad \hbar^{2}\eta_{k}<\epsilon^{-1}
\end{align}
 The latter follows from
\begin{align*}
	\hbar^{2}\eta_{k}\leq \frac{\lambda^{2}(\sigma_{k})}{\lambda^{2}(\sigma_{k+1})}\hbar^{2}\eta_{k}<\epsilon^{-1}.
\end{align*}
For the former one, we note that $\xi$ is the input variable of another application of the operator $\calK_{\hbar,\epsilon}^{(0),d,2}$. Let $\tilde{\xi}$ be the corresponding output variable, then we have, due to the diagonal restriction,
\begin{align*}
	\epsilon\leq \hbar^{2}\tilde{\xi}\leq \left(1+\frac{1}{\ell}\right)\hbar^{2}\xi,\quad \Rightarrow\quad \hbar^{2}\xi\gtrsim\epsilon.
\end{align*}
Plugging \eqref{xi etak epsilon} back to \eqref{compare xi etak}, we obtain
\begin{align*}
\sigma_{k}^{1+\frac{1}{\nu}}>\frac{1}{\left(1+\frac{1}{\ell}\right)^{\frac{k}{2}}}\tau^{1+\frac{1}{\nu}}\epsilon,\quad \Rightarrow\quad \left(\frac{\sigma_{k}}{\tau}\right)^{1+\frac{1}{\nu}}>\frac{\epsilon}{\left(1+\frac{1}{\ell}\right)^{\frac{k}{2}}}.
\end{align*}
For sufficiently large $\ell$ and $k$, $\frac{1}{\left(1+\frac{1}{\ell}\right)^{\frac{k}{2}}}$ has a positive lower bound. Therefore for fixed but sufficiently small $\epsilon$, we have
\begin{align*}
	\frac{1}{\left(1+\frac{1}{\ell}\right)^{\frac{k}{2}}}>\epsilon^{\frac{1}{\nu}},\quad \Rightarrow\quad \left(\frac{\sigma_{k}}{\tau}\right)^{1+\frac{1}{\nu}}>\epsilon^{1+\frac{1}{\nu}},\quad \Rightarrow\quad \sigma_{k}>\epsilon\tau,\quad 1\leq k\leq j-1.
\end{align*}
With all above preparation, we are ready  to show how to gain smallness in \eqref{reiteration main} for fixed $\epsilon>0$, provided $j\in\bbN$ sufficiently large. In view of Propositions \ref{prop: F hbar boundedness}, \ref{prop:transferenceonsingularngeq2}, and the proof of Proposition \ref{prop:sourcetermswithtransferenceopngeq2specialregular}, we know that the transference operators preserves the structures of a  ``good"  function. So we only need to take care of the temporal integrals in \eqref{reiteration main}. We have
\begin{align*}
	&\left\|\left(\beta_{\nu}(\tau)\calK_{\hbar,\epsilon}^{(0),d,2}\calD_{\tau}U\right)^{j}\left(\beta_{\nu}(\tau)\calK_{\hbar,\epsilon}^{(0),d,2}f\right)\right\|_{good}\\\lesssim &\beta_{\nu}(\tau)\int_{\epsilon\tau}^{\tau}\beta_{\nu}(\sigma_{1})\,\int_{\epsilon\tau}^{\sigma_{1}}\beta_{\nu}(\sigma_{2})\,...\,\int_{\epsilon\tau}^{\sigma_{j-1}}\beta_{\nu}(\sigma_{j})\left\|f(\sigma_{j},\cdot)\right\|_{good}\,d\sigma_{j}\,...\,d\sigma_{1}.
\end{align*} 
Pulling out the norm $\sup_{\sigma_{j}\simeq\tau}\left\|f(\sigma_{j},\cdot)\right\|_{good}$ of the temporal integrals, we need to bound
\begin{align*}
	&\beta_{\nu}(\tau)\int_{\epsilon\tau}^{\tau}\beta_{\nu}(\sigma_{1})\,\int_{\epsilon\tau}^{\sigma_{1}}\beta_{\nu}(\sigma_{2})\,...\,\int_{\epsilon\tau}^{\sigma_{j-1}}\beta_{\nu}(\sigma_{j})\,d\sigma_{j}\,...\,d\sigma_{1}\\
	\lesssim &\epsilon^{-j}\tau^{-j-1}\int_{\epsilon\tau}^{\tau}\,\int_{\epsilon\tau}^{\sigma_{1}}\,...\,\int_{\epsilon\tau}^{\sigma_{j-1}}\,d\sigma_{j}\,...\,d\sigma_{1}\\
	\lesssim&\epsilon^{-j}\tau^{-j-1}\int_{\epsilon\tau}^{\tau}\,\int_{\epsilon\tau}^{\sigma_{1}}\,...\,\int_{\epsilon\tau}^{\sigma_{j-2}}(\sigma_{j-1}-\epsilon\tau)d\sigma_{j-1}\,...\,d\sigma_{1}\\
	\lesssim &\frac12\epsilon^{-j}\tau^{-j-1}\int_{\epsilon\tau}^{\tau}\,\int_{\epsilon\tau}^{\sigma_{1}}\,...\,\int_{\epsilon\tau}^{\sigma_{j-3}}(\sigma_{j-2}-\epsilon\tau)^{2}d\sigma_{j-2}\,...\,d\sigma_{1}\\
	\lesssim&\,...\\
	\lesssim &\frac{\epsilon^{-j}}{j!}\cdot\frac{(1-\epsilon)^{j}}{\tau}\lesssim \tau^{-1}\frac{\epsilon^{-j}}{j!},
\end{align*}
which is desired. Here the factor $\tau^{-1}$ is not a smallness gain, since the reiteration starts with a term like $\tau^{-1}\calK_{\hbar}^{(0)}f$.
\end{proof}

The preceding finally entails the desired
\begin{proposition}\label{prop:ngeq2fourierwavesoln} There exists $\tau_0$ large enough, independently of $\hbar$, such that if $\mathcal{F}^{\hbar}(F_{\pm}(n))$ is a good source function, then the problem \eqref{eq:ngeq2Fourier1} admits a solution on $[\tau_0,\infty)$ with vanishing data 
	\[
	\left(\overline{x}^{\hbar}(\tau_0,\cdot),\,\mathcal{D}_{\tau}\overline{x}^{\hbar}(\tau_0,\cdot)\right) = \left(0,0\right)
	\]
	at time $\tau = \tau_0$, and such that 
	\[
	\left\| \overline{x}^{\hbar}\right\|_{\text{good}}\lesssim \left\| \mathcal{F}^{\hbar}(F_{\pm}(n))\right\|_{\text{goodsource}}.
	\]
	The implied constant is uniform in $\hbar$. 
\end{proposition}
\begin{proof} Pick $\epsilon = \frac12$ and choose $k$ large enough and then $\tau_0$ large enough such that $\epsilon^{\gamma k}\cdot e^{\epsilon^{-2}}<\frac12$. 
	Then set up the iterative scheme 
	\begin{align*}
		&-\left(\mathcal{D}_{\tau}^2 + \frac{\lambda'(\tau)}{\lambda(\tau)}\mathcal{D}_{\tau} + \xi \right)\overline{x}_{l+1}^{\hbar}\\
		& =  \mathcal{F}^{\hbar}\left(F_\pm(n)\right) +  2\frac{\lambda'(\tau)}{\lambda(\tau)}\mathcal{K}_{\hbar}^{(0)}\mathcal{D}_{\tau}\overline{x}_l^{\hbar} + \left(\frac{\lambda'(\tau)}{\lambda(\tau)}\right)'\mathcal{K}_{\hbar}^{(0)}\overline{x}_l^{\hbar} + \frac{\lambda'(\tau)}{\lambda(\tau)}\left[\mathcal{D}_{\tau},\mathcal{K}_{\hbar}^{(0)}\right]\overline{x}_l^{\hbar} + \left(\frac{\lambda'(\tau)}{\lambda(\tau)}\right)^2\left(\left(\mathcal{K}_{\hbar}^{(0)}\right)^2 + \mathcal{K}_{\hbar}^{(0)}\right)\overline{x}_l^{\hbar},\quad l\geq 0, 
	\end{align*}
	where $\overline{x}_0^{\hbar}$ solves the inhomogeneous problem 
	\[
	-\left(\mathcal{D}_{\tau}^2 + \frac{\lambda'(\tau)}{\lambda(\tau)}\mathcal{D}_{\tau} + \xi \right)\overline{x}_0^{\hbar} =  \mathcal{F}^{\hbar}\left(F_\pm(n)\right) 
	\]
	with vanishing initial data throughout. Then writing explicitly the iterated Duhamel for $\overline{x}_k^{\hbar}$, one gains smallness for all arising expressions upon choosing $\tau_0$ larger, if necessary, except for the expression $\Phi^k \mathcal{D}_{\tau}\overline{x}_0^{\hbar}$, for which smallness follows from the previous proposition. 
\end{proof}
\subsection{The final estimates for all the source terms near the light cone in the case $|n|\geq 2$.}
In order to finally wrap things up for the $|n|\geq 2$ modes, we of course also need to define the concept of a good (Fourier representation) function for the exceptional modes, and here we also need to take the coefficients $c_j(\tau)$ of the instabilities $\phi_j(R)$ in mind, recalling \eqref{eq:levelup-rep}, \eqref{eq:levelup0rep}, \eqref{eq:levelup1rep}. Recalling from Proposition~\ref{prop:bilin2}, Proposition~\ref{prop:bilin3} that we have already defined norms incorporating both the coefficients $c_n$ and the Fourier coefficients of $\mathcal{D}_jf$, we generalize things naturally as follows:
\begin{definition}\label{defi:goodfourierrepnless2} Let $n\in \{0,\pm1\}$, and assume that the angular momentum $n$ function $\phi(\tau, R)$ admits the representation 
	\[
	\phi(\tau,R) = c_n(\tau)\phi_n(R) + \phi_n(R)\cdot\int_0^R\left[\phi_n(s)\right]^{-1}\mathcal{D}_n\phi(\tau, s)\,ds,\quad \mathcal{D}_n\phi(\tau,R) = \int_0^\infty \xb(\tau,\xi)\phi_{n}(R,\xi)\tilde{\rho}_{n}(\xi)\,d\xi
	\]
	Then we say that the function $\phi$ is a good function, or alternatively the pair $\left(c_n(\tau),  \xb(\tau,\xi)\right)$ is good, provided 
	\begin{align*}
		\left|c_n(\tau)\right| + \tau\cdot \left|c_n'(\tau)\right|\lesssim \tau^{-2+10\nu},
	\end{align*}
	and the function $\xb(\tau,\xi)$ is 'good' in the sense of Def.~\ref{defi:xsingulartermsnless2smooth}.
	
		We then introduce the corresponding `good norm' for the pair of functions $\left(c_n(\tau),  x(\tau,\xi)\right)$ by 
	\begin{align*}
		\left\|\left(c_n(\tau),  \xb(\tau,\xi)\right)\right\|_{\text{good}}: &= \left\|\tau^{2-10\nu}\cdot\left[\left|c_n(\tau)\right| + \tau\cdot \left|c_n'(\tau)\right|\right]\right\|_{L_{\tau}^{\infty}([\tau_0,\infty)} + \big\| \xb\big\|_{good}
		\end{align*}
	where we refer to Def.~\ref{defi:xsingulartermsnless2smooth} for the last norm on the right.  We also have the natural analogues of the concept of `restricted principal singular part' with correspondingly modified norms, in analogy to Def.~\ref{defy:goodfunctionsnorm}.
\end{definition}

Combining Definition~\ref{defi:goodfourierrepnless2}, Definition~\ref{defi:goodfourierrepngeq2}, we finally have the tools that are sufficiently precise to derive the multilinear estimates to handle all the source terms arising for the $|n|\geq 2$ modes. Specifically, we strive to obtain analogues to Prop.~\ref{prop:smoothlinearsource}, Prop.~\ref{prop:bilinwithUregular2}, Prop.~\ref{prop:bilinwithUregular3}, but here we shall have to refer to the more sophisticated functional framework developed in the preceding. Recall the formulae
\begin{align*}
	\varphi_1 = \frac12\left[\varepsilon_{+} + \varepsilon_{-}\right],\quad \varphi_2 = \frac{1}{2i}\left[\varepsilon_{-} -\varepsilon_{+}\right].
\end{align*}
as well as the decompositions
\begin{equation*}
	\varepsilon_{+} = \sum_{n\in\Z}\varepsilon_{+}(n)e^{in\theta},\quad \varepsilon_{-} = \sum_{n\in\Z}\varepsilon_{-}(n)e^{in\theta}. 
\end{equation*}
and where we have $\varepsilon_{-}(-n) = \overline{\varepsilon_{+}(n)}$ since the solutions we consider are real valued. We shall assume that for $|n|\geq 2$ the function $\varepsilon_{+}(n)$ is a good angular momentum $n$ function in the sense that 
\begin{align*}
	\varepsilon_{+}(n)(\tau, R) = \int_0^\infty \phi_{n}(R,\xi)\xb_n(\tau,\xi)\rho_{n}(\xi)\,d\xi, 
\end{align*}
where $\xb_n(\tau,\xi)$ is a good angular momentum $n$ function in the sense of Definition~\ref{defi:goodfourierrepngeq2}. Similarly, if $n \in \{0,\pm 1\}$, we assume that 
\[
\varepsilon_{+}(n)(\tau, R) = c_n(\tau)\phi_n(R) + \phi_n(R)\cdot\int_0^R \left[\phi_n(s)\right]^{-1}\cdot \mathcal{D}_n\epsilon_{+}(n)(\tau, s)\,ds,
\]
where 
\begin{align*}
	\mathcal{D}_n\varepsilon_{+}(n)(\tau, R) = \int_0^\infty \phi_n(R,\xi)\xb_n(\tau,\xi)\tilde{\rho}_{n}(\xi)\,d\xi, 
\end{align*}
and the pair $\left(c_n(\tau),\,\xb_n(\tau,\xi)\right)$ is good in the sense of Definition~\ref{defi:goodfourierrepnless2}. Finally, set 
\begin{equation}\label{eq:LambdaDef}
	\Lambda: = \sum_{n\in\Z,\,|n|\geq 2}\left\langle n\right\rangle^{C}\left\| \xb_n\right\|_{\text{good}} + \sum_{n\in \{0,\pm1\}}\left\|\left(c_n(\tau),\,\xb_n(\tau,\xi)\right)\right\|_{\text{good}}. 
\end{equation}
We have the following analogue of Proposition~\ref{prop:generalsourcetermsfirstboundnearlightcone}:
\begin{proposition}\label{prop:ngeq2finalsourcetermestimatesingoodspaces} Assume that $\Lambda\ll 1$, and that $\tau_0\gg 1$ is sufficiently large. Then for each $n\in \Z, |n|\geq 2$, $j = 1,\ldots$, there exist angular momentum $n$ functions 
	\[
	\psi^{\pm}_j(n) = \int_0^\infty \xb_j^{(n)}(\tau,\xi)\phi_{n}(R,\xi)\rho_{n}(\xi)\,d\xi,
	\]
	with 
	\[
	\sum_{|n|\geq 2}n^{C}\left\| \xb_j^{(n)}\right\|_{\text{good}}\lesssim \Lambda^{3}, 
	\]
	and such that if $F_j^{(\pm)}$ represents any one of the functions occurring in Prop.~\ref{prop:smoothlinearsource}, Prop.~\ref{prop:bilinwithUregular2}, Prop.~\ref{prop:bilinwithUregular3}, and we write 
	\[
	F_j^{(\pm)} = \sum_{n\in \Z}F_j^{(\pm)}(n)e^{in\theta}, 
	\]
	then we have for each $|n|\geq 2$
	\begin{align*}
		\left(F_j^{(\pm)}(n) - \Box_n\psi^{\pm}_j(n)\right)|_{R<\nu\tau} = G_j^{(\pm)}(n)|_{R<\nu\tau}, 
	\end{align*}
	where $G_j^{(\pm)}(n)$ is a good angular momentum $n$ source function and more quantitatively setting 
	\[
	G_j^{(\pm)}(n) = \int_0^\infty \phi_{n}(R,\xi)\yb_n(\tau,\xi)\rho_{n}(\xi)\,d\xi,\,
	\]
	we have 
	\begin{align*}
		\sum_{|n|\geq 2}n^{C}\left\| \yb_n(\tau,\xi)\right\|_{\text{goodsource}}\lesssim \left(\tau_0^{-1}+\Lambda\right)\cdot\Lambda\ll\Lambda. 
	\end{align*}
	If we restrict to functions $\xb_n$ with restricted principal singular part, and correspondingly use $\left\| \xb_n\right\|_{\text{good}(r)}$, then $\yb_n$ also has restricted singular principal part, and we may replace 
	$\left\| \yb_n(\tau,\xi)\right\|_{\text{goodsource}}$ by $\left\| \yb_n(\tau,\xi)\right\|_{\text{goodsource}(r)}$. 
\end{proposition}

The proof of this proposition is completely analogous to the one of Proposition~\ref{prop:generalsourcetermsfirstboundnearlightcone}, the only difference being that instead of Lemma~\ref{lem:singFouriertiphysicalngeq2adm}, Lemma~\ref{lem:singFouriertiphysicalngeq2admDeriv}, Lemma~\ref{lem:derLinfty}, we take advantage of Lemma~\ref{lem:structuredlowfreqsmoothbasic} to describe functions at angular momentum $|n|\geq 2$ and with structured smooth but poorly temporally decaying distorted Fourier transform $\xb_{1,smooth}$. 
\section{The exceptional angular momenta: $n\in \{0,\pm 1\}$ and modulation theory}\label{sec:exceptional modes}
\subsection{A revisit to the background solution}
In order to perform the modulation argument, we need a more precise description than the one provided by Theorem \ref{thm:KSTGaoKBase} of the background solution near the light cone. More specifically we have (see \cite{GaoK, KST})
\begin{theorem}\label{thm:KSTGao precise}
	Let $U=Q(\lambda(t)r)+\epsilon(t,r), \epsilon=\epsilon_{1}+\epsilon_{2}$ be given as in Theorem \ref{thm:KSTGaoKBase}. Then $\epsilon_{1}$ admits the following decomposition:
	\begin{align*}
		\epsilon_{1}(t,r)=\epsilon_{1a}(t,r)+\epsilon_{1b}(t,r).
	\end{align*}
Here the leading order part $\epsilon_{1a}(t,r)$ is given by, with $R:=\lambda(t)r, \tau=\nu^{-1}t^{-\nu}$,
\begin{align*}
\epsilon_{1a}(t,r)=&\frac{1}{\nu\tau}\left[\left(g_{0}(1)+\frac{g_{1}(1)}{(\nu\tau)^{\nu+\frac12}}(\nu\tau-R)^{\nu+\frac12}\right)\log R\right]\\
&+\frac{1}{\nu\tau}\left[h_{0}(1)+\frac{h_{1}(1)}{(\nu\tau)^{\nu+\frac12}}(\nu\tau-R)^{\nu+\frac12}+\frac{h_{2}(1)}{(\nu\tau)^{\nu+\frac12}}(\nu\tau-R)^{\nu+\frac12}\left(\log(\nu\tau-R)-\log\nu\tau\right)\right].
\end{align*}
The functions $g_{0}(a), g_{1}(a)$ are given by 
\begin{align*}
g_{0}(a)=&2(\nu+1)(\nu^{3}-\nu^{2}-2)k^{-1}\phi_{1}(a)\int_{0}^{a}\phi_{2}^{0}(a')(1+a')^{-\nu-\frac12}a^{\prime2}\,da'\\
&+2(\nu+1)(\nu^{3}-\nu^{2}-2)k^{-1}\phi_{2}^{0}(a)(1-a)^{\nu+\frac12}\int_{a}^{1}\phi_{1}(a')\left(1-a^{\prime2}\right)^{-\nu-\frac12}a^{\prime2}\,da'\\
g_{1}(a)=&-2(\nu+1)(\nu^{3}-v^{2}-2)k^{-1}\phi_{2}^{0}(a)\int_{0}^{1}\phi_{1}(a')\left(1-a^{\prime2}\right)^{-\nu-\frac12}a^{\prime2}\,da'
\end{align*}
The constants $h_{0}(1), h_{1}(1), h_{2}(1)$ are given by
\begin{align*}
h_{0}(1)=&k^{-1}\int_{0}^{1}\phi_{2}^{0}(a')(1+a')^{-\nu-\frac12}a^{\prime}\,f^{2}_{0}(a')\,da'\\
&+k^{-1}\int_{0}^{1}\phi_{2}^{0}(a')(1+a')^{-\nu-\frac12}a^{\prime}\,f^{1}_{0}(a')(1-a')^{\nu-\frac12}\,da'\\
&+k^{-1}\int_{0}^{1}\phi_{2}^{0}(a')(1+a')^{-\nu-\frac12}a^{\prime}\,f^{3}_{0}(a')(1-a')^{\nu+\frac12}\,da'\\
&+k^{-1}\int_{0}^{1}\phi_{2}^{0}(a')(1+a')^{-\nu-\frac12}a^{\prime}\,f^{4}_{0}(a')(1-a')^{\nu+\frac32}\,da'\\
h_{1}(1)=&-k^{-1}\frac{f_{0}^{1}(1)}{2^{\nu+\frac12}\left(\nu+\frac12\right)}\\
&+k^{-1}\int_{0}^{1}\phi_{1}(a')\left(1-a^{\prime2}\right)^{-\nu-\frac12}a^{\prime}\,\left(f^{3}_{0}(a')(1-a')^{\nu+\frac12}+f^{4}_{0}(a')(1-a')^{\nu+\frac32}\right)\,da'\\
&-k^{-1}\int_{0}^{1}\phi_{1}(a')\left(1-a^{\prime2}\right)^{-\nu-\frac12}a'\,f^{2}_{0}(a')\,da'.\\
h_{2}(1)=&k^{-1}\frac{f_{0}^{1}(1)}{2^{\nu+\frac12}}.
\end{align*}
Here $\phi_{1}(a):=1+\sum_{\ell=1}^{\infty}\mu_{\ell}(1-a)^{\ell}$ and $(1-a)^{\nu+\frac12}\phi_{2}^{0}(a)=(1-a)^{\nu+\frac12}\left[1+\sum_{\ell=1}\tilde{\mu}_{\ell}(1-a)^{\ell}\right]:=\phi_{2}(a)$ constitute a fundamental system of the operator
\begin{align*}
	\calL_{\nu}:=\left(1-a^{2}\right)\partial_{a}^{2}+\left(a^{-1}+2a\nu-2a\right)\partial_{a}+\left(-\nu^{2}+\nu-a^{-2}\right),
\end{align*}
as shown in \cite{KST}. The non-zero constant $k$ is defined such that the Wronskian $W(a):=\phi_{1}(a)\phi^{\prime}_{2}(a)-\phi^{\prime}_{1}(a)\phi_{2}(a)$ is $k\left(1-a^{2}\right)^{\nu-\frac12}a^{-1}$. The functions $f_{0}^{1}, f_{0}^{2}, f_{0}^{3}, f_{0}^{4}$ are given by
\begin{align*}
	f_{0}^{1}(a):=&-\nu\left(\nu+\frac12\right)a\,g_{1}(a)(1-a)^{\nu-\frac12}\\
	f_{0}^{2}(a):=&\nu(\nu-2)(1-\nu^{2})a-\left(a^{-1}-(1+\nu)a\right)g^{\prime}_{0}(a)-\left((\nu+1)\nu+a^{-2}\right)g_{0}(a)\\
	f_{0}^{3}(a):=&\left(\nu+\frac12\right)\frac{1+a}{a}g_{1}(a)-\left((\nu+1)\nu+a^{-2}\right)g_{1}(a)+\nu ag_{1}^{\prime}(a)\\
	f_{0}^{4}(a):=&-\frac{1+a}{a}g_{1}^{\prime}(a).
\end{align*}
\vspace{3mm}
The lower order part $\epsilon_{1b}(t,r)$ enjoys the following trichotomy:
\begin{itemize}
\item $\epsilon_{1b}$ decays as fast as $\epsilon_{1a}$ when $t\rightarrow0_{+}$ but smoother than $\epsilon_{1a}$ near the cone $\{(t,r)|t=r\}$,
\item  $\epsilon_{1b}$ is as smooth as $\epsilon_{1a}$ near the light cone $\{(t,r)|t=r\}$ but decays faster than $\epsilon_{1a}$ when $t\rightarrow0_{+}$,
\item $\epsilon_{1b}$ decays faster than $\epsilon_{1a}$ when $t\rightarrow0_{+}$ and is smoother than $\epsilon_{1a}$ near the cone $\{(t,r)|t=r\}$.
 \end{itemize}
\end{theorem}
\begin{proof}
According to the analysis in \cite{GaoK}, $U(t,r)$ is given by
	\begin{align*}
		U(t,r)=u_{k-1}(t,r)+\epsilon_{k}(t,r).
	\end{align*}
Here $u_{k-1}(t,r)$ is an approximation and $e_{k-1}(t,r)$ is its corresponding error. The $k$-th order approximation is in term given by
\begin{align*}
	u_{k}(t,r)=Q(\lambda(t)r)+\sum_{j=1}^{k}v_{j}(t,r),\quad u_{0}(t,r)=Q(\lambda(t)r)
\end{align*}
The error $e_{k}$ are defined by
\begin{align}\label{def ek}
	e_{k}:=\partial_{t}^{2}u_{k}-\partial_{r}^{2}u_{k}-\frac1r\partial_{r}u_{k}+\frac{\sin(2u_{k})}{2r^{2}},
\end{align}
and it is decomposed into $e_{k}=e^{0}_{k}+e^{1}_{k}$. The leading order contribution $e_{k}^{0}$ satisfies
\begin{align}\label{def ek0}
	\begin{split}
		&\left(-\partial_{t}^{2}+\partial_{r}^{2}+\frac1r\partial_{r}-\frac{1}{r^{2}}\right)v_{2k}=e_{2k-1}^{0},\quad\left(\partial_{r}^{2}+\frac1r\partial_{r}-\frac{\cos(2u_{0})}{r^{2}}\right)v_{2k+1}=e_{2k}^{0}.
	\end{split}
\end{align}
Then a direct computation gives the successive errors:
\begin{align*}
	e_{2k}=e_{2k-1}^{1}+N_{2k}(v_{2k}),\quad e_{2k+1}=e_{2k}^{1}+\partial_{t}^{2}v_{2k+1}+N_{2k+1}(v_{2k+1}),
\end{align*}
where
\begin{align*}
	-N_{2k+1}(v)=&\frac{\cos(2u_{0})-\cos(2u_{2k})}{r^{2}}v+\frac{\sin(2u_{2k})}{2r^{2}}(1-\cos(2v))+\frac{\cos(2u_{2k})}{2r^{2}}(2v-\sin(2v)),
\end{align*}
and 
\begin{align*}
	-N_{2k}(v)=&\frac{1-\cos(2u_{2k-1})}{r^{2}}v+\frac{\sin(2u_{2k-1})}{2r^{2}}(1-\cos(2v))+\frac{\cos(2u_{2k-1})}{2r^{2}}(2v-\sin(2v)).
\end{align*}
The approximation pieces $v_{k}$ in the above construction satisfy
\begin{align*}
	v_{2k-1}\in \frac{1}{(t\lambda)^{2k}}IS^{3}\left(R(\log R)^{2k-1},\calQ_{k-1}\right),\quad v_{2k}\in \frac{1}{(t\lambda)^{2k+2}}IS^{3}\left(R^{3}(\log R)^{2k-1},\calQ_{k}\right).
\end{align*} 
Since $t\lambda(t)=\nu\tau$, as $k$ increases, $v_{k}$ decays faster in $\tau$. On the other hand, recall the reason we modulate the solution: upon taking the distorted Fourier transforms, certain singular profiles in the approximation solution $u_{k}$ do not have sufficient decay in $\tau$. So we need to modulate the solution to annihilate the contribution from these singular profiles. It turns out that we only need to annihilate the contribution from $v_{2}$ and starting from $v_{3}$, there is sufficient decay in $\tau$. As we shall see later, for $k>1$, the distorted Fourier transform of the singular profile decays faster in large frequency, which in turn gives extra temporal decay when it gets mapped by the inhomogeneous parametrix. This implies that we only need to compute $v_{2}$ explicitly. Recall from \cite{GaoK}, $v_{2}$ is given by $v_{2}=w_{2}+\tilde{w}_{2}$ where
\begin{align}\label{w2 profile}
	w_{2}=\frac{1}{t\lambda}\left(W_{2}^{1}(a)\,\log R+W_{2}^{0}(a)\right),\quad \tilde{w}_{2}=\frac{1}{(t\lambda)^{2}}\left(\tilde{W}_{2}^{1}(a)\,\log R+\tilde{W}_{2}^{0}(a)\right).
\end{align}
Here $a=\frac{r}{t}=\frac{R}{t\lambda}$ and the functions $W_{2}^{i}, \tilde{W}_{2}^{i}, i=0,1$ are to be determined. Note that $\tilde{w}_{2}$ has better decay in $\tau$, so we only focus on $w_{2}$. To find an explicit formula for $w_{2}$, the starting point is the equation
\begin{align*}
	t^{2}\tilde{\Box}v_{2}=t^{2}e_{1}^{0},\quad \tilde{\Box}:=-\partial_{t}^{2}+\partial_{r}^{2}+\frac{1}{r}\partial_{r}-\frac{1}{r^{2}},
\end{align*}
and by \cite{GaoK}
\begin{align*}
	t^{2}e_{1}^{0}=\frac{1}{t\lambda}\left(c_{1}\,R\log R+c_{2}\,R+c_{3}\,\log R+c_{4}\right),
\end{align*}
and we need to determine $c_{i}, i=1,2,3,4$ explicitly. To this end, we start with $v_{1}$ which together with $u_{0}$ gives an explicit expression for $e_{1}$. Following \cite{GaoK}, $v_{1}$ satisfies 
\begin{align*}
	(t\lambda)^{2}\tilde{\calL}v_{1}=t^{2}e_{0},\quad \Rightarrow\quad -\calL(R^{\frac12}v_{1})=-R^{\frac12}t^{2}e_{0},
\end{align*}
where we denote
\begin{align*}
	\tilde{\calL}:=\partial_{R}^{2}+\frac{1}{R}\partial_{R}-\frac{\cos(2u_{0})}{R^{2}}=\partial_{R}^{2}+\frac{1}{R}\partial_{R}-\frac{1}{R^{2}}\frac{1-6R^{2}+R^{4}}{(1+R^{2})^{2}},\quad -\calL:=\partial_{R}^{2}-\frac{3}{4R^{2}}+\frac{8}{(1+R^{2})^{2}}.
\end{align*}
A fundamental system for $\calL$ is given by
\begin{align*}
	\phi(R)=\frac{R^{\frac32}}{1+R^{2}},\quad \theta(R)=\frac{-1+4R^{2}\log R+R^{4}}{R^{\frac12}(1+R^{2})},
\end{align*}
which in turn gives
\begin{align*}
(t\lambda)^{2}v_{1}=\frac12 R^{-\frac12}\theta(R)\int_{0}^{R}\phi(R')R^{\prime\frac12}t^{2}e_{0}(R')\,dR'-\frac12R^{-\frac12}\phi(R)\int_{0}^{R}\theta(R')R^{\prime\frac12}t^{2}e_{0}(R')\,dR'.
\end{align*}
Recall that 
\begin{align*}
	t^{2}e_{0}=t^{2}\partial_{t}^{2}u_{0}=\left(-\nu(\nu+1)\frac{2R}{1+R^{2}}+(\nu+1)^{2}\frac{4R}{(1+R^{2})^{2}}\right).
\end{align*}
We then have (we only keep the terms up to constant order)
\begin{align*}
	&R^{-\frac12}\theta(R)\int_{0}^{R}\phi(R^{\prime})R^{\prime\frac12} t^{2}e_{0}(R^{\prime})\,dR'\\
	=&\frac{-1+4R^{2}\log R+R^{4}}{R(1+R^{2})}\int_{0}^{R}\left(1-\frac{1}{1+R^{\prime2}}\right)\left(-\nu(\nu+1)\frac{1}{1+R^{\prime2}}+(\nu+1)^{2}\frac{2}{(1+R^{\prime2})^{2}}\right)\,d(R^{\prime2}+1)\\
	=&\left(R+\frac{4\log R}{R}-\frac{1}{R}-\frac{4\log R}{R(R^{2}+1)}\right)\left(-\nu(\nu+1)\log(R^{2}+1)+\nu(\nu+1)\right)+\textrm{l.o.t.}\\
	=&-2\nu(\nu+1)R\log R+\nu(\nu+1)R+\textrm{l.o.t.}
\end{align*}
and
\begin{align*}
	&-R^{-\frac12}\phi(R)\int_{0}^{R}\theta(R')R^{\prime\frac12}t^{2}e_{0}(R')\,dR'\\
	=&\frac{R}{1+R^{2}}\int_{0}^{R}\frac{R^{\prime4}+4R^{\prime2}\log R^{\prime}-1}{1+R^{\prime2}}\left(\nu(\nu+1)\frac{1}{1+R^{\prime2}}-(\nu+1)^{2}\frac{2R^{\prime2}}{(1+R^{\prime2})^{2}}\right)\,d(R^{\prime2}+1)\\
	=&\frac{R}{1+R^{2}}\int_{0}^{R}\left(R^{\prime2}+4\log R^{\prime}-1-\frac{4\log R^{\prime}}{R^{\prime2}+1}\right)\left(\nu(\nu+1)\frac{1}{1+R^{\prime2}}-(\nu+1)^{2}\frac{2R^{\prime2}}{(1+R^{\prime2})^{2}}\right)\,d(R^{\prime2}+1)\\
	=&\frac{R}{1+R^{2}}\int_{0}^{R}\left(R^{\prime2}+1+4\log R^{\prime}-2-\frac{4\log R'}{R^{\prime2}+1}\right)\left((\nu(\nu+1)-2(\nu+1)^{2})\frac{1}{1+R^{\prime2}}+\frac{2(\nu+1)^{2}}{(1+R^{\prime2})^{2}}\right)\,d(R^{\prime2}+1)\\
	=&(\nu(\nu+1)-2(\nu+1)^{2})R+\textrm{l.o.t.}
\end{align*}
Summarizing, we have
\begin{align}\label{v1 asymp}
	\begin{split}
		v_{1}=-\frac{1}{(t\lambda)^{2}}\left(\nu(\nu+1)R\log R+(\nu+1)R\right)+\textrm{l.o.t.}
	\end{split}
\end{align}
Recall from \cite{GaoK} and using the fact $|v_{1}|\ll1$,
\begin{align*}
	t^{2}e_{1}=\,&t^{2}\partial_{t}^{2}v_{1}-\frac{\sin(2u_{0})}{2R^{2}}(t\lambda)^{2}\left(1-\cos(2v_{1})\right)-\frac{\cos(2u_{0})}{2R^{2}}(t\lambda)^{2}\left(2v_{1}-\sin(2v_{1})\right)\\
	\simeq\,&t^{2}\partial_{t}^{2}v_{1}-\frac{\sin(2u_{0})}{2R^{2}}(t\lambda)^{2}\cdot\frac{v_{1}^{2}}{2}-\frac{\cos(2u_{0})}{2R^{2}}(t\lambda)^{2}\cdot\frac{v_{1}^{3}}{6}.
\end{align*}
We drop out the higher order powers in $v_{1}$ because if we have one more power in $v_{1}$, then we at least gain a power of 
\begin{align*}
	\frac{R\log R}{(t\lambda)^{2}}\lesssim \frac{\log R}{R},
\end{align*}
which makes the corresponding contribution lower order.
For the same reason we can further drop out the cubic contribution involving $v_{1}^{3}$. Now we compute $t^{2}\partial_{t}^{2}v_{1}$ explicitly. We shall use the fact
\begin{align*}
	v_{1}=-\nu(\nu+1)t^{-1+\nu}r\left[-(1+\nu)\log t+\log r\right]-(\nu+1)t^{-1+\nu}r+\textrm{l.o.t.},
\end{align*}
which implies
\begin{align*}
	&\partial_{t}v_{1}=\nu(1-\nu^{2})t^{-2+\nu}r\left[-(1+\nu)\log t+\log r\right]+\nu(\nu+1)^{2}t^{-2+\nu}r+(1-\nu^{2})t^{-2+\nu}r+\textrm{l.o.t.},
\end{align*}
which in turn implies
\begin{align*}
	\partial_{t}^{2}v_{1}=&\nu(\nu-2)(1-\nu^{2})t^{-3+\nu}r\left[-(1+\nu)\log t+\log r\right]-\nu(1-\nu^{2})(1+\nu)t^{-3+\nu}r\\
	&+\nu(\nu+1)^{2}(\nu-2)t^{-3+\nu}r+(1-\nu^{2})(\nu-2)t^{-3+\nu}r+\textrm{l.o.t.}
\end{align*}
Therefore we have
\begin{align*}
	t^{2}\partial_{t}^{2}v_{2}=&\frac{1}{(t\lambda)^{2}}\left[\nu(\nu-2)(1-\nu^{2})R\log R+2(\nu+1)(v^{3}-v^{2}-2)R\right]+\textrm{l.o.t.}
\end{align*}
For the contribution from $-\frac{\sin(2u_{0})}{2R^{2}}(t\lambda)^{2}\cdot\frac{v_{1}^{2}}{2}$, we have
\begin{align*}
	-\frac{\sin(2u_{0})}{2R^{2}}(t\lambda)^{2}\cdot v_{1}\simeq&\frac{4(1-R^{2})}{R(R^{2}+1)^{2}}\cdot \left(\nu(\nu+1)R\log R+(\nu+1)R\right)\\
	\simeq &\frac{1}{(t\lambda)^{3}}\cdot \left(\nu(\nu+1)R\log R+(\nu+1)R\right),
\end{align*}
which already decays faster than $v_{1}$ itself. Therefore $-\frac{\sin(2u_{0})}{2R^{2}}(t\lambda)^{2}\cdot\frac{v_{1}^{2}}{2}$ is lower order.
This finally shows that we can even drop out the contribution from the term involving $v_{1}^{2}$ and simply write
\begin{align}\label{t2e1}
	t^{2}e_{1}=&\frac{1}{(t\lambda)^{2}}\left[\nu(\nu-2)(1-\nu^{2})R\log R+2(\nu+1)(v^{3}-v^{2}-2)R\right]+\textrm{l.o.t.}
\end{align}
Now we are ready to find an explicit profile for $w_{2}$. Recall from \eqref{w2 profile}, we need to find out $W_{2}^{1}(a)$ and $W_{2}^{0}(a)$. Following \cite{GaoK}, we have
\begin{align*}
	\calL_{\nu}W_{2}^{i}(a)=a\,c_{i+1}-F_{i}(a),\quad i=0,1.
\end{align*}
Here 
\begin{align*}
	c_{1}=\nu(\nu-2)(1-\nu^{2}),\quad c_{2}=2(\nu+1)(\nu^{3}-\nu^{2}-2),
\end{align*}
and 
\begin{align*}
	F_{1}(a)=0,\quad F_{0}(a)=\left((\nu+1)\nu+a^{-2}\right)W_{2}^{1}(a)+\left(a^{-1}-(1+\nu)a\right)\partial_{a}W_{2}^{1}(a),
\end{align*}
as well as
\begin{align}\label{def L nu}
	\calL_{\nu}:=\left(1-a^{2}\right)\partial_{a}^{2}+\left(a^{-1}+2a\nu-2a\right)\partial_{a}+\left(-\nu^{2}+\nu-a^{-2}\right)
\end{align}
According to Lemma 3.9 in \cite{KST}, for sufficiently small $\nu>0$, $W_{2}^{1}(a)$ has the following ansatz:
\begin{align*}
	W_{2}^{1}(a)=g_{0}(a)+g_{1}(a)(1-a)^{\nu+\frac12}.
\end{align*}
Following the argument in \cite{KST}, the operator \eqref{def L nu} admits a fundamental system of the form
\begin{align}\label{funda L nu}
	\begin{split}
	\phi_{1}(a)=1+\sum_{\ell=1}^{\infty}\mu_{\ell}(1-a)^{\ell},\quad \phi_{2}(a)=&(1-a)^{\nu+\frac12}\left[1+\sum_{\ell=1}^{\infty}\tilde{\mu}_{\ell}(1-a)^{\ell}\right]\\
	=:&(1-a)^{\nu+\frac12}\phi_{2}^{0}(a).
	\end{split}
\end{align}
A direct calculation shows that the Wronskian $W(a):=\phi_{1}(a)\phi^{\prime}_{2}(a)-\phi^{\prime}_{1}(a)\phi_{2}(a)$ is given by
\begin{align*}
	W(a)=k\left(1-a^{2}\right)^{\nu-\frac12}a^{-1},\quad k\slashed{=}0.
\end{align*}
Define the Green's function
\begin{align*}
	G_{\nu}(a,a'):=\phi_{1}(a)\phi_{2}(a')-\phi_{1}(a')\phi_{2}(a).
\end{align*}
Then $W_{2}^{1}(a)$, which satisfies the equation $\calL_{\nu}W_{2}^{1}(a)=2(\nu+1)(\nu^{3}-\nu^{2}-2)a$, is given by
\begin{align*}
	W_{2}^{1}(a)=&2(\nu+1)(\nu^{3}-\nu^{2}-2)\int_{0}^{a}G_{\nu}(a,a')\,W(a')^{-1}\left(1-a^{\prime2}\right)^{-1}a'\,da'\\
	=&2(\nu+1)(\nu^{3}-\nu^{2}-2)k^{-1}\int_{0}^{a}G_{\nu}(a,a^{\prime})\left(1-a^{\prime2}\right)^{-\nu-\frac12}a^{\prime2}\,da'\\
	=&2(\nu+1)(\nu^{3}-\nu^{2}-2)k^{-1}\phi_{1}(a)\int_{0}^{a}\phi_{2}^{0}(a')(1+a^{\prime})^{-\nu-\frac12}a^{\prime2}\,da'\\
	&-2(\nu+1)(\nu^{3}-\nu^{2}-2)k^{-1}\phi_{2}^{0}(a)(1-a)^{\nu+\frac12}\int_{0}^{1}\phi_{1}(a')\left(1-a^{\prime2}\right)^{-\nu-\frac12}a^{\prime2}\,da^{\prime}\\
	&+2(\nu+1)(\nu^{3}-\nu^{2}-2)k^{-1}\phi_{2}^{0}(a)(1-a)^{\nu+\frac12}\int_{a}^{1}\phi_{1}(a')\left(1-a^{\prime2}\right)^{-\nu-\frac12}a^{\prime2}\,da^{\prime}.
\end{align*}
Note  that the first and third terms on the RHS above are both analytic near $a=1$ and the third term vanishes at $a=1$ of order $O(1-a)$. Therefore we have
\begin{align}\label{coe W21}
	\begin{split}
		g_{0}(a)=&2(\nu+1)(\nu^{3}-\nu^{2}-2)k^{-1}\phi_{1}(a)\int_{0}^{a}\phi_{2}^{0}(a')(1+a^{\prime})^{-\nu-\frac12}a^{\prime2}\,da'\\
		&+2(\nu+1)(\nu^{3}-\nu^{2}-2)k^{-1}\phi_{2}^{0}(a)(1-a)^{\nu+\frac12}\int_{a}^{1}\phi_{1}(a')\left(1-a^{\prime2}\right)^{-\nu-\frac12}a^{\prime2}\,da^{\prime},\\
		g_{1}(a)=&-2(\nu+1)(\nu^{3}-\nu^{2}-2)k^{-1}\phi_{2}^{0}(a)\int_{0}^{1}\phi_{1}(a')\left(1-a^{\prime2}\right)^{-\nu-\frac12}a^{\prime2}\,da^{\prime},
	\end{split}
\end{align}
and
\begin{align}\label{coe W21 at 1}
	\begin{split}
		g_{0}(1)=&2(\nu+1)(\nu^{3}-\nu^{2}-2)k^{-1}\int_{0}^{1}\phi_{2}^{0}(a')(1+a^{\prime})^{-\nu-\frac12}a^{\prime2}\,da',\\
		g_{1}(1)=&-2(\nu+1)(\nu^{3}-\nu^{2}-2)k^{-1}\int_{0}^{1}\phi_{1}(a')\left(1-a^{\prime2}\right)^{-\nu-\frac12}a^{\prime2}\,da^{\prime}.
	\end{split}
\end{align}
Next we turn to $W_{2}^{0}(a)$. We start by calculating its source:
\begin{align*}
	F_{0}(a)=&\left((\nu+1)\nu+a^{-2}\right)\left(g_{0}(a)+g_{1}(a)(1-a)^{\nu+\frac12}\right)\\
	&+\left(a^{-1}-(1+\nu)a\right)\left(g^{\prime}_{0}(a)+g^{\prime}_{1}(a)(1-a)^{\nu+\frac12}-\left(\nu+\frac12\right)g_{1}(a)(1-a)^{\nu-\frac12}\right).
\end{align*}
Therefore the source for $W_{2}^{0}(a)$ is given by
\begin{align}\label{W20 source}
	\begin{split}
		f_{0}(a):=&\left(\nu+\frac12\right)\left(a^{-1}-(1+\nu)a\right)g_{1}(a)(1-a)^{\nu-\frac12}\\
		&+\nu(\nu-2)(1-\nu^{2})a-\left(a^{-1}-(1+\nu)a\right)g^{\prime}_{0}(a)-\left((\nu+1)\nu+a^{-2}\right)g_{0}(a)\\
		&-\left((\nu+1)\nu+a^{-2}\right)g_{1}(a)(1-a)^{\nu+\frac12}-\left(a^{-1}-(1+\nu)a\right)g^{\prime}_{1}(a)(1-a)^{\nu+\frac12}\\
		=&-\nu\left(\nu+\frac12\right)a\,g_{1}(a)(1-a)^{\nu-\frac12}\\
		&+\nu(\nu-2)(1-\nu^{2})a-\left(a^{-1}-(1+\nu)a\right)g^{\prime}_{0}(a)-\left((\nu+1)\nu+a^{-2}\right)g_{0}(a)\\
		&+\left(\nu+\frac12\right)\frac{1+a}{a}g_{1}(a)(1-a)^{\nu+\frac12}-\left((\nu+1)\nu+a^{-2}\right)g_{1}(a)(1-a)^{\nu+\frac12}+\nu a\,g_{1}^{\prime}(a)(1-a)^{\nu+\frac12}\\
		&-\frac{1+a}{a}g^{\prime}_{1}(a)(1-a)^{\nu+\frac32}\\
		=:&f_{0}^{1}(a)\,(1-a)^{\nu-\frac12}+f_{0}^{2}(a)+f_{0}^{3}(a)\,(1-a)^{\nu+\frac12}+f_{0}^{4}(a)\,(1-a)^{\nu+\frac32}.
	\end{split}
\end{align}
Here $f_{0}^{i}(a), i=1,2,3,4$ are analytic near $a=1$ and integrable on $[0,1]$ due to the fact that  $g_{0}(a)$ and $g_{1}(a)$ vanish at $a=0$ of order $O\left(a^{3}\right)$. Using Green's function, we have
\begin{align}\label{W20 pre}
	\begin{split}
W_{2}^{0}(a)=&k^{-1}\int_{0}^{a}G_{\nu}(a,a')\left(1-a^{\prime2}\right)^{-\nu-\frac12}a^{\prime}\,f_{0}(a^{\prime})\,da^{\prime}\\
=&k^{-1}\phi_{1}(a)\int_{0}^{a}\phi_{2}^{0}(a')(1+a')^{-\nu-\frac12}a^{\prime}\,f_{0}(a')\,da'\\
&-k^{-1}\phi_{2}^{0}(a)(1-a)^{\nu+\frac12}\int_{0}^{a}\phi_{1}(a')\left(1-a^{\prime2}\right)^{-\nu-\frac12}\,a^{\prime}\,f_{0}(a')\,da'.
\end{split}
\end{align}
The first term on the RHS of \eqref{W20 pre} can be written as
\begin{align*}
	h^{1}_{0}(a)+h^{1}_{1}(a)(1-a)^{\nu+\frac12},
\end{align*}
where 
\begin{align}\label{coe W20 a}
	\begin{split}
	h^{1}_{0}(a):=&k^{-1}\phi_{1}(a)\int_{0}^{a}\phi_{2}^{0}(a')(1+a')^{-\nu-\frac12}a^{\prime}\,f^{2}_{0}(a')\,da',\\
	&+k^{-1}\phi_{1}(a)\int_{0}^{1}\phi_{2}^{0}(a')(1+a')^{-\nu-\frac12}a^{\prime}\,f^{1}_{0}(a')(1-a')^{\nu-\frac12}\,da'\\
	&+k^{-1}\phi_{1}(a)\int_{0}^{1}\phi_{2}^{0}(a')(1+a')^{-\nu-\frac12}a^{\prime}\,f^{3}_{0}(a')(1-a')^{\nu+\frac12}\,da'\\
	&+k^{-1}\phi_{1}(a)\int_{0}^{1}\phi_{2}^{0}(a')(1+a')^{-\nu-\frac12}a^{\prime}\,f^{4}_{0}(a')(1-a')^{\nu+\frac32}\,da'.\\
	h^{1}_{1}(a):=&-k^{-1}\phi_{1}(a)(1-a)^{-\nu-\frac12}\int_{a}^{1}\phi_{2}^{0}(a')(1+a')^{-\nu-\frac12}a^{\prime}\,f^{1}_{0}(a')(1-a')^{\nu-\frac12}\,da'\\
	&-k^{-1}\phi_{1}(a)(1-a)^{-\nu-\frac12}\int_{a}^{1}\phi_{2}^{0}(a')(1+a')^{-\nu-\frac12}a^{\prime}\,f^{3}_{0}(a')(1-a')^{\nu+\frac12}\,da'\\
	&-k^{-1}\phi_{1}(a)(1-a)^{-\nu-\frac12}\int_{a}^{1}\phi_{2}^{0}(a')(1+a')^{-\nu-\frac12}a^{\prime}\,f^{4}_{0}(a')(1-a')^{\nu+\frac32}\,da'
	\end{split}
\end{align}
and 
\begin{align}\label{coe W20 a at 1}
	\begin{split}
		h_{0}^{1}(1)=&k^{-1}\int_{0}^{1}\phi_{2}^{0}(a')(1+a')^{-\nu-\frac12}a^{\prime}\,f^{2}_{0}(a')\,da'\\
		&+k^{-1}\int_{0}^{1}\phi_{2}^{0}(a')(1+a')^{-\nu-\frac12}a^{\prime}\,f^{1}_{0}(a')(1-a')^{\nu-\frac12}\,da'\\
		&+k^{-1}\int_{0}^{1}\phi_{2}^{0}(a')(1+a')^{-\nu-\frac12}a^{\prime}\,f^{3}_{0}(a')(1-a')^{\nu+\frac12}\,da'\\
		&+k^{-1}\int_{0}^{1}\phi_{2}^{0}(a')(1+a')^{-\nu-\frac12}a^{\prime}\,f^{4}_{0}(a')(1-a')^{\nu+\frac32}\,da'.\\
		h_{1}^{1}(1)=&-k^{-1}\frac{f_{0}^{1}(1)}{2^{\nu+\frac12}\left(\nu+\frac12\right)}.
	\end{split}
\end{align}
The second term on the RHS of \eqref{W20 pre} can be written as
\begin{align*}
	h_{0}^{2}(a)+h_{1}^{2}(a)(1-a)^{\nu+\frac12}+h_{2}(a)(1-a)^{\nu+\frac12}\log(1-a)
\end{align*}
where
\begin{align}\label{coe W20 b}
\begin{split}
h_{0}^{2}(a):=&k^{-1}\phi_{2}^{0}(a)(1-a)^{\nu+\frac12}\int_{a}^{1}\phi_{1}(a')\left(1-a^{\prime2}\right)^{-\nu-\frac12}a'\,f_{0}^{2}(a')\,da'\\
h_{1}^{2}(a):=&-k^{-1}\phi_{2}^{0}(a)\int_{0}^{a}\phi_{1}(a')\left(1-a^{\prime2}\right)^{-\nu-\frac12}a^{\prime}\,\left(f^{3}_{0}(a')(1-a')^{\nu+\frac12}+f^{4}_{0}(a')(1-a')^{\nu+\frac32}\right)\,da'\\
&-k^{-1}\phi_{2}^{0}(a)\int_{0}^{1}\phi_{1}(a')\left(1-a^{\prime2}\right)^{-\nu-\frac12}a'\,f^{2}_{0}(a')\,da'.\\
h_{2}(a):=&-k^{-1}\frac{\phi_{2}^{0}(a)}{\log(1-a)}\int_{0}^{a}\phi_{1}(a')\left(1-a^{\prime2}\right)^{-\nu-\frac12}a^{\prime}\,f_{0}^{1}(a')(1-a')^{\nu-\frac12}\,da'
\end{split}
\end{align}
and
\begin{align}\label{coe W20 b at 1}
	\begin{split}
		h_{0}^{2}(1)=&0,\\
		h_{1}^{2}(1)=&-k^{-1}\int_{0}^{1}\phi_{1}(a')\left(1-a^{\prime2}\right)^{-\nu-\frac12}a^{\prime}\,\left(f^{3}_{0}(a')(1-a')^{\nu+\frac12}+f^{4}_{0}(a')(1-a')^{\nu+\frac32}\right)\,da'\\
		&-k^{-1}\int_{0}^{1}\phi_{1}(a')\left(1-a^{\prime2}\right)^{-\nu-\frac12}a'\,f^{2}_{0}(a')\,da',\\
		h_{2}(1)=&k^{-1}\frac{f_{0}^{1}(1)}{2^{\nu+\frac12}}.
	\end{split}
\end{align}
Finally we obtain
\begin{align}\label{W20 final}
	\begin{split}
		&W_{2}^{0}(a)=h_{0}(a)+h_{1}(a)(1-a)^{\nu+\frac12}+h_{2}(a)(1-a)^{\nu+\frac12}\log(1-a),\\& h_{0}(a)=h_{0}^{1}(a)+h_{0}^{2}(a),\quad h_{1}(a)=h_{1}^{1}(a)+h_{1}^{2}(a),
	\end{split}
\end{align}
and the proof completes.
\end{proof}
\subsection{The equations for the modes $n=0,\pm1$}

\subsubsection{The angular mode $n=1$, equations on the physical side.}
Recall the equation 
\begin{align}\label{eq:n=+1modequation}
	\begin{split}
	&\ddtau \veps^{1}_{+}+H_{1}^{+}\veps_{+}^{1}=F_{+}(1)
	\end{split}
\end{align}
where we may as well specialize to the $+$-case by conjugation symmetry. We also recall the representation 
\begin{align*}
	\veps_{+}^{1}=\phi_{+}\left(\calD_{+}\veps_{+}^{1}\right)+c_{+}(\tau)\phi_{1}(R).
\end{align*}
where we have $\calD_{+}=\partial_{R}+\frac{2R}{1+R^{2}}$, as well as
\begin{align*}
	\phi_{1}(R):=\frac{1}{1+R^{2}}, \quad \phi_{+}(g):=\phi_{1}(R)\int_{0}^{R}\left(\phi_{1}(s)\right)^{-1}g(s)\,ds.
\end{align*}
From \eqref{eq:n=+1modequation} we infer two equations, one for the `better variable' $\mathcal{D}_+\varepsilon_+^1$ by applying $\mathcal{D}_+$ to it, and an ODE for $c_+(\tau)$ by analyzing the vanishing behavior of the original variable $\varepsilon_+^1$ at $R = 0$. To begin with, straightforward differentiation of \eqref{eq:n=+1modequation} and computation of some commutators leads to the following equation for $\mathcal{D}_+\varepsilon_+^1$: 
\begin{align}\label{eq:vareps1eqn}
	\begin{split}
&\ddDtau\calD_{+}\veps_{+}^{1}+\tH_{1}^{+}\calD_{+}\veps_{+}^{1}-\left(2\left(\frac{\lambda'}{\lambda}\right)^{2}+\left(\frac{\lambda'}{\lambda}\right)'\right)\calD_{+}\veps_{+}^{1}\\
=&-2\frac{\lambda'}{\lambda}\frac{4R}{(1+R^{2})^{2}}\left(\partial_{\tau}+\frac{\lambda'}{\lambda}R\partial_{R}\right)\veps_{+}^{1}+\left(\frac{\lambda'}{\lambda}\right)^{2}\left(\frac{4R}{(1+R^{2})^{2}}-\frac{16R}{(1+R^{2})^{3}}\right)\veps_{+}^{1}\\
&-\frac{4R}{(1+R^{2})^{2}}\left(\frac{\lambda'}{\lambda}\right)'\veps_{+}^{1}+\calD_{+}\left(F_{+}(1)\right)\\
=:&\calR_{+}(\veps_{+}^{1},\calD_{+}\veps_{+}^{1})+\calD_{+}\left(F_{+}(1)\right).
	\end{split}
\end{align}
Here we recall that the `super-symmetrical' operator $\tilde{H}_1^{+} = -\mathcal{D}_+\mathcal{D}_{+}^*$, i.e., the factors have been switched compared to $H_1^{+}$. 
The preceding equation gets complemented by one for $c_+(\tau)$, obtained by analyzing the first equation \eqref{eq:n=+1modequation} at $R = 0$, and tracking the value of both sides at $R = 0$: 
\begin{equation}\label{eq:c+evolutioneqn}
	-\partial_{\tau}^2c_+ -  \frac{\lambda_{\tau}}{\lambda}\partial_{\tau}c_+ + \lim_{R\rightarrow 0}H_1^{+}\varepsilon_+^1 =  \lim_{R\rightarrow 0}F_+(1)
\end{equation}
At this point it is crucial to observe that while the source term $F_+(1)$ does have components which depend linearly on $c(\tau)$ and are otherwise independent on the perturbation, namely those terms involving interactions of $\varepsilon_+(1)$ and $\epsilon$, the latter quantity referring to the original blow up profile $Q\left(\lambda(t)r\right) + \epsilon(t,r) = U(t,r)$. These terms are explicitly given by 
\begin{align*}
	&\frac{2\sin\left(2Q+\epsilon\right)\sin\epsilon}{R^2}\varepsilon_+^1,\quad \frac{4\sin\left(Q+\frac{\epsilon}{2}\right)\sin\left(\frac{\epsilon}{2}\right)}{R^2}i\varepsilon_+^1,\\&\left(\frac{2\partial_R\epsilon}{1+R^2} + \left(\partial_R\epsilon\right)^2 - \frac{\lambda'}{\lambda}\cdot\frac{8R}{1+R^2}\left(\partial_{\tau}\epsilon + \frac{\lambda'}{\lambda}R\partial_R\epsilon\right) - \left(\partial_{\tau}\epsilon + \frac{\lambda'}{\lambda}R\partial_R\epsilon\right)^2\right)i\varepsilon_+^1
\end{align*}
Since by their very construction the function $\epsilon$ vanishes to third order at $R = 0$, it is seen that each of these terms vanishes at the origin. This implies that the operator 
\[
\partial_{\tau}^2 + \frac{\lambda_{\tau}}{\lambda}\partial_{\tau}
\]
is responsible to leading order for the evolution of $c_+$. We observe that a fundamental system for this operator consists of the functions $1,\,\tau^{-\nu^{-1}}$. The presence of the function $1$ here distinguishes the mode $n = +1$ from the modes $n =0, n = -1$ treated below, and for which the fundamental system consists of two rapidly decaying functions. Our way to deal with this will invoke modulation theory, by applying a carefully chosen rotation to the bulk part. For now, we relegate this issue for later, and deal with the challenges occurring in the multilinear estimates as well as the iterative step, after translating things to the Fourier side, and making the hypothesis that $c_+(\tau)$ does decay sufficiently toward $\tau = +\infty$, in accordance with the setup in subsection~\ref{subsec:nexcepadmfns}. Forcing this decay assumption will be the role of the final modulation step. 
\\

For later reference, we also observe a solution of the inhomogeneous problem associated to the above differential operator, namely
\[
\partial_{\tau}^2 c + \frac{\lambda_{\tau}}{\lambda}\partial_{\tau} c= h
\]
is solved by 
\begin{equation}\label{eq:cplusinhomsoln}
	c(\tau) = \nu\int_{\tau_0}^{\tau}\sigma h(\sigma)\,d\sigma - \nu\tau^{-\nu^{-1}}\int_{\tau_0}^{\tau}\sigma^{1+\nu^{-1}}h(\sigma)\,d\sigma.
\end{equation}

\subsubsection{Translation of the equation \eqref{eq:vareps1eqn} to the Fourier side}
The case $n = +1$ plays a somewhat special role here as well, since the spectral measure for the operator $\tilde{H}_1^{+}$ is particularly simple, namely a constant multiple of $\xi$, which simplifies both the dilation operator $\mathcal{D}_{\tau,1}$ as well as the formula for the wave parametrix. Even more importantly, the transference operator completely vanishes here, due to the extremely special structure of the Fourier basis. Specifically, we find that 
\begin{align}\label{def Dtau n1}
	\calD_{\tau,1}:=\partial_{\tau}-2\frac{\lambda'}{\lambda}\xi\partial_{\xi}-2\frac{\lambda'}{\lambda}.
\end{align}
Then we translate \eqref{eq:vareps1eqn} into 
\begin{align}\label{eq:Fouriervareps1eqn}
	-\left(\calD_{\tau,1}^{2}+\frac{\lambda'}{\lambda}\calD_{\tau,1}+\xi\right)\xb^{1}=\calF^{(1)}\left(\calR_{+}\left(\veps_{+}^{1},\calD_{+}\veps_{+}^{1}\right)\right)+\calF^{(1)}\left(\calD_{+}\left(F_{+}(\veps)(1)\right)\right).
\end{align}
For the solution of the homogeneous wave equation corresponding to the operator on the left hand side, we have the following 
\begin{lemma}\label{lem:n=1Fourierwavehomogeneous}
	The homogeneous initial value problem 
	\begin{align}\label{parametrix n1 homo eq}
		\left(\calD^{2}_{\tau,1}+\dftau\calD_{\tau,1}+\xi\right)\xb(\tau,\xi)=0;\quad \xb(\tau_{0},\xi)=\xb_{0}(\xi),\quad \calD_{\tau,1}\xb(\tau_{0},\xi)=\xb_{1}(\xi),
	\end{align}
	is solved by
	\begin{align}\label{parametrix n1 homo}
	\begin{split}
	\xb(\tau,\xi)=&\frac{\lambda^{2}(\tau)}{\lambda^{2}(\tau_{0})}\cos\left(\lambda(\tau)\xi^{\frac12}\int_{\tau_{0}}^{\tau}\lambda(u)^{-1}du\right)\xb_{0}\left(\frac{\lambda(\tau)^{2}}{\lambda(\tau_{0})^{2}}\xi\right)\\
	&+\xi^{-\frac12}\frac{\lambda(\tau)}{\lambda(\tau_{0})}\sin\left(\lambda(\tau)\xi^{\frac12}\int_{\tau_{0}}^{\tau}\lambda(u)^{-1}du\right)\xb_{1}\left(\frac{\lambda(\tau)^{2}}{\lambda(\tau_{0})^{2}}\xi\right).
	\end{split}
	\end{align}
This implies the fundamental $S^{(+1)}$-space propagation bounds
\begin{align*}
	\left\| \xb^{1}(\tau,\xi)\right\|_{S_0^{(1)}} + \left\|\mathcal{D}_{\tau,1} \xb^{1}(\tau,\xi)\right\|_{S_1^{(1)}}
	\lesssim \left(\frac{\lambda(\tau_0)}{\lambda(\tau)}\right)^{1-\delta}\cdot \left[\left\|\xb_0\right\|_{S_0^{(1)}} + \left\|\xb_1\right\|_{S_1^{(1)}}\right].
\end{align*}
In particular, choosing $\nu$ small enough, this quantity decays faster than any prescribed negative power of $\tau$. 
	If $\xb(\tau,\xi)$ satisfies
	\begin{align}\label{parametrix n1 inhomo eq}
		\left(\calD_{\tau,1}+\dftau\calD_{\tau,1}+\xi\right)\xb(\tau,\xi)=f(\tau,\xi);\quad \xb(\tau_{0},\xi)=\calD_{\tau,1}\xb(\tau_{0},\xi)=0,
	\end{align}
	then 
	\begin{align}\label{parametrix n1 inhomo}
		\begin{split}
		\xb(\tau,\xi)=&\xi^{-\frac12}\int_{\tau_{0}}^{\tau}\frac{\lambda(\tau)}{\lambda(\sigma)}\sin\left(\lambda(\tau)\xi^{\frac12}\int_{\sigma}^{\tau}\lambda(u)^{-1}du\right)f\left(\sigma,\frac{\lambda(\tau)^{2}}{\lambda(\sigma)^{2}}\xi\right)\,d\sigma.
		\end{split}
	\end{align}
\end{lemma}

 On the other hand, for the inhomogeneous problem, we have a direct analogue of Prop.~\ref{prop:parametrixongoodsource}, but with a subtle difference: 

\begin{proposition}\label{prop:parametrixongoodsourcenless2}
	Let the function $\yb(\tau,\xi)$  be a 'good source function' in accordance with Def.~\ref{defi:xsingulartermsnless2smooth}, Def.~\ref{defi:xsingulartermsnless2proto}. Assume that for the admissibly singular source part the terms with $k = 1, 2, 3$ and $l = 0$ vanish. Denoting by $U^{(1)}(\tau, \sigma,\xi)$ Duhamel propagator for the inhomogeneous problem associated to $ \mathcal{D}_{\tau,1}^2 + \frac{\lambda'}{\lambda} \mathcal{D}_{\tau,1} + \xi$, we have 
	\begin{align*}
		\left\|\int_{\tau_0}^\tau U^{(1)}\left(\tau, \sigma,\xi\right)\cdot \yb\left(\sigma,\frac{\lambda^2(\tau)}{\lambda^2(\sigma)}\xi\right)\,d\sigma\right\|_{\text{good}}\lesssim \left\|\yb\right\|_{\text{goodsource}},
	\end{align*}
	where the norms are of course in the sense of angular momentum $n = 1$ functions. If $\yb$ has principal singular part of restricted type, so does the left hand side, and the norms can be adjusted accordingly. 
\end{proposition}
\begin{proof} The reason for excluding the cases $k = 1, 2, 3$ comes from the fact that we shall gain factors $\big(\frac{\lambda(\sigma)}{\lambda(\tau)}\big)^{k\nu}$, and if $k\geq 4$, these suffice to recover the decay rate $\tau^{-3\nu+}$ postulated in Def.~\ref{defi:xsingulartermsnless2proto}.
	Recall the kernel of the parametrix 
	\begin{align*}
		U^{(1)}(\tau, \sigma,\xi):=\xi^{-\frac12}\frac{\lambda(\tau)}{\lambda(\sigma)}\sin\left(\lambda(\tau)\xi^{\frac12}\int_{\sigma}^{\tau}\lambda(u)^{-1}\,du\right).
	\end{align*}
We first consider the ``smooth part". Assume that the source function $\yb(\tau,\xi)$ is given by
\begin{align*}
	\yb(\tau,\xi)=\tau^{-1}\xi^{\frac12}\zb(\tau,\xi)
\end{align*}
where $\zb(\tau,\cdot)\in S_{0}^{+}$ with appropriate temporal decay. We want to prove that the function
\begin{align*}
\xb(\tau,\xi)=\xi^{-\frac12}\int_{\tau_{0}}^{\tau}\frac{\lambda(\tau)}{\lambda(\sigma)}\sin\left(\lambda(\tau)\xi^{\frac12}\int_{\sigma}^{\tau}\lambda(u)^{-1}\,du\right)\,\yb\left(\sigma,\frac{\lambda^{2}(\tau)}{\lambda^{2}(\sigma)}\xi\right)\,d\sigma
\end{align*}
belongs to $S_{0}^{+}$. To this end, we write $\xb(\tau,\xi)$ as
\begin{align*}
	\xb(\tau,\xi)=\int_{\tau_{0}}^{\tau}\frac{\lambda(\tau)}{\lambda(\sigma)}\sin\left(\lambda(\tau)\xi^{\frac12}\int_{\sigma}^{\tau}\lambda(u)^{-1}\,du\right)\,\sigma^{-1}\frac{\lambda(\tau)}{\lambda(\sigma)}\zb\left(\sigma,\frac{\lambda^{2}(\tau)}{\lambda^{2}(\sigma)}\xi\right)\,d\sigma. 
\end{align*}
So we have
\begin{align*}
	\left\|\xb(\tau,\cdot)\right\|_{S_{0}^{+}}\lesssim& \int_{\tau_{0}}^{\tau}\sigma^{-1}\left(\frac{\lambda(\tau)}{\lambda(\sigma)}\right)^{2}\left\|\xi^{1-\frac{\delta}{2}}\langle\xi\rangle^{\frac32+\delta}\zb\left(\sigma,\frac{\lambda^{2}(\tau)}{\lambda^{2}(\sigma)}\xi\right)\right\|_{L^{2}_{d\xi}}\,d\sigma\\
	\lesssim &\int_{\tau_{0}}^{\tau}\sigma^{-1}\left(\frac{\lambda(\sigma)}{\lambda(\tau)}\right)^{1-\delta}\left\|\zb(\sigma,\cdot)\right\|_{S_{0}^{+}}\,d\sigma\lesssim \sup_{\tau_0\leq\sigma\leq\tau}\left\|\zb(\sigma,\cdot)\right\|_{S_{0}^{+}}.
\end{align*}
Note that the temporal decay of the $S_{0}^{+}$-norm follows from the above estimate. 

Next we consider  the singular part. Assume that $\zb(\tau,\xi)$ is of the form in Definition \ref{defi:xsingulartermsnless2proto}. We first consider the case when $\zb$ is of \emph{prototype}. Let
\begin{align}\label{zb 1 source}
	\begin{split}
	\zb\left(\sigma,\frac{\lambda^{2}(\tau)}{\lambda^{2}(\sigma)}\xi\right)=&\sum_{\pm}\sum_{k=1}^{N}\sum_{i=0}^{N_{1}}\chi_{\frac{\lambda^{2}(\tau)}{\lambda^{2}(\sigma)}\xi\geq1}\frac{e^{\pm i\nu\sigma\frac{\lambda(\tau)}{\lambda(\sigma)}\xi^{\frac12}}}{\left(\frac{\lambda^{2}(\tau)}{\lambda^{2}(\sigma)}\xi\right)^{1+k\frac{\nu}{2}}}\left(\log\frac{\lambda^{2}(\tau)}{\lambda^{2}(\sigma)}\xi\right)^{i}\cdot a_{k,i}(\sigma)\\
	&+\sum_{\pm}\sum_{l=1}^{7}\sum_{k=1}^{N}\sum_{i=0}^{N_{1}}\chi_{\frac{\lambda^{2}(\tau)}{\lambda^{2}(\sigma)}\xi\geq1}\frac{e^{\pm i\nu\sigma\frac{\lambda(\tau)}{\lambda(\sigma)}\xi^{\frac12}}}{\left(\frac{\lambda^{2}(\tau)}{\lambda^{2}(\sigma)}\xi\right)^{1+k\frac{\nu}{2}+\frac{l}{4}}}\left(\log\frac{\lambda^{2}(\tau)}{\lambda^{2}(\sigma)}\xi\right)^{i}\cdot F_{l,k,i}\left(\sigma,\frac{\lambda^{2}(\tau)}{\lambda^{2}(\sigma)}\xi\right)
	\end{split}
\end{align}
with the coefficients $a_{k,i}$ and $F_{l,k,i}$ satisfying the pointwise bounds in Definition \ref{defi:xsingulartermsnless2proto}. To consider the contribution of \eqref{zb 1 source} to the inhomogeneous parametrix, we record the resulting  phase as
\begin{align*}
	e^{\pm i\nu\sigma\frac{\lambda(\tau)}{\lambda(\sigma)}\xi^{\frac12}}\cdot \sin\left(\lambda(\tau)\xi^{\frac12}\int_{\sigma}^{\tau}\lambda(u)^{-1}\,du\right)=\frac{\pm 1}{2i}\left(e^{-i\left(\nu\tau-2\nu\sigma\frac{\lambda(\tau)}{\lambda(\sigma)}\right)\xi^{\frac12}}-e^{i\nu\tau\xi^{\frac12}}\right)
\end{align*}
We start with the 1st term in \eqref{zb 1 source}. In what follows, we omit all the constant coefficients but only focus on  the structures. The contribution from the input phase $e^{\pm i\nu\tau\xi^{\frac12}}$ is given by
\begin{align}\label{inhomo 1 a}
	\begin{split}
	\sum_{\pm}\sum_{k=1}^{N}\sum_{i=0}^{N_{1}}\frac{e^{\pm i\nu\tau\xi^{\frac12}}}{\xi^{1+k\frac{\nu}{2}}}	\int_{\tau_{0}}^{\tau}\sigma^{-1}\frac{\lambda(\tau)}{\lambda(\sigma)}\cdot\left(\frac{\lambda(\sigma)}{\lambda(\tau)}\right)^{1+k\nu}\chi_{\frac{\lambda^{2}(\tau)}{\lambda^{2}(\sigma)}\xi\geq1}\left(\log\frac{\lambda^{2}(\tau)}{\lambda^{2}(\sigma)}\xi\right)^{i}\cdot a_{k,i}(\sigma)\,d\sigma,
	\end{split}
\end{align}
which is decomposed as
\begin{align}\label{inhomo 1 a decom}
	\begin{split}
	&\chi_{\xi\geq1}\sum_{\pm}\sum_{k=1}^{N}\sum_{i=0}^{N_{1}}\frac{e^{\pm i\nu\tau\xi^{\frac12}}}{\xi^{1+k\frac{\nu}{2}}}	\int_{\tau_{0}}^{\tau}\sigma^{-1}\frac{\lambda(\tau)}{\lambda(\sigma)}\cdot\left(\frac{\lambda(\sigma)}{\lambda(\tau)}\right)^{1+k\nu}\chi_{\frac{\lambda^{2}(\tau)}{\lambda^{2}(\sigma)}\xi\geq1}\left(\log\frac{\lambda^{2}(\tau)}{\lambda^{2}(\sigma)}\xi\right)^{i}\cdot a_{k,i}(\sigma)\,d\sigma\\
	&+\left(1-\chi_{\xi\geq1}\right)\sum_{\pm}\sum_{k=1}^{N}\sum_{i=0}^{N_{1}}\frac{e^{\pm i\nu\tau\xi^{\frac12}}}{\xi^{1+k\frac{\nu}{2}}}	\int_{\tau_{0}}^{\tau}\sigma^{-1}\frac{\lambda(\tau)}{\lambda(\sigma)}\cdot\left(\frac{\lambda(\sigma)}{\lambda(\tau)}\right)^{1+k\nu}\chi_{\frac{\lambda^{2}(\tau)}{\lambda^{2}(\sigma)}\xi\geq1}\left(\log\frac{\lambda^{2}(\tau)}{\lambda^{2}(\sigma)}\xi\right)^{i}\cdot a_{k,i}(\sigma)\,d\sigma\\
	=&\chi_{\xi\geq1}\sum_{\pm}\sum_{k=1}^{N}\sum_{i=0}^{N_{1}}\frac{e^{\pm i\nu\tau\xi^{\frac12}}}{\xi^{1+k\frac{\nu}{2}}}	\int_{\tau_{0}}^{\tau}\sigma^{-1}\frac{\lambda(\tau)}{\lambda(\sigma)}\cdot\left(\frac{\lambda(\sigma)}{\lambda(\tau)}\right)^{1+k\nu}\left(\log\frac{\lambda^{2}(\tau)}{\lambda^{2}(\sigma)}\xi\right)^{i}\cdot a_{k,i}(\sigma)\,d\sigma\\
	&-\chi_{\xi\geq1}\sum_{\pm}\sum_{k=1}^{N}\sum_{i=0}^{N_{1}}\frac{e^{\pm i\nu\tau\xi^{\frac12}}}{\xi^{1+k\frac{\nu}{2}}}	\int_{\tau_{0}}^{\tau}\sigma^{-1}\frac{\lambda(\tau)}{\lambda(\sigma)}\cdot\left(\frac{\lambda(\sigma)}{\lambda(\tau)}\right)^{1+k\nu}\left(1-\chi_{\frac{\lambda^{2}(\tau)}{\lambda^{2}(\sigma)}\xi\geq1}\right)\left(\log\frac{\lambda^{2}(\tau)}{\lambda^{2}(\sigma)}\xi\right)^{i}\cdot a_{k,i}(\sigma)\,d\sigma\\
	&+\left(1-\chi_{\xi\geq1}\right)\sum_{\pm}\sum_{k=1}^{N}\sum_{i=0}^{N_{1}}\frac{e^{\pm i\nu\tau\xi^{\frac12}}}{\xi^{1+k\frac{\nu}{2}}}	\int_{\tau_{0}}^{\tau}\sigma^{-1}\frac{\lambda(\tau)}{\lambda(\sigma)}\cdot\left(\frac{\lambda(\sigma)}{\lambda(\tau)}\right)^{1+k\nu}\chi_{\frac{\lambda^{2}(\tau)}{\lambda^{2}(\sigma)}\xi\geq1}\left(\log\frac{\lambda^{2}(\tau)}{\lambda^{2}(\sigma)}\xi\right)^{i}\cdot a_{k,i}(\sigma)\,d\sigma\\
	=:&\xb^{p}_{1}+\xb^{p}_{2}+\xb^{p}_{3}.
	\end{split}
\end{align}
Recall that $\chi$ is  a smooth cutoff, so in the support of $\xb^{p}_{2}$ we have $\xi\simeq 1, \tau\simeq\sigma$, which implies $\xb^{p}_{2}\in S_{0}^{+}$ with desired temporal decay. For $\xb^{p}_{3}$, we rewrite it as:
\begin{align*}
	\left(1-\chi_{\xi\geq1}\right)\sum_{\pm}\sum_{k=1}^{N}\sum_{i=0}^{N_{1}}\frac{e^{\pm i\nu\tau\xi^{\frac12}}}{\xi^{\frac12}}\int_{\tau_{0}}^{\tau}\chi_{\frac{\lambda^{2}(\tau)}{\lambda^{2}(\sigma)}\xi\geq1}\frac{\lambda(\tau)}{\lambda(\sigma)}\frac{1}{\left(\frac{\lambda^{2}(\tau)}{\lambda^{2}(\sigma)}\xi\right)^{\frac12+k\frac{\nu}{2}}}\left(\log\frac{\lambda^{2}(\tau)}{\lambda^{2}(\sigma)}\xi\right)^{i}\cdot\sigma^{-1} a_{k,i}(\sigma)\,d\sigma
\end{align*}
Since the above integral $\int_{\tau_{0}}^{\tau}...$ is finite, we have $\xb^{p}_{3}\in S_{0}^{+}$ at each fixed time $\tau\geq \tau_0$, and is of smooth structured type $\xb_{1,smooth}$ as defined in Def.~\ref{defi:xsingulartermsnless2smooth}. For the term $\xb^{p}_{1}$, we first expand
\begin{align*}
	\left(\log\left(\frac{\lambda^{2}(\tau)}{\lambda^{2}(\sigma)}\xi\right)\right)^{i}=\sum_{k'+l = i}C_{k',l}\left(\log\left(\frac{\lambda^{2}(\tau)}{\lambda^{2}(\sigma)}\right)\right)^{k'}\cdot\left(\log\xi\right)^{l},
\end{align*}
then we define
\begin{align*}
	a_{k,p}(\tau,\sigma):=\sigma^{-1}\left(\frac{\lambda(\sigma)}{\lambda(\tau)}\right)^{k\nu}C_{k',l}\left(\log\frac{\lambda^{2}(\tau)}{\lambda^{2}(\sigma)}\right)^{i-p}\cdot a_{k,i}(\sigma).
\end{align*}

For the output phase $e^{\pm i\nu\left(\tau-2\sigma\frac{\lambda(\tau)}{\lambda(\sigma)}\right)\xi^{\frac12}}$, we obtain an outgoing singular term and we omit the similar details. For the contribution from the other profile of a prototype singular term we again obtain an incoming singular term if the output phase  is $e^{\pm i\nu\tau\xi^{\frac12}}$, and an outgoing singular term if  the output phase  is $e^{\pm i\nu\left(\tau-2\sigma\frac{\lambda(\tau)}{\lambda(\sigma)}\right)\xi^{\frac12}}$. We omit the similar details.

Next we consider the contribution from a \emph{incoming} singular part defined in Definition \ref{defi:xsingulartermsnless2proto}. Assume  that $\zb(\tau,\xi)$ is given by
\begin{align}\label{zb 2 source}
	\begin{split}
	\zb\left(\sigma,\frac{\lambda^{2}(\tau)}{\lambda^{2}(\sigma)}\xi\right):=&\sum_{\pm}\sum_{k=1}^{N}\sum_{j=0}^{N_{1}}\chi_{\frac{\lambda^{2}(\tau)}{\lambda^{2}(\sigma)}\xi\geq1}\frac{e^{\pm i\nu\sigma\frac{\lambda(\tau)}{\lambda(\sigma)}\xi^{\frac12}}}{\left(\frac{\lambda^{2}(\tau)}{\lambda^{2}(\sigma)}\xi\right)^{1+k\frac{\nu}{2}}}\left(\log\frac{\lambda^{2}(\tau)}{\lambda^{2}(\sigma)}\xi\right)^{j}\int_{\tau_{0}}^{\sigma}a^{(\pm)}_{k,j}(\sigma,\sigma_{1})\,d\sigma_{1}\\
	&+\sum_{\pm}\sum_{l=1}^{7}\sum_{k=1}^{N}\sum_{i=0}^{N_{1}}\chi_{\frac{\lambda^{2}(\tau)}{\lambda^{2}(\sigma)}\xi\geq 1}\left\langle\frac{\lambda^{2}(\tau)}{\lambda^{2}(\sigma)}\xi\right\rangle^{-\frac{l}{4}}\frac{e^{\pm i\nu\sigma\frac{\lambda(\tau)}{\lambda(\sigma)}\xi^{\frac12}}}{\left(\frac{\lambda^{2}(\tau)}{\lambda^{2}(\sigma)}\xi\right)^{1+k\frac{\nu}{2}}}\left(\log\frac{\lambda^{2}(\tau)}{\lambda^{2}(\sigma)}\xi\right)^{i}\\
	&\cdot\int_{\tau_{0}}^{\sigma}F^{(\pm)}_{l,k,i}\left(\sigma,\sigma_{1},\frac{\lambda^{2}(\tau)}{\lambda^{2}(\sigma_{1})}\xi\right)\,d\sigma_{1}
	\end{split}
\end{align}
The contribution from the 1st term in \eqref{zb 2 source} can be handled  identically as in the proof of Proposition \ref{prop:parametrixonadmissiblesource} except that now there is  no extra phase $\rho$. Next we give more details for the contribution from the 2nd term in \eqref{zb 2 source}. If the output phase is $e^{\pm i\nu\tau\xi^{\frac12}}$, we have
\begin{align*}
	&\sum_{\pm}\sum_{l=1}^{7}\sum_{k=1}^{N}\sum_{i=0}^{N_{1}}\frac{e^{\pm i\nu\tau\xi^{\frac12}}}{\xi^{\frac12}}\int_{\tau_{0}}^{\tau}\chi_{\frac{\lambda^{2}(\tau)}{\lambda^{2}(\sigma)}\xi\geq 1}\frac{\lambda(\tau)}{\lambda(\sigma)}\left\langle\frac{\lambda^{2}(\tau)}{\lambda^{2}(\sigma)}\xi\right\rangle^{-\frac{l}{4}}\frac{\left(\log\frac{\lambda^{2}(\tau)}{\lambda^{2}(\sigma)}\xi\right)^{i}}{\left(\frac{\lambda^{2}(\tau)}{\lambda^{2}(\sigma)}\xi\right)^{\frac12+k\frac{\nu}{2}}}\cdot\int_{\tau_{0}}^{\sigma}\sigma^{-1}F^{(\pm)}_{l,k,i}\left(\sigma,\sigma_{1},\frac{\lambda^{2}(\tau)}{\lambda^{2}(\sigma_{1})}\xi\right)\,d\sigma_{1}\,d\sigma\\
	=&\sum_{\pm}\sum_{l=1}^{7}\sum_{k=1}^{N}\sum_{i=0}^{N_{1}}\frac{e^{\pm i\nu\tau\xi^{\frac12}}}{\xi^{\frac12}}\int_{\tau_{0}}^{\tau}\left(\int_{\sigma_{1}}^{\tau}\chi_{\frac{\lambda^{2}(\tau)}{\lambda^{2}(\sigma)}\xi\geq 1}\frac{\lambda(\tau)}{\lambda(\sigma)}\left\langle\frac{\lambda^{2}(\tau)}{\lambda^{2}(\sigma)}\xi\right\rangle^{-\frac{l}{4}}\frac{\left(\log\frac{\lambda^{2}(\tau)}{\lambda^{2}(\sigma)}\xi\right)^{i}}{\left(\frac{\lambda^{2}(\tau)}{\lambda^{2}(\sigma)}\xi\right)^{\frac12+k\frac{\nu}{2}}}\cdot\sigma^{-1}F^{(\pm)}_{l,k,i}\left(\sigma,\sigma_{1},\frac{\lambda^{2}(\tau)}{\lambda^{2}(\sigma_{1})}\xi\right)\,d\sigma\right)\,d\sigma_{1}
\end{align*}
As in \eqref{inhomo 1 a decom}, we decompose the above expression into
\begin{align*}
	\chi_{\xi\geq 1}\int_{\tau_{0}}^{\tau}\int_{\sigma_{1}}^{\tau}...\,d\sigma\,d\sigma_{1}-\chi_{\xi\geq1}\int_{\tau_{0}}^{\tau}\int_{\sigma_{1}}^{\tau}\left(1-\chi_{\frac{\lambda^{2}(\tau)}{\lambda^{2}(\sigma)}\xi\geq1}\right)...\,d\sigma\,d\sigma_{1}+\left(1-\chi_{\xi\geq1}\right)\int_{\tau_{0}}^{\tau}\int_{\sigma_{1}}^{\tau}\chi_{\frac{\lambda^{2}(\tau)}{\lambda^{2}(\sigma)}\xi\geq1}...\,d\sigma\,d\sigma_{1}.
\end{align*}
Using identical arguments as for \eqref{inhomo 1 a decom}, we can show that the 2nd and the 3rd terms above belong to $S_{0}^{+}$. The 1st term is given by
\begin{align}\label{inhomo 2 b main}
	\begin{split}
	&\sum_{\pm}\sum_{l=1}^{7}\sum_{k=1}^{N}\sum_{i=0}^{N_{1}}\chi_{\xi\geq1}\frac{e^{\pm i\nu\tau\xi^{\frac12}}}{\xi^{\frac12}}\int_{\tau_{0}}^{\tau}\left(\int_{\sigma_{1}}^{\tau}\frac{\lambda(\tau)}{\lambda(\sigma)}\left\langle\frac{\lambda^{2}(\tau)}{\lambda^{2}(\sigma)}\xi\right\rangle^{-\frac{l}{4}}\frac{\left(\log\frac{\lambda^{2}(\tau)}{\lambda^{2}(\sigma)}\xi\right)^{i}}{\left(\frac{\lambda^{2}(\tau)}{\lambda^{2}(\sigma)}\xi\right)^{\frac12+k\frac{\nu}{2}}}\cdot\sigma^{-1}F^{(\pm)}_{l,k,i}\left(\sigma,\sigma_{1},\frac{\lambda^{2}(\tau)}{\lambda^{2}(\sigma_{1})}\xi\right)\,d\sigma\right)\,d\sigma_{1}\\
	=&\sum_{\pm}\sum_{l=1}^{7}\sum_{k=1}^{N}\sum_{i=0}^{N_{1}}\chi_{\xi\geq1}\frac{e^{\pm i\nu\tau\xi^{\frac12}}}{\xi^{1+k\frac{\nu}{2}}}\left\langle\xi\right\rangle^{-\frac{l}{4}}\\&\cdot\int_{\tau_{0}}^{\tau}\left(\int_{\sigma_{1}}^{\tau}\left\langle\xi\right\rangle^{\frac{l}{4}}\left\langle\frac{\lambda^{2}(\tau)}{\lambda^{2}(\sigma)}\xi\right\rangle^{-\frac{l}{4}}\left(\frac{\lambda(\sigma)}{\lambda(\tau)}\right)^{k\nu}\left(\log\frac{\lambda^{2}(\tau)}{\lambda^{2}(\sigma)}\xi\right)^{i}\cdot\sigma^{-1}F^{(\pm)}_{l,k,i}\left(\sigma,\sigma_{1},\frac{\lambda^{2}(\tau)}{\lambda^{2}(\sigma_{1})}\xi\right)\,d\sigma\right)\,d\sigma_{1}.
	\end{split}
\end{align}
We expand
\begin{align*}
	\left(\log\frac{\lambda^{2}(\tau)}{\lambda^{2}(\sigma)}\xi\right)^{i}=&\sum_{k'+l=i}C_{k',l}\left(\log\left(\frac{\lambda^{2}(\tau)}{\lambda^{2}(\sigma)}\right)\right)^{k'}\cdot\left(\log\xi\right)^{l},
\end{align*}
and define
\begin{align*}
	\tilde{F}^{(\pm)}_{l,k,p}\left(\tau,\sigma_{1},\xi\right):=&\chi_{\xi\geq1}\int_{\sigma_{1}}^{\tau}\left\langle\xi\right\rangle^{\frac{l}{4}}\left\langle\frac{\lambda^{2}(\tau)}{\lambda^{2}(\sigma)}\xi\right\rangle^{-\frac{l}{4}}\left(\frac{\lambda(\sigma)}{\lambda(\tau)}\right)^{k\nu}\left(\log\frac{\lambda^{2}(\tau)}{\lambda^{2}(\sigma)}\right)^{i-p}\cdot\sigma^{-1}F^{(\pm)}_{l,k,i}\left(\sigma,\sigma_{1},\frac{\lambda^{2}(\tau)}{\lambda^{2}(\sigma_{1})}\xi\right)\,d\sigma,
\end{align*}
then we obtain the desired output incoming singular profile. If the output source is $e^{\pm i\nu\left(\tau-2\sigma\frac{\lambda(\tau)}{\lambda(\sigma)}\xi^{\frac12}\right)}$, the principal contribution in consideration is given by
\begin{align*}
	&\sum_{\pm}\sum_{l=1}^{7}\sum_{k=1}^{N}\sum_{i=0}^{N_{1}}\chi_{\xi\geq1}\frac{e^{\pm i\nu\tau\xi^{\frac12}}}{\xi^{1+k\frac{\nu}{2}}}\left\langle\xi\right\rangle^{-\frac{l}{4}}\\&\cdot\int_{\tau_{0}}^{\tau}\left(\int_{\sigma_{1}}^{\tau}e^{\mp 2i\nu\sigma\frac{\lambda(\tau)}{\lambda(\sigma)}\xi^{\frac12}}\left\langle\xi\right\rangle^{\frac{l}{4}}\left\langle\frac{\lambda^{2}(\tau)}{\lambda^{2}(\sigma)}\xi\right\rangle^{-\frac{l}{4}}\left(\frac{\lambda(\sigma)}{\lambda(\tau)}\right)^{k\nu}\left(\log\frac{\lambda^{2}(\tau)}{\lambda^{2}(\sigma)}\xi\right)^{i}\cdot\sigma^{-1}F^{(\pm)}_{l,k,i}\left(\sigma,\sigma_{1},\frac{\lambda^{2}(\tau)}{\lambda^{2}(\sigma_{1})}\xi\right)\,d\sigma\right)\,d\sigma_{1}.
\end{align*}
As in the proof for Proposition \ref{prop:parametrixonadmissiblesource}, we write 
\begin{align*}
	e^{\pm i\nu\tau\xi^{\frac12}}\cdot e^{\mp 2i\nu\sigma\frac{\lambda(\tau)}{\lambda(\sigma)}\xi^{\frac12}}=e^{\mp i\nu\left(\tau+2\sigma\frac{\lambda(\tau)}{\lambda(\sigma)}-2\tau\right)\xi^{\frac12}},
\end{align*}
and introduce the variable
\begin{align*}
	x:=\frac{2\sigma\frac{\lambda(\tau)}{\lambda(\sigma)}-2\tau}{\frac{\lambda(\tau)}{\lambda(\sigma_{1})}}.
\end{align*}
Then set, with $\sigma=\sigma(x,\tau,\sigma_{1})$,
\begin{align*}
	G^{(\pm)}_{l,k,p}(\tau,\sigma_{1},x,\xi):=&\chi_{\xi\geq1}\left\langle\xi\right\rangle^{\frac{l}{4}}\left\langle\frac{\lambda^{2}(\tau)}{\lambda^{2}(\sigma)}\xi\right\rangle^{-\frac{l}{4}}\left(\frac{\lambda(\sigma)}{\lambda(\tau)}\right)^{k\nu}C_{i-p,l}\left(\log\left(\frac{\lambda^{2}(\tau)}{\lambda^{2}(\sigma)}\right)\right)^{i-p}\\
	&\cdot \sigma^{-1}F^{(\pm)}_{l,k,i}\left(\sigma,\sigma_{1},\frac{\lambda^{2}(\tau)}{\lambda^{2}(\sigma_{1})}\xi\right)\cdot\frac{\partial\sigma}{\partial x}\cdot\chi^{(\text{sharp})}_{\sigma\in[\sigma_{1},\tau]},
\end{align*}
we obtain the desired outgoing singular term.

Now we consider the contribution from an \emph{outgoing} singular term. Suppose the source term  $\zb(\tau,\xi)$ is given by
\begin{align}\label{zb 3 source}
	\begin{split}
		\zb\left(\sigma,\frac{\lambda^{2}(\tau)}{\lambda^{2}(\sigma)}\xi\right)=&\sum_{\pm}\sum_{l=0}^{7}\sum_{k=1}^{N}\sum_{i=0}^{N_{1}}\chi_{\frac{\lambda^{2}(\tau)}{\lambda^{2}(\sigma)}\xi\geq1}\left\langle\frac{\lambda^{2}(\tau)}{\lambda^{2}(\sigma)}\xi\right\rangle^{-\frac{l}{4}}\frac{\left(\log\frac{\lambda^{2}(\tau)}{\lambda^{2}(\sigma)}\xi\right)^{i}}{\left(\frac{\lambda^{2}(\tau)}{\lambda^{2}(\sigma)}\xi\right)^{1+k\frac{\nu }{2}}}\\
		&\cdot \int_{\tau_{0}}^{\sigma}e^{\pm i\nu\left(\sigma\frac{\lambda(\tau)}{\lambda(\sigma)}-2\sigma_{1}\frac{\lambda(\tau)}{\lambda(\sigma_{1})}\right)\xi^{\frac12}}\cdot  F^{\pm}_{l,k,i}\left(\sigma,\sigma_{1},\frac{\lambda^{2}(\tau)}{\lambda^{2}(\sigma_{1})}\xi\right)\,d\sigma_{1}.
	\end{split}
\end{align}
Combing with the oscillatory phase, the output phase is of the form:
\begin{align}\label{output phase outgoing}
	e^{\pm i\nu\left(\tau-2\sigma_{1}\frac{\lambda(\tau)}{\lambda(\sigma_{1})}\right)\xi^{\frac12}}\quad \text{and}\quad e^{\pm i\nu\left(\tau-2\sigma\frac{\lambda(\tau)}{\lambda(\sigma)}+2\sigma_{1}\frac{\lambda(\tau)}{\lambda(\sigma_{1})}\right)\xi^{\frac12}}.
\end{align}
Applying the inhomogeneous parametrix and interchanging the order of the integral variables $\sigma$ and $\sigma_{1}$, we again decompose the output expression into
\begin{align*}
	\chi_{\xi\geq1}\int_{\tau_{0}}^{\tau}\int_{\sigma_{1}}^{\tau}...\,d\sigma\,d\sigma_{1}-\chi_{\xi\geq1}\int_{\tau_{0}}^{\tau}\int_{\sigma_{1}}^{\tau}\left(1-\chi_{\frac{\lambda^{2}(\tau)}{\lambda^{2}(\sigma)}\xi\geq1}\right)...\,d\sigma\,d\sigma_{1}+\left(1-\chi_{\xi\geq1}\right)\int_{\tau_{0}}^{\tau}\int_{\sigma_{1}}^{\tau}\chi_{\frac{\lambda^{2}(\tau)}{\lambda^{2}(\sigma)}\xi\geq1}...\,d\sigma\,d\sigma_{1},
\end{align*}
and the 2nd and the 3rd terms are smooth so we only focus on the 1st term which is the principal contribution. For the 1st output phase in \eqref{output phase outgoing}, we have
\begin{align}\label{inhomo 3 b}
	\begin{split}
	&\sum_{\pm}\sum_{l=0}^{7}\sum_{k=1}^{N}\sum_{i=0}^{N_{1}}\chi_{\xi\geq1}\frac{1}{\xi^{1+k\frac{\nu}{2}}}\left\langle\xi\right\rangle^{-\frac{l}{4}}\\
	&\cdot \int_{\tau_{0}}^{\tau}e^{\pm i\nu\left(\tau-2\sigma_{1}\frac{\lambda(\tau)}{\lambda(\sigma_{1})}\right)\xi^{\frac12}}\left(\int_{\sigma_{1}}^{\tau}\langle\xi\rangle^{\frac{l}{4}}\left\langle\frac{\lambda^{2}(\tau)}{\lambda^{2}(\sigma)}\xi\right\rangle^{-\frac{l}{4}}\left(\frac{\lambda(\sigma)}{\lambda(\tau)}\right)^{k\nu}\left(\log\frac{\lambda^{2}(\tau)}{\lambda^{2}(\sigma)}\xi\right)^{i}\cdot\sigma^{-1}F^{(\pm)}\left(\sigma,\sigma_{1},\frac{\lambda^{2}(\tau)}{\lambda^{2}(\sigma_{1})}\xi\right)\,d\sigma\right)\,d\sigma_{1}.
	\end{split}
\end{align}
This can be handled in the same way as \eqref{inhomo 2 b main} and we obtain an outgoing singular term of the form $\xb_{\text{out},1}$. For the 2nd output phase in \eqref{output phase outgoing}, the principal output in consideration is given by
\begin{align}\label{inhomo 3 c main}
	\begin{split}
	\sum_{\pm}\sum_{l=0}^{7}\sum_{k=1}^{N}\sum_{i=0}^{N_{1}}\chi_{\xi\geq1}\frac{1}{\xi^{1+k\frac{\nu}{2}}}\left\langle\xi\right\rangle^{-\frac{l}{4}}\cdot \int_{\tau_{0}}^{\tau}\left(\int_{\sigma_{1}}^{\tau}e^{\pm i\nu\left(\tau-2\sigma\frac{\lambda(\tau)}{\lambda(\sigma)}+2\sigma_{1}\frac{\lambda(\tau)}{\lambda(\sigma_{1})}\right)\xi^{\frac12}}\,...\,d\sigma\right)\,d\sigma_{1}.
	\end{split}
\end{align}
Here the ``$...$" are identical as the integrand in \eqref{inhomo 3 b}. Introducing the new variable
\begin{align*}
	x:=-\frac{2\sigma\frac{\lambda(\tau)}{\lambda(\sigma)}-2\tau}{\frac{\lambda(\tau)}{\lambda(\sigma_{1})}}+2\sigma_{1},
\end{align*}
the output phase function above becomes
\begin{align*}
	\tau-2\sigma\frac{\lambda(\tau)}{\lambda(\sigma)}+2\sigma_{1}\frac{\lambda(\tau)}{\lambda(\sigma_{1})}=\tau+x\frac{\lambda(\tau)}{\lambda(\sigma_{1})}.
\end{align*}
When $\sigma\in[\sigma_{1},\tau]$, $x$ remains non-negative. With this new variable, \eqref{inhomo 3 c main} becomes
\begin{align*}
	\sum_{\pm}\sum_{l=0}^{7}\sum_{k=1}^{N}\sum_{p=0}^{N_{1}}\chi_{\xi\geq1}\langle\xi\rangle^{-\frac{l}{4}}\frac{\left(\log\xi\right)^{p}}{\xi^{1+k\frac{\nu}{2}}}\cdot\int_{0}^{\infty}\int_{\tau_{0}}^{\tau}e^{\pm i\nu\left(\tau+\frac{\lambda(\tau)}{\lambda(\sigma_{1})}x\right)\xi^{\frac12}}\cdot G_{l,k,p}^{\pm}\left(\tau,\sigma_{1},x,\frac{\lambda^{2}(\tau)}{\lambda^{2}(\sigma_{1})}\xi\right)\,d\sigma_{1}\,dx,
\end{align*}
where
\begin{align*}
	G_{l,k,p}^{\pm}\left(\tau,\sigma_{1},x,\frac{\lambda^{2}(\tau)}{\lambda^{2}(\sigma_{1})}\xi\right):=&\left\langle\xi\right\rangle^{\frac{l}{4}}\left\langle\frac{\lambda^{2}(\tau)}{\lambda^{2}(\sigma)}\xi\right\rangle^{-\frac{l}{4}}\left(\frac{\lambda(\sigma)}{\lambda(\tau)}\right)^{k\nu}C_{i-p,l}\left(\log\left(\frac{\lambda^{2}(\tau)}{\lambda^{2}(\sigma)}\right)\right)^{i-p}\\
	&\cdot \sigma^{-1}F^{(\pm)}_{l,k,i}\left(\sigma,\sigma_{1},\frac{\lambda^{2}(\tau)}{\lambda^{2}(\sigma_{1})}\xi\right)\cdot\frac{\partial\sigma}{\partial x}\cdot\chi^{(\text{sharp})}_{\sigma\in[\sigma_{1},\tau]}.
\end{align*}
Here we again think of $\sigma$ as a function of $x,\tau$ and $\sigma_{1}$. Therefore we obtain the desired outgoing singular function.

Finally, the contribution from an outgoing singular source of the form $\tau^{-1}\xi^{\frac12}\xb_{out,2}(\tau,\xi)$ is handled in the same way as Proposition \ref{prop:parametrixonadmissiblesource} and we  omit the details.
\end{proof}

In order to solve the inhomogeneous linear problem \eqref{eq:vareps1eqn}, we have the following analogue of Prop.~\ref{prop:ngeq2fourierwavesoln}:
\begin{proposition}\label{prop:neq1fourierwavesoln} There exists $\tau_0$ large enough, such that if $\mathcal{F}^{(1)}(\mathcal{D}_+F_{\pm}(1))$ is a good source function according to Def.~\ref{defi:xsingulartermsnless2smooth} whose principal singular part is of restricted type and with the $k = 1,2, 3$-components vanishing, with $\mathcal{F}^{(1)}$ denoting the distorted Fourier transform at angular momentum $n = 1$ in the sense of subsection~\ref{subsec:nexcspaces}, then the problem \eqref{eq:vareps1eqn} admits a solution on $[\tau_0,\infty)$ with (distorted Fourier) vanishing data 
	\[
	\left(\overline{x}^{(1)}(\tau_0,\cdot),\,\mathcal{D}_{\tau}\overline{x}^{(1)}(\tau_0,\cdot)\right) = \left(0,0\right)
	\]
	at time $\tau = \tau_0$, and such that (in terms of the distorted Fourier transform of $\varepsilon_+^1$)
	\[
	\left\| \overline{x}\right\|_{\text{good}}\lesssim \left\| \mathcal{F}^{(1)}(F_{\pm}(1))\right\|_{\text{goodsource}} + \big\|\big(c_1(\cdot),\,0\big)\big\|_{good}.
	\]
\end{proposition}
\begin{remark}\label{rem:prop:neq1fourierwavesoln} In order to handle the contribution of the principal singular part with $k = 1, 2, 3$, we shall have to force a certain vanishing condition below by means of suitable modulations. This shall ensure both good bounds for $c_1$ as well as the for $\|\xb\|_{good}$.  
\end{remark}
The proof of this is the same as the one of Prop.~\ref{prop:ngeq2fourierwavesoln}, relying on the fact that the terms linear in $\varepsilon_+$ on the right hand side in \eqref{eq:vareps1eqn} result in terms linear in $c_+(\tau)$ as well as linear terms involving $\mathcal{D}_+\varepsilon$, and which admit a description analogous to the transference operator when expressed on the Fourier side. 

The following is then the direct analogue of Prop.~\ref{prop:ngeq2finalsourcetermestimatesingoodspaces}, and proved the same way: 
\begin{proposition}\label{prop:neq12finalsourcetermestimatesingoodspaces} Assume that $\Lambda\ll 1$ is defined via \eqref{eq:LambdaDef}, and that $\tau_0\gg 1$ is sufficiently large. Then there exists an angular momentum $1$ function
	\[
	\psi^{\pm}_j(1) = \int_0^\infty \xb_j^{(1)}(\tau,\xi)\phi_{1}(R,\xi)\tilde{\rho}_{1}(\xi)\,d\xi,
	\]
	see subsection~\ref{subsec:nexcspaces}, with 
	\[
	\left\| \xb_j^{(1)}\right\|_{\text{good}}\lesssim \Lambda^{3}, 
	\]
	and such that if $F_j^{(\pm)}$ represents any one of the functions occurring in Prop.~\ref{prop:smoothlinearsource}, Prop.~\ref{prop:bilinwithUregular2}, Prop.~\ref{prop:bilinwithUregular3}, and we write 
	\[
	F_j^{(\pm)} = \sum_{n\in \Z}F_j^{(\pm)}(n)e^{in\theta}, 
	\]
	then we have (here $\Box_1$ denotes the wave operator on the left of \eqref{eq:vareps1eqn})
	\begin{align*}
		\left(\mathcal{D}_+F_j^{(\pm)}(1) - \Box_1\psi^{\pm}_j(1)\right)|_{R<\nu\tau} = G_j^{(\pm)}(1)|_{R<\nu\tau}, 
	\end{align*}
	where $G_j^{(\pm)}(1)$ is a good angular momentum $1$ source function and more quantitatively setting 
	\[
	G_j^{(\pm)}(1) = \int_0^\infty \phi_{n}(R,\xi)\yb_1(\tau,\xi)\tilde{\rho}_{1}(\xi)\,d\xi,\,
	\]
	we have 
	\begin{align*}
	\left\| \yb_1(\tau,\xi)\right\|_{\text{goodsource}}\lesssim \left(\tau_0^{-1}+\Lambda\right)\cdot\Lambda\ll\Lambda. 
	\end{align*}
	If we restrict to functions $\xb_n$ with restricted principal singular part, and correspondingly use $\left\| \xb_n\right\|_{\text{good}(r)}$, then $\yb_1$ also has restricted singular principal part, and we may replace 
	$\left\| \yb_1(\tau,\xi)\right\|_{\text{goodsource}}$ by $\left\| \yb_1(\tau,\xi)\right\|_{\text{goodsource}(r)}$. 
\end{proposition}

We next analyze more closely the equation \eqref{eq:c+evolutioneqn} for the coefficient $c_+(\tau)$ of the instability, and specifically the delicate source term 
\begin{equation}\label{eq:delicatesourcetermforcplusode}
	\lim_{R\rightarrow 0}H_1^{+}\varepsilon_+^1
\end{equation}
there. Observe that this term is expected to decay at most like $\mathcal{D}_+\varepsilon_+^1$ 
which means solving the ordinary differential equation \eqref{eq:c+evolutioneqn} will lead to a loss of two additional powers of $\tau$-decay relative to the decay of $\mathcal{D}_+\varepsilon_+^1$. 
\\
To begin with, we render explicit the precise structure of the Fourier transform of the source term $\mathcal{D}_+\left(F_+(1)\right)$ as far as the top order singular terms are concerned, in order to prepare the modulation step. This proposition takes advantage of the fact that all bilinear null-form expressions in the non-linearity take real values, and also of the fact that the structure simplifies since high angular momentum terms which contribute to low angular momentum expressions have to come in pairs. 

\begin{proposition}\label{prop:fourierboundsforn=1termlocalizedaway} Denote by $\mathcal{F}^{(1)}$ the distorted Fourier transform at angular momentum $n = 1$ in the sense of subsection~\ref{subsec:nexcspaces}. Assume that we have \eqref{eq:varepsilonpmangulardecomp} and define $\Lambda$ as before, with $\Lambda\ll 1$. Finally, assume that the distorted Fourier transforms of the $\varepsilon_{\pm}(n)$ (both for the $|n|\geq 2$ and the exceptional modes) have restricted principal singular type (and in particular, we have to define $\Lambda$ by using $\left\|\cdot\right\|_{\text{good}(r)}$). 
	Then we can write
	\begin{equation}\label{eq:obstructionfor=1mode}
		\mathcal{F}^{(1)}\left(\mathcal{D}_+\left(F_+(1)\right)(\sigma, \cdot)\right)(\xi) = \sum_{\pm}\sum_{k=1,2,3}\sum_{j\leq N_1}\frac{e^{\pm i\nu\sigma\xi^{\frac12}}}{\xi^{\frac12+\frac{k\nu}{2}}}\left(\log\xi\right)^j\cdot\beta_{\pm}^{(k,j)}(\sigma) + \yb(\sigma, \xi) + \tilde{\yb}(\sigma, \xi),
	\end{equation}
	where $\yb(\tau, \xi)$ is a good source term at angular momentum $n = 1$ and such that all terms with $k\in \{1,2,3\}, l = 0$ in the expansion of its singular part (according to Definition~\ref{defi:xsingulartermsnless2proto}) vanish. Furthermore, $\tilde{\yb}(\sigma,\cdot)\in \langle \xi^{\frac14}\rangle\cdot S_1^{(1)}$, and with 
	\begin{align*}
	\big\|\tilde{\yb}(\sigma,\cdot)\big\|_{\langle \xi^{\frac14}\rangle\cdot S_1^{(1)}} \lesssim \sigma^{-4+}\cdot \Lambda\cdot(\tau_0^{-1}+\Lambda),\,\sigma\geq \tau_0. 
	\end{align*}
		Finally, we have the bound
	\begin{align*}
		\left|\beta^{(k,j)}(\sigma)\right|\lesssim \sigma^{-4-}\cdot \Lambda\cdot(\tau_0^{-1}+\Lambda),\,\sigma\geq \tau_0 ,
	\end{align*}
	as well as 
	\begin{align*}
	\big\| \yb\big\|_{goodsource}\lesssim  \Lambda\cdot(\tau_0^{-1}+\Lambda).
	\end{align*}
	provided we have the bound $\left|c_0(\tau)\right| + \tau\cdot \left|c_0'(\tau)\right|\lesssim \Lambda\cdot\tau^{-2-}$, which hence improves the generic assumption for the $c_n(\tau)$ in the definition of $\Lambda$. 
	\end{proposition}
\begin{proof} The structure of the leading term (first term on the right) follows from Def.~\ref{defi:xsingulartermsnless2proto}, as well as Prop.~\ref{prop:neq12finalsourcetermestimatesingoodspaces}. Note that for the source functions at exceptional angular momenta, the correction $-\Box_n\psi(n)$ is only required due to the quadratic null-form terms where a factor in $S_0^{(\hbar)}$ is hit by $\partial_\tau + \frac{\lambda_{\tau}}{\lambda}R\partial_R - \partial_R$ while an admissibly singula factor gets hit by a 'bad' derivative $\partial_\tau + \frac{\lambda_{\tau}}{\lambda}R\partial_R + \partial_R$, which leads to terms neither structured nor in $S_1^{(n)}$. For this recall the discussion at the beginning of subsection~\ref{subset:hardnullformestimates}, and also Lemma~\ref{lem:N02bbilin2}. However, the correction $-\Box_n\psi(n)$ is not needed if we accept source terms in $\langle \xi^{\frac14}\rangle\cdot S_1^{(1)}$. We give a more detailed proof of these error estimates in Lemma~\ref{lem:tildec1control}.
\end{proof}
\begin{remark}\label{rem:prop:fourierboundsforn=1termlocalizedaway} The stronger assumption on $c_0(\tau)$ will be ensured via a secondary modulation step below, which will in effect allow us to annihilate $c_0(\tau)$ completely. 
\end{remark}

To get a more detailed description of the coefficients $\beta_{\pm}^{(k,j)}(\sigma) $ of the principal part, we observe that according to Lemma~\ref{lem:singFouriertiphysicaln0adm}, considering a leading singular term in $\mathcal{D}_+\left(F_+(1)\right)(\tau, \cdot)$ of the form 
\begin{equation}\label{eq:leadsing2}
\partial_{R}\left(\left(\nu\sigma-R\right)^{-\frac12+k\nu}\cdot \big[\log(\nu\tau - R)\big]^2\right)\cdot \frac{g_2(\tau)}{\tau^{\frac12}}, 
\end{equation}
which we interpret as 
\begin{align*}
\lim_{\delta\downarrow 0}\partial_{R}\left(\chi_{\nu\tau - R>\delta}\left(\nu\sigma-R\right)^{-\frac12+k\nu}\cdot \big[\log(\nu\tau - R)\big]^2\right)\cdot \frac{g_2(\tau)}{\tau^{\frac12}}
\end{align*}
for a smooth cutoff $\chi_{\nu\tau - R>\delta}$, we obtain in light of \eqref{modified funda sys} (which needs to be multiplied by $R^{-\frac12}$) and integration by parts (which is legitimate for each $\delta>0$) the leading symbol for $\mathcal{F}^{(1)}\left(\mathcal{D}_+\left(F_+(1)\right)(\tau, \cdot)\right)$
\begin{equation}\label{eq:leadsing2Fourier}
\sum_{\pm}\sum_{j=0,1,2}\chi_{\xi\geq 1}\frac{e^{\pm i\nu\tau\xi^{\frac12}}}{\xi^{\frac12+k\nu}}\cdot (\log\xi)^j\cdot \gamma_{k,2-j}^{\pm}\cdot g_2(\tau), 
\end{equation}
where we set 
\begin{align*}
\gamma_{k,2-j}^{\pm} = C_{\pm,k,2-j}\cdot\int_0^\infty e^{\pm ix}\cdot \big(\log x\big)^{2-j}\cdot x^{-\frac12+k\nu}\,dx. 
\end{align*}
for certain non-vanishing constants $C_{\pm,k,2-j}$ independent of $\nu$, and which satisfy $C_{+,k,2-j} = \overline{C_{-,k,2-j}}$.

\begin{remark}\label{rem:fourierboundsforn=1termlocalizedaway} This proposition gives one precise obstruction to obtaining the desired decay for the $c_+(\tau)$-coefficient, as far as the contributions from the source are 
	concerned. In fact, it is precisely the first term on the right in \eqref{eq:obstructionfor=1mode} which a priori leads to a poorly decaying contribution to $c_+(\tau)$, upon application of the wave parametrix stated in Lemma~\ref{lem:n=1Fourierwavehomogeneous}, as will become clear from the next lemma. However, this obstruction lives in a finite dimensional vector space, and suitable choice of modulation parameters below will allow us to eliminate it. 
\end{remark}
\begin{remark}\label{rem:fourierboundsforn=1termlocalizedaway1} The power two of the logarithm has been chosen since in effect for the $k = 1,2,3$ and $l = 0$ singular part the solutions constructed in \cite{KST} we can set $N_1 = 2$, i. e. there is at most a square of a logarithm. This plays no essential role in our method but simplifies explicit calculations. 
\end{remark}

To continue, we now observe the following simple identity which expresses the delicate source term \eqref{eq:delicatesourcetermforcplusode} via the Fourier transform: setting 
\[
\mathcal{D}_+\varepsilon_+^1(\tau, R) = \int_0^\infty \phi_1(R,\xi)\cdot \xb(\tau, \xi)\cdot \tilde{\rho}_{1}(\xi)\,d\xi,
\]
we infer that (see \eqref{D star plus} and \eqref{phi1 small R})
\begin{equation}\label{eq:hardsourceforcplus}
	- \lim_{R\rightarrow 0}\mathcal{D}_+^*\mathcal{D}_+\varepsilon_+^1 = -\frac{\pi}{2}\cdot \int_0^\infty \xb(\tau, \xi)\cdot \tilde{\rho}_{1}(\xi)\,d\xi.
\end{equation}
Let us apply the parametrix from Proposition~\ref{prop:parametrixongoodsourcenless2} to the term \eqref{eq:obstructionfor=1mode}, which is of course only the first approximation to the true Fourier transform of $ \mathcal{D}_+\varepsilon_+^1$, into the preceding integral expression. We first make the following basic observation: 
\begin{lemma}\label{lem:n=1maintermafterparametrix} Denote by $U^{(1)}\left(\tau,\sigma, \xi\right)$ the Duhamel propagator according to Proposition~\ref{prop:parametrixongoodsourcenless2}.
	 Then setting 
	\[
	\xb(\tau, \xi) = \int_{\tau_0}^{\tau}U^{(1)}\left(\tau,\sigma, \xi\right)\cdot  \mathcal{F}^{(1)}\left(\mathcal{D}_+\left(F_+(1)\right)(\sigma, \cdot)\right)\left(\frac{\lambda^2(\tau)}{\lambda^2(\sigma)}\xi\right)\,d\sigma,  
	\]
	we have (for an arbitrary large but finite number $M$ and a smooth cutoff $\chi_{\lambda^2(\tau)\xi<M}$)
	\begin{equation}\label{eq:c_+sourcetermobstruction}\begin{split}
		&\int_0^\infty \chi_{\lambda^2(\tau)\xi<M}\xb(\tau, \xi)\tilde{\rho}_{1}(\xi)\,d\xi\\& = \sum_{\pm}\sum_{k=1,2,3}\sum_{0\leq j_1\leq N_1}\int_0^\infty \chi_{\lambda^2(\tau)\xi<M}\frac{e^{\pm i\nu\tau\xi^{\frac12}}}{\xi^{\frac{k\nu}{2}}}\left[\log\left(\xi\lambda^2(\tau)\right)\right]^{j_1}\,d\xi\\&\hspace{6cm}\cdot \sum_{j_2\leq N_1-j_1} \int_{\tau_0}^{\tau}\left(\log\sigma\right)^{j_2}\left[\frac{\lambda(\sigma)}{\lambda(\tau)}\right]^{\nu k}\beta_{\pm}^{(k,j_1,j_2)}(\sigma)\,d\sigma\\
		& + \tilde{c}_1(\tau), 
	\end{split}\end{equation}
	where the error satisfies 
	\[
	\left|\tilde{c}_1(\tau)\right|\lesssim \tau^{-4-}\cdot  \Lambda\cdot(\tau_0^{-1}+\Lambda),
	\]
	provided we impose the same strengthening on the decay rate of $c_0(\tau)$ as in the preceding proposition. 
\end{lemma}
\begin{proof}
According to the preceding lemma, as well as Definition~\ref{defi:xsingulartermsnless2proto}, Definition~\ref{defi:xsingulartermsnless2smooth}  and the structure of $U^{(1)}\left(\tau,\sigma, \xi\right)$ (see \eqref{parametrix n1 inhomo}), it suffices to show that 
\begin{align*}
\Big|\int_0^\infty \chi_{\lambda^2(\tau)\xi<M}\xb(\tau, \xi)\tilde{\rho}_{1}(\xi)\,d\xi\Big|\lesssim  \tau^{-4-}\cdot  \Lambda\cdot(\tau_0^{-1}+\Lambda)
\end{align*}
provided $\xb$ is either of incoming singular type with $k\geq 4$, of outgoing singular type, of structured smooth type $\xb_{1,smooth}$ or in $\langle\xi^{\frac14}\rangle S_0^{(1)}$, all arising from the source term $\mathcal{D}_+\left(F_+(1)\right)$. For the first three of these, we observe that the $\xi$-integral gains $\tau^{-(2-)}$ since we can essentially restrict to $\tau\xi^{\frac12}\lesssim 1$ by using integration by parts in $\xi$. For the remaining situation where $\xb(\tau,\cdot)\in \langle\xi^{\frac14}\rangle S_0^{(1)}$, we note that the $\xi$-integral converges since 
\begin{align*}
\Big|\int_0^\infty \chi_{\lambda^2(\tau)\xi<M}\xb(\tau, \xi)\tilde{\rho}_{1}(\xi)\,d\xi\Big|\lesssim \big\|\xb(\tau, \cdot)\big\|_{\langle\xi^{\frac14}\rangle S_0^{(1)}}
\end{align*}
uniformly in $M$, in light of \eqref{eq:S01}. To get the temporal decay $ \tau^{-4-}$, we refer to Lemma~\ref{lem:tildec1control}. 
\end{proof}
The preceding lemma shows that the obstruction to obtaining a $\tau^{-4-}$-decay rate for the source term \eqref{eq:hardsourceforcplus} for the evolution equation \eqref{eq:c+evolutioneqn} comes from the first term on the right in \eqref{eq:c_+sourcetermobstruction}. However, assuming for each $j_1\in \{0,\ldots, N_1\}$ the vanishing relation 
\begin{equation}\label{eq:n=1vanishingrelations}
\sum_{j_2\leq N_1-j_1}\int_{\tau_0}^{\infty}\left(\log\sigma\right)^{j_2}\left[\lambda(\sigma)\right]^{\nu k}\beta_{\pm}^{(k,j_1,j_2)}(\sigma)\,d\sigma = 0,
\end{equation}
we infer the decay rate $\tau^{-5+}$ for this term. 
We shall achieve the desired vanishing relations by exploiting suitable rotations on the target sphere and the effect such rotations have on the singular part on the light cone. 
\\
We can make the structure of the functions $\beta_{\pm}^{(k,j_1,j_2)}(\sigma)$ more transparent for the specific source term \eqref{eq:leadsing2}. In fact, using the same notation as for \eqref{eq:leadsing2Fourier}, we find (with $N_1 = 2$ compared to the above)
\begin{equation}\label{eq:leadsing2Duhamelspecialcoefficients}
\beta_{\pm}^{(k,j_1,j_2)}(\sigma) = C_{\pm,k,j_1,j_2}\cdot \gamma^{\pm}_{k, 2-j_1 - j_2}\cdot g_2(\sigma). 
\end{equation}
where $C_{\pm,k,j_1,j_2}$ are suitable non-vanishing (complex) constants independent of $\nu$. 

 \subsubsection{The angular mode $n = 0$, equations on the physical side}. 

Recall the equation
\begin{equation}\label{eq:n=+0modequation}
	-\left(\partial_{\tau}+\frac{\lambda_{\tau}}{\lambda}R\partial_R\right)^2\varepsilon_+^0 - \frac{\lambda_{\tau}}{\lambda}\left(\partial_{\tau}+\frac{\lambda_{\tau}}{\lambda}R\partial_R\right)\varepsilon_+^0 + H_0^{+}\varepsilon_+^0 = F_+(0), 
\end{equation}
where we also focus on the $+$-case. Recall the representation 
\[
\varepsilon_+^0(\tau, R) = c_0(\tau)\phi_0(R) + \phi_0(R)\cdot\int_0^R[\phi_0(s)]^{-1}\mathcal{D}_0\varepsilon_+^0(\tau, s)\,ds, 
\]
where we set $\mathcal{D}_0 = \partial_R + \frac{1}{R}\cdot\frac{R^2-1}{R^2+1}$, and we recall $\phi_0(R) = \frac{R}{1+R^2}$. Commuting $\mathcal{D}_0$ past the equation results in 
\begin{align}\label{eq:vareps0eqn}
	\begin{split}
		&-\left(\left(\partial_{\tau}+\frac{\lambda'}{\lambda}R\partial_{R}\right)^{2}+3\frac{\lambda'}{\lambda}\left(\partial_{\tau}+\frac{\lambda'}{\lambda}R\partial_{R}\right)\right)\calD_0\veps_+^0 + \tilde{H}_0^{+}\calD_0\veps_+^0\\
		=&\calD_0\left(F_+(0)\right)-\frac{4R}{(R^{2}+1)^{2}}\left(2\left(\frac{\lambda'}{\lambda}\right)^{2}+\left(\frac{\lambda'}{\lambda}\right)'\right)\veps_+^0\\
		&-\frac{\lambda'}{\lambda}\frac{4R}{(R^{2}+1)^{2}}\left(\partial_{\tau}+\frac{\lambda'}{\lambda}R\partial_{R}\right)\veps_+^0-\frac{\lambda'}{\lambda}\left(\partial_{\tau}+\frac{\lambda'}{\lambda}R\partial_{R}\right)\left(\frac{4R}{(R^{2}+1)^{2}}\veps_+^0\right)\\
		&+\left(2\left(\frac{\lambda'}{\lambda}\right)^{2}+\left(\frac{\lambda'}{\lambda}\right)'\right)\calD\veps_+^0\\
		=:&\calD_0\left(F_+(0)\right)+\calR^0(\veps_+^0,\calD\veps_+^0)+\left(2\left(\frac{\lambda'}{\lambda}\right)^{2}+\left(\frac{\lambda'}{\lambda}\right)'\right)\calD\veps_+^0.
	\end{split}
\end{align}
The operator $\tilde{H}_0^{+}$ is the super-symmetric cousin of the operator $H_0^{+}$, i.e., $\tilde{H}_0^{+} = -\calD_0\calD_0^*$. To complete things, we analyze \eqref{eq:n=+0modequation} around $R = 0$ to extract the evolution law governing the coefficient $c_0(\tau)$ of the instability, which becomes the following:
\begin{equation}\label{eq:c0evolutioneqn}
	-\left(\partial_{\tau} + \frac{\lambda'}{\lambda}\right)^2c_+ -  \frac{\lambda_{\tau}}{\lambda}\left(\partial_{\tau}+ \frac{\lambda'}{\lambda}\right)c_+ + \lim_{R\rightarrow 0}R^{-1}H_0^{+}\varepsilon_+^0 =  \lim_{R\rightarrow 0}R^{-1}F_+(0)
\end{equation}
This time the key operator governing the evolution of the nonlinearity is given by 
\[
\left(\partial_{\tau} + \frac{\lambda'}{\lambda}\right)^2 + \frac{\lambda_{\tau}}{\lambda}\left(\partial_{\tau}+ \frac{\lambda'}{\lambda}\right)
\]
with fundamental system given by $\tau^{-1-\nu^{-1}},\,\tau^{-1 - 2\nu^{-1}}$. Moreover, we can solve the equation 
\[
\left(\partial_{\tau} + \frac{\lambda'}{\lambda}\right)^2c + \frac{\lambda_{\tau}}{\lambda}\left(\partial_{\tau}+ \frac{\lambda'}{\lambda}\right)c = h
\]
by means of 
\begin{equation}\label{eq:czeroinhomsoln}
	c(\tau) = \nu\left(\tau^{-1-\nu^{-1}}\int_{\tau_0}^{\tau}\sigma^{2+\nu^{-1}} h(\sigma)\,d\sigma - \tau^{-1-2\nu^{-1}}\int_{\tau_0}^{\tau}\sigma^{2+2\nu^{-1}}h(\sigma)\,d\sigma\right). 
\end{equation}

\subsubsection{Translation of the equation \eqref{eq:vareps0eqn} to the Fourier side} 

From an algebraic standpoint, the situation here is more complicated than the preceding case, since the transference operator is non-vanishing for $n = 0$, and the spectral measure is only implicit. Introduce the auxiliary operator 
\[
\mathcal{D}_{\tau}^{(0)} = \partial_{\tau} - 2\frac{\lambda_{\tau}}{\lambda}\xi\partial_{\xi} - \frac{\lambda_{\tau}}{\lambda}\frac{\left(\tilde{\rho}_{0}(\xi)\right)'\xi}{\tilde{\rho}_{0}(\xi)} - \frac{\lambda_{\tau}}{\lambda}. 
\]
By means of it, we translate the equation \eqref{eq:vareps0eqn} to the Fourier variables $\xb(\tau, \xi)$, which satisfies 
\[
\calD_{0}\veps_+^0(\tau, R) = \int_0^\infty \phi_0(R,\xi)\cdot \xb(\tau, \xi)\cdot \tilde{\rho}_{0}(\xi)\,d\xi,
\]
in perfect analogy to \eqref{eq:ngeq2Fourier1} as follows: 
\begin{equation}\label{eq:neq0Fourier1}\begin{split}
		&-\left(\left(\mathcal{D}_{\tau}^{(0)}\right)^2 + \frac{\lambda'(\tau)}{\lambda(\tau)}\mathcal{D}_{\tau}^{(0)} + \xi \right)\xb\\
		& =  \mathcal{F}^{(0)}\left(F_\pm(0)\right) +  2\frac{\lambda'(\tau)}{\lambda(\tau)}\mathcal{K}_{0}^{(0)}\mathcal{D}_{\tau}^{(0)}\xb + \left(\frac{\lambda'(\tau)}{\lambda(\tau)}\right)'\mathcal{K}_{0}^{(0)}\xb + \frac{\lambda'(\tau)}{\lambda(\tau)}\left[\mathcal{D}_{\tau}^{(0)},\mathcal{K}_{0}^{(0)}\right]\xb + \left(\frac{\lambda'(\tau)}{\lambda(\tau)}\right)^2\left(\left(\mathcal{K}_{0}^{(0)}\right)^2 + \mathcal{K}_{0}^{(0)}\right)\xb \\
		& +  \mathcal{F}^{(0)}\left(\calR^0\left(\veps_+^0,\calD\veps_+^0\right)\right)
\end{split}\end{equation}
where $\mathcal{K}_0^{(0)}$ is the off-diagonal part of the transference operator associated with $\phi_0(R,\xi)$. 
\\
In analogy to Lemma~\ref{lem:n=1Fourierwavehomogeneous}, we have  
\begin{lemma}\label{lem:n=0Fourierwavehomogeneous} The homogeneous initial value problem 
	\[
	\left(\left(\mathcal{D}_{\tau}^{(0)}\right)^2 + \frac{\lambda'}{\lambda} \mathcal{D}_{\tau}^{(0)} + \xi\right)\overline{x}(\tau,\xi) = 0,\quad \overline{x}(\tau_0,\xi) = \xb_0(\xi),\quad \mathcal{D}_{\tau}^{(0)}\overline{x}(\tau_0,\xi) = \xb_1(\xi),
	\]
	is solved by the function 
	\begin{align*}
		\overline{x}(\tau,\xi) &= \frac{\lambda(\tau)}{\lambda(\tau_0)}\cdot\frac{\left(\tilde{\rho}_{0}\right)^{\frac12}\left(\frac{\lambda^2(\tau)}{\lambda^2(\tau_0)}\xi\right)}{\left(\tilde{\rho}_{0}\right)^{\frac12}(\xi)}\cdot \cos\left[\lambda(\tau)\xi^{\frac12}\int_{\tau_0}^{\tau}\lambda^{-1}(u)\,du\right]\cdot \xb_0\left(\frac{\lambda^2(\tau)}{\lambda^2(\tau_0)}\xi\right)\\
		& +\frac{\left(\tilde{\rho}_{0}\right)^{\frac12}\left(\frac{\lambda^2(\tau)}{\lambda^2(\tau_0)}\xi\right)}{\left(\tilde{\rho}_{0}\right)^{\frac12}(\xi)}\cdot \xi^{-\frac12}\sin\left[\lambda(\tau)\xi^{\frac12}\int_{\tau_0}^{\tau}\lambda^{-1}(u)\,du\right]\cdot \xb_1\left(\frac{\lambda^2(\tau)}{\lambda^2(\tau_0)}\xi\right)\\
	\end{align*}
	This implies the fundamental $S^{(0)}$-space propagation bounds(recall the definition \eqref{eq:S00} of this norm)
	\begin{align*}
		\left\| \overline{x}(\tau,\xi)\right\|_{S_0^{(0)}} + \left\|\mathcal{D}_{\tau}^{(0)} \overline{x}(\tau,\xi)\right\|_{S_1^{(0)}}
		\lesssim \left(\frac{\lambda(\tau_0)}{\lambda(\tau)}\right)^{1-\delta}\cdot \left[\left\|\xb_0\right\|_{S_0^{(0)}} + \left\|\xb_1\right\|_{S_1^{(0)}}\right]
	\end{align*}
	In particular, choosing $\nu$ small enough, this quantity decays faster than any prescribed negative power of $\tau$. 
\end{lemma}
\begin{lemma}\label{lem:n=0Fourierwaveinhomogeneous}
	The inhomogeneous initial value problem
		\[
	\left(\left(\mathcal{D}_{\tau}^{(0)}\right)^2 + \frac{\lambda'}{\lambda} \mathcal{D}_{\tau}^{(0)} + \xi\right)\overline{x}(\tau,\xi) = f(\tau,\xi),\quad \overline{x}(\tau_0,\xi) = 0,\quad \mathcal{D}_{\tau}^{(0)}\overline{x}(\tau_0,\xi) = 0,
	\]
	is solved by the function
	\begin{align*}
		\xb(\tau,\xi)=&\int_{\tau_{0}}^{\tau}\frac{\left(\tilde{\rho}_{0}\right)^{\frac12}\left(\frac{\lambda^2(\tau)}{\lambda^2(\sigma)}\xi\right)}{\left(\tilde{\rho}_{0}\right)^{\frac12}(\xi)}\cdot\xi^{-\frac12}\sin\left[\lambda(\tau)\xi^{\frac12}\int_{\sigma}^{\tau}\lambda^{-1}(u)\,du\right]\cdot f\left(\sigma,\frac{\lambda^{2}(\tau)}{\lambda^{2}(\sigma)}\xi\right)\,d\sigma.
	\end{align*}
\end{lemma}
Furthermore, recalling Definition~\ref{defi:goodfourierrepnless2} , we have the following analogue of Proposition \ref{prop:parametrixongoodsourcenless2}:
\begin{proposition}\label{prop:parametrixongoodsourcenless20} 
	Given a ``good source" function $\yb(\tau,\xi)$  such that for the admissibly singular source part the terms with $l=0$ and $k = 1, 2, 3$ vanish, and letting $U^{(0)}\left(\tau, \sigma,\xi\right)$ denote the Duhamel propagator for the inhomogeneous problem associated to $ \left(\mathcal{D}_{\tau}^{(0)}\right)^2 + \frac{\lambda'}{\lambda} \mathcal{D}_{\tau}^{(0)} + \xi$, we have 
	\begin{align*}
		\left\|\int_{\tau_0}^\tau U^{(0)}\left(\tau, \sigma,\xi\right)\cdot \yb\left(\sigma,\frac{\lambda^2(\tau)}{\lambda^2(\sigma)}\xi\right)\,d\sigma\right\|_{\text{good}}\lesssim \left\|\yb\right\|_{\text{goodsource}},
	\end{align*}
	where the norms are of course in the sense of angular momentum $n = 0$ functions. If $\yb$ has principal singular part of restricted type, so does the left hand side, and the norms can be adjusted accordingly. 
\end{proposition}
\begin{proof}
	Recall the expression of the kernel $U^{(0)}(\tau,\sigma,\xi)$:
	\begin{align*}
		U^{(0)}(\tau,\sigma,\xi)=\frac{\left(\tilde{\rho}_{0}\right)^{\frac12}\left(\frac{\lambda^2(\tau)}{\lambda^2(\sigma)}\xi\right)}{\left(\tilde{\rho}_{0}\right)^{\frac12}(\xi)}\cdot\xi^{-\frac12}\sin\left[\lambda(\tau)\xi^{\frac12}\int_{\sigma}^{\tau}\lambda^{-1}(u)\,du\right].
	\end{align*}
	We start with the smooth part. Assume the source function $\yb(\tau,\xi)$ is of the form $\yb(\tau,\xi)=\tau^{-1}\xi^{\frac12}\zb(\tau,\xi)$ where $\zb(\tau,\cdot)\in S_{0}^{0}$, and we want to prove that
	\begin{align*}
		\xb(\tau,\xi)=\int_{\tau_{0}}^{\tau}U^{(0)}(\tau,\sigma,\xi)\yb\left(\sigma,\frac{\lambda^{2}(\tau)}{\lambda^{2}(\sigma)}\xi\right)\,d\sigma
	\end{align*}
belongs to $S_{0}^{0}$, for which we recall its definition:
\begin{align*}
	\left\|\xb(\xi)\right\|_{S_{0}^{0}}:=\left\|\left(\frac{\langle\log\langle\xi\rangle\rangle}{\langle\log\xi\rangle}\right)^{1+\delta}\xi^{\frac12}\langle\xi\rangle^{\frac52+\frac{\delta}{2}}\xb(\xi)\right\|_{L^{2}_{d\xi}}.
\end{align*}
	We first investigate the behavior of $\frac{\left(\tilde{\rho}_{0}\right)^{\frac12}\left(\frac{\lambda^2(\tau)}{\lambda^2(\sigma)}\xi\right)}{\left(\tilde{\rho}_{0}\right)^{\frac12}(\xi)}$. If $\underline{\frac{\lambda^{2}(\tau)}{\lambda^{2}(\sigma)}\xi\lesssim 1}$, we have
	\begin{align*}
		\frac{\left(\tilde{\rho}_{0}\right)^{\frac12}\left(\frac{\lambda^2(\tau)}{\lambda^2(\sigma)}\xi\right)}{\left(\tilde{\rho}_{0}\right)^{\frac12}(\xi)}\simeq\frac{\langle\log\xi\rangle}{\left\langle\log\frac{\lambda^{2}(\tau)}{\lambda^{2}(\sigma)}\xi\right\rangle}
	\end{align*}
and in this case we have
\begin{align*}
	\frac{\left(\tilde{\rho}_{0}\right)^{\frac12}\left(\frac{\lambda^2(\tau)}{\lambda^2(\sigma)}\xi\right)}{\left(\tilde{\rho}_{0}\right)^{\frac12}(\xi)}\cdot\left(\frac{\langle\log\langle\xi\rangle\rangle}{\langle\log\xi\rangle}\right)^{1+\delta}\lesssim \frac{\langle\log\langle\xi\rangle\rangle}{\left\langle\log\frac{\lambda^{2}(\tau)}{\lambda^{2}(\sigma)}\xi\right\rangle}\lesssim 1,
\end{align*}
and
\begin{align*}
	\left\|\xb(\tau,\cdot)\right\|_{S_{0}^{0}}\lesssim&\int_{\tau_{0}}^{\tau}\sigma^{-1}\frac{\lambda(\tau)}{\lambda(\sigma)}\left\|\xi^{\frac12}\zb\left(\sigma,\frac{\lambda^{2}(\tau)}{\lambda^{2}(\sigma)}\xi\right)\right\|_{L^{2}_{d\xi}}\,d\sigma\\
	\lesssim &\int_{\tau_{0}}^{\tau}\sigma^{-1}\frac{\lambda(\sigma)}{\lambda(\tau)}\cdot\left\|\zb(\sigma,\cdot)\right\|_{S_{0}^{0}}\,d\sigma\lesssim \sup_{\sigma\in[\tau_{0},\tau]}\left\|\zb(\sigma,\cdot)\right\|_{S_{0}^{0}}.
\end{align*}
Next we turn to the case when \underline{$\xi\lesssim 1$ and $\frac{\lambda^{2}(\tau)}{\lambda^{2}(\sigma)}\xi\gg1$}. In this case we have
\begin{align*}
	\frac{\left(\tilde{\rho}_{0}\right)^{\frac12}\left(\frac{\lambda^2(\tau)}{\lambda^2(\sigma)}\xi\right)}{\left(\tilde{\rho}_{0}\right)^{\frac12}(\xi)}\simeq\frac{\lambda^{2}(\tau)}{\lambda^{2}(\sigma)}\xi\cdot\langle\log\xi\rangle,
\end{align*}
and 
\begin{align*}
	\frac{\left(\tilde{\rho}_{0}\right)^{\frac12}\left(\frac{\lambda^2(\tau)}{\lambda^2(\sigma)}\xi\right)}{\left(\tilde{\rho}_{0}\right)^{\frac12}(\xi)}\cdot\left(\frac{\langle\log\langle\xi\rangle\rangle}{\langle\log\xi\rangle}\right)^{1+\delta}\lesssim  \frac{\lambda^{2}(\tau)}{\lambda^{2}(\sigma)}\xi\cdot\langle\log\xi\rangle^{-\delta},
\end{align*}
which implies
\begin{align*}
	\|\xb(\tau,\cdot)\|_{S_{0}^{0}}\lesssim&\int_{\tau_{0}}^{\tau}\sigma^{-1}\frac{\lambda^{2}(\tau)}{\lambda^{2}(\sigma)}\cdot\frac{\lambda(\tau)}{\lambda(\sigma)}\left\|\xi\langle\log\xi\rangle^{-\delta}\cdot \xi^{\frac12}\zb\left(\sigma,\frac{\lambda^{2}(\tau)}{\lambda^{2}(\sigma)}\xi\right)\right\|_{L^{2}_{d\xi}}\,d\sigma\\
	\lesssim&\int_{\tau_{0}}^{\tau}\sigma^{-1}\frac{\lambda(\sigma)}{\lambda(\tau)}\left\|\zb(\sigma,\cdot)\right\|_{S_{0}^{0}}\,d\sigma\lesssim \sup_{\sigma\in[\tau_{0},\tau]}\left\|\zb(\sigma,\cdot)\right\|_{S_{0}^{0}}.
\end{align*}
Finally we  consider the case  \underline{$\xi\gg1$}. In this case we have
\begin{align*}
	\frac{\left(\tilde{\rho}_{0}\right)^{\frac12}\left(\frac{\lambda^2(\tau)}{\lambda^2(\sigma)}\xi\right)}{\left(\tilde{\rho}_{0}\right)^{\frac12}(\xi)}\simeq\frac{\lambda^{2}(\tau)}{\lambda^{2}(\sigma)}\quad \Rightarrow\quad \frac{\left(\tilde{\rho}_{0}\right)^{\frac12}\left(\frac{\lambda^2(\tau)}{\lambda^2(\sigma)}\xi\right)}{\left(\tilde{\rho}_{0}\right)^{\frac12}(\xi)}\cdot\left(\frac{\langle\log\langle\xi\rangle\rangle}{\langle\log\xi\rangle}\right)^{1+\delta}\simeq \frac{\lambda^{2}(\tau)}{\lambda^{2}(\sigma)},
\end{align*}
and
\begin{align*}
		\|\xb(\tau,\cdot)\|_{S_{0}^{0}}\lesssim&\int_{\tau_{0}}^{\tau}\sigma^{-1}\frac{\lambda^{2}(\tau)}{\lambda^{2}(\sigma)}\cdot\frac{\lambda(\tau)}{\lambda(\sigma)}\left\|\xi^{3+\frac{\delta}{2}}\zb\left(\sigma,\frac{\lambda^{2}(\tau)}{\lambda^{2}(\sigma)}\xi\right)\right\|_{L^{2}_{d\xi}}\,d\sigma\\
		\lesssim&\int_{\tau_{0}}^{\tau}\sigma^{-1}\left(\frac{\lambda(\sigma)}{\lambda(\tau)}\right)^{4+\delta}\left\|\zb(\sigma,\cdot)\right\|_{S_{0}^{0}}\,d\sigma\lesssim \sup_{\sigma\in[\tau_{0},\tau]}\left\|\zb(\sigma,\cdot)\right\|_{S_{0}^{0}}.
\end{align*}
This  completes the discussion for the smooth part.

The discussion for the singular part is similar to that of $U^{(1)}(\tau,\sigma,\xi)$. We only need to pay attention on how to extract the principal contribution from the output. Let us consider a prototypical singular input for example. Suppose the source is given by $\tau^{-1}\xi^{\frac12}\zb(\tau,\xi)$ where $\zb(\tau,\xi)$ as source is given by
\begin{align*}
	\zb\left(\sigma,\frac{\lambda^{2}(\tau)}{\lambda^{2}(\sigma)}\xi\right)=\sum_{\pm}\sum_{k=1}^{N}\sum_{i=0}^{N_{1}}\chi_{\frac{\lambda^{2}(\tau)}{\lambda^{2}(\sigma)}\xi\geq1}\frac{e^{\pm i\nu\sigma\frac{\lambda(\tau)}{\lambda(\sigma)}\xi^{\frac12}}}{\left(\frac{\lambda^{2}(\tau)}{\lambda^{2}(\sigma)}\xi\right)^{\frac32+k\frac{\nu}{2}}}\left(\log\frac{\lambda^{2}(\tau)}{\lambda^{2}(\sigma)}\xi\right)^{i}\cdot a_{k,i}(\sigma).
\end{align*}
Applying the parametrix, we obtain (for incoming output, for instance)
\begin{align}\label{zb 0 source}
	\xb(\tau,\xi)=\xi^{-\frac12}\int_{\tau_{0}}^{\tau}\frac{\left(\tilde{\rho}_{0}\right)^{\frac12}\left(\frac{\lambda^2(\tau)}{\lambda^2(\sigma)}\xi\right)}{\left(\tilde{\rho}_{0}\right)^{\frac12}(\xi)}\sum_{\pm}\sum_{k=1}^{N}\sum_{i=0}^{N_{1}}\chi_{\frac{\lambda^{2}(\tau)}{\lambda^{2}(\sigma)}\xi\geq1}\frac{e^{\pm i\nu\tau\xi^{\frac12}}}{\left(\frac{\lambda^{2}(\tau)}{\lambda^{2}(\sigma)}\xi\right)^{\frac 32+k\frac{\nu}{2}}}\left(\log\frac{\lambda^{2}(\tau)}{\lambda^{2}(\sigma)}\xi\right)^{i}\cdot a_{k,i}(\sigma)\,d\sigma.
\end{align}
We again decompose the above integral according to
\begin{align}\label{zb 0 source decom}
	\chi_{\xi\geq1}\int_{\tau_{0}}^{\tau}...\,d\sigma-\chi_{\xi\geq1}\int_{\tau_{0}}^{\tau}\left(1-\chi_{\frac{\lambda^{2}(\tau)}{\lambda^{2}(\sigma)}\xi\geq1}\right)...\,d\sigma+\left(1-\chi_{\xi\geq1}\right)\int_{\tau_{0}}^{\tau}\chi_{\frac{\lambda^{2}(\tau)}{\lambda^{2}(\sigma)}\xi\geq1}...\,d\sigma
\end{align}
The 1st integral  in \eqref{zb 0 source decom} is the principal part and handled in the same way the principal part in the case $n=1$. For the 2nd and the 3rd terms, to show that they belong to $S_{0}^{0}$, we need to take care of the temporal decay. To this end, we only need to show that the factor $\frac{\left(\tilde{\rho}_{0}\right)^{\frac12}\left(\frac{\lambda^2(\tau)}{\lambda^2(\sigma)}\xi\right)}{\left(\tilde{\rho}_{0}\right)^{\frac12}(\xi)}$ gives a temporal growth of $\left(\frac{\lambda(\tau)}{\lambda(\sigma)}\right)^{2}$ at most. For the 2nd term in \eqref{zb 0 source decom},  we have $\xi\simeq1$ in this regime, so
\begin{align*}
	\frac{\left(\tilde{\rho}_{0}\right)^{\frac12}\left(\frac{\lambda^2(\tau)}{\lambda^2(\sigma)}\xi\right)}{\left(\tilde{\rho}_{0}\right)^{\frac12}(\xi)}\simeq\frac{\left(\log\xi\right)^{2}}{\left(\log\frac{\lambda^{2}(\tau)}{\lambda^{2}(\sigma)}\xi\right)^{2}}\lesssim 1. 
\end{align*}
For the 3rd term in \eqref{zb 0 source decom}, we have, since $\xi\leq 1$,
\begin{align*}
	\frac{\left(\tilde{\rho}_{0}\right)^{\frac12}\left(\frac{\lambda^2(\tau)}{\lambda^2(\sigma)}\xi\right)}{\left(\tilde{\rho}_{0}\right)^{\frac12}(\xi)}\simeq\frac{\lambda^{2}(\tau)}{\lambda^{2}(\sigma)}\xi\left(\log\xi\right)\lesssim \frac{\lambda^{2}(\tau)}{\lambda^{2}(\sigma)},
\end{align*}
which is as desired.
\end{proof}
The following proposition is analogous to Prop.~\ref{prop:neq1fourierwavesoln} and the proof in \cite{KMiao}(which is similar to the one of Prop.~\ref{prop:ngeq2fourierwavesoln}) can be directly adjusted to our setting:
\begin{proposition}\label{prop:neq0fourierwavesoln} There exists $\tau_0$ large enough, such that if $\mathcal{F}^{(0)}(\mathcal{D}_0F_{\pm}(0))$ is a good source function according to Def.~\ref{defi:xsingulartermsnless2smooth} whose principal singular part is of restricted type and with the $k = 1,2, 3$-components for $l = 0$ vanishing, with $\mathcal{F}^{(0)}$ denoting the distorted Fourier transform at angular momentum $n = 0$ in the sense of subsection~\ref{subsec:nexcspaces}, then the problem \eqref{eq:vareps0eqn} admits a solution on $[\tau_0,\infty)$ with (distorted Fourier) vanishing data 
	\[
	\left(\overline{x}^{(0)}(\tau_0,\cdot),\,\mathcal{D}_{\tau}\overline{x}^{(0)}(\tau_0,\cdot)\right) = \left(0,0\right)
	\]
	at time $\tau = \tau_0$, and such that (in terms of the distorted Fourier transform of $\varepsilon_+^0$)
	\[
	\left\| \overline{x}^{(0)}\right\|_{\text{good}}\lesssim \left\| \mathcal{F}^{(0)}(F_{\pm}(0))\right\|_{\text{goodsource}} + \left\|\big(c_0(\cdot),\,0\big)\right\|_{good}.
	\]
\end{proposition}

We have the following analogue of Prop.~\ref{prop:neq12finalsourcetermestimatesingoodspaces} 
\begin{proposition}\label{prop:neq22finalsourcetermestimatesingoodspaces} Assume that $\Lambda\ll 1$ is defined via \eqref{eq:LambdaDef}, and that $\tau_0\gg 1$ is sufficiently large. Then there exists an angular momentum $0$ function
	\[
	\psi^{\pm}_j(0) = \int_0^\infty \xb_j^{(0)}(\tau,\xi)\phi_{0}(R,\xi)\tilde{\rho}_{1}(\xi)\,d\xi,
	\]
	see subsection~\ref{subsec:nexcspaces}, with 
	\[
	\left\| \xb_j^{(0)}\right\|_{\text{good}}\lesssim \Lambda^{3}, 
	\]
	and such that if $F_j^{(\pm)}$ represents any one of the functions occurring in Prop.~\ref{prop:smoothlinearsource}, Prop.~\ref{prop:bilinwithUregular2}, Prop.~\ref{prop:bilinwithUregular3}, and we write 
	\[
	F_j^{(\pm)} = \sum_{n\in \Z}F_j^{(\pm)}(n)e^{in\theta}, 
	\]
	then we have (here $\Box_1$ denotes the wave operator on the left of \eqref{eq:vareps1eqn})
	\begin{align*}
		\left(\mathcal{D}_0F_j^{(\pm)}(0) - \Box_1\psi^{\pm}_j(0)\right)|_{R<\nu\tau} = G_j^{(\pm)}(0)|_{R<\nu\tau}, 
	\end{align*}
	where $G_j^{(\pm)}(0)$ is a good angular momentum $1$ source function and more quantitatively setting 
	\[
	G_j^{(\pm)}(0) = \int_0^\infty \phi_{n}(R,\xi)\yb_0(\tau,\xi)\tilde{\rho}_{1}(\xi)\,d\xi,\,
	\]
	we have 
	\begin{align*}
	\left\| \yb_0(\tau,\xi)\right\|_{\text{goodsource}}\lesssim \left(\tau_0^{-1}+\Lambda\right)\cdot\Lambda\ll\Lambda. 
	\end{align*}
	If we restrict to functions $\xb_n$ with restricted principal singular part, and correspondingly use $\left\| \xb_n\right\|_{\text{good}(r)}$, then $\yb_0$ also has restricted singular principal part, and we may replace 
	$\left\| \yb_0(\tau,\xi)\right\|_{\text{goodsource}}$ by $\left\| \yb_0(\tau,\xi)\right\|_{\text{goodsource}(r)}$. 
\end{proposition}
 Continuing in the vein of the case $n = 1$, we next consider the analogue of \eqref{eq:delicatesourcetermforcplusode}, which is the source term 
\begin{equation}\label{eq:delicatesourcetermforcode}
	\lim_{R\rightarrow 0}R^{-1}H_0^{+}\varepsilon_+^0 =  \lim_{R\rightarrow 0}R^{-1}\mathcal{D}_0^{*}\mathcal{D}_0\varepsilon_+^0 = c\int_0^\infty \xb(\tau, \xi)\cdot\tilde{\rho}_{0}(\xi)\,d\xi
\end{equation}
for a suitable real number $c\neq 0$. To proceed, we use a proposition analogous to Prop.~\ref{prop:fourierboundsforn=1termlocalizedaway}, which differs subtly in the low frequency regime:
\begin{proposition}\label{prop:fourierboundsforn=0termlocalizedaway} Denote by $\mathcal{F}^{(0)}$ the distorted Fourier transform at angular momentum $n = 0$ in the sense of subsection~\ref{subsec:nexcspaces}. Assume that we have \eqref{eq:varepsilonpmangulardecomp} and define $\Lambda$ as before, with $\Lambda\ll 1$. Finally, assume that the distorted Fourier transforms of the $\varepsilon_{\pm}(n)$ (both for the $|n|\geq 2$ and the exceptional modes) have restricted principal singular type (and in particular, we have to define $\Lambda$ by using $\left\|\cdot\right\|_{\text{good}(r)}$). 
	Then we can write
	\begin{equation}\label{eq:obstructionfor=0mode}
		\mathcal{F}^{(0)}\left(\mathcal{D}_0\left(F_+(0)\right)(\sigma, \cdot)\right)(\xi) = \sum_{\pm}\sum_{k=1,2,3}\sum_{j\leq N_1}\frac{e^{\pm i\nu\sigma\xi^{\frac12}}}{\xi^{1+\frac{k\nu}{2}}}\left(\log\xi\right)^j\cdot\beta_{\pm}^{(k,j)}(\sigma) + \yb(\sigma, \xi) + \tilde{\yb}(\sigma, \xi),
	\end{equation}
	where $\yb(\tau, \xi)$ is a good source term at angular momentum $n = 0$ and such that all terms with $k\in \{1,2,3\}, l = 0$ in the expansion of its singular part (according to Definition~\ref{defi:xsingulartermsnless2proto}) vanish. Furthermore, $\tilde{\yb}(\sigma,\cdot)\in \langle \xi^{\frac14}\rangle\cdot S_1^{(0)}$, and with 
	\begin{align*}
	\big\|\tilde{\yb}(\sigma,\cdot)\big\|_{\langle \xi^{\frac14}\rangle\cdot S_1^{(0)}} \lesssim \sigma^{-4+}\cdot \Lambda\cdot(\tau_0^{-1}+\Lambda),\,\sigma\geq \tau_0. 
	\end{align*}
		Finally, we have the bound
	\begin{align*}
		\left|\beta^{(k,j)}(\sigma)\right|\lesssim \sigma^{-4-}\cdot \Lambda\cdot(\tau_0^{-1}+\Lambda),\,\sigma\geq \tau_0 ,
	\end{align*}
	as well as 
	\begin{align*}
	\big\| \yb\big\|_{goodsource}\lesssim  \Lambda\cdot(\tau_0^{-1}+\Lambda).
	\end{align*}
	provided we have the bound $\left|c_0(\tau)\right| + \tau\cdot \left|c_0'(\tau)\right|\lesssim \Lambda\cdot\tau^{-2-}$, which hence improves the generic assumption for the $c_n(\tau)$ in the definition of $\Lambda$. 
	\end{proposition}

\begin{proof}
The proof is similar to that of Proposition \ref{prop:fourierboundsforn=1termlocalizedaway}. The only difference is that now for $R\xi^{\frac12}\geq 1$ the principal part of the Fourier basis is $\xi^{-\frac54}R^{-\frac12}e^{iR\xi^{\frac12}}$. We omit the analogous details.
\end{proof}
From here we can again infer the main obstruction to obtaining a good bound for the source term \eqref{eq:delicatesourcetermforcode} in the ODE \eqref{eq:c0evolutioneqn}:
\begin{lemma}\label{lem:n=0maintermafterparametrix} Denote by $U^{(0)}\left(\tau,\sigma, \xi\right)$ the Duhamel propagator according to Proposition~\ref{prop:parametrixongoodsourcenless2}.
	 Then setting 
	\[
	\xb(\tau, \xi) = \int_{\tau_0}^{\tau}U^{(0)}\left(\tau,\sigma, \xi\right)\cdot  \mathcal{F}^{(0)}\left(\mathcal{D}_0\left(F_+(0)\right)(\sigma, \cdot)\right)\left(\frac{\lambda^2(\tau)}{\lambda^2(\sigma)}\xi\right)\,d\sigma,  
	\]
	we have (for an arbitrary large but finite number $M$ and a smooth cutoff $\chi_{\lambda^2(\tau)\xi<M}$)
	\begin{equation}\label{eq:c_+sourcetermobstruction2}\begin{split}
		&\int_0^\infty \chi_{\lambda^2(\tau)\xi<M}\xb(\tau, \xi)\tilde{\rho}_{0}(\xi)\,d\xi\\& = \sum_{\pm}\sum_{k=1,2,3}\sum_{0\leq j_1\leq N_1}\int_0^\infty \chi_{\lambda^2(\tau)\xi<M}\frac{e^{\pm i\nu\tau\xi^{\frac12}}}{\xi^{\frac12+\frac{k\nu}{2}}}\left[\log\left(\xi\lambda^2(\tau)\right)\right]^{j_1}\tilde{\rho}_{0}^{\frac12}(\xi)\,d\xi\\&\hspace{6cm}\cdot \sum_{j_2\leq N_1-j_1} \int_{\tau_0}^{\tau}\left(\log\sigma\right)^{j_2}\left[\frac{\lambda(\sigma)}{\lambda(\tau)}\right]^{\nu k}\beta_{\pm}^{(k,j_1,j_2)}(\sigma)\,d\sigma\\
		& + \tilde{c}_0(\tau), 
	\end{split}\end{equation}
	where the error satisfies 
	\[
	\left|\tilde{c}_0(\tau)\right|\lesssim \tau^{-4-}\cdot  \Lambda\cdot(\tau_0^{-1}+\Lambda),
	\]
	provided we impose the same strengthening on the decay rate of $c_0(\tau)$ as in the preceding proposition. 
\end{lemma}
\begin{proof}
	The  proof is similar to that of Lemma \ref{lem:n=1maintermafterparametrix} and we omit the details.
\end{proof}
We note that the coefficients $\beta_{\pm}^{(k,j_1,j_2)}(\sigma)$ admit a description like \eqref{eq:leadsing2Duhamelspecialcoefficients}. 

Observe that the smooth cutoff $\chi_{\lambda^2(\tau)\xi<M}$ ensures that the first integral in fact has a uniform bound in $M$, and in fact we have 
\begin{align*}
	\left| \int_0^\infty \chi_{\lambda^2(\tau)\xi<M}\frac{e^{\pm i\nu\tau\xi^{\frac12}}}{\xi^{\frac12+\frac{k\nu}{2}}}\left[\log\left(\xi\lambda^2(\tau)\right)\right]^{j_1}\left(\tilde{\rho}_{0}\right)^{\frac12}(\xi)\,d\xi\right|\lesssim \tau^{-(1-k\nu)}\left(\log\tau\right)^{N_1}, 
\end{align*}
as one sees by inserting the smooth cutoffs $\chi_{\tau\xi^{\frac12}\lesssim 1}, \chi_{\tau\xi^{\frac12}\gtrsim 1}$, and performing integration by parts in case of inserting the latter cutoff. As in the case $n = 1$, we could infer the good bound 
\[
\left| \int_0^\infty \chi_{\lambda^2(\tau)\xi<M}\xb(\tau, \xi)\tilde{\rho}_{0}(\xi)\,d\xi\right|\lesssim \tau^{-4-},
\]
in case we could enforce the vanishing conditions 
\begin{equation}\label{eq:n=0vanishingrelations}
\sum_{j_2\leq N_1-j_1} \int_{\tau_0}^{\infty}\left(\log\sigma\right)^{j_2}\left[\lambda(\sigma)\right]^{\nu k}\beta_{\pm}^{(k,j_1,j_2)}(\sigma)\,d\sigma = 0
\end{equation}
in case $j_1 = 0,1,\ldots, N_1$. 
\\

The required cancellations here shall be enforced by exploiting scaling invariance, as well as the one remaining rotation on the target.

\subsubsection{The angular mode $n=-1$, equations on the physical side}
Finally, we turn to the equation for the $n= -1$ angular momentum mode, which we recall is given by 
\begin{equation}\label{eq:n=-1modequation}
	-\left(\partial_{\tau}+\frac{\lambda_{\tau}}{\lambda}R\partial_R\right)^2\varepsilon_+^{-1} - \frac{\lambda_{\tau}}{\lambda}\left(\partial_{\tau}+\frac{\lambda_{\tau}}{\lambda}R\partial_R\right)\varepsilon_+^{-1} + H_{-1}^{+}\varepsilon_+^{-1} = F_+(-1), 
\end{equation}
where we also focus on the $+$-case. Recall the representation 
\[
\varepsilon_+^{-1}(\tau, R) = c_{-1}(\tau)\phi_{-1}(R) + \phi_{-1}(R)\cdot\int_0^R\left[\phi_{-1}(s)\right]^{-1}\mathcal{D}_{-}\varepsilon_+^{-1}(\tau, s)\,ds, 
\]
where we set $\mathcal{D}_- = \partial_R  - \frac{2}{R} + \frac{2R}{1+R^2}$, and we recall $\phi_{-1}(R) = \frac{R^2}{1+R^2}$. Commuting $\mathcal{D}_-$ past the equation results in 
\begin{equation}\label{eq:vareps-1eqn}\begin{split}
		&-\left(\partial_{\tau}+\frac{\lambda_{\tau}}{\lambda}R\partial_R\right)^2\mathcal{D}_{-}\varepsilon_+^{-1} - 3\frac{\lambda_{\tau}}{\lambda}\left(\partial_{\tau}+\frac{\lambda_{\tau}}{\lambda}R\partial_R\right)\mathcal{D}_-\varepsilon_+^{-1} + \tilde{H}_{-1}^{+}\varepsilon_+^{-1} - \left(2\left(\frac{\lambda'}{\lambda}\right)^2 + \left(\frac{\lambda'}{\lambda}\right)'\right)\mathcal{D}_-\varepsilon_+^{-1}\\
		& = -\frac{\lambda_{\tau}}{\lambda}\frac{8R}{(1+R^2)^2}\left(\partial_{\tau}+\frac{\lambda_{\tau}}{\lambda}R\partial_R\right)\varepsilon_+^{-1} + \left(\frac{\lambda'}{\lambda}\right)^2\cdot\left(\frac{4R}{(1+R^2)^2} - \frac{16R}{(1+R^2)^3}\right)\varepsilon_+^{-1}\\
		&- \frac{4R}{(1+R^2)^2} \left(\frac{\lambda'}{\lambda}\right)'\varepsilon_+^{-1} + \mathcal{D}_-\left(F_+(-1)\right)\\
		&=:\mathcal{R}_+^{-1}\left(\varepsilon_+^{-1},\mathcal{D}_-\varepsilon_+^{-1}\right) + { \mathcal{D}_-\left(F_+(-1)\right)}.
\end{split}\end{equation}
Here we recall that the `super-symmetrical' operator $\tilde{H}_{-1}^{+}\varepsilon_+^1 = -\mathcal{D}_-\mathcal{D}_{-}^*\veps_{+}^{1}$, i.e., the factors have been switched compared to $H_{-1}^{+}$. 
The preceding equation gets complemented by the one for $c_-(\tau)$ which arises by analyzing the terms vanishing to lowest order at the origin in \eqref{eq:n=-1modequation}, namely those vanishing to second order: 
\begin{equation}\label{eq:c-1evolutioneqn}
	-\left(\partial_{\tau} + 2\frac{\lambda'}{\lambda}\right)^2 c_- -  \frac{\lambda'}{\lambda}\left(\partial_{\tau}+ 2\frac{\lambda'}{\lambda}\right)c_- + \lim_{R\rightarrow 0}R^{-2}H_{-1}^{+}\varepsilon_+^{-1} =  \lim_{R\rightarrow 0}R^{-2}F_+(-1)
\end{equation}
The operator here occurring on the left 
\[
\left(\partial_{\tau} + 2\frac{\lambda'}{\lambda}\right)^2 +  \frac{\lambda'}{\lambda}\left(\partial_{\tau}+ 2\frac{\lambda'}{\lambda}\right)
\]
admits the fundamental system $\{\tau^{-2-2\nu^{-1}},\,\tau^{-2-3\nu^{-1}}\}$,  and the corresponding inhomogeneous problem 
\[
\left(\partial_{\tau} + 2\frac{\lambda'}{\lambda}\right)^2c_{-1} +  \frac{\lambda'}{\lambda}\left(\partial_{\tau}+ 2\frac{\lambda'}{\lambda}\right)c_{-1} = h
\]
is solved by the explicit expression 
\begin{equation}\label{eq:cminusinhomsoln}
	c_{-1}(\tau) = \nu\left(\tau^{-2-2\nu^{-1}}\int_{\tau_0}^{\tau}\sigma^{3+2\nu^{-1}} h(\sigma)\,d\sigma - \tau^{-2-3\nu^{-1}}\int_{\tau_0}^{\tau}\sigma^{3+3\nu^{-1}}h(\sigma)\,d\sigma\right). 
\end{equation}

\subsubsection{Translation of the equation \eqref{eq:vareps-1eqn} to the Fourier side} 

The situation here is formally quite analogous to the one in the case $n = 0$, except the asymptotics of the spectral measure are quite different. 
Introduce the auxiliary operator 
\[
\mathcal{D}_{\tau}^{(-1)} = \partial_{\tau} - 2\frac{\lambda_{\tau}}{\lambda}\xi\partial_{\xi} - \frac{\lambda_{\tau}}{\lambda}\frac{\left(\tilde{\rho}_{-1}(\xi)\right)'\xi}{\tilde{\rho}_{-1}(\xi)} - \frac{\lambda_{\tau}}{\lambda}. 
\]
By means of it, we translate the equation \eqref{eq:vareps0eqn} to the Fourier variables $\xb(\tau, \xi)$, which satisfies 
\[
\calD_-\veps_+^{-1}(\tau, R) = \int_0^\infty \phi_{-1}(R,\xi)\cdot \xb(\tau, \xi)\cdot \tilde{\rho}_{-1}(\xi)\,d\xi,
\]
in perfect analogy to \eqref{eq:ngeq2Fourier1} as follows: 
\begin{equation}\label{eq:neq-1Fourier1}\begin{split}
		&-\left(\left(\mathcal{D}_{\tau}^{(-1)}\right)^2 + \frac{\lambda'(\tau)}{\lambda(\tau)}\mathcal{D}_{\tau}^{(-1)} + \xi \right)\xb\\
		& =  \mathcal{F}^{(-1)}\left(F_+(-1)\right) +  2\frac{\lambda'(\tau)}{\lambda(\tau)}\mathcal{K}_{0}^{(-1)}\mathcal{D}_{\tau}^{(-1)}\xb + \left(\frac{\lambda'(\tau)}{\lambda(\tau)}\right)'\mathcal{K}_{0}^{(-1)}\xb + \frac{\lambda'(\tau)}{\lambda(\tau)}\left[\mathcal{D}_{\tau}^{(-1)},\mathcal{K}_{0}^{(-1)}\right]\xb + \left(\frac{\lambda'(\tau)}{\lambda(\tau)}\right)^2\left(\left(\mathcal{K}_{0}^{(-1)}\right)^2 + \mathcal{K}_{0}^{(-1)}\right)\xb \\
		& +  \mathcal{F}^{(-1)}\left(\calR^{-1}\left(\veps_+^{-1},\calD_{-}\veps_+^{-1}\right)\right)
\end{split}\end{equation}
where $\mathcal{K}_0^{(-1)}$ is the off-diagonal part of the transference operator associated with $\phi_{-1}(R,\xi)$. 
\\
The evolution under the linear operator on the left is described by the following 
\begin{lemma}\label{lem:n=-1Fourierwavehomogeneous} The homogeneous initial value problem 
	\[
	\left(\left(\mathcal{D}_{\tau}^{(-1)}\right)^2 + \frac{\lambda'}{\lambda} \mathcal{D}_{\tau}^{(-1)} + \xi\right)\overline{x}(\tau,\xi) = 0,\quad\overline{x}(\tau_0,\xi) = \xb_0(\xi),\quad \mathcal{D}_{\tau}^{(-1)}\overline{x}(\tau_0,\xi) = \xb_1(\xi)
	\]
	is solved by the function 
	\begin{align*}
		\overline{x}(\tau,\xi) &= \frac{\lambda(\tau)}{\lambda(\tau_0)}\cdot\frac{\left(\tilde{\rho}_{-1}\right)^{\frac12}\left(\frac{\lambda^2(\tau)}{\lambda^2(\tau_0)}\xi\right)}{\left(\tilde{\rho}_{-1}\right)^{\frac12}(\xi)}\cdot \cos\left[\lambda(\tau)\xi^{\frac12}\int_{\tau_0}^{\tau}\lambda^{-1}(u)\,du\right]\cdot \xb_0\left(\frac{\lambda^2(\tau)}{\lambda^2(\tau_0)}\xi\right)\\
		& +\frac{\left(\tilde{\rho}_{-1}\right)^{\frac12}\left(\frac{\lambda^2(\tau)}{\lambda^2(\tau_0)}\xi\right)}{\left(\tilde{\rho}_{-1}\right)^{\frac12}(\xi)}\cdot \xi^{-\frac12}\sin\left[\lambda(\tau)\xi^{\frac12}\int_{\tau_0}^{\tau}\lambda^{-1}(u)\,du\right]\cdot \xb_1\left(\frac{\lambda^2(\tau)}{\lambda^2(\tau_0)}\xi\right)\\
	\end{align*}
	This implies the fundamental $S^{(0)}$-space propagation bounds (recall the definition \eqref{eq:S0-1} of this norm)
	\begin{align*}
		\left\| \overline{x}(\tau,\xi)\right\|_{S_0^{(-1)}} + \left\|\mathcal{D}_{\tau}^{(-1)} \overline{x}(\tau,\xi)\right\|_{S_1^{(-1)}}
		\lesssim \left(\frac{\lambda(\tau_0)}{\lambda(\tau)}\right)^{1-\delta}\cdot \left[\left\|\xb_0\right\|_{S_0^{(-1)}} + \left\|\xb_1\right\|_{S_1^{(-1)}}\right]
	\end{align*}
	In particular, choosing $\nu$ small enough, this quantity decays faster than any prescribed negative power of $\tau$. 
\end{lemma}
\begin{lemma}\label{lem:n=-1Fourierwaveinhomogeneous} 
	The inhomogeneous initial value problem 
	\[
	\left(\left(\mathcal{D}_{\tau}^{(-1)}\right)^2 + \frac{\lambda'}{\lambda} \mathcal{D}_{\tau}^{(-1)} + \xi\right)\overline{x}(\tau,\xi) = f(\tau,\xi),\quad\overline{x}(\tau_0,\xi) = 0,\quad \mathcal{D}_{\tau}^{(-1)}\overline{x}(\tau_0,\xi) = 0
	\]
	is solved by the function 
	\begin{align*}
		\xb(\tau,\xi)=\xi^{-\frac12}\int_{\tau_{0}}^{\tau}\frac{\left(\tilde{\rho}_{-1}\right)^{\frac12}\left(\frac{\lambda^2(\tau)}{\lambda^2(\sigma)}\xi\right)}{\left(\tilde{\rho}_{-1}\right)^{\frac12}(\xi)}\sin\left[\lambda(\tau)\xi^{\frac12}\int_{\sigma}^{\tau}\lambda^{-1}(u)\,du\right]\cdot f\left(\sigma,\frac{\lambda^2(\tau)}{\lambda^2(\sigma)}\xi\right)\,d\sigma
	\end{align*}
\end{lemma}
Again in light of  Definition~\ref{defi:goodfourierrepnless2} , we have the following analogue of Prop.~\ref{prop:parametrixongoodsourcenless2}:
\begin{proposition}\label{prop:parametrixongoodsourcenless21}
	Let the function $\yb(\tau,\xi)$  be a 'good source function' at angular momentum $n = -1$, in accordance with Def.~\ref{defi:xsingulartermsnless2smooth}, Def.~\ref{defi:xsingulartermsnless2proto}. Assume that for the admissibly singular source part the terms with $k = 1, 2, 3$ and $l = 0$ vanish. Denoting by $U^{(-1)}(\tau, \sigma,\xi)$ Duhamel propagator for the inhomogeneous problem associated to $ \mathcal{D}_{\tau,-1}^2 + \frac{\lambda'}{\lambda} \mathcal{D}_{\tau,-1} + \xi$, we have 
	\begin{align*}
		\left\|\int_{\tau_0}^\tau U^{(-1)}\left(\tau, \sigma,\xi\right)\cdot \yb\left(\sigma,\frac{\lambda^2(\tau)}{\lambda^2(\sigma)}\xi\right)\,d\sigma\right\|_{\text{good}}\lesssim \left\|\yb\right\|_{\text{goodsource}},
	\end{align*}
	where the norms are of course in the sense of angular momentum $n = 1$ functions. If $\yb$ has principal singular part of restricted type, so does the left hand side, and the norms can be adjusted accordingly. 
\end{proposition}
\begin{proof}
	Recall the expression of the kernel $U^{(-1)}(\tau,\sigma,\xi)$:
	\begin{align*}
	U^{(-1)}(\tau,\sigma,\xi)=\xi^{-\frac12}\frac{\left(\tilde{\rho}_{-1}\right)^{\frac12}\left(\frac{\lambda^2(\tau)}{\lambda^2(\sigma)}\xi\right)}{\left(\tilde{\rho}_{-1}\right)^{\frac12}(\xi)}\sin\left[\lambda(\tau)\xi^{\frac12}\int_{\sigma}^{\tau}\lambda^{-1}(u)\,du\right].
	\end{align*}
We assume the source $\yb(\tau,\xi)$ is given by $\yb(\tau,\xi)=\tau^{-1}\xi^{\frac12}\zb(\tau,\xi)$. We start with the smooth part by assuming $\zb(\tau,\cdot)\in S_{0}^{-}$. Then we  want to prove that $\xb(\tau,\cdot)\in S_{0}^{-}$. We also recall the definition of $S_{0}^{-}$:
\begin{align*}
	\left\|\xb(\xi)\right\|_{S_{0}^{-}}:=\left\|\xi^{1-\frac{\delta}{2}}\langle\xi\rangle^{\frac52+\delta}\xb(\xi)\right\|_{L^{2}_{d\xi}}.
\end{align*}
We start with the case \underline{$\frac{\lambda^{2}(\tau)}{\lambda^{2}(\sigma)}\xi\lesssim 1$}. Now we have
\begin{align*}
	\frac{\left(\tilde{\rho}_{-1}\right)^{\frac12}\left(\frac{\lambda^2(\tau)}{\lambda^2(\sigma)}\xi\right)}{\left(\tilde{\rho}_{-1}\right)^{\frac12}(\xi)}\simeq \frac{\lambda(\tau)}{\lambda(\sigma)},
\end{align*}
which gives
\begin{align*}
\left\|\xb(\tau,\cdot)\right\|_{S_{0}^{-}}\lesssim& \int_{\tau_{0}}^{\tau}\sigma^{-1}\frac{\lambda(\tau)}{\lambda(\sigma)}\cdot\frac{\lambda(\tau)}{\lambda(\sigma)} \left\|\xi^{1-\frac{\delta}{2}}\zb\left(\sigma,\frac{\lambda^{2}(\tau)}{\lambda^{2}(\sigma)}\xi\right)\right\|_{L^{2}_{d\xi}}\,d\sigma\\
\lesssim &\int_{\tau_{0}}^{\tau}\sigma^{-1}\left(\frac{\lambda(\sigma)}{\lambda(\tau)}\right)^{1-\delta}\|\zb(\sigma,\cdot)\|_{S_{0}^{-}}\,d\sigma\lesssim\sup_{\sigma\in[\tau_{0},\tau]}\|\zb(\sigma,\cdot)\|_{S_{0}^{-}}.
\end{align*}
Next we  turn to the case \underline{$\xi\lesssim1$ and $\frac{\lambda^{2}(\tau)}{\lambda^{2}(\sigma)}\xi\gg1$}. In this case we have
\begin{align*}
	\frac{\left(\tilde{\rho}_{-1}\right)^{\frac12}\left(\frac{\lambda^2(\tau)}{\lambda^2(\sigma)}\xi\right)}{\left(\tilde{\rho}_{-1}\right)^{\frac12}(\xi)}\simeq\frac{\lambda^{3}(\tau)}{\lambda^{3}(\sigma)}\xi,
\end{align*}
and 
\begin{align*}
	\|\xb(\tau,\cdot)\|_{S_{0}^{-}}\lesssim&\int_{\tau_{0}}^{\tau}\sigma^{-1}\frac{\lambda^{3}(\tau)}{\lambda^{3}(\sigma)}\cdot\frac{\lambda(\tau)}{\lambda(\sigma)}\left\|\xi\cdot\xi^{1-\frac{\delta}{2}}\zb\left(\sigma,\frac{\lambda^{2}(\tau)}{\lambda^{2}(\sigma)}\xi\right)\right\|_{L^{2}_{d\xi}}\,d\sigma\\
	\lesssim& \int_{\tau_{0}}^{\tau}\sigma^{-1}\left(\frac{\lambda(\sigma)}{\lambda(\tau)}\right)^{1-\delta}\left\|\zb(\sigma,\cdot)\right\|_{S_{0}^{-}}\,d\sigma\lesssim\sup_{\sigma\in[\tau_{0},\tau]}\|\zb(\sigma,\cdot)\|_{S_{0}^{-}}.
\end{align*}
Finally we consider the case \underline{$\xi\gg1$}. In this case we have
\begin{align*}
	\frac{\left(\tilde{\rho}_{-1}\right)^{\frac12}\left(\frac{\lambda^2(\tau)}{\lambda^2(\sigma)}\xi\right)}{\left(\tilde{\rho}_{-1}\right)^{\frac12}(\xi)}\simeq\frac{\lambda^{3}(\tau)}{\lambda^{3}(\sigma)},
\end{align*}
and 
\begin{align*}
	\|\xb(\tau,\cdot)\|_{S_{0}^{-}}\lesssim &\int_{\tau_{0}}^{\tau}\sigma^{-1}\frac{\lambda^{3}(\tau)}{\lambda^{3}(\sigma)}\cdot\frac{\lambda(\tau)}{\lambda(\sigma)}\left\|\xi^{\frac72+\frac{\delta}{2}}\zb\left(\sigma,\frac{\lambda^{2}(\tau)}{\lambda^{2}(\sigma)}\xi\right)\right\|_{L^{2}_{d\xi}}\,d\sigma\\
	\lesssim&\int_{\tau_{0}}^{\tau}\sigma^{-1}\left(\frac{\lambda(\sigma)}{\lambda(\tau)}\right)^{4+\delta}\left\|\zb(\sigma,\cdot)\right\|_{S_{0}^{-}}\,d\sigma\lesssim \sup_{\sigma\in[\tau_{0},\tau]}\left\|\zb(\sigma,\cdot)\right\|_{S_{0}^{-}}.
\end{align*}
The singular part can be handled similar as in the proof for Propositions \ref{prop:parametrixongoodsourcenless2} and \ref{prop:parametrixongoodsourcenless20}. We omit  the analogous details.
\end{proof}

In analogy to the cases of angular momenta $n = 0,\,n = 1$, we have 
\begin{proposition}\label{prop:neq-12finalsourcetermestimatesingoodspaces} Assume that $\Lambda\ll 1$ is defined via \eqref{eq:LambdaDef}, and that $\tau_0\gg 1$ is sufficiently large. Then there exists an angular momentum $0$ function
	\[
	\psi^{+}_j(-1) = \int_0^\infty \xb_j^{(-1)}(\tau,\xi)\phi_{-1}(R,\xi)\tilde{\rho}_{-1}(\xi)\,d\xi,
	\]
	see subsection~\ref{subsec:nexcspaces}, with 
	\[
	\left\| \xb_j^{(-1)}\right\|_{\text{good}}\lesssim \Lambda^{3}, 
	\]
	and such that if $F_j^{(\pm)}$ represents any one of the functions occurring in Prop.~\ref{prop:smoothlinearsource}, Prop.~\ref{prop:bilinwithUregular2}, Prop.~\ref{prop:bilinwithUregular3}, and we write 
	\[
	F_j^{(\pm)} = \sum_{n\in \Z}F_j^{(\pm)}(n)e^{in\theta}, 
	\]
	then we have (here $\Box_{-1}$ denotes the wave operator on the left of \eqref{eq:vareps-1eqn})
	\begin{align*}
		\left(\mathcal{D}_{-}F_j^{(+)}(-1) - \Box_{-1}\psi^{+}_j(-1)\right)|_{R<\nu\tau} = G_j^{(+)}(-1)|_{R<\nu\tau}, 
	\end{align*}
	where $G_j^{(+)}(-1)$ is a good angular momentum $-1$ source function and more quantitatively setting 
	\[
	G_j^{(+)}(-1) = \int_0^\infty \phi_{n}(R,\xi)\yb_{-1}(\tau,\xi)\tilde{\rho}_{-1}(\xi)\,d\xi,\,
	\]
	we have 
	\begin{align*}
	\left\| \yb_{-1}(\tau,\xi)\right\|_{\text{goodsource}}\lesssim \left(\tau_0^{-1}+\Lambda\right)\cdot\Lambda\ll\Lambda. 
	\end{align*}
	If we restrict to functions $\xb_n$ with restricted principal singular part, and correspondingly use $\left\| \xb_n\right\|_{\text{good}(r)}$, then $\yb_{-1}$ also has restricted singular principal part, and we may replace 
	$\left\| \yb_{-1}(\tau,\xi)\right\|_{\text{goodsource}}$ by $\left\| \yb_{-1}(\tau,\xi)\right\|_{\text{goodsource}(r)}$. 
\end{proposition}

Finally, we analyze the delicate term $\lim_{R\rightarrow 0}R^{-2}H_{-1}^{+}\varepsilon_+^{-1}$ in \eqref{eq:c-1evolutioneqn}, for which we have the customary relation 
\[
\lim_{R\rightarrow 0}R^{-2}H_{-1}^{+}\varepsilon_+^{-1} = d\int_0^\infty  \xb(\tau, \xi)\cdot \tilde{\rho}_{-1}(\xi)\,d\xi
\]
for suitable $d\neq 0$, provided $\xb(\tau, \xi)$ represents the distorted Fourier transform of $\mathcal{D}_{-}\varepsilon_+^{-1}$. In order to control this, we again need to understand the most singular terms in the source, which result in poorly temporally decaying contributions to the preceding expression: 
\begin{proposition}\label{prop:fourierboundsforn=-1termlocalizedaway} Denote by $\mathcal{F}^{(-1)}$ the distorted Fourier transform at angular momentum $n = -1$ in the sense of subsection~\ref{subsec:nexcspaces}. Assume that we have \eqref{eq:varepsilonpmangulardecomp} and define $\Lambda$ as before, with $\Lambda\ll 1$. Finally, assume that the distorted Fourier transforms of the $\varepsilon_{\pm}(n)$ (both for the $|n|\geq 2$ and the exceptional modes) have restricted principal singular type (and in particular, we have to define $\Lambda$ by using $\left\|\cdot\right\|_{\text{good}(r)}$). 
	Then we can write
	\begin{equation}\label{eq:obstructionfor=-1mode}
		\mathcal{F}^{(0)}\left(\mathcal{D}_-\left(F_+(-1)\right)(\sigma, \cdot)\right)(\xi) = \sum_{\pm}\sum_{k=1,2,3}\sum_{j\leq N_1}\frac{e^{\pm i\nu\sigma\xi^{\frac12}}}{\xi^{\frac32+\frac{k\nu}{2}}}\left(\log\xi\right)^j\cdot\beta_{\pm}^{(k,j)}(\sigma) + \yb(\sigma, \xi) + \tilde{\yb}(\sigma, \xi),
	\end{equation}
	where $\yb(\tau, \xi)$ is a good source term at angular momentum $n = -1$ and such that all terms with $k\in \{1,2,3\}, l = 0$ in the expansion of its singular part (according to Definition~\ref{defi:xsingulartermsnless2proto}) vanish. Furthermore, $\tilde{\yb}(\sigma,\cdot)\in \langle \xi^{\frac14}\rangle\cdot S_1^{(-1)}$, and with 
	\begin{align*}
	\big\|\tilde{\yb}(\sigma,\cdot)\big\|_{\langle \xi^{\frac14}\rangle\cdot S_1^{(-1)}} \lesssim \sigma^{-4+}\cdot \Lambda\cdot(\tau_0^{-1}+\Lambda),\,\sigma\geq \tau_0. 
	\end{align*}
		Finally, we have the bound
	\begin{align*}
		\left|\beta^{(k,j)}(\sigma)\right|\lesssim \sigma^{-4-}\cdot \Lambda\cdot(\tau_0^{-1}+\Lambda),\,\sigma\geq \tau_0 ,
	\end{align*}
	as well as 
	\begin{align*}
	\big\| \yb\big\|_{goodsource}\lesssim  \Lambda\cdot(\tau_0^{-1}+\Lambda).
	\end{align*}
	provided we have the bound $\left|c_0(\tau)\right| + \tau\cdot \left|c_0'(\tau)\right|\lesssim \Lambda\cdot\tau^{-2-}$, which hence improves the generic assumption for the $c_n(\tau)$ in the definition of $\Lambda$. 
	\end{proposition}

\begin{proof}
	The proof is again similar to that of Proposition \ref{prop:fourierboundsforn=1termlocalizedaway} and the only difference is that now for $R\xi^{\frac12}\geq 1$ the principal part of  the Fourier basis is $\xi^{-\frac74}R^{-\frac12}e^{iR\xi^{\frac12}}$. We omit the analogous details.
\end{proof}
As in the preceding cases, we then infer the following structure result for the truncated integral expressing $\lim_{R\rightarrow 0}R^{-2}H_{-1}^{+}\varepsilon_+^{-1}$: 
\begin{lemma}\label{lem:n=-1maintermafterparametrix} Denote by $U^{(-1)}\left(\tau,\sigma, \xi\right)$ the Duhamel propagator according to Proposition~\ref{prop:parametrixongoodsourcenless2}.
	 Then setting 
	\[
	\xb(\tau, \xi) = \int_{\tau_0}^{\tau}U^{(-1)}\left(\tau,\sigma, \xi\right)\cdot  \mathcal{F}^{(-1)}\left(\mathcal{D}_-\left(F_+(-1)\right)(\sigma, \cdot)\right)\left(\frac{\lambda^2(\tau)}{\lambda^2(\sigma)}\xi\right)\,d\sigma,  
	\]
	we have (for an arbitrary large but finite number $M$ and a smooth cutoff $\chi_{\lambda^2(\tau)\xi<M}$)
	\begin{equation}\label{eq:c_-sourcetermobstruction}\begin{split}
		&\int_0^\infty \chi_{\lambda^2(\tau)\xi<M}\xb(\tau, \xi)\tilde{\rho}_{-1}(\xi)\,d\xi\\& = \sum_{\pm}\sum_{k=1,2,3}\sum_{0\leq j_1\leq N_1}\int_0^\infty \chi_{\lambda^2(\tau)\xi<M}\frac{e^{\pm i\nu\tau\xi^{\frac12}}}{\xi^{\frac12+\frac{k\nu}{2}}}\left[\log\left(\xi\lambda^2(\tau)\right)\right]^{j_1}\,d\xi\\&\hspace{6cm}\cdot \sum_{j_2\leq N_1-j_1} \int_{\tau_0}^{\tau}\left(\log\sigma\right)^{j_2}\left[\frac{\lambda(\sigma)}{\lambda(\tau)}\right]^{\nu k}\beta_{\pm}^{(k,j_1,j_2)}(\sigma)\,d\sigma\\
		& + \tilde{c}_{-1}(\tau), 
	\end{split}\end{equation}
	where the error satisfies 
	\[
	\left|\tilde{c}_{-1}(\tau)\right|\lesssim \tau^{-4-}\cdot  \Lambda\cdot(\tau_0^{-1}+\Lambda),
	\]
	provided we impose the same strengthening on the decay rate of $c_0(\tau)$ as in the preceding proposition. 
\end{lemma}

As before the preceding implies the bound
\[
\left| \int_0^\infty \chi_{\lambda^2(\tau)\xi<M}\xb(\tau, \xi)\tilde{\rho}_{-1}(\xi)\,d\xi\right|\lesssim \tau^{-4-},
\]
in case we could enforce the vanishing conditions 
\begin{equation}\label{eq:n=-1vanishingrelations}
\sum_{j_2\leq N_1-j_1} \int_{\tau_0}^{\infty}\left(\log\sigma\right)^{j_2}\left[\lambda(\sigma)\right]^{\nu k}\beta_{\pm}^{(k,j_1,j_2)}(\sigma)\,d\sigma = 0
\end{equation}
in case $j_1 = 0,1,\ldots, N_1$. 
\begin{proof}
	The proof is again similar to that of Lemma \ref{lem:n=1maintermafterparametrix} and we omit the details.
\end{proof}
\subsection{The modulation step for the exceptional modes; forcing the vanishing conditions}\label{subsec:modulation}
In this subsection, we introduce the final tool which will allow us to close all the estimates introduced gradually in the preceding. Specifically, recalling the decomposition of our Wave Map
\begin{align*}
	\Psi &= \Phi + \Pi_{\Phi^{\perp}}\varphi + a\left(\Pi_{\Phi^{\perp}}\varphi\right)\Phi\\
	& = \Phi + \vphi_1 E_1 + \vphi_2 E_2 + a\left(\Pi_{\Phi^{\perp}}\varphi\right)\Phi,\\
\end{align*}
we modify this by applying  a re-scaling $\mathcal{S}_{c(t)}$, and a rotation $\mathcal{R}_{h(t)}^{\alpha(t),\beta(t)}$ as well as a Lorentz transform $\mathcal{L}_{v(t)}$ to the full expression. For technical reasons, the Lorentz transform as well as the scaling will be chosen to be constant from some time $t_1<t_0$, where $t_0$ is the initial time for the perturbation. In fact, there will be a lot of flexibility in choosing these modulation parameters, since their role will be to enforce certain moment conditions, which are in effect time integrals of these functions multiplied against certain weights\footnote{However, their terminal values will be uniquely determined, of course}. In total, we shall then pass to the following ansatz
\begin{equation}\label{eq:modulatedphi}
	\Psi = \mathcal{L}_{v(t)}\mathcal{R}_{h(t)}^{\alpha(t),\beta(t)}\mathcal{S}_{c(t)}\left( \Phi + \vphi_1 E_1 + \vphi_2 E_2 + a\left(\Pi_{\Phi^{\perp}}\varphi\right)\Phi\right),\quad \Phi =  \left(\begin{array}{c}\sin U\cos\theta\\ \sin  U \sin\theta\\ \cos U\end{array}\right),
\end{equation}
whence a complete description of the Wave Map evolution shall consist of the tuple of functions 
\[
\{\vphi_1,\,\vphi_2,\,\alpha(t),\beta(t),\,h(t),\,c(t),\,v_1(t),\,v_2(t)\},
\]
where in turn $\vphi_{1,2}$ are described in terms of $\varepsilon_{\pm} = \vphi_1 \mp i\vphi_2$, each in turn decomposed into angular momentum $n$ pieces, which are described and measured as in the preceding. 
While the equation for $\vphi_{1,2}$ is exactly identical in the case of constant coefficients $\alpha$ etc, the time dependence of these will introduce additional source terms. In the sequel, we shall analyze the different modulations parameters and the effect they have to leading order on the equations. Since the leading order effect will be linear, we can treat the contribution of each modulation parameter separately, which will clarify the analysis. We emphasize that our version of modulation theory differs from the more standard kind, where usually the action of the symmetries on the bulk part $Q$ are used to enforce various vanishing conditions on fixed time slices. In our version, the action of the symmetries on the singular part of the profile $U$ (and thus, a lower order term) are used to counteract singular terms arising via interactions between the perturbation and $U$, and the vanishing conditions refer to suitable time integrals. 

\subsubsection{The contribution of the rotations $\alpha(t), \beta(t)$.}\label{subsubsec:alphabetamodulation}

Our convention shall be that the angle $\alpha(t)$ corresponds to the rotation 
\begin{align*}
	\left(\begin{array}{ccc}\cos\alpha(t)&0&\sin\alpha(t)\\0&1&0\\-\sin\alpha(t)&0&\cos\alpha(t)\end{array}\right),
\end{align*}
while the angle $\beta(t)$ corresponds to 
\begin{align*}
	\left(\begin{array}{ccc}1&0&0\\0&\cos\beta(t)&\sin\beta(t)\\0&-\sin\beta(t)&\cos\beta(t)\end{array}\right).
\end{align*}

It turns out that modulating on $\alpha$ contributes both in terms of the effect on the bulk profile $Q$ as well as on the singular terms in $U$, where we recall that the polar angle $U$ of the unperturbed blow up solution is given by 
\[
U = Q(R) + \epsilon,\quad Q(R) = 2\arctan R.
\]
In order to simplify the computations at first, we start by considering the effect of $\mathcal{R}_{0}^{\alpha(t),0}$. Denoting by $(\partial_t^2)'$ only those terms where at least one derivative falls on $\alpha(t)$, we compute 
\begin{align*}
	\left(\partial_t^2\right)'\left(\mathcal{R}_{0}^{\alpha(t),0}\Phi\right)&= \mathcal{R}^{\alpha(t),0}_{0}\alpha''(t)\left(\begin{array}{ccc}\cos\alpha(t)&0&-\sin\alpha(t)\\0&1&0\\\sin\alpha(t)&0&\cos\alpha(t)\end{array}\right)\left(\begin{array}{c}-\cos\theta\sin U\sin\alpha + \cos U\cos\alpha\\0\\-\cos\theta\sin U\cos\alpha - \cos U\sin\alpha\end{array}\right)\\
	& + 2\alpha'(t)\lambda'(t)rQ'(R) \mathcal{R}^{\alpha(t),0}_{0}\left(\begin{array}{ccc}\cos\alpha(t)&0&-\sin\alpha(t)\\0&1&0\\\sin\alpha(t)&0&\cos\alpha(t)\end{array}\right)\left(\begin{array}{c}-\cos\theta\cos U\sin\alpha - \sin U\cos\alpha\\0\\-\cos\theta\cos U\cos\alpha + \sin U\sin\alpha\end{array}\right)\\
	&+2\alpha'(t)\lambda\cdot\left[\epsilon_{\tau} + \epsilon_X\right]\cdot  \mathcal{R}^{\alpha(t),0}_{0}\left(\begin{array}{ccc}\cos\alpha(t)&0&-\sin\alpha(t)\\0&1&0\\\sin\alpha(t)&0&\cos\alpha(t)\end{array}\right)\left(\begin{array}{c}-\cos\theta\cos U\sin\alpha - \sin U\cos\alpha\\0\\-\cos\theta\cos U\cos\alpha + \sin U\sin\alpha\end{array}\right)
\end{align*}
Here we have introduce the new variable $X: = \nu\tau - R = \lambda\cdot(t-r)$, whence the singularity of $\epsilon$ across the light cone gets expressed in terms of $X$. 
We simplify the terms as follows: 
\begin{align*}
	&\left(\begin{array}{ccc}\cos\alpha(t)&0&-\sin\alpha(t)\\0&1&0\\\sin\alpha(t)&0&\cos\alpha(t)\end{array}\right)\left(\begin{array}{c}-\cos\theta\sin U\sin\alpha + \cos U\cos\alpha\\0\\-\cos\theta\sin U\cos\alpha - \cos U\sin\alpha\end{array}\right)
	= \left(\begin{array}{c}\cos U\\ 0\\ -\cos\theta\sin U\end{array}\right)\\
	&\left(\begin{array}{ccc}\cos\alpha(t)&0&-\sin\alpha(t)\\0&1&0\\\sin\alpha(t)&0&\cos\alpha(t)\end{array}\right)\left(\begin{array}{c}-\cos\theta\cos U\sin\alpha - \sin U\cos\alpha\\0\\-\cos\theta\cos U\cos\alpha + \sin U\sin\alpha\end{array}\right)
	=  \left(\begin{array}{c}-\sin U\\ 0\\ -\cos\theta\cos U\end{array}\right)
\end{align*}
Also recall the customary change of temporal variables 
\[
\tau = \int_t^\infty \lambda(s)\,ds\longrightarrow \frac{d\tau}{dt} = -\lambda, 
\]
whence 
\begin{align*}
	&\alpha'(t) = -\lambda\alpha'(\tau),\quad \alpha''(t) = \lambda^2\alpha''(\tau) + \lambda\lambda_{\tau}\alpha'(\tau),\\
	& 2\alpha'(t)\lambda'(t)rQ'(R) =  2\lambda^2\alpha'(\tau)\lambda_{\tau}rQ'(R) = 2\lambda^2\alpha'(\tau)\frac{\lambda_{\tau}}{\lambda}RQ'(R) 
\end{align*}
In order to translate things to the $(\varepsilon_{1,2})$-coordinates, we expand things in terms of the $(E_{1,2},\Phi)$-frame, where we recall the formulae 
\begin{align*}
	E_1 = \left(\begin{array}{c}\cos\theta\cos U\\\sin\theta\cos U\\-\sin U\end{array}\right),\quad E_2 = \left(\begin{array}{c}-\sin\theta\\ \cos\theta\\ 0\end{array}\right).
\end{align*}
\begin{align*}
	&\left(\begin{array}{c}\cos U\\ 0\\ -\cos\theta\sin U\end{array}\right) = \cos\theta\cdot E_1 + (-\sin\theta\cos U)\cdot E_2 + (*)\cdot \Phi,\\
	&  \left(\begin{array}{c}-\sin U\\ 0\\ -\cos\theta\cos U\end{array}\right) = 0\cdot E_1 +\sin\theta\sin U\cdot  E_2 + (*)\cdot \Phi.\\
\end{align*}
This implies that the $n = 1$ component in the source term of the equation for $\varepsilon_+^1 = \vphi_1^{(1)} - i\vphi_2^{(1)}$ due to modulation in $\alpha$ is given by (after division by $\lambda^2$)
\begin{align}\label{eq:alphamodulationn=1sourcecontribution}
	&\lambda^{-2}\alpha''(t)\cdot \left(\frac12 - i\cdot \frac{i}{2}\cos U\right) + 2\alpha'(\tau)\frac{\lambda_{\tau}}{\lambda}RQ'(R)\cdot\left(0 +i\cdot \frac{i}{2}\sin U\right)\\
	& +  \alpha'(\tau)\cdot \left[\epsilon_{\tau} + \epsilon_X\right]\cdot \sin U\nonumber\\
	& = \frac{\alpha''(\tau)+\frac{\lambda_\tau}{\lambda}\alpha'(\tau)}{2}\cdot(1+\cos U) - \alpha'(\tau)\frac{\lambda_{\tau}}{\lambda}RQ'(R)\sin U\label{eq:alphamodulationn=1A}\\
	&+ \alpha'(\tau)\cdot [\epsilon_{\tau} + \epsilon_X]\cdot \sin U.\label{eq:alphamodulationn=1B}
\end{align}
Recalling \eqref{eq:c+evolutioneqn}, we now analyze the principal contributions of the preceding expression to the term $\lim_{R\rightarrow 0}F_+(1)$, as well as the term $\lim_{R\rightarrow 0}H_1^+\varepsilon_+^1$, where the latter term requires evaluation of the wave parametrix (for the $n = +1$ mode) on $\mathcal{D_+}\left(\eqref{eq:alphamodulationn=1sourcecontribution}\right)$. We shall distinguish between the contribution of the first line \eqref{eq:alphamodulationn=1A} in the last expression in \eqref{eq:alphamodulationn=1sourcecontribution}, i.e., the smooth part, and the last term \eqref{eq:alphamodulationn=1B} (involving $\epsilon_{\tau} + \epsilon_X$), which is the non-smooth part. 
\\

{\it{(i): Contribution of \eqref{eq:alphamodulationn=1A}  to $\lim_{R\rightarrow 0}F_+(1)$.}} This follows directly from the fact that $U(\tau, 0) = 0$, and so this contribution equals
\[
\alpha''(\tau)+\frac{\lambda_\tau}{\lambda}\alpha'(\tau).
\]
This expression is that one source term for the ODE governing the evolution of $c_+(\tau)$, according to  \eqref{eq:c+evolutioneqn}.
\\

{\it{(ii): Contribution of \eqref{eq:alphamodulationn=1A}  to $\lim_{R\rightarrow 0}H_1^+\varepsilon_+^1$.}} Here we apply $\mathcal{D}_+$ to the extra source term and then the wave parametrix (after translating things to the Fourier side), and then we extract the leading behavior. Since the contribution can be handled the same way as the lower order terms on the RHS of \eqref{eq:obstructionfor=1mode}, we replace $U$ by $Q(R)$, the bulk profile. 
Note that 
\[
\cos Q(R) = \cos^2(\arctan R) -\sin^2(\arctan R) = \frac{1-R^2}{1+R^2},
\]
and so $\frac{1+\cos U}{2} = \frac{1}{1+R^2}$, which is of course killed by $\mathcal{D}_+$. Next, we have 
\begin{align*}
	RQ'(R)\sin Q(R) = \frac{2R}{1+R^2}\cdot \frac{2R}{1+R^2} = \frac{4R^2}{(1+R^2)^2} = \frac{4}{1+R^2} - \frac{4}{(1+R^2)^2}. 
\end{align*}
We conclude that to leading order we have 
\begin{align*}
	\mathcal{D}_+\left( - \alpha'(\tau)\frac{\lambda_{\tau}}{\lambda}RQ'(R)\sin U\right) = + \alpha'(\tau)\frac{\lambda_{\tau}}{\lambda}\mathcal{D}_+\left(\frac{4}{(1+R^2)^2}\right),
\end{align*}
which is a new source term in the wave equation for $\mathcal{D}_+\varepsilon_+^1$ due to modulating in $\alpha$. 
\\
Apply now the wave parametrix according to Lemma~\ref{lem:n=1Fourierwavehomogeneous}, which results in the integral  
\begin{align*}
	\int_{\tau_0}^{\tau} \xi^{-\frac12}\sin\left[\lambda(\tau)\xi^{\frac12}\int_{\sigma}^{\tau}\lambda^{-1}(u)\,du\right]\cdot \mathcal{F}^{(1)}\left(\alpha'(\sigma)\frac{\lambda_{\sigma}}{\lambda}\mathcal{D}_+\left(\frac{4}{(1+R^2)^2}\right)\right)\left(\frac{\lambda^2(\tau)}{\lambda^2(\sigma)}\xi\right)\,d\sigma
\end{align*}
To extract the main contribution from this term, we perform integration by parts with respect to $\sigma$, which produces a boundary term at $\sigma = \tau$, as well as a negligible boundary term at $\sigma = \tau_0$ and a better term with an extra $\sigma$-derivative falling on the Fourier transform. Relegating the treatment of all these errors to later, we then wind up with the principal term
\[
-\frac{1}{\xi}\cdot \mathcal{F}^{(1)}\left(\alpha'(\tau)\frac{\lambda_{\tau}(\tau)}{\lambda(\tau)}\mathcal{D}_+\left(\frac{4}{(1+R^2)^2}\right)\right)(\xi),
\]
which contributes the term 
\[
-\int_0^\infty\frac{1}{\xi}\cdot \mathcal{F}^{(1)}\left(\alpha'(\tau)\frac{\lambda_{\tau}(\tau)}{\lambda(\tau)}\mathcal{D}_+\left(\frac{4}{(1+R^2)^2}\right)\right)(\xi)\tilde{\rho}_{1}(\xi)\,d\xi
\]
to the term $\lim_{R\rightarrow 0}H_1^+\varepsilon_+^1$. We claim that this can be rendered more explicit. To begin with, we can replace $\frac{4}{(1+R^2)^2}$ by $\frac{4}{(1+R^2)^2} - \frac{c}{1+R^2}$ for arbitrary $c$.
Then we pick $c$ in such a way that 
\[
\frac{4}{(1+R^2)^2} - \frac{c}{1+R^2} = \mathcal{D}_+^*g
\]
where $g$ is such that $\int_0^\infty \left(\mathcal{F}^{(1)}g\right)(\xi)\tilde{\rho}_{1}(\xi)\,d\xi$ converges. Specifically, recall that 
\[
\mathcal{D}_+^* = -\partial_R - \frac{1}{R} + \frac{2R}{1+R^2},
\]
which annihilates $\psi(R) = \frac{1+R^2}{R}$, and so we can set 
\begin{align*}
	g(R) = -\frac{R^2+1}{R}\cdot\int_0^R \frac{s}{s^2+1}\cdot \left(\frac{4}{(1+s^2)^2} -\frac{c}{1+s^2}\right) \,ds &= \frac{R^2+1}{R}\cdot \left(\frac{1}{(1+R^2)^2} - \frac{c}{2(1+R^2)}\right)\big|_0^R\\
	& = \frac{R^2+1}{R}\cdot \frac{1-(1+R^2)}{(1+R^2)^2}\\
	& = -\frac{R}{1+R^2},
\end{align*}
provided we set $c = 2$. It follows that 
\begin{align*}
	&-\int_0^\infty\frac{1}{\xi}\cdot \mathcal{F}^{(1)}\left(\alpha'(\tau)\frac{\lambda_{\tau}(\tau)}{\lambda(\tau)}\mathcal{D}_+\left(\frac{4}{(1+R^2)^2}\right)\right)(\xi)\tilde{\rho}_{1}(\xi)\,d\xi\\
	& = -\alpha'(\tau)\frac{\lambda_{\tau}(\tau)}{\lambda(\tau)}\cdot \int_0^\infty \frac{1}{\xi}\cdot \mathcal{F}^{(1)}\left(\mathcal{D}_+\mathcal{D}_+^*g\right)(\xi)\tilde{\rho}_{1}(\xi)\,d\xi\\
	& = -\alpha'(\tau)\frac{\lambda_{\tau}(\tau)}{\lambda(\tau)}\cdot \int_0^\infty \mathcal{F}^{(1)}\left(g\right)(\xi)\tilde{\rho}_{1}(\xi)\,d\xi\\
	& = c_{0}\alpha'(\tau)\frac{\lambda_{\tau}(\tau)}{\lambda(\tau)}
\end{align*}
for certain constant $c_{0}$. This is the leading order contribution to $\lim_{R\rightarrow 0}H_1^+\varepsilon_+^1$ arising from modulating in $\alpha(t)$. 
\\

{\it{(iii): Combined leading order contribution of \eqref{eq:alphamodulationn=1A}  to \eqref{eq:c+evolutioneqn}}}. From {\it{(i)}} and {\it{(ii)}}, this is seen to be the equation 
\[
c_+''(\tau) + \frac{\lambda_{\tau}}{\lambda}c_+'(\tau) = -\alpha''(\tau)+(c_{0}-1)\frac{\lambda'(\tau)}{\lambda(\tau)}\alpha'(\tau).
\]
which in light of \eqref{eq:cplusinhomsoln}
\begin{align*}
	c_{+}(\tau)=&-\nu\int_{\tau_{0}}^{\tau}\sigma\left(\alpha''(\sigma)+(1-c_{0})(1+\nu^{-1})\sigma^{-1}\alpha'(\sigma)\right)\,d\sigma\\
	&+\nu\tau^{-\nu^{-1}}\int_{\tau_{0}}^{\tau}\sigma^{1+\nu^{-1}}\left(\alpha''(\sigma)+(1-c_{0})(1+\nu^{-1})\sigma^{-1}\alpha'(\sigma)\right)\,d\sigma.
\end{align*}
The parameter $\alpha(\tau)$ will be chosen such that $\alpha'(\tau)$ is compactly supported in $(\tau_{0},\infty)$. So for sufficiently large $\tau$ we have $\alpha(\tau)=\alpha(\infty)$. Therefore if we denote by $\alpha_{\infty}:=\alpha(\infty)-\alpha(\tau_{0})$, we have for sufficiently large $\tau$,
\begin{align*}
	c_{+}(\tau)=&\left(c_{0}(1+\nu)-1\right)\alpha_{\infty}+\nu\tau^{-\nu^{-1}}\int_{\tau_{0}}^{\infty}\sigma^{1+\nu^{-1}}\alpha''(\sigma)\,d\sigma+(c_{0}-1)(1+\nu)\tau^{-\nu^{-1}}\int_{\tau_{0}}^{\infty}\sigma^{\nu^{-1}}\alpha'(\sigma)\,d\sigma\\
	&+(1+\nu)\tau^{-\nu^{-1}}\int_{\tau}^{\infty}\sigma^{\nu^{-1}}\alpha'(\sigma)\,d\sigma+(1-c_{0})(1+\nu)\tau^{-\nu^{-1}}\int_{\tau}^{\infty}\sigma^{\nu^{-1}}\alpha'(\sigma)\,d\sigma.
\end{align*}
Except the term involving $\alpha_{\infty}$, all the other terms on the RHS above decays rapidly in $\tau$ as $\tau\rightarrow\infty$. On the other hand, since the operator $\frac{d^{2}}{d\tau^{2}}+\frac{\lambda_{\tau}}{\lambda}\frac{d}{d\tau}$ admits a constant fundamental solution. Therefore we choose $\alpha_{\infty}$ appropriately (using the flexibility of the moment condition satisfied by $\alpha'$, which is weaker than the ones for other modulation parameters except $\beta$) to cancel the constant fundamental solution such that $c_{+}(\tau)$ decays rapidly as $\tau\rightarrow\infty$.
\\

{\it{(iv): Contribution of \eqref{eq:alphamodulationn=1B}  to $\lim_{R\rightarrow 0}F_+(1)$.}} This is negligible and in fact of order $\alpha'(\tau)\tau^{-N}$ due to the structure of $\epsilon$. 
\\

{\it{(v): Contribution of \eqref{eq:alphamodulationn=1B}  to $\lim_{R\rightarrow 0}H_1^+\varepsilon_+^1$.}} Precisely, we shall want to use this term to cancel the troublesome terms on the right in \eqref{eq:obstructionfor=1mode}. 
For this, we shall use 
\begin{lemma}\label{lem:alphamodulationn=1singularterms} There is a function $H(\tau, R)$ coinciding with the function  $\calD_{+}\left(\chi_{R\gtrsim \tau}\cdot \eqref{eq:alphamodulationn=1B}\right)$ near the light cone $R=\nu\tau$, and such that we have the formula 
\begin{equation}\label{eq:modulationFourierstructureforn=1mode}\begin{split}
\mathcal{F}^{(1)}\left(\left(H\right)(\sigma, \cdot)\right)(\xi) &= \sum_{\pm}\sum_{k=1,2,3}\sum_{j_1+j_2\leq r,\,r\leq 2}\frac{e^{\pm i\nu\sigma\xi^{\frac12}}}{\xi^{\frac12+\frac{k\nu}{2}}}\left(\log\xi\right)^{j_1}\cdot\gamma_{k,r-j_1 - j_2}^{\pm}\cdot\frac{\zeta_{k,r,j_{1,2}}\cdot\log^{j_2}\sigma\cdot\alpha'(\sigma)}{\sigma^{1+k\nu}}\cdot \sin U(\sigma, \nu\sigma)\\&  + \yb(\sigma, \xi),
	\end{split}\end{equation}
	with $\gamma_{k,r-j_1-j_2}^{\pm}$ defined as in \eqref{eq:leadsing2Fourier}, while $\zeta_{k,r,j_{1,2}}$ is a suitable non-vanishing real constant, and
	where $\yb(\tau, \xi)$ is a good source term at angular momentum $n = 1$ and such that all terms with $k\in \{1,2,3\}, l = 0$ in the expansion of its singular part (according to Definition~\ref{defi:xsingulartermsnless2proto}) vanish. Furthermore,  we have the bound 
	\begin{align*}
	\big\|\yb\big\|_{goodsource} \lesssim \tau_0^{-2+}\cdot \big\|\alpha'\big\|_{L^\infty}. 
	\end{align*}
\end{lemma}
\begin{remark}\label{rem:lem:alphamodulationn=1singularterms} We shall later choose $\alpha$ with a bound $\big\|\tau^{2+}\cdot\alpha'\big\|_{L^\infty}\lesssim (\tau_0^{-1}+\Lambda)\cdot \Lambda$, whence the term $\yb$ will be perturbative. 

\end{remark}
\begin{proof} This is analogous to \eqref{eq:leadsing2}, \eqref{eq:leadsing2Fourier}, except that we now set $g_2(\tau) = \tau^{-1}\cdot\alpha'(\tau)$ up to a constant, and we also allow lower powers of the logarithm $\log(\nu\tau - R)$. The fact that we only need to sum to $r = 2$ follows from Remark~\ref{rem:fourierboundsforn=1termlocalizedaway1}.

\end{proof}
We can now deduce the following direct analogue of Lemma~\ref{lem:n=1maintermafterparametrix}:
\begin{lemma}\label{lem:alphamodulationn=1singulartermspropagated} 
 Denote by $U^{(1)}\left(\tau,\sigma, \xi\right)$ the Duhamel propagator according to Proposition~\ref{prop:parametrixongoodsourcenless2}.
	 Then setting 
	\[
	\xb(\tau, \xi) = \int_{\tau_0}^{\tau}U^{(1)}\left(\tau,\sigma, \xi\right)\cdot  \mathcal{F}^{(1)}\left(H(\sigma, \cdot)\right)\left(\frac{\lambda^2(\tau)}{\lambda^2(\sigma)}\xi\right)\,d\sigma,  
	\]
	we have (for an arbitrary large but finite number $M$ and a smooth cutoff $\chi_{\lambda^2(\tau)\xi<M}$)
	\begin{equation}\label{eq:c_+sourcetermobstruction3}\begin{split}
		&\int_0^\infty \chi_{\lambda^2(\tau)\xi<M}\xb(\tau, \xi)\tilde{\rho}_{1}(\xi)\,d\xi\\& = \sum_{\pm}\sum_{k=1,2,3}\sum_{0\leq r\leq 2}\sum_{0\leq j_1\leq r}\int_0^\infty \chi_{\lambda^2(\tau)\xi<M}\frac{e^{\pm i\nu\tau\xi^{\frac12}}}{\xi^{\frac{k\nu}{2}}}\left[\log\left(\xi\lambda^2(\tau)\right)\right]^{j_1}\,d\xi\\&\hspace{6cm}\cdot \sum_{0\leq j_2\leq r-j_1} \int_{\tau_0}^{\tau}\left(\log\sigma\right)^{j_2}\left[\frac{\lambda(\sigma)}{\lambda(\tau)}\right]^{\nu k}\tilde{\beta}_{\pm}^{(k,r,j_1,j_2)}(\sigma)\,d\sigma\\
		& + \tilde{\tilde{c}}_1(\tau), 
	\end{split}\end{equation}
	where the error satisfies under the assumption $\left\|\tau^{2+}\cdot\alpha'\right\|_{L^\infty}\lesssim (\tau_0^{-1}+\Lambda)\cdot \Lambda$ the bound
	\[
	\left|\tilde{\tilde{c}}_1(\tau)\right|\lesssim \tau^{-4-}\cdot  \Lambda\cdot(\tau_0^{-1}+\Lambda),
	\]
	and the functions $\tilde{\beta}_{\pm}^{(k,j_1,j_2)}(\sigma)$ are explicitly given by 
	\begin{align*}
	\tilde{\beta}_{\pm}^{(k,r,j_1,j_2)}(\sigma) = \sigma^{-(k-1)\nu}\cdot\lambda^{-\nu}(\sigma)\cdot \gamma_{k,r-j_1-j_2}^{\pm}\cdot\zeta_{k,r,j_{1,2}}\cdot \alpha'(\sigma)\cdot \sin U(\sigma, \nu\sigma). 
	\end{align*}
\end{lemma}
\begin{proof} This is a consequence of the preceding lemma, relying on a straightforward analogue of Lemma~\ref{lem:alphamodulationn=1singularterms}. 
\end{proof}

The preceding lemma is of course still not enough to force the vanishing relations \eqref{eq:n=1vanishingrelations}, since $\alpha'(\sigma)$ is real valued, but the other coefficients are mostly complex valued. This has to do with the fact that we have neglected modulations in the angle $\beta(t)$ up to now, which we do below. \\
However, let us briefly verify that under the assumption that the we replace \eqref{eq:leadsing2} by
\begin{equation}\label{eq:leadsingleq2}
\sum_{0\leq r\leq 2}\partial_{R}\left(\left(\nu\sigma-R\right)^{-\frac12+k\nu}\cdot \big[\log(\nu\tau - R)\big]^r\right)\cdot \frac{g_r(\tau)}{\tau^{\frac12}}, 
\end{equation}
with {\it{real valued $g_r$}}, appropriate choice of $\alpha(\tau)$ can force the vanishing relations \eqref{eq:n=1vanishingrelations}: in fact, in light of \eqref{eq:c_+sourcetermobstruction}, Remark~\ref{rem:fourierboundsforn=1termlocalizedaway1} as well as \eqref{eq:leadsing2Fourier} (and simple analogues for $r$ replacing $2$), the vanishing relations amount to the following equations which need to be satisfied for each $k = 1, 2, 3$ and each $j_1 = 0, 1, 2$: 
\begin{equation}\label{eq:alphamodulationconditionsexplicit1}\begin{split}
&\sum_{0\leq r\leq 2}\sum_{0\leq j_2\leq r-j_1} \gamma_{k,r-j_1-j_2}^{\pm}\\&\hspace{2cm}\cdot \int_{\tau_0}^\infty \big(\log\sigma\big)^{j_2}\cdot \Big[\zeta_{k,r,j_{1,2}}\cdot\alpha'(\sigma)\cdot \sin U(\sigma, \nu\sigma)\cdot \big(\frac{\lambda(\sigma)}{\sigma}\big)^{(k-1)\nu} + \big(\lambda(\sigma)\big)^{\nu k}\cdot g_r(\sigma)\Big]\,d\sigma\\
&=0.
\end{split}\end{equation}
Then we have 
\begin{lemma}\label{lem:alphamodulationconditionssolved} Assuming that $\big|g_r(\sigma)\big|\lesssim \sigma^{-4-}\cdot \Lambda$, there exists a linear map $g_r\rightarrow T\big(g_r\big)\in C_0^1\big([\tau_0, 2\tau_0]\big)$ such that 
\begin{align*}
\alpha'(\sigma) = T\big(g_r\big)
\end{align*}
satisfies the relations \eqref{eq:alphamodulationconditionsexplicit1} for $k = 1, 2, 3$, $j_1 = 0, 1, 2$. Moreover, we have the estimate 
\begin{align*}
\big\|\sigma^{2+}\cdot\alpha'\big\|_{L^\infty}\lesssim \Lambda.
\end{align*}
\end{lemma}
\begin{proof} We represent the singular source term due to modulating in $\alpha(\tau)$ by the following analogue of \eqref{eq:leadsingleq2}, namely 
\begin{equation}\label{eq:leadsingleqalphasource}
\sum_{0\leq r\leq 2}c_r\cdot \partial_{R}\left(\left(\nu\sigma-R\right)^{-\frac12+k\nu}\cdot \big[\log(\nu\tau - R)-\log(\nu\tau)\big]^r\right)\cdot \frac{\alpha'(\tau)}{\tau^{\frac32+k\nu}}, 
\end{equation}
where the coefficients $c_r\in \mathbb{R}\backslash\{0\}$. 
Applying the Duhamel propagator to the distorted Fourier transform (at $n = 1$) and reverting to the representation \eqref{eq:c_+sourcetermobstruction}, we infer a vanishing condition for each $j_1\in \{0, 1, 2\}$, starting with $j_1 = 2$, the top order singularity:
throughout set $\tilde{\alpha}(\sigma): =  \sin U(\sigma, \nu\sigma)\cdot \big(\frac{\lambda(\sigma)}{\sigma}\big)^{(k-1)\nu}\cdot\alpha'(\sigma)$, $\tilde{g}_r(\sigma) = \lambda^{k\nu}(\sigma)\cdot g_r(\sigma)$.
\\

{\it{$j_1 = 2$}}. We infer the simple relation
\begin{align*}
\gamma_{k,0}^{\pm}\cdot \int_{\tau_0}^\infty \big(c_2\cdot \tilde{\alpha}(\sigma) + \tilde{g}_2(\sigma)\big)\,d\sigma = 0. 
\end{align*}
Note that since we assume for now that $g_r$ is real-valued, this can be accomplished by means of real-valued $\alpha'(\sigma)$ for both signs. 
\\

{\it{$j_1 = 1$}}. This is more complicated but there is some simplification due to the first relation above: 
\begin{align*}
&-2\gamma^{\pm}_{k,1}\cdot \int_{\tau_0}^\infty \big(c_2\cdot \tilde{\alpha}(\sigma) + \tilde{g}_2(\sigma)\big)\,d\sigma\\
&+ 2\gamma^{\pm}_{k,0}\cdot \int_{\tau_0}^\infty\big(\log(\nu\sigma) - \log(\lambda(\sigma))\big)\cdot \big(c_2\cdot \tilde{\alpha}(\sigma) + \tilde{g}_2(\sigma)\big)\,d\sigma\\
&-\gamma^{\pm}_{k,0}\cdot \int_{\tau_0}^\infty \big(c_1\cdot \tilde{\alpha}(\sigma) + \tilde{g}_1(\sigma)\big)\,d\sigma\\
& = 0. 
\end{align*}
The first line here vanishes due to the condition for $j_1 = 2$, while the rest can be accomplished via real-valued $\alpha'(\sigma)$ if $g_1$ is real-valued; we shall verify this formally below. 
\\

{\it{$j_1 = 0$}}. Again things simplify on account of the preceding cases. We obtain 
\begin{align*}
&\gamma^{\pm}_{k,2}\cdot \int_{\tau_0}^\infty \big(c_2\cdot \tilde{\alpha}(\sigma) + \tilde{g}_2(\sigma)\big)\,d\sigma\\
& + 2\gamma^{\pm}_{k,1}\cdot \int_{\tau_0}^\infty \big(-\log(\nu\sigma) + \log(\lambda(\sigma))\big)\cdot  \big(c_2\cdot \tilde{\alpha}(\sigma) + \tilde{g}_2(\sigma)\big)\,d\sigma\\
&+c_1\gamma^{\pm}_{k,1}\cdot \int_{\tau_0}^\infty \big(c_1\cdot \tilde{\alpha}(\sigma) + \tilde{g}_1(\sigma)\big)\,d\sigma\\
& + \gamma^{\pm}_{k,0}\cdot \int_{\tau_0}^\infty\big(\log(\lambda(\sigma)) - \log(\nu\sigma)\big)^2\cdot \big(c_2\cdot \tilde{\alpha}(\sigma) + \tilde{g}_2(\sigma)\big)\,d\sigma\\
& + \gamma^{\pm}_{k,0}\cdot \int_{\tau_0}^\infty\big(\log(\lambda(\sigma)) - \log(\nu\sigma)\big)\cdot \big(c_1\cdot \tilde{\alpha}(\sigma) + \tilde{g}_1(\sigma)\big)\,d\sigma\\
& + \gamma^{\pm}_{k,0}\cdot \int_{\tau_0}^\infty\big(c_3\cdot \tilde{\alpha}(\sigma) + \tilde{g}_0(\sigma)\big)\,d\sigma.
\end{align*}
The first line again vanishes due to the condition for $j_1 = 2$, while the sum of the second and third line vanish due to the condition for $j_1 = 1$. The sum of the last three lines can be made to vanish by picking $\alpha'(\tau)$ suitably as a real function, provided all $\tilde{g}_r$ are real valued. 
To prove the existence of $\alpha'(\tau)$ as asserted in the lemma rigorously, we make the following observations: 
For $\tau_0$ sufficiently large, we have $U(\tau,\nu\tau)\sim \tau^{-1}\cdot\log\tau$ for $\tau\geq \tau_0$, and the map 
\begin{align*}
Z: C_0^1\big([\tau_0, 2\tau_0]\big)\longrightarrow \mathbb{R}^9
\end{align*}
given by $Z = \left(Z_{kj}\right)_{\substack{1\leq k\leq 3\\0\leq j\leq 2}}$ with 
\[
Z_{kj}(\phi) = \tau_0^{-k}\cdot\int_{\tau_0}^{2\tau_0}\phi(\sigma)\cdot \sin U(\sigma, \nu\sigma)\cdot \big(\frac{\lambda(\sigma)}{\sigma}\big)^{(k-1)\nu}\cdot \log^j\sigma\,d\sigma
\]
is boundedly invertible for $\tau_0$ large enough, with 
\begin{align*}
\big\|Z^{-1}({\bf{p}})\big\|_{L^\infty}\lesssim \big\|{\bf{p}}\big\|. 
\end{align*}
Determining the components of ${\bf{p}}$ as before inductively, and inverting $Z$, the bound for $\phi = \alpha'$ follows easily from the assumption on the $g_r$. 
\end{proof}
\begin{remark}\label{rem:lem:alphamodulationconditionssolved} The fact that the integrals defining the components $Z_{kj}$ in the preceding proof all involve the factor $\sin U(\sigma, \nu\sigma)$ implies that we can always independently prescribe the terminal value of $\alpha(\tau)$, i .e. 
\[
\alpha_{\infty} = \int_{\tau_0}^{2\tau_0}\alpha'(\sigma)\,d\sigma. 
\]
\end{remark}

Let us now see how picking suitable $\beta$ in addition to $\alpha$  allows us to force the full vanishing condition \eqref{eq:n=1vanishingrelations}, without assuming the functions $g_r$ in \eqref{eq:leadsingleq2} to be real valued.
Since we are in the process of describing a number of moment conditions on $\alpha, \beta$, we may as well assume that $\alpha$ is constant on the support of $\beta$ in order to simplify the computation, and equals its limiting value there. 
Denoting by $\left(\partial_t^2\right)'$ only those terms where at least one derivative falls on $\beta(t)$, 
\begin{align*}
	&\left(\partial_t^2\right)'\left(\mathcal{R}^{\alpha,\beta(t)}_{0}\Phi\right)\\=
	&\mathcal{R}^{\alpha,\beta(t)}_{0}\beta''(t)\left(\begin{array}{ccc}\cos\alpha&0&-\sin\alpha\\0&1&0\\\sin\alpha&0&\cos\alpha\end{array}\right)\cdot\left(\begin{array}{ccc}1&0&0\\0&\cos\beta(t)&-\sin\beta(t)\\0&\sin\beta(t)&\cos\beta(t)\end{array}\right)\\&\hspace{3cm}\cdot \left(\begin{array}{c}0\\-\sin\beta\sin\theta\sin U - \cos\beta\sin\alpha\cos\theta\sin U+\cos U\cos\alpha\cos\beta\\-\cos\beta\sin\theta\sin U + \sin\beta\sin\alpha\sin U\cos\theta - \cos U\cos\alpha\sin\beta\end{array}\right)\\
	& + \mathcal{R}^{\alpha,\beta(t)}_{0}2\beta'(t)\lambda'(t)rQ'(R)\left(\begin{array}{ccc}\cos\alpha&0&-\sin\alpha\\0&1&0\\\sin\alpha&0&\cos\alpha\end{array}\right)\cdot\left(\begin{array}{ccc}1&0&0\\0&\cos\beta(t)&-\sin\beta(t)\\0&\sin\beta(t)&\cos\beta(t)\end{array}\right)\\&\hspace{3cm}\cdot \left(\begin{array}{c}0\\-\sin\beta\sin\theta\cos U - \cos\beta\sin\alpha\cos\theta\cos U-\sin U\cos\alpha\cos\beta\\-\cos\beta\sin\theta\cos U + \sin\beta\sin\alpha\cos U\cos\theta +\sin U\cos\alpha\sin\beta\end{array}\right)\\
	& +  \mathcal{R}^{\alpha,\beta(t)}_{0}2\beta'(t)\lambda(t)\cdot\left[\epsilon_{\tau} + \epsilon_X\right]\cdot\left(\begin{array}{ccc}\cos\alpha&0&-\sin\alpha\\0&1&0\\\sin\alpha&0&\cos\alpha\end{array}\right)\cdot\left(\begin{array}{ccc}1&0&0\\0&\cos\beta(t)&-\sin\beta(t)\\0&\sin\beta(t)&\cos\beta(t)\end{array}\right)\\&\hspace{3cm}\cdot \left(\begin{array}{c}0\\-\sin\beta\sin\theta\cos U - \cos\beta\sin\alpha\cos\theta\cos U-\sin U\cos\alpha\cos\beta\\-\cos\beta\sin\theta\cos U + \sin\beta\sin\alpha\cos U\cos\theta +\sin U\cos\alpha\sin\beta\end{array}\right)
\end{align*}
We compute 
\begin{align*}
	&\left(\begin{array}{ccc}\cos\alpha&0&-\sin\alpha\\0&1&0\\\sin\alpha&0&\cos\alpha\end{array}\right)\cdot\left(\begin{array}{ccc}1&0&0\\0&\cos\beta(t)&-\sin\beta(t)\\0&\sin\beta(t)&\cos\beta(t)\end{array}\right)\\&\hspace{3cm}\cdot \left(\begin{array}{c}0\\-\sin\beta\sin\theta\sin U - \cos\beta\sin\alpha\cos\theta\sin U+\cos U\cos\alpha\cos\beta\\-\cos\beta\sin\theta\sin U + \sin\beta\sin\alpha\sin U\cos\theta - \cos U\cos\alpha\sin\beta\end{array}\right)\\
	& = \left(\begin{array}{c}\sin\alpha\sin\theta\sin U\\ -\sin\alpha\cos\theta\sin U + \cos\alpha\cos U\\-\cos\alpha\sin\theta\sin U\end{array}\right)=\cos\alpha\sin\theta E_{1}+\cos\alpha\cos\theta\cos U E_{2}+(*)\Phi
\end{align*}
as well as 
\begin{align*}
	&\left(\begin{array}{ccc}\cos\alpha&0&-\sin\alpha\\0&1&0\\\sin\alpha&0&\cos\alpha\end{array}\right)\cdot\left(\begin{array}{ccc}1&0&0\\0&\cos\beta(t)&-\sin\beta(t)\\0&\sin\beta(t)&\cos\beta(t)\end{array}\right)\\&\hspace{3cm}\cdot \left(\begin{array}{c}0\\-\sin\beta\sin\theta\cos U - \cos\beta\sin\alpha\cos\theta\cos U-\sin U\cos\alpha\cos\beta\\-\cos\beta\sin\theta\cos U + \sin\beta\sin\alpha\cos U\cos\theta +\sin U\cos\alpha\sin\beta\end{array}\right)\\
	& = \left(\begin{array}{c}\sin\alpha\sin\theta\cos U\\ -\sin\alpha\cos\theta\cos U - \cos\alpha\sin U\\ -\cos\alpha\sin\theta\cos U\end{array}\right)=0\cdot E_{1}-\left(\sin\alpha\cos U+\cos\alpha\sin U\cos\theta\right)E_{2}+(*)\Phi.
\end{align*}
We conclude that the contribution to the source term for $\varepsilon_+^1$ coming from modulating in $\beta$ is given by the expression (after division by $\lambda^2$)
\begin{equation}\label{eq:n=1sourcetermduetobetamodulation}
	-i\cos\alpha\cdot\left[\frac{\beta''(\tau)+\frac{\lambda_\tau}{\lambda}\beta'(\tau)}{2}\cdot\left(1+\cos U\right) - \beta'(\tau)\frac{\lambda_{\tau}}{\lambda}RQ'(R)\sin U - \beta'(t)\lambda(t)^{-1}\cdot\left[\epsilon_{\tau} + \epsilon_X\right]\sin U\right]
\end{equation}
As a consequence, proceeding in exact analogy to the angle $\alpha$, we infer 
\\

{\it{(vi)  Combined leading order contribution of smooth source terms generated by modulating in $\alpha, \beta$ to \eqref{eq:c+evolutioneqn}}}. In analogy to {\it{(iii)}} before, this is given by the formula (setting $\beta_{\infty}:=\beta(\infty)-\beta(\tau_{0})$)
\begin{align*}
	c_{+}(\tau)=&\left(c_{0}(1+\nu)-1\right)\left(\alpha_{\infty}-i\beta_{\infty}\right)\\
	&+\nu\tau^{-\nu^{-1}}\int_{\tau_{0}}^{\infty}\sigma^{1+\nu^{-1}}\left(\alpha''(\sigma)-i\beta''(\sigma)\right)\,d\sigma+(c_{0}-1)(1+\nu)\tau^{-\nu^{-1}}\int_{\tau_{0}}^{\infty}\sigma^{\nu^{-1}}\left(\alpha'(\sigma)-i\beta'(\sigma)\right)\,d\sigma\\
	&+(1+\nu)\tau^{-\nu^{-1}}\int_{\tau}^{\infty}\sigma^{\nu^{-1}}\left(\alpha'(\sigma)-i\beta'(\sigma)\right)\,d\sigma+(1-c_{0})(1+\nu)\tau^{-\nu^{-1}}\int_{\tau}^{\infty}\sigma^{\nu^{-1}}\left(\alpha'(\sigma)-i\beta'(\sigma)\right)\,d\sigma.
\end{align*}
A similar analysis as in {\it{(iii)}} on the terminating value $\alpha_{\infty}-i\beta_{\infty}$ applies here. 
\\

{\it{(vii): The combined effect of modulating on $\alpha, \beta$ on  $\lim_{R\rightarrow 0}H_1^+\varepsilon_+^1$.}} Here by proceeding analogously for $\beta$ as for $\alpha$, we conclude 
\begin{lemma}\label{lem:alphabetamodulationn=1singularterms} The same conclusion as in Lemma~\ref{lem:alphamodulationn=1singulartermspropagated} obtains but with $\alpha$ replaced by $\alpha - i\beta$. 
\end{lemma}

Up to this stage, we have only considered the effect that modulating in the angles $\alpha, \beta$ has on the $n = 1$ mode. However, there are also leading order effects on the $n = -1$ mode. Here only the effect of the singular term matters, due to the rapid decay of the fundamental system describing the $c_{-1}(\tau)$ evolution. Specifically, in light of \eqref{eq:c-1evolutioneqn}, we need to determine the leading order effect on $\lim_{R\rightarrow 0}R^{-2}H_{-1}^{+}\varepsilon_+^{-1}$ arising after applying the $n= -1$ parametrix to the singular source term generated by modulating in $\alpha, \beta$:
\\

{\it{(viii): The combined effect of modulating on $\alpha, \beta$ on the evolution on $c_{-1}(\tau)$ via $\lim_{R\rightarrow 0}R^{-2}H_{-1}^{+}\varepsilon_+^{-1}$}}. We formulate this directly in analogy to Lemma~\ref{lem:alphamodulationn=1singulartermspropagated}:  
\begin{lemma}\label{lem:alphabetaeffectonsingularn=-1} There exists a function $H(\tau, R)$ coinciding with 
	\[
	\mathcal{D}_-\left(\left(i\beta'(\tau)+\alpha'(\tau)\right)\cdot\left[\epsilon_{\tau} + \epsilon_X\right]\cdot \sin U\right)
	\]
	near the light cone $R<\nu\tau$, and such that setting 
	\[
	\xb(\tau, \xi): = \int_{\tau_0}^{\tau}U^{(-1)}\left(\tau, \sigma,\xi\right)\cdot\mathcal{F}^{(-1)}\left(H\left(\sigma, \cdot\right)\right)\left(\frac{\lambda^2(\tau)}{\lambda^2(\sigma)}\xi\right)\,d\sigma, 
	\]
	we have (where as usual $ \chi_{\lambda^2(\tau)\xi<M}$ is a smooth cutoff)
	\begin{equation}\label{eq:c_-sourcetermobstruction1}\begin{split}
		&\int_0^\infty \chi_{\lambda^2(\tau)\xi<M}\xb(\tau, \xi)\tilde{\rho}_{-1}(\xi)\,d\xi\\& = \sum_{\pm}\sum_{k=1,2,3}\sum_{0\leq r\leq 2}\sum_{0\leq j_1\leq r}\int_0^\infty \chi_{\lambda^2(\tau)\xi<M}\frac{e^{\pm i\nu\tau\xi^{\frac12}}}{\xi^{\frac12+\frac{k\nu}{2}}}\left[\log\left(\xi\lambda^2(\tau)\right)\right]^{j_1}\,d\xi\\&\hspace{6cm}\cdot \sum_{0\leq j_2\leq r-j_1} \int_{\tau_0}^{\tau}\left(\log\sigma\right)^{j_2}\left[\frac{\lambda(\sigma)}{\lambda(\tau)}\right]^{\nu k}\tilde{\tilde{\beta}}_{\pm}^{(k,r,j_1,j_2)}(\sigma)\,d\sigma\\
		& + \tilde{\tilde{c}}_{-1}(\tau), 
	\end{split}\end{equation}
	where the error satisfies 
	\[
	\left|\tilde{\tilde{c}}_{-1}(\tau)\right|\lesssim \tau^{-4-}\cdot  \Lambda\cdot(\tau_0^{-1}+\Lambda),
	\]
	and the functions $\tilde{\tilde{\beta}}_{\pm}^{(k,j_1,j_2)}(\sigma)$ are explicitly given by 
	\begin{align*}
	\tilde{\beta}_{\pm}^{(k,r,j_1,j_2)}(\sigma) = \sigma^{-(k-1)\nu}\cdot\lambda^{-\nu}(\sigma)\cdot \gamma_{k,r-j_1-j_2}^{\pm}\cdot\zeta_{k,r,j_{1,2}}\cdot \big(\alpha'(\sigma)+i\beta'(\tau)\big)\cdot \sin U(\sigma, \nu\sigma). 
	\end{align*}
\end{lemma}
\begin{proof}
	The proof  is similar to that of Lemma \ref{lem:alphamodulationn=1singulartermspropagated}.
\end{proof}
Lemma~\ref{lem:alphabetaeffectonsingularn=-1} and Lemma~\ref{lem:alphabetamodulationn=1singularterms} give the principal effect of modulating on the angles $\alpha,\beta$ on the principal ingoing singular part of the $n = \mp 1$ modes and from there to the delicate source terms 
\[
\lim_{R\rightarrow 0}R^{-2}H_{-1}^{+}\varepsilon_+^{-1},\quad  \lim_{R\rightarrow 0}H_1^+\varepsilon_+^1
\]
in the ODEs for $c_{-1}(\tau), c_{+1}(\tau)$. In order to force the vanishing conditions \eqref{eq:n=1vanishingrelations}, \eqref{eq:n=-1vanishingrelations}, we have to complement the effect of modulating on $\alpha, \beta$ by the effect of modulating on the Lorentz transform parameters $v_1, v_2$. 

\subsubsection{The contribution of the Lorentz parameters $v_1(t), v_2(t)$.} Since the parameters will be very small (their size depending on the initial perturbation), we shall neglect terms quadratic in the $v_j$ in the ensuing discussion. This means that for some estimates, it is permissible to replace the Lorentz transform by the simpler Galilean transform $\bf{x}\longrightarrow {\bf{x}} - t\cdot{\bf{v}}$, ${\bf{v}} = \left(\begin{array}{c}v_1\\ v_2\end{array}\right)$. Of course as far as the effect of the Lorentz transform on the singular terms is concerned, the precise structure will be of crucial importance. Observe that modulating on the Lorentz parameters will affect both the $n = +1$ and the $n = -1$ modes. We first analyze the effect on the bulk term $Q$ and thereby directly on the evolution of $c_+(\tau)$ via $\lim_{R\rightarrow 0}F_+(1)$ (compare with \eqref{eq:c+evolutioneqn}), and which will turn out to be of small order $O(|\bfv|^2)$, and then we shall analyze the contribution to the source terms 
\[
\lim_{R\rightarrow 0}R^{-2}H_{-1}^{+}\varepsilon_+^{-1},\quad  \lim_{R\rightarrow 0}H_1^+\varepsilon_+^1
\]
via the effect on the singular part of $U$. 
\\

{\it{Effect of Lorentz modulating on the $n = +1$ mode via the bulk part}}. Here we may replace the bulk part $Q(R)$ by $2R$, and neglecting terms which vanish at the origin $R = 0$, the bulk part $\Phi$ gets replaced by 
\[
\left(\begin{array}{c}\lambda(t)x_1\\ \lambda(t)x_2\\ 1\end{array}\right)
\]
which gets mapped into 
\[
\left(\begin{array}{c}\lambda(t - {\bf{x}}\cdot{\bf{v}})(x_1 - tv_1)\\ \lambda(t - {\bf{x}}\cdot{\bf{v}})(x_2 -tv_2)\\ 1\end{array}\right).
\] 
Using  $\left(\partial_t^2\right)'$ as usual the operator where at least one derivative falls on a component $v_j(t)$, we infer 
\begin{align*}
	\left(\partial_t^2\right)'\left(\begin{array}{c}\lambda(t - {\bf{x}}\cdot{\bf{v}})(x_1 - tv_1)\\ \lambda(t - {\bf{x}}\cdot{\bf{v}})(x_2 -tv_2)\\ 1\end{array}\right) = \left(\begin{array}{c}-2\lambda'(t)\cdot tv_1'(t) - \lambda(t)\cdot (tv_1)''\\ -2\lambda'(t)\cdot tv_2'(t)- \lambda(t)\cdot (tv_2)''\\0\end{array}\right)+\text{error}, 
\end{align*}
where we again neglect terms quadratic in ${\bf{v}}$ and the terms which vanishes at $R=0$. Projecting the preceding onto the frame $\{E_1, E_2\}$ and extracting the $\vphi_1^{(1)} - i\vphi_2^{(1)}$-component leads after division by $\lambda^{-2}$ and conversion to the variable $\tau$ to 
\begin{align*}
	-\nu\left[\tau v_{1,\tau\tau} + 2v_{1,\tau} - i(\tau v_{2,\tau\tau} + 2v_{2,\tau})\right]-(1+\nu)\left[v_{1,\tau}-iv_{2,\tau}\right]
\end{align*}
The contribution from the first bracket above to $c_{+}$ is trivial in the leading order. In fact, for any twice continuously differentiable function $f(\tau)$ with $f_{\tau} = 0$ for large enough $\tau$ and $f_{\tau}(\tau_0) = 0$, we have 
\begin{align*}
	\int_{\tau_0}^{\tau}\sigma\cdot\left[\sigma f_{\sigma\sigma} + 2f_{\sigma}\right]\,d\sigma& =\int_{\tau_0}^{\tau}\sigma\cdot (\sigma f)_{\sigma\sigma}\,d\sigma\\
	& = \sigma(\sigma f)_{\sigma}|_{\tau_0}^{\tau} - \int_{\tau_0}^{\tau}(\sigma f)_{\sigma}\,d\sigma\\
	& = \sigma^2 f_{\sigma}|_{\tau_0}^{\tau} + \sigma f(\sigma)|_{\tau_0}^{\tau} -  \sigma f(\sigma)|_{\tau_0}^{\tau} = 0
\end{align*}
for large enough $\tau$, and since we can impose these requirements on $v_{1,2}(\sigma)$.  The contribution from the second bracket is 
\begin{align*}
	&-(1+\nu)\int_{\tau_{0}}^{\tau}\sigma\left(v_{1,\sigma}(\sigma)-iv_{2,\sigma}(\sigma)\right)\,d\sigma\\
	=&-(1+\nu)\left(\tau\left(v_{1}(\tau)-iv_{2}(\tau)\right)-\tau_{0}\left(v_{1}(\tau_{0})-iv_{2}(\tau_{0})\right)\right)+(1+\nu)\int_{\tau_{0}}^{\tau}\left(v_{1}(\sigma)-iv_{2}(\sigma)\right)\,d\sigma.
\end{align*}
In view of the moment condition satisfied by $v_{1}, v_{2}$, we don't force the expression on the RHS above to vanish, but instead we choose $\alpha_{\infty}-i\beta_{\infty}$ appropriately (since the moment conditions satisfied by $\alpha,\beta$ are weaker), to cancel both the non-trivial contribution of $v_{1},v_{2}$ and the contribution from the constant fundamental solution of the operator $\frac{d^{2}}{d\tau^{2}}+\frac{\lambda'}{\lambda}\frac{d}{d\tau}$.
\\

{\it{Effect of Lorentz modulating on the singular parts of the $n = \pm 1$ modes}}. Here we finally complete forcing the vanishing conditions \eqref{eq:n=1vanishingrelations}, \eqref{eq:n=-1vanishingrelations} to leading order by combining the effects of modulating on the angles $\alpha, \beta$ as well as on the Lorentz parameters. For this we have to analyze the leading order effect of Lorentz transforming the singular term $\epsilon$ inherent in the bulk term $U = Q + \epsilon$. Recall the definitions in subsection~\ref{subsec:introshock}, and write $\epsilon = \epsilon(t, t-r) =:\epsilon(t, \tilde{X})$, where we have 
\[
\tilde{X} = t - r = \frac{t^2 - r^2}{t+r}
\]
Then observe that setting $X = \frac{t^2 - r^2}{2t}$, we have 
\[
\epsilon(t, \tilde{X}) = \epsilon(t, X) + \text{error},
\]
where the error term is one degree smoother than the principal term $\epsilon(t, X)$, and thus in terms of the top order singularity, it suffices to work with $ \epsilon(t, X)$. We now determine the maximally singular error terms generated by Lorentz modulating. Recalling that 
\begin{align*}
	\Phi = \left(\begin{array}{c} \sin U \cos\theta\\ \sin U\sin\theta\\ \cos U\end{array}\right) = \left(\begin{array}{c} \sin U \frac{x_1}{\sqrt{x_1^2+x_2^2}}\\ \sin U\frac{x_2}{\sqrt{x_1^2+x_2^2}}\\ \cos U\end{array}\right),\quad U = Q+\epsilon,
\end{align*}
and (neglecting terms quadratic in $v$) replacing $x_1\longrightarrow x_1 - tv_1(t),\,x_2\longrightarrow x_2 - tv_2(t),\,t\longrightarrow t_{\bf{v}}: =  t - \textbf{x}\cdot\textbf{v} (t)$, while $\epsilon = \epsilon(t, X)$, and denoting by $(\partial_t^2)'$ the operator where at least one derivative falls on a factor $v_j(t)$, we compute 
\begin{align}\label{Lorentz singular}
	\begin{split}
	(\partial_t^2)'\left(\mathcal{L}_{v(t)}\Phi\right) &= 2\left(\begin{array}{c}\cos(Q+\epsilon)\epsilon_X\cdot \frac{\partial X}{\partial t}\cdot \frac{-tv_1'}{\sqrt{(x_1 - tv_1(t))^2 + (x_2 - tv_2(t))^2}}\\ \cos(Q+\epsilon)\epsilon_X\cdot \frac{\partial X}{\partial t}\cdot \frac{-tv_2'}{\sqrt{(x_1 - tv_1(t))^2 + (x_2 - tv_2(t))^2}}\\ 0\end{array}\right)\\
	& +  2\left(\begin{array}{c}\cos(Q+\epsilon)\epsilon_X\cdot \frac{\partial X}{\partial t}\cdot \frac{-(-tv_1')(x_1 - tv_1(t))^2 - (-tv_2')(x_1 - tv_1(t))(x_2 - tv_2(t))}{\left[\sqrt{(x_1 - tv_1(t))^2 + (x_2 - tv_2(t))^2}\right]^{3}}\\ \cos(Q+\epsilon)\epsilon_X\cdot \frac{\partial X}{\partial t}\cdot \frac{-(-tv_2')(x_2 - tv_2(t))^2 - (-tv_1')(x_1 - tv_1(t))(x_2 - tv_2(t))}{\left[\sqrt{(x_1 - tv_1(t))^2 + (x_2 - tv_2(t))^2}\right]^3}\\ 0\end{array}\right)\\
	& + \left(\begin{array}{c}\cos(Q+\epsilon){\bf{v}}'(t)\cdot\nabla_{\bf{v}}(t_{\bf{v}})\epsilon_{tX} \frac{x_1-tv_1}{\sqrt{(x_1 - tv_1(t))^2 + (x_2 - tv_2(t))^2}}\\ \cos(Q+\epsilon){\bf{v}}'(t)\cdot\nabla_{\bf{v}}(t_{\bf{v}})\epsilon_{tX} \frac{x_2-tv_2}{\sqrt{(x_1 - tv_1(t))^2 + (x_2 - tv_2(t))^2}}\\ -\sin(Q+\epsilon){\bf{v}}'(t)\cdot\nabla_{\bf{v}}(t_{\bf{v}})\epsilon_{tX}\end{array}\right)\\
	& +  \left(\begin{array}{c}\cos(Q+\epsilon)\left[(\partial_t^2)'X\right]\epsilon_{X} \frac{x_1-tv_1}{\sqrt{(x_1 - tv_1(t))^2 + (x_2 - tv_2(t))^2}}\\ \cos(Q+\epsilon)\left[(\partial_t^2)'X\right]\epsilon_{X} \frac{x_2-tv_2}{\sqrt{(x_1 - tv_1(t))^2 + (x_2 - tv_2(t))^2}}\\ -\sin(Q+\epsilon)\left[(\partial_t^2)'X\right]\epsilon_{X}\end{array}\right)\\
	& +  \left(\begin{array}{c}\cos(Q+\epsilon)(\partial_t)'X\cdot\epsilon_{XX} \frac{x_1-tv_1}{\sqrt{(x_1 - tv_1(t))^2 + (x_2 - tv_2(t))^2}}\\ \cos(Q+\epsilon)(\partial_t)'X\cdot\epsilon_{XX} \frac{x_2-tv_2}{\sqrt{(x_1 - tv_1(t))^2 + (x_2 - tv_2(t))^2}}\\ -\sin(Q+\epsilon)(\partial_t)'X\epsilon_{XX}\end{array}\right)\\
	& + \text{error},
	\end{split}
\end{align}
where the last term denotes expressions that are either quadratic in $v$ or one degree smoother than $\epsilon_X$. 
Observe that we can simplify the sum of the first two vectors on the right to 
\begin{align*}
	\left(\begin{array}{c}\cos(Q+\epsilon)\epsilon_X\cdot \frac{\partial X}{\partial t}\cdot \frac{(-tv_1')(x_2 - tv_2(t))^2 - (-tv_2')(x_1 - tv_1(t))(x_2 - tv_2(t))}{\left[\sqrt{(x_1 - tv_1(t))^2 + (x_2 - tv_2(t))^2}\right]^{3}}\\ \cos(Q+\epsilon)\epsilon_X\cdot \frac{\partial X}{\partial t}\cdot \frac{(-tv_2')(x_1 - tv_1(t))^2 - (-tv_1')(x_1 - tv_1(t))(x_2 - tv_2(t))}{[\sqrt{(x_1 - tv_1(t))^2 + (x_2 - tv_2(t))^2}]^3}\\ 0\end{array}\right),
\end{align*}
which, up to more regular terms, is equal to
\begin{align}\label{Lorentz singular 1}
	\left(\begin{array}{c}
		\cos U\eps_{X}\cdot\left(\sin\theta\cos\theta v'_{2}-\sin^{2}\theta v'_{1}\right)\\
		\cos U\eps_{X}\cdot\left(\sin\theta\cos\theta v'_{1}-\cos^{2}\theta v'_{2}\right)\\
		0
	\end{array}\right)=\cos U\eps_{X}\left(-\cos\theta v'_{2}(t)+\sin\theta v'_{1}(t)\right)E_{2}.
\end{align} 
Here we also used the fact that up to more regular terms and quadratic terms in $\bfv$, we have $\frac{\partial X}{\partial t}\simeq 1$. The last three terms in \eqref{Lorentz singular} can be written as
\begin{align*}
	\left({\bf{v}}'(t)\cdot\nabla_{\bf{v}}(t_{\bf{v}})\epsilon_{tX} +\left[(\partial_t^2)'X\right]\epsilon_{X} +(\partial_t)'X\cdot\epsilon_{XX} \right)E_{1}.
\end{align*}
Again up to the more regular terms as well as the quadratic terms in $|\bfv|$, we have
\begin{align*}
	(\partial_{t})'X\simeq &(t-r)\cdot \left(\cos\theta v'_{1}(t)+\sin\theta v'_{2}(t)\right),\\(\partial^{2}_{t})'X\simeq &\left(\cos\theta v'_{1}(t)+\sin\theta v'_{2}(t)\right),\\
	\bfv'(t)\cdot\nabla_{\bfv}(t_{\bfv})\simeq& -t\left(\cos\theta v'_{1}(t)+\sin\theta v'_{2}(t)\right)
\end{align*}
Therefore the last three terms in \eqref{Lorentz singular} is equal to, up to the more regular terms and the quadratic terms in $|\bfv|$,
\begin{align}\label{Lorentz singular 2}
	\begin{split}
		\left(-t\left(\cos\theta v'_{1}(t)+\sin\theta v'_{2}(t)\right)\eps_{tX}+\left(\cos\theta v'_{1}(t)+\sin\theta v'_{2}(t)\right)\eps_{X}+(t-r)\cdot \left(\cos\theta v'_{1}(t)+\sin\theta v'_{2}(t)\right)\eps_{XX}\right)E_{1}
	\end{split}
\end{align}
Moreover, a direct calculation shows
\begin{align*}
\frac{1}{2\pi}\int_{0}^{2\pi}\left(\cos\theta v'_{1}(t)+\sin\theta v'_{2}(t)\right)e^{\pm i\theta}\,d\theta=\frac12\left(v'_{1}(t)\pm iv'_{2}(t)\right),
\end{align*}
and 
\begin{align*}
	\frac{1}{2\pi}\int_{0}^{2\pi}\left(-\cos\theta v'_{2}(t)+\sin\theta v'_{1}(t)\right)e^{\pm i\theta}\,d\theta=-\frac12 v'_{2}(t)\pm \frac{i}{2}v'_{1}(t).
\end{align*}
Interpreting the terms in \eqref{Lorentz singular} as $\mathcal{L}_{v(t)}\tilde{\Phi}$, and writing 
\[
\tilde{\Phi} = \tilde{\vphi}_1 E_1 + \tilde{\vphi}_2 E_2 + (*)\Phi, 
\]
we easily infer the following formula for the contribution to $\tilde{\vphi}_1^{(1)} - i\tilde{\vphi}_2^{(1)}$: 
\begin{equation}\label{eq:n=1contributionfromLorentzmodulationsingular}
	\tilde{\vphi}_1^{(1)} - i\tilde{\vphi}_2^{(1)} = \frac12\left((t-r)\epsilon _{XX} - t\epsilon_{tX} + \epsilon_X\right)(v_1'(t) - iv_2'(t)) - (\cos U) \epsilon_X(v_1'(t) - iv_2'(t)). 
\end{equation}
Similarly, as for the contribution to the $n = -1$ mode, we find 
\begin{equation}\label{eq:n=-1contributionfromLorentzmodulationsingular}
	\tilde{\vphi}_1^{(-1)} - i\tilde{\vphi}_2^{(-1)} = \frac12\left((t-r)\epsilon _{XX} - t\epsilon_{tX} + \epsilon_X\right)(v_1'(t) + iv_2'(t)) - (\cos U) \epsilon_X(v_1'(t) + iv_2'(t)). 
\end{equation}
Observe here that the expression $(t-r)\epsilon _{XX} - t\epsilon_{tX} + \epsilon_X$ is again a type of singularity as described in subsection~\ref{subsec:introshock}, specifically of type $\mathcal{Q}'$. Normalizing things by dividing by $\lambda^{-2}$, we can then formulate an analogue of Lemma~\ref{lem:alphabetaeffectonsingularn=-1} , Lemma~\ref{lem:alphabetamodulationn=1singularterms}, which takes into account the combined effect of modulating on the angles $\alpha,\beta$ as well as the Lorentz parameters $v_1, v_2$, in terms of the effect on the top order singular terms and their contribution to the source terms for the ODEs governing $c_{+}(\tau), c_-(\tau)$. 
\begin{lemma}\label{lem:alphabetavonetwoeffectonsingularn=plusminus1} There exists a function $H(\tau, R)$ coinciding with 
	\[
	\mathcal{D}_+\left(\eqref{eq:alphamodulationn=1B} + \eqref{eq:n=1sourcetermduetobetamodulation} + \lambda^{-2}\cdot\eqref{eq:n=1contributionfromLorentzmodulationsingular}
	\right)
	\]
	on the light cone $R<\nu\tau$, and setting
	\[
	\xb(\tau, \xi): = \int_{\tau_0}^{\tau}U^{(1)}\left(\tau, \sigma,\xi\right)\cdot\mathcal{F}^{(1)}\left(H(\sigma, \cdot)\right)\left(\frac{\lambda^2(\tau)}{\lambda^2(\sigma)}\xi\right)\,d\sigma, 
	\]
	we have (where as usual $ \chi_{\lambda^2(\tau)\xi<M}$ is a smooth cutoff)
	\begin{align*}
		&\int_0^\infty \chi_{\lambda^2(\tau)\xi<M}\xb(\tau, \xi)\tilde{\rho}_{1}(\xi)\,d\xi\\& =  \sum_{\pm}\sum_{k=1,2,3}\sum_{0\leq r\leq 2}\sum_{0\leq j_1\leq r}\int_0^\infty \chi_{\lambda^2(\tau)\xi<M}\frac{e^{\pm i\nu\tau\xi^{\frac12}}}{\xi^{\frac{k\nu}{2}}}\left[\log\left(\xi\lambda^2(\tau)\right)\right]^{j_1}\,d\xi\\&\hspace{6cm}\cdot \sum_{0\leq j_2\leq r-j_1} \int_{\tau_0}^{\tau}\left(\log\sigma\right)^{j_2}\left[\frac{\lambda(\sigma)}{\lambda(\tau)}\right]^{\nu k} B_{\pm}^{(k,r,j_1,j_2)}(\sigma)\,d\sigma\\
		& + C_1(\tau), 
	\end{align*}
	where the error satisfies 
	\[
	\left|C_1(\tau)\right|\lesssim \tau^{-4-}\cdot  \Lambda\cdot(\tau_0^{-1}+\Lambda),
	\]
	and the functions $B_{\pm}^{(k,j_1,j_2)}(\sigma)$ are explicitly given by 
	\begin{align*}
	B_{\pm}^{(k,r,j_1,j_2)}(\sigma) &= \sigma^{-(k-1)\nu}\cdot\lambda^{-\nu}(\sigma)\cdot \gamma_{k,r-j_1-j_2}^{\pm}\cdot\zeta_{k,r,j_{1,2}}\cdot (\alpha'(\sigma) - i\beta'(\sigma))\cdot \sin U(\sigma, \nu\sigma)\\
	& + \sigma^{-(k-1)\nu}\cdot\lambda^{-\nu}(\sigma)\cdot \gamma_{k,r-j_1-j_2}^{\pm}\cdot\left(\kappa_{k,r,j_{1,2}}+O(\sigma^{-1})\right)\cdot (\nu_1'(\sigma) - i\nu_2'(\sigma)).
	\end{align*}
        for non-vanishing real constants $\zeta_{k,r,j_{1,2}}, \kappa_{k,r,j_{1,2}}$, the error $O(\sigma^{-1})$ being real-valued as well.  
		
	Replacing $n=+1$ by $n = -1$ and proceeding analogously to the preceding, we infer mutatis mutandis the following formula (where we think of all parameters $\alpha, \beta, v_{1,2}$ as functions of $\tau$, i.e., the re-scaled time):
	\begin{align*}
		&\int_0^\infty \chi_{\lambda^2(\tau)\xi<M}\xb(\tau, \xi)\tilde{\rho}_{-1}(\xi)\,d\xi\\& =  \sum_{\pm}\sum_{k=1,2,3}\sum_{0\leq r\leq 2}\sum_{0\leq j_1\leq r}\int_0^\infty \chi_{\lambda^2(\tau)\xi<M}\frac{e^{\pm i\nu\tau\xi^{\frac12}}}{\frac12\xi^{\frac{k\nu}{2}}}\left[\log\left(\xi\lambda^2(\tau)\right)\right]^{j_1}\,d\xi\\&\hspace{6cm}\cdot \sum_{0\leq j_2\leq r-j_1} \int_{\tau_0}^{\tau}\left(\log\sigma\right)^{j_2}\left[\frac{\lambda(\sigma)}{\lambda(\tau)}\right]^{\nu k} \tilde{B}_{\pm}^{(k,r,j_1,j_2)}(\sigma)\,d\sigma\\
		& + \tilde{C}_{-1}(\tau),
		\end{align*}
	where the error satisfies 
	\[
	\left|\tilde{C}_{-1}(\tau)\right|\lesssim \tau^{-4-}. 
	\]
	The functions $ \tilde{B}_{\pm}^{(k,r,j_1,j_2)}$ admit an expansion similar to the one above: 
	\begin{align*}
	\tilde{B}_{\pm}^{(k,r,j_1,j_2)}(\sigma) &= -\sigma^{-(k-1)\nu}\cdot\lambda^{-\nu}(\sigma)\cdot \gamma_{k,r-j_1-j_2}^{\pm}\cdot\zeta_{k,r,j_{1,2}}\cdot (\alpha'(\sigma) + i\beta'(\sigma))\cdot \sin U(\sigma, \nu\sigma)\\
	& + \sigma^{-(k-1)\nu}\cdot\lambda^{-\nu}(\sigma)\cdot \gamma_{k,r-j_1-j_2}^{\pm}\cdot\big(\kappa_{k,r,j_{1,2}}+O(\sigma^{-1})\big)\cdot (\nu_1'(\sigma) + i\nu_2'(\sigma)).
	\end{align*}
	for 
	\end{lemma}
\begin{proof} The proof is analogous to the one of Lemma~\ref{lem:alphamodulationn=1singulartermspropagated} by relying on an analogue of Lemma~\ref{lem:alphamodulationn=1singularterms} which in turn follows from \eqref{eq:n=1contributionfromLorentzmodulationsingular}, \eqref{eq:n=-1contributionfromLorentzmodulationsingular}.
\end{proof}
By taking advantage of the freedom to specify the quadruple of real valued functions $\alpha, \beta, v_1, v_2$, we can finally enforce the required vanishing conditions \eqref{eq:n=-1vanishingrelations}, \eqref{eq:n=1vanishingrelations}:

\begin{proposition}\label{prop:forcedvanishingatn=plusminus1} 
Assume the functions $\beta_{\pm}^{k,j_1,j_2}$ are as in \eqref{eq:leadsing2Duhamelspecialcoefficients}, and similarly for $\tilde{\beta}_{\pm}^{k,j_1,j_2}$, with the latter defined in terms of $\tilde{g}_2(\sigma)$. Assume we have the bound 
\begin{align*}
\big|\alpha_{\infty}\big| + \big|\beta_{\infty}\big|+\big|g_2(\sigma)\big| + \big|\tilde{g}_2(\sigma)\big|\lesssim \sigma^{-4-}\cdot (\tau_0^{-1}+\Lambda)\cdot \Lambda\,\forall \sigma\geq \tau_0. 
\end{align*}
Here the functions $g_2, \tilde{g}_2$ are complex valued and continuous. Then there is a $\mathbb{R}$-linear map 
\[
S: \big(\alpha_{\infty}, \beta_{\infty}, g_2,\,\tilde{g}_2\big)\longrightarrow \Big(C_0^0\big([\tau_0, 2\tau_0];\mathbb{R}\big)\Big)^4
\]
 such that if we set 
\[
\big(\alpha',\,\beta',\,\nu_1',\,\nu_2'\big) = S\big(\alpha_{\infty}, \beta_{\infty}, g_2,\,\tilde{g}_2\big)
\]
then defining $B_{\pm}^{k,r,j_1,j_2},\,\tilde{B}_{\pm}^{k,r,j_1,j_2}$ as in Lemma~\ref{lem:alphabetavonetwoeffectonsingularn=plusminus1}, we have the vanishing relations 
\begin{align*}
&\sum_{j_2\leq 2-j_1}\int_{\tau_0}^{\infty}\left(\log\sigma\right)^{j_2}\left[\lambda(\sigma)\right]^{\nu k}\big[\beta_{\pm}^{(k,j_1,j_2)}(\sigma) + \sum_{j_1\leq r\leq 2}B_{\pm}^{k,r,j_1,j_2}(\sigma)\big]\,d\sigma = 0,\\
&\sum_{j_2\leq 2-j_1}\int_{\tau_0}^{\infty}\left(\log\sigma\right)^{j_2}\left[\lambda(\sigma)\right]^{\nu k}\big[\tilde{\beta}_{\pm}^{(k,j_1,j_2)}(\sigma) + \sum_{j_1\leq r\leq 2}\tilde{B}_{\pm}^{k,r,j_1,j_2}(\sigma)\big]\,d\sigma = 0\\
&\int_{\tau_0}^{\infty}\alpha'(\sigma)\,d\sigma - \alpha_{\infty} = 0,\,\int_{\tau_0}^{\infty}\beta'(\sigma)\,d\sigma - \beta_{\infty} = 0.
\end{align*}
Moreover we have the bound 
\begin{align*}
\left\|\sigma^{2+}\alpha'(\sigma)\right\|_{L^\infty} + \left\|\sigma^{2+}\beta'(\sigma)\right\|_{L^\infty} + \sum_{j=1,2}\left\|\sigma^{3+}\nu_j'(\sigma)\right\|_{L^\infty}\lesssim (\tau_0^{-1}+\Lambda)\cdot \Lambda.
\end{align*}

\end{proposition}
\begin{proof} This is analogous to the proof of Lemma~\ref{lem:alphamodulationconditionssolved}, also keeping in mind Remark~\ref{rem:lem:alphamodulationconditionssolved}, and considering the real and imaginary parts separately.

\end{proof}

At this stage, the only obstruction left to force sufficient decay for all parameters at time $\tau = \infty$ comes from \eqref{eq:n=0vanishingrelations}, arising for the $n = 0$ mode. We shall force these by modulating in the remaining angle $h(t)$ as well as the scaling parameter, meaning replacing $\lambda$ by $\tilde{\lambda}$. 

\subsubsection{The contribution of the remaining angle $h(t)$ and scaling}

At this point we have not yet exploited the final remaining rotation which acts in terms of the angle $h(t)$ via 
\[
\left(\begin{array}{ccc}\cos h(t)& -\sin h(t)&0\\\sin h(t)& \cos h(t)& 0\\ 0&0&1\end{array}\right)
\]
Alternatively, the effect on the bulk part will be that the angle $\theta$ in the representation $\Phi = \left(\begin{array}{c}\cos\theta\sin U\\ \sin \theta\sin U\\ \cos U\end{array}\right)$ gets replaced by $\theta + h(t)$. We compute $(\partial_t^2)'\mathcal{R}_{h(t)}\Phi$ to leading order, where $(\partial_t^2)'$ means at least one derivative hits $h(t)$, and we only keep track of terms linear in $h$, and which are maximally singular at the light cone, i.e., which contribute to the principal incoming part. This immediately leads to the term 
\begin{align*}
	&h'(t)\mathcal{R}_{h(t)}\circ \mathcal{R}_{-h(t)}\left(\begin{array}{ccc}-\sin h(t)& -\cos h(t)&0\\\cos h(t)& -\sin h(t)& 0\\ 0&0&0\end{array}\right)\cdot \left(\begin{array}{c}\cos\theta\cos U\epsilon_X\\ \sin \theta\cos U \epsilon_X\\ -\sin U\epsilon_X\end{array}\right)\\
	& = h'(t)\mathcal{R}_{h(t)}\left(\begin{array}{c}-\sin\theta\cos U\\ \cos\theta\cos U\\ 0\end{array}\right)\cdot\epsilon_X + \text{error}. 
\end{align*}
and in terms of the frame $\{E_1, E_2\}$, we have 
\[
h'(t)\left(\begin{array}{c}-\sin\theta\cos U\\ \cos\theta\cos U\\ 0\end{array}\right)\cdot\epsilon_X = h'(t)\cos U\epsilon_X E_2, 
\]
whence modulating in $h$ contributes to leading order a term $h'(t)\cos U\epsilon_X$ to $\vphi_2^{(0)}$, which constitutes the imaginary part of the quantity 
\[
\varepsilon_+^{(0)} = \vphi_1^{(0)} - i\vphi_2^{(0)}. 
\]
In order to handle the real part, we have to use our final modulation parameter, namely scaling. Precisely, let us set $\tilde{\lambda} = c(t)\cdot \lambda(t)$ where $\lambda(t) = t^{-1-\nu}$, and where $c(t_0) =1$ and $c'(t)$ is compactly supported on some interval $[t_1, t_0]$ with $t_1>0$.  Now computing 
\[
\left(\partial_t^2\right)'\Phi_{c(t)}, 
\]
where the subscript denotes a space-time rescaling of the function, and $\left(\partial_t^2\right)'$ indicates that at least one derivative falls on $c(t)$, we compute (again up to leading order both in terms of its dependence on $c'(t)$ as well as the singularity at the light cone), we get 
\begin{align*}
	\left(\partial_t^2\right)'\Phi_{c(t)} = \left(\begin{array}{c}\cos\theta\cos U\\ \sin\theta\cos U\\ -\sin U\end{array}\right)_{c(t)}\cdot \left(2c'\epsilon_X + 2tcc'\epsilon_{tX} + 2\epsilon_{XX}(t-r)cc'\right) + \text{error},  
\end{align*}
where $\text{error}$ also comprises top order singular terms depending linearly on $c'$ but which have an extra $\tau^{-1}$ smallness factor. Then the main term on the right will clearly only contribute to $\vphi_1^{(0)}$, namely the term 
\[
2c'\epsilon_X + 2tcc'\epsilon_{tX} + 2\epsilon_{XX}(t-r)cc',
\]
which up to smaller order terms can be equated with 
\[
c'\left(2\epsilon_X + 2t\epsilon_{tX} +  2(t-r)\epsilon_{XX}\right).
\]
Combining the effects of modulating in $h$ and in $c$, we find the following principal contribution to $\vphi_1^{(0)} - i\vphi_2^{(0)}$: 
\begin{align}\label{eq:n=0effectofmodulatingincandh}
c'\left(2\epsilon_X + 2t\epsilon_{tX} +  2(t-r)\epsilon_{XX}\right) - ih'(t)\cos U\epsilon_X
\end{align}
We can then formulate the final lemma giving the remaining contribution to the principal singular part at angular momentum $n = 0$: 
\begin{lemma}\label{lem:alphabetavonetwoeffectonsingularn=0} There exists a function $H(\tau, R)$ coinciding with 
	\[
	\mathcal{D}_0\left(\lambda^{-2}\left[c'\left(2\epsilon_X + 2t\epsilon_{tX} +  2(t-r)\epsilon_{XX}\right) - ih'(t)\cos U\epsilon_X\right]\right)
	\]
	near the light cone $R=\nu\tau$, and setting
	\[
	\xb(\tau, \xi): = \int_{\tau_0}^{\tau}U^{(0)}\left(\tau, \sigma,\xi\right)\cdot\mathcal{F}^{(0)}\left(H(\sigma, \cdot)\right)\left(\frac{\lambda^2(\tau)}{\lambda^2(\sigma)}\xi\right)\,d\sigma, 
	\]
	we have (where as usual $ \chi_{\xi<M}$ is a smooth cutoff, and by abuse of notation, we interpret $c, h$ as functions of the re-scaled variable $\tau$)
\begin{equation}\label{eq:c_0sourcetermobstruction}\begin{split}
		&\int_0^\infty \chi_{\lambda^2(\tau)\xi<M}\xb(\tau, \xi)\tilde{\rho}_{0}(\xi)\,d\xi\\& = \sum_{\pm}\sum_{k=1,2,3}\sum_{0\leq r\leq 2}\sum_{0\leq j_1\leq r}\int_0^\infty \chi_{\lambda^2(\tau)\xi<M}\frac{e^{\pm i\nu\tau\xi^{\frac12}}}{\xi^{\frac12+\frac{k\nu}{2}}}\left[\log\left(\xi\lambda^2(\tau)\right)\right]^{j_1}\tilde{\rho}_{0}^{\frac12}(\xi)\,d\xi\\&\hspace{6cm}\cdot \sum_{0\leq j_2\leq r-j_1} \int_{\tau_0}^{\tau}\left(\log\sigma\right)^{j_2}\left[\frac{\lambda(\sigma)}{\lambda(\tau)}\right]^{\nu k}\nu_{\pm}^{(k,r,j_1,j_2)}(\sigma)\,d\sigma\\
		& + \tilde{\tilde{c}}_0(\tau), 
	\end{split}\end{equation}
	where the error satisfies under the assumption $\big\|\tau^{2+}\cdot(c'(\tau)-ih'(\tau))\big\|_{L^\infty}\lesssim (\tau_0^{-1}+\Lambda)\cdot \Lambda$ the bound
	\[
	\left|\tilde{\tilde{c}}_0(\tau)\right|\lesssim \tau^{-4-}\cdot  \Lambda\cdot(\tau_0^{-1}+\Lambda),
	\]
	and the functions $\nu_{\pm}^{(k,j_1,j_2)}(\sigma)$ are explicitly given by 
	\begin{align*}
	\nu_{\pm}^{(k,r,j_1,j_2)}(\sigma) = \sigma^{-(k-1)\nu}\cdot\lambda^{-\nu}(\sigma)\cdot \gamma_{k,r-j_1-j_2}^{\pm}\cdot \left(o_{k,r,j_{1,2}}+O(\sigma^{-1})\right)\cdot \left(c'(\sigma)-ih'(\sigma)\right),
	\end{align*}
	where the coefficients $o_{k,r,j_{1,2}}$ and error terms $O(\sigma^{-1})$ are real-valued. 
\end{lemma}
Then we can formulate the analogue of Prop.~\ref{prop:forcedvanishingatn=plusminus1} for the angular momentum $n = 0$ case, which is proved analogously to Proposition \ref{prop:forcedvanishingatn=plusminus1}:
\begin{proposition}\label{prop:forcedvanishingatn=0}  Let the coefficient functions $\beta_{\pm}^{k,j_1,j_2}$  be as in Lemma~\ref{lem:n=0maintermafterparametrix} with $N_1 = 2$,  $k = 1, 2, 3$,satisfying the bound 
\[
\sum_{k,j_1,j_2}\left|\beta_{\pm}^{k,j_1,j_2}(\sigma)\right|\lesssim \sigma^{-4-}\cdot (\tau_0^{-1}+\Lambda)\cdot \Lambda. 
\]
Then there is a $\mathbb{R}$-linear map $\tilde{S}: \{\beta_{\pm}^{k,j_1,j_2}\}\longrightarrow \left(C_0\left([\tau_0, 2\tau_0];\mathbb{R}\right)\right)^2$ such that if we set 
\[
\left(c', h'\right) = \tilde{S}\left(\{\beta_{\pm}^{k,j_1,j_2}\}\right),
\]
then the vanishing relations 
\[
\sum_{j_2\leq 2-j_1}\int_{\tau_0}^{\infty}\left(\log\sigma\right)^{j_2}\left[\lambda(\sigma)\right]^{\nu k}\left[\beta_{\pm}^{(k,j_1,j_2)}(\sigma) + \sum_{j_1\leq r\leq 2}\nu_{\pm}^{k,r,j_1,j_2}(\sigma)\right]\,d\sigma = 0
\]
hold for each $k = 1, 2, 3,, j_1=0, 1, 2$. 
Moreover we have the bound 
\begin{align*}
\left\|\sigma^{2+}c'(\sigma)\right\|_{L^\infty} + \left\|\sigma^{2+}h'(\sigma)\right\|_{L^\infty}\lesssim (\tau_0^{-1}+\Lambda)\cdot \Lambda.
\end{align*}

\end{proposition}
The proof is again the same as for Lemma~\ref{lem:alphamodulationconditionssolved}.
\\

Taking advantage of  Prop.~\ref{prop:forcedvanishingatn=plusminus1} and Prop.~\ref{prop:forcedvanishingatn=0}, we can finally formulate the analogue of the key result Prop.~\ref{prop:ngeq2finalsourcetermestimatesingoodspaces} in the context of the exceptional modes $n = 0,\pm 1$. 
\begin{proposition}\label{prop:nless23finalsourcetermestimatesingoodspaceswithvanishingcondition} Let $\Lambda\ll 1$ be defined as in \eqref{eq:LambdaDef}, and $\tau_0\gg 1$. Then there exist modulation parameters 
	\[
	\left(\alpha'(\tau), \beta'(\tau), v_1'(\tau), v_2'(\tau), c'(\tau), h'(\tau)\right)\in C_0^\infty\left([\tau_0, 2\tau_0]; \mathbb{R}\right)
	\]
	satisfying the bounds stated in  Prop.~\ref{prop:forcedvanishingatn=0}, Prop.~\ref{prop:forcedvanishingatn=plusminus1}, and such that the following conclusion holds: there exist functions $\yb^{(n)}$, $n = 0,\,\pm 1$, such that setting 
	\[
	\psi^{+}(n) = \int_0^\infty \yb^{(n)}(\tau,\xi)\phi_{n}(R,\xi)\tilde{\rho}_{n}(\xi)\,d\xi,\,n = 0,\,\pm1,
	\]
	and letting $F^{\pm}(n)$ be the sum of all the source terms in Prop.~\ref{prop:smoothlinearsource}, Prop.~\ref{prop:bilinwithUregular2}, Prop.~\ref{prop:bilinwithUregular3}, as well as the terms generated by modulating on the above parameters, i.e. \eqref{eq:alphamodulationn=1sourcecontribution}, \eqref{eq:n=1sourcetermduetobetamodulation}, \eqref{eq:n=1contributionfromLorentzmodulationsingular},  \eqref{eq:n=-1contributionfromLorentzmodulationsingular}, \eqref{eq:n=0effectofmodulatingincandh}, and defining $G^{(+)}(n),\,n = 0,\,\pm 1$, in analogy to Prop.~\ref{prop:neq12finalsourcetermestimatesingoodspaces}, Prop.~\ref{prop:neq22finalsourcetermestimatesingoodspaces}, Prop.~\ref{prop:neq-12finalsourcetermestimatesingoodspaces}, the same bounds as in these propositions obtain. Furthermore, introducing 
	\[
	\xb^{(n)}(\tau,\xi): = \int_{\tau_0}^{\tau}U^{(n)}\left(\tau, \sigma,\xi\right)\cdot\mathcal{F}^{(n)}\left(G^{(+)}(n)(\sigma, \cdot)\right)\left(\frac{\lambda^2(\tau)}{\lambda^2(\sigma)}\xi\right)\,d\sigma
	\]
	we have the good source bounds 
	\begin{align*}
		\left|\int_0^\infty \chi_{\lambda^2(\tau)\xi<M}\xb^{(n)}(\tau, \xi)\tilde{\rho}_{n}(\xi)\,d\xi\right|\lesssim \tau^{-4-}\cdot \left(\tau_0^{-1}+\Lambda\right)\cdot \Lambda
	\end{align*}
	provided  we make the somewhat stronger assumption $\left|c_0(\tau)\right| + \tau\cdot \left|c_0'(\tau)\right|\lesssim \Lambda\cdot\tau^{-2-}$. Furthermore, the solution to \eqref{eq:c+evolutioneqn} satisfies the bound 
	\begin{align*}
		\left|c_+(\tau)\right| + \tau\cdot\left|c_+'(\tau)\right|\lesssim \tau^{-2-}\cdot\left(\tau_0^{-1}+\Lambda\right)\cdot \Lambda. 
	\end{align*}
\end{proposition}
We recall that the evolution equation for $c_+$ is exceptional due to the presence of a non-decaying mode for the operator $\partial_{\tau}^2 + \frac{\lambda_{\tau}}{\lambda}\partial_{\tau}$, and that we force vanishing toward $\tau = +\infty$ for $c_+$ by choosing $\alpha_{\infty}$ suitably, see the discussion in {\it{(iii)}} of subsection~\ref{subsubsec:alphabetamodulation}. This can be done independently of the vanishing conditions `neutralizing' the leading order singular part on the light cone according to Prop.~\ref{prop:forcedvanishingatn=plusminus1}. 
\subsection{Secondary modulations}\label{subsec:SecondarymodDetails}

In the previous subsection, we use modulation to eliminate the top order singular terms on the light cone. As we have seen, after modulation, the source terms in the ODE for $c_{0}(\tau), c_{\pm 1}(\tau)$ decay like $\tau^{-4+}$, which implies that $c_{0}(\tau), c_{\pm1}(\tau)$ decay $\tau^{-2+}$. On the other hand, $c_{0}$ contributes quadratic source to the ODE satisfied by $c_{i}, i=0,\pm1$, and  the $\tau^{-2+}$ decay for $c_0$ would not be enough to close things, see e.g. the statement of Prop.~\ref{prop:fourierboundsforn=-1termlocalizedaway}. Therefore we need a secondary modulation to eliminate this kind of quadratic source. Based on \eqref{eq:modulatedphi}, we consider the following ansatz:
\begin{align}\label{eq:precisemodulation}
	\Psi = \mathcal{L}_{v(t)}\mathcal{R}_{h(t)}^{\alpha(t),\beta(t)}\mathcal{S}_{c(t)}\left(\Phi^{(\text{mod}_2(t))} + \varphi_1E_1 + \varphi_2E_2 + a(\Pi_{\Phi^{\perp}}\varphi)\right).
\end{align}
Here $v(t), h(t),\alpha(t), \beta(t), c(t)$ are as in the previous subsection, and in particular they are chosen to be constant after a finite time interval whose length is chosen such that the perturbation part $\varphi$ has not yet interacted with the shock on the light cone. The secondary modulation parameters act on $\Phi$, resulting in $\Phi^{(\text{mod}_2(t))}$, and are given explicitly by 
\begin{equation}\label{eq:secondarymod}
	\Phi^{(\text{mod}_2(t))} =\mathcal{R}_{\bar{h}(t)}^{0, 0}\mathcal{S}_{\bar{c}(t)}\big(\Phi_{\text{smooth}}\big) + \Phi - \Phi_{\text{smooth}},
\end{equation}
where we carefully observe that the Lorentz modulations are omitted. We shall pick these secondary modulation parameters in such a way that they vanish towards $\tau = +\infty$. Furthermore, we let $\Phi_{\text{smooth}}$ be defined as a partially truncated version of $\Phi$, namely 
\begin{equation}\label{eqLPhismoothcone}
	\Phi_{\text{smooth}} = \left(\begin{array}{c}\sin U_{\text{smooth}} \cos\theta\\ \sin U_{\text{smooth}} \sin \theta\\ \cos U_{\text{smooth}}\end{array}\right),
\end{equation}
and where 
\begin{align}\label{Usmoothcone}
	U_{\text{smooth}} = Q(\lambda(t)r) + \chi_{R<\delta\tau}\epsilon(t,r),.
\end{align}
Here $ \chi_{R<\delta\tau}$ smoothly localizes to the indicated region, and we shall pick $\delta\ll \nu$. The point here is that the secondary modulation parameters do not act on the singular portion of $\epsilon$ near the light cone. Moreover, while the primary modulation parameters also act on the bulk part (away from the light cone), their effect on the evolution of the instabilities via this action on the bulk part will be decaying very fast toward $\tau = +\infty$, while the main effect comes from the action of the secondary modulation parameters on the bulk part. Also, carefully observe that while the primary modulation parameters act on the perturbation $\varphi_{1,2}$, which in principle leads to a derivative loss, this is not the case for the secondary modulation parameters. This derivative loss is harmless since it only concerns a finite time interval, where $\varphi_{1,2}$ are very smooth $H^{C_1}$ since the perturbation has not yet interacted with the shock. 
\\

We now analyze the effect of the secondary modulation: this will reveal that the parameters $\bar{c}(t), \bar{h}(t)$ will have their primary effect via the modification of $\lim_{R\rightarrow 0}R^{-j(n)}F_+(n)$, $n = 0, 1$, see \eqref{eq:c+evolutioneqn}, \eqref{eq:c0evolutioneqn}, \eqref{eq:c-1evolutioneqn}, while the more indirect effect via inhomogeneous parametrix is weaker. 
\\
\subsubsection{Secondary scaling modulation in $\bar{c}$.}

From the bulk term $Q(\lambda(t)r)$, we generate the error terms(which we re-scale right away for the equation in terms of $\tau, R$)
\begin{equation}\label{eq:Qc2modulationterms}
	\lambda^{-2}\bar{c}''(t)RQ'(R) + 2\lambda^{-2}\bar{c}'(t)\lambda'(t)(R\partial_R)^2Q(R) +  \lambda^{-2}(\bar{c}'(t))^2(R\partial_R)^2Q(R).
\end{equation}
The effect from modulating on the correction term $\varepsilon$ is more complicated due to the 'gluing procedure' connecting the inner solution in the regime $R<\delta\tau$ to the outer region. First, as before we encounter the terms
\begin{equation}\label{eq:Qc2modulationtermsepsilon0}
	\chi\lambda^{-2}\bar{c}''(t)R\epsilon'(R) + 2\chi\lambda^{-2}\bar{c}'(t)\lambda'(t)(R\partial_R)^2\epsilon(R) +  \lambda^{-2}\chi(\bar{c}'(t))^2(R\partial_R)^2\epsilon(R).
\end{equation}
Next, from the 'gluing procedure', we generate the additional error terms 
\begin{equation}\label{eq:Qc2modulationtermsepsilon1}\begin{split}
		&\frac{1}{R}\partial_R\left(\chi_{R<\delta\tau}\right)\left(\mathcal{S}_{\bar{c}} - I\right)\epsilon,\,\partial_R^2\left(\chi_{R<\delta\tau}\right)\left(\mathcal{S}_{\bar{c}} - I\right)\epsilon,\,\lambda^{-2}\bar{c}'\partial_t\left(\chi_{R<\delta\tau}\right)R\epsilon'(R),\\
		&\frac{\sin\left(2\left(Q_{\bar{c}} + \chi\epsilon_{\bar{c}} + (1-\chi)\epsilon\right)\right)}{2R^2} - \chi\frac{\sin\left(2\left(Q_{\bar{c}} + \epsilon_{\bar{c}}\right)\right)}{2R^2} - (1-\chi)\frac{\sin\left(2\left(Q_{\bar{c}}+\epsilon\right)\right)}{2R^2}\\
		& - \lambda^{-2}(1-\chi)\Box^{(0)}\left(Q_{\bar{c}} - Q\right) - (1-\chi)\left[\frac{\sin\left(2\left(Q+\epsilon\right)\right)}{2R^2} - \frac{\sin\left(2\left(Q_{\bar{c}}+\epsilon\right)\right)}{2R^2}\right],\\
\end{split}
\end{equation}
where $\Box^{(0)}$ is the wave operator with respect to the original variables $t, r$ acting on the $n = 0$ mode. There is a subtlety here which comes from the logarithmic degeneracy both in $\epsilon$ as well as in the distorted Fourier basis. 
\\

To begin with, we note that only the terms \eqref{eq:Qc2modulationterms} contribute directly via $\lim_{R\rightarrow 0}R^{-1}F_+(0)$. We shall show that this contribution is dominant, the other ones being perturbative. 
To see this, consider first the delicate term
\begin{align*}
	\frac{1}{R}\partial_R\left(\chi_{R<\delta\tau}\right)\left(\mathcal{S}_{\bar{c}} - I\right)\epsilon.
\end{align*}
Note that we can replace $\varepsilon$ by its leading order term $\log R\cdot\frac{R}{(\lambda t)^2}$. We need to analyze the contribution of this term to the expression 
\begin{align*}
	\int_0^\infty x(\tau,\xi)\tilde{\rho}_0(\xi)\,d\xi,
\end{align*}
where $x(\tau,\xi)$ is the Duhamel propagator of 
\[
\mathcal{D}\left(\frac{1}{R}\partial_R\left(\chi_{R<\delta\tau}\right)\left(\mathcal{S}_{\bar{c}} - I\right)\epsilon\right).
\]
Since this term is smooth, integration by parts with respect to $\sigma$ in the Duhamel propagator reveals that we may assume that $\xi^{\frac12}\lesssim \tau^{-1}$ for the output frequency. Then observe that 
\begin{align*}
	&\left\langle \mathcal{D}\left(\frac{1}{R}\partial_R\left(\chi_{R<\delta\tau}\right)\left(\mathcal{S}_{\bar{c}} - I\right)\epsilon\right), \phi(R;\xi)\right\rangle_{L^2_{R\,dR}}\\
	& = -\left\langle \left(\frac{1}{R}\partial_R\left(\chi_{R<\delta\tau}\right)\left(\mathcal{S}_{\bar{c}} - I\right)\epsilon\right), \left(\partial_R +O\left(R^{-3}\right)\right)\phi(R;\xi)\right\rangle_{L^2_{R\,dR}}\\
\end{align*}
Here observe that 
\[
\left(\partial_R+O\left(R^{-3}\right)\right)\phi(R;\xi) = R^{-1} + O(\log R\cdot R\xi).
\]
For the contribution of the first term on the right we have gained a $\log R$, which means the corresponding contribution gains a $\log^{-1}\tau$. For the second term on the right since 
\[
R^2\xi\lesssim \delta^2, 
\]
we get 
\begin{align*}
	&\left \langle \left(\frac{1}{R}\partial_R\left(\chi_{R<\delta\tau}\right)\left(\mathcal{S}_{\bar{c}} - I\right)\epsilon\right), O(\log R\cdot R\xi)\right\rangle_{L^2_{R\,dR}}\\
	&\sim \left\langle \tilde{\chi}_{R\sim\delta\tau}R^{-2}\cdot \log R\cdot \frac{\bar{c}}{\tau^2}, \log R\cdot \delta^2\right\rangle_{L^2_{R\,dR}}.
\end{align*}
To see how this affects the contribution to the source term 
\[
\int_0^\infty x(\tau,\xi)\tilde{\rho}_0(\xi)\,d\xi,
\]
we need a lemma according to the following lines:
\begin{lemma}\label{lem:c2dispsource} Let 
	\begin{align*}
		U(\tau)f(\cdot, \xi) = \int_{\tau_0}^{\tau}\frac{\sin\left(\lambda(\tau)\xi^{\frac12}\int_{\sigma}^{\tau}\lambda^{-1}(s)\,ds\right)}{\xi^{\frac12}}f\left(\sigma,\frac{\lambda^2(\tau)}{\lambda^2(\sigma)}\xi\right)\,d\sigma, 
	\end{align*}
	the inhomogeneous Schrodinger propagator corresponding to $n = 0$. Also assume $f(\tau_0,\cdot) = 0$ as well as $\int_0^\infty \xi^{-1}f(\tau,\xi)\tilde{\rho}_0(\xi)\,d\xi = 0$ for all $\tau>\tau_0$. Then for $0<p<p_*(\nu)$ we have the bound 
	\begin{align*}
		\tau^p\cdot \left|\int_0^\infty  U(\tau)f(\cdot, \xi)\tilde{\rho}_0(\xi)\,d\xi\right|&\lesssim\left\|\frac{\tau^p}{\log^2\tau}\cdot f(\tau,\cdot)\right\|_{L_\tau^\infty l^1L^\infty_{d\xi}} + \left\| \frac{\tau^p}{\log^2\tau}\cdot(\tau\partial_{\tau})f(\tau,\cdot)\right\|_{L_\tau^\infty l^1L^\infty_{d\xi}}\\
		& + \left\|\frac{\tau^p}{\log^2\tau}\cdot (\xi\partial_{\xi})f(\tau,\cdot)\right\|_{L_\tau^\infty l^1L^\infty_{d\xi}}
	\end{align*}
	Here we use the norm 
	\begin{align*}
		\left\|f(\cdot)\right\|_{l^1 L^\infty_{d\xi}}: = \sum_{\lambda\in \Z}\sum_{\xi\sim 2^\lambda}\left\|f\right\|_{L^\infty(\xi\sim 2^{\lambda})}.
	\end{align*}
	Furthermore, for given $\delta_*>0$, we also have the bound 
	\begin{align*}
		\tau^p\cdot \left|\int_0^\infty  U(\tau)f(\cdot, \xi)\tilde{\rho}_0(\xi)\,d\xi\right|&\lesssim\delta_*^{-1}\left\|\frac{\tau^p}{\log^2\tau}\cdot f(\tau,\cdot)\right\|_{L_\tau^\infty l^1L^\infty_{d\xi}(\xi<\delta_*^{-2}\tau^{-2})}\\& +\delta_*\left\| \frac{\tau^p}{\log^2\tau}\cdot(\tau\partial_{\tau})^2f(\tau,\cdot)\right\|_{L_\tau^\infty l^1L^\infty_{d\xi}(\xi>\delta_*^{-2}\tau^{-2})}\\
		& + \delta_*\left\|\frac{\tau^p}{\log^2\tau}\cdot (\xi\partial_{\xi})^2f(\tau,\cdot)\right\|_{L_\tau^\infty l^1L^\infty_{d\xi}(\xi>\delta_*^{-2}\tau^{-2})}
	\end{align*}
	
\end{lemma}
\begin{proof} {\it{First estimate}}. We decompose the integral
	\begin{align*}
		\int_0^\infty  U(\tau)f(\cdot, \xi)\tilde{\rho}_0(\xi)\,d\xi = \sum_{j=1,2}\int_0^\infty  U_j(\tau)f(\cdot, \xi)\tilde{\rho}_0(\xi)\,d\xi,
	\end{align*}
	where we define 
	\begin{align*}
		U_1(\tau)f(\cdot, \xi) = \int_{\tau_0}^{\tau}\chi_{\sigma\frac{\lambda(\tau)}{\lambda(\sigma)}\xi^{\frac12}\lesssim 1}\frac{\sin\left(\lambda(\tau)\xi^{\frac12}\int_{\sigma}^{\tau}\lambda^{-1}(s)\,ds\right)}{\xi^{\frac12}}f\left(\sigma,\frac{\lambda^2(\tau)}{\lambda^2(\sigma)}\xi\right)\,d\sigma
	\end{align*}
	as well as 
	\[
	U_2(\tau)f(\cdot, \xi) = U(\tau)f(\cdot, \xi) - U_1(\tau)f(\cdot, \xi).
	\]
	{\it{Estimate for $U_1f$.}} We can bound 
	\begin{align*}
		\left|\int_0^\infty  U_1(\tau)f(\cdot, \xi)\tilde{\rho}_0(\xi)\,d\xi\right|\lesssim  \int_{\tau_0}^{\tau}\log^{-2}\sigma\cdot \left(\sigma\frac{\lambda(\tau)}{\lambda(\sigma)}\right)^{-1}\cdot \left\|f(\sigma,\cdot)\right\|_{L^\infty_{d\xi}}\,d\sigma. 
	\end{align*}
	It follows that for $p<p_*(\nu)$ we have 
	\begin{align*}
		\tau^p\cdot \left|\int_0^\infty  U_1(\tau)f(\cdot, \xi)\tilde{\rho}_0(\xi)\,d\xi\right|\lesssim \left\|\frac{\tau^p}{\log^2\tau}\cdot f(\tau,\cdot)\right\|_{L_\tau^\infty L^\infty_{d\xi}}
	\end{align*}
	which is better than what we need. 
	\\
	
	{\it{Estimate for $U_2f$.}} Here we perform integration by parts with respect to $\sigma$. This results in the boundary term 
	\begin{align*}
		B_2(\tau): = \int_0^\infty \chi_{\xi^{\frac12}\tau\gtrsim 1}\cdot \xi^{-1}f(\tau,\xi)\cdot \tilde{\rho}_0(\xi)\,d\xi,
	\end{align*}
	which according to our assumption also equals 
	\begin{align*}
		B_2(\tau) = -\int_0^\infty \chi_{\xi^{\frac12}\tau\lesssim 1}\cdot \xi^{-1}f(\tau,\xi)\cdot \tilde{\rho}_0(\xi)\,d\xi,
	\end{align*}
	But then we have 
	\begin{align*}
		\left|\tau^p\cdot B_2(\tau)\right|\lesssim \left\|\frac{\tau^p}{\log^2\tau}\cdot f(\tau,\cdot)\right\|_{L_\tau^\infty l^1L^\infty_{d\xi}}.  
	\end{align*}
	We next arrive at the double integral 
	\begin{align*}
		\int_0^\infty  \tilde{U}_2(\tau)f(\tau, \xi)\tilde{\rho}_0(\xi)\,d\xi,
	\end{align*}
	where we set 
	\begin{align*}
		\tilde{U}_2(\tau)f(\tau, \xi): =  \int_{\tau_0}^{\tau}\frac{\lambda(\sigma)}{\lambda(\tau)}\cdot \frac{\cos\left(\lambda(\tau)\xi^{\frac12}\int_{\sigma}^{\tau}\lambda^{-1}(s)\,ds\right)}{\xi}\partial_{\sigma}\left(\chi_{\sigma\frac{\lambda(\tau)}{\lambda(\sigma)}\xi^{\frac12}\lesssim 1}f\left(\sigma,\frac{\lambda^2(\tau)}{\lambda^2(\sigma)}\xi\right)\right)\,d\sigma
	\end{align*}
	This can then be bounded by 
	\begin{align*}
		\left|\tau^p\cdot \int_0^\infty  \tilde{U}_2(\tau)f(\tau, \xi)\tilde{\rho}_0(\xi)\,d\xi\right|\lesssim \left[\left\|\frac{\tau^p}{\log^2\tau}\cdot f(\tau,\cdot)\right\|_{L_\tau^\infty l^1L^\infty_{d\xi}} + \left\| \frac{\tau^p}{\log^2\tau}\cdot(\tau\partial_{\tau})f(\tau,\cdot)\right\|_{L_\tau^\infty l^1L^\infty_{d\xi}}\right]
	\end{align*}
	provided $p<p_*(\nu)$. 
	\\
	
	{\it{Second estimate}}. Here one uses the cutoff $\chi_{\sigma\frac{\lambda(\tau)}{\lambda(\sigma)}\xi^{\frac12}<\delta_*^{-1}}$ instead and performs one more integration by parts. 
	
\end{proof}
To apply the preceding lemma in the context of the error terms \eqref{eq:Qc2modulationtermsepsilon1}\, we state 
\begin{lemma}\label{lem:Qc2varepstermsbounds} Denoting the first three terms in \eqref{eq:Qc2modulationtermsepsilon1} as $E_j,\,j = 1, 2, 3$, and letting 
	\begin{align*}
		\left\|f\right\|_{N_{\delta_*}}
	\end{align*}
	denote the norm on the right at the end of the preceding lemma, we have the bounds ($j = 1, 2, 3$)
	\begin{align*}
		\left\|\mathcal{F}^{(0)}\left(\mathcal{D}(N_j)\right)\right\|_{N_{\sqrt{\delta}}}\lesssim \delta^{\frac14}\cdot \left(\sup_{\tau>\tau_0}\left\|\tau^p\cdot \frac{\bar{c}(\tau)-1}{\tau^2}\right\|_{L^\infty_{d\tau}} + \sup_{\tau>\tau_0}\left\|\tau^p\cdot \frac{\bar{c}'(\tau)}{\tau}\right\|_{L^\infty_{d\tau}}\right).
	\end{align*}
\end{lemma}
\begin{proof} These terms are all similar, and the first of them is representative. Using the considerations before the last lemma, we reduce to bounding the three expressions at the end of Lemma~\ref{lem:c2dispsource} with 
	\[
	f(\tau,\xi) = \left\langle \left(\frac{1}{R}\partial_R\left(\chi_{R<\delta\tau}\right)\left(\mathcal{S}_{\bar{c}} - I\right)\epsilon\right), O(\log R\cdot R\xi)\right\rangle_{L^2_{R\,dR}}
	\]
	Using $\delta_* = \delta^{\frac12}$, and recalling 
	\begin{equation}\label{eq:varepsasymptotics}
		\epsilon(\tau, R) = \log R\cdot \frac{R}{(\lambda t)^2} + O\left(\frac{R}{(\lambda t)^2}\right), 
	\end{equation}
	we get 
	\begin{align*}
		\delta_*^{-1}\left\|\frac{\tau^p}{\log^2\tau}\cdot f(\tau,\cdot)\right\|_{L_\tau^\infty l^1L^\infty_{d\xi}(\xi<\delta_*^{-2}\tau^{-2})}\lesssim \delta^{\frac12-}\cdot \sup_{\tau>\tau_0}\left\|\tau^p\cdot \frac{\bar{c}(\tau)-1}{\tau^2}\right\|_{L^\infty_{d\tau}}, 
	\end{align*}
	which is as desired. To bound the expressions 
	\begin{align*}
		\delta_*\left\| \frac{\tau^p}{\log^2\tau}\cdot(\tau\partial_{\tau})^2f(\tau,\cdot)\right\|_{L_\tau^\infty l^1L^\infty_{d\xi}(\xi>\delta_*^{-2}\tau^{-2})},\, \delta_*\left\|\frac{\tau^p}{\log^2\tau}\cdot (\xi\partial_{\xi})^2f(\tau,\cdot)\right\|_{L_\tau^\infty l^1L^\infty_{d\xi}(\xi>\delta_*^{-2}\tau^{-2})},
	\end{align*}
	we proceed as for the previous term when $R^2\xi\lesssim 1$, while we use the asymptotic
	\[
	\left|\phi(R;\xi)\right|\lesssim |\log\xi|\cdot (R\xi^{\frac12})^{-\frac12}
	\]
	in the regime $R^2\xi\gtrsim 1$. Also, observe that 
	\begin{align*}
		(\xi\partial_\xi)\phi(R;\xi) = (R\partial_R)\phi(R;\xi) + O\left(\log R(R\xi^{\frac12})^{-1}\right),
	\end{align*}
	which allows us to reduce the term $(\xi\partial_{\xi})f$ to one without the dilation operator via integration by parts with respect to $R$. 
\end{proof}

It remains to deal with the terms in the second and third line of \eqref{eq:Qc2modulationtermsepsilon1}. Commence with the terms on the third line. We expand 
\begin{align*}
	(1-\chi)\left[\frac{\sin\left(2(Q+\epsilon)\right)}{2R^2} - \frac{\sin\left(2(Q_{\bar{c}}+\epsilon)\right)}{2R^2}\right] &= (1-\chi)\cos\left(2\epsilon\right)\cdot \left(\frac{\sin(2Q)}{2R^2} - \frac{\sin(2Q_{\bar{c}})}{2R^2}\right)\\
	& + (1-\chi)\sin\left(2\epsilon\right)\cdot \left(\frac{\cos(2Q)}{2R^2} - \frac{\cos(2Q_{\bar{c}})}{2R^2}\right)\\
\end{align*}
The first term on the right partially cancels against the first term on the third line in \eqref{eq:Qc2modulationtermsepsilon1}. In fact, write 
\begin{align*}
	&N_4: = (1-\chi)\left(\partial_{RR} + \frac{1}{R}\partial_R\right)\left(Q - Q_{\bar{c}}\right) - (1-\chi)\cos\left(2\epsilon\right)\cdot \left(\frac{\sin(2Q)}{2R^2} - \frac{\sin(2Q_{\bar{c}})}{2R^2}\right),\\
	&N_5: = (1-\chi)\sin\left(2\varepsilon\right)\cdot \left(\frac{\cos(2Q)}{2R^2} - \frac{\cos(2Q_{\bar{c}})}{2R^2}\right)\\
	&N_6: =  \lambda^{-2}(1-\chi)(\lambda\partial_{\tau})^2\left(Q - Q_{\bar{c}}\right).
\end{align*}
Then we have 
\begin{lemma}\label{lem:N45} We have the bounds (for $j = 4, 5$)
	\begin{align*}
		&\tau^{p+}\cdot \left\|\mathcal{F}^{(0)}\big(\mathcal{D}(N_j)\big)\right\|_{S_1^0\cap L^\infty_{d\xi}} + \tau^{p+}\cdot \left\|(\tau\partial_{\tau})\mathcal{F}^{(0)}\left(\mathcal{D}(N_j)\right)\right\|_{S_1^0\cap L^\infty_{d\xi}} \\&\lesssim \sup_{\tau>\tau_0}\left\|\tau^p\cdot \frac{\bar{c}(\tau)-1}{\tau^2}\right\|_{L^\infty_{d\tau}} + \sup_{\tau>\tau_0}\left\|\tau^p\cdot \frac{\bar{c}'(\tau)}{\tau}\right\|_{L^\infty_{d\tau}} + \sup_{\tau>\tau_0}\left\|\tau^p\cdot \bar{c}''(\tau)\right\|_{L^\infty_{d\tau}}
	\end{align*}
	In particular, given $\delta>0$ there is $\tau_0 = \tau_0(\delta)$ such that for $\tau\geq \tau_0$ we have 
	\begin{align*}
		&\tau^{p+}\cdot \left\|\mathcal{F}^{(0)}\left(\mathcal{D}(N_j)\right)\right\|_{S_1^0\cap L^\infty_{d\xi}}+ \tau^{p+}\cdot \left\|(\tau\partial_{\tau})\mathcal{F}^{(0)}\left(\mathcal{D}(N_j)\right)\right\|_{S_1^0\cap L^\infty_{d\xi}}\\&\lesssim \delta\cdot\left(\sup_{\tau>\tau_0}\left\|\tau^p\cdot \frac{\bar{c}(\tau)-1}{\tau^2}\right\|_{L^\infty_{d\tau}} + \sup_{\tau>\tau_0}\left\|\tau^p\cdot \frac{\bar{c}'(\tau)}{\tau}\right\|_{L^\infty_{d\tau}}+ \sup_{\tau>\tau_0}\left\|\tau^p\cdot \bar{c}''(\tau)\right\|_{L^\infty_{d\tau}}
		\right). 
	\end{align*}
	The same bound but without the term involving $\tau\partial_{\tau}$ and with replacing $\tau^{p+}$ by $\frac{\tau^p}{\log\tau}$ applies to $N_6$. Furthermore, 
	\[
	\mathcal{F}^{(0)}\left(\mathcal{D} N_6\right)(\tau, \xi)
	\]
	decays rapidly beyond $\xi_* := (\delta\tau)^{-2}$. 
\end{lemma}
\begin{proof} To begin with, observe that 
	\begin{align*}
		N_4 = (1-\chi)\cdot O(\varepsilon^2)\cdot \frac{\left(Q - Q_{\bar{c}}\right)}{R^2}.
	\end{align*}
	Using \eqref{eq:varepsasymptotics} it is straightforward to check that 
	\begin{align*}
		\tau^{p}\cdot \left\|(\tau\partial_{\tau})^{\kappa}\mathcal{F}^{(0)}\left(\mathcal{D}(N_4)\right)\right\|_{S_1^0\cap L^\infty_{d\xi}}\lesssim \tau^{-(2-)}\cdot \sup_{\tau>\tau_0}\left\|\tau^p\cdot \frac{\bar{c}(\tau)-1}{\tau^2}\right\|_{L^\infty_{d\tau}},\quad \kappa = 0,\,1, 
	\end{align*}
	which gives the desired bound for this term with $p+ = p+2-$. 
	\\
	The term $N_5$ is handled similarly by using 
	\begin{align*}
		\frac{\cos(2Q)}{2R^2} - \frac{\cos(2Q_{\bar{c}})}{2R^2} = \frac{\sin(Q - Q_{\bar{c}})\cdot \sin(Q + Q_{\bar{c}})}{R^2}. 
	\end{align*}
	Next we consider the term $N_6$. Here while the bulk term 
	\[
	Q - Q_{\bar{c}}
	\]
	does not lose an extra logarithm in $\tau$, we do lose one factor $\log R$ from the low frequency asymptotics of $\phi(R;\xi)$. Specifically, note that 
	\begin{align*}
		\mathcal{D}\left((1 - \chi)\left(Q - Q_{\bar{c}}\right)\right) = -\chi'(\bar{c} - 1)\Lambda Q(R) - \chi'\frac{(\bar{c}-1)^2}{2}\Lambda^2Q(R) + O\left((1-\chi)(\bar{c} - 1)R^{-4}\right),
	\end{align*}
	which implies that 
	\begin{align*}
		\left\langle \lambda^{-2}\mathcal{D}\left((1 - \chi)(\lambda\partial_{\tau})^2\left(Q - Q_{\bar{c}}\right)\right), \phi(R;\xi)\right\rangle_{L^2_{R\,dR}} &= \log\xi \cdot \zeta_1(\tau;\xi)\cdot \frac{\bar{c} - 1}{\tau^2} +  \log\xi \cdot \zeta_2(\tau;\xi)\cdot \frac{\bar{c}'}{\tau}\\
		& +  \log\xi \cdot \zeta_3(\tau;\xi)\cdot \bar{c}''(\tau),
	\end{align*}
	where the functions $\zeta_j(\tau;\xi)$ are bounded with symbol type bounds and decay rapidly beyond $\xi_* = (\delta\tau)^{-2}$, as one sees via integration by parts with respect to $R$.
\end{proof}

In order to complete the control of the contribution of $N_6$ to 
\[
\int_0^\infty x(\tau;\xi)\tilde{\rho}_0(\xi)\,d\xi,
\]
we now mention
\begin{lemma}\label{lem:N6duhamel} Assume that the source term $f(\tau,\xi)$ satisfies 
	\[
	\left|f(\tau,\xi)\right|\leq A\cdot \tau^{-p}\log\tau\cdot (1+\delta\tau\xi^{\frac12})^{-N}. 
	\]
	Then provided $p<p_*(\nu)$ we have the bound 
	\begin{align*}
		\tau^p\cdot \left|\int_0^\infty Uf(\tau;\xi)\tilde{\rho}_0(\xi)\,d\xi\right|\lesssim \log^{-1}\tau\cdot A. 
	\end{align*}
	where $Uf$ is defined as in Lemma~\ref{lem:c2dispsource}.
\end{lemma}
\begin{proof}
	Write 
	\begin{align*}
		\int_0^\infty Uf(\tau;\xi)\tilde{\rho}_0(\xi)\,d\xi = \int_{\tau_0}^{\tau}\int_0^\infty \frac{\sin\left(\lambda(\tau)\xi^{\frac12}\int_{\sigma}^{\tau}\lambda^{-1}(s)\,ds\right)}{\xi^{\frac12}}f\left(\sigma,\frac{\lambda^2(\tau)}{\lambda^2(\sigma)}\xi\right)\tilde{\rho}_0(\xi)\,d\xi d\sigma
	\end{align*}
	Then using the low frequency asymptotics $0<\tilde{\rho}_0(\xi)\lesssim \log^{-2}\xi$, $0<\xi\ll 1$, we infer that 
	\begin{align*}
		\left|\int_0^\infty \frac{\sin\left(\lambda(\tau)\xi^{\frac12}\int_{\sigma}^{\tau}\lambda^{-1}(s)\,ds\right)}{\xi^{\frac12}}f\left(\sigma,\frac{\lambda^2(\tau)}{\lambda^2(\sigma)}\xi\right)\tilde{\rho}_0(\xi)\,d\xi \right|\lesssim \log^{-1}\sigma\cdot \sigma^{-1-p}\cdot \frac{\lambda(\sigma)}{\lambda(\tau)}\cdot A,
	\end{align*}
	and we have 
	\begin{align*}
		\tau^p\cdot \int_{\tau_0}^{\tau}\log^{-1}\sigma\cdot \sigma^{-1-p}\cdot \frac{\lambda(\sigma)}{\lambda(\tau)}\cdot A\,d\sigma\lesssim A\cdot \log^{-1}(\tau), 
	\end{align*}
	provided $p<p_*(\nu)$. 
\end{proof}

Turning next to the terms \eqref{eq:Qc2modulationtermsepsilon0}, we note 
\begin{lemma}\label{lem:simplerQc2varepsilonerrors} For each of the three terms $N_j,\,j = 7, 8, 9$ in  \eqref{eq:Qc2modulationtermsepsilon0} we can write 
	\begin{align*}
		\mathcal{F}^{(0)}\left(\mathcal{D} N_j\right) = \log^2\tau\cdot \zeta_j^{(1)}(\tau;\xi)\cdot \bar{c}''(\tau) + \log^2\tau\cdot \zeta_j^{(2)}(\tau;\xi)\cdot \frac{\bar{c}'(\tau)}{\tau},
	\end{align*}
	where the functions $\zeta_j(\tau;\xi)$ are bounded, have symbol behavior, and satisfy the bound
	\[
	\left|\zeta_j^{(\kappa)}(\tau;\xi)\right|\lesssim_{\nu} \delta^2
	\]
\end{lemma}
\begin{proof} This follows from \eqref{eq:varepsasymptotics} as well as the fact that $R\leq \delta\tau$ on the support of the expressions. 
\end{proof}

Arguing as for the proof of Lemma~\ref{lem:N6duhamel}, we now infer 
\begin{lemma}\label{lem:N789} Letting $f_j = \mathcal{F}^{(0)}\left(\mathcal{D} N_j\right)$, $j = 7, 8, 9$, we have the bound 
	\begin{align*}
		\left|\tau^p\cdot \int_0^\infty U(\tau)f_j(\cdot,\xi)\tilde{\rho}_0(\xi)\,d\xi\right|\lesssim \delta^2\cdot \left( \sup_{\tau>\tau_0}\left\|\tau^p\cdot \frac{\bar{c}'(\tau)}{\tau}\right\|_{L^\infty_{d\tau}}+ \sup_{\tau>\tau_0}\left\|\tau^p\cdot \bar{c}''(\tau)\right\|_{L^\infty_{d\tau}}\right). 
	\end{align*}
\end{lemma}

It remains to deal with the terms \eqref{eq:Qc2modulationterms}. Here the key is that applying $\mathcal{D}$ results in much better spatial decay (and elimination of the first of the terms), and in particular no logarithmic growth factor. In fact, we observe that letting $N_j,\,j = 10, 11, 12$ be those terms, we have 
\begin{align*}
	\left|(\tau\partial_{\tau})^{\kappa}\mathcal{F}^{(0)}\left(\mathcal{D}N_j\right)(\tau,\cdot)\right|\lesssim \left|\frac{\bar{c}'(\tau)}{\tau}\right| + \left|\bar{c}''(\tau)\right|,\,\kappa = 0, 1,
\end{align*}
and also 
\begin{align*}
	\left|(\xi\partial_{\xi})\mathcal{F}^{(0)}\left(\mathcal{D}N_j\right)(\tau,\xi)\right|\lesssim \left|\frac{\bar{c}'(\tau)}{\tau}\right| + \left|\bar{c}''(\tau)\right|.
\end{align*}
Now a simple modification of the proof of Lemma~\ref{lem:c2dispsource} allows us to replace the norm $l^1L^\infty_{d\xi}$ there at the expense of one factor $\log\tau$. We thereby infer the 
\begin{lemma}\label{lem:N101112} We have the bound, for $j = 10, 11, 12$ and $f_j = \mathcal{F}^{(0)}\left(\mathcal{D} N_j\right)$,
	\begin{align*}
		\left|\tau^p\cdot \int_0^\infty U(\tau)f_j(\cdot,\xi)\tilde{\rho}_0(\xi)\,d\xi\right|\lesssim \log^{-1}\tau\cdot \left(\sup_{\tau>\tau_0}\tau^p\left|\frac{\bar{c}'(\tau)}{\tau}\right| + \sup_{\tau>\tau_0}\tau^p\left|\bar{c}''(\tau)\right|\right). 
	\end{align*}
\end{lemma}

Combining the Lemma~\ref{lem:N101112}, Lemma~\ref{lem:N789}, Lemma~\ref{lem:N45} and Lemma~\ref{lem:N6duhamel} as well as Lemma~\ref{lem:Qc2varepstermsbounds}, we can finally infer the following 
\begin{lemma}\label{lem:c2secondarycontribution} Given $\gamma>0$, there is $\tau_0 = \tau_0(\gamma)$ and $\delta = \delta(\gamma)$ such that the following hold: Letting $N$ be any one of the terms in \eqref{eq:Qc2modulationterms}, \eqref{eq:Qc2modulationtermsepsilon0}, \eqref{eq:Qc2modulationtermsepsilon1}. Then we have the estimate 
	\begin{align*}
		\tau^p\cdot \left|\int_0^\infty U(\tau)\left(\mathcal{F}^{(0)}(\mathcal{D} N)\right)(\xi)\tilde{\rho}_0(\xi)\,d\xi\right|\leq \gamma\cdot\left(\sup_{\tau>\tau_0}\tau^p\left|\frac{\bar{c}(\tau)-1}{\tau^2}\right|+\sup_{\tau>\tau_0}\tau^p\left|\frac{\bar{c}'(\tau)}{\tau}\right| + \sup_{\tau>\tau_0}\tau^p\left|\bar{c}''(\tau)\right|\right)
	\end{align*}
	provided $\tau\geq\tau_0$. 
\end{lemma}

The preceding lemma gives the first of an infinite number of contributions to 
\[
\int_0^\infty x(\tau;\xi)\tilde{\rho}_0(\xi)\,d\xi
\]
where $x(\tau;\xi)$ is the exact Fourier transform of the wave propagator (at angular momentum $n =0$) of the source terms treated before, the remaining contributions being due to the action of the transference operator. To handle these latter terms, which all turn out to be perturbative as well, we require the following two lemmas:
\begin{lemma}\label{lem:crudeboundc2disperrors} All of the source terms $N_j,\,j = 1, 2,\ldots, 12$ satisfy either the bound 
	\begin{align*}
		\tau^p\cdot \left\|\langle \xi\tau\rangle^A\mathcal{F}^{(0)}\left(\mathcal{D}(N_j)\right)(\tau,\cdot)\right\|_{L^M_{d\xi}}\lesssim_A \log^2\tau\cdot \left(\sup_{\tau>\tau_0}\tau^p\left|\frac{\bar{c}(\tau)-1}{\tau^2}\right|+\sup_{\tau>\tau_0}\tau^p\left|\frac{\bar{c}'(\tau)}{\tau}\right| + \sup_{\tau>\tau_0}\tau^p\left|\bar{c}''(\tau)\right|\right)
	\end{align*}
	for any $A\geq 1$, or else the bounds ($\iota, \kappa = 0, 1, 2$)
	\begin{align*}
		&\tau^p\cdot \left\|(\xi\partial_j)^{\iota}(\tau\partial_{\tau})^{\kappa}\mathcal{F}^{(0)}\left(\mathcal{D}(N_j)\right)(\tau,\cdot)\right\|_{L^M_{d\xi}}\lesssim\cdot \left(\sup_{\tau>\tau_0}\tau^p\left|\frac{\bar{c}(\tau)-1}{\tau^2}\right|+\sup_{\tau>\tau_0}\tau^p\left|\frac{\bar{c}'(\tau)}{\tau}\right| + \sup_{\tau>\tau_0}\tau^p\left|\bar{c}''(\tau)\right|\right)
	\end{align*}
\end{lemma}
The proof follows from the preceding considerations. The first kind of estimate obtains for the source terms $N_j,\,j = 1,\ldots, 9$, which are located in the regime $R\gtrsim\tau$, while the second bound applies to $N_{j}, j = 10, 11, 12$. 

\begin{lemma}\label{lem:transferenceonc2dispersiveerrors} Let $\mathcal{K}_0^{(0)}$ the transference operator at level $n = 0$, and let $U$ as usual the the Duhamel propagator (at level $n = 0$). Then we have for $A\geq j\geq 1$ the relation 
	\begin{align*}
		\xi^{\frac12}\cdot U\left(\frac{\lambda_{\tau}}{\lambda}\mathcal{K}_0^{(0)}\mathcal{D}_{\tau}U\right)^j f = g_j + h_j
	\end{align*}
	where we have the bounds (for some $M\gg 1$ and $C = C(M,\nu,p)$ and $0<\nu<\nu_*(M,p)$)
	\begin{align*}
		&\tau^p\cdot \left\|\langle \xi^{\frac12}\tau\rangle^{A - j} g_j(\tau,\xi)\right\|_{L^{M}_{d\xi}}\leq \frac{C^j}{\log^{(2-)\cdot j}\tau}\cdot \sup_{\sigma>\tau_0}\sigma^p\cdot\left\|\langle \xi^{\frac12}\sigma\rangle^A f(\sigma,\xi)\right\|_{L^M_{d\xi}}\\
		&\tau^{p+}\cdot \left\|h_j(\tau,\xi)\right\|_{L^{M}_{d\xi}}\leq \gamma^j\cdot \sup_{\sigma>\tau_0}\sigma^p\cdot\left\| f(\sigma,\xi)\right\|_{L^M_{d\xi}}\\
	\end{align*}
	provided $0<\gamma$ and $\tau_0\geq \tau_0(\gamma, \nu, p)$. We also have the bound (under the same assumptions on $\nu, \tau_0$)
	\begin{align*}
		\tau^p\cdot\left\|\xi^{\frac12}\cdot U\left(\frac{\lambda_{\tau}}{\lambda}\mathcal{K}_0^{(0)}\mathcal{D}_{\tau}U\right)^j f\right\|_{L^M_{d\xi}}&\leq C_{\gamma}\cdot\gamma^j\cdot \log^{-1}\tau\cdot  \sum_{\iota,\kappa = 0, 1}\sup_{\sigma>\tau_0}\sigma^p\cdot\left\|(\xi\partial_{\xi})^{\iota}(\sigma\partial_\sigma)^{\kappa} f(\sigma,\xi)\right\|_{L^M_{d\xi}}\\
	\end{align*}
	
\end{lemma}
\begin{proof}(sketch) For the first part of the lemma, recalling that $\mathcal{K}_0^{(0)}$ is given in terms of an integral kernel 
	\[
	\frac{F(\xi,\eta)\tilde{\rho}_0(\eta)}{\xi - \eta},
	\]
	define the modified operator $\mathcal{K}_{0, d}^{(0)}$ by inclusion of an extra smooth localizer $\chi_{\xi\sim \eta}$ and further
	\begin{align*}
		\mathcal{K}_{0, nd}^{(0)} = \mathcal{K}_0^{(0)} -\mathcal{K}_{0, d}^{(0)}.
	\end{align*}
	Then we set 
	\begin{align*}
		g_j: = U\left(\frac{\lambda_{\tau}}{\lambda}\mathcal{K}_{0,d}^{(0)}\mathcal{D}_{\tau}U\right)^j f 
	\end{align*}
	Then we note that the factor $\langle \xi\tau\rangle$ can be passed through the expression due to the fact that the operator $\mathcal{K}_{0,d}^{(0)}$ 'essentially diagonalizes' the frequencies, and in fact one may include a factor $\langle\eta^{\frac12}\sigma\rangle$ in front of each $\mathcal{K}_{0, d}^{(0)}$. This essentially localizes these terms to frequency $\eta\lesssim \sigma^{-2}$, and the weight $\tilde{\rho}_0(\eta)\sim \frac{1}{\log^2\eta}$ for $0<\eta\ll 1$ in the kernel for $\mathcal{K}_{0, d}^{(0)}$ as well as the $L^M_{d\xi}$-boundedness of $\mathcal{K}_{0, d}^{(0)}$ and an argument as in \cite{KMiao} implies the first bound. The second bound follows by writing 
	\begin{align*}
		\mathcal{K}_0^{(0)} -\mathcal{K}_{0, d}^{(0)} 
	\end{align*}
	as a sum of compositions of the form 
	\begin{align*}
		U\circ\prod_{l=1}^j T_l
	\end{align*}
	where 
	\[
	T_l = \frac{\lambda_{\tau}}{\lambda}\mathcal{K}_{0,d}^{(0)}\mathcal{D}_{\tau}U,\,T_l = \frac{\lambda_{\tau}}{\lambda}\mathcal{K}_{0,nd}^{(0)}\mathcal{D}_{\tau}U
	\]
	with at least one operator $\mathcal{K}_{0,nd}^{(0)}$, and one gains in time by performing integration by parts with respect to the temporal variable as in the proof of Lemma 12.0.32 in \cite{BKmemo}. 
	\\
	For the last estimate of the lemma, one may again replace all operators $\mathcal{K}_{0}^{(0)}$ by operators $\mathcal{K}_{0,d}^{(0)}$ as one gains a power of $\tau$ otherwise via temporal integration by parts. Then one  performs integration by parts with respect to $\sigma$ in $\mathcal{D}_{\tau}Uf$, which does not result in a boundary term at $\tau = \sigma$. Further integration by parts with respect to the time variable $\sigma$ in the leftmost propagator $U$ allows us to reduce $\xi^{\frac12}$ to the region $\xi^{\frac12}\lesssim \frac{\lambda(\sigma)}{\lambda(\tau)}\sigma^{-1}$ which will propagate to all other spatial frequencies due to the frequency localization, and the gain of $\log^{-1}\tau$ is then a consequence of the presence of at least one factor $\tilde{\rho}_0$. The gains $\gamma^j$ are again obtained as in the proof of Prop. 11.2 in \cite{KMiao}. 
\end{proof}

If we use the first part of the preceding lemma with $A = 4$ in conjunction with the terms obeying the first inequality in Lemma~\ref{lem:crudeboundc2disperrors} and the second part for those terms obeying the second inequality in Lemma~\ref{lem:crudeboundc2disperrors}, we now infer the following 
\begin{lemma}\label{lem:c2transferencedispcontribu} Let $\gamma>0$ be given, as well as $p>0$. There exists $\tau_0 = \tau_0(\gamma,\nu,p)\gg1 $ such that for any $\tau>\tau_0$ we have the bound  
	\begin{align*}
		&\sum_{j\geq 1}\tau^p\cdot \left|\int_0^\infty U\left(\frac{\lambda_{\tau}}{\lambda}\mathcal{K}_0^{(0)}\mathcal{D}_{\tau}U\right)^j \left(\mathcal{F}^{(0)}\left(\mathcal{D}N\right)\right)(\tau,\xi)\tilde{\rho}_0(\xi)\,d\xi\right|\\&\leq \gamma\cdot \left(\sup_{\tau>\tau_0}\tau^p\left|\frac{\bar{c}(\tau)-1}{\tau^2}\right|+\sup_{\tau>\tau_0}\tau^p\left|\frac{\bar{c}'(\tau)}{\tau}\right| + \sup_{\tau>\tau_0}\tau^p\left|\bar{c}''(\tau)\right|\right)
	\end{align*}
	where $N$ denotes any one of the source terms displayed in \eqref{eq:Qc2modulationterms}, \eqref{eq:Qc2modulationtermsepsilon0}, \eqref{eq:Qc2modulationtermsepsilon1}.
\end{lemma}

Coming back to the ansatz \eqref{eq:precisemodulation}, we still need to account for errors due to modulating in $c_2$ and involving the interaction of the modified bulk term $\Phi^{(\text{mod}_2(t))}$ with the perturbation 
\[
\phi_1E_1 + \phi_2E_2 + a(\Pi_{\Phi^{\perp}}\phi),
\]
since the leading order potential term is now slightly different. The key is that the $n = 0$-components of $\varepsilon_{\pm}$ decay well, namely 
\begin{align*}
	\left|\langle R\rangle^{-1}\varepsilon_{\pm}(0)(\tau, R)\right|\lesssim \tau^{-3+\kappa},
\end{align*}
as follows from Definition~\ref{defi:xsingulartermsnless2smooth}.
We now analyze the additional errors generated by replacing $U = Q_{\lambda} + \epsilon$ by $\chi_{R<\delta\tau}\cdot\mathcal{S}_{\bar{c}}U + (1-\chi)\cdot\left(\mathcal{S}_{\bar{c}}Q + \epsilon\right)$ and involving interactions with the perturbation. 
To leading order, these errors are linear in the perturbation, and of the form (see \eqref{coe eq more precise 1}, \eqref{coe eq more precise 2})
\begin{align*}
	\frac{\cos(2Q_{\bar{c}}) - \cos(2Q)}{R^2}\cdot \varepsilon_{\pm},\,\chi\left(\frac{2\sin\left(2(Q_{\bar{c}}+\epsilon_{\bar{c}})\right)\sin\epsilon_{\bar{c}}}{R^2} - \frac{2\sin(2Q+\epsilon)\sin\epsilon}{R^2}\right)\cdot \varepsilon_{\pm}.
\end{align*}
We can summarize the contribution of these to the source of the instability at $n = 0$ via the following 
\begin{lemma}\label{lem:c2potentialchange} Letting $N$ be any one of the immediately preceding error terms, given $\gamma>0$, there is $\tau_0 = \tau_0(\gamma,\nu,p)\gg 1$ such that we have 
	\begin{align*}
		\tau^p\cdot \left|\int_0^1 U\left(\mathcal{F}^{(0)}\left(\mathcal{D}N\right)\right)(\tau,\xi)\tilde{\rho}_0(\xi)\,d\xi\right|\leq \gamma\cdot \left(\sup_{\tau>\tau_0}\tau^p\left|\frac{\bar{c}(\tau)-1}{\tau^2}\right|+\sup_{\tau>\tau_0}\tau^p\left|\frac{\bar{c}'(\tau)}{\tau}\right| + \sup_{\tau>\tau_0}\tau^p\left|\bar{c}''(\tau)\right|\right).
	\end{align*}
	Moreover, we can write 
	\begin{align*}
		\int_1^\infty U\left(\mathcal{F}^{(0)}\left(\mathcal{D}N\right)\right)(\tau,\xi)\tilde{\rho}_0(\xi)\,d\xi = \partial_{\tau}A(\tau),
	\end{align*}
	where we have the bound 
	\begin{align*}
		\left|\tau^p\cdot A(\tau)\right|\leq \gamma\cdot \left(\sup_{\tau>\tau_0}\tau^p\left|\frac{\bar{c}(\tau)-1}{\tau^2}\right|+\sup_{\tau>\tau_0}\tau^p\left|\frac{\bar{c}'(\tau)}{\tau}\right| + \sup_{\tau>\tau_0}\tau^p\left|\bar{c}''(\tau)\right|\right).
	\end{align*}
	Similarly, we have the bound 
	\begin{align*}
		&\sum_{j\geq 1}\tau^p\cdot \left|\int_0^\infty U\left(\frac{\lambda_{\tau}}{\lambda}\mathcal{K}_0^{(0)}\mathcal{D}_{\tau}U\right)^j \left(\mathcal{F}^{(0)}\left(\mathcal{D}N\right)\right)(\tau,\xi)\tilde{\rho}_0(\xi)\,d\xi\right|\\&\leq \gamma^j\cdot \left(\sup_{\tau>\tau_0}\tau^p\left|\frac{\bar{c}(\tau)-1}{\tau^2}\right|+\sup_{\tau>\tau_0}\tau^p\left|\frac{\bar{c}'(\tau)}{\tau}\right| + \sup_{\tau>\tau_0}\tau^p\left|\bar{c}''(\tau)\right|\right).
	\end{align*}
\end{lemma}
\begin{proof} We explain this for the first term, the other ones being handled similarly. Write 
	\begin{align*}
		\frac{\cos(2Q_{\bar{c}}) - \cos(2Q)}{R^2}\cdot \varepsilon_{\pm} = \frac{2\sin\left(Q_{\bar{c}} - Q\right)\cdot \left(Q_{\bar{c}}+Q\right)}{R^2}\cdot \varepsilon_{\pm}\sim \frac{1}{\langle R\rangle^4}\cdot (\bar{c} - 1)\cdot \varepsilon_{\pm}.
	\end{align*}
	This implies that calling the preceding term $N$, we have ($\kappa = 0, 1$)
	\begin{align*}
		\left|\partial_{\xi^{\frac12}}^{\kappa}\mathcal{F}^{(0)}\left(\mathcal{D}(N)\right)(\tau,\cdot)\right|\lesssim \left|\bar{c} - 1\right|\cdot \left(\left\|\frac{\varepsilon_{\pm}(\tau,\cdot)}{\langle R\rangle}\right\|_{L^\infty_{dR}} + \left\|\mathcal{D}\varepsilon_{\pm}(\tau,\cdot)\right\|_{L^\infty_{dR}}\right).
	\end{align*}
	To handle the large frequencies ($\xi\gg 1$), we complement this with 
	\begin{align*}
		\left|\langle\xi\rangle^{2+}\partial_{\xi^{\frac12}}^{\kappa}\mathcal{F}^{(0)}\left(\mathcal{D}(N)\right)(\tau,\cdot)\right|\lesssim \tau^{-3+\kappa}\cdot \left|\bar{c} - 1\right|\cdot\left\|\mathcal{F}^{(0)}(\varepsilon_+)\right\|_{good}
	\end{align*}
	To deduce the first bound of the lemma, switch the orders of integration between $\xi$ and $\sigma$(inside the $U$-propagator), and perform integration by parts with respect to $\xi^{\frac12}$, thereby gaining a factor
	\[
	\lambda(\tau)\int_{\sigma}^{\tau}\lambda^{-1}(s)\,ds;
	\]
	this helps us avoid a loss with respect to $\tau$-decay up to a factor $\log\tau$, and we infer the bound 
	\begin{align*}
		\tau^p\cdot \left|\int_0^1 U\left(\mathcal{F}^{(0)}\left(\mathcal{D}N\right)\right)(\tau,\xi)\tilde{\rho}_0(\xi)\,d\xi\right|\leq \tau^{-1+\kappa}\log\tau\cdot \sup_{\tau>\tau_0}\tau^{p-2}\left|\bar{c}(\tau) - 1\right|\cdot \left\|\mathcal{F}^{(0)}(\varepsilon_+)\right\|_{good}.
	\end{align*}
	which implies the desired inequality. 
\end{proof}

We can now formulate the key result for how to determine the evolution of $\bar{c}$, given a source term $h$:
\begin{proposition}\label{prop:determiningc2} Given $h:[\tau_0,\infty)\rightarrow \R$ with $\sup_{\tau>\tau_0}\left|\tau^p\cdot h(\tau)\right|<\gamma\ll 1$, and assuming $2<p<p_*(\nu)$,  there is a function $\bar{c}: [\tau_0,\infty)\rightarrow \R$ with 
	\[
	\sup_{\tau>\tau_0}\tau^p\left|\frac{\bar{c}(\tau)-1}{\tau^2}\right|+\sup_{\tau>\tau_0}\tau^p\left|\frac{\bar{c}'(\tau)}{\tau}\right| + \sup_{\tau>\tau_0}\tau^p\left|\bar{c}''(\tau)\right|\lesssim \gamma, 
	\]
	such that the following holds: let ${\bf{\mathcal{K}}_0}\overline{x}$ be the sum of the last five terms in \eqref{eq:neq0Fourier1}, where $\varepsilon_+^0$ is determined via 
	\begin{align*}
		\varepsilon_+^0(\tau, R) = \phi_0(R)\cdot\int_0^R [\phi_0(s)]^{-1}\cdot \mathcal{D}_0\varepsilon_+^0(\tau, s)\,ds,\,\mathcal{D}_0\varepsilon_+^0(\tau, R) = \int_0^\infty \phi_0(R;\xi)\cdot \overline{x}(\tau,\xi)\cdot \tilde{\rho}_0(\xi)\,d\xi. 
	\end{align*}
	If $N$ is the sum of all the source terms generated by modulating in $\bar{c}$, and if 
	\[
	n(\tau) = \lim_{R\rightarrow 0}R^{-1}N(\tau, R), 
	\]
	then we have the identity
	\begin{align*}
		h(\tau) = n(\tau) + \sum_{j\geq 0}\int_0^\infty U(\tau)\left({\bf{\mathcal{K}}_0} U\right)^j\mathcal{F}^{(0)}\left(\mathcal{D}N\right)\,\tilde{\rho}_0(\xi)\,d\xi
	\end{align*}
\end{proposition}
\begin{proof} We first determine the precise form of $n(\tau)$: observe that due to the higher order vanishing of $\epsilon(\tau, R)$ at $R = 0$, only the terms \eqref{eq:Qc2modulationterms} contribute. This leads to 
	\begin{align*}
		n(\tau) = 2\cdot\left(\bar{c}''(\tau) + 3\cdot (1+\nu^{-1})\cdot \frac{\bar{c}'(\tau)}{\tau} + \gamma\cdot [\bar{c}'(\tau)]^2\right)
	\end{align*}
	for a suitable constant $\gamma = \gamma(\nu)$. Then we observe that the equation 
	\begin{align*}
		\bar{c}''(\tau) + 3\cdot (1+\nu^{-1})\cdot \frac{\bar{c}'(\tau)}{\tau} = \frac{h(\tau)}{2}
	\end{align*}
	can be solved by means of the variation of constants formula 
	\begin{align*}
		\bar{c}(\tau) = 1+(3+2\nu^{-1})^{-1}\left(\tau^{-2-3\nu^{-1}}\cdot\int_{\tau_0}^{\tau}s^{3+3\nu^{-1}}h(s)\,ds - \int_{\tau_0}^{\tau}s h(s)\,ds\right)
	\end{align*}
	and this satisfies the desired bound provided $p\neq 4+3\nu^{-1}$. The assumption $p>2$ and the smallness of $\gamma$ imply that 
	\begin{align*}
		\sup_{\tau>\tau_0}\left|\tau^p\cdot [c'(\tau)]^2\right|\ll \gamma, 
	\end{align*}
	whence the final term in the equation for $n(\tau)$ and quadratic in $c'(\tau)$ is perturbative. Furthermore, invoking Lemma~\ref{lem:c2potentialchange}, we find that 
	\begin{align*}
		\sup_{\tau>\tau_0}\tau^p\cdot\left|\sum_{j\geq 0}\int_0^\infty U(\tau)\left({\bf{\mathcal{K}}_0} U\right)^j\mathcal{F}^{(0)}\left(\mathcal{D}N\right)\,\tilde{\rho}_0(\xi)\,d\xi\right|\ll \sup_{\tau>\tau_0}\tau^p\left|\frac{\bar{c}(\tau)-1}{\tau^2}\right|+\sup_{\tau>\tau_0}\tau^p\left|\frac{\bar{c}'(\tau)}{\tau}\right| + \sup_{\tau>\tau_0}\tau^p\left|\bar{c}''(\tau)\right|.
	\end{align*}
	The claim of the proposition now follows from a simple fixed point argument. 
	\end{proof}

\subsection{Controlling the source term $h$}\label{subsec:hcontrol}

In order to take advantage of Proposition~\ref{prop:determiningc2}, we now need to analyze the $h$ which needs to be canceled via appropriate choice of $\bar{c}$. From \eqref{eq:c0evolutioneqn}, we know that $h$ consists of the contributions of the terms in $F_+(0)$ on the right in \eqref{eq diag phys} to 
\[
\lim_{R\rightarrow 0}R^{-1}F_+(0),\quad \lim_{R\rightarrow 0}R^{-1}H_0^+\varepsilon_+^0. 
\]
We observe here that as far as determining $\bar{c}$ via Proposition~\ref{prop:determiningc2} is concerned, it is only the real part of these source terms which matters, while the imaginary part will serve as a source term for the equation determining $\bar{h}$ later on.
Dealing with the first type of term is straightforward, we can state 
\begin{proposition}\label{prop:F+(0)directsource} Let $N(\tau, R)$ be any one of the terms in \eqref{coe eq precise 1}, \eqref{coe eq precise 2},\eqref{nonlinearity tau R}, projected onto the $n = 0$ mode, and assume that it is $k$-linear in terms of its dependence on $\phi_{1,2}$. Then assuming the factors to be $|n|\geq 2$ angular momentum functions with distorted Fourier transforms $\overline{x_j}(\tau,\xi)$ as angular momentum $n_j$ functions (with $|n_j|\geq 2,\,\sum_j n_j = 0$), we have the bound 
	\begin{align*}
		\left|\lim_{R\rightarrow 0}R^{-1}N(\tau, R)\right|\lesssim_{\max\{|n_j|\}} \tau^{-6+2\nu}\cdot \prod_j \left\|\overline{x_j}\right\|_{good_r}.
	\end{align*}
\end{proposition}
\begin{proof} From the construction of $\epsilon$, we observe that the explicitly displayed terms in \eqref{nonlinearity tau R} which depend linearly on $\phi_j$, $j = 1, 2$, vanish to order three at $R = 0$ and hence do not contribute to the source term $\lim_{R\rightarrow 0}R^{-1}N(\tau, R)$. It hence suffices to consider the contributions of $\lambda^{-2}N(\phi_j),\,j = 1, 2$.  
	\\
	
	We consider here the contribution of the cubic null-forms (see subsection~\ref{subset:hardnullformestimates}) as well as the quadratic terms. Thus, to begin with, consider a term of the form 
	\begin{align*}
		\mathcal{N}_0(\varepsilon_1, \varepsilon_2, \varepsilon_3) = \varepsilon_1\cdot \left(\partial_{\tau} + \frac{\lambda_{\tau}}{\lambda}R\partial_R - \partial_R\right)\varepsilon_2\cdot  \left(\partial_{\tau} + \frac{\lambda_{\tau}}{\lambda}R\partial_R + \partial_R\right)\varepsilon_3,
	\end{align*}
	where the functions $\varepsilon_j$ stand for $\varepsilon_{\pm}(n_j) = \hat{\varphi}_1(n_j) \mp i \hat{\varphi}_2(n_j)$ with $|n_j|\geq 2$, and $\sum_j n_j = 0$. Then from Lemma~\ref{lem:singFouriertiphysicalngeq2admDeriv}, Lemma~\ref{lem:derLinfty}, Lemma~\ref{lem:structuredlowfreqsmoothbasic},  we infer the bounds
	\begin{align*}
		\left|R^{-1}\varepsilon_1\right|\lesssim_{\max\{|n_j|\}} \tau^{-3+\nu}\cdot \left\|\overline{x_1}\right\|_{good_r},\quad \left|\left(\partial_{\tau} + \frac{\lambda_{\tau}}{\lambda}R\partial_R - \partial_R\right)\varepsilon_2\right|\lesssim_{\max\{|n_j|\}} \tau^{-3+\nu}\cdot \left\|\overline{x_2}\right\|_{good_r},
	\end{align*}
	and similarly for the third factor $\left(\partial_{\tau} + \frac{\lambda_{\tau}}{\lambda}R\partial_R + \partial_R\right)\varepsilon_3$, and so the desired bound follows easily, in fact with a temporal decay of $\tau^{-9+3\nu}$. 
	\\
	As far as the quadratic terms are concerned, these are of the form (recalling that we project onto the $n = 0$ angular mode)
	\begin{align*}
		\varepsilon_1\cdot \left(\partial_{\tau} + \frac{\lambda_{\tau}}{\lambda}R\partial_R - \partial_R\right)\varepsilon_2\cdot  \left(\partial_{\tau} + \frac{\lambda_{\tau}}{\lambda}R\partial_R + \partial_R\right)U,\,\varepsilon_1\cdot \left(\partial_{\tau} + \frac{\lambda_{\tau}}{\lambda}R\partial_R - \partial_R\right)U\cdot  \left(\partial_{\tau} + \frac{\lambda_{\tau}}{\lambda}R\partial_R + \partial_R\right)\varepsilon_2,
	\end{align*}
	contributed by the second line in \eqref{coe eq precise 1} as well as by a term coming from \eqref{non linear 1 2} in interaction with $\varphi_1$ in  \eqref{coe eq precise 1}, and further the terms
	\begin{align*}
		\frac{2\sin U}{R^2}\varepsilon_{1,\theta}\cdot \varepsilon_2,\quad \frac{\sin 2U}{R^2}\varepsilon_1\cdot \varepsilon_2. 
	\end{align*}
	For each of these terms the required bound for the contribution to $\lim_{R\rightarrow 0}R^{-1}N(\tau, R)$ follows directly from the three lemmas cited above. The remaining terms of higher degree of multilinearity are handled similarly. 
\end{proof}

We now turn to the more delicate task of bounding the contribution of the source terms \eqref{coe eq precise 1}, \eqref{coe eq precise 2},\eqref{nonlinearity tau R} to $\lim_{R\rightarrow 0}R^{-1}H_0^+\varepsilon_+^0$. We recall that this will be of the form 
\begin{align*}
	\int_0^\infty U\left(\mathcal{F}^{(0)}\left(\mathcal{D}N\right)\right)\tilde{\rho}_0(\xi)\,d\xi
\end{align*}
for the primary contribution, as well as 
\begin{align*}
	\sum_{j\geq 1}\int_0^\infty U(\tau)\left({\bf{\mathcal{K}}_0} U\right)^j\mathcal{F}^{(0)}\left(\mathcal{D}N\right)\,\tilde{\rho}_0(\xi)\,d\xi\
\end{align*} for the secondary contributions via the transference operator, 
where $N$ denotes the original source terms (not those arising from modulating over $c_2$). We shall have to control the contributions from all terms in \eqref{coe eq precise 1}, \eqref{coe eq precise 2},\eqref{nonlinearity tau R} also keeping in mind (3.49). We start with the most delicate type of terms, namely the cubic singular ones, where we have taken advantage of the primary modulation step to ensure better temporal decay.  The following bounds are then useful for the determination of both $\bar{c}$ and $\bar{h}$. 
\\

{\it{Contribution to $h$ from cubic singular terms}}. As in the preceding, we shall have to control the {\it{direct effect}} via the Duhamel parametrix and the {\it{indirect effect}} via the (iterated) transference operator contributions. The key point shall be that except for the leading order singular terms and incoming singularities, we are only dealing with perturbative contributions. We begin with 
\begin{lemma}\label{lem:hcubicmainsingular} Let 
	\[
	N = \varepsilon_1\cdot\left(\partial_{\tau} + \frac{\lambda_{\tau}}{\lambda}R\partial_R - \partial_R\right)\varepsilon_2\cdot \left(\partial_{\tau} + \frac{\lambda_{\tau}}{\lambda}R\partial_R + \partial_R\right)\varepsilon_3,
	\]
	at angular momentum $n = 0$ (meaning only angular momenta $m_j$, $j = 1, 2, 3$ components for $\varepsilon_j$ with $\sum_j m_j = 0$ contribute). Introduce the smoothly truncated\footnote{Observe that this corresponds to $M\cdot(t-r)\gtrsim 1$, i. e. a smaller light cone.} 
	\begin{align*}
		\chi_{M(\nu\tau-R)\gtrsim \lambda(\tau)}N =: N_{M}.
	\end{align*}
	Then we have the estimate 
	\begin{align*}
		\left|\int_0^\infty U(\tau)\left(\mathcal{F}^{(0)}\left(\mathcal{D} N_M\right)\right)\tilde{\rho}^{(0)}(\xi)\,d\xi\right|\lesssim \tau^{-(4-10\nu)}\cdot \prod_{j=1}^3\left\|\overline{x}_j\right\|_{good_r},
	\end{align*}
	provided the vanishing conditions \eqref{eq:n=0vanishingrelations} are satisfied for this source term. 
\end{lemma}
\begin{proof} We  consider the situation where all factors (`inputs') are angular momentum $|m_j|\geq 2$ functions, the case of exceptional angular momenta only requiring small adjustmens. We may neglect the dependence on $|m_j|$ in the bounds, since the largest such angular momentum has to cancel against the remaining two, and so there is at least another comparable angular momentum. The main singularity only comes from the third factor 
	\[
	\left(\partial_{\tau} + \frac{\lambda_{\tau}}{\lambda}R\partial_R + \partial_R\right)\varepsilon_3.
	\]
	We can expand it by taking advantage of Lemma~\ref{lem:DerivonPrincSing}. Furthermore use Lemma~\ref{lem:singFouriertiphysicalngeq2admDeriv}  for the middle factor, and we can take advantage of a first order Taylor expansion around $R = \nu\tau$ for the first factor. We then decompose the Fourier coefficient into two parts:
	\begin{equation}\label{eq:FzeroDNsingular}
		\mathcal{F}^{(0)}\left(\mathcal{D} N_M\right) = \left\langle \mathcal{D} N_M,\,\chi_{R\xi^{\frac12}\gtrsim 1}\phi_0(R;\xi)\right\rangle_{L^2_{R\,dR}} +  \left\langle \mathcal{D} N_M,\,\chi_{R\xi^{\frac12}\lesssim 1}\phi_0(R;\xi)\right\rangle_{L^2_{R\,dR}},
	\end{equation}
	where we focus for now on the {\it{high-frequency regime}} $\xi\gg 1$. Then the second term on the right leads to a straightforward contribution, since we necessarily have $R\ll 1$ in the inner product. Taking advantage of Lemma \ref{lem:singFouriertiphysicalngeq2adm}, we directly infer the bound 
	\begin{align*}
		\left|\left\langle \mathcal{D} N_M,\,\chi_{R\xi^{\frac12}\lesssim 1}\phi_0(R;\xi)\right\rangle_{L^2_{R\,dR}}\right|\lesssim \tau^{-5}\cdot \prod_{j=1}^3\left\|\overline{x}_j\right\|_{good_r}, 
	\end{align*}
	which is much better than what we require. For this estimate observe that only the term $f_1$ in Lemma \ref{lem:singFouriertiphysicalngeq2adm} occurs for the singular contribution of $\varepsilon_3$ (should include a lemma somewhere which deduces this bound for general good function in regime $R\ll\tau$). 
	We hence reduce to the contribution in the singular region, i. e. the term 
	\[
	\left\langle \mathcal{D} N_M,\,\chi_{R\xi^{\frac12}\gtrsim 1}\phi_0(R;\xi)\right\rangle_{L^2_{R\,dR}} 
	\]
	We can immediately introduce a cutoff $\chi_{R\gtrsim\tau}$, since else the same bound as for the previous term is obtained.  We then use the asymptotic description (in the regime $R\xi^{\frac12}\gtrsim 1$)
	\begin{align*}
		\phi_0(R;\xi) = R^{-\frac12}\cdot\sum_{\pm}\frac{e^{\pm iR\xi^{\frac12}}}{\xi^{\frac54}}\cdot\sigma_{\pm}(R\xi^{\frac12}, R),
	\end{align*}
	see \cite{KMiao} as well as \cite{KST}. We then take advantage of the fine structure of $\left(\partial_{\tau} + \frac{\lambda_{\tau}}{\lambda}R\partial_R + \partial_R\right)\varepsilon_3$ via Lemma \ref{lem:singFouriertiphysicalngeq2admDeriv} and replace the remaining terms in terms of a first order Taylor expansion around $R = \nu\tau$. Observe that the integration over $\nu\tau - R$ `costs' nothing thanks to the decay with respect to this variable for the singular part of $\varepsilon_3$. We hence arrive at the asymptotic formula (assuming $\xi\gtrsim 1$)
	\begin{equation}\label{eq:preciseasympto1forlhighfreq}\begin{split}
			&\left\langle \mathcal{D} N_M,\,\chi_{R\xi^{\frac12}\gtrsim 1}\phi_0(R;\xi)\right\rangle_{L^2_{R\,dR}}\\&  =\sum_{k = 1}^3 \chi_{\lambda(\tau)\xi^{\frac12} M^{-1}\lesssim 1}^{(1,k)}\sum_{\pm} \log\xi\cdot \frac{e^{\pm i\nu\tau\xi^{\frac12}}}{\xi^{1+k\cdot\frac{\nu}{2}}}\cdot a(\tau) + \sum_{k = 1}^3 \chi_{\lambda(\tau)\xi^{\frac12} M^{-1}\lesssim 1}^{(2,k)}\sum_{\pm} \frac{e^{\pm i\nu\tau\xi^{\frac12}}}{\xi^{1+k\cdot\frac{\nu}{2}}}\cdot b(\tau) + O\left(\log\tau\cdot \tau^{-4-\nu}\cdot \log\xi\cdot\xi^{-\frac32}\right),\\
			&\left|a(\tau)\right| + \left|b(\tau)\right|\lesssim \tau^{-4-\nu}\cdot\log\tau,
	\end{split}\end{equation}
	where $ \chi_{\lambda(\tau)\xi^{\frac12} M^{-1}\lesssim 1}^{(j,k)}$ are smooth functions which vanish rapidly beyond scale $\frac{M}{\lambda(\tau)}$, and with the property that (locally uniformly in $\xi$)
	\begin{align*}
		\lim_{M\rightarrow+\infty}\chi_{\lambda(\tau)\xi^{\frac12} M^{-1}\lesssim 1}^{(j,k)} = c_*^{(j,k)}
	\end{align*}
	where $c_*^{(j,k)}$ are constants. Specifically, the functions $\chi_{\lambda(\tau)\xi^{\frac12} M^{-1}\lesssim 1}^{(j,k)}$ are given by 
	\begin{align*}
		\chi_{\lambda(\tau)\xi^{\frac12} M^{-1}\lesssim 1}^{(1,k)} =\sum_{\pm}\int_0^\infty \chi_{M(\nu\tau - R)\gtrsim \lambda(\tau)}\cdot\frac{c_{\pm}\cdot e^{\pm i(R-\nu\tau)\xi^{\frac12}}}{[\xi^{\frac12}(\nu\tau - R)]^{\frac12-k\nu}}d\left(\xi^{\frac12}(\nu\tau - R)\right)
	\end{align*}
	for suitable complex constants $c_{\pm}$, while we have 
	\begin{align*}
		\chi_{\lambda(\tau)\xi^{\frac12} M^{-1}\lesssim 1}^{(2,k)} =\sum_{\pm}\int_0^\infty \chi_{M(\nu\tau - R)\gtrsim \lambda(\tau)}\cdot\frac{c_{\pm}\cdot e^{\pm i(R-\nu\tau)\xi^{\frac12}}}{[\xi^{\frac12}(\nu\tau - R)]^{\frac12-k\nu}}\log\left(\xi^{\frac12}(\nu\tau - R)\right)\cdot d\left(\xi^{\frac12}(\nu\tau - R)\right).
	\end{align*}
	We note that by means of the substitution $y = \xi^{\frac12}(\nu\tau - R)$ we can write the preceding integrals in the form 
	\begin{align*}
		\chi_{\lambda(\tau)\xi^{\frac12} M^{-1}\lesssim 1}^{(1,k)} &=\sum_{\pm}\int_0^\infty \chi_{M y\xi^{-\frac12}\gtrsim \lambda(\tau)}\cdot\frac{c_{\pm}\cdot e^{\pm i y}}{y^{\frac12-k\nu}}d y\\
		& = \sum_{\pm}\int_0^\infty\frac{c_{\pm}\cdot e^{\pm i y}}{y^{\frac12-k\nu}}d y - \sum_{\pm}\int_0^\infty \chi_{M y\xi^{-\frac12}\lesssim \lambda(\tau)}\cdot\frac{c_{\pm}\cdot e^{\pm i y}}{y^{\frac12-k\nu}}d y\\
	\end{align*}
	and similarly for $ \chi_{\xi^{\frac12} M^{-1}\lesssim 1}^{(2,k)}$. In particular we infer that 
	\begin{align*}
		c_*^{(1,k)} = \sum_{\pm}\int_0^\infty \frac{c_{\pm}\cdot e^{\pm i y}}{y^{\frac12-k\nu}}d y,
	\end{align*}
	and analogously for $ c_*^{(2,k)}$. 
	Observe that the presence of the operator $\mathcal{D}$ 'costs' an additional factor $\xi^{\frac12}$. We also note that the error term $O\left(\ldots\right)$ admits a finer structure 
	\begin{align*}
		O\left(\ldots\right) = \sum_{\pm}e^{\pm i\nu\tau\xi^{\frac12}}\cdot h_{\pm}(\tau,\xi) + \tilde{h}(\tau,\xi)
	\end{align*}
	where we have 
	\begin{align*}
		&\left|\partial_{\xi^{\frac12}}^{\kappa}h_{\pm}(\tau,\xi)\right|\lesssim\log\tau\cdot \tau^{-4-\nu}\cdot \log\xi\cdot\xi^{-\frac32-\frac{\kappa}{2}},\,\kappa\in \{0, 1\}, \\
		&\left|\tilde{h}(\tau,\xi)\right|\lesssim \log\tau\cdot \tau^{-5-\nu}\cdot \log\xi\cdot\xi^{-2-}.
	\end{align*}
	\\
	The preceding concludes our analysis of the large frequency regime $\xi\gg 1$. 
	\\
	
	Next, we turn to the small frequency case $\xi\lesssim 1$. Here we use the asymptotics 
	\begin{align*}
		\phi_0(R;\xi) = R^{-\frac12}\cdot \sum_{\pm}\log\xi\cdot \frac{e^{\pm iR\xi^{\frac12}}}{\xi^{\frac14}}\cdot \sigma_{\pm}(R\xi^{\frac12}, R)
	\end{align*}
	provided $R\xi^{\frac12}\gtrsim 1$, while we have the more crude bound 
	\begin{align*}
		\left|\phi_0(R;\xi)\right|\lesssim \langle \log\langle R\rangle\rangle
	\end{align*}
	in the regime $R\xi^{\frac12}\lesssim 1$. If the third factor $\varepsilon_3$ is of singular type, then we perform integration by parts in order to move the operator $\mathcal{D}$ onto $\phi_0(R;\xi)$, which gains at least a factor $\xi^{\frac12}$, and we note that 
	\begin{align*}
		\left|\xi^{\frac12}\cdot\phi_0(R;\xi)\right|\lesssim R^{-\frac12}\xi^{\frac14}\cdot \left|\log\xi\right|. 
	\end{align*}
	Keeping in mind the fast decay with respect to $\nu\tau - R$ for the singular terms according to Lemma 7.10, we infer the bound (throughout $\xi\lesssim 1$)
	\begin{align*}
		\left|\left\langle \mathcal{D} N_M,\,\phi_0(R;\xi)\right\rangle_{L^2_{R\,dR}}\right|\lesssim \xi^{\frac14}\left|\log\xi\right|\cdot \tau^{-4-\nu}\log^3\tau
	\end{align*}
	More precisely, by using the cutoffs $\chi_{R\xi^{\frac12}\gtrsim 1}$ and $\chi_{R\xi^{\frac12}\lesssim 1}$ and replacing $e^{\pm iR\xi^{\frac12}}$ by $e^{\pm i\nu\sigma\xi^{\frac12}}$, we obtain 
	\begin{equation}\label{eq:preciseasympto1forlowfreq}
		\left\langle \mathcal{D} N_M,\,\phi_0(R;\xi)\right\rangle_{L^2_{R\,dR}} = \sum_{\pm}e^{\pm i\nu\tau\xi^{\frac12}}\cdot n_{\pm}(\tau,\xi) + n_0(\tau,\xi), 
	\end{equation}
	where we have the derivative bounds 
	\begin{align*}
		\left|\partial_{\xi^{\frac12}}^{\kappa}n_j(\tau,\xi)\right|\lesssim \xi^{\frac14-\frac{\kappa}{2}}\cdot \left|\log\xi\right|\cdot \tau^{-4-\nu}\log^3\tau, 
	\end{align*}
	provided $\kappa\in \{0, 1\}$ and $j = \pm, 0$ (this being due to the fact that we replace the phases $e^{\pm i R\xi^{\frac12}}$ in the oscillatory regime by $e^{\pm i(R-\nu\tau)\xi^{\frac12}}$, and a loss of a factor $R-\nu\tau)$ due to differentiation with respect to $\xi^{\frac12}$ is counteracted by the decay from $\varepsilon_3$) . 
	\\
	In the other case that $\varepsilon_3$ is of smooth type, the additional derivative allows us to gain an additional decay with respect to $\nu\tau - R$, see Lemma~\ref{lem:structuredlowfreqsmoothbasic}. The extra derivative $\mathcal{D}$ can then be let to fall on the triple product, which at least gains a factor (Lemma~\ref{lem:structuredlowfreqsmoothbasic})
	\[
	\max\{(\nu\tau - R)^{-1}, R^{-1}\}.
	\]
	We then conclude that the same representation as before obtains, i. e. \eqref{eq:preciseasympto1forlowfreq}.
	We can now infer the effect of the Duhamel parametrix and deduce the desired bound for this contribution, based on \eqref{eq:preciseasympto1forlhighfreq} together with \eqref{eq:preciseasympto1forlowfreq}. 
	Write
	\begin{align*}
		\int_0^\infty U(\tau)\left(\mathcal{F}^{(0)}\left(\mathcal{D} N_M\right)\right)\tilde{\rho}^{(0)}(\xi)\,d\xi &= \int_0^\infty\int_{\tau_0}^{\tau} \chi_{\frac{\lambda^2(\tau)}{\lambda^2(\sigma)}\xi\lesssim 1}U(\tau,\sigma;\xi))\left(\mathcal{F}^{(0)}\left(\mathcal{D} N_M\right)\right)\left(\sigma, \frac{\lambda^2(\tau)}{\lambda^2(\sigma)}\xi\right)\tilde{\rho}^{(0)}(\xi)\,d\xi\\
		& +   \int_0^\infty\int_{\tau_0}^{\tau} \chi_{\frac{\lambda^2(\tau)}{\lambda^2(\sigma)}\xi\gtrsim 1}U(\tau,\sigma;\xi))\left(\mathcal{F}^{(0)}\left(\mathcal{D} N_M\right)\right)\left(\sigma, \frac{\lambda^2(\tau)}{\lambda^2(\sigma)}\xi\right)\tilde{\rho}^{(0)}(\xi)\,d\xi\\
		& =: I + II, 
	\end{align*}
	where we set 
	\begin{align*}
		U(\tau,\sigma;\xi)) = \frac{\tilde{\rho}_0^{\frac12}\left(\frac{\lambda^2(\tau)}{\lambda^2(\sigma)}\xi\right)}{\tilde{\rho}_{0}^{\frac12}(\xi)}\cdot \frac{\sin\left(\lambda(\tau)\xi^{\frac12}\int_{\sigma}^{\tau}\lambda^{-1}(s)\,ds\right)}{\xi^{\frac12}}.
	\end{align*}
	
	{\it{The bound for $I$.}} Here we use the asymptotics \eqref{eq:preciseasympto1forlowfreq}. Writing the propagator $U(\tau,\sigma;\xi))$ in terms of exponentials, we need to bound 
	\begin{align*}
		\int_0^\infty\int_{\tau_0}^{\tau} \chi_{\frac{\lambda^2(\tau)}{\lambda^2(\sigma)}\xi\gtrsim 1}U(\tau,\sigma;\xi))\cdot \left(\sum_{\pm}e^{\pm i\nu\sigma\frac{\lambda(\tau)}{\lambda(\sigma)}\cdot\xi^{\frac12}}\cdot n_{\pm}\left(\sigma,\frac{\lambda(\tau)}{\lambda(\sigma)}\cdot\xi\right) + n_0\left(\sigma,\frac{\lambda(\tau)}{\lambda(\sigma)}\cdot\xi\right)\right)\tilde{\rho}^{(0)}(\xi)\,d\sigma d\xi,
	\end{align*}
	where combining the oscillatory phases we encounter the following exponentials:
	\[
	e^{\pm i\nu\tau\xi^{\frac12}},\,e^{\pm i\left(\nu\tau - \nu\sigma\frac{\lambda(\tau)}{\lambda(\sigma)}\right)\xi^{\frac12}},\,e^{\pm i\left(\nu\tau - 2\nu\sigma\frac{\lambda(\tau)}{\lambda(\sigma)}\right)\xi^{\frac12}}.
	\]
	In each case, we interpolate between a bound obtained via integration by parts with respect to $\xi^{\frac12}$ and a bound obtained by exploiting the localisation $\xi^{\frac12}\lesssim \frac{\lambda(\sigma)}{\lambda(\tau)}$. Letting $\kappa=10\nu>0$ a small constant, this results in 
	\begin{align*}
		&\left|I\right|\lesssim \int_0^\infty\int_{\tau_0}^{\tau}(\tau-\sigma)^{-(1-\kappa)}\cdot \left(\frac{\lambda(\sigma)}{\lambda(\tau)}\right)^{\kappa}\cdot\left(\left|\left\langle\xi\partial_{\xi}\right\rangle n_{\pm}\left(\sigma,\frac{\lambda(\tau)}{\lambda(\sigma)}\cdot\xi\right)\right| +  \left|\left\langle\xi\partial_{\xi}\right\rangle n_0\left(\sigma,\frac{\lambda(\tau)}{\lambda(\sigma)}\cdot\xi\right)\right|\right)\xi^{-1}\tilde{\rho}^{(0)}(\xi)\,d\sigma d\xi\\
		&\ll_{\tau_0}\tau^{-(4-10\nu)}. 
	\end{align*}
	
	Next we turn to the contribution of the term $II$, where we shall rely on \eqref{eq:preciseasympto1forlhighfreq}, including the fine structure of the error term. We start with the contribution of the principal singular part, which consists of the first two terms on the right hand side in \eqref{eq:preciseasympto1forlhighfreq}. To begin with, note that the re-scaled cutoff
	\begin{align*}
		\chi_{\lambda(\sigma)\left(\frac{\lambda(\tau)}{\lambda(\sigma)}\xi\right)^{\frac12} M^{-1}\lesssim 1}^{(j,k)} = \chi_{\lambda(\tau)\xi^{\frac12} M^{-1}\lesssim 1}^{(j,k)}
	\end{align*}
	is independent of the time integration variable $\sigma$ and hence can be moved outside of $\int_{\tau_0}^{\tau}\ldots\,d\sigma$. Combining the oscillatory phase
	\[
	e^{\pm i\nu\sigma\cdot\frac{\lambda(\tau)}{\lambda(\sigma)}\xi^{\frac12}}
	\]
	with the oscillatory kernel 
	\[
	U(\tau,\sigma;\xi)
	\]
	results in the oscillatory phases 
	\[
	e^{\pm i\nu\tau\xi^{\frac12}},\,e^{\pm i\left(\nu\tau - 2\nu\sigma\cdot\frac{\lambda(\tau)}{\lambda(\sigma)}\right)\xi^{\frac12}}.
	\]
	We first show that the contribution of the {\it{second oscillatory term}} leads to a good bound. Keeping in mind that for situation $II$ we include a cutoff $\chi_{\frac{\lambda^2(\tau)}{\lambda^2(\sigma)}\xi\gtrsim 1}$, we reduce to bounding (for $\kappa\in \{0, 1\}$)
	\begin{align*}
		II_1: = \int_0^\infty\chi_{\lambda(\tau)\xi^{\frac12} M^{-1}\lesssim 1}^{(j,k)}\left(\int_{\tau_0}^{\tau}\chi_{\frac{\lambda^2(\tau)}{\lambda^2(\sigma)}\xi\gtrsim 1}\frac{\tilde{\rho}_0^{\frac12}\left(\frac{\lambda^2(\tau)}{\lambda^2(\sigma)}\xi\right)}{\tilde{\rho}_{0}^{\frac12}(\xi)}\cdot \frac{e^{\pm i\left(\nu\tau - 2\nu\sigma\cdot\frac{\lambda(\tau)}{\lambda(\sigma)}\right)\xi^{\frac12}}}{\xi^{\frac12}\cdot\left(\frac{\lambda^2(\tau)}{\lambda^2(\sigma)}\xi\right)^{1+\frac{k\nu}{2}}}\cdot\log^{\kappa}\left(\frac{\lambda^2(\tau)}{\lambda^2(\sigma)}\xi\right)\cdot a(\sigma)\,d\sigma\right) \tilde{\rho}_{0}(\xi)\,d\xi,
	\end{align*}
	Here we can use the bound
	\begin{align*}
		\left|\chi_{\frac{\lambda^2(\tau)}{\lambda^2(\sigma)}\xi\gtrsim 1}\frac{\tilde{\rho}_0^{\frac12}\left(\frac{\lambda^2(\tau)}{\lambda^2(\sigma)}\xi\right)}{\tilde{\rho}_{0}^{\frac12}(\xi)}\cdot \frac{1}{\xi^{\frac12}\cdot\left(\frac{\lambda^2(\tau)}{\lambda^2(\sigma)}\xi\right)^{1+\frac{k\nu}{2}}}\right|\lesssim \frac{\left|\log\xi\right|}{\left(\frac{\lambda^2(\tau)}{\lambda^2(\sigma)}\xi\right)^{\frac{k\nu}{2}}\cdot \xi^{\frac12}},\,\xi\in (0, 1),
	\end{align*}
	as well as 
	\begin{align*}
		\left|\chi_{\frac{\lambda^2(\tau)}{\lambda^2(\sigma)}\xi\gtrsim 1}\frac{\tilde{\rho}_0^{\frac12}\left(\frac{\lambda^2(\tau)}{\lambda^2(\sigma)}\xi\right)}{\tilde{\rho}_{0}^{\frac12}(\xi)}\cdot \frac{1}{\xi^{\frac12}\cdot\left(\frac{\lambda^2(\tau)}{\lambda^2(\sigma)}\xi\right)^{1+\frac{k\nu}{2}}}\right|\lesssim \frac{\left|\log\xi\right|}{\left(\frac{\lambda^2(\tau)}{\lambda^2(\sigma)}\xi\right)^{\frac{k\nu}{2}}\cdot \xi^{\frac32}},\,\xi\in [1,\infty).
	\end{align*}
	Then we take advantage of the fact that 
	\[
	\left|\nu\tau - 2\nu\sigma\cdot\frac{\lambda(\tau)}{\lambda(\sigma)}\right|\gtrsim (\tau-\sigma)\cdot \frac{\lambda(\tau)}{\lambda(\sigma)},\,\tau_0<\sigma<\tau. 
	\]
	Dividing the $\xi$-integral into the regions 
	\[
	\left|\nu\tau - 2\nu\sigma\cdot\frac{\lambda(\tau)}{\lambda(\sigma)}\right|\xi^{\frac12}\gtrsim 1
	\]
	and its complement via a smooth partition of unity and performing integration by parts with respect to $\xi^{\frac12}$ in the latter (and exploiting the absence of boundary terms due to the large frequency cutoff), we easily deduce the bound 
	\begin{align*}
		\left|II_1\right|\lesssim \tau^{-4+10\nu}
	\end{align*}
	in light of the bounds stated in \eqref{eq:preciseasympto1forlhighfreq}.
	\\
	We continue with the contribution of the less oscillatory phase $e^{\pm i\nu\tau\xi^{\frac12}}$, which in analogy to the term $II_1$ leads to integrals of the form 
	\begin{align*}
		II_2 = \int_0^\infty\chi_{\lambda(\tau)\xi^{\frac12} M^{-1}\lesssim 1}^{(j,k)}\Big(\int_{\tau_0}^{\tau}\chi_{\frac{\lambda^2(\tau)}{\lambda^2(\sigma)}\xi\gtrsim 1}\frac{\tilde{\rho}_0^{\frac12}(\frac{\lambda^2(\tau)}{\lambda^2(\sigma)}\xi)}{\tilde{\rho}_{0}^{\frac12}(\xi)}\cdot \frac{e^{\pm i\nu\tau\xi^{\frac12}}}{\xi^{\frac12}\cdot\big(\frac{\lambda^2(\tau)}{\lambda^2(\sigma)}\xi\big)^{1+\frac{k\nu}{2}}}\cdot\log^{\kappa}\big(\frac{\lambda^2(\tau)}{\lambda^2(\sigma)}\xi\big)\cdot a(\sigma)\,d\sigma\Big) \tilde{\rho}_{0}(\xi)\,d\xi.
	\end{align*}
	Here we expand
	\begin{align*}
		\chi_{\frac{\lambda^2(\tau)}{\lambda^2(\sigma)}\xi\gtrsim 1}\tilde{\rho}_0^{\frac12}(\frac{\lambda^2(\tau)}{\lambda^2(\sigma)}\xi) = c\chi_{\frac{\lambda^2(\tau)}{\lambda^2(\sigma)}\xi\gtrsim 1}\cdot \frac{\lambda^2(\tau)}{\lambda^2(\sigma)}\xi + O(1)
	\end{align*}
	The contribution of the term $O(1)$ can be handled in analogy to case $I$. 
	Next splitting 
	\begin{align*}
		\chi_{\frac{\lambda^2(\tau)}{\lambda^2(\sigma)}\xi\gtrsim 1} = 1 - \chi_{\frac{\lambda^2(\tau)}{\lambda^2(\sigma)}\xi\lesssim 1}
	\end{align*}
	and treating the contribution of the low frequency term involving $\chi_{\frac{\lambda^2(\tau)}{\lambda^2(\sigma)}\xi\lesssim 1}$ in analogy to term $I$ above, we reduce to bounding 
	\begin{align*}
		\tilde{II}_2 &= \int_0^\infty\chi_{\lambda(\tau)\xi^{\frac12} M^{-1}\lesssim 1}^{(j,k)}\left(\int_{\tau_0}^{\tau}\frac{\frac{\lambda^2(\tau)}{\lambda^2(\sigma)}\xi}{\tilde{\rho}_{0}^{\frac12}(\xi)}\cdot \frac{e^{\pm i\nu\tau\xi^{\frac12}}}{\xi^{\frac12}\cdot\left(\frac{\lambda^2(\tau)}{\lambda^2(\sigma)}\xi\right)^{1+\frac{k\nu}{2}}}\cdot\log^{\kappa}\left(\frac{\lambda^2(\tau)}{\lambda^2(\sigma)}\xi\right)\cdot a(\sigma)\,d\sigma\right) \tilde{\rho}_{0}(\xi)\,d\xi\\
		& = \int_0^\infty\chi_{\lambda(\tau)\xi^{\frac12} M^{-1}\lesssim 1}^{(j,k)}\left(\int_{\tau_0}^{\tau}\frac{1}{\tilde{\rho}_{0}^{\frac12}(\xi)}\cdot \frac{e^{\pm i\nu\tau\xi^{\frac12}}}{\xi^{\frac12}\cdot\left(\frac{\lambda^2(\tau)}{\lambda^2(\sigma)}\xi\right)^{\frac{k\nu}{2}}}\cdot\log^{\kappa}\left(\frac{\lambda^2(\tau)}{\lambda^2(\sigma)}\xi\right)\cdot a(\sigma)\,d\sigma\right) \tilde{\rho}_{0}(\xi)\,d\xi\\
	\end{align*}
	Thanks to the primary modulation step we have that 
	\begin{align*}
		\tilde{II}_2 = \int_0^\infty\chi_{\lambda(\tau)\xi^{\frac12} M^{-1}\lesssim 1}^{(j,k)}\left(\int_{\tau}^{\infty}\frac{1}{\tilde{\rho}_{0}^{\frac12}(\xi)}\cdot \frac{e^{\pm i\nu\tau\xi^{\frac12}}}{\xi^{\frac12}\cdot\left(\frac{\lambda^2(\tau)}{\lambda^2(\sigma)}\xi\right)^{\frac{k\nu}{2}}}\cdot\log^{\kappa}\left(\frac{\lambda^2(\tau)}{\lambda^2(\sigma)}\xi\right)\cdot a(\sigma)\,d\sigma\right) \tilde{\rho}_{0}(\xi)\,d\xi
	\end{align*}
	provided $k = 1, 2$. For $k = 3$ we do not use this cancellation, but instead estimate directly 
	\begin{align*}
		&\left|\int_{\tau_0}^{\tau}\frac{1}{\tilde{\rho}_{0}^{\frac12}(\xi)}\cdot \frac{1}{\xi^{\frac12}\cdot\left(\frac{\lambda^2(\tau)}{\lambda^2(\sigma)}\xi\right)^{\frac{k\nu}{2}}}\cdot\log^{\kappa}\left(\frac{\lambda^2(\tau)}{\lambda^2(\sigma)}\xi\right)\cdot a(\sigma)\,d\sigma\right|\\
		&\lesssim \tau^{-3-2\nu}\cdot\log^2 \tau\cdot \frac{1}{\tilde{\rho}_{0}^{\frac12}(\xi)}\cdot \xi^{-\frac12-\frac{k\nu}{2}}\cdot\left|\log\xi\right|. 
	\end{align*}
	In each case we split the $\xi$-integral smoothly into the cases $\xi^{\frac12}\lesssim 1$, $\xi^{\frac12}\gtrsim 1$ and perform integration by parts in the latter, taking advantage of the presence of the oscillatory factor $e^{\pm i\nu\tau\xi^{\frac12}}$. We then easily infer the same bound as for the term $II_1$. This concludes the proof of the lemma. 
\end{proof}

In analogy to the preceding lemma, we formulate the following result concerning the effect of the (iterated) transference operator:
\begin{lemma}\label{lem:hcubicmainsingulartransference} Letting $N_M$ be defined as in the preceding lemma, and letting $\delta>0$ a small constant, there is $\tau_0 = \tau_0(\delta)\ll 1$ such that for $\tau\geq \tau_0$ we have the bounds
	\begin{align*} 
		\left|\int_0^\infty U(\tau)\left(\left[\frac{\lambda_{\tau}}{\lambda}\mathcal{K}_0\mathcal{D}_{\tau}U\right]^j\mathcal{F}^{(0)}\left(\mathcal{D} N_M\right)\right)\tilde{\rho}^{(0)}(\xi)\,d\xi\right|\lesssim \delta^j\cdot\tau^{-(4-10\nu)}\cdot \prod_{j=1}^3\left\|\overline{x}_j\right\|_{good_r},\,j\geq 1. 
	\end{align*}
	provided the vanishing conditions \eqref{eq:n=0vanishingrelations} are satisfied for this source term. 
\end{lemma}
\begin{proof} We show that application of the operator $\frac{\lambda_{\tau}}{\lambda}\mathcal{K}_0\mathcal{D}_{\tau}U$ to $\mathcal{F}^{(0)}\left(\mathcal{D} N_M\right)$ leads to an improved bound which can then be reiterated and improved upon by subsequent applications of $\frac{\lambda_{\tau}}{\lambda}\mathcal{K}_0\mathcal{D}_{\tau}U$. Thus to begin with we make the 
	\\
	{\bf{Claim}}: {\it{We have the following decomposition}}
	\begin{equation}\label{eq:transferenceforcubicmainstructure1}\begin{split}
			\frac{\lambda_{\tau}}{\lambda}\mathcal{K}_0\mathcal{D}_{\tau}U\left(\mathcal{F}^{(0)}\left(\mathcal{D} N_M\right)\right) &=e^{\pm i\nu\tau\xi^{\frac12}}G_1^{\pm}(\tau,\xi)\\
			& +\sum_{\pm} \int_{\tau_0}^{\tau}e^{\pm i\left(\nu\tau - 2\nu\sigma\frac{\lambda(\tau)}{\lambda(\sigma)}\right)\xi^{\frac12}}\cdot G_2^{\pm}(\tau,\sigma,\xi)\,d\sigma\\
			& + G_3(\tau,\xi)
	\end{split}\end{equation}
	{\it{where we have the bounds}}
	\begin{equation}\label{eq:transferenceforcubicmainstructure1bounds}\begin{split}
			&\left\|\xi^{-\gamma}\cdot \langle\xi\rangle^{1+2\gamma}\langle\xi\partial_{\xi}\rangle^{\kappa}G_1(\tau,\xi)\right\|_{L^P_{\zeta(\tau,\xi)d\xi}\cap L^M_{d\xi}} + \left\|\xi^{-\gamma}\cdot \langle\xi\rangle^{3+2\gamma}\langle\xi\partial_{\xi}\rangle^{\kappa}G_3(\tau,\xi)\right\|_{L^P_{\zeta(\tau,\xi)d\xi}\cap L^M_{d\xi}} 
			\\&+\int_{\tau_0}^{\tau}\left(\frac{\lambda(\sigma)}{\lambda(\tau)}\right)^{5\nu}\cdot \left\|\xi^{-\frac12-\gamma}\cdot \langle\xi\rangle^{1+2\gamma}\langle\xi\partial_{\xi}\rangle^{\kappa}G_2(\tau,\sigma,\xi)\right\|_{L^P_{\zeta(\tau,\xi)d\xi}\cap L^M_{d\xi}}\,d\sigma\\&
			\lesssim \tau^{-4-\nu}\cdot\log\tau\cdot \prod_{j=1}^3\left\|\overline{x}_j\right\|_{good_r},\,\kappa\in \{0, 1,2,3\}. 
	\end{split}\end{equation}
	and we have introduced the weight function $\zeta(\tau,\xi): = \min\{\langle\lambda^2(\tau)\xi M^{-2}\rangle^{10},\langle\xi^2\rangle\}$. 
	{\it{for suitable $\gamma>0, \infty>P>1$, $P-1\ll1 , M\gg 1$.}}
	\begin{proof}(Claim, sketch) This relies on using the relations \eqref{eq:preciseasympto1forlhighfreq}, \eqref{eq:preciseasympto1forlowfreq} in the high- and low-frequency regimes, respectively. Consider for example the contribution of the first term on the right in \eqref{eq:preciseasympto1forlowfreq}. Then we can write this contribution to \eqref{eq:transferenceforcubicmainstructure1} explicitly as
		\begin{align*}
			\sum_{j=1,2}N_j^{\pm}(\tau, \xi),
		\end{align*}
		where we set 
		\begin{align*}
			&N_1^{\pm}(\tau, \xi): = e^{\pm i\nu\tau\xi^{\frac12}}\cdot \frac{\lambda_{\tau}}{\lambda}\cdot\mathcal{K}_0\left( e^{\pm i\nu\tau\left(\eta^{\frac12}-\xi^{\frac12}\right)}\cdot\int_{\tau_0}^{\tau}\frac{\tilde{\rho}_0\left(\frac{\lambda^2(\tau)}{\lambda^2(\sigma)}\eta\right)\cdot n_{\pm}\left(\sigma, \frac{\lambda^2(\tau)}{\lambda^2(\sigma)}\eta\right)}{\tilde{\rho}_0(\eta)}\,d\sigma\right) = :e^{\pm i\nu\tau\xi^{\frac12}}\cdot G_{11}^{\pm}(\tau,\xi),\\
			&N_2^{\pm}(\tau, \xi): =  \int_{\tau_0}^{\tau}e^{\pm i\left(\nu\tau - 2\nu\sigma\frac{\lambda(\tau)}{\lambda(\sigma)}\right)\xi^{\frac12}}\cdot \frac{\lambda_{\tau}}{\lambda}\cdot\mathcal{K}_0\left(e^{\pm i\left(\nu\tau - 2\nu\sigma\frac{\lambda(\tau)}{\lambda(\sigma)}\right)\left(\eta^{\frac12}-\xi^{\frac12}\right)}\cdot\frac{\tilde{\rho}_0\left(\frac{\lambda^2(\tau)}{\lambda^2(\sigma)}\eta\right)\cdot n_{\pm}\left(\sigma, \frac{\lambda^2(\tau)}{\lambda^2(\sigma)}\eta\right)}{\tilde{\rho}_0(\eta)}\right)\,d\sigma
		\end{align*}
		In these expressions we may assume that $n_{\pm}(\sigma,\xi)$ are supported at $\xi\lesssim 1$. Then use for the first term $N_1^{\pm}(\tau, \xi)$ and $\kappa = 0$ the bounds 
		\begin{align*}
			\left\|\mathcal{K}_0(\chi_{1\gtrsim \eta\gg\xi}f)\right\|_{L^\infty_{d\xi}}\lesssim \left\|\eta^{-\gamma}f\right\|_{L^P_{d\eta}}
		\end{align*}
		provided $P-1$ is sufficiently small depending on $\gamma>0$, and similarly 
		\begin{align*}
			\left\|\mathcal{K}_0(\chi_{1\gtrsim \eta\sim\xi}f)\right\|_{L^P_{d\xi}}\lesssim \left\|\eta^{-\gamma}f\right\|_{L^P_{d\eta}},\,\left\|\langle\xi\rangle^{10}\mathcal{K}_0(\chi_{\eta\ll\min\{\xi, 1\}}f)\right\|_{L^P_{d\xi}}\lesssim \left\|\eta^{-\gamma}f\right\|_{L^P_{d\eta}}
		\end{align*}
		Using that (taking advantage of \eqref{eq:preciseasympto1forlowfreq} )
		\begin{align*}
			\frac{\lambda_{\tau}}{\lambda}\cdot\Big\|\eta^{-\frac12-\gamma}\cdot e^{\pm i\nu\tau(\eta^{\frac12}-\xi^{\frac12})}\cdot\int_{\tau_0}^{\tau}\frac{\tilde{\rho}_0(\frac{\lambda^2(\tau)}{\lambda^2(\sigma)}\eta)\cdot n_{\pm}(\sigma, \frac{\lambda^2(\tau)}{\lambda^2(\sigma)}\eta)}{\tilde{\rho}_0(\eta)}\,d\sigma\Big\|_{L^p_{d\eta}}\lesssim \tau^{-4-\nu}\log \tau
		\end{align*}
		provided $\gamma\leq \frac14$ (say) and $P-1$ sufficiently small, we get the better than needed estimate 
		\begin{align*}
			\left\|\xi^{-\frac12-\gamma}\cdot \langle\xi\rangle^{1+2\gamma}\langle\xi\partial_{\xi}\rangle^{\kappa}G_{11}^{\pm}(\tau,\xi)\right\|_{L^P_{\langle\xi^2\rangle d\xi}}\lesssim \tau^{-4-\nu}\log \tau.
		\end{align*}
		The bounds for the differentiated expression (i. e. applying $(\xi\partial_{\xi})^{\kappa})$ follow similarly by expressing the operator $\mathcal{K}_0$ as integral operator with respect to the variable $\tilde{\eta}: = \eta^{\frac12}-\xi^{\frac12}$ as well as using the assumption \eqref{eq:preciseasympto1forlowfreq}.
		\\
		The contribution of $N_2^{\pm}(\tau, \xi)$ is handled similarly, and setting 
		\begin{align*}
			G_{22}^{\pm}(\tau,\sigma;\xi): = \frac{\lambda_{\tau}}{\lambda}\cdot\mathcal{K}_0\left(e^{\pm i\left(\nu\tau - 2\nu\sigma\frac{\lambda(\tau)}{\lambda(\sigma)}\right)\left(\eta^{\frac12}-\xi^{\frac12}\right)}\cdot\frac{\tilde{\rho}_0\left(\frac{\lambda^2(\tau)}{\lambda^2(\sigma)}\eta\right)\cdot n_{\pm}\left(\sigma, \frac{\lambda^2(\tau)}{\lambda^2(\sigma)}\eta\right)}{\tilde{\rho}_0(\eta)}\right),
		\end{align*}
		we in fact have the better bound 
		\begin{align*}
			\int_{\tau_0}^{\tau}\left\|\xi^{-\gamma}\cdot \langle\xi\rangle^{1+2\gamma}\langle\xi\partial_{\xi}\rangle^{\kappa}G_{22}(\tau,\sigma,\xi)\right\|_{L^P_{\langle\xi^2\rangle d\xi}}\,d\sigma\lesssim \tau^{-4-\nu}\log \tau,\,\kappa = 0, 1, 2, 3. 
		\end{align*}
		We next verify the claim for the contribution of the first two terms on the right in \eqref{eq:preciseasympto1forlhighfreq}. This leads to the contribution of incoming principal singular terms (we only consider the contribution involving the coefficient $a(\sigma)$, the one with $b(\sigma)$ being similar)
		\begin{align*}
			N_3^{\pm}(\tau,\xi): = e^{\pm i\nu\tau\xi^{\frac12}}\cdot  \frac{\lambda_{\tau}}{\lambda}\cdot\mathcal{K}_0\left(\chi_{\lambda(\tau)\eta^{\frac12} M^{-1}\lesssim 1}^{(j,k)} \frac{e^{\pm i\nu\tau\left(\eta^{\frac12}-\xi^{\frac12}\right)}}{\tilde{\rho}_0^{\frac12}(\eta)}\int_{\tau_0}^{\tau}\frac{a(\sigma)}{\left(\frac{\lambda^2(\tau)}{\lambda^2(\sigma)}\eta\right)^{\frac{k\nu}{2}}}\cdot \log^{\kappa}\left(\frac{\lambda^2(\tau)}{\lambda^2(\sigma)}\eta\right)\,d\sigma\right),\,\kappa\in\{0,1\},
		\end{align*}
		where we have generated the errors\footnote{The second error arises by replacing $\tilde{\rho}_0^{\frac12}\left(\frac{\lambda^2(\tau)}{\lambda^2(\sigma)}\eta\right)$ by $c\frac{\lambda^2(\tau)}{\lambda^2(\sigma)}\eta$.} (with $\kappa\in\{0,1\}$)
		\begin{align*}
			&N_4^{\pm}(\tau,\xi): = e^{\pm i\nu\tau\xi^{\frac12}}\cdot  \frac{\lambda_{\tau}}{\lambda}\cdot\mathcal{K}_0\left(\chi_{\lambda(\tau)\eta^{\frac12} M^{-1}\lesssim 1}^{(j,k)}\cdot\frac{e^{\pm i\nu\tau\left(\eta^{\frac12}-\xi^{\frac12}\right)}}{\tilde{\rho}_0^{\frac12}(\eta)}\int_{\tau_0}^{\tau}\chi_{\frac{\lambda^2(\tau)}{\lambda^2(\sigma)}\eta\lesssim 1}\cdot\frac{a(\sigma)}{\left(\frac{\lambda^2(\tau)}{\lambda^2(\sigma)}\eta\right)^{\frac{k\nu}{2}}}\cdot \log^{\kappa}\left(\frac{\lambda^2(\tau)}{\lambda^2(\sigma)}\eta\right)\,d\sigma\right),\\
			&N_5^{\pm}(\tau,\xi): = e^{\pm i\nu\tau\xi^{\frac12}}\cdot  \frac{\lambda_{\tau}}{\lambda}\cdot\mathcal{K}_0\left(\chi_{\lambda(\tau)\eta^{\frac12} M^{-1}\lesssim 1}^{(j,k)}\cdot\frac{e^{\pm i\nu\tau\left(\eta^{\frac12}-\xi^{\frac12}\right)}}{\tilde{\rho}_0^{\frac12}(\eta)}\int_{\tau_0}^{\tau}\frac{O(1)\cdot a(\sigma)}{\left(\frac{\lambda^2(\tau)}{\lambda^2(\sigma)}\eta\right)^{1+\frac{k\nu}{2}}}\cdot \log^{\kappa}\left(\frac{\lambda^2(\tau)}{\lambda^2(\sigma)}\eta\right)\,d\sigma\right),\\
		\end{align*}
		where $O(1)$ is a bounded function of $\tau,\sigma,\eta$ with symbol type behavior with respect to all variables. Here the last two terms can be treated in analogy to $N_1^{\pm}$, since we gain a factor $\left(\frac{\lambda(\sigma)}{\lambda(\tau)}\right)^{\gamma}$ due to the restriction on $\eta$ or the better large frequency decay. On the other hand, for the first term $N_3^{\pm}$ we can exploit the vanishing condition for $a(\sigma)$ to replace the integral $\int_{\tau_0}^{\tau}\ldots $ by the integral $\int_{\tau}^\infty\ldots$. If we then set 
		\begin{align*}
			N_3^{\pm}(\tau,\xi)=:  e^{\pm i\nu\tau\xi^{\frac12}}\cdot G^{\pm}_{31}(\tau,\xi), 
		\end{align*}
		the bound 
		\begin{align*}
			\left\|\xi^{-\gamma}\cdot \langle\xi\rangle^{1+2\gamma}\left\langle\xi\partial_{\xi}\right\rangle^{\kappa}G^{\pm}_{31}(\tau,\xi)\right\|_{L^P_{\zeta(\tau,\xi)d\xi}}\lesssim \tau^{-4-\nu}\log \tau,\,\kappa \in\{0, 1, 2, 3\}
		\end{align*}
		follows. 
		\\
		Next, we consider the outgoing singular contribution which can be written as 
		\begin{align*}
			&N_6^{\pm}(\tau,\xi):\\& = \int_{\tau_0}^{\tau}e^{\pm i\left(\nu\tau - 2\nu\sigma\frac{\lambda(\tau)}{\lambda(\sigma)}\right)\xi^{\frac12}}\cdot  \frac{\lambda_{\tau}}{\lambda}\cdot\mathcal{K}_0\left(\chi_{\lambda(\tau)\eta^{\frac12} M^{-1}\lesssim 1}^{(j,k)} \frac{e^{\pm i\left(\nu\tau- 2\nu\sigma\frac{\lambda(\tau)}{\lambda(\sigma)}\right)\left(\eta^{\frac12}-\xi^{\frac12}\right)}}{\tilde{\rho}_0^{\frac12}(\eta)}\frac{\chi_{\frac{\lambda^2(\tau)}{\lambda^2(\sigma)}\eta\gtrsim 1}a(\sigma)}{\big(\frac{\lambda^2(\tau)}{\lambda^2(\sigma)}\eta\big)^{\frac{k\nu}{2}}}\cdot \log^{\kappa}\left(\frac{\lambda^2(\tau)}{\lambda^2(\sigma)}\eta\right)\,d\sigma\right),\\&\kappa\in\{0,1\},
		\end{align*}
		Interpreting the preceding expression as 
		\begin{align*}
			\sum_{\pm} \int_{\tau_0}^{\tau}e^{i\left(\nu\tau - 2\nu\sigma\frac{\lambda(\tau)}{\lambda(\sigma)}\right)\xi^{\frac12}}\cdot G_{62}^{\pm}(\tau,\sigma,\xi)\,d\sigma,
		\end{align*}
		we have the bound
		\begin{align*}
			\int_{\tau_0}^{\tau}\left(\frac{\lambda(\sigma)}{\lambda(\tau)}\right)^{\gamma}\cdot \left\|\xi^{-\gamma}\cdot \langle\xi\rangle^{1+2\gamma}\langle\xi\partial_{\xi}\rangle^{\kappa}G_{62}(\tau,\sigma,\xi)\right\|_{L^P_{\zeta(\tau,\xi)d\xi}}\,d\sigma\lesssim \tau^{-4-\nu}\log \tau,\,\kappa \in\{0, 1, 2, 3\}.
		\end{align*}
		The claim for the remaining contributions of \eqref{eq:preciseasympto1forlhighfreq}, \eqref{eq:preciseasympto1forlowfreq} follows similarly. 
	\end{proof}
	
	Continuing with the proof of the lemma, we first consider the case $j = 1$. Thus we need to show that 
	\begin{align*}
		\left|\int_0^\infty U(\tau)H(\tau,\xi)\tilde{\rho}_0(\xi)\,d\xi\right|\lesssim \tau^{-4-\nu}\log\tau
	\end{align*}
	where $H(\tau,\xi)$ stands for one of the three terms on the right of \eqref{eq:transferenceforcubicmainstructure1}. Observe that the propagator $U(\tau)$ involves another time integration, which needs to be compensated for by the oscillatory nature of the $\xi$-integral. We treat the three types of terms separately.
	\\
	
	{\it{Contribution of $G_3$.}} The propagator $U(\tau)G_3$ can be rendered explicitly as 
	\begin{align*}
		\int_{\tau_0}^{\tau}\frac{\tilde{\rho}_0^{\frac12}\left(\frac{\lambda^2(\tau)}{\lambda^2(\sigma)}\xi\right)}{\tilde{\rho}_0^{\frac12}(\xi)}\cdot \frac{\sin\left[\lambda(\tau)\xi^{\frac12}\int_{\sigma}^{\tau}\lambda^{-1}(s)\,ds\right]}{\xi^{\frac12}}\cdot G_3\left(\sigma, \frac{\lambda^2(\tau)}{\lambda^2(\sigma)}\xi\right)\,d\sigma. 
	\end{align*}
	Observe that we have 
	\begin{align*}
		\left|\tilde{\rho}_0^{\frac12}\left(\frac{\lambda^2(\tau)}{\lambda^2(\sigma)}\xi\right)\cdot G_3\left(\sigma, \frac{\lambda^2(\tau)}{\lambda^2(\sigma)}\xi\right)\right|\lesssim \left(\frac{\lambda^2(\tau)}{\lambda^2(\sigma)}\xi\right)^{-2-\gamma}\cdot \sigma^{-4-\nu}\log\sigma
	\end{align*}
	provided $\frac{\lambda^2(\tau)}{\lambda^2(\sigma)}\xi\gtrsim 1$. This implies that the integral 
	\begin{align*}
		\int_0^\infty U(\tau) G_3\tilde{\rho}_0(\xi)\,d\xi
	\end{align*}
	converges absolutely. If we then switch the orders of integration if $\sigma, \xi$, and smoothly partition the integration region for $\xi$ for fixed $\sigma$ into $\left(\nu\tau - \nu\sigma\frac{\lambda(\tau)}{\lambda(\sigma)}\right)\xi^{\frac12}\lesssim 1, \left(\nu\tau - \nu\sigma\frac{\lambda(\tau)}{\lambda(\sigma)}\right)\xi^{\frac12}\gtrsim 1$ and perform integration by parts with respect to $\xi^{\frac12}$ in the latter, the estimate 
	\begin{align*}
		\left|\int_0^\infty U(\tau) G_3\tilde{\rho}_0(\xi)\,d\xi\right|\lesssim \tau^{-(4-10\nu)}
	\end{align*}
	easily follows. 
	\\
	
	{\it{Contribution of the middle term in \eqref{eq:transferenceforcubicmainstructure1}.}} For this term there are two potentially problematic time integrations (namely over both $\tau$ and $\sigma$) where we have to avoid losing powers of $\tau$. For the $\sigma$-integral this will have to exploit a faster oscillation of the $\xi^{\frac12}$-phase than for the previous contribution. Replacing $\sigma$ by $\sigma_1$ and $\tau$ in the formula in \eqref{eq:transferenceforcubicmainstructure1} by $\sigma$, application of the propagator $U(\tau)$ results in the expressions
	\begin{equation}\label{eq:middleterm}\begin{split}
			&\int_{\tau_0}^{\tau}\int_{\tau_0}^{\sigma}\frac{\tilde{\rho}_0^{\frac12}\left(\frac{\lambda^2(\tau)}{\lambda^2(\sigma)}\xi\right)}{\tilde{\rho}_0^{\frac12}(\xi)}\cdot \frac{e^{\pm i\left(\nu\tau - 2\nu\sigma_1\frac{\lambda(\tau)}{\lambda(\sigma_1)}\right)\xi^{\frac12}}}{\xi^{\frac12}}\cdot G_2^{\pm}\left(\sigma,\sigma_1,\frac{\lambda^2(\tau)}{\lambda^2(\sigma)}\xi\right)\,d\sigma_1 d\sigma\\
			&\int_{\tau_0}^{\tau}\int_{\tau_0}^{\sigma}\frac{\tilde{\rho}_0^{\frac12}\left(\frac{\lambda^2(\tau)}{\lambda^2(\sigma)}\xi\right)}{\tilde{\rho}_0^{\frac12}(\xi)}\cdot \frac{e^{\pm i\left(\nu\tau -2\nu\sigma\frac{\lambda(\tau)}{\lambda(\sigma)}+ 2\nu\sigma_1\frac{\lambda(\tau)}{\lambda(\sigma_1)}\right)\xi^{\frac12}}}{\xi^{\frac12}}\cdot G_2^{\pm}\left(\sigma,\sigma_1,\frac{\lambda^2(\tau)}{\lambda^2(\sigma)}\xi\right)\,d\sigma_1 d\sigma\\
			& =: \eqref{eq:middleterm}_1 + \eqref{eq:middleterm}_2,
	\end{split}\end{equation}
	which then need to be integrated against $\tilde{\rho}_0(\xi)\,d\xi$. For the first type of expression, we exploit that 
	\[
	\left|\nu\tau - 2\nu\sigma_1\frac{\lambda(\tau)}{\lambda(\sigma_1)}\right|\geq \nu\sigma_1\frac{\lambda(\tau)}{\lambda(\sigma_1)}\geq \nu\tau, 
	\]
	whence also (for $\nu$ small enough depending on $\gamma$)
	\begin{align*}
		\left|\nu\tau - 2\nu\sigma_1\frac{\lambda(\tau)}{\lambda(\sigma_1)}\right|\gtrsim\left(\frac{\lambda(\tau)}{\lambda(\sigma_1)}\right)^{2\gamma}\cdot \sigma.
	\end{align*}
	Distinguishing between the regions 
	\begin{align*}
		\left|\nu\tau - 2\nu\sigma_1\frac{\lambda(\tau)}{\lambda(\sigma_1)}\right|\xi^{\frac12}\gtrsim 1,\,\left|\nu\tau - 2\nu\sigma_1\frac{\lambda(\tau)}{\lambda(\sigma_1)}\right|\xi^{\frac12}\lesssim 1
	\end{align*}
	and performing integration by parts with respect to $\xi^{\frac12}$ in the latter, we easily gain a factor 
	\begin{align*}
		\left(\frac{\lambda(\tau)}{\lambda(\sigma_1)}\right)^{-2\gamma}\cdot \sigma^{-1}\lesssim \left(\frac{\lambda(\tau)}{\lambda(\sigma_1)}\right)^{-\gamma}\cdot \left(\frac{\lambda(\tau)}{\lambda(\sigma)}\right)^{-\gamma}\sigma^{-1}
	\end{align*}
	allowing to bound the double time integral without a loss and repeated integration by parts in the region $\xi\gtrsim 1$ ensures convergence for the $\xi$-integral. 
	\\
	Next, consider the second type of double integral $\eqref{eq:middleterm}_2$ with the more complicated phase function. Observe that if 
	\[
	\sigma\sim\sigma_1,
	\]
	then $\frac{\lambda(\tau)}{\lambda(\sigma)}\sim \frac{\lambda(\tau)}{\lambda(\sigma_1)}\geq \frac{\lambda(\sigma)}{\lambda(\sigma_1)}$, and furthermore that (for $\tau\geq\sigma\geq\sigma_1$)
	\[
	\left|\nu\tau -2\nu\sigma\frac{\lambda(\tau)}{\lambda(\sigma)}+ 2\nu\sigma_1\frac{\lambda(\tau)}{\lambda(\sigma_1)}\right|\geq \nu\tau.
	\]
	If we then use the bounds (for $\kappa = 0, 1, 2, 3$)
	\begin{align*}
		\left\|\frac{\tilde{\rho}_0^{\frac12}(\frac{\lambda^2(\tau)}{\lambda^2(\sigma)}\xi)}{\tilde{\rho}_0^{\frac12}(\xi)}\cdot \xi^{-\frac12}\cdot \langle\xi\partial_{\xi}\rangle^{\kappa}G_2^{\pm}\left(\sigma,\sigma_1,\frac{\lambda^2(\tau)}{\lambda^2(\sigma)}\xi\right)\right\|_{L^{1}_{\tilde{\rho}_0 d\xi}}\lesssim \left(\frac{\lambda(\sigma)}{\lambda(\tau)}\right)^{1-}\cdot \left\|\xi^{-\frac12}\langle\xi\partial_{\xi}\rangle^{\kappa} G_2(\sigma,\sigma_1,\xi)\right\|_{L^{P}_{d\xi}}
	\end{align*}
	as well as
	\begin{align*}
		&\left\|\frac{\tilde{\rho}_0^{\frac12}\left(\frac{\lambda^2(\tau)}{\lambda^2(\sigma)}\xi\right)}{\tilde{\rho}_0^{\frac12}(\xi)}\cdot \xi^{-\frac12}\cdot G_2^{\pm}\left(\sigma,\sigma_1,\frac{\lambda^2(\tau)}{\lambda^2(\sigma)}\xi\right)\right\|_{L^{1}_{\tilde{\rho}_0d\xi}(\xi^{\frac12}\lesssim \tau^{-1})}\lesssim \tau^{-(1-)}\cdot \left\|\langle\xi\rangle G_2^{\pm}(\sigma,\sigma_1,\xi)\right\|_{L^M_{d\xi}},\\
		&\left\|\frac{\tilde{\rho}_0^{\frac12}\left(\frac{\lambda^2(\tau)}{\lambda^2(\sigma)}\xi\right)}{\tilde{\rho}_0^{\frac12}(\xi)}\cdot \xi^{-\frac12}\cdot \langle\xi\rangle^{-2}\langle\xi\partial_{\xi}\rangle^{\kappa}G_2^{\pm}\left(\sigma,\sigma_1,\frac{\lambda^2(\tau)}{\lambda^2(\sigma)}\xi\right)\right\|_{L^{1}_{\tilde{\rho}_0 d\xi}(\xi^{\frac12}\gtrsim \tau^{-1})}\lesssim \tau^{-(1-)}\cdot \left\|\langle\xi\rangle \langle\xi\partial_{\xi}\rangle^{\kappa}G_2^{\pm}(\sigma,\sigma_1,\xi)\right\|_{L^M_{d\xi}},
	\end{align*}
	the second bound in \eqref{eq:middleterm} under the restriction $\sigma\sim\sigma_1$ follows by integrating by parts with respect to $\xi^{\frac12}$ in the region $\xi^{\frac12}\gtrsim \tau^{-1}$ (after using a smooth partition of unity) and interpolating between the preceding two bounds: 
	\begin{align*}
		&\left|\int_0^\infty \eqref{eq:middleterm}_2\tilde{\rho}_0(\xi)\,d\xi\right|\\&\lesssim \tau^{-(1-10\nu)}\cdot \int_{\tau_0}^{\tau}\left(\frac{\lambda(\sigma)}{\lambda(\tau)}\right)^{10\nu}\int_{\tau_0}^{\sigma}\left(\frac{\lambda(\sigma_1)}{\lambda(\sigma)}\right)^{10\nu}\left[\left\|\xi^{-\frac12}\langle\xi\partial_{\xi}\rangle^{\kappa} G_2(\sigma,\sigma_1,\xi)\right\|_{L^{P}_{d\xi}}
		+\left\|\langle\xi\rangle \langle\xi\partial_{\xi}\rangle^{\kappa}G_2^{\pm}(\sigma,\sigma_1,\xi)\right\|_{L^M_{d\xi}}\right]\,d\sigma_1 d\sigma\\
		&\lesssim  \tau^{-(1-10\nu)}\cdot \int_{\tau_0}^{\tau}\left(\frac{\lambda(\sigma)}{\lambda(\tau)}\right)^{5\nu} \sigma^{-4-\nu}\log\tau\,d\sigma\ll \tau^{-(4-10\nu)}. 
	\end{align*}
	It remains to deal with the situation $\sigma\gg\sigma_1$ in $\eqref{eq:middleterm}_2$. Here we have 
	\begin{align*}
		\left|\nu\tau -2\nu\sigma\frac{\lambda(\tau)}{\lambda(\sigma)}+ 2\nu\sigma_1\frac{\lambda(\tau)}{\lambda(\sigma_1)}\right|\gtrsim \sigma_1\frac{\lambda(\tau)}{\lambda(\sigma_1)}, 
	\end{align*}
	and so the required bound follows by the same argument as for the term $\eqref{eq:middleterm}_1$.
	\\
	
	It remains to deal with the cases $j\geq 2$. Here we shall require an inductive step which shows that application of the operator 
	\[
	\frac{\lambda_{\tau}}{\lambda}\mathcal{K}_0\mathcal{D}_{\tau}U 
	\]
	improves the estimate as far as the high frequency regime $\xi\gg 1$ is concerned. In fact, we make the 
	\\
	{\bf{Claim}}: {\it{Assuming that $H(\tau,\xi)$ can be written as sum of three terms as on the right of \eqref{eq:transferenceforcubicmainstructure1}, then}} 
	\[
	\frac{\lambda_{\tau}}{\lambda}\mathcal{K}_0\mathcal{D}_{\tau}U(H)
	\]
	{\it{can be written as a sum of three terms admitting an analogous description but with estimates better estimates, in addition to a fourth term: in fact, assuming \eqref{eq:transferenceforcubicmainstructure1bounds} for the 'input' but with $\prod_{j=1}^3\left\|\overline{x}_j\right\|_{good_r}$ replaced by $A$, we can write}}
	\begin{align*}
		\frac{\lambda_{\tau}}{\lambda}\mathcal{K}_0\mathcal{D}_{\tau}U\left(\mathcal{F}^{(0)}(H)\right) &=e^{\pm i\nu\tau\xi^{\frac12}}H_1^{\pm}(\tau,\xi)\\
		& +\sum_{\pm} \int_{\tau_0}^{\tau}e^{\pm i\left(\nu\tau - 2\nu\sigma\frac{\lambda(\tau)}{\lambda(\sigma)}\right)\xi^{\frac12}}\cdot H_2^{\pm}(\tau,\sigma,\xi)\,d\sigma\\
		& + H_3(\tau,\xi) + \sum_{\pm}\int_0^\infty e^{\pm i\nu(\tau + x)\xi^{\frac12}}\cdot H_4(\tau,\xi;x)\,dx,
	\end{align*}
	{\it{with}}
	\begin{align*}
		&\left\|\xi^{-\gamma}\cdot \langle\xi\rangle^{\frac32+2\gamma}\langle\xi\partial_{\xi}\rangle^{\kappa}H_1(\tau,\xi)\right\|_{L^P_{\zeta(\tau,\xi)d\xi}\cap L^M_{d\xi}} + \left\|\xi^{-\gamma}\cdot \langle\xi\rangle^{\frac72+2\gamma}\langle\xi\partial_{\xi}\rangle^{\kappa}H_3(\tau,\xi)\right\|_{L^P_{\zeta(\tau,\xi)d\xi}\cap L^M_{d\xi}} 
		\\&+\int_{\tau_0}^{\tau}\left(\frac{\lambda(\sigma)}{\lambda(\tau)}\right)^{5\nu}\cdot \left\|\xi^{-\gamma}\cdot \langle\xi\rangle^{\frac32+2\gamma}\langle\xi\partial_{\xi}\rangle^{\kappa}H_2(\tau,\sigma,\xi)\right\|_{L^P_{\zeta(\tau,\xi)d\xi}\cap L^M_{d\xi}}\,d\sigma\\&\hspace{3cm} + \left\|\langle x\rangle^{-\gamma}\xi^{-\gamma}\cdot \langle\xi\rangle^{\frac32+2\gamma}\langle\xi\partial_{\xi}\rangle^{\kappa}H_4(\tau,\xi;x)\right\|_{(L^P_{\zeta(\tau,\xi)d\xi}\cap L^M_{d\xi})L_x^1(0,\infty)}\\&
		\lesssim \tau^{-4-\nu}\cdot\log\tau\cdot A,\,\kappa\in \{0, 1,2,3\}. 
	\end{align*}
	This follows by an argument similar to the one for the preceding claim. Let us explain this for the contribution of the middle term in \eqref{eq:transferenceforcubicmainstructure1}. Applying the propagator $\mathcal{D}_{\tau}U$ to it, we obtain the double integrals 
	\begin{align*}
		&(1)_{\pm}: = \int_{\tau_0}^{\tau}\int_{\tau_0}^{\sigma}\frac{\tilde{\rho}_0^{\frac12}\left(\frac{\lambda^2(\tau)}{\lambda^2(\sigma)}\xi\right)}{\tilde{\rho}_0^{\frac12}(\xi)}\cdot e^{\pm i\left(\nu\tau - 2\nu\sigma_1\frac{\lambda(\tau)}{\lambda(\sigma_1)}\right)\xi^{\frac12}}\cdot G_2^{\pm}\left(\sigma,\sigma_1,\frac{\lambda^2(\tau)}{\lambda^2(\sigma)}\xi\right)\,d\sigma_1 d\sigma,\\
		&(2)_{\pm}: = \int_{\tau_0}^{\tau}\int_{\tau_0}^{\sigma}\frac{\tilde{\rho}_0^{\frac12}\left(\frac{\lambda^2(\tau)}{\lambda^2(\sigma)}\xi\right)}{\tilde{\rho}_0^{\frac12}(\xi)}\cdot e^{\pm i\left(\nu\tau - 2\nu\sigma\frac{\lambda(\tau)}{\lambda(\sigma)} + 2\nu\sigma_1\frac{\lambda(\tau)}{\lambda(\sigma_1)}\right)\xi^{\frac12}}\cdot G_2^{\pm}\left(\sigma,\sigma_1,\frac{\lambda^2(\tau)}{\lambda^2(\sigma)}\xi\right)\,d\sigma_1 d\sigma.
	\end{align*}
	Then we see that we can interpret the first terms in the form 
	\begin{align*}
		(1)_{\pm}: = \int_{\tau_0}^{\tau}e^{\pm i\left(\nu\tau - 2\nu\sigma_1\frac{\lambda(\tau)}{\lambda(\sigma_1)}\right)\xi^{\frac12}}\cdot \tilde{G}_2^{\pm}(\tau,\sigma_1, \xi)\,d\sigma_1,
	\end{align*}
	where we set (for $\tau_0\leq\sigma_1\leq \tau$)
	\begin{align*}
		\tilde{G}_2^{\pm}(\tau,\sigma_1, \xi) = \int_{\sigma_1}^{\tau}\frac{\tilde{\rho}_0^{\frac12}\left(\frac{\lambda^2(\tau)}{\lambda^2(\sigma)}\xi\right)}{\tilde{\rho}_0^{\frac12}(\xi)}\cdot G_2^{\pm}\left(\sigma,\sigma_1,\frac{\lambda^2(\tau)}{\lambda^2(\sigma)}\xi\right)\,d\sigma.
	\end{align*}
	Then we can bound (for $\kappa\in \{0, 1,2,3\}$)
	\begin{align*}
		&\int_{\tau_0}^{\tau}\left(\frac{\lambda(\sigma_1)}{\lambda(\tau)}\right)^{5\nu}\cdot \left\|\xi^{-\gamma}\cdot \langle\xi\rangle^{2\gamma}\langle\xi\partial_{\xi}\rangle^{\kappa}\tilde{G}_2^{\pm}(\tau,\sigma_1, \xi)\right\|_{L^P_{\zeta(\tau,\xi)d\xi}\cap L^M_{d\xi}}\,d\sigma_1\\
		&\lesssim \int_{\tau_0}^{\tau}\left(\frac{\lambda(\sigma)}{\lambda(\tau)}\right)^{5\nu}\cdot \left(\int_{\tau_0}^{\sigma}\left\|\xi^{-\gamma}\cdot \langle\xi\rangle^{1+2\gamma}\langle\xi\partial_{\xi}\rangle^{\kappa}\tilde{G}_2^{\pm}(\tau,\sigma_1, \xi)\right\|_{L^P_{\zeta(\tau,\xi)d\xi}\cap L^M_{d\xi}}\,d\sigma_1\right) d\sigma \lesssim \tau^{-3-\nu}\cdot\log\tau\cdot A.
	\end{align*}
	Arguing in analogy to the previous claim, we then see that we can write 
	\begin{align*}
		\frac{\lambda_{\tau}}{\lambda}\mathcal{K}_0((1)_{\pm}) = \sum_{\pm} \int_{\tau_0}^{\tau}e^{\pm i\left(\nu\tau - 2\nu\sigma_1\frac{\lambda(\tau)}{\lambda(\sigma_1)}\right)\xi^{\frac12}}\cdot H_2^{\pm}(\tau,\sigma_1,\xi)\,d\sigma
	\end{align*}
	where 
	\begin{align*}
		H_2^{\pm}(\tau,\sigma_1,\xi) = \frac{\lambda_{\tau}}{\lambda}\mathcal{K}_0\left(e^{\pm i\left(\nu\tau - 2\nu\sigma_1\frac{\lambda(\tau)}{\lambda(\sigma_1)}\right)\left(\eta^{\frac12} - \xi^{\frac12}\right)}\cdot \tilde{G}_2^{\pm}(\tau,\sigma_1, \eta)\right)
	\end{align*}
	satisfies the bound asserted in the claim. 
	\\
	
	Next, concerning the term $(2)_{\pm}$, we identify its contribution with a term of the type 
	\[
	\sum_{\pm}\int_0^\infty e^{\pm i\nu(\tau + x)\xi^{\frac12}}\cdot H_4(\tau,\xi;x)\,dx.
	\]
	To begin with, set 
	\begin{align*}
		x: = - 2\nu\sigma\frac{\lambda(\tau)}{\lambda(\sigma)} + 2\nu\sigma_1\frac{\lambda(\tau)}{\lambda(\sigma_1)},
	\end{align*}
	which replaces $\sigma_1$. Interpreting $\sigma_1$ as function of $\sigma\geq\tau_0, x\geq 0$, we write 
	\begin{align*}
		(2)_{\pm} =  \int_0^\infty e^{\pm i(\nu\tau + x)\xi^{\frac12}}\cdot \tilde{G}_2^{\pm}(\tau, \xi; x)\, dx,
	\end{align*}
	where we define 
	\begin{align*}
		\tilde{G}_2^{\pm}(\tau, \xi; x): = \int_{\tau_0}^{\tau}\frac{\partial\sigma_1}{\partial x}\cdot \frac{\tilde{\rho}_0^{\frac12}\left(\frac{\lambda^2(\tau)}{\lambda^2(\sigma)}\xi\right)}{\tilde{\rho}_0^{\frac12}(\xi)}\cdot \chi_{x\leq x_*(\tau,\sigma,\tau_0)}\cdot G_2^{\pm}\left(\sigma,\sigma_1(x),\frac{\lambda^2(\tau)}{\lambda^2(\sigma)}\xi\right)\,d\sigma,
	\end{align*}
	and where $x_*(\tau,\sigma,\tau_0) = - 2\nu\sigma\frac{\lambda(\tau)}{\lambda(\sigma)} + 2\nu\tau_0\frac{\lambda(\tau)}{\lambda(\tau_0)}$, the cutoff $\chi$ being the sharp cutoff. Since 
	\[
	x\sim \sigma_1\cdot\frac{\lambda(\tau)}{\lambda(\sigma_1)},\,\sigma_1\ll\sigma, 
	\]
	we have that ($\kappa\in\{0, 1, 2, 3\}$)
	\begin{align*}
		&\left\|\langle x\rangle^{-\gamma}\langle\xi\partial_{\xi}\rangle^{\kappa}\tilde{G}_2^{\pm}(\tau, \xi; x)\right\|_{(L^P_{\zeta(\tau,\xi)d\xi}\cap L^M_{d\xi})L^1_{dx}}\\&\lesssim \int_{\tau_0}^{\tau}\left(\frac{\lambda(\sigma)}{\lambda(\tau)}\right)^{5\nu}\cdot \int_{\tau_0}^{\tau}\left(\frac{\lambda(\sigma_1)}{\lambda(\sigma)}\right)^{5\nu}\cdot \left\|\langle\xi\rangle\langle\xi\partial_{\xi}\rangle^{\kappa}G_2^{\pm}(\sigma,\sigma_1, \xi)\right\|_{L^P_{\zeta(\tau,\xi)d\xi}\cap L^M_{d\xi}}\,d\sigma_1 d\sigma\\
		&\lesssim \tau^{-3-\nu}\log\tau\cdot A. 
	\end{align*}
	To conclude the estimate for this contribution after applying $\frac{\lambda_{\tau}}{\lambda}\mathcal{K}_0$, we then decompose 
	\begin{align*}
		\frac{\lambda_{\tau}}{\lambda}\mathcal{K}_0((2)_{\pm}) &= \int_0^\infty e^{\pm i(\nu\tau + x)\xi^{\frac12}}\cdot \frac{\lambda_{\tau}}{\lambda}\mathcal{K}_0\left(e^{\pm i(\nu\tau + x)\left(\eta^{\frac12} - \xi^{\frac12}\right)}\cdot  \tilde{G}_2^{\pm}(\tau, \xi; x)\right)\,dx\\
		& =:  \int_0^\infty e^{\pm i(\nu\tau + x)\xi^{\frac12}}\cdot H_4^{\pm}(\tau,\xi;x)\,dx
	\end{align*}
	where $H_4$ satisfies the bound stated before. In fact, this follows from the kernel estimates for $\mathcal{K}_0$ in. Proposition 4.2 of \cite{KMiao}. 
	\\
	This concludes the treatment for the contribution of the middle term in \eqref{eq:transferenceforcubicmainstructure1}, and the other ones are treated similarly. To conclude we show that starting with a term of the fourth type 
	\begin{equation}\label{eq:fourthtype}
		\sum_{\pm}\int_0^\infty e^{\pm i(\nu\tau + x)\xi^{\frac12}}\cdot H_4^{\pm}(\tau,\xi;x)\,dx
	\end{equation}
	satisfying the bound stated in the preceding claim, and applying $\frac{\lambda_{\tau}}{\lambda}\mathcal{K}_0\mathcal{D}_{\tau}U$, leads to another such term, which then implies that all the iterates 
	\[
	\left(\left[\frac{\lambda_{\tau}}{\lambda}\mathcal{K}_0\mathcal{D}_{\tau}U\right]^j\mathcal{F}^{(0)}\left(\mathcal{D} N_M\right)\right)
	\]
	can be described in this manner. In fact, applying $\mathcal{D}_{\tau}U$ to \eqref{eq:fourthtype}, we arrive at either of the following terms:
	\begin{align*}
		&(1)_{\pm} = \sum_{\pm}\int_0^\infty \int_{\tau_0}^{\tau} \frac{\tilde{\rho}_0^{\frac12}\left(\frac{\lambda^2(\tau)}{\lambda^2(\sigma)}\xi\right)}{\tilde{\rho}_0^{\frac12}(\xi)}\cdot e^{\pm i\left(\nu\tau + x\cdot\frac{\lambda(\tau)}{\lambda(\sigma)}\right)\xi^{\frac12}}\cdot H_4^{\pm}\left(\tau,\frac{\lambda^2(\tau)}{\lambda^2(\sigma)}\xi;x\right)\,d\sigma dx\\
		&(2)_{\pm} = \sum_{\pm}\int_0^\infty \int_{\tau_0}^{\tau} \frac{\tilde{\rho}_0^{\frac12}\left(\frac{\lambda^2(\tau)}{\lambda^2(\sigma)}\xi\right)}{\tilde{\rho}_0^{\frac12}(\xi)}\cdot e^{\pm i\left(\nu\tau - 2\nu\sigma\frac{\lambda(\tau)}{\lambda(\sigma)} - x\cdot\frac{\lambda(\tau)}{\lambda(\sigma)}\right)\xi^{\frac12}}\cdot H_4^{\pm}\left(\tau,\frac{\lambda^2(\tau)}{\lambda^2(\sigma)}\xi;x\right)\,d\sigma dx.\\
	\end{align*}
	We can reformulate each of these integrals as 
	\begin{align*}
		\sum_{\pm}\int_0^\infty e^{\pm i(\nu\tau + x)\xi^{\frac12}}\cdot \tilde{H}_4^{\pm}(\tau,\xi;x)\,dx,
	\end{align*}
	where for the terms $(1)_{\pm} $ we set  
	\begin{align*}
		\tilde{H}_4^{\pm}(\tau,\xi;x) =  \int_{\tau_0}^{\tau} \frac{\tilde{\rho}_0^{\frac12}\left(\frac{\lambda^2(\tau)}{\lambda^2(\sigma)}\xi\right)}{\tilde{\rho}_0^{\frac12}(\xi)}\cdot \frac{\lambda(\sigma)}{\lambda(\tau)}H_4^{\pm}\left(\tau,\frac{\lambda^2(\tau)}{\lambda^2(\sigma)}\xi;\frac{\lambda(\sigma)}{\lambda(\tau)}x\right)\,d\sigma.
	\end{align*}
	For the more complicated term $(2)_{\pm} $ we can set 
	\begin{align*}
		\tilde{H}_4^{\pm}(\tau,\xi;y) = \int_{\tau_0}^{\tau} \frac{\lambda(\sigma)}{\lambda(\tau)}\cdot \chi_{y\geq 2\nu\sigma\frac{\lambda(\tau)}{\lambda(\sigma)} - 2\nu\tau}\cdot H^{\mp}_4\left(\tau, \frac{\lambda^2(\tau)}{\lambda^2(\sigma)}\xi; \frac{\lambda(\sigma)}{\lambda(\tau)}\left[y + 2\nu\tau - 2\nu\sigma\frac{\lambda(\tau)}{\lambda(\sigma)} \right]\right)\,d\sigma.
	\end{align*}
	Then it is straightforward to verify that we have the bounds (for $\kappa\in \{0, 1, 2, 3\}$)
	\begin{align*}
		\left\|\langle x\rangle^{-\gamma}\xi^{-\gamma}\cdot \langle\xi\rangle^{2\gamma}\langle\xi\partial_{\xi}\rangle^{\kappa}\tilde{H}_4(\tau,\xi;x)\right\|_{(L^P_{\zeta(\tau,\xi)d\xi}\cap L^M_{d\xi})L_x^1(0,\infty)}\lesssim \tau^{3-\nu}\log\tau\cdot A. 
	\end{align*}
	In fact, for the second term $(2)_{\pm} $ we observe that 
	\begin{align*}
		&\langle y\rangle^{-\gamma}\cdot \left|\frac{\lambda(\sigma)}{\lambda(\tau)}\cdot \chi_{y\geq 2\nu\sigma\frac{\lambda(\tau)}{\lambda(\sigma)} - 2\nu\tau}\cdot H^{\mp}_4\left(\tau, \frac{\lambda^2(\tau)}{\lambda^2(\sigma)}\xi; \frac{\lambda(\sigma)}{\lambda(\tau)}\left[y + 2\nu\tau - 2\nu\sigma\frac{\lambda(\tau)}{\lambda(\sigma)} \right]\right)\right|\\
		&\lesssim \left(\frac{\lambda(\sigma)}{\lambda(\tau)}\right)^{\gamma}\cdot \left|\frac{\lambda(\sigma)}{\lambda(\tau)}\cdot \chi_{y\geq 2\nu\sigma\frac{\lambda(\tau)}{\lambda(\sigma)} - 2\nu\tau}\cdot\left(\langle y\rangle^{-\gamma}H^{\mp}_4\right)\left(\tau, \frac{\lambda^2(\tau)}{\lambda^2(\sigma)}\xi; \frac{\lambda(\sigma)}{\lambda(\tau)}\left[y + 2\nu\tau - 2\nu\sigma\frac{\lambda(\tau)}{\lambda(\sigma)} \right]\right)\right|,
	\end{align*}
	and the extra factor $\left(\frac{\lambda(\sigma)}{\lambda(\tau)}\right)^{\gamma}$ allows us to translate the decay of $\sigma^{-4-\nu}\log\sigma$ into $\tau^{-3-\nu}\log \tau$, after integration over $\sigma$ accounting for one power less of decay. Applying the operator $\frac{\lambda_{\tau}}{\lambda}\mathcal{K}_0$ leads to a term of the form \eqref{eq:fourthtype} with the bound stated in the claim, i. e. of the same type.  In fact, we gain $\langle \xi\rangle^{-\frac12}$ decay for large frequencies. 
	\\
	
	Using the customary Volterra re-iteration argument as in Proposition \ref{prop:keyreiteratengeq2}, we infer that 
	\begin{align*}
		\left[\frac{\lambda_{\tau}}{\lambda}\mathcal{K}_0\mathcal{D}_{\tau}U\right]^j\mathcal{F}^{(0)}\left(\mathcal{D} N_M\right)
	\end{align*}
	admits a representation as in the preceding claim but with the bound $\lesssim \tau^{-4-\nu}\log\tau\cdot A$ replaced by 
	\[
	\lesssim \delta^j\cdot\tau^{-4-\nu}\log\tau\cdot A,
	\]
	provided $\tau_0 \geq \tau_0(\nu,\delta)$. The lemma is now a straightforward consequence of this observation and integration by parts with respect to $\xi^{\frac12}$. 
\end{proof}

Lemma~\ref{lem:hcubicmainsingular} together with Lemma~\ref{lem:hcubicmainsingulartransference} complete controlling the contribution of the (regularized) source term $\mathcal{D} N_M$ to $h$. 
\\

{\it{Contribution to $h$ from linear correction term in \eqref{coe eq tau R} as well as the linear in $\varphi$ interactions with $\epsilon$ in $\frak{N}(\varphi_j)$ in \eqref{nonlinearity tau R}.}} We will show that these contributions, as well as all the remaining ones, will be of smaller order than the ones treated before. 
\\

{\it{(1) Contribution of the term $L_1: = \frac{(1+\nu)^2}{\nu^2\tau^2}\cdot\frac{4R^2}{1+2R^2 + R^4}\varphi_2$.}} For this term, we have 
\begin{lemma}\label{lem:linearsource} Denoting by means of the subscript $M$ the regularisation of the term as in the preceding (provided $\varphi_2$ is of admissibly singular type), we have the estimate
	\begin{align*}
		\left|\int_0^\infty U(\tau)\left(\mathcal{F}^{(0)}\left(\mathcal{D} L_{1,M}\right)\right)\tilde{\rho}^{(0)}(\xi)\,d\xi\right|\lesssim \tau^{-4+\nu}\cdot\left\|\bar{x}\right\|_{good_r},
	\end{align*}
	provided $\bar{x}$ is the distorted Fourier transform at level $n = 0$ of the angular momentum $n= 0$ function $\varphi_2$. 
\end{lemma}

\begin{proof}
	Observing that angular momentum $n = 0$ functions satisfy better estimates (see Definition~\ref{defi:xsingulartermsnless2proto}, Definition~\ref{defi:xsingulartermsnless2smooth}), we infer the following analogues of \eqref{eq:preciseasympto1forlowfreq}, \eqref{eq:preciseasympto1forlhighfreq}: 
	\\
	
	{\it{(1): High frequency regime $\xi\gtrsim 1$}}. 
	\begin{equation}\label{eq:preciseasympto2forlhighfreq}\begin{split}
			&\left\langle \mathcal{D} L_{1,M},\,\chi_{R\xi^{\frac12}\gtrsim 1}\phi_0(R;\xi)\right\rangle_{L^2_{R\,dR}}\\&  =\sum_{k = 1}^3 \chi_{\lambda(\tau)\xi^{\frac12} M^{-1}\lesssim 1}^{(1,k)}\sum_{\pm} \log\xi\cdot \frac{e^{\pm i\nu\tau\xi^{\frac12}}}{\xi^{\frac32+k\cdot\frac{\nu}{2}}}\cdot a(\tau) + \sum_{k = 1}^3 \chi_{\lambda(\tau)\xi^{\frac12} M^{-1}\lesssim 1}^{(2,k)}\sum_{\pm} \frac{e^{\pm i\nu\tau\xi^{\frac12}}}{\xi^{\frac32+k\cdot\frac{\nu}{2}}}\cdot b(\tau) + O\left(\log\tau\cdot \tau^{-5+\nu}\cdot \log\xi\cdot\xi^{-2}\right),\\
			&|a(\tau)| + |b(\tau)|\lesssim \tau^{-5+\nu}\cdot\log\tau,
	\end{split}\end{equation}
	where the error term $O(\ldots)$ admits a fine structure analogous to the one of the error term $O(\ldots)$ in \eqref{eq:preciseasympto1forlhighfreq}, and the 'damping terms' $\chi_{\lambda(\tau)\xi^{\frac12} M^{-1}\lesssim 1}^{(j,k)}$ are also defined as before. 
	\\
	
	{\it{(2): Low frequency regime $\xi\lesssim 1$}}. Here we have the representation 
	\begin{equation}\label{eq:preciseasympto2forlowfreq}
		\left\langle \mathcal{D} L_{1,M},\,\phi_0(R;\xi)\right\rangle_{L^2_{R\,dR}} = \sum_{\pm}e^{\pm i\nu\tau\xi^{\frac12}}\cdot n_{\pm}(\tau,\xi) + n_0(\tau,\xi), 
	\end{equation}
	where we have the derivative bounds 
	\begin{align*}
		\left|\partial_{\xi^{\frac12}}^{\kappa}n_j(\tau,\xi)\right|\lesssim \xi^{\frac14-\frac{\kappa}{2}}\cdot |\log\xi|\cdot \tau^{-5-\nu}\log^3\tau, 
	\end{align*}
	provided $\kappa\in \{0, 1\}$ and $j = \pm, 0$.  
	\\
	Applying the propagator $U$ (as before at angular momentum $n = 0$) we find (without having to exploit any vanishing properties) 
	\begin{align*}
		&U\left(\mathcal{F}^{(0)}\left( \mathcal{D} L_{1,M}\right)\right)(\tau,\xi)\\& = \sum_{\pm} \chi_{\lambda(\tau)\xi^{\frac12} M^{-1}\lesssim 1} e^{\pm i\nu\tau\xi^{\frac12}}\cdot g_{\pm}(\tau,\xi) + \int_{\tau_0}^{\tau} \chi_{\lambda(\tau)\xi^{\frac12} M^{-1}\lesssim 1} e^{\pm i\left(\nu\tau - 2\nu\sigma\frac{\lambda(\tau)}{\lambda(\sigma)}\right)\xi^{\frac12}}\cdot h_{\pm}(\tau,\sigma,\xi)\,d\sigma,\\
		& +H(\tau,\xi), 
	\end{align*}
	where we have the bounds (for $\kappa\in \{0, 1, 2, 3\}$)
	\begin{align*}
		&\left|\langle\xi\partial_{\xi}\rangle^{\kappa} g_{\pm}(\tau,\xi)\right|\lesssim \tau^{-4-\nu}\log\tau\cdot \min\{\xi^{-(1-)}, \xi^{-2}\},\\
		&\left|\langle\xi\partial_{\xi}\rangle^{\kappa} \int_{\tau_0}^{\tau}h_{\pm}(\tau,\sigma, \xi)\right|\,d\sigma \lesssim \tau^{-4-\nu}\log\tau\cdot \min\{\xi^{-(1-)}, \xi^{-2}\},\\
		&|H(\tau,\xi)|\lesssim \tau^{-4-\nu}\log\tau\cdot \min\{\xi^{-\frac12}, \langle\xi\rangle^{-3-}\}. 
	\end{align*}
	The lemma is then an easy consequence upon integration over $\xi$ and integration by parts in the large $\xi$-regime where necessary. 
\end{proof}
To conclude the contribution of the source term $L_{1,M}$, we also need to control the terms arising upon repeatedly applying the transference operator. This is done in analogy to Lemma~\ref{lem:hcubicmainsingulartransference} and leads to the bound 
\begin{align*} 
	\left|\int_0^\infty U(\tau)\left(\left[\frac{\lambda_{\tau}}{\lambda}\mathcal{K}_0\mathcal{D}_{\tau}U\right]^j\mathcal{F}^{(0)}\left(\mathcal{D} L_{1,M}\right)\right)\tilde{\rho}^{(0)}(\xi)\,d\xi\right|\lesssim \delta^j\cdot\tau^{-(4-\nu)}\log\tau\cdot \prod_{j=1}^3\left\|\bar{x}\right\|_{good_r}. 
\end{align*}
provided $\tau\geq \tau_0(\delta, \nu)$. We omit the similar details. 
\\

{\it{(2) Contribution of the terms on the third line in \eqref{nonlinearity tau R}.}} We observe that the expression
\begin{align*}
	L_{2,M}: = \chi_{\nu\tau-R\gtrsim M^{-1}\lambda(\tau)}\cdot \left[(\partial_R\epsilon)^2 - \left(\partial_{\tau}\epsilon + \frac{\lambda_{\tau}}{\lambda}R\partial_R\epsilon\right)^2\right]\phi_2
\end{align*}
is a null-form just like the one treated in Lemma~\ref{lem:hcubicmainsingular}. Here the function $\phi_2$ is indeed an angular momentum $n = 0$ function. To treat this contribution, we use 
\begin{lemma}\label{lem:linearinphitrilinearnull} We have the estimate 
	\begin{align*}
		\left|\int_0^\infty U(\tau)\mathcal{F}^{(0)}\left(\mathcal{D}L_{2,M}\right)(\xi)\tilde{\rho}_0(\xi)\,d\xi\right|\lesssim \tau^{-4+\nu}\log\tau\cdot \left\|\bar{x}\right\|_{good_r},
	\end{align*}
	where $\bar{x}$ is the distorted Fourier transform of $\varphi_2$ as an angular momentum $n = 0$ function, provided the vanishing conditions  \eqref{eq:n=0vanishingrelations} are satisfied for this source term(reference to the main manuscript); in particular, this source term is source admissible. 
\end{lemma}
\begin{proof} To begin with, we use Lemma~\ref{lem:singFouriertiphysicaln0adm}
	to infer the representation 
	\begin{align*}
		\varphi_2(\tau, R) = \chi_{|\nu\tau - R|\lesssim 1}\sum_{k=1}^N(\nu\tau - R)^{\frac12+k\nu}\left(\log(\nu\tau - R)\right)^j\cdot \frac{G_k(\tau, \nu\tau - R)}{\tau^{\frac12}} + f(\tau, R),
	\end{align*}
	where we have 
	\begin{equation}\label{eq:linearinphitrilinearnullinput}
		\sum_{j=1}^5 \left((\nu\tau - R)\partial_{\nu\tau - R}\right)^j G_k(\tau, \nu\tau - R) + \left|\langle R\rangle^{-1}f(\tau, R)\right| + \sum_{k=1}^5\left|\partial_R^kf(\tau, R)\right|\lesssim \tau^{-3+\nu}\log\tau. 
	\end{equation}
	If we next expand $\left[(\partial_R\epsilon)^2 - \left(\partial_{\tau}\epsilon + \frac{\lambda_{\tau}}{\lambda}R\partial_R\epsilon\right)^2\right]$ using Theorem \ref{thm:KSTGao precise}, we deduce that 
	\begin{align*}
		&\chi_{\nu\tau-R\gtrsim M^{-1}\lambda(\tau)}\cdot \mathcal{D}\left(\left[(\partial_R\epsilon)^2 - \left(\partial_{\tau}\epsilon + \frac{\lambda_{\tau}}{\lambda}R\partial_R\epsilon\right)^2\right]\varphi_2\right)\\& = 
		\chi_{M^{-1}\lambda(\tau)\lesssim|\nu\tau - R|\lesssim 1}\sum_{k=1}^N(\nu\tau - R)^{-\frac32+k\nu}\left(\log(\nu\tau - R)\right)^j\cdot \frac{\tilde{G}_k(\tau, \nu\tau - R)}{\tau^{\frac12}} + \tilde{f}(\tau, R),
	\end{align*}
	where the functions $\tilde{G}_k, \tilde{f}$ satisfy a bound analogous to \eqref{eq:linearinphitrilinearnullinput} but with $j$ ranging from $1$ to $4$ and the right hand side replaced by $\tau^{-6+\nu}\log\tau$. This implies the following analogue of  \eqref{eq:preciseasympto1forlhighfreq} in the high frequency regime$\xi\gtrsim 1$: 
	\begin{equation}\label{eq:preciseasympto3forlhighfreq}\begin{split}
			&\mathcal{F}^{(0)}\left(\mathcal{D}L_{2,M}\right)(\tau,\xi)\\&
			=\sum_{k = 1}^N \chi_{\lambda(\tau)\xi^{\frac12} M^{-1}\lesssim 1}^{(1,k)}\sum_{\pm} \log\xi\cdot \frac{e^{\pm i\nu\tau\xi^{\frac12}}}{\xi^{1+k\cdot\frac{\nu}{2}}}\cdot a(\tau) + \sum_{k = 1}^N \chi_{\lambda(\tau)\xi^{\frac12} M^{-1}\lesssim 1}^{(2,k)}\sum_{\pm} \frac{e^{\pm i\nu\tau\xi^{\frac12}}}{\xi^{1+k\cdot\frac{\nu}{2}}}\cdot b(\tau) + O\left(\log\tau\cdot \tau^{-6-\nu}\cdot \log\xi\cdot\xi^{-\frac54}\right),\\
			&|a(\tau)| + |b(\tau)|\lesssim \tau^{-6-\nu}\cdot\log\tau,
	\end{split}\end{equation}
	where the error term has a fine structure as specified after \eqref{eq:preciseasympto1forlhighfreq}  but with the improved temporal decay. As for the low frequency regime $\xi\lesssim 1$, we can afford to use a more crude estimate 
	\begin{align*}
		\left\|\chi_{\xi\lesssim 1}\mathcal{F}^{(0)}\left(\mathcal{D}L_{2,M}\right)(\tau,\xi)\right\|_{S_1^{(0)}}\lesssim \tau^{-6-\nu}\log\tau. 
	\end{align*}
	The assertion of the lemma then follows by using the argument in the proof of Lemma~\ref{lem:hcubicmainsingular} for the high frequency contribution and crude estimates for the contribution of the low frequency term, namely the fact that 
	\begin{align*}
		&\left|\int_0^\infty \int_{\tau_0}^{\tau}U(\tau,\sigma;\xi)\cdot f(\sigma, \xi)\tilde{\rho}_0(\xi)\,d\xi\right|\lesssim \tau^{-p+1}\cdot \sup_{\sigma\geq \tau_0}\sigma^{p}\cdot \left\|f(\sigma, \cdot)\right\|_{S_1^{(0)}},\\
		&U(\tau,\sigma;\xi) = \frac{\tilde{\rho}_0^{\frac12}\left(\frac{\lambda^2(\tau)}{\lambda^2(\sigma)}\xi\right)}{\tilde{\rho}_0^{\frac12}(\xi)}\cdot \frac{\sin\left[\lambda(\tau)\xi^{\frac12}\int_{\sigma}^{\tau}\lambda^{-1}(s)\,ds\right]}{\xi^{\frac12}}, 
	\end{align*}
	provided $\nu\leq \nu_*(p)$. 
\end{proof}

In analogy to Lemma~\ref{lem:hcubicmainsingulartransference}
\begin{lemma}\label{lem:L2Mtransference} Letting $L_{2,M}$ be defined as in the preceding lemma, $\nu<\nu_*$, and letting $\delta>0$ a small constant, there is $\tau_0 = \tau_0(\delta,\nu)\ll 1$ such that for $\tau\geq \tau_0$ we have the bounds
	\begin{align*} 
		\left|\int_0^\infty U(\tau)\left(\left[\frac{\lambda_{\tau}}{\lambda}\mathcal{K}_0\mathcal{D}_{\tau}U\right]^j\mathcal{F}^{(0)}\left(\mathcal{D} L_{2,M}\right)\right)\tilde{\rho}^{(0)}(\xi)\,d\xi\right|\lesssim \delta^j\cdot\tau^{-(4-\nu)}\cdot\left\|\overline{x}\right\|_{good_r},\,j\geq 1. 
	\end{align*}
	provided the vanishing conditions \eqref{eq:n=0vanishingrelations} are satisfied for this source term. 
\end{lemma}

To conclude the contributions of the terms on the third line of \eqref{nonlinearity tau R}, we also have the following completely analogous lemma for the remaining terms: 
\begin{lemma}\label{lem:360remainingthirdline} Letting 
	\[
	L_{3,M} = \chi_{\nu\tau-R\gtrsim M^{-1}\lambda(\tau)}\cdot \left(\frac{\partial\epsilon_R}{1+R^2} - \frac{\lambda_{\tau}}{\lambda}\cdot\left(\epsilon_{\tau} + \frac{\lambda_{\tau}}{\lambda}R\partial_R\epsilon\right)\right)\varphi_2
	\]
	where $\varphi_2$ is an angular momentum $n = 0$ function with distorted Fourier transform $\overline{x}$, then the conclusions of the preceding two lemmas obtain with $L_{2,M}$ replaced by $L_{3,M}$. The same conclusion applies to the terms 
	\begin{align*}
		L_{4,M}: = \frac{\sin(2Q+\epsilon)\sin\epsilon}{R^2}\varphi_2.
	\end{align*}
\end{lemma}

{\it{(3) Contribution of the terms contained in $N(\varphi_1), N(\varphi_2)$.}} The main new feature here concerns certain quadratic terms where we have to ensure that their contribution is indeed small enough. 
\\

{\it{Quadratic terms in $N(\varphi_1), N(\varphi_2)$.}} Recalling (3.49) and the fact that $a\left(\Pi_{\Phi^{\perp}}\varphi\right)$ is quadratic in the perturbation, we arrive at the following list of quadratic terms:
\begin{align*}
	&\lambda^{-2}\left(U_t\varphi_{1,t} - U_r\varphi_{1,r} - \frac{\sin U}{r^2}\cos U\varphi_1\right)\varphi,\, -\frac{2\sin U}{R^2}\varphi_{2,\theta}\varphi,\,\lambda^{-2}\left(\sum_{j=1,2}(U_r\varphi_j\varphi_{j,r} - U_t\varphi_j\varphi_{j,t})\right),\\
	&\frac{2\sin U}{R^2}\left(\sum_{j=1,2}\varphi_j\varphi_{j,\theta}\right), 
\end{align*}
where a term $\varphi$ without subscript stands for one of $\varphi_{1,2}$. Throughout we recall that $U = Q + \epsilon$. Labelling the preceding list of terms $I - IV$, we need 
\begin{lemma}\label{lem:ItoIVhcontribution} Letting $F$ stand for any one of $I$ to $IV$, and assuming that the factors $\varphi$ are angular momentum $|n_j|\geq 2$ functions compatible with angular momentum zero output, i. e. $\sum n_j = 0$, then we have 
	\begin{align*}
		\left|\int_{0}^{\infty}U(\tau)\left(\mathcal{F}^{(0)}\left(\chi_{\nu\tau - R\gtrsim M^{-1}\lambda(\tau)}\mathcal{D}F\right)\right)(\tau,\xi)\tilde{\rho}_0(\xi)\,d\xi\right|\lesssim \tau^{-(4-\nu)}\cdot \prod_{j=1}^2\left\|\overline{x}_j\right\|_{good_r}. 
	\end{align*}
	provided the vanishing conditions \eqref{eq:n=0vanishingrelations} are satisfied for the top order singular contribution of the terms involving a derivative falling on $\varphi$. Here $U$ is the Duhamel propagator at angular momentum $n = 0$. 
\end{lemma}
\begin{proof} Let us start with a schematic term of the form
	\begin{align*}
		\chi_{\nu\tau - R\gtrsim M^{-1}\lambda(\tau)}\frac{\sin U}{R^2}\cdot \varphi_1\cdot\varphi_2 =: F_M
	\end{align*}
	where both factors $\varphi_j$ are angular momentum $|n_j|\geq 2$ with good distorted Fourier transform. Since these terms are one degree smoother than the top order singular terms, no vanishing condition shall be required for their contribution. Now we distinguish between different scenarios:
	\\
	
	{\it{(1): both factors have distorted Fourier transform of type $\overline{x}_{1l,smooth}$.}} This case as well as the corresponding high high frequency contributions with Fourier transform of type $\overline{x}_{1h,smooth}$ are delicate due to the very weak temporal decay of order $\tau^{-1}$. 
	Thus we have to take crucial advantage of the fine structure in order to eke out the necessary gains. Let us denote by $\xi_1, \xi_2$ the frequencies in the Fourier representations of $\varphi_1, \varphi_2$. Observe that 
	\begin{align*}
		U(\tau)\left(\mathcal{F}^{(0)}(\mathcal{D}F_M)\right)(\xi) = \int_{\tau_0}^{\tau}\frac{\tilde{\rho}_0^{\frac12}(\frac{\lambda^2(\tau)}{\lambda^2(\sigma)}\xi)}{\tilde{\rho}_0^{\frac12}(\xi)}\cdot \frac{\sin\left[\lambda(\tau)\xi^{\frac12}\int_{\sigma}^{\tau}\lambda^{-1}(s)\,ds\right]}{\xi^{\frac12}}\cdot \left(\mathcal{F}^{(0)}(\mathcal{D}F_M)\right)\left(\frac{\lambda^2(\tau)}{\lambda^2(\sigma)}\xi\right)\,d\sigma,
	\end{align*}
	and distinguish between the following two cases:
	\\
	
	{\it{(1.a): $\frac{\lambda^2(\tau)}{\lambda^2(\sigma)}\xi\lesssim \max\{\xi_1,\xi_2\}$.}} This implies that $\xi\lesssim_{\hbar}1$, and we can exploit that integration over $\xi$ gains $\lesssim \max\{\xi_1^{\frac12}, \xi_2^{\frac12}\}$. Assuming for example that $\xi_1\gtrsim \frac{\lambda^2(\tau)}{\lambda^2(\sigma)}\xi$, we have that\footnote{The subscript indicates the corresponding localization of the frequency in the Fourier representation} 
	\begin{align*}
		&\left|(\partial_R\varphi_1)_{\xi_1\gtrsim \frac{\lambda^2(\tau)}{\lambda^2(\sigma)}\xi}\right|\lesssim R^{-\frac12}\cdot \left|\nu\tau - R\right|^{-1}\cdot \left(\frac{\lambda^2(\tau)}{\lambda^2(\sigma)}\xi\right)^{-\frac14}\cdot \tau^{-(1-\delta)}\cdot \left\|\overline{x}_1\right\|_{good_r}, \\
		&|\varphi_2|\lesssim \tau^{-(1-\delta)}\cdot \left\|\overline{x}_2\right\|_{good_r},
	\end{align*}
	and using that $|\phi(R;\xi)|\lesssim (\langle \log R\rangle + \langle \log\xi\rangle) \cdot \xi^{-\frac14}R^{-\frac12}$ where $\phi$ here denotes the angular momentum $n = 0$ Fourier basis, we easily infer the bound\footnote{The fast decaying factor $\left\langle\frac{\lambda^2(\tau)}{\lambda^2(\sigma)}\xi\right\rangle^{-N}$ can be included after repeated integrations by parts, as the product of low frequency factors is also low frequency up to rapidly decaying tails.}
	\begin{align*}
		\left|\left(\mathcal{F}^{(0)}\left(\chi_{R\sim\tau}\mathcal{D}F_M\right)\right)\left(\frac{\lambda^2(\tau)}{\lambda^2(\sigma)}\xi\right)\right|\lesssim_{\hbar, N} \tau^{-(5-2\delta)}\cdot \left(\frac{\lambda^2(\tau)}{\lambda^2(\sigma)}\xi\right)^{-\frac12}\cdot\left\langle\frac{\lambda^2(\tau)}{\lambda^2(\sigma)}\xi\right\rangle^{-N}\cdot \prod_{j=1,2}\left\|\overline{x}_j\right\|_{good_r}
	\end{align*}
	In order to treat the contribution of the source term in the regime $R\ll\tau$ (but under the same restriction on the frequencies), one exploits that the Fourier basis at angular momentum $|n|\geq 2$ satisfies
	\[
	\left|R^{-1}\phi(R;\xi)\right|\lesssim \hbar\cdot\xi^{\frac12}
	\]
	and performs integrations by parts with respect to $\xi_1$, $\xi_2$ in the Fourier integrals representing $\varphi_{1,2}$, respectively, which leads to the same bound (with an extra logarithmic loss). We then infer that under the current frequency localizations implicitly enforced, we have the bound 
	\begin{align*}
		\int_0^\infty \left| U(\tau)\left(\mathcal{F}^{(0)}(\mathcal{D}F_M)\right)(\xi) \right|\tilde{\rho}_0(\xi)\,d\xi\lesssim \tau^{-4+2\delta}\cdot \prod_{j=1,2}\left\|\overline{x}_j\right\|_{good_r}, 
	\end{align*}
	which is effectively better than what we need. 
	\\
	
	{\it{(1.b): $\frac{\lambda^2(\tau)}{\lambda^2(\sigma)}\xi\gg\max\{\xi_1,\xi_2\}$.}} Here one performs integration by parts with respect to $R$ to trade powers of $frac{\lambda^2(\tau)}{\lambda^2(\sigma)}\xi$ for powers of $\xi_{1,2}$, otherwise one uses the same estimates as in the preceding case. We omit the details. 
	\\
	
	{\it{(2): one factor $\varphi_1$ has distorted Fourier transform of type $\overline{x}_{1,smooth}$, the other $\varphi_2$ of admissibly singular type.}} Relying on Lemma \ref{lem:singFouriertiphysicalngeq2adm} and using Lemma~\ref{lem:structuredlowfreqsmoothbasic}, we can control the contribution of the inner region $R\ll\tau$: 
	\begin{align*}
		\left|\left(\mathcal{F}^{(0)}\left(\chi_{R\ll\tau}\mathcal{D} F_M\right)\right)\left(\frac{\lambda^2(\tau)}{\lambda^2(\sigma)}\xi\right)\right|\lesssim _N\left\langle \left(\frac{\lambda^2(\tau)}{\lambda^2(\sigma)}\xi\right)\right\rangle^{-5}\cdot \tau^{-7+\delta}\cdot \prod_{j=1,2}\left\|\overline{x}_j\right\|_{good_r}. 
	\end{align*}
	We then easily infer a bound which is much better than needed for this contribution: 
	\begin{align*}
		\int_0^\infty \left|U(\tau)\left(\mathcal{F}^{(0)}(\chi_{R\ll\tau}\mathcal{D} F_M)\right)(\xi)\right|\tilde{\rho}_0(\xi)\,d\xi\lesssim \tau^{-6+\delta}\cdot \prod_{j=1,2}\left\|\overline{x}_j\right\|_{good_r}. 
	\end{align*}
	Next, we restrict the source term to the singular region $R\gtrsim \tau$.  Here we again derive an analogue of \eqref{eq:preciseasympto1forlhighfreq}. Specifically, expanding $\varphi_1$ into a Taylor expansion of order $4$ around $R = \nu\tau$, and using Lemma \ref{lem:singFouriertiphysicalngeq2adm}, Lemma \ref{lem:singFouriertiphysicalngeq2}, we can write for $\xi\gtrsim 1$
	\begin{align*}
		&\left(\mathcal{F}^{(0)}(\chi_{R\ll\tau}\mathcal{D} F_M)\right)(\xi)\\& = \sum_{k = 1}^3 \chi_{\lambda(\tau)\xi^{\frac12} M^{-1}\lesssim 1}^{(1,k)}\sum_{\pm} \log\xi\cdot \frac{e^{\pm i\nu\tau\xi^{\frac12}}}{\xi^{\frac32+k\cdot\frac{\nu}{2}}}\cdot a(\tau) + \sum_{k = 1}^3 \chi_{\lambda(\tau)\xi^{\frac12} M^{-1}\lesssim 1}^{(2,k)}\sum_{\pm} \frac{e^{\pm i\nu\tau\xi^{\frac12}}}{\xi^{\frac32+k\cdot\frac{\nu}{2}}}\cdot b(\tau) + O\left(\log\tau\cdot \tau^{-5-\nu+\delta}\cdot \log\xi\cdot\xi^{-2}\right),\\
		&|a(\tau)| + |b(\tau)|\lesssim \tau^{-5-\nu+\delta}\cdot\log\tau,
	\end{align*}
	The desired estimate 
	\begin{align*}
		\int_0^\infty \left|U(\tau)\left(\chi_{\xi\gtrsim 1}\mathcal{F}^{(0)}(\chi_{R\gtrsim\tau}\mathcal{D} F_M)\right)(\xi)\right|\tilde{\rho}_0(\xi)\,d\xi\lesssim  \tau^{-4-\nu+\delta}\cdot\log\tau \cdot \prod_{j=1,2}\left\|\overline{x}_j\right\|_{good_r}
	\end{align*}
	is then a straightforward consequence of the properties of the propagator $U(\tau)$. Note in particular that no cancellation property is required in order to recover the bound, since this term is one degree smoother than the null-form terms. 
	\\
	Furthermore, for the low-frequency contribution corresponding to $\xi\lesssim 1$, it suffices to use the bound 
	\begin{align*}
		\left|\left(\mathcal{F}^{(0)}(\chi_{R\ll\tau}\mathcal{D} F_M)\right)(\tau, \xi)\right|\lesssim \log\tau\cdot \tau^{-5-\nu+\delta}\xi^{-\frac14}|\log\xi|\cdot \prod_{j=1,2}\left\|\overline{x}_j\right\|_{good_r},\,\xi\lesssim 1, 
	\end{align*}
	and this leads directly to the same bound as before 
	\begin{align*}
		\int_0^\infty \left|U(\tau)\left(\chi_{\xi\lesssim 1}\mathcal{F}^{(0)}(\chi_{R\gtrsim\tau}\mathcal{D} F_M)\right)(\xi)\right|\tilde{\rho}_0(\xi)\,d\xi\lesssim  \tau^{-4-\nu+\delta}\cdot\log\tau \cdot \prod_{j=1,2}\left\|\overline{x}_j\right\|_{good_r}
	\end{align*}
	
	{\it{(3): one factor $\varphi_1$ has distorted Fourier transform of type $\overline{x}_{1,smooth}$, the other $\varphi_2$ has distorted Fourier transform of smooth type $\overline{x}_{2,smooth}$.}} Here we use the fact that 
	\begin{align*}
		\left|\langle R\rangle^{-1}\mathcal{D}^{\kappa}\varphi_1(\tau, R)\right|\lesssim \tau^{-2+\delta}\cdot \left\|\overline{x}_1\right\|_{good_r},\,\kappa\in\{0,1\},
	\end{align*}
	see Lemma~\ref{lem:structuredlowfreqsmoothbasic}. Further, we have (with $\kappa\in\{0,1\}$)
	\begin{align*}
		\left|\langle R\rangle^{-\frac12}\mathcal{D}^{\kappa}\varphi_2(\tau, R)\right|\lesssim_{\hbar}\left\|\overline{x}_{2,smooth}\right\|_{S_0^{\hbar}}\lesssim  \tau^{-3+\nu}\cdot \left\|\overline{x}_{2,smooth}\right\|_{good_r}. 
	\end{align*}
	Taking advantage of the bound $\left|\phi(R;\xi)\right|\lesssim \xi^{-\frac14}\cdot |\log\xi|\cdot R^{-\frac12}$ for the $n = 0$ Fourier basis in the low frequency regime $\xi\lesssim 1$, and the smoothness of the source term, we infer the estimate 
	\begin{align*}
		\left\langle \mathcal{F}\left(\mathcal{D}F_M\right), \phi(R;\xi)\right\rangle_{L^2_{R\,dR}}\lesssim \langle\xi\rangle^{-\frac72-\frac{\nu}{2}}\cdot \xi^{-\frac14}|\log\xi|\cdot \tau^{-5+\nu+\delta}\cdot \prod_{j=1,2}\left\|\overline{x}_j\right\|_{good_r}.
	\end{align*}
	From here we easily infer that 
	\begin{align*}
		\int_0^\infty \left|U(\tau)\left(\mathcal{F}^{(0)}(\mathcal{D} F_M)\right)(\xi)\right|\tilde{\rho}_0(\xi)\,d\xi\lesssim  \tau^{-4-\nu+\delta}\cdot \prod_{j=1,2}\left\|\overline{x}_j\right\|_{good_r}
	\end{align*}
	The remaining interactions are handled similarly and omitted here, which concludes the treatment of the terms 
	\[
	\chi_{\nu\tau - R\gtrsim M^{-1}\lambda(\tau)}\frac{\sin U}{R^2}\cdot \varphi_1\cdot\varphi_2.
	\]
	We note that the terms of type II and IV are covered by these estimates. It then remains to deal with the expressions
	\begin{align*}
		&\lambda^{-2}\left(U_t\varphi_{1,t} - U_r\varphi_{1,r}\right)\varphi,\,\lambda^{-2}\left(\sum_{j=1,2}(U_r\varphi_j\varphi_{j,r} - U_t\varphi_j\varphi_{j,t})\right),
	\end{align*}
	These are standard null-forms which can be re-formulated in the form of the term $N$ in Lemma~\ref{lem:hcubicmainsingular} but with either $\varepsilon_2$ or $\varepsilon_3$ replaced by $U = Q+\epsilon$. Substituting $\epsilon$ for $U$, we arrive at an expression for which one easily verifies the relations \eqref{eq:preciseasympto1forlhighfreq}, \eqref{eq:preciseasympto1forlowfreq}, and hence arrives at the desired conclusion. If we substitute $U = Q$, we in addition arrive at the terms
	\begin{equation}\label{eq:delicateremainingquadratic}
		\frac{\lambda_{\tau}}{\lambda}R\partial_RQ\cdot \left(\partial_{\tau} + \frac{\lambda_{\tau}}{\lambda}R\partial_R\right)\varphi_1\cdot \varphi,\,Q_R\cdot \varphi_{1,R}\cdot \varphi,
	\end{equation}
	where we recall that $Q_R\sim \frac{1}{1+R^2}$. We claim that for these terms in the high frequency regime $\xi\gg 1$ we also have the asymptotics  \eqref{eq:preciseasympto1forlhighfreq}. To see this we have to again consider the various possiblities for the 'inputs' $\varphi_1, \varphi$. Consider for example the case that $\varphi_1$ is of admissibly singular type, while $\varphi$ has distorted Fourier transform of type $\overline{x}_{1l,smooth}$. Taking advantage of Lemma 7.10 for $\varphi_{1,R} $, and developing $\varphi$ into a Taylor series around $R = \nu\tau$, we deduce that (with the right most factor $\phi(R;\xi)$ the $n = 0$ Fourier basis)
	\begin{align*}
		\left\langle \chi_{M(\nu\tau - R)\gtrsim \lambda(\tau)}\cdot \mathcal{D}\left(Q_R\cdot \varphi_{1,R}\cdot \varphi\right), \phi(R;\xi)\right\rangle_{L^2_{R\,dR}},\,\xi\gg 1
	\end{align*}
	can be written in the form \eqref{eq:preciseasympto1forlhighfreq}, while for the low frequency regime $\xi\lesssim 1$ the representation \eqref{eq:preciseasympto1forlowfreq} applies.  
	\\
	If both `inputs' $\varphi_1, \varphi$ have distorted Fourier transform of type $\overline{x}_{1l,smooth}$, we take advantage of Lemma~\ref{lem:structuredlowfreqsmoothbasic}, and the argument is a little more delicate. To begin with, in the non-oscillatory regime $R\xi^{\frac12}\lesssim 1$, we use the bound 
	\begin{align*}
		\left|\langle R\rangle^{-1}\varphi_{1,R}\right| + \left|\varphi_{1,RR}\right|\lesssim \tau^{-3}\cdot \left\|\overline{x}_1\right\|_{good_r},\,R\ll\tau, 
	\end{align*}
	which implies the bound 
	\begin{equation}\label{eq:largeRnutua-Rbbounds0}\begin{split}
			&\left|\left\langle\chi_{R\ll\tau}\mathcal{D}\left(Q_R\cdot \varphi_{1,R}\cdot \varphi\right),\,\chi_{R\xi^{\frac12}\lesssim 1}\phi(R;\xi)\right\rangle\right|\lesssim \tau^{-4+\nu}\prod_{j=1,2}\left\|\overline{x}_j\right\|_{good_r},\\
			&\left|(\xi\partial_{\xi})\left\langle\chi_{R\ll\tau}\mathcal{D}\left(Q_R\cdot \varphi_{1,R}\cdot \varphi\right),\,\chi_{R\xi^{\frac12}\lesssim 1}\phi(R;\xi)\right\rangle\right|\lesssim \tau^{-4+\nu}\prod_{j=1,2}\left\|\overline{x}_j\right\|_{good_r},\\
	\end{split}\end{equation}
	the latter bound being a consequence of the fact that $\left|(\xi\partial_{\xi})\chi_{R\xi^{\frac12}\lesssim 1}\phi(R;\xi)\right|\lesssim |\log\xi|$, $\xi\ll1$ and\\ $\left|(\xi\partial_{\xi})\chi_{R\xi^{\frac12}\lesssim 1}\phi(R;\xi)\right|\lesssim 1$, $\xi\gtrsim 1$. For the oscillatory regime $R\xi^{\frac12}\gtrsim 1$, we have 
	the same bounds, where we have to integrate by parts with respect to $R$ once to account for the effect of the operator $\xi\partial_{\xi}$ in the second inequality, and we have to use the asymptotics of $\phi(R;\xi)$ in the oscillatory regime. 
	\\
	Next, in the regime $R\gtrsim \tau$, we use the bounds 
	\begin{equation}\label{eq:largeRnutua-Rbbounds}\begin{split}
			&\left|\chi_{R\gtrsim \tau}\cdot Q_R\cdot \varphi_{1,R}\cdot\varphi\right|\lesssim \tau^{-\frac92}\cdot \left|\nu\tau - R\right|^{-\frac12}\cdot \prod_{j=1,2}\left\|\overline{x}_j\right\|_{good_r},\\
			&\left|\chi_{R\gtrsim \tau}\cdot Q_R\cdot \varphi_{1,RR}\cdot\varphi\right|\lesssim \tau^{-\frac92}\cdot \left|\nu\tau - R\right|^{-\frac32}\cdot \prod_{j=1,2}\left\|\overline{x}_j\right\|_{good_r}.\\
	\end{split}\end{equation}
	The first of these bounds suffices to yield the following estimates:
	\begin{equation}\label{eq::largeRnutua-Rbbounds1}\begin{split}
			&\left|\left\langle \mathcal{D}\left(\chi_{R\gtrsim \tau}\cdot Q_R\cdot \varphi_{1,R}\cdot\varphi\right), \chi_{R\xi^{\frac12}\lesssim1}\phi(R;\xi)\right\rangle_{L^2_{R\,dR}}\right|\lesssim \tau^{-4+\nu}\cdot \prod_{j=1,2}\left\|\overline{x}_j\right\|_{good_r},\\
			&\left|(\xi\partial_{\xi})\left\langle \mathcal{D}\left(\chi_{R\gtrsim \tau}\cdot Q_R\cdot \varphi_{1,R}\cdot\varphi\right), \chi_{R\xi^{\frac12}\lesssim1}\phi(R;\xi)\right\rangle_{L^2_{R\,dR}}\right|\lesssim \tau^{-4+\nu}\cdot \prod_{j=1,2}\left\|\overline{x}_j\right\|_{good_r}.
	\end{split}\end{equation}
	In fact, replacing $\mathcal{D}$ by $\langle R\rangle^{-1}$ this is a direct consequence of the first of the preceding two bounds, while replacing $\mathcal{D}$ by $\partial_R$ this follows by integrating by parts with respect to $R$ and observing that 
	\begin{align*}
		\left|\partial_R\left(\chi_{R\xi^{\frac12}\lesssim1}\phi(R;\xi)\right)\right|\lesssim R^{-1}.
	\end{align*}
	For the regime $R\xi^{\frac12}\gtrsim 1$, we claim that we can write in analogy to \eqref{eq:preciseasympto1forlowfreq}
	\begin{equation}\label{eq::largeRnutua-Rbbounds2}
		\left\langle \mathcal{D}\left(\chi_{R\gtrsim \tau}\cdot Q_R\cdot \varphi_{1,R}\cdot\varphi\right), \chi_{R\xi^{\frac12}\gtrsim 1}\phi(R;\xi)\right\rangle_{L^2_{R\,dR}} = \sum_{\pm}e^{\pm i\nu\tau\xi^{\frac12}}\cdot h_{\pm}(\tau,\xi), 
	\end{equation}
	where we have the bounds 
	\begin{equation}\label{eq::largeRnutua-Rbbounds3}
		\left|(\xi\partial_{\xi})^{\kappa}h_{\pm}(\tau,\xi) \right| \lesssim \tau^{-4+\nu}\cdot \prod_{j=1,2}\left\|\overline{x}_j\right\|_{good_r}, \kappa\in \{0, 1\}. 
	\end{equation}
	For this one uses the asymptotic expansion of $\chi_{R\xi^{\frac12}\gtrsim 1}\phi(R;\xi)$ and replaces the exponent $e^{\pm iR\xi^{\frac12}}$ by $e^{\pm i\nu\tau\xi^{\frac12}}$: write
	\begin{align*}
		\phi(R;\xi) =\sum_{\pm}e^{\pm iR\xi^{\frac12}}\cdot a_{\pm}(R;\xi),  
	\end{align*}
	and so we can expand
	\begin{align*}
		&\left\langle \mathcal{D}\left(\chi_{R\gtrsim \tau}\cdot Q_R\cdot \varphi_{1,R}\cdot\varphi\right), \chi_{R\xi^{\frac12}\gtrsim 1}\phi(R;\xi)\right\rangle_{L^2_{R\,dR}}\\& = \sum_{\pm}e^{\pm i\nu\tau\xi^{\frac12}}\left\langle \mathcal{D}\left(\chi_{R\gtrsim \tau}\cdot Q_R\cdot \varphi_{1,R}\cdot\varphi\right), \chi_{R\xi^{\frac12}\gtrsim 1}a_{\pm}(R;\xi)\cdot \frac{e^{\mp i(\nu\tau-R)\xi^{\frac12}}}{\left(R\xi^{\frac12}\right)^{\frac12}}\right\rangle_{L^2_{R\,dR}},  
	\end{align*}
	where the functions $a_{\pm}(R;\xi)$ obey symbol type bounds with respect to $R, \xi$, as well as the description
	\[
	a_{\pm}(R;\xi) = a_{\pm}(\xi)\left(1+O\left(R^{-1}\xi^{-\frac12}\right)\right),
	\]
	and further 
	\[
	|a(\xi)|\lesssim |\log\xi|,\,\xi\ll 1,\,|a(\xi)|\lesssim \xi^{-1},\,\xi\gg 1. 
	\]
	Writing 
	\begin{align*}
		\sum_{\pm}e^{\pm i\nu\tau\xi^{\frac12}}\langle\ldots\rangle_{\pm} =  \sum_{\pm}e^{\pm i\nu\tau\xi^{\frac12}}\left\langle\chi_{(\nu\tau - R)\xi^{\frac12}\lesssim 1}\ldots\right\rangle_{\pm} +  \sum_{\pm}e^{\pm i\nu\tau\xi^{\frac12}}\left\langle\chi_{(\nu\tau - R)\xi^{\frac12}\gtrsim 1}\ldots\right\rangle_{\pm}
	\end{align*}
	we can use the first bound in \eqref{eq:largeRnutua-Rbbounds} to infer that 
	\begin{align*}
		\left|(\xi\partial_{\xi})^{\kappa}\left\langle\chi_{(\nu\tau - R)\xi^{\frac12}\lesssim 1}\ldots\right\rangle_{\pm}\right|\lesssim \tau^{-4}\cdot \prod_{j=1,2}\left\|\overline{x}_j\right\|_{good_r},\,\kappa\in \{0, 1\}. 
	\end{align*}
	In fact, replacing $\mathcal{D}$ by $R^{-1}$ this bound is a direct consequence of  \eqref{eq:largeRnutua-Rbbounds}, while if we replace $\mathcal{D}$ by $\partial_R$, this is a consequence of integration by parts with respect to $R$, whence moving the $\partial_R$ onto 
	\[
	\chi_{(\nu\tau - R)\xi^{\frac12}\lesssim 1}\cdot \chi_{R\xi^{\frac12}\gtrsim 1}a_{\pm}(R;\xi)\cdot \frac{e^{\mp i(\nu\tau-R)\xi^{\frac12}}}{(R\xi^{\frac12})^{\frac12}}
	\]
	which can then be bounded by 
	\[
	\lesssim |\log\xi|\cdot R^{-\frac12}\cdot \xi^{\frac14}\lesssim |\log\xi| R^{-\frac12}\cdot(\nu\tau - R)^{-\frac12}. 
	\]
	In order to bound 
	\[
	(\xi\partial_{\xi})^{\kappa}\left\langle\chi_{(\nu\tau - R)\xi^{\frac12}\gtrsim 1}\ldots\right\rangle_{\pm},\,\kappa\in \{0, 1\}, 
	\]
	we take advantage of the second bound in \eqref{eq:largeRnutua-Rbbounds}, as well as another integration by parts with respect to $R$ in case $\kappa = 1$, this time in order to gain a factor 
	\[
	(\nu\tau - R)^{-\frac12}\xi^{-\frac12}. 
	\]
	The estimates \eqref{eq:largeRnutua-Rbbounds0}, \eqref{eq::largeRnutua-Rbbounds1}, \eqref{eq::largeRnutua-Rbbounds2},  \eqref{eq::largeRnutua-Rbbounds2} and simple integration by parts with respect to $\xi^{\frac12}$ in the region $\xi^{\frac12}(\nu - \tau)\gtrsim 1$ then imply the desired estimate\footnote{The cutoff $\chi_{\nu\tau - R\gtrsim M^{-1}\lambda(\tau)}$ plays no role here.} 
	\begin{align*}
		\left|\int_{0}^{\infty}U(\tau)\left(\mathcal{F}^{(0)}\left(\chi_{\nu\tau - R\gtrsim M^{-1}\lambda(\tau)}\mathcal{D}F\right)\right)(\tau,\xi)\tilde{\rho}_0(\xi)\,d\xi\right|\lesssim \tau^{-(4-\nu)}\cdot \prod_{j=1}^2\left\|\overline{x}_j\right\|_{good_r},
	\end{align*}
	where $F$ denotes one of \eqref{eq:delicateremainingquadratic}, under the assumption that both $\varphi_1, \varphi$ have distorted Fourier transform of type $\overline{x}_{1l,smooth}$. The remaining possible interactions, depending on the type of $\varphi_1, \varphi$, are handled analogously.  \end{proof}

Arguing as for the proof of Lemma~\ref{lem:hcubicmainsingulartransference}, we can also infer the following lemma covering the effect of the transference operator: 
\begin{lemma}\label{lem:ItoIVhtransferencecontribution} Letting $F$ be as in the previous lemma, and given $\delta>0$, there is $\tau_0 = \tau_0(\delta,\nu)\gg1 $ such that for $\tau\geq\tau_0$ we have 
	\begin{align*}
		\left|\int_{0}^{\infty}U(\tau)\left(\left[\frac{\lambda_{\tau}}{\lambda}\mathcal{K}_0\mathcal{D}_{\tau}U\right]^j\mathcal{F}^{(0)}\left(\chi_{\nu\tau - R\gtrsim M^{-1}\lambda(\tau)}\mathcal{D}F\right)\right)(\tau,\xi)\tilde{\rho}_0(\xi)\,d\xi\right|\lesssim \delta^j\cdot \tau^{-(4-\nu)}\cdot \prod_{j=1}^2\left\|\overline{x}_j\right\|_{good_r},\,j\geq 1. 
	\end{align*}
\end{lemma}

\subsubsection{Other secondary modulations}
Next we consider the errors generated by modulating with $\bar{h}(t)$. Recall the decomposition \eqref{eqLPhismoothcone}-\eqref{Usmoothcone}, we have
\begin{align*}
	\Phi^{(\text{mod}_{2}(t))}=&\calR_{\bar{h}(t)}\left(\Phi_{\text{smooth}}\right)+\Phi_{\text{cone}}\\
	=&\left(\begin{array}{c}
		\sin U_{\text{smooth}}\cos(\theta+\bar{h}(t))\\
		\sin U_{\text{smooth}}\sin(\theta+\bar{h}(t))\\
		\cos U_{\text{smooth}}
	\end{array}\right)+\left(\begin{array}{c}
		\sin U_{\text{cone}}\cos\theta\\
		\sin U_{\text{cone}}\sin\theta\\
		\cos U_{\text{cone}}
	\end{array}\right).
\end{align*}
This modulation generates the terms (omitting the lower order terms),
\begin{align*}
	\lambda^{-2}\left(\begin{array}{c}
		-\sin U_{\text{smooth}}\sin\theta\\
		\sin U_{\text{smooth}}\cos\theta\\
		0
	\end{array}\right)\cdot \bar{h}^{\prime\prime}(t),
\end{align*}
\begin{align*}
	\left(\partial_{R}^{2}\chi_{R\leq \delta\tau}+\frac{1}{R}\partial_{R}\chi_{R\leq\delta\tau}\right)\cdot \epsilon(R)\cdot\left(\left(\begin{array}{c}
		\cos U_{\text{smooth}}\cos(\theta+\bar{h}(t))\\
		\cos U_{\text{smooth}}\sin(\theta+\bar{h}(t))\\
		-\sin U_{\text{smooth}}
	\end{array}\right)-\left(\begin{array}{c}
		\cos U_{\text{smooth}}\cos\theta\\
		\cos U_{\text{smooth}}\sin\theta\\
		-\sin U_{\text{smooth}}
	\end{array}\right) \right),
\end{align*}
and
\begin{align*}
	\partial_{R}\chi_{R\leq\delta\tau}\cdot\epsilon^{\prime}(R)\cdot\left(\left(\begin{array}{c}
		\cos U_{\text{smooth}}\cos(\theta+\bar{h}(t))\\
		\cos U_{\text{smooth}}\sin(\theta+\bar{h}(t))\\
		-\sin U_{\text{smooth}}
	\end{array}\right)-\left(\begin{array}{c}
		\cos U_{\text{smooth}}\cos\theta\\
		\cos U_{\text{smooth}}\sin\theta\\
		-\sin U_{\text{smooth}}
	\end{array}\right) \right) \end{align*}
Only
\begin{align*}
	\lambda^{-2}\left(\begin{array}{c}
		-\sin Q(R)\sin\theta\\
		\sin Q(R)\cos\theta\\
		0
	\end{array}\right)\cdot \bar{h}^{\prime\prime}(t) \end{align*}
contributes directly to $\lim_{R\rightarrow0^{+}}R^{-1}F_{+}(0)$.
Following the arguments in the preceding subsection, one derives the following analogue of Prop.~\ref{prop:determiningc2}: 
\begin{proposition}\label{prop:determiningh2} Given $h:[\tau_0,\infty)\rightarrow \R$ with $\sup_{\tau>\tau_0}\left|\tau^p\cdot h(\tau)\right|<\gamma\ll 1$, and assuming $2<p<p_*(\nu)$,  there is a function $\bar{h}: [\tau_0,\infty)\rightarrow \R$ with 
	\[
	\sup_{\tau>\tau_0}\tau^p\left|\frac{\bar{h}(\tau)}{\tau^2}\right|+\sup_{\tau>\tau_0}\tau^p\left|\frac{\bar{h}'(\tau)}{\tau}\right| + \sup_{\tau>\tau_0}\tau^p\left|\bar{h}''(\tau)\right|\lesssim \gamma, 
	\]
	such that the following holds: let ${\bf{\mathcal{K}}_0}\overline{x}$ be the sum of the last five terms in \eqref{eq:neq0Fourier1}, where $\varepsilon_+^0$ is determined via 
	\begin{align*}
		\varepsilon_+^0(\tau, R) = \phi_0(R)\cdot\int_0^R [\phi_0(s)]^{-1}\cdot \mathcal{D}_0\varepsilon_+^0(\tau, s)\,ds,\,\mathcal{D}_0\varepsilon_+^0(\tau, R) = \int_0^\infty \phi_0(R;\xi)\cdot \overline{x}(\tau,\xi)\cdot \tilde{\rho}_0(\xi)\,d\xi. 
	\end{align*}
	If $N$ is the sum of all the source terms generated by modulating in $\bar{h}$, and if 
	\[
	n(\tau) = \lim_{R\rightarrow 0}R^{-1}N(\tau, R), 
	\]
	then we have the identity
	\begin{align*}
		h(\tau) = n(\tau) + \sum_{j\geq 0}\int_0^\infty U(\tau)\left({\bf{\mathcal{K}}_0} U\right)^j\mathcal{F}^{(0)}\left(\mathcal{D}N\right)\,\tilde{\rho}_0(\xi)\,d\xi
	\end{align*}
\end{proposition}

As for the source term $h$, one can rely here on Prop.~\ref{prop:F+(0)directsource} - Prop.~\ref{lem:ItoIVhtransferencecontribution}.

\subsection{The final estimates for the source terms arising for the exceptional modes, both away from and near the light cone} 

Recalling the tools and methods from section~\ref{subsec:allsourcetermsngeq2} as well as section~\ref{sec:multilinestimatesnearcone}, and the preceding modulation techniques, we can finally formulate the analogues of Prop.~\ref{prop:smoothlinearsource}, Prop.~\ref{prop:bilinwithUregular1}, Prop.~\ref{prop:bilinwithUregular3} as well as Prop.~\ref{prop:ngeq2finalsourcetermestimatesingoodspaces}, and more specifically, the sharpened version of Prop.~\ref{prop:nless23finalsourcetermestimatesingoodspaceswithvanishingcondition} without a stronger condition on $c_0$. This is of course at the expense of having to use secondary modulations.  Recall the ansatz \eqref{eq:modulatedphi}. Then in analogy to \eqref{eq:RegularFinestructure11} (which is in fact valid for all angular momenta), we have modified equations with our new representation of the solution, which take into account the time dependence of the modulation parameters. In case of angular momenta $|n|\geq 2$, where \eqref{eq:RegularFinestructure11} is indeed the relevant equation, let us write the new equations in the final form 
\begin{equation}\label{eq:RegularFinestructurWithmodulation}\begin{split}
	&-\left(\left(\partial_{\tau} + \frac{\lambda_{\tau}}{\lambda}R\partial_R\right)^2 + \frac{\lambda_{\tau}}{\lambda}\left(\partial_{\tau} + \frac{\lambda_{\tau}}{\lambda}R\partial_R\right)\right)\varepsilon_{\pm}(n)+ H_n^{\pm}\varepsilon_{\pm}(n) = F_{\pm}(n)
\end{split}\end{equation}
where the term $F_{\pm}(n)$ takes into account the additional terms generated by modulating (both primary and secondary modulations) and at angular momentum $|n|\geq 2$. For the exceptional modes, we write the corresponding equations in the form 
\begin{equation}\label{eq:varepsnp0plusminus1eqnWithmodulation}\begin{split}
	&-\left(\partial_{\tau}+\frac{\lambda_{\tau}}{\lambda}R\partial_R\right)^2\mathcal{D}_n\varepsilon_+^n - 3\frac{\lambda_{\tau}}{\lambda}\left(\partial_{\tau}+\frac{\lambda_{\tau}}{\lambda}R\partial_R\right)\mathcal{D}_n\varepsilon_+^n + \tilde{H}_n^{+}\mathcal{D}_n\varepsilon_+^n - \left(2\left(\frac{\lambda'}{\lambda}\right)^2 + \left(\frac{\lambda'}{\lambda}\right)'\right)\mathcal{D}_n\varepsilon_+^n\\
	&=:\mathcal{R}_+^n\left(\varepsilon_+^n,\mathcal{D}_n\varepsilon_+^n\right) + \mathcal{D}_n\left(F_+(n)\right),\\
\end{split}\end{equation}
to be complemented by the evolution equations for the unstable modes \eqref{eq:c+evolutioneqn}, \eqref{eq:c0evolutioneqn}, \eqref{eq:c-1evolutioneqn}, where of course $F_+(n)$ is now interpreted to also involve the source terms generated by modulation. Also introduce the notation 
\begin{align*}
-\left(\partial_{\tau}+\frac{\lambda_{\tau}}{\lambda}R\partial_R\right)^2\mathcal{D}_n\varepsilon_+^n - 3\frac{\lambda_{\tau}}{\lambda}\left(\partial_{\tau}+\frac{\lambda_{\tau}}{\lambda}R\partial_R\right)\mathcal{D}_n\varepsilon_+^n + \tilde{H}_{n}^{+}\mathcal{D}_n\varepsilon_+^n - \left(2\left(\frac{\lambda'}{\lambda}\right)^2 + \left(\frac{\lambda'}{\lambda}\right)'\right)\mathcal{D}_n\varepsilon_+^n
= \tilde{\Box}_n\left(\mathcal{D}_n\varepsilon_+^n\right)
\end{align*}
Throughout we now refer to the refined ansatz \eqref{eq:precisemodulation}:
\begin{proposition}\label{prop:exceptionalmodesfinalsourceboundswithmodulation} Let $\Lambda\ll 1$ be defined as in \eqref{eq:LambdaDef}, and $\tau_0\gg 1$. Then there exist primary modulation parameters 
\[
\left(\alpha'(\tau), \beta'(\tau), v_1'(\tau), v_2'(\tau), c'(\tau), h'(\tau)\right)\in C_0^\infty\left([\tau_0, 2\tau_0]; \mathbb{R}\right)
\]
satisfying the bounds stated in  Prop.~\ref{prop:forcedvanishingatn=0}, Prop.~\ref{prop:forcedvanishingatn=plusminus1}, as well as secondary modulation parameters $\bar{c}, \bar{h}\in C^2([\tau_0,\infty)$ satisfying the bounds
\begin{align*}
	\sup_{\tau\geq \tau_0}\sum_{i=0}^2 \tau^{i+2-}\cdot \left|(\bar{c}-1)^{(i)}(\tau)\right| + \sup_{\tau\geq \tau_0}\sum_{i=0}^2 \tau^{i+2-}\cdot \left|\bar{h}^{(i)}(\tau)\right|\lesssim \left(\tau_0^{-1}+\Lambda\right)\cdot \Lambda,
\end{align*}
and such that the following conclusion holds: there exist functions $\xb^{(n)}$, $n = 0,\,\pm 1$, such that setting 
\[
\psi^{+}(n) = \int_0^\infty \xb^{(n)}(\tau,\xi)\phi_{n}(R,\xi)\tilde{\rho}_{n}(\xi)\,d\xi,\,n = 0,\,\pm1,
\]
and letting $F^{\pm}(n)$ be the sum of all the source terms in Prop.~\ref{prop:smoothlinearsource}, Prop.~\ref{prop:bilinwithUregular2}, Prop.~\ref{prop:bilinwithUregular3}, as well as the terms generated by modulating on the above parameters, i.e. \eqref{eq:alphamodulationn=1sourcecontribution}, \eqref{eq:n=1sourcetermduetobetamodulation}, \eqref{eq:n=1contributionfromLorentzmodulationsingular},  \eqref{eq:n=-1contributionfromLorentzmodulationsingular}, \eqref{eq:n=0effectofmodulatingincandh}, \eqref{eq:Qc2modulationterms}, \eqref{eq:Qc2modulationtermsepsilon0}, \eqref{eq:Qc2modulationtermsepsilon1} as well as the analogous terms associated to $\bar{h}$,  and defining $G^{(+)}(n),\,n = 0,\,\pm 1$, in analogy to Prop.~\ref{prop:neq12finalsourcetermestimatesingoodspaces}, Prop.~\ref{prop:neq22finalsourcetermestimatesingoodspaces}, Prop.~\ref{prop:neq-12finalsourcetermestimatesingoodspaces}, the same bounds as in these propositions obtain. Introducing 
\[
{\bf{\xb}}^{(n)}(\tau,\xi): = \int_{\tau_0}^{\tau}U^{(n)}\left(\tau, \sigma,\xi\right)\cdot\mathcal{F}^{(n)}\left(G^{(+)}(n)(\sigma, \cdot)\right)\left(\frac{\lambda^2(\tau)}{\lambda^2(\sigma)}\xi\right)\,d\sigma,\,n\in \{0,\,\pm 1\},
\]
we have the good source bounds 
\begin{align*}
	\left|\int_0^\infty \chi_{\lambda^2(\tau)\xi<M}\xb^{(n)}(\tau, \xi)\tilde{\rho}_{n}(\xi)\,d\xi\right|\lesssim \tau^{-4-}\cdot \left(\tau_0^{-1}+\Lambda\right)\cdot \Lambda.
\end{align*}
Furthermore, the solution to \eqref{eq:c+evolutioneqn}, \eqref{eq:c-1evolutioneqn} satisfy the bounds 
\begin{align*}
	&\left|c_+(\tau)\right| + \tau\cdot\left|c_+'(\tau)\right|\lesssim \tau^{-2-}\cdot\left(\tau_0^{-1}+\Lambda\right)\cdot \Lambda,\\
	&\left|c_-(\tau)\right| + \tau\cdot\left|c_-'(\tau)\right|\lesssim \tau^{-2-}\cdot\left(\tau_0^{-1}+\Lambda\right)\cdot \Lambda,
\end{align*}
while we have $c_0 = 0$. \\
For the Fourier modes $|n|\geq 2$, the conclusions of Prop.~\ref{prop:ngeq2finalsourcetermestimatesingoodspaces} apply. 
\end{proposition}

This proposition results from choosing $\bar{c}, \bar{h}$ in accordance with Prop.~\ref{prop:determiningc2}, Prop.~\ref{prop:determiningh2} with source term $h$ controlled via Prop.~\ref{prop:F+(0)directsource} - Prop.~\ref{lem:ItoIVhtransferencecontribution}, and finally invoking Prop.~\ref{prop:nless23finalsourcetermestimatesingoodspaceswithvanishingcondition}.

Proposition~\ref{prop:exceptionalmodesfinalsourceboundswithmodulation}, finally furnishes us with the tools to control an iterative process 
leading to the desired solution of the initial value problem. 
\section{The iterative process and construction of the solution}\label{sec:final}
\subsection{Basic setup and zeroth iterate}
Let $\Phi = \left(\begin{array}{c}\cos\theta\sin U\\ \sin\theta\sin U\\ \cos U\end{array}\right)$ be the unperturbed co-rotational blow up solution, and at time $t = t_0$, corresponding to $\tau = \tau_0 = \int_{t_0}^\infty \lambda(s)\,ds$ consider a perturbation of its data $\left(\Phi, \Phi_t\right)|_{t = t_0}$, which we write in the form 
\begin{equation}\label{eq:initialperturbation}
\left(\Psi, \Psi_t\right) = \left(\Phi|_{t=t_0} + \sum_{j=1,2}\varphi_j^{ini} E_j|_{t=t_0} + (*)\Phi,\,\sum_{j=1,2}\tilde{\varphi}_j^{ini} E_j|_{t=t_0}\right).
\end{equation}
Here recall that the frame $\{E_1, E_2\}$ is given by $\left\{\left(\begin{array}{c}\cos\theta\sin U\\ \sin\theta\sin U\\ \cos U\end{array}\right),\,\left(\begin{array}{c}-\sin\theta\\ \cos\theta\\ 0\end{array}\right)\right\}$.
We shall make very conservative assumptions on the functions $\varphi_j, \tilde{\varphi}_j$, certainly far from optimal: first, we assume that all $\vphi_j$ are supported on $r<\frac{t_0}{2}$, i.e., away from the light cone\footnote{This assumption is purely technical to simplify the details of the modulation step and can certainly be eliminated.}. Further, we assume the bound
\begin{equation}\label{eq:perturbationsmallness}
\sum_{j=1,2}\left\|\varphi_j^{ini}\right\|_{H^{C_1}} + \sum_{j=1,2}\left\|\tilde{\varphi}_j^{ini}\right\|_{H^{C_1 - 1}}\leq \tau_0^{-C_*},\,C_*, C_1\gg 1.
\end{equation}
Expanding each component into angular modes, we write (using now the re-scaled variable $R = r\lambda(t_0)$)
\begin{align*}
\varphi_j^{ini} = \sum_{n\in \Z}\varphi_j^{ini}(n)(R) e^{in\theta},\quad \tilde{\varphi}_j = \sum_{n\in \Z}\tilde{\varphi}_j^{ini}(n)(R) e^{in\theta}
\end{align*}
and our assumption \eqref{eq:perturbationsmallness} easily implies, passing to the variables(which are each functions of $R$) 
\[
\varepsilon_{\pm}^{ini}(n) = \varphi_1^{ini}(n) \mp i\varphi_2^{ini}(n),\quad \tilde{\varepsilon}_{\pm}^{ini}(n) = \lambda^{-1}(t_0)\left[\tilde{\varphi}_1^{ini}(n) \mp i\tilde{\varphi}_2^{ini}(n)\right], 
\]
and further defining 
\begin{align*}
\tilde{\tilde{\varepsilon}}_{\pm}^{ini}(n) = \tilde{\varepsilon}_{\pm}^{ini}(n) + \frac{\lambda_{\tau}}{\lambda}|_{\tau = \tau_0}\cdot \left(R\partial_R\right)\varepsilon_{\pm}^{ini}(n),
\end{align*}
we invoke the ansatz \eqref{eq:precisemodulation} and consider the family of problems (recall \eqref{eq:RegularFinestructure1})
\begin{equation}\label{eq:RegularFinestructure11again}\begin{split}
	&-\left(\left(\partial_{\tau} + \frac{\lambda_{\tau}}{\lambda}R\partial_R\right)^2 + \frac{\lambda_{\tau}}{\lambda}\left(\partial_{\tau} + \frac{\lambda_{\tau}}{\lambda}R\partial_R\right)\right)\varepsilon_{\pm}(n)+ H_n^{\pm}\varepsilon_{\pm}(n) = F_{\pm}(n) + \tilde{F}_{\pm}(n),\\
	&\varepsilon_{\pm}(n)|_{\tau = \tau_0} = \varepsilon_{\pm}^{ini}(n),\quad \left(\partial_{\tau} + \frac{\lambda_{\tau}}{\lambda}R\partial_R\right)\varepsilon_{\pm}(n)|_{\tau = \tau_0} = \tilde{\tilde{\varepsilon}}_{\pm}^{ini}(n).
\end{split}\end{equation}
The source terms $ F_{\pm}(n)$ in turn are as in \eqref{eq:RegularFinestructure2}, while the term $\tilde{F}_{\pm}(n)$ reflects the additional source terms generated by modulating in the symmetries. These are the equations which we solve for the $|n|\geq 2$ angular momenta, the exceptional angular momenta $n\in \{0,\pm1\}$ requiring passage to the supersymmetric formulation and application of symmetry modulations. 
In order to apply the methods from subsection~\ref{subsec:inhomfouriersolution}, everything has to be converted to the (distorted) Fourier side, for which we use the simple 
\begin{lemma}\label{lem:ngeq2initialdata} Writing 
\[
\varepsilon_{+}(n)(\tau, R) = \int_0^\infty \xb_n(\tau,\xi)\phi_{n}(R,\xi)\rho_{n}(\xi)\,d\xi,\,\varepsilon_{-}(n)(\tau, R) = \overline{\varepsilon_{+}(-n)(\tau, R)},
\]
we have (recall \eqref{eq:Dtaungeq2definition})
\begin{align*}
	\xb_n(\tau_0,\xi) = \left\langle \varepsilon_{+}^{ini}(n), \phi_{n}(R,\xi)\right\rangle_{L^2_{R\,dR}},\quad \mathcal{D}_{\tau} \xb_n(\tau,\xi)|_{\tau = \tau_0} &=  \left\langle \tilde{\tilde{\varepsilon}}_{+}^{ini}(n), \phi_{n}(R,\xi)\right\rangle_{L^2_{R\,dR}} - \frac{\lambda'(\tau)}{\lambda(\tau)}\left(\frac{\rho_n'(\xi)\xi}{\rho_n(\xi)} + 2\right) \xb_n(\tau_0,\xi)\\& - \frac{\lambda'(\tau)}{\lambda(\tau)}\mathcal{K}_{\hbar}^{(0)} \xb_n(\tau_0,\xi)
\end{align*}
Moreover, we have the bound (recall the assumption for Prop.~\ref{prop:ngeq2finalsourcetermestimatesingoodspaces})
\[
\sum_{|n|\geq 2}|n|^{C_2}\left[\left\| \xb_n(\tau_0,\xi)\right\|_{S_0^{(n)}} + \left\| \mathcal{D}_{\tau}x_n(\tau_0,\xi)\right\|_{S_1^{(n)}}\right]\lesssim \tau_0^{-C_*}\ll 1, \,C_2 = C_1 - 10, 
\]
provided $\tau_0\gg 1$. In particular, letting $\varepsilon_{\pm}^{(0)}(n)$ be the solution of the auxiliary linear equation
\begin{align*}
	&-\left(\left(\partial_{\tau} + \frac{\lambda_{\tau}}{\lambda}R\partial_R\right)^2 + \frac{\lambda_{\tau}}{\lambda}\left(\partial_{\tau} + \frac{\lambda_{\tau}}{\lambda}R\partial_R\right)\right)\varepsilon_{\pm}^{(0)}(n)+ H_n^{\pm}\varepsilon_{\pm}^{(0)}(n) = 0,\\
	&\varepsilon_{\pm}(n)|_{\tau = \tau_0} = \varepsilon_{\pm}^{ini}(n),\quad \left(\partial_{\tau} + \frac{\lambda_{\tau}}{\lambda}R\partial_R\right)\varepsilon_{\pm}(n)|_{\tau = \tau_0} = \tilde{\tilde{\varepsilon}}_{\pm}^{ini}(n),
\end{align*}
and writing $\varepsilon_{+}^{(0)}(n)(\tau, R) = \int_0^\infty \xb_n^{(0)}(\tau,\xi)\phi_{n}(R,\xi)\rho_{n}(\xi)\,d\xi$, $\varepsilon_{-}^{(0)}(n)(\tau, R) = \overline{\varepsilon_{+}^{(0)}(-n)(\tau, R)}$,  
we have (always assuming $|n|\geq 2$) the bound 
\[
\sum_{|n|\geq 2}n^{C_2}\left\| \xb_n^{(0)}(\tau,\xi)\right\|_{\text{good}}\lesssim \tau_0^{-C_*+10},
\]
where we keep in mind Definition~\ref{defy:goodfunctionsnorm} for the norm. In fact $\xb_n^{(0)}(\tau,\xi)$ only consists of a smooth part which decays (much faster than) $\tau^{-3}$. 
\end{lemma}

This lemma is a consequence of the Plancherel's theorem for the distorted Fourier transform and the mapping bounds in Prop.~\ref{prop: F hbar boundedness}, as well as of Prop.~\ref{prop: inhomo para decay hbar}.
\\

Dealing with the exceptional modes is more complicated, but we still define the zeroth iterate in the simplest nontrivial possible way: recalling \eqref{eq:varepsnp0plusminus1eqnWithmodulation}, as well as the representation of $\varepsilon_\pm(n)$, $n\in \{0,\pm1\}$, in terms of $\left(c_n, \mathcal{D}_n\varepsilon_\pm(n)\right)$, we define the first iterate, as pair of functions $\left(c_n^{(0)}, \mathcal{D}_n\varepsilon_\pm^{(0)}(n)\right)$, via 
\begin{align*}
\mathcal{D}_n\varepsilon_{-}^{(0)}(n) = \overline{ \mathcal{D}_{-n}\varepsilon_{+}^{(0)}(-n)}
\end{align*}
\begin{equation}\label{eq:varepsnp0plusminus1eqnWithmidulationzerothiterate}\begin{split}
	&-\left(\partial_{\tau}+\frac{\lambda_{\tau}}{\lambda}R\partial_R\right)^2\mathcal{D}_n\varepsilon_+^{(0)}(n) - 3\frac{\lambda_{\tau}}{\lambda}\left(\partial_{\tau}+\frac{\lambda_{\tau}}{\lambda}R\partial_R\right)\mathcal{D}_n\varepsilon_+^{(0)}(n)\\
	& + \tilde{H}_n^{+}\mathcal{D}_n\varepsilon_+^{(0)}(n) - \left(2\left(\frac{\lambda'}{\lambda}\right)^2 + \left(\frac{\lambda'}{\lambda}\right)'\right)\mathcal{D}_n\varepsilon_+^{(0)}(n)=0,\\
	&\mathcal{D}_n\varepsilon_{+}^{(0)}(n)|_{\tau = \tau_0} = \mathcal{D}_n\varepsilon_{+}^{ini}(n),\quad \left(\partial_{\tau} + \frac{\lambda_{\tau}}{\lambda}R\partial_R\right)\mathcal{D}_n\varepsilon_{+}(n)|_{\tau = \tau_0} = \mathcal{D}_n\tilde{\tilde{\varepsilon}}_{+}^{ini}(n) + \frac{\lambda_{\tau}}{\lambda}\left[R\partial_R, \mathcal{D}_n\right]\varepsilon_{+}^{ini}(n)\\
	&\hspace{8.8cm} =:\widetilde{ \varepsilon_{+}^{ini}}(n)\\
	&c_n^{(0)} = 0. 
\end{split}\end{equation}
Invoking the slightly different dilation type operator on the Fourier side 
\[
\mathcal{D}_{\tau} = \mathcal{D}_{\tau}^{(n)} = \partial_{\tau} - 2\frac{\lambda_{\tau}}{\lambda}\xi\partial_{\xi} - \frac{\lambda_{\tau}}{\lambda}\frac{\left(\tilde{\rho}_{n}(\xi)\right)'\xi}{\tilde{\rho}_{n}(\xi)} - \frac{\lambda_{\tau}}{\lambda}, 
\]
we can then formulate the analogue of the preceding lemma for the exceptional modes and the corresponding zeroth iterate. The proof uses the analogous mapping properties of the transference operators at the exceptional angular momenta. Note that at this stage modulation theory does not come in yet, since we force the vanishing of the instabilities: 
\begin{lemma}\label{lem:nless2initialdataandzerothiterate} Let $n\in \{0, \pm1\}$. Writing 
\[
\mathcal{D}_n\varepsilon_{\pm}^{(0)}(n)(\tau, R) = \int_0^\infty \xb_n^{(0)}(\tau,\xi)\phi_n(R,\xi)\tilde{\rho}_{n}(\xi)\,d\xi,
\]
we have 
\begin{align*}
	\xb_n^{(0)}(\tau_0,\xi) = \left\langle \mathcal{D}_n\varepsilon_{\pm}^{ini}(n), \phi_n(R,\xi)\right\rangle_{L^2_{R\,dR}},\quad \mathcal{D}_{\tau} \xb_n^{(0)}(\tau,\xi)|_{\tau = \tau_0} &=  \left\langle \widetilde{\varepsilon_{\pm}^{ini}}(n), \phi_n(R,\xi)\right\rangle_{L^2_{R\,dR}}\\& - \frac{\lambda'(\tau)}{\lambda(\tau)}\left(\frac{\left(\tilde{\rho}_{n}\right)'(\xi)\xi}{\tilde{\rho}_{n}(\xi)} + 1\right) \xb_n(\tau_0,\xi)\\& - \frac{\lambda'(\tau)}{\lambda(\tau)}\mathcal{K}_{n}^{(0)} \xb_n(\tau_0,\xi).
\end{align*}
Moreover, we have the bound (recall subsection~\ref{subsec:nexcspaces} for the corresponding norms)
\begin{align*}
	\left\| \xb_n^{(0)}(\tau_0,\xi)\right\|_{S_0^{(n)}} + \left\|\mathcal{D}_{\tau} \xb_n^{(0)}(\tau_{0},\xi)\right\|_{S_1^{(n)}}\lesssim \tau_0^{-C_*}\ll 1, 
\end{align*}
provided $\tau_0\gg 1$. Moreover, writing the solution of the first part of \eqref{eq:varepsnp0plusminus1eqnWithmidulationzerothiterate} in the form 
\[
\mathcal{D}_n\varepsilon_+^{(0)}(n) = \int_0^\infty \xb_n^{(0)}(\tau,\xi)\phi_n(R,\xi)\tilde{\rho}_{n}(\xi)\,d\xi,
\]
we have the bound (recalling Definition~\ref{defi:goodfourierrepnless2})
\[
\left\|\left(0, \xb_n^{(0)}(\tau,\xi)\right)\right\|_{\text{good}}\lesssim \tau_0^{-C_*+10}. 
\]
In fact $\xb_n^{(0)}(\tau,\xi)$ only consists of a smooth part which decays (much faster than) $\tau^{-3}$. 
\end{lemma}
To complete the information for the zeroth step, we also add the (trivial) information on the modulation parameters, namely we set 
\begin{equation}\label{eq:zerothiteratemodulationparameters}
\left(\alpha^{(0)}(t),\beta^{(0)}(t),\,h^{(0)}(t),\,c^{(0)}(t),\,v_1^{(0)}(t),\,v_2^{(0)}(t),\,\underline{c}^{(0)}(t),\,\underline{h}^{(0)}(t)\right) = \left(0,0,0,1,0,0,1,0\right).
\end{equation}

\subsection{The higher iterates} For the iterative step, consider a tuple of functions 
\begin{equation}\label{eq:inputtuple}
\left(\vphi_1^{(l-1)},\,\vphi_2^{(l-1)},\,\alpha^{(l-1)},\beta^{(l-1)},\,h^{(l-1)},\,(c^{(l-1)}),\,v_1^{(l-1)},\,v_2^{(l-1)},\,\underline{c}^{(l-1)},\,\underline{h}^{(l-1)}\right),\quad l\geq 1, 
\end{equation}
where 
\begin{align*}
\left(\big(\alpha^{(l-1)}\big)',\big(\beta^{(l-1)}\big)',\,\big(h^{(l-1)}\big)',\,\big(c^{(l-1)}\big)',\,\big(v_1^{(l-1)}\big)',\,\big(v_2^{(l-1)}\big)'\right)\in \times_{j=1}^6 C_0^\infty([\tau_0, 2\tau_0]),\,\Big(\underline{c}^{(l-1)},\,\underline{h}^{(l-1)}\Big)\in C^2([\tau_0,\infty))^2.
\end{align*}
\begin{equation}\label{eq:inductionhypmodulationbounds}\begin{split}
	&	\left\|\big(\alpha^{(l-1)}\big)'\right\|_{L^\infty} + \left\|\big(\beta^{(l-1)}\big)'\right\|_{L^\infty}  + \left\|\big(v_1^{(l-1)}\big)'\right\|_{L^\infty}  + \left\|\big(v_2^{(l-1)}\big)'\right\|_{L^\infty}  + \left\|\big(h^{(l-1)}\big)'\right\|_{L^\infty} + \left|\left(c^{(l-1)}\right)'\right\|_{L^\infty}\leq \tau_0^{-C_*+10},\\
	&\sup_{\tau>\tau_0}\tau^p\big|\frac{\bar{c}^{(l-1)}(\tau)-1}{\tau^2}\big|+\sup_{\tau>\tau_0}\tau^p\big|\frac{\big(\bar{c}^{(l-1)}\big)'(\tau)}{\tau}\big| + \sup_{\tau>\tau_0}\tau^p\big|\big(\bar{c}^{(l-1)}\big)''\tau)\Big)\leq \tau_0^{-C_*+10},\,p = 4-\nu,\\
	&\sup_{\tau>\tau_0}\tau^p\big|\frac{\bar{h}^{(l-1)}(\tau)}{\tau^2}\big|+\sup_{\tau>\tau_0}\tau^p\big|\frac{\big(\bar{h}^{(l-1)}\big)'(\tau)}{\tau}\big| + \sup_{\tau>\tau_0}\tau^p\big|\big(\bar{h}^{(l-1)}\big)''(\tau)\Big)\leq \tau_0^{-C_*+10},\,p = 4-\nu.
\end{split}\end{equation}
We shall impose initial conditions on these functions consistent with \eqref{eq:zerothiteratemodulationparameters}.
Furthermore, setting 
\[
\varepsilon_{\pm}^{(l-1)} = \vphi_1^{(l-1)} \mp i\vphi_2^{(l-1)},
\]
we use the decompositions \eqref{eq:varepsilonpmangulardecomp}, i.e., write (using the coordinates $(\tau,R,\theta)$)
\begin{align*}
\varepsilon_{\pm}^{(l-1)} = \sum_{n\in\Z}\varepsilon_{\pm}^{(l-1)}(n)e^{in\theta},\quad \varepsilon_{\pm}^{(l-1)}(n) = \varepsilon_{\pm}^{(l-1)}(n)(\tau,R),
\end{align*}
where for each $|n|\geq 2$ the function $\varepsilon_{\pm}^{(l-1)}(n)$ admits the representation 
\[
\varepsilon_{+}^{(l-1)}(n)(\tau, R) = \int_0^\infty \phi_{n}(R,\xi)\xb_n^{(l-1)}(\tau,\xi)\rho_{n}(\xi)\,d\xi,
\]
with $\xb_n^{(l-1)}(\tau,\xi)$ being restricted singular part in the sense of Definition~\ref{defi:goodfourierrepngeq2}, and with the bound 
\begin{equation}\label{eq:ngeq2inductionhypbounds}
\sum_{|n|\geq 2} n^{C_2}\left\|\xb_n^{(l-1)}(\tau,\xi)\right\|_{\text{good}(r)}\leq \tau_0^{-C_*+10},\quad |n|\geq 2. 
\end{equation}
Here $C_2$ is as in Lemma~\ref{lem:ngeq2initialdata}.
We also assume that the functions $\varepsilon_{\pm}^{(l-1)}(n)$ satisfy the boundary conditions \eqref{eq:RegularFinestructure11again} at time $\tau = \tau_0$. 
\\
We recall that the preceding suffices for a description of all the functions $\varepsilon_{\pm}^{(l-1)}(n)$, $|n|\geq 2$, due to the conjugation symmetry $\varepsilon_{-}^{(l-1)}(n) = \overline{\varepsilon_{+}^{(l-1)}(-n)}$.\\ 
As for the exceptional modes $n\in \{0,\pm1\}$, we write them in the form 
\[
\varepsilon_{\pm}^{(l-1)}(n)(\tau, R) = c_{n,\pm}^{(l-1)}(\tau)\phi_n(R) + \phi_n(R)\int_0^R\left[\phi_n(s)\right]^{-1}\mathcal{D}_n\varepsilon_{\pm}^{(l-1)}(n)(\tau, s)\,ds,\quad n\in\{0,\pm1\}, 
\]
and we impose the initial conditions 
\[
c_{n,\pm}^{(l-1)}(\tau_0) = \lim_{R\rightarrow 0}\left[\phi_n(R)\right]^{-1}\varepsilon_{\pm}^{ini}(n)(R),\quad  \left(c_{n,\pm}^{(l-1)}\right)'(\tau_0) = \lim_{R\rightarrow 0}\left[\phi_n(R)\right]^{-1}\tilde{\varepsilon}_{\pm}^{ini}(n)(R),
\]
as well as the same initial conditions for $\mathcal{D}_n\varepsilon_{\pm}^{(l-1)}(n)$ as those stated for $\mathcal{D}_n\varepsilon_{\pm}^{(0)}(n)$ in \eqref{eq:varepsnp0plusminus1eqnWithmidulationzerothiterate}. Further, writing 
\[
\mathcal{D}_n\varepsilon_{\pm}^{(l-1)}(n)(\tau, R) = \int_0^\infty \xb_n^{(l-1)}(\tau,\xi)\phi_n(R,\xi)\tilde{\rho}_{n}(\xi)\,d\xi, 
\]
we assume that $\left(c_{n,\pm}^{(l-1)}(\tau), \xb_n^{(l-1)}(\tau,\xi)\right)$ is good in the sense of Definition~\ref{defi:goodfourierrepnless2} and with $\xb_n^{(l-1)}(\tau,\xi)$ of restricted principal singular type, and we assume the bound 
\begin{equation}\label{eq:nless2inductionhypbounds}
\left\|\left(c_{n,\pm}^{(l-1)}(\tau), \xb_n^{(l-1)}(\tau,\xi)\right)\right\|_{\text{good}(r)}\leq \tau_0^{-C_*}. 
\end{equation}
Finally, recalling \eqref{eq:RegularFinestructure11again}, let $F_{\pm}^{(l-1)}(n)$ be defined in terms of $\varphi_j^{(l-1)},\,j = 1,2$, using the expressions in \eqref{eq:RegularFinestructure2}, and similarly we let $G_{\pm}(n)^{(l-1)}$ be the errors generated by modulating with respect to the parameters and inputs at stage $l-1$.  We can now formulate the final theorem giving the induction step and convergence to the desired solution: 

\begin{theorem}\label{thm:inductionstepandconvergence} Let $\nu_0>0$ sufficiently small. For $0<\nu\leq \nu_0$, there are positive constants $C_* = C_*(\nu), C_2 = C_2(\nu)$, such that letting $\tau_0 =\tau_0(\nu, C_*, C_2)\gg 1$ sufficiently large, the following holds. Given a tuple \eqref{eq:inputtuple} satisfying the bounds \eqref{eq:inductionhypmodulationbounds}, \eqref{eq:ngeq2inductionhypbounds},\eqref{eq:nless2inductionhypbounds}, there exist 
\begin{equation}\label{eq:inputtuple}
	\left(\vphi_1^{(l)},\,\vphi_2^{(l)},\,\alpha^{(l)},\beta^{(l)},\,h^{(l)},\,(c^{(l)}),\,v_1^{(l)},\,v_2^{(l)},\,\underline{c}^{(l)},\,\underline{h}^{(l)}\right),\quad l\geq 1, 
\end{equation}
where 
\begin{align*}
	\left(\big(\alpha^{(l)}\big)',\big(\beta^{(l)}\big)',\,\big(h^{(l)}\big)',\,\big(c^{(l)}\big)',\,\big(v_1^{(l)}\big)',\,\big(v_2^{(l)}\big)'\right)\in \times_{j=1}^6 C_0^\infty([\tau_0, 2\tau_0]),\,\Big(\underline{c}^{(l)},\,\underline{h}^{(l)}\Big)\in C^2([\tau_0,\infty))^2.
\end{align*}
satisfying the bounds
\begin{equation}\label{eq:inductionhypmodulationbounds}\begin{split}
		&	\left\|\big(\alpha^{(l)}\big)'\right\|_{L^\infty} + \left\|\big(\beta^{(l)}\big)'\right\|_{L^\infty}  + \left\|\big(v_1^{(l)}\big)'\right\|_{L^\infty}  + \left\|\big(v_2^{(l)}\big)'\right\|_{L^\infty}  + \left\|\big(h^{(l)}\big)'\right\|_{L^\infty} + \left|\left(c^{(l)}\right)'\right\|_{L^\infty}\leq \tau_0^{-C_*+10},\\
		&\sup_{\tau>\tau_0}\tau^p\big|\frac{\bar{c}^{(l)}(\tau)-1}{\tau^2}\big|+\sup_{\tau>\tau_0}\tau^p\big|\frac{\big(\bar{c}^{(l)}\big)'(\tau)}{\tau}\big| + \sup_{\tau>\tau_0}\tau^p\big|\big(\bar{c}^{(l)}\big)''\tau)\Big)\leq \tau_0^{-C_*+10},\,p = 4-\nu,\\
		&\sup_{\tau>\tau_0}\tau^p\big|\frac{\bar{h}^{(l)}(\tau)}{\tau^2}\big|+\sup_{\tau>\tau_0}\tau^p\big|\frac{\big(\bar{h}^{(l)}\big)'(\tau)}{\tau}\big| + \sup_{\tau>\tau_0}\tau^p\big|\big(\bar{h}^{(l)}\big)''(\tau)\Big)\leq \tau_0^{-C_*+10},\,p = 4-\nu,
\end{split}\end{equation}
and such that for all $|n|\geq 2$ the equations (where $\tilde{F}_{\pm}^{(l)}$ denotes the additional source terms generated by the primary and secondary modulations)
\begin{equation}\label{eq:RegularFinestructure111}\begin{split}
		&-\left(\left(\partial_{\tau} + \frac{\lambda_{\tau}}{\lambda}R\partial_R\right)^2 + \frac{\lambda_{\tau}}{\lambda}\left(\partial_{\tau} + \frac{\lambda_{\tau}}{\lambda}R\partial_R\right)\right)\varepsilon_{\pm}^{(l)}(n)+ H_n^{\pm}\varepsilon_{\pm}^{(l)}(n) = F_{\pm}^{(l-1)}(n) + \tilde{F}_{\pm}^{(l)}(n),\\
		&\varepsilon_{\pm}^{(l)}(n)|_{\tau = \tau_0} = \varepsilon_{\pm}^{ini}(n),\quad \left(\partial_{\tau} + \frac{\lambda_{\tau}}{\lambda}R\partial_R\right)\varepsilon_{\pm}^{(l)}(n)|_{\tau = \tau_0} = \tilde{\tilde{\varepsilon}}_{\pm}^{ini}(n).
\end{split}\end{equation}
admit solutions on the light cone $R<\nu\tau$ (meaning global solutions upon choice of suitable extensions of $\left(F_{\pm}^{(l-1)}(n) + \tilde{F}_{\pm}^{(l)}(n)\right)\big|_{R<\nu\tau}$ beyond the light cone), which, using the representations
\[
\varepsilon_{\pm}^{(l)}(n)(\tau, R) = \int_0^\infty \xb_n^{(l)}(\tau,\xi)\phi_{n}(R,\xi)\rho_{n}(\xi)\,d\xi,
\]
satisfy the bounds 
\begin{equation}\label{eq:ngeq2inductionhypbounds1}
	\sum_{|n|\geq 2} n^{C_2}\left\|\xb_n^{(l)}(\tau,\xi)\right\|_{\text{good}(r)}\leq \tau_0^{-C_*+10}. 
\end{equation}
Moreover, for the exceptional modes $n\in\{0,\pm 1\}$, the equations 
\begin{equation}\label{eq:varepsnp0plusminus1eqnWithmidulation}\begin{split}
		&-\left(\partial_{\tau}+\frac{\lambda_{\tau}}{\lambda}R\partial_R\right)^2\mathcal{D}_n\varepsilon_+^{(l)}(n) - 3\frac{\lambda_{\tau}}{\lambda}\left(\partial_{\tau}+\frac{\lambda_{\tau}}{\lambda}R\partial_R\right)\mathcal{D}_n\varepsilon_+^{(l)}(n) + \tilde{H}_n^{+}\mathcal{D}_n\varepsilon_+^{(l)}(n) - \left(2\left(\frac{\lambda'}{\lambda}\right)^2 + \left(\frac{\lambda'}{\lambda}\right)'\right)\mathcal{D}_n\varepsilon_+^{(l)}(n)\\
		&=:\mathcal{R}_+^n\left(\varepsilon_+^{(l-1)}(n),\mathcal{D}_n\varepsilon_+^{(l-1)}(n)\right) + \mathcal{D}_n\left(F_+^{(l-1)}(n)\right) + \mathcal{D}_n\left(G_+^{(l-1)}(n)\right),\\
		&\mathcal{D}_n\varepsilon_{\pm}^{(l)}(n)|_{\tau = \tau_0} = \mathcal{D}_n\varepsilon_{\pm}^{ini}(n),\quad \left(\partial_{\tau} + \frac{\lambda_{\tau}}{\lambda}R\partial_R\right)\mathcal{D}_n\varepsilon_{\pm}^{(l)}(n)|_{\tau = \tau_0} = \mathcal{D}_n\tilde{\tilde{\varepsilon}}_{\pm}^{ini}(n) + \frac{\lambda_{\tau}}{\lambda}\left[R\partial_R, \mathcal{D}_n\right]\varepsilon_{\pm}^{ini}(n)\\
		&\hspace{8.8cm} =\widetilde{ \varepsilon_{\pm}^{ini}}(n)\\
\end{split}\end{equation}
with suitable extensions of the right hand sides outside the light cone, 
in conjunction with the evolution for the unstable modes given by 
\begin{align*}
	&-\left(\partial_{\tau} + 2\frac{\lambda'}{\lambda}\right)^2c_-^{(l)} -  \frac{\lambda_{\tau}}{\lambda}\left(\partial_{\tau}+ 2\frac{\lambda'}{\lambda}\right)c_-^{(l)} + \lim_{R\rightarrow 0}R^{-2}H_{-1}^{+}\varepsilon_+^{(l)}(-1) =  \lim_{R\rightarrow 0}R^{-2}\left[F_+(-1) + \tilde{F}_+(-1)\right]\\
	& c_{-}^{(l)}(\tau_0) = \lim_{R\rightarrow 0}\left[\phi_{-1}(R)\right]^{-1}\varepsilon_{+}^{ini}(-1)(R),\quad  \left(c_{-}^{(l)}\right)'(\tau_0) = \lim_{R\rightarrow 0}\left[\phi_{-1}(R)\right]^{-1}\tilde{\varepsilon}_{+}^{ini}(-1)(R)\\
	&-\left(\partial_{\tau} + \frac{\lambda'}{\lambda}\right)^2c_0^{(l)} -  \frac{\lambda_{\tau}}{\lambda}\left(\partial_{\tau}+ \frac{\lambda'}{\lambda}\right)c_0^{(l)} + \lim_{R\rightarrow 0}R^{-2}H_{0}^{+}\varepsilon_+^{(l)}(0) =  \lim_{R\rightarrow 0}R^{-1}\left[F_+(0) + \tilde{F}_+(0)\right]\\
	& c_{0}^{(l)}(\tau_0) = \lim_{R\rightarrow 0}\left[\phi_{0}(R)\right]^{-1}\varepsilon_{+}^{ini}(0)(R),\quad \left(c_{0}^{(l)}\right)'(\tau_0) = \lim_{R\rightarrow 0}\left[\phi_{0}(R)\right]^{-1}\tilde{\varepsilon}_{+}^{ini}(0)(R)\\
	&-\left(\partial_{\tau} \right)^2c_+^{(l)} -  \frac{\lambda_{\tau}}{\lambda}\left(\partial_{\tau}\right)c_+^{(l)} + \lim_{R\rightarrow 0}R^{-2}H_{1}^{+}\varepsilon_+^{(l)}(0) =  \lim_{R\rightarrow 0}\left[F_+(1) + \tilde{F}_+(1)\right]\\
	& c_{+}^{(l)}(\tau_0) = \lim_{R\rightarrow 0}\left[\phi_{1}(R)\right]^{-1}\varepsilon_{+}^{ini}(1)(R),\quad \left(c_{+}^{(l)}\right)'(\tau_0) = \lim_{R\rightarrow 0}\left(\phi_{1}(R)\right)^{-1}\tilde{\varepsilon}_{+}^{ini}(1)(R)\\
\end{align*}
admit solutions on $(\tau_0,\infty)$, such that 
\begin{equation}\label{eq:ngeq2inductionhypbounds2}
	\sum_{|n|<2}\left\|\left(c_n^{(l)}, \xb_n^{(l)}(\tau,\xi)\right)\right\|_{\text{good}(r)}\leq \tau_0^{-C_*+10}. 
\end{equation}
In fact, we have $c_0^{(l)} = 0$. 
\\
Defining $\Lambda^{(l)}$ to be the sum of 
\begin{align*}
	\sum_{n\in\Z,\,|n|\geq 2}\left\langle n\right\rangle^{C_2}\left\| \xb_n^{(l)}\right\|_{\text{good}(r)} + \sum_{n\in \{0,\pm1\}}\left\|\left(c_n^{(l)}(\tau),\,\xb_n^{(l)}(\tau,\xi)\right)\right\|_{\text{good}(r)}
\end{align*}
and the norms on the left  in \eqref{eq:inductionhypmodulationbounds}, we have the estimate 
\begin{align*}
	\Lambda^{(l)}\leq \tau_0^{-C_*} + \Lambda^{(l-1)}\cdot \big(\tau_0^{-1} + \Lambda^{(l-1)}\big).
\end{align*}
Defining the zeroth iterate as in Lemma~\ref{lem:ngeq2initialdata}, Lemma~\ref{lem:nless2initialdataandzerothiterate}, in particular $c_n^{(0)} = 0$, the preceding inductively constructed sequence of iterates converges, in the sense that the (obviously defined) difference norms are $\leq \delta^l$ for some $\delta<1$. 
\\
Defining for $|n|\geq 2$
\[
\varepsilon_{\pm}(n) = \lim_{l\rightarrow\infty}\varepsilon_+^{(l)}(n),\,
\]
and for $|n|<2$ setting $c_n = \lim_{l\rightarrow\infty} c_n^{(l)}$, $\calD_{n}\veps_{+}(n)=\lim_{l\rightarrow\infty}\mathcal{D}_n\epsilon_+^{(l)}(n)$, and 
\[
\varepsilon_+(n) = c_n\phi_n(R) + \phi_n(R)\int_0^R\left[\phi_n(s)\right]^{-1}\mathcal{D}_n\varepsilon_+(n)(\tau,s)\,ds,
\]
and finally $\varepsilon_-(n) = \overline{\varepsilon_+(-n)}$, defining 
\[
\varphi_1(n) = \frac12\sum_{\pm}\varepsilon_{\pm}(n),\quad \varphi_2(n) = \frac{1}{2i}\left[\varepsilon_-(n) - \varepsilon_+(n)\right], 
\]
and further 
\[
\vphi_j = \sum_{n\in Z}\vphi_j(n)e^{in\theta},\quad j = 1,2, 
\]
the function 
\[
\Psi = \mathcal{L}_{v(t)}\mathcal{R}_{h(t)}^{\alpha(t),\beta(t)}\mathcal{S}_{c(t)}\left( \Phi^{\left(\text{mod}_2(t)\right)} + \vphi_1 E_1 + \vphi_2 E_2 + a\left(\Pi_{\Phi^{\perp}}\varphi\right)\Phi\right),\quad \Phi =  \left(\begin{array}{c}\sin U\cos\theta\\ \sin  U \sin\theta\\ \cos U\end{array}\right),
\]
where the modulation parameters are the limits of $\alpha^{(l)}$ etc, and $\Phi^{\left(\text{mod}_2(t)\right)} $ is defined in \eqref{eq:secondarymod}, solves the Wave Maps equation with initial data \eqref{eq:initialperturbation}. 

\end{theorem}
The proof of this theorem follows from Proposition~\ref{prop:exceptionalmodesfinalsourceboundswithmodulation} and passing to difference estimates between subsequent iterates. 

\section{Appendix}

We gather here some technical lemmas which are needed at specific places in the text. 

\begin{lemma}\label{lem:FtoGimprovement} Let $f\in L^\infty(\hbar^{-2},\infty)$. Then if $\nu, l$ are as in the paper, the function
\begin{align*}
G(x): = \int_0^\infty \chi_{x\xi\gtrsim 1}\cdot e^{ix\xi}\cdot \frac{f(\xi)}{[x\xi]^{p + \frac{l}{2} + k\nu}}\cdot\xi^{-1}\,d\xi,\,p>0,\,x>0
\end{align*}
satisfies 
\begin{align*}
H: =\left(x\partial_x\right)^{\lfloor\alpha\rfloor}G\in L^\infty,\,x^{\alpha - \lfloor\alpha\rfloor  +\delta}\cdot  H\in \dot{C}^{\alpha  - \lfloor\alpha\rfloor + \delta}, 
\end{align*}
provided $0<\delta\ll 1$, where we set $\alpha: =  \frac{l}{2} + k\nu$. 
\end{lemma}
\begin{proof} One checks directly that 
\begin{align*}
\left(x\partial_x\right)^{\lfloor\alpha\rfloor}G =  \int_0^\infty \chi_{x\xi\gtrsim 1}\cdot e^{ix\xi}\cdot \frac{f(\xi)}{[x\xi]^{p +\alpha  - \lfloor\alpha\rfloor}}\cdot\xi^{-1}\,d\xi
\end{align*}
Then it suffices to observe that for $\delta<p +\alpha  - \lfloor\alpha\rfloor$ we have (with $x\neq y$)
\begin{align*}
|x-y|^{-\delta}\cdot \left|\frac{x^\delta\cdot  \chi_{x\xi\gtrsim 1}e^{ix\xi}}{[x\xi]^{p +\alpha  - \lfloor\alpha\rfloor}} - \frac{y^\delta\cdot  \chi_{y\xi\gtrsim 1}e^{iy\xi}}{[y\xi]^{p +\alpha  - \lfloor\alpha\rfloor}}\right|\lesssim [x\xi]^{-\delta'} +  [y\xi]^{-\delta'}
\end{align*}
for $\delta' = p +\alpha  - \lfloor\alpha\rfloor - \delta$. 
\end{proof}

\begin{lemma}\label{lem:fromGtoF} Assume the function $G$ is supported on $0\leq x\lesssim 1$ and satisfies 
\[
H: = \left(x\partial_x\right)^{k}G\in L^\infty,\,0\leq k\leq  \lfloor\alpha\rfloor,\,x^{\delta+\alpha - \lfloor\alpha\rfloor}\cdot H\in \dot{C}^{\delta+\alpha - \lfloor\alpha\rfloor},
\]
where it is assumed that $0<\delta+\alpha - \lfloor\alpha\rfloor<1$. Then we have that 
\begin{align*}
F(\xi): = \int_0^\infty\chi_{x\xi\gtrsim 1} e^{ix\cdot\xi}\cdot \left(x\xi\right)^{\alpha}\cdot G(x)\,dx\in \left(\xi^{-\delta}\cdot\dot{C}^{\delta}\cap L^\infty\right)(0,\infty). 
\end{align*}
\end{lemma} 
\begin{proof}
Integration by parts with respect to $x$ results in 
\begin{align*}
F(\xi) = \sum_{0\leq k\leq \lfloor\alpha\rfloor}\int_0^\infty e^{ix\cdot\xi}\cdot \xi^{-\lfloor\alpha\rfloor}\Phi^{\lfloor\alpha\rfloor;k}\left(\chi_{x\xi\gtrsim 1}\cdot\left(x\xi\right)^{\alpha}\right)\cdot \left(x\partial_x\right)^kG(x)\,dx,
\end{align*}
where the operator $\Phi^{\lfloor\alpha\rfloor;k}$ is a finite linear combination of `monomials' of the form $\prod_{j=1}^{\lfloor\alpha\rfloor}T_j$ where $k$ of the $T_j$'s equal $\frac{1}{x}$ and $\lfloor\alpha\rfloor - k$ of the $T_j$'s equal $\partial_x$. 
Performing one more integration by parts, we infer 
\begin{align*}
&\sum_{0\leq k\leq \lfloor\alpha\rfloor-1} \int_0^\infty e^{ix\cdot\xi}\cdot \xi^{-\lfloor\alpha\rfloor}\Phi^{\lfloor\alpha\rfloor;k}\left(\chi_{x\xi\gtrsim 1}\cdot\big(x\xi\big)^{\alpha}\right)\cdot \left(x\partial_x\right)^kG(x)\,dx\\
& = \sum_{0\leq k\leq \lfloor\alpha\rfloor} \int_0^\infty e^{ix\cdot\xi}\cdot \xi^{-\lfloor\alpha\rfloor-1}\Psi^{\lfloor\alpha\rfloor+1;k}\left(\chi_{x\xi\gtrsim 1}\cdot\big(x\xi\big)^{\alpha}\right)\cdot \left(x\partial_x\right)^kG(x)\,dx\\
\end{align*}
where the $\Psi^{\lfloor\alpha\rfloor+1;k}$ are defined analogously to the $\Phi^{\lfloor\alpha\rfloor;k}$  whence
\begin{align*}
\left|\xi^{-\lfloor\alpha\rfloor-1}\Psi^{\lfloor\alpha\rfloor+1;k}\left(\chi_{x\xi\gtrsim 1}\cdot\left(x\xi\right)^{\alpha}\right)\right|\lesssim (x\xi)^{\alpha - \lfloor\alpha\rfloor-1}
\end{align*}
Since furthermore we have 
\begin{align*}
\left\|\xi^{\delta}\cdot \xi^{-\lfloor\alpha\rfloor-1}\Psi^{\lfloor\alpha\rfloor+1;k}\left(\chi_{x\xi\gtrsim 1}\cdot\left(x\xi\right)^{\alpha}\right)\right\|_{\dot{C}^{\delta}}\lesssim 1,\,\left\|\xi^{\delta}\cdot e^{ix\xi}\right\|_{\dot{C}^{\delta}(\xi\sim\mu)}\lesssim (x\mu)^{\delta}, 
\end{align*}
we find that 
\begin{align*}
\left\|\xi^{\delta}\cdot\sum_{0\leq k\leq \lfloor\alpha\rfloor} \int_0^\infty e^{ix\cdot\xi}\cdot \xi^{-\lfloor\alpha\rfloor-1}\Psi^{\lfloor\alpha\rfloor+1;k}\left(\chi_{x\xi\gtrsim 1}\cdot\left(x\xi\right)^{\alpha}\right)\cdot \left(x\partial_x\right)^kG(x)\,dx\right\|_{\dot{C}^{\delta}}\lesssim 1
\end{align*}
provided $\delta+\alpha - \lfloor\alpha\rfloor<1$. Finally, considering the case $k = \lfloor\alpha\rfloor$, we have 
\begin{align*}
&\int_0^\infty\chi_{x\xi\gtrsim 1} e^{ix\cdot\xi}\cdot \left(x\xi\right)^{\alpha-\lfloor\alpha\rfloor}\cdot H(x)\,dx\\
& = -\frac12\cdot\int_0^\infty\chi_{x\xi\gtrsim 1} e^{ix\cdot\xi}\cdot \left(\left[x + \frac{\pi}{\xi}\right]\xi\right)^{\alpha-\lfloor\alpha\rfloor}\cdot H\left(\left[x + \frac{\pi}{\xi}\right]\right)\,dx\\
& + \frac12\cdot \int_0^\infty\chi_{x\xi\gtrsim 1} e^{ix\cdot\xi}\cdot \left(x\xi\right)^{\alpha-\lfloor\alpha\rfloor}\cdot H(x)\,dx,\\
\end{align*}
and we have the estimate 
\begin{align*}
&\left|\left(\left[x + \frac{\pi}{\xi}\right]\xi\right)^{\alpha-\lfloor\alpha\rfloor}\cdot H\left(\left[x + \frac{\pi}{\xi}\right]\right) -  \left(x\xi\right)^{\alpha-\lfloor\alpha\rfloor}\cdot H(x)\right|\\
&\lesssim (x\xi)^{-\delta}\cdot \left(\left\|H\right\|_{L^\infty} + \left\|x^{\delta + \alpha - \lfloor\alpha\rfloor}H\right\|_{\dot{C}^{\delta + \alpha - \lfloor\alpha\rfloor}}\right). 
\end{align*}
The desired Holder bound follows easily from this. 
\end{proof}

\begin{lemma}\label{lem:transferenceonadmsingoffdiag} Let $\chi_{\gtrsim 1}$ be a smooth cutoff localising to the indicated region. Let $F_n(\xi,\eta)$ be the numerator occurring in the $|n|\geq 2$ transference operator(continuous part), see Prop.~\ref{prop: K operator}. Then if $x(\tau,\xi)$ is admissibly singular and $\delta>0$, the function 
\begin{align*}
y(\tau,\xi): = \int_0^\infty \chi_{\gtrsim 1}\left(\frac{\left|\xi^{\frac12} - \eta^{\frac12}\right|}{\xi^{\delta}}\right)\frac{F_n(\xi,\eta)\tilde{\rho}^{(\hbar)}(\eta)}{\xi - \eta}\cdot x(\tau,\eta)\,d\eta
\end{align*}
is in $S_0^{(\hbar)}$, and decaying rapidly in time $\tau$. Precisely, we have the bound 
\begin{align*}
\tau^5\cdot\left\|y(\tau,\cdot)\right\|_{S_0^{(\hbar)}}\lesssim \left\|\overline{x}\right\|_{adm}. 
\end{align*}
If $\overline{x}$ is source admissibly singular, then a similar conclusion applies upon replacing $S_0^{(\hbar)}$ by $S_1^{(\hbar)}$. 
\end{lemma}

We note that this lemma justifies inclusion of an extra $\chi_{\lesssim 1}\left(\frac{\left|\xi^{\frac12} - \eta^{\frac12}\right|}{\xi^{\delta}}\right)$, and since $\eta\gtrsim\hbar^{-2}$ on the support of the admissibly singular factor $x(\tau,\eta)$, this means 
\[
\left|\frac{\left|\xi^{\frac12} - \eta^{\frac12}\right|}{\eta^{\frac12}}\right|\lesssim \hbar^{1-\frac{\delta}{2}}, 
\]
and hence for $\delta<\frac23$, say, we have that $\left|\frac{\left|\xi^{\frac12} - \eta^{\frac12}\right|}{\eta^{\frac12}}\right|\lesssim \hbar^{\frac23}$, which is a necessary condition to derive good trace derivative bounds on the kernel $F(\xi,\eta)$. 
\begin{proof}(lemma) 

We consider the case that $\overline{x}$ is of prototypical singular type, of the form 
\begin{align*}
\overline{x}(\tau, \xi) = \hbar^{-1}\cdot \chi_{\hbar^2\xi\gtrsim 1}\frac{e^{\pm i\left(\nu\tau\xi^{\frac12} + \hbar^{-1}\rho_n(x_{\tau};\alpha,\hbar)\right)}}{\xi^{1+\frac{k\nu}{2}}}\left(\log\xi\right)^j\cdot F_{k,j}^{\pm}(\tau,\xi), 
\end{align*}
where $F_{k,j}^{\pm} = F_{0,k,j}^{\pm}$ satisfies the bounds in Definition~\ref{defi:xsingulartermsngeq2proto}. The contributions of the formally more complicated terms in an admissibly singular function are handled analogously. 
Recalling the proof of Prop.~\ref{prop: K operator}, and specifically \eqref{off diagonal F pre}, we write 
\begin{align*}
F_n(\xi,\eta) = \int_0^\infty W_n^+(R)\phi_n(R;\xi)\cdot \phi_n(R;\eta)R\,dR.
\end{align*}
Denoting by $x_t^{(\eta)}$ the turning point in the sense of Lemma~\ref{lem: monotonicity of root in alpha}
 and for frequency $\eta$ and at fixed angular momentum $n$, we decompose 
\begin{align*}
F_n(\xi,\eta) &=  \int_0^\infty W_n^+(R)\chi_{x^{(\eta)}<\frac{x_t^{(\eta)}}{2}}\phi_n(R;\xi)\cdot \phi_n(R;\eta)R\,dR\\
& + \int_0^\infty W_n^+(R)\chi_{x^{(\eta)}\in[\frac{x_t^{(\eta)}}{2}, x_t^{(\eta)} + \hbar^{\frac23}]}\phi_n(R;\xi)\cdot \phi_n(R;\eta)R\,dR\\
& +  \int_0^\infty W_n^+(R)\chi_{x^{(\eta)}\in[x_t^{(\eta)} + \hbar^{\frac23}, (1+\gamma)x_t^{(\eta)}]}\phi_n(R;\xi)\cdot \phi_n(R;\eta)R\,dR\\
& + \int_0^\infty W_n^+(R)\chi_{x^{(\eta)}\in[(1+\gamma)x_t^{(\eta)},\infty)}\phi_n(R;\xi)\cdot \phi_n(R;\eta)R\,dR\\
& =: F_{n,I} + F_{n,II} + F_{n,III} + F_{n,IV}. 
\end{align*}
for some $0<\gamma\ll 1$. Here the cutoff $\chi_{x^{(\eta)}<\frac{x_t^{(\eta)}}{2}}$ may be realized by a smooth function, while the cutoff $\chi_{x^{(\eta)}\in[\frac{x_t^{(\eta)}}{2}, x_t^{(\eta)} + \hbar^{\frac23}]}$ is to be thought of as a sum of smooth cutoffs localizing to the intervals
\[
[x_t^{(\eta)} - \hbar^{\frac23}, x_t^{(\eta)} + \hbar^{\frac23}]], \{[x_t - 2^j\cdot \hbar^{\frac23}, x_t - 2^{j+2}\cdot \hbar^{\frac23}]\}_{j=1}^{c|\log\hbar|},
\]
and analogously for $\chi_{x^{(\eta)}\in[x_t^{(\eta)} + \hbar^{\frac23}, (1+\gamma)x_t^{(\eta)}]}$, while $\chi_{x^{(\eta)}\in[(1+\gamma)x_t^{(\eta)},\infty)}$ is again realized via a smooth cutoff. Each of the smooth cutoffs has symbol behavior with respect to differentiation (with respect to either $R$ or $\eta$), and summing over the cutoffs costs as most a factor $\sim|\log\hbar|$, which will be overwhelmed by gains later on. We shall hence commit abuse of notation and assume that all cutoffs in the terms $I - IV$ have symbol behavior. We now derive the required bound independently for the terms $I - IV$ by exploiting the asymptotics for $\phi_n(R;\eta), \phi_n(R;\xi)$ furnished by Prop.~\ref{prop:DFT nlarge}. 
\\

{\it{The estimate for the contribution of the term $F_{n,I}$.}} Observe that under the condition $\eta\gtrsim \hbar^{-2}$ we have (with $\beta = \eta^{\frac12}\cdot\hbar$)
\begin{align*}
\partial_{\eta^{\frac12}}\left[\nu\tau\eta^{\frac12} + \hbar^{-1}\rho_n(x_{\tau}^{(\eta)};\beta,\hbar)\right] = \nu\tau + \hbar^{-1}\cdot \tau\hbar\cdot \rho_{n,x}(x_{\tau}^{(\eta)};\beta,\hbar) + \rho_{n,\beta}(x_{\tau}^{(\eta)};\beta,\hbar) = \nu\tau + O(1), 
\end{align*}
since due to Lemma~\ref{lem: Lemma 3.4 CDST} we have 
\begin{align*}
\left| \hbar^{-1}\cdot \tau\hbar\cdot \rho_{n,x}(x_{\tau}^{(\eta)};\beta,\hbar)\right| + \left|\rho_{n,\beta}(x_{\tau}^{(\eta)};\beta,\hbar)\right| = O(1)
\end{align*}
uniformly in all parameters, including $\hbar$. We shall then perform integration by parts with respect to $\eta^{\frac12}$ in the integral defining $y(\tau,\xi)$. Noting that 
\begin{align*}
\left|\partial_{\eta^{\frac12}}^k\left(\chi_{\gtrsim 1}\left(\frac{|\xi^{\frac12} - \eta^{\frac12}|}{\xi^{\delta}}\right)\cdot \frac{1}{\xi^{\frac12} - \eta^{\frac12}}\right)\right|\lesssim_k \xi^{-\delta k},\, 
\end{align*}
and further, from the proof of Prop.~\ref{prop: K operator} and Prop.~\ref{prop:DFT nlarge} the bounds 
\begin{align*}
\left|\partial_{\eta^{\frac12}}^k F_{n,I}\right|\lesssim c^{\hbar^{-1}}\cdot \eta^{-\frac{k}{2}}\cdot \left\langle \frac{\xi}{\eta}\right\rangle^{-10},\,0\leq k\ll\hbar^{-1},\,\left|\partial_{\eta^{\frac12}}^k\tilde{\rho}^{(\hbar)}(\eta)\right|\lesssim_k\eta^{-\frac{k}{2}},
\end{align*}
where $c\in (0, 1)$ is a suitable constant. It follows that upon setting $\psi(\eta,\tau;\hbar): = \nu\tau\eta^{\frac12} + \hbar^{-1}\rho_n(x_{\tau}^{(\eta)};\beta,\hbar)$, we have 
\begin{align*}
&\left|\int_0^\infty \chi_{\gtrsim 1}\left(\frac{|\xi^{\frac12} - \eta^{\frac12}|}{\xi^{\delta}}\right)\frac{F_{n,I}(\xi,\eta)\tilde{\rho}^{(\hbar)}(\eta)}{\xi - \eta}\cdot \overline{x}(\tau,\eta)\,d\eta\right|\\
 &= \left|\hbar^{-1}\int_0^\infty \chi_{\gtrsim 1}\left(\frac{|\xi^{\frac12} - \eta^{\frac12}|}{\xi^{\delta}}\right)\frac{F_{n,I}(\xi,\eta)\tilde{\rho}^{(\hbar)}(\eta)}{\xi - \eta}\cdot \left((\partial_{\eta^{\frac12}}\psi)^{-1}(\eta,\tau;\hbar)\cdot\partial_{\eta}^{\frac12}\right)^p\left(e^{i\psi}\right)\cdot  \frac{\chi_{\hbar^2\eta\gtrsim 1}\cdot \left(\log\eta\right)^j\cdot F_{k,j}^{\pm}(\tau,\eta)}{\eta^{1+\frac{k\nu}{2}}}\,d\eta\right|\\
 &\lesssim  \hbar^{-1}\cdot c^{\hbar^{-1}}\cdot \langle\xi\rangle^{-10}\cdot \tau^{-10}\cdot \sum_{0\leq k_1\leq 10}\left|\eta^{k_1}\partial_{\eta}^{k_1} F_{k,j}^{\pm}(\tau,\eta)\right|.
\end{align*}
Noting that 
\begin{align*}
 \sum_{0\leq k_1\leq 10}\left|\eta^{k_1}\partial_{\eta}^{k_1} F_{k,j}^{\pm}(\tau,\eta)\right|\lesssim \tau^{-1-}\cdot \left\|\overline{x}\right\|_{adm}, 
\end{align*}
we infer the bound 
\begin{align*}
\left\|\int_0^\infty \chi_{\gtrsim 1}\left(\frac{|\xi^{\frac12} - \eta^{\frac12}|}{\xi^{\delta}}\right)\frac{F_{n,I}(\xi,\eta)\tilde{\rho}^{(\hbar)}(\eta)}{\xi - \eta}\cdot x(\tau,\eta)\,d\eta\right\|_{S_0^{\hbar}}\lesssim \tilde{c}^{\hbar^{-1}}\cdot\tau^{-10}\cdot \left\|\overline{x}\right\|_{adm}, 
\end{align*}
which is much better than needed. 
\\

{\it{The estimate for the contribution of the term $F_{n,II}$.}} Using Prop.~\ref{prop:DFT nlarge}, we infer the bounds 
\begin{align*}
\left|\partial_{\eta^{\frac12}}^k\phi_n(R;\eta)\right|\lesssim_k \hbar^{\frac{k}{3}}\cdot R^k\cdot e^{-\mu^{\frac32}}\lesssim \hbar^{\frac{k}{3}}\cdot (\hbar\eta^{\frac12})^{-k}\cdot e^{-\mu^{\frac32}}
\end{align*}
provided (with $\mu\geq 1$)
\[
x-x_t\sim -\mu\cdot\hbar^{\frac23},\,\eta\gtrsim \hbar^{-2}.
\]
Since $\hbar\eta^{\frac12}\gtrsim 1$, we can also deduce
\begin{align*}
\left|\partial_{\eta^{\frac12}}^k\phi_n(R;\eta)\right|\lesssim_k \hbar^{\frac{k}{6}}\cdot\eta^{-\frac{k}{6}}\cdot e^{-\mu^{\frac32}}.
\end{align*}
Furthermore we have the bound 
\begin{align*}
\left|\partial_{\eta^{\frac12}}^k\partial_R^l\phi_n(R;\eta)\right|\lesssim \hbar^{\frac{k+2l}{6}}\cdot \eta^{\frac{3l-k}{6}}\cdot e^{-\mu^{\frac32}}
\end{align*}
If $\xi\gg \eta$, then since $x^{(\eta)}\sim1$ on the support of $F_{n,II}$, we have $x^{(\xi)}\gg 1$, whence $\phi(R;\xi)$ will be in the oscillatory regime. Using integration by parts with respect to $R$ in this latter case, in conjunction with the preceding bound, furnishes the estimate (after summing over $\mu$)
\begin{align*}
\left|\partial_{\eta^{\frac12}}^k F_{n,II}\right|\lesssim_{k} \hbar^{\frac{k}{6}}\cdot \left\langle \frac{\xi}{\eta}\right\rangle^{-10}\cdot \eta^{-\frac{k}{6}}. 
\end{align*}
One then completes the bound as in the preceding case, except that this time we only gain polynomially in $\hbar$ .
\\

{\it{The estimate for the contribution of the term $F_{n,III}$.}} This is similar to the preceding case, except one does not gain in $\mu$ for $x = x_t + \hbar^{\frac23}\cdot\mu\in [x_t + \hbar^{\frac23}, (1+\gamma)x_t]$. This results in a loss of $|\log\hbar|$ from summing over dyadic $\mu$, but this is absorbed by the polynomial gains in $\hbar$. 
\\

{\it{The estimate for the contribution of the term $F_{n,IV}$.}} Here the function $\phi_n(R;\eta)$ is in the oscillatory regime, and can be written in the form 
\[
\phi_n(R;\eta) = \sum_{\pm}e^{\pm i\left(R\eta^{\frac12} - \hbar^{-1}y(\beta;\hbar) + \hbar^{-1}\rho_n(x^{(\eta)};\beta,\hbar)\right)}\cdot\left(1 + \hbar\tilde{a}_{1,\pm}(\zeta,\beta;\hbar)\right), 
\]
where we refer to section~\ref{sec:LAN} and Prop.~\ref{prop:DFT nlarge} for the notation and properties of the displayed functions. Now if $R\ll\tau$, whence the phases $e^{\pm i\left(\nu\tau\eta^{\frac12} +\hbar^{-1}\rho_n(x_{\tau}^{(\eta)};\beta,\hbar)\right)}$ and $e^{\pm i(R  \hbar^{-1}\rho_n(x^{(\eta)};\beta,\hbar))}$ are non-resonant, then repeated integration by parts with respect to $\eta^{\frac12}$ and with respect to $R$ if $\xi\gg \eta$ results in the bound 
\begin{align*}
\left|\int_0^\infty \chi_{\gtrsim 1}\left(\frac{\left|\xi^{\frac12} - \eta^{\frac12}\right|}{\xi^{\delta}}\right)\frac{F_{n,IV}(\xi,\eta)\tilde{\rho}^{(\hbar)}(\eta)}{\xi - \eta}\cdot x(\tau,\eta)\,d\eta\right|\lesssim \langle\xi\rangle^{-10}\cdot \tau^{-10}\cdot \left\|\overline{x}\right\|_{adm}, 
\end{align*}
which again amply suffices to recover the bound claimed in the lemma. \\
If we impose the additional localization $R\gtrsim \tau$ in the integral defining $F_{n,IV}$, resulting in $F_{n,IV}^{R\gtrsim\tau}$, then integration by parts with respect to $R$ and the fact that $|\xi^{\frac12}-\eta^{\frac12}|\gtrsim \xi^{\delta}$ result in a similar bound, since we then have 
\begin{align*}
\left| F_{n,IV}^{R\gtrsim\tau}\right|\lesssim \langle\xi\rangle^{-10}\cdot \tau^{-10}, 
\end{align*}
and from here the bound 
\begin{align*}
\left|\int_0^\infty \chi_{\gtrsim 1}\left(\frac{\left|\xi^{\frac12} - \eta^{\frac12}\right|}{\xi^{\delta}}\right)\frac{F_{n,IV}^{R\gtrsim \tau}(\xi,\eta)\tilde{\rho}^{(\hbar)}(\eta)}{\xi - \eta}\cdot x(\tau,\eta)\,d\eta\right|\lesssim \langle\xi\rangle^{-10}\cdot \tau^{-10}\cdot \left\|\overline{x}\right\|_{adm}
\end{align*}
easily follows. This again confirms the claim of the lemma for this contribution. 
\end{proof}

\begin{lemma}\label{lem:lem:singFouriertiphysicalngeq2admrefined} Let $\lambda\gtrsim \hbar^{-1}$. Then if $\phi$ is admissibly singular, with distorted Fourier coefficient $\overline{x}$, the functions $P_{\geq\lambda}\phi, P_{\lambda}\phi$, defined as in subsection~\ref{subsubsec:Fourierlocal}, satisfy the conclusion of Lemma~\ref{lem:singFouriertiphysicalngeq2adm}, and in the inequalities for $f_1, f_3, G_{k,l,j}$, we gain an additional factor $\langle\lambda (\nu\tau-R)\rangle^{-\gamma}$ for some $\gamma>0$. We also have the bounds
\begin{align*}
\sum_{j=1,3}\left(\left\|f_j^{(\geq\lambda)}\right\|_{\tilde{S}_0^{(\hbar)}} + \left\|f_j^{(\lambda)}\right\|_{\tilde{S}_0^{(\hbar)}}\right)\lesssim \lambda^{-(0+)}\cdot \tau^{-1-\nu}\cdot\left(\log\tau\right)^{N_1}\cdot \left\|\overline{x}\right\|_{adm},
\end{align*}
where $f_j^{(\geq\lambda)},\,f_j^{(\lambda)}$ refers to the components of the decomposition of $P_{\geq\lambda}\phi, P_{\lambda}\phi$, in accordance with Lemma~\ref{lem:singFouriertiphysicalngeq2adm}.
Finally, we have the improved estimates
\begin{align*}
\hbar^{1+\alpha}\left\|\partial_R^{\alpha}f_j^{(\geq \lambda)}(\tau, \cdot)\right\|_{L^\infty}\lesssim \lambda^{-\frac12}\cdot \tau^{-1-\nu}\cdot(\log\tau)^{N_1}\cdot \left\|\overline{x}\right\|_{adm},\,j = 1, 3. 
\end{align*}
and similarly for $f_j^{(\lambda)},\,j = 1, 3$. 
\end{lemma}
\begin{proof} This follows since one gains additional factors $[\xi^{\frac12}(\nu\tau - R)]^{-1}$ in the proof of Lemma~\ref{lem:singFouriertiphysicalngeq2adm} for $f_1, f_3$ by integration by parts with respect to $\xi^{\frac12}$, or by simply examining the structure of the $\xi$-integrals giving rise to the components of $f_2$. 
\end{proof}

\begin{lemma}\label{lem:basicproductsource} Let $\phi_1, \phi_2$ be angular momentum $n_j,\,j = 1, 2$ functions with $|n_j|\geq 2$, and such that their distorted Fourier transform (at angular momentum $n_j$) satisfies 
\[
\left\langle \phi_{n_j}(R;\xi),\,\phi_j\right\rangle_{L^2_{R\,dR}} = \overline{y}_j + \overline{z}_j
\]
with 
\[
\left\|\overline{y}_j\right\|_{adm} + \sup_{\tau\geq \tau_0}\tau^3\left\|\overline{z}_j\right\|_{S_0^{\hbar_j}} =:\Lambda_j<\infty. 
\]
Then letting $\{m, n_1, n_2\}$ be admissible, the product 
\[
\Phi: = \chi_{1\lesssim R\lesssim\nu\tau}\cdot \langle R\rangle^{-2}\phi_1\cdot\phi_2,
\]
satisfies 
\[
\left\langle \phi_m(R;\xi),\,\Phi\right\rangle_{L^2_{R\,dR}} = \overline{y} + \overline{z},
\]
where 
\[
\left\|\overline{y}\right\|_{sourceadm} +  \sup_{\tau\geq \tau_0}\tau^4\left\| \overline{z}\right\|_{S_1^{\hbar}}\lesssim \min\{|n_{1,2}|\}^C\cdot \prod_j\Lambda_j,\,\hbar = |m|^{-1}. 
\]
The same conclusion applies to the expression 
\[
\Phi: = \chi_{1\lesssim R\lesssim\nu\tau}\cdot \langle R\rangle^{-2}\partial_{\theta}\phi_1\cdot\partial_{\theta}\phi_2,
\]
\end{lemma}
\begin{remark}\label{rem::basicproductsource} The lemma can be refined for certain more special situations. For example, when including a cutoff $\chi_{1\lesssim R\ll \nu\tau}$, it suffices to assume a weaker bound 
\[
\sup_{\tau\geq \tau_0}\tau^2\left\|\overline{z}_j\right\|_{S_0^{\hbar_j}}
\]
for one of the two factors. 
\end{remark}
\begin{proof} We may assume that $|n_1|\leq |n_2|$. Then we distinguish between the following cases:
\\

{\it{(1): Both $\phi_{1,2}$ are in $S_0^{\hbar_j}$.}} In this case we show that $\Phi\in \tilde{S}_1^{\hbar}$, with suitable temporal decay. In fact, from Prop.~\ref{prop:bilin1}, slightly modified with $R^{-1}$ instead of $\partial_R$, we obtain the even stronger bound 
\begin{align*}
 \sup_{\tau\geq \tau_0}\left\|\tau^6\cdot \left\langle \phi_m(R;\xi),\,\Phi\right\rangle_{L^2_{R\,dR}}\right\|_{S_1^{\hbar}}\lesssim |n_1|^C\cdot \prod_j\Lambda_j.
\end{align*}
{\it{{(2)}: $\phi_1$ admissibly singular, $\phi_2\in S_0^{\hbar_2}$.} } We decompose $\phi_1$ according to Lemma~\ref{lem:singFouriertiphysicalngeq2adm}, 
\[
\phi_1 = \sum_{j=1}^3\phi_{1j}
\]
where the functions $\phi_{1j},\,j =1,\,3$ are in $S_0^{\hbar_1}$ as detailed in Lemma~\ref{lem:singFouriertiphysicalngeq2adm}, while $\phi_{12}$ is of the explicit singular kind. By direct inspection of Lemma~\ref{lem:singFouriertiphysicalngeq2adm}, we have the bounds 
\begin{align*}
\tau^{2+}\cdot\sum_{j=1,3}\left\|\langle R\rangle^{-1}\phi_{1j}(\tau,\cdot)\right\|_{S_1^{\hbar_1}}\lesssim \Lambda_1.
\end{align*}
Again referring to Prop.~\ref{prop:bilin1}, we find
\begin{align*}
\tau^{5+}\cdot \left\|\left\langle \phi_{m}(R;\xi),\,\left(\sum_{j=1,3}\phi_{1j}\right)\cdot \phi_2\right\rangle_{L^2_{R\,dR}}\right\|_{S_1^{\hbar}}\lesssim \prod_{j=1,2}\Lambda_j. 
\end{align*}
On the other hand, replacing $\phi_1$ by $\phi_{12}$, we decompose the product $\phi_{12}\cdot \phi_2$ into a structured singular term and a better but less structured term in $S_1^{\hbar}$. Precisely, write as in Lemma~\ref{lem:singFouriertiphysicalngeq2adm}
\begin{align*}
\phi_{12}(\tau, R) = \chi_{|\nu\tau - R|\lesssim \hbar_1}\cdot \sum \frac{G_{k,l,j}(\tau, \nu\tau - R)}{\tau^{\frac12}}\cdot \hbar_1^{-\frac{l+1}{2}}\cdot [\nu\tau - R]^{\frac{l}{2}+k\nu}\cdot\left(\log(\nu\tau - R)\right)^j.
\end{align*}
It is easily seen that this function, when localized to $|\nu\tau - R|\gtrsim \hbar_2$, is in $\tilde{S}_0^{\hbar_2}$, and so the product of this contribution and $\phi_2$ can be handled as in the previous case. Henceforth we restrict it to $|\nu\tau - R|\lesssim \hbar_2$
Further, decompose $\phi_2$ into a fourth order Taylor polynomial and an error term around $R = \nu\tau$: 
\begin{align*}
\phi_2(\tau, R) = \mathbb{P}_4\left(\phi_2(\tau,\cdot)\right)(R) + \mathbb{E}_4\left(\phi_2(\tau,\cdot)\right)(R),
\end{align*}
where we have 
\begin{align*}
\mathbb{P}_4\left(\phi_2(\tau,\cdot)\right)(R)  = \sum_{j=0}^4 a_j(\tau)\cdot(\nu\tau - R)^j
\end{align*}
with the bounds 
\begin{align*}
\left| a_j(\tau)\right|\lesssim \hbar_2^{-1-j}\cdot \tau\cdot\left\|\phi_2\right\|_{S_0^{\hbar_2}},\,j = 0,\ldots, 4.
\end{align*}
see Lemma~\ref{lem:derLinfty}. It follows that 
\begin{align*}
&\phi_{12}(\tau, R)\cdot \mathbb{P}_4\left(\phi_2\right)(\tau, R)\\&
 =  \chi_{|\nu\tau - R|\lesssim \hbar_2}\cdot \sum a_r(\tau)\cdot \frac{G_{k,l,j}(\tau, \nu\tau - R)}{\tau^{\frac12}}\cdot \hbar_1^{-\frac{l+1}{2}}\cdot [\nu\tau - R]^{\frac{l}{2}+r+k\nu}\cdot\big(\log(\nu\tau - R)\big)^j\\
 & =   \chi_{|\nu\tau - R|\lesssim \hbar_2}\cdot \sum \tilde{a}_r(\tau)\cdot \frac{G_{k,l,j}(\tau, \nu\tau - R)}{\tau^{\frac12}}\cdot \hbar_2^{-\frac{l+2r+1}{2}}\cdot [\nu\tau - R]^{\frac{l}{2}+r+k\nu}\cdot\left(\log(\nu\tau - R)\right)^j\,\
\end{align*}
where we have 
\begin{align*}
\left|\tilde{a}_r(\tau)\right|\lesssim \hbar_1^{-\frac{l+1}{2}}\cdot \tau\cdot\left\|\phi_2\right\|_{S_0^{\hbar_2}},\,r = 0,\ldots, 4. 
\end{align*}
The preceding term is then seen to be a singular source admissible term, in light of Lemma~\ref{lem:singPhysicaltoFourierngeq2}, but with a factor $\hbar_1^{-C}\sim |n_1|^C$. \\
The same applies to the product 
\[
\mathbb{P}_4\left(\langle R\rangle^{-2}\right)\cdot \phi_{12}(\tau, R)\cdot \mathbb{P}_4\left(\phi_2\right)(\tau, R)
\]
where we have replaced the factor $\langle R\rangle^{-2}$ by its 4th order Taylor polynomial around $R = \nu\tau$. \\
As for the contribution of $\mathbb{E}_4\left(\phi_2(\tau,\cdot)\right)$, this will be in $S_1^{\hbar}$. In fact, we have for $r\in \{0,1,\ldots,4\}$
\begin{equation}\label{eq:E4derivativestructure}
\left|\partial_R^r\mathbb{E}_4\left(\phi_2\right)(\tau, R)\right|\lesssim \hbar_2^{-5-}\cdot (\nu\tau - R)^{\frac92 - r+}\cdot \left\|\phi_2\right\|_{S_0^{\hbar_2}}.
\end{equation}
We then estimate 
\begin{align*}
\left\|\langle R\rangle^{-2}\cdot \phi_{12}(\tau, R)\cdot \mathbb{E}_4\left(\phi_2\right)\right\|_{\tilde{S}_1^{\hbar}} \lesssim &\hbar^{\frac12}\cdot \left\|\langle \phi_m(R;\xi),\,\langle R\rangle^{-2}\cdot \phi_{12}(\tau, R)\cdot \mathbb{E}_4\left(\phi_2\right)\rangle_{L^2_{R\,dR}}\right\|_{L^2_{d\xi}(\xi^{\frac12}\hbar\lesssim 1)}\\
& + \hbar^{5+}\left\|\xi^{2+}\langle \phi_m(R;\xi),\,\langle R\rangle^{-2}\cdot \phi_{12}(\tau, R)\cdot \mathbb{E}_4\left(\phi_2\right)\rangle_{L^2_{R\,dR}}\right\|_{L^2_{d\xi}(\xi^{\frac12}\hbar\gtrsim 1)},
\end{align*}
where the first term on the right is bounded by (using Plancherel's theorem for the distorted Fourier transform)
\begin{align*}
&\hbar^{\frac12}\cdot \left\|\langle \phi_m(R;\xi),\,\langle R\rangle^{-2}\cdot \phi_{12}(\tau, R)\cdot \mathbb{E}_4\left(\phi_2\right)\rangle_{L^2_{R\,dR}}\right\|_{L^2_{d\xi}(\xi^{\frac12}\hbar\lesssim 1)}\\
&\lesssim \hbar^{\frac12}\cdot \tau^{-2}\cdot \left\| \phi_{12}(\tau, R)\right\|_{L^\infty_{dR}}\cdot \left\|\chi_{|\nu\tau - R|\lesssim \hbar_2} \mathbb{E}_4\left(\phi_2\right)\right\|_{L^2_{R\,dR}}\\
&\lesssim  \tau^{-\frac32}\cdot\hbar_1^{-C}\cdot \left\|\overline{y}_1\right\|_{adm}\cdot \left\|\phi_2\right\|_{\tilde{S}_0^{\hbar_2}}. 
\end{align*}
Here we have used Lemma~\ref{lem:singFouriertiphysicalngeq2adm} to bound the coefficients of the explicit expansion for $\phi_{12}$, and we have used the fact that $R\sim \tau$ on the support of $\phi_{12}$. \\
For the second term above in the high frequency regime $\xi^{\frac12}\hbar\gtrsim 1$, we can bound this by using Plancherel's theorem and the Leibniz rule to distribute derivatives:
\begin{align*}
&\hbar^{5+}\left\|\xi^{2+}\langle \phi_m(R;\xi),\,\langle R\rangle^{-2}\cdot \phi_{12}(\tau, R)\cdot \mathbb{E}_4\left(\phi_2\right)\rangle_{L^2_{R\,dR}}\right\|_{L^2_{d\xi}(\xi^{\frac12}\hbar\gtrsim 1)}\\
&\lesssim \sum_{j+k = 4}\hbar^{5+}\left\|\xi^{0+}\langle \phi_m(R;\xi),\,\partial_R^j\left(\langle R\rangle^{-2}\cdot \phi_{12}(\tau, R)\right)\cdot \partial_R^k\mathbb{E}_4\left(\phi_2\right)\rangle_{L^2_{R\,dR}}\right\|_{L^2_{d\xi}(\xi^{\frac12}\hbar\gtrsim 1)}\\
\end{align*}
Taking advantage of \eqref{eq:E4derivativestructure} as well as the fine structure of $\phi_{12}$ displayed earlier (via Lemma~\ref{lem:singFouriertiphysicalngeq2adm}), we can bound the preceding by 
\begin{align*}
&\sum_{j+k = 4}\hbar^{5+}\left\|\xi^{0+}\langle \phi_m(R;\xi),\,\partial_R^j\left(\langle R\rangle^{-2}\cdot \phi_{12}(\tau, R)\right)\cdot \partial_R^k\mathbb{E}_4\left(\phi_2\right)\rangle_{L^2_{R\,dR}}\right\|_{L^2_{d\xi}(\xi^{\frac12}\hbar\gtrsim 1)}\\
&\lesssim \hbar^{5+}\cdot \sum_{j = 0,\ldots,4}\left\|\partial_R^j\left(\langle R\rangle^{-2}\cdot \phi_{12}(\tau, R)\right)\cdot (\nu\tau - R)^{\frac12+j}\right\|_{H^{0+}_{R\,dR}}\\
&\lesssim \tau^{-1}\cdot \hbar_1^{-C}\cdot \left\|\overline{y}_1\right\|_{adm}\cdot \left\|\phi_2\right\|_{\tilde{S}_0^{\hbar_2}}. 
\end{align*}
Here we have also used that $\hbar^{5+}\cdot \hbar_2^{-5-}\lesssim \hbar_1^{-5-}$ since $|n_1|\leq |n_2|$ and either $|m|\sim |n_2|$ or $|n_1|\sim |n_2|$. 
\\

{\it{{(3)}: $\phi_2$ admissibly singular, $\phi_1\in S_0^{\hbar_1}$.} } This is similar to the preceding case, we omit the details. 
\\

{\it{(4): both $\phi_j,\,j = 1, 2$ of admissibly singular type.}} Here we implement the decomposition 
\[
\phi_j = \sum_{k=1}^3\phi_{jk}
\]
according to Lemma~\ref{lem:singFouriertiphysicalngeq2adm} for both factors. Then the estimate 
\begin{align*}
\tau^4\left\|\chi_{1\lesssim R\lesssim \nu\tau}\langle R\rangle^{-2}\prod_{j=1,2}\left(\sum_{k=1,3}\phi_{jk}\right)\right\|_{S_1^{\hbar}}\lesssim |n_1|^C\cdot\prod_{j=1,2}\left\|\overline{y}_j\right\|_{adm}
\end{align*}
follows directly from the pointwise bounds in Lemma~\ref{lem:singFouriertiphysicalngeq2adm} as well as Plancherel's theorem. 
\\
If we substitute $\phi_{j2}$ for at least one of the factors $\phi_j$, the desired bound follows as in {\it{(2)}}, using a Taylor decomposition around $R = \nu\tau$ for the regular factor.  
\\

The last statement of the lemma is proved in analogous fashion, observing that we can absorb the loss of factors $n_j$ due to the presence of the $\partial_{\theta}$ as the resulting singular terms will be of level $l\geq 1$ as source admissible terms, as well as Prop.~\ref{prop:singularmultiplier}  and Prop.~\ref{prop:bilin1}. 

\end{proof}

The preceding lemma dealt with terms localized away from the origin, due to the cutoff $\chi_{1\lesssim R}$. The following lemma, is analogous to Prop.~\ref{prop:bilin4}, and proved similarly by invoking the explicit structure away from $R = \nu\tau$ of the admissibly singular terms as exemplified by Lemma~\ref{lem:singFouriertiphysicalngeq2adm}: 
\begin{lemma}\label{lem:basicproductsourcenear origin} Let $\phi_j,\,j = 1, 2$ as well as $m, n_1, n_2, \Lambda_{1,2}$ be as in the preceding lemma. Then there exist constants $\gamma_k,\,k = 0, 1, 2, 3$, such that 
\[
\chi_{R\lesssim 1}\cdot\left(\phi_1\cdot\phi_2 - \sum_{k=0}^3\gamma_k\cdot R^k\right)\in \tilde{S}_1^{\hbar}, 
\]
and we have the estimate 
\begin{align*}
\sup_{\tau\geq\tau_0}\tau^4\left\|\left\langle \phi_m(R;\xi),\,\chi_{R\lesssim 1}\cdot\left(\phi_1\cdot\phi_2 - \sum_{k=0}^3\gamma_k\cdot R^k\right)\right\rangle_{L^2_{R^3\,dR}}\right\|_{S_1^{\hbar}}\lesssim \min\{|n_{1,2}|\}^C\cdot \prod_j\Lambda_j.
\end{align*}

\end{lemma} 

The next lemma is also proved analogously to the preceding two, using the precise structure of $U$ in Theorem \ref{thm:KSTGao precise}:
\begin{lemma}\label{lem:termswithU1} Let $\phi_{1,2}, n_{1,2}, m,\Lambda_{1,2}$ be as in Lemma~\ref{lem:basicproductsource}. Then the expression 
\[
 \Phi: = \chi_{1\lesssim R\lesssim\nu\tau}\cdot \sin^{l_1}U\cdot \cos^{l_2}U\langle R\rangle^{-2}\phi_1\cdot\phi_2,\,l_{1,2}\in \mathbb{N}_{\geq 0}
 \]
 satisfies the same conclusion as $\Phi$ in Lemma~\ref{lem:basicproductsource}.\\
 Similarly, setting 
 \[
 \Phi: = \chi_{R\lesssim 1}\cdot \sin^{l_1}U\cdot \cos^{l_2}U\langle R\rangle^{-2}\phi_1\cdot\phi_2,\,l_{1,2}\in \mathbb{N}_{\geq 0},
 \]
 the same conclusion as in Lemma~\ref{lem:basicproductsourcenear origin} obtains. 

\end{lemma}

The proof is again a simple variation of the preceding. 

\begin{lemma}\label{lem:transferenceonadmsingoffdiaglowfreq} Assume that $\overline{x} = \overline{x}_{1,smooth}^{(l)}$. Then setting 
\begin{align*}
y(\tau,\xi): = \int_0^\infty \chi_{\gtrsim \hbar^{\frac23}}\left(\frac{\left|\xi^{\frac12} - \eta^{\frac12}\right|}{\xi^{\frac12}}\right)\frac{F_n(\xi,\eta)\tilde{\rho}^{(\hbar)}(\eta)}{\xi - \eta}\cdot \overline{x}(\tau,\eta)\,d\eta,
\end{align*}
we have the estimate 
\begin{align*}
\tau^{3+}\cdot \left\|y(\tau,\cdot)\right\|_{S_0^{\hbar}}\lesssim \left\|\overline{x}\right\|_{smooth}. 
\end{align*}
\end{lemma}
\begin{proof} Recall from the definition of functions of type $\overline{x}_{1,smooth}^{(l)}$ that $\eta\lesssim \hbar^{-2}$ on their support. If $\xi\gg \hbar^{-2}$, then using the gain in $\left(\frac{\eta}{\xi}\right)^k$ coming from Prop.~\ref{prop:DFT nlarge}, one can easily repeat the proof of the preceding lemma in this case. We shall henceforth assume $\xi\lesssim \hbar^{-2}$. Furthermore, assume 
\begin{align*}
\overline{x}(\tau, \xi) = \hbar^{-1}\cdot \chi_{\hbar^2\xi\lesssim 1}\frac{e^{\pm i\left(\nu\tau\xi^{\frac12} + \hbar^{-1}\rho_n(x_{\tau};\alpha,\hbar)\right)}}{\xi^{1+\frac{\nu}{2}}}\left(\log\xi\right)^j\cdot \chi_{x_{\tau}\gg1}F_{j}^{\pm}(\tau,\xi), 
\end{align*}
We need to show that 
\begin{align*}
\tau^{3+}\cdot \left\|\chi_{\xi\lesssim \hbar^{-2}}\cdot y(\tau,\cdot)\right\|_{S_0^{\hbar}}\lesssim \left\|\overline{x}\right\|_{smooth}. 
\end{align*}
In order to obtain the added temporal decay, we have to exploit integration by parts with respect to $\eta^{\frac12}$. We note that similarly to an inequality in the proof of the preceding lemma, we have 
\[
\partial_{\eta^{\frac12}}\left(\nu\tau + \hbar^{-1}\rho_n(x_{\tau};\alpha,\hbar)\right)\sim \tau,
\]
provided $x_{\tau}\gg 1$. Taking advantage of Prop.~\ref{prop:DFT nlarge}, we infer the estimates 
\begin{align*}
\left|\partial_{\eta^{\frac12}}^k\left(\frac{F_n(\xi,\eta)\tilde{\rho}^{(\hbar)}(\eta)}{\xi^{\frac12} - \eta^{\frac12}}\right)\right|\lesssim \left(\hbar^{\frac23}\xi^{\frac12}\right)^{-1-k}\cdot \left(\xi^{\frac12}\hbar\right)^2,\,0\leq k\leq 2,\,|\xi^{\frac12} - \eta^{\frac12}|\gtrsim \hbar^{\frac23}\cdot\xi^{\frac12}, 
\end{align*}
and in particular
\begin{align*}
\left|\partial_{\eta^{\frac12}}^2\left(\frac{F_n(\xi,\eta)\tilde{\rho}^{(\hbar)}(\eta)}{\xi^{\frac12} - \eta^{\frac12}}\right)\right|\lesssim \left(\hbar^{\frac23}\xi^{\frac12}\right)^{-3}\cdot \left(\xi^{\frac12}\hbar\right)^2\lesssim \xi^{-\frac12}. 
\end{align*}
Moreover, if $\xi\ll\eta$ or $\eta\ll\xi$, we gain additional factors $\left(\frac{\xi}{\eta}\right)^{\pm k}$, $k\gg 1$. It follows after twofold integration by parts with respect to $\eta^{\frac12}$ that we can write 
\begin{align*}
\chi_{\xi\lesssim \hbar^{-2}}\cdot y(\tau,\cdot) = \chi_{\xi\lesssim \hbar^{-2}}\cdot \int_0^\infty \chi_{\eta^2\hbar\lesssim 1}H(\xi,\eta;\tau)\,d(\eta^{\frac12}), 
\end{align*}
where we have the bound  (with $k\gg 1$)
\begin{align*}
\left|H(\xi,\eta;\tau)\right|\lesssim \hbar^{-1}\cdot \tau^{-3-\nu+}\cdot \xi^{-\frac12}\cdot \eta^{-1-\frac{\nu}{2}}\cdot \frac{1}{\left(\frac{\xi}{\eta}\right)^k +\left( \frac{\eta}{\xi}\right)^k}\cdot \left\|\overline{x}\right\|_{smooth}. 
\end{align*}
This inequality then implies the desired estimate 
\begin{align*}
\left\|\chi_{\xi\lesssim \hbar^{-2}}\cdot y(\tau,\cdot)\right\|_{S_0^{\hbar}}&\lesssim \tau^{-3-\nu+}\cdot \left\|\left(\hbar^2\xi\right)^{1-\frac{\delta}{2}}\cdot \chi_{\xi\lesssim \hbar^{-2}}\cdot\int_0^\infty \chi_{\eta^2\hbar\lesssim 1}H(\xi,\eta;\tau)\,d(\eta^{\frac12})\right\|_{L^2_{d\xi}}\\
&\lesssim  \tau^{-3-\nu+}\cdot \left\|\overline{x}\right\|_{smooth}.
\end{align*}
\end{proof}

Noting the restriction to $x_{\tau}\gg 1$ in the definition of the first of the structured smooth terms (with poor temporal decay) in Definition~\ref{defi:goodfourierrepngeq2}, we also need an additional lemma allowing the introduction of such a cutoff after application of the transference operator. This follows from 
\begin{lemma}\label{lem:transferenceonadmsingoffdiaglowfreqverylowoutputfreq} Assume that  $\overline{x} = \overline{x}_{1,smooth}^{(l)}$ and set 
\begin{align*}
y(\tau,\xi): = \chi_{x_{\tau}\lesssim 1}\int_0^\infty \frac{F_n(\xi,\eta)\tilde{\rho}^{(\hbar)}(\eta)}{\xi - \eta}\cdot \overline{x}(\tau,\eta)\,d\eta. 
\end{align*}
Then we have the estimate 
\begin{align*}
\left\|y(\tau,\cdot)\right\|_{S_0^{\hbar}}\lesssim \tau^{-4+}\cdot \left\|\overline{x}\right\|_{smooth}. 
\end{align*}
\end{lemma}
\begin{proof} Note that $x_{\tau}\lesssim 1$ implies $\xi^{\frac12}\hbar\lesssim \tau^{-1}$. Taking advantage of Prop.~\ref{prop: K operator}, Prop.~\ref{prop:DFT nlarge} and \eqref{def S0 hbar}, we easily infer the bound 
\begin{align*}
\left\|y(\tau,\cdot)\right\|_{S_0^{\hbar}}\lesssim\left\|(\xi\hbar^2)^{1-\frac{\delta}{2}}\cdot y(\tau,\cdot)\right\|_{L^2_{d\xi}}\lesssim \tau^{-4+}\cdot \left\|\overline{x}\right\|_{smooth}. 
\end{align*}

\end{proof}

The following definition is the analogue of Definition~\ref{defi:goodfourierrepngeq2}, but adjusted to the case of exceptional angular momenta. It gives the final concept of `good function' to be used for the distorted Fourier transforms of $\mathcal{D}_+\varepsilon_+[n],\,n\in \{0,\pm 1\}$. Again we note that the decay rates here are better than for the angular modes $|n|\geq 2$, which comes from application of modulations to improve the features of the exceptional modes: 

\begin{definition}\label{defi:xsingulartermsnless2smooth} For functions of angular momentum $n =1$, we say that $\xb(\tau, \xi)$, the distorted Fourier transform of $\mathcal{D}_1\varepsilon_+[1]$  is a good function, provided we can write 
\[
\xb = \xb_{smooth} + \xb_{singular},
\]
where $ \xb_{singular}$ is admissibly singular according to Def.~\ref{defi:xsingulartermsnless2proto}, while the smooth part $\xb_{smooth}$ can be decomposed as
\[
\xb_{smooth} = \xb_{1,smooth} + \xb_{2,smooth},
\]  
with the following properties. For each $\tau\geq \tau_0$, we have the estimate 
\begin{align*}
\left\| \xb_{2,smooth}\right\|_{S_0^{(1)}} + \left\| \mathcal{D}_{\tau}\xb_{2,smooth}\right\|_{S_1^{(1)}}\lesssim \tau^{-(3-10\nu)}. 
\end{align*}
Furthermore we can write 
\begin{align*}
 \xb_{1,smooth} =& \left\langle\xi\right\rangle^{-2}\frac{e^{\pm i\nu\tau\xi^{\frac12}}}{\xi^{1+\delta}}\cdot \int_{\tau_0}^{\tau} \chi_{\xi^{\frac12}\sigma\cdot \frac{\lambda(\tau)}{\lambda(\sigma)}\gtrsim 1}F^{(\pm)}\left(\tau, \sigma, \frac{\lambda^2(\tau)}{\lambda^2(\sigma)}\xi\right)\,d\sigma\\
 & + \left\langle\xi\right\rangle^{-2}\frac{1}{\xi^{1+\delta}}\cdot\int_{\tau_0}^{\tau} \chi_{\xi^{\frac12}\sigma\cdot \frac{\lambda(\tau)}{\lambda(\sigma)}\gtrsim 1}e^{\pm i\left[\left(\nu\tau - 2\frac{\lambda(\tau)}{\lambda(\sigma)}\nu\sigma\right)\xi^{\frac12} \right]}\cdot \tilde{F}^{\pm}\left(\tau,\sigma, \frac{\lambda^2(\tau)}{\lambda^2(\sigma)}\xi\right)\,d\sigma\\
 & + \left\langle\xi\right\rangle^{-2}\frac{1}{\xi^{1+\delta}}\cdot\int_0^\infty \int_{\tau_0}^{\tau} \chi_{\xi^{\frac12}\sigma\cdot \frac{\lambda(\tau)}{\lambda(\sigma)}\gtrsim 1}e^{\pm i\left[\nu\left(\frac{\lambda(\tau)}{\lambda(\sigma)}x + \tau\right)\xi^{\frac12} \right]}\cdot G^{\pm}\left(\tau, \sigma,x,\frac{\lambda^2(\tau)}{\lambda^2(\sigma)}\xi\right)\,d\sigma dx
\end{align*}
where we have the bounds
\begin{align*}
		&\left|\xi^{k_2}\partial_{\tau}^{\iota}\partial_{\xi}^{k_2} F^{(\pm)}(\tau, \sigma, \xi)\right|\lesssim  \left(\log\tau\right)^{N_1}\tau^{-3-\iota - \nu}\cdot\sigma^{-1}\cdot \left[\sigma^{-2} + \kappa\left(\xi^{\frac12}\right)\right],\quad 0\leq k_2\leq 10,\quad\iota\in \{0,1\}\\
		&\left|\xi^{k_2}\partial_{\tau}^{\iota}\partial_{\xi}^{k_2} \tilde{F}^{(\pm)}(\tau, \sigma, \xi)\right|\lesssim  \left(\log\tau\right)^{N_1}\tau^{-3-\iota-\nu}\cdot\sigma^{-1}\cdot \left[\sigma^{-2} + \kappa\left(\xi^{\frac12}\right)\right],\quad 0\leq k_2\leq 10,\quad \iota\in \{0,1\}\\
		&\left\|\xi^{10+\delta}\partial_{\tau}^{\iota}\partial_{\xi}^{10} F^{(\pm)}(\tau, \sigma, \xi)\right\|_{\dot{C}^{\delta}_{\xi}(\xi\simeq\lambda)}\lesssim \left(\log\tau\right)^{N_1}\tau^{-3-\iota-\nu}\cdot\sigma^{-1}\cdot \left[\sigma^{-2} + \kappa\left(\lambda^{\frac12}\right)\right],\quad \iota\in \{0,1\}\\
		&\left\|\xi^{10+\delta}\partial_{\tau}^{\iota}\partial_{\xi}^{10} \tilde{F}^{(\pm)}(\tau, \sigma, \xi)\right\|_{\dot{C}^{\delta}_{\xi}(\xi\simeq\lambda)}\lesssim \left(\log\tau\right)^{N_1}\tau^{-3-\iota-\nu}\cdot\sigma^{-1}\cdot \left[\sigma^{-2} + \kappa\left(\lambda^{\frac12}\right)\right],\quad \iota\in \{0,1\}\\
		&\left\|\xi^{k_2}\partial_{\xi}^{k_2} \partial_{\tau}^{\iota}G^{(\pm)}(\tau, \sigma, x, \xi)\right\|_{L_x^1}\lesssim  \left(\log\tau\right)^{N_1}\tau^{-3-\iota-\nu}\cdot\sigma^{-1}\cdot \left[\sigma^{-2} + \kappa\left(\xi^{\frac12}\right)\right],\quad 0\leq k_2\leq 10,\\
		&\left\|\left\|\xi^{10+\delta}\partial_{\xi}^{10} \partial_{\tau}^{\iota}G^{(\pm)}(\tau, \sigma, x, \xi)\right\|_{\dot{C}^{\delta}_{\xi}(\xi\simeq\lambda)}\right\|_{L_x^1}\lesssim  \left(\log\tau\right)^{N_1}\tau^{-3-\iota-\nu}\cdot\sigma^{-1}\cdot \left[\sigma^{-2} + \kappa\left(\lambda^{\frac12}\right)\right],
	\end{align*}
The norm $\|\cdot\|_{good}$ is then defined in analogy to Def.~\ref{defy:goodfunctionsnorm}. The concept of a 'good source function' at angular momentum $n = 1$ is then defined analogously to the case $|n|\geq 2$ in Definition~\ref{defi:goodfourierrepngeq2}, as is the norm $\|\cdot\|_{goodsource}$. 
The cases $n = 0,\,-1$ lead to similar defnitions by adjusting the weights in analogy to  Def.~\ref{defi:xsingulartermsnless2proto}.
\end{definition}

The following lemma is analogous to Lemma~\ref{lem:transferenceongoodngeq2}:
\begin{lemma}\label{lem:transferenceongoodnless2} The conclusion of Lemma~\ref{lem:transferenceongoodngeq2} holds upon replacing $\mathcal{K}_{\hbar}^{(0)}$ by $\mathcal{K}_{-1}^0$, $\mathcal{K}_0$, see subsections~\ref{subsec:n=-1transference}, ~\ref{subsec:linear end}, provided we use the preceding definition for the norms
\end{lemma}
The proof is analogous to the one of Prop.~\ref{prop:transferenceonsingularngeq2}.
\\

In order to work with functions whose distorted Fourier transform conforms to Def.~\ref{defi:goodfourierrepngeq2}, we need a lemma that allows us to translate functions which are smooth but with weak temporal decay but added structure, namely those with distorted Fourier transform of type $x_{1,smooth}^{l}$, see Remark~\ref{rem:defi:goodfourierrepngeq2}: 
\begin{lemma}\label{lem:structuredlowfreqsmoothbasic} Let $f = f(\tau, R)$ be an angular momentum $|n|\geq 2$ function whose Fourier transform is of type $\overline{x} = \overline{x}_{1,smooth}^{(l)}$ as in Def.~\ref{defi:goodfourierrepngeq2} . Then the function $f$ satisfies
\[
\left|f|_{R<\nu\tau}\right|\lesssim  \hbar^{-\frac12}\cdot \tau^{-1}\cdot\left\|\overline{x}\right\|_{smooth}
\]
where the norm is as in Def.~\ref{defy:goodfunctionsnorm}. For the derivatives we have the better bounds 
\begin{align*}
&\left|\left(\partial_{\tau} + \frac{\lambda_{\tau}}{\lambda}R\partial_R - \partial_R\right)f|_{R<\nu\tau}\right|\lesssim \hbar^{-\frac32}\cdot \tau^{-2}\cdot\left\|\overline{x}\right\|_{smooth},\\&\left|\left(\partial_{\tau} + \frac{\lambda_{\tau}}{\lambda}R\partial_R + \partial_R\right)f|_{R<\nu\tau}\right|\lesssim \hbar^{-\frac32}\cdot \tau^{-1}\cdot\langle R-\nu\tau\rangle^{-1}\cdot\left\|\overline{x}\right\|_{smooth}\\
&\left|\left(\partial_{\tau} + \frac{\lambda_{\tau}}{\lambda}R\partial_R + \partial_R\right)f|_{R<\nu\tau}\right|\lesssim \hbar^{-\frac32}\cdot \tau^{-1}\cdot\langle R-\nu\tau\rangle^{-\frac12}\cdot \langle R\rangle^{-\frac12}\cdot\left\|\overline{x}\right\|_{smooth}
\end{align*}
and we have the improved weighted bound 
\begin{align*}
\left|R^{-1}\cdot f\right|\lesssim \hbar^{-\frac12}\cdot  \tau^{-2}\cdot\left\|\overline{x}\right\|_{smooth}. 
\end{align*}
For the higher derivatives, we have the bounds ($l\geq 1$)
\begin{align*}
&\left|\nabla_{R}^l\left(\partial_{\tau} + \frac{\lambda_{\tau}}{\lambda}R\partial_R - \partial_R\right)f|_{R<\nu\tau}\right|\lesssim \hbar^{-\frac32-l}\cdot \langle \nu\tau - R\rangle^{-1-l}\cdot\tau^{-2}\cdot\left\|\overline{x}\right\|_{smooth},\\&\left|\nabla_{R}^l\left(\partial_{\tau} + \frac{\lambda_{\tau}}{\lambda}R\partial_R + \partial_R\right)f|_{R<\nu\tau}\right|\lesssim \hbar^{-\frac32-l}\cdot \tau^{-1}\cdot \langle \nu\tau - R\rangle^{-2-l}\cdot\left\|\overline{x}\right\|_{smooth},\\
&\left|R^{-1}\cdot \nabla_{R}^lf\right|\lesssim \hbar^{-\frac12-l}\cdot \langle \nu\tau - R\rangle^{-1-l}\cdot  \tau^{-2}\cdot\left\|\overline{x}\right\|_{smooth}. 
\end{align*}
\end{lemma}
\begin{proof} This follows from the detailed asymptotics in Prop.~\ref{prop:DFT nlarge}. For the first estimate, it suffices to invoke the bound $\left|\phi_n(R;\xi)\right|\lesssim \hbar^{\frac12}$, $\left|\xi^{-\delta-}\cdot \phi_n(R;\xi)\right|\lesssim \hbar^{\frac12}\cdot R^{\delta+}$. Note that from Def.~\ref{defi:goodfourierrepngeq2}  we have for $\overline{x} = \overline{x}_{1,smooth}^{(l)}$ the bound 
\begin{align*}
\left\|\min\{\xi^{\delta+}, 1\}\cdot \overline{x}(\tau,\xi)\right\|_{L^1_{d\xi}}\lesssim \hbar^{-1}\cdot \tau^{-1-\nu+}\cdot \left\|\overline{x}\right\|_{smooth}. 
\end{align*}
It follows that 
\begin{align*}
\left|\int_0^\infty \phi_n(R;\xi)\cdot \overline{x}(\tau,\xi)\cdot\rho_n(\xi)\,d\xi|_{R<\nu\tau}\right|\lesssim \hbar^{-\frac12}\cdot\tau^{-1-\nu+\delta+}\cdot \left\|\overline{x}\right\|_{smooth}\lesssim  \hbar^{-\frac12}\cdot\tau^{-1}\cdot \left\|\overline{x}\right\|_{smooth}
\end{align*}
provided $\delta\ll\nu$. 
\\
For the estimate involving the 'good' derivative $\partial_{\tau} + \frac{\lambda_{\tau}}{\lambda}R\partial_R - \partial_R$, we need to recall the fine structure of $\phi_n(R;\xi)$. To begin with, we recall the identity
\begin{equation}\label{eq:goodderidentity}
\partial_{\tau} + \frac{\lambda_{\tau}}{\lambda}R\partial_R - \partial_R = \partial_{\tau} + \nu\partial_R + (1+\nu^{-1})\cdot \frac{R-\nu\tau}{\tau}\partial_R. 
\end{equation}
Then we decompose the Fourier representation of $f$ into three regions, depending on the oscillatory character of $\phi_n(R;\xi)$. Following the proof of Prop.~\ref{prop:derivative1}, we distinguish between the cases $R\hbar\xi^{\frac12}<\frac{x_t}{2},\,R\hbar\xi^{\frac12}\in [\frac{x_t}{2}, 2x_t],\,R\hbar\xi^{\frac12}>2x_t$, where $x_t$ denotes the turning point. Thus we write
\[
f = f_I + f_{II} + f_{III}
\]
where the terms on the right are obtained by including a smooth cutoff $\chi_{R\hbar\xi^{\frac12}<\frac{x_t}{2}}, \chi_{R\hbar\xi^{\frac12}\in [\frac{x_t}{2}, (1+c)x_t]}, \chi_{R\hbar\xi^{\frac12}>(1+c)x_t},\,0<c\ll 1$, in the Fourier integral. We estimate these terms individually by using the asymptotics for $\phi_n(R;\xi)$ in each sub-interval:
\\

{\it{The estimate for $f_I$; non-oscillatory regime.}} Write 
\begin{align*}
f_I = \int_0^\infty\chi_{R\xi^{\frac12}\hbar<\frac{x_t}{2}}\cdot \phi_n(R;\xi)\cdot  \overline{x}(\tau,\xi)\cdot\rho_n(\xi)\,d\xi. 
\end{align*}
From Prop.~\ref{prop:DFT nlarge} we can write
\begin{align*}
\partial_R\left(\chi_{R\xi^{\frac12}\hbar<\frac{x_t}{2}}\cdot\phi_n(R;\xi)\right) = \xi^{\frac12}\cdot \psi_n(R;\xi), 
\end{align*}
where we have the bounds 
\begin{equation}\label{eq:psinbounds}
\left|\partial_{\xi^{\frac12}}^{\kappa}\left(\psi_n(R;\xi)\right)\right|\lesssim_{\kappa}\min\{\xi^{\frac12}R\hbar,1\}(c_1x_t)^{c_2n}\xi^{-\frac{\kappa}{2}},\,0\leq\kappa\ll \hbar^{-1}
\end{equation}
for $c_1<1,\,c_2\ll 1$. Now assuming that $ \overline{x}(\tau,\xi) =  \overline{x}_{1,smooth}^{(l)}(\tau,\xi) $ is of the first type in Def.~\ref{defi:goodfourierrepngeq2}, the other contributions being treated similarly,  we perform one integration by parts with respect to $\xi^{\frac12}$ to infer that 
\begin{align*}
\partial_R f_I &=  \int_0^\infty\xi^{\frac12}\cdot\psi_n(R;\xi)\cdot  \overline{x}(\tau,\xi)\cdot\rho_n(\xi)\,d\xi\\
& = \int_0^\infty\partial_{\xi^{\frac12}}\left(\xi^{\frac12}\cdot\psi_n(R;\xi)\cdot\rho_n(\xi)\right)\cdot  \overline{x}_1(\tau,\xi)\,d\xi\\
& +  \int_0^\infty\xi^{\frac12}\cdot\psi_n(R;\xi)\cdot  \overline{x}_2(\tau,\xi)\cdot\rho_n(\xi)\,d\xi,\\
\end{align*}
where we set 
\begin{align*}
&\overline{x}_1(\tau,\xi) = \chi_{\hbar^2\xi\lesssim 1}\hbar^{-1}\cdot \frac{e^{\pm i\nu\tau\xi^{\frac12}}}{\mp i\nu\tau\xi^{1+\delta}}\cdot\int_{\tau_0}^{\tau}e^{\pm i\hbar^{-1}\rho_n\left(x_{\sigma\cdot\frac{\lambda(\tau)}{\lambda(\sigma)}};\alpha\cdot\frac{\lambda(\tau)}{\lambda(\sigma)},\hbar\right)}\cdot F^{(\pm)}\left(\tau,\sigma, \frac{\lambda^2(\tau)}{\lambda^2(\sigma)}\xi\right)\,d\sigma,\\
&\overline{x}_2(\tau,\xi) = e^{\pm i\nu\tau\xi^{\frac12}}\cdot\partial_{\xi^{\frac12}}\left(\chi_{\hbar^2\xi\lesssim 1}\hbar^{-1}\cdot \frac{1}{\mp i\nu\tau\xi^{1+\delta}}\cdot\int_{\tau_0}^{\tau}e^{\pm i\hbar^{-1}\rho_n\left(x_{\sigma\cdot\frac{\lambda(\tau)}{\lambda(\sigma)}};\alpha\cdot\frac{\lambda(\tau)}{\lambda(\sigma)},\hbar\right)}\cdot F^{(\pm)}\left(\tau,\sigma, \frac{\lambda^2(\tau)}{\lambda^2(\sigma)}\xi\right)\,d\sigma\right).
\end{align*}
Taking advantage of Lemma~\ref{lem: Lemma 3.4 CDST} and Def.~\ref{defi:goodfourierrepngeq2}, we conclude that 
\begin{align*}
\xi^{1+\delta}\cdot\tau^{2+\nu-}\cdot\left(\left|\overline{x}_1(\tau,\xi)\right| + \xi^{\frac12}\left|\overline{x}_2(\tau,\xi)\right|\right)\lesssim \left\|\overline{x}\right\|_{smooth}. 
\end{align*}
Also taking advantage of \eqref{eq:psinbounds} and the preceding decomposition of $\partial_R f_I $ as well as Prop.~\ref{prop:DFT nlarge}, we conclude that 
\begin{align*}
\left|\partial_R f_I|_{R<\nu\tau}\right|\lesssim \tau^{-2}\cdot  \left\|\overline{x}\right\|_{smooth}. 
\end{align*}
The corresponding estimate for $\partial_{\tau}f_I$ is analogous, and the desired estimate for 
\[
\left(\partial_{\tau} + \frac{\lambda_{\tau}}{\lambda}R\partial_R - \partial_R\right)f_I|_{R<\nu\tau}
\]
follows. In particular, there is no power loss in $\hbar^{-1}$ for this contribution.
\\

{\it{The estimate for $f_{II}$; transition regime.}} We use arguments from the proof of Prop.~\ref{prop:derivative1}. We need to distinguish between the non-oscillatory regime $x - x_t\in [-\frac{x_t}{2}, \hbar^{\frac23}]$ and the (weakly) oscillatory regime $x - x_t\in [\hbar^{\frac23}, (1+c)x_t],\,0<c\ll1$, where we recall subsections~\ref{sec: ngeq2 linear}, ~\ref{sec:transfAiry}  for the properties of $x, x_t$. In the former (non-oscillatory) situation, we may localize smoothly to $x - x_t\sim -\mu\hbar^{\frac23}$, $\mu\geq 1$, or $|x-x_t|\lesssim \hbar^{\frac23}$, and exploit the bound (as in the proof of Prop.~\ref{prop:derivative1})
\begin{align*}
\left|\partial_R\phi_n(R;\xi)\right|\lesssim \hbar^{\frac23}\xi^{\frac12}\mu^{-\frac14}e^{-\frac23\mu^{\frac32}},\,x - x_t\sim -\mu\hbar^{\frac23},
\end{align*}
together with the double derivative bound 
\begin{align*}
\left|\partial^2_{\xi^{\frac12}R}\phi_n(R;\xi)\right|\lesssim\mu^{-\frac14}e^{-\frac23\mu^{\frac32}},\,x - x_t\sim -\mu\hbar^{\frac23},
\end{align*}
Since for fixed $R$ the variable $\xi$ is restricted to an interval $I$ of length $\lesssim \hbar^{\frac23}\cdot \xi\cdot \mu$ provided $x - x_t\sim -\mu\hbar^{\frac23}$, we deduce upon integrating by parts once with respect to $\xi^{\frac12}$ the bounds\footnote{The additional loss of $\hbar^{-1}$ comes from the phase function $e^{\pm i\hbar^{-1}\rho\big(x_{\sigma\cdot\frac{\lambda(\tau)}{\lambda(\sigma)}};\alpha\cdot\frac{\lambda(\tau)}{\lambda(\sigma)},\hbar\big)}$ implicit in the definition of $\overline{x}$.} (with $R\lesssim \tau$)
\begin{align*}
&\left|\int_0^\infty \chi_{x - x_t\sim -\mu\hbar^{\frac23}}\cdot \partial_R\phi_n(R;\xi)\cdot \overline{x}(\tau,\xi)\cdot \rho_n(\xi)\,d\xi\right|\lesssim \hbar^{-\frac43}\cdot \mu^{\frac34}e^{-\frac23\mu^{\frac32}}\cdot \tau^{-2}\cdot\left\|\overline{x}\right\|_{smooth},\\
&\left|\int_0^\infty \chi_{|x - x_t|\lesssim \hbar^{\frac23}}\cdot \partial_R\phi_n(R;\xi)\cdot \overline{x}(\tau,\xi)\cdot \rho_n(\xi)\,d\xi\right|\lesssim \hbar^{-\frac43}\cdot \tau^{-2}\cdot\left\|\overline{x}\right\|_{smooth}.
\end{align*}
Summing over $\mu\geq 1$ gives the desired bound for the contribution of $\partial_R f_{II}$, and the function $\partial_{\tau}f_{II}$ is handled similarly. 
\\
In the weakly oscillatory regime $x - x_t\in [\hbar^{\frac23}, (1+c)x_t],\,0<c\ll1$, we have to combine the oscillatory phases in $\phi(R;\xi)$ and $\overline{x}(\tau,\xi)$ (which we assume is of type $\overline{x}_{1,smooth}^{(l)}$) for the integration by parts with respect to $\xi^{\frac12}$. We observe that we can write
\begin{align*}
\partial_R\phi_n(R;\xi) = \hbar^{\frac23}\xi^{\frac12}\cdot \psi_n(R,\xi)\cdot \Ai\left(-\hbar^{-\frac23}\tilde{\tau}\right),\,x - x_t\in [\hbar^{\frac23}, (1+c)x_t],
\end{align*}
where the variable $\tilde{\tau}$ is defined like the variable $\tau$ in Prop.~\ref{prop:DFT nlarge}, and the function $\psi_n(R,\xi)$ is bounded and has symbol type bounds with respect to the variables $R,\xi$ uniformly in $\hbar$. Expanding the Airy function as in the proof of Prop.~\ref{prop:DFT nlarge}, we encounter the oscillatory function 
\[
e^{i\hbar^{-1}\cdot \tilde{\tau}^{\frac32}},
\]
and we have 
\[
\partial_{\xi^{\frac12}}\left(\hbar^{-1}\cdot \tilde{\tau}^{\frac32}\right) \sim R\cdot \tilde{\tau}^{\frac12}. 
\]
Since $\tilde{\tau}\sim x - x_t\ll1$, we infer that the oscillatory phase arising by combining the oscillatory part  $e^{\pm i\nu\tau\xi^{\frac12}}$ from $\overline{x}$ with the preceding satisfies 
\begin{align*}
\partial_{\xi^{\frac12}}\left(\pm i\nu\tau\xi^{\frac12} + \hbar^{-1}\cdot \tilde{\tau}^{\frac32}\right)\sim \pm i\tau
\end{align*}
provided $R<\nu\tau$. If we then use integration by parts with respect to $\xi^{\frac12}$, we infer the same bounds as before for this contribution:
\begin{align*}
\left|\int_0^\infty \chi_{x - x_t\in [\hbar^{\frac23}, (1+c)x_t]}\cdot \partial_R\phi_n(R;\xi)\cdot \overline{x}(\tau,\xi)\cdot \rho_n(\xi)\,d\xi\right|\lesssim \hbar^{-\frac43}\cdot \tau^{-2}\cdot\left\|\overline{x}\right\|_{smooth}.
\end{align*}
Next we consider the effect of applying $\partial_{\tau}$ instead of $\partial_R$. Recalling the fine structure of $\overline{x} = \overline{x}_{1,smooth}^{(l)}$, and assuming that $\overline{x} $ is of the first of the three possible types, we infer 
\begin{align*}
\partial_{\tau}\overline{x} = \overline{x}_1 + \overline{x}_2 + \overline{x}_3,  
\end{align*}
where 
\begin{align*}
&\overline{x}_1 =  \chi_{\hbar^2\xi\lesssim 1}\hbar^{-1}\cdot \frac{e^{\pm i\nu\tau\xi^{\frac12}}}{\mp i\nu\xi^{\frac12+\delta}}\cdot\int_{\tau_0}^{\tau}e^{\pm i\hbar^{-1}\rho_n\left(x_{\sigma\cdot\frac{\lambda(\tau)}{\lambda(\sigma)}};\alpha\cdot\frac{\lambda(\tau)}{\lambda(\sigma)},\hbar\right)}\cdot F^{(\pm)}\left(\tau,\sigma, \frac{\lambda^2(\tau)}{\lambda^2(\sigma)}\xi\right)\,d\sigma,\\
&\overline{x}_2 =  \chi_{\hbar^2\xi\lesssim 1}\hbar^{-1}\cdot \frac{e^{\pm i\nu\tau\xi^{\frac12}}}{\xi^{1+\delta}}\cdot\int_{\tau_0}^{\tau}\partial_{\tau}\left(e^{\pm i\hbar^{-1}\rho_n\left(x_{\sigma\cdot\frac{\lambda(\tau)}{\lambda(\sigma)}};\alpha\cdot\frac{\lambda(\tau)}{\lambda(\sigma)},\hbar\right)}\right)\cdot F^{(\pm)}\left(\tau,\sigma, \frac{\lambda^2(\tau)}{\lambda^2(\sigma)}\xi\right)\,d\sigma,\\
\end{align*}
and finally
\begin{align*}
&\overline{x}_3 =  \chi_{\hbar^2\xi\lesssim 1}\hbar^{-1}\cdot \frac{e^{\pm i\nu\tau\xi^{\frac12}}}{\xi^{1+\delta}}\cdot\int_{\tau_0}^{\tau}e^{\pm i\hbar^{-1}\rho_n\left(x_{\sigma\cdot\frac{\lambda(\tau)}{\lambda(\sigma)}};\alpha\cdot\frac{\lambda(\tau)}{\lambda(\sigma)},\hbar\right)}\cdot \partial_{\tau}\left(F^{(\pm)}\left(\tau,\sigma, \frac{\lambda^2(\tau)}{\lambda^2(\sigma)}\xi\right)\right)\,d\sigma
\end{align*}
Since using Lemma~\ref{lem: Lemma 3.2 CDST} we have 
\begin{align*}
\left|\partial_{\tau}\left(e^{\pm i\hbar^{-1}\rho_n\left(x_{\sigma\cdot\frac{\lambda(\tau)}{\lambda(\sigma)}};\alpha\cdot\frac{\lambda(\tau)}{\lambda(\sigma)},\hbar\right)}\right)\right|\lesssim \hbar^{-1}\tau^{-1}\lesssim \xi^{\frac12},\, \int_{\tau_0}^{\tau}\left|\partial_{\tau}\left(F^{(\pm)}\left(\tau,\sigma, \frac{\lambda^2(\tau)}{\lambda^2(\sigma)}\xi\right)\right)\right|\,d\sigma&\lesssim \tau^{-1}\hbar^{-1}\cdot\left\|\overline{x}\right\|_{smooth}\\
&\lesssim \xi^{\frac12}\cdot \left\|\overline{x}\right\|_{smooth},\\
\end{align*}
it suffices to treat the contribution of the term $\overline{x}_1$: note that since $x\sim x_t\sim 1$, we have $R\hbar\xi^{\frac12}\sim 1$. Since $R\leq \tau$, we have 
\[
\tau \hbar\xi^{\frac12}\gtrsim 1\Rightarrow \xi^{-\delta}\lesssim (\tau\hbar)^{2\delta}. 
\]
Assuming for now that $R\sim \tau$ and using the decomposition 
\[
\phi_n(R;\xi) = \hbar^{\frac13}\cdot \tilde{\psi}_n(R,\xi)\cdot \Ai\left(-\hbar^{-\frac23}\tilde{\tau}\right),\,x - x_t\in [\hbar^{\frac23}, (1+c)x_t],
\]
where $ \tilde{\psi}_n(R,\xi)$ satisfies analogous properties as $\psi_n(R;\xi)$ in the decomposition of $\partial_R\phi_n(R;\xi)$ (see Prop.~\ref{prop:DFT nlarge}), we infer from the asymptotics of Airy functions (see \eqref{eq:Ai}) that 
\begin{align*}
&\left|\int_0^\infty \chi_{x - x_t\in [\hbar^{\frac23}, (1+c)x_t]}\cdot\overline{x}_1\cdot \phi_n(R;\xi)\cdot \rho_n(\xi)\,d\xi\right|\\
&\lesssim \int_0^\infty \hbar^{\frac12}\cdot \chi_{x - x_t\in [\hbar^{\frac23}, (1+c)x_t]}\cdot\tilde{\tau}^{-\frac14}\cdot \left|\overline{x}_1\right|\cdot \rho_n(\xi)\,d\xi\\
&\lesssim \hbar^{-\frac32}\tau^{-2}\cdot \left\|\overline{x}\right\|_{smooth}\cdot \int_0^{O(1)}\tilde{\tau}^{-\frac14}\,d\tilde{\tau}\\
&\lesssim \hbar^{-\frac32}\tau^{-2}\cdot \left\|\overline{x}\right\|_{smooth}.
\end{align*}
Here we have used that $\xi^{-\frac12}\cdot d\xi = 2d(\xi^{\frac12}) = \kappa(\xi)\cdot (\hbar R)^{-1}d\tilde{\tau}$, where the function $\kappa(\xi)$ is bounded for $x - x_t\in [\hbar^{\frac23}, (1+c)x_t]$, and furthermore, $\tilde{\tau}\in [0, O(1)]$ under the latter condition. 
The case $R\ll \tau$ requires a further integration by parts with respect to $\xi^{\frac12}$, taking into account the oscillatory character of $\Ai\left(-\hbar^{-\frac23}\tilde{\tau}\right)$. 
The preceding bounds in conjunction with \eqref{eq:goodderidentity} give the desired bound for the 'good derivative' in the transition regime, i. e. for the term $f_{II}$. 
\\

{\it{The estimate for $f_{III}$; oscillatory regime.}} Here the Airy function $\Ai\left(-\hbar^{-\frac23}\tilde{\tau}\right)$ is in the oscillatory regime, and we can take advantage of Lemma~\ref{lem: Lemma 3.4 CDST} in order to infer
\begin{align*}
\frac23\cdot \tilde{\tau}^{\frac32} = x - y(\alpha;\hbar) + \rho(x;\alpha,\hbar),\,\alpha = \xi^{\frac12}\hbar,\,x = R\alpha, 
\end{align*}
and where $y,\,\rho$ satisfy the bounds in loc. cit. Using Prop.~\ref{prop:DFT nlarge}, we can then expand (for $x - x_t\in [[1+c)x_t,\,\infty)$)
\begin{align*}
\phi_n(R;\xi) = \sum_{\pm} R^{-\frac12}\xi^{-\frac14}\cdot e^{\pm iR\xi^{\frac12}}\cdot e^{\pm i\hbar^{-1}\left(- y(\alpha;\hbar) + \rho_n(x;\alpha,\hbar)\right)}\cdot \tilde{\Xi}(R,\xi;\hbar), 
\end{align*}
where the function $\tilde{\Xi}(R,\xi;\hbar)$ is smooth and bounded and satisfies symbol type bounds with respect to the variables $R, \xi$. Also, note that under our restriction on $x$, we have 
\[
R^{-\frac12}\xi^{-\frac14}\lesssim \hbar^{\frac12}. 
\]
Consider then the integral 
\begin{align*}
\int_0^\infty \chi_{x - x_t\in [[1+c)x_t,\,\infty)}\cdot \overline{x}(\tau,\xi)\cdot \phi_n(R;\xi)\cdot \rho_n(\xi)\,d\xi. 
\end{align*}
The idea then is to combine the oscillatory term $e^{\pm i\nu\tau\xi^{\frac12}}$ occuring in $\overline{x}$ (according to Def.~\ref{defi:goodfourierrepngeq2}  and the phase $e^{\pm iR\xi^{\frac12}}$. Observe that the resonant phases 
\[
e^{\pm i(R-\nu\tau)\xi^{\frac12}}
\]
are annihilated by the operator $\partial_{\tau} + \nu\partial_R$ in \eqref{eq:goodderidentity}. On the other hand, using integration by parts with respect to $\xi^{\frac12}$, we infer (for $R\leq\tau$)
\begin{align*}
&\left|\frac{R-\nu\tau}{\tau}\partial_R\left(\int_0^\infty \chi_{x - x_t\in [[1+c)x_t,\,\infty)}\cdot \overline{x}(\tau,\xi)\cdot \phi_n(R;\xi)\cdot \rho_n(\xi)\,d\xi\right)\right|\\
&\lesssim \hbar^{-\frac32}\tau^{-2}\cdot \left\|\overline{x}\right\|_{smooth}.
\end{align*}
We note that the extra loss of a factor $\hbar^{-1}$ arises when $\partial_{\xi^{\frac12}}$ hits $e^{\pm i\hbar^{-1}\left(\pm\rho_n(x;\alpha,\hbar)\pm \rho_n\left(x_{\sigma\cdot\frac{\lambda(\tau)}{\lambda(\sigma)}};\alpha\cdot\frac{\lambda(\tau)}{\lambda(\sigma)},\hbar\right)\right)}$. 
This concludes the argument for the second estimate of the lemma. The remaining ones follow very similar reasoning. For the fourth bound, one uses the fact that 
\begin{align*}
\left|R^{-1}\phi_n(R;\xi)\right|\lesssim \hbar^{\frac32}\xi^{\frac12},
\end{align*}
as well as integration by parts with respect to $\xi^{\frac12}$ in the region $R\ll\tau$. 
\end{proof}

\begin{lemma}\label{lem:tildec1control} Assume that $\big|c_0(\tau)\big| + \tau\cdot \big|c_0'(\tau)\big|\lesssim \tau^{-2-}$. Recalling Lemma~\ref{prop:fourierboundsforn=1termlocalizedaway} , we have the bound 
\begin{align*}
\left\|\tilde{\overline{y}}\right\|_{\langle\xi^{\frac14}\rangle S_1^{(1)}}\lesssim \tau^{-4+}\cdot \Lambda\cdot \left(\tau_0^{-1}+\Lambda\right). 
\end{align*}
Recalling Lemma~\ref{lem:n=1maintermafterparametrix}, the following estimate holds
\begin{align*}
\left||\tilde{c}_1(\tau)\right|\lesssim \tau^{-4-}\cdot \Lambda\cdot \left(\tau_0^{-1}+\Lambda\right). 
\end{align*}
Here $\Lambda$ is defined in \eqref{eq:LambdaDef}. 
\end{lemma}
\begin{proof} We shall outline the bounds for the various nonlinear source terms consituting $\mathcal{D}_+\left(F_+(1)\right)$ which are bounded in $\langle\xi^{\frac14}\rangle S_1^{(1)}$. Throughout we work at angular momentum $n = 1$, whence large angular momenta in the non-linearity need to cancel. The fast decay of angular momentum $n$-functions with $|n|\gg1 $ as implied by the definition of $\Lambda$ allows us to completely neglect all issues related to potential losses with respect to $n$ in our nonlinear estimates. We assume throughout that the angular momenta of the factors sum up to $n = 1$ output angular momentum. 
\\

{\it{(1): Cubic null-form}} 
\[
X_1: = \mathcal{D}_+\left(\chi_{R\lesssim\tau}\cdot\phi_1\cdot \left(\partial_{\tau} + \frac{\lambda_{\tau}}{\lambda}R\partial_R - \partial_R\right)\phi_2\cdot \left(\partial_{\tau} + \frac{\lambda_{\tau}}{\lambda}R\partial_R + \partial_R\right)\phi_3\right).
\]
We need to consider the various possibilities for the factors leading to a (non-structured) term in $\langle\xi^{\frac14}\rangle S_1^{(1)}$:
\\

{\it{(1.a): all factors with distorted Fourier transform of type $\xb_{1,smooth}$.}} This is in some sense the worst case since these factors have the weakest temporal decay. We shall take advantage of the pointwise bounds established in Lemma~\ref{lem:structuredlowfreqsmoothbasic}. 
Then the first bound is a direct consequence of the fact that 
\begin{align*}
\left\|\mathcal{F}^{(1)}\left(X_1\right)\right\|_{\langle\xi^{\frac14}\rangle S_1^{(1)}}\lesssim \left\|\langle \nabla_R\rangle^4 X_1\right\|_{L^2_{R\,dR}}\lesssim_{\hbar_j}\tau^{-4+}\cdot \prod_{j=1}^3\left\|\xb_j\right\|_{smooth}
\end{align*}
where $\xb_j$ denotes the distorted Fourier transform of the $j$-th factor. In fact, this estimate follows by using the pointwise bounds of the factors in light of Lemma~\ref{lem:structuredlowfreqsmoothbasic}. 
\\
For the better bound for $\tilde{c}_1(\tau)$, we need to take advantage of the $\xi$-integral. Precisely, recalling the expression for $U^{(1)}(\tau,\sigma,\xi)$, and combining the oscillatory phase there with the one coming from $R^{-\frac12}\cdot\Phi_1(R;\frac{\lambda^2(\tau)}{\lambda^2(\sigma)}\xi)$ involved in $\mathcal{F}^{(1)}\left(X_1\right)\left(\sigma, \frac{\lambda^2(\tau)}{\lambda^2(\sigma)}\xi\right)$, we arrive at the oscillatory phase $e^{\pm i(\nu\tau - \nu\sigma\cdot \frac{\lambda(\tau)}{\lambda(\sigma)}\pm R\cdot \frac{\lambda(\tau)}{\lambda(\sigma)})\xi^{\frac12}}$. Then we either have 
\begin{align*}
\left|\nu\tau - \nu\sigma\cdot \frac{\lambda(\tau)}{\lambda(\sigma)}\pm R\cdot \frac{\lambda(\tau)}{\lambda(\sigma)}\right|\gtrsim \tau,
\end{align*}
or else
\begin{align*}
\left|\left(R - \nu\sigma\right)\cdot \frac{\lambda(\tau)}{\lambda(\sigma)}\right|\gtrsim \tau. 
\end{align*}
Recalling that $\tilde{c}_1(\tau)$ is a contribution to 
\begin{align*}
\int_0^\infty \chi_{\xi\lambda^2(\tau)\lesssim M}\cdot \xb(\tau,\xi)\cdot\tilde{\rho}_1(\xi)\,d\xi,
\end{align*}
where $ \xb(\tau,\cdot)$ is the Duhamel propagator (with kernel $U^{(1)}(\tau,\sigma,\xi)$) of the preceding source term, in the {\it{former case}}, we perform integration by parts with respect to $\xi^{\frac12}$ twice, which gains a $\tau^{-2}$ at the expense of a factor $\xi^{-1}$, which gets compensated by the spectral density  $\tilde{\rho}_1(\xi)\sim \xi$ at small frequency $\xi\ll 1$. 
In the latter case, we perform integration by parts with respect to $R$ twice in the integral defining $\mathcal{F}^{(1)}\left(X_1\right)\left(\sigma, \frac{\lambda^2(\tau)}{\lambda^2(\sigma)}\xi\right)$, which in light of Lemma~\ref{lem:structuredlowfreqsmoothbasic} gains $(\nu\tau - R)^{-2}$ at the expense of a factor $\left(\frac{\lambda^2(\tau)}{\lambda^2(\sigma)}\xi\right)^{-1}$. One then again gains a factor $\tau^{-2}$, and one factor $\tau$ is spent in bounding the Duhamel time integral over $\sigma$. It follows that this contribution to $\tilde{c}_1(\tau)$ is bounded by $\lesssim\tau^{-\frac92}$, which is as desired. 
\\

{\it{(1.b): all factors with distorted Fourier transform either of type $\xb_{1,smooth}$ or of the form $c_n\cdot \phi_n(R)$, $n = 0, 1, 2$.}} For the first estimate, we observe that by assumption $|c_n(\tau)|\lesssim \tau^{-2+}$, and we have 
\[
 \left\|\chi_{R\lesssim\tau}\cdot \left(\partial_{\tau} + \frac{\lambda_{\tau}}{\lambda}R\partial_R \mp \partial_R\right)\phi_n(R)\right\|_{L^{1+}_{R\,dR}}\lesssim 1.
  \]
  One then easily infers the better bound 
  \begin{align*}
  \left\|\mathcal{F}^{(1)}\left(X_1(\tau,\cdot)\right)\right\|_{\langle\xi^{\frac14}\rangle S_1^{(1)}}\lesssim\tau^{-4+}\cdot \Lambda^3
  \end{align*}
  if at least one factor is of the form $c_n\cdot \phi_n(R)$, $n\in \{0, 1, 2\}$. As for the estimate for $\tilde{c}_1$, we distinguish between the following cases: if there is precisely one factor of type $c_n\cdot \phi_n(R)$, assume first that it is in the un-differentiated position. if $\mathcal{D}_+$ hits $c_n\cdot \phi_n(R)$, we gain a factor $R^{-1}$, if it hits the other factors, we gain $|R-\nu\tau|^{-1}$ (recalling Lemma~\ref{lem:structuredlowfreqsmoothbasic}). In the former case, we can perform integration by parts in the $\xi$-integral giving $\tilde{c}_1$, gaining $\left(\nu\tau - \nu\sigma\cdot \frac{\lambda(\tau)}{\lambda(\sigma)}\right)^{-1}$ but losing a factor $R$ if $\partial_{\xi}$ hits $\mathcal{F}^{(1)}\left(X_1\right)\left(\sigma, \frac{\lambda^2(\tau)}{\lambda^2(\sigma)}\xi\right)$. which is compensated by the gain of $R^{-1}$. Since we have the estimate 
  \begin{align*}
  \left|\mathcal{F}^{(1)}\left(X_1\right)\left(\sigma, \frac{\lambda^2(\tau)}{\lambda^2(\sigma)}\xi\right)\right|\lesssim \sigma^{-5+}
  \end{align*}
   in this case, we can compensate for the $\sigma$-integral by the factor $\left(\nu\tau - \nu\sigma\cdot \frac{\lambda(\tau)}{\lambda(\sigma)}\right)^{-1}$ and obtain a contribution to $\tilde{c}_1$ of size $\tau^{-5+}$. 
   \\
   If $\mathcal{D}_+$ hits one of the other two factors (with distorted Fourier transform of type $\xb_{1,smooth}$, one can argue precisely as in the case {\it{(1.a)}}, by combining the oscillatory phases and distinguishing between a non-oscillatory and an oscillatory case. 
   \\
   Similar arguments apply to the cases when two or more factors of the form $c_n\cdot \phi_n(R)$ are present. 
   \\
   
   {\it{(1.c): at least one factor has distorted Fourier transform in $S_0^{(\hbar)}$ (or in $S_0^{(n)}$ for the exceptional modes).}} This is very similar to the preceding case, as one can invoke the decay rate $\tau^{-3+}$ for these terms (according to the definition of $\|\cdot\|_{good}$), as well as the bounds in Lemma~\ref{lem:derLinfty}, Prop.~\ref{prop:exceptionalnderivative1}.
   \\
   
   {\it{(1.d): at least one admissibly singular factor, but product is placed in $\langle\xi^{\frac14}\rangle S_1^{(1)}$.}} Recalling the estimates in subsection~\ref{subset:hardnullformestimates}, this situation occurs if we apply a finite Taylor expansion to a regular factor (with distorted Fourier transform in $S_0^{(\hbar)}$ or of type $\xb_{1,smooth}$), and replace the regular factor by the Cauchy error in the Taylor expansion. Given the fact that admissibly singular terms enjoy similar spatial and temporal decay estimates (away from the light cone), this case is handled similarly to the preceding two situations. 
\\

In the following we only concentrate on the more difficult estimate for $\tilde{c}_1(\tau)$, as the other bound will follow easily from the estimates below. 
\\

{\it{(2): Quadratic terms:}} We recall from \eqref{non linear 1 2}, \eqref{coe eq precise 1}, \eqref{coe eq precise 2} that these terms are schematically of the form 
\begin{align*}
\mathcal{D}_+\left(\chi_{R\lesssim\tau}\cdot\frac{\sin U}{R^2}\cdot \phi_1\cdot \phi_{2,\theta}\right),\,\mathcal{D}_+\left(\chi_{R\lesssim\tau}\cdot U_R\cdot \phi_1\cdot \phi_{2,R}\right). 
\end{align*}

{\it{(2.a): both factors $\phi_j$ have distorted Fourier transform of type $\xb_{1,smooth}$.}} This case is potentially delicate due to the poor pointwise decay bounds for these in light of Lemma~\ref{lem:structuredlowfreqsmoothbasic}. Here we can exploit the fact that since we restrict the source term to angular momentum $n = 1$, if neither factor is of exceptional angular momentum $n_j\in\{0, 1, 2\}$ (for which we have better pointwise decay due to Def.~\ref{defi:xsingulartermsnless2smooth}), at least one $|n_j|\geq 3$. For such $n_j$, the vanishing properties of $\phi_{n_j}(R;\xi)$ at  $\xi = 0$ easily imply a sharpening of Lemma~\ref{lem:structuredlowfreqsmoothbasic}, namely 
\begin{align*}
\left|\langle R\rangle^{-2}\cdot \phi_j\right| + \left|\langle R\rangle^{-1}\cdot \mathcal{D}_+\phi_j\right|\lesssim_{n_j}\tau^{-3}\cdot \left\|\xb_j\right\|_{smooth}. 
\end{align*}
We then conclude for $\tau\geq \sigma\geq \tau_0$ that we have the estimate 
\begin{align*}
&\left|\langle \zeta_{R,\sigma,\tau}^{\pm}\cdot \mathcal{D}_+\left(\chi_{R\lesssim\tau}\cdot\frac{\sin U}{R^2}\cdot \phi_1\cdot \phi_{2,\theta}\right),\,\Phi_1\left(R;\frac{\lambda^2(\tau)}{\lambda^2(\sigma)}\xi\right)\rangle_{L^2_{R\,dR}}\right|\\
&\lesssim_{n_j}\left(\frac{\lambda(\tau)}{\lambda(\sigma)}\right)^{-\frac32}\cdot\xi^{-\frac34}\cdot \prod_{j=1,2}\left[\left\|\langle R\rangle^{-j-1}\cdot \phi_j\right\|_{L^\infty_{dR}} + \left\|\langle R\rangle^{-j}\cdot \mathcal{D}_+\phi_j\right\|_{L^\infty_{dR}}\right]\\
&\lesssim \left(\frac{\lambda(\tau)}{\lambda(\sigma)}\right)^{-\frac32}\cdot\xi^{-\frac34}\cdot\tau^{-5+}\cdot \prod_{j=1,2}\left\|\xb_j\right\|_{smooth}. 
\end{align*}
where we set $\zeta_{R,\sigma,\tau}^{\pm} = \langle \nu\tau - \nu\sigma\cdot\frac{\lambda(\tau)}{\lambda(\sigma)}\pm R\rangle^{-1}\cdot R^{\frac12}$. For this also recall \eqref{modified funda sys}. Using twofold integration by parts with respect to $\xi$, it then follows that we have the estimate 
\begin{align*}
&\left|\int_0^\infty \chi_{\xi\lambda^2(\tau)\lesssim M}\cdot \int_{\tau_0}^{\tau}U^{(1)}\left(\tau,\sigma,\xi\right)\cdot \mathcal{F}^{(1)}\left(\mathcal{D}_+\left(\chi_{R\lesssim\tau}\cdot\frac{\sin U}{R^2}\cdot \phi_1\cdot \phi_{2,\theta}\right)\right)\left(\sigma,\frac{\lambda^2(\tau)}{\lambda^2(\sigma)}\xi\right)\tilde{\rho}_1(\xi)\,d\sigma d\xi\right|\\
&\lesssim_{n_j} \tau^{-\frac{9}{2}+}\cdot \prod_{j=1,2}\left\|\xb_j\right\|_{smooth}\lesssim \tau^{-4-}\cdot \Lambda\cdot \left(\tau_0^{-1}+\Lambda\right). 
\end{align*}
The second quadratic expression above is handled similarly. 
\\

{\it{(2.b): one factor of type $c_{n_1}\cdot \phi_{n_1}(R)$, the other with distorted Fourier transform of type $\xb_{1,smooth}$.}} Since the source term is restricted to angular momentum $n = 1$, if $n_1 = -1$ then the other factor needs to have $n_2 = 2$. Here we have to use integration by parts with respect to $\sigma$ in order to arrive at an acceptable bound. For this we note that the boundary term at $\sigma = \tau$ generated by such integration by parts vanishes asymptotically as $M\rightarrow +\infty$ as 
\begin{align*}
\lim_{M\rightarrow+\infty}\int_0^\infty \chi_{\xi\lambda^2(\tau)\lesssim M}\cdot\mathcal{F}^{(1)}\left(\mathcal{D}_+\left(\chi_{R\lesssim\sigma}\cdot\frac{\sin U}{R^2}\cdot \phi_1\cdot \phi_{2,\theta}\right)\right)(\tau,\xi)\tilde{\rho}_1(\xi)\,d\sigma d\xi = 0, 
\end{align*}
see e. g. (10.5) in \cite{KMiao}. On the other hand, the other boundary term at $\tau = \tau_0$ is easily seen to lead to an expression bounded by 
\[
\lesssim_{\tau_0}\lambda^{-1}(\tau)\cdot \Lambda\cdot \left(\tau_0^{-1}+\Lambda\right)\ll \tau^{-4-} \Lambda\cdot \left(\tau_0^{-1}+\Lambda\right)
\]
for $\nu$ sufficiently small. We conclude that at the expense of a factor $\left(\xi^{\frac12}\cdot \frac{\lambda(\tau)}{\lambda(\sigma)}\right)^{-1}$ and replacing $\sin\left[\lambda(\tau)\xi^{\frac12}\cdot \int_{\sigma}^{\tau}\lambda^{-1}(s)\,ds\right]$ by $\cos\left[\ldots\right]$, which we indicate by replacing $U^{(1)}(\ldots)$ by $\tilde{U}^{(1)}(\ldots)$, we reduce this case to bounding 
\begin{align*}
\int_0^\infty \chi_{\xi\lambda^2(\tau)\lesssim M}\cdot \int_{\tau_0}^{\tau}\tilde{U}^{(1)}\left(\tau,\sigma,\xi\right)\cdot \partial_{\sigma}\mathcal{F}^{(1)}\left(\mathcal{D}_+\left(\chi_{R\lesssim\sigma}\cdot\frac{\sin U}{R^2}\cdot \phi_1\cdot \phi_{2,\theta}\right)\right)\left(\sigma,\frac{\lambda^2(\tau)}{\lambda^2(\sigma)}\xi\right)\tilde{\rho}_1(\xi)\,d\sigma d\xi
\end{align*}
The operator $\partial_{\sigma}$ gains an extra $\sigma^{-1}$ of decay, and another integration by parts with respect to $\xi$ in order to avoid a loss in time decay due to the integration over $\sigma$ show that this contribution to $\tilde{c}_1(\tau)$ is bounded by $\lesssim \tau^{-5+}\cdot  \Lambda\cdot \left(\tau_0^{-1}+\Lambda\right)$. 
\\
If $n_1 = 0$, then the other factor is also exceptional (namely $n_2 = 1$), and the desired estimate follows more easily due to Def.~\ref{defi:xsingulartermsnless2smooth} and integration by parts with respect to $\xi$ to avoid a loss due to integration over $\sigma$, resulting in a bound as before for this contribution. 
The same argument applies if $n_1 = 1$. 
\\

{\it{(2.c): both factors of type $c_{n_j}\cdot \phi_{n_j}(R)$.}} In this case one of the $n_j = 0$, and we use the assumption of better temporal decay for this in the statement of the lemma, as well as integration by parts with respect to $\xi$ to avoid a loss due to integration over $\sigma$. 
\\

{\it{(2.d): at least one factor with distorted Fourier transform of (unstructured but better decaying) type $S_0^{\hbar}$ or $S_0^{(n)}$.}} This is simpler since on account of Lemma~\ref{lem:derLinfty}, Prop.~\ref{prop:exceptionalnderivative1}, we always get the bound 
\begin{align*}
\left|\langle\partial_{\xi}\rangle\mathcal{F}^{(1)}\left(\mathcal{D}_+\left(\chi_{R\lesssim\sigma}\cdot\frac{\sin U}{R^2}\cdot \phi_1\cdot \phi_{2,\theta}\right)\right)(\tau,\cdot)\right|\lesssim \tau^{-5+}\cdot \Lambda\cdot \left(\tau_0^{-1}+\Lambda\right), 
\end{align*}
and the desired bound for $\tilde{c}_1$ follows by integration by parts with respect to $\xi$ to avoid a loss due to integration over $\sigma$. Again the second quadratic form is treated in the same way.
 \\
 
 {\it{(3) Terms linear in $\phi$.}} We treat here the term 
 \[
 \mathcal{D}_+\left(\frac{\partial_R\epsilon}{1+R^2}\cdot \phi\right)
 \]
 from \eqref{coe eq more precise 2}, the other linear terms being treated similarly. Since the source function is at angular momentum $n = 1$, so is $\phi$. As usual one decomposes $\phi$ into a finite Taylor expansion around $R = \nu\tau$ and includes the contribution of the Cauchy error to this source term into a function whose distorted Fourier transform lives in $\langle\xi^{\frac14}\rangle\cdot S_1^{(1)}$. Observe that 
 \begin{align*}
 \left| \mathcal{D}_+\left(\frac{\partial_R\epsilon}{1+R^2}\right)\right|\lesssim \langle R\rangle^{-3}\cdot \tau^{-2}\log\tau. 
 \end{align*}
 If $\phi = \phi_1(R)\cdot \int_0^R \phi_1^{-1}(s)\mathcal{D}_+\phi(\cdot,s)\,ds$ as in Prop.~\ref{prop:exceptionalnderivative1}, and the distorted Fourier support is of type $\xb_{1,smooth}$ as in  Def.~\ref{defi:xsingulartermsnless2smooth}, then we easily obtain the estimate 
 \begin{align*}
 \left|\langle \partial_{\xi}\rangle\mathcal{F}^{(1)}\left( \mathcal{D}_+\left(\frac{\partial_R\epsilon}{1+R^2}\cdot \phi\right)\right)\right|\lesssim \tau^{-5+}\cdot \Lambda.
 \end{align*}
 In fact this follows from Lemma~\ref{lem:structuredlowfreqsmoothbasic}.The same bound obtains when the distorted Fourier transform of $\mathcal{D}_+\phi$ is in $S_0^{(1)}$, as follows from Prop.~\ref{prop:exceptionalnderivative1} and the definition of $\|\xb\|_{good}$. 
 One then obtains a bound of $\left|\tilde{c}_1(\tau)\right|\lesssim \tau^{-5+}\cdot \Lambda$ by exploiting integration by parts in the $\xi$-integral to counteract the $\sigma$-integral in the Duhamel propagator. 
 \\
 The situation when $\phi = c_1(\tau)\cdot\phi_1(R)$  is handled like the case {\it{(2.b)}} before via integration by parts with respect to $\sigma$, exploiting integration by parts with respect to $\sigma$ as well as the high order of vanishing of $\epsilon$ in $R = 0$. 
 
\end{proof}

\bibliographystyle{abbrv}
\bibliography{kst.bib}

\bigskip

\centerline{\scshape Joachim Krieger}
\smallskip
{\footnotesize
	\centerline{B\^atiment des Math\'ematiques, \'Ecole Polytechnique F\'ed\'erale de Lausanne (EPFL)}
	\centerline{Lausanne 1015, Switzerland}
	\centerline{\email{joachim.krieger@epfl.ch}}
} 

\medskip

\centerline{\scshape Shuang Miao}
\smallskip
{\footnotesize
	\centerline{School of Mathematics and Statistics, Wuhan University}
	\centerline{Wuhan, Hubei 430072, China}
	\centerline{\email{shuang.m@whu.edu.cn}}
} 

\medskip

\centerline{\scshape Wilhelm Schlag}
\smallskip
{\footnotesize
	\centerline{Department of Mathematics, Yale University}
	\centerline{New Haven, CT 06511, U.S.A.}
	\centerline{\email{wilhelm.schlag@yale.edu}}
} 
\end{document}